\documentclass[a4paper,10pt]{amsart}

\usepackage{amssymb}
\usepackage{fullpage}
\usepackage[OT4]{fontenc}
\usepackage{texdraw}

\theoremstyle{definition}
\newtheorem{definition}{Definition}
\theoremstyle{plain}
\newtheorem{theorem}[definition]{Theorem}
\newtheorem{proposition}[definition]{Proposition}
\newtheorem{lemma}[definition]{Lemma}

\newtheorem{criterion}{Criterion}
\theoremstyle{remark}
\newtheorem{remark}[definition]{Remark}
\newtheorem{example}[definition]{Example}

\DeclareMathOperator{\edim}{edim}
\DeclareMathOperator{\vdim}{vdim}
\DeclareMathOperator{\rank}{rank}
\DeclareMathOperator{\mult}{mult}
\DeclareMathOperator{\mdeg}{mdeg}
\DeclareMathOperator{\rows}{rows}
\DeclareMathOperator{\cols}{cols}

\DeclareMathOperator{\Cr}{Cr}
\DeclareMathOperator{\Hp}{Hp}

\def\ass{\longleftarrow}
\def\field{\mathbb{K}}
\def\N{\mathbb{N}}
\def\Z{\mathbb{Z}}

\def\PP{\mathbb{P}}

\def\sys{\mathcal{L}}

\def\rdf{:=}
\def\dff{\it}
\def\birto{\longmapsto}
\def\stbun{\mathcal O_{\PP^2}(1)}

\let\limto\to
\let\to\longrightarrow

\newenvironment{algorithm}[1]{\medskip\noindent{\bf Algorithm} #1\\}{\medskip}

\begin{document}

\title{Regularity and non-emptyness of linear systems in $\mathbb P^n$}

\author{Marcin Dumnicki}

\dedicatory{
Institute of Mathematics, Jagiellonian University, \\
Reymonta 4, 30-059 Krak\'ow, Poland \\
}

\thanks{Email address: Marcin.Dumnicki@im.uj.edu.pl}





\begin{abstract}
The main goal of this paper is to present a new algorithm bounding the regularity
and ``alpha'' (the lowest degree of existing hypersurface) of a linear system of 
hypersurfaces (in $\mathbb P^n$) passing through 
multiple points in general position.
To do the above we formulate and prove new theorem, which allows
to show non-specialty of linear system by splitting it into non-special
(and simpler) systems.
As a result we give new bounds for
multiple point Seshadri constants on $\PP^2$. 
\end{abstract}

\maketitle

\section{Introduction}

In what follows we assume that the ground field $\field$ is of characteristic zero.
Let $n \geq 2$, let $d \in \Z$, let $m_1,\dots,m_r \in \N$.
By $\sys_n(d;m_1,\dots,m_r)$ we denote the linear system of hypersurfaces
(in $\PP^n \rdf \PP^n(\field)$)
of degree $d$ passing through $r$ points $p_1,\dots,p_r$ in general position with
multiplicities at least $m_1,\dots,m_r$. The dimension of such system is denoted by
$\dim \sys_n(d;m_1,\dots,m_r)$.
Define the {\dff virtual dimension of $L = \sys_n(d;m_1,\dots,m_r)$}
$$\vdim L \rdf \binom{d+n}{n} - \sum_{j=1}^{r} \binom{m_j+n-1}{n} - 1$$
and the {\dff expected dimension of $L$}
$$\edim L  \rdf \max \{ \vdim L, -1 \}.$$
Observe that $\dim L \geq \edim L$. If this inequality is strict then
$L$ is called {\dff special}, {\dff non-special} otherwise.
The system $L$ is called {\dff non-empty} if $\dim L \geq 0$, {\dff empty} otherwise.


The problem of classifying special systems has been widely studied
by many authors, especially for $n=2$ (most of 
papers included in reference section deal with the case $n=2$; there
are many others). For low multiplicities (i.e. bounded by some constant)
the classification is known, it begun with \cite{HirF}, where the case
$m_1=\dots=m_r \leq 3$ was solved;
special systems with (arbitrary) multiplicities
up to $11$ has been classified in \cite{mdwj};
the homogeneous case ($m_1=\dots=m_r$)
with multiplicities up to $42$ in \cite{md};
quasi-homogeneous case ($m_1=\dots=m_{r-1}$, $m_r$ arbitrary)
for $m_1=5$ (for $m_1 < 5$ the classification is also done) in \cite{base5}.
For $r \leq 9$ the classical result is due to Nagata (\cite{nagata}).
In the general case ($n \geq 3$) we note that very few is known, 
results concern either low multiplicities ($m_1=\dots=m_r=2$, $n$ arbitrary) (see \cite{AH2}, \cite{chandler}), or few points
($r \leq 8$, $n=3$ in \cite{p3}).

For a system of multiplicities $m_1,\dots,m_r$ and $n \geq 2$ we define (after Harbourne, \cite{surv})
\begin{align*}
\alpha_n(m_1,\dots,m_r) & \rdf  \min \{ d \in \N : \sys_n(d;m_1,\dots,m_r) \text{ is non-empty}\},\\
\tau_n(m_1,\dots,m_r) & \rdf  \min \{ d \in \N : \sys_n(d;m_1,\dots,m_r) \text{ is non-special and non-empty}\}.
\end{align*}

The second number is often called the {\dff regularity}.
Observe that for $d < \alpha_n(m_1,\dots,m_r)$ the system
$\sys_n(d;m_1,\dots,m_r)$ is empty, similarly,
for $d \geq \tau_n(m_1,\dots,m_r)$ the system $\sys_n(d;m_1,\dots,m_r)$ is
non-special.
More information about $\alpha_2$ and $\tau_2$ can be found in \cite{surv} and \cite{knowalpha}.
Note that usually (for small values of $n$) the conjectured values of $\alpha_n$ and $\tau_n$ are known
(for $n=2$ by Harbourne--Hirschowitz conjecture (see e.g. \cite{HCON}); for $n = 3$ the conjecture
can be found in \cite{conj3}; for $n \geq 4$ in many cases we can guess), which
gives upper bound for $\alpha_n$ and lower bound for $\tau_n$.
In this paper we present a new algorithm to bound $\alpha_n$ (from below)
and $\tau_n$ (from above). We study the behavior of this algorithm for
$\tau_2$ and $\alpha_n$ for $n=2,3,4$. Additionally we give new bounds
on multiple point Seshadri constants on $\PP^2$.

The paper is organized as follows: Next section is devoted to present
main facts used in the paper. Theorem \ref{maintool} is the main tool
and (as author believes) is interesting on its own. Other theorems
are (more or less) known. In section 3 we present two algorithms
(called {\sc NSsplit} and {\sc NSglue}). Both are used to show
non-specialty of a given system. The first one tries to split
system into many smaller (i.e. lower degree, less multiplicities) systems
(it is more probable to find a suitable criterion of non-specialty for a small system;
some of such criterions
are given in section 2). The second one uses Thm. \ref{maintool} to
,,glue'' multiplicities, which allows usage of birational isomorphism (Thm. \ref{cremona}).
In section 4 we present examples of bounds on $\tau_2$,
$\alpha_2$, $\alpha_3$ and $\alpha_4$. We focus on
quasi-homogeneous systems (for explanation see Thm. \ref{quasi}).
The bounds (especially for regularity)
are (as far as author knows) much better than bounds known before.
In the last section we show how our method, together with
Thm. \ref{eckl} (proposed by Eckl), can be used to produce bounds
on multiple point Seshadri constant on $\PP^2$.

\section{Main tools}

In \cite{CMirdeg} Ciliberto and Miranda proposed a method of
computing the dimension of a linear system in $\PP^2$ by splitting the problem
into several ones (possibly easier).
Similar method (based on splitting) was used by Biran (\cite{Biran}) to prove
ampleness and nefness of divisors on blow-ups of $\PP^2$.
Our theorem uses the same concept, but the proof is completely different
and works in any dimension ($n \geq 2$).

\begin{theorem}
\label{maintool}
Let $n \geq 2$, let $d,k,m_1,\dots,m_r,m_{r+1},\dots,m_s \in \N$.
If
\begin{itemize}
\item
$L_1 = \sys_n(k;m_{1},\dots,m_s)$ is non-special,
\item
$L_2 = \sys_n(d;m_{s+1},\dots,m_r,k+1)$ is non-special,
\item
$(\vdim L_1+1)(\vdim L_2+1) \geq 0$,
\end{itemize}
then the system $L=\sys_n(d;m_1,\dots,m_r)$ is non-special.
\end{theorem}

In the proof of Theorem \ref{maintool} we will use the 
reduction method introduced in \cite{md}.
We must adapt this method  to the $n$-dimensional case.
The following notations are used only in the proof of Theorem \ref{maintool}.

Every system $\sys_n(d;m_1,\dots,m_r)$ can be considered
as a vector space $V_n(d;m_1,\dots,m_r)$ (over $\field$) of
polynomials $f \in \field[x] = \field[x_1,\dots,x_n]$ of degree at most $d$, such that
$\mult_{p_j}f \geq m_j$, for $p_1,\dots,p_r \in \field^n$ in general position.
We have
$$\dim_{\field} V_n(d;m_1,\dots,m_r) = \dim \sys_n(d;m_1,\dots,m_r) + 1.$$
Let
\begin{align*}
\vdim V_n(d;m_1,\dots,m_r) & \rdf  \vdim \sys_n(d;m_1,\dots,m_r) + 1,\\
\edim V_n(d;m_1,\dots,m_r) & \rdf  \edim \sys_n(d;m_1,\dots,m_r) + 1.
\end{align*}
The space $V = V_n(d;m_1,\dots,m_r)$ is said to be {\dff non-special}
if $\dim_{\field} V=\edim V$, which is equivalent to the non-specialty of
$\sys_n(d;m_1,\dots,m_r)$. 
 
By $V_n(d)$ we understand the space of all polynomials of degree at most $d$
with no conditions imposed.
Let
$$\varphi_j : V_n(d) \to \field^{\binom{m_j+n-1}{n}}$$
be a linear
function which maps $f$ to a collection of all partial derivatives
of $f$ up to degree $m_j-1$ evaluated at $p_j$. Put
$$\Phi : V_n(d) \ni f \longmapsto (\varphi_1(f),\dots,\varphi_r(f)) \in \field^{c},$$
where
$$c = \sum_{j=1}^{r}\binom{m_j+n-1}{n}.$$
It is easy to observe that $V_n(d;m_1,\dots,m_r) = \ker \Phi$.
Let $M_n(d;m_1,\dots,m_r)$ be the matrix of $\Phi$ in monomial basis
of $V_n(d)$ and canonical basis of $\field^{c}$. 
Columns of $M_n(d;m_1,\dots,m_r)$ are indexed by monomials of degree
at most $d$, rows are indexed by conditions (i.e. points and symbols of
partial derivatives).
This matrix depends on
coordinates of points $p_1,\dots,p_r$, hence we consider each entry of
$M_n(d;m_1,\dots,m_r)$ as a polynomial (in fact, it is a monomial with coefficient) 
in $P \rdf \field[p_1^1,\dots,p_1^n,\dots,p_r^1,\dots,p_r^n]$,
where $p_j^k$ is the $k$-th coordinate of $j$-th point.
Observe that system $\sys_n(d;m_1,\dots,m_r)$ is non-special if and only if
the matrix $M_n(d;m_1,\dots,m_r)$ has maximal rank (as a matrix over $P$;
in other words, there exists a square submatrix $M$ of $M_n(d;m_1,\dots,m_r)$
of maximal size such that $\det M$ is a non-zero polynomial in $P$).

Define the multidegree function $\mdeg$ on $P$ by setting
$$\mdeg(p_j^k) \rdf (0,\dots,0,1,0,\dots,0) \in \N^n$$
with $1$ on $k$-th position.

\begin{lemma}
\label{ismono}
If $M$ is a square submatrix
of $M_n(d;m_1,\dots,m_r)$ then $\det M$ is multi-homogeneous with respect
to $\mdeg$. Moreover, 
$$\mdeg(\det M) = \beta - \gamma,$$
where $\beta$ is equal to the standard multidegree (in $\field[x_1,\dots,x_n]$) of 
product of monomials indexing columns of $M$, and
$\gamma$ depends only on conditions indexing rows of $M$.
\end{lemma}

\begin{proof}
The proof can be done easily by examining the entries of $M$, and
is esentially the same as the proof of Proposition 11 in \cite{md}.
\end{proof}

\begin{proof}[Proof of Theorem \ref{maintool}]
Consider the space $V_n(d;m_{s+1},\dots,m_r,k+1)$.
Since a translation $T : \field^n \to \field^n$ produces an isomorphism of
$V_n(d;m_{s+1},\dots,m_r,k+1)$ based on $p_{s+1},\dots,p_r,p_{r+1}$ with the same space based
on $T(p_{s+1}),\dots,T(p_r),T(p_{r+1})$, we may assume that 
$p_{r+1} = (0,\dots,0)$.
Now the space $V_n(d;m_{s+1},\dots,m_r,k+1)$ is nothing else than the space
$V_n(d,k;m_{s+1},\dots,m_r)$ of polynomials generated (over $\field$) by
the set
$$\{ x^\beta : k < \deg x^\beta \leq d \}$$
with respective multiplicities at $p_{s+1},\dots,p_{r}$.

Let us consider the matrix $M=M_n(d;m_1,\dots,m_r)$. If necessary, we may
reorder columns and rows of $M$ such that the first
$\binom{k+n}{n}$ columns are indexed by monomials of degree at most $k$,
first $\sum_{j=1}^{s}\binom{m_j+n-1}{n}$ rows indexed by
points $p_1,\dots,p_s$. Let $M_1$ be the left upper submatrix of $M$
consisting of $\binom{k+n}{n}$ columns and $\sum_{j=1}^{s}\binom{m_j+n-1}{n}$ rows,
$$
M = \left[ \begin{array}{c|c}
M_1 & K_1 \\ \hline
K_2 & M_2 \\ \end{array} \right].
$$
Observe that $M_1 = M_n(k;m_1,\dots,m_s)$,
$M_2 = M_n(d,k;m_{s+1},\dots,m_r)$, so both $M_1$ and $M_2$ have maximal rank.

Now, let us consider the case $\vdim L_1 \geq -1$ and $\vdim L_2 \geq -1$.
Then $\rows(M_1) = \rank(M_1) \leq \cols(M_1)$,
$\rows(M_2) = \rank(M_2) \leq \cols(M_2)$. Take
a square submatrix $M_1'$ of $M_1$ (resp. $M_2'$ of $M_2$) such that
$\rows(M_1') = \rows(M_1)$, $\det M_1' \neq 0$ (resp. 
$\rows(M_2') = \rows(M_2)$, $\det M_2' \neq 0$). Consider the matrix
$$
M' = \left[ \begin{array}{c|c}
M_1' & K_1' \\ \hline
K_2' & M_2' \\ \end{array} \right].
$$
Observe that $M'$ is a square submatrix of $M$ of size $\rows(M)$. It is enough
to show that $\det M' \neq 0$ to complete the proof.

Let $D_\ell$ be the set of monomials indexing columns of $M_\ell'$, $\ell=1,2$,
put $D \rdf D_1 \cup D_2$.
Let $U = [M_1' \mid K_1']$, $L = [K_2' \mid M_2']$ be submatrices
of $M'$, let 
$$\mathcal{C} \rdf \{ C \subset D \mid \#C = \#D_1 \}.$$
For $C \subset D$ define $U_{C}$ (resp. $L_{C}$)
as the submatrix of $U$ (resp. $L$) given by taking the columns indexed by
elements of $C$.
Now we can compute
$$\det M = \sum_{C \in \mathcal{C}} \epsilon(C) \det U_{C} \det L_{D \setminus C},$$
for $\epsilon(C) = \pm 1$.
Observe that 
\begin{align*}
\det U_{C} & \in \field[p_1^1,\dots,p_1^n,\dots,p_s^1,\dots,p_s^n], \\
\det L_{D \setminus C} & \in \field[p_{s+1}^1,\dots,p_{s+1}^n,\dots,p_{r}^1,\dots,p_{r}^n].
\end{align*}
We have
$$\det M' = \pm \det M_1' \det M_2' + f.$$
Assume that $\det M'=0$. As $\det M_1' \neq 0$ and $\det M_2' \neq 0$,
there exists $C \in \mathcal{C}$, $C \neq D_1$,
such that $\mdeg(\det U_C) = \mdeg(\det M_1')$.
By Lemma \ref{ismono}
$$
\mdeg(\det U_C) = \big(\sum_{x^\beta \in C} \beta \big) - \gamma, 
\quad
\mdeg(\det M_1') = \big(\sum_{x^\beta \in D_1} \beta \big) - \gamma,
$$
so
$$\sum_{x^\beta \in C} |\beta| = \sum_{x^\beta \in D_1} |\beta|.$$
However, for any $x^{\beta_1} \in D_1$ and $x^{\beta_2} \in D \setminus D_1$ we have
$|\beta_1| < |\beta_2|$, which leads to a contradiction.
The case 
$\vdim L_1 \leq -1$ and $\vdim L_2 \leq -1$ can be done similarly.
\end{proof}

\begin{theorem}
\label{cremona}
Let $n \geq 2$, let $d,m_1,\dots,m_r \in \N$, let
$k = (n-1)d-(m_1+\dots+m_{n+1})$.
If $m_j+k \geq 0$ for $j=1,\dots,n+1$ then
$$\dim \sys_n(d;m_1,\dots,m_r) = \dim \Cr(\sys_n(d;m_1,\dots,m_r)),$$
where
$$\Cr(\sys_n(d;m_1,\dots,m_r)) \rdf \sys_n(d+k;m_1+k,\dots,m_{n+1}+k,m_{n+2},\dots,m_{r}).$$
\end{theorem}

\begin{proof}
The proof for $n=2$ is well-known.
In \cite{conj3} we can find the proof for $n=3$
and that idea can be applied to arbitrary $n$. Namely, using projective
change of coordinates assume $p_1,\dots,p_{n+1}$ to be fundamental ones.
The system $\sys_n(d;m_1,\dots,m_r)$ is equivalent to the system of
hypersurfaces $\sys_n(M;m_{n+2},\dots,m_r)$ generated by the set
$$M = \{x_1^{\beta_1}\cdot \ldots \cdot x_{n+1}^{\beta_{n+1}} : |\beta| = d, \, \beta_j \leq d-m_j, \, j=1,\dots,n+1\}.$$
Similarly, system $\sys_n(d+k;m_1+k,\dots,m_{n+1}+k,m_{n+2},\dots,m_r)$
is equivalent to $\sys_n(M';m_{n+2},\dots,m_r)$ for
$$M' = \{x_1^{\beta_1}\cdot \ldots \cdot x_{n+1}^{\beta_{n+1}} : |\beta| = d+k, \, \beta_j \leq (d+k)-(m_j+k), \, j=1,\dots,n+1\}.$$
The standard birational transformation
$$\PP^n : (x_1:\ldots:x_{n+1}) \birto (x_1^{-1}:\ldots:x_{n+1}^{-1}) \in \PP^n$$
induces a bijection
$$ M \ni x_1^{\beta_1}\cdot\ldots\cdot x_{n+1}^{\beta_{n+1}} \longmapsto
x_1^{d-m_1-\beta_1}\cdot\ldots\cdot x_{n+1}^{d-m_{n+1}-\beta_{n+1}} \in M',$$
which completes the proof.
\end{proof}

Geometrically speaking, the system $\Cr(L)$ is an image of
$L$ by a birational transformation. Such transformation (in $\PP^2$)
is often reffered as to {\dff Cremona transformation}. In what follows
we use the name ``birational transformation'' to denote
$\Cr$ operation.

\begin{theorem}
\label{hyper}
Let $n \geq 2$, let $d,m_1,\dots,m_r \in \Z$.
If $(n-1)d-\sum_{j=1}^{n}m_j < 0$, $m_j > 0$ for $j=1,\dots,n$ then
$$\dim \sys_n(d;m_1,\dots,m_r) = \dim \Hp(\sys_n(d;m_1,\dots,m_r)),$$
where
$$\Hp(\sys_n(d;m_1,\dots,m_r)) \rdf \sys_n(d-1;m_1-1,\dots,m_n-1,m_{n+1},\dots,m_r).$$
\end{theorem}

\begin{proof}
The proof for $n=3$ can be found in \cite{conj3}, we use the same reasoning.
Consider the situation as in the previous proof. Since
$\beta_1+\dots+\beta_n \leq nd-(m_1+\dots+m_n) < d$, the map
$$
M \ni x_1^{\beta_1}\cdot\ldots\cdot x_{n+1}^{\beta_{n+1}} \longmapsto
x_1^{\beta_1}\cdot\ldots\cdot x_n^{\beta_n} \cdot x_{n+1}^{\beta_{n+1}-1}$$
is one-to-one and its image corresponds to the system
$\sys_n(d-1;m_1-1,\dots,m_n-1,m_{n+1},\dots,m_r)$.
Geometrically speaking, the hyperplane passing through the first $n$ points is in the
base locus of $\sys_n(d;m_1,\dots,m_r)$ and can be ``taken out''.
\end{proof}

\begin{definition}
We say that {\dff $\sys_n(d;m_1,\dots,m_r)$ is in standard form}
if $d < 0$ or the following holds:
\begin{itemize}
\item
$m_1,\dots,m_r$ are non-increasing,
\item
$(n-1)d - \sum_{j=1}^{n+1} m_j \geq 0$.
\end{itemize}
\end{definition}

\begin{proposition}
\label{standard}
Each system $L$ has its standard form $L'$ such that $\dim L = \dim L'$.
\end{proposition}

\begin{proof}
The following set of instructions allows to transform any system $L$ into
a standard form.
\begin{itemize}
\item
sort multiplicities in non-increasing order,
\item
if $(n-1)d-\sum_{j=1}^n m_j < 0$ then take $L \ass \Hp(L)$ and go back
to the first step,
\item
if $(n-1)d-\sum_{j=1}^{n+1} m_j < 0$ then take $L \ass \Cr(L)$ and go back
to the first step.
\end{itemize}
Observe that each time (in steps 2 and 3) the degree of $L$ decreases
and the dimension does not change.
\end{proof}

The end of this section is devoted to prepare a family of non-special
systems (which will be used later) by presenting criterions of non-specialty.
We will use the following notation: $m^{\times k}$ denotes the
sequence of $m$'s taken $k$ times,
$$m^{\times k} = \underbrace{(m,\dots,m)}_{k}.$$

\begin{criterion}
\label{crendcr}
Let $\sys_n(d;m_1,\dots,m_r)$ be a standard form of a system $L$.
If either $d < 0$ or $m_2 \leq 1$ then $L$ is non-special.
\end{criterion}

\begin{proof}
The case $d<0$ is obvious. If $m_2 \leq 1$ and $\sys_n(d;m_1,\dots,m_r)$
is in standard form then $m_j \leq 1$ for $j=2,\dots,r$. Each point
of muliplicity one (or zero) imposes an independent condition.
\end{proof}

\begin{remark}
\label{negcremona}
It is known that we can perform $\Cr$ operation on the system with
$(n-1)d-\sum_{j=1}^{n} m_j < 0$. This leads to negative multiplicities, which
can be properly understood using divisors on the blow-up of $\PP^n$
in $r$ points. However,
it is easy to show (using Thms. \ref{cremona} and \ref{hyper}) that
if $L=\sys_n(d;m_1,\dots,m_r)$ can be transformed (by birational transformations) to
$\sys_n(d';m_1',\dots,m_r')$, for $m_1',\dots,m_r' \in \Z$ and $d < 0$, then
$L$ is empty.
\end{remark}

\begin{criterion}
\label{reductionplus}
Let $L=\sys_2(d;m_1,\dots,m_r)$. If $\vdim L \geq -1$ and $L$ can be 
reduced to a system $L'$ (using reductions decribed in \cite{mdwj}; the system
$L'$ can be generated by some finite set of monomials) such that
$\vdim L'=\vdim L$ and $L'$ has no conditions (multiplicities imposed)
then $L$ is non-special.
\end{criterion}

\begin{proof}
By Thm. 17 in \cite{mdwj} $\dim L \leq \dim L'$. Since $L'$ is condition-free,
we have $\dim L'=\vdim L'=\vdim L$.
\end{proof}

\begin{criterion}
\label{reductionminus}
Let $L=\sys_2(d;m_1,\dots,m_r)$. If $\vdim L \leq -1$ and $L$ can be 
reduced to a system $L'$ (using weak reductions decribed in \cite{mdwj}) such that
$L'$ has no monomials, then $L$ is empty (and hence non-special).
\end{criterion}

\begin{proof}
Again, by Thm. 17 in \cite{mdwj} we have
$\dim L \leq \dim L' = -1$.
\end{proof}

\begin{criterion}
\label{4points}
Let $m \geq 0$. The systems $\sys_2(2m-1;m^{\times 4})$ and $\sys_2(2m;m^{\times 4})$ are non-special.
\end{criterion}

\begin{proof}
Use birational transformation.
\end{proof}

\begin{criterion}
\label{squarepoints}
Let $r \geq 3$, let $m > (r-2)/4$. Then
$\tau_2(m^{\times r^2}) = rm+\lceil(r-3)/2\rceil$.
\end{criterion}

\begin{proof}
See \cite{tausquare}, Lemma 5.3.
\end{proof}

\section{Algorithms}

To bound $\alpha_n$ and $\tau_n$ it is sufficient to find
an algorithm, which, given $n,d,m_1,\dots,m_r$, returns
either {\sc non-special} (and then the system $\sys_n(d;m_1,\dots,m_r)$
is non-special), or {\sc not-decided}. We will focus on such algorithms.
The answer {\sc special} is also allowed (of course only
if $\sys_n(d;m_1,\dots,m_r)$ is special).

Let us assume we are given a collection (called {\sc Eclass}) of non-special systems.
We will assume that if $L \in \text{{\sc Eclass}}$ then
also $\Cr(L)$ and $\Hp(L)$ are in {\sc Eclass}. The first algorithm
({\sc NSsplit}) makes use of Theorem \ref{maintool} and tries to
show non-specialty of a given system by splitting it into
systems belonging to {\sc Eclass}.

\begin{algorithm}{{\sc NSsplit}}

\noindent 
\begin{tabular}{rl}
{\bf Input:} & $n \geq 2$, $d,m_1,\dots,m_r \in \N$ \\
{\bf Output:} & {\sc non-special} or {\sc not-decided} or {\sc special}.
\end{tabular}
\\

\noindent
$v_1 \ass \vdim \sys_n(d;m_1,\dots,m_r)$;\\
change $\sys_n(d;m_1,\dots,m_r)$ into standard form;\\
{\bf if} $\vdim \sys_n(d;m_1,\dots,m_r) \geq \max \{0, v_1\}$ {\bf then} {\bf return} {\sc special};\\
{\bf if} $\sys_n(d;m_1,\dots,m_r) \in \text{{\sc Eclass}}$ {\bf then} {\bf return} {\sc non-special};\\
{\bf for each} $\{i_1,\dots,i_s\} \subset \{1,\dots,r\}$ and $0 \leq k < d$ {\bf do}\\
\hspace*{0.5cm} $\{j_1,\dots,j_{r-s}\} \ass \{1,\dots,r\} \setminus \{i_1,\dots,i_s\}$;\\
\hspace*{0.5cm} $v \ass (\vdim \sys_n(k;m_{i_1},\dots,m_{i_s})+1)(\vdim \sys_n(d;k+1,m_{j_1},\dots,m_{j_{r-s}})+1)$;\\
\hspace*{0.5cm} $\text{{\sc a1}} \ass \text{{\sc NSsplit}}(n,k,m_{i_1},\dots,m_{i_s})$;\\
\hspace*{0.5cm} $\text{{\sc a2}} \ass \text{{\sc NSsplit}}(n,d,k+1,m_{j_1},\dots,m_{j_{r-s}})$;\\
\hspace*{0.5cm} {\bf if} $v \geq 0$ {\bf and} $\text{{\sc a1}} = \text{{\sc non-special}}$ {\bf and} $\text{{\sc a2}} = \text{{\sc non-special}}$ {\bf then} {\bf return} {\sc non-special};\\
{\bf end for each}\\
{\bf return} {\sc not-decided};
\end{algorithm}

The above algorithm can be implemented in a more subtle way, for example
for $v < 0$ we can skip the next two steps (running {\sc NSsplit}).
Also one can use criterions of specialty to avoid running {\sc NSsplit}
on special systems.

\begin{example}
We will show that $\sys_2(72;10^{\times 50})$ is non-special (i.e. empty).
Applying {\sc NSsplit} with
$$\text{{\sc Eclass}} \rdf \{ L=\sys_2(d;m_1,\dots,m_r) : L \text{ is empty by Crit. \ref{crendcr}}\}$$
we can find the following solution (presented in 5 steps):
\begin{enumerate}
\item
split $\sys_2(72;10^{\times 50})$ into $\sys_2(72;61,10^{\times 15})$, which
is empty by Crit. \ref{crendcr}, and $\sys_2(60;10^{\times 35})$;
\item
split $\sys_2(60;10^{\times 35})$ into $\sys_2(45;10^{\times 20})$ and
$\sys_2(60;46,10^{\times 15})$, which, by Cremona transformation,
is equivalent to the system $\sys_2(18;10,4^{\times 15})$;
\item
split $\sys_2(18;10,4^{\times 15})$ into $\sys_2(10;4^{\times 7})$ and
$\sys_2(18;11,10,4^{\times 8})$; both systems are empty by Crit. \ref{crendcr},
\item
split $\sys_2(45;10^{\times 20})$ into $\sys_2(28;10^{\times 8})$ and
$\sys_2(45;29,10^{\times 12})$; the first system is empty by Crit. \ref{crendcr},
the second one can be transformed into $\sys_2(21;6^{\times 15},5)$;
\item
split $\sys_2(21;6^{\times 15},5)$ into $\sys_2(11;6^{\times 4})$ and
$\sys_2(21;12,6^{\times 8},5)$; both are empty by Crit. \ref{crendcr}.
\end{enumerate}
\end{example}

Algorithm {\sc NSsplit} is very slow for systems with many multiplicities.
Therefore we will present another algorithm, which does not search through
the tree of all possibilities, but tries to ``glue'' some multiplicities.
The aim is to obtain multiplicities relatively big with respect to the degree,
i.e. $(n-1)d < \sum_{j=1}^{n+1} m_j$, and use birational transformation(s).

\begin{theorem}
\label{glue}
Let $\sys_n(k;m^s)$ be non-special,
let 
\begin{align*}
L_1 & = \sys_n(d;m_1,\dots,m_r,m^s),\\
L_2 & = \sys_n(d;m_1,\dots,m_r,k+1).
\end{align*}
If either $-1 \leq \vdim L_2 \leq \vdim L_1$ or $\vdim L_1 \leq \vdim L_2 \leq -1$
then in order to show non-specialty of $L_1$ it is enough to show non-specialty
of $L_2$.
\end{theorem}

\begin{proof}
Follows from Theorem \ref{maintool}.
\end{proof}

The algorithm {\sc NSglue} depends on two 
sets ({\sc Eclass} and {\sc Gclass}) of non-special systems.
The set {\sc Gclass} containing systems of the form $\sys_n(k;m^{\times s})$
should be relatively small
and should contain systems with virtual dimensions close to $-1$.

\begin{algorithm}{{\sc NSglue}}

\noindent 
\begin{tabular}{rl}
{\bf Input:} & $n \geq 2$, $d,m_1,\dots,m_r \in \N$ \\
{\bf Output:} & {\sc non-special} or {\sc not-decided} or {\sc special}. 
\end{tabular}
\\

\medskip

\noindent
$v_1 \ass \vdim \sys_n(d;m_1,\dots,m_r)$;\\
change $\sys_n(d;m_1,\dots,m_r)$ into standard form;\\
{\bf if} $\vdim \sys_n(d;m_1,\dots,m_r) \geq \max \{0, v_1\}$ {\bf then} {\bf return} {\sc special};\\
{\bf if} $\sys_n(d;m_1,\dots,m_r) \in \text{{\sc Eclass}}$ {\bf then} {\bf return} {\sc non-special};\\
$v_2 \ass \vdim \sys_n(d;m_1,\dots,m_r)$;\\
{\bf for each} $\sys_n(k;m^{\times s}) \in \text{{\sc Gclass}}$ such that
$m_{i_1}=\dots=m_{i_s}=m$ for $\{i_1,\dots,i_s\} \subset \{1,\dots,r\}$ {\bf do}\\
\hspace*{0.5cm} $\{j_1,\dots,j_{r-s}\} \ass \{1,\dots,r\} \setminus \{i_1,\dots,i_s\}$;\\
\hspace*{0.5cm} $v_3 \ass \vdim \sys_n(d;k+1,m_{j_1},\dots,m_{j_{r-s}})$;\\
\hspace*{0.5cm} $\text{{\sc a}} \ass \text{{\sc NSglue}}(n,d,k+1,m_{j_1},\dots,m_{j_{r-s}})$;\\
\hspace*{0.5cm} {\bf if} ($-1 \leq v_3 \leq v_2$ {\bf or} $v_2 \leq v_3 \leq -1$) {\bf and} $\text{{\sc a}} = \text{{\sc non-special}}$ {\bf then} {\bf return} {\sc non-special};\\
{\bf end for each}\\
{\bf return} {\sc not-decided};
\end{algorithm}

Observe that
$\vdim \sys_n(d;m_1,\dots,m_r,k+1) = \vdim \sys_n(d;m_1,\dots,m_r,m^{\times s})
- \vdim \sys_n(k;m^{\times s})$.
Therefore in order to show non-specialty of a system with non-negative (resp. non-positive)
virtual dimension we must use ``glueing'' systems (from {\sc Gclass}) with
non-negative (resp. non-positive) virtual dimension.
In both cases the virtual dimension after ``glueing'' will go closer
to $-1$.

\begin{example}
Let us show how {\sc NSglue} works with the system
$\sys_2(5918;4000,1000^{\times 19})$. Put
\begin{align*}
\text{{\sc Eclass}} & \rdf \{ L : L \text{ is non-special by Crit. \ref{crendcr} or \ref{reductionplus}} \},\\
\text{{\sc Gclass}} & \rdf \{ \sys_2(2m;m^{\times 4}) \}.
\end{align*}
By Crit. \ref{4points} the set {\sc Gclass} contains non-special systems
of dimension $m+1$.

The standard form of
$\sys_2(5918;4000,1000^{\times 19})$ is
$\sys_2(5180;3262,1000,918^{\times 18})$. Now we glue the four points
($918^{\times 4} \longrightarrow 2\cdot 918+1=1837$) obtaining $\sys_2(5180;3262,1837,1000,918^{\times 14})$, which can
be transformed into its standard form
$\sys_2(4261;2343,918^{\times 15},81)$.
Again, we glue the points to consider the system
$\sys_2(4261;2343,1837,918^{\times 11},81)$.
This can be transformed into
$\sys_2(3424;1506,1000,918^{\times 10},81^{\times 2})$.
The last glueing gives
$\sys_2(3424;1837,1506,1000,918^{\times 6},81^{\times 2})$.
The standard form is $\sys_2(112;22,9^{\times 3},4^{\times 7})$.
By Crit. \ref{reductionplus} the last system is non-special.
Alternatively, we can show that $\tau_2(22,9^{\times 3},4^{\times 7}) \geq 31$.
Glueing $4$ points of multiplicity $4$ we get the system
$\sys_2(31;22,9^{\times 4},4^{\times 3})$. The standard form is equal to
$\sys_2(13;4^{\times 4})$, which is non-special.
\end{example}

\begin{example}
\label{197}
Let us consider the system $L=\sys_2(d;m^{\times 21},h)$. We will show
that if $197d-42(21m+h)<0$ then $L$ is empty.
By Thm. \ref{glue} it is enough to show that
$L'=\sys_2(d;(2m)^{\times 4},m^{\times 5},h)$ is empty. The system $L'$
has $10$ base points, say $p_1,\dots,p_4$ with multiplicity $2m$,
$p_5,\dots,p_9$ with multiplicity $m$ and $p_{10}$ with multiplicity $h$.
We perform birational transformations based on triples of points
(allowing negative multiplicities, see Rem. \ref{negcremona}), according to
the following sequence (of numbers of points):
$(1,2,3)$, $(4,5,6)$, $(4,7,8)$, $(4,9,10)$, $(1,2,3)$, $(5,6,7)$,
$(8,9,10)$, $(5,6,7)$, $(1,2,3)$, $(4,5,6)$, $(4,7,8)$, $(4,9,10)$, $(1,2,3)$.
After all computations (done by hand or with Singular procedure 
available at \cite{MYWWW}) we get the system
$\sys_2(197d-882m-42h;(84d-376m-18h)^{\times 4},(42d-188m-9h)^{\times 5},42d-189m-8h)$. 
This result allows to bound the Seshadri constant of $\stbun$
for $22$ points in general position (see Prop. \ref{22}).
\end{example}

\section{Results for $n=2,3,4$}

\subsection{Case $n=2$}

For $n=2$ many algorithms bounding $\alpha_2$ and $\tau_2$ are known (see e.g.
\cite{harsesh}, \cite{roe1} and \cite{roe2} for Ro\'e's unloading method, 
\cite{harroe} for modified unloading, \cite{M}, and finally \cite{surv} for
survey through many algorithms). We will use algorithm
{\sc NSglue} with
\begin{align*}
\text{{\sc Eclass}} & = \{ L : L \text{ is non-special by Criterion \ref{crendcr}, \ref{reductionplus}} \text{ or } \ref{reductionminus}\},\\
\text{{\sc Gclass}} & = \{ \sys_2(2m;m^{\times 4}) \} \cup \{ \sys_2((2s+1)m+s-1;m^{\times (2s+1)^2}), \, s=1,\dots,10 \}.
\end{align*}
By Crit. \ref{squarepoints} the set {\sc Gclass} contains non-special systems
(this criterion does not work for all $m$, but even in the worst case
it works for $m \geq 5$; for $m \leq 4$ everything is known, see e.g. \cite{CMir}).
The implementation of {\sc NSglue} with the above classes can be found and
downloaded from \cite{MYWWW}. The work is done in FreePascal.

For a very large family of quasi-homogeneous systems our algorithm
gives better bounds than any other one. 
In \cite{harroe2} we can find a list of homogeneous systems for $r=10,\dots,99$
being not a square (one system for each $r$). They are conjectured to be empty,
but the authors of \cite{harroe2} were unable to show this using known
algorithms. For $64$ of these systems the algorithm {\sc NSglue} gives the
expected answer, for $18$ we must use {\sc NSsplit}. Only for $r=10,11$
our algorithms do not work.

We will present only few examples,
first two of them appearing in \cite{harroe2},
four of them proposed in \cite{M}. 
In the Tab. \ref{tab1} one can find results
the auhor was able to find using
known algorithms (for convenience, we present the difference
between bound and conjectured value). For $\alpha_2$ they were all obtained by
one of Harbourne/Ro\'e's algorithms (using \cite{harwww} or by
author own implementation for quasi-homogeneous case; the author is not sure that they are the best),
for $\tau_2$ the bounds were either computed with \cite{harwww} (homogeneous case), 
or obtained by Monserrat (quasi-homogeneous case; \cite{M}), who checked
other algorithms and claimed to find the best bounds.

\begin{table}[ht!]
$$
\begin{array}{|c|c|c|c|c|c|}
\hline
\text{system} & \text{conj. } \alpha_2=\tau_2 & \text{Har--Ro\'e } \alpha_2 & \text{Mon/Har } \tau_2 & \text{{\sc NSglue} } \alpha_2 & \text{{\sc NSglue} } \tau_2 \\
\hline
(24^{\times 12})        & 84   & -1   & +1    & 0  & +1 \\
\hline
(173^{\times 96})       & 1699 & -19  & +4     & -3 & +1 \\
\hline
(4000,1000^{\times 19}) & 5917 & -255 & +92 & -1 & +1 \\
\hline
(6000,1500^{\times 19}) & 8875 & -382 & +138 & -1 & +1 \\
\hline
(500^{\times 1000}) & 15826 & -25 & +16  & -3 & +0 \\
\hline
(1200^{\times 1000}) & 37962 & -41 & +17 & -6 & +0 \\
\hline
\end{array}
$$
\caption{Examples of bounds for $\alpha_2$ and $\tau_2$}\label{tab1}
\end{table}

Following Harbourne (\cite{surv}) we present the results of
{\sc NSglue} for homogeneous systems with $1 \leq m \leq 150$,
$10 \leq r \leq 450$. On the graph presented on Fig. \ref{pic1} each dot
denotes the success --- the bound of $\tau_2$ is equal to its conjectured value,
the same rule is used to present results for $\alpha_2$ on Fig. \ref{pic2}.
The graphs presenting partial successes (the difference is not greater that some
constant) and those situations, when bound for $\alpha_2$ is not less
that conjectured by Nagata (see the last section) can be found at \cite{MYWWW}.

\begin{figure}[ht!]
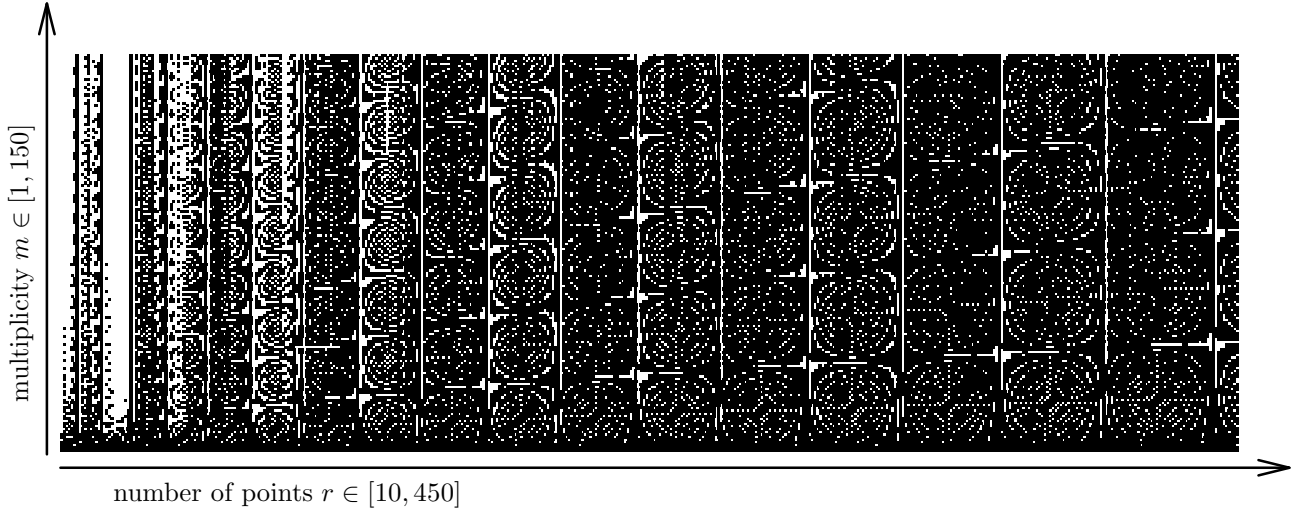

\centertexdraw{
\drawdim pt
\move(10 1)
\lvec(451 1)
\lvec(451 2)
\lvec(10 2)
\ifill f:0
\move(10 2)
\lvec(451 2)
\lvec(451 3)
\lvec(10 3)
\ifill f:0
\move(10 3)
\lvec(15 3)
\lvec(15 4)
\lvec(10 4)
\ifill f:0
\move(16 3)
\lvec(42 3)
\lvec(42 4)
\lvec(16 4)
\ifill f:0
\move(43 3)
\lvec(54 3)
\lvec(54 4)
\lvec(43 4)
\ifill f:0
\move(55 3)
\lvec(99 3)
\lvec(99 4)
\lvec(55 4)
\ifill f:0
\move(100 3)
\lvec(117 3)
\lvec(117 4)
\lvec(100 4)
\ifill f:0
\move(118 3)
\lvec(180 3)
\lvec(180 4)
\lvec(118 4)
\ifill f:0
\move(181 3)
\lvec(204 3)
\lvec(204 4)
\lvec(181 4)
\ifill f:0
\move(205 3)
\lvec(285 3)
\lvec(285 4)
\lvec(205 4)
\ifill f:0
\move(286 3)
\lvec(315 3)
\lvec(315 4)
\lvec(286 4)
\ifill f:0
\move(316 3)
\lvec(414 3)
\lvec(414 4)
\lvec(316 4)
\ifill f:0
\move(415 3)
\lvec(450 3)
\lvec(450 4)
\lvec(415 4)
\ifill f:0
\move(10 4)
\lvec(17 4)
\lvec(17 5)
\lvec(10 5)
\ifill f:0
\move(18 4)
\lvec(23 4)
\lvec(23 5)
\lvec(18 5)
\ifill f:0
\move(24 4)
\lvec(35 4)
\lvec(35 5)
\lvec(24 5)
\ifill f:0
\move(36 4)
\lvec(56 4)
\lvec(56 5)
\lvec(36 5)
\ifill f:0
\move(57 4)
\lvec(70 4)
\lvec(70 5)
\lvec(57 5)
\ifill f:0
\move(71 4)
\lvec(74 4)
\lvec(74 5)
\lvec(71 5)
\ifill f:0
\move(75 4)
\lvec(86 4)
\lvec(86 5)
\lvec(75 5)
\ifill f:0
\move(87 4)
\lvec(90 4)
\lvec(90 5)
\lvec(87 5)
\ifill f:0
\move(91 4)
\lvec(108 4)
\lvec(108 5)
\lvec(91 5)
\ifill f:0
\move(109 4)
\lvec(143 4)
\lvec(143 5)
\lvec(109 5)
\ifill f:0
\move(144 4)
\lvec(165 4)
\lvec(165 5)
\lvec(144 5)
\ifill f:0
\move(166 4)
\lvec(171 4)
\lvec(171 5)
\lvec(166 5)
\ifill f:0
\move(172 4)
\lvec(189 4)
\lvec(189 5)
\lvec(172 5)
\ifill f:0
\move(190 4)
\lvec(195 4)
\lvec(195 5)
\lvec(190 5)
\ifill f:0
\move(196 4)
\lvec(221 4)
\lvec(221 5)
\lvec(196 5)
\ifill f:0
\move(222 4)
\lvec(270 4)
\lvec(270 5)
\lvec(222 5)
\ifill f:0
\move(271 4)
\lvec(300 4)
\lvec(300 5)
\lvec(271 5)
\ifill f:0
\move(301 4)
\lvec(308 4)
\lvec(308 5)
\lvec(301 5)
\ifill f:0
\move(309 4)
\lvec(332 4)
\lvec(332 5)
\lvec(309 5)
\ifill f:0
\move(333 4)
\lvec(340 4)
\lvec(340 5)
\lvec(333 5)
\ifill f:0
\move(341 4)
\lvec(374 4)
\lvec(374 5)
\lvec(341 5)
\ifill f:0
\move(375 4)
\lvec(437 4)
\lvec(437 5)
\lvec(375 5)
\ifill f:0
\move(438 4)
\lvec(451 4)
\lvec(451 5)
\lvec(438 5)
\ifill f:0
\move(11 5)
\lvec(27 5)
\lvec(27 6)
\lvec(11 6)
\ifill f:0
\move(28 5)
\lvec(33 5)
\lvec(33 6)
\lvec(28 6)
\ifill f:0
\move(34 5)
\lvec(37 5)
\lvec(37 6)
\lvec(34 6)
\ifill f:0
\move(38 5)
\lvec(57 5)
\lvec(57 6)
\lvec(38 6)
\ifill f:0
\move(58 5)
\lvec(60 5)
\lvec(60 6)
\lvec(58 6)
\ifill f:0
\move(61 5)
\lvec(63 5)
\lvec(63 6)
\lvec(61 6)
\ifill f:0
\move(64 5)
\lvec(72 5)
\lvec(72 6)
\lvec(64 6)
\ifill f:0
\move(73 5)
\lvec(75 5)
\lvec(75 6)
\lvec(73 6)
\ifill f:0
\move(76 5)
\lvec(78 5)
\lvec(78 6)
\lvec(76 6)
\ifill f:0
\move(79 5)
\lvec(88 5)
\lvec(88 6)
\lvec(79 6)
\ifill f:0
\move(89 5)
\lvec(95 5)
\lvec(95 6)
\lvec(89 6)
\ifill f:0
\move(96 5)
\lvec(106 5)
\lvec(106 6)
\lvec(96 6)
\ifill f:0
\move(107 5)
\lvec(110 5)
\lvec(110 6)
\lvec(107 6)
\ifill f:0
\move(111 5)
\lvec(114 5)
\lvec(114 6)
\lvec(111 6)
\ifill f:0
\move(115 5)
\lvec(126 5)
\lvec(126 6)
\lvec(115 6)
\ifill f:0
\move(127 5)
\lvec(130 5)
\lvec(130 6)
\lvec(127 6)
\ifill f:0
\move(131 5)
\lvec(134 5)
\lvec(134 6)
\lvec(131 6)
\ifill f:0
\move(135 5)
\lvec(147 5)
\lvec(147 6)
\lvec(135 6)
\ifill f:0
\move(148 5)
\lvec(156 5)
\lvec(156 6)
\lvec(148 6)
\ifill f:0
\move(157 5)
\lvec(170 5)
\lvec(170 6)
\lvec(157 6)
\ifill f:0
\move(171 5)
\lvec(175 5)
\lvec(175 6)
\lvec(171 6)
\ifill f:0
\move(176 5)
\lvec(180 5)
\lvec(180 6)
\lvec(176 6)
\ifill f:0
\move(181 5)
\lvec(195 5)
\lvec(195 6)
\lvec(181 6)
\ifill f:0
\move(196 5)
\lvec(200 5)
\lvec(200 6)
\lvec(196 6)
\ifill f:0
\move(201 5)
\lvec(205 5)
\lvec(205 6)
\lvec(201 6)
\ifill f:0
\move(206 5)
\lvec(221 5)
\lvec(221 6)
\lvec(206 6)
\ifill f:0
\move(222 5)
\lvec(232 5)
\lvec(232 6)
\lvec(222 6)
\ifill f:0
\move(233 5)
\lvec(249 5)
\lvec(249 6)
\lvec(233 6)
\ifill f:0
\move(250 5)
\lvec(255 5)
\lvec(255 6)
\lvec(250 6)
\ifill f:0
\move(256 5)
\lvec(261 5)
\lvec(261 6)
\lvec(256 6)
\ifill f:0
\move(262 5)
\lvec(279 5)
\lvec(279 6)
\lvec(262 6)
\ifill f:0
\move(280 5)
\lvec(285 5)
\lvec(285 6)
\lvec(280 6)
\ifill f:0
\move(286 5)
\lvec(291 5)
\lvec(291 6)
\lvec(286 6)
\ifill f:0
\move(292 5)
\lvec(310 5)
\lvec(310 6)
\lvec(292 6)
\ifill f:0
\move(311 5)
\lvec(323 5)
\lvec(323 6)
\lvec(311 6)
\ifill f:0
\move(324 5)
\lvec(343 5)
\lvec(343 6)
\lvec(324 6)
\ifill f:0
\move(344 5)
\lvec(350 5)
\lvec(350 6)
\lvec(344 6)
\ifill f:0
\move(351 5)
\lvec(357 5)
\lvec(357 6)
\lvec(351 6)
\ifill f:0
\move(358 5)
\lvec(378 5)
\lvec(378 6)
\lvec(358 6)
\ifill f:0
\move(379 5)
\lvec(385 5)
\lvec(385 6)
\lvec(379 6)
\ifill f:0
\move(386 5)
\lvec(392 5)
\lvec(392 6)
\lvec(386 6)
\ifill f:0
\move(393 5)
\lvec(414 5)
\lvec(414 6)
\lvec(393 6)
\ifill f:0
\move(415 5)
\lvec(429 5)
\lvec(429 6)
\lvec(415 6)
\ifill f:0
\move(430 5)
\lvec(451 5)
\lvec(451 6)
\lvec(430 6)
\ifill f:0
\move(10 6)
\lvec(12 6)
\lvec(12 7)
\lvec(10 7)
\ifill f:0
\move(15 6)
\lvec(22 6)
\lvec(22 7)
\lvec(15 7)
\ifill f:0
\move(23 6)
\lvec(28 6)
\lvec(28 7)
\lvec(23 7)
\ifill f:0
\move(29 6)
\lvec(33 6)
\lvec(33 7)
\lvec(29 7)
\ifill f:0
\move(34 6)
\lvec(37 6)
\lvec(37 7)
\lvec(34 7)
\ifill f:0
\move(38 6)
\lvec(39 6)
\lvec(39 7)
\lvec(38 7)
\ifill f:0
\move(40 6)
\lvec(45 6)
\lvec(45 7)
\lvec(40 7)
\ifill f:0
\move(46 6)
\lvec(47 6)
\lvec(47 7)
\lvec(46 7)
\ifill f:0
\move(48 6)
\lvec(51 6)
\lvec(51 7)
\lvec(48 7)
\ifill f:0
\move(52 6)
\lvec(63 6)
\lvec(63 7)
\lvec(52 7)
\ifill f:0
\move(64 6)
\lvec(68 6)
\lvec(68 7)
\lvec(64 7)
\ifill f:0
\move(69 6)
\lvec(73 6)
\lvec(73 7)
\lvec(69 7)
\ifill f:0
\move(74 6)
\lvec(84 6)
\lvec(84 7)
\lvec(74 7)
\ifill f:0
\move(85 6)
\lvec(87 6)
\lvec(87 7)
\lvec(85 7)
\ifill f:0
\move(88 6)
\lvec(90 6)
\lvec(90 7)
\lvec(88 7)
\ifill f:0
\move(91 6)
\lvec(99 6)
\lvec(99 7)
\lvec(91 7)
\ifill f:0
\move(100 6)
\lvec(102 6)
\lvec(102 7)
\lvec(100 7)
\ifill f:0
\move(103 6)
\lvec(105 6)
\lvec(105 7)
\lvec(103 7)
\ifill f:0
\move(106 6)
\lvec(118 6)
\lvec(118 7)
\lvec(106 7)
\ifill f:0
\move(119 6)
\lvec(125 6)
\lvec(125 7)
\lvec(119 7)
\ifill f:0
\move(126 6)
\lvec(132 6)
\lvec(132 7)
\lvec(126 7)
\ifill f:0
\move(133 6)
\lvec(139 6)
\lvec(139 7)
\lvec(133 7)
\ifill f:0
\move(140 6)
\lvec(150 6)
\lvec(150 7)
\lvec(140 7)
\ifill f:0
\move(151 6)
\lvec(154 6)
\lvec(154 7)
\lvec(151 7)
\ifill f:0
\move(155 6)
\lvec(158 6)
\lvec(158 7)
\lvec(155 7)
\ifill f:0
\move(159 6)
\lvec(162 6)
\lvec(162 7)
\lvec(159 7)
\ifill f:0
\move(163 6)
\lvec(174 6)
\lvec(174 7)
\lvec(163 7)
\ifill f:0
\move(175 6)
\lvec(178 6)
\lvec(178 7)
\lvec(175 7)
\ifill f:0
\move(179 6)
\lvec(182 6)
\lvec(182 7)
\lvec(179 7)
\ifill f:0
\move(183 6)
\lvec(186 6)
\lvec(186 7)
\lvec(183 7)
\ifill f:0
\move(187 6)
\lvec(199 6)
\lvec(199 7)
\lvec(187 7)
\ifill f:0
\move(200 6)
\lvec(208 6)
\lvec(208 7)
\lvec(200 7)
\ifill f:0
\move(209 6)
\lvec(217 6)
\lvec(217 7)
\lvec(209 7)
\ifill f:0
\move(218 6)
\lvec(226 6)
\lvec(226 7)
\lvec(218 7)
\ifill f:0
\move(227 6)
\lvec(240 6)
\lvec(240 7)
\lvec(227 7)
\ifill f:0
\move(241 6)
\lvec(245 6)
\lvec(245 7)
\lvec(241 7)
\ifill f:0
\move(246 6)
\lvec(250 6)
\lvec(250 7)
\lvec(246 7)
\ifill f:0
\move(251 6)
\lvec(255 6)
\lvec(255 7)
\lvec(251 7)
\ifill f:0
\move(256 6)
\lvec(270 6)
\lvec(270 7)
\lvec(256 7)
\ifill f:0
\move(271 6)
\lvec(275 6)
\lvec(275 7)
\lvec(271 7)
\ifill f:0
\move(276 6)
\lvec(280 6)
\lvec(280 7)
\lvec(276 7)
\ifill f:0
\move(281 6)
\lvec(285 6)
\lvec(285 7)
\lvec(281 7)
\ifill f:0
\move(286 6)
\lvec(301 6)
\lvec(301 7)
\lvec(286 7)
\ifill f:0
\move(302 6)
\lvec(312 6)
\lvec(312 7)
\lvec(302 7)
\ifill f:0
\move(313 6)
\lvec(323 6)
\lvec(323 7)
\lvec(313 7)
\ifill f:0
\move(324 6)
\lvec(334 6)
\lvec(334 7)
\lvec(324 7)
\ifill f:0
\move(335 6)
\lvec(351 6)
\lvec(351 7)
\lvec(335 7)
\ifill f:0
\move(352 6)
\lvec(357 6)
\lvec(357 7)
\lvec(352 7)
\ifill f:0
\move(358 6)
\lvec(363 6)
\lvec(363 7)
\lvec(358 7)
\ifill f:0
\move(364 6)
\lvec(369 6)
\lvec(369 7)
\lvec(364 7)
\ifill f:0
\move(370 6)
\lvec(387 6)
\lvec(387 7)
\lvec(370 7)
\ifill f:0
\move(388 6)
\lvec(393 6)
\lvec(393 7)
\lvec(388 7)
\ifill f:0
\move(394 6)
\lvec(399 6)
\lvec(399 7)
\lvec(394 7)
\ifill f:0
\move(400 6)
\lvec(405 6)
\lvec(405 7)
\lvec(400 7)
\ifill f:0
\move(406 6)
\lvec(424 6)
\lvec(424 7)
\lvec(406 7)
\ifill f:0
\move(425 6)
\lvec(437 6)
\lvec(437 7)
\lvec(425 7)
\ifill f:0
\move(438 6)
\lvec(450 6)
\lvec(450 7)
\lvec(438 7)
\ifill f:0
\move(10 7)
\lvec(29 7)
\lvec(29 8)
\lvec(10 8)
\ifill f:0
\move(30 7)
\lvec(32 7)
\lvec(32 8)
\lvec(30 8)
\ifill f:0
\move(33 7)
\lvec(47 7)
\lvec(47 8)
\lvec(33 8)
\ifill f:0
\move(48 7)
\lvec(51 7)
\lvec(51 8)
\lvec(48 8)
\ifill f:0
\move(52 7)
\lvec(53 7)
\lvec(53 8)
\lvec(52 8)
\ifill f:0
\move(54 7)
\lvec(59 7)
\lvec(59 8)
\lvec(54 8)
\ifill f:0
\move(60 7)
\lvec(61 7)
\lvec(61 8)
\lvec(60 8)
\ifill f:0
\move(62 7)
\lvec(63 7)
\lvec(63 8)
\lvec(62 8)
\ifill f:0
\move(64 7)
\lvec(65 7)
\lvec(65 8)
\lvec(64 8)
\ifill f:0
\move(66 7)
\lvec(74 7)
\lvec(74 8)
\lvec(66 8)
\ifill f:0
\move(75 7)
\lvec(86 7)
\lvec(86 8)
\lvec(75 8)
\ifill f:0
\move(87 7)
\lvec(91 7)
\lvec(91 8)
\lvec(87 8)
\ifill f:0
\move(92 7)
\lvec(99 7)
\lvec(99 8)
\lvec(92 8)
\ifill f:0
\move(100 7)
\lvec(107 7)
\lvec(107 8)
\lvec(100 8)
\ifill f:0
\move(108 7)
\lvec(110 7)
\lvec(110 8)
\lvec(108 8)
\ifill f:0
\move(111 7)
\lvec(143 7)
\lvec(143 8)
\lvec(111 8)
\ifill f:0
\move(144 7)
\lvec(146 7)
\lvec(146 8)
\lvec(144 8)
\ifill f:0
\move(147 7)
\lvec(156 7)
\lvec(156 8)
\lvec(147 8)
\ifill f:0
\move(157 7)
\lvec(166 7)
\lvec(166 8)
\lvec(157 8)
\ifill f:0
\move(167 7)
\lvec(173 7)
\lvec(173 8)
\lvec(167 8)
\ifill f:0
\move(174 7)
\lvec(180 7)
\lvec(180 8)
\lvec(174 8)
\ifill f:0
\move(181 7)
\lvec(191 7)
\lvec(191 8)
\lvec(181 8)
\ifill f:0
\move(192 7)
\lvec(206 7)
\lvec(206 8)
\lvec(192 8)
\ifill f:0
\move(207 7)
\lvec(210 7)
\lvec(210 8)
\lvec(207 8)
\ifill f:0
\move(211 7)
\lvec(214 7)
\lvec(214 8)
\lvec(211 8)
\ifill f:0
\move(215 7)
\lvec(218 7)
\lvec(218 8)
\lvec(215 8)
\ifill f:0
\move(219 7)
\lvec(230 7)
\lvec(230 8)
\lvec(219 8)
\ifill f:0
\move(231 7)
\lvec(234 7)
\lvec(234 8)
\lvec(231 8)
\ifill f:0
\move(235 7)
\lvec(238 7)
\lvec(238 8)
\lvec(235 8)
\ifill f:0
\move(239 7)
\lvec(242 7)
\lvec(242 8)
\lvec(239 8)
\ifill f:0
\move(243 7)
\lvec(259 7)
\lvec(259 8)
\lvec(243 8)
\ifill f:0
\move(260 7)
\lvec(272 7)
\lvec(272 8)
\lvec(260 8)
\ifill f:0
\move(273 7)
\lvec(281 7)
\lvec(281 8)
\lvec(273 8)
\ifill f:0
\move(282 7)
\lvec(290 7)
\lvec(290 8)
\lvec(282 8)
\ifill f:0
\move(291 7)
\lvec(304 7)
\lvec(304 8)
\lvec(291 8)
\ifill f:0
\move(305 7)
\lvec(318 7)
\lvec(318 8)
\lvec(305 8)
\ifill f:0
\move(319 7)
\lvec(323 7)
\lvec(323 8)
\lvec(319 8)
\ifill f:0
\move(324 7)
\lvec(378 7)
\lvec(378 8)
\lvec(324 8)
\ifill f:0
\move(379 7)
\lvec(383 7)
\lvec(383 8)
\lvec(379 8)
\ifill f:0
\move(384 7)
\lvec(399 7)
\lvec(399 8)
\lvec(384 8)
\ifill f:0
\move(400 7)
\lvec(415 7)
\lvec(415 8)
\lvec(400 8)
\ifill f:0
\move(416 7)
\lvec(426 7)
\lvec(426 8)
\lvec(416 8)
\ifill f:0
\move(427 7)
\lvec(437 7)
\lvec(437 8)
\lvec(427 8)
\ifill f:0
\move(438 7)
\lvec(451 7)
\lvec(451 8)
\lvec(438 8)
\ifill f:0
\move(11 8)
\lvec(12 8)
\lvec(12 9)
\lvec(11 9)
\ifill f:0
\move(13 8)
\lvec(17 8)
\lvec(17 9)
\lvec(13 9)
\ifill f:0
\move(18 8)
\lvec(26 8)
\lvec(26 9)
\lvec(18 9)
\ifill f:0
\move(28 8)
\lvec(30 8)
\lvec(30 9)
\lvec(28 9)
\ifill f:0
\move(32 8)
\lvec(34 8)
\lvec(34 9)
\lvec(32 9)
\ifill f:0
\move(35 8)
\lvec(52 8)
\lvec(52 9)
\lvec(35 9)
\ifill f:0
\move(53 8)
\lvec(54 8)
\lvec(54 9)
\lvec(53 9)
\ifill f:0
\move(55 8)
\lvec(63 8)
\lvec(63 9)
\lvec(55 9)
\ifill f:0
\move(64 8)
\lvec(65 8)
\lvec(65 9)
\lvec(64 9)
\ifill f:0
\move(66 8)
\lvec(67 8)
\lvec(67 9)
\lvec(66 9)
\ifill f:0
\move(68 8)
\lvec(69 8)
\lvec(69 9)
\lvec(68 9)
\ifill f:0
\move(70 8)
\lvec(75 8)
\lvec(75 9)
\lvec(70 9)
\ifill f:0
\move(76 8)
\lvec(77 8)
\lvec(77 9)
\lvec(76 9)
\ifill f:0
\move(78 8)
\lvec(79 8)
\lvec(79 9)
\lvec(78 9)
\ifill f:0
\move(80 8)
\lvec(83 8)
\lvec(83 9)
\lvec(80 9)
\ifill f:0
\move(84 8)
\lvec(92 8)
\lvec(92 9)
\lvec(84 9)
\ifill f:0
\move(93 8)
\lvec(94 8)
\lvec(94 9)
\lvec(93 9)
\ifill f:0
\move(95 8)
\lvec(99 8)
\lvec(99 9)
\lvec(95 9)
\ifill f:0
\move(100 8)
\lvec(101 8)
\lvec(101 9)
\lvec(100 9)
\ifill f:0
\move(102 8)
\lvec(111 8)
\lvec(111 9)
\lvec(102 9)
\ifill f:0
\move(112 8)
\lvec(116 8)
\lvec(116 9)
\lvec(112 9)
\ifill f:0
\move(117 8)
\lvec(124 8)
\lvec(124 9)
\lvec(117 9)
\ifill f:0
\move(125 8)
\lvec(129 8)
\lvec(129 9)
\lvec(125 9)
\ifill f:0
\move(130 8)
\lvec(132 8)
\lvec(132 9)
\lvec(130 9)
\ifill f:0
\move(133 8)
\lvec(137 8)
\lvec(137 9)
\lvec(133 9)
\ifill f:0
\move(138 8)
\lvec(140 8)
\lvec(140 9)
\lvec(138 9)
\ifill f:0
\move(141 8)
\lvec(143 8)
\lvec(143 9)
\lvec(141 9)
\ifill f:0
\move(144 8)
\lvec(160 8)
\lvec(160 9)
\lvec(144 9)
\ifill f:0
\move(161 8)
\lvec(166 8)
\lvec(166 9)
\lvec(161 9)
\ifill f:0
\move(167 8)
\lvec(182 8)
\lvec(182 9)
\lvec(167 9)
\ifill f:0
\move(183 8)
\lvec(185 8)
\lvec(185 9)
\lvec(183 9)
\ifill f:0
\move(186 8)
\lvec(188 8)
\lvec(188 9)
\lvec(186 9)
\ifill f:0
\move(189 8)
\lvec(195 8)
\lvec(195 9)
\lvec(189 9)
\ifill f:0
\move(196 8)
\lvec(198 8)
\lvec(198 9)
\lvec(196 9)
\ifill f:0
\move(199 8)
\lvec(205 8)
\lvec(205 9)
\lvec(199 9)
\ifill f:0
\move(206 8)
\lvec(208 8)
\lvec(208 9)
\lvec(206 9)
\ifill f:0
\move(209 8)
\lvec(215 8)
\lvec(215 9)
\lvec(209 9)
\ifill f:0
\move(216 8)
\lvec(222 8)
\lvec(222 9)
\lvec(216 9)
\ifill f:0
\move(223 8)
\lvec(229 8)
\lvec(229 9)
\lvec(223 9)
\ifill f:0
\move(230 8)
\lvec(236 8)
\lvec(236 9)
\lvec(230 9)
\ifill f:0
\move(237 8)
\lvec(240 8)
\lvec(240 9)
\lvec(237 9)
\ifill f:0
\move(241 8)
\lvec(247 8)
\lvec(247 9)
\lvec(241 9)
\ifill f:0
\move(248 8)
\lvec(251 8)
\lvec(251 9)
\lvec(248 9)
\ifill f:0
\move(252 8)
\lvec(266 8)
\lvec(266 9)
\lvec(252 9)
\ifill f:0
\move(267 8)
\lvec(270 8)
\lvec(270 9)
\lvec(267 9)
\ifill f:0
\move(271 8)
\lvec(274 8)
\lvec(274 9)
\lvec(271 9)
\ifill f:0
\move(275 8)
\lvec(278 8)
\lvec(278 9)
\lvec(275 9)
\ifill f:0
\move(279 8)
\lvec(282 8)
\lvec(282 9)
\lvec(279 9)
\ifill f:0
\move(283 8)
\lvec(294 8)
\lvec(294 9)
\lvec(283 9)
\ifill f:0
\move(295 8)
\lvec(298 8)
\lvec(298 9)
\lvec(295 9)
\ifill f:0
\move(299 8)
\lvec(302 8)
\lvec(302 9)
\lvec(299 9)
\ifill f:0
\move(303 8)
\lvec(306 8)
\lvec(306 9)
\lvec(303 9)
\ifill f:0
\move(307 8)
\lvec(310 8)
\lvec(310 9)
\lvec(307 9)
\ifill f:0
\move(311 8)
\lvec(327 8)
\lvec(327 9)
\lvec(311 9)
\ifill f:0
\move(328 8)
\lvec(331 8)
\lvec(331 9)
\lvec(328 9)
\ifill f:0
\move(332 8)
\lvec(340 8)
\lvec(340 9)
\lvec(332 9)
\ifill f:0
\move(341 8)
\lvec(344 8)
\lvec(344 9)
\lvec(341 9)
\ifill f:0
\move(345 8)
\lvec(353 8)
\lvec(353 9)
\lvec(345 9)
\ifill f:0
\move(354 8)
\lvec(362 8)
\lvec(362 9)
\lvec(354 9)
\ifill f:0
\move(363 8)
\lvec(371 8)
\lvec(371 9)
\lvec(363 9)
\ifill f:0
\move(372 8)
\lvec(380 8)
\lvec(380 9)
\lvec(372 9)
\ifill f:0
\move(381 8)
\lvec(385 8)
\lvec(385 9)
\lvec(381 9)
\ifill f:0
\move(386 8)
\lvec(394 8)
\lvec(394 9)
\lvec(386 9)
\ifill f:0
\move(395 8)
\lvec(399 8)
\lvec(399 9)
\lvec(395 9)
\ifill f:0
\move(400 8)
\lvec(408 8)
\lvec(408 9)
\lvec(400 9)
\ifill f:0
\move(409 8)
\lvec(413 8)
\lvec(413 9)
\lvec(409 9)
\ifill f:0
\move(414 8)
\lvec(418 8)
\lvec(418 9)
\lvec(414 9)
\ifill f:0
\move(419 8)
\lvec(442 8)
\lvec(442 9)
\lvec(419 9)
\ifill f:0
\move(443 8)
\lvec(447 8)
\lvec(447 9)
\lvec(443 9)
\ifill f:0
\move(448 8)
\lvec(451 8)
\lvec(451 9)
\lvec(448 9)
\ifill f:0
\move(12 9)
\lvec(13 9)
\lvec(13 10)
\lvec(12 10)
\ifill f:0
\move(14 9)
\lvec(17 9)
\lvec(17 10)
\lvec(14 10)
\ifill f:0
\move(18 9)
\lvec(19 9)
\lvec(19 10)
\lvec(18 10)
\ifill f:0
\move(20 9)
\lvec(21 9)
\lvec(21 10)
\lvec(20 10)
\ifill f:0
\move(22 9)
\lvec(24 9)
\lvec(24 10)
\lvec(22 10)
\ifill f:0
\move(25 9)
\lvec(26 9)
\lvec(26 10)
\lvec(25 10)
\ifill f:0
\move(31 9)
\lvec(34 9)
\lvec(34 10)
\lvec(31 10)
\ifill f:0
\move(35 9)
\lvec(38 9)
\lvec(38 10)
\lvec(35 10)
\ifill f:0
\move(39 9)
\lvec(42 9)
\lvec(42 10)
\lvec(39 10)
\ifill f:0
\move(43 9)
\lvec(46 9)
\lvec(46 10)
\lvec(43 10)
\ifill f:0
\move(47 9)
\lvec(50 9)
\lvec(50 10)
\lvec(47 10)
\ifill f:0
\move(51 9)
\lvec(52 9)
\lvec(52 10)
\lvec(51 10)
\ifill f:0
\move(53 9)
\lvec(55 9)
\lvec(55 10)
\lvec(53 10)
\ifill f:0
\move(56 9)
\lvec(60 9)
\lvec(60 10)
\lvec(56 10)
\ifill f:0
\move(61 9)
\lvec(65 9)
\lvec(65 10)
\lvec(61 10)
\ifill f:0
\move(66 9)
\lvec(79 9)
\lvec(79 10)
\lvec(66 10)
\ifill f:0
\move(80 9)
\lvec(83 9)
\lvec(83 10)
\lvec(80 10)
\ifill f:0
\move(84 9)
\lvec(85 9)
\lvec(85 10)
\lvec(84 10)
\ifill f:0
\move(86 9)
\lvec(87 9)
\lvec(87 10)
\lvec(86 10)
\ifill f:0
\move(88 9)
\lvec(93 9)
\lvec(93 10)
\lvec(88 10)
\ifill f:0
\move(94 9)
\lvec(95 9)
\lvec(95 10)
\lvec(94 10)
\ifill f:0
\move(96 9)
\lvec(97 9)
\lvec(97 10)
\lvec(96 10)
\ifill f:0
\move(98 9)
\lvec(99 9)
\lvec(99 10)
\lvec(98 10)
\ifill f:0
\move(100 9)
\lvec(101 9)
\lvec(101 10)
\lvec(100 10)
\ifill f:0
\move(102 9)
\lvec(112 9)
\lvec(112 10)
\lvec(102 10)
\ifill f:0
\move(113 9)
\lvec(119 9)
\lvec(119 10)
\lvec(113 10)
\ifill f:0
\move(120 9)
\lvec(126 9)
\lvec(126 10)
\lvec(120 10)
\ifill f:0
\move(127 9)
\lvec(128 9)
\lvec(128 10)
\lvec(127 10)
\ifill f:0
\move(129 9)
\lvec(133 9)
\lvec(133 10)
\lvec(129 10)
\ifill f:0
\move(134 9)
\lvec(138 9)
\lvec(138 10)
\lvec(134 10)
\ifill f:0
\move(139 9)
\lvec(143 9)
\lvec(143 10)
\lvec(139 10)
\ifill f:0
\move(144 9)
\lvec(148 9)
\lvec(148 10)
\lvec(144 10)
\ifill f:0
\move(149 9)
\lvec(156 9)
\lvec(156 10)
\lvec(149 10)
\ifill f:0
\move(157 9)
\lvec(161 9)
\lvec(161 10)
\lvec(157 10)
\ifill f:0
\move(162 9)
\lvec(164 9)
\lvec(164 10)
\lvec(162 10)
\ifill f:0
\move(165 9)
\lvec(172 9)
\lvec(172 10)
\lvec(165 10)
\ifill f:0
\move(173 9)
\lvec(180 9)
\lvec(180 10)
\lvec(173 10)
\ifill f:0
\move(181 9)
\lvec(183 9)
\lvec(183 10)
\lvec(181 10)
\ifill f:0
\move(184 9)
\lvec(186 9)
\lvec(186 10)
\lvec(184 10)
\ifill f:0
\move(187 9)
\lvec(189 9)
\lvec(189 10)
\lvec(187 10)
\ifill f:0
\move(190 9)
\lvec(192 9)
\lvec(192 10)
\lvec(190 10)
\ifill f:0
\move(193 9)
\lvec(195 9)
\lvec(195 10)
\lvec(193 10)
\ifill f:0
\move(196 9)
\lvec(198 9)
\lvec(198 10)
\lvec(196 10)
\ifill f:0
\move(199 9)
\lvec(207 9)
\lvec(207 10)
\lvec(199 10)
\ifill f:0
\move(208 9)
\lvec(210 9)
\lvec(210 10)
\lvec(208 10)
\ifill f:0
\move(211 9)
\lvec(213 9)
\lvec(213 10)
\lvec(211 10)
\ifill f:0
\move(214 9)
\lvec(216 9)
\lvec(216 10)
\lvec(214 10)
\ifill f:0
\move(217 9)
\lvec(219 9)
\lvec(219 10)
\lvec(217 10)
\ifill f:0
\move(220 9)
\lvec(222 9)
\lvec(222 10)
\lvec(220 10)
\ifill f:0
\move(223 9)
\lvec(235 9)
\lvec(235 10)
\lvec(223 10)
\ifill f:0
\move(236 9)
\lvec(238 9)
\lvec(238 10)
\lvec(236 10)
\ifill f:0
\move(239 9)
\lvec(245 9)
\lvec(245 10)
\lvec(239 10)
\ifill f:0
\move(246 9)
\lvec(248 9)
\lvec(248 10)
\lvec(246 10)
\ifill f:0
\move(249 9)
\lvec(255 9)
\lvec(255 10)
\lvec(249 10)
\ifill f:0
\move(256 9)
\lvec(258 9)
\lvec(258 10)
\lvec(256 10)
\ifill f:0
\move(259 9)
\lvec(265 9)
\lvec(265 10)
\lvec(259 10)
\ifill f:0
\move(266 9)
\lvec(272 9)
\lvec(272 10)
\lvec(266 10)
\ifill f:0
\move(273 9)
\lvec(275 9)
\lvec(275 10)
\lvec(273 10)
\ifill f:0
\move(276 9)
\lvec(279 9)
\lvec(279 10)
\lvec(276 10)
\ifill f:0
\move(280 9)
\lvec(286 9)
\lvec(286 10)
\lvec(280 10)
\ifill f:0
\move(287 9)
\lvec(293 9)
\lvec(293 10)
\lvec(287 10)
\ifill f:0
\move(294 9)
\lvec(297 9)
\lvec(297 10)
\lvec(294 10)
\ifill f:0
\move(298 9)
\lvec(304 9)
\lvec(304 10)
\lvec(298 10)
\ifill f:0
\move(305 9)
\lvec(308 9)
\lvec(308 10)
\lvec(305 10)
\ifill f:0
\move(309 9)
\lvec(315 9)
\lvec(315 10)
\lvec(309 10)
\ifill f:0
\move(316 9)
\lvec(319 9)
\lvec(319 10)
\lvec(316 10)
\ifill f:0
\move(320 9)
\lvec(330 9)
\lvec(330 10)
\lvec(320 10)
\ifill f:0
\move(331 9)
\lvec(334 9)
\lvec(334 10)
\lvec(331 10)
\ifill f:0
\move(335 9)
\lvec(338 9)
\lvec(338 10)
\lvec(335 10)
\ifill f:0
\move(339 9)
\lvec(342 9)
\lvec(342 10)
\lvec(339 10)
\ifill f:0
\move(343 9)
\lvec(346 9)
\lvec(346 10)
\lvec(343 10)
\ifill f:0
\move(347 9)
\lvec(350 9)
\lvec(350 10)
\lvec(347 10)
\ifill f:0
\move(351 9)
\lvec(354 9)
\lvec(354 10)
\lvec(351 10)
\ifill f:0
\move(355 9)
\lvec(366 9)
\lvec(366 10)
\lvec(355 10)
\ifill f:0
\move(367 9)
\lvec(370 9)
\lvec(370 10)
\lvec(367 10)
\ifill f:0
\move(371 9)
\lvec(374 9)
\lvec(374 10)
\lvec(371 10)
\ifill f:0
\move(375 9)
\lvec(378 9)
\lvec(378 10)
\lvec(375 10)
\ifill f:0
\move(379 9)
\lvec(382 9)
\lvec(382 10)
\lvec(379 10)
\ifill f:0
\move(383 9)
\lvec(386 9)
\lvec(386 10)
\lvec(383 10)
\ifill f:0
\move(387 9)
\lvec(390 9)
\lvec(390 10)
\lvec(387 10)
\ifill f:0
\move(391 9)
\lvec(403 9)
\lvec(403 10)
\lvec(391 10)
\ifill f:0
\move(404 9)
\lvec(407 9)
\lvec(407 10)
\lvec(404 10)
\ifill f:0
\move(408 9)
\lvec(416 9)
\lvec(416 10)
\lvec(408 10)
\ifill f:0
\move(417 9)
\lvec(420 9)
\lvec(420 10)
\lvec(417 10)
\ifill f:0
\move(421 9)
\lvec(429 9)
\lvec(429 10)
\lvec(421 10)
\ifill f:0
\move(430 9)
\lvec(433 9)
\lvec(433 10)
\lvec(430 10)
\ifill f:0
\move(434 9)
\lvec(442 9)
\lvec(442 10)
\lvec(434 10)
\ifill f:0
\move(443 9)
\lvec(451 9)
\lvec(451 10)
\lvec(443 10)
\ifill f:0
\move(11 10)
\lvec(12 10)
\lvec(12 11)
\lvec(11 11)
\ifill f:0
\move(13 10)
\lvec(14 10)
\lvec(14 11)
\lvec(13 11)
\ifill f:0
\move(15 10)
\lvec(17 10)
\lvec(17 11)
\lvec(15 11)
\ifill f:0
\move(18 10)
\lvec(23 10)
\lvec(23 11)
\lvec(18 11)
\ifill f:0
\move(25 10)
\lvec(26 10)
\lvec(26 11)
\lvec(25 11)
\ifill f:0
\move(27 10)
\lvec(29 10)
\lvec(29 11)
\lvec(27 11)
\ifill f:0
\move(35 10)
\lvec(50 10)
\lvec(50 11)
\lvec(35 11)
\ifill f:0
\move(51 10)
\lvec(53 10)
\lvec(53 11)
\lvec(51 11)
\ifill f:0
\move(55 10)
\lvec(56 10)
\lvec(56 11)
\lvec(55 11)
\ifill f:0
\move(57 10)
\lvec(68 10)
\lvec(68 11)
\lvec(57 11)
\ifill f:0
\move(69 10)
\lvec(71 10)
\lvec(71 11)
\lvec(69 11)
\ifill f:0
\move(72 10)
\lvec(76 10)
\lvec(76 11)
\lvec(72 11)
\ifill f:0
\move(77 10)
\lvec(86 10)
\lvec(86 11)
\lvec(77 11)
\ifill f:0
\move(87 10)
\lvec(88 10)
\lvec(88 11)
\lvec(87 11)
\ifill f:0
\move(89 10)
\lvec(95 10)
\lvec(95 11)
\lvec(89 11)
\ifill f:0
\move(96 10)
\lvec(99 10)
\lvec(99 11)
\lvec(96 11)
\ifill f:0
\move(100 10)
\lvec(101 10)
\lvec(101 11)
\lvec(100 11)
\ifill f:0
\move(102 10)
\lvec(103 10)
\lvec(103 11)
\lvec(102 11)
\ifill f:0
\move(104 10)
\lvec(105 10)
\lvec(105 11)
\lvec(104 11)
\ifill f:0
\move(106 10)
\lvec(107 10)
\lvec(107 11)
\lvec(106 11)
\ifill f:0
\move(108 10)
\lvec(113 10)
\lvec(113 11)
\lvec(108 11)
\ifill f:0
\move(114 10)
\lvec(115 10)
\lvec(115 11)
\lvec(114 11)
\ifill f:0
\move(116 10)
\lvec(117 10)
\lvec(117 11)
\lvec(116 11)
\ifill f:0
\move(118 10)
\lvec(119 10)
\lvec(119 11)
\lvec(118 11)
\ifill f:0
\move(120 10)
\lvec(123 10)
\lvec(123 11)
\lvec(120 11)
\ifill f:0
\move(124 10)
\lvec(125 10)
\lvec(125 11)
\lvec(124 11)
\ifill f:0
\move(126 10)
\lvec(134 10)
\lvec(134 11)
\lvec(126 11)
\ifill f:0
\move(135 10)
\lvec(143 10)
\lvec(143 11)
\lvec(135 11)
\ifill f:0
\move(144 10)
\lvec(145 10)
\lvec(145 11)
\lvec(144 11)
\ifill f:0
\move(146 10)
\lvec(150 10)
\lvec(150 11)
\lvec(146 11)
\ifill f:0
\move(151 10)
\lvec(157 10)
\lvec(157 11)
\lvec(151 11)
\ifill f:0
\move(158 10)
\lvec(162 10)
\lvec(162 11)
\lvec(158 11)
\ifill f:0
\move(163 10)
\lvec(174 10)
\lvec(174 11)
\lvec(163 11)
\ifill f:0
\move(175 10)
\lvec(182 10)
\lvec(182 11)
\lvec(175 11)
\ifill f:0
\move(183 10)
\lvec(187 10)
\lvec(187 11)
\lvec(183 11)
\ifill f:0
\move(188 10)
\lvec(192 10)
\lvec(192 11)
\lvec(188 11)
\ifill f:0
\move(193 10)
\lvec(195 10)
\lvec(195 11)
\lvec(193 11)
\ifill f:0
\move(196 10)
\lvec(203 10)
\lvec(203 11)
\lvec(196 11)
\ifill f:0
\move(204 10)
\lvec(211 10)
\lvec(211 11)
\lvec(204 11)
\ifill f:0
\move(212 10)
\lvec(214 10)
\lvec(214 11)
\lvec(212 11)
\ifill f:0
\move(215 10)
\lvec(228 10)
\lvec(228 11)
\lvec(215 11)
\ifill f:0
\move(229 10)
\lvec(231 10)
\lvec(231 11)
\lvec(229 11)
\ifill f:0
\move(232 10)
\lvec(234 10)
\lvec(234 11)
\lvec(232 11)
\ifill f:0
\move(235 10)
\lvec(237 10)
\lvec(237 11)
\lvec(235 11)
\ifill f:0
\move(238 10)
\lvec(240 10)
\lvec(240 11)
\lvec(238 11)
\ifill f:0
\move(241 10)
\lvec(243 10)
\lvec(243 11)
\lvec(241 11)
\ifill f:0
\move(244 10)
\lvec(252 10)
\lvec(252 11)
\lvec(244 11)
\ifill f:0
\move(253 10)
\lvec(255 10)
\lvec(255 11)
\lvec(253 11)
\ifill f:0
\move(256 10)
\lvec(258 10)
\lvec(258 11)
\lvec(256 11)
\ifill f:0
\move(259 10)
\lvec(261 10)
\lvec(261 11)
\lvec(259 11)
\ifill f:0
\move(262 10)
\lvec(264 10)
\lvec(264 11)
\lvec(262 11)
\ifill f:0
\move(265 10)
\lvec(267 10)
\lvec(267 11)
\lvec(265 11)
\ifill f:0
\move(268 10)
\lvec(270 10)
\lvec(270 11)
\lvec(268 11)
\ifill f:0
\move(271 10)
\lvec(283 10)
\lvec(283 11)
\lvec(271 11)
\ifill f:0
\move(284 10)
\lvec(286 10)
\lvec(286 11)
\lvec(284 11)
\ifill f:0
\move(287 10)
\lvec(296 10)
\lvec(296 11)
\lvec(287 11)
\ifill f:0
\move(297 10)
\lvec(299 10)
\lvec(299 11)
\lvec(297 11)
\ifill f:0
\move(300 10)
\lvec(306 10)
\lvec(306 11)
\lvec(300 11)
\ifill f:0
\move(307 10)
\lvec(309 10)
\lvec(309 11)
\lvec(307 11)
\ifill f:0
\move(310 10)
\lvec(316 10)
\lvec(316 11)
\lvec(310 11)
\ifill f:0
\move(317 10)
\lvec(319 10)
\lvec(319 11)
\lvec(317 11)
\ifill f:0
\move(320 10)
\lvec(323 10)
\lvec(323 11)
\lvec(320 11)
\ifill f:0
\move(324 10)
\lvec(326 10)
\lvec(326 11)
\lvec(324 11)
\ifill f:0
\move(327 10)
\lvec(333 10)
\lvec(333 11)
\lvec(327 11)
\ifill f:0
\move(334 10)
\lvec(340 10)
\lvec(340 11)
\lvec(334 11)
\ifill f:0
\move(341 10)
\lvec(347 10)
\lvec(347 11)
\lvec(341 11)
\ifill f:0
\move(348 10)
\lvec(351 10)
\lvec(351 11)
\lvec(348 11)
\ifill f:0
\move(352 10)
\lvec(354 10)
\lvec(354 11)
\lvec(352 11)
\ifill f:0
\move(355 10)
\lvec(358 10)
\lvec(358 11)
\lvec(355 11)
\ifill f:0
\move(359 10)
\lvec(365 10)
\lvec(365 11)
\lvec(359 11)
\ifill f:0
\move(366 10)
\lvec(369 10)
\lvec(369 11)
\lvec(366 11)
\ifill f:0
\move(370 10)
\lvec(376 10)
\lvec(376 11)
\lvec(370 11)
\ifill f:0
\move(377 10)
\lvec(380 10)
\lvec(380 11)
\lvec(377 11)
\ifill f:0
\move(381 10)
\lvec(391 10)
\lvec(391 11)
\lvec(381 11)
\ifill f:0
\move(392 10)
\lvec(395 10)
\lvec(395 11)
\lvec(392 11)
\ifill f:0
\move(396 10)
\lvec(410 10)
\lvec(410 11)
\lvec(396 11)
\ifill f:0
\move(411 10)
\lvec(414 10)
\lvec(414 11)
\lvec(411 11)
\ifill f:0
\move(415 10)
\lvec(418 10)
\lvec(418 11)
\lvec(415 11)
\ifill f:0
\move(419 10)
\lvec(422 10)
\lvec(422 11)
\lvec(419 11)
\ifill f:0
\move(423 10)
\lvec(426 10)
\lvec(426 11)
\lvec(423 11)
\ifill f:0
\move(427 10)
\lvec(430 10)
\lvec(430 11)
\lvec(427 11)
\ifill f:0
\move(431 10)
\lvec(434 10)
\lvec(434 11)
\lvec(431 11)
\ifill f:0
\move(435 10)
\lvec(446 10)
\lvec(446 11)
\lvec(435 11)
\ifill f:0
\move(447 10)
\lvec(450 10)
\lvec(450 11)
\lvec(447 11)
\ifill f:0
\move(11 11)
\lvec(13 11)
\lvec(13 12)
\lvec(11 12)
\ifill f:0
\move(14 11)
\lvec(17 11)
\lvec(17 12)
\lvec(14 12)
\ifill f:0
\move(18 11)
\lvec(24 11)
\lvec(24 12)
\lvec(18 12)
\ifill f:0
\move(25 11)
\lvec(30 11)
\lvec(30 12)
\lvec(25 12)
\ifill f:0
\move(35 11)
\lvec(42 11)
\lvec(42 12)
\lvec(35 12)
\ifill f:0
\move(43 11)
\lvec(50 11)
\lvec(50 12)
\lvec(43 12)
\ifill f:0
\move(52 11)
\lvec(54 11)
\lvec(54 12)
\lvec(52 12)
\ifill f:0
\move(56 11)
\lvec(62 11)
\lvec(62 12)
\lvec(56 12)
\ifill f:0
\move(63 11)
\lvec(69 11)
\lvec(69 12)
\lvec(63 12)
\ifill f:0
\move(70 11)
\lvec(72 11)
\lvec(72 12)
\lvec(70 12)
\ifill f:0
\move(73 11)
\lvec(78 11)
\lvec(78 12)
\lvec(73 12)
\ifill f:0
\move(79 11)
\lvec(84 11)
\lvec(84 12)
\lvec(79 12)
\ifill f:0
\move(85 11)
\lvec(89 11)
\lvec(89 12)
\lvec(85 12)
\ifill f:0
\move(90 11)
\lvec(94 11)
\lvec(94 12)
\lvec(90 12)
\ifill f:0
\move(95 11)
\lvec(99 11)
\lvec(99 12)
\lvec(95 12)
\ifill f:0
\move(100 11)
\lvec(101 11)
\lvec(101 12)
\lvec(100 12)
\ifill f:0
\move(102 11)
\lvec(119 11)
\lvec(119 12)
\lvec(102 12)
\ifill f:0
\move(120 11)
\lvec(123 11)
\lvec(123 12)
\lvec(120 12)
\ifill f:0
\move(124 11)
\lvec(125 11)
\lvec(125 12)
\lvec(124 12)
\ifill f:0
\move(126 11)
\lvec(127 11)
\lvec(127 12)
\lvec(126 12)
\ifill f:0
\move(128 11)
\lvec(129 11)
\lvec(129 12)
\lvec(128 12)
\ifill f:0
\move(130 11)
\lvec(135 11)
\lvec(135 12)
\lvec(130 12)
\ifill f:0
\move(136 11)
\lvec(137 11)
\lvec(137 12)
\lvec(136 12)
\ifill f:0
\move(138 11)
\lvec(139 11)
\lvec(139 12)
\lvec(138 12)
\ifill f:0
\move(140 11)
\lvec(141 11)
\lvec(141 12)
\lvec(140 12)
\ifill f:0
\move(142 11)
\lvec(143 11)
\lvec(143 12)
\lvec(142 12)
\ifill f:0
\move(144 11)
\lvec(145 11)
\lvec(145 12)
\lvec(144 12)
\ifill f:0
\move(146 11)
\lvec(147 11)
\lvec(147 12)
\lvec(146 12)
\ifill f:0
\move(148 11)
\lvec(158 11)
\lvec(158 12)
\lvec(148 12)
\ifill f:0
\move(159 11)
\lvec(160 11)
\lvec(160 12)
\lvec(159 12)
\ifill f:0
\move(161 11)
\lvec(167 11)
\lvec(167 12)
\lvec(161 12)
\ifill f:0
\move(168 11)
\lvec(171 11)
\lvec(171 12)
\lvec(168 12)
\ifill f:0
\move(172 11)
\lvec(176 11)
\lvec(176 12)
\lvec(172 12)
\ifill f:0
\move(177 11)
\lvec(183 11)
\lvec(183 12)
\lvec(177 12)
\ifill f:0
\move(184 11)
\lvec(190 11)
\lvec(190 12)
\lvec(184 12)
\ifill f:0
\move(191 11)
\lvec(195 11)
\lvec(195 12)
\lvec(191 12)
\ifill f:0
\move(196 11)
\lvec(197 11)
\lvec(197 12)
\lvec(196 12)
\ifill f:0
\move(198 11)
\lvec(200 11)
\lvec(200 12)
\lvec(198 12)
\ifill f:0
\move(201 11)
\lvec(210 11)
\lvec(210 12)
\lvec(201 12)
\ifill f:0
\move(211 11)
\lvec(215 11)
\lvec(215 12)
\lvec(211 12)
\ifill f:0
\move(216 11)
\lvec(220 11)
\lvec(220 12)
\lvec(216 12)
\ifill f:0
\move(221 11)
\lvec(228 11)
\lvec(228 12)
\lvec(221 12)
\ifill f:0
\move(229 11)
\lvec(233 11)
\lvec(233 12)
\lvec(229 12)
\ifill f:0
\move(234 11)
\lvec(241 11)
\lvec(241 12)
\lvec(234 12)
\ifill f:0
\move(242 11)
\lvec(244 11)
\lvec(244 12)
\lvec(242 12)
\ifill f:0
\move(245 11)
\lvec(252 11)
\lvec(252 12)
\lvec(245 12)
\ifill f:0
\move(253 11)
\lvec(255 11)
\lvec(255 12)
\lvec(253 12)
\ifill f:0
\move(256 11)
\lvec(266 11)
\lvec(266 12)
\lvec(256 12)
\ifill f:0
\move(267 11)
\lvec(269 11)
\lvec(269 12)
\lvec(267 12)
\ifill f:0
\move(270 11)
\lvec(272 11)
\lvec(272 12)
\lvec(270 12)
\ifill f:0
\move(273 11)
\lvec(286 11)
\lvec(286 12)
\lvec(273 12)
\ifill f:0
\move(287 11)
\lvec(292 11)
\lvec(292 12)
\lvec(287 12)
\ifill f:0
\move(293 11)
\lvec(295 11)
\lvec(295 12)
\lvec(293 12)
\ifill f:0
\move(296 11)
\lvec(298 11)
\lvec(298 12)
\lvec(296 12)
\ifill f:0
\move(299 11)
\lvec(301 11)
\lvec(301 12)
\lvec(299 12)
\ifill f:0
\move(302 11)
\lvec(304 11)
\lvec(304 12)
\lvec(302 12)
\ifill f:0
\move(305 11)
\lvec(307 11)
\lvec(307 12)
\lvec(305 12)
\ifill f:0
\move(308 11)
\lvec(323 11)
\lvec(323 12)
\lvec(308 12)
\ifill f:0
\move(324 11)
\lvec(326 11)
\lvec(326 12)
\lvec(324 12)
\ifill f:0
\move(327 11)
\lvec(329 11)
\lvec(329 12)
\lvec(327 12)
\ifill f:0
\move(330 11)
\lvec(332 11)
\lvec(332 12)
\lvec(330 12)
\ifill f:0
\move(333 11)
\lvec(342 11)
\lvec(342 12)
\lvec(333 12)
\ifill f:0
\move(343 11)
\lvec(345 11)
\lvec(345 12)
\lvec(343 12)
\ifill f:0
\move(346 11)
\lvec(348 11)
\lvec(348 12)
\lvec(346 12)
\ifill f:0
\move(349 11)
\lvec(355 11)
\lvec(355 12)
\lvec(349 12)
\ifill f:0
\move(356 11)
\lvec(358 11)
\lvec(358 12)
\lvec(356 12)
\ifill f:0
\move(359 11)
\lvec(368 11)
\lvec(368 12)
\lvec(359 12)
\ifill f:0
\move(369 11)
\lvec(375 11)
\lvec(375 12)
\lvec(369 12)
\ifill f:0
\move(376 11)
\lvec(378 11)
\lvec(378 12)
\lvec(376 12)
\ifill f:0
\move(379 11)
\lvec(385 11)
\lvec(385 12)
\lvec(379 12)
\ifill f:0
\move(386 11)
\lvec(392 11)
\lvec(392 12)
\lvec(386 12)
\ifill f:0
\move(393 11)
\lvec(395 11)
\lvec(395 12)
\lvec(393 12)
\ifill f:0
\move(396 11)
\lvec(399 11)
\lvec(399 12)
\lvec(396 12)
\ifill f:0
\move(400 11)
\lvec(402 11)
\lvec(402 12)
\lvec(400 12)
\ifill f:0
\move(403 11)
\lvec(409 11)
\lvec(409 12)
\lvec(403 12)
\ifill f:0
\move(410 11)
\lvec(413 11)
\lvec(413 12)
\lvec(410 12)
\ifill f:0
\move(414 11)
\lvec(416 11)
\lvec(416 12)
\lvec(414 12)
\ifill f:0
\move(417 11)
\lvec(420 11)
\lvec(420 12)
\lvec(417 12)
\ifill f:0
\move(421 11)
\lvec(427 11)
\lvec(427 12)
\lvec(421 12)
\ifill f:0
\move(428 11)
\lvec(434 11)
\lvec(434 12)
\lvec(428 12)
\ifill f:0
\move(435 11)
\lvec(438 11)
\lvec(438 12)
\lvec(435 12)
\ifill f:0
\move(439 11)
\lvec(445 11)
\lvec(445 12)
\lvec(439 12)
\ifill f:0
\move(446 11)
\lvec(449 11)
\lvec(449 12)
\lvec(446 12)
\ifill f:0
\move(450 11)
\lvec(451 11)
\lvec(451 12)
\lvec(450 12)
\ifill f:0
\move(13 12)
\lvec(15 12)
\lvec(15 13)
\lvec(13 13)
\ifill f:0
\move(16 12)
\lvec(17 12)
\lvec(17 13)
\lvec(16 13)
\ifill f:0
\move(18 12)
\lvec(19 12)
\lvec(19 13)
\lvec(18 13)
\ifill f:0
\move(20 12)
\lvec(24 12)
\lvec(24 13)
\lvec(20 13)
\ifill f:0
\move(25 12)
\lvec(26 12)
\lvec(26 13)
\lvec(25 13)
\ifill f:0
\move(31 12)
\lvec(34 12)
\lvec(34 13)
\lvec(31 13)
\ifill f:0
\move(35 12)
\lvec(37 12)
\lvec(37 13)
\lvec(35 13)
\ifill f:0
\move(38 12)
\lvec(50 12)
\lvec(50 13)
\lvec(38 13)
\ifill f:0
\move(54 12)
\lvec(56 12)
\lvec(56 13)
\lvec(54 13)
\ifill f:0
\move(58 12)
\lvec(62 12)
\lvec(62 13)
\lvec(58 13)
\ifill f:0
\move(63 12)
\lvec(66 12)
\lvec(66 13)
\lvec(63 13)
\ifill f:0
\move(67 12)
\lvec(74 12)
\lvec(74 13)
\lvec(67 13)
\ifill f:0
\move(75 12)
\lvec(82 12)
\lvec(82 13)
\lvec(75 13)
\ifill f:0
\move(83 12)
\lvec(84 12)
\lvec(84 13)
\lvec(83 13)
\ifill f:0
\move(86 12)
\lvec(90 12)
\lvec(90 13)
\lvec(86 13)
\ifill f:0
\move(91 12)
\lvec(93 12)
\lvec(93 13)
\lvec(91 13)
\ifill f:0
\move(94 12)
\lvec(96 12)
\lvec(96 13)
\lvec(94 13)
\ifill f:0
\move(97 12)
\lvec(99 12)
\lvec(99 13)
\lvec(97 13)
\ifill f:0
\move(100 12)
\lvec(109 12)
\lvec(109 13)
\lvec(100 13)
\ifill f:0
\move(110 12)
\lvec(128 12)
\lvec(128 13)
\lvec(110 13)
\ifill f:0
\move(129 12)
\lvec(130 12)
\lvec(130 13)
\lvec(129 13)
\ifill f:0
\move(131 12)
\lvec(143 12)
\lvec(143 13)
\lvec(131 13)
\ifill f:0
\move(144 12)
\lvec(145 12)
\lvec(145 13)
\lvec(144 13)
\ifill f:0
\move(146 12)
\lvec(147 12)
\lvec(147 13)
\lvec(146 13)
\ifill f:0
\move(148 12)
\lvec(149 12)
\lvec(149 13)
\lvec(148 13)
\ifill f:0
\move(150 12)
\lvec(151 12)
\lvec(151 13)
\lvec(150 13)
\ifill f:0
\move(152 12)
\lvec(153 12)
\lvec(153 13)
\lvec(152 13)
\ifill f:0
\move(154 12)
\lvec(159 12)
\lvec(159 13)
\lvec(154 13)
\ifill f:0
\move(160 12)
\lvec(161 12)
\lvec(161 13)
\lvec(160 13)
\ifill f:0
\move(162 12)
\lvec(163 12)
\lvec(163 13)
\lvec(162 13)
\ifill f:0
\move(164 12)
\lvec(165 12)
\lvec(165 13)
\lvec(164 13)
\ifill f:0
\move(166 12)
\lvec(167 12)
\lvec(167 13)
\lvec(166 13)
\ifill f:0
\move(168 12)
\lvec(171 12)
\lvec(171 13)
\lvec(168 13)
\ifill f:0
\move(172 12)
\lvec(173 12)
\lvec(173 13)
\lvec(172 13)
\ifill f:0
\move(174 12)
\lvec(184 12)
\lvec(184 13)
\lvec(174 13)
\ifill f:0
\move(185 12)
\lvec(195 12)
\lvec(195 13)
\lvec(185 13)
\ifill f:0
\move(196 12)
\lvec(197 12)
\lvec(197 13)
\lvec(196 13)
\ifill f:0
\move(198 12)
\lvec(204 12)
\lvec(204 13)
\lvec(198 13)
\ifill f:0
\move(205 12)
\lvec(211 12)
\lvec(211 13)
\lvec(205 13)
\ifill f:0
\move(212 12)
\lvec(218 12)
\lvec(218 13)
\lvec(212 13)
\ifill f:0
\move(219 12)
\lvec(230 12)
\lvec(230 13)
\lvec(219 13)
\ifill f:0
\move(231 12)
\lvec(232 12)
\lvec(232 13)
\lvec(231 13)
\ifill f:0
\move(233 12)
\lvec(235 12)
\lvec(235 13)
\lvec(233 13)
\ifill f:0
\move(236 12)
\lvec(240 12)
\lvec(240 13)
\lvec(236 13)
\ifill f:0
\move(241 12)
\lvec(247 12)
\lvec(247 13)
\lvec(241 13)
\ifill f:0
\move(248 12)
\lvec(250 12)
\lvec(250 13)
\lvec(248 13)
\ifill f:0
\move(251 12)
\lvec(255 12)
\lvec(255 13)
\lvec(251 13)
\ifill f:0
\move(256 12)
\lvec(257 12)
\lvec(257 13)
\lvec(256 13)
\ifill f:0
\move(258 12)
\lvec(260 12)
\lvec(260 13)
\lvec(258 13)
\ifill f:0
\move(261 12)
\lvec(268 12)
\lvec(268 13)
\lvec(261 13)
\ifill f:0
\move(269 12)
\lvec(273 12)
\lvec(273 13)
\lvec(269 13)
\ifill f:0
\move(274 12)
\lvec(281 12)
\lvec(281 13)
\lvec(274 13)
\ifill f:0
\move(282 12)
\lvec(284 12)
\lvec(284 13)
\lvec(282 13)
\ifill f:0
\move(285 12)
\lvec(292 12)
\lvec(292 13)
\lvec(285 13)
\ifill f:0
\move(293 12)
\lvec(303 12)
\lvec(303 13)
\lvec(293 13)
\ifill f:0
\move(304 12)
\lvec(306 12)
\lvec(306 13)
\lvec(304 13)
\ifill f:0
\move(307 12)
\lvec(311 12)
\lvec(311 13)
\lvec(307 13)
\ifill f:0
\move(312 12)
\lvec(314 12)
\lvec(314 13)
\lvec(312 13)
\ifill f:0
\move(315 12)
\lvec(317 12)
\lvec(317 13)
\lvec(315 13)
\ifill f:0
\move(318 12)
\lvec(320 12)
\lvec(320 13)
\lvec(318 13)
\ifill f:0
\move(321 12)
\lvec(323 12)
\lvec(323 13)
\lvec(321 13)
\ifill f:0
\move(324 12)
\lvec(340 12)
\lvec(340 13)
\lvec(324 13)
\ifill f:0
\move(341 12)
\lvec(349 12)
\lvec(349 13)
\lvec(341 13)
\ifill f:0
\move(350 12)
\lvec(352 12)
\lvec(352 13)
\lvec(350 13)
\ifill f:0
\move(353 12)
\lvec(355 12)
\lvec(355 13)
\lvec(353 13)
\ifill f:0
\move(356 12)
\lvec(358 12)
\lvec(358 13)
\lvec(356 13)
\ifill f:0
\move(359 12)
\lvec(364 12)
\lvec(364 13)
\lvec(359 13)
\ifill f:0
\move(365 12)
\lvec(367 12)
\lvec(367 13)
\lvec(365 13)
\ifill f:0
\move(368 12)
\lvec(380 12)
\lvec(380 13)
\lvec(368 13)
\ifill f:0
\move(381 12)
\lvec(383 12)
\lvec(383 13)
\lvec(381 13)
\ifill f:0
\move(384 12)
\lvec(386 12)
\lvec(386 13)
\lvec(384 13)
\ifill f:0
\move(387 12)
\lvec(389 12)
\lvec(389 13)
\lvec(387 13)
\ifill f:0
\move(390 12)
\lvec(392 12)
\lvec(392 13)
\lvec(390 13)
\ifill f:0
\move(393 12)
\lvec(399 12)
\lvec(399 13)
\lvec(393 13)
\ifill f:0
\move(400 12)
\lvec(402 12)
\lvec(402 13)
\lvec(400 13)
\ifill f:0
\move(403 12)
\lvec(405 12)
\lvec(405 13)
\lvec(403 13)
\ifill f:0
\move(406 12)
\lvec(415 12)
\lvec(415 13)
\lvec(406 13)
\ifill f:0
\move(416 12)
\lvec(418 12)
\lvec(418 13)
\lvec(416 13)
\ifill f:0
\move(419 12)
\lvec(425 12)
\lvec(425 13)
\lvec(419 13)
\ifill f:0
\move(426 12)
\lvec(428 12)
\lvec(428 13)
\lvec(426 13)
\ifill f:0
\move(429 12)
\lvec(431 12)
\lvec(431 13)
\lvec(429 13)
\ifill f:0
\move(432 12)
\lvec(438 12)
\lvec(438 13)
\lvec(432 13)
\ifill f:0
\move(439 12)
\lvec(445 12)
\lvec(445 13)
\lvec(439 13)
\ifill f:0
\move(446 12)
\lvec(448 12)
\lvec(448 13)
\lvec(446 13)
\ifill f:0
\move(449 12)
\lvec(451 12)
\lvec(451 13)
\lvec(449 13)
\ifill f:0
\move(11 13)
\lvec(14 13)
\lvec(14 14)
\lvec(11 14)
\ifill f:0
\move(16 13)
\lvec(17 13)
\lvec(17 14)
\lvec(16 14)
\ifill f:0
\move(18 13)
\lvec(22 13)
\lvec(22 14)
\lvec(18 14)
\ifill f:0
\move(23 13)
\lvec(26 13)
\lvec(26 14)
\lvec(23 14)
\ifill f:0
\move(29 13)
\lvec(30 13)
\lvec(30 14)
\lvec(29 14)
\ifill f:0
\move(33 13)
\lvec(34 13)
\lvec(34 14)
\lvec(33 14)
\ifill f:0
\move(35 13)
\lvec(37 13)
\lvec(37 14)
\lvec(35 14)
\ifill f:0
\move(38 13)
\lvec(39 13)
\lvec(39 14)
\lvec(38 14)
\ifill f:0
\move(40 13)
\lvec(42 13)
\lvec(42 14)
\lvec(40 14)
\ifill f:0
\move(45 13)
\lvec(47 13)
\lvec(47 14)
\lvec(45 14)
\ifill f:0
\move(49 13)
\lvec(50 13)
\lvec(50 14)
\lvec(49 14)
\ifill f:0
\move(56 13)
\lvec(61 13)
\lvec(61 14)
\lvec(56 14)
\ifill f:0
\move(62 13)
\lvec(67 13)
\lvec(67 14)
\lvec(62 14)
\ifill f:0
\move(68 13)
\lvec(72 13)
\lvec(72 14)
\lvec(68 14)
\ifill f:0
\move(73 13)
\lvec(82 13)
\lvec(82 14)
\lvec(73 14)
\ifill f:0
\move(83 13)
\lvec(85 13)
\lvec(85 14)
\lvec(83 14)
\ifill f:0
\move(87 13)
\lvec(89 13)
\lvec(89 14)
\lvec(87 14)
\ifill f:0
\move(90 13)
\lvec(92 13)
\lvec(92 14)
\lvec(90 14)
\ifill f:0
\move(93 13)
\lvec(95 13)
\lvec(95 14)
\lvec(93 14)
\ifill f:0
\move(97 13)
\lvec(110 13)
\lvec(110 14)
\lvec(97 14)
\ifill f:0
\move(111 13)
\lvec(124 13)
\lvec(124 14)
\lvec(111 14)
\ifill f:0
\move(125 13)
\lvec(126 13)
\lvec(126 14)
\lvec(125 14)
\ifill f:0
\move(127 13)
\lvec(131 13)
\lvec(131 14)
\lvec(127 14)
\ifill f:0
\move(132 13)
\lvec(138 13)
\lvec(138 14)
\lvec(132 14)
\ifill f:0
\move(139 13)
\lvec(143 13)
\lvec(143 14)
\lvec(139 14)
\ifill f:0
\move(144 13)
\lvec(145 13)
\lvec(145 14)
\lvec(144 14)
\ifill f:0
\move(146 13)
\lvec(154 13)
\lvec(154 14)
\lvec(146 14)
\ifill f:0
\move(155 13)
\lvec(167 13)
\lvec(167 14)
\lvec(155 14)
\ifill f:0
\move(168 13)
\lvec(171 13)
\lvec(171 14)
\lvec(168 14)
\ifill f:0
\move(172 13)
\lvec(173 13)
\lvec(173 14)
\lvec(172 14)
\ifill f:0
\move(174 13)
\lvec(175 13)
\lvec(175 14)
\lvec(174 14)
\ifill f:0
\move(176 13)
\lvec(177 13)
\lvec(177 14)
\lvec(176 14)
\ifill f:0
\move(178 13)
\lvec(179 13)
\lvec(179 14)
\lvec(178 14)
\ifill f:0
\move(180 13)
\lvec(185 13)
\lvec(185 14)
\lvec(180 14)
\ifill f:0
\move(186 13)
\lvec(187 13)
\lvec(187 14)
\lvec(186 14)
\ifill f:0
\move(188 13)
\lvec(189 13)
\lvec(189 14)
\lvec(188 14)
\ifill f:0
\move(190 13)
\lvec(191 13)
\lvec(191 14)
\lvec(190 14)
\ifill f:0
\move(192 13)
\lvec(193 13)
\lvec(193 14)
\lvec(192 14)
\ifill f:0
\move(194 13)
\lvec(195 13)
\lvec(195 14)
\lvec(194 14)
\ifill f:0
\move(196 13)
\lvec(197 13)
\lvec(197 14)
\lvec(196 14)
\ifill f:0
\move(198 13)
\lvec(212 13)
\lvec(212 14)
\lvec(198 14)
\ifill f:0
\move(213 13)
\lvec(223 13)
\lvec(223 14)
\lvec(213 14)
\ifill f:0
\move(224 13)
\lvec(227 13)
\lvec(227 14)
\lvec(224 14)
\ifill f:0
\move(228 13)
\lvec(232 13)
\lvec(232 14)
\lvec(228 14)
\ifill f:0
\move(233 13)
\lvec(234 13)
\lvec(234 14)
\lvec(233 14)
\ifill f:0
\move(235 13)
\lvec(241 13)
\lvec(241 14)
\lvec(235 14)
\ifill f:0
\move(242 13)
\lvec(248 13)
\lvec(248 14)
\lvec(242 14)
\ifill f:0
\move(249 13)
\lvec(250 13)
\lvec(250 14)
\lvec(249 14)
\ifill f:0
\move(251 13)
\lvec(255 13)
\lvec(255 14)
\lvec(251 14)
\ifill f:0
\move(256 13)
\lvec(257 13)
\lvec(257 14)
\lvec(256 14)
\ifill f:0
\move(258 13)
\lvec(267 13)
\lvec(267 14)
\lvec(258 14)
\ifill f:0
\move(268 13)
\lvec(272 13)
\lvec(272 14)
\lvec(268 14)
\ifill f:0
\move(273 13)
\lvec(274 13)
\lvec(274 14)
\lvec(273 14)
\ifill f:0
\move(275 13)
\lvec(279 13)
\lvec(279 14)
\lvec(275 14)
\ifill f:0
\move(280 13)
\lvec(294 13)
\lvec(294 14)
\lvec(280 14)
\ifill f:0
\move(295 13)
\lvec(297 13)
\lvec(297 14)
\lvec(295 14)
\ifill f:0
\move(298 13)
\lvec(302 13)
\lvec(302 14)
\lvec(298 14)
\ifill f:0
\move(303 13)
\lvec(307 13)
\lvec(307 14)
\lvec(303 14)
\ifill f:0
\move(308 13)
\lvec(315 13)
\lvec(315 14)
\lvec(308 14)
\ifill f:0
\move(316 13)
\lvec(320 13)
\lvec(320 14)
\lvec(316 14)
\ifill f:0
\move(321 13)
\lvec(323 13)
\lvec(323 14)
\lvec(321 14)
\ifill f:0
\move(324 13)
\lvec(331 13)
\lvec(331 14)
\lvec(324 14)
\ifill f:0
\move(332 13)
\lvec(339 13)
\lvec(339 14)
\lvec(332 14)
\ifill f:0
\move(340 13)
\lvec(342 13)
\lvec(342 14)
\lvec(340 14)
\ifill f:0
\move(343 13)
\lvec(350 13)
\lvec(350 14)
\lvec(343 14)
\ifill f:0
\move(351 13)
\lvec(353 13)
\lvec(353 14)
\lvec(351 14)
\ifill f:0
\move(354 13)
\lvec(364 13)
\lvec(364 14)
\lvec(354 14)
\ifill f:0
\move(365 13)
\lvec(367 13)
\lvec(367 14)
\lvec(365 14)
\ifill f:0
\move(368 13)
\lvec(381 13)
\lvec(381 14)
\lvec(368 14)
\ifill f:0
\move(382 13)
\lvec(384 13)
\lvec(384 14)
\lvec(382 14)
\ifill f:0
\move(385 13)
\lvec(387 13)
\lvec(387 14)
\lvec(385 14)
\ifill f:0
\move(388 13)
\lvec(390 13)
\lvec(390 14)
\lvec(388 14)
\ifill f:0
\move(391 13)
\lvec(393 13)
\lvec(393 14)
\lvec(391 14)
\ifill f:0
\move(394 13)
\lvec(396 13)
\lvec(396 14)
\lvec(394 14)
\ifill f:0
\move(397 13)
\lvec(399 13)
\lvec(399 14)
\lvec(397 14)
\ifill f:0
\move(400 13)
\lvec(402 13)
\lvec(402 14)
\lvec(400 14)
\ifill f:0
\move(403 13)
\lvec(405 13)
\lvec(405 14)
\lvec(403 14)
\ifill f:0
\move(406 13)
\lvec(414 13)
\lvec(414 14)
\lvec(406 14)
\ifill f:0
\move(415 13)
\lvec(417 13)
\lvec(417 14)
\lvec(415 14)
\ifill f:0
\move(418 13)
\lvec(420 13)
\lvec(420 14)
\lvec(418 14)
\ifill f:0
\move(421 13)
\lvec(423 13)
\lvec(423 14)
\lvec(421 14)
\ifill f:0
\move(424 13)
\lvec(426 13)
\lvec(426 14)
\lvec(424 14)
\ifill f:0
\move(427 13)
\lvec(429 13)
\lvec(429 14)
\lvec(427 14)
\ifill f:0
\move(430 13)
\lvec(432 13)
\lvec(432 14)
\lvec(430 14)
\ifill f:0
\move(433 13)
\lvec(435 13)
\lvec(435 14)
\lvec(433 14)
\ifill f:0
\move(436 13)
\lvec(438 13)
\lvec(438 14)
\lvec(436 14)
\ifill f:0
\move(439 13)
\lvec(451 13)
\lvec(451 14)
\lvec(439 14)
\ifill f:0
\move(11 14)
\lvec(12 14)
\lvec(12 15)
\lvec(11 15)
\ifill f:0
\move(16 14)
\lvec(17 14)
\lvec(17 15)
\lvec(16 15)
\ifill f:0
\move(20 14)
\lvec(21 14)
\lvec(21 15)
\lvec(20 15)
\ifill f:0
\move(22 14)
\lvec(26 14)
\lvec(26 15)
\lvec(22 15)
\ifill f:0
\move(28 14)
\lvec(29 14)
\lvec(29 15)
\lvec(28 15)
\ifill f:0
\move(34 14)
\lvec(37 14)
\lvec(37 15)
\lvec(34 15)
\ifill f:0
\move(38 14)
\lvec(47 14)
\lvec(47 15)
\lvec(38 15)
\ifill f:0
\move(49 14)
\lvec(50 14)
\lvec(50 15)
\lvec(49 15)
\ifill f:0
\move(52 14)
\lvec(54 14)
\lvec(54 15)
\lvec(52 15)
\ifill f:0
\move(58 14)
\lvec(59 14)
\lvec(59 15)
\lvec(58 15)
\ifill f:0
\move(60 14)
\lvec(82 14)
\lvec(82 15)
\lvec(60 15)
\ifill f:0
\move(83 14)
\lvec(86 14)
\lvec(86 15)
\lvec(83 15)
\ifill f:0
\move(87 14)
\lvec(90 14)
\lvec(90 15)
\lvec(87 15)
\ifill f:0
\move(91 14)
\lvec(95 14)
\lvec(95 15)
\lvec(91 15)
\ifill f:0
\move(96 14)
\lvec(98 14)
\lvec(98 15)
\lvec(96 15)
\ifill f:0
\move(99 14)
\lvec(105 14)
\lvec(105 15)
\lvec(99 15)
\ifill f:0
\move(106 14)
\lvec(118 14)
\lvec(118 15)
\lvec(106 15)
\ifill f:0
\move(119 14)
\lvec(122 14)
\lvec(122 15)
\lvec(119 15)
\ifill f:0
\move(123 14)
\lvec(124 14)
\lvec(124 15)
\lvec(123 15)
\ifill f:0
\move(125 14)
\lvec(127 14)
\lvec(127 15)
\lvec(125 15)
\ifill f:0
\move(128 14)
\lvec(132 14)
\lvec(132 15)
\lvec(128 15)
\ifill f:0
\move(133 14)
\lvec(135 14)
\lvec(135 15)
\lvec(133 15)
\ifill f:0
\move(136 14)
\lvec(138 14)
\lvec(138 15)
\lvec(136 15)
\ifill f:0
\move(139 14)
\lvec(140 14)
\lvec(140 15)
\lvec(139 15)
\ifill f:0
\move(141 14)
\lvec(143 14)
\lvec(143 15)
\lvec(141 15)
\ifill f:0
\move(144 14)
\lvec(150 14)
\lvec(150 15)
\lvec(144 15)
\ifill f:0
\move(151 14)
\lvec(155 14)
\lvec(155 15)
\lvec(151 15)
\ifill f:0
\move(156 14)
\lvec(162 14)
\lvec(162 15)
\lvec(156 15)
\ifill f:0
\move(163 14)
\lvec(176 14)
\lvec(176 15)
\lvec(163 15)
\ifill f:0
\move(177 14)
\lvec(180 14)
\lvec(180 15)
\lvec(177 15)
\ifill f:0
\move(181 14)
\lvec(191 14)
\lvec(191 15)
\lvec(181 15)
\ifill f:0
\move(192 14)
\lvec(195 14)
\lvec(195 15)
\lvec(192 15)
\ifill f:0
\move(196 14)
\lvec(197 14)
\lvec(197 15)
\lvec(196 15)
\ifill f:0
\move(198 14)
\lvec(199 14)
\lvec(199 15)
\lvec(198 15)
\ifill f:0
\move(200 14)
\lvec(201 14)
\lvec(201 15)
\lvec(200 15)
\ifill f:0
\move(202 14)
\lvec(203 14)
\lvec(203 15)
\lvec(202 15)
\ifill f:0
\move(204 14)
\lvec(205 14)
\lvec(205 15)
\lvec(204 15)
\ifill f:0
\move(206 14)
\lvec(207 14)
\lvec(207 15)
\lvec(206 15)
\ifill f:0
\move(208 14)
\lvec(213 14)
\lvec(213 15)
\lvec(208 15)
\ifill f:0
\move(214 14)
\lvec(215 14)
\lvec(215 15)
\lvec(214 15)
\ifill f:0
\move(216 14)
\lvec(217 14)
\lvec(217 15)
\lvec(216 15)
\ifill f:0
\move(218 14)
\lvec(219 14)
\lvec(219 15)
\lvec(218 15)
\ifill f:0
\move(220 14)
\lvec(221 14)
\lvec(221 15)
\lvec(220 15)
\ifill f:0
\move(222 14)
\lvec(223 14)
\lvec(223 15)
\lvec(222 15)
\ifill f:0
\move(224 14)
\lvec(227 14)
\lvec(227 15)
\lvec(224 15)
\ifill f:0
\move(228 14)
\lvec(229 14)
\lvec(229 15)
\lvec(228 15)
\ifill f:0
\move(230 14)
\lvec(242 14)
\lvec(242 15)
\lvec(230 15)
\ifill f:0
\move(243 14)
\lvec(244 14)
\lvec(244 15)
\lvec(243 15)
\ifill f:0
\move(245 14)
\lvec(255 14)
\lvec(255 15)
\lvec(245 15)
\ifill f:0
\move(256 14)
\lvec(257 14)
\lvec(257 15)
\lvec(256 15)
\ifill f:0
\move(258 14)
\lvec(264 14)
\lvec(264 15)
\lvec(258 15)
\ifill f:0
\move(265 14)
\lvec(273 14)
\lvec(273 15)
\lvec(265 15)
\ifill f:0
\move(274 14)
\lvec(280 14)
\lvec(280 15)
\lvec(274 15)
\ifill f:0
\move(281 14)
\lvec(294 14)
\lvec(294 15)
\lvec(281 15)
\ifill f:0
\move(295 14)
\lvec(296 14)
\lvec(296 15)
\lvec(295 15)
\ifill f:0
\move(297 14)
\lvec(301 14)
\lvec(301 15)
\lvec(297 15)
\ifill f:0
\move(302 14)
\lvec(306 14)
\lvec(306 15)
\lvec(302 15)
\ifill f:0
\move(307 14)
\lvec(308 14)
\lvec(308 15)
\lvec(307 15)
\ifill f:0
\move(309 14)
\lvec(313 14)
\lvec(313 15)
\lvec(309 15)
\ifill f:0
\move(314 14)
\lvec(318 14)
\lvec(318 15)
\lvec(314 15)
\ifill f:0
\move(319 14)
\lvec(320 14)
\lvec(320 15)
\lvec(319 15)
\ifill f:0
\move(321 14)
\lvec(323 14)
\lvec(323 15)
\lvec(321 15)
\ifill f:0
\move(324 14)
\lvec(325 14)
\lvec(325 15)
\lvec(324 15)
\ifill f:0
\move(326 14)
\lvec(328 14)
\lvec(328 15)
\lvec(326 15)
\ifill f:0
\move(329 14)
\lvec(333 14)
\lvec(333 15)
\lvec(329 15)
\ifill f:0
\move(334 14)
\lvec(338 14)
\lvec(338 15)
\lvec(334 15)
\ifill f:0
\move(339 14)
\lvec(343 14)
\lvec(343 15)
\lvec(339 15)
\ifill f:0
\move(344 14)
\lvec(351 14)
\lvec(351 15)
\lvec(344 15)
\ifill f:0
\move(352 14)
\lvec(356 14)
\lvec(356 15)
\lvec(352 15)
\ifill f:0
\move(357 14)
\lvec(364 14)
\lvec(364 15)
\lvec(357 15)
\ifill f:0
\move(365 14)
\lvec(369 14)
\lvec(369 15)
\lvec(365 15)
\ifill f:0
\move(370 14)
\lvec(377 14)
\lvec(377 15)
\lvec(370 15)
\ifill f:0
\move(378 14)
\lvec(380 14)
\lvec(380 15)
\lvec(378 15)
\ifill f:0
\move(381 14)
\lvec(385 14)
\lvec(385 15)
\lvec(381 15)
\ifill f:0
\move(386 14)
\lvec(388 14)
\lvec(388 15)
\lvec(386 15)
\ifill f:0
\move(389 14)
\lvec(396 14)
\lvec(396 15)
\lvec(389 15)
\ifill f:0
\move(397 14)
\lvec(399 14)
\lvec(399 15)
\lvec(397 15)
\ifill f:0
\move(400 14)
\lvec(407 14)
\lvec(407 15)
\lvec(400 15)
\ifill f:0
\move(408 14)
\lvec(410 14)
\lvec(410 15)
\lvec(408 15)
\ifill f:0
\move(411 14)
\lvec(421 14)
\lvec(421 15)
\lvec(411 15)
\ifill f:0
\move(422 14)
\lvec(424 14)
\lvec(424 15)
\lvec(422 15)
\ifill f:0
\move(425 14)
\lvec(427 14)
\lvec(427 15)
\lvec(425 15)
\ifill f:0
\move(428 14)
\lvec(444 14)
\lvec(444 15)
\lvec(428 15)
\ifill f:0
\move(445 14)
\lvec(447 14)
\lvec(447 15)
\lvec(445 15)
\ifill f:0
\move(448 14)
\lvec(450 14)
\lvec(450 15)
\lvec(448 15)
\ifill f:0
\move(12 15)
\lvec(17 15)
\lvec(17 16)
\lvec(12 16)
\ifill f:0
\move(18 15)
\lvec(23 15)
\lvec(23 16)
\lvec(18 16)
\ifill f:0
\move(24 15)
\lvec(26 15)
\lvec(26 16)
\lvec(24 16)
\ifill f:0
\move(27 15)
\lvec(28 15)
\lvec(28 16)
\lvec(27 16)
\ifill f:0
\move(35 15)
\lvec(37 15)
\lvec(37 16)
\lvec(35 16)
\ifill f:0
\move(38 15)
\lvec(42 15)
\lvec(42 16)
\lvec(38 16)
\ifill f:0
\move(43 15)
\lvec(48 15)
\lvec(48 16)
\lvec(43 16)
\ifill f:0
\move(49 15)
\lvec(53 15)
\lvec(53 16)
\lvec(49 16)
\ifill f:0
\move(54 15)
\lvec(55 15)
\lvec(55 16)
\lvec(54 16)
\ifill f:0
\move(56 15)
\lvec(65 15)
\lvec(65 16)
\lvec(56 16)
\ifill f:0
\move(66 15)
\lvec(72 15)
\lvec(72 16)
\lvec(66 16)
\ifill f:0
\move(73 15)
\lvec(74 15)
\lvec(74 16)
\lvec(73 16)
\ifill f:0
\move(75 15)
\lvec(82 15)
\lvec(82 16)
\lvec(75 16)
\ifill f:0
\move(84 15)
\lvec(88 15)
\lvec(88 16)
\lvec(84 16)
\ifill f:0
\move(89 15)
\lvec(93 15)
\lvec(93 16)
\lvec(89 16)
\ifill f:0
\move(95 15)
\lvec(98 15)
\lvec(98 16)
\lvec(95 16)
\ifill f:0
\move(99 15)
\lvec(102 15)
\lvec(102 16)
\lvec(99 16)
\ifill f:0
\move(103 15)
\lvec(110 15)
\lvec(110 16)
\lvec(103 16)
\ifill f:0
\move(111 15)
\lvec(122 15)
\lvec(122 16)
\lvec(111 16)
\ifill f:0
\move(123 15)
\lvec(125 15)
\lvec(125 16)
\lvec(123 16)
\ifill f:0
\move(126 15)
\lvec(128 15)
\lvec(128 16)
\lvec(126 16)
\ifill f:0
\move(129 15)
\lvec(131 15)
\lvec(131 16)
\lvec(129 16)
\ifill f:0
\move(132 15)
\lvec(134 15)
\lvec(134 16)
\lvec(132 16)
\ifill f:0
\move(135 15)
\lvec(137 15)
\lvec(137 16)
\lvec(135 16)
\ifill f:0
\move(138 15)
\lvec(148 15)
\lvec(148 16)
\lvec(138 16)
\ifill f:0
\move(149 15)
\lvec(156 15)
\lvec(156 16)
\lvec(149 16)
\ifill f:0
\move(157 15)
\lvec(159 15)
\lvec(159 16)
\lvec(157 16)
\ifill f:0
\move(160 15)
\lvec(172 15)
\lvec(172 16)
\lvec(160 16)
\ifill f:0
\move(173 15)
\lvec(174 15)
\lvec(174 16)
\lvec(173 16)
\ifill f:0
\move(175 15)
\lvec(181 15)
\lvec(181 16)
\lvec(175 16)
\ifill f:0
\move(182 15)
\lvec(188 15)
\lvec(188 16)
\lvec(182 16)
\ifill f:0
\move(189 15)
\lvec(197 15)
\lvec(197 16)
\lvec(189 16)
\ifill f:0
\move(198 15)
\lvec(208 15)
\lvec(208 16)
\lvec(198 16)
\ifill f:0
\move(209 15)
\lvec(223 15)
\lvec(223 16)
\lvec(209 16)
\ifill f:0
\move(224 15)
\lvec(227 15)
\lvec(227 16)
\lvec(224 16)
\ifill f:0
\move(228 15)
\lvec(229 15)
\lvec(229 16)
\lvec(228 16)
\ifill f:0
\move(230 15)
\lvec(231 15)
\lvec(231 16)
\lvec(230 16)
\ifill f:0
\move(232 15)
\lvec(233 15)
\lvec(233 16)
\lvec(232 16)
\ifill f:0
\move(234 15)
\lvec(235 15)
\lvec(235 16)
\lvec(234 16)
\ifill f:0
\move(236 15)
\lvec(237 15)
\lvec(237 16)
\lvec(236 16)
\ifill f:0
\move(238 15)
\lvec(243 15)
\lvec(243 16)
\lvec(238 16)
\ifill f:0
\move(244 15)
\lvec(245 15)
\lvec(245 16)
\lvec(244 16)
\ifill f:0
\move(246 15)
\lvec(247 15)
\lvec(247 16)
\lvec(246 16)
\ifill f:0
\move(248 15)
\lvec(249 15)
\lvec(249 16)
\lvec(248 16)
\ifill f:0
\move(250 15)
\lvec(251 15)
\lvec(251 16)
\lvec(250 16)
\ifill f:0
\move(252 15)
\lvec(253 15)
\lvec(253 16)
\lvec(252 16)
\ifill f:0
\move(254 15)
\lvec(255 15)
\lvec(255 16)
\lvec(254 16)
\ifill f:0
\move(256 15)
\lvec(257 15)
\lvec(257 16)
\lvec(256 16)
\ifill f:0
\move(258 15)
\lvec(274 15)
\lvec(274 16)
\lvec(258 16)
\ifill f:0
\move(275 15)
\lvec(287 15)
\lvec(287 16)
\lvec(275 16)
\ifill f:0
\move(288 15)
\lvec(291 15)
\lvec(291 16)
\lvec(288 16)
\ifill f:0
\move(292 15)
\lvec(293 15)
\lvec(293 16)
\lvec(292 16)
\ifill f:0
\move(294 15)
\lvec(298 15)
\lvec(298 16)
\lvec(294 16)
\ifill f:0
\move(299 15)
\lvec(307 15)
\lvec(307 16)
\lvec(299 16)
\ifill f:0
\move(308 15)
\lvec(309 15)
\lvec(309 16)
\lvec(308 16)
\ifill f:0
\move(310 15)
\lvec(316 15)
\lvec(316 16)
\lvec(310 16)
\ifill f:0
\move(317 15)
\lvec(323 15)
\lvec(323 16)
\lvec(317 16)
\ifill f:0
\move(324 15)
\lvec(325 15)
\lvec(325 16)
\lvec(324 16)
\ifill f:0
\move(326 15)
\lvec(330 15)
\lvec(330 16)
\lvec(326 16)
\ifill f:0
\move(331 15)
\lvec(337 15)
\lvec(337 16)
\lvec(331 16)
\ifill f:0
\move(338 15)
\lvec(342 15)
\lvec(342 16)
\lvec(338 16)
\ifill f:0
\move(343 15)
\lvec(349 15)
\lvec(349 16)
\lvec(343 16)
\ifill f:0
\move(350 15)
\lvec(354 15)
\lvec(354 16)
\lvec(350 16)
\ifill f:0
\move(355 15)
\lvec(366 15)
\lvec(366 16)
\lvec(355 16)
\ifill f:0
\move(367 15)
\lvec(371 15)
\lvec(371 16)
\lvec(367 16)
\ifill f:0
\move(372 15)
\lvec(376 15)
\lvec(376 16)
\lvec(372 16)
\ifill f:0
\move(377 15)
\lvec(381 15)
\lvec(381 16)
\lvec(377 16)
\ifill f:0
\move(382 15)
\lvec(386 15)
\lvec(386 16)
\lvec(382 16)
\ifill f:0
\move(387 15)
\lvec(391 15)
\lvec(391 16)
\lvec(387 16)
\ifill f:0
\move(392 15)
\lvec(399 15)
\lvec(399 16)
\lvec(392 16)
\ifill f:0
\move(400 15)
\lvec(401 15)
\lvec(401 16)
\lvec(400 16)
\ifill f:0
\move(402 15)
\lvec(404 15)
\lvec(404 16)
\lvec(402 16)
\ifill f:0
\move(405 15)
\lvec(412 15)
\lvec(412 16)
\lvec(405 16)
\ifill f:0
\move(413 15)
\lvec(417 15)
\lvec(417 16)
\lvec(413 16)
\ifill f:0
\move(418 15)
\lvec(420 15)
\lvec(420 16)
\lvec(418 16)
\ifill f:0
\move(421 15)
\lvec(425 15)
\lvec(425 16)
\lvec(421 16)
\ifill f:0
\move(426 15)
\lvec(433 15)
\lvec(433 16)
\lvec(426 16)
\ifill f:0
\move(434 15)
\lvec(436 15)
\lvec(436 16)
\lvec(434 16)
\ifill f:0
\move(437 15)
\lvec(444 15)
\lvec(444 16)
\lvec(437 16)
\ifill f:0
\move(445 15)
\lvec(451 15)
\lvec(451 16)
\lvec(445 16)
\ifill f:0
\move(11 16)
\lvec(12 16)
\lvec(12 17)
\lvec(11 17)
\ifill f:0
\move(14 16)
\lvec(17 16)
\lvec(17 17)
\lvec(14 17)
\ifill f:0
\move(18 16)
\lvec(22 16)
\lvec(22 17)
\lvec(18 17)
\ifill f:0
\move(23 16)
\lvec(26 16)
\lvec(26 17)
\lvec(23 17)
\ifill f:0
\move(27 16)
\lvec(28 16)
\lvec(28 17)
\lvec(27 17)
\ifill f:0
\move(35 16)
\lvec(37 16)
\lvec(37 17)
\lvec(35 17)
\ifill f:0
\move(38 16)
\lvec(48 16)
\lvec(48 17)
\lvec(38 17)
\ifill f:0
\move(49 16)
\lvec(53 16)
\lvec(53 17)
\lvec(49 17)
\ifill f:0
\move(57 16)
\lvec(65 16)
\lvec(65 17)
\lvec(57 17)
\ifill f:0
\move(66 16)
\lvec(82 16)
\lvec(82 17)
\lvec(66 17)
\ifill f:0
\move(86 16)
\lvec(90 16)
\lvec(90 17)
\lvec(86 17)
\ifill f:0
\move(93 16)
\lvec(98 16)
\lvec(98 17)
\lvec(93 17)
\ifill f:0
\move(99 16)
\lvec(122 16)
\lvec(122 17)
\lvec(99 17)
\ifill f:0
\move(123 16)
\lvec(125 16)
\lvec(125 17)
\lvec(123 17)
\ifill f:0
\move(126 16)
\lvec(129 16)
\lvec(129 17)
\lvec(126 17)
\ifill f:0
\move(130 16)
\lvec(132 16)
\lvec(132 17)
\lvec(130 17)
\ifill f:0
\move(133 16)
\lvec(139 16)
\lvec(139 17)
\lvec(133 17)
\ifill f:0
\move(140 16)
\lvec(166 16)
\lvec(166 17)
\lvec(140 17)
\ifill f:0
\move(167 16)
\lvec(172 16)
\lvec(172 17)
\lvec(167 17)
\ifill f:0
\move(173 16)
\lvec(177 16)
\lvec(177 17)
\lvec(173 17)
\ifill f:0
\move(178 16)
\lvec(182 16)
\lvec(182 17)
\lvec(178 17)
\ifill f:0
\move(183 16)
\lvec(190 16)
\lvec(190 17)
\lvec(183 17)
\ifill f:0
\move(191 16)
\lvec(197 16)
\lvec(197 17)
\lvec(191 17)
\ifill f:0
\move(198 16)
\lvec(202 16)
\lvec(202 17)
\lvec(198 17)
\ifill f:0
\move(203 16)
\lvec(209 16)
\lvec(209 17)
\lvec(203 17)
\ifill f:0
\move(210 16)
\lvec(216 16)
\lvec(216 17)
\lvec(210 17)
\ifill f:0
\move(217 16)
\lvec(227 16)
\lvec(227 17)
\lvec(217 17)
\ifill f:0
\move(228 16)
\lvec(238 16)
\lvec(238 17)
\lvec(228 17)
\ifill f:0
\move(239 16)
\lvec(255 16)
\lvec(255 17)
\lvec(239 17)
\ifill f:0
\move(256 16)
\lvec(257 16)
\lvec(257 17)
\lvec(256 17)
\ifill f:0
\move(258 16)
\lvec(259 16)
\lvec(259 17)
\lvec(258 17)
\ifill f:0
\move(260 16)
\lvec(261 16)
\lvec(261 17)
\lvec(260 17)
\ifill f:0
\move(262 16)
\lvec(263 16)
\lvec(263 17)
\lvec(262 17)
\ifill f:0
\move(264 16)
\lvec(265 16)
\lvec(265 17)
\lvec(264 17)
\ifill f:0
\move(266 16)
\lvec(267 16)
\lvec(267 17)
\lvec(266 17)
\ifill f:0
\move(268 16)
\lvec(269 16)
\lvec(269 17)
\lvec(268 17)
\ifill f:0
\move(270 16)
\lvec(275 16)
\lvec(275 17)
\lvec(270 17)
\ifill f:0
\move(276 16)
\lvec(277 16)
\lvec(277 17)
\lvec(276 17)
\ifill f:0
\move(278 16)
\lvec(279 16)
\lvec(279 17)
\lvec(278 17)
\ifill f:0
\move(280 16)
\lvec(281 16)
\lvec(281 17)
\lvec(280 17)
\ifill f:0
\move(282 16)
\lvec(283 16)
\lvec(283 17)
\lvec(282 17)
\ifill f:0
\move(284 16)
\lvec(285 16)
\lvec(285 17)
\lvec(284 17)
\ifill f:0
\move(286 16)
\lvec(287 16)
\lvec(287 17)
\lvec(286 17)
\ifill f:0
\move(288 16)
\lvec(291 16)
\lvec(291 17)
\lvec(288 17)
\ifill f:0
\move(292 16)
\lvec(293 16)
\lvec(293 17)
\lvec(292 17)
\ifill f:0
\move(294 16)
\lvec(308 16)
\lvec(308 17)
\lvec(294 17)
\ifill f:0
\move(309 16)
\lvec(310 16)
\lvec(310 17)
\lvec(309 17)
\ifill f:0
\move(311 16)
\lvec(321 16)
\lvec(321 17)
\lvec(311 17)
\ifill f:0
\move(322 16)
\lvec(323 16)
\lvec(323 17)
\lvec(322 17)
\ifill f:0
\move(324 16)
\lvec(325 16)
\lvec(325 17)
\lvec(324 17)
\ifill f:0
\move(326 16)
\lvec(334 16)
\lvec(334 17)
\lvec(326 17)
\ifill f:0
\move(335 16)
\lvec(343 16)
\lvec(343 17)
\lvec(335 17)
\ifill f:0
\move(344 16)
\lvec(352 16)
\lvec(352 17)
\lvec(344 17)
\ifill f:0
\move(353 16)
\lvec(359 16)
\lvec(359 17)
\lvec(353 17)
\ifill f:0
\move(360 16)
\lvec(363 16)
\lvec(363 17)
\lvec(360 17)
\ifill f:0
\move(364 16)
\lvec(368 16)
\lvec(368 17)
\lvec(364 17)
\ifill f:0
\move(369 16)
\lvec(375 16)
\lvec(375 17)
\lvec(369 17)
\ifill f:0
\move(376 16)
\lvec(380 16)
\lvec(380 17)
\lvec(376 17)
\ifill f:0
\move(381 16)
\lvec(382 16)
\lvec(382 17)
\lvec(381 17)
\ifill f:0
\move(383 16)
\lvec(387 16)
\lvec(387 17)
\lvec(383 17)
\ifill f:0
\move(388 16)
\lvec(394 16)
\lvec(394 17)
\lvec(388 17)
\ifill f:0
\move(395 16)
\lvec(399 16)
\lvec(399 17)
\lvec(395 17)
\ifill f:0
\move(400 16)
\lvec(401 16)
\lvec(401 17)
\lvec(400 17)
\ifill f:0
\move(402 16)
\lvec(404 16)
\lvec(404 17)
\lvec(402 17)
\ifill f:0
\move(405 16)
\lvec(411 16)
\lvec(411 17)
\lvec(405 17)
\ifill f:0
\move(412 16)
\lvec(416 16)
\lvec(416 17)
\lvec(412 17)
\ifill f:0
\move(417 16)
\lvec(421 16)
\lvec(421 17)
\lvec(417 17)
\ifill f:0
\move(422 16)
\lvec(426 16)
\lvec(426 17)
\lvec(422 17)
\ifill f:0
\move(427 16)
\lvec(431 16)
\lvec(431 17)
\lvec(427 17)
\ifill f:0
\move(432 16)
\lvec(436 16)
\lvec(436 17)
\lvec(432 17)
\ifill f:0
\move(437 16)
\lvec(446 16)
\lvec(446 17)
\lvec(437 17)
\ifill f:0
\move(447 16)
\lvec(449 16)
\lvec(449 17)
\lvec(447 17)
\ifill f:0
\move(450 16)
\lvec(451 16)
\lvec(451 17)
\lvec(450 17)
\ifill f:0
\move(11 17)
\lvec(14 17)
\lvec(14 18)
\lvec(11 18)
\ifill f:0
\move(15 17)
\lvec(17 17)
\lvec(17 18)
\lvec(15 18)
\ifill f:0
\move(18 17)
\lvec(19 17)
\lvec(19 18)
\lvec(18 18)
\ifill f:0
\move(20 17)
\lvec(22 17)
\lvec(22 18)
\lvec(20 18)
\ifill f:0
\move(23 17)
\lvec(26 17)
\lvec(26 18)
\lvec(23 18)
\ifill f:0
\move(28 17)
\lvec(29 17)
\lvec(29 18)
\lvec(28 18)
\ifill f:0
\move(36 17)
\lvec(37 17)
\lvec(37 18)
\lvec(36 18)
\ifill f:0
\move(38 17)
\lvec(39 17)
\lvec(39 18)
\lvec(38 18)
\ifill f:0
\move(40 17)
\lvec(42 17)
\lvec(42 18)
\lvec(40 18)
\ifill f:0
\move(43 17)
\lvec(50 17)
\lvec(50 18)
\lvec(43 18)
\ifill f:0
\move(51 17)
\lvec(52 17)
\lvec(52 18)
\lvec(51 18)
\ifill f:0
\move(54 17)
\lvec(59 17)
\lvec(59 18)
\lvec(54 18)
\ifill f:0
\move(60 17)
\lvec(65 17)
\lvec(65 18)
\lvec(60 18)
\ifill f:0
\move(66 17)
\lvec(71 17)
\lvec(71 18)
\lvec(66 18)
\ifill f:0
\move(76 17)
\lvec(78 17)
\lvec(78 18)
\lvec(76 18)
\ifill f:0
\move(81 17)
\lvec(82 17)
\lvec(82 18)
\lvec(81 18)
\ifill f:0
\move(87 17)
\lvec(96 17)
\lvec(96 18)
\lvec(87 18)
\ifill f:0
\move(97 17)
\lvec(98 17)
\lvec(98 18)
\lvec(97 18)
\ifill f:0
\move(99 17)
\lvec(104 17)
\lvec(104 18)
\lvec(99 18)
\ifill f:0
\move(105 17)
\lvec(110 17)
\lvec(110 18)
\lvec(105 18)
\ifill f:0
\move(111 17)
\lvec(122 17)
\lvec(122 18)
\lvec(111 18)
\ifill f:0
\move(123 17)
\lvec(126 17)
\lvec(126 18)
\lvec(123 18)
\ifill f:0
\move(127 17)
\lvec(130 17)
\lvec(130 18)
\lvec(127 18)
\ifill f:0
\move(131 17)
\lvec(134 17)
\lvec(134 18)
\lvec(131 18)
\ifill f:0
\move(135 17)
\lvec(138 17)
\lvec(138 18)
\lvec(135 18)
\ifill f:0
\move(139 17)
\lvec(142 17)
\lvec(142 18)
\lvec(139 18)
\ifill f:0
\move(143 17)
\lvec(156 17)
\lvec(156 18)
\lvec(143 18)
\ifill f:0
\move(157 17)
\lvec(163 17)
\lvec(163 18)
\lvec(157 18)
\ifill f:0
\move(164 17)
\lvec(166 17)
\lvec(166 18)
\lvec(164 18)
\ifill f:0
\move(167 17)
\lvec(170 17)
\lvec(170 18)
\lvec(167 18)
\ifill f:0
\move(171 17)
\lvec(172 17)
\lvec(172 18)
\lvec(171 18)
\ifill f:0
\move(174 17)
\lvec(175 17)
\lvec(175 18)
\lvec(174 18)
\ifill f:0
\move(176 17)
\lvec(189 17)
\lvec(189 18)
\lvec(176 18)
\ifill f:0
\move(190 17)
\lvec(192 17)
\lvec(192 18)
\lvec(190 18)
\ifill f:0
\move(193 17)
\lvec(200 17)
\lvec(200 18)
\lvec(193 18)
\ifill f:0
\move(201 17)
\lvec(205 17)
\lvec(205 18)
\lvec(201 18)
\ifill f:0
\move(206 17)
\lvec(210 17)
\lvec(210 18)
\lvec(206 18)
\ifill f:0
\move(211 17)
\lvec(215 17)
\lvec(215 18)
\lvec(211 18)
\ifill f:0
\move(216 17)
\lvec(220 17)
\lvec(220 18)
\lvec(216 18)
\ifill f:0
\move(221 17)
\lvec(232 17)
\lvec(232 18)
\lvec(221 18)
\ifill f:0
\move(233 17)
\lvec(239 17)
\lvec(239 18)
\lvec(233 18)
\ifill f:0
\move(240 17)
\lvec(246 17)
\lvec(246 18)
\lvec(240 18)
\ifill f:0
\move(247 17)
\lvec(248 17)
\lvec(248 18)
\lvec(247 18)
\ifill f:0
\move(249 17)
\lvec(255 17)
\lvec(255 18)
\lvec(249 18)
\ifill f:0
\move(256 17)
\lvec(257 17)
\lvec(257 18)
\lvec(256 18)
\ifill f:0
\move(258 17)
\lvec(270 17)
\lvec(270 18)
\lvec(258 18)
\ifill f:0
\move(271 17)
\lvec(287 17)
\lvec(287 18)
\lvec(271 18)
\ifill f:0
\move(288 17)
\lvec(291 17)
\lvec(291 18)
\lvec(288 18)
\ifill f:0
\move(292 17)
\lvec(293 17)
\lvec(293 18)
\lvec(292 18)
\ifill f:0
\move(294 17)
\lvec(295 17)
\lvec(295 18)
\lvec(294 18)
\ifill f:0
\move(296 17)
\lvec(297 17)
\lvec(297 18)
\lvec(296 18)
\ifill f:0
\move(298 17)
\lvec(299 17)
\lvec(299 18)
\lvec(298 18)
\ifill f:0
\move(300 17)
\lvec(301 17)
\lvec(301 18)
\lvec(300 18)
\ifill f:0
\move(302 17)
\lvec(303 17)
\lvec(303 18)
\lvec(302 18)
\ifill f:0
\move(304 17)
\lvec(309 17)
\lvec(309 18)
\lvec(304 18)
\ifill f:0
\move(310 17)
\lvec(311 17)
\lvec(311 18)
\lvec(310 18)
\ifill f:0
\move(312 17)
\lvec(313 17)
\lvec(313 18)
\lvec(312 18)
\ifill f:0
\move(314 17)
\lvec(315 17)
\lvec(315 18)
\lvec(314 18)
\ifill f:0
\move(316 17)
\lvec(317 17)
\lvec(317 18)
\lvec(316 18)
\ifill f:0
\move(318 17)
\lvec(319 17)
\lvec(319 18)
\lvec(318 18)
\ifill f:0
\move(320 17)
\lvec(321 17)
\lvec(321 18)
\lvec(320 18)
\ifill f:0
\move(322 17)
\lvec(323 17)
\lvec(323 18)
\lvec(322 18)
\ifill f:0
\move(324 17)
\lvec(325 17)
\lvec(325 18)
\lvec(324 18)
\ifill f:0
\move(326 17)
\lvec(344 17)
\lvec(344 18)
\lvec(326 18)
\ifill f:0
\move(345 17)
\lvec(346 17)
\lvec(346 18)
\lvec(345 18)
\ifill f:0
\move(347 17)
\lvec(359 17)
\lvec(359 18)
\lvec(347 18)
\ifill f:0
\move(360 17)
\lvec(363 17)
\lvec(363 18)
\lvec(360 18)
\ifill f:0
\move(364 17)
\lvec(365 17)
\lvec(365 18)
\lvec(364 18)
\ifill f:0
\move(366 17)
\lvec(370 17)
\lvec(370 18)
\lvec(366 18)
\ifill f:0
\move(371 17)
\lvec(372 17)
\lvec(372 18)
\lvec(371 18)
\ifill f:0
\move(373 17)
\lvec(381 17)
\lvec(381 18)
\lvec(373 18)
\ifill f:0
\move(382 17)
\lvec(383 17)
\lvec(383 18)
\lvec(382 18)
\ifill f:0
\move(384 17)
\lvec(390 17)
\lvec(390 18)
\lvec(384 18)
\ifill f:0
\move(391 17)
\lvec(392 17)
\lvec(392 18)
\lvec(391 18)
\ifill f:0
\move(393 17)
\lvec(399 17)
\lvec(399 18)
\lvec(393 18)
\ifill f:0
\move(400 17)
\lvec(401 17)
\lvec(401 18)
\lvec(400 18)
\ifill f:0
\move(402 17)
\lvec(406 17)
\lvec(406 18)
\lvec(402 18)
\ifill f:0
\move(407 17)
\lvec(413 17)
\lvec(413 18)
\lvec(407 18)
\ifill f:0
\move(414 17)
\lvec(420 17)
\lvec(420 18)
\lvec(414 18)
\ifill f:0
\move(421 17)
\lvec(427 17)
\lvec(427 18)
\lvec(421 18)
\ifill f:0
\move(428 17)
\lvec(434 17)
\lvec(434 18)
\lvec(428 18)
\ifill f:0
\move(435 17)
\lvec(446 17)
\lvec(446 18)
\lvec(435 18)
\ifill f:0
\move(447 17)
\lvec(451 17)
\lvec(451 18)
\lvec(447 18)
\ifill f:0
\move(13 18)
\lvec(14 18)
\lvec(14 19)
\lvec(13 19)
\ifill f:0
\move(15 18)
\lvec(17 18)
\lvec(17 19)
\lvec(15 19)
\ifill f:0
\move(19 18)
\lvec(26 18)
\lvec(26 19)
\lvec(19 19)
\ifill f:0
\move(28 18)
\lvec(30 18)
\lvec(30 19)
\lvec(28 19)
\ifill f:0
\move(36 18)
\lvec(37 18)
\lvec(37 19)
\lvec(36 19)
\ifill f:0
\move(38 18)
\lvec(39 18)
\lvec(39 19)
\lvec(38 19)
\ifill f:0
\move(40 18)
\lvec(46 18)
\lvec(46 19)
\lvec(40 19)
\ifill f:0
\move(47 18)
\lvec(50 18)
\lvec(50 19)
\lvec(47 19)
\ifill f:0
\move(53 18)
\lvec(55 18)
\lvec(55 19)
\lvec(53 19)
\ifill f:0
\move(57 18)
\lvec(61 18)
\lvec(61 19)
\lvec(57 19)
\ifill f:0
\move(62 18)
\lvec(65 18)
\lvec(65 19)
\lvec(62 19)
\ifill f:0
\move(66 18)
\lvec(79 18)
\lvec(79 19)
\lvec(66 19)
\ifill f:0
\move(81 18)
\lvec(82 18)
\lvec(82 19)
\lvec(81 19)
\ifill f:0
\move(85 18)
\lvec(87 18)
\lvec(87 19)
\lvec(85 19)
\ifill f:0
\move(97 18)
\lvec(98 18)
\lvec(98 19)
\lvec(97 19)
\ifill f:0
\move(99 18)
\lvec(101 18)
\lvec(101 19)
\lvec(99 19)
\ifill f:0
\move(102 18)
\lvec(106 18)
\lvec(106 19)
\lvec(102 19)
\ifill f:0
\move(107 18)
\lvec(114 18)
\lvec(114 19)
\lvec(107 19)
\ifill f:0
\move(115 18)
\lvec(122 18)
\lvec(122 19)
\lvec(115 19)
\ifill f:0
\move(123 18)
\lvec(127 18)
\lvec(127 19)
\lvec(123 19)
\ifill f:0
\move(128 18)
\lvec(132 18)
\lvec(132 19)
\lvec(128 19)
\ifill f:0
\move(133 18)
\lvec(137 18)
\lvec(137 19)
\lvec(133 19)
\ifill f:0
\move(139 18)
\lvec(142 18)
\lvec(142 19)
\lvec(139 19)
\ifill f:0
\move(143 18)
\lvec(146 18)
\lvec(146 19)
\lvec(143 19)
\ifill f:0
\move(147 18)
\lvec(150 18)
\lvec(150 19)
\lvec(147 19)
\ifill f:0
\move(151 18)
\lvec(154 18)
\lvec(154 19)
\lvec(151 19)
\ifill f:0
\move(155 18)
\lvec(158 18)
\lvec(158 19)
\lvec(155 19)
\ifill f:0
\move(159 18)
\lvec(170 18)
\lvec(170 19)
\lvec(159 19)
\ifill f:0
\move(171 18)
\lvec(173 18)
\lvec(173 19)
\lvec(171 19)
\ifill f:0
\move(174 18)
\lvec(176 18)
\lvec(176 19)
\lvec(174 19)
\ifill f:0
\move(177 18)
\lvec(179 18)
\lvec(179 19)
\lvec(177 19)
\ifill f:0
\move(180 18)
\lvec(182 18)
\lvec(182 19)
\lvec(180 19)
\ifill f:0
\move(183 18)
\lvec(192 18)
\lvec(192 19)
\lvec(183 19)
\ifill f:0
\move(193 18)
\lvec(203 18)
\lvec(203 19)
\lvec(193 19)
\ifill f:0
\move(204 18)
\lvec(206 18)
\lvec(206 19)
\lvec(204 19)
\ifill f:0
\move(207 18)
\lvec(217 18)
\lvec(217 19)
\lvec(207 19)
\ifill f:0
\move(218 18)
\lvec(228 18)
\lvec(228 19)
\lvec(218 19)
\ifill f:0
\move(229 18)
\lvec(230 18)
\lvec(230 19)
\lvec(229 19)
\ifill f:0
\move(231 18)
\lvec(233 18)
\lvec(233 19)
\lvec(231 19)
\ifill f:0
\move(234 18)
\lvec(235 18)
\lvec(235 19)
\lvec(234 19)
\ifill f:0
\move(236 18)
\lvec(240 18)
\lvec(240 19)
\lvec(236 19)
\ifill f:0
\move(241 18)
\lvec(245 18)
\lvec(245 19)
\lvec(241 19)
\ifill f:0
\move(246 18)
\lvec(250 18)
\lvec(250 19)
\lvec(246 19)
\ifill f:0
\move(251 18)
\lvec(257 18)
\lvec(257 19)
\lvec(251 19)
\ifill f:0
\move(258 18)
\lvec(264 18)
\lvec(264 19)
\lvec(258 19)
\ifill f:0
\move(265 18)
\lvec(271 18)
\lvec(271 19)
\lvec(265 19)
\ifill f:0
\move(272 18)
\lvec(280 18)
\lvec(280 19)
\lvec(272 19)
\ifill f:0
\move(281 18)
\lvec(291 18)
\lvec(291 19)
\lvec(281 19)
\ifill f:0
\move(292 18)
\lvec(298 18)
\lvec(298 19)
\lvec(292 19)
\ifill f:0
\move(299 18)
\lvec(302 18)
\lvec(302 19)
\lvec(299 19)
\ifill f:0
\move(303 18)
\lvec(304 18)
\lvec(304 19)
\lvec(303 19)
\ifill f:0
\move(305 18)
\lvec(319 18)
\lvec(319 19)
\lvec(305 19)
\ifill f:0
\move(320 18)
\lvec(323 18)
\lvec(323 19)
\lvec(320 19)
\ifill f:0
\move(324 18)
\lvec(325 18)
\lvec(325 19)
\lvec(324 19)
\ifill f:0
\move(326 18)
\lvec(327 18)
\lvec(327 19)
\lvec(326 19)
\ifill f:0
\move(328 18)
\lvec(329 18)
\lvec(329 19)
\lvec(328 19)
\ifill f:0
\move(330 18)
\lvec(331 18)
\lvec(331 19)
\lvec(330 19)
\ifill f:0
\move(332 18)
\lvec(333 18)
\lvec(333 19)
\lvec(332 19)
\ifill f:0
\move(334 18)
\lvec(335 18)
\lvec(335 19)
\lvec(334 19)
\ifill f:0
\move(336 18)
\lvec(337 18)
\lvec(337 19)
\lvec(336 19)
\ifill f:0
\move(338 18)
\lvec(339 18)
\lvec(339 19)
\lvec(338 19)
\ifill f:0
\move(340 18)
\lvec(345 18)
\lvec(345 19)
\lvec(340 19)
\ifill f:0
\move(346 18)
\lvec(347 18)
\lvec(347 19)
\lvec(346 19)
\ifill f:0
\move(348 18)
\lvec(349 18)
\lvec(349 19)
\lvec(348 19)
\ifill f:0
\move(350 18)
\lvec(351 18)
\lvec(351 19)
\lvec(350 19)
\ifill f:0
\move(352 18)
\lvec(353 18)
\lvec(353 19)
\lvec(352 19)
\ifill f:0
\move(354 18)
\lvec(355 18)
\lvec(355 19)
\lvec(354 19)
\ifill f:0
\move(356 18)
\lvec(357 18)
\lvec(357 19)
\lvec(356 19)
\ifill f:0
\move(358 18)
\lvec(359 18)
\lvec(359 19)
\lvec(358 19)
\ifill f:0
\move(360 18)
\lvec(363 18)
\lvec(363 19)
\lvec(360 19)
\ifill f:0
\move(364 18)
\lvec(365 18)
\lvec(365 19)
\lvec(364 19)
\ifill f:0
\move(366 18)
\lvec(382 18)
\lvec(382 19)
\lvec(366 19)
\ifill f:0
\move(383 18)
\lvec(384 18)
\lvec(384 19)
\lvec(383 19)
\ifill f:0
\move(385 18)
\lvec(386 18)
\lvec(386 19)
\lvec(385 19)
\ifill f:0
\move(387 18)
\lvec(397 18)
\lvec(397 19)
\lvec(387 19)
\ifill f:0
\move(398 18)
\lvec(399 18)
\lvec(399 19)
\lvec(398 19)
\ifill f:0
\move(400 18)
\lvec(401 18)
\lvec(401 19)
\lvec(400 19)
\ifill f:0
\move(402 18)
\lvec(410 18)
\lvec(410 19)
\lvec(402 19)
\ifill f:0
\move(411 18)
\lvec(412 18)
\lvec(412 19)
\lvec(411 19)
\ifill f:0
\move(413 18)
\lvec(421 18)
\lvec(421 19)
\lvec(413 19)
\ifill f:0
\move(422 18)
\lvec(423 18)
\lvec(423 19)
\lvec(422 19)
\ifill f:0
\move(424 18)
\lvec(430 18)
\lvec(430 19)
\lvec(424 19)
\ifill f:0
\move(431 18)
\lvec(432 18)
\lvec(432 19)
\lvec(431 19)
\ifill f:0
\move(433 18)
\lvec(439 18)
\lvec(439 19)
\lvec(433 19)
\ifill f:0
\move(440 18)
\lvec(443 18)
\lvec(443 19)
\lvec(440 19)
\ifill f:0
\move(444 18)
\lvec(448 18)
\lvec(448 19)
\lvec(444 19)
\ifill f:0
\move(449 18)
\lvec(450 18)
\lvec(450 19)
\lvec(449 19)
\ifill f:0
\move(11 19)
\lvec(13 19)
\lvec(13 20)
\lvec(11 20)
\ifill f:0
\move(14 19)
\lvec(17 19)
\lvec(17 20)
\lvec(14 20)
\ifill f:0
\move(18 19)
\lvec(19 19)
\lvec(19 20)
\lvec(18 20)
\ifill f:0
\move(20 19)
\lvec(22 19)
\lvec(22 20)
\lvec(20 20)
\ifill f:0
\move(24 19)
\lvec(26 19)
\lvec(26 20)
\lvec(24 20)
\ifill f:0
\move(36 19)
\lvec(37 19)
\lvec(37 20)
\lvec(36 20)
\ifill f:0
\move(38 19)
\lvec(42 19)
\lvec(42 20)
\lvec(38 20)
\ifill f:0
\move(43 19)
\lvec(46 19)
\lvec(46 20)
\lvec(43 20)
\ifill f:0
\move(47 19)
\lvec(50 19)
\lvec(50 20)
\lvec(47 20)
\ifill f:0
\move(52 19)
\lvec(53 19)
\lvec(53 20)
\lvec(52 20)
\ifill f:0
\move(55 19)
\lvec(58 19)
\lvec(58 20)
\lvec(55 20)
\ifill f:0
\move(59 19)
\lvec(62 19)
\lvec(62 20)
\lvec(59 20)
\ifill f:0
\move(63 19)
\lvec(65 19)
\lvec(65 20)
\lvec(63 20)
\ifill f:0
\move(66 19)
\lvec(71 19)
\lvec(71 20)
\lvec(66 20)
\ifill f:0
\move(73 19)
\lvec(79 19)
\lvec(79 20)
\lvec(73 20)
\ifill f:0
\move(81 19)
\lvec(93 19)
\lvec(93 20)
\lvec(81 20)
\ifill f:0
\move(95 19)
\lvec(101 19)
\lvec(101 20)
\lvec(95 20)
\ifill f:0
\move(102 19)
\lvec(110 19)
\lvec(110 20)
\lvec(102 20)
\ifill f:0
\move(112 19)
\lvec(122 19)
\lvec(122 20)
\lvec(112 20)
\ifill f:0
\move(124 19)
\lvec(129 19)
\lvec(129 20)
\lvec(124 20)
\ifill f:0
\move(130 19)
\lvec(135 19)
\lvec(135 20)
\lvec(130 20)
\ifill f:0
\move(136 19)
\lvec(141 19)
\lvec(141 20)
\lvec(136 20)
\ifill f:0
\move(142 19)
\lvec(147 19)
\lvec(147 20)
\lvec(142 20)
\ifill f:0
\move(148 19)
\lvec(156 19)
\lvec(156 20)
\lvec(148 20)
\ifill f:0
\move(157 19)
\lvec(161 19)
\lvec(161 20)
\lvec(157 20)
\ifill f:0
\move(162 19)
\lvec(170 19)
\lvec(170 20)
\lvec(162 20)
\ifill f:0
\move(171 19)
\lvec(173 19)
\lvec(173 20)
\lvec(171 20)
\ifill f:0
\move(174 19)
\lvec(177 19)
\lvec(177 20)
\lvec(174 20)
\ifill f:0
\move(178 19)
\lvec(184 19)
\lvec(184 20)
\lvec(178 20)
\ifill f:0
\move(185 19)
\lvec(191 19)
\lvec(191 20)
\lvec(185 20)
\ifill f:0
\move(192 19)
\lvec(204 19)
\lvec(204 20)
\lvec(192 20)
\ifill f:0
\move(205 19)
\lvec(207 19)
\lvec(207 20)
\lvec(205 20)
\ifill f:0
\move(208 19)
\lvec(210 19)
\lvec(210 20)
\lvec(208 20)
\ifill f:0
\move(211 19)
\lvec(216 19)
\lvec(216 20)
\lvec(211 20)
\ifill f:0
\move(217 19)
\lvec(219 19)
\lvec(219 20)
\lvec(217 20)
\ifill f:0
\move(220 19)
\lvec(222 19)
\lvec(222 20)
\lvec(220 20)
\ifill f:0
\move(223 19)
\lvec(226 19)
\lvec(226 20)
\lvec(223 20)
\ifill f:0
\move(227 19)
\lvec(228 19)
\lvec(228 20)
\lvec(227 20)
\ifill f:0
\move(229 19)
\lvec(236 19)
\lvec(236 20)
\lvec(229 20)
\ifill f:0
\move(237 19)
\lvec(239 19)
\lvec(239 20)
\lvec(237 20)
\ifill f:0
\move(240 19)
\lvec(247 19)
\lvec(247 20)
\lvec(240 20)
\ifill f:0
\move(248 19)
\lvec(252 19)
\lvec(252 20)
\lvec(248 20)
\ifill f:0
\move(253 19)
\lvec(272 19)
\lvec(272 20)
\lvec(253 20)
\ifill f:0
\move(273 19)
\lvec(277 19)
\lvec(277 20)
\lvec(273 20)
\ifill f:0
\move(278 19)
\lvec(296 19)
\lvec(296 20)
\lvec(278 20)
\ifill f:0
\move(297 19)
\lvec(305 19)
\lvec(305 20)
\lvec(297 20)
\ifill f:0
\move(306 19)
\lvec(314 19)
\lvec(314 20)
\lvec(306 20)
\ifill f:0
\move(315 19)
\lvec(323 19)
\lvec(323 20)
\lvec(315 20)
\ifill f:0
\move(324 19)
\lvec(325 19)
\lvec(325 20)
\lvec(324 20)
\ifill f:0
\move(326 19)
\lvec(338 19)
\lvec(338 20)
\lvec(326 20)
\ifill f:0
\move(339 19)
\lvec(340 19)
\lvec(340 20)
\lvec(339 20)
\ifill f:0
\move(341 19)
\lvec(359 19)
\lvec(359 20)
\lvec(341 20)
\ifill f:0
\move(360 19)
\lvec(363 19)
\lvec(363 20)
\lvec(360 20)
\ifill f:0
\move(364 19)
\lvec(365 19)
\lvec(365 20)
\lvec(364 20)
\ifill f:0
\move(366 19)
\lvec(367 19)
\lvec(367 20)
\lvec(366 20)
\ifill f:0
\move(368 19)
\lvec(369 19)
\lvec(369 20)
\lvec(368 20)
\ifill f:0
\move(370 19)
\lvec(371 19)
\lvec(371 20)
\lvec(370 20)
\ifill f:0
\move(372 19)
\lvec(373 19)
\lvec(373 20)
\lvec(372 20)
\ifill f:0
\move(374 19)
\lvec(375 19)
\lvec(375 20)
\lvec(374 20)
\ifill f:0
\move(376 19)
\lvec(377 19)
\lvec(377 20)
\lvec(376 20)
\ifill f:0
\move(378 19)
\lvec(383 19)
\lvec(383 20)
\lvec(378 20)
\ifill f:0
\move(384 19)
\lvec(385 19)
\lvec(385 20)
\lvec(384 20)
\ifill f:0
\move(386 19)
\lvec(387 19)
\lvec(387 20)
\lvec(386 20)
\ifill f:0
\move(388 19)
\lvec(389 19)
\lvec(389 20)
\lvec(388 20)
\ifill f:0
\move(390 19)
\lvec(391 19)
\lvec(391 20)
\lvec(390 20)
\ifill f:0
\move(392 19)
\lvec(393 19)
\lvec(393 20)
\lvec(392 20)
\ifill f:0
\move(394 19)
\lvec(395 19)
\lvec(395 20)
\lvec(394 20)
\ifill f:0
\move(396 19)
\lvec(397 19)
\lvec(397 20)
\lvec(396 20)
\ifill f:0
\move(398 19)
\lvec(399 19)
\lvec(399 20)
\lvec(398 20)
\ifill f:0
\move(400 19)
\lvec(401 19)
\lvec(401 20)
\lvec(400 20)
\ifill f:0
\move(402 19)
\lvec(403 19)
\lvec(403 20)
\lvec(402 20)
\ifill f:0
\move(404 19)
\lvec(422 19)
\lvec(422 20)
\lvec(404 20)
\ifill f:0
\move(423 19)
\lvec(424 19)
\lvec(424 20)
\lvec(423 20)
\ifill f:0
\move(425 19)
\lvec(439 19)
\lvec(439 20)
\lvec(425 20)
\ifill f:0
\move(440 19)
\lvec(443 19)
\lvec(443 20)
\lvec(440 20)
\ifill f:0
\move(444 19)
\lvec(445 19)
\lvec(445 20)
\lvec(444 20)
\ifill f:0
\move(446 19)
\lvec(451 19)
\lvec(451 20)
\lvec(446 20)
\ifill f:0
\move(11 20)
\lvec(12 20)
\lvec(12 21)
\lvec(11 21)
\ifill f:0
\move(14 20)
\lvec(15 20)
\lvec(15 21)
\lvec(14 21)
\ifill f:0
\move(16 20)
\lvec(17 20)
\lvec(17 21)
\lvec(16 21)
\ifill f:0
\move(18 20)
\lvec(19 20)
\lvec(19 21)
\lvec(18 21)
\ifill f:0
\move(20 20)
\lvec(22 20)
\lvec(22 21)
\lvec(20 21)
\ifill f:0
\move(25 20)
\lvec(26 20)
\lvec(26 21)
\lvec(25 21)
\ifill f:0
\move(34 20)
\lvec(35 20)
\lvec(35 21)
\lvec(34 21)
\ifill f:0
\move(36 20)
\lvec(37 20)
\lvec(37 21)
\lvec(36 21)
\ifill f:0
\move(38 20)
\lvec(43 20)
\lvec(43 21)
\lvec(38 21)
\ifill f:0
\move(44 20)
\lvec(45 20)
\lvec(45 21)
\lvec(44 21)
\ifill f:0
\move(46 20)
\lvec(50 20)
\lvec(50 21)
\lvec(46 21)
\ifill f:0
\move(52 20)
\lvec(53 20)
\lvec(53 21)
\lvec(52 21)
\ifill f:0
\move(54 20)
\lvec(55 20)
\lvec(55 21)
\lvec(54 21)
\ifill f:0
\move(57 20)
\lvec(59 20)
\lvec(59 21)
\lvec(57 21)
\ifill f:0
\move(60 20)
\lvec(62 20)
\lvec(62 21)
\lvec(60 21)
\ifill f:0
\move(63 20)
\lvec(65 20)
\lvec(65 21)
\lvec(63 21)
\ifill f:0
\move(67 20)
\lvec(70 20)
\lvec(70 21)
\lvec(67 21)
\ifill f:0
\move(71 20)
\lvec(74 20)
\lvec(74 21)
\lvec(71 21)
\ifill f:0
\move(76 20)
\lvec(80 20)
\lvec(80 21)
\lvec(76 21)
\ifill f:0
\move(81 20)
\lvec(86 20)
\lvec(86 21)
\lvec(81 21)
\ifill f:0
\move(87 20)
\lvec(88 20)
\lvec(88 21)
\lvec(87 21)
\ifill f:0
\move(91 20)
\lvec(101 20)
\lvec(101 21)
\lvec(91 21)
\ifill f:0
\move(102 20)
\lvec(122 20)
\lvec(122 21)
\lvec(102 21)
\ifill f:0
\move(126 20)
\lvec(132 20)
\lvec(132 21)
\lvec(126 21)
\ifill f:0
\move(134 20)
\lvec(140 20)
\lvec(140 21)
\lvec(134 21)
\ifill f:0
\move(141 20)
\lvec(147 20)
\lvec(147 21)
\lvec(141 21)
\ifill f:0
\move(148 20)
\lvec(153 20)
\lvec(153 21)
\lvec(148 21)
\ifill f:0
\move(154 20)
\lvec(159 20)
\lvec(159 21)
\lvec(154 21)
\ifill f:0
\move(160 20)
\lvec(164 20)
\lvec(164 21)
\lvec(160 21)
\ifill f:0
\move(165 20)
\lvec(170 20)
\lvec(170 21)
\lvec(165 21)
\ifill f:0
\move(171 20)
\lvec(182 20)
\lvec(182 21)
\lvec(171 21)
\ifill f:0
\move(183 20)
\lvec(190 20)
\lvec(190 21)
\lvec(183 21)
\ifill f:0
\move(192 20)
\lvec(194 20)
\lvec(194 21)
\lvec(192 21)
\ifill f:0
\move(195 20)
\lvec(205 20)
\lvec(205 21)
\lvec(195 21)
\ifill f:0
\move(206 20)
\lvec(226 20)
\lvec(226 21)
\lvec(206 21)
\ifill f:0
\move(227 20)
\lvec(228 20)
\lvec(228 21)
\lvec(227 21)
\ifill f:0
\move(229 20)
\lvec(231 20)
\lvec(231 21)
\lvec(229 21)
\ifill f:0
\move(232 20)
\lvec(234 20)
\lvec(234 21)
\lvec(232 21)
\ifill f:0
\move(235 20)
\lvec(240 20)
\lvec(240 21)
\lvec(235 21)
\ifill f:0
\move(241 20)
\lvec(243 20)
\lvec(243 21)
\lvec(241 21)
\ifill f:0
\move(244 20)
\lvec(246 20)
\lvec(246 21)
\lvec(244 21)
\ifill f:0
\move(247 20)
\lvec(260 20)
\lvec(260 21)
\lvec(247 21)
\ifill f:0
\move(261 20)
\lvec(271 20)
\lvec(271 21)
\lvec(261 21)
\ifill f:0
\move(272 20)
\lvec(276 20)
\lvec(276 21)
\lvec(272 21)
\ifill f:0
\move(277 20)
\lvec(294 20)
\lvec(294 21)
\lvec(277 21)
\ifill f:0
\move(295 20)
\lvec(299 20)
\lvec(299 21)
\lvec(295 21)
\ifill f:0
\move(300 20)
\lvec(306 20)
\lvec(306 21)
\lvec(300 21)
\ifill f:0
\move(307 20)
\lvec(311 20)
\lvec(311 21)
\lvec(307 21)
\ifill f:0
\move(312 20)
\lvec(318 20)
\lvec(318 21)
\lvec(312 21)
\ifill f:0
\move(319 20)
\lvec(323 20)
\lvec(323 21)
\lvec(319 21)
\ifill f:0
\move(324 20)
\lvec(325 20)
\lvec(325 21)
\lvec(324 21)
\ifill f:0
\move(326 20)
\lvec(332 20)
\lvec(332 21)
\lvec(326 21)
\ifill f:0
\move(333 20)
\lvec(341 20)
\lvec(341 21)
\lvec(333 21)
\ifill f:0
\move(342 20)
\lvec(350 20)
\lvec(350 21)
\lvec(342 21)
\ifill f:0
\move(351 20)
\lvec(363 20)
\lvec(363 21)
\lvec(351 21)
\ifill f:0
\move(364 20)
\lvec(376 20)
\lvec(376 21)
\lvec(364 21)
\ifill f:0
\move(377 20)
\lvec(378 20)
\lvec(378 21)
\lvec(377 21)
\ifill f:0
\move(379 20)
\lvec(399 20)
\lvec(399 21)
\lvec(379 21)
\ifill f:0
\move(400 20)
\lvec(401 20)
\lvec(401 21)
\lvec(400 21)
\ifill f:0
\move(402 20)
\lvec(403 20)
\lvec(403 21)
\lvec(402 21)
\ifill f:0
\move(404 20)
\lvec(405 20)
\lvec(405 21)
\lvec(404 21)
\ifill f:0
\move(406 20)
\lvec(407 20)
\lvec(407 21)
\lvec(406 21)
\ifill f:0
\move(408 20)
\lvec(409 20)
\lvec(409 21)
\lvec(408 21)
\ifill f:0
\move(410 20)
\lvec(411 20)
\lvec(411 21)
\lvec(410 21)
\ifill f:0
\move(412 20)
\lvec(413 20)
\lvec(413 21)
\lvec(412 21)
\ifill f:0
\move(414 20)
\lvec(415 20)
\lvec(415 21)
\lvec(414 21)
\ifill f:0
\move(416 20)
\lvec(417 20)
\lvec(417 21)
\lvec(416 21)
\ifill f:0
\move(418 20)
\lvec(423 20)
\lvec(423 21)
\lvec(418 21)
\ifill f:0
\move(424 20)
\lvec(425 20)
\lvec(425 21)
\lvec(424 21)
\ifill f:0
\move(426 20)
\lvec(427 20)
\lvec(427 21)
\lvec(426 21)
\ifill f:0
\move(428 20)
\lvec(429 20)
\lvec(429 21)
\lvec(428 21)
\ifill f:0
\move(430 20)
\lvec(431 20)
\lvec(431 21)
\lvec(430 21)
\ifill f:0
\move(432 20)
\lvec(433 20)
\lvec(433 21)
\lvec(432 21)
\ifill f:0
\move(434 20)
\lvec(435 20)
\lvec(435 21)
\lvec(434 21)
\ifill f:0
\move(436 20)
\lvec(437 20)
\lvec(437 21)
\lvec(436 21)
\ifill f:0
\move(438 20)
\lvec(439 20)
\lvec(439 21)
\lvec(438 21)
\ifill f:0
\move(440 20)
\lvec(443 20)
\lvec(443 21)
\lvec(440 21)
\ifill f:0
\move(444 20)
\lvec(445 20)
\lvec(445 21)
\lvec(444 21)
\ifill f:0
\move(446 20)
\lvec(451 20)
\lvec(451 21)
\lvec(446 21)
\ifill f:0
\move(13 21)
\lvec(14 21)
\lvec(14 22)
\lvec(13 22)
\ifill f:0
\move(16 21)
\lvec(17 21)
\lvec(17 22)
\lvec(16 22)
\ifill f:0
\move(18 21)
\lvec(22 21)
\lvec(22 22)
\lvec(18 22)
\ifill f:0
\move(23 21)
\lvec(24 21)
\lvec(24 22)
\lvec(23 22)
\ifill f:0
\move(25 21)
\lvec(26 21)
\lvec(26 22)
\lvec(25 22)
\ifill f:0
\move(27 21)
\lvec(29 21)
\lvec(29 22)
\lvec(27 22)
\ifill f:0
\move(36 21)
\lvec(37 21)
\lvec(37 22)
\lvec(36 22)
\ifill f:0
\move(38 21)
\lvec(39 21)
\lvec(39 22)
\lvec(38 22)
\ifill f:0
\move(40 21)
\lvec(42 21)
\lvec(42 22)
\lvec(40 22)
\ifill f:0
\move(43 21)
\lvec(47 21)
\lvec(47 22)
\lvec(43 22)
\ifill f:0
\move(48 21)
\lvec(50 21)
\lvec(50 22)
\lvec(48 22)
\ifill f:0
\move(51 21)
\lvec(52 21)
\lvec(52 22)
\lvec(51 22)
\ifill f:0
\move(56 21)
\lvec(57 21)
\lvec(57 22)
\lvec(56 22)
\ifill f:0
\move(58 21)
\lvec(60 21)
\lvec(60 22)
\lvec(58 22)
\ifill f:0
\move(61 21)
\lvec(65 21)
\lvec(65 22)
\lvec(61 22)
\ifill f:0
\move(66 21)
\lvec(72 21)
\lvec(72 22)
\lvec(66 22)
\ifill f:0
\move(73 21)
\lvec(76 21)
\lvec(76 22)
\lvec(73 22)
\ifill f:0
\move(77 21)
\lvec(82 21)
\lvec(82 22)
\lvec(77 22)
\ifill f:0
\move(83 21)
\lvec(85 21)
\lvec(85 22)
\lvec(83 22)
\ifill f:0
\move(88 21)
\lvec(92 21)
\lvec(92 22)
\lvec(88 22)
\ifill f:0
\move(93 21)
\lvec(94 21)
\lvec(94 22)
\lvec(93 22)
\ifill f:0
\move(95 21)
\lvec(101 21)
\lvec(101 22)
\lvec(95 22)
\ifill f:0
\move(102 21)
\lvec(109 21)
\lvec(109 22)
\lvec(102 22)
\ifill f:0
\move(115 21)
\lvec(117 21)
\lvec(117 22)
\lvec(115 22)
\ifill f:0
\move(121 21)
\lvec(122 21)
\lvec(122 22)
\lvec(121 22)
\ifill f:0
\move(127 21)
\lvec(145 21)
\lvec(145 22)
\lvec(127 22)
\ifill f:0
\move(146 21)
\lvec(148 21)
\lvec(148 22)
\lvec(146 22)
\ifill f:0
\move(150 21)
\lvec(156 21)
\lvec(156 22)
\lvec(150 22)
\ifill f:0
\move(157 21)
\lvec(163 21)
\lvec(163 22)
\lvec(157 22)
\ifill f:0
\move(164 21)
\lvec(170 21)
\lvec(170 22)
\lvec(164 22)
\ifill f:0
\move(171 21)
\lvec(175 21)
\lvec(175 22)
\lvec(171 22)
\ifill f:0
\move(176 21)
\lvec(180 21)
\lvec(180 22)
\lvec(176 22)
\ifill f:0
\move(181 21)
\lvec(185 21)
\lvec(185 22)
\lvec(181 22)
\ifill f:0
\move(186 21)
\lvec(189 21)
\lvec(189 22)
\lvec(186 22)
\ifill f:0
\move(190 21)
\lvec(194 21)
\lvec(194 22)
\lvec(190 22)
\ifill f:0
\move(195 21)
\lvec(198 21)
\lvec(198 22)
\lvec(195 22)
\ifill f:0
\move(199 21)
\lvec(202 21)
\lvec(202 22)
\lvec(199 22)
\ifill f:0
\move(203 21)
\lvec(210 21)
\lvec(210 22)
\lvec(203 22)
\ifill f:0
\move(211 21)
\lvec(214 21)
\lvec(214 22)
\lvec(211 22)
\ifill f:0
\move(215 21)
\lvec(226 21)
\lvec(226 22)
\lvec(215 22)
\ifill f:0
\move(227 21)
\lvec(229 21)
\lvec(229 22)
\lvec(227 22)
\ifill f:0
\move(230 21)
\lvec(232 21)
\lvec(232 22)
\lvec(230 22)
\ifill f:0
\move(233 21)
\lvec(242 21)
\lvec(242 22)
\lvec(233 22)
\ifill f:0
\move(243 21)
\lvec(245 21)
\lvec(245 22)
\lvec(243 22)
\ifill f:0
\move(246 21)
\lvec(248 21)
\lvec(248 22)
\lvec(246 22)
\ifill f:0
\move(249 21)
\lvec(272 21)
\lvec(272 22)
\lvec(249 22)
\ifill f:0
\move(273 21)
\lvec(275 21)
\lvec(275 22)
\lvec(273 22)
\ifill f:0
\move(276 21)
\lvec(290 21)
\lvec(290 22)
\lvec(276 22)
\ifill f:0
\move(291 21)
\lvec(292 21)
\lvec(292 22)
\lvec(291 22)
\ifill f:0
\move(293 21)
\lvec(297 21)
\lvec(297 22)
\lvec(293 22)
\ifill f:0
\move(298 21)
\lvec(305 21)
\lvec(305 22)
\lvec(298 22)
\ifill f:0
\move(306 21)
\lvec(310 21)
\lvec(310 22)
\lvec(306 22)
\ifill f:0
\move(311 21)
\lvec(315 21)
\lvec(315 22)
\lvec(311 22)
\ifill f:0
\move(316 21)
\lvec(325 21)
\lvec(325 22)
\lvec(316 22)
\ifill f:0
\move(326 21)
\lvec(330 21)
\lvec(330 22)
\lvec(326 22)
\ifill f:0
\move(331 21)
\lvec(335 21)
\lvec(335 22)
\lvec(331 22)
\ifill f:0
\move(336 21)
\lvec(342 21)
\lvec(342 22)
\lvec(336 22)
\ifill f:0
\move(343 21)
\lvec(354 21)
\lvec(354 22)
\lvec(343 22)
\ifill f:0
\move(355 21)
\lvec(368 21)
\lvec(368 22)
\lvec(355 22)
\ifill f:0
\move(369 21)
\lvec(370 21)
\lvec(370 22)
\lvec(369 22)
\ifill f:0
\move(371 21)
\lvec(377 21)
\lvec(377 22)
\lvec(371 22)
\ifill f:0
\move(378 21)
\lvec(379 21)
\lvec(379 22)
\lvec(378 22)
\ifill f:0
\move(380 21)
\lvec(388 21)
\lvec(388 22)
\lvec(380 22)
\ifill f:0
\move(389 21)
\lvec(390 21)
\lvec(390 22)
\lvec(389 22)
\ifill f:0
\move(391 21)
\lvec(401 21)
\lvec(401 22)
\lvec(391 22)
\ifill f:0
\move(402 21)
\lvec(403 21)
\lvec(403 22)
\lvec(402 22)
\ifill f:0
\move(404 21)
\lvec(416 21)
\lvec(416 22)
\lvec(404 22)
\ifill f:0
\move(417 21)
\lvec(418 21)
\lvec(418 22)
\lvec(417 22)
\ifill f:0
\move(419 21)
\lvec(439 21)
\lvec(439 22)
\lvec(419 22)
\ifill f:0
\move(440 21)
\lvec(443 21)
\lvec(443 22)
\lvec(440 22)
\ifill f:0
\move(444 21)
\lvec(445 21)
\lvec(445 22)
\lvec(444 22)
\ifill f:0
\move(446 21)
\lvec(447 21)
\lvec(447 22)
\lvec(446 22)
\ifill f:0
\move(448 21)
\lvec(449 21)
\lvec(449 22)
\lvec(448 22)
\ifill f:0
\move(450 21)
\lvec(451 21)
\lvec(451 22)
\lvec(450 22)
\ifill f:0
\move(11 22)
\lvec(12 22)
\lvec(12 23)
\lvec(11 23)
\ifill f:0
\move(16 22)
\lvec(17 22)
\lvec(17 23)
\lvec(16 23)
\ifill f:0
\move(20 22)
\lvec(24 22)
\lvec(24 23)
\lvec(20 23)
\ifill f:0
\move(25 22)
\lvec(26 22)
\lvec(26 23)
\lvec(25 23)
\ifill f:0
\move(36 22)
\lvec(37 22)
\lvec(37 23)
\lvec(36 23)
\ifill f:0
\move(38 22)
\lvec(39 22)
\lvec(39 23)
\lvec(38 23)
\ifill f:0
\move(40 22)
\lvec(50 22)
\lvec(50 23)
\lvec(40 23)
\ifill f:0
\move(51 22)
\lvec(52 22)
\lvec(52 23)
\lvec(51 23)
\ifill f:0
\move(57 22)
\lvec(65 22)
\lvec(65 23)
\lvec(57 23)
\ifill f:0
\move(66 22)
\lvec(68 22)
\lvec(68 23)
\lvec(66 23)
\ifill f:0
\move(69 22)
\lvec(82 22)
\lvec(82 23)
\lvec(69 23)
\ifill f:0
\move(83 22)
\lvec(84 22)
\lvec(84 23)
\lvec(83 23)
\ifill f:0
\move(86 22)
\lvec(89 22)
\lvec(89 23)
\lvec(86 23)
\ifill f:0
\move(91 22)
\lvec(93 22)
\lvec(93 23)
\lvec(91 23)
\ifill f:0
\move(97 22)
\lvec(101 22)
\lvec(101 23)
\lvec(97 23)
\ifill f:0
\move(102 22)
\lvec(104 22)
\lvec(104 23)
\lvec(102 23)
\ifill f:0
\move(105 22)
\lvec(119 22)
\lvec(119 23)
\lvec(105 23)
\ifill f:0
\move(121 22)
\lvec(122 22)
\lvec(122 23)
\lvec(121 23)
\ifill f:0
\move(126 22)
\lvec(128 22)
\lvec(128 23)
\lvec(126 23)
\ifill f:0
\move(136 22)
\lvec(137 22)
\lvec(137 23)
\lvec(136 23)
\ifill f:0
\move(139 22)
\lvec(145 22)
\lvec(145 23)
\lvec(139 23)
\ifill f:0
\move(146 22)
\lvec(151 22)
\lvec(151 23)
\lvec(146 23)
\ifill f:0
\move(152 22)
\lvec(161 22)
\lvec(161 23)
\lvec(152 23)
\ifill f:0
\move(162 22)
\lvec(170 22)
\lvec(170 23)
\lvec(162 23)
\ifill f:0
\move(171 22)
\lvec(176 22)
\lvec(176 23)
\lvec(171 23)
\ifill f:0
\move(177 22)
\lvec(182 22)
\lvec(182 23)
\lvec(177 23)
\ifill f:0
\move(183 22)
\lvec(188 22)
\lvec(188 23)
\lvec(183 23)
\ifill f:0
\move(189 22)
\lvec(194 22)
\lvec(194 23)
\lvec(189 23)
\ifill f:0
\move(195 22)
\lvec(203 22)
\lvec(203 23)
\lvec(195 23)
\ifill f:0
\move(204 22)
\lvec(226 22)
\lvec(226 23)
\lvec(204 23)
\ifill f:0
\move(227 22)
\lvec(229 22)
\lvec(229 23)
\lvec(227 23)
\ifill f:0
\move(230 22)
\lvec(233 22)
\lvec(233 23)
\lvec(230 23)
\ifill f:0
\move(234 22)
\lvec(240 22)
\lvec(240 23)
\lvec(234 23)
\ifill f:0
\move(241 22)
\lvec(244 22)
\lvec(244 23)
\lvec(241 23)
\ifill f:0
\move(245 22)
\lvec(264 22)
\lvec(264 23)
\lvec(245 23)
\ifill f:0
\move(265 22)
\lvec(274 22)
\lvec(274 23)
\lvec(265 23)
\ifill f:0
\move(275 22)
\lvec(277 22)
\lvec(277 23)
\lvec(275 23)
\ifill f:0
\move(278 22)
\lvec(280 22)
\lvec(280 23)
\lvec(278 23)
\ifill f:0
\move(281 22)
\lvec(283 22)
\lvec(283 23)
\lvec(281 23)
\ifill f:0
\move(284 22)
\lvec(286 22)
\lvec(286 23)
\lvec(284 23)
\ifill f:0
\move(287 22)
\lvec(290 22)
\lvec(290 23)
\lvec(287 23)
\ifill f:0
\move(291 22)
\lvec(292 22)
\lvec(292 23)
\lvec(291 23)
\ifill f:0
\move(293 22)
\lvec(295 22)
\lvec(295 23)
\lvec(293 23)
\ifill f:0
\move(296 22)
\lvec(298 22)
\lvec(298 23)
\lvec(296 23)
\ifill f:0
\move(299 22)
\lvec(306 22)
\lvec(306 23)
\lvec(299 23)
\ifill f:0
\move(307 22)
\lvec(309 22)
\lvec(309 23)
\lvec(307 23)
\ifill f:0
\move(310 22)
\lvec(320 22)
\lvec(320 23)
\lvec(310 23)
\ifill f:0
\move(321 22)
\lvec(328 22)
\lvec(328 23)
\lvec(321 23)
\ifill f:0
\move(329 22)
\lvec(333 22)
\lvec(333 23)
\lvec(329 23)
\ifill f:0
\move(334 22)
\lvec(341 22)
\lvec(341 23)
\lvec(334 23)
\ifill f:0
\move(342 22)
\lvec(346 22)
\lvec(346 23)
\lvec(342 23)
\ifill f:0
\move(347 22)
\lvec(351 22)
\lvec(351 23)
\lvec(347 23)
\ifill f:0
\move(352 22)
\lvec(356 22)
\lvec(356 23)
\lvec(352 23)
\ifill f:0
\move(357 22)
\lvec(366 22)
\lvec(366 23)
\lvec(357 23)
\ifill f:0
\move(367 22)
\lvec(373 22)
\lvec(373 23)
\lvec(367 23)
\ifill f:0
\move(374 22)
\lvec(380 22)
\lvec(380 23)
\lvec(374 23)
\ifill f:0
\move(381 22)
\lvec(385 22)
\lvec(385 23)
\lvec(381 23)
\ifill f:0
\move(386 22)
\lvec(392 22)
\lvec(392 23)
\lvec(386 23)
\ifill f:0
\move(393 22)
\lvec(401 22)
\lvec(401 23)
\lvec(393 23)
\ifill f:0
\move(402 22)
\lvec(410 22)
\lvec(410 23)
\lvec(402 23)
\ifill f:0
\move(411 22)
\lvec(419 22)
\lvec(419 23)
\lvec(411 23)
\ifill f:0
\move(420 22)
\lvec(430 22)
\lvec(430 23)
\lvec(420 23)
\ifill f:0
\move(431 22)
\lvec(443 22)
\lvec(443 23)
\lvec(431 23)
\ifill f:0
\move(444 22)
\lvec(451 22)
\lvec(451 23)
\lvec(444 23)
\ifill f:0
\move(11 23)
\lvec(12 23)
\lvec(12 24)
\lvec(11 24)
\ifill f:0
\move(14 23)
\lvec(17 23)
\lvec(17 24)
\lvec(14 24)
\ifill f:0
\move(18 23)
\lvec(19 23)
\lvec(19 24)
\lvec(18 24)
\ifill f:0
\move(20 23)
\lvec(22 23)
\lvec(22 24)
\lvec(20 24)
\ifill f:0
\move(23 23)
\lvec(26 23)
\lvec(26 24)
\lvec(23 24)
\ifill f:0
\move(29 23)
\lvec(30 23)
\lvec(30 24)
\lvec(29 24)
\ifill f:0
\move(36 23)
\lvec(37 23)
\lvec(37 24)
\lvec(36 24)
\ifill f:0
\move(38 23)
\lvec(41 23)
\lvec(41 24)
\lvec(38 24)
\ifill f:0
\move(43 23)
\lvec(46 23)
\lvec(46 24)
\lvec(43 24)
\ifill f:0
\move(47 23)
\lvec(50 23)
\lvec(50 24)
\lvec(47 24)
\ifill f:0
\move(51 23)
\lvec(52 23)
\lvec(52 24)
\lvec(51 24)
\ifill f:0
\move(54 23)
\lvec(55 23)
\lvec(55 24)
\lvec(54 24)
\ifill f:0
\move(56 23)
\lvec(57 23)
\lvec(57 24)
\lvec(56 24)
\ifill f:0
\move(58 23)
\lvec(59 23)
\lvec(59 24)
\lvec(58 24)
\ifill f:0
\move(60 23)
\lvec(63 23)
\lvec(63 24)
\lvec(60 24)
\ifill f:0
\move(64 23)
\lvec(65 23)
\lvec(65 24)
\lvec(64 24)
\ifill f:0
\move(66 23)
\lvec(72 23)
\lvec(72 24)
\lvec(66 24)
\ifill f:0
\move(73 23)
\lvec(75 23)
\lvec(75 24)
\lvec(73 24)
\ifill f:0
\move(76 23)
\lvec(82 23)
\lvec(82 24)
\lvec(76 24)
\ifill f:0
\move(83 23)
\lvec(84 23)
\lvec(84 24)
\lvec(83 24)
\ifill f:0
\move(85 23)
\lvec(88 23)
\lvec(88 24)
\lvec(85 24)
\ifill f:0
\move(89 23)
\lvec(96 23)
\lvec(96 24)
\lvec(89 24)
\ifill f:0
\move(97 23)
\lvec(101 23)
\lvec(101 24)
\lvec(97 24)
\ifill f:0
\move(102 23)
\lvec(103 23)
\lvec(103 24)
\lvec(102 24)
\ifill f:0
\move(104 23)
\lvec(110 23)
\lvec(110 24)
\lvec(104 24)
\ifill f:0
\move(111 23)
\lvec(120 23)
\lvec(120 24)
\lvec(111 24)
\ifill f:0
\move(121 23)
\lvec(145 23)
\lvec(145 24)
\lvec(121 24)
\ifill f:0
\move(146 23)
\lvec(156 23)
\lvec(156 24)
\lvec(146 24)
\ifill f:0
\move(157 23)
\lvec(170 23)
\lvec(170 24)
\lvec(157 24)
\ifill f:0
\move(172 23)
\lvec(178 23)
\lvec(178 24)
\lvec(172 24)
\ifill f:0
\move(179 23)
\lvec(186 23)
\lvec(186 24)
\lvec(179 24)
\ifill f:0
\move(188 23)
\lvec(194 23)
\lvec(194 24)
\lvec(188 24)
\ifill f:0
\move(195 23)
\lvec(199 23)
\lvec(199 24)
\lvec(195 24)
\ifill f:0
\move(200 23)
\lvec(210 23)
\lvec(210 24)
\lvec(200 24)
\ifill f:0
\move(211 23)
\lvec(220 23)
\lvec(220 24)
\lvec(211 24)
\ifill f:0
\move(221 23)
\lvec(226 23)
\lvec(226 24)
\lvec(221 24)
\ifill f:0
\move(227 23)
\lvec(234 23)
\lvec(234 24)
\lvec(227 24)
\ifill f:0
\move(235 23)
\lvec(238 23)
\lvec(238 24)
\lvec(235 24)
\ifill f:0
\move(239 23)
\lvec(242 23)
\lvec(242 24)
\lvec(239 24)
\ifill f:0
\move(243 23)
\lvec(250 23)
\lvec(250 24)
\lvec(243 24)
\ifill f:0
\move(251 23)
\lvec(254 23)
\lvec(254 24)
\lvec(251 24)
\ifill f:0
\move(255 23)
\lvec(265 23)
\lvec(265 24)
\lvec(255 24)
\ifill f:0
\move(266 23)
\lvec(272 23)
\lvec(272 24)
\lvec(266 24)
\ifill f:0
\move(273 23)
\lvec(279 23)
\lvec(279 24)
\lvec(273 24)
\ifill f:0
\move(280 23)
\lvec(290 23)
\lvec(290 24)
\lvec(280 24)
\ifill f:0
\move(291 23)
\lvec(292 23)
\lvec(292 24)
\lvec(291 24)
\ifill f:0
\move(293 23)
\lvec(305 23)
\lvec(305 24)
\lvec(293 24)
\ifill f:0
\move(306 23)
\lvec(311 23)
\lvec(311 24)
\lvec(306 24)
\ifill f:0
\move(312 23)
\lvec(314 23)
\lvec(314 24)
\lvec(312 24)
\ifill f:0
\move(315 23)
\lvec(320 23)
\lvec(320 24)
\lvec(315 24)
\ifill f:0
\move(321 23)
\lvec(331 23)
\lvec(331 24)
\lvec(321 24)
\ifill f:0
\move(332 23)
\lvec(342 23)
\lvec(342 24)
\lvec(332 24)
\ifill f:0
\move(343 23)
\lvec(353 23)
\lvec(353 24)
\lvec(343 24)
\ifill f:0
\move(354 23)
\lvec(364 23)
\lvec(364 24)
\lvec(354 24)
\ifill f:0
\move(365 23)
\lvec(366 23)
\lvec(366 24)
\lvec(365 24)
\ifill f:0
\move(367 23)
\lvec(379 23)
\lvec(379 24)
\lvec(367 24)
\ifill f:0
\move(380 23)
\lvec(384 23)
\lvec(384 24)
\lvec(380 24)
\ifill f:0
\move(385 23)
\lvec(389 23)
\lvec(389 24)
\lvec(385 24)
\ifill f:0
\move(390 23)
\lvec(394 23)
\lvec(394 24)
\lvec(390 24)
\ifill f:0
\move(395 23)
\lvec(401 23)
\lvec(401 24)
\lvec(395 24)
\ifill f:0
\move(402 23)
\lvec(406 23)
\lvec(406 24)
\lvec(402 24)
\ifill f:0
\move(407 23)
\lvec(413 23)
\lvec(413 24)
\lvec(407 24)
\ifill f:0
\move(414 23)
\lvec(420 23)
\lvec(420 24)
\lvec(414 24)
\ifill f:0
\move(421 23)
\lvec(427 23)
\lvec(427 24)
\lvec(421 24)
\ifill f:0
\move(428 23)
\lvec(434 23)
\lvec(434 24)
\lvec(428 24)
\ifill f:0
\move(435 23)
\lvec(450 23)
\lvec(450 24)
\lvec(435 24)
\ifill f:0
\move(14 24)
\lvec(17 24)
\lvec(17 25)
\lvec(14 25)
\ifill f:0
\move(18 24)
\lvec(22 24)
\lvec(22 25)
\lvec(18 25)
\ifill f:0
\move(24 24)
\lvec(26 24)
\lvec(26 25)
\lvec(24 25)
\ifill f:0
\move(28 24)
\lvec(29 24)
\lvec(29 25)
\lvec(28 25)
\ifill f:0
\move(36 24)
\lvec(37 24)
\lvec(37 25)
\lvec(36 25)
\ifill f:0
\move(38 24)
\lvec(44 24)
\lvec(44 25)
\lvec(38 25)
\ifill f:0
\move(45 24)
\lvec(50 24)
\lvec(50 25)
\lvec(45 25)
\ifill f:0
\move(51 24)
\lvec(53 24)
\lvec(53 25)
\lvec(51 25)
\ifill f:0
\move(55 24)
\lvec(56 24)
\lvec(56 25)
\lvec(55 25)
\ifill f:0
\move(57 24)
\lvec(59 24)
\lvec(59 25)
\lvec(57 25)
\ifill f:0
\move(60 24)
\lvec(63 24)
\lvec(63 25)
\lvec(60 25)
\ifill f:0
\move(64 24)
\lvec(65 24)
\lvec(65 25)
\lvec(64 25)
\ifill f:0
\move(66 24)
\lvec(69 24)
\lvec(69 25)
\lvec(66 25)
\ifill f:0
\move(70 24)
\lvec(71 24)
\lvec(71 25)
\lvec(70 25)
\ifill f:0
\move(72 24)
\lvec(82 24)
\lvec(82 25)
\lvec(72 25)
\ifill f:0
\move(84 24)
\lvec(86 24)
\lvec(86 25)
\lvec(84 25)
\ifill f:0
\move(88 24)
\lvec(90 24)
\lvec(90 25)
\lvec(88 25)
\ifill f:0
\move(91 24)
\lvec(94 24)
\lvec(94 25)
\lvec(91 25)
\ifill f:0
\move(95 24)
\lvec(98 24)
\lvec(98 25)
\lvec(95 25)
\ifill f:0
\move(99 24)
\lvec(101 24)
\lvec(101 25)
\lvec(99 25)
\ifill f:0
\move(102 24)
\lvec(107 24)
\lvec(107 25)
\lvec(102 25)
\ifill f:0
\move(108 24)
\lvec(120 24)
\lvec(120 25)
\lvec(108 25)
\ifill f:0
\move(121 24)
\lvec(131 24)
\lvec(131 25)
\lvec(121 25)
\ifill f:0
\move(134 24)
\lvec(145 24)
\lvec(145 25)
\lvec(134 25)
\ifill f:0
\move(146 24)
\lvec(168 24)
\lvec(168 25)
\lvec(146 25)
\ifill f:0
\move(169 24)
\lvec(170 24)
\lvec(170 25)
\lvec(169 25)
\ifill f:0
\move(174 24)
\lvec(182 24)
\lvec(182 25)
\lvec(174 25)
\ifill f:0
\move(183 24)
\lvec(192 24)
\lvec(192 25)
\lvec(183 25)
\ifill f:0
\move(193 24)
\lvec(194 24)
\lvec(194 25)
\lvec(193 25)
\ifill f:0
\move(195 24)
\lvec(197 24)
\lvec(197 25)
\lvec(195 25)
\ifill f:0
\move(198 24)
\lvec(200 24)
\lvec(200 25)
\lvec(198 25)
\ifill f:0
\move(201 24)
\lvec(207 24)
\lvec(207 25)
\lvec(201 25)
\ifill f:0
\move(208 24)
\lvec(213 24)
\lvec(213 25)
\lvec(208 25)
\ifill f:0
\move(214 24)
\lvec(219 24)
\lvec(219 25)
\lvec(214 25)
\ifill f:0
\move(220 24)
\lvec(226 24)
\lvec(226 25)
\lvec(220 25)
\ifill f:0
\move(227 24)
\lvec(230 24)
\lvec(230 25)
\lvec(227 25)
\ifill f:0
\move(231 24)
\lvec(240 24)
\lvec(240 25)
\lvec(231 25)
\ifill f:0
\move(241 24)
\lvec(250 24)
\lvec(250 25)
\lvec(241 25)
\ifill f:0
\move(251 24)
\lvec(254 24)
\lvec(254 25)
\lvec(251 25)
\ifill f:0
\move(255 24)
\lvec(258 24)
\lvec(258 25)
\lvec(255 25)
\ifill f:0
\move(259 24)
\lvec(262 24)
\lvec(262 25)
\lvec(259 25)
\ifill f:0
\move(263 24)
\lvec(270 24)
\lvec(270 25)
\lvec(263 25)
\ifill f:0
\move(271 24)
\lvec(274 24)
\lvec(274 25)
\lvec(271 25)
\ifill f:0
\move(275 24)
\lvec(290 24)
\lvec(290 25)
\lvec(275 25)
\ifill f:0
\move(291 24)
\lvec(293 24)
\lvec(293 25)
\lvec(291 25)
\ifill f:0
\move(294 24)
\lvec(296 24)
\lvec(296 25)
\lvec(294 25)
\ifill f:0
\move(297 24)
\lvec(306 24)
\lvec(306 25)
\lvec(297 25)
\ifill f:0
\move(307 24)
\lvec(313 24)
\lvec(313 25)
\lvec(307 25)
\ifill f:0
\move(314 24)
\lvec(316 24)
\lvec(316 25)
\lvec(314 25)
\ifill f:0
\move(317 24)
\lvec(319 24)
\lvec(319 25)
\lvec(317 25)
\ifill f:0
\move(320 24)
\lvec(358 24)
\lvec(358 25)
\lvec(320 25)
\ifill f:0
\move(359 24)
\lvec(362 24)
\lvec(362 25)
\lvec(359 25)
\ifill f:0
\move(363 24)
\lvec(364 24)
\lvec(364 25)
\lvec(363 25)
\ifill f:0
\move(365 24)
\lvec(372 24)
\lvec(372 25)
\lvec(365 25)
\ifill f:0
\move(373 24)
\lvec(375 24)
\lvec(375 25)
\lvec(373 25)
\ifill f:0
\move(376 24)
\lvec(380 24)
\lvec(380 25)
\lvec(376 25)
\ifill f:0
\move(381 24)
\lvec(383 24)
\lvec(383 25)
\lvec(381 25)
\ifill f:0
\move(384 24)
\lvec(388 24)
\lvec(388 25)
\lvec(384 25)
\ifill f:0
\move(389 24)
\lvec(391 24)
\lvec(391 25)
\lvec(389 25)
\ifill f:0
\move(392 24)
\lvec(396 24)
\lvec(396 25)
\lvec(392 25)
\ifill f:0
\move(397 24)
\lvec(409 24)
\lvec(409 25)
\lvec(397 25)
\ifill f:0
\move(410 24)
\lvec(414 24)
\lvec(414 25)
\lvec(410 25)
\ifill f:0
\move(415 24)
\lvec(419 24)
\lvec(419 25)
\lvec(415 25)
\ifill f:0
\move(420 24)
\lvec(436 24)
\lvec(436 25)
\lvec(420 25)
\ifill f:0
\move(437 24)
\lvec(448 24)
\lvec(448 25)
\lvec(437 25)
\ifill f:0
\move(449 24)
\lvec(451 24)
\lvec(451 25)
\lvec(449 25)
\ifill f:0
\move(11 25)
\lvec(12 25)
\lvec(12 26)
\lvec(11 26)
\ifill f:0
\move(15 25)
\lvec(17 25)
\lvec(17 26)
\lvec(15 26)
\ifill f:0
\move(18 25)
\lvec(19 25)
\lvec(19 26)
\lvec(18 26)
\ifill f:0
\move(20 25)
\lvec(23 25)
\lvec(23 26)
\lvec(20 26)
\ifill f:0
\move(24 25)
\lvec(26 25)
\lvec(26 26)
\lvec(24 26)
\ifill f:0
\move(36 25)
\lvec(37 25)
\lvec(37 26)
\lvec(36 26)
\ifill f:0
\move(38 25)
\lvec(39 25)
\lvec(39 26)
\lvec(38 26)
\ifill f:0
\move(40 25)
\lvec(50 25)
\lvec(50 26)
\lvec(40 26)
\ifill f:0
\move(52 25)
\lvec(54 25)
\lvec(54 26)
\lvec(52 26)
\ifill f:0
\move(56 25)
\lvec(57 25)
\lvec(57 26)
\lvec(56 26)
\ifill f:0
\move(58 25)
\lvec(60 25)
\lvec(60 26)
\lvec(58 26)
\ifill f:0
\move(61 25)
\lvec(63 25)
\lvec(63 26)
\lvec(61 26)
\ifill f:0
\move(64 25)
\lvec(65 25)
\lvec(65 26)
\lvec(64 26)
\ifill f:0
\move(66 25)
\lvec(74 25)
\lvec(74 26)
\lvec(66 26)
\ifill f:0
\move(75 25)
\lvec(82 25)
\lvec(82 26)
\lvec(75 26)
\ifill f:0
\move(84 25)
\lvec(85 25)
\lvec(85 26)
\lvec(84 26)
\ifill f:0
\move(87 25)
\lvec(93 25)
\lvec(93 26)
\lvec(87 26)
\ifill f:0
\move(96 25)
\lvec(98 25)
\lvec(98 26)
\lvec(96 26)
\ifill f:0
\move(99 25)
\lvec(101 25)
\lvec(101 26)
\lvec(99 26)
\ifill f:0
\move(103 25)
\lvec(106 25)
\lvec(106 26)
\lvec(103 26)
\ifill f:0
\move(107 25)
\lvec(122 25)
\lvec(122 26)
\lvec(107 26)
\ifill f:0
\move(123 25)
\lvec(126 25)
\lvec(126 26)
\lvec(123 26)
\ifill f:0
\move(129 25)
\lvec(137 25)
\lvec(137 26)
\lvec(129 26)
\ifill f:0
\move(138 25)
\lvec(145 25)
\lvec(145 26)
\lvec(138 26)
\ifill f:0
\move(146 25)
\lvec(154 25)
\lvec(154 26)
\lvec(146 26)
\ifill f:0
\move(162 25)
\lvec(164 25)
\lvec(164 26)
\lvec(162 26)
\ifill f:0
\move(169 25)
\lvec(170 25)
\lvec(170 26)
\lvec(169 26)
\ifill f:0
\move(174 25)
\lvec(176 25)
\lvec(176 26)
\lvec(174 26)
\ifill f:0
\move(177 25)
\lvec(190 25)
\lvec(190 26)
\lvec(177 26)
\ifill f:0
\move(193 25)
\lvec(197 25)
\lvec(197 26)
\lvec(193 26)
\ifill f:0
\move(198 25)
\lvec(201 25)
\lvec(201 26)
\lvec(198 26)
\ifill f:0
\move(202 25)
\lvec(210 25)
\lvec(210 26)
\lvec(202 26)
\ifill f:0
\move(212 25)
\lvec(218 25)
\lvec(218 26)
\lvec(212 26)
\ifill f:0
\move(219 25)
\lvec(226 25)
\lvec(226 26)
\lvec(219 26)
\ifill f:0
\move(227 25)
\lvec(231 25)
\lvec(231 26)
\lvec(227 26)
\ifill f:0
\move(233 25)
\lvec(237 25)
\lvec(237 26)
\lvec(233 26)
\ifill f:0
\move(238 25)
\lvec(243 25)
\lvec(243 26)
\lvec(238 26)
\ifill f:0
\move(244 25)
\lvec(248 25)
\lvec(248 26)
\lvec(244 26)
\ifill f:0
\move(249 25)
\lvec(263 25)
\lvec(263 26)
\lvec(249 26)
\ifill f:0
\move(264 25)
\lvec(272 25)
\lvec(272 26)
\lvec(264 26)
\ifill f:0
\move(273 25)
\lvec(281 25)
\lvec(281 26)
\lvec(273 26)
\ifill f:0
\move(282 25)
\lvec(290 25)
\lvec(290 26)
\lvec(282 26)
\ifill f:0
\move(291 25)
\lvec(293 25)
\lvec(293 26)
\lvec(291 26)
\ifill f:0
\move(294 25)
\lvec(297 25)
\lvec(297 26)
\lvec(294 26)
\ifill f:0
\move(298 25)
\lvec(308 25)
\lvec(308 26)
\lvec(298 26)
\ifill f:0
\move(309 25)
\lvec(315 25)
\lvec(315 26)
\lvec(309 26)
\ifill f:0
\move(316 25)
\lvec(319 25)
\lvec(319 26)
\lvec(316 26)
\ifill f:0
\move(320 25)
\lvec(322 25)
\lvec(322 26)
\lvec(320 26)
\ifill f:0
\move(323 25)
\lvec(329 25)
\lvec(329 26)
\lvec(323 26)
\ifill f:0
\move(330 25)
\lvec(339 25)
\lvec(339 26)
\lvec(330 26)
\ifill f:0
\move(340 25)
\lvec(342 25)
\lvec(342 26)
\lvec(340 26)
\ifill f:0
\move(343 25)
\lvec(358 25)
\lvec(358 26)
\lvec(343 26)
\ifill f:0
\move(359 25)
\lvec(362 25)
\lvec(362 26)
\lvec(359 26)
\ifill f:0
\move(363 25)
\lvec(364 25)
\lvec(364 26)
\lvec(363 26)
\ifill f:0
\move(365 25)
\lvec(367 25)
\lvec(367 26)
\lvec(365 26)
\ifill f:0
\move(368 25)
\lvec(370 25)
\lvec(370 26)
\lvec(368 26)
\ifill f:0
\move(371 25)
\lvec(373 25)
\lvec(373 26)
\lvec(371 26)
\ifill f:0
\move(374 25)
\lvec(376 25)
\lvec(376 26)
\lvec(374 26)
\ifill f:0
\move(377 25)
\lvec(379 25)
\lvec(379 26)
\lvec(377 26)
\ifill f:0
\move(380 25)
\lvec(390 25)
\lvec(390 26)
\lvec(380 26)
\ifill f:0
\move(391 25)
\lvec(393 25)
\lvec(393 26)
\lvec(391 26)
\ifill f:0
\move(394 25)
\lvec(404 25)
\lvec(404 26)
\lvec(394 26)
\ifill f:0
\move(405 25)
\lvec(412 25)
\lvec(412 26)
\lvec(405 26)
\ifill f:0
\move(413 25)
\lvec(420 25)
\lvec(420 26)
\lvec(413 26)
\ifill f:0
\move(421 25)
\lvec(423 25)
\lvec(423 26)
\lvec(421 26)
\ifill f:0
\move(424 25)
\lvec(428 25)
\lvec(428 26)
\lvec(424 26)
\ifill f:0
\move(429 25)
\lvec(436 25)
\lvec(436 26)
\lvec(429 26)
\ifill f:0
\move(437 25)
\lvec(442 25)
\lvec(442 26)
\lvec(437 26)
\ifill f:0
\move(443 25)
\lvec(446 25)
\lvec(446 26)
\lvec(443 26)
\ifill f:0
\move(447 25)
\lvec(451 25)
\lvec(451 26)
\lvec(447 26)
\ifill f:0
\move(11 26)
\lvec(12 26)
\lvec(12 27)
\lvec(11 27)
\ifill f:0
\move(15 26)
\lvec(17 26)
\lvec(17 27)
\lvec(15 27)
\ifill f:0
\move(23 26)
\lvec(26 26)
\lvec(26 27)
\lvec(23 27)
\ifill f:0
\move(27 26)
\lvec(28 26)
\lvec(28 27)
\lvec(27 27)
\ifill f:0
\move(36 26)
\lvec(37 26)
\lvec(37 27)
\lvec(36 27)
\ifill f:0
\move(38 26)
\lvec(39 26)
\lvec(39 27)
\lvec(38 27)
\ifill f:0
\move(40 26)
\lvec(50 26)
\lvec(50 27)
\lvec(40 27)
\ifill f:0
\move(54 26)
\lvec(55 26)
\lvec(55 27)
\lvec(54 27)
\ifill f:0
\move(57 26)
\lvec(63 26)
\lvec(63 27)
\lvec(57 27)
\ifill f:0
\move(64 26)
\lvec(65 26)
\lvec(65 27)
\lvec(64 27)
\ifill f:0
\move(66 26)
\lvec(68 26)
\lvec(68 27)
\lvec(66 27)
\ifill f:0
\move(69 26)
\lvec(71 26)
\lvec(71 27)
\lvec(69 27)
\ifill f:0
\move(72 26)
\lvec(73 26)
\lvec(73 27)
\lvec(72 27)
\ifill f:0
\move(74 26)
\lvec(75 26)
\lvec(75 27)
\lvec(74 27)
\ifill f:0
\move(76 26)
\lvec(77 26)
\lvec(77 27)
\lvec(76 27)
\ifill f:0
\move(78 26)
\lvec(79 26)
\lvec(79 27)
\lvec(78 27)
\ifill f:0
\move(80 26)
\lvec(82 26)
\lvec(82 27)
\lvec(80 27)
\ifill f:0
\move(84 26)
\lvec(85 26)
\lvec(85 27)
\lvec(84 27)
\ifill f:0
\move(86 26)
\lvec(90 26)
\lvec(90 27)
\lvec(86 27)
\ifill f:0
\move(91 26)
\lvec(92 26)
\lvec(92 27)
\lvec(91 27)
\ifill f:0
\move(93 26)
\lvec(95 26)
\lvec(95 27)
\lvec(93 27)
\ifill f:0
\move(96 26)
\lvec(98 26)
\lvec(98 27)
\lvec(96 27)
\ifill f:0
\move(99 26)
\lvec(101 26)
\lvec(101 27)
\lvec(99 27)
\ifill f:0
\move(102 26)
\lvec(105 26)
\lvec(105 27)
\lvec(102 27)
\ifill f:0
\move(106 26)
\lvec(122 26)
\lvec(122 27)
\lvec(106 27)
\ifill f:0
\move(123 26)
\lvec(125 26)
\lvec(125 27)
\lvec(123 27)
\ifill f:0
\move(127 26)
\lvec(131 26)
\lvec(131 27)
\lvec(127 27)
\ifill f:0
\move(133 26)
\lvec(138 26)
\lvec(138 27)
\lvec(133 27)
\ifill f:0
\move(139 26)
\lvec(140 26)
\lvec(140 27)
\lvec(139 27)
\ifill f:0
\move(141 26)
\lvec(145 26)
\lvec(145 27)
\lvec(141 27)
\ifill f:0
\move(146 26)
\lvec(149 26)
\lvec(149 27)
\lvec(146 27)
\ifill f:0
\move(150 26)
\lvec(167 26)
\lvec(167 27)
\lvec(150 27)
\ifill f:0
\move(169 26)
\lvec(170 26)
\lvec(170 27)
\lvec(169 27)
\ifill f:0
\move(175 26)
\lvec(177 26)
\lvec(177 27)
\lvec(175 27)
\ifill f:0
\move(186 26)
\lvec(187 26)
\lvec(187 27)
\lvec(186 27)
\ifill f:0
\move(188 26)
\lvec(197 26)
\lvec(197 27)
\lvec(188 27)
\ifill f:0
\move(198 26)
\lvec(204 26)
\lvec(204 27)
\lvec(198 27)
\ifill f:0
\move(205 26)
\lvec(215 26)
\lvec(215 27)
\lvec(205 27)
\ifill f:0
\move(217 26)
\lvec(226 26)
\lvec(226 27)
\lvec(217 27)
\ifill f:0
\move(227 26)
\lvec(233 26)
\lvec(233 27)
\lvec(227 27)
\ifill f:0
\move(234 26)
\lvec(240 26)
\lvec(240 27)
\lvec(234 27)
\ifill f:0
\move(241 26)
\lvec(247 26)
\lvec(247 27)
\lvec(241 27)
\ifill f:0
\move(248 26)
\lvec(253 26)
\lvec(253 27)
\lvec(248 27)
\ifill f:0
\move(254 26)
\lvec(259 26)
\lvec(259 27)
\lvec(254 27)
\ifill f:0
\move(260 26)
\lvec(264 26)
\lvec(264 27)
\lvec(260 27)
\ifill f:0
\move(265 26)
\lvec(275 26)
\lvec(275 27)
\lvec(265 27)
\ifill f:0
\move(276 26)
\lvec(280 26)
\lvec(280 27)
\lvec(276 27)
\ifill f:0
\move(281 26)
\lvec(290 26)
\lvec(290 27)
\lvec(281 27)
\ifill f:0
\move(291 26)
\lvec(298 26)
\lvec(298 27)
\lvec(291 27)
\ifill f:0
\move(299 26)
\lvec(302 26)
\lvec(302 27)
\lvec(299 27)
\ifill f:0
\move(303 26)
\lvec(306 26)
\lvec(306 27)
\lvec(303 27)
\ifill f:0
\move(307 26)
\lvec(310 26)
\lvec(310 27)
\lvec(307 27)
\ifill f:0
\move(311 26)
\lvec(314 26)
\lvec(314 27)
\lvec(311 27)
\ifill f:0
\move(315 26)
\lvec(318 26)
\lvec(318 27)
\lvec(315 27)
\ifill f:0
\move(319 26)
\lvec(322 26)
\lvec(322 27)
\lvec(319 27)
\ifill f:0
\move(323 26)
\lvec(333 26)
\lvec(333 27)
\lvec(323 27)
\ifill f:0
\move(334 26)
\lvec(337 26)
\lvec(337 27)
\lvec(334 27)
\ifill f:0
\move(338 26)
\lvec(344 26)
\lvec(344 27)
\lvec(338 27)
\ifill f:0
\move(345 26)
\lvec(362 26)
\lvec(362 27)
\lvec(345 27)
\ifill f:0
\move(363 26)
\lvec(364 26)
\lvec(364 27)
\lvec(363 27)
\ifill f:0
\move(365 26)
\lvec(377 26)
\lvec(377 27)
\lvec(365 27)
\ifill f:0
\move(378 26)
\lvec(380 26)
\lvec(380 27)
\lvec(378 27)
\ifill f:0
\move(381 26)
\lvec(410 26)
\lvec(410 27)
\lvec(381 27)
\ifill f:0
\move(411 26)
\lvec(413 26)
\lvec(413 27)
\lvec(411 27)
\ifill f:0
\move(414 26)
\lvec(416 26)
\lvec(416 27)
\lvec(414 27)
\ifill f:0
\move(417 26)
\lvec(427 26)
\lvec(427 27)
\lvec(417 27)
\ifill f:0
\move(428 26)
\lvec(430 26)
\lvec(430 27)
\lvec(428 27)
\ifill f:0
\move(431 26)
\lvec(442 26)
\lvec(442 27)
\lvec(431 27)
\ifill f:0
\move(443 26)
\lvec(444 26)
\lvec(444 27)
\lvec(443 27)
\ifill f:0
\move(445 26)
\lvec(449 26)
\lvec(449 27)
\lvec(445 27)
\ifill f:0
\move(450 26)
\lvec(451 26)
\lvec(451 27)
\lvec(450 27)
\ifill f:0
\move(14 27)
\lvec(15 27)
\lvec(15 28)
\lvec(14 28)
\ifill f:0
\move(16 27)
\lvec(17 27)
\lvec(17 28)
\lvec(16 28)
\ifill f:0
\move(18 27)
\lvec(19 27)
\lvec(19 28)
\lvec(18 28)
\ifill f:0
\move(20 27)
\lvec(22 27)
\lvec(22 28)
\lvec(20 28)
\ifill f:0
\move(23 27)
\lvec(26 27)
\lvec(26 28)
\lvec(23 28)
\ifill f:0
\move(36 27)
\lvec(37 27)
\lvec(37 28)
\lvec(36 28)
\ifill f:0
\move(39 27)
\lvec(43 27)
\lvec(43 28)
\lvec(39 28)
\ifill f:0
\move(44 27)
\lvec(45 27)
\lvec(45 28)
\lvec(44 28)
\ifill f:0
\move(47 27)
\lvec(48 27)
\lvec(48 28)
\lvec(47 28)
\ifill f:0
\move(49 27)
\lvec(50 27)
\lvec(50 28)
\lvec(49 28)
\ifill f:0
\move(56 27)
\lvec(58 27)
\lvec(58 28)
\lvec(56 28)
\ifill f:0
\move(60 27)
\lvec(63 27)
\lvec(63 28)
\lvec(60 28)
\ifill f:0
\move(64 27)
\lvec(65 27)
\lvec(65 28)
\lvec(64 28)
\ifill f:0
\move(66 27)
\lvec(79 27)
\lvec(79 28)
\lvec(66 28)
\ifill f:0
\move(80 27)
\lvec(82 27)
\lvec(82 28)
\lvec(80 28)
\ifill f:0
\move(84 27)
\lvec(85 27)
\lvec(85 28)
\lvec(84 28)
\ifill f:0
\move(88 27)
\lvec(89 27)
\lvec(89 28)
\lvec(88 28)
\ifill f:0
\move(90 27)
\lvec(91 27)
\lvec(91 28)
\lvec(90 28)
\ifill f:0
\move(92 27)
\lvec(96 27)
\lvec(96 28)
\lvec(92 28)
\ifill f:0
\move(97 27)
\lvec(98 27)
\lvec(98 28)
\lvec(97 28)
\ifill f:0
\move(99 27)
\lvec(101 27)
\lvec(101 28)
\lvec(99 28)
\ifill f:0
\move(102 27)
\lvec(104 27)
\lvec(104 28)
\lvec(102 28)
\ifill f:0
\move(105 27)
\lvec(107 27)
\lvec(107 28)
\lvec(105 28)
\ifill f:0
\move(108 27)
\lvec(122 27)
\lvec(122 28)
\lvec(108 28)
\ifill f:0
\move(123 27)
\lvec(124 27)
\lvec(124 28)
\lvec(123 28)
\ifill f:0
\move(126 27)
\lvec(129 27)
\lvec(129 28)
\lvec(126 28)
\ifill f:0
\move(131 27)
\lvec(134 27)
\lvec(134 28)
\lvec(131 28)
\ifill f:0
\move(136 27)
\lvec(141 27)
\lvec(141 28)
\lvec(136 28)
\ifill f:0
\move(142 27)
\lvec(145 27)
\lvec(145 28)
\lvec(142 28)
\ifill f:0
\move(146 27)
\lvec(147 27)
\lvec(147 28)
\lvec(146 28)
\ifill f:0
\move(149 27)
\lvec(155 27)
\lvec(155 28)
\lvec(149 28)
\ifill f:0
\move(156 27)
\lvec(167 27)
\lvec(167 28)
\lvec(156 28)
\ifill f:0
\move(169 27)
\lvec(170 27)
\lvec(170 28)
\lvec(169 28)
\ifill f:0
\move(171 27)
\lvec(197 27)
\lvec(197 28)
\lvec(171 28)
\ifill f:0
\move(198 27)
\lvec(209 27)
\lvec(209 28)
\lvec(198 28)
\ifill f:0
\move(212 27)
\lvec(226 27)
\lvec(226 28)
\lvec(212 28)
\ifill f:0
\move(228 27)
\lvec(235 27)
\lvec(235 28)
\lvec(228 28)
\ifill f:0
\move(236 27)
\lvec(244 27)
\lvec(244 28)
\lvec(236 28)
\ifill f:0
\move(245 27)
\lvec(252 27)
\lvec(252 28)
\lvec(245 28)
\ifill f:0
\move(253 27)
\lvec(259 27)
\lvec(259 28)
\lvec(253 28)
\ifill f:0
\move(260 27)
\lvec(266 27)
\lvec(266 28)
\lvec(260 28)
\ifill f:0
\move(267 27)
\lvec(272 27)
\lvec(272 28)
\lvec(267 28)
\ifill f:0
\move(273 27)
\lvec(278 27)
\lvec(278 28)
\lvec(273 28)
\ifill f:0
\move(279 27)
\lvec(290 27)
\lvec(290 28)
\lvec(279 28)
\ifill f:0
\move(291 27)
\lvec(294 27)
\lvec(294 28)
\lvec(291 28)
\ifill f:0
\move(295 27)
\lvec(299 27)
\lvec(299 28)
\lvec(295 28)
\ifill f:0
\move(300 27)
\lvec(304 27)
\lvec(304 28)
\lvec(300 28)
\ifill f:0
\move(305 27)
\lvec(309 27)
\lvec(309 28)
\lvec(305 28)
\ifill f:0
\move(310 27)
\lvec(313 27)
\lvec(313 28)
\lvec(310 28)
\ifill f:0
\move(314 27)
\lvec(322 27)
\lvec(322 28)
\lvec(314 28)
\ifill f:0
\move(323 27)
\lvec(326 27)
\lvec(326 28)
\lvec(323 28)
\ifill f:0
\move(327 27)
\lvec(330 27)
\lvec(330 28)
\lvec(327 28)
\ifill f:0
\move(331 27)
\lvec(334 27)
\lvec(334 28)
\lvec(331 28)
\ifill f:0
\move(335 27)
\lvec(338 27)
\lvec(338 28)
\lvec(335 28)
\ifill f:0
\move(339 27)
\lvec(342 27)
\lvec(342 28)
\lvec(339 28)
\ifill f:0
\move(343 27)
\lvec(346 27)
\lvec(346 28)
\lvec(343 28)
\ifill f:0
\move(347 27)
\lvec(362 27)
\lvec(362 28)
\lvec(347 28)
\ifill f:0
\move(363 27)
\lvec(365 27)
\lvec(365 28)
\lvec(363 28)
\ifill f:0
\move(366 27)
\lvec(368 27)
\lvec(368 28)
\lvec(366 28)
\ifill f:0
\move(369 27)
\lvec(375 27)
\lvec(375 28)
\lvec(369 28)
\ifill f:0
\move(376 27)
\lvec(379 27)
\lvec(379 28)
\lvec(376 28)
\ifill f:0
\move(380 27)
\lvec(382 27)
\lvec(382 28)
\lvec(380 28)
\ifill f:0
\move(383 27)
\lvec(392 27)
\lvec(392 28)
\lvec(383 28)
\ifill f:0
\move(393 27)
\lvec(395 27)
\lvec(395 28)
\lvec(393 28)
\ifill f:0
\move(396 27)
\lvec(411 27)
\lvec(411 28)
\lvec(396 28)
\ifill f:0
\move(412 27)
\lvec(414 27)
\lvec(414 28)
\lvec(412 28)
\ifill f:0
\move(415 27)
\lvec(417 27)
\lvec(417 28)
\lvec(415 28)
\ifill f:0
\move(418 27)
\lvec(420 27)
\lvec(420 28)
\lvec(418 28)
\ifill f:0
\move(421 27)
\lvec(423 27)
\lvec(423 28)
\lvec(421 28)
\ifill f:0
\move(424 27)
\lvec(429 27)
\lvec(429 28)
\lvec(424 28)
\ifill f:0
\move(430 27)
\lvec(432 27)
\lvec(432 28)
\lvec(430 28)
\ifill f:0
\move(433 27)
\lvec(435 27)
\lvec(435 28)
\lvec(433 28)
\ifill f:0
\move(436 27)
\lvec(438 27)
\lvec(438 28)
\lvec(436 28)
\ifill f:0
\move(439 27)
\lvec(442 27)
\lvec(442 28)
\lvec(439 28)
\ifill f:0
\move(443 27)
\lvec(444 27)
\lvec(444 28)
\lvec(443 28)
\ifill f:0
\move(445 27)
\lvec(450 27)
\lvec(450 28)
\lvec(445 28)
\ifill f:0
\move(11 28)
\lvec(12 28)
\lvec(12 29)
\lvec(11 29)
\ifill f:0
\move(14 28)
\lvec(15 28)
\lvec(15 29)
\lvec(14 29)
\ifill f:0
\move(16 28)
\lvec(17 28)
\lvec(17 29)
\lvec(16 29)
\ifill f:0
\move(18 28)
\lvec(19 28)
\lvec(19 29)
\lvec(18 29)
\ifill f:0
\move(21 28)
\lvec(26 28)
\lvec(26 29)
\lvec(21 29)
\ifill f:0
\move(28 28)
\lvec(29 28)
\lvec(29 29)
\lvec(28 29)
\ifill f:0
\move(36 28)
\lvec(37 28)
\lvec(37 29)
\lvec(36 29)
\ifill f:0
\move(38 28)
\lvec(41 28)
\lvec(41 29)
\lvec(38 29)
\ifill f:0
\move(42 28)
\lvec(46 28)
\lvec(46 29)
\lvec(42 29)
\ifill f:0
\move(49 28)
\lvec(50 28)
\lvec(50 29)
\lvec(49 29)
\ifill f:0
\move(51 28)
\lvec(52 28)
\lvec(52 29)
\lvec(51 29)
\ifill f:0
\move(57 28)
\lvec(62 28)
\lvec(62 29)
\lvec(57 29)
\ifill f:0
\move(63 28)
\lvec(65 28)
\lvec(65 29)
\lvec(63 29)
\ifill f:0
\move(66 28)
\lvec(70 28)
\lvec(70 29)
\lvec(66 29)
\ifill f:0
\move(71 28)
\lvec(73 28)
\lvec(73 29)
\lvec(71 29)
\ifill f:0
\move(75 28)
\lvec(76 28)
\lvec(76 29)
\lvec(75 29)
\ifill f:0
\move(77 28)
\lvec(79 28)
\lvec(79 29)
\lvec(77 29)
\ifill f:0
\move(80 28)
\lvec(82 28)
\lvec(82 29)
\lvec(80 29)
\ifill f:0
\move(83 28)
\lvec(84 28)
\lvec(84 29)
\lvec(83 29)
\ifill f:0
\move(85 28)
\lvec(86 28)
\lvec(86 29)
\lvec(85 29)
\ifill f:0
\move(87 28)
\lvec(88 28)
\lvec(88 29)
\lvec(87 29)
\ifill f:0
\move(89 28)
\lvec(90 28)
\lvec(90 29)
\lvec(89 29)
\ifill f:0
\move(91 28)
\lvec(93 28)
\lvec(93 29)
\lvec(91 29)
\ifill f:0
\move(95 28)
\lvec(96 28)
\lvec(96 29)
\lvec(95 29)
\ifill f:0
\move(97 28)
\lvec(98 28)
\lvec(98 29)
\lvec(97 29)
\ifill f:0
\move(99 28)
\lvec(101 28)
\lvec(101 29)
\lvec(99 29)
\ifill f:0
\move(102 28)
\lvec(106 28)
\lvec(106 29)
\lvec(102 29)
\ifill f:0
\move(107 28)
\lvec(122 28)
\lvec(122 29)
\lvec(107 29)
\ifill f:0
\move(123 28)
\lvec(124 28)
\lvec(124 29)
\lvec(123 29)
\ifill f:0
\move(125 28)
\lvec(128 28)
\lvec(128 29)
\lvec(125 29)
\ifill f:0
\move(129 28)
\lvec(132 28)
\lvec(132 29)
\lvec(129 29)
\ifill f:0
\move(133 28)
\lvec(145 28)
\lvec(145 29)
\lvec(133 29)
\ifill f:0
\move(146 28)
\lvec(147 28)
\lvec(147 29)
\lvec(146 29)
\ifill f:0
\move(148 28)
\lvec(153 28)
\lvec(153 29)
\lvec(148 29)
\ifill f:0
\move(154 28)
\lvec(159 28)
\lvec(159 29)
\lvec(154 29)
\ifill f:0
\move(161 28)
\lvec(168 28)
\lvec(168 29)
\lvec(161 29)
\ifill f:0
\move(169 28)
\lvec(180 28)
\lvec(180 29)
\lvec(169 29)
\ifill f:0
\move(181 28)
\lvec(182 28)
\lvec(182 29)
\lvec(181 29)
\ifill f:0
\move(184 28)
\lvec(197 28)
\lvec(197 29)
\lvec(184 29)
\ifill f:0
\move(198 28)
\lvec(224 28)
\lvec(224 29)
\lvec(198 29)
\ifill f:0
\move(225 28)
\lvec(226 28)
\lvec(226 29)
\lvec(225 29)
\ifill f:0
\move(230 28)
\lvec(240 28)
\lvec(240 29)
\lvec(230 29)
\ifill f:0
\move(242 28)
\lvec(251 28)
\lvec(251 29)
\lvec(242 29)
\ifill f:0
\move(252 28)
\lvec(257 28)
\lvec(257 29)
\lvec(252 29)
\ifill f:0
\move(258 28)
\lvec(260 28)
\lvec(260 29)
\lvec(258 29)
\ifill f:0
\move(261 28)
\lvec(268 28)
\lvec(268 29)
\lvec(261 29)
\ifill f:0
\move(270 28)
\lvec(283 28)
\lvec(283 29)
\lvec(270 29)
\ifill f:0
\move(284 28)
\lvec(290 28)
\lvec(290 29)
\lvec(284 29)
\ifill f:0
\move(291 28)
\lvec(295 28)
\lvec(295 29)
\lvec(291 29)
\ifill f:0
\move(296 28)
\lvec(306 28)
\lvec(306 29)
\lvec(296 29)
\ifill f:0
\move(307 28)
\lvec(322 28)
\lvec(322 29)
\lvec(307 29)
\ifill f:0
\move(323 28)
\lvec(331 28)
\lvec(331 29)
\lvec(323 29)
\ifill f:0
\move(332 28)
\lvec(340 28)
\lvec(340 29)
\lvec(332 29)
\ifill f:0
\move(341 28)
\lvec(362 28)
\lvec(362 29)
\lvec(341 29)
\ifill f:0
\move(363 28)
\lvec(365 28)
\lvec(365 29)
\lvec(363 29)
\ifill f:0
\move(366 28)
\lvec(369 28)
\lvec(369 29)
\lvec(366 29)
\ifill f:0
\move(370 28)
\lvec(380 28)
\lvec(380 29)
\lvec(370 29)
\ifill f:0
\move(381 28)
\lvec(384 28)
\lvec(384 29)
\lvec(381 29)
\ifill f:0
\move(385 28)
\lvec(391 28)
\lvec(391 29)
\lvec(385 29)
\ifill f:0
\move(392 28)
\lvec(395 28)
\lvec(395 29)
\lvec(392 29)
\ifill f:0
\move(396 28)
\lvec(398 28)
\lvec(398 29)
\lvec(396 29)
\ifill f:0
\move(399 28)
\lvec(405 28)
\lvec(405 29)
\lvec(399 29)
\ifill f:0
\move(406 28)
\lvec(412 28)
\lvec(412 29)
\lvec(406 29)
\ifill f:0
\move(413 28)
\lvec(425 28)
\lvec(425 29)
\lvec(413 29)
\ifill f:0
\move(426 28)
\lvec(442 28)
\lvec(442 29)
\lvec(426 29)
\ifill f:0
\move(443 28)
\lvec(444 28)
\lvec(444 29)
\lvec(443 29)
\ifill f:0
\move(445 28)
\lvec(447 28)
\lvec(447 29)
\lvec(445 29)
\ifill f:0
\move(448 28)
\lvec(450 28)
\lvec(450 29)
\lvec(448 29)
\ifill f:0
\move(11 29)
\lvec(12 29)
\lvec(12 30)
\lvec(11 30)
\ifill f:0
\move(16 29)
\lvec(17 29)
\lvec(17 30)
\lvec(16 30)
\ifill f:0
\move(18 29)
\lvec(22 29)
\lvec(22 30)
\lvec(18 30)
\ifill f:0
\move(24 29)
\lvec(26 29)
\lvec(26 30)
\lvec(24 30)
\ifill f:0
\move(36 29)
\lvec(37 29)
\lvec(37 30)
\lvec(36 30)
\ifill f:0
\move(38 29)
\lvec(39 29)
\lvec(39 30)
\lvec(38 30)
\ifill f:0
\move(40 29)
\lvec(43 29)
\lvec(43 30)
\lvec(40 30)
\ifill f:0
\move(44 29)
\lvec(48 29)
\lvec(48 30)
\lvec(44 30)
\ifill f:0
\move(49 29)
\lvec(50 29)
\lvec(50 30)
\lvec(49 30)
\ifill f:0
\move(51 29)
\lvec(53 29)
\lvec(53 30)
\lvec(51 30)
\ifill f:0
\move(54 29)
\lvec(55 29)
\lvec(55 30)
\lvec(54 30)
\ifill f:0
\move(56 29)
\lvec(60 29)
\lvec(60 30)
\lvec(56 30)
\ifill f:0
\move(62 29)
\lvec(65 29)
\lvec(65 30)
\lvec(62 30)
\ifill f:0
\move(66 29)
\lvec(67 29)
\lvec(67 30)
\lvec(66 30)
\ifill f:0
\move(68 29)
\lvec(75 29)
\lvec(75 30)
\lvec(68 30)
\ifill f:0
\move(76 29)
\lvec(82 29)
\lvec(82 30)
\lvec(76 30)
\ifill f:0
\move(83 29)
\lvec(84 29)
\lvec(84 30)
\lvec(83 30)
\ifill f:0
\move(85 29)
\lvec(89 29)
\lvec(89 30)
\lvec(85 30)
\ifill f:0
\move(90 29)
\lvec(91 29)
\lvec(91 30)
\lvec(90 30)
\ifill f:0
\move(92 29)
\lvec(93 29)
\lvec(93 30)
\lvec(92 30)
\ifill f:0
\move(96 29)
\lvec(98 29)
\lvec(98 30)
\lvec(96 30)
\ifill f:0
\move(100 29)
\lvec(101 29)
\lvec(101 30)
\lvec(100 30)
\ifill f:0
\move(102 29)
\lvec(118 29)
\lvec(118 30)
\lvec(102 30)
\ifill f:0
\move(119 29)
\lvec(122 29)
\lvec(122 30)
\lvec(119 30)
\ifill f:0
\move(123 29)
\lvec(124 29)
\lvec(124 30)
\lvec(123 30)
\ifill f:0
\move(125 29)
\lvec(127 29)
\lvec(127 30)
\lvec(125 30)
\ifill f:0
\move(128 29)
\lvec(134 29)
\lvec(134 30)
\lvec(128 30)
\ifill f:0
\move(135 29)
\lvec(138 29)
\lvec(138 30)
\lvec(135 30)
\ifill f:0
\move(139 29)
\lvec(142 29)
\lvec(142 30)
\lvec(139 30)
\ifill f:0
\move(143 29)
\lvec(145 29)
\lvec(145 30)
\lvec(143 30)
\ifill f:0
\move(146 29)
\lvec(151 29)
\lvec(151 30)
\lvec(146 30)
\ifill f:0
\move(152 29)
\lvec(170 29)
\lvec(170 30)
\lvec(152 30)
\ifill f:0
\move(171 29)
\lvec(176 29)
\lvec(176 30)
\lvec(171 30)
\ifill f:0
\move(178 29)
\lvec(187 29)
\lvec(187 30)
\lvec(178 30)
\ifill f:0
\move(189 29)
\lvec(197 29)
\lvec(197 30)
\lvec(189 30)
\ifill f:0
\move(198 29)
\lvec(210 29)
\lvec(210 30)
\lvec(198 30)
\ifill f:0
\move(217 29)
\lvec(219 29)
\lvec(219 30)
\lvec(217 30)
\ifill f:0
\move(225 29)
\lvec(226 29)
\lvec(226 30)
\lvec(225 30)
\ifill f:0
\move(230 29)
\lvec(249 29)
\lvec(249 30)
\lvec(230 30)
\ifill f:0
\move(251 29)
\lvec(257 29)
\lvec(257 30)
\lvec(251 30)
\ifill f:0
\move(258 29)
\lvec(262 29)
\lvec(262 30)
\lvec(258 30)
\ifill f:0
\move(263 29)
\lvec(272 29)
\lvec(272 30)
\lvec(263 30)
\ifill f:0
\move(273 29)
\lvec(281 29)
\lvec(281 30)
\lvec(273 30)
\ifill f:0
\move(282 29)
\lvec(290 29)
\lvec(290 30)
\lvec(282 30)
\ifill f:0
\move(291 29)
\lvec(296 29)
\lvec(296 30)
\lvec(291 30)
\ifill f:0
\move(297 29)
\lvec(298 29)
\lvec(298 30)
\lvec(297 30)
\ifill f:0
\move(299 29)
\lvec(303 29)
\lvec(303 30)
\lvec(299 30)
\ifill f:0
\move(304 29)
\lvec(309 29)
\lvec(309 30)
\lvec(304 30)
\ifill f:0
\move(310 29)
\lvec(315 29)
\lvec(315 30)
\lvec(310 30)
\ifill f:0
\move(316 29)
\lvec(321 29)
\lvec(321 30)
\lvec(316 30)
\ifill f:0
\move(323 29)
\lvec(332 29)
\lvec(332 30)
\lvec(323 30)
\ifill f:0
\move(333 29)
\lvec(337 29)
\lvec(337 30)
\lvec(333 30)
\ifill f:0
\move(338 29)
\lvec(342 29)
\lvec(342 30)
\lvec(338 30)
\ifill f:0
\move(343 29)
\lvec(347 29)
\lvec(347 30)
\lvec(343 30)
\ifill f:0
\move(348 29)
\lvec(362 29)
\lvec(362 30)
\lvec(348 30)
\ifill f:0
\move(363 29)
\lvec(370 29)
\lvec(370 30)
\lvec(363 30)
\ifill f:0
\move(371 29)
\lvec(378 29)
\lvec(378 30)
\lvec(371 30)
\ifill f:0
\move(379 29)
\lvec(382 29)
\lvec(382 30)
\lvec(379 30)
\ifill f:0
\move(383 29)
\lvec(390 29)
\lvec(390 30)
\lvec(383 30)
\ifill f:0
\move(391 29)
\lvec(394 29)
\lvec(394 30)
\lvec(391 30)
\ifill f:0
\move(395 29)
\lvec(398 29)
\lvec(398 30)
\lvec(395 30)
\ifill f:0
\move(399 29)
\lvec(413 29)
\lvec(413 30)
\lvec(399 30)
\ifill f:0
\move(414 29)
\lvec(420 29)
\lvec(420 30)
\lvec(414 30)
\ifill f:0
\move(421 29)
\lvec(434 29)
\lvec(434 30)
\lvec(421 30)
\ifill f:0
\move(435 29)
\lvec(442 29)
\lvec(442 30)
\lvec(435 30)
\ifill f:0
\move(443 29)
\lvec(445 29)
\lvec(445 30)
\lvec(443 30)
\ifill f:0
\move(446 29)
\lvec(451 29)
\lvec(451 30)
\lvec(446 30)
\ifill f:0
\move(16 30)
\lvec(17 30)
\lvec(17 31)
\lvec(16 31)
\ifill f:0
\move(20 30)
\lvec(22 30)
\lvec(22 31)
\lvec(20 31)
\ifill f:0
\move(25 30)
\lvec(26 30)
\lvec(26 31)
\lvec(25 31)
\ifill f:0
\move(36 30)
\lvec(37 30)
\lvec(37 31)
\lvec(36 31)
\ifill f:0
\move(38 30)
\lvec(39 30)
\lvec(39 31)
\lvec(38 31)
\ifill f:0
\move(40 30)
\lvec(48 30)
\lvec(48 31)
\lvec(40 31)
\ifill f:0
\move(49 30)
\lvec(50 30)
\lvec(50 31)
\lvec(49 31)
\ifill f:0
\move(51 30)
\lvec(53 30)
\lvec(53 31)
\lvec(51 31)
\ifill f:0
\move(56 30)
\lvec(57 30)
\lvec(57 31)
\lvec(56 31)
\ifill f:0
\move(58 30)
\lvec(59 30)
\lvec(59 31)
\lvec(58 31)
\ifill f:0
\move(60 30)
\lvec(65 30)
\lvec(65 31)
\lvec(60 31)
\ifill f:0
\move(66 30)
\lvec(68 30)
\lvec(68 31)
\lvec(66 31)
\ifill f:0
\move(69 30)
\lvec(74 30)
\lvec(74 31)
\lvec(69 31)
\ifill f:0
\move(75 30)
\lvec(82 30)
\lvec(82 31)
\lvec(75 31)
\ifill f:0
\move(83 30)
\lvec(85 30)
\lvec(85 31)
\lvec(83 31)
\ifill f:0
\move(86 30)
\lvec(87 30)
\lvec(87 31)
\lvec(86 31)
\ifill f:0
\move(88 30)
\lvec(90 30)
\lvec(90 31)
\lvec(88 31)
\ifill f:0
\move(91 30)
\lvec(95 30)
\lvec(95 31)
\lvec(91 31)
\ifill f:0
\move(96 30)
\lvec(99 30)
\lvec(99 31)
\lvec(96 31)
\ifill f:0
\move(100 30)
\lvec(101 30)
\lvec(101 31)
\lvec(100 31)
\ifill f:0
\move(102 30)
\lvec(105 30)
\lvec(105 31)
\lvec(102 31)
\ifill f:0
\move(106 30)
\lvec(107 30)
\lvec(107 31)
\lvec(106 31)
\ifill f:0
\move(108 30)
\lvec(109 30)
\lvec(109 31)
\lvec(108 31)
\ifill f:0
\move(110 30)
\lvec(111 30)
\lvec(111 31)
\lvec(110 31)
\ifill f:0
\move(112 30)
\lvec(116 30)
\lvec(116 31)
\lvec(112 31)
\ifill f:0
\move(117 30)
\lvec(122 30)
\lvec(122 31)
\lvec(117 31)
\ifill f:0
\move(124 30)
\lvec(126 30)
\lvec(126 31)
\lvec(124 31)
\ifill f:0
\move(127 30)
\lvec(129 30)
\lvec(129 31)
\lvec(127 31)
\ifill f:0
\move(130 30)
\lvec(132 30)
\lvec(132 31)
\lvec(130 31)
\ifill f:0
\move(133 30)
\lvec(135 30)
\lvec(135 31)
\lvec(133 31)
\ifill f:0
\move(136 30)
\lvec(138 30)
\lvec(138 31)
\lvec(136 31)
\ifill f:0
\move(139 30)
\lvec(142 30)
\lvec(142 31)
\lvec(139 31)
\ifill f:0
\move(143 30)
\lvec(145 30)
\lvec(145 31)
\lvec(143 31)
\ifill f:0
\move(147 30)
\lvec(150 30)
\lvec(150 31)
\lvec(147 31)
\ifill f:0
\move(151 30)
\lvec(154 30)
\lvec(154 31)
\lvec(151 31)
\ifill f:0
\move(155 30)
\lvec(163 30)
\lvec(163 31)
\lvec(155 31)
\ifill f:0
\move(164 30)
\lvec(170 30)
\lvec(170 31)
\lvec(164 31)
\ifill f:0
\move(171 30)
\lvec(173 30)
\lvec(173 31)
\lvec(171 31)
\ifill f:0
\move(174 30)
\lvec(175 30)
\lvec(175 31)
\lvec(174 31)
\ifill f:0
\move(176 30)
\lvec(190 30)
\lvec(190 31)
\lvec(176 31)
\ifill f:0
\move(192 30)
\lvec(197 30)
\lvec(197 31)
\lvec(192 31)
\ifill f:0
\move(198 30)
\lvec(202 30)
\lvec(202 31)
\lvec(198 31)
\ifill f:0
\move(203 30)
\lvec(223 30)
\lvec(223 31)
\lvec(203 31)
\ifill f:0
\move(225 30)
\lvec(226 30)
\lvec(226 31)
\lvec(225 31)
\ifill f:0
\move(232 30)
\lvec(234 30)
\lvec(234 31)
\lvec(232 31)
\ifill f:0
\move(243 30)
\lvec(257 30)
\lvec(257 31)
\lvec(243 31)
\ifill f:0
\move(258 30)
\lvec(265 30)
\lvec(265 31)
\lvec(258 31)
\ifill f:0
\move(267 30)
\lvec(278 30)
\lvec(278 31)
\lvec(267 31)
\ifill f:0
\move(280 30)
\lvec(290 30)
\lvec(290 31)
\lvec(280 31)
\ifill f:0
\move(291 30)
\lvec(298 30)
\lvec(298 31)
\lvec(291 31)
\ifill f:0
\move(300 30)
\lvec(306 30)
\lvec(306 31)
\lvec(300 31)
\ifill f:0
\move(307 30)
\lvec(314 30)
\lvec(314 31)
\lvec(307 31)
\ifill f:0
\move(315 30)
\lvec(321 30)
\lvec(321 31)
\lvec(315 31)
\ifill f:0
\move(322 30)
\lvec(327 30)
\lvec(327 31)
\lvec(322 31)
\ifill f:0
\move(328 30)
\lvec(333 30)
\lvec(333 31)
\lvec(328 31)
\ifill f:0
\move(334 30)
\lvec(339 30)
\lvec(339 31)
\lvec(334 31)
\ifill f:0
\move(340 30)
\lvec(345 30)
\lvec(345 31)
\lvec(340 31)
\ifill f:0
\move(346 30)
\lvec(356 30)
\lvec(356 31)
\lvec(346 31)
\ifill f:0
\move(357 30)
\lvec(362 30)
\lvec(362 31)
\lvec(357 31)
\ifill f:0
\move(363 30)
\lvec(366 30)
\lvec(366 31)
\lvec(363 31)
\ifill f:0
\move(367 30)
\lvec(380 30)
\lvec(380 31)
\lvec(367 31)
\ifill f:0
\move(381 30)
\lvec(385 30)
\lvec(385 31)
\lvec(381 31)
\ifill f:0
\move(386 30)
\lvec(389 30)
\lvec(389 31)
\lvec(386 31)
\ifill f:0
\move(390 30)
\lvec(394 30)
\lvec(394 31)
\lvec(390 31)
\ifill f:0
\move(395 30)
\lvec(398 30)
\lvec(398 31)
\lvec(395 31)
\ifill f:0
\move(399 30)
\lvec(402 30)
\lvec(402 31)
\lvec(399 31)
\ifill f:0
\move(403 30)
\lvec(406 30)
\lvec(406 31)
\lvec(403 31)
\ifill f:0
\move(407 30)
\lvec(410 30)
\lvec(410 31)
\lvec(407 31)
\ifill f:0
\move(411 30)
\lvec(418 30)
\lvec(418 31)
\lvec(411 31)
\ifill f:0
\move(419 30)
\lvec(422 30)
\lvec(422 31)
\lvec(419 31)
\ifill f:0
\move(423 30)
\lvec(426 30)
\lvec(426 31)
\lvec(423 31)
\ifill f:0
\move(427 30)
\lvec(442 30)
\lvec(442 31)
\lvec(427 31)
\ifill f:0
\move(443 30)
\lvec(445 30)
\lvec(445 31)
\lvec(443 31)
\ifill f:0
\move(446 30)
\lvec(448 30)
\lvec(448 31)
\lvec(446 31)
\ifill f:0
\move(449 30)
\lvec(451 30)
\lvec(451 31)
\lvec(449 31)
\ifill f:0
\move(14 31)
\lvec(17 31)
\lvec(17 32)
\lvec(14 32)
\ifill f:0
\move(18 31)
\lvec(19 31)
\lvec(19 32)
\lvec(18 32)
\ifill f:0
\move(20 31)
\lvec(21 31)
\lvec(21 32)
\lvec(20 32)
\ifill f:0
\move(23 31)
\lvec(24 31)
\lvec(24 32)
\lvec(23 32)
\ifill f:0
\move(25 31)
\lvec(26 31)
\lvec(26 32)
\lvec(25 32)
\ifill f:0
\move(27 31)
\lvec(29 31)
\lvec(29 32)
\lvec(27 32)
\ifill f:0
\move(36 31)
\lvec(37 31)
\lvec(37 32)
\lvec(36 32)
\ifill f:0
\move(38 31)
\lvec(50 31)
\lvec(50 32)
\lvec(38 32)
\ifill f:0
\move(51 31)
\lvec(52 31)
\lvec(52 32)
\lvec(51 32)
\ifill f:0
\move(55 31)
\lvec(65 31)
\lvec(65 32)
\lvec(55 32)
\ifill f:0
\move(66 31)
\lvec(71 31)
\lvec(71 32)
\lvec(66 32)
\ifill f:0
\move(72 31)
\lvec(82 31)
\lvec(82 32)
\lvec(72 32)
\ifill f:0
\move(83 31)
\lvec(85 31)
\lvec(85 32)
\lvec(83 32)
\ifill f:0
\move(87 31)
\lvec(88 31)
\lvec(88 32)
\lvec(87 32)
\ifill f:0
\move(89 31)
\lvec(91 31)
\lvec(91 32)
\lvec(89 32)
\ifill f:0
\move(95 31)
\lvec(99 31)
\lvec(99 32)
\lvec(95 32)
\ifill f:0
\move(100 31)
\lvec(101 31)
\lvec(101 32)
\lvec(100 32)
\ifill f:0
\move(102 31)
\lvec(122 31)
\lvec(122 32)
\lvec(102 32)
\ifill f:0
\move(124 31)
\lvec(125 31)
\lvec(125 32)
\lvec(124 32)
\ifill f:0
\move(126 31)
\lvec(128 31)
\lvec(128 32)
\lvec(126 32)
\ifill f:0
\move(129 31)
\lvec(131 31)
\lvec(131 32)
\lvec(129 32)
\ifill f:0
\move(132 31)
\lvec(145 31)
\lvec(145 32)
\lvec(132 32)
\ifill f:0
\move(146 31)
\lvec(149 31)
\lvec(149 32)
\lvec(146 32)
\ifill f:0
\move(150 31)
\lvec(160 31)
\lvec(160 32)
\lvec(150 32)
\ifill f:0
\move(161 31)
\lvec(170 31)
\lvec(170 32)
\lvec(161 32)
\ifill f:0
\move(171 31)
\lvec(173 31)
\lvec(173 32)
\lvec(171 32)
\ifill f:0
\move(175 31)
\lvec(179 31)
\lvec(179 32)
\lvec(175 32)
\ifill f:0
\move(180 31)
\lvec(185 31)
\lvec(185 32)
\lvec(180 32)
\ifill f:0
\move(186 31)
\lvec(191 31)
\lvec(191 32)
\lvec(186 32)
\ifill f:0
\move(193 31)
\lvec(197 31)
\lvec(197 32)
\lvec(193 32)
\ifill f:0
\move(198 31)
\lvec(200 31)
\lvec(200 32)
\lvec(198 32)
\ifill f:0
\move(201 31)
\lvec(210 31)
\lvec(210 32)
\lvec(201 32)
\ifill f:0
\move(211 31)
\lvec(223 31)
\lvec(223 32)
\lvec(211 32)
\ifill f:0
\move(225 31)
\lvec(226 31)
\lvec(226 32)
\lvec(225 32)
\ifill f:0
\move(227 31)
\lvec(257 31)
\lvec(257 32)
\lvec(227 32)
\ifill f:0
\move(258 31)
\lvec(271 31)
\lvec(271 32)
\lvec(258 32)
\ifill f:0
\move(273 31)
\lvec(290 31)
\lvec(290 32)
\lvec(273 32)
\ifill f:0
\move(292 31)
\lvec(301 31)
\lvec(301 32)
\lvec(292 32)
\ifill f:0
\move(303 31)
\lvec(311 31)
\lvec(311 32)
\lvec(303 32)
\ifill f:0
\move(312 31)
\lvec(320 31)
\lvec(320 32)
\lvec(312 32)
\ifill f:0
\move(321 31)
\lvec(325 31)
\lvec(325 32)
\lvec(321 32)
\ifill f:0
\move(326 31)
\lvec(328 31)
\lvec(328 32)
\lvec(326 32)
\ifill f:0
\move(329 31)
\lvec(335 31)
\lvec(335 32)
\lvec(329 32)
\ifill f:0
\move(336 31)
\lvec(342 31)
\lvec(342 32)
\lvec(336 32)
\ifill f:0
\move(343 31)
\lvec(355 31)
\lvec(355 32)
\lvec(343 32)
\ifill f:0
\move(356 31)
\lvec(362 31)
\lvec(362 32)
\lvec(356 32)
\ifill f:0
\move(363 31)
\lvec(367 31)
\lvec(367 32)
\lvec(363 32)
\ifill f:0
\move(368 31)
\lvec(372 31)
\lvec(372 32)
\lvec(368 32)
\ifill f:0
\move(373 31)
\lvec(383 31)
\lvec(383 32)
\lvec(373 32)
\ifill f:0
\move(384 31)
\lvec(388 31)
\lvec(388 32)
\lvec(384 32)
\ifill f:0
\move(389 31)
\lvec(402 31)
\lvec(402 32)
\lvec(389 32)
\ifill f:0
\move(403 31)
\lvec(420 31)
\lvec(420 32)
\lvec(403 32)
\ifill f:0
\move(421 31)
\lvec(433 31)
\lvec(433 32)
\lvec(421 32)
\ifill f:0
\move(434 31)
\lvec(442 31)
\lvec(442 32)
\lvec(434 32)
\ifill f:0
\move(443 31)
\lvec(445 31)
\lvec(445 32)
\lvec(443 32)
\ifill f:0
\move(446 31)
\lvec(449 31)
\lvec(449 32)
\lvec(446 32)
\ifill f:0
\move(450 31)
\lvec(451 31)
\lvec(451 32)
\lvec(450 32)
\ifill f:0
\move(11 32)
\lvec(12 32)
\lvec(12 33)
\lvec(11 33)
\ifill f:0
\move(14 32)
\lvec(17 32)
\lvec(17 33)
\lvec(14 33)
\ifill f:0
\move(18 32)
\lvec(21 32)
\lvec(21 33)
\lvec(18 33)
\ifill f:0
\move(23 32)
\lvec(24 32)
\lvec(24 33)
\lvec(23 33)
\ifill f:0
\move(25 32)
\lvec(26 32)
\lvec(26 33)
\lvec(25 33)
\ifill f:0
\move(36 32)
\lvec(37 32)
\lvec(37 33)
\lvec(36 33)
\ifill f:0
\move(38 32)
\lvec(46 32)
\lvec(46 33)
\lvec(38 33)
\ifill f:0
\move(47 32)
\lvec(50 32)
\lvec(50 33)
\lvec(47 33)
\ifill f:0
\move(53 32)
\lvec(55 32)
\lvec(55 33)
\lvec(53 33)
\ifill f:0
\move(56 32)
\lvec(57 32)
\lvec(57 33)
\lvec(56 33)
\ifill f:0
\move(59 32)
\lvec(65 32)
\lvec(65 33)
\lvec(59 33)
\ifill f:0
\move(66 32)
\lvec(75 32)
\lvec(75 33)
\lvec(66 33)
\ifill f:0
\move(76 32)
\lvec(82 32)
\lvec(82 33)
\lvec(76 33)
\ifill f:0
\move(83 32)
\lvec(86 32)
\lvec(86 33)
\lvec(83 33)
\ifill f:0
\move(87 32)
\lvec(90 32)
\lvec(90 33)
\lvec(87 33)
\ifill f:0
\move(91 32)
\lvec(93 32)
\lvec(93 33)
\lvec(91 33)
\ifill f:0
\move(95 32)
\lvec(96 32)
\lvec(96 33)
\lvec(95 33)
\ifill f:0
\move(97 32)
\lvec(99 32)
\lvec(99 33)
\lvec(97 33)
\ifill f:0
\move(100 32)
\lvec(101 32)
\lvec(101 33)
\lvec(100 33)
\ifill f:0
\move(102 32)
\lvec(104 32)
\lvec(104 33)
\lvec(102 33)
\ifill f:0
\move(105 32)
\lvec(111 32)
\lvec(111 33)
\lvec(105 33)
\ifill f:0
\move(112 32)
\lvec(113 32)
\lvec(113 33)
\lvec(112 33)
\ifill f:0
\move(114 32)
\lvec(115 32)
\lvec(115 33)
\lvec(114 33)
\ifill f:0
\move(116 32)
\lvec(119 32)
\lvec(119 33)
\lvec(116 33)
\ifill f:0
\move(120 32)
\lvec(122 32)
\lvec(122 33)
\lvec(120 33)
\ifill f:0
\move(124 32)
\lvec(125 32)
\lvec(125 33)
\lvec(124 33)
\ifill f:0
\move(126 32)
\lvec(132 32)
\lvec(132 33)
\lvec(126 33)
\ifill f:0
\move(133 32)
\lvec(137 32)
\lvec(137 33)
\lvec(133 33)
\ifill f:0
\move(138 32)
\lvec(145 32)
\lvec(145 33)
\lvec(138 33)
\ifill f:0
\move(146 32)
\lvec(158 32)
\lvec(158 33)
\lvec(146 33)
\ifill f:0
\move(159 32)
\lvec(170 32)
\lvec(170 33)
\lvec(159 33)
\ifill f:0
\move(171 32)
\lvec(173 32)
\lvec(173 33)
\lvec(171 33)
\ifill f:0
\move(174 32)
\lvec(177 32)
\lvec(177 33)
\lvec(174 33)
\ifill f:0
\move(178 32)
\lvec(187 32)
\lvec(187 33)
\lvec(178 33)
\ifill f:0
\move(188 32)
\lvec(192 32)
\lvec(192 33)
\lvec(188 33)
\ifill f:0
\move(194 32)
\lvec(197 32)
\lvec(197 33)
\lvec(194 33)
\ifill f:0
\move(198 32)
\lvec(199 32)
\lvec(199 33)
\lvec(198 33)
\ifill f:0
\move(200 32)
\lvec(215 32)
\lvec(215 33)
\lvec(200 33)
\ifill f:0
\move(216 32)
\lvec(224 32)
\lvec(224 33)
\lvec(216 33)
\ifill f:0
\move(225 32)
\lvec(240 32)
\lvec(240 33)
\lvec(225 33)
\ifill f:0
\move(242 32)
\lvec(257 32)
\lvec(257 33)
\lvec(242 33)
\ifill f:0
\move(258 32)
\lvec(288 32)
\lvec(288 33)
\lvec(258 33)
\ifill f:0
\move(289 32)
\lvec(290 32)
\lvec(290 33)
\lvec(289 33)
\ifill f:0
\move(294 32)
\lvec(306 32)
\lvec(306 33)
\lvec(294 33)
\ifill f:0
\move(308 32)
\lvec(318 32)
\lvec(318 33)
\lvec(308 33)
\ifill f:0
\move(321 32)
\lvec(325 32)
\lvec(325 33)
\lvec(321 33)
\ifill f:0
\move(326 32)
\lvec(329 32)
\lvec(329 33)
\lvec(326 33)
\ifill f:0
\move(330 32)
\lvec(338 32)
\lvec(338 33)
\lvec(330 33)
\ifill f:0
\move(339 32)
\lvec(346 32)
\lvec(346 33)
\lvec(339 33)
\ifill f:0
\move(347 32)
\lvec(362 32)
\lvec(362 33)
\lvec(347 33)
\ifill f:0
\move(363 32)
\lvec(368 32)
\lvec(368 33)
\lvec(363 33)
\ifill f:0
\move(369 32)
\lvec(374 32)
\lvec(374 33)
\lvec(369 33)
\ifill f:0
\move(375 32)
\lvec(380 32)
\lvec(380 33)
\lvec(375 33)
\ifill f:0
\move(381 32)
\lvec(386 32)
\lvec(386 33)
\lvec(381 33)
\ifill f:0
\move(387 32)
\lvec(392 32)
\lvec(392 33)
\lvec(387 33)
\ifill f:0
\move(393 32)
\lvec(397 32)
\lvec(397 33)
\lvec(393 33)
\ifill f:0
\move(398 32)
\lvec(432 32)
\lvec(432 33)
\lvec(398 33)
\ifill f:0
\move(433 32)
\lvec(442 32)
\lvec(442 33)
\lvec(433 33)
\ifill f:0
\move(443 32)
\lvec(450 32)
\lvec(450 33)
\lvec(443 33)
\ifill f:0
\move(15 33)
\lvec(17 33)
\lvec(17 34)
\lvec(15 34)
\ifill f:0
\move(18 33)
\lvec(19 33)
\lvec(19 34)
\lvec(18 34)
\ifill f:0
\move(20 33)
\lvec(22 33)
\lvec(22 34)
\lvec(20 34)
\ifill f:0
\move(23 33)
\lvec(26 33)
\lvec(26 34)
\lvec(23 34)
\ifill f:0
\move(36 33)
\lvec(37 33)
\lvec(37 34)
\lvec(36 34)
\ifill f:0
\move(40 33)
\lvec(41 33)
\lvec(41 34)
\lvec(40 34)
\ifill f:0
\move(42 33)
\lvec(46 33)
\lvec(46 34)
\lvec(42 34)
\ifill f:0
\move(47 33)
\lvec(50 33)
\lvec(50 34)
\lvec(47 34)
\ifill f:0
\move(52 33)
\lvec(53 33)
\lvec(53 34)
\lvec(52 34)
\ifill f:0
\move(56 33)
\lvec(59 33)
\lvec(59 34)
\lvec(56 34)
\ifill f:0
\move(61 33)
\lvec(65 33)
\lvec(65 34)
\lvec(61 34)
\ifill f:0
\move(66 33)
\lvec(71 33)
\lvec(71 34)
\lvec(66 34)
\ifill f:0
\move(72 33)
\lvec(74 33)
\lvec(74 34)
\lvec(72 34)
\ifill f:0
\move(75 33)
\lvec(82 33)
\lvec(82 34)
\lvec(75 34)
\ifill f:0
\move(84 33)
\lvec(87 33)
\lvec(87 34)
\lvec(84 34)
\ifill f:0
\move(88 33)
\lvec(95 33)
\lvec(95 34)
\lvec(88 34)
\ifill f:0
\move(97 33)
\lvec(99 33)
\lvec(99 34)
\lvec(97 34)
\ifill f:0
\move(100 33)
\lvec(101 33)
\lvec(101 34)
\lvec(100 34)
\ifill f:0
\move(102 33)
\lvec(107 33)
\lvec(107 34)
\lvec(102 34)
\ifill f:0
\move(108 33)
\lvec(112 33)
\lvec(112 34)
\lvec(108 34)
\ifill f:0
\move(113 33)
\lvec(119 33)
\lvec(119 34)
\lvec(113 34)
\ifill f:0
\move(120 33)
\lvec(122 33)
\lvec(122 34)
\lvec(120 34)
\ifill f:0
\move(124 33)
\lvec(127 33)
\lvec(127 34)
\lvec(124 34)
\ifill f:0
\move(128 33)
\lvec(129 33)
\lvec(129 34)
\lvec(128 34)
\ifill f:0
\move(130 33)
\lvec(131 33)
\lvec(131 34)
\lvec(130 34)
\ifill f:0
\move(132 33)
\lvec(133 33)
\lvec(133 34)
\lvec(132 34)
\ifill f:0
\move(134 33)
\lvec(138 33)
\lvec(138 34)
\lvec(134 34)
\ifill f:0
\move(139 33)
\lvec(145 33)
\lvec(145 34)
\lvec(139 34)
\ifill f:0
\move(146 33)
\lvec(159 33)
\lvec(159 34)
\lvec(146 34)
\ifill f:0
\move(160 33)
\lvec(170 33)
\lvec(170 34)
\lvec(160 34)
\ifill f:0
\move(171 33)
\lvec(172 33)
\lvec(172 34)
\lvec(171 34)
\ifill f:0
\move(173 33)
\lvec(176 33)
\lvec(176 34)
\lvec(173 34)
\ifill f:0
\move(177 33)
\lvec(180 33)
\lvec(180 34)
\lvec(177 34)
\ifill f:0
\move(181 33)
\lvec(189 33)
\lvec(189 34)
\lvec(181 34)
\ifill f:0
\move(190 33)
\lvec(197 33)
\lvec(197 34)
\lvec(190 34)
\ifill f:0
\move(198 33)
\lvec(204 33)
\lvec(204 34)
\lvec(198 34)
\ifill f:0
\move(205 33)
\lvec(210 33)
\lvec(210 34)
\lvec(205 34)
\ifill f:0
\move(211 33)
\lvec(233 33)
\lvec(233 34)
\lvec(211 34)
\ifill f:0
\move(236 33)
\lvec(247 33)
\lvec(247 34)
\lvec(236 34)
\ifill f:0
\move(248 33)
\lvec(257 33)
\lvec(257 34)
\lvec(248 34)
\ifill f:0
\move(258 33)
\lvec(269 33)
\lvec(269 34)
\lvec(258 34)
\ifill f:0
\move(270 33)
\lvec(271 33)
\lvec(271 34)
\lvec(270 34)
\ifill f:0
\move(280 33)
\lvec(282 33)
\lvec(282 34)
\lvec(280 34)
\ifill f:0
\move(289 33)
\lvec(290 33)
\lvec(290 34)
\lvec(289 34)
\ifill f:0
\move(294 33)
\lvec(295 33)
\lvec(295 34)
\lvec(294 34)
\ifill f:0
\move(296 33)
\lvec(298 33)
\lvec(298 34)
\lvec(296 34)
\ifill f:0
\move(299 33)
\lvec(316 33)
\lvec(316 34)
\lvec(299 34)
\ifill f:0
\move(318 33)
\lvec(325 33)
\lvec(325 34)
\lvec(318 34)
\ifill f:0
\move(326 33)
\lvec(330 33)
\lvec(330 34)
\lvec(326 34)
\ifill f:0
\move(331 33)
\lvec(342 33)
\lvec(342 34)
\lvec(331 34)
\ifill f:0
\move(343 33)
\lvec(352 33)
\lvec(352 34)
\lvec(343 34)
\ifill f:0
\move(353 33)
\lvec(362 33)
\lvec(362 34)
\lvec(353 34)
\ifill f:0
\move(363 33)
\lvec(369 33)
\lvec(369 34)
\lvec(363 34)
\ifill f:0
\move(371 33)
\lvec(377 33)
\lvec(377 34)
\lvec(371 34)
\ifill f:0
\move(378 33)
\lvec(390 33)
\lvec(390 34)
\lvec(378 34)
\ifill f:0
\move(391 33)
\lvec(397 33)
\lvec(397 34)
\lvec(391 34)
\ifill f:0
\move(398 33)
\lvec(403 33)
\lvec(403 34)
\lvec(398 34)
\ifill f:0
\move(404 33)
\lvec(409 33)
\lvec(409 34)
\lvec(404 34)
\ifill f:0
\move(410 33)
\lvec(420 33)
\lvec(420 34)
\lvec(410 34)
\ifill f:0
\move(421 33)
\lvec(436 33)
\lvec(436 34)
\lvec(421 34)
\ifill f:0
\move(437 33)
\lvec(442 33)
\lvec(442 34)
\lvec(437 34)
\ifill f:0
\move(443 33)
\lvec(446 33)
\lvec(446 34)
\lvec(443 34)
\ifill f:0
\move(447 33)
\lvec(451 33)
\lvec(451 34)
\lvec(447 34)
\ifill f:0
\move(15 34)
\lvec(17 34)
\lvec(17 35)
\lvec(15 35)
\ifill f:0
\move(20 34)
\lvec(22 34)
\lvec(22 35)
\lvec(20 35)
\ifill f:0
\move(24 34)
\lvec(26 34)
\lvec(26 35)
\lvec(24 35)
\ifill f:0
\move(28 34)
\lvec(29 34)
\lvec(29 35)
\lvec(28 35)
\ifill f:0
\move(36 34)
\lvec(37 34)
\lvec(37 35)
\lvec(36 35)
\ifill f:0
\move(38 34)
\lvec(39 34)
\lvec(39 35)
\lvec(38 35)
\ifill f:0
\move(40 34)
\lvec(47 34)
\lvec(47 35)
\lvec(40 35)
\ifill f:0
\move(48 34)
\lvec(50 34)
\lvec(50 35)
\lvec(48 35)
\ifill f:0
\move(52 34)
\lvec(53 34)
\lvec(53 35)
\lvec(52 35)
\ifill f:0
\move(56 34)
\lvec(57 34)
\lvec(57 35)
\lvec(56 35)
\ifill f:0
\move(58 34)
\lvec(61 34)
\lvec(61 35)
\lvec(58 35)
\ifill f:0
\move(62 34)
\lvec(65 34)
\lvec(65 35)
\lvec(62 35)
\ifill f:0
\move(66 34)
\lvec(68 34)
\lvec(68 35)
\lvec(66 35)
\ifill f:0
\move(69 34)
\lvec(82 34)
\lvec(82 35)
\lvec(69 35)
\ifill f:0
\move(86 34)
\lvec(89 34)
\lvec(89 35)
\lvec(86 35)
\ifill f:0
\move(92 34)
\lvec(93 34)
\lvec(93 35)
\lvec(92 35)
\ifill f:0
\move(96 34)
\lvec(98 34)
\lvec(98 35)
\lvec(96 35)
\ifill f:0
\move(100 34)
\lvec(101 34)
\lvec(101 35)
\lvec(100 35)
\ifill f:0
\move(102 34)
\lvec(105 34)
\lvec(105 35)
\lvec(102 35)
\ifill f:0
\move(106 34)
\lvec(108 34)
\lvec(108 35)
\lvec(106 35)
\ifill f:0
\move(109 34)
\lvec(111 34)
\lvec(111 35)
\lvec(109 35)
\ifill f:0
\move(112 34)
\lvec(122 34)
\lvec(122 35)
\lvec(112 35)
\ifill f:0
\move(123 34)
\lvec(124 34)
\lvec(124 35)
\lvec(123 35)
\ifill f:0
\move(125 34)
\lvec(126 34)
\lvec(126 35)
\lvec(125 35)
\ifill f:0
\move(127 34)
\lvec(128 34)
\lvec(128 35)
\lvec(127 35)
\ifill f:0
\move(129 34)
\lvec(130 34)
\lvec(130 35)
\lvec(129 35)
\ifill f:0
\move(131 34)
\lvec(132 34)
\lvec(132 35)
\lvec(131 35)
\ifill f:0
\move(133 34)
\lvec(134 34)
\lvec(134 35)
\lvec(133 35)
\ifill f:0
\move(135 34)
\lvec(138 34)
\lvec(138 35)
\lvec(135 35)
\ifill f:0
\move(139 34)
\lvec(145 34)
\lvec(145 35)
\lvec(139 35)
\ifill f:0
\move(146 34)
\lvec(150 34)
\lvec(150 35)
\lvec(146 35)
\ifill f:0
\move(151 34)
\lvec(155 34)
\lvec(155 35)
\lvec(151 35)
\ifill f:0
\move(156 34)
\lvec(160 34)
\lvec(160 35)
\lvec(156 35)
\ifill f:0
\move(161 34)
\lvec(163 34)
\lvec(163 35)
\lvec(161 35)
\ifill f:0
\move(164 34)
\lvec(166 34)
\lvec(166 35)
\lvec(164 35)
\ifill f:0
\move(167 34)
\lvec(170 34)
\lvec(170 35)
\lvec(167 35)
\ifill f:0
\move(171 34)
\lvec(172 34)
\lvec(172 35)
\lvec(171 35)
\ifill f:0
\move(173 34)
\lvec(175 34)
\lvec(175 35)
\lvec(173 35)
\ifill f:0
\move(176 34)
\lvec(186 34)
\lvec(186 35)
\lvec(176 35)
\ifill f:0
\move(187 34)
\lvec(194 34)
\lvec(194 35)
\lvec(187 35)
\ifill f:0
\move(195 34)
\lvec(197 34)
\lvec(197 35)
\lvec(195 35)
\ifill f:0
\move(198 34)
\lvec(203 34)
\lvec(203 35)
\lvec(198 35)
\ifill f:0
\move(204 34)
\lvec(213 34)
\lvec(213 35)
\lvec(204 35)
\ifill f:0
\move(214 34)
\lvec(226 34)
\lvec(226 35)
\lvec(214 35)
\ifill f:0
\move(227 34)
\lvec(232 34)
\lvec(232 35)
\lvec(227 35)
\ifill f:0
\move(233 34)
\lvec(240 34)
\lvec(240 35)
\lvec(233 35)
\ifill f:0
\move(241 34)
\lvec(250 34)
\lvec(250 35)
\lvec(241 35)
\ifill f:0
\move(251 34)
\lvec(257 34)
\lvec(257 35)
\lvec(251 35)
\ifill f:0
\move(258 34)
\lvec(287 34)
\lvec(287 35)
\lvec(258 35)
\ifill f:0
\move(289 34)
\lvec(290 34)
\lvec(290 35)
\lvec(289 35)
\ifill f:0
\move(297 34)
\lvec(299 34)
\lvec(299 35)
\lvec(297 35)
\ifill f:0
\move(312 34)
\lvec(325 34)
\lvec(325 35)
\lvec(312 35)
\ifill f:0
\move(326 34)
\lvec(334 34)
\lvec(334 35)
\lvec(326 35)
\ifill f:0
\move(336 34)
\lvec(349 34)
\lvec(349 35)
\lvec(336 35)
\ifill f:0
\move(350 34)
\lvec(362 34)
\lvec(362 35)
\lvec(350 35)
\ifill f:0
\move(363 34)
\lvec(371 34)
\lvec(371 35)
\lvec(363 35)
\ifill f:0
\move(372 34)
\lvec(380 34)
\lvec(380 35)
\lvec(372 35)
\ifill f:0
\move(381 34)
\lvec(389 34)
\lvec(389 35)
\lvec(381 35)
\ifill f:0
\move(390 34)
\lvec(396 34)
\lvec(396 35)
\lvec(390 35)
\ifill f:0
\move(397 34)
\lvec(401 34)
\lvec(401 35)
\lvec(397 35)
\ifill f:0
\move(402 34)
\lvec(410 34)
\lvec(410 35)
\lvec(402 35)
\ifill f:0
\move(411 34)
\lvec(417 34)
\lvec(417 35)
\lvec(411 35)
\ifill f:0
\move(418 34)
\lvec(423 34)
\lvec(423 35)
\lvec(418 35)
\ifill f:0
\move(424 34)
\lvec(435 34)
\lvec(435 35)
\lvec(424 35)
\ifill f:0
\move(436 34)
\lvec(442 34)
\lvec(442 35)
\lvec(436 35)
\ifill f:0
\move(443 34)
\lvec(451 34)
\lvec(451 35)
\lvec(443 35)
\ifill f:0
\move(11 35)
\lvec(12 35)
\lvec(12 36)
\lvec(11 36)
\ifill f:0
\move(14 35)
\lvec(15 35)
\lvec(15 36)
\lvec(14 36)
\ifill f:0
\move(16 35)
\lvec(17 35)
\lvec(17 36)
\lvec(16 36)
\ifill f:0
\move(18 35)
\lvec(22 35)
\lvec(22 36)
\lvec(18 36)
\ifill f:0
\move(24 35)
\lvec(26 35)
\lvec(26 36)
\lvec(24 36)
\ifill f:0
\move(36 35)
\lvec(37 35)
\lvec(37 36)
\lvec(36 36)
\ifill f:0
\move(38 35)
\lvec(47 35)
\lvec(47 36)
\lvec(38 36)
\ifill f:0
\move(48 35)
\lvec(50 35)
\lvec(50 36)
\lvec(48 36)
\ifill f:0
\move(54 35)
\lvec(55 35)
\lvec(55 36)
\lvec(54 36)
\ifill f:0
\move(56 35)
\lvec(58 35)
\lvec(58 36)
\lvec(56 36)
\ifill f:0
\move(59 35)
\lvec(65 35)
\lvec(65 36)
\lvec(59 36)
\ifill f:0
\move(66 35)
\lvec(74 35)
\lvec(74 36)
\lvec(66 36)
\ifill f:0
\move(79 35)
\lvec(80 35)
\lvec(80 36)
\lvec(79 36)
\ifill f:0
\move(81 35)
\lvec(82 35)
\lvec(82 36)
\lvec(81 36)
\ifill f:0
\move(87 35)
\lvec(92 35)
\lvec(92 36)
\lvec(87 36)
\ifill f:0
\move(95 35)
\lvec(98 35)
\lvec(98 36)
\lvec(95 36)
\ifill f:0
\move(100 35)
\lvec(101 35)
\lvec(101 36)
\lvec(100 36)
\ifill f:0
\move(103 35)
\lvec(106 35)
\lvec(106 36)
\lvec(103 36)
\ifill f:0
\move(107 35)
\lvec(122 35)
\lvec(122 36)
\lvec(107 36)
\ifill f:0
\move(123 35)
\lvec(124 35)
\lvec(124 36)
\lvec(123 36)
\ifill f:0
\move(125 35)
\lvec(126 35)
\lvec(126 36)
\lvec(125 36)
\ifill f:0
\move(127 35)
\lvec(131 35)
\lvec(131 36)
\lvec(127 36)
\ifill f:0
\move(132 35)
\lvec(133 35)
\lvec(133 36)
\lvec(132 36)
\ifill f:0
\move(134 35)
\lvec(135 35)
\lvec(135 36)
\lvec(134 36)
\ifill f:0
\move(136 35)
\lvec(137 35)
\lvec(137 36)
\lvec(136 36)
\ifill f:0
\move(138 35)
\lvec(139 35)
\lvec(139 36)
\lvec(138 36)
\ifill f:0
\move(140 35)
\lvec(141 35)
\lvec(141 36)
\lvec(140 36)
\ifill f:0
\move(142 35)
\lvec(143 35)
\lvec(143 36)
\lvec(142 36)
\ifill f:0
\move(144 35)
\lvec(145 35)
\lvec(145 36)
\lvec(144 36)
\ifill f:0
\move(146 35)
\lvec(147 35)
\lvec(147 36)
\lvec(146 36)
\ifill f:0
\move(148 35)
\lvec(161 35)
\lvec(161 36)
\lvec(148 36)
\ifill f:0
\move(162 35)
\lvec(170 35)
\lvec(170 36)
\lvec(162 36)
\ifill f:0
\move(171 35)
\lvec(184 35)
\lvec(184 36)
\lvec(171 36)
\ifill f:0
\move(185 35)
\lvec(187 35)
\lvec(187 36)
\lvec(185 36)
\ifill f:0
\move(188 35)
\lvec(194 35)
\lvec(194 36)
\lvec(188 36)
\ifill f:0
\move(195 35)
\lvec(197 35)
\lvec(197 36)
\lvec(195 36)
\ifill f:0
\move(199 35)
\lvec(206 35)
\lvec(206 36)
\lvec(199 36)
\ifill f:0
\move(207 35)
\lvec(210 35)
\lvec(210 36)
\lvec(207 36)
\ifill f:0
\move(211 35)
\lvec(215 35)
\lvec(215 36)
\lvec(211 36)
\ifill f:0
\move(216 35)
\lvec(220 35)
\lvec(220 36)
\lvec(216 36)
\ifill f:0
\move(221 35)
\lvec(226 35)
\lvec(226 36)
\lvec(221 36)
\ifill f:0
\move(227 35)
\lvec(237 35)
\lvec(237 36)
\lvec(227 36)
\ifill f:0
\move(238 35)
\lvec(244 35)
\lvec(244 36)
\lvec(238 36)
\ifill f:0
\move(245 35)
\lvec(250 35)
\lvec(250 36)
\lvec(245 36)
\ifill f:0
\move(251 35)
\lvec(252 35)
\lvec(252 36)
\lvec(251 36)
\ifill f:0
\move(253 35)
\lvec(257 35)
\lvec(257 36)
\lvec(253 36)
\ifill f:0
\move(258 35)
\lvec(261 35)
\lvec(261 36)
\lvec(258 36)
\ifill f:0
\move(262 35)
\lvec(271 35)
\lvec(271 36)
\lvec(262 36)
\ifill f:0
\move(274 35)
\lvec(288 35)
\lvec(288 36)
\lvec(274 36)
\ifill f:0
\move(289 35)
\lvec(290 35)
\lvec(290 36)
\lvec(289 36)
\ifill f:0
\move(291 35)
\lvec(325 35)
\lvec(325 36)
\lvec(291 36)
\ifill f:0
\move(326 35)
\lvec(341 35)
\lvec(341 36)
\lvec(326 36)
\ifill f:0
\move(344 35)
\lvec(360 35)
\lvec(360 36)
\lvec(344 36)
\ifill f:0
\move(361 35)
\lvec(362 35)
\lvec(362 36)
\lvec(361 36)
\ifill f:0
\move(364 35)
\lvec(374 35)
\lvec(374 36)
\lvec(364 36)
\ifill f:0
\move(376 35)
\lvec(385 35)
\lvec(385 36)
\lvec(376 36)
\ifill f:0
\move(387 35)
\lvec(395 35)
\lvec(395 36)
\lvec(387 36)
\ifill f:0
\move(396 35)
\lvec(401 35)
\lvec(401 36)
\lvec(396 36)
\ifill f:0
\move(402 35)
\lvec(404 35)
\lvec(404 36)
\lvec(402 36)
\ifill f:0
\move(405 35)
\lvec(420 35)
\lvec(420 36)
\lvec(405 36)
\ifill f:0
\move(421 35)
\lvec(427 35)
\lvec(427 36)
\lvec(421 36)
\ifill f:0
\move(428 35)
\lvec(434 35)
\lvec(434 36)
\lvec(428 36)
\ifill f:0
\move(435 35)
\lvec(442 35)
\lvec(442 36)
\lvec(435 36)
\ifill f:0
\move(443 35)
\lvec(447 35)
\lvec(447 36)
\lvec(443 36)
\ifill f:0
\move(448 35)
\lvec(451 35)
\lvec(451 36)
\lvec(448 36)
\ifill f:0
\move(14 36)
\lvec(15 36)
\lvec(15 37)
\lvec(14 37)
\ifill f:0
\move(16 36)
\lvec(17 36)
\lvec(17 37)
\lvec(16 37)
\ifill f:0
\move(18 36)
\lvec(19 36)
\lvec(19 37)
\lvec(18 37)
\ifill f:0
\move(20 36)
\lvec(22 36)
\lvec(22 37)
\lvec(20 37)
\ifill f:0
\move(23 36)
\lvec(26 36)
\lvec(26 37)
\lvec(23 37)
\ifill f:0
\move(27 36)
\lvec(28 36)
\lvec(28 37)
\lvec(27 37)
\ifill f:0
\move(36 36)
\lvec(37 36)
\lvec(37 37)
\lvec(36 37)
\ifill f:0
\move(38 36)
\lvec(43 36)
\lvec(43 37)
\lvec(38 37)
\ifill f:0
\move(44 36)
\lvec(50 36)
\lvec(50 37)
\lvec(44 37)
\ifill f:0
\move(51 36)
\lvec(52 36)
\lvec(52 37)
\lvec(51 37)
\ifill f:0
\move(56 36)
\lvec(59 36)
\lvec(59 37)
\lvec(56 37)
\ifill f:0
\move(60 36)
\lvec(62 36)
\lvec(62 37)
\lvec(60 37)
\ifill f:0
\move(63 36)
\lvec(65 36)
\lvec(65 37)
\lvec(63 37)
\ifill f:0
\move(66 36)
\lvec(70 36)
\lvec(70 37)
\lvec(66 37)
\ifill f:0
\move(71 36)
\lvec(79 36)
\lvec(79 37)
\lvec(71 37)
\ifill f:0
\move(81 36)
\lvec(82 36)
\lvec(82 37)
\lvec(81 37)
\ifill f:0
\move(83 36)
\lvec(84 36)
\lvec(84 37)
\lvec(83 37)
\ifill f:0
\move(92 36)
\lvec(98 36)
\lvec(98 37)
\lvec(92 37)
\ifill f:0
\move(100 36)
\lvec(101 36)
\lvec(101 37)
\lvec(100 37)
\ifill f:0
\move(102 36)
\lvec(103 36)
\lvec(103 37)
\lvec(102 37)
\ifill f:0
\move(104 36)
\lvec(118 36)
\lvec(118 37)
\lvec(104 37)
\ifill f:0
\move(119 36)
\lvec(122 36)
\lvec(122 37)
\lvec(119 37)
\ifill f:0
\move(123 36)
\lvec(124 36)
\lvec(124 37)
\lvec(123 37)
\ifill f:0
\move(125 36)
\lvec(127 36)
\lvec(127 37)
\lvec(125 37)
\ifill f:0
\move(128 36)
\lvec(132 36)
\lvec(132 37)
\lvec(128 37)
\ifill f:0
\move(133 36)
\lvec(134 36)
\lvec(134 37)
\lvec(133 37)
\ifill f:0
\move(135 36)
\lvec(143 36)
\lvec(143 37)
\lvec(135 37)
\ifill f:0
\move(144 36)
\lvec(145 36)
\lvec(145 37)
\lvec(144 37)
\ifill f:0
\move(146 36)
\lvec(147 36)
\lvec(147 37)
\lvec(146 37)
\ifill f:0
\move(148 36)
\lvec(149 36)
\lvec(149 37)
\lvec(148 37)
\ifill f:0
\move(150 36)
\lvec(151 36)
\lvec(151 37)
\lvec(150 37)
\ifill f:0
\move(152 36)
\lvec(153 36)
\lvec(153 37)
\lvec(152 37)
\ifill f:0
\move(154 36)
\lvec(155 36)
\lvec(155 37)
\lvec(154 37)
\ifill f:0
\move(156 36)
\lvec(170 36)
\lvec(170 37)
\lvec(156 37)
\ifill f:0
\move(172 36)
\lvec(174 36)
\lvec(174 37)
\lvec(172 37)
\ifill f:0
\move(175 36)
\lvec(185 36)
\lvec(185 37)
\lvec(175 37)
\ifill f:0
\move(186 36)
\lvec(188 36)
\lvec(188 37)
\lvec(186 37)
\ifill f:0
\move(189 36)
\lvec(191 36)
\lvec(191 37)
\lvec(189 37)
\ifill f:0
\move(192 36)
\lvec(194 36)
\lvec(194 37)
\lvec(192 37)
\ifill f:0
\move(195 36)
\lvec(197 36)
\lvec(197 37)
\lvec(195 37)
\ifill f:0
\move(198 36)
\lvec(201 36)
\lvec(201 37)
\lvec(198 37)
\ifill f:0
\move(202 36)
\lvec(226 36)
\lvec(226 37)
\lvec(202 37)
\ifill f:0
\move(227 36)
\lvec(229 36)
\lvec(229 37)
\lvec(227 37)
\ifill f:0
\move(231 36)
\lvec(235 36)
\lvec(235 37)
\lvec(231 37)
\ifill f:0
\move(236 36)
\lvec(240 36)
\lvec(240 37)
\lvec(236 37)
\ifill f:0
\move(241 36)
\lvec(246 36)
\lvec(246 37)
\lvec(241 37)
\ifill f:0
\move(247 36)
\lvec(257 36)
\lvec(257 37)
\lvec(247 37)
\ifill f:0
\move(258 36)
\lvec(260 36)
\lvec(260 37)
\lvec(258 37)
\ifill f:0
\move(261 36)
\lvec(268 36)
\lvec(268 37)
\lvec(261 37)
\ifill f:0
\move(269 36)
\lvec(277 36)
\lvec(277 37)
\lvec(269 37)
\ifill f:0
\move(278 36)
\lvec(288 36)
\lvec(288 37)
\lvec(278 37)
\ifill f:0
\move(289 36)
\lvec(305 36)
\lvec(305 37)
\lvec(289 37)
\ifill f:0
\move(308 36)
\lvec(325 36)
\lvec(325 37)
\lvec(308 37)
\ifill f:0
\move(326 36)
\lvec(360 36)
\lvec(360 37)
\lvec(326 37)
\ifill f:0
\move(361 36)
\lvec(362 36)
\lvec(362 37)
\lvec(361 37)
\ifill f:0
\move(366 36)
\lvec(380 36)
\lvec(380 37)
\lvec(366 37)
\ifill f:0
\move(381 36)
\lvec(394 36)
\lvec(394 37)
\lvec(381 37)
\ifill f:0
\move(396 36)
\lvec(401 36)
\lvec(401 37)
\lvec(396 37)
\ifill f:0
\move(402 36)
\lvec(405 36)
\lvec(405 37)
\lvec(402 37)
\ifill f:0
\move(406 36)
\lvec(415 36)
\lvec(415 37)
\lvec(406 37)
\ifill f:0
\move(416 36)
\lvec(433 36)
\lvec(433 37)
\lvec(416 37)
\ifill f:0
\move(434 36)
\lvec(442 36)
\lvec(442 37)
\lvec(434 37)
\ifill f:0
\move(443 36)
\lvec(448 36)
\lvec(448 37)
\lvec(443 37)
\ifill f:0
\move(450 36)
\lvec(451 36)
\lvec(451 37)
\lvec(450 37)
\ifill f:0
\move(16 37)
\lvec(17 37)
\lvec(17 38)
\lvec(16 38)
\ifill f:0
\move(18 37)
\lvec(19 37)
\lvec(19 38)
\lvec(18 38)
\ifill f:0
\move(20 37)
\lvec(21 37)
\lvec(21 38)
\lvec(20 38)
\ifill f:0
\move(23 37)
\lvec(26 37)
\lvec(26 38)
\lvec(23 38)
\ifill f:0
\move(36 37)
\lvec(37 37)
\lvec(37 38)
\lvec(36 38)
\ifill f:0
\move(38 37)
\lvec(39 37)
\lvec(39 38)
\lvec(38 38)
\ifill f:0
\move(40 37)
\lvec(46 37)
\lvec(46 38)
\lvec(40 38)
\ifill f:0
\move(47 37)
\lvec(50 37)
\lvec(50 38)
\lvec(47 38)
\ifill f:0
\move(51 37)
\lvec(52 37)
\lvec(52 38)
\lvec(51 38)
\ifill f:0
\move(54 37)
\lvec(55 37)
\lvec(55 38)
\lvec(54 38)
\ifill f:0
\move(56 37)
\lvec(60 37)
\lvec(60 38)
\lvec(56 38)
\ifill f:0
\move(61 37)
\lvec(65 37)
\lvec(65 38)
\lvec(61 38)
\ifill f:0
\move(66 37)
\lvec(69 37)
\lvec(69 38)
\lvec(66 38)
\ifill f:0
\move(70 37)
\lvec(73 37)
\lvec(73 38)
\lvec(70 38)
\ifill f:0
\move(74 37)
\lvec(79 37)
\lvec(79 38)
\lvec(74 38)
\ifill f:0
\move(81 37)
\lvec(82 37)
\lvec(82 38)
\lvec(81 38)
\ifill f:0
\move(83 37)
\lvec(93 37)
\lvec(93 38)
\lvec(83 38)
\ifill f:0
\move(94 37)
\lvec(95 37)
\lvec(95 38)
\lvec(94 38)
\ifill f:0
\move(97 37)
\lvec(98 37)
\lvec(98 38)
\lvec(97 38)
\ifill f:0
\move(100 37)
\lvec(101 37)
\lvec(101 38)
\lvec(100 38)
\ifill f:0
\move(102 37)
\lvec(109 37)
\lvec(109 38)
\lvec(102 38)
\ifill f:0
\move(110 37)
\lvec(122 37)
\lvec(122 38)
\lvec(110 38)
\ifill f:0
\move(123 37)
\lvec(125 37)
\lvec(125 38)
\lvec(123 38)
\ifill f:0
\move(126 37)
\lvec(131 37)
\lvec(131 38)
\lvec(126 38)
\ifill f:0
\move(132 37)
\lvec(133 37)
\lvec(133 38)
\lvec(132 38)
\ifill f:0
\move(134 37)
\lvec(138 37)
\lvec(138 38)
\lvec(134 38)
\ifill f:0
\move(139 37)
\lvec(143 37)
\lvec(143 38)
\lvec(139 38)
\ifill f:0
\move(144 37)
\lvec(145 37)
\lvec(145 38)
\lvec(144 38)
\ifill f:0
\move(146 37)
\lvec(160 37)
\lvec(160 38)
\lvec(146 38)
\ifill f:0
\move(161 37)
\lvec(170 37)
\lvec(170 38)
\lvec(161 38)
\ifill f:0
\move(172 37)
\lvec(173 37)
\lvec(173 38)
\lvec(172 38)
\ifill f:0
\move(174 37)
\lvec(176 37)
\lvec(176 38)
\lvec(174 38)
\ifill f:0
\move(177 37)
\lvec(181 37)
\lvec(181 38)
\lvec(177 38)
\ifill f:0
\move(182 37)
\lvec(186 37)
\lvec(186 38)
\lvec(182 38)
\ifill f:0
\move(187 37)
\lvec(189 37)
\lvec(189 38)
\lvec(187 38)
\ifill f:0
\move(190 37)
\lvec(191 37)
\lvec(191 38)
\lvec(190 38)
\ifill f:0
\move(192 37)
\lvec(194 37)
\lvec(194 38)
\lvec(192 38)
\ifill f:0
\move(195 37)
\lvec(197 37)
\lvec(197 38)
\lvec(195 38)
\ifill f:0
\move(198 37)
\lvec(207 37)
\lvec(207 38)
\lvec(198 38)
\ifill f:0
\move(208 37)
\lvec(226 37)
\lvec(226 38)
\lvec(208 38)
\ifill f:0
\move(227 37)
\lvec(229 37)
\lvec(229 38)
\lvec(227 38)
\ifill f:0
\move(230 37)
\lvec(233 37)
\lvec(233 38)
\lvec(230 38)
\ifill f:0
\move(234 37)
\lvec(238 37)
\lvec(238 38)
\lvec(234 38)
\ifill f:0
\move(239 37)
\lvec(242 37)
\lvec(242 38)
\lvec(239 38)
\ifill f:0
\move(243 37)
\lvec(248 37)
\lvec(248 38)
\lvec(243 38)
\ifill f:0
\move(249 37)
\lvec(253 37)
\lvec(253 38)
\lvec(249 38)
\ifill f:0
\move(254 37)
\lvec(257 37)
\lvec(257 38)
\lvec(254 38)
\ifill f:0
\move(258 37)
\lvec(259 37)
\lvec(259 38)
\lvec(258 38)
\ifill f:0
\move(260 37)
\lvec(280 37)
\lvec(280 38)
\lvec(260 38)
\ifill f:0
\move(281 37)
\lvec(298 37)
\lvec(298 38)
\lvec(281 38)
\ifill f:0
\move(301 37)
\lvec(313 37)
\lvec(313 38)
\lvec(301 38)
\ifill f:0
\move(315 37)
\lvec(325 37)
\lvec(325 38)
\lvec(315 38)
\ifill f:0
\move(326 37)
\lvec(341 37)
\lvec(341 38)
\lvec(326 38)
\ifill f:0
\move(351 37)
\lvec(353 37)
\lvec(353 38)
\lvec(351 38)
\ifill f:0
\move(361 37)
\lvec(362 37)
\lvec(362 38)
\lvec(361 38)
\ifill f:0
\move(366 37)
\lvec(370 37)
\lvec(370 38)
\lvec(366 38)
\ifill f:0
\move(371 37)
\lvec(391 37)
\lvec(391 38)
\lvec(371 38)
\ifill f:0
\move(393 37)
\lvec(401 37)
\lvec(401 38)
\lvec(393 38)
\ifill f:0
\move(402 37)
\lvec(407 37)
\lvec(407 38)
\lvec(402 38)
\ifill f:0
\move(408 37)
\lvec(420 37)
\lvec(420 38)
\lvec(408 38)
\ifill f:0
\move(421 37)
\lvec(431 37)
\lvec(431 38)
\lvec(421 38)
\ifill f:0
\move(432 37)
\lvec(442 37)
\lvec(442 38)
\lvec(432 38)
\ifill f:0
\move(443 37)
\lvec(450 37)
\lvec(450 38)
\lvec(443 38)
\ifill f:0
\move(11 38)
\lvec(12 38)
\lvec(12 39)
\lvec(11 39)
\ifill f:0
\move(16 38)
\lvec(17 38)
\lvec(17 39)
\lvec(16 39)
\ifill f:0
\move(19 38)
\lvec(21 38)
\lvec(21 39)
\lvec(19 39)
\ifill f:0
\move(22 38)
\lvec(26 38)
\lvec(26 39)
\lvec(22 39)
\ifill f:0
\move(28 38)
\lvec(29 38)
\lvec(29 39)
\lvec(28 39)
\ifill f:0
\move(36 38)
\lvec(37 38)
\lvec(37 39)
\lvec(36 39)
\ifill f:0
\move(38 38)
\lvec(39 38)
\lvec(39 39)
\lvec(38 39)
\ifill f:0
\move(41 38)
\lvec(42 38)
\lvec(42 39)
\lvec(41 39)
\ifill f:0
\move(44 38)
\lvec(46 38)
\lvec(46 39)
\lvec(44 39)
\ifill f:0
\move(47 38)
\lvec(50 38)
\lvec(50 39)
\lvec(47 39)
\ifill f:0
\move(51 38)
\lvec(53 38)
\lvec(53 39)
\lvec(51 39)
\ifill f:0
\move(54 38)
\lvec(65 38)
\lvec(65 39)
\lvec(54 39)
\ifill f:0
\move(66 38)
\lvec(68 38)
\lvec(68 39)
\lvec(66 39)
\ifill f:0
\move(69 38)
\lvec(71 38)
\lvec(71 39)
\lvec(69 39)
\ifill f:0
\move(72 38)
\lvec(74 38)
\lvec(74 39)
\lvec(72 39)
\ifill f:0
\move(76 38)
\lvec(80 38)
\lvec(80 39)
\lvec(76 39)
\ifill f:0
\move(81 38)
\lvec(82 38)
\lvec(82 39)
\lvec(81 39)
\ifill f:0
\move(83 38)
\lvec(85 38)
\lvec(85 39)
\lvec(83 39)
\ifill f:0
\move(87 38)
\lvec(90 38)
\lvec(90 39)
\lvec(87 39)
\ifill f:0
\move(92 38)
\lvec(94 38)
\lvec(94 39)
\lvec(92 39)
\ifill f:0
\move(97 38)
\lvec(98 38)
\lvec(98 39)
\lvec(97 39)
\ifill f:0
\move(99 38)
\lvec(101 38)
\lvec(101 39)
\lvec(99 39)
\ifill f:0
\move(102 38)
\lvec(105 38)
\lvec(105 39)
\lvec(102 39)
\ifill f:0
\move(107 38)
\lvec(122 38)
\lvec(122 39)
\lvec(107 39)
\ifill f:0
\move(123 38)
\lvec(125 38)
\lvec(125 39)
\lvec(123 39)
\ifill f:0
\move(126 38)
\lvec(132 38)
\lvec(132 39)
\lvec(126 39)
\ifill f:0
\move(133 38)
\lvec(135 38)
\lvec(135 39)
\lvec(133 39)
\ifill f:0
\move(136 38)
\lvec(138 38)
\lvec(138 39)
\lvec(136 39)
\ifill f:0
\move(139 38)
\lvec(143 38)
\lvec(143 39)
\lvec(139 39)
\ifill f:0
\move(144 38)
\lvec(145 38)
\lvec(145 39)
\lvec(144 39)
\ifill f:0
\move(146 38)
\lvec(150 38)
\lvec(150 39)
\lvec(146 39)
\ifill f:0
\move(151 38)
\lvec(159 38)
\lvec(159 39)
\lvec(151 39)
\ifill f:0
\move(160 38)
\lvec(161 38)
\lvec(161 39)
\lvec(160 39)
\ifill f:0
\move(162 38)
\lvec(163 38)
\lvec(163 39)
\lvec(162 39)
\ifill f:0
\move(164 38)
\lvec(167 38)
\lvec(167 39)
\lvec(164 39)
\ifill f:0
\move(168 38)
\lvec(170 38)
\lvec(170 39)
\lvec(168 39)
\ifill f:0
\move(172 38)
\lvec(173 38)
\lvec(173 39)
\lvec(172 39)
\ifill f:0
\move(174 38)
\lvec(180 38)
\lvec(180 39)
\lvec(174 39)
\ifill f:0
\move(181 38)
\lvec(187 38)
\lvec(187 39)
\lvec(181 39)
\ifill f:0
\move(188 38)
\lvec(192 38)
\lvec(192 39)
\lvec(188 39)
\ifill f:0
\move(193 38)
\lvec(194 38)
\lvec(194 39)
\lvec(193 39)
\ifill f:0
\move(195 38)
\lvec(197 38)
\lvec(197 39)
\lvec(195 39)
\ifill f:0
\move(198 38)
\lvec(200 38)
\lvec(200 39)
\lvec(198 39)
\ifill f:0
\move(201 38)
\lvec(203 38)
\lvec(203 39)
\lvec(201 39)
\ifill f:0
\move(204 38)
\lvec(215 38)
\lvec(215 39)
\lvec(204 39)
\ifill f:0
\move(216 38)
\lvec(226 38)
\lvec(226 39)
\lvec(216 39)
\ifill f:0
\move(227 38)
\lvec(228 38)
\lvec(228 39)
\lvec(227 39)
\ifill f:0
\move(229 38)
\lvec(232 38)
\lvec(232 39)
\lvec(229 39)
\ifill f:0
\move(233 38)
\lvec(240 38)
\lvec(240 39)
\lvec(233 39)
\ifill f:0
\move(241 38)
\lvec(249 38)
\lvec(249 39)
\lvec(241 39)
\ifill f:0
\move(250 38)
\lvec(257 38)
\lvec(257 39)
\lvec(250 39)
\ifill f:0
\move(258 38)
\lvec(264 38)
\lvec(264 39)
\lvec(258 39)
\ifill f:0
\move(265 38)
\lvec(269 38)
\lvec(269 39)
\lvec(265 39)
\ifill f:0
\move(270 38)
\lvec(275 38)
\lvec(275 39)
\lvec(270 39)
\ifill f:0
\move(276 38)
\lvec(282 38)
\lvec(282 39)
\lvec(276 39)
\ifill f:0
\move(283 38)
\lvec(290 38)
\lvec(290 39)
\lvec(283 39)
\ifill f:0
\move(291 38)
\lvec(296 38)
\lvec(296 39)
\lvec(291 39)
\ifill f:0
\move(298 38)
\lvec(305 38)
\lvec(305 39)
\lvec(298 39)
\ifill f:0
\move(307 38)
\lvec(317 38)
\lvec(317 39)
\lvec(307 39)
\ifill f:0
\move(319 38)
\lvec(325 38)
\lvec(325 39)
\lvec(319 39)
\ifill f:0
\move(326 38)
\lvec(332 38)
\lvec(332 39)
\lvec(326 39)
\ifill f:0
\move(333 38)
\lvec(359 38)
\lvec(359 39)
\lvec(333 39)
\ifill f:0
\move(361 38)
\lvec(362 38)
\lvec(362 39)
\lvec(361 39)
\ifill f:0
\move(370 38)
\lvec(372 38)
\lvec(372 39)
\lvec(370 39)
\ifill f:0
\move(384 38)
\lvec(387 38)
\lvec(387 39)
\lvec(384 39)
\ifill f:0
\move(388 38)
\lvec(401 38)
\lvec(401 39)
\lvec(388 39)
\ifill f:0
\move(402 38)
\lvec(411 38)
\lvec(411 39)
\lvec(402 39)
\ifill f:0
\move(412 38)
\lvec(427 38)
\lvec(427 39)
\lvec(412 39)
\ifill f:0
\move(428 38)
\lvec(442 38)
\lvec(442 39)
\lvec(428 39)
\ifill f:0
\move(443 38)
\lvec(451 38)
\lvec(451 39)
\lvec(443 39)
\ifill f:0
\move(14 39)
\lvec(17 39)
\lvec(17 40)
\lvec(14 40)
\ifill f:0
\move(18 39)
\lvec(19 39)
\lvec(19 40)
\lvec(18 40)
\ifill f:0
\move(20 39)
\lvec(22 39)
\lvec(22 40)
\lvec(20 40)
\ifill f:0
\move(25 39)
\lvec(26 39)
\lvec(26 40)
\lvec(25 40)
\ifill f:0
\move(36 39)
\lvec(37 39)
\lvec(37 40)
\lvec(36 40)
\ifill f:0
\move(39 39)
\lvec(50 39)
\lvec(50 40)
\lvec(39 40)
\ifill f:0
\move(52 39)
\lvec(54 39)
\lvec(54 40)
\lvec(52 40)
\ifill f:0
\move(56 39)
\lvec(57 39)
\lvec(57 40)
\lvec(56 40)
\ifill f:0
\move(58 39)
\lvec(59 39)
\lvec(59 40)
\lvec(58 40)
\ifill f:0
\move(60 39)
\lvec(61 39)
\lvec(61 40)
\lvec(60 40)
\ifill f:0
\move(62 39)
\lvec(63 39)
\lvec(63 40)
\lvec(62 40)
\ifill f:0
\move(64 39)
\lvec(65 39)
\lvec(65 40)
\lvec(64 40)
\ifill f:0
\move(66 39)
\lvec(70 39)
\lvec(70 40)
\lvec(66 40)
\ifill f:0
\move(71 39)
\lvec(73 39)
\lvec(73 40)
\lvec(71 40)
\ifill f:0
\move(74 39)
\lvec(80 39)
\lvec(80 40)
\lvec(74 40)
\ifill f:0
\move(81 39)
\lvec(82 39)
\lvec(82 40)
\lvec(81 40)
\ifill f:0
\move(83 39)
\lvec(85 39)
\lvec(85 40)
\lvec(83 40)
\ifill f:0
\move(88 39)
\lvec(94 39)
\lvec(94 40)
\lvec(88 40)
\ifill f:0
\move(95 39)
\lvec(98 39)
\lvec(98 40)
\lvec(95 40)
\ifill f:0
\move(99 39)
\lvec(101 39)
\lvec(101 40)
\lvec(99 40)
\ifill f:0
\move(102 39)
\lvec(109 39)
\lvec(109 40)
\lvec(102 40)
\ifill f:0
\move(110 39)
\lvec(122 39)
\lvec(122 40)
\lvec(110 40)
\ifill f:0
\move(123 39)
\lvec(126 39)
\lvec(126 40)
\lvec(123 40)
\ifill f:0
\move(127 39)
\lvec(130 39)
\lvec(130 40)
\lvec(127 40)
\ifill f:0
\move(131 39)
\lvec(134 39)
\lvec(134 40)
\lvec(131 40)
\ifill f:0
\move(135 39)
\lvec(143 39)
\lvec(143 40)
\lvec(135 40)
\ifill f:0
\move(144 39)
\lvec(145 39)
\lvec(145 40)
\lvec(144 40)
\ifill f:0
\move(146 39)
\lvec(158 39)
\lvec(158 40)
\lvec(146 40)
\ifill f:0
\move(159 39)
\lvec(163 39)
\lvec(163 40)
\lvec(159 40)
\ifill f:0
\move(164 39)
\lvec(167 39)
\lvec(167 40)
\lvec(164 40)
\ifill f:0
\move(168 39)
\lvec(170 39)
\lvec(170 40)
\lvec(168 40)
\ifill f:0
\move(172 39)
\lvec(175 39)
\lvec(175 40)
\lvec(172 40)
\ifill f:0
\move(176 39)
\lvec(177 39)
\lvec(177 40)
\lvec(176 40)
\ifill f:0
\move(178 39)
\lvec(179 39)
\lvec(179 40)
\lvec(178 40)
\ifill f:0
\move(180 39)
\lvec(190 39)
\lvec(190 40)
\lvec(180 40)
\ifill f:0
\move(191 39)
\lvec(192 39)
\lvec(192 40)
\lvec(191 40)
\ifill f:0
\move(193 39)
\lvec(197 39)
\lvec(197 40)
\lvec(193 40)
\ifill f:0
\move(198 39)
\lvec(216 39)
\lvec(216 40)
\lvec(198 40)
\ifill f:0
\move(217 39)
\lvec(222 39)
\lvec(222 40)
\lvec(217 40)
\ifill f:0
\move(223 39)
\lvec(226 39)
\lvec(226 40)
\lvec(223 40)
\ifill f:0
\move(227 39)
\lvec(228 39)
\lvec(228 40)
\lvec(227 40)
\ifill f:0
\move(229 39)
\lvec(235 39)
\lvec(235 40)
\lvec(229 40)
\ifill f:0
\move(236 39)
\lvec(242 39)
\lvec(242 40)
\lvec(236 40)
\ifill f:0
\move(243 39)
\lvec(246 39)
\lvec(246 40)
\lvec(243 40)
\ifill f:0
\move(247 39)
\lvec(254 39)
\lvec(254 40)
\lvec(247 40)
\ifill f:0
\move(255 39)
\lvec(257 39)
\lvec(257 40)
\lvec(255 40)
\ifill f:0
\move(258 39)
\lvec(283 39)
\lvec(283 40)
\lvec(258 40)
\ifill f:0
\move(284 39)
\lvec(290 39)
\lvec(290 40)
\lvec(284 40)
\ifill f:0
\move(291 39)
\lvec(295 39)
\lvec(295 40)
\lvec(291 40)
\ifill f:0
\move(296 39)
\lvec(302 39)
\lvec(302 40)
\lvec(296 40)
\ifill f:0
\move(303 39)
\lvec(310 39)
\lvec(310 40)
\lvec(303 40)
\ifill f:0
\move(311 39)
\lvec(319 39)
\lvec(319 40)
\lvec(311 40)
\ifill f:0
\move(320 39)
\lvec(325 39)
\lvec(325 40)
\lvec(320 40)
\ifill f:0
\move(326 39)
\lvec(329 39)
\lvec(329 40)
\lvec(326 40)
\ifill f:0
\move(331 39)
\lvec(342 39)
\lvec(342 40)
\lvec(331 40)
\ifill f:0
\move(344 39)
\lvec(359 39)
\lvec(359 40)
\lvec(344 40)
\ifill f:0
\move(361 39)
\lvec(362 39)
\lvec(362 40)
\lvec(361 40)
\ifill f:0
\move(363 39)
\lvec(401 39)
\lvec(401 40)
\lvec(363 40)
\ifill f:0
\move(402 39)
\lvec(419 39)
\lvec(419 40)
\lvec(402 40)
\ifill f:0
\move(422 39)
\lvec(440 39)
\lvec(440 40)
\lvec(422 40)
\ifill f:0
\move(441 39)
\lvec(442 39)
\lvec(442 40)
\lvec(441 40)
\ifill f:0
\move(444 39)
\lvec(451 39)
\lvec(451 40)
\lvec(444 40)
\ifill f:0
\move(14 40)
\lvec(17 40)
\lvec(17 41)
\lvec(14 41)
\ifill f:0
\move(18 40)
\lvec(22 40)
\lvec(22 41)
\lvec(18 41)
\ifill f:0
\move(25 40)
\lvec(26 40)
\lvec(26 41)
\lvec(25 41)
\ifill f:0
\move(36 40)
\lvec(37 40)
\lvec(37 41)
\lvec(36 41)
\ifill f:0
\move(38 40)
\lvec(41 40)
\lvec(41 41)
\lvec(38 41)
\ifill f:0
\move(43 40)
\lvec(50 40)
\lvec(50 41)
\lvec(43 41)
\ifill f:0
\move(54 40)
\lvec(55 40)
\lvec(55 41)
\lvec(54 41)
\ifill f:0
\move(56 40)
\lvec(58 40)
\lvec(58 41)
\lvec(56 41)
\ifill f:0
\move(59 40)
\lvec(63 40)
\lvec(63 41)
\lvec(59 41)
\ifill f:0
\move(64 40)
\lvec(65 40)
\lvec(65 41)
\lvec(64 41)
\ifill f:0
\move(66 40)
\lvec(67 40)
\lvec(67 41)
\lvec(66 41)
\ifill f:0
\move(68 40)
\lvec(71 40)
\lvec(71 41)
\lvec(68 41)
\ifill f:0
\move(72 40)
\lvec(74 40)
\lvec(74 41)
\lvec(72 41)
\ifill f:0
\move(75 40)
\lvec(82 40)
\lvec(82 41)
\lvec(75 41)
\ifill f:0
\move(83 40)
\lvec(84 40)
\lvec(84 41)
\lvec(83 41)
\ifill f:0
\move(86 40)
\lvec(90 40)
\lvec(90 41)
\lvec(86 41)
\ifill f:0
\move(94 40)
\lvec(98 40)
\lvec(98 41)
\lvec(94 41)
\ifill f:0
\move(115 40)
\lvec(122 40)
\lvec(122 41)
\lvec(115 41)
\ifill f:0
\move(123 40)
\lvec(127 40)
\lvec(127 41)
\lvec(123 41)
\ifill f:0
\move(128 40)
\lvec(132 40)
\lvec(132 41)
\lvec(128 41)
\ifill f:0
\move(133 40)
\lvec(139 40)
\lvec(139 41)
\lvec(133 41)
\ifill f:0
\move(140 40)
\lvec(145 40)
\lvec(145 41)
\lvec(140 41)
\ifill f:0
\move(146 40)
\lvec(149 40)
\lvec(149 41)
\lvec(146 41)
\ifill f:0
\move(150 40)
\lvec(157 40)
\lvec(157 41)
\lvec(150 41)
\ifill f:0
\move(158 40)
\lvec(162 40)
\lvec(162 41)
\lvec(158 41)
\ifill f:0
\move(163 40)
\lvec(167 40)
\lvec(167 41)
\lvec(163 41)
\ifill f:0
\move(168 40)
\lvec(170 40)
\lvec(170 41)
\lvec(168 41)
\ifill f:0
\move(171 40)
\lvec(176 40)
\lvec(176 41)
\lvec(171 41)
\ifill f:0
\move(177 40)
\lvec(180 40)
\lvec(180 41)
\lvec(177 41)
\ifill f:0
\move(181 40)
\lvec(197 40)
\lvec(197 41)
\lvec(181 41)
\ifill f:0
\move(198 40)
\lvec(204 40)
\lvec(204 41)
\lvec(198 41)
\ifill f:0
\move(205 40)
\lvec(214 40)
\lvec(214 41)
\lvec(205 41)
\ifill f:0
\move(215 40)
\lvec(226 40)
\lvec(226 41)
\lvec(215 41)
\ifill f:0
\move(227 40)
\lvec(228 40)
\lvec(228 41)
\lvec(227 41)
\ifill f:0
\move(229 40)
\lvec(234 40)
\lvec(234 41)
\lvec(229 41)
\ifill f:0
\move(235 40)
\lvec(237 40)
\lvec(237 41)
\lvec(235 41)
\ifill f:0
\move(238 40)
\lvec(240 40)
\lvec(240 41)
\lvec(238 41)
\ifill f:0
\move(241 40)
\lvec(247 40)
\lvec(247 41)
\lvec(241 41)
\ifill f:0
\move(248 40)
\lvec(250 40)
\lvec(250 41)
\lvec(248 41)
\ifill f:0
\move(251 40)
\lvec(254 40)
\lvec(254 41)
\lvec(251 41)
\ifill f:0
\move(255 40)
\lvec(257 40)
\lvec(257 41)
\lvec(255 41)
\ifill f:0
\move(259 40)
\lvec(266 40)
\lvec(266 41)
\lvec(259 41)
\ifill f:0
\move(267 40)
\lvec(279 40)
\lvec(279 41)
\lvec(267 41)
\ifill f:0
\move(280 40)
\lvec(284 40)
\lvec(284 41)
\lvec(280 41)
\ifill f:0
\move(285 40)
\lvec(290 40)
\lvec(290 41)
\lvec(285 41)
\ifill f:0
\move(291 40)
\lvec(293 40)
\lvec(293 41)
\lvec(291 41)
\ifill f:0
\move(295 40)
\lvec(299 40)
\lvec(299 41)
\lvec(295 41)
\ifill f:0
\move(300 40)
\lvec(306 40)
\lvec(306 41)
\lvec(300 41)
\ifill f:0
\move(307 40)
\lvec(313 40)
\lvec(313 41)
\lvec(307 41)
\ifill f:0
\move(314 40)
\lvec(319 40)
\lvec(319 41)
\lvec(314 41)
\ifill f:0
\move(321 40)
\lvec(325 40)
\lvec(325 41)
\lvec(321 41)
\ifill f:0
\move(326 40)
\lvec(328 40)
\lvec(328 41)
\lvec(326 41)
\ifill f:0
\move(329 40)
\lvec(337 40)
\lvec(337 41)
\lvec(329 41)
\ifill f:0
\move(338 40)
\lvec(348 40)
\lvec(348 41)
\lvec(338 41)
\ifill f:0
\move(349 40)
\lvec(360 40)
\lvec(360 41)
\lvec(349 41)
\ifill f:0
\move(361 40)
\lvec(379 40)
\lvec(379 41)
\lvec(361 41)
\ifill f:0
\move(382 40)
\lvec(401 40)
\lvec(401 41)
\lvec(382 41)
\ifill f:0
\move(402 40)
\lvec(440 40)
\lvec(440 41)
\lvec(402 41)
\ifill f:0
\move(441 40)
\lvec(442 40)
\lvec(442 41)
\lvec(441 41)
\ifill f:0
\move(446 40)
\lvec(451 40)
\lvec(451 41)
\lvec(446 41)
\ifill f:0
\move(11 41)
\lvec(12 41)
\lvec(12 42)
\lvec(11 42)
\ifill f:0
\move(15 41)
\lvec(17 41)
\lvec(17 42)
\lvec(15 42)
\ifill f:0
\move(18 41)
\lvec(19 41)
\lvec(19 42)
\lvec(18 42)
\ifill f:0
\move(20 41)
\lvec(24 41)
\lvec(24 42)
\lvec(20 42)
\ifill f:0
\move(25 41)
\lvec(26 41)
\lvec(26 42)
\lvec(25 42)
\ifill f:0
\move(27 41)
\lvec(29 41)
\lvec(29 42)
\lvec(27 42)
\ifill f:0
\move(36 41)
\lvec(37 41)
\lvec(37 42)
\lvec(36 42)
\ifill f:0
\move(38 41)
\lvec(39 41)
\lvec(39 42)
\lvec(38 42)
\ifill f:0
\move(41 41)
\lvec(46 41)
\lvec(46 42)
\lvec(41 42)
\ifill f:0
\move(48 41)
\lvec(50 41)
\lvec(50 42)
\lvec(48 42)
\ifill f:0
\move(56 41)
\lvec(57 41)
\lvec(57 42)
\lvec(56 42)
\ifill f:0
\move(58 41)
\lvec(60 41)
\lvec(60 42)
\lvec(58 42)
\ifill f:0
\move(61 41)
\lvec(63 41)
\lvec(63 42)
\lvec(61 42)
\ifill f:0
\move(64 41)
\lvec(65 41)
\lvec(65 42)
\lvec(64 42)
\ifill f:0
\move(66 41)
\lvec(75 41)
\lvec(75 42)
\lvec(66 42)
\ifill f:0
\move(76 41)
\lvec(78 41)
\lvec(78 42)
\lvec(76 42)
\ifill f:0
\move(79 41)
\lvec(82 41)
\lvec(82 42)
\lvec(79 42)
\ifill f:0
\move(83 41)
\lvec(84 41)
\lvec(84 42)
\lvec(83 42)
\ifill f:0
\move(85 41)
\lvec(88 41)
\lvec(88 42)
\lvec(85 42)
\ifill f:0
\move(90 41)
\lvec(94 41)
\lvec(94 42)
\lvec(90 42)
\ifill f:0
\move(96 41)
\lvec(101 41)
\lvec(101 42)
\lvec(96 42)
\ifill f:0
\move(102 41)
\lvec(111 41)
\lvec(111 42)
\lvec(102 42)
\ifill f:0
\move(112 41)
\lvec(122 41)
\lvec(122 42)
\lvec(112 42)
\ifill f:0
\move(125 41)
\lvec(128 41)
\lvec(128 42)
\lvec(125 42)
\ifill f:0
\move(130 41)
\lvec(134 41)
\lvec(134 42)
\lvec(130 42)
\ifill f:0
\move(135 41)
\lvec(138 41)
\lvec(138 42)
\lvec(135 42)
\ifill f:0
\move(139 41)
\lvec(142 41)
\lvec(142 42)
\lvec(139 42)
\ifill f:0
\move(143 41)
\lvec(145 41)
\lvec(145 42)
\lvec(143 42)
\ifill f:0
\move(146 41)
\lvec(159 41)
\lvec(159 42)
\lvec(146 42)
\ifill f:0
\move(160 41)
\lvec(170 41)
\lvec(170 42)
\lvec(160 42)
\ifill f:0
\move(171 41)
\lvec(172 41)
\lvec(172 42)
\lvec(171 42)
\ifill f:0
\move(173 41)
\lvec(174 41)
\lvec(174 42)
\lvec(173 42)
\ifill f:0
\move(175 41)
\lvec(183 41)
\lvec(183 42)
\lvec(175 42)
\ifill f:0
\move(184 41)
\lvec(185 41)
\lvec(185 42)
\lvec(184 42)
\ifill f:0
\move(186 41)
\lvec(187 41)
\lvec(187 42)
\lvec(186 42)
\ifill f:0
\move(188 41)
\lvec(189 41)
\lvec(189 42)
\lvec(188 42)
\ifill f:0
\move(190 41)
\lvec(193 41)
\lvec(193 42)
\lvec(190 42)
\ifill f:0
\move(194 41)
\lvec(195 41)
\lvec(195 42)
\lvec(194 42)
\ifill f:0
\move(196 41)
\lvec(197 41)
\lvec(197 42)
\lvec(196 42)
\ifill f:0
\move(198 41)
\lvec(208 41)
\lvec(208 42)
\lvec(198 42)
\ifill f:0
\move(209 41)
\lvec(215 41)
\lvec(215 42)
\lvec(209 42)
\ifill f:0
\move(216 41)
\lvec(226 41)
\lvec(226 42)
\lvec(216 42)
\ifill f:0
\move(227 41)
\lvec(233 41)
\lvec(233 42)
\lvec(227 42)
\ifill f:0
\move(234 41)
\lvec(236 41)
\lvec(236 42)
\lvec(234 42)
\ifill f:0
\move(237 41)
\lvec(248 41)
\lvec(248 42)
\lvec(237 42)
\ifill f:0
\move(249 41)
\lvec(251 41)
\lvec(251 42)
\lvec(249 42)
\ifill f:0
\move(252 41)
\lvec(257 41)
\lvec(257 42)
\lvec(252 42)
\ifill f:0
\move(258 41)
\lvec(280 41)
\lvec(280 42)
\lvec(258 42)
\ifill f:0
\move(281 41)
\lvec(290 41)
\lvec(290 42)
\lvec(281 42)
\ifill f:0
\move(291 41)
\lvec(293 41)
\lvec(293 42)
\lvec(291 42)
\ifill f:0
\move(294 41)
\lvec(298 41)
\lvec(298 42)
\lvec(294 42)
\ifill f:0
\move(299 41)
\lvec(309 41)
\lvec(309 42)
\lvec(299 42)
\ifill f:0
\move(310 41)
\lvec(315 41)
\lvec(315 42)
\lvec(310 42)
\ifill f:0
\move(316 41)
\lvec(320 41)
\lvec(320 42)
\lvec(316 42)
\ifill f:0
\move(322 41)
\lvec(325 41)
\lvec(325 42)
\lvec(322 42)
\ifill f:0
\move(326 41)
\lvec(342 41)
\lvec(342 42)
\lvec(326 42)
\ifill f:0
\move(343 41)
\lvec(351 41)
\lvec(351 42)
\lvec(343 42)
\ifill f:0
\move(352 41)
\lvec(373 41)
\lvec(373 42)
\lvec(352 42)
\ifill f:0
\move(374 41)
\lvec(387 41)
\lvec(387 42)
\lvec(374 42)
\ifill f:0
\move(388 41)
\lvec(389 41)
\lvec(389 42)
\lvec(388 42)
\ifill f:0
\move(390 41)
\lvec(401 41)
\lvec(401 42)
\lvec(390 42)
\ifill f:0
\move(402 41)
\lvec(418 41)
\lvec(418 42)
\lvec(402 42)
\ifill f:0
\move(430 41)
\lvec(432 41)
\lvec(432 42)
\lvec(430 42)
\ifill f:0
\move(441 41)
\lvec(442 41)
\lvec(442 42)
\lvec(441 42)
\ifill f:0
\move(446 41)
\lvec(451 41)
\lvec(451 42)
\lvec(446 42)
\ifill f:0
\move(15 42)
\lvec(17 42)
\lvec(17 43)
\lvec(15 43)
\ifill f:0
\move(20 42)
\lvec(21 42)
\lvec(21 43)
\lvec(20 43)
\ifill f:0
\move(23 42)
\lvec(24 42)
\lvec(24 43)
\lvec(23 43)
\ifill f:0
\move(25 42)
\lvec(26 42)
\lvec(26 43)
\lvec(25 43)
\ifill f:0
\move(36 42)
\lvec(37 42)
\lvec(37 43)
\lvec(36 43)
\ifill f:0
\move(38 42)
\lvec(39 42)
\lvec(39 43)
\lvec(38 43)
\ifill f:0
\move(40 42)
\lvec(42 42)
\lvec(42 43)
\lvec(40 43)
\ifill f:0
\move(43 42)
\lvec(46 42)
\lvec(46 43)
\lvec(43 43)
\ifill f:0
\move(49 42)
\lvec(50 42)
\lvec(50 43)
\lvec(49 43)
\ifill f:0
\move(56 42)
\lvec(59 42)
\lvec(59 43)
\lvec(56 43)
\ifill f:0
\move(60 42)
\lvec(63 42)
\lvec(63 43)
\lvec(60 43)
\ifill f:0
\move(64 42)
\lvec(65 42)
\lvec(65 43)
\lvec(64 43)
\ifill f:0
\move(66 42)
\lvec(68 42)
\lvec(68 43)
\lvec(66 43)
\ifill f:0
\move(69 42)
\lvec(76 42)
\lvec(76 43)
\lvec(69 43)
\ifill f:0
\move(77 42)
\lvec(82 42)
\lvec(82 43)
\lvec(77 43)
\ifill f:0
\move(85 42)
\lvec(86 42)
\lvec(86 43)
\lvec(85 43)
\ifill f:0
\move(88 42)
\lvec(90 42)
\lvec(90 43)
\lvec(88 43)
\ifill f:0
\move(92 42)
\lvec(94 42)
\lvec(94 43)
\lvec(92 43)
\ifill f:0
\move(97 42)
\lvec(101 42)
\lvec(101 43)
\lvec(97 43)
\ifill f:0
\move(102 42)
\lvec(105 42)
\lvec(105 43)
\lvec(102 43)
\ifill f:0
\move(106 42)
\lvec(122 42)
\lvec(122 43)
\lvec(106 43)
\ifill f:0
\move(125 42)
\lvec(131 42)
\lvec(131 43)
\lvec(125 43)
\ifill f:0
\move(133 42)
\lvec(137 42)
\lvec(137 43)
\lvec(133 43)
\ifill f:0
\move(139 42)
\lvec(142 42)
\lvec(142 43)
\lvec(139 43)
\ifill f:0
\move(143 42)
\lvec(145 42)
\lvec(145 43)
\lvec(143 43)
\ifill f:0
\move(147 42)
\lvec(150 42)
\lvec(150 43)
\lvec(147 43)
\ifill f:0
\move(151 42)
\lvec(157 42)
\lvec(157 43)
\lvec(151 43)
\ifill f:0
\move(158 42)
\lvec(170 42)
\lvec(170 43)
\lvec(158 43)
\ifill f:0
\move(171 42)
\lvec(172 42)
\lvec(172 43)
\lvec(171 43)
\ifill f:0
\move(173 42)
\lvec(177 42)
\lvec(177 43)
\lvec(173 43)
\ifill f:0
\move(178 42)
\lvec(184 42)
\lvec(184 43)
\lvec(178 43)
\ifill f:0
\move(185 42)
\lvec(191 42)
\lvec(191 43)
\lvec(185 43)
\ifill f:0
\move(192 42)
\lvec(195 42)
\lvec(195 43)
\lvec(192 43)
\ifill f:0
\move(196 42)
\lvec(197 42)
\lvec(197 43)
\lvec(196 43)
\ifill f:0
\move(198 42)
\lvec(199 42)
\lvec(199 43)
\lvec(198 43)
\ifill f:0
\move(200 42)
\lvec(203 42)
\lvec(203 43)
\lvec(200 43)
\ifill f:0
\move(204 42)
\lvec(205 42)
\lvec(205 43)
\lvec(204 43)
\ifill f:0
\move(206 42)
\lvec(207 42)
\lvec(207 43)
\lvec(206 43)
\ifill f:0
\move(208 42)
\lvec(211 42)
\lvec(211 43)
\lvec(208 43)
\ifill f:0
\move(212 42)
\lvec(218 42)
\lvec(218 43)
\lvec(212 43)
\ifill f:0
\move(219 42)
\lvec(226 42)
\lvec(226 43)
\lvec(219 43)
\ifill f:0
\move(228 42)
\lvec(235 42)
\lvec(235 43)
\lvec(228 43)
\ifill f:0
\move(236 42)
\lvec(240 42)
\lvec(240 43)
\lvec(236 43)
\ifill f:0
\move(241 42)
\lvec(243 42)
\lvec(243 43)
\lvec(241 43)
\ifill f:0
\move(244 42)
\lvec(257 42)
\lvec(257 43)
\lvec(244 43)
\ifill f:0
\move(258 42)
\lvec(264 42)
\lvec(264 43)
\lvec(258 43)
\ifill f:0
\move(265 42)
\lvec(267 42)
\lvec(267 43)
\lvec(265 43)
\ifill f:0
\move(268 42)
\lvec(274 42)
\lvec(274 43)
\lvec(268 43)
\ifill f:0
\move(275 42)
\lvec(281 42)
\lvec(281 43)
\lvec(275 43)
\ifill f:0
\move(282 42)
\lvec(290 42)
\lvec(290 43)
\lvec(282 43)
\ifill f:0
\move(291 42)
\lvec(293 42)
\lvec(293 43)
\lvec(291 43)
\ifill f:0
\move(294 42)
\lvec(297 42)
\lvec(297 43)
\lvec(294 43)
\ifill f:0
\move(298 42)
\lvec(306 42)
\lvec(306 43)
\lvec(298 43)
\ifill f:0
\move(307 42)
\lvec(311 42)
\lvec(311 43)
\lvec(307 43)
\ifill f:0
\move(312 42)
\lvec(316 42)
\lvec(316 43)
\lvec(312 43)
\ifill f:0
\move(317 42)
\lvec(321 42)
\lvec(321 43)
\lvec(317 43)
\ifill f:0
\move(322 42)
\lvec(325 42)
\lvec(325 43)
\lvec(322 43)
\ifill f:0
\move(326 42)
\lvec(327 42)
\lvec(327 43)
\lvec(326 43)
\ifill f:0
\move(328 42)
\lvec(333 42)
\lvec(333 43)
\lvec(328 43)
\ifill f:0
\move(334 42)
\lvec(339 42)
\lvec(339 43)
\lvec(334 43)
\ifill f:0
\move(340 42)
\lvec(346 42)
\lvec(346 43)
\lvec(340 43)
\ifill f:0
\move(347 42)
\lvec(353 42)
\lvec(353 43)
\lvec(347 43)
\ifill f:0
\move(354 42)
\lvec(362 42)
\lvec(362 43)
\lvec(354 43)
\ifill f:0
\move(363 42)
\lvec(370 42)
\lvec(370 43)
\lvec(363 43)
\ifill f:0
\move(371 42)
\lvec(379 42)
\lvec(379 43)
\lvec(371 43)
\ifill f:0
\move(381 42)
\lvec(392 42)
\lvec(392 43)
\lvec(381 43)
\ifill f:0
\move(394 42)
\lvec(401 42)
\lvec(401 43)
\lvec(394 43)
\ifill f:0
\move(402 42)
\lvec(409 42)
\lvec(409 43)
\lvec(402 43)
\ifill f:0
\move(410 42)
\lvec(436 42)
\lvec(436 43)
\lvec(410 43)
\ifill f:0
\move(437 42)
\lvec(439 42)
\lvec(439 43)
\lvec(437 43)
\ifill f:0
\move(441 42)
\lvec(442 42)
\lvec(442 43)
\lvec(441 43)
\ifill f:0
\move(14 43)
\lvec(15 43)
\lvec(15 44)
\lvec(14 44)
\ifill f:0
\move(16 43)
\lvec(17 43)
\lvec(17 44)
\lvec(16 44)
\ifill f:0
\move(18 43)
\lvec(21 43)
\lvec(21 44)
\lvec(18 44)
\ifill f:0
\move(23 43)
\lvec(26 43)
\lvec(26 44)
\lvec(23 44)
\ifill f:0
\move(36 43)
\lvec(37 43)
\lvec(37 44)
\lvec(36 44)
\ifill f:0
\move(38 43)
\lvec(43 43)
\lvec(43 44)
\lvec(38 44)
\ifill f:0
\move(44 43)
\lvec(46 43)
\lvec(46 44)
\lvec(44 44)
\ifill f:0
\move(47 43)
\lvec(48 43)
\lvec(48 44)
\lvec(47 44)
\ifill f:0
\move(49 43)
\lvec(50 43)
\lvec(50 44)
\lvec(49 44)
\ifill f:0
\move(51 43)
\lvec(53 43)
\lvec(53 44)
\lvec(51 44)
\ifill f:0
\move(54 43)
\lvec(57 43)
\lvec(57 44)
\lvec(54 44)
\ifill f:0
\move(60 43)
\lvec(62 43)
\lvec(62 44)
\lvec(60 44)
\ifill f:0
\move(64 43)
\lvec(65 43)
\lvec(65 44)
\lvec(64 44)
\ifill f:0
\move(66 43)
\lvec(74 43)
\lvec(74 44)
\lvec(66 44)
\ifill f:0
\move(75 43)
\lvec(82 43)
\lvec(82 44)
\lvec(75 44)
\ifill f:0
\move(84 43)
\lvec(85 43)
\lvec(85 44)
\lvec(84 44)
\ifill f:0
\move(87 43)
\lvec(89 43)
\lvec(89 44)
\lvec(87 44)
\ifill f:0
\move(90 43)
\lvec(93 43)
\lvec(93 44)
\lvec(90 44)
\ifill f:0
\move(94 43)
\lvec(96 43)
\lvec(96 44)
\lvec(94 44)
\ifill f:0
\move(97 43)
\lvec(101 43)
\lvec(101 44)
\lvec(97 44)
\ifill f:0
\move(102 43)
\lvec(103 43)
\lvec(103 44)
\lvec(102 44)
\ifill f:0
\move(104 43)
\lvec(111 43)
\lvec(111 44)
\lvec(104 44)
\ifill f:0
\move(112 43)
\lvec(113 43)
\lvec(113 44)
\lvec(112 44)
\ifill f:0
\move(114 43)
\lvec(115 43)
\lvec(115 44)
\lvec(114 44)
\ifill f:0
\move(118 43)
\lvec(119 43)
\lvec(119 44)
\lvec(118 44)
\ifill f:0
\move(121 43)
\lvec(122 43)
\lvec(122 44)
\lvec(121 44)
\ifill f:0
\move(127 43)
\lvec(135 43)
\lvec(135 44)
\lvec(127 44)
\ifill f:0
\move(136 43)
\lvec(145 43)
\lvec(145 44)
\lvec(136 44)
\ifill f:0
\move(146 43)
\lvec(147 43)
\lvec(147 44)
\lvec(146 44)
\ifill f:0
\move(148 43)
\lvec(159 43)
\lvec(159 44)
\lvec(148 44)
\ifill f:0
\move(160 43)
\lvec(163 43)
\lvec(163 44)
\lvec(160 44)
\ifill f:0
\move(164 43)
\lvec(166 43)
\lvec(166 44)
\lvec(164 44)
\ifill f:0
\move(167 43)
\lvec(170 43)
\lvec(170 44)
\lvec(167 44)
\ifill f:0
\move(171 43)
\lvec(172 43)
\lvec(172 44)
\lvec(171 44)
\ifill f:0
\move(174 43)
\lvec(175 43)
\lvec(175 44)
\lvec(174 44)
\ifill f:0
\move(176 43)
\lvec(183 43)
\lvec(183 44)
\lvec(176 44)
\ifill f:0
\move(184 43)
\lvec(195 43)
\lvec(195 44)
\lvec(184 44)
\ifill f:0
\move(196 43)
\lvec(197 43)
\lvec(197 44)
\lvec(196 44)
\ifill f:0
\move(198 43)
\lvec(226 43)
\lvec(226 44)
\lvec(198 44)
\ifill f:0
\move(228 43)
\lvec(229 43)
\lvec(229 44)
\lvec(228 44)
\ifill f:0
\move(230 43)
\lvec(234 43)
\lvec(234 44)
\lvec(230 44)
\ifill f:0
\move(235 43)
\lvec(239 43)
\lvec(239 44)
\lvec(235 44)
\ifill f:0
\move(240 43)
\lvec(244 43)
\lvec(244 44)
\lvec(240 44)
\ifill f:0
\move(245 43)
\lvec(252 43)
\lvec(252 44)
\lvec(245 44)
\ifill f:0
\move(253 43)
\lvec(257 43)
\lvec(257 44)
\lvec(253 44)
\ifill f:0
\move(258 43)
\lvec(266 43)
\lvec(266 44)
\lvec(258 44)
\ifill f:0
\move(267 43)
\lvec(269 43)
\lvec(269 44)
\lvec(267 44)
\ifill f:0
\move(270 43)
\lvec(282 43)
\lvec(282 44)
\lvec(270 44)
\ifill f:0
\move(283 43)
\lvec(290 43)
\lvec(290 44)
\lvec(283 44)
\ifill f:0
\move(291 43)
\lvec(292 43)
\lvec(292 44)
\lvec(291 44)
\ifill f:0
\move(293 43)
\lvec(296 43)
\lvec(296 44)
\lvec(293 44)
\ifill f:0
\move(297 43)
\lvec(308 43)
\lvec(308 44)
\lvec(297 44)
\ifill f:0
\move(309 43)
\lvec(317 43)
\lvec(317 44)
\lvec(309 44)
\ifill f:0
\move(318 43)
\lvec(325 43)
\lvec(325 44)
\lvec(318 44)
\ifill f:0
\move(326 43)
\lvec(342 43)
\lvec(342 44)
\lvec(326 44)
\ifill f:0
\move(343 43)
\lvec(362 43)
\lvec(362 44)
\lvec(343 44)
\ifill f:0
\move(363 43)
\lvec(368 43)
\lvec(368 44)
\lvec(363 44)
\ifill f:0
\move(369 43)
\lvec(376 43)
\lvec(376 44)
\lvec(369 44)
\ifill f:0
\move(377 43)
\lvec(385 43)
\lvec(385 44)
\lvec(377 44)
\ifill f:0
\move(386 43)
\lvec(394 43)
\lvec(394 44)
\lvec(386 44)
\ifill f:0
\move(396 43)
\lvec(401 43)
\lvec(401 44)
\lvec(396 44)
\ifill f:0
\move(402 43)
\lvec(419 43)
\lvec(419 44)
\lvec(402 44)
\ifill f:0
\move(420 43)
\lvec(421 43)
\lvec(421 44)
\lvec(420 44)
\ifill f:0
\move(422 43)
\lvec(439 43)
\lvec(439 44)
\lvec(422 44)
\ifill f:0
\move(441 43)
\lvec(442 43)
\lvec(442 44)
\lvec(441 44)
\ifill f:0
\move(443 43)
\lvec(451 43)
\lvec(451 44)
\lvec(443 44)
\ifill f:0
\move(11 44)
\lvec(12 44)
\lvec(12 45)
\lvec(11 45)
\ifill f:0
\move(16 44)
\lvec(17 44)
\lvec(17 45)
\lvec(16 45)
\ifill f:0
\move(18 44)
\lvec(19 44)
\lvec(19 45)
\lvec(18 45)
\ifill f:0
\move(20 44)
\lvec(22 44)
\lvec(22 45)
\lvec(20 45)
\ifill f:0
\move(24 44)
\lvec(26 44)
\lvec(26 45)
\lvec(24 45)
\ifill f:0
\move(36 44)
\lvec(37 44)
\lvec(37 45)
\lvec(36 45)
\ifill f:0
\move(38 44)
\lvec(40 44)
\lvec(40 45)
\lvec(38 45)
\ifill f:0
\move(41 44)
\lvec(42 44)
\lvec(42 45)
\lvec(41 45)
\ifill f:0
\move(43 44)
\lvec(45 44)
\lvec(45 45)
\lvec(43 45)
\ifill f:0
\move(46 44)
\lvec(48 44)
\lvec(48 45)
\lvec(46 45)
\ifill f:0
\move(49 44)
\lvec(50 44)
\lvec(50 45)
\lvec(49 45)
\ifill f:0
\move(51 44)
\lvec(53 44)
\lvec(53 45)
\lvec(51 45)
\ifill f:0
\move(56 44)
\lvec(62 44)
\lvec(62 45)
\lvec(56 45)
\ifill f:0
\move(64 44)
\lvec(65 44)
\lvec(65 45)
\lvec(64 45)
\ifill f:0
\move(67 44)
\lvec(75 44)
\lvec(75 45)
\lvec(67 45)
\ifill f:0
\move(76 44)
\lvec(77 44)
\lvec(77 45)
\lvec(76 45)
\ifill f:0
\move(78 44)
\lvec(82 44)
\lvec(82 45)
\lvec(78 45)
\ifill f:0
\move(84 44)
\lvec(85 44)
\lvec(85 45)
\lvec(84 45)
\ifill f:0
\move(86 44)
\lvec(88 44)
\lvec(88 45)
\lvec(86 45)
\ifill f:0
\move(89 44)
\lvec(90 44)
\lvec(90 45)
\lvec(89 45)
\ifill f:0
\move(92 44)
\lvec(94 44)
\lvec(94 45)
\lvec(92 45)
\ifill f:0
\move(95 44)
\lvec(101 44)
\lvec(101 45)
\lvec(95 45)
\ifill f:0
\move(102 44)
\lvec(119 44)
\lvec(119 45)
\lvec(102 45)
\ifill f:0
\move(121 44)
\lvec(122 44)
\lvec(122 45)
\lvec(121 45)
\ifill f:0
\move(124 44)
\lvec(125 44)
\lvec(125 45)
\lvec(124 45)
\ifill f:0
\move(133 44)
\lvec(140 44)
\lvec(140 45)
\lvec(133 45)
\ifill f:0
\move(141 44)
\lvec(145 44)
\lvec(145 45)
\lvec(141 45)
\ifill f:0
\move(146 44)
\lvec(153 44)
\lvec(153 45)
\lvec(146 45)
\ifill f:0
\move(154 44)
\lvec(170 44)
\lvec(170 45)
\lvec(154 45)
\ifill f:0
\move(171 44)
\lvec(173 44)
\lvec(173 45)
\lvec(171 45)
\ifill f:0
\move(174 44)
\lvec(176 44)
\lvec(176 45)
\lvec(174 45)
\ifill f:0
\move(177 44)
\lvec(179 44)
\lvec(179 45)
\lvec(177 45)
\ifill f:0
\move(180 44)
\lvec(187 44)
\lvec(187 45)
\lvec(180 45)
\ifill f:0
\move(188 44)
\lvec(195 44)
\lvec(195 45)
\lvec(188 45)
\ifill f:0
\move(196 44)
\lvec(197 44)
\lvec(197 45)
\lvec(196 45)
\ifill f:0
\move(198 44)
\lvec(211 44)
\lvec(211 45)
\lvec(198 45)
\ifill f:0
\move(212 44)
\lvec(215 44)
\lvec(215 45)
\lvec(212 45)
\ifill f:0
\move(216 44)
\lvec(217 44)
\lvec(217 45)
\lvec(216 45)
\ifill f:0
\move(218 44)
\lvec(219 44)
\lvec(219 45)
\lvec(218 45)
\ifill f:0
\move(220 44)
\lvec(221 44)
\lvec(221 45)
\lvec(220 45)
\ifill f:0
\move(222 44)
\lvec(223 44)
\lvec(223 45)
\lvec(222 45)
\ifill f:0
\move(224 44)
\lvec(226 44)
\lvec(226 45)
\lvec(224 45)
\ifill f:0
\move(228 44)
\lvec(229 44)
\lvec(229 45)
\lvec(228 45)
\ifill f:0
\move(230 44)
\lvec(238 44)
\lvec(238 45)
\lvec(230 45)
\ifill f:0
\move(239 44)
\lvec(240 44)
\lvec(240 45)
\lvec(239 45)
\ifill f:0
\move(241 44)
\lvec(245 44)
\lvec(245 45)
\lvec(241 45)
\ifill f:0
\move(246 44)
\lvec(257 44)
\lvec(257 45)
\lvec(246 45)
\ifill f:0
\move(258 44)
\lvec(260 44)
\lvec(260 45)
\lvec(258 45)
\ifill f:0
\move(261 44)
\lvec(290 44)
\lvec(290 45)
\lvec(261 45)
\ifill f:0
\move(291 44)
\lvec(292 44)
\lvec(292 45)
\lvec(291 45)
\ifill f:0
\move(293 44)
\lvec(298 44)
\lvec(298 45)
\lvec(293 45)
\ifill f:0
\move(300 44)
\lvec(306 44)
\lvec(306 45)
\lvec(300 45)
\ifill f:0
\move(307 44)
\lvec(310 44)
\lvec(310 45)
\lvec(307 45)
\ifill f:0
\move(311 44)
\lvec(314 44)
\lvec(314 45)
\lvec(311 45)
\ifill f:0
\move(315 44)
\lvec(318 44)
\lvec(318 45)
\lvec(315 45)
\ifill f:0
\move(319 44)
\lvec(322 44)
\lvec(322 45)
\lvec(319 45)
\ifill f:0
\move(323 44)
\lvec(325 44)
\lvec(325 45)
\lvec(323 45)
\ifill f:0
\move(326 44)
\lvec(335 44)
\lvec(335 45)
\lvec(326 45)
\ifill f:0
\move(336 44)
\lvec(340 44)
\lvec(340 45)
\lvec(336 45)
\ifill f:0
\move(341 44)
\lvec(345 44)
\lvec(345 45)
\lvec(341 45)
\ifill f:0
\move(346 44)
\lvec(350 44)
\lvec(350 45)
\lvec(346 45)
\ifill f:0
\move(351 44)
\lvec(362 44)
\lvec(362 45)
\lvec(351 45)
\ifill f:0
\move(363 44)
\lvec(367 44)
\lvec(367 45)
\lvec(363 45)
\ifill f:0
\move(368 44)
\lvec(380 44)
\lvec(380 45)
\lvec(368 45)
\ifill f:0
\move(381 44)
\lvec(387 44)
\lvec(387 45)
\lvec(381 45)
\ifill f:0
\move(389 44)
\lvec(395 44)
\lvec(395 45)
\lvec(389 45)
\ifill f:0
\move(397 44)
\lvec(401 44)
\lvec(401 45)
\lvec(397 45)
\ifill f:0
\move(402 44)
\lvec(405 44)
\lvec(405 45)
\lvec(402 45)
\ifill f:0
\move(406 44)
\lvec(415 44)
\lvec(415 45)
\lvec(406 45)
\ifill f:0
\move(416 44)
\lvec(426 44)
\lvec(426 45)
\lvec(416 45)
\ifill f:0
\move(428 44)
\lvec(440 44)
\lvec(440 45)
\lvec(428 45)
\ifill f:0
\move(441 44)
\lvec(451 44)
\lvec(451 45)
\lvec(441 45)
\ifill f:0
\move(16 45)
\lvec(17 45)
\lvec(17 46)
\lvec(16 46)
\ifill f:0
\move(18 45)
\lvec(19 45)
\lvec(19 46)
\lvec(18 46)
\ifill f:0
\move(24 45)
\lvec(26 45)
\lvec(26 46)
\lvec(24 46)
\ifill f:0
\move(36 45)
\lvec(37 45)
\lvec(37 46)
\lvec(36 46)
\ifill f:0
\move(40 45)
\lvec(41 45)
\lvec(41 46)
\lvec(40 46)
\ifill f:0
\move(42 45)
\lvec(43 45)
\lvec(43 46)
\lvec(42 46)
\ifill f:0
\move(44 45)
\lvec(45 45)
\lvec(45 46)
\lvec(44 46)
\ifill f:0
\move(46 45)
\lvec(50 45)
\lvec(50 46)
\lvec(46 46)
\ifill f:0
\move(51 45)
\lvec(52 45)
\lvec(52 46)
\lvec(51 46)
\ifill f:0
\move(62 45)
\lvec(65 45)
\lvec(65 46)
\lvec(62 46)
\ifill f:0
\move(66 45)
\lvec(67 45)
\lvec(67 46)
\lvec(66 46)
\ifill f:0
\move(68 45)
\lvec(71 45)
\lvec(71 46)
\lvec(68 46)
\ifill f:0
\move(72 45)
\lvec(74 45)
\lvec(74 46)
\lvec(72 46)
\ifill f:0
\move(75 45)
\lvec(79 45)
\lvec(79 46)
\lvec(75 46)
\ifill f:0
\move(80 45)
\lvec(82 45)
\lvec(82 46)
\lvec(80 46)
\ifill f:0
\move(83 45)
\lvec(85 45)
\lvec(85 46)
\lvec(83 46)
\ifill f:0
\move(86 45)
\lvec(87 45)
\lvec(87 46)
\lvec(86 46)
\ifill f:0
\move(88 45)
\lvec(92 45)
\lvec(92 46)
\lvec(88 46)
\ifill f:0
\move(93 45)
\lvec(94 45)
\lvec(94 46)
\lvec(93 46)
\ifill f:0
\move(96 45)
\lvec(98 45)
\lvec(98 46)
\lvec(96 46)
\ifill f:0
\move(99 45)
\lvec(101 45)
\lvec(101 46)
\lvec(99 46)
\ifill f:0
\move(102 45)
\lvec(106 45)
\lvec(106 46)
\lvec(102 46)
\ifill f:0
\move(107 45)
\lvec(111 45)
\lvec(111 46)
\lvec(107 46)
\ifill f:0
\move(113 45)
\lvec(120 45)
\lvec(120 46)
\lvec(113 46)
\ifill f:0
\move(121 45)
\lvec(122 45)
\lvec(122 46)
\lvec(121 46)
\ifill f:0
\move(123 45)
\lvec(138 45)
\lvec(138 46)
\lvec(123 46)
\ifill f:0
\move(139 45)
\lvec(145 45)
\lvec(145 46)
\lvec(139 46)
\ifill f:0
\move(146 45)
\lvec(148 45)
\lvec(148 46)
\lvec(146 46)
\ifill f:0
\move(150 45)
\lvec(155 45)
\lvec(155 46)
\lvec(150 46)
\ifill f:0
\move(156 45)
\lvec(170 45)
\lvec(170 46)
\lvec(156 46)
\ifill f:0
\move(171 45)
\lvec(173 45)
\lvec(173 46)
\lvec(171 46)
\ifill f:0
\move(174 45)
\lvec(177 45)
\lvec(177 46)
\lvec(174 46)
\ifill f:0
\move(178 45)
\lvec(180 45)
\lvec(180 46)
\lvec(178 46)
\ifill f:0
\move(181 45)
\lvec(183 45)
\lvec(183 46)
\lvec(181 46)
\ifill f:0
\move(184 45)
\lvec(189 45)
\lvec(189 46)
\lvec(184 46)
\ifill f:0
\move(190 45)
\lvec(192 45)
\lvec(192 46)
\lvec(190 46)
\ifill f:0
\move(193 45)
\lvec(195 45)
\lvec(195 46)
\lvec(193 46)
\ifill f:0
\move(196 45)
\lvec(197 45)
\lvec(197 46)
\lvec(196 46)
\ifill f:0
\move(198 45)
\lvec(200 45)
\lvec(200 46)
\lvec(198 46)
\ifill f:0
\move(201 45)
\lvec(212 45)
\lvec(212 46)
\lvec(201 46)
\ifill f:0
\move(213 45)
\lvec(223 45)
\lvec(223 46)
\lvec(213 46)
\ifill f:0
\move(224 45)
\lvec(226 45)
\lvec(226 46)
\lvec(224 46)
\ifill f:0
\move(228 45)
\lvec(229 45)
\lvec(229 46)
\lvec(228 46)
\ifill f:0
\move(230 45)
\lvec(231 45)
\lvec(231 46)
\lvec(230 46)
\ifill f:0
\move(232 45)
\lvec(233 45)
\lvec(233 46)
\lvec(232 46)
\ifill f:0
\move(234 45)
\lvec(235 45)
\lvec(235 46)
\lvec(234 46)
\ifill f:0
\move(236 45)
\lvec(237 45)
\lvec(237 46)
\lvec(236 46)
\ifill f:0
\move(238 45)
\lvec(248 45)
\lvec(248 46)
\lvec(238 46)
\ifill f:0
\move(249 45)
\lvec(250 45)
\lvec(250 46)
\lvec(249 46)
\ifill f:0
\move(251 45)
\lvec(257 45)
\lvec(257 46)
\lvec(251 46)
\ifill f:0
\move(258 45)
\lvec(262 45)
\lvec(262 46)
\lvec(258 46)
\ifill f:0
\move(263 45)
\lvec(267 45)
\lvec(267 46)
\lvec(263 46)
\ifill f:0
\move(268 45)
\lvec(275 45)
\lvec(275 46)
\lvec(268 46)
\ifill f:0
\move(276 45)
\lvec(280 45)
\lvec(280 46)
\lvec(276 46)
\ifill f:0
\move(281 45)
\lvec(283 45)
\lvec(283 46)
\lvec(281 46)
\ifill f:0
\move(284 45)
\lvec(286 45)
\lvec(286 46)
\lvec(284 46)
\ifill f:0
\move(287 45)
\lvec(290 45)
\lvec(290 46)
\lvec(287 46)
\ifill f:0
\move(291 45)
\lvec(292 45)
\lvec(292 46)
\lvec(291 46)
\ifill f:0
\move(293 45)
\lvec(295 45)
\lvec(295 46)
\lvec(293 46)
\ifill f:0
\move(296 45)
\lvec(298 45)
\lvec(298 46)
\lvec(296 46)
\ifill f:0
\move(299 45)
\lvec(308 45)
\lvec(308 46)
\lvec(299 46)
\ifill f:0
\move(309 45)
\lvec(311 45)
\lvec(311 46)
\lvec(309 46)
\ifill f:0
\move(312 45)
\lvec(315 45)
\lvec(315 46)
\lvec(312 46)
\ifill f:0
\move(316 45)
\lvec(322 45)
\lvec(322 46)
\lvec(316 46)
\ifill f:0
\move(323 45)
\lvec(325 45)
\lvec(325 46)
\lvec(323 46)
\ifill f:0
\move(327 45)
\lvec(330 45)
\lvec(330 46)
\lvec(327 46)
\ifill f:0
\move(331 45)
\lvec(334 45)
\lvec(334 46)
\lvec(331 46)
\ifill f:0
\move(335 45)
\lvec(338 45)
\lvec(338 46)
\lvec(335 46)
\ifill f:0
\move(339 45)
\lvec(342 45)
\lvec(342 46)
\lvec(339 46)
\ifill f:0
\move(343 45)
\lvec(356 45)
\lvec(356 46)
\lvec(343 46)
\ifill f:0
\move(357 45)
\lvec(362 45)
\lvec(362 46)
\lvec(357 46)
\ifill f:0
\move(363 45)
\lvec(366 45)
\lvec(366 46)
\lvec(363 46)
\ifill f:0
\move(367 45)
\lvec(377 45)
\lvec(377 46)
\lvec(367 46)
\ifill f:0
\move(378 45)
\lvec(383 45)
\lvec(383 46)
\lvec(378 46)
\ifill f:0
\move(384 45)
\lvec(390 45)
\lvec(390 46)
\lvec(384 46)
\ifill f:0
\move(391 45)
\lvec(401 45)
\lvec(401 46)
\lvec(391 46)
\ifill f:0
\move(402 45)
\lvec(404 45)
\lvec(404 46)
\lvec(402 46)
\ifill f:0
\move(405 45)
\lvec(412 45)
\lvec(412 46)
\lvec(405 46)
\ifill f:0
\move(413 45)
\lvec(430 45)
\lvec(430 46)
\lvec(413 46)
\ifill f:0
\move(431 45)
\lvec(440 45)
\lvec(440 46)
\lvec(431 46)
\ifill f:0
\move(441 45)
\lvec(451 45)
\lvec(451 46)
\lvec(441 46)
\ifill f:0
\move(16 46)
\lvec(17 46)
\lvec(17 47)
\lvec(16 47)
\ifill f:0
\move(19 46)
\lvec(22 46)
\lvec(22 47)
\lvec(19 47)
\ifill f:0
\move(23 46)
\lvec(26 46)
\lvec(26 47)
\lvec(23 47)
\ifill f:0
\move(27 46)
\lvec(28 46)
\lvec(28 47)
\lvec(27 47)
\ifill f:0
\move(36 46)
\lvec(37 46)
\lvec(37 47)
\lvec(36 47)
\ifill f:0
\move(38 46)
\lvec(39 46)
\lvec(39 47)
\lvec(38 47)
\ifill f:0
\move(40 46)
\lvec(42 46)
\lvec(42 47)
\lvec(40 47)
\ifill f:0
\move(43 46)
\lvec(46 46)
\lvec(46 47)
\lvec(43 47)
\ifill f:0
\move(47 46)
\lvec(50 46)
\lvec(50 47)
\lvec(47 47)
\ifill f:0
\move(53 46)
\lvec(55 46)
\lvec(55 47)
\lvec(53 47)
\ifill f:0
\move(56 46)
\lvec(59 46)
\lvec(59 47)
\lvec(56 47)
\ifill f:0
\move(60 46)
\lvec(65 46)
\lvec(65 47)
\lvec(60 47)
\ifill f:0
\move(66 46)
\lvec(68 46)
\lvec(68 47)
\lvec(66 47)
\ifill f:0
\move(69 46)
\lvec(76 46)
\lvec(76 47)
\lvec(69 47)
\ifill f:0
\move(77 46)
\lvec(79 46)
\lvec(79 47)
\lvec(77 47)
\ifill f:0
\move(80 46)
\lvec(82 46)
\lvec(82 47)
\lvec(80 47)
\ifill f:0
\move(83 46)
\lvec(84 46)
\lvec(84 47)
\lvec(83 47)
\ifill f:0
\move(85 46)
\lvec(86 46)
\lvec(86 47)
\lvec(85 47)
\ifill f:0
\move(87 46)
\lvec(90 46)
\lvec(90 47)
\lvec(87 47)
\ifill f:0
\move(91 46)
\lvec(93 46)
\lvec(93 47)
\lvec(91 47)
\ifill f:0
\move(94 46)
\lvec(95 46)
\lvec(95 47)
\lvec(94 47)
\ifill f:0
\move(96 46)
\lvec(98 46)
\lvec(98 47)
\lvec(96 47)
\ifill f:0
\move(99 46)
\lvec(101 46)
\lvec(101 47)
\lvec(99 47)
\ifill f:0
\move(102 46)
\lvec(105 46)
\lvec(105 47)
\lvec(102 47)
\ifill f:0
\move(106 46)
\lvec(109 46)
\lvec(109 47)
\lvec(106 47)
\ifill f:0
\move(110 46)
\lvec(120 46)
\lvec(120 47)
\lvec(110 47)
\ifill f:0
\move(121 46)
\lvec(122 46)
\lvec(122 47)
\lvec(121 47)
\ifill f:0
\move(123 46)
\lvec(131 46)
\lvec(131 47)
\lvec(123 47)
\ifill f:0
\move(132 46)
\lvec(133 46)
\lvec(133 47)
\lvec(132 47)
\ifill f:0
\move(136 46)
\lvec(138 46)
\lvec(138 47)
\lvec(136 47)
\ifill f:0
\move(139 46)
\lvec(145 46)
\lvec(145 47)
\lvec(139 47)
\ifill f:0
\move(146 46)
\lvec(150 46)
\lvec(150 47)
\lvec(146 47)
\ifill f:0
\move(152 46)
\lvec(158 46)
\lvec(158 47)
\lvec(152 47)
\ifill f:0
\move(160 46)
\lvec(164 46)
\lvec(164 47)
\lvec(160 47)
\ifill f:0
\move(165 46)
\lvec(170 46)
\lvec(170 47)
\lvec(165 47)
\ifill f:0
\move(171 46)
\lvec(174 46)
\lvec(174 47)
\lvec(171 47)
\ifill f:0
\move(175 46)
\lvec(185 46)
\lvec(185 47)
\lvec(175 47)
\ifill f:0
\move(186 46)
\lvec(192 46)
\lvec(192 47)
\lvec(186 47)
\ifill f:0
\move(193 46)
\lvec(195 46)
\lvec(195 47)
\lvec(193 47)
\ifill f:0
\move(196 46)
\lvec(197 46)
\lvec(197 47)
\lvec(196 47)
\ifill f:0
\move(198 46)
\lvec(203 46)
\lvec(203 47)
\lvec(198 47)
\ifill f:0
\move(204 46)
\lvec(211 46)
\lvec(211 47)
\lvec(204 47)
\ifill f:0
\move(212 46)
\lvec(223 46)
\lvec(223 47)
\lvec(212 47)
\ifill f:0
\move(224 46)
\lvec(226 46)
\lvec(226 47)
\lvec(224 47)
\ifill f:0
\move(227 46)
\lvec(240 46)
\lvec(240 47)
\lvec(227 47)
\ifill f:0
\move(241 46)
\lvec(242 46)
\lvec(242 47)
\lvec(241 47)
\ifill f:0
\move(243 46)
\lvec(257 46)
\lvec(257 47)
\lvec(243 47)
\ifill f:0
\move(258 46)
\lvec(264 46)
\lvec(264 47)
\lvec(258 47)
\ifill f:0
\move(265 46)
\lvec(266 46)
\lvec(266 47)
\lvec(265 47)
\ifill f:0
\move(267 46)
\lvec(271 46)
\lvec(271 47)
\lvec(267 47)
\ifill f:0
\move(272 46)
\lvec(276 46)
\lvec(276 47)
\lvec(272 47)
\ifill f:0
\move(277 46)
\lvec(281 46)
\lvec(281 47)
\lvec(277 47)
\ifill f:0
\move(282 46)
\lvec(290 46)
\lvec(290 47)
\lvec(282 47)
\ifill f:0
\move(291 46)
\lvec(297 46)
\lvec(297 47)
\lvec(291 47)
\ifill f:0
\move(298 46)
\lvec(306 46)
\lvec(306 47)
\lvec(298 47)
\ifill f:0
\move(307 46)
\lvec(309 46)
\lvec(309 47)
\lvec(307 47)
\ifill f:0
\move(310 46)
\lvec(319 46)
\lvec(319 47)
\lvec(310 47)
\ifill f:0
\move(320 46)
\lvec(322 46)
\lvec(322 47)
\lvec(320 47)
\ifill f:0
\move(323 46)
\lvec(325 46)
\lvec(325 47)
\lvec(323 47)
\ifill f:0
\move(326 46)
\lvec(329 46)
\lvec(329 47)
\lvec(326 47)
\ifill f:0
\move(330 46)
\lvec(333 46)
\lvec(333 47)
\lvec(330 47)
\ifill f:0
\move(334 46)
\lvec(344 46)
\lvec(344 47)
\lvec(334 47)
\ifill f:0
\move(345 46)
\lvec(362 46)
\lvec(362 47)
\lvec(345 47)
\ifill f:0
\move(363 46)
\lvec(365 46)
\lvec(365 47)
\lvec(363 47)
\ifill f:0
\move(366 46)
\lvec(370 46)
\lvec(370 47)
\lvec(366 47)
\ifill f:0
\move(371 46)
\lvec(375 46)
\lvec(375 47)
\lvec(371 47)
\ifill f:0
\move(376 46)
\lvec(391 46)
\lvec(391 47)
\lvec(376 47)
\ifill f:0
\move(392 46)
\lvec(397 46)
\lvec(397 47)
\lvec(392 47)
\ifill f:0
\move(398 46)
\lvec(401 46)
\lvec(401 47)
\lvec(398 47)
\ifill f:0
\move(402 46)
\lvec(403 46)
\lvec(403 47)
\lvec(402 47)
\ifill f:0
\move(404 46)
\lvec(410 46)
\lvec(410 47)
\lvec(404 47)
\ifill f:0
\move(411 46)
\lvec(424 46)
\lvec(424 47)
\lvec(411 47)
\ifill f:0
\move(425 46)
\lvec(432 46)
\lvec(432 47)
\lvec(425 47)
\ifill f:0
\move(433 46)
\lvec(450 46)
\lvec(450 47)
\lvec(433 47)
\ifill f:0
\move(11 47)
\lvec(12 47)
\lvec(12 48)
\lvec(11 48)
\ifill f:0
\move(14 47)
\lvec(17 47)
\lvec(17 48)
\lvec(14 48)
\ifill f:0
\move(18 47)
\lvec(19 47)
\lvec(19 48)
\lvec(18 48)
\ifill f:0
\move(23 47)
\lvec(26 47)
\lvec(26 48)
\lvec(23 48)
\ifill f:0
\move(36 47)
\lvec(37 47)
\lvec(37 48)
\lvec(36 48)
\ifill f:0
\move(38 47)
\lvec(40 47)
\lvec(40 48)
\lvec(38 48)
\ifill f:0
\move(41 47)
\lvec(43 47)
\lvec(43 48)
\lvec(41 48)
\ifill f:0
\move(44 47)
\lvec(46 47)
\lvec(46 48)
\lvec(44 48)
\ifill f:0
\move(47 47)
\lvec(50 47)
\lvec(50 48)
\lvec(47 48)
\ifill f:0
\move(52 47)
\lvec(53 47)
\lvec(53 48)
\lvec(52 48)
\ifill f:0
\move(57 47)
\lvec(65 47)
\lvec(65 48)
\lvec(57 48)
\ifill f:0
\move(66 47)
\lvec(70 47)
\lvec(70 48)
\lvec(66 48)
\ifill f:0
\move(72 47)
\lvec(75 47)
\lvec(75 48)
\lvec(72 48)
\ifill f:0
\move(76 47)
\lvec(78 47)
\lvec(78 48)
\lvec(76 48)
\ifill f:0
\move(79 47)
\lvec(82 47)
\lvec(82 48)
\lvec(79 48)
\ifill f:0
\move(83 47)
\lvec(84 47)
\lvec(84 48)
\lvec(83 48)
\ifill f:0
\move(87 47)
\lvec(91 47)
\lvec(91 48)
\lvec(87 48)
\ifill f:0
\move(92 47)
\lvec(93 47)
\lvec(93 48)
\lvec(92 48)
\ifill f:0
\move(94 47)
\lvec(96 47)
\lvec(96 48)
\lvec(94 48)
\ifill f:0
\move(97 47)
\lvec(98 47)
\lvec(98 48)
\lvec(97 48)
\ifill f:0
\move(99 47)
\lvec(101 47)
\lvec(101 48)
\lvec(99 48)
\ifill f:0
\move(102 47)
\lvec(111 47)
\lvec(111 48)
\lvec(102 48)
\ifill f:0
\move(112 47)
\lvec(122 47)
\lvec(122 48)
\lvec(112 48)
\ifill f:0
\move(123 47)
\lvec(126 47)
\lvec(126 48)
\lvec(123 48)
\ifill f:0
\move(130 47)
\lvec(145 47)
\lvec(145 48)
\lvec(130 48)
\ifill f:0
\move(146 47)
\lvec(154 47)
\lvec(154 48)
\lvec(146 48)
\ifill f:0
\move(156 47)
\lvec(163 47)
\lvec(163 48)
\lvec(156 48)
\ifill f:0
\move(164 47)
\lvec(170 47)
\lvec(170 48)
\lvec(164 48)
\ifill f:0
\move(171 47)
\lvec(175 47)
\lvec(175 48)
\lvec(171 48)
\ifill f:0
\move(176 47)
\lvec(179 47)
\lvec(179 48)
\lvec(176 48)
\ifill f:0
\move(181 47)
\lvec(191 47)
\lvec(191 48)
\lvec(181 48)
\ifill f:0
\move(192 47)
\lvec(195 47)
\lvec(195 48)
\lvec(192 48)
\ifill f:0
\move(196 47)
\lvec(197 47)
\lvec(197 48)
\lvec(196 48)
\ifill f:0
\move(198 47)
\lvec(201 47)
\lvec(201 48)
\lvec(198 48)
\ifill f:0
\move(202 47)
\lvec(215 47)
\lvec(215 48)
\lvec(202 48)
\ifill f:0
\move(216 47)
\lvec(226 47)
\lvec(226 48)
\lvec(216 48)
\ifill f:0
\move(227 47)
\lvec(232 47)
\lvec(232 48)
\lvec(227 48)
\ifill f:0
\move(233 47)
\lvec(241 47)
\lvec(241 48)
\lvec(233 48)
\ifill f:0
\move(242 47)
\lvec(243 47)
\lvec(243 48)
\lvec(242 48)
\ifill f:0
\move(244 47)
\lvec(247 47)
\lvec(247 48)
\lvec(244 48)
\ifill f:0
\move(248 47)
\lvec(249 47)
\lvec(249 48)
\lvec(248 48)
\ifill f:0
\move(250 47)
\lvec(253 47)
\lvec(253 48)
\lvec(250 48)
\ifill f:0
\move(254 47)
\lvec(255 47)
\lvec(255 48)
\lvec(254 48)
\ifill f:0
\move(256 47)
\lvec(257 47)
\lvec(257 48)
\lvec(256 48)
\ifill f:0
\move(258 47)
\lvec(270 47)
\lvec(270 48)
\lvec(258 48)
\ifill f:0
\move(271 47)
\lvec(279 47)
\lvec(279 48)
\lvec(271 48)
\ifill f:0
\move(280 47)
\lvec(284 47)
\lvec(284 48)
\lvec(280 48)
\ifill f:0
\move(285 47)
\lvec(290 47)
\lvec(290 48)
\lvec(285 48)
\ifill f:0
\move(292 47)
\lvec(294 47)
\lvec(294 48)
\lvec(292 48)
\ifill f:0
\move(295 47)
\lvec(299 47)
\lvec(299 48)
\lvec(295 48)
\ifill f:0
\move(300 47)
\lvec(302 47)
\lvec(302 48)
\lvec(300 48)
\ifill f:0
\move(303 47)
\lvec(319 47)
\lvec(319 48)
\lvec(303 48)
\ifill f:0
\move(320 47)
\lvec(322 47)
\lvec(322 48)
\lvec(320 48)
\ifill f:0
\move(323 47)
\lvec(325 47)
\lvec(325 48)
\lvec(323 48)
\ifill f:0
\move(326 47)
\lvec(342 47)
\lvec(342 48)
\lvec(326 48)
\ifill f:0
\move(343 47)
\lvec(362 47)
\lvec(362 48)
\lvec(343 48)
\ifill f:0
\move(363 47)
\lvec(365 47)
\lvec(365 48)
\lvec(363 48)
\ifill f:0
\move(366 47)
\lvec(369 47)
\lvec(369 48)
\lvec(366 48)
\ifill f:0
\move(370 47)
\lvec(378 47)
\lvec(378 48)
\lvec(370 48)
\ifill f:0
\move(379 47)
\lvec(387 47)
\lvec(387 48)
\lvec(379 48)
\ifill f:0
\move(388 47)
\lvec(392 47)
\lvec(392 48)
\lvec(388 48)
\ifill f:0
\move(393 47)
\lvec(401 47)
\lvec(401 48)
\lvec(393 48)
\ifill f:0
\move(402 47)
\lvec(414 47)
\lvec(414 48)
\lvec(402 48)
\ifill f:0
\move(415 47)
\lvec(427 47)
\lvec(427 48)
\lvec(415 48)
\ifill f:0
\move(428 47)
\lvec(434 47)
\lvec(434 48)
\lvec(428 48)
\ifill f:0
\move(435 47)
\lvec(449 47)
\lvec(449 48)
\lvec(435 48)
\ifill f:0
\move(450 47)
\lvec(451 47)
\lvec(451 48)
\lvec(450 48)
\ifill f:0
\move(15 48)
\lvec(17 48)
\lvec(17 49)
\lvec(15 49)
\ifill f:0
\move(18 48)
\lvec(19 48)
\lvec(19 49)
\lvec(18 49)
\ifill f:0
\move(20 48)
\lvec(21 48)
\lvec(21 49)
\lvec(20 49)
\ifill f:0
\move(23 48)
\lvec(26 48)
\lvec(26 49)
\lvec(23 49)
\ifill f:0
\move(28 48)
\lvec(29 48)
\lvec(29 49)
\lvec(28 49)
\ifill f:0
\move(36 48)
\lvec(37 48)
\lvec(37 49)
\lvec(36 49)
\ifill f:0
\move(38 48)
\lvec(42 48)
\lvec(42 49)
\lvec(38 49)
\ifill f:0
\move(43 48)
\lvec(46 48)
\lvec(46 49)
\lvec(43 49)
\ifill f:0
\move(48 48)
\lvec(50 48)
\lvec(50 49)
\lvec(48 49)
\ifill f:0
\move(52 48)
\lvec(53 48)
\lvec(53 49)
\lvec(52 49)
\ifill f:0
\move(56 48)
\lvec(57 48)
\lvec(57 49)
\lvec(56 49)
\ifill f:0
\move(59 48)
\lvec(65 48)
\lvec(65 49)
\lvec(59 49)
\ifill f:0
\move(66 48)
\lvec(73 48)
\lvec(73 49)
\lvec(66 49)
\ifill f:0
\move(75 48)
\lvec(78 48)
\lvec(78 49)
\lvec(75 49)
\ifill f:0
\move(79 48)
\lvec(82 48)
\lvec(82 49)
\lvec(79 49)
\ifill f:0
\move(83 48)
\lvec(84 48)
\lvec(84 49)
\lvec(83 49)
\ifill f:0
\move(86 48)
\lvec(87 48)
\lvec(87 49)
\lvec(86 49)
\ifill f:0
\move(88 48)
\lvec(90 48)
\lvec(90 49)
\lvec(88 49)
\ifill f:0
\move(91 48)
\lvec(93 48)
\lvec(93 49)
\lvec(91 49)
\ifill f:0
\move(95 48)
\lvec(96 48)
\lvec(96 49)
\lvec(95 49)
\ifill f:0
\move(97 48)
\lvec(98 48)
\lvec(98 49)
\lvec(97 49)
\ifill f:0
\move(99 48)
\lvec(101 48)
\lvec(101 49)
\lvec(99 49)
\ifill f:0
\move(102 48)
\lvec(106 48)
\lvec(106 49)
\lvec(102 49)
\ifill f:0
\move(107 48)
\lvec(109 48)
\lvec(109 49)
\lvec(107 49)
\ifill f:0
\move(110 48)
\lvec(122 48)
\lvec(122 49)
\lvec(110 49)
\ifill f:0
\move(123 48)
\lvec(125 48)
\lvec(125 49)
\lvec(123 49)
\ifill f:0
\move(127 48)
\lvec(133 48)
\lvec(133 49)
\lvec(127 49)
\ifill f:0
\move(136 48)
\lvec(145 48)
\lvec(145 49)
\lvec(136 49)
\ifill f:0
\move(146 48)
\lvec(160 48)
\lvec(160 49)
\lvec(146 49)
\ifill f:0
\move(162 48)
\lvec(170 48)
\lvec(170 49)
\lvec(162 49)
\ifill f:0
\move(171 48)
\lvec(176 48)
\lvec(176 49)
\lvec(171 49)
\ifill f:0
\move(177 48)
\lvec(181 48)
\lvec(181 49)
\lvec(177 49)
\ifill f:0
\move(182 48)
\lvec(186 48)
\lvec(186 49)
\lvec(182 49)
\ifill f:0
\move(187 48)
\lvec(190 48)
\lvec(190 49)
\lvec(187 49)
\ifill f:0
\move(192 48)
\lvec(194 48)
\lvec(194 49)
\lvec(192 49)
\ifill f:0
\move(196 48)
\lvec(197 48)
\lvec(197 49)
\lvec(196 49)
\ifill f:0
\move(198 48)
\lvec(208 48)
\lvec(208 49)
\lvec(198 49)
\ifill f:0
\move(209 48)
\lvec(211 48)
\lvec(211 49)
\lvec(209 49)
\ifill f:0
\move(212 48)
\lvec(214 48)
\lvec(214 49)
\lvec(212 49)
\ifill f:0
\move(215 48)
\lvec(226 48)
\lvec(226 49)
\lvec(215 49)
\ifill f:0
\move(227 48)
\lvec(228 48)
\lvec(228 49)
\lvec(227 49)
\ifill f:0
\move(229 48)
\lvec(235 48)
\lvec(235 49)
\lvec(229 49)
\ifill f:0
\move(236 48)
\lvec(242 48)
\lvec(242 49)
\lvec(236 49)
\ifill f:0
\move(243 48)
\lvec(253 48)
\lvec(253 49)
\lvec(243 49)
\ifill f:0
\move(254 48)
\lvec(255 48)
\lvec(255 49)
\lvec(254 49)
\ifill f:0
\move(256 48)
\lvec(257 48)
\lvec(257 49)
\lvec(256 49)
\ifill f:0
\move(258 48)
\lvec(259 48)
\lvec(259 49)
\lvec(258 49)
\ifill f:0
\move(260 48)
\lvec(263 48)
\lvec(263 49)
\lvec(260 49)
\ifill f:0
\move(264 48)
\lvec(265 48)
\lvec(265 49)
\lvec(264 49)
\ifill f:0
\move(266 48)
\lvec(267 48)
\lvec(267 49)
\lvec(266 49)
\ifill f:0
\move(268 48)
\lvec(269 48)
\lvec(269 49)
\lvec(268 49)
\ifill f:0
\move(270 48)
\lvec(271 48)
\lvec(271 49)
\lvec(270 49)
\ifill f:0
\move(272 48)
\lvec(290 48)
\lvec(290 49)
\lvec(272 49)
\ifill f:0
\move(292 48)
\lvec(296 48)
\lvec(296 49)
\lvec(292 49)
\ifill f:0
\move(297 48)
\lvec(298 48)
\lvec(298 49)
\lvec(297 49)
\ifill f:0
\move(299 48)
\lvec(306 48)
\lvec(306 49)
\lvec(299 49)
\ifill f:0
\move(307 48)
\lvec(314 48)
\lvec(314 49)
\lvec(307 49)
\ifill f:0
\move(315 48)
\lvec(317 48)
\lvec(317 49)
\lvec(315 49)
\ifill f:0
\move(318 48)
\lvec(322 48)
\lvec(322 49)
\lvec(318 49)
\ifill f:0
\move(323 48)
\lvec(325 48)
\lvec(325 49)
\lvec(323 49)
\ifill f:0
\move(326 48)
\lvec(344 48)
\lvec(344 49)
\lvec(326 49)
\ifill f:0
\move(345 48)
\lvec(354 48)
\lvec(354 49)
\lvec(345 49)
\ifill f:0
\move(355 48)
\lvec(362 48)
\lvec(362 49)
\lvec(355 49)
\ifill f:0
\move(363 48)
\lvec(368 48)
\lvec(368 49)
\lvec(363 49)
\ifill f:0
\move(369 48)
\lvec(401 48)
\lvec(401 49)
\lvec(369 49)
\ifill f:0
\move(402 48)
\lvec(423 48)
\lvec(423 49)
\lvec(402 49)
\ifill f:0
\move(424 48)
\lvec(429 48)
\lvec(429 49)
\lvec(424 49)
\ifill f:0
\move(430 48)
\lvec(435 48)
\lvec(435 49)
\lvec(430 49)
\ifill f:0
\move(436 48)
\lvec(442 48)
\lvec(442 49)
\lvec(436 49)
\ifill f:0
\move(443 48)
\lvec(448 48)
\lvec(448 49)
\lvec(443 49)
\ifill f:0
\move(449 48)
\lvec(451 48)
\lvec(451 49)
\lvec(449 49)
\ifill f:0
\move(15 49)
\lvec(17 49)
\lvec(17 50)
\lvec(15 50)
\ifill f:0
\move(18 49)
\lvec(21 49)
\lvec(21 50)
\lvec(18 50)
\ifill f:0
\move(25 49)
\lvec(26 49)
\lvec(26 50)
\lvec(25 50)
\ifill f:0
\move(36 49)
\lvec(37 49)
\lvec(37 50)
\lvec(36 50)
\ifill f:0
\move(38 49)
\lvec(39 49)
\lvec(39 50)
\lvec(38 50)
\ifill f:0
\move(40 49)
\lvec(47 49)
\lvec(47 50)
\lvec(40 50)
\ifill f:0
\move(48 49)
\lvec(50 49)
\lvec(50 50)
\lvec(48 50)
\ifill f:0
\move(54 49)
\lvec(55 49)
\lvec(55 50)
\lvec(54 50)
\ifill f:0
\move(57 49)
\lvec(59 49)
\lvec(59 50)
\lvec(57 50)
\ifill f:0
\move(61 49)
\lvec(65 49)
\lvec(65 50)
\lvec(61 50)
\ifill f:0
\move(66 49)
\lvec(82 49)
\lvec(82 50)
\lvec(66 50)
\ifill f:0
\move(83 49)
\lvec(85 49)
\lvec(85 50)
\lvec(83 50)
\ifill f:0
\move(87 49)
\lvec(88 49)
\lvec(88 50)
\lvec(87 50)
\ifill f:0
\move(89 49)
\lvec(91 49)
\lvec(91 50)
\lvec(89 50)
\ifill f:0
\move(92 49)
\lvec(93 49)
\lvec(93 50)
\lvec(92 50)
\ifill f:0
\move(96 49)
\lvec(97 49)
\lvec(97 50)
\lvec(96 50)
\ifill f:0
\move(100 49)
\lvec(101 49)
\lvec(101 50)
\lvec(100 50)
\ifill f:0
\move(102 49)
\lvec(103 49)
\lvec(103 50)
\lvec(102 50)
\ifill f:0
\move(104 49)
\lvec(111 49)
\lvec(111 50)
\lvec(104 50)
\ifill f:0
\move(112 49)
\lvec(114 49)
\lvec(114 50)
\lvec(112 50)
\ifill f:0
\move(115 49)
\lvec(122 49)
\lvec(122 50)
\lvec(115 50)
\ifill f:0
\move(123 49)
\lvec(124 49)
\lvec(124 50)
\lvec(123 50)
\ifill f:0
\move(126 49)
\lvec(130 49)
\lvec(130 50)
\lvec(126 50)
\ifill f:0
\move(131 49)
\lvec(137 49)
\lvec(137 50)
\lvec(131 50)
\ifill f:0
\move(139 49)
\lvec(145 49)
\lvec(145 50)
\lvec(139 50)
\ifill f:0
\move(146 49)
\lvec(170 49)
\lvec(170 50)
\lvec(146 50)
\ifill f:0
\move(172 49)
\lvec(184 49)
\lvec(184 50)
\lvec(172 50)
\ifill f:0
\move(185 49)
\lvec(189 49)
\lvec(189 50)
\lvec(185 50)
\ifill f:0
\move(190 49)
\lvec(194 49)
\lvec(194 50)
\lvec(190 50)
\ifill f:0
\move(196 49)
\lvec(197 49)
\lvec(197 50)
\lvec(196 50)
\ifill f:0
\move(199 49)
\lvec(202 49)
\lvec(202 50)
\lvec(199 50)
\ifill f:0
\move(203 49)
\lvec(216 49)
\lvec(216 50)
\lvec(203 50)
\ifill f:0
\move(217 49)
\lvec(219 49)
\lvec(219 50)
\lvec(217 50)
\ifill f:0
\move(220 49)
\lvec(222 49)
\lvec(222 50)
\lvec(220 50)
\ifill f:0
\move(223 49)
\lvec(226 49)
\lvec(226 50)
\lvec(223 50)
\ifill f:0
\move(227 49)
\lvec(228 49)
\lvec(228 50)
\lvec(227 50)
\ifill f:0
\move(229 49)
\lvec(241 49)
\lvec(241 50)
\lvec(229 50)
\ifill f:0
\move(242 49)
\lvec(248 49)
\lvec(248 50)
\lvec(242 50)
\ifill f:0
\move(249 49)
\lvec(255 49)
\lvec(255 50)
\lvec(249 50)
\ifill f:0
\move(256 49)
\lvec(257 49)
\lvec(257 50)
\lvec(256 50)
\ifill f:0
\move(258 49)
\lvec(290 49)
\lvec(290 50)
\lvec(258 50)
\ifill f:0
\move(292 49)
\lvec(293 49)
\lvec(293 50)
\lvec(292 50)
\ifill f:0
\move(294 49)
\lvec(298 49)
\lvec(298 50)
\lvec(294 50)
\ifill f:0
\move(299 49)
\lvec(305 49)
\lvec(305 50)
\lvec(299 50)
\ifill f:0
\move(306 49)
\lvec(310 49)
\lvec(310 50)
\lvec(306 50)
\ifill f:0
\move(311 49)
\lvec(315 49)
\lvec(315 50)
\lvec(311 50)
\ifill f:0
\move(316 49)
\lvec(320 49)
\lvec(320 50)
\lvec(316 50)
\ifill f:0
\move(321 49)
\lvec(325 49)
\lvec(325 50)
\lvec(321 50)
\ifill f:0
\move(326 49)
\lvec(328 49)
\lvec(328 50)
\lvec(326 50)
\ifill f:0
\move(329 49)
\lvec(362 49)
\lvec(362 50)
\lvec(329 50)
\ifill f:0
\move(363 49)
\lvec(364 49)
\lvec(364 50)
\lvec(363 50)
\ifill f:0
\move(365 49)
\lvec(382 49)
\lvec(382 50)
\lvec(365 50)
\ifill f:0
\move(383 49)
\lvec(386 49)
\lvec(386 50)
\lvec(383 50)
\ifill f:0
\move(387 49)
\lvec(390 49)
\lvec(390 50)
\lvec(387 50)
\ifill f:0
\move(391 49)
\lvec(394 49)
\lvec(394 50)
\lvec(391 50)
\ifill f:0
\move(395 49)
\lvec(398 49)
\lvec(398 50)
\lvec(395 50)
\ifill f:0
\move(399 49)
\lvec(401 49)
\lvec(401 50)
\lvec(399 50)
\ifill f:0
\move(402 49)
\lvec(411 49)
\lvec(411 50)
\lvec(402 50)
\ifill f:0
\move(412 49)
\lvec(425 49)
\lvec(425 50)
\lvec(412 50)
\ifill f:0
\move(426 49)
\lvec(430 49)
\lvec(430 50)
\lvec(426 50)
\ifill f:0
\move(431 49)
\lvec(442 49)
\lvec(442 50)
\lvec(431 50)
\ifill f:0
\move(443 49)
\lvec(447 49)
\lvec(447 50)
\lvec(443 50)
\ifill f:0
\move(448 49)
\lvec(451 49)
\lvec(451 50)
\lvec(448 50)
\ifill f:0
\move(15 50)
\lvec(17 50)
\lvec(17 51)
\lvec(15 51)
\ifill f:0
\move(20 50)
\lvec(21 50)
\lvec(21 51)
\lvec(20 51)
\ifill f:0
\move(25 50)
\lvec(26 50)
\lvec(26 51)
\lvec(25 51)
\ifill f:0
\move(36 50)
\lvec(37 50)
\lvec(37 51)
\lvec(36 51)
\ifill f:0
\move(38 50)
\lvec(39 50)
\lvec(39 51)
\lvec(38 51)
\ifill f:0
\move(40 50)
\lvec(41 50)
\lvec(41 51)
\lvec(40 51)
\ifill f:0
\move(42 50)
\lvec(46 50)
\lvec(46 51)
\lvec(42 51)
\ifill f:0
\move(47 50)
\lvec(50 50)
\lvec(50 51)
\lvec(47 51)
\ifill f:0
\move(51 50)
\lvec(52 50)
\lvec(52 51)
\lvec(51 51)
\ifill f:0
\move(55 50)
\lvec(57 50)
\lvec(57 51)
\lvec(55 51)
\ifill f:0
\move(58 50)
\lvec(61 50)
\lvec(61 51)
\lvec(58 51)
\ifill f:0
\move(62 50)
\lvec(65 50)
\lvec(65 51)
\lvec(62 51)
\ifill f:0
\move(66 50)
\lvec(68 50)
\lvec(68 51)
\lvec(66 51)
\ifill f:0
\move(69 50)
\lvec(70 50)
\lvec(70 51)
\lvec(69 51)
\ifill f:0
\move(71 50)
\lvec(72 50)
\lvec(72 51)
\lvec(71 51)
\ifill f:0
\move(76 50)
\lvec(82 50)
\lvec(82 51)
\lvec(76 51)
\ifill f:0
\move(83 50)
\lvec(86 50)
\lvec(86 51)
\lvec(83 51)
\ifill f:0
\move(87 50)
\lvec(89 50)
\lvec(89 51)
\lvec(87 51)
\ifill f:0
\move(91 50)
\lvec(92 50)
\lvec(92 51)
\lvec(91 51)
\ifill f:0
\move(94 50)
\lvec(95 50)
\lvec(95 51)
\lvec(94 51)
\ifill f:0
\move(96 50)
\lvec(98 50)
\lvec(98 51)
\lvec(96 51)
\ifill f:0
\move(100 50)
\lvec(101 50)
\lvec(101 51)
\lvec(100 51)
\ifill f:0
\move(102 50)
\lvec(103 50)
\lvec(103 51)
\lvec(102 51)
\ifill f:0
\move(104 50)
\lvec(105 50)
\lvec(105 51)
\lvec(104 51)
\ifill f:0
\move(106 50)
\lvec(122 50)
\lvec(122 51)
\lvec(106 51)
\ifill f:0
\move(123 50)
\lvec(124 50)
\lvec(124 51)
\lvec(123 51)
\ifill f:0
\move(125 50)
\lvec(128 50)
\lvec(128 51)
\lvec(125 51)
\ifill f:0
\move(129 50)
\lvec(133 50)
\lvec(133 51)
\lvec(129 51)
\ifill f:0
\move(134 50)
\lvec(138 50)
\lvec(138 51)
\lvec(134 51)
\ifill f:0
\move(139 50)
\lvec(140 50)
\lvec(140 51)
\lvec(139 51)
\ifill f:0
\move(141 50)
\lvec(145 50)
\lvec(145 51)
\lvec(141 51)
\ifill f:0
\move(146 50)
\lvec(149 50)
\lvec(149 51)
\lvec(146 51)
\ifill f:0
\move(151 50)
\lvec(170 50)
\lvec(170 51)
\lvec(151 51)
\ifill f:0
\move(174 50)
\lvec(180 50)
\lvec(180 51)
\lvec(174 51)
\ifill f:0
\move(181 50)
\lvec(194 50)
\lvec(194 51)
\lvec(181 51)
\ifill f:0
\move(195 50)
\lvec(197 50)
\lvec(197 51)
\lvec(195 51)
\ifill f:0
\move(198 50)
\lvec(203 50)
\lvec(203 51)
\lvec(198 51)
\ifill f:0
\move(204 50)
\lvec(215 50)
\lvec(215 51)
\lvec(204 51)
\ifill f:0
\move(216 50)
\lvec(226 50)
\lvec(226 51)
\lvec(216 51)
\ifill f:0
\move(227 50)
\lvec(228 50)
\lvec(228 51)
\lvec(227 51)
\ifill f:0
\move(229 50)
\lvec(231 50)
\lvec(231 51)
\lvec(229 51)
\ifill f:0
\move(232 50)
\lvec(234 50)
\lvec(234 51)
\lvec(232 51)
\ifill f:0
\move(235 50)
\lvec(250 50)
\lvec(250 51)
\lvec(235 51)
\ifill f:0
\move(251 50)
\lvec(255 50)
\lvec(255 51)
\lvec(251 51)
\ifill f:0
\move(256 50)
\lvec(257 50)
\lvec(257 51)
\lvec(256 51)
\ifill f:0
\move(258 50)
\lvec(264 50)
\lvec(264 51)
\lvec(258 51)
\ifill f:0
\move(265 50)
\lvec(271 50)
\lvec(271 51)
\lvec(265 51)
\ifill f:0
\move(272 50)
\lvec(275 50)
\lvec(275 51)
\lvec(272 51)
\ifill f:0
\move(276 50)
\lvec(277 50)
\lvec(277 51)
\lvec(276 51)
\ifill f:0
\move(278 50)
\lvec(279 50)
\lvec(279 51)
\lvec(278 51)
\ifill f:0
\move(280 50)
\lvec(281 50)
\lvec(281 51)
\lvec(280 51)
\ifill f:0
\move(282 50)
\lvec(285 50)
\lvec(285 51)
\lvec(282 51)
\ifill f:0
\move(286 50)
\lvec(287 50)
\lvec(287 51)
\lvec(286 51)
\ifill f:0
\move(288 50)
\lvec(290 50)
\lvec(290 51)
\lvec(288 51)
\ifill f:0
\move(292 50)
\lvec(293 50)
\lvec(293 51)
\lvec(292 51)
\ifill f:0
\move(294 50)
\lvec(295 50)
\lvec(295 51)
\lvec(294 51)
\ifill f:0
\move(296 50)
\lvec(304 50)
\lvec(304 51)
\lvec(296 51)
\ifill f:0
\move(305 50)
\lvec(313 50)
\lvec(313 51)
\lvec(305 51)
\ifill f:0
\move(314 50)
\lvec(320 50)
\lvec(320 51)
\lvec(314 51)
\ifill f:0
\move(321 50)
\lvec(325 50)
\lvec(325 51)
\lvec(321 51)
\ifill f:0
\move(326 50)
\lvec(333 50)
\lvec(333 51)
\lvec(326 51)
\ifill f:0
\move(334 50)
\lvec(338 50)
\lvec(338 51)
\lvec(334 51)
\ifill f:0
\move(339 50)
\lvec(341 50)
\lvec(341 51)
\lvec(339 51)
\ifill f:0
\move(342 50)
\lvec(352 50)
\lvec(352 51)
\lvec(342 51)
\ifill f:0
\move(353 50)
\lvec(355 50)
\lvec(355 51)
\lvec(353 51)
\ifill f:0
\move(356 50)
\lvec(358 50)
\lvec(358 51)
\lvec(356 51)
\ifill f:0
\move(359 50)
\lvec(362 50)
\lvec(362 51)
\lvec(359 51)
\ifill f:0
\move(363 50)
\lvec(364 50)
\lvec(364 51)
\lvec(363 51)
\ifill f:0
\move(365 50)
\lvec(367 50)
\lvec(367 51)
\lvec(365 51)
\ifill f:0
\move(368 50)
\lvec(370 50)
\lvec(370 51)
\lvec(368 51)
\ifill f:0
\move(371 50)
\lvec(377 50)
\lvec(377 51)
\lvec(371 51)
\ifill f:0
\move(378 50)
\lvec(387 50)
\lvec(387 51)
\lvec(378 51)
\ifill f:0
\move(388 50)
\lvec(394 50)
\lvec(394 51)
\lvec(388 51)
\ifill f:0
\move(395 50)
\lvec(398 50)
\lvec(398 51)
\lvec(395 51)
\ifill f:0
\move(399 50)
\lvec(401 50)
\lvec(401 51)
\lvec(399 51)
\ifill f:0
\move(403 50)
\lvec(406 50)
\lvec(406 51)
\lvec(403 51)
\ifill f:0
\move(407 50)
\lvec(410 50)
\lvec(410 51)
\lvec(407 51)
\ifill f:0
\move(411 50)
\lvec(414 50)
\lvec(414 51)
\lvec(411 51)
\ifill f:0
\move(415 50)
\lvec(418 50)
\lvec(418 51)
\lvec(415 51)
\ifill f:0
\move(419 50)
\lvec(427 50)
\lvec(427 51)
\lvec(419 51)
\ifill f:0
\move(428 50)
\lvec(436 50)
\lvec(436 51)
\lvec(428 51)
\ifill f:0
\move(437 50)
\lvec(442 50)
\lvec(442 51)
\lvec(437 51)
\ifill f:0
\move(443 50)
\lvec(445 50)
\lvec(445 51)
\lvec(443 51)
\ifill f:0
\move(447 50)
\lvec(451 50)
\lvec(451 51)
\lvec(447 51)
\ifill f:0
\move(14 51)
\lvec(15 51)
\lvec(15 52)
\lvec(14 52)
\ifill f:0
\move(16 51)
\lvec(17 51)
\lvec(17 52)
\lvec(16 52)
\ifill f:0
\move(20 51)
\lvec(24 51)
\lvec(24 52)
\lvec(20 52)
\ifill f:0
\move(25 51)
\lvec(26 51)
\lvec(26 52)
\lvec(25 52)
\ifill f:0
\move(27 51)
\lvec(29 51)
\lvec(29 52)
\lvec(27 52)
\ifill f:0
\move(36 51)
\lvec(37 51)
\lvec(37 52)
\lvec(36 52)
\ifill f:0
\move(39 51)
\lvec(45 51)
\lvec(45 52)
\lvec(39 52)
\ifill f:0
\move(47 51)
\lvec(50 51)
\lvec(50 52)
\lvec(47 52)
\ifill f:0
\move(51 51)
\lvec(52 51)
\lvec(52 52)
\lvec(51 52)
\ifill f:0
\move(57 51)
\lvec(58 51)
\lvec(58 52)
\lvec(57 52)
\ifill f:0
\move(59 51)
\lvec(62 51)
\lvec(62 52)
\lvec(59 52)
\ifill f:0
\move(63 51)
\lvec(65 51)
\lvec(65 52)
\lvec(63 52)
\ifill f:0
\move(66 51)
\lvec(82 51)
\lvec(82 52)
\lvec(66 52)
\ifill f:0
\move(84 51)
\lvec(87 51)
\lvec(87 52)
\lvec(84 52)
\ifill f:0
\move(88 51)
\lvec(91 51)
\lvec(91 52)
\lvec(88 52)
\ifill f:0
\move(92 51)
\lvec(93 51)
\lvec(93 52)
\lvec(92 52)
\ifill f:0
\move(95 51)
\lvec(99 51)
\lvec(99 52)
\lvec(95 52)
\ifill f:0
\move(100 51)
\lvec(101 51)
\lvec(101 52)
\lvec(100 52)
\ifill f:0
\move(102 51)
\lvec(113 51)
\lvec(113 52)
\lvec(102 52)
\ifill f:0
\move(114 51)
\lvec(122 51)
\lvec(122 52)
\lvec(114 52)
\ifill f:0
\move(123 51)
\lvec(124 51)
\lvec(124 52)
\lvec(123 52)
\ifill f:0
\move(125 51)
\lvec(127 51)
\lvec(127 52)
\lvec(125 52)
\ifill f:0
\move(128 51)
\lvec(131 51)
\lvec(131 52)
\lvec(128 52)
\ifill f:0
\move(132 51)
\lvec(134 51)
\lvec(134 52)
\lvec(132 52)
\ifill f:0
\move(136 51)
\lvec(145 51)
\lvec(145 52)
\lvec(136 52)
\ifill f:0
\move(146 51)
\lvec(147 51)
\lvec(147 52)
\lvec(146 52)
\ifill f:0
\move(149 51)
\lvec(159 51)
\lvec(159 52)
\lvec(149 52)
\ifill f:0
\move(166 51)
\lvec(167 51)
\lvec(167 52)
\lvec(166 52)
\ifill f:0
\move(169 51)
\lvec(170 51)
\lvec(170 52)
\lvec(169 52)
\ifill f:0
\move(175 51)
\lvec(176 51)
\lvec(176 52)
\lvec(175 52)
\ifill f:0
\move(177 51)
\lvec(185 51)
\lvec(185 52)
\lvec(177 52)
\ifill f:0
\move(187 51)
\lvec(193 51)
\lvec(193 52)
\lvec(187 52)
\ifill f:0
\move(195 51)
\lvec(197 51)
\lvec(197 52)
\lvec(195 52)
\ifill f:0
\move(198 51)
\lvec(199 51)
\lvec(199 52)
\lvec(198 52)
\ifill f:0
\move(200 51)
\lvec(209 51)
\lvec(209 52)
\lvec(200 52)
\ifill f:0
\move(210 51)
\lvec(226 51)
\lvec(226 52)
\lvec(210 52)
\ifill f:0
\move(227 51)
\lvec(229 51)
\lvec(229 52)
\lvec(227 52)
\ifill f:0
\move(230 51)
\lvec(232 51)
\lvec(232 52)
\lvec(230 52)
\ifill f:0
\move(233 51)
\lvec(235 51)
\lvec(235 52)
\lvec(233 52)
\ifill f:0
\move(236 51)
\lvec(238 51)
\lvec(238 52)
\lvec(236 52)
\ifill f:0
\move(239 51)
\lvec(241 51)
\lvec(241 52)
\lvec(239 52)
\ifill f:0
\move(242 51)
\lvec(244 51)
\lvec(244 52)
\lvec(242 52)
\ifill f:0
\move(245 51)
\lvec(250 51)
\lvec(250 52)
\lvec(245 52)
\ifill f:0
\move(251 51)
\lvec(252 51)
\lvec(252 52)
\lvec(251 52)
\ifill f:0
\move(253 51)
\lvec(255 51)
\lvec(255 52)
\lvec(253 52)
\ifill f:0
\move(256 51)
\lvec(257 51)
\lvec(257 52)
\lvec(256 52)
\ifill f:0
\move(258 51)
\lvec(267 51)
\lvec(267 52)
\lvec(258 52)
\ifill f:0
\move(268 51)
\lvec(274 51)
\lvec(274 52)
\lvec(268 52)
\ifill f:0
\move(275 51)
\lvec(287 51)
\lvec(287 52)
\lvec(275 52)
\ifill f:0
\move(288 51)
\lvec(290 51)
\lvec(290 52)
\lvec(288 52)
\ifill f:0
\move(292 51)
\lvec(293 51)
\lvec(293 52)
\lvec(292 52)
\ifill f:0
\move(294 51)
\lvec(297 51)
\lvec(297 52)
\lvec(294 52)
\ifill f:0
\move(298 51)
\lvec(299 51)
\lvec(299 52)
\lvec(298 52)
\ifill f:0
\move(300 51)
\lvec(301 51)
\lvec(301 52)
\lvec(300 52)
\ifill f:0
\move(302 51)
\lvec(305 51)
\lvec(305 52)
\lvec(302 52)
\ifill f:0
\move(306 51)
\lvec(316 51)
\lvec(316 52)
\lvec(306 52)
\ifill f:0
\move(317 51)
\lvec(325 51)
\lvec(325 52)
\lvec(317 52)
\ifill f:0
\move(326 51)
\lvec(362 51)
\lvec(362 52)
\lvec(326 52)
\ifill f:0
\move(363 51)
\lvec(385 51)
\lvec(385 52)
\lvec(363 52)
\ifill f:0
\move(386 51)
\lvec(395 51)
\lvec(395 52)
\lvec(386 52)
\ifill f:0
\move(396 51)
\lvec(401 51)
\lvec(401 52)
\lvec(396 52)
\ifill f:0
\move(402 51)
\lvec(409 51)
\lvec(409 52)
\lvec(402 52)
\ifill f:0
\move(410 51)
\lvec(442 51)
\lvec(442 52)
\lvec(410 52)
\ifill f:0
\move(443 51)
\lvec(445 51)
\lvec(445 52)
\lvec(443 52)
\ifill f:0
\move(446 51)
\lvec(450 51)
\lvec(450 52)
\lvec(446 52)
\ifill f:0
\move(16 52)
\lvec(17 52)
\lvec(17 53)
\lvec(16 53)
\ifill f:0
\move(18 52)
\lvec(19 52)
\lvec(19 53)
\lvec(18 53)
\ifill f:0
\move(20 52)
\lvec(21 52)
\lvec(21 53)
\lvec(20 53)
\ifill f:0
\move(23 52)
\lvec(24 52)
\lvec(24 53)
\lvec(23 53)
\ifill f:0
\move(25 52)
\lvec(26 52)
\lvec(26 53)
\lvec(25 53)
\ifill f:0
\move(36 52)
\lvec(37 52)
\lvec(37 53)
\lvec(36 53)
\ifill f:0
\move(38 52)
\lvec(43 52)
\lvec(43 53)
\lvec(38 53)
\ifill f:0
\move(44 52)
\lvec(45 52)
\lvec(45 53)
\lvec(44 53)
\ifill f:0
\move(47 52)
\lvec(50 52)
\lvec(50 53)
\lvec(47 53)
\ifill f:0
\move(51 52)
\lvec(53 52)
\lvec(53 53)
\lvec(51 53)
\ifill f:0
\move(54 52)
\lvec(55 52)
\lvec(55 53)
\lvec(54 53)
\ifill f:0
\move(56 52)
\lvec(57 52)
\lvec(57 53)
\lvec(56 53)
\ifill f:0
\move(58 52)
\lvec(59 52)
\lvec(59 53)
\lvec(58 53)
\ifill f:0
\move(60 52)
\lvec(62 52)
\lvec(62 53)
\lvec(60 53)
\ifill f:0
\move(63 52)
\lvec(65 52)
\lvec(65 53)
\lvec(63 53)
\ifill f:0
\move(67 52)
\lvec(71 52)
\lvec(71 53)
\lvec(67 53)
\ifill f:0
\move(73 52)
\lvec(82 52)
\lvec(82 53)
\lvec(73 53)
\ifill f:0
\move(86 52)
\lvec(88 52)
\lvec(88 53)
\lvec(86 53)
\ifill f:0
\move(90 52)
\lvec(93 52)
\lvec(93 53)
\lvec(90 53)
\ifill f:0
\move(95 52)
\lvec(96 52)
\lvec(96 53)
\lvec(95 53)
\ifill f:0
\move(97 52)
\lvec(99 52)
\lvec(99 53)
\lvec(97 53)
\ifill f:0
\move(100 52)
\lvec(101 52)
\lvec(101 53)
\lvec(100 53)
\ifill f:0
\move(102 52)
\lvec(106 52)
\lvec(106 53)
\lvec(102 53)
\ifill f:0
\move(107 52)
\lvec(116 52)
\lvec(116 53)
\lvec(107 53)
\ifill f:0
\move(117 52)
\lvec(122 52)
\lvec(122 53)
\lvec(117 53)
\ifill f:0
\move(124 52)
\lvec(126 52)
\lvec(126 53)
\lvec(124 53)
\ifill f:0
\move(127 52)
\lvec(129 52)
\lvec(129 53)
\lvec(127 53)
\ifill f:0
\move(130 52)
\lvec(133 52)
\lvec(133 53)
\lvec(130 53)
\ifill f:0
\move(134 52)
\lvec(137 52)
\lvec(137 53)
\lvec(134 53)
\ifill f:0
\move(138 52)
\lvec(145 52)
\lvec(145 53)
\lvec(138 53)
\ifill f:0
\move(146 52)
\lvec(153 52)
\lvec(153 53)
\lvec(146 53)
\ifill f:0
\move(155 52)
\lvec(163 52)
\lvec(163 53)
\lvec(155 53)
\ifill f:0
\move(164 52)
\lvec(167 52)
\lvec(167 53)
\lvec(164 53)
\ifill f:0
\move(169 52)
\lvec(170 52)
\lvec(170 53)
\lvec(169 53)
\ifill f:0
\move(172 52)
\lvec(173 52)
\lvec(173 53)
\lvec(172 53)
\ifill f:0
\move(178 52)
\lvec(179 52)
\lvec(179 53)
\lvec(178 53)
\ifill f:0
\move(181 52)
\lvec(192 52)
\lvec(192 53)
\lvec(181 53)
\ifill f:0
\move(193 52)
\lvec(194 52)
\lvec(194 53)
\lvec(193 53)
\ifill f:0
\move(195 52)
\lvec(197 52)
\lvec(197 53)
\lvec(195 53)
\ifill f:0
\move(198 52)
\lvec(200 52)
\lvec(200 53)
\lvec(198 53)
\ifill f:0
\move(201 52)
\lvec(206 52)
\lvec(206 53)
\lvec(201 53)
\ifill f:0
\move(207 52)
\lvec(226 52)
\lvec(226 53)
\lvec(207 53)
\ifill f:0
\move(227 52)
\lvec(229 52)
\lvec(229 53)
\lvec(227 53)
\ifill f:0
\move(230 52)
\lvec(233 52)
\lvec(233 53)
\lvec(230 53)
\ifill f:0
\move(234 52)
\lvec(236 52)
\lvec(236 53)
\lvec(234 53)
\ifill f:0
\move(237 52)
\lvec(257 52)
\lvec(257 53)
\lvec(237 53)
\ifill f:0
\move(258 52)
\lvec(260 52)
\lvec(260 53)
\lvec(258 53)
\ifill f:0
\move(261 52)
\lvec(280 52)
\lvec(280 53)
\lvec(261 53)
\ifill f:0
\move(281 52)
\lvec(287 52)
\lvec(287 53)
\lvec(281 53)
\ifill f:0
\move(288 52)
\lvec(290 52)
\lvec(290 53)
\lvec(288 53)
\ifill f:0
\move(291 52)
\lvec(298 52)
\lvec(298 53)
\lvec(291 53)
\ifill f:0
\move(299 52)
\lvec(310 52)
\lvec(310 53)
\lvec(299 53)
\ifill f:0
\move(311 52)
\lvec(323 52)
\lvec(323 53)
\lvec(311 53)
\ifill f:0
\move(324 52)
\lvec(325 52)
\lvec(325 53)
\lvec(324 53)
\ifill f:0
\move(326 52)
\lvec(341 52)
\lvec(341 53)
\lvec(326 53)
\ifill f:0
\move(342 52)
\lvec(362 52)
\lvec(362 53)
\lvec(342 53)
\ifill f:0
\move(363 52)
\lvec(369 52)
\lvec(369 53)
\lvec(363 53)
\ifill f:0
\move(370 52)
\lvec(383 52)
\lvec(383 53)
\lvec(370 53)
\ifill f:0
\move(384 52)
\lvec(386 52)
\lvec(386 53)
\lvec(384 53)
\ifill f:0
\move(387 52)
\lvec(389 52)
\lvec(389 53)
\lvec(387 53)
\ifill f:0
\move(390 52)
\lvec(392 52)
\lvec(392 53)
\lvec(390 53)
\ifill f:0
\move(393 52)
\lvec(395 52)
\lvec(395 53)
\lvec(393 53)
\ifill f:0
\move(396 52)
\lvec(401 52)
\lvec(401 53)
\lvec(396 53)
\ifill f:0
\move(402 52)
\lvec(405 52)
\lvec(405 53)
\lvec(402 53)
\ifill f:0
\move(406 52)
\lvec(415 52)
\lvec(415 53)
\lvec(406 53)
\ifill f:0
\move(416 52)
\lvec(433 52)
\lvec(433 53)
\lvec(416 53)
\ifill f:0
\move(434 52)
\lvec(442 52)
\lvec(442 53)
\lvec(434 53)
\ifill f:0
\move(443 52)
\lvec(445 52)
\lvec(445 53)
\lvec(443 53)
\ifill f:0
\move(446 52)
\lvec(449 52)
\lvec(449 53)
\lvec(446 53)
\ifill f:0
\move(450 52)
\lvec(451 52)
\lvec(451 53)
\lvec(450 53)
\ifill f:0
\move(16 53)
\lvec(17 53)
\lvec(17 54)
\lvec(16 54)
\ifill f:0
\move(18 53)
\lvec(19 53)
\lvec(19 54)
\lvec(18 54)
\ifill f:0
\move(20 53)
\lvec(21 53)
\lvec(21 54)
\lvec(20 54)
\ifill f:0
\move(23 53)
\lvec(26 53)
\lvec(26 54)
\lvec(23 54)
\ifill f:0
\move(36 53)
\lvec(37 53)
\lvec(37 54)
\lvec(36 54)
\ifill f:0
\move(38 53)
\lvec(39 53)
\lvec(39 54)
\lvec(38 54)
\ifill f:0
\move(40 53)
\lvec(46 53)
\lvec(46 54)
\lvec(40 54)
\ifill f:0
\move(47 53)
\lvec(50 53)
\lvec(50 54)
\lvec(47 54)
\ifill f:0
\move(52 53)
\lvec(53 53)
\lvec(53 54)
\lvec(52 54)
\ifill f:0
\move(57 53)
\lvec(58 53)
\lvec(58 54)
\lvec(57 54)
\ifill f:0
\move(59 53)
\lvec(60 53)
\lvec(60 54)
\lvec(59 54)
\ifill f:0
\move(61 53)
\lvec(65 53)
\lvec(65 54)
\lvec(61 54)
\ifill f:0
\move(66 53)
\lvec(74 53)
\lvec(74 54)
\lvec(66 54)
\ifill f:0
\move(77 53)
\lvec(78 53)
\lvec(78 54)
\lvec(77 54)
\ifill f:0
\move(81 53)
\lvec(82 53)
\lvec(82 54)
\lvec(81 54)
\ifill f:0
\move(87 53)
\lvec(91 53)
\lvec(91 54)
\lvec(87 54)
\ifill f:0
\move(92 53)
\lvec(95 53)
\lvec(95 54)
\lvec(92 54)
\ifill f:0
\move(97 53)
\lvec(99 53)
\lvec(99 54)
\lvec(97 54)
\ifill f:0
\move(100 53)
\lvec(101 53)
\lvec(101 54)
\lvec(100 54)
\ifill f:0
\move(102 53)
\lvec(107 53)
\lvec(107 54)
\lvec(102 54)
\ifill f:0
\move(108 53)
\lvec(122 53)
\lvec(122 54)
\lvec(108 54)
\ifill f:0
\move(124 53)
\lvec(131 53)
\lvec(131 54)
\lvec(124 54)
\ifill f:0
\move(132 53)
\lvec(134 53)
\lvec(134 54)
\lvec(132 54)
\ifill f:0
\move(135 53)
\lvec(138 53)
\lvec(138 54)
\lvec(135 54)
\ifill f:0
\move(139 53)
\lvec(142 53)
\lvec(142 54)
\lvec(139 54)
\ifill f:0
\move(143 53)
\lvec(145 53)
\lvec(145 54)
\lvec(143 54)
\ifill f:0
\move(146 53)
\lvec(151 53)
\lvec(151 54)
\lvec(146 54)
\ifill f:0
\move(152 53)
\lvec(167 53)
\lvec(167 54)
\lvec(152 54)
\ifill f:0
\move(169 53)
\lvec(170 53)
\lvec(170 54)
\lvec(169 54)
\ifill f:0
\move(171 53)
\lvec(189 53)
\lvec(189 54)
\lvec(171 54)
\ifill f:0
\move(193 53)
\lvec(194 53)
\lvec(194 54)
\lvec(193 54)
\ifill f:0
\move(195 53)
\lvec(197 53)
\lvec(197 54)
\lvec(195 54)
\ifill f:0
\move(198 53)
\lvec(201 53)
\lvec(201 54)
\lvec(198 54)
\ifill f:0
\move(202 53)
\lvec(215 53)
\lvec(215 54)
\lvec(202 54)
\ifill f:0
\move(216 53)
\lvec(220 53)
\lvec(220 54)
\lvec(216 54)
\ifill f:0
\move(221 53)
\lvec(226 53)
\lvec(226 54)
\lvec(221 54)
\ifill f:0
\move(227 53)
\lvec(234 53)
\lvec(234 54)
\lvec(227 54)
\ifill f:0
\move(235 53)
\lvec(242 53)
\lvec(242 54)
\lvec(235 54)
\ifill f:0
\move(243 53)
\lvec(248 53)
\lvec(248 54)
\lvec(243 54)
\ifill f:0
\move(249 53)
\lvec(257 53)
\lvec(257 54)
\lvec(249 54)
\ifill f:0
\move(258 53)
\lvec(269 53)
\lvec(269 54)
\lvec(258 54)
\ifill f:0
\move(270 53)
\lvec(290 53)
\lvec(290 54)
\lvec(270 54)
\ifill f:0
\move(291 53)
\lvec(296 53)
\lvec(296 54)
\lvec(291 54)
\ifill f:0
\move(297 53)
\lvec(307 53)
\lvec(307 54)
\lvec(297 54)
\ifill f:0
\move(308 53)
\lvec(309 53)
\lvec(309 54)
\lvec(308 54)
\ifill f:0
\move(310 53)
\lvec(311 53)
\lvec(311 54)
\lvec(310 54)
\ifill f:0
\move(312 53)
\lvec(313 53)
\lvec(313 54)
\lvec(312 54)
\ifill f:0
\move(314 53)
\lvec(315 53)
\lvec(315 54)
\lvec(314 54)
\ifill f:0
\move(316 53)
\lvec(317 53)
\lvec(317 54)
\lvec(316 54)
\ifill f:0
\move(318 53)
\lvec(319 53)
\lvec(319 54)
\lvec(318 54)
\ifill f:0
\move(320 53)
\lvec(321 53)
\lvec(321 54)
\lvec(320 54)
\ifill f:0
\move(322 53)
\lvec(323 53)
\lvec(323 54)
\lvec(322 54)
\ifill f:0
\move(324 53)
\lvec(325 53)
\lvec(325 54)
\lvec(324 54)
\ifill f:0
\move(326 53)
\lvec(340 53)
\lvec(340 54)
\lvec(326 54)
\ifill f:0
\move(341 53)
\lvec(356 53)
\lvec(356 54)
\lvec(341 54)
\ifill f:0
\move(357 53)
\lvec(362 53)
\lvec(362 54)
\lvec(357 54)
\ifill f:0
\move(363 53)
\lvec(366 53)
\lvec(366 54)
\lvec(363 54)
\ifill f:0
\move(367 53)
\lvec(371 53)
\lvec(371 54)
\lvec(367 54)
\ifill f:0
\move(372 53)
\lvec(379 53)
\lvec(379 54)
\lvec(372 54)
\ifill f:0
\move(380 53)
\lvec(387 53)
\lvec(387 54)
\lvec(380 54)
\ifill f:0
\move(388 53)
\lvec(390 53)
\lvec(390 54)
\lvec(388 54)
\ifill f:0
\move(391 53)
\lvec(401 53)
\lvec(401 54)
\lvec(391 54)
\ifill f:0
\move(402 53)
\lvec(417 53)
\lvec(417 54)
\lvec(402 54)
\ifill f:0
\move(418 53)
\lvec(427 53)
\lvec(427 54)
\lvec(418 54)
\ifill f:0
\move(428 53)
\lvec(434 53)
\lvec(434 54)
\lvec(428 54)
\ifill f:0
\move(435 53)
\lvec(442 53)
\lvec(442 54)
\lvec(435 54)
\ifill f:0
\move(443 53)
\lvec(451 53)
\lvec(451 54)
\lvec(443 54)
\ifill f:0
\move(16 54)
\lvec(17 54)
\lvec(17 55)
\lvec(16 55)
\ifill f:0
\move(19 54)
\lvec(21 54)
\lvec(21 55)
\lvec(19 55)
\ifill f:0
\move(22 54)
\lvec(23 54)
\lvec(23 55)
\lvec(22 55)
\ifill f:0
\move(24 54)
\lvec(26 54)
\lvec(26 55)
\lvec(24 55)
\ifill f:0
\move(36 54)
\lvec(37 54)
\lvec(37 55)
\lvec(36 55)
\ifill f:0
\move(38 54)
\lvec(39 54)
\lvec(39 55)
\lvec(38 55)
\ifill f:0
\move(40 54)
\lvec(43 54)
\lvec(43 55)
\lvec(40 55)
\ifill f:0
\move(45 54)
\lvec(50 54)
\lvec(50 55)
\lvec(45 55)
\ifill f:0
\move(56 54)
\lvec(57 54)
\lvec(57 55)
\lvec(56 55)
\ifill f:0
\move(58 54)
\lvec(65 54)
\lvec(65 55)
\lvec(58 55)
\ifill f:0
\move(66 54)
\lvec(71 54)
\lvec(71 55)
\lvec(66 55)
\ifill f:0
\move(73 54)
\lvec(78 54)
\lvec(78 55)
\lvec(73 55)
\ifill f:0
\move(81 54)
\lvec(82 54)
\lvec(82 55)
\lvec(81 55)
\ifill f:0
\move(89 54)
\lvec(93 54)
\lvec(93 55)
\lvec(89 55)
\ifill f:0
\move(96 54)
\lvec(99 54)
\lvec(99 55)
\lvec(96 55)
\ifill f:0
\move(100 54)
\lvec(101 54)
\lvec(101 55)
\lvec(100 55)
\ifill f:0
\move(102 54)
\lvec(105 54)
\lvec(105 55)
\lvec(102 55)
\ifill f:0
\move(106 54)
\lvec(111 54)
\lvec(111 55)
\lvec(106 55)
\ifill f:0
\move(112 54)
\lvec(115 54)
\lvec(115 55)
\lvec(112 55)
\ifill f:0
\move(116 54)
\lvec(117 54)
\lvec(117 55)
\lvec(116 55)
\ifill f:0
\move(118 54)
\lvec(122 54)
\lvec(122 55)
\lvec(118 55)
\ifill f:0
\move(124 54)
\lvec(125 54)
\lvec(125 55)
\lvec(124 55)
\ifill f:0
\move(126 54)
\lvec(130 54)
\lvec(130 55)
\lvec(126 55)
\ifill f:0
\move(131 54)
\lvec(133 54)
\lvec(133 55)
\lvec(131 55)
\ifill f:0
\move(134 54)
\lvec(135 54)
\lvec(135 55)
\lvec(134 55)
\ifill f:0
\move(136 54)
\lvec(138 54)
\lvec(138 55)
\lvec(136 55)
\ifill f:0
\move(139 54)
\lvec(142 54)
\lvec(142 55)
\lvec(139 55)
\ifill f:0
\move(143 54)
\lvec(145 54)
\lvec(145 55)
\lvec(143 55)
\ifill f:0
\move(147 54)
\lvec(150 54)
\lvec(150 55)
\lvec(147 55)
\ifill f:0
\move(151 54)
\lvec(155 54)
\lvec(155 55)
\lvec(151 55)
\ifill f:0
\move(156 54)
\lvec(159 54)
\lvec(159 55)
\lvec(156 55)
\ifill f:0
\move(162 54)
\lvec(168 54)
\lvec(168 55)
\lvec(162 55)
\ifill f:0
\move(169 54)
\lvec(170 54)
\lvec(170 55)
\lvec(169 55)
\ifill f:0
\move(171 54)
\lvec(191 54)
\lvec(191 55)
\lvec(171 55)
\ifill f:0
\move(192 54)
\lvec(197 54)
\lvec(197 55)
\lvec(192 55)
\ifill f:0
\move(198 54)
\lvec(203 54)
\lvec(203 55)
\lvec(198 55)
\ifill f:0
\move(204 54)
\lvec(212 54)
\lvec(212 55)
\lvec(204 55)
\ifill f:0
\move(213 54)
\lvec(219 54)
\lvec(219 55)
\lvec(213 55)
\ifill f:0
\move(220 54)
\lvec(226 54)
\lvec(226 55)
\lvec(220 55)
\ifill f:0
\move(227 54)
\lvec(230 54)
\lvec(230 55)
\lvec(227 55)
\ifill f:0
\move(231 54)
\lvec(235 54)
\lvec(235 55)
\lvec(231 55)
\ifill f:0
\move(236 54)
\lvec(247 54)
\lvec(247 55)
\lvec(236 55)
\ifill f:0
\move(248 54)
\lvec(251 54)
\lvec(251 55)
\lvec(248 55)
\ifill f:0
\move(252 54)
\lvec(257 54)
\lvec(257 55)
\lvec(252 55)
\ifill f:0
\move(258 54)
\lvec(264 54)
\lvec(264 55)
\lvec(258 55)
\ifill f:0
\move(265 54)
\lvec(270 54)
\lvec(270 55)
\lvec(265 55)
\ifill f:0
\move(271 54)
\lvec(273 54)
\lvec(273 55)
\lvec(271 55)
\ifill f:0
\move(274 54)
\lvec(284 54)
\lvec(284 55)
\lvec(274 55)
\ifill f:0
\move(285 54)
\lvec(290 54)
\lvec(290 55)
\lvec(285 55)
\ifill f:0
\move(291 54)
\lvec(294 54)
\lvec(294 55)
\lvec(291 55)
\ifill f:0
\move(295 54)
\lvec(299 54)
\lvec(299 55)
\lvec(295 55)
\ifill f:0
\move(300 54)
\lvec(308 54)
\lvec(308 55)
\lvec(300 55)
\ifill f:0
\move(309 54)
\lvec(319 54)
\lvec(319 55)
\lvec(309 55)
\ifill f:0
\move(320 54)
\lvec(321 54)
\lvec(321 55)
\lvec(320 55)
\ifill f:0
\move(322 54)
\lvec(323 54)
\lvec(323 55)
\lvec(322 55)
\ifill f:0
\move(324 54)
\lvec(325 54)
\lvec(325 55)
\lvec(324 55)
\ifill f:0
\move(326 54)
\lvec(327 54)
\lvec(327 55)
\lvec(326 55)
\ifill f:0
\move(328 54)
\lvec(329 54)
\lvec(329 55)
\lvec(328 55)
\ifill f:0
\move(330 54)
\lvec(331 54)
\lvec(331 55)
\lvec(330 55)
\ifill f:0
\move(332 54)
\lvec(333 54)
\lvec(333 55)
\lvec(332 55)
\ifill f:0
\move(334 54)
\lvec(335 54)
\lvec(335 55)
\lvec(334 55)
\ifill f:0
\move(336 54)
\lvec(337 54)
\lvec(337 55)
\lvec(336 55)
\ifill f:0
\move(338 54)
\lvec(339 54)
\lvec(339 55)
\lvec(338 55)
\ifill f:0
\move(340 54)
\lvec(352 54)
\lvec(352 55)
\lvec(340 55)
\ifill f:0
\move(353 54)
\lvec(362 54)
\lvec(362 55)
\lvec(353 55)
\ifill f:0
\move(364 54)
\lvec(368 54)
\lvec(368 55)
\lvec(364 55)
\ifill f:0
\move(369 54)
\lvec(370 54)
\lvec(370 55)
\lvec(369 55)
\ifill f:0
\move(371 54)
\lvec(388 54)
\lvec(388 55)
\lvec(371 55)
\ifill f:0
\move(389 54)
\lvec(396 54)
\lvec(396 55)
\lvec(389 55)
\ifill f:0
\move(397 54)
\lvec(401 54)
\lvec(401 55)
\lvec(397 55)
\ifill f:0
\move(402 54)
\lvec(407 54)
\lvec(407 55)
\lvec(402 55)
\ifill f:0
\move(408 54)
\lvec(410 54)
\lvec(410 55)
\lvec(408 55)
\ifill f:0
\move(411 54)
\lvec(425 54)
\lvec(425 55)
\lvec(411 55)
\ifill f:0
\move(426 54)
\lvec(428 54)
\lvec(428 55)
\lvec(426 55)
\ifill f:0
\move(429 54)
\lvec(442 54)
\lvec(442 55)
\lvec(429 55)
\ifill f:0
\move(443 54)
\lvec(444 54)
\lvec(444 55)
\lvec(443 55)
\ifill f:0
\move(445 54)
\lvec(448 54)
\lvec(448 55)
\lvec(445 55)
\ifill f:0
\move(449 54)
\lvec(451 54)
\lvec(451 55)
\lvec(449 55)
\ifill f:0
\move(14 55)
\lvec(17 55)
\lvec(17 56)
\lvec(14 56)
\ifill f:0
\move(20 55)
\lvec(21 55)
\lvec(21 56)
\lvec(20 56)
\ifill f:0
\move(24 55)
\lvec(26 55)
\lvec(26 56)
\lvec(24 56)
\ifill f:0
\move(36 55)
\lvec(37 55)
\lvec(37 56)
\lvec(36 56)
\ifill f:0
\move(38 55)
\lvec(41 55)
\lvec(41 56)
\lvec(38 56)
\ifill f:0
\move(42 55)
\lvec(46 55)
\lvec(46 56)
\lvec(42 56)
\ifill f:0
\move(48 55)
\lvec(50 55)
\lvec(50 56)
\lvec(48 56)
\ifill f:0
\move(54 55)
\lvec(56 55)
\lvec(56 56)
\lvec(54 56)
\ifill f:0
\move(57 55)
\lvec(59 55)
\lvec(59 56)
\lvec(57 56)
\ifill f:0
\move(60 55)
\lvec(61 55)
\lvec(61 56)
\lvec(60 56)
\ifill f:0
\move(62 55)
\lvec(65 55)
\lvec(65 56)
\lvec(62 56)
\ifill f:0
\move(66 55)
\lvec(70 55)
\lvec(70 56)
\lvec(66 56)
\ifill f:0
\move(71 55)
\lvec(74 55)
\lvec(74 56)
\lvec(71 56)
\ifill f:0
\move(75 55)
\lvec(79 55)
\lvec(79 56)
\lvec(75 56)
\ifill f:0
\move(81 55)
\lvec(82 55)
\lvec(82 56)
\lvec(81 56)
\ifill f:0
\move(83 55)
\lvec(90 55)
\lvec(90 56)
\lvec(83 56)
\ifill f:0
\move(92 55)
\lvec(94 55)
\lvec(94 56)
\lvec(92 56)
\ifill f:0
\move(95 55)
\lvec(98 55)
\lvec(98 56)
\lvec(95 56)
\ifill f:0
\move(100 55)
\lvec(101 55)
\lvec(101 56)
\lvec(100 56)
\ifill f:0
\move(102 55)
\lvec(106 55)
\lvec(106 56)
\lvec(102 56)
\ifill f:0
\move(107 55)
\lvec(109 55)
\lvec(109 56)
\lvec(107 56)
\ifill f:0
\move(110 55)
\lvec(119 55)
\lvec(119 56)
\lvec(110 56)
\ifill f:0
\move(120 55)
\lvec(122 55)
\lvec(122 56)
\lvec(120 56)
\ifill f:0
\move(123 55)
\lvec(125 55)
\lvec(125 56)
\lvec(123 56)
\ifill f:0
\move(126 55)
\lvec(127 55)
\lvec(127 56)
\lvec(126 56)
\ifill f:0
\move(128 55)
\lvec(129 55)
\lvec(129 56)
\lvec(128 56)
\ifill f:0
\move(130 55)
\lvec(131 55)
\lvec(131 56)
\lvec(130 56)
\ifill f:0
\move(132 55)
\lvec(134 55)
\lvec(134 56)
\lvec(132 56)
\ifill f:0
\move(135 55)
\lvec(145 55)
\lvec(145 56)
\lvec(135 56)
\ifill f:0
\move(146 55)
\lvec(149 55)
\lvec(149 56)
\lvec(146 56)
\ifill f:0
\move(150 55)
\lvec(153 55)
\lvec(153 56)
\lvec(150 56)
\ifill f:0
\move(154 55)
\lvec(168 55)
\lvec(168 56)
\lvec(154 56)
\ifill f:0
\move(169 55)
\lvec(170 55)
\lvec(170 56)
\lvec(169 56)
\ifill f:0
\move(171 55)
\lvec(176 55)
\lvec(176 56)
\lvec(171 56)
\ifill f:0
\move(178 55)
\lvec(197 55)
\lvec(197 56)
\lvec(178 56)
\ifill f:0
\move(198 55)
\lvec(207 55)
\lvec(207 56)
\lvec(198 56)
\ifill f:0
\move(208 55)
\lvec(226 55)
\lvec(226 56)
\lvec(208 56)
\ifill f:0
\move(227 55)
\lvec(231 55)
\lvec(231 56)
\lvec(227 56)
\ifill f:0
\move(232 55)
\lvec(242 55)
\lvec(242 56)
\lvec(232 56)
\ifill f:0
\move(243 55)
\lvec(246 55)
\lvec(246 56)
\lvec(243 56)
\ifill f:0
\move(247 55)
\lvec(250 55)
\lvec(250 56)
\lvec(247 56)
\ifill f:0
\move(251 55)
\lvec(254 55)
\lvec(254 56)
\lvec(251 56)
\ifill f:0
\move(255 55)
\lvec(257 55)
\lvec(257 56)
\lvec(255 56)
\ifill f:0
\move(258 55)
\lvec(265 55)
\lvec(265 56)
\lvec(258 56)
\ifill f:0
\move(266 55)
\lvec(290 55)
\lvec(290 56)
\lvec(266 56)
\ifill f:0
\move(291 55)
\lvec(292 55)
\lvec(292 56)
\lvec(291 56)
\ifill f:0
\move(293 55)
\lvec(297 55)
\lvec(297 56)
\lvec(293 56)
\ifill f:0
\move(298 55)
\lvec(302 55)
\lvec(302 56)
\lvec(298 56)
\ifill f:0
\move(303 55)
\lvec(307 55)
\lvec(307 56)
\lvec(303 56)
\ifill f:0
\move(308 55)
\lvec(314 55)
\lvec(314 56)
\lvec(308 56)
\ifill f:0
\move(315 55)
\lvec(321 55)
\lvec(321 56)
\lvec(315 56)
\ifill f:0
\move(322 55)
\lvec(323 55)
\lvec(323 56)
\lvec(322 56)
\ifill f:0
\move(324 55)
\lvec(325 55)
\lvec(325 56)
\lvec(324 56)
\ifill f:0
\move(326 55)
\lvec(362 55)
\lvec(362 56)
\lvec(326 56)
\ifill f:0
\move(364 55)
\lvec(365 55)
\lvec(365 56)
\lvec(364 56)
\ifill f:0
\move(366 55)
\lvec(379 55)
\lvec(379 56)
\lvec(366 56)
\ifill f:0
\move(380 55)
\lvec(391 55)
\lvec(391 56)
\lvec(380 56)
\ifill f:0
\move(392 55)
\lvec(401 55)
\lvec(401 56)
\lvec(392 56)
\ifill f:0
\move(402 55)
\lvec(404 55)
\lvec(404 56)
\lvec(402 56)
\ifill f:0
\move(405 55)
\lvec(412 55)
\lvec(412 56)
\lvec(405 56)
\ifill f:0
\move(413 55)
\lvec(423 55)
\lvec(423 56)
\lvec(413 56)
\ifill f:0
\move(424 55)
\lvec(426 55)
\lvec(426 56)
\lvec(424 56)
\ifill f:0
\move(427 55)
\lvec(429 55)
\lvec(429 56)
\lvec(427 56)
\ifill f:0
\move(430 55)
\lvec(432 55)
\lvec(432 56)
\lvec(430 56)
\ifill f:0
\move(433 55)
\lvec(435 55)
\lvec(435 56)
\lvec(433 56)
\ifill f:0
\move(436 55)
\lvec(438 55)
\lvec(438 56)
\lvec(436 56)
\ifill f:0
\move(439 55)
\lvec(442 55)
\lvec(442 56)
\lvec(439 56)
\ifill f:0
\move(443 55)
\lvec(444 55)
\lvec(444 56)
\lvec(443 56)
\ifill f:0
\move(445 55)
\lvec(451 55)
\lvec(451 56)
\lvec(445 56)
\ifill f:0
\move(15 56)
\lvec(17 56)
\lvec(17 57)
\lvec(15 57)
\ifill f:0
\move(18 56)
\lvec(19 56)
\lvec(19 57)
\lvec(18 57)
\ifill f:0
\move(20 56)
\lvec(22 56)
\lvec(22 57)
\lvec(20 57)
\ifill f:0
\move(23 56)
\lvec(26 56)
\lvec(26 57)
\lvec(23 57)
\ifill f:0
\move(27 56)
\lvec(28 56)
\lvec(28 57)
\lvec(27 57)
\ifill f:0
\move(36 56)
\lvec(37 56)
\lvec(37 57)
\lvec(36 57)
\ifill f:0
\move(38 56)
\lvec(43 56)
\lvec(43 57)
\lvec(38 57)
\ifill f:0
\move(44 56)
\lvec(46 56)
\lvec(46 57)
\lvec(44 57)
\ifill f:0
\move(49 56)
\lvec(50 56)
\lvec(50 57)
\lvec(49 57)
\ifill f:0
\move(56 56)
\lvec(58 56)
\lvec(58 57)
\lvec(56 57)
\ifill f:0
\move(59 56)
\lvec(63 56)
\lvec(63 57)
\lvec(59 57)
\ifill f:0
\move(64 56)
\lvec(65 56)
\lvec(65 57)
\lvec(64 57)
\ifill f:0
\move(66 56)
\lvec(67 56)
\lvec(67 57)
\lvec(66 57)
\ifill f:0
\move(68 56)
\lvec(72 56)
\lvec(72 57)
\lvec(68 57)
\ifill f:0
\move(73 56)
\lvec(74 56)
\lvec(74 57)
\lvec(73 57)
\ifill f:0
\move(76 56)
\lvec(80 56)
\lvec(80 57)
\lvec(76 57)
\ifill f:0
\move(81 56)
\lvec(82 56)
\lvec(82 57)
\lvec(81 57)
\ifill f:0
\move(83 56)
\lvec(85 56)
\lvec(85 57)
\lvec(83 57)
\ifill f:0
\move(87 56)
\lvec(98 56)
\lvec(98 57)
\lvec(87 57)
\ifill f:0
\move(100 56)
\lvec(101 56)
\lvec(101 57)
\lvec(100 57)
\ifill f:0
\move(102 56)
\lvec(103 56)
\lvec(103 57)
\lvec(102 57)
\ifill f:0
\move(104 56)
\lvec(107 56)
\lvec(107 57)
\lvec(104 57)
\ifill f:0
\move(108 56)
\lvec(111 56)
\lvec(111 57)
\lvec(108 57)
\ifill f:0
\move(112 56)
\lvec(122 56)
\lvec(122 57)
\lvec(112 57)
\ifill f:0
\move(123 56)
\lvec(124 56)
\lvec(124 57)
\lvec(123 57)
\ifill f:0
\move(125 56)
\lvec(126 56)
\lvec(126 57)
\lvec(125 57)
\ifill f:0
\move(127 56)
\lvec(135 56)
\lvec(135 57)
\lvec(127 57)
\ifill f:0
\move(136 56)
\lvec(140 56)
\lvec(140 57)
\lvec(136 57)
\ifill f:0
\move(141 56)
\lvec(145 56)
\lvec(145 57)
\lvec(141 57)
\ifill f:0
\move(146 56)
\lvec(155 56)
\lvec(155 57)
\lvec(146 57)
\ifill f:0
\move(156 56)
\lvec(159 56)
\lvec(159 57)
\lvec(156 57)
\ifill f:0
\move(160 56)
\lvec(163 56)
\lvec(163 57)
\lvec(160 57)
\ifill f:0
\move(164 56)
\lvec(170 56)
\lvec(170 57)
\lvec(164 57)
\ifill f:0
\move(171 56)
\lvec(173 56)
\lvec(173 57)
\lvec(171 57)
\ifill f:0
\move(176 56)
\lvec(184 56)
\lvec(184 57)
\lvec(176 57)
\ifill f:0
\move(186 56)
\lvec(197 56)
\lvec(197 57)
\lvec(186 57)
\ifill f:0
\move(198 56)
\lvec(214 56)
\lvec(214 57)
\lvec(198 57)
\ifill f:0
\move(216 56)
\lvec(226 56)
\lvec(226 57)
\lvec(216 57)
\ifill f:0
\move(227 56)
\lvec(233 56)
\lvec(233 57)
\lvec(227 57)
\ifill f:0
\move(235 56)
\lvec(239 56)
\lvec(239 57)
\lvec(235 57)
\ifill f:0
\move(240 56)
\lvec(250 56)
\lvec(250 57)
\lvec(240 57)
\ifill f:0
\move(251 56)
\lvec(254 56)
\lvec(254 57)
\lvec(251 57)
\ifill f:0
\move(255 56)
\lvec(257 56)
\lvec(257 57)
\lvec(255 57)
\ifill f:0
\move(259 56)
\lvec(262 56)
\lvec(262 57)
\lvec(259 57)
\ifill f:0
\move(263 56)
\lvec(266 56)
\lvec(266 57)
\lvec(263 57)
\ifill f:0
\move(267 56)
\lvec(273 56)
\lvec(273 57)
\lvec(267 57)
\ifill f:0
\move(274 56)
\lvec(280 56)
\lvec(280 57)
\lvec(274 57)
\ifill f:0
\move(281 56)
\lvec(283 56)
\lvec(283 57)
\lvec(281 57)
\ifill f:0
\move(284 56)
\lvec(286 56)
\lvec(286 57)
\lvec(284 57)
\ifill f:0
\move(287 56)
\lvec(290 56)
\lvec(290 57)
\lvec(287 57)
\ifill f:0
\move(291 56)
\lvec(292 56)
\lvec(292 57)
\lvec(291 57)
\ifill f:0
\move(293 56)
\lvec(295 56)
\lvec(295 57)
\lvec(293 57)
\ifill f:0
\move(296 56)
\lvec(298 56)
\lvec(298 57)
\lvec(296 57)
\ifill f:0
\move(299 56)
\lvec(303 56)
\lvec(303 57)
\lvec(299 57)
\ifill f:0
\move(304 56)
\lvec(311 56)
\lvec(311 57)
\lvec(304 57)
\ifill f:0
\move(312 56)
\lvec(323 56)
\lvec(323 57)
\lvec(312 57)
\ifill f:0
\move(324 56)
\lvec(325 56)
\lvec(325 57)
\lvec(324 57)
\ifill f:0
\move(326 56)
\lvec(347 56)
\lvec(347 57)
\lvec(326 57)
\ifill f:0
\move(348 56)
\lvec(349 56)
\lvec(349 57)
\lvec(348 57)
\ifill f:0
\move(350 56)
\lvec(351 56)
\lvec(351 57)
\lvec(350 57)
\ifill f:0
\move(352 56)
\lvec(353 56)
\lvec(353 57)
\lvec(352 57)
\ifill f:0
\move(354 56)
\lvec(357 56)
\lvec(357 57)
\lvec(354 57)
\ifill f:0
\move(358 56)
\lvec(359 56)
\lvec(359 57)
\lvec(358 57)
\ifill f:0
\move(360 56)
\lvec(362 56)
\lvec(362 57)
\lvec(360 57)
\ifill f:0
\move(364 56)
\lvec(365 56)
\lvec(365 57)
\lvec(364 57)
\ifill f:0
\move(366 56)
\lvec(378 56)
\lvec(378 57)
\lvec(366 57)
\ifill f:0
\move(379 56)
\lvec(387 56)
\lvec(387 57)
\lvec(379 57)
\ifill f:0
\move(388 56)
\lvec(394 56)
\lvec(394 57)
\lvec(388 57)
\ifill f:0
\move(395 56)
\lvec(401 56)
\lvec(401 57)
\lvec(395 57)
\ifill f:0
\move(402 56)
\lvec(419 56)
\lvec(419 57)
\lvec(402 57)
\ifill f:0
\move(420 56)
\lvec(427 56)
\lvec(427 57)
\lvec(420 57)
\ifill f:0
\move(428 56)
\lvec(442 56)
\lvec(442 57)
\lvec(428 57)
\ifill f:0
\move(443 56)
\lvec(444 56)
\lvec(444 57)
\lvec(443 57)
\ifill f:0
\move(445 56)
\lvec(447 56)
\lvec(447 57)
\lvec(445 57)
\ifill f:0
\move(448 56)
\lvec(450 56)
\lvec(450 57)
\lvec(448 57)
\ifill f:0
\move(15 57)
\lvec(17 57)
\lvec(17 58)
\lvec(15 58)
\ifill f:0
\move(18 57)
\lvec(21 57)
\lvec(21 58)
\lvec(18 58)
\ifill f:0
\move(23 57)
\lvec(26 57)
\lvec(26 58)
\lvec(23 58)
\ifill f:0
\move(36 57)
\lvec(37 57)
\lvec(37 58)
\lvec(36 58)
\ifill f:0
\move(40 57)
\lvec(46 57)
\lvec(46 58)
\lvec(40 58)
\ifill f:0
\move(47 57)
\lvec(48 57)
\lvec(48 58)
\lvec(47 58)
\ifill f:0
\move(49 57)
\lvec(50 57)
\lvec(50 58)
\lvec(49 58)
\ifill f:0
\move(51 57)
\lvec(53 57)
\lvec(53 58)
\lvec(51 58)
\ifill f:0
\move(54 57)
\lvec(55 57)
\lvec(55 58)
\lvec(54 58)
\ifill f:0
\move(57 57)
\lvec(60 57)
\lvec(60 58)
\lvec(57 58)
\ifill f:0
\move(61 57)
\lvec(63 57)
\lvec(63 58)
\lvec(61 58)
\ifill f:0
\move(64 57)
\lvec(65 57)
\lvec(65 58)
\lvec(64 58)
\ifill f:0
\move(66 57)
\lvec(71 57)
\lvec(71 58)
\lvec(66 58)
\ifill f:0
\move(72 57)
\lvec(80 57)
\lvec(80 58)
\lvec(72 58)
\ifill f:0
\move(81 57)
\lvec(82 57)
\lvec(82 58)
\lvec(81 58)
\ifill f:0
\move(83 57)
\lvec(85 57)
\lvec(85 58)
\lvec(83 58)
\ifill f:0
\move(92 57)
\lvec(93 57)
\lvec(93 58)
\lvec(92 58)
\ifill f:0
\move(97 57)
\lvec(98 57)
\lvec(98 58)
\lvec(97 58)
\ifill f:0
\move(100 57)
\lvec(101 57)
\lvec(101 58)
\lvec(100 58)
\ifill f:0
\move(102 57)
\lvec(109 57)
\lvec(109 58)
\lvec(102 58)
\ifill f:0
\move(110 57)
\lvec(122 57)
\lvec(122 58)
\lvec(110 58)
\ifill f:0
\move(123 57)
\lvec(124 57)
\lvec(124 58)
\lvec(123 58)
\ifill f:0
\move(125 57)
\lvec(126 57)
\lvec(126 58)
\lvec(125 58)
\ifill f:0
\move(127 57)
\lvec(131 57)
\lvec(131 58)
\lvec(127 58)
\ifill f:0
\move(132 57)
\lvec(138 57)
\lvec(138 58)
\lvec(132 58)
\ifill f:0
\move(139 57)
\lvec(145 57)
\lvec(145 58)
\lvec(139 58)
\ifill f:0
\move(146 57)
\lvec(148 57)
\lvec(148 58)
\lvec(146 58)
\ifill f:0
\move(149 57)
\lvec(157 57)
\lvec(157 58)
\lvec(149 58)
\ifill f:0
\move(158 57)
\lvec(170 57)
\lvec(170 58)
\lvec(158 58)
\ifill f:0
\move(171 57)
\lvec(173 57)
\lvec(173 58)
\lvec(171 58)
\ifill f:0
\move(175 57)
\lvec(179 57)
\lvec(179 58)
\lvec(175 58)
\ifill f:0
\move(181 57)
\lvec(187 57)
\lvec(187 58)
\lvec(181 58)
\ifill f:0
\move(188 57)
\lvec(189 57)
\lvec(189 58)
\lvec(188 58)
\ifill f:0
\move(190 57)
\lvec(197 57)
\lvec(197 58)
\lvec(190 58)
\ifill f:0
\move(198 57)
\lvec(226 57)
\lvec(226 58)
\lvec(198 58)
\ifill f:0
\move(228 57)
\lvec(235 57)
\lvec(235 58)
\lvec(228 58)
\ifill f:0
\move(236 57)
\lvec(242 57)
\lvec(242 58)
\lvec(236 58)
\ifill f:0
\move(243 57)
\lvec(248 57)
\lvec(248 58)
\lvec(243 58)
\ifill f:0
\move(249 57)
\lvec(257 57)
\lvec(257 58)
\lvec(249 58)
\ifill f:0
\move(258 57)
\lvec(263 57)
\lvec(263 58)
\lvec(258 58)
\ifill f:0
\move(264 57)
\lvec(275 57)
\lvec(275 58)
\lvec(264 58)
\ifill f:0
\move(276 57)
\lvec(279 57)
\lvec(279 58)
\lvec(276 58)
\ifill f:0
\move(280 57)
\lvec(290 57)
\lvec(290 58)
\lvec(280 58)
\ifill f:0
\move(291 57)
\lvec(293 57)
\lvec(293 58)
\lvec(291 58)
\ifill f:0
\move(294 57)
\lvec(307 57)
\lvec(307 58)
\lvec(294 58)
\ifill f:0
\move(308 57)
\lvec(310 57)
\lvec(310 58)
\lvec(308 58)
\ifill f:0
\move(311 57)
\lvec(315 57)
\lvec(315 58)
\lvec(311 58)
\ifill f:0
\move(316 57)
\lvec(323 57)
\lvec(323 58)
\lvec(316 58)
\ifill f:0
\move(324 57)
\lvec(325 57)
\lvec(325 58)
\lvec(324 58)
\ifill f:0
\move(326 57)
\lvec(335 57)
\lvec(335 58)
\lvec(326 58)
\ifill f:0
\move(336 57)
\lvec(344 57)
\lvec(344 58)
\lvec(336 58)
\ifill f:0
\move(345 57)
\lvec(359 57)
\lvec(359 58)
\lvec(345 58)
\ifill f:0
\move(360 57)
\lvec(362 57)
\lvec(362 58)
\lvec(360 58)
\ifill f:0
\move(364 57)
\lvec(365 57)
\lvec(365 58)
\lvec(364 58)
\ifill f:0
\move(366 57)
\lvec(369 57)
\lvec(369 58)
\lvec(366 58)
\ifill f:0
\move(370 57)
\lvec(371 57)
\lvec(371 58)
\lvec(370 58)
\ifill f:0
\move(372 57)
\lvec(373 57)
\lvec(373 58)
\lvec(372 58)
\ifill f:0
\move(374 57)
\lvec(375 57)
\lvec(375 58)
\lvec(374 58)
\ifill f:0
\move(376 57)
\lvec(377 57)
\lvec(377 58)
\lvec(376 58)
\ifill f:0
\move(378 57)
\lvec(379 57)
\lvec(379 58)
\lvec(378 58)
\ifill f:0
\move(380 57)
\lvec(390 57)
\lvec(390 58)
\lvec(380 58)
\ifill f:0
\move(391 57)
\lvec(392 57)
\lvec(392 58)
\lvec(391 58)
\ifill f:0
\move(393 57)
\lvec(401 57)
\lvec(401 58)
\lvec(393 58)
\ifill f:0
\move(402 57)
\lvec(413 57)
\lvec(413 58)
\lvec(402 58)
\ifill f:0
\move(414 57)
\lvec(433 57)
\lvec(433 58)
\lvec(414 58)
\ifill f:0
\move(434 57)
\lvec(442 57)
\lvec(442 58)
\lvec(434 58)
\ifill f:0
\move(443 57)
\lvec(449 57)
\lvec(449 58)
\lvec(443 58)
\ifill f:0
\move(450 57)
\lvec(451 57)
\lvec(451 58)
\lvec(450 58)
\ifill f:0
\move(15 58)
\lvec(17 58)
\lvec(17 59)
\lvec(15 59)
\ifill f:0
\move(20 58)
\lvec(21 58)
\lvec(21 59)
\lvec(20 59)
\ifill f:0
\move(23 58)
\lvec(26 58)
\lvec(26 59)
\lvec(23 59)
\ifill f:0
\move(28 58)
\lvec(29 58)
\lvec(29 59)
\lvec(28 59)
\ifill f:0
\move(36 58)
\lvec(37 58)
\lvec(37 59)
\lvec(36 59)
\ifill f:0
\move(38 58)
\lvec(39 58)
\lvec(39 59)
\lvec(38 59)
\ifill f:0
\move(40 58)
\lvec(45 58)
\lvec(45 59)
\lvec(40 59)
\ifill f:0
\move(46 58)
\lvec(48 58)
\lvec(48 59)
\lvec(46 59)
\ifill f:0
\move(49 58)
\lvec(50 58)
\lvec(50 59)
\lvec(49 59)
\ifill f:0
\move(51 58)
\lvec(53 58)
\lvec(53 59)
\lvec(51 59)
\ifill f:0
\move(54 58)
\lvec(55 58)
\lvec(55 59)
\lvec(54 59)
\ifill f:0
\move(56 58)
\lvec(59 58)
\lvec(59 59)
\lvec(56 59)
\ifill f:0
\move(60 58)
\lvec(63 58)
\lvec(63 59)
\lvec(60 59)
\ifill f:0
\move(64 58)
\lvec(65 58)
\lvec(65 59)
\lvec(64 59)
\ifill f:0
\move(66 58)
\lvec(70 58)
\lvec(70 59)
\lvec(66 59)
\ifill f:0
\move(71 58)
\lvec(74 58)
\lvec(74 59)
\lvec(71 59)
\ifill f:0
\move(75 58)
\lvec(77 58)
\lvec(77 59)
\lvec(75 59)
\ifill f:0
\move(78 58)
\lvec(82 58)
\lvec(82 59)
\lvec(78 59)
\ifill f:0
\move(83 58)
\lvec(84 58)
\lvec(84 59)
\lvec(83 59)
\ifill f:0
\move(87 58)
\lvec(93 58)
\lvec(93 59)
\lvec(87 59)
\ifill f:0
\move(97 58)
\lvec(98 58)
\lvec(98 59)
\lvec(97 59)
\ifill f:0
\move(100 58)
\lvec(101 58)
\lvec(101 59)
\lvec(100 59)
\ifill f:0
\move(102 58)
\lvec(105 58)
\lvec(105 59)
\lvec(102 59)
\ifill f:0
\move(107 58)
\lvec(111 58)
\lvec(111 59)
\lvec(107 59)
\ifill f:0
\move(112 58)
\lvec(118 58)
\lvec(118 59)
\lvec(112 59)
\ifill f:0
\move(119 58)
\lvec(122 58)
\lvec(122 59)
\lvec(119 59)
\ifill f:0
\move(123 58)
\lvec(124 58)
\lvec(124 59)
\lvec(123 59)
\ifill f:0
\move(125 58)
\lvec(127 58)
\lvec(127 59)
\lvec(125 59)
\ifill f:0
\move(128 58)
\lvec(129 58)
\lvec(129 59)
\lvec(128 59)
\ifill f:0
\move(130 58)
\lvec(134 58)
\lvec(134 59)
\lvec(130 59)
\ifill f:0
\move(135 58)
\lvec(138 58)
\lvec(138 59)
\lvec(135 59)
\ifill f:0
\move(139 58)
\lvec(145 58)
\lvec(145 59)
\lvec(139 59)
\ifill f:0
\move(146 58)
\lvec(150 58)
\lvec(150 59)
\lvec(146 59)
\ifill f:0
\move(151 58)
\lvec(155 58)
\lvec(155 59)
\lvec(151 59)
\ifill f:0
\move(156 58)
\lvec(170 58)
\lvec(170 59)
\lvec(156 59)
\ifill f:0
\move(171 58)
\lvec(172 58)
\lvec(172 59)
\lvec(171 59)
\ifill f:0
\move(174 58)
\lvec(191 58)
\lvec(191 59)
\lvec(174 59)
\ifill f:0
\move(192 58)
\lvec(197 58)
\lvec(197 59)
\lvec(192 59)
\ifill f:0
\move(198 58)
\lvec(201 58)
\lvec(201 59)
\lvec(198 59)
\ifill f:0
\move(202 58)
\lvec(203 58)
\lvec(203 59)
\lvec(202 59)
\ifill f:0
\move(204 58)
\lvec(226 58)
\lvec(226 59)
\lvec(204 59)
\ifill f:0
\move(230 58)
\lvec(238 58)
\lvec(238 59)
\lvec(230 59)
\ifill f:0
\move(239 58)
\lvec(247 58)
\lvec(247 59)
\lvec(239 59)
\ifill f:0
\move(248 58)
\lvec(253 58)
\lvec(253 59)
\lvec(248 59)
\ifill f:0
\move(254 58)
\lvec(257 58)
\lvec(257 59)
\lvec(254 59)
\ifill f:0
\move(258 58)
\lvec(259 58)
\lvec(259 59)
\lvec(258 59)
\ifill f:0
\move(260 58)
\lvec(264 58)
\lvec(264 59)
\lvec(260 59)
\ifill f:0
\move(265 58)
\lvec(269 58)
\lvec(269 59)
\lvec(265 59)
\ifill f:0
\move(270 58)
\lvec(290 58)
\lvec(290 59)
\lvec(270 59)
\ifill f:0
\move(291 58)
\lvec(293 58)
\lvec(293 59)
\lvec(291 59)
\ifill f:0
\move(294 58)
\lvec(296 58)
\lvec(296 59)
\lvec(294 59)
\ifill f:0
\move(297 58)
\lvec(320 58)
\lvec(320 59)
\lvec(297 59)
\ifill f:0
\move(321 58)
\lvec(323 58)
\lvec(323 59)
\lvec(321 59)
\ifill f:0
\move(324 58)
\lvec(325 58)
\lvec(325 59)
\lvec(324 59)
\ifill f:0
\move(326 58)
\lvec(328 58)
\lvec(328 59)
\lvec(326 59)
\ifill f:0
\move(329 58)
\lvec(333 58)
\lvec(333 59)
\lvec(329 59)
\ifill f:0
\move(334 58)
\lvec(338 58)
\lvec(338 59)
\lvec(334 59)
\ifill f:0
\move(339 58)
\lvec(359 58)
\lvec(359 59)
\lvec(339 59)
\ifill f:0
\move(360 58)
\lvec(362 58)
\lvec(362 59)
\lvec(360 59)
\ifill f:0
\move(363 58)
\lvec(399 58)
\lvec(399 59)
\lvec(363 59)
\ifill f:0
\move(400 58)
\lvec(401 58)
\lvec(401 59)
\lvec(400 59)
\ifill f:0
\move(402 58)
\lvec(410 58)
\lvec(410 59)
\lvec(402 59)
\ifill f:0
\move(411 58)
\lvec(419 58)
\lvec(419 59)
\lvec(411 59)
\ifill f:0
\move(420 58)
\lvec(431 58)
\lvec(431 59)
\lvec(420 59)
\ifill f:0
\move(432 58)
\lvec(436 58)
\lvec(436 59)
\lvec(432 59)
\ifill f:0
\move(437 58)
\lvec(442 58)
\lvec(442 59)
\lvec(437 59)
\ifill f:0
\move(443 58)
\lvec(446 58)
\lvec(446 59)
\lvec(443 59)
\ifill f:0
\move(447 58)
\lvec(451 58)
\lvec(451 59)
\lvec(447 59)
\ifill f:0
\move(14 59)
\lvec(15 59)
\lvec(15 60)
\lvec(14 60)
\ifill f:0
\move(16 59)
\lvec(17 59)
\lvec(17 60)
\lvec(16 60)
\ifill f:0
\move(20 59)
\lvec(21 59)
\lvec(21 60)
\lvec(20 60)
\ifill f:0
\move(25 59)
\lvec(26 59)
\lvec(26 60)
\lvec(25 60)
\ifill f:0
\move(36 59)
\lvec(37 59)
\lvec(37 60)
\lvec(36 60)
\ifill f:0
\move(38 59)
\lvec(45 59)
\lvec(45 60)
\lvec(38 60)
\ifill f:0
\move(46 59)
\lvec(50 59)
\lvec(50 60)
\lvec(46 60)
\ifill f:0
\move(51 59)
\lvec(52 59)
\lvec(52 60)
\lvec(51 60)
\ifill f:0
\move(59 59)
\lvec(62 59)
\lvec(62 60)
\lvec(59 60)
\ifill f:0
\move(64 59)
\lvec(65 59)
\lvec(65 60)
\lvec(64 60)
\ifill f:0
\move(66 59)
\lvec(71 59)
\lvec(71 60)
\lvec(66 60)
\ifill f:0
\move(72 59)
\lvec(73 59)
\lvec(73 60)
\lvec(72 60)
\ifill f:0
\move(74 59)
\lvec(75 59)
\lvec(75 60)
\lvec(74 60)
\ifill f:0
\move(76 59)
\lvec(78 59)
\lvec(78 60)
\lvec(76 60)
\ifill f:0
\move(79 59)
\lvec(82 59)
\lvec(82 60)
\lvec(79 60)
\ifill f:0
\move(83 59)
\lvec(84 59)
\lvec(84 60)
\lvec(83 60)
\ifill f:0
\move(85 59)
\lvec(89 59)
\lvec(89 60)
\lvec(85 60)
\ifill f:0
\move(91 59)
\lvec(94 59)
\lvec(94 60)
\lvec(91 60)
\ifill f:0
\move(95 59)
\lvec(98 59)
\lvec(98 60)
\lvec(95 60)
\ifill f:0
\move(99 59)
\lvec(101 59)
\lvec(101 60)
\lvec(99 60)
\ifill f:0
\move(102 59)
\lvec(107 59)
\lvec(107 60)
\lvec(102 60)
\ifill f:0
\move(108 59)
\lvec(113 59)
\lvec(113 60)
\lvec(108 60)
\ifill f:0
\move(114 59)
\lvec(122 59)
\lvec(122 60)
\lvec(114 60)
\ifill f:0
\move(123 59)
\lvec(125 59)
\lvec(125 60)
\lvec(123 60)
\ifill f:0
\move(126 59)
\lvec(130 59)
\lvec(130 60)
\lvec(126 60)
\ifill f:0
\move(131 59)
\lvec(133 59)
\lvec(133 60)
\lvec(131 60)
\ifill f:0
\move(134 59)
\lvec(135 59)
\lvec(135 60)
\lvec(134 60)
\ifill f:0
\move(136 59)
\lvec(137 59)
\lvec(137 60)
\lvec(136 60)
\ifill f:0
\move(138 59)
\lvec(139 59)
\lvec(139 60)
\lvec(138 60)
\ifill f:0
\move(140 59)
\lvec(141 59)
\lvec(141 60)
\lvec(140 60)
\ifill f:0
\move(142 59)
\lvec(143 59)
\lvec(143 60)
\lvec(142 60)
\ifill f:0
\move(144 59)
\lvec(145 59)
\lvec(145 60)
\lvec(144 60)
\ifill f:0
\move(146 59)
\lvec(147 59)
\lvec(147 60)
\lvec(146 60)
\ifill f:0
\move(148 59)
\lvec(159 59)
\lvec(159 60)
\lvec(148 60)
\ifill f:0
\move(160 59)
\lvec(170 59)
\lvec(170 60)
\lvec(160 60)
\ifill f:0
\move(171 59)
\lvec(172 59)
\lvec(172 60)
\lvec(171 60)
\ifill f:0
\move(173 59)
\lvec(176 59)
\lvec(176 60)
\lvec(173 60)
\ifill f:0
\move(177 59)
\lvec(180 59)
\lvec(180 60)
\lvec(177 60)
\ifill f:0
\move(182 59)
\lvec(186 59)
\lvec(186 60)
\lvec(182 60)
\ifill f:0
\move(187 59)
\lvec(191 59)
\lvec(191 60)
\lvec(187 60)
\ifill f:0
\move(193 59)
\lvec(197 59)
\lvec(197 60)
\lvec(193 60)
\ifill f:0
\move(198 59)
\lvec(200 59)
\lvec(200 60)
\lvec(198 60)
\ifill f:0
\move(201 59)
\lvec(216 59)
\lvec(216 60)
\lvec(201 60)
\ifill f:0
\move(221 59)
\lvec(222 59)
\lvec(222 60)
\lvec(221 60)
\ifill f:0
\move(225 59)
\lvec(226 59)
\lvec(226 60)
\lvec(225 60)
\ifill f:0
\move(230 59)
\lvec(243 59)
\lvec(243 60)
\lvec(230 60)
\ifill f:0
\move(245 59)
\lvec(257 59)
\lvec(257 60)
\lvec(245 60)
\ifill f:0
\move(258 59)
\lvec(266 59)
\lvec(266 60)
\lvec(258 60)
\ifill f:0
\move(267 59)
\lvec(271 59)
\lvec(271 60)
\lvec(267 60)
\ifill f:0
\move(272 59)
\lvec(290 59)
\lvec(290 60)
\lvec(272 60)
\ifill f:0
\move(291 59)
\lvec(293 59)
\lvec(293 60)
\lvec(291 60)
\ifill f:0
\move(294 59)
\lvec(297 59)
\lvec(297 60)
\lvec(294 60)
\ifill f:0
\move(298 59)
\lvec(307 59)
\lvec(307 60)
\lvec(298 60)
\ifill f:0
\move(308 59)
\lvec(320 59)
\lvec(320 60)
\lvec(308 60)
\ifill f:0
\move(321 59)
\lvec(323 59)
\lvec(323 60)
\lvec(321 60)
\ifill f:0
\move(324 59)
\lvec(325 59)
\lvec(325 60)
\lvec(324 60)
\ifill f:0
\move(326 59)
\lvec(331 59)
\lvec(331 60)
\lvec(326 60)
\ifill f:0
\move(332 59)
\lvec(339 59)
\lvec(339 60)
\lvec(332 60)
\ifill f:0
\move(340 59)
\lvec(362 59)
\lvec(362 60)
\lvec(340 60)
\ifill f:0
\move(363 59)
\lvec(368 59)
\lvec(368 60)
\lvec(363 60)
\ifill f:0
\move(369 59)
\lvec(370 59)
\lvec(370 60)
\lvec(369 60)
\ifill f:0
\move(371 59)
\lvec(379 59)
\lvec(379 60)
\lvec(371 60)
\ifill f:0
\move(380 59)
\lvec(383 59)
\lvec(383 60)
\lvec(380 60)
\ifill f:0
\move(384 59)
\lvec(385 59)
\lvec(385 60)
\lvec(384 60)
\ifill f:0
\move(386 59)
\lvec(387 59)
\lvec(387 60)
\lvec(386 60)
\ifill f:0
\move(388 59)
\lvec(389 59)
\lvec(389 60)
\lvec(388 60)
\ifill f:0
\move(390 59)
\lvec(391 59)
\lvec(391 60)
\lvec(390 60)
\ifill f:0
\move(392 59)
\lvec(395 59)
\lvec(395 60)
\lvec(392 60)
\ifill f:0
\move(396 59)
\lvec(397 59)
\lvec(397 60)
\lvec(396 60)
\ifill f:0
\move(398 59)
\lvec(399 59)
\lvec(399 60)
\lvec(398 60)
\ifill f:0
\move(400 59)
\lvec(401 59)
\lvec(401 60)
\lvec(400 60)
\ifill f:0
\move(402 59)
\lvec(403 59)
\lvec(403 60)
\lvec(402 60)
\ifill f:0
\move(404 59)
\lvec(418 59)
\lvec(418 60)
\lvec(404 60)
\ifill f:0
\move(419 59)
\lvec(427 59)
\lvec(427 60)
\lvec(419 60)
\ifill f:0
\move(428 59)
\lvec(436 59)
\lvec(436 60)
\lvec(428 60)
\ifill f:0
\move(437 59)
\lvec(442 59)
\lvec(442 60)
\lvec(437 60)
\ifill f:0
\move(443 59)
\lvec(451 59)
\lvec(451 60)
\lvec(443 60)
\ifill f:0
\move(16 60)
\lvec(17 60)
\lvec(17 61)
\lvec(16 61)
\ifill f:0
\move(18 60)
\lvec(21 60)
\lvec(21 61)
\lvec(18 61)
\ifill f:0
\move(25 60)
\lvec(26 60)
\lvec(26 61)
\lvec(25 61)
\ifill f:0
\move(36 60)
\lvec(37 60)
\lvec(37 61)
\lvec(36 61)
\ifill f:0
\move(38 60)
\lvec(41 60)
\lvec(41 61)
\lvec(38 61)
\ifill f:0
\move(42 60)
\lvec(46 60)
\lvec(46 61)
\lvec(42 61)
\ifill f:0
\move(47 60)
\lvec(50 60)
\lvec(50 61)
\lvec(47 61)
\ifill f:0
\move(54 60)
\lvec(55 60)
\lvec(55 61)
\lvec(54 61)
\ifill f:0
\move(56 60)
\lvec(62 60)
\lvec(62 61)
\lvec(56 61)
\ifill f:0
\move(64 60)
\lvec(65 60)
\lvec(65 61)
\lvec(64 61)
\ifill f:0
\move(67 60)
\lvec(74 60)
\lvec(74 61)
\lvec(67 61)
\ifill f:0
\move(75 60)
\lvec(76 60)
\lvec(76 61)
\lvec(75 61)
\ifill f:0
\move(77 60)
\lvec(78 60)
\lvec(78 61)
\lvec(77 61)
\ifill f:0
\move(79 60)
\lvec(82 60)
\lvec(82 61)
\lvec(79 61)
\ifill f:0
\move(85 60)
\lvec(86 60)
\lvec(86 61)
\lvec(85 61)
\ifill f:0
\move(88 60)
\lvec(91 60)
\lvec(91 61)
\lvec(88 61)
\ifill f:0
\move(92 60)
\lvec(93 60)
\lvec(93 61)
\lvec(92 61)
\ifill f:0
\move(94 60)
\lvec(98 60)
\lvec(98 61)
\lvec(94 61)
\ifill f:0
\move(99 60)
\lvec(101 60)
\lvec(101 61)
\lvec(99 61)
\ifill f:0
\move(102 60)
\lvec(111 60)
\lvec(111 61)
\lvec(102 61)
\ifill f:0
\move(112 60)
\lvec(122 60)
\lvec(122 61)
\lvec(112 61)
\ifill f:0
\move(123 60)
\lvec(125 60)
\lvec(125 61)
\lvec(123 61)
\ifill f:0
\move(126 60)
\lvec(128 60)
\lvec(128 61)
\lvec(126 61)
\ifill f:0
\move(129 60)
\lvec(131 60)
\lvec(131 61)
\lvec(129 61)
\ifill f:0
\move(132 60)
\lvec(134 60)
\lvec(134 61)
\lvec(132 61)
\ifill f:0
\move(135 60)
\lvec(143 60)
\lvec(143 61)
\lvec(135 61)
\ifill f:0
\move(144 60)
\lvec(145 60)
\lvec(145 61)
\lvec(144 61)
\ifill f:0
\move(146 60)
\lvec(147 60)
\lvec(147 61)
\lvec(146 61)
\ifill f:0
\move(148 60)
\lvec(149 60)
\lvec(149 61)
\lvec(148 61)
\ifill f:0
\move(150 60)
\lvec(151 60)
\lvec(151 61)
\lvec(150 61)
\ifill f:0
\move(152 60)
\lvec(163 60)
\lvec(163 61)
\lvec(152 61)
\ifill f:0
\move(164 60)
\lvec(166 60)
\lvec(166 61)
\lvec(164 61)
\ifill f:0
\move(167 60)
\lvec(170 60)
\lvec(170 61)
\lvec(167 61)
\ifill f:0
\move(171 60)
\lvec(172 60)
\lvec(172 61)
\lvec(171 61)
\ifill f:0
\move(173 60)
\lvec(175 60)
\lvec(175 61)
\lvec(173 61)
\ifill f:0
\move(176 60)
\lvec(179 60)
\lvec(179 61)
\lvec(176 61)
\ifill f:0
\move(180 60)
\lvec(183 60)
\lvec(183 61)
\lvec(180 61)
\ifill f:0
\move(184 60)
\lvec(187 60)
\lvec(187 61)
\lvec(184 61)
\ifill f:0
\move(188 60)
\lvec(192 60)
\lvec(192 61)
\lvec(188 61)
\ifill f:0
\move(194 60)
\lvec(197 60)
\lvec(197 61)
\lvec(194 61)
\ifill f:0
\move(198 60)
\lvec(199 60)
\lvec(199 61)
\lvec(198 61)
\ifill f:0
\move(200 60)
\lvec(223 60)
\lvec(223 61)
\lvec(200 61)
\ifill f:0
\move(225 60)
\lvec(226 60)
\lvec(226 61)
\lvec(225 61)
\ifill f:0
\move(229 60)
\lvec(230 60)
\lvec(230 61)
\lvec(229 61)
\ifill f:0
\move(236 60)
\lvec(251 60)
\lvec(251 61)
\lvec(236 61)
\ifill f:0
\move(252 60)
\lvec(257 60)
\lvec(257 61)
\lvec(252 61)
\ifill f:0
\move(258 60)
\lvec(274 60)
\lvec(274 61)
\lvec(258 61)
\ifill f:0
\move(275 60)
\lvec(280 60)
\lvec(280 61)
\lvec(275 61)
\ifill f:0
\move(281 60)
\lvec(290 60)
\lvec(290 61)
\lvec(281 61)
\ifill f:0
\move(291 60)
\lvec(298 60)
\lvec(298 61)
\lvec(291 61)
\ifill f:0
\move(299 60)
\lvec(302 60)
\lvec(302 61)
\lvec(299 61)
\ifill f:0
\move(303 60)
\lvec(309 60)
\lvec(309 61)
\lvec(303 61)
\ifill f:0
\move(310 60)
\lvec(323 60)
\lvec(323 61)
\lvec(310 61)
\ifill f:0
\move(324 60)
\lvec(325 60)
\lvec(325 61)
\lvec(324 61)
\ifill f:0
\move(326 60)
\lvec(343 60)
\lvec(343 61)
\lvec(326 61)
\ifill f:0
\move(344 60)
\lvec(356 60)
\lvec(356 61)
\lvec(344 61)
\ifill f:0
\move(357 60)
\lvec(362 60)
\lvec(362 61)
\lvec(357 61)
\ifill f:0
\move(363 60)
\lvec(366 60)
\lvec(366 61)
\lvec(363 61)
\ifill f:0
\move(367 60)
\lvec(382 60)
\lvec(382 61)
\lvec(367 61)
\ifill f:0
\move(383 60)
\lvec(395 60)
\lvec(395 61)
\lvec(383 61)
\ifill f:0
\move(396 60)
\lvec(397 60)
\lvec(397 61)
\lvec(396 61)
\ifill f:0
\move(398 60)
\lvec(399 60)
\lvec(399 61)
\lvec(398 61)
\ifill f:0
\move(400 60)
\lvec(401 60)
\lvec(401 61)
\lvec(400 61)
\ifill f:0
\move(402 60)
\lvec(403 60)
\lvec(403 61)
\lvec(402 61)
\ifill f:0
\move(404 60)
\lvec(405 60)
\lvec(405 61)
\lvec(404 61)
\ifill f:0
\move(406 60)
\lvec(409 60)
\lvec(409 61)
\lvec(406 61)
\ifill f:0
\move(410 60)
\lvec(411 60)
\lvec(411 61)
\lvec(410 61)
\ifill f:0
\move(412 60)
\lvec(413 60)
\lvec(413 61)
\lvec(412 61)
\ifill f:0
\move(414 60)
\lvec(415 60)
\lvec(415 61)
\lvec(414 61)
\ifill f:0
\move(416 60)
\lvec(432 60)
\lvec(432 61)
\lvec(416 61)
\ifill f:0
\move(433 60)
\lvec(442 60)
\lvec(442 61)
\lvec(433 61)
\ifill f:0
\move(444 60)
\lvec(450 60)
\lvec(450 61)
\lvec(444 61)
\ifill f:0
\move(16 61)
\lvec(17 61)
\lvec(17 62)
\lvec(16 62)
\ifill f:0
\move(18 61)
\lvec(19 61)
\lvec(19 62)
\lvec(18 62)
\ifill f:0
\move(20 61)
\lvec(21 61)
\lvec(21 62)
\lvec(20 62)
\ifill f:0
\move(23 61)
\lvec(24 61)
\lvec(24 62)
\lvec(23 62)
\ifill f:0
\move(25 61)
\lvec(26 61)
\lvec(26 62)
\lvec(25 62)
\ifill f:0
\move(27 61)
\lvec(28 61)
\lvec(28 62)
\lvec(27 62)
\ifill f:0
\move(36 61)
\lvec(37 61)
\lvec(37 62)
\lvec(36 62)
\ifill f:0
\move(38 61)
\lvec(39 61)
\lvec(39 62)
\lvec(38 62)
\ifill f:0
\move(40 61)
\lvec(43 61)
\lvec(43 62)
\lvec(40 62)
\ifill f:0
\move(44 61)
\lvec(46 61)
\lvec(46 62)
\lvec(44 62)
\ifill f:0
\move(47 61)
\lvec(50 61)
\lvec(50 62)
\lvec(47 62)
\ifill f:0
\move(54 61)
\lvec(55 61)
\lvec(55 62)
\lvec(54 62)
\ifill f:0
\move(62 61)
\lvec(63 61)
\lvec(63 62)
\lvec(62 62)
\ifill f:0
\move(64 61)
\lvec(65 61)
\lvec(65 62)
\lvec(64 62)
\ifill f:0
\move(66 61)
\lvec(67 61)
\lvec(67 62)
\lvec(66 62)
\ifill f:0
\move(68 61)
\lvec(71 61)
\lvec(71 62)
\lvec(68 62)
\ifill f:0
\move(72 61)
\lvec(73 61)
\lvec(73 62)
\lvec(72 62)
\ifill f:0
\move(75 61)
\lvec(82 61)
\lvec(82 62)
\lvec(75 62)
\ifill f:0
\move(84 61)
\lvec(85 61)
\lvec(85 62)
\lvec(84 62)
\ifill f:0
\move(87 61)
\lvec(90 61)
\lvec(90 62)
\lvec(87 62)
\ifill f:0
\move(91 61)
\lvec(93 61)
\lvec(93 62)
\lvec(91 62)
\ifill f:0
\move(96 61)
\lvec(101 61)
\lvec(101 62)
\lvec(96 62)
\ifill f:0
\move(102 61)
\lvec(115 61)
\lvec(115 62)
\lvec(102 62)
\ifill f:0
\move(116 61)
\lvec(122 61)
\lvec(122 62)
\lvec(116 62)
\ifill f:0
\move(123 61)
\lvec(126 61)
\lvec(126 62)
\lvec(123 62)
\ifill f:0
\move(127 61)
\lvec(129 61)
\lvec(129 62)
\lvec(127 62)
\ifill f:0
\move(131 61)
\lvec(133 61)
\lvec(133 62)
\lvec(131 62)
\ifill f:0
\move(134 61)
\lvec(138 61)
\lvec(138 62)
\lvec(134 62)
\ifill f:0
\move(139 61)
\lvec(143 61)
\lvec(143 62)
\lvec(139 62)
\ifill f:0
\move(144 61)
\lvec(145 61)
\lvec(145 62)
\lvec(144 62)
\ifill f:0
\move(146 61)
\lvec(159 61)
\lvec(159 62)
\lvec(146 62)
\ifill f:0
\move(160 61)
\lvec(170 61)
\lvec(170 62)
\lvec(160 62)
\ifill f:0
\move(171 61)
\lvec(174 61)
\lvec(174 62)
\lvec(171 62)
\ifill f:0
\move(175 61)
\lvec(178 61)
\lvec(178 62)
\lvec(175 62)
\ifill f:0
\move(179 61)
\lvec(185 61)
\lvec(185 62)
\lvec(179 62)
\ifill f:0
\move(186 61)
\lvec(189 61)
\lvec(189 62)
\lvec(186 62)
\ifill f:0
\move(190 61)
\lvec(197 61)
\lvec(197 62)
\lvec(190 62)
\ifill f:0
\move(198 61)
\lvec(211 61)
\lvec(211 62)
\lvec(198 62)
\ifill f:0
\move(212 61)
\lvec(223 61)
\lvec(223 62)
\lvec(212 62)
\ifill f:0
\move(225 61)
\lvec(226 61)
\lvec(226 62)
\lvec(225 62)
\ifill f:0
\move(227 61)
\lvec(248 61)
\lvec(248 62)
\lvec(227 62)
\ifill f:0
\move(251 61)
\lvec(257 61)
\lvec(257 62)
\lvec(251 62)
\ifill f:0
\move(258 61)
\lvec(262 61)
\lvec(262 62)
\lvec(258 62)
\ifill f:0
\move(263 61)
\lvec(271 61)
\lvec(271 62)
\lvec(263 62)
\ifill f:0
\move(272 61)
\lvec(284 61)
\lvec(284 62)
\lvec(272 62)
\ifill f:0
\move(285 61)
\lvec(290 61)
\lvec(290 62)
\lvec(285 62)
\ifill f:0
\move(291 61)
\lvec(294 61)
\lvec(294 62)
\lvec(291 62)
\ifill f:0
\move(295 61)
\lvec(299 61)
\lvec(299 62)
\lvec(295 62)
\ifill f:0
\move(300 61)
\lvec(303 61)
\lvec(303 62)
\lvec(300 62)
\ifill f:0
\move(304 61)
\lvec(315 61)
\lvec(315 62)
\lvec(304 62)
\ifill f:0
\move(316 61)
\lvec(319 61)
\lvec(319 62)
\lvec(316 62)
\ifill f:0
\move(320 61)
\lvec(325 61)
\lvec(325 62)
\lvec(320 62)
\ifill f:0
\move(326 61)
\lvec(329 61)
\lvec(329 62)
\lvec(326 62)
\ifill f:0
\move(330 61)
\lvec(362 61)
\lvec(362 62)
\lvec(330 62)
\ifill f:0
\move(363 61)
\lvec(364 61)
\lvec(364 62)
\lvec(363 62)
\ifill f:0
\move(365 61)
\lvec(376 61)
\lvec(376 62)
\lvec(365 62)
\ifill f:0
\move(377 61)
\lvec(390 61)
\lvec(390 62)
\lvec(377 62)
\ifill f:0
\move(391 61)
\lvec(399 61)
\lvec(399 62)
\lvec(391 62)
\ifill f:0
\move(400 61)
\lvec(401 61)
\lvec(401 62)
\lvec(400 62)
\ifill f:0
\move(402 61)
\lvec(442 61)
\lvec(442 62)
\lvec(402 62)
\ifill f:0
\move(444 61)
\lvec(445 61)
\lvec(445 62)
\lvec(444 62)
\ifill f:0
\move(446 61)
\lvec(451 61)
\lvec(451 62)
\lvec(446 62)
\ifill f:0
\move(16 62)
\lvec(17 62)
\lvec(17 63)
\lvec(16 63)
\ifill f:0
\move(20 62)
\lvec(21 62)
\lvec(21 63)
\lvec(20 63)
\ifill f:0
\move(23 62)
\lvec(24 62)
\lvec(24 63)
\lvec(23 63)
\ifill f:0
\move(25 62)
\lvec(26 62)
\lvec(26 63)
\lvec(25 63)
\ifill f:0
\move(36 62)
\lvec(37 62)
\lvec(37 63)
\lvec(36 63)
\ifill f:0
\move(38 62)
\lvec(39 62)
\lvec(39 63)
\lvec(38 63)
\ifill f:0
\move(40 62)
\lvec(41 62)
\lvec(41 63)
\lvec(40 63)
\ifill f:0
\move(42 62)
\lvec(46 62)
\lvec(46 63)
\lvec(42 63)
\ifill f:0
\move(48 62)
\lvec(50 62)
\lvec(50 63)
\lvec(48 63)
\ifill f:0
\move(52 62)
\lvec(53 62)
\lvec(53 63)
\lvec(52 63)
\ifill f:0
\move(56 62)
\lvec(65 62)
\lvec(65 63)
\lvec(56 63)
\ifill f:0
\move(66 62)
\lvec(75 62)
\lvec(75 63)
\lvec(66 63)
\ifill f:0
\move(76 62)
\lvec(79 62)
\lvec(79 63)
\lvec(76 63)
\ifill f:0
\move(80 62)
\lvec(82 62)
\lvec(82 63)
\lvec(80 63)
\ifill f:0
\move(84 62)
\lvec(85 62)
\lvec(85 63)
\lvec(84 63)
\ifill f:0
\move(86 62)
\lvec(88 62)
\lvec(88 63)
\lvec(86 63)
\ifill f:0
\move(89 62)
\lvec(91 62)
\lvec(91 63)
\lvec(89 63)
\ifill f:0
\move(93 62)
\lvec(94 62)
\lvec(94 63)
\lvec(93 63)
\ifill f:0
\move(97 62)
\lvec(101 62)
\lvec(101 63)
\lvec(97 63)
\ifill f:0
\move(102 62)
\lvec(105 62)
\lvec(105 63)
\lvec(102 63)
\ifill f:0
\move(112 62)
\lvec(113 62)
\lvec(113 63)
\lvec(112 63)
\ifill f:0
\move(114 62)
\lvec(122 62)
\lvec(122 63)
\lvec(114 63)
\ifill f:0
\move(123 62)
\lvec(127 62)
\lvec(127 63)
\lvec(123 63)
\ifill f:0
\move(128 62)
\lvec(131 62)
\lvec(131 63)
\lvec(128 63)
\ifill f:0
\move(132 62)
\lvec(134 62)
\lvec(134 63)
\lvec(132 63)
\ifill f:0
\move(136 62)
\lvec(138 62)
\lvec(138 63)
\lvec(136 63)
\ifill f:0
\move(139 62)
\lvec(140 62)
\lvec(140 63)
\lvec(139 63)
\ifill f:0
\move(141 62)
\lvec(143 62)
\lvec(143 63)
\lvec(141 63)
\ifill f:0
\move(144 62)
\lvec(145 62)
\lvec(145 63)
\lvec(144 63)
\ifill f:0
\move(146 62)
\lvec(150 62)
\lvec(150 63)
\lvec(146 63)
\ifill f:0
\move(151 62)
\lvec(162 62)
\lvec(162 63)
\lvec(151 63)
\ifill f:0
\move(163 62)
\lvec(170 62)
\lvec(170 63)
\lvec(163 63)
\ifill f:0
\move(172 62)
\lvec(174 62)
\lvec(174 63)
\lvec(172 63)
\ifill f:0
\move(175 62)
\lvec(177 62)
\lvec(177 63)
\lvec(175 63)
\ifill f:0
\move(178 62)
\lvec(180 62)
\lvec(180 63)
\lvec(178 63)
\ifill f:0
\move(181 62)
\lvec(183 62)
\lvec(183 63)
\lvec(181 63)
\ifill f:0
\move(184 62)
\lvec(186 62)
\lvec(186 63)
\lvec(184 63)
\ifill f:0
\move(187 62)
\lvec(190 62)
\lvec(190 63)
\lvec(187 63)
\ifill f:0
\move(191 62)
\lvec(194 62)
\lvec(194 63)
\lvec(191 63)
\ifill f:0
\move(195 62)
\lvec(197 62)
\lvec(197 63)
\lvec(195 63)
\ifill f:0
\move(198 62)
\lvec(203 62)
\lvec(203 63)
\lvec(198 63)
\ifill f:0
\move(204 62)
\lvec(224 62)
\lvec(224 63)
\lvec(204 63)
\ifill f:0
\move(225 62)
\lvec(226 62)
\lvec(226 63)
\lvec(225 63)
\ifill f:0
\move(227 62)
\lvec(242 62)
\lvec(242 63)
\lvec(227 63)
\ifill f:0
\move(243 62)
\lvec(257 62)
\lvec(257 63)
\lvec(243 63)
\ifill f:0
\move(258 62)
\lvec(264 62)
\lvec(264 63)
\lvec(258 63)
\ifill f:0
\move(265 62)
\lvec(266 62)
\lvec(266 63)
\lvec(265 63)
\ifill f:0
\move(267 62)
\lvec(275 62)
\lvec(275 63)
\lvec(267 63)
\ifill f:0
\move(276 62)
\lvec(283 62)
\lvec(283 63)
\lvec(276 63)
\ifill f:0
\move(284 62)
\lvec(290 62)
\lvec(290 63)
\lvec(284 63)
\ifill f:0
\move(291 62)
\lvec(295 62)
\lvec(295 63)
\lvec(291 63)
\ifill f:0
\move(296 62)
\lvec(300 62)
\lvec(300 63)
\lvec(296 63)
\ifill f:0
\move(301 62)
\lvec(305 62)
\lvec(305 63)
\lvec(301 63)
\ifill f:0
\move(306 62)
\lvec(310 62)
\lvec(310 63)
\lvec(306 63)
\ifill f:0
\move(311 62)
\lvec(314 62)
\lvec(314 63)
\lvec(311 63)
\ifill f:0
\move(315 62)
\lvec(318 62)
\lvec(318 63)
\lvec(315 63)
\ifill f:0
\move(319 62)
\lvec(322 62)
\lvec(322 63)
\lvec(319 63)
\ifill f:0
\move(323 62)
\lvec(325 62)
\lvec(325 63)
\lvec(323 63)
\ifill f:0
\move(326 62)
\lvec(333 62)
\lvec(333 63)
\lvec(326 63)
\ifill f:0
\move(334 62)
\lvec(343 62)
\lvec(343 63)
\lvec(334 63)
\ifill f:0
\move(344 62)
\lvec(362 62)
\lvec(362 63)
\lvec(344 63)
\ifill f:0
\move(363 62)
\lvec(364 62)
\lvec(364 63)
\lvec(363 63)
\ifill f:0
\move(365 62)
\lvec(377 62)
\lvec(377 63)
\lvec(365 63)
\ifill f:0
\move(378 62)
\lvec(387 62)
\lvec(387 63)
\lvec(378 63)
\ifill f:0
\move(388 62)
\lvec(392 62)
\lvec(392 63)
\lvec(388 63)
\ifill f:0
\move(393 62)
\lvec(399 62)
\lvec(399 63)
\lvec(393 63)
\ifill f:0
\move(400 62)
\lvec(401 62)
\lvec(401 63)
\lvec(400 63)
\ifill f:0
\move(402 62)
\lvec(410 62)
\lvec(410 63)
\lvec(402 63)
\ifill f:0
\move(411 62)
\lvec(423 62)
\lvec(423 63)
\lvec(411 63)
\ifill f:0
\move(424 62)
\lvec(425 62)
\lvec(425 63)
\lvec(424 63)
\ifill f:0
\move(426 62)
\lvec(427 62)
\lvec(427 63)
\lvec(426 63)
\ifill f:0
\move(428 62)
\lvec(429 62)
\lvec(429 63)
\lvec(428 63)
\ifill f:0
\move(430 62)
\lvec(431 62)
\lvec(431 63)
\lvec(430 63)
\ifill f:0
\move(432 62)
\lvec(433 62)
\lvec(433 63)
\lvec(432 63)
\ifill f:0
\move(434 62)
\lvec(435 62)
\lvec(435 63)
\lvec(434 63)
\ifill f:0
\move(436 62)
\lvec(437 62)
\lvec(437 63)
\lvec(436 63)
\ifill f:0
\move(438 62)
\lvec(439 62)
\lvec(439 63)
\lvec(438 63)
\ifill f:0
\move(440 62)
\lvec(442 62)
\lvec(442 63)
\lvec(440 63)
\ifill f:0
\move(444 62)
\lvec(445 62)
\lvec(445 63)
\lvec(444 63)
\ifill f:0
\move(446 62)
\lvec(451 62)
\lvec(451 63)
\lvec(446 63)
\ifill f:0
\move(14 63)
\lvec(17 63)
\lvec(17 64)
\lvec(14 64)
\ifill f:0
\move(19 63)
\lvec(22 63)
\lvec(22 64)
\lvec(19 64)
\ifill f:0
\move(23 63)
\lvec(26 63)
\lvec(26 64)
\lvec(23 64)
\ifill f:0
\move(36 63)
\lvec(37 63)
\lvec(37 64)
\lvec(36 64)
\ifill f:0
\move(39 63)
\lvec(43 63)
\lvec(43 64)
\lvec(39 64)
\ifill f:0
\move(44 63)
\lvec(47 63)
\lvec(47 64)
\lvec(44 64)
\ifill f:0
\move(48 63)
\lvec(50 63)
\lvec(50 64)
\lvec(48 64)
\ifill f:0
\move(54 63)
\lvec(55 63)
\lvec(55 64)
\lvec(54 64)
\ifill f:0
\move(58 63)
\lvec(65 63)
\lvec(65 64)
\lvec(58 64)
\ifill f:0
\move(66 63)
\lvec(74 63)
\lvec(74 64)
\lvec(66 64)
\ifill f:0
\move(75 63)
\lvec(79 63)
\lvec(79 64)
\lvec(75 64)
\ifill f:0
\move(80 63)
\lvec(82 63)
\lvec(82 64)
\lvec(80 64)
\ifill f:0
\move(84 63)
\lvec(85 63)
\lvec(85 64)
\lvec(84 64)
\ifill f:0
\move(86 63)
\lvec(87 63)
\lvec(87 64)
\lvec(86 64)
\ifill f:0
\move(88 63)
\lvec(90 63)
\lvec(90 64)
\lvec(88 64)
\ifill f:0
\move(91 63)
\lvec(93 63)
\lvec(93 64)
\lvec(91 64)
\ifill f:0
\move(94 63)
\lvec(96 63)
\lvec(96 64)
\lvec(94 64)
\ifill f:0
\move(97 63)
\lvec(101 63)
\lvec(101 64)
\lvec(97 64)
\ifill f:0
\move(102 63)
\lvec(103 63)
\lvec(103 64)
\lvec(102 64)
\ifill f:0
\move(104 63)
\lvec(122 63)
\lvec(122 64)
\lvec(104 64)
\ifill f:0
\move(125 63)
\lvec(128 63)
\lvec(128 64)
\lvec(125 64)
\ifill f:0
\move(129 63)
\lvec(133 63)
\lvec(133 64)
\lvec(129 64)
\ifill f:0
\move(134 63)
\lvec(143 63)
\lvec(143 64)
\lvec(134 64)
\ifill f:0
\move(144 63)
\lvec(145 63)
\lvec(145 64)
\lvec(144 64)
\ifill f:0
\move(146 63)
\lvec(148 63)
\lvec(148 64)
\lvec(146 64)
\ifill f:0
\move(149 63)
\lvec(153 63)
\lvec(153 64)
\lvec(149 64)
\ifill f:0
\move(154 63)
\lvec(170 63)
\lvec(170 64)
\lvec(154 64)
\ifill f:0
\move(172 63)
\lvec(173 63)
\lvec(173 64)
\lvec(172 64)
\ifill f:0
\move(174 63)
\lvec(176 63)
\lvec(176 64)
\lvec(174 64)
\ifill f:0
\move(177 63)
\lvec(187 63)
\lvec(187 64)
\lvec(177 64)
\ifill f:0
\move(188 63)
\lvec(194 63)
\lvec(194 64)
\lvec(188 64)
\ifill f:0
\move(195 63)
\lvec(197 63)
\lvec(197 64)
\lvec(195 64)
\ifill f:0
\move(199 63)
\lvec(202 63)
\lvec(202 64)
\lvec(199 64)
\ifill f:0
\move(203 63)
\lvec(211 63)
\lvec(211 64)
\lvec(203 64)
\ifill f:0
\move(212 63)
\lvec(226 63)
\lvec(226 64)
\lvec(212 64)
\ifill f:0
\move(227 63)
\lvec(234 63)
\lvec(234 64)
\lvec(227 64)
\ifill f:0
\move(237 63)
\lvec(257 63)
\lvec(257 64)
\lvec(237 64)
\ifill f:0
\move(258 63)
\lvec(269 63)
\lvec(269 64)
\lvec(258 64)
\ifill f:0
\move(270 63)
\lvec(271 63)
\lvec(271 64)
\lvec(270 64)
\ifill f:0
\move(272 63)
\lvec(281 63)
\lvec(281 64)
\lvec(272 64)
\ifill f:0
\move(282 63)
\lvec(290 63)
\lvec(290 64)
\lvec(282 64)
\ifill f:0
\move(291 63)
\lvec(296 63)
\lvec(296 64)
\lvec(291 64)
\ifill f:0
\move(297 63)
\lvec(298 63)
\lvec(298 64)
\lvec(297 64)
\ifill f:0
\move(299 63)
\lvec(302 63)
\lvec(302 64)
\lvec(299 64)
\ifill f:0
\move(303 63)
\lvec(308 63)
\lvec(308 64)
\lvec(303 64)
\ifill f:0
\move(309 63)
\lvec(313 63)
\lvec(313 64)
\lvec(309 64)
\ifill f:0
\move(314 63)
\lvec(322 63)
\lvec(322 64)
\lvec(314 64)
\ifill f:0
\move(323 63)
\lvec(325 63)
\lvec(325 64)
\lvec(323 64)
\ifill f:0
\move(327 63)
\lvec(330 63)
\lvec(330 64)
\lvec(327 64)
\ifill f:0
\move(331 63)
\lvec(334 63)
\lvec(334 64)
\lvec(331 64)
\ifill f:0
\move(335 63)
\lvec(338 63)
\lvec(338 64)
\lvec(335 64)
\ifill f:0
\move(339 63)
\lvec(358 63)
\lvec(358 64)
\lvec(339 64)
\ifill f:0
\move(359 63)
\lvec(362 63)
\lvec(362 64)
\lvec(359 64)
\ifill f:0
\move(363 63)
\lvec(364 63)
\lvec(364 64)
\lvec(363 64)
\ifill f:0
\move(365 63)
\lvec(367 63)
\lvec(367 64)
\lvec(365 64)
\ifill f:0
\move(368 63)
\lvec(370 63)
\lvec(370 64)
\lvec(368 64)
\ifill f:0
\move(371 63)
\lvec(381 63)
\lvec(381 64)
\lvec(371 64)
\ifill f:0
\move(382 63)
\lvec(394 63)
\lvec(394 64)
\lvec(382 64)
\ifill f:0
\move(395 63)
\lvec(399 63)
\lvec(399 64)
\lvec(395 64)
\ifill f:0
\move(400 63)
\lvec(401 63)
\lvec(401 64)
\lvec(400 64)
\ifill f:0
\move(402 63)
\lvec(439 63)
\lvec(439 64)
\lvec(402 64)
\ifill f:0
\move(440 63)
\lvec(442 63)
\lvec(442 64)
\lvec(440 64)
\ifill f:0
\move(444 63)
\lvec(445 63)
\lvec(445 64)
\lvec(444 64)
\ifill f:0
\move(446 63)
\lvec(447 63)
\lvec(447 64)
\lvec(446 64)
\ifill f:0
\move(448 63)
\lvec(449 63)
\lvec(449 64)
\lvec(448 64)
\ifill f:0
\move(450 63)
\lvec(451 63)
\lvec(451 64)
\lvec(450 64)
\ifill f:0
\move(15 64)
\lvec(17 64)
\lvec(17 65)
\lvec(15 65)
\ifill f:0
\move(18 64)
\lvec(19 64)
\lvec(19 65)
\lvec(18 65)
\ifill f:0
\move(22 64)
\lvec(23 64)
\lvec(23 65)
\lvec(22 65)
\ifill f:0
\move(24 64)
\lvec(26 64)
\lvec(26 65)
\lvec(24 65)
\ifill f:0
\move(36 64)
\lvec(37 64)
\lvec(37 65)
\lvec(36 65)
\ifill f:0
\move(38 64)
\lvec(46 64)
\lvec(46 65)
\lvec(38 65)
\ifill f:0
\move(47 64)
\lvec(50 64)
\lvec(50 65)
\lvec(47 65)
\ifill f:0
\move(51 64)
\lvec(52 64)
\lvec(52 65)
\lvec(51 65)
\ifill f:0
\move(56 64)
\lvec(58 64)
\lvec(58 65)
\lvec(56 65)
\ifill f:0
\move(60 64)
\lvec(65 64)
\lvec(65 65)
\lvec(60 65)
\ifill f:0
\move(66 64)
\lvec(71 64)
\lvec(71 65)
\lvec(66 65)
\ifill f:0
\move(72 64)
\lvec(79 64)
\lvec(79 65)
\lvec(72 65)
\ifill f:0
\move(80 64)
\lvec(82 64)
\lvec(82 65)
\lvec(80 65)
\ifill f:0
\move(83 64)
\lvec(84 64)
\lvec(84 65)
\lvec(83 65)
\ifill f:0
\move(85 64)
\lvec(86 64)
\lvec(86 65)
\lvec(85 65)
\ifill f:0
\move(87 64)
\lvec(89 64)
\lvec(89 65)
\lvec(87 65)
\ifill f:0
\move(90 64)
\lvec(91 64)
\lvec(91 65)
\lvec(90 65)
\ifill f:0
\move(92 64)
\lvec(93 64)
\lvec(93 65)
\lvec(92 65)
\ifill f:0
\move(95 64)
\lvec(98 64)
\lvec(98 65)
\lvec(95 65)
\ifill f:0
\move(99 64)
\lvec(101 64)
\lvec(101 65)
\lvec(99 65)
\ifill f:0
\move(102 64)
\lvec(109 64)
\lvec(109 65)
\lvec(102 65)
\ifill f:0
\move(111 64)
\lvec(122 64)
\lvec(122 65)
\lvec(111 65)
\ifill f:0
\move(126 64)
\lvec(130 64)
\lvec(130 65)
\lvec(126 65)
\ifill f:0
\move(132 64)
\lvec(135 64)
\lvec(135 65)
\lvec(132 65)
\ifill f:0
\move(136 64)
\lvec(139 64)
\lvec(139 65)
\lvec(136 65)
\ifill f:0
\move(140 64)
\lvec(143 64)
\lvec(143 65)
\lvec(140 65)
\ifill f:0
\move(144 64)
\lvec(145 64)
\lvec(145 65)
\lvec(144 65)
\ifill f:0
\move(146 64)
\lvec(149 64)
\lvec(149 65)
\lvec(146 65)
\ifill f:0
\move(150 64)
\lvec(154 64)
\lvec(154 65)
\lvec(150 65)
\ifill f:0
\move(155 64)
\lvec(159 64)
\lvec(159 65)
\lvec(155 65)
\ifill f:0
\move(160 64)
\lvec(161 64)
\lvec(161 65)
\lvec(160 65)
\ifill f:0
\move(162 64)
\lvec(163 64)
\lvec(163 65)
\lvec(162 65)
\ifill f:0
\move(164 64)
\lvec(165 64)
\lvec(165 65)
\lvec(164 65)
\ifill f:0
\move(166 64)
\lvec(167 64)
\lvec(167 65)
\lvec(166 65)
\ifill f:0
\move(168 64)
\lvec(170 64)
\lvec(170 65)
\lvec(168 65)
\ifill f:0
\move(172 64)
\lvec(173 64)
\lvec(173 65)
\lvec(172 65)
\ifill f:0
\move(174 64)
\lvec(180 64)
\lvec(180 65)
\lvec(174 65)
\ifill f:0
\move(181 64)
\lvec(191 64)
\lvec(191 65)
\lvec(181 65)
\ifill f:0
\move(192 64)
\lvec(194 64)
\lvec(194 65)
\lvec(192 65)
\ifill f:0
\move(195 64)
\lvec(197 64)
\lvec(197 65)
\lvec(195 65)
\ifill f:0
\move(198 64)
\lvec(201 64)
\lvec(201 65)
\lvec(198 65)
\ifill f:0
\move(202 64)
\lvec(219 64)
\lvec(219 65)
\lvec(202 65)
\ifill f:0
\move(220 64)
\lvec(226 64)
\lvec(226 65)
\lvec(220 65)
\ifill f:0
\move(227 64)
\lvec(232 64)
\lvec(232 65)
\lvec(227 65)
\ifill f:0
\move(233 64)
\lvec(242 64)
\lvec(242 65)
\lvec(233 65)
\ifill f:0
\move(245 64)
\lvec(257 64)
\lvec(257 65)
\lvec(245 65)
\ifill f:0
\move(258 64)
\lvec(277 64)
\lvec(277 65)
\lvec(258 65)
\ifill f:0
\move(278 64)
\lvec(279 64)
\lvec(279 65)
\lvec(278 65)
\ifill f:0
\move(281 64)
\lvec(290 64)
\lvec(290 65)
\lvec(281 65)
\ifill f:0
\move(292 64)
\lvec(298 64)
\lvec(298 65)
\lvec(292 65)
\ifill f:0
\move(300 64)
\lvec(305 64)
\lvec(305 65)
\lvec(300 65)
\ifill f:0
\move(306 64)
\lvec(311 64)
\lvec(311 65)
\lvec(306 65)
\ifill f:0
\move(312 64)
\lvec(322 64)
\lvec(322 65)
\lvec(312 65)
\ifill f:0
\move(323 64)
\lvec(325 64)
\lvec(325 65)
\lvec(323 65)
\ifill f:0
\move(326 64)
\lvec(347 64)
\lvec(347 65)
\lvec(326 65)
\ifill f:0
\move(348 64)
\lvec(362 64)
\lvec(362 65)
\lvec(348 65)
\ifill f:0
\move(363 64)
\lvec(388 64)
\lvec(388 65)
\lvec(363 65)
\ifill f:0
\move(389 64)
\lvec(391 64)
\lvec(391 65)
\lvec(389 65)
\ifill f:0
\move(392 64)
\lvec(396 64)
\lvec(396 65)
\lvec(392 65)
\ifill f:0
\move(397 64)
\lvec(401 64)
\lvec(401 65)
\lvec(397 65)
\ifill f:0
\move(402 64)
\lvec(430 64)
\lvec(430 65)
\lvec(402 65)
\ifill f:0
\move(431 64)
\lvec(439 64)
\lvec(439 65)
\lvec(431 65)
\ifill f:0
\move(440 64)
\lvec(442 64)
\lvec(442 65)
\lvec(440 65)
\ifill f:0
\move(443 64)
\lvec(451 64)
\lvec(451 65)
\lvec(443 65)
\ifill f:0
\move(15 65)
\lvec(17 65)
\lvec(17 66)
\lvec(15 66)
\ifill f:0
\move(18 65)
\lvec(21 65)
\lvec(21 66)
\lvec(18 66)
\ifill f:0
\move(24 65)
\lvec(26 65)
\lvec(26 66)
\lvec(24 66)
\ifill f:0
\move(36 65)
\lvec(37 65)
\lvec(37 66)
\lvec(36 66)
\ifill f:0
\move(38 65)
\lvec(39 65)
\lvec(39 66)
\lvec(38 66)
\ifill f:0
\move(40 65)
\lvec(41 65)
\lvec(41 66)
\lvec(40 66)
\ifill f:0
\move(43 65)
\lvec(45 65)
\lvec(45 66)
\lvec(43 66)
\ifill f:0
\move(47 65)
\lvec(50 65)
\lvec(50 66)
\lvec(47 66)
\ifill f:0
\move(51 65)
\lvec(52 65)
\lvec(52 66)
\lvec(51 66)
\ifill f:0
\move(57 65)
\lvec(59 65)
\lvec(59 66)
\lvec(57 66)
\ifill f:0
\move(61 65)
\lvec(65 65)
\lvec(65 66)
\lvec(61 66)
\ifill f:0
\move(66 65)
\lvec(74 65)
\lvec(74 66)
\lvec(66 66)
\ifill f:0
\move(76 65)
\lvec(78 65)
\lvec(78 66)
\lvec(76 66)
\ifill f:0
\move(79 65)
\lvec(82 65)
\lvec(82 66)
\lvec(79 66)
\ifill f:0
\move(83 65)
\lvec(84 65)
\lvec(84 66)
\lvec(83 66)
\ifill f:0
\move(87 65)
\lvec(88 65)
\lvec(88 66)
\lvec(87 66)
\ifill f:0
\move(89 65)
\lvec(90 65)
\lvec(90 66)
\lvec(89 66)
\ifill f:0
\move(91 65)
\lvec(94 65)
\lvec(94 66)
\lvec(91 66)
\ifill f:0
\move(96 65)
\lvec(98 65)
\lvec(98 66)
\lvec(96 66)
\ifill f:0
\move(99 65)
\lvec(101 65)
\lvec(101 66)
\lvec(99 66)
\ifill f:0
\move(103 65)
\lvec(106 65)
\lvec(106 66)
\lvec(103 66)
\ifill f:0
\move(107 65)
\lvec(115 65)
\lvec(115 66)
\lvec(107 66)
\ifill f:0
\move(119 65)
\lvec(120 65)
\lvec(120 66)
\lvec(119 66)
\ifill f:0
\move(121 65)
\lvec(122 65)
\lvec(122 66)
\lvec(121 66)
\ifill f:0
\move(127 65)
\lvec(133 65)
\lvec(133 66)
\lvec(127 66)
\ifill f:0
\move(135 65)
\lvec(138 65)
\lvec(138 66)
\lvec(135 66)
\ifill f:0
\move(139 65)
\lvec(145 65)
\lvec(145 66)
\lvec(139 66)
\ifill f:0
\move(146 65)
\lvec(155 65)
\lvec(155 66)
\lvec(146 66)
\ifill f:0
\move(156 65)
\lvec(160 65)
\lvec(160 66)
\lvec(156 66)
\ifill f:0
\move(161 65)
\lvec(163 65)
\lvec(163 66)
\lvec(161 66)
\ifill f:0
\move(164 65)
\lvec(167 65)
\lvec(167 66)
\lvec(164 66)
\ifill f:0
\move(168 65)
\lvec(170 65)
\lvec(170 66)
\lvec(168 66)
\ifill f:0
\move(172 65)
\lvec(173 65)
\lvec(173 66)
\lvec(172 66)
\ifill f:0
\move(174 65)
\lvec(175 65)
\lvec(175 66)
\lvec(174 66)
\ifill f:0
\move(176 65)
\lvec(177 65)
\lvec(177 66)
\lvec(176 66)
\ifill f:0
\move(178 65)
\lvec(179 65)
\lvec(179 66)
\lvec(178 66)
\ifill f:0
\move(180 65)
\lvec(184 65)
\lvec(184 66)
\lvec(180 66)
\ifill f:0
\move(185 65)
\lvec(186 65)
\lvec(186 66)
\lvec(185 66)
\ifill f:0
\move(187 65)
\lvec(189 65)
\lvec(189 66)
\lvec(187 66)
\ifill f:0
\move(190 65)
\lvec(191 65)
\lvec(191 66)
\lvec(190 66)
\ifill f:0
\move(192 65)
\lvec(194 65)
\lvec(194 66)
\lvec(192 66)
\ifill f:0
\move(195 65)
\lvec(197 65)
\lvec(197 66)
\lvec(195 66)
\ifill f:0
\move(198 65)
\lvec(211 65)
\lvec(211 66)
\lvec(198 66)
\ifill f:0
\move(212 65)
\lvec(220 65)
\lvec(220 66)
\lvec(212 66)
\ifill f:0
\move(221 65)
\lvec(226 65)
\lvec(226 66)
\lvec(221 66)
\ifill f:0
\move(227 65)
\lvec(231 65)
\lvec(231 66)
\lvec(227 66)
\ifill f:0
\move(232 65)
\lvec(238 65)
\lvec(238 66)
\lvec(232 66)
\ifill f:0
\move(239 65)
\lvec(247 65)
\lvec(247 66)
\lvec(239 66)
\ifill f:0
\move(249 65)
\lvec(257 65)
\lvec(257 66)
\lvec(249 66)
\ifill f:0
\move(258 65)
\lvec(290 65)
\lvec(290 66)
\lvec(258 66)
\ifill f:0
\move(292 65)
\lvec(308 65)
\lvec(308 66)
\lvec(292 66)
\ifill f:0
\move(310 65)
\lvec(315 65)
\lvec(315 66)
\lvec(310 66)
\ifill f:0
\move(316 65)
\lvec(321 65)
\lvec(321 66)
\lvec(316 66)
\ifill f:0
\move(323 65)
\lvec(325 65)
\lvec(325 66)
\lvec(323 66)
\ifill f:0
\move(326 65)
\lvec(362 65)
\lvec(362 66)
\lvec(326 66)
\ifill f:0
\move(363 65)
\lvec(365 65)
\lvec(365 66)
\lvec(363 66)
\ifill f:0
\move(366 65)
\lvec(368 65)
\lvec(368 66)
\lvec(366 66)
\ifill f:0
\move(369 65)
\lvec(381 65)
\lvec(381 66)
\lvec(369 66)
\ifill f:0
\move(382 65)
\lvec(384 65)
\lvec(384 66)
\lvec(382 66)
\ifill f:0
\move(385 65)
\lvec(387 65)
\lvec(387 66)
\lvec(385 66)
\ifill f:0
\move(388 65)
\lvec(390 65)
\lvec(390 66)
\lvec(388 66)
\ifill f:0
\move(391 65)
\lvec(393 65)
\lvec(393 66)
\lvec(391 66)
\ifill f:0
\move(394 65)
\lvec(401 65)
\lvec(401 66)
\lvec(394 66)
\ifill f:0
\move(402 65)
\lvec(404 65)
\lvec(404 66)
\lvec(402 66)
\ifill f:0
\move(405 65)
\lvec(412 65)
\lvec(412 66)
\lvec(405 66)
\ifill f:0
\move(413 65)
\lvec(427 65)
\lvec(427 66)
\lvec(413 66)
\ifill f:0
\move(428 65)
\lvec(442 65)
\lvec(442 66)
\lvec(428 66)
\ifill f:0
\move(443 65)
\lvec(450 65)
\lvec(450 66)
\lvec(443 66)
\ifill f:0
\move(15 66)
\lvec(17 66)
\lvec(17 67)
\lvec(15 67)
\ifill f:0
\move(23 66)
\lvec(26 66)
\lvec(26 67)
\lvec(23 67)
\ifill f:0
\move(27 66)
\lvec(28 66)
\lvec(28 67)
\lvec(27 67)
\ifill f:0
\move(36 66)
\lvec(37 66)
\lvec(37 67)
\lvec(36 67)
\ifill f:0
\move(38 66)
\lvec(39 66)
\lvec(39 67)
\lvec(38 67)
\ifill f:0
\move(40 66)
\lvec(45 66)
\lvec(45 67)
\lvec(40 67)
\ifill f:0
\move(47 66)
\lvec(50 66)
\lvec(50 67)
\lvec(47 67)
\ifill f:0
\move(51 66)
\lvec(53 66)
\lvec(53 67)
\lvec(51 67)
\ifill f:0
\move(54 66)
\lvec(55 66)
\lvec(55 67)
\lvec(54 67)
\ifill f:0
\move(56 66)
\lvec(58 66)
\lvec(58 67)
\lvec(56 67)
\ifill f:0
\move(59 66)
\lvec(61 66)
\lvec(61 67)
\lvec(59 67)
\ifill f:0
\move(62 66)
\lvec(65 66)
\lvec(65 67)
\lvec(62 67)
\ifill f:0
\move(66 66)
\lvec(70 66)
\lvec(70 67)
\lvec(66 67)
\ifill f:0
\move(72 66)
\lvec(74 66)
\lvec(74 67)
\lvec(72 67)
\ifill f:0
\move(75 66)
\lvec(78 66)
\lvec(78 67)
\lvec(75 67)
\ifill f:0
\move(79 66)
\lvec(82 66)
\lvec(82 67)
\lvec(79 67)
\ifill f:0
\move(83 66)
\lvec(84 66)
\lvec(84 67)
\lvec(83 67)
\ifill f:0
\move(86 66)
\lvec(87 66)
\lvec(87 67)
\lvec(86 67)
\ifill f:0
\move(88 66)
\lvec(89 66)
\lvec(89 67)
\lvec(88 67)
\ifill f:0
\move(90 66)
\lvec(91 66)
\lvec(91 67)
\lvec(90 67)
\ifill f:0
\move(92 66)
\lvec(93 66)
\lvec(93 67)
\lvec(92 67)
\ifill f:0
\move(94 66)
\lvec(95 66)
\lvec(95 67)
\lvec(94 67)
\ifill f:0
\move(96 66)
\lvec(98 66)
\lvec(98 67)
\lvec(96 67)
\ifill f:0
\move(99 66)
\lvec(101 66)
\lvec(101 67)
\lvec(99 67)
\ifill f:0
\move(102 66)
\lvec(105 66)
\lvec(105 67)
\lvec(102 67)
\ifill f:0
\move(106 66)
\lvec(109 66)
\lvec(109 67)
\lvec(106 67)
\ifill f:0
\move(110 66)
\lvec(119 66)
\lvec(119 67)
\lvec(110 67)
\ifill f:0
\move(121 66)
\lvec(122 66)
\lvec(122 67)
\lvec(121 67)
\ifill f:0
\move(123 66)
\lvec(124 66)
\lvec(124 67)
\lvec(123 67)
\ifill f:0
\move(132 66)
\lvec(136 66)
\lvec(136 67)
\lvec(132 67)
\ifill f:0
\move(139 66)
\lvec(142 66)
\lvec(142 67)
\lvec(139 67)
\ifill f:0
\move(143 66)
\lvec(145 66)
\lvec(145 67)
\lvec(143 67)
\ifill f:0
\move(146 66)
\lvec(150 66)
\lvec(150 67)
\lvec(146 67)
\ifill f:0
\move(151 66)
\lvec(167 66)
\lvec(167 67)
\lvec(151 67)
\ifill f:0
\move(168 66)
\lvec(170 66)
\lvec(170 67)
\lvec(168 67)
\ifill f:0
\move(171 66)
\lvec(174 66)
\lvec(174 67)
\lvec(171 67)
\ifill f:0
\move(175 66)
\lvec(176 66)
\lvec(176 67)
\lvec(175 67)
\ifill f:0
\move(177 66)
\lvec(185 66)
\lvec(185 67)
\lvec(177 67)
\ifill f:0
\move(186 66)
\lvec(187 66)
\lvec(187 67)
\lvec(186 67)
\ifill f:0
\move(188 66)
\lvec(192 66)
\lvec(192 67)
\lvec(188 67)
\ifill f:0
\move(193 66)
\lvec(194 66)
\lvec(194 67)
\lvec(193 67)
\ifill f:0
\move(195 66)
\lvec(197 66)
\lvec(197 67)
\lvec(195 67)
\ifill f:0
\move(198 66)
\lvec(200 66)
\lvec(200 67)
\lvec(198 67)
\ifill f:0
\move(201 66)
\lvec(203 66)
\lvec(203 67)
\lvec(201 67)
\ifill f:0
\move(204 66)
\lvec(226 66)
\lvec(226 67)
\lvec(204 67)
\ifill f:0
\move(227 66)
\lvec(229 66)
\lvec(229 67)
\lvec(227 67)
\ifill f:0
\move(231 66)
\lvec(235 66)
\lvec(235 67)
\lvec(231 67)
\ifill f:0
\move(236 66)
\lvec(241 66)
\lvec(241 67)
\lvec(236 67)
\ifill f:0
\move(243 66)
\lvec(250 66)
\lvec(250 67)
\lvec(243 67)
\ifill f:0
\move(251 66)
\lvec(257 66)
\lvec(257 67)
\lvec(251 67)
\ifill f:0
\move(258 66)
\lvec(264 66)
\lvec(264 67)
\lvec(258 67)
\ifill f:0
\move(265 66)
\lvec(290 66)
\lvec(290 67)
\lvec(265 67)
\ifill f:0
\move(294 66)
\lvec(304 66)
\lvec(304 67)
\lvec(294 67)
\ifill f:0
\move(306 66)
\lvec(313 66)
\lvec(313 67)
\lvec(306 67)
\ifill f:0
\move(315 66)
\lvec(321 66)
\lvec(321 67)
\lvec(315 67)
\ifill f:0
\move(323 66)
\lvec(325 66)
\lvec(325 67)
\lvec(323 67)
\ifill f:0
\move(326 66)
\lvec(327 66)
\lvec(327 67)
\lvec(326 67)
\ifill f:0
\move(328 66)
\lvec(333 66)
\lvec(333 67)
\lvec(328 67)
\ifill f:0
\move(334 66)
\lvec(362 66)
\lvec(362 67)
\lvec(334 67)
\ifill f:0
\move(363 66)
\lvec(365 66)
\lvec(365 67)
\lvec(363 67)
\ifill f:0
\move(366 66)
\lvec(369 66)
\lvec(369 67)
\lvec(366 67)
\ifill f:0
\move(370 66)
\lvec(376 66)
\lvec(376 67)
\lvec(370 67)
\ifill f:0
\move(377 66)
\lvec(383 66)
\lvec(383 67)
\lvec(377 67)
\ifill f:0
\move(384 66)
\lvec(386 66)
\lvec(386 67)
\lvec(384 67)
\ifill f:0
\move(387 66)
\lvec(401 66)
\lvec(401 67)
\lvec(387 67)
\ifill f:0
\move(402 66)
\lvec(410 66)
\lvec(410 67)
\lvec(402 67)
\ifill f:0
\move(411 66)
\lvec(421 66)
\lvec(421 67)
\lvec(411 67)
\ifill f:0
\move(422 66)
\lvec(442 66)
\lvec(442 67)
\lvec(422 67)
\ifill f:0
\move(443 66)
\lvec(451 66)
\lvec(451 67)
\lvec(443 67)
\ifill f:0
\move(14 67)
\lvec(15 67)
\lvec(15 68)
\lvec(14 68)
\ifill f:0
\move(16 67)
\lvec(17 67)
\lvec(17 68)
\lvec(16 68)
\ifill f:0
\move(20 67)
\lvec(21 67)
\lvec(21 68)
\lvec(20 68)
\ifill f:0
\move(22 67)
\lvec(26 67)
\lvec(26 68)
\lvec(22 68)
\ifill f:0
\move(36 67)
\lvec(37 67)
\lvec(37 68)
\lvec(36 68)
\ifill f:0
\move(38 67)
\lvec(41 67)
\lvec(41 68)
\lvec(38 68)
\ifill f:0
\move(43 67)
\lvec(46 67)
\lvec(46 68)
\lvec(43 68)
\ifill f:0
\move(47 67)
\lvec(50 67)
\lvec(50 68)
\lvec(47 68)
\ifill f:0
\move(52 67)
\lvec(53 67)
\lvec(53 68)
\lvec(52 68)
\ifill f:0
\move(55 67)
\lvec(56 67)
\lvec(56 68)
\lvec(55 68)
\ifill f:0
\move(57 67)
\lvec(58 67)
\lvec(58 68)
\lvec(57 68)
\ifill f:0
\move(60 67)
\lvec(62 67)
\lvec(62 68)
\lvec(60 68)
\ifill f:0
\move(63 67)
\lvec(65 67)
\lvec(65 68)
\lvec(63 68)
\ifill f:0
\move(66 67)
\lvec(67 67)
\lvec(67 68)
\lvec(66 68)
\ifill f:0
\move(68 67)
\lvec(76 67)
\lvec(76 68)
\lvec(68 68)
\ifill f:0
\move(78 67)
\lvec(82 67)
\lvec(82 68)
\lvec(78 68)
\ifill f:0
\move(83 67)
\lvec(85 67)
\lvec(85 68)
\lvec(83 68)
\ifill f:0
\move(87 67)
\lvec(90 67)
\lvec(90 68)
\lvec(87 68)
\ifill f:0
\move(91 67)
\lvec(92 67)
\lvec(92 68)
\lvec(91 68)
\ifill f:0
\move(95 67)
\lvec(96 67)
\lvec(96 68)
\lvec(95 68)
\ifill f:0
\move(97 67)
\lvec(98 67)
\lvec(98 68)
\lvec(97 68)
\ifill f:0
\move(99 67)
\lvec(101 67)
\lvec(101 68)
\lvec(99 68)
\ifill f:0
\move(102 67)
\lvec(108 67)
\lvec(108 68)
\lvec(102 68)
\ifill f:0
\move(109 67)
\lvec(111 67)
\lvec(111 68)
\lvec(109 68)
\ifill f:0
\move(112 67)
\lvec(120 67)
\lvec(120 68)
\lvec(112 68)
\ifill f:0
\move(121 67)
\lvec(122 67)
\lvec(122 68)
\lvec(121 68)
\ifill f:0
\move(123 67)
\lvec(132 67)
\lvec(132 68)
\lvec(123 68)
\ifill f:0
\move(136 67)
\lvec(141 67)
\lvec(141 68)
\lvec(136 68)
\ifill f:0
\move(142 67)
\lvec(145 67)
\lvec(145 68)
\lvec(142 68)
\ifill f:0
\move(146 67)
\lvec(147 67)
\lvec(147 68)
\lvec(146 68)
\ifill f:0
\move(148 67)
\lvec(151 67)
\lvec(151 68)
\lvec(148 68)
\ifill f:0
\move(152 67)
\lvec(155 67)
\lvec(155 68)
\lvec(152 68)
\ifill f:0
\move(156 67)
\lvec(161 67)
\lvec(161 68)
\lvec(156 68)
\ifill f:0
\move(162 67)
\lvec(170 67)
\lvec(170 68)
\lvec(162 68)
\ifill f:0
\move(171 67)
\lvec(172 67)
\lvec(172 68)
\lvec(171 68)
\ifill f:0
\move(173 67)
\lvec(174 67)
\lvec(174 68)
\lvec(173 68)
\ifill f:0
\move(175 67)
\lvec(176 67)
\lvec(176 68)
\lvec(175 68)
\ifill f:0
\move(177 67)
\lvec(188 67)
\lvec(188 68)
\lvec(177 68)
\ifill f:0
\move(189 67)
\lvec(190 67)
\lvec(190 68)
\lvec(189 68)
\ifill f:0
\move(191 67)
\lvec(192 67)
\lvec(192 68)
\lvec(191 68)
\ifill f:0
\move(193 67)
\lvec(194 67)
\lvec(194 68)
\lvec(193 68)
\ifill f:0
\move(195 67)
\lvec(197 67)
\lvec(197 68)
\lvec(195 68)
\ifill f:0
\move(198 67)
\lvec(208 67)
\lvec(208 68)
\lvec(198 68)
\ifill f:0
\move(209 67)
\lvec(211 67)
\lvec(211 68)
\lvec(209 68)
\ifill f:0
\move(212 67)
\lvec(214 67)
\lvec(214 68)
\lvec(212 68)
\ifill f:0
\move(215 67)
\lvec(226 67)
\lvec(226 68)
\lvec(215 68)
\ifill f:0
\move(227 67)
\lvec(229 67)
\lvec(229 68)
\lvec(227 68)
\ifill f:0
\move(230 67)
\lvec(233 67)
\lvec(233 68)
\lvec(230 68)
\ifill f:0
\move(234 67)
\lvec(239 67)
\lvec(239 68)
\lvec(234 68)
\ifill f:0
\move(240 67)
\lvec(250 67)
\lvec(250 68)
\lvec(240 68)
\ifill f:0
\move(252 67)
\lvec(257 67)
\lvec(257 68)
\lvec(252 68)
\ifill f:0
\move(258 67)
\lvec(261 67)
\lvec(261 68)
\lvec(258 68)
\ifill f:0
\move(262 67)
\lvec(277 67)
\lvec(277 68)
\lvec(262 68)
\ifill f:0
\move(278 67)
\lvec(279 67)
\lvec(279 68)
\lvec(278 68)
\ifill f:0
\move(285 67)
\lvec(286 67)
\lvec(286 68)
\lvec(285 68)
\ifill f:0
\move(289 67)
\lvec(290 67)
\lvec(290 68)
\lvec(289 68)
\ifill f:0
\move(294 67)
\lvec(295 67)
\lvec(295 68)
\lvec(294 68)
\ifill f:0
\move(296 67)
\lvec(298 67)
\lvec(298 68)
\lvec(296 68)
\ifill f:0
\move(299 67)
\lvec(310 67)
\lvec(310 68)
\lvec(299 68)
\ifill f:0
\move(312 67)
\lvec(320 67)
\lvec(320 68)
\lvec(312 68)
\ifill f:0
\move(322 67)
\lvec(325 67)
\lvec(325 68)
\lvec(322 68)
\ifill f:0
\move(326 67)
\lvec(328 67)
\lvec(328 68)
\lvec(326 68)
\ifill f:0
\move(329 67)
\lvec(335 67)
\lvec(335 68)
\lvec(329 68)
\ifill f:0
\move(336 67)
\lvec(341 67)
\lvec(341 68)
\lvec(336 68)
\ifill f:0
\move(342 67)
\lvec(362 67)
\lvec(362 68)
\lvec(342 68)
\ifill f:0
\move(363 67)
\lvec(370 67)
\lvec(370 68)
\lvec(363 68)
\ifill f:0
\move(371 67)
\lvec(392 67)
\lvec(392 68)
\lvec(371 68)
\ifill f:0
\move(393 67)
\lvec(395 67)
\lvec(395 68)
\lvec(393 68)
\ifill f:0
\move(396 67)
\lvec(401 67)
\lvec(401 68)
\lvec(396 68)
\ifill f:0
\move(402 67)
\lvec(436 67)
\lvec(436 68)
\lvec(402 68)
\ifill f:0
\move(437 67)
\lvec(442 67)
\lvec(442 68)
\lvec(437 68)
\ifill f:0
\move(443 67)
\lvec(446 67)
\lvec(446 68)
\lvec(443 68)
\ifill f:0
\move(447 67)
\lvec(451 67)
\lvec(451 68)
\lvec(447 68)
\ifill f:0
\move(16 68)
\lvec(17 68)
\lvec(17 69)
\lvec(16 69)
\ifill f:0
\move(18 68)
\lvec(22 68)
\lvec(22 69)
\lvec(18 69)
\ifill f:0
\move(23 68)
\lvec(26 68)
\lvec(26 69)
\lvec(23 69)
\ifill f:0
\move(36 68)
\lvec(37 68)
\lvec(37 69)
\lvec(36 69)
\ifill f:0
\move(38 68)
\lvec(50 68)
\lvec(50 69)
\lvec(38 69)
\ifill f:0
\move(56 68)
\lvec(57 68)
\lvec(57 69)
\lvec(56 69)
\ifill f:0
\move(58 68)
\lvec(59 68)
\lvec(59 69)
\lvec(58 69)
\ifill f:0
\move(60 68)
\lvec(62 68)
\lvec(62 69)
\lvec(60 69)
\ifill f:0
\move(63 68)
\lvec(65 68)
\lvec(65 69)
\lvec(63 69)
\ifill f:0
\move(67 68)
\lvec(71 68)
\lvec(71 69)
\lvec(67 69)
\ifill f:0
\move(72 68)
\lvec(74 68)
\lvec(74 69)
\lvec(72 69)
\ifill f:0
\move(76 68)
\lvec(82 68)
\lvec(82 69)
\lvec(76 69)
\ifill f:0
\move(83 68)
\lvec(86 68)
\lvec(86 69)
\lvec(83 69)
\ifill f:0
\move(87 68)
\lvec(91 68)
\lvec(91 69)
\lvec(87 69)
\ifill f:0
\move(92 68)
\lvec(94 68)
\lvec(94 69)
\lvec(92 69)
\ifill f:0
\move(95 68)
\lvec(96 68)
\lvec(96 69)
\lvec(95 69)
\ifill f:0
\move(97 68)
\lvec(98 68)
\lvec(98 69)
\lvec(97 69)
\ifill f:0
\move(99 68)
\lvec(101 68)
\lvec(101 69)
\lvec(99 69)
\ifill f:0
\move(102 68)
\lvec(106 68)
\lvec(106 69)
\lvec(102 69)
\ifill f:0
\move(107 68)
\lvec(120 68)
\lvec(120 69)
\lvec(107 69)
\ifill f:0
\move(121 68)
\lvec(122 68)
\lvec(122 69)
\lvec(121 69)
\ifill f:0
\move(123 68)
\lvec(140 68)
\lvec(140 69)
\lvec(123 69)
\ifill f:0
\move(141 68)
\lvec(145 68)
\lvec(145 69)
\lvec(141 69)
\ifill f:0
\move(146 68)
\lvec(147 68)
\lvec(147 69)
\lvec(146 69)
\ifill f:0
\move(148 68)
\lvec(160 68)
\lvec(160 69)
\lvec(148 69)
\ifill f:0
\move(161 68)
\lvec(170 68)
\lvec(170 69)
\lvec(161 69)
\ifill f:0
\move(171 68)
\lvec(172 68)
\lvec(172 69)
\lvec(171 69)
\ifill f:0
\move(173 68)
\lvec(177 68)
\lvec(177 69)
\lvec(173 69)
\ifill f:0
\move(178 68)
\lvec(179 68)
\lvec(179 69)
\lvec(178 69)
\ifill f:0
\move(180 68)
\lvec(197 68)
\lvec(197 69)
\lvec(180 69)
\ifill f:0
\move(198 68)
\lvec(226 68)
\lvec(226 69)
\lvec(198 69)
\ifill f:0
\move(227 68)
\lvec(228 68)
\lvec(228 69)
\lvec(227 69)
\ifill f:0
\move(229 68)
\lvec(232 68)
\lvec(232 69)
\lvec(229 69)
\ifill f:0
\move(233 68)
\lvec(241 68)
\lvec(241 69)
\lvec(233 69)
\ifill f:0
\move(242 68)
\lvec(257 68)
\lvec(257 69)
\lvec(242 69)
\ifill f:0
\move(258 68)
\lvec(269 68)
\lvec(269 69)
\lvec(258 69)
\ifill f:0
\move(270 68)
\lvec(283 68)
\lvec(283 69)
\lvec(270 69)
\ifill f:0
\move(285 68)
\lvec(287 68)
\lvec(287 69)
\lvec(285 69)
\ifill f:0
\move(289 68)
\lvec(290 68)
\lvec(290 69)
\lvec(289 69)
\ifill f:0
\move(293 68)
\lvec(294 68)
\lvec(294 69)
\lvec(293 69)
\ifill f:0
\move(300 68)
\lvec(302 68)
\lvec(302 69)
\lvec(300 69)
\ifill f:0
\move(303 68)
\lvec(318 68)
\lvec(318 69)
\lvec(303 69)
\ifill f:0
\move(321 68)
\lvec(325 68)
\lvec(325 69)
\lvec(321 69)
\ifill f:0
\move(326 68)
\lvec(329 68)
\lvec(329 69)
\lvec(326 69)
\ifill f:0
\move(330 68)
\lvec(337 68)
\lvec(337 69)
\lvec(330 69)
\ifill f:0
\move(338 68)
\lvec(344 68)
\lvec(344 69)
\lvec(338 69)
\ifill f:0
\move(345 68)
\lvec(350 68)
\lvec(350 69)
\lvec(345 69)
\ifill f:0
\move(351 68)
\lvec(356 68)
\lvec(356 69)
\lvec(351 69)
\ifill f:0
\move(357 68)
\lvec(362 68)
\lvec(362 69)
\lvec(357 69)
\ifill f:0
\move(363 68)
\lvec(366 68)
\lvec(366 69)
\lvec(363 69)
\ifill f:0
\move(367 68)
\lvec(375 68)
\lvec(375 69)
\lvec(367 69)
\ifill f:0
\move(376 68)
\lvec(391 68)
\lvec(391 69)
\lvec(376 69)
\ifill f:0
\move(392 68)
\lvec(395 68)
\lvec(395 69)
\lvec(392 69)
\ifill f:0
\move(396 68)
\lvec(401 68)
\lvec(401 69)
\lvec(396 69)
\ifill f:0
\move(402 68)
\lvec(405 68)
\lvec(405 69)
\lvec(402 69)
\ifill f:0
\move(406 68)
\lvec(418 68)
\lvec(418 69)
\lvec(406 69)
\ifill f:0
\move(419 68)
\lvec(421 68)
\lvec(421 69)
\lvec(419 69)
\ifill f:0
\move(422 68)
\lvec(424 68)
\lvec(424 69)
\lvec(422 69)
\ifill f:0
\move(425 68)
\lvec(427 68)
\lvec(427 69)
\lvec(425 69)
\ifill f:0
\move(428 68)
\lvec(430 68)
\lvec(430 69)
\lvec(428 69)
\ifill f:0
\move(431 68)
\lvec(442 68)
\lvec(442 69)
\lvec(431 69)
\ifill f:0
\move(443 68)
\lvec(444 68)
\lvec(444 69)
\lvec(443 69)
\ifill f:0
\move(445 68)
\lvec(449 68)
\lvec(449 69)
\lvec(445 69)
\ifill f:0
\move(450 68)
\lvec(451 68)
\lvec(451 69)
\lvec(450 69)
\ifill f:0
\move(16 69)
\lvec(17 69)
\lvec(17 70)
\lvec(16 70)
\ifill f:0
\move(18 69)
\lvec(19 69)
\lvec(19 70)
\lvec(18 70)
\ifill f:0
\move(20 69)
\lvec(21 69)
\lvec(21 70)
\lvec(20 70)
\ifill f:0
\move(25 69)
\lvec(26 69)
\lvec(26 70)
\lvec(25 70)
\ifill f:0
\move(36 69)
\lvec(37 69)
\lvec(37 70)
\lvec(36 70)
\ifill f:0
\move(40 69)
\lvec(42 69)
\lvec(42 70)
\lvec(40 70)
\ifill f:0
\move(43 69)
\lvec(45 69)
\lvec(45 70)
\lvec(43 70)
\ifill f:0
\move(49 69)
\lvec(50 69)
\lvec(50 70)
\lvec(49 70)
\ifill f:0
\move(54 69)
\lvec(55 69)
\lvec(55 70)
\lvec(54 70)
\ifill f:0
\move(57 69)
\lvec(58 69)
\lvec(58 70)
\lvec(57 70)
\ifill f:0
\move(59 69)
\lvec(60 69)
\lvec(60 70)
\lvec(59 70)
\ifill f:0
\move(61 69)
\lvec(65 69)
\lvec(65 70)
\lvec(61 70)
\ifill f:0
\move(66 69)
\lvec(70 69)
\lvec(70 70)
\lvec(66 70)
\ifill f:0
\move(71 69)
\lvec(82 69)
\lvec(82 70)
\lvec(71 70)
\ifill f:0
\move(84 69)
\lvec(87 69)
\lvec(87 70)
\lvec(84 70)
\ifill f:0
\move(88 69)
\lvec(90 69)
\lvec(90 70)
\lvec(88 70)
\ifill f:0
\move(91 69)
\lvec(94 69)
\lvec(94 70)
\lvec(91 70)
\ifill f:0
\move(96 69)
\lvec(97 69)
\lvec(97 70)
\lvec(96 70)
\ifill f:0
\move(100 69)
\lvec(101 69)
\lvec(101 70)
\lvec(100 70)
\ifill f:0
\move(102 69)
\lvec(103 69)
\lvec(103 70)
\lvec(102 70)
\ifill f:0
\move(104 69)
\lvec(111 69)
\lvec(111 70)
\lvec(104 70)
\ifill f:0
\move(112 69)
\lvec(122 69)
\lvec(122 70)
\lvec(112 70)
\ifill f:0
\move(123 69)
\lvec(126 69)
\lvec(126 70)
\lvec(123 70)
\ifill f:0
\move(134 69)
\lvec(135 69)
\lvec(135 70)
\lvec(134 70)
\ifill f:0
\move(139 69)
\lvec(145 69)
\lvec(145 70)
\lvec(139 70)
\ifill f:0
\move(146 69)
\lvec(148 69)
\lvec(148 70)
\lvec(146 70)
\ifill f:0
\move(150 69)
\lvec(154 69)
\lvec(154 70)
\lvec(150 70)
\ifill f:0
\move(155 69)
\lvec(159 69)
\lvec(159 70)
\lvec(155 70)
\ifill f:0
\move(160 69)
\lvec(163 69)
\lvec(163 70)
\lvec(160 70)
\ifill f:0
\move(164 69)
\lvec(170 69)
\lvec(170 70)
\lvec(164 70)
\ifill f:0
\move(171 69)
\lvec(172 69)
\lvec(172 70)
\lvec(171 70)
\ifill f:0
\move(173 69)
\lvec(175 69)
\lvec(175 70)
\lvec(173 70)
\ifill f:0
\move(176 69)
\lvec(180 69)
\lvec(180 70)
\lvec(176 70)
\ifill f:0
\move(181 69)
\lvec(187 69)
\lvec(187 70)
\lvec(181 70)
\ifill f:0
\move(188 69)
\lvec(189 69)
\lvec(189 70)
\lvec(188 70)
\ifill f:0
\move(190 69)
\lvec(191 69)
\lvec(191 70)
\lvec(190 70)
\ifill f:0
\move(192 69)
\lvec(195 69)
\lvec(195 70)
\lvec(192 70)
\ifill f:0
\move(196 69)
\lvec(197 69)
\lvec(197 70)
\lvec(196 70)
\ifill f:0
\move(198 69)
\lvec(201 69)
\lvec(201 70)
\lvec(198 70)
\ifill f:0
\move(202 69)
\lvec(206 69)
\lvec(206 70)
\lvec(202 70)
\ifill f:0
\move(207 69)
\lvec(219 69)
\lvec(219 70)
\lvec(207 70)
\ifill f:0
\move(220 69)
\lvec(226 69)
\lvec(226 70)
\lvec(220 70)
\ifill f:0
\move(227 69)
\lvec(228 69)
\lvec(228 70)
\lvec(227 70)
\ifill f:0
\move(229 69)
\lvec(235 69)
\lvec(235 70)
\lvec(229 70)
\ifill f:0
\move(236 69)
\lvec(239 69)
\lvec(239 70)
\lvec(236 70)
\ifill f:0
\move(240 69)
\lvec(248 69)
\lvec(248 70)
\lvec(240 70)
\ifill f:0
\move(249 69)
\lvec(253 69)
\lvec(253 70)
\lvec(249 70)
\ifill f:0
\move(254 69)
\lvec(257 69)
\lvec(257 70)
\lvec(254 70)
\ifill f:0
\move(258 69)
\lvec(259 69)
\lvec(259 70)
\lvec(258 70)
\ifill f:0
\move(260 69)
\lvec(266 69)
\lvec(266 70)
\lvec(260 70)
\ifill f:0
\move(267 69)
\lvec(274 69)
\lvec(274 70)
\lvec(267 70)
\ifill f:0
\move(276 69)
\lvec(287 69)
\lvec(287 70)
\lvec(276 70)
\ifill f:0
\move(289 69)
\lvec(290 69)
\lvec(290 70)
\lvec(289 70)
\ifill f:0
\move(291 69)
\lvec(315 69)
\lvec(315 70)
\lvec(291 70)
\ifill f:0
\move(316 69)
\lvec(317 69)
\lvec(317 70)
\lvec(316 70)
\ifill f:0
\move(318 69)
\lvec(325 69)
\lvec(325 70)
\lvec(318 70)
\ifill f:0
\move(326 69)
\lvec(330 69)
\lvec(330 70)
\lvec(326 70)
\ifill f:0
\move(331 69)
\lvec(340 69)
\lvec(340 70)
\lvec(331 70)
\ifill f:0
\move(341 69)
\lvec(348 69)
\lvec(348 70)
\lvec(341 70)
\ifill f:0
\move(349 69)
\lvec(355 69)
\lvec(355 70)
\lvec(349 70)
\ifill f:0
\move(356 69)
\lvec(362 69)
\lvec(362 70)
\lvec(356 70)
\ifill f:0
\move(363 69)
\lvec(367 69)
\lvec(367 70)
\lvec(363 70)
\ifill f:0
\move(368 69)
\lvec(377 69)
\lvec(377 70)
\lvec(368 70)
\ifill f:0
\move(378 69)
\lvec(386 69)
\lvec(386 70)
\lvec(378 70)
\ifill f:0
\move(387 69)
\lvec(390 69)
\lvec(390 70)
\lvec(387 70)
\ifill f:0
\move(391 69)
\lvec(398 69)
\lvec(398 70)
\lvec(391 70)
\ifill f:0
\move(399 69)
\lvec(401 69)
\lvec(401 70)
\lvec(399 70)
\ifill f:0
\move(402 69)
\lvec(416 69)
\lvec(416 70)
\lvec(402 70)
\ifill f:0
\move(417 69)
\lvec(423 69)
\lvec(423 70)
\lvec(417 70)
\ifill f:0
\move(424 69)
\lvec(426 69)
\lvec(426 70)
\lvec(424 70)
\ifill f:0
\move(427 69)
\lvec(429 69)
\lvec(429 70)
\lvec(427 70)
\ifill f:0
\move(430 69)
\lvec(432 69)
\lvec(432 70)
\lvec(430 70)
\ifill f:0
\move(433 69)
\lvec(438 69)
\lvec(438 70)
\lvec(433 70)
\ifill f:0
\move(439 69)
\lvec(442 69)
\lvec(442 70)
\lvec(439 70)
\ifill f:0
\move(443 69)
\lvec(444 69)
\lvec(444 70)
\lvec(443 70)
\ifill f:0
\move(445 69)
\lvec(451 69)
\lvec(451 70)
\lvec(445 70)
\ifill f:0
\move(16 70)
\lvec(17 70)
\lvec(17 71)
\lvec(16 71)
\ifill f:0
\move(20 70)
\lvec(21 70)
\lvec(21 71)
\lvec(20 71)
\ifill f:0
\move(25 70)
\lvec(26 70)
\lvec(26 71)
\lvec(25 71)
\ifill f:0
\move(36 70)
\lvec(37 70)
\lvec(37 71)
\lvec(36 71)
\ifill f:0
\move(38 70)
\lvec(39 70)
\lvec(39 71)
\lvec(38 71)
\ifill f:0
\move(40 70)
\lvec(41 70)
\lvec(41 71)
\lvec(40 71)
\ifill f:0
\move(42 70)
\lvec(43 70)
\lvec(43 71)
\lvec(42 71)
\ifill f:0
\move(44 70)
\lvec(46 70)
\lvec(46 71)
\lvec(44 71)
\ifill f:0
\move(49 70)
\lvec(50 70)
\lvec(50 71)
\lvec(49 71)
\ifill f:0
\move(56 70)
\lvec(57 70)
\lvec(57 71)
\lvec(56 71)
\ifill f:0
\move(58 70)
\lvec(65 70)
\lvec(65 71)
\lvec(58 71)
\ifill f:0
\move(66 70)
\lvec(73 70)
\lvec(73 71)
\lvec(66 71)
\ifill f:0
\move(75 70)
\lvec(82 70)
\lvec(82 71)
\lvec(75 71)
\ifill f:0
\move(86 70)
\lvec(88 70)
\lvec(88 71)
\lvec(86 71)
\ifill f:0
\move(89 70)
\lvec(93 70)
\lvec(93 71)
\lvec(89 71)
\ifill f:0
\move(94 70)
\lvec(95 70)
\lvec(95 71)
\lvec(94 71)
\ifill f:0
\move(96 70)
\lvec(98 70)
\lvec(98 71)
\lvec(96 71)
\ifill f:0
\move(100 70)
\lvec(101 70)
\lvec(101 71)
\lvec(100 71)
\ifill f:0
\move(102 70)
\lvec(103 70)
\lvec(103 71)
\lvec(102 71)
\ifill f:0
\move(104 70)
\lvec(105 70)
\lvec(105 71)
\lvec(104 71)
\ifill f:0
\move(106 70)
\lvec(113 70)
\lvec(113 71)
\lvec(106 71)
\ifill f:0
\move(114 70)
\lvec(122 70)
\lvec(122 71)
\lvec(114 71)
\ifill f:0
\move(123 70)
\lvec(125 70)
\lvec(125 71)
\lvec(123 71)
\ifill f:0
\move(128 70)
\lvec(134 70)
\lvec(134 71)
\lvec(128 71)
\ifill f:0
\move(135 70)
\lvec(138 70)
\lvec(138 71)
\lvec(135 71)
\ifill f:0
\move(139 70)
\lvec(145 70)
\lvec(145 71)
\lvec(139 71)
\ifill f:0
\move(146 70)
\lvec(150 70)
\lvec(150 71)
\lvec(146 71)
\ifill f:0
\move(151 70)
\lvec(170 70)
\lvec(170 71)
\lvec(151 71)
\ifill f:0
\move(171 70)
\lvec(173 70)
\lvec(173 71)
\lvec(171 71)
\ifill f:0
\move(174 70)
\lvec(176 70)
\lvec(176 71)
\lvec(174 71)
\ifill f:0
\move(177 70)
\lvec(181 70)
\lvec(181 71)
\lvec(177 71)
\ifill f:0
\move(182 70)
\lvec(186 70)
\lvec(186 71)
\lvec(182 71)
\ifill f:0
\move(187 70)
\lvec(191 70)
\lvec(191 71)
\lvec(187 71)
\ifill f:0
\move(192 70)
\lvec(193 70)
\lvec(193 71)
\lvec(192 71)
\ifill f:0
\move(194 70)
\lvec(195 70)
\lvec(195 71)
\lvec(194 71)
\ifill f:0
\move(196 70)
\lvec(197 70)
\lvec(197 71)
\lvec(196 71)
\ifill f:0
\move(198 70)
\lvec(201 70)
\lvec(201 71)
\lvec(198 71)
\ifill f:0
\move(202 70)
\lvec(203 70)
\lvec(203 71)
\lvec(202 71)
\ifill f:0
\move(204 70)
\lvec(212 70)
\lvec(212 71)
\lvec(204 71)
\ifill f:0
\move(213 70)
\lvec(226 70)
\lvec(226 71)
\lvec(213 71)
\ifill f:0
\move(227 70)
\lvec(231 70)
\lvec(231 71)
\lvec(227 71)
\ifill f:0
\move(232 70)
\lvec(234 70)
\lvec(234 71)
\lvec(232 71)
\ifill f:0
\move(235 70)
\lvec(241 70)
\lvec(241 71)
\lvec(235 71)
\ifill f:0
\move(242 70)
\lvec(249 70)
\lvec(249 71)
\lvec(242 71)
\ifill f:0
\move(250 70)
\lvec(257 70)
\lvec(257 71)
\lvec(250 71)
\ifill f:0
\move(258 70)
\lvec(264 70)
\lvec(264 71)
\lvec(258 71)
\ifill f:0
\move(265 70)
\lvec(288 70)
\lvec(288 71)
\lvec(265 71)
\ifill f:0
\move(289 70)
\lvec(290 70)
\lvec(290 71)
\lvec(289 71)
\ifill f:0
\move(291 70)
\lvec(310 70)
\lvec(310 71)
\lvec(291 71)
\ifill f:0
\move(311 70)
\lvec(325 70)
\lvec(325 71)
\lvec(311 71)
\ifill f:0
\move(326 70)
\lvec(333 70)
\lvec(333 71)
\lvec(326 71)
\ifill f:0
\move(334 70)
\lvec(362 70)
\lvec(362 71)
\lvec(334 71)
\ifill f:0
\move(363 70)
\lvec(368 70)
\lvec(368 71)
\lvec(363 71)
\ifill f:0
\move(369 70)
\lvec(379 70)
\lvec(379 71)
\lvec(369 71)
\ifill f:0
\move(380 70)
\lvec(384 70)
\lvec(384 71)
\lvec(380 71)
\ifill f:0
\move(385 70)
\lvec(389 70)
\lvec(389 71)
\lvec(385 71)
\ifill f:0
\move(390 70)
\lvec(394 70)
\lvec(394 71)
\lvec(390 71)
\ifill f:0
\move(395 70)
\lvec(398 70)
\lvec(398 71)
\lvec(395 71)
\ifill f:0
\move(399 70)
\lvec(401 70)
\lvec(401 71)
\lvec(399 71)
\ifill f:0
\move(403 70)
\lvec(410 70)
\lvec(410 71)
\lvec(403 71)
\ifill f:0
\move(411 70)
\lvec(428 70)
\lvec(428 71)
\lvec(411 71)
\ifill f:0
\move(429 70)
\lvec(442 70)
\lvec(442 71)
\lvec(429 71)
\ifill f:0
\move(443 70)
\lvec(444 70)
\lvec(444 71)
\lvec(443 71)
\ifill f:0
\move(445 70)
\lvec(450 70)
\lvec(450 71)
\lvec(445 71)
\ifill f:0
\move(14 71)
\lvec(17 71)
\lvec(17 72)
\lvec(14 72)
\ifill f:0
\move(19 71)
\lvec(21 71)
\lvec(21 72)
\lvec(19 72)
\ifill f:0
\move(23 71)
\lvec(24 71)
\lvec(24 72)
\lvec(23 72)
\ifill f:0
\move(25 71)
\lvec(26 71)
\lvec(26 72)
\lvec(25 72)
\ifill f:0
\move(27 71)
\lvec(28 71)
\lvec(28 72)
\lvec(27 72)
\ifill f:0
\move(36 71)
\lvec(37 71)
\lvec(37 72)
\lvec(36 72)
\ifill f:0
\move(38 71)
\lvec(42 71)
\lvec(42 72)
\lvec(38 72)
\ifill f:0
\move(43 71)
\lvec(46 71)
\lvec(46 72)
\lvec(43 72)
\ifill f:0
\move(47 71)
\lvec(48 71)
\lvec(48 72)
\lvec(47 72)
\ifill f:0
\move(49 71)
\lvec(50 71)
\lvec(50 72)
\lvec(49 72)
\ifill f:0
\move(51 71)
\lvec(53 71)
\lvec(53 72)
\lvec(51 72)
\ifill f:0
\move(57 71)
\lvec(59 71)
\lvec(59 72)
\lvec(57 72)
\ifill f:0
\move(60 71)
\lvec(61 71)
\lvec(61 72)
\lvec(60 72)
\ifill f:0
\move(62 71)
\lvec(65 71)
\lvec(65 72)
\lvec(62 72)
\ifill f:0
\move(66 71)
\lvec(71 71)
\lvec(71 72)
\lvec(66 72)
\ifill f:0
\move(72 71)
\lvec(74 71)
\lvec(74 72)
\lvec(72 72)
\ifill f:0
\move(76 71)
\lvec(78 71)
\lvec(78 72)
\lvec(76 72)
\ifill f:0
\move(80 71)
\lvec(82 71)
\lvec(82 72)
\lvec(80 72)
\ifill f:0
\move(87 71)
\lvec(90 71)
\lvec(90 72)
\lvec(87 72)
\ifill f:0
\move(92 71)
\lvec(93 71)
\lvec(93 72)
\lvec(92 72)
\ifill f:0
\move(95 71)
\lvec(99 71)
\lvec(99 72)
\lvec(95 72)
\ifill f:0
\move(100 71)
\lvec(101 71)
\lvec(101 72)
\lvec(100 72)
\ifill f:0
\move(102 71)
\lvec(109 71)
\lvec(109 72)
\lvec(102 72)
\ifill f:0
\move(110 71)
\lvec(111 71)
\lvec(111 72)
\lvec(110 72)
\ifill f:0
\move(112 71)
\lvec(117 71)
\lvec(117 72)
\lvec(112 72)
\ifill f:0
\move(118 71)
\lvec(122 71)
\lvec(122 72)
\lvec(118 72)
\ifill f:0
\move(123 71)
\lvec(124 71)
\lvec(124 72)
\lvec(123 72)
\ifill f:0
\move(126 71)
\lvec(130 71)
\lvec(130 72)
\lvec(126 72)
\ifill f:0
\move(133 71)
\lvec(145 71)
\lvec(145 72)
\lvec(133 72)
\ifill f:0
\move(146 71)
\lvec(152 71)
\lvec(152 72)
\lvec(146 72)
\ifill f:0
\move(154 71)
\lvec(160 71)
\lvec(160 72)
\lvec(154 72)
\ifill f:0
\move(161 71)
\lvec(170 71)
\lvec(170 72)
\lvec(161 72)
\ifill f:0
\move(171 71)
\lvec(173 71)
\lvec(173 72)
\lvec(171 72)
\ifill f:0
\move(175 71)
\lvec(180 71)
\lvec(180 72)
\lvec(175 72)
\ifill f:0
\move(181 71)
\lvec(185 71)
\lvec(185 72)
\lvec(181 72)
\ifill f:0
\move(186 71)
\lvec(195 71)
\lvec(195 72)
\lvec(186 72)
\ifill f:0
\move(196 71)
\lvec(197 71)
\lvec(197 72)
\lvec(196 72)
\ifill f:0
\move(198 71)
\lvec(213 71)
\lvec(213 72)
\lvec(198 72)
\ifill f:0
\move(214 71)
\lvec(220 71)
\lvec(220 72)
\lvec(214 72)
\ifill f:0
\move(221 71)
\lvec(226 71)
\lvec(226 72)
\lvec(221 72)
\ifill f:0
\move(227 71)
\lvec(233 71)
\lvec(233 72)
\lvec(227 72)
\ifill f:0
\move(234 71)
\lvec(242 71)
\lvec(242 72)
\lvec(234 72)
\ifill f:0
\move(243 71)
\lvec(250 71)
\lvec(250 72)
\lvec(243 72)
\ifill f:0
\move(251 71)
\lvec(254 71)
\lvec(254 72)
\lvec(251 72)
\ifill f:0
\move(255 71)
\lvec(257 71)
\lvec(257 72)
\lvec(255 72)
\ifill f:0
\move(258 71)
\lvec(268 71)
\lvec(268 72)
\lvec(258 72)
\ifill f:0
\move(269 71)
\lvec(280 71)
\lvec(280 72)
\lvec(269 72)
\ifill f:0
\move(281 71)
\lvec(288 71)
\lvec(288 72)
\lvec(281 72)
\ifill f:0
\move(289 71)
\lvec(290 71)
\lvec(290 72)
\lvec(289 72)
\ifill f:0
\move(291 71)
\lvec(298 71)
\lvec(298 72)
\lvec(291 72)
\ifill f:0
\move(300 71)
\lvec(301 71)
\lvec(301 72)
\lvec(300 72)
\ifill f:0
\move(302 71)
\lvec(325 71)
\lvec(325 72)
\lvec(302 72)
\ifill f:0
\move(326 71)
\lvec(338 71)
\lvec(338 72)
\lvec(326 72)
\ifill f:0
\move(340 71)
\lvec(352 71)
\lvec(352 72)
\lvec(340 72)
\ifill f:0
\move(353 71)
\lvec(362 71)
\lvec(362 72)
\lvec(353 72)
\ifill f:0
\move(363 71)
\lvec(369 71)
\lvec(369 72)
\lvec(363 72)
\ifill f:0
\move(371 71)
\lvec(382 71)
\lvec(382 72)
\lvec(371 72)
\ifill f:0
\move(383 71)
\lvec(401 71)
\lvec(401 72)
\lvec(383 72)
\ifill f:0
\move(402 71)
\lvec(411 71)
\lvec(411 72)
\lvec(402 72)
\ifill f:0
\move(412 71)
\lvec(427 71)
\lvec(427 72)
\lvec(412 72)
\ifill f:0
\move(428 71)
\lvec(434 71)
\lvec(434 72)
\lvec(428 72)
\ifill f:0
\move(435 71)
\lvec(442 71)
\lvec(442 72)
\lvec(435 72)
\ifill f:0
\move(443 71)
\lvec(445 71)
\lvec(445 72)
\lvec(443 72)
\ifill f:0
\move(446 71)
\lvec(451 71)
\lvec(451 72)
\lvec(446 72)
\ifill f:0
\move(15 72)
\lvec(17 72)
\lvec(17 73)
\lvec(15 73)
\ifill f:0
\move(18 72)
\lvec(19 72)
\lvec(19 73)
\lvec(18 73)
\ifill f:0
\move(20 72)
\lvec(21 72)
\lvec(21 73)
\lvec(20 73)
\ifill f:0
\move(23 72)
\lvec(24 72)
\lvec(24 73)
\lvec(23 73)
\ifill f:0
\move(25 72)
\lvec(26 72)
\lvec(26 73)
\lvec(25 73)
\ifill f:0
\move(36 72)
\lvec(37 72)
\lvec(37 73)
\lvec(36 73)
\ifill f:0
\move(38 72)
\lvec(41 72)
\lvec(41 73)
\lvec(38 73)
\ifill f:0
\move(42 72)
\lvec(43 72)
\lvec(43 73)
\lvec(42 73)
\ifill f:0
\move(44 72)
\lvec(45 72)
\lvec(45 73)
\lvec(44 73)
\ifill f:0
\move(46 72)
\lvec(48 72)
\lvec(48 73)
\lvec(46 73)
\ifill f:0
\move(49 72)
\lvec(50 72)
\lvec(50 73)
\lvec(49 73)
\ifill f:0
\move(51 72)
\lvec(53 72)
\lvec(53 73)
\lvec(51 73)
\ifill f:0
\move(54 72)
\lvec(55 72)
\lvec(55 73)
\lvec(54 73)
\ifill f:0
\move(56 72)
\lvec(57 72)
\lvec(57 73)
\lvec(56 73)
\ifill f:0
\move(59 72)
\lvec(63 72)
\lvec(63 73)
\lvec(59 73)
\ifill f:0
\move(64 72)
\lvec(65 72)
\lvec(65 73)
\lvec(64 73)
\ifill f:0
\move(66 72)
\lvec(67 72)
\lvec(67 73)
\lvec(66 73)
\ifill f:0
\move(68 72)
\lvec(73 72)
\lvec(73 73)
\lvec(68 73)
\ifill f:0
\move(74 72)
\lvec(78 72)
\lvec(78 73)
\lvec(74 73)
\ifill f:0
\move(81 72)
\lvec(82 72)
\lvec(82 73)
\lvec(81 73)
\ifill f:0
\move(88 72)
\lvec(92 72)
\lvec(92 73)
\lvec(88 73)
\ifill f:0
\move(95 72)
\lvec(96 72)
\lvec(96 73)
\lvec(95 73)
\ifill f:0
\move(97 72)
\lvec(99 72)
\lvec(99 73)
\lvec(97 73)
\ifill f:0
\move(100 72)
\lvec(101 72)
\lvec(101 73)
\lvec(100 73)
\ifill f:0
\move(102 72)
\lvec(106 72)
\lvec(106 73)
\lvec(102 73)
\ifill f:0
\move(107 72)
\lvec(115 72)
\lvec(115 73)
\lvec(107 73)
\ifill f:0
\move(116 72)
\lvec(122 72)
\lvec(122 73)
\lvec(116 73)
\ifill f:0
\move(123 72)
\lvec(124 72)
\lvec(124 73)
\lvec(123 73)
\ifill f:0
\move(125 72)
\lvec(129 72)
\lvec(129 73)
\lvec(125 73)
\ifill f:0
\move(130 72)
\lvec(134 72)
\lvec(134 73)
\lvec(130 73)
\ifill f:0
\move(137 72)
\lvec(145 72)
\lvec(145 73)
\lvec(137 73)
\ifill f:0
\move(146 72)
\lvec(164 72)
\lvec(164 73)
\lvec(146 73)
\ifill f:0
\move(165 72)
\lvec(170 72)
\lvec(170 73)
\lvec(165 73)
\ifill f:0
\move(171 72)
\lvec(174 72)
\lvec(174 73)
\lvec(171 73)
\ifill f:0
\move(175 72)
\lvec(184 72)
\lvec(184 73)
\lvec(175 73)
\ifill f:0
\move(185 72)
\lvec(187 72)
\lvec(187 73)
\lvec(185 73)
\ifill f:0
\move(188 72)
\lvec(195 72)
\lvec(195 73)
\lvec(188 73)
\ifill f:0
\move(196 72)
\lvec(197 72)
\lvec(197 73)
\lvec(196 73)
\ifill f:0
\move(198 72)
\lvec(204 72)
\lvec(204 73)
\lvec(198 73)
\ifill f:0
\move(205 72)
\lvec(226 72)
\lvec(226 73)
\lvec(205 73)
\ifill f:0
\move(228 72)
\lvec(230 72)
\lvec(230 73)
\lvec(228 73)
\ifill f:0
\move(231 72)
\lvec(235 72)
\lvec(235 73)
\lvec(231 73)
\ifill f:0
\move(236 72)
\lvec(238 72)
\lvec(238 73)
\lvec(236 73)
\ifill f:0
\move(239 72)
\lvec(241 72)
\lvec(241 73)
\lvec(239 73)
\ifill f:0
\move(242 72)
\lvec(244 72)
\lvec(244 73)
\lvec(242 73)
\ifill f:0
\move(245 72)
\lvec(247 72)
\lvec(247 73)
\lvec(245 73)
\ifill f:0
\move(248 72)
\lvec(250 72)
\lvec(250 73)
\lvec(248 73)
\ifill f:0
\move(251 72)
\lvec(254 72)
\lvec(254 73)
\lvec(251 73)
\ifill f:0
\move(255 72)
\lvec(257 72)
\lvec(257 73)
\lvec(255 73)
\ifill f:0
\move(259 72)
\lvec(262 72)
\lvec(262 73)
\lvec(259 73)
\ifill f:0
\move(263 72)
\lvec(266 72)
\lvec(266 73)
\lvec(263 73)
\ifill f:0
\move(267 72)
\lvec(271 72)
\lvec(271 73)
\lvec(267 73)
\ifill f:0
\move(272 72)
\lvec(282 72)
\lvec(282 73)
\lvec(272 73)
\ifill f:0
\move(283 72)
\lvec(290 72)
\lvec(290 73)
\lvec(283 73)
\ifill f:0
\move(291 72)
\lvec(297 72)
\lvec(297 73)
\lvec(291 73)
\ifill f:0
\move(298 72)
\lvec(309 72)
\lvec(309 73)
\lvec(298 73)
\ifill f:0
\move(311 72)
\lvec(325 72)
\lvec(325 73)
\lvec(311 73)
\ifill f:0
\move(326 72)
\lvec(347 72)
\lvec(347 73)
\lvec(326 73)
\ifill f:0
\move(348 72)
\lvec(362 72)
\lvec(362 73)
\lvec(348 73)
\ifill f:0
\move(363 72)
\lvec(371 72)
\lvec(371 73)
\lvec(363 73)
\ifill f:0
\move(372 72)
\lvec(379 72)
\lvec(379 73)
\lvec(372 73)
\ifill f:0
\move(380 72)
\lvec(392 72)
\lvec(392 73)
\lvec(380 73)
\ifill f:0
\move(393 72)
\lvec(401 72)
\lvec(401 73)
\lvec(393 73)
\ifill f:0
\move(402 72)
\lvec(417 72)
\lvec(417 73)
\lvec(402 73)
\ifill f:0
\move(418 72)
\lvec(442 72)
\lvec(442 73)
\lvec(418 73)
\ifill f:0
\move(443 72)
\lvec(445 72)
\lvec(445 73)
\lvec(443 73)
\ifill f:0
\move(446 72)
\lvec(448 72)
\lvec(448 73)
\lvec(446 73)
\ifill f:0
\move(449 72)
\lvec(451 72)
\lvec(451 73)
\lvec(449 73)
\ifill f:0
\move(15 73)
\lvec(17 73)
\lvec(17 74)
\lvec(15 74)
\ifill f:0
\move(18 73)
\lvec(19 73)
\lvec(19 74)
\lvec(18 74)
\ifill f:0
\move(20 73)
\lvec(22 73)
\lvec(22 74)
\lvec(20 74)
\ifill f:0
\move(24 73)
\lvec(26 73)
\lvec(26 74)
\lvec(24 74)
\ifill f:0
\move(36 73)
\lvec(37 73)
\lvec(37 74)
\lvec(36 74)
\ifill f:0
\move(38 73)
\lvec(39 73)
\lvec(39 74)
\lvec(38 74)
\ifill f:0
\move(40 73)
\lvec(42 73)
\lvec(42 74)
\lvec(40 74)
\ifill f:0
\move(43 73)
\lvec(45 73)
\lvec(45 74)
\lvec(43 74)
\ifill f:0
\move(47 73)
\lvec(50 73)
\lvec(50 74)
\lvec(47 74)
\ifill f:0
\move(51 73)
\lvec(52 73)
\lvec(52 74)
\lvec(51 74)
\ifill f:0
\move(57 73)
\lvec(60 73)
\lvec(60 74)
\lvec(57 74)
\ifill f:0
\move(61 73)
\lvec(63 73)
\lvec(63 74)
\lvec(61 74)
\ifill f:0
\move(64 73)
\lvec(65 73)
\lvec(65 74)
\lvec(64 74)
\ifill f:0
\move(66 73)
\lvec(71 73)
\lvec(71 74)
\lvec(66 74)
\ifill f:0
\move(72 73)
\lvec(74 73)
\lvec(74 74)
\lvec(72 74)
\ifill f:0
\move(75 73)
\lvec(79 73)
\lvec(79 74)
\lvec(75 74)
\ifill f:0
\move(81 73)
\lvec(82 73)
\lvec(82 74)
\lvec(81 74)
\ifill f:0
\move(83 73)
\lvec(88 73)
\lvec(88 74)
\lvec(83 74)
\ifill f:0
\move(92 73)
\lvec(93 73)
\lvec(93 74)
\lvec(92 74)
\ifill f:0
\move(94 73)
\lvec(95 73)
\lvec(95 74)
\lvec(94 74)
\ifill f:0
\move(97 73)
\lvec(99 73)
\lvec(99 74)
\lvec(97 74)
\ifill f:0
\move(100 73)
\lvec(101 73)
\lvec(101 74)
\lvec(100 74)
\ifill f:0
\move(102 73)
\lvec(109 73)
\lvec(109 74)
\lvec(102 74)
\ifill f:0
\move(110 73)
\lvec(111 73)
\lvec(111 74)
\lvec(110 74)
\ifill f:0
\move(112 73)
\lvec(122 73)
\lvec(122 74)
\lvec(112 74)
\ifill f:0
\move(123 73)
\lvec(124 73)
\lvec(124 74)
\lvec(123 74)
\ifill f:0
\move(125 73)
\lvec(127 73)
\lvec(127 74)
\lvec(125 74)
\ifill f:0
\move(128 73)
\lvec(131 73)
\lvec(131 74)
\lvec(128 74)
\ifill f:0
\move(133 73)
\lvec(138 73)
\lvec(138 74)
\lvec(133 74)
\ifill f:0
\move(139 73)
\lvec(145 73)
\lvec(145 74)
\lvec(139 74)
\ifill f:0
\move(146 73)
\lvec(162 73)
\lvec(162 74)
\lvec(146 74)
\ifill f:0
\move(164 73)
\lvec(170 73)
\lvec(170 74)
\lvec(164 74)
\ifill f:0
\move(171 73)
\lvec(175 73)
\lvec(175 74)
\lvec(171 74)
\ifill f:0
\move(176 73)
\lvec(179 73)
\lvec(179 74)
\lvec(176 74)
\ifill f:0
\move(180 73)
\lvec(186 73)
\lvec(186 74)
\lvec(180 74)
\ifill f:0
\move(187 73)
\lvec(189 73)
\lvec(189 74)
\lvec(187 74)
\ifill f:0
\move(190 73)
\lvec(192 73)
\lvec(192 74)
\lvec(190 74)
\ifill f:0
\move(193 73)
\lvec(195 73)
\lvec(195 74)
\lvec(193 74)
\ifill f:0
\move(196 73)
\lvec(197 73)
\lvec(197 74)
\lvec(196 74)
\ifill f:0
\move(198 73)
\lvec(200 73)
\lvec(200 74)
\lvec(198 74)
\ifill f:0
\move(201 73)
\lvec(207 73)
\lvec(207 74)
\lvec(201 74)
\ifill f:0
\move(208 73)
\lvec(226 73)
\lvec(226 74)
\lvec(208 74)
\ifill f:0
\move(228 73)
\lvec(229 73)
\lvec(229 74)
\lvec(228 74)
\ifill f:0
\move(230 73)
\lvec(232 73)
\lvec(232 74)
\lvec(230 74)
\ifill f:0
\move(233 73)
\lvec(234 73)
\lvec(234 74)
\lvec(233 74)
\ifill f:0
\move(235 73)
\lvec(237 73)
\lvec(237 74)
\lvec(235 74)
\ifill f:0
\move(238 73)
\lvec(242 73)
\lvec(242 74)
\lvec(238 74)
\ifill f:0
\move(243 73)
\lvec(248 73)
\lvec(248 74)
\lvec(243 74)
\ifill f:0
\move(249 73)
\lvec(251 73)
\lvec(251 74)
\lvec(249 74)
\ifill f:0
\move(252 73)
\lvec(257 73)
\lvec(257 74)
\lvec(252 74)
\ifill f:0
\move(258 73)
\lvec(265 73)
\lvec(265 74)
\lvec(258 74)
\ifill f:0
\move(266 73)
\lvec(269 73)
\lvec(269 74)
\lvec(266 74)
\ifill f:0
\move(270 73)
\lvec(283 73)
\lvec(283 74)
\lvec(270 74)
\ifill f:0
\move(284 73)
\lvec(290 73)
\lvec(290 74)
\lvec(284 74)
\ifill f:0
\move(291 73)
\lvec(295 73)
\lvec(295 74)
\lvec(291 74)
\ifill f:0
\move(296 73)
\lvec(304 73)
\lvec(304 74)
\lvec(296 74)
\ifill f:0
\move(305 73)
\lvec(314 73)
\lvec(314 74)
\lvec(305 74)
\ifill f:0
\move(316 73)
\lvec(325 73)
\lvec(325 74)
\lvec(316 74)
\ifill f:0
\move(326 73)
\lvec(362 73)
\lvec(362 74)
\lvec(326 74)
\ifill f:0
\move(364 73)
\lvec(373 73)
\lvec(373 74)
\lvec(364 74)
\ifill f:0
\move(374 73)
\lvec(383 73)
\lvec(383 74)
\lvec(374 74)
\ifill f:0
\move(384 73)
\lvec(390 73)
\lvec(390 74)
\lvec(384 74)
\ifill f:0
\move(391 73)
\lvec(397 73)
\lvec(397 74)
\lvec(391 74)
\ifill f:0
\move(398 73)
\lvec(401 73)
\lvec(401 74)
\lvec(398 74)
\ifill f:0
\move(402 73)
\lvec(403 73)
\lvec(403 74)
\lvec(402 74)
\ifill f:0
\move(404 73)
\lvec(414 73)
\lvec(414 74)
\lvec(404 74)
\ifill f:0
\move(415 73)
\lvec(419 73)
\lvec(419 74)
\lvec(415 74)
\ifill f:0
\move(420 73)
\lvec(424 73)
\lvec(424 74)
\lvec(420 74)
\ifill f:0
\move(425 73)
\lvec(442 73)
\lvec(442 74)
\lvec(425 74)
\ifill f:0
\move(443 73)
\lvec(445 73)
\lvec(445 74)
\lvec(443 74)
\ifill f:0
\move(446 73)
\lvec(449 73)
\lvec(449 74)
\lvec(446 74)
\ifill f:0
\move(450 73)
\lvec(451 73)
\lvec(451 74)
\lvec(450 74)
\ifill f:0
\move(15 74)
\lvec(17 74)
\lvec(17 75)
\lvec(15 75)
\ifill f:0
\move(19 74)
\lvec(21 74)
\lvec(21 75)
\lvec(19 75)
\ifill f:0
\move(24 74)
\lvec(26 74)
\lvec(26 75)
\lvec(24 75)
\ifill f:0
\move(36 74)
\lvec(37 74)
\lvec(37 75)
\lvec(36 75)
\ifill f:0
\move(38 74)
\lvec(39 74)
\lvec(39 75)
\lvec(38 75)
\ifill f:0
\move(40 74)
\lvec(46 74)
\lvec(46 75)
\lvec(40 75)
\ifill f:0
\move(47 74)
\lvec(50 74)
\lvec(50 75)
\lvec(47 75)
\ifill f:0
\move(60 74)
\lvec(63 74)
\lvec(63 75)
\lvec(60 75)
\ifill f:0
\move(64 74)
\lvec(65 74)
\lvec(65 75)
\lvec(64 75)
\ifill f:0
\move(66 74)
\lvec(75 74)
\lvec(75 75)
\lvec(66 75)
\ifill f:0
\move(77 74)
\lvec(80 74)
\lvec(80 75)
\lvec(77 75)
\ifill f:0
\move(81 74)
\lvec(82 74)
\lvec(82 75)
\lvec(81 75)
\ifill f:0
\move(83 74)
\lvec(85 74)
\lvec(85 75)
\lvec(83 75)
\ifill f:0
\move(86 74)
\lvec(93 74)
\lvec(93 75)
\lvec(86 75)
\ifill f:0
\move(96 74)
\lvec(98 74)
\lvec(98 75)
\lvec(96 75)
\ifill f:0
\move(100 74)
\lvec(101 74)
\lvec(101 75)
\lvec(100 75)
\ifill f:0
\move(102 74)
\lvec(105 74)
\lvec(105 75)
\lvec(102 75)
\ifill f:0
\move(106 74)
\lvec(116 74)
\lvec(116 75)
\lvec(106 75)
\ifill f:0
\move(117 74)
\lvec(122 74)
\lvec(122 75)
\lvec(117 75)
\ifill f:0
\move(124 74)
\lvec(126 74)
\lvec(126 75)
\lvec(124 75)
\ifill f:0
\move(127 74)
\lvec(130 74)
\lvec(130 75)
\lvec(127 75)
\ifill f:0
\move(131 74)
\lvec(133 74)
\lvec(133 75)
\lvec(131 75)
\ifill f:0
\move(135 74)
\lvec(138 74)
\lvec(138 75)
\lvec(135 75)
\ifill f:0
\move(139 74)
\lvec(140 74)
\lvec(140 75)
\lvec(139 75)
\ifill f:0
\move(141 74)
\lvec(145 74)
\lvec(145 75)
\lvec(141 75)
\ifill f:0
\move(146 74)
\lvec(149 74)
\lvec(149 75)
\lvec(146 75)
\ifill f:0
\move(154 74)
\lvec(155 74)
\lvec(155 75)
\lvec(154 75)
\ifill f:0
\move(158 74)
\lvec(159 74)
\lvec(159 75)
\lvec(158 75)
\ifill f:0
\move(161 74)
\lvec(170 74)
\lvec(170 75)
\lvec(161 75)
\ifill f:0
\move(171 74)
\lvec(176 74)
\lvec(176 75)
\lvec(171 75)
\ifill f:0
\move(177 74)
\lvec(188 74)
\lvec(188 75)
\lvec(177 75)
\ifill f:0
\move(189 74)
\lvec(192 74)
\lvec(192 75)
\lvec(189 75)
\ifill f:0
\move(193 74)
\lvec(195 74)
\lvec(195 75)
\lvec(193 75)
\ifill f:0
\move(196 74)
\lvec(197 74)
\lvec(197 75)
\lvec(196 75)
\ifill f:0
\move(198 74)
\lvec(203 74)
\lvec(203 75)
\lvec(198 75)
\ifill f:0
\move(204 74)
\lvec(208 74)
\lvec(208 75)
\lvec(204 75)
\ifill f:0
\move(209 74)
\lvec(219 74)
\lvec(219 75)
\lvec(209 75)
\ifill f:0
\move(220 74)
\lvec(221 74)
\lvec(221 75)
\lvec(220 75)
\ifill f:0
\move(222 74)
\lvec(223 74)
\lvec(223 75)
\lvec(222 75)
\ifill f:0
\move(224 74)
\lvec(226 74)
\lvec(226 75)
\lvec(224 75)
\ifill f:0
\move(228 74)
\lvec(229 74)
\lvec(229 75)
\lvec(228 75)
\ifill f:0
\move(230 74)
\lvec(246 74)
\lvec(246 75)
\lvec(230 75)
\ifill f:0
\move(247 74)
\lvec(257 74)
\lvec(257 75)
\lvec(247 75)
\ifill f:0
\move(258 74)
\lvec(264 74)
\lvec(264 75)
\lvec(258 75)
\ifill f:0
\move(265 74)
\lvec(271 74)
\lvec(271 75)
\lvec(265 75)
\ifill f:0
\move(272 74)
\lvec(275 74)
\lvec(275 75)
\lvec(272 75)
\ifill f:0
\move(276 74)
\lvec(279 74)
\lvec(279 75)
\lvec(276 75)
\ifill f:0
\move(280 74)
\lvec(284 74)
\lvec(284 75)
\lvec(280 75)
\ifill f:0
\move(285 74)
\lvec(290 74)
\lvec(290 75)
\lvec(285 75)
\ifill f:0
\move(291 74)
\lvec(293 74)
\lvec(293 75)
\lvec(291 75)
\ifill f:0
\move(294 74)
\lvec(301 74)
\lvec(301 75)
\lvec(294 75)
\ifill f:0
\move(302 74)
\lvec(308 74)
\lvec(308 75)
\lvec(302 75)
\ifill f:0
\move(309 74)
\lvec(317 74)
\lvec(317 75)
\lvec(309 75)
\ifill f:0
\move(319 74)
\lvec(325 74)
\lvec(325 75)
\lvec(319 75)
\ifill f:0
\move(326 74)
\lvec(333 74)
\lvec(333 75)
\lvec(326 75)
\ifill f:0
\move(334 74)
\lvec(362 74)
\lvec(362 75)
\lvec(334 75)
\ifill f:0
\move(366 74)
\lvec(377 74)
\lvec(377 75)
\lvec(366 75)
\ifill f:0
\move(378 74)
\lvec(379 74)
\lvec(379 75)
\lvec(378 75)
\ifill f:0
\move(380 74)
\lvec(388 74)
\lvec(388 75)
\lvec(380 75)
\ifill f:0
\move(389 74)
\lvec(396 74)
\lvec(396 75)
\lvec(389 75)
\ifill f:0
\move(397 74)
\lvec(401 74)
\lvec(401 75)
\lvec(397 75)
\ifill f:0
\move(402 74)
\lvec(410 74)
\lvec(410 75)
\lvec(402 75)
\ifill f:0
\move(411 74)
\lvec(416 74)
\lvec(416 75)
\lvec(411 75)
\ifill f:0
\move(417 74)
\lvec(427 74)
\lvec(427 75)
\lvec(417 75)
\ifill f:0
\move(428 74)
\lvec(432 74)
\lvec(432 75)
\lvec(428 75)
\ifill f:0
\move(433 74)
\lvec(442 74)
\lvec(442 75)
\lvec(433 75)
\ifill f:0
\move(443 74)
\lvec(450 74)
\lvec(450 75)
\lvec(443 75)
\ifill f:0
\move(16 75)
\lvec(17 75)
\lvec(17 76)
\lvec(16 76)
\ifill f:0
\move(20 75)
\lvec(22 75)
\lvec(22 76)
\lvec(20 76)
\ifill f:0
\move(24 75)
\lvec(26 75)
\lvec(26 76)
\lvec(24 76)
\ifill f:0
\move(36 75)
\lvec(37 75)
\lvec(37 76)
\lvec(36 76)
\ifill f:0
\move(43 75)
\lvec(45 75)
\lvec(45 76)
\lvec(43 76)
\ifill f:0
\move(47 75)
\lvec(50 75)
\lvec(50 76)
\lvec(47 76)
\ifill f:0
\move(54 75)
\lvec(55 75)
\lvec(55 76)
\lvec(54 76)
\ifill f:0
\move(57 75)
\lvec(58 75)
\lvec(58 76)
\lvec(57 76)
\ifill f:0
\move(59 75)
\lvec(62 75)
\lvec(62 76)
\lvec(59 76)
\ifill f:0
\move(64 75)
\lvec(65 75)
\lvec(65 76)
\lvec(64 76)
\ifill f:0
\move(66 75)
\lvec(71 75)
\lvec(71 76)
\lvec(66 76)
\ifill f:0
\move(72 75)
\lvec(74 75)
\lvec(74 76)
\lvec(72 76)
\ifill f:0
\move(75 75)
\lvec(77 75)
\lvec(77 76)
\lvec(75 76)
\ifill f:0
\move(78 75)
\lvec(80 75)
\lvec(80 76)
\lvec(78 76)
\ifill f:0
\move(81 75)
\lvec(82 75)
\lvec(82 76)
\lvec(81 76)
\ifill f:0
\move(83 75)
\lvec(85 75)
\lvec(85 76)
\lvec(83 76)
\ifill f:0
\move(92 75)
\lvec(93 75)
\lvec(93 76)
\lvec(92 76)
\ifill f:0
\move(95 75)
\lvec(98 75)
\lvec(98 76)
\lvec(95 76)
\ifill f:0
\move(100 75)
\lvec(101 75)
\lvec(101 76)
\lvec(100 76)
\ifill f:0
\move(103 75)
\lvec(106 75)
\lvec(106 76)
\lvec(103 76)
\ifill f:0
\move(107 75)
\lvec(109 75)
\lvec(109 76)
\lvec(107 76)
\ifill f:0
\move(110 75)
\lvec(111 75)
\lvec(111 76)
\lvec(110 76)
\ifill f:0
\move(112 75)
\lvec(122 75)
\lvec(122 76)
\lvec(112 76)
\ifill f:0
\move(124 75)
\lvec(125 75)
\lvec(125 76)
\lvec(124 76)
\ifill f:0
\move(127 75)
\lvec(129 75)
\lvec(129 76)
\lvec(127 76)
\ifill f:0
\move(130 75)
\lvec(131 75)
\lvec(131 76)
\lvec(130 76)
\ifill f:0
\move(133 75)
\lvec(135 75)
\lvec(135 76)
\lvec(133 76)
\ifill f:0
\move(137 75)
\lvec(141 75)
\lvec(141 76)
\lvec(137 76)
\ifill f:0
\move(142 75)
\lvec(145 75)
\lvec(145 76)
\lvec(142 76)
\ifill f:0
\move(146 75)
\lvec(147 75)
\lvec(147 76)
\lvec(146 76)
\ifill f:0
\move(149 75)
\lvec(170 75)
\lvec(170 76)
\lvec(149 76)
\ifill f:0
\move(172 75)
\lvec(177 75)
\lvec(177 76)
\lvec(172 76)
\ifill f:0
\move(178 75)
\lvec(187 75)
\lvec(187 76)
\lvec(178 76)
\ifill f:0
\move(188 75)
\lvec(191 75)
\lvec(191 76)
\lvec(188 76)
\ifill f:0
\move(192 75)
\lvec(195 75)
\lvec(195 76)
\lvec(192 76)
\ifill f:0
\move(196 75)
\lvec(197 75)
\lvec(197 76)
\lvec(196 76)
\ifill f:0
\move(198 75)
\lvec(201 75)
\lvec(201 76)
\lvec(198 76)
\ifill f:0
\move(202 75)
\lvec(214 75)
\lvec(214 76)
\lvec(202 76)
\ifill f:0
\move(215 75)
\lvec(223 75)
\lvec(223 76)
\lvec(215 76)
\ifill f:0
\move(224 75)
\lvec(226 75)
\lvec(226 76)
\lvec(224 76)
\ifill f:0
\move(228 75)
\lvec(229 75)
\lvec(229 76)
\lvec(228 76)
\ifill f:0
\move(230 75)
\lvec(231 75)
\lvec(231 76)
\lvec(230 76)
\ifill f:0
\move(232 75)
\lvec(233 75)
\lvec(233 76)
\lvec(232 76)
\ifill f:0
\move(234 75)
\lvec(235 75)
\lvec(235 76)
\lvec(234 76)
\ifill f:0
\move(236 75)
\lvec(242 75)
\lvec(242 76)
\lvec(236 76)
\ifill f:0
\move(243 75)
\lvec(252 75)
\lvec(252 76)
\lvec(243 76)
\ifill f:0
\move(253 75)
\lvec(257 75)
\lvec(257 76)
\lvec(253 76)
\ifill f:0
\move(258 75)
\lvec(266 75)
\lvec(266 76)
\lvec(258 76)
\ifill f:0
\move(267 75)
\lvec(273 75)
\lvec(273 76)
\lvec(267 76)
\ifill f:0
\move(274 75)
\lvec(280 75)
\lvec(280 76)
\lvec(274 76)
\ifill f:0
\move(281 75)
\lvec(284 75)
\lvec(284 76)
\lvec(281 76)
\ifill f:0
\move(285 75)
\lvec(290 75)
\lvec(290 76)
\lvec(285 76)
\ifill f:0
\move(291 75)
\lvec(293 75)
\lvec(293 76)
\lvec(291 76)
\ifill f:0
\move(294 75)
\lvec(298 75)
\lvec(298 76)
\lvec(294 76)
\ifill f:0
\move(300 75)
\lvec(310 75)
\lvec(310 76)
\lvec(300 76)
\ifill f:0
\move(312 75)
\lvec(319 75)
\lvec(319 76)
\lvec(312 76)
\ifill f:0
\move(320 75)
\lvec(325 75)
\lvec(325 76)
\lvec(320 76)
\ifill f:0
\move(326 75)
\lvec(329 75)
\lvec(329 76)
\lvec(326 76)
\ifill f:0
\move(331 75)
\lvec(351 75)
\lvec(351 76)
\lvec(331 76)
\ifill f:0
\move(356 75)
\lvec(357 75)
\lvec(357 76)
\lvec(356 76)
\ifill f:0
\move(361 75)
\lvec(362 75)
\lvec(362 76)
\lvec(361 76)
\ifill f:0
\move(366 75)
\lvec(384 75)
\lvec(384 76)
\lvec(366 76)
\ifill f:0
\move(385 75)
\lvec(401 75)
\lvec(401 76)
\lvec(385 76)
\ifill f:0
\move(402 75)
\lvec(404 75)
\lvec(404 76)
\lvec(402 76)
\ifill f:0
\move(405 75)
\lvec(412 75)
\lvec(412 76)
\lvec(405 76)
\ifill f:0
\move(413 75)
\lvec(419 75)
\lvec(419 76)
\lvec(413 76)
\ifill f:0
\move(420 75)
\lvec(425 75)
\lvec(425 76)
\lvec(420 76)
\ifill f:0
\move(426 75)
\lvec(436 75)
\lvec(436 76)
\lvec(426 76)
\ifill f:0
\move(437 75)
\lvec(442 75)
\lvec(442 76)
\lvec(437 76)
\ifill f:0
\move(443 75)
\lvec(446 75)
\lvec(446 76)
\lvec(443 76)
\ifill f:0
\move(447 75)
\lvec(451 75)
\lvec(451 76)
\lvec(447 76)
\ifill f:0
\move(16 76)
\lvec(17 76)
\lvec(17 77)
\lvec(16 77)
\ifill f:0
\move(18 76)
\lvec(19 76)
\lvec(19 77)
\lvec(18 77)
\ifill f:0
\move(20 76)
\lvec(21 76)
\lvec(21 77)
\lvec(20 77)
\ifill f:0
\move(23 76)
\lvec(26 76)
\lvec(26 77)
\lvec(23 77)
\ifill f:0
\move(36 76)
\lvec(37 76)
\lvec(37 77)
\lvec(36 77)
\ifill f:0
\move(38 76)
\lvec(45 76)
\lvec(45 77)
\lvec(38 77)
\ifill f:0
\move(48 76)
\lvec(50 76)
\lvec(50 77)
\lvec(48 77)
\ifill f:0
\move(52 76)
\lvec(53 76)
\lvec(53 77)
\lvec(52 77)
\ifill f:0
\move(57 76)
\lvec(61 76)
\lvec(61 77)
\lvec(57 77)
\ifill f:0
\move(64 76)
\lvec(65 76)
\lvec(65 77)
\lvec(64 77)
\ifill f:0
\move(67 76)
\lvec(74 76)
\lvec(74 77)
\lvec(67 77)
\ifill f:0
\move(76 76)
\lvec(77 76)
\lvec(77 77)
\lvec(76 77)
\ifill f:0
\move(78 76)
\lvec(80 76)
\lvec(80 77)
\lvec(78 77)
\ifill f:0
\move(81 76)
\lvec(82 76)
\lvec(82 77)
\lvec(81 77)
\ifill f:0
\move(83 76)
\lvec(84 76)
\lvec(84 77)
\lvec(83 77)
\ifill f:0
\move(87 76)
\lvec(98 76)
\lvec(98 77)
\lvec(87 77)
\ifill f:0
\move(100 76)
\lvec(101 76)
\lvec(101 77)
\lvec(100 77)
\ifill f:0
\move(102 76)
\lvec(103 76)
\lvec(103 77)
\lvec(102 77)
\ifill f:0
\move(104 76)
\lvec(115 76)
\lvec(115 77)
\lvec(104 77)
\ifill f:0
\move(116 76)
\lvec(119 76)
\lvec(119 77)
\lvec(116 77)
\ifill f:0
\move(120 76)
\lvec(122 76)
\lvec(122 77)
\lvec(120 77)
\ifill f:0
\move(124 76)
\lvec(125 76)
\lvec(125 77)
\lvec(124 77)
\ifill f:0
\move(126 76)
\lvec(128 76)
\lvec(128 77)
\lvec(126 77)
\ifill f:0
\move(129 76)
\lvec(130 76)
\lvec(130 77)
\lvec(129 77)
\ifill f:0
\move(131 76)
\lvec(133 76)
\lvec(133 77)
\lvec(131 77)
\ifill f:0
\move(134 76)
\lvec(137 76)
\lvec(137 77)
\lvec(134 77)
\ifill f:0
\move(138 76)
\lvec(145 76)
\lvec(145 77)
\lvec(138 77)
\ifill f:0
\move(146 76)
\lvec(147 76)
\lvec(147 77)
\lvec(146 77)
\ifill f:0
\move(148 76)
\lvec(155 76)
\lvec(155 77)
\lvec(148 77)
\ifill f:0
\move(157 76)
\lvec(170 76)
\lvec(170 77)
\lvec(157 77)
\ifill f:0
\move(175 76)
\lvec(179 76)
\lvec(179 77)
\lvec(175 77)
\ifill f:0
\move(181 76)
\lvec(185 76)
\lvec(185 77)
\lvec(181 77)
\ifill f:0
\move(187 76)
\lvec(190 76)
\lvec(190 77)
\lvec(187 77)
\ifill f:0
\move(192 76)
\lvec(194 76)
\lvec(194 77)
\lvec(192 77)
\ifill f:0
\move(196 76)
\lvec(197 76)
\lvec(197 77)
\lvec(196 77)
\ifill f:0
\move(198 76)
\lvec(218 76)
\lvec(218 77)
\lvec(198 77)
\ifill f:0
\move(219 76)
\lvec(223 76)
\lvec(223 77)
\lvec(219 77)
\ifill f:0
\move(224 76)
\lvec(226 76)
\lvec(226 77)
\lvec(224 77)
\ifill f:0
\move(227 76)
\lvec(243 76)
\lvec(243 77)
\lvec(227 77)
\ifill f:0
\move(244 76)
\lvec(257 76)
\lvec(257 77)
\lvec(244 77)
\ifill f:0
\move(258 76)
\lvec(260 76)
\lvec(260 77)
\lvec(258 77)
\ifill f:0
\move(261 76)
\lvec(271 76)
\lvec(271 77)
\lvec(261 77)
\ifill f:0
\move(272 76)
\lvec(278 76)
\lvec(278 77)
\lvec(272 77)
\ifill f:0
\move(279 76)
\lvec(290 76)
\lvec(290 77)
\lvec(279 77)
\ifill f:0
\move(291 76)
\lvec(293 76)
\lvec(293 77)
\lvec(291 77)
\ifill f:0
\move(294 76)
\lvec(297 76)
\lvec(297 77)
\lvec(294 77)
\ifill f:0
\move(298 76)
\lvec(302 76)
\lvec(302 77)
\lvec(298 77)
\ifill f:0
\move(303 76)
\lvec(307 76)
\lvec(307 77)
\lvec(303 77)
\ifill f:0
\move(308 76)
\lvec(313 76)
\lvec(313 77)
\lvec(308 77)
\ifill f:0
\move(314 76)
\lvec(319 76)
\lvec(319 77)
\lvec(314 77)
\ifill f:0
\move(321 76)
\lvec(325 76)
\lvec(325 77)
\lvec(321 77)
\ifill f:0
\move(326 76)
\lvec(328 76)
\lvec(328 77)
\lvec(326 77)
\ifill f:0
\move(329 76)
\lvec(338 76)
\lvec(338 77)
\lvec(329 77)
\ifill f:0
\move(340 76)
\lvec(356 76)
\lvec(356 77)
\lvec(340 77)
\ifill f:0
\move(357 76)
\lvec(359 76)
\lvec(359 77)
\lvec(357 77)
\ifill f:0
\move(361 76)
\lvec(362 76)
\lvec(362 77)
\lvec(361 77)
\ifill f:0
\move(366 76)
\lvec(367 76)
\lvec(367 77)
\lvec(366 77)
\ifill f:0
\move(373 76)
\lvec(375 76)
\lvec(375 77)
\lvec(373 77)
\ifill f:0
\move(376 76)
\lvec(394 76)
\lvec(394 77)
\lvec(376 77)
\ifill f:0
\move(396 76)
\lvec(401 76)
\lvec(401 77)
\lvec(396 77)
\ifill f:0
\move(402 76)
\lvec(422 76)
\lvec(422 77)
\lvec(402 77)
\ifill f:0
\move(423 76)
\lvec(442 76)
\lvec(442 77)
\lvec(423 77)
\ifill f:0
\move(443 76)
\lvec(451 76)
\lvec(451 77)
\lvec(443 77)
\ifill f:0
\move(16 77)
\lvec(17 77)
\lvec(17 78)
\lvec(16 78)
\ifill f:0
\move(18 77)
\lvec(21 77)
\lvec(21 78)
\lvec(18 78)
\ifill f:0
\move(22 77)
\lvec(26 77)
\lvec(26 78)
\lvec(22 78)
\ifill f:0
\move(36 77)
\lvec(37 77)
\lvec(37 78)
\lvec(36 78)
\ifill f:0
\move(38 77)
\lvec(39 77)
\lvec(39 78)
\lvec(38 78)
\ifill f:0
\move(40 77)
\lvec(41 77)
\lvec(41 78)
\lvec(40 78)
\ifill f:0
\move(42 77)
\lvec(47 77)
\lvec(47 78)
\lvec(42 78)
\ifill f:0
\move(48 77)
\lvec(50 77)
\lvec(50 78)
\lvec(48 78)
\ifill f:0
\move(62 77)
\lvec(63 77)
\lvec(63 78)
\lvec(62 78)
\ifill f:0
\move(64 77)
\lvec(65 77)
\lvec(65 78)
\lvec(64 78)
\ifill f:0
\move(66 77)
\lvec(67 77)
\lvec(67 78)
\lvec(66 78)
\ifill f:0
\move(68 77)
\lvec(70 77)
\lvec(70 78)
\lvec(68 78)
\ifill f:0
\move(71 77)
\lvec(75 77)
\lvec(75 78)
\lvec(71 78)
\ifill f:0
\move(76 77)
\lvec(78 77)
\lvec(78 78)
\lvec(76 78)
\ifill f:0
\move(79 77)
\lvec(82 77)
\lvec(82 78)
\lvec(79 78)
\ifill f:0
\move(83 77)
\lvec(84 77)
\lvec(84 78)
\lvec(83 78)
\ifill f:0
\move(85 77)
\lvec(90 77)
\lvec(90 78)
\lvec(85 78)
\ifill f:0
\move(97 77)
\lvec(98 77)
\lvec(98 78)
\lvec(97 78)
\ifill f:0
\move(100 77)
\lvec(101 77)
\lvec(101 78)
\lvec(100 78)
\ifill f:0
\move(102 77)
\lvec(111 77)
\lvec(111 78)
\lvec(102 78)
\ifill f:0
\move(112 77)
\lvec(119 77)
\lvec(119 78)
\lvec(112 78)
\ifill f:0
\move(120 77)
\lvec(122 77)
\lvec(122 78)
\lvec(120 78)
\ifill f:0
\move(124 77)
\lvec(125 77)
\lvec(125 78)
\lvec(124 78)
\ifill f:0
\move(126 77)
\lvec(127 77)
\lvec(127 78)
\lvec(126 78)
\ifill f:0
\move(128 77)
\lvec(129 77)
\lvec(129 78)
\lvec(128 78)
\ifill f:0
\move(130 77)
\lvec(132 77)
\lvec(132 78)
\lvec(130 78)
\ifill f:0
\move(133 77)
\lvec(134 77)
\lvec(134 78)
\lvec(133 78)
\ifill f:0
\move(136 77)
\lvec(138 77)
\lvec(138 78)
\lvec(136 78)
\ifill f:0
\move(139 77)
\lvec(145 77)
\lvec(145 78)
\lvec(139 78)
\ifill f:0
\move(146 77)
\lvec(151 77)
\lvec(151 78)
\lvec(146 78)
\ifill f:0
\move(153 77)
\lvec(159 77)
\lvec(159 78)
\lvec(153 78)
\ifill f:0
\move(160 77)
\lvec(163 77)
\lvec(163 78)
\lvec(160 78)
\ifill f:0
\move(167 77)
\lvec(168 77)
\lvec(168 78)
\lvec(167 78)
\ifill f:0
\move(169 77)
\lvec(170 77)
\lvec(170 78)
\lvec(169 78)
\ifill f:0
\move(175 77)
\lvec(176 77)
\lvec(176 78)
\lvec(175 78)
\ifill f:0
\move(177 77)
\lvec(183 77)
\lvec(183 78)
\lvec(177 78)
\ifill f:0
\move(184 77)
\lvec(189 77)
\lvec(189 78)
\lvec(184 78)
\ifill f:0
\move(190 77)
\lvec(194 77)
\lvec(194 78)
\lvec(190 78)
\ifill f:0
\move(196 77)
\lvec(197 77)
\lvec(197 78)
\lvec(196 78)
\ifill f:0
\move(199 77)
\lvec(202 77)
\lvec(202 78)
\lvec(199 78)
\ifill f:0
\move(203 77)
\lvec(220 77)
\lvec(220 78)
\lvec(203 78)
\ifill f:0
\move(221 77)
\lvec(226 77)
\lvec(226 78)
\lvec(221 78)
\ifill f:0
\move(227 77)
\lvec(232 77)
\lvec(232 78)
\lvec(227 78)
\ifill f:0
\move(233 77)
\lvec(234 77)
\lvec(234 78)
\lvec(233 78)
\ifill f:0
\move(235 77)
\lvec(248 77)
\lvec(248 78)
\lvec(235 78)
\ifill f:0
\move(249 77)
\lvec(250 77)
\lvec(250 78)
\lvec(249 78)
\ifill f:0
\move(251 77)
\lvec(257 77)
\lvec(257 78)
\lvec(251 78)
\ifill f:0
\move(258 77)
\lvec(273 77)
\lvec(273 78)
\lvec(258 78)
\ifill f:0
\move(274 77)
\lvec(279 77)
\lvec(279 78)
\lvec(274 78)
\ifill f:0
\move(280 77)
\lvec(282 77)
\lvec(282 78)
\lvec(280 78)
\ifill f:0
\move(283 77)
\lvec(290 77)
\lvec(290 78)
\lvec(283 78)
\ifill f:0
\move(291 77)
\lvec(292 77)
\lvec(292 78)
\lvec(291 78)
\ifill f:0
\move(293 77)
\lvec(305 77)
\lvec(305 78)
\lvec(293 78)
\ifill f:0
\move(306 77)
\lvec(315 77)
\lvec(315 78)
\lvec(306 78)
\ifill f:0
\move(316 77)
\lvec(320 77)
\lvec(320 78)
\lvec(316 78)
\ifill f:0
\move(322 77)
\lvec(325 77)
\lvec(325 78)
\lvec(322 78)
\ifill f:0
\move(326 77)
\lvec(335 77)
\lvec(335 78)
\lvec(326 78)
\ifill f:0
\move(336 77)
\lvec(345 77)
\lvec(345 78)
\lvec(336 78)
\ifill f:0
\move(346 77)
\lvec(359 77)
\lvec(359 78)
\lvec(346 78)
\ifill f:0
\move(361 77)
\lvec(362 77)
\lvec(362 78)
\lvec(361 78)
\ifill f:0
\move(363 77)
\lvec(390 77)
\lvec(390 78)
\lvec(363 78)
\ifill f:0
\move(393 77)
\lvec(401 77)
\lvec(401 78)
\lvec(393 78)
\ifill f:0
\move(402 77)
\lvec(407 77)
\lvec(407 78)
\lvec(402 78)
\ifill f:0
\move(408 77)
\lvec(418 77)
\lvec(418 78)
\lvec(408 78)
\ifill f:0
\move(419 77)
\lvec(427 77)
\lvec(427 78)
\lvec(419 78)
\ifill f:0
\move(428 77)
\lvec(434 77)
\lvec(434 78)
\lvec(428 78)
\ifill f:0
\move(435 77)
\lvec(442 77)
\lvec(442 78)
\lvec(435 78)
\ifill f:0
\move(443 77)
\lvec(451 77)
\lvec(451 78)
\lvec(443 78)
\ifill f:0
\move(16 78)
\lvec(17 78)
\lvec(17 79)
\lvec(16 79)
\ifill f:0
\move(20 78)
\lvec(22 78)
\lvec(22 79)
\lvec(20 79)
\ifill f:0
\move(24 78)
\lvec(26 78)
\lvec(26 79)
\lvec(24 79)
\ifill f:0
\move(36 78)
\lvec(37 78)
\lvec(37 79)
\lvec(36 79)
\ifill f:0
\move(38 78)
\lvec(39 78)
\lvec(39 79)
\lvec(38 79)
\ifill f:0
\move(41 78)
\lvec(46 78)
\lvec(46 79)
\lvec(41 79)
\ifill f:0
\move(47 78)
\lvec(50 78)
\lvec(50 79)
\lvec(47 79)
\ifill f:0
\move(51 78)
\lvec(52 78)
\lvec(52 79)
\lvec(51 79)
\ifill f:0
\move(54 78)
\lvec(55 78)
\lvec(55 79)
\lvec(54 79)
\ifill f:0
\move(57 78)
\lvec(63 78)
\lvec(63 79)
\lvec(57 79)
\ifill f:0
\move(64 78)
\lvec(65 78)
\lvec(65 79)
\lvec(64 79)
\ifill f:0
\move(66 78)
\lvec(71 78)
\lvec(71 79)
\lvec(66 79)
\ifill f:0
\move(72 78)
\lvec(74 78)
\lvec(74 79)
\lvec(72 79)
\ifill f:0
\move(75 78)
\lvec(76 78)
\lvec(76 79)
\lvec(75 79)
\ifill f:0
\move(77 78)
\lvec(78 78)
\lvec(78 79)
\lvec(77 79)
\ifill f:0
\move(79 78)
\lvec(82 78)
\lvec(82 79)
\lvec(79 79)
\ifill f:0
\move(85 78)
\lvec(86 78)
\lvec(86 79)
\lvec(85 79)
\ifill f:0
\move(89 78)
\lvec(93 78)
\lvec(93 79)
\lvec(89 79)
\ifill f:0
\move(97 78)
\lvec(98 78)
\lvec(98 79)
\lvec(97 79)
\ifill f:0
\move(100 78)
\lvec(101 78)
\lvec(101 79)
\lvec(100 79)
\ifill f:0
\move(102 78)
\lvec(105 78)
\lvec(105 79)
\lvec(102 79)
\ifill f:0
\move(107 78)
\lvec(113 78)
\lvec(113 79)
\lvec(107 79)
\ifill f:0
\move(114 78)
\lvec(116 78)
\lvec(116 79)
\lvec(114 79)
\ifill f:0
\move(117 78)
\lvec(122 78)
\lvec(122 79)
\lvec(117 79)
\ifill f:0
\move(123 78)
\lvec(124 78)
\lvec(124 79)
\lvec(123 79)
\ifill f:0
\move(125 78)
\lvec(126 78)
\lvec(126 79)
\lvec(125 79)
\ifill f:0
\move(127 78)
\lvec(131 78)
\lvec(131 79)
\lvec(127 79)
\ifill f:0
\move(132 78)
\lvec(133 78)
\lvec(133 79)
\lvec(132 79)
\ifill f:0
\move(134 78)
\lvec(136 78)
\lvec(136 79)
\lvec(134 79)
\ifill f:0
\move(137 78)
\lvec(138 78)
\lvec(138 79)
\lvec(137 79)
\ifill f:0
\move(140 78)
\lvec(142 78)
\lvec(142 79)
\lvec(140 79)
\ifill f:0
\move(143 78)
\lvec(145 78)
\lvec(145 79)
\lvec(143 79)
\ifill f:0
\move(146 78)
\lvec(150 78)
\lvec(150 79)
\lvec(146 79)
\ifill f:0
\move(151 78)
\lvec(155 78)
\lvec(155 79)
\lvec(151 79)
\ifill f:0
\move(156 78)
\lvec(163 78)
\lvec(163 79)
\lvec(156 79)
\ifill f:0
\move(164 78)
\lvec(167 78)
\lvec(167 79)
\lvec(164 79)
\ifill f:0
\move(169 78)
\lvec(170 78)
\lvec(170 79)
\lvec(169 79)
\ifill f:0
\move(171 78)
\lvec(172 78)
\lvec(172 79)
\lvec(171 79)
\ifill f:0
\move(178 78)
\lvec(179 78)
\lvec(179 79)
\lvec(178 79)
\ifill f:0
\move(180 78)
\lvec(187 78)
\lvec(187 79)
\lvec(180 79)
\ifill f:0
\move(188 78)
\lvec(194 78)
\lvec(194 79)
\lvec(188 79)
\ifill f:0
\move(195 78)
\lvec(197 78)
\lvec(197 79)
\lvec(195 79)
\ifill f:0
\move(198 78)
\lvec(203 78)
\lvec(203 79)
\lvec(198 79)
\ifill f:0
\move(204 78)
\lvec(226 78)
\lvec(226 79)
\lvec(204 79)
\ifill f:0
\move(227 78)
\lvec(228 78)
\lvec(228 79)
\lvec(227 79)
\ifill f:0
\move(229 78)
\lvec(230 78)
\lvec(230 79)
\lvec(229 79)
\ifill f:0
\move(231 78)
\lvec(235 78)
\lvec(235 79)
\lvec(231 79)
\ifill f:0
\move(236 78)
\lvec(237 78)
\lvec(237 79)
\lvec(236 79)
\ifill f:0
\move(238 78)
\lvec(257 78)
\lvec(257 79)
\lvec(238 79)
\ifill f:0
\move(258 78)
\lvec(264 78)
\lvec(264 79)
\lvec(258 79)
\ifill f:0
\move(265 78)
\lvec(266 78)
\lvec(266 79)
\lvec(265 79)
\ifill f:0
\move(267 78)
\lvec(269 78)
\lvec(269 79)
\lvec(267 79)
\ifill f:0
\move(270 78)
\lvec(277 78)
\lvec(277 79)
\lvec(270 79)
\ifill f:0
\move(278 78)
\lvec(290 78)
\lvec(290 79)
\lvec(278 79)
\ifill f:0
\move(291 78)
\lvec(292 78)
\lvec(292 79)
\lvec(291 79)
\ifill f:0
\move(293 78)
\lvec(299 78)
\lvec(299 79)
\lvec(293 79)
\ifill f:0
\move(300 78)
\lvec(307 78)
\lvec(307 79)
\lvec(300 79)
\ifill f:0
\move(308 78)
\lvec(311 78)
\lvec(311 79)
\lvec(308 79)
\ifill f:0
\move(312 78)
\lvec(321 78)
\lvec(321 79)
\lvec(312 79)
\ifill f:0
\move(322 78)
\lvec(325 78)
\lvec(325 79)
\lvec(322 79)
\ifill f:0
\move(326 78)
\lvec(327 78)
\lvec(327 79)
\lvec(326 79)
\ifill f:0
\move(328 78)
\lvec(333 78)
\lvec(333 79)
\lvec(328 79)
\ifill f:0
\move(334 78)
\lvec(340 78)
\lvec(340 79)
\lvec(334 79)
\ifill f:0
\move(341 78)
\lvec(349 78)
\lvec(349 79)
\lvec(341 79)
\ifill f:0
\move(350 78)
\lvec(360 78)
\lvec(360 79)
\lvec(350 79)
\ifill f:0
\move(361 78)
\lvec(362 78)
\lvec(362 79)
\lvec(361 79)
\ifill f:0
\move(363 78)
\lvec(401 78)
\lvec(401 79)
\lvec(363 79)
\ifill f:0
\move(402 78)
\lvec(410 78)
\lvec(410 79)
\lvec(402 79)
\ifill f:0
\move(411 78)
\lvec(423 78)
\lvec(423 79)
\lvec(411 79)
\ifill f:0
\move(424 78)
\lvec(433 78)
\lvec(433 79)
\lvec(424 79)
\ifill f:0
\move(434 78)
\lvec(442 78)
\lvec(442 79)
\lvec(434 79)
\ifill f:0
\move(443 78)
\lvec(448 78)
\lvec(448 79)
\lvec(443 79)
\ifill f:0
\move(449 78)
\lvec(451 78)
\lvec(451 79)
\lvec(449 79)
\ifill f:0
\move(15 79)
\lvec(17 79)
\lvec(17 80)
\lvec(15 80)
\ifill f:0
\move(19 79)
\lvec(21 79)
\lvec(21 80)
\lvec(19 80)
\ifill f:0
\move(25 79)
\lvec(26 79)
\lvec(26 80)
\lvec(25 80)
\ifill f:0
\move(36 79)
\lvec(37 79)
\lvec(37 80)
\lvec(36 80)
\ifill f:0
\move(38 79)
\lvec(39 79)
\lvec(39 80)
\lvec(38 80)
\ifill f:0
\move(40 79)
\lvec(43 79)
\lvec(43 80)
\lvec(40 80)
\ifill f:0
\move(44 79)
\lvec(45 79)
\lvec(45 80)
\lvec(44 80)
\ifill f:0
\move(47 79)
\lvec(50 79)
\lvec(50 80)
\lvec(47 80)
\ifill f:0
\move(51 79)
\lvec(52 79)
\lvec(52 80)
\lvec(51 80)
\ifill f:0
\move(55 79)
\lvec(56 79)
\lvec(56 80)
\lvec(55 80)
\ifill f:0
\move(59 79)
\lvec(65 79)
\lvec(65 80)
\lvec(59 80)
\ifill f:0
\move(66 79)
\lvec(73 79)
\lvec(73 80)
\lvec(66 80)
\ifill f:0
\move(74 79)
\lvec(75 79)
\lvec(75 80)
\lvec(74 80)
\ifill f:0
\move(76 79)
\lvec(82 79)
\lvec(82 80)
\lvec(76 80)
\ifill f:0
\move(84 79)
\lvec(85 79)
\lvec(85 80)
\lvec(84 80)
\ifill f:0
\move(87 79)
\lvec(90 79)
\lvec(90 80)
\lvec(87 80)
\ifill f:0
\move(92 79)
\lvec(94 79)
\lvec(94 80)
\lvec(92 80)
\ifill f:0
\move(95 79)
\lvec(98 79)
\lvec(98 80)
\lvec(95 80)
\ifill f:0
\move(100 79)
\lvec(101 79)
\lvec(101 80)
\lvec(100 80)
\ifill f:0
\move(102 79)
\lvec(106 79)
\lvec(106 80)
\lvec(102 80)
\ifill f:0
\move(107 79)
\lvec(122 79)
\lvec(122 80)
\lvec(107 80)
\ifill f:0
\move(123 79)
\lvec(124 79)
\lvec(124 80)
\lvec(123 80)
\ifill f:0
\move(125 79)
\lvec(126 79)
\lvec(126 80)
\lvec(125 80)
\ifill f:0
\move(127 79)
\lvec(128 79)
\lvec(128 80)
\lvec(127 80)
\ifill f:0
\move(129 79)
\lvec(130 79)
\lvec(130 80)
\lvec(129 80)
\ifill f:0
\move(131 79)
\lvec(132 79)
\lvec(132 80)
\lvec(131 80)
\ifill f:0
\move(133 79)
\lvec(134 79)
\lvec(134 80)
\lvec(133 80)
\ifill f:0
\move(135 79)
\lvec(145 79)
\lvec(145 80)
\lvec(135 80)
\ifill f:0
\move(146 79)
\lvec(149 79)
\lvec(149 80)
\lvec(146 80)
\ifill f:0
\move(150 79)
\lvec(153 79)
\lvec(153 80)
\lvec(150 80)
\ifill f:0
\move(154 79)
\lvec(159 79)
\lvec(159 80)
\lvec(154 80)
\ifill f:0
\move(160 79)
\lvec(167 79)
\lvec(167 80)
\lvec(160 80)
\ifill f:0
\move(169 79)
\lvec(170 79)
\lvec(170 80)
\lvec(169 80)
\ifill f:0
\move(171 79)
\lvec(182 79)
\lvec(182 80)
\lvec(171 80)
\ifill f:0
\move(185 79)
\lvec(186 79)
\lvec(186 80)
\lvec(185 80)
\ifill f:0
\move(187 79)
\lvec(193 79)
\lvec(193 80)
\lvec(187 80)
\ifill f:0
\move(195 79)
\lvec(197 79)
\lvec(197 80)
\lvec(195 80)
\ifill f:0
\move(198 79)
\lvec(199 79)
\lvec(199 80)
\lvec(198 80)
\ifill f:0
\move(200 79)
\lvec(204 79)
\lvec(204 80)
\lvec(200 80)
\ifill f:0
\move(205 79)
\lvec(219 79)
\lvec(219 80)
\lvec(205 80)
\ifill f:0
\move(220 79)
\lvec(226 79)
\lvec(226 80)
\lvec(220 80)
\ifill f:0
\move(227 79)
\lvec(228 79)
\lvec(228 80)
\lvec(227 80)
\ifill f:0
\move(229 79)
\lvec(233 79)
\lvec(233 80)
\lvec(229 80)
\ifill f:0
\move(234 79)
\lvec(238 79)
\lvec(238 80)
\lvec(234 80)
\ifill f:0
\move(239 79)
\lvec(247 79)
\lvec(247 80)
\lvec(239 80)
\ifill f:0
\move(248 79)
\lvec(249 79)
\lvec(249 80)
\lvec(248 80)
\ifill f:0
\move(250 79)
\lvec(251 79)
\lvec(251 80)
\lvec(250 80)
\ifill f:0
\move(252 79)
\lvec(255 79)
\lvec(255 80)
\lvec(252 80)
\ifill f:0
\move(256 79)
\lvec(257 79)
\lvec(257 80)
\lvec(256 80)
\ifill f:0
\move(258 79)
\lvec(259 79)
\lvec(259 80)
\lvec(258 80)
\ifill f:0
\move(260 79)
\lvec(280 79)
\lvec(280 80)
\lvec(260 80)
\ifill f:0
\move(281 79)
\lvec(283 79)
\lvec(283 80)
\lvec(281 80)
\ifill f:0
\move(284 79)
\lvec(286 79)
\lvec(286 80)
\lvec(284 80)
\ifill f:0
\move(287 79)
\lvec(290 79)
\lvec(290 80)
\lvec(287 80)
\ifill f:0
\move(291 79)
\lvec(292 79)
\lvec(292 80)
\lvec(291 80)
\ifill f:0
\move(293 79)
\lvec(295 79)
\lvec(295 80)
\lvec(293 80)
\ifill f:0
\move(296 79)
\lvec(298 79)
\lvec(298 80)
\lvec(296 80)
\ifill f:0
\move(299 79)
\lvec(305 79)
\lvec(305 80)
\lvec(299 80)
\ifill f:0
\move(306 79)
\lvec(309 79)
\lvec(309 80)
\lvec(306 80)
\ifill f:0
\move(310 79)
\lvec(313 79)
\lvec(313 80)
\lvec(310 80)
\ifill f:0
\move(314 79)
\lvec(317 79)
\lvec(317 80)
\lvec(314 80)
\ifill f:0
\move(318 79)
\lvec(321 79)
\lvec(321 80)
\lvec(318 80)
\ifill f:0
\move(322 79)
\lvec(325 79)
\lvec(325 80)
\lvec(322 80)
\ifill f:0
\move(326 79)
\lvec(337 79)
\lvec(337 80)
\lvec(326 80)
\ifill f:0
\move(338 79)
\lvec(344 79)
\lvec(344 80)
\lvec(338 80)
\ifill f:0
\move(345 79)
\lvec(362 79)
\lvec(362 80)
\lvec(345 80)
\ifill f:0
\move(363 79)
\lvec(375 79)
\lvec(375 80)
\lvec(363 80)
\ifill f:0
\move(376 79)
\lvec(401 79)
\lvec(401 80)
\lvec(376 80)
\ifill f:0
\move(402 79)
\lvec(416 79)
\lvec(416 80)
\lvec(402 80)
\ifill f:0
\move(417 79)
\lvec(442 79)
\lvec(442 80)
\lvec(417 80)
\ifill f:0
\move(443 79)
\lvec(450 79)
\lvec(450 80)
\lvec(443 80)
\ifill f:0
\move(15 80)
\lvec(17 80)
\lvec(17 81)
\lvec(15 81)
\ifill f:0
\move(18 80)
\lvec(19 80)
\lvec(19 81)
\lvec(18 81)
\ifill f:0
\move(20 80)
\lvec(23 80)
\lvec(23 81)
\lvec(20 81)
\ifill f:0
\move(25 80)
\lvec(26 80)
\lvec(26 81)
\lvec(25 81)
\ifill f:0
\move(36 80)
\lvec(37 80)
\lvec(37 81)
\lvec(36 81)
\ifill f:0
\move(38 80)
\lvec(41 80)
\lvec(41 81)
\lvec(38 81)
\ifill f:0
\move(42 80)
\lvec(45 80)
\lvec(45 81)
\lvec(42 81)
\ifill f:0
\move(47 80)
\lvec(50 80)
\lvec(50 81)
\lvec(47 81)
\ifill f:0
\move(51 80)
\lvec(52 80)
\lvec(52 81)
\lvec(51 81)
\ifill f:0
\move(54 80)
\lvec(55 80)
\lvec(55 81)
\lvec(54 81)
\ifill f:0
\move(57 80)
\lvec(58 80)
\lvec(58 81)
\lvec(57 81)
\ifill f:0
\move(60 80)
\lvec(65 80)
\lvec(65 81)
\lvec(60 81)
\ifill f:0
\move(66 80)
\lvec(70 80)
\lvec(70 81)
\lvec(66 81)
\ifill f:0
\move(72 80)
\lvec(74 80)
\lvec(74 81)
\lvec(72 81)
\ifill f:0
\move(76 80)
\lvec(79 80)
\lvec(79 81)
\lvec(76 81)
\ifill f:0
\move(80 80)
\lvec(82 80)
\lvec(82 81)
\lvec(80 81)
\ifill f:0
\move(84 80)
\lvec(85 80)
\lvec(85 81)
\lvec(84 81)
\ifill f:0
\move(86 80)
\lvec(89 80)
\lvec(89 81)
\lvec(86 81)
\ifill f:0
\move(90 80)
\lvec(93 80)
\lvec(93 81)
\lvec(90 81)
\ifill f:0
\move(95 80)
\lvec(98 80)
\lvec(98 81)
\lvec(95 81)
\ifill f:0
\move(99 80)
\lvec(101 80)
\lvec(101 81)
\lvec(99 81)
\ifill f:0
\move(102 80)
\lvec(109 80)
\lvec(109 81)
\lvec(102 81)
\ifill f:0
\move(110 80)
\lvec(114 80)
\lvec(114 81)
\lvec(110 81)
\ifill f:0
\move(115 80)
\lvec(118 80)
\lvec(118 81)
\lvec(115 81)
\ifill f:0
\move(119 80)
\lvec(122 80)
\lvec(122 81)
\lvec(119 81)
\ifill f:0
\move(123 80)
\lvec(124 80)
\lvec(124 81)
\lvec(123 81)
\ifill f:0
\move(125 80)
\lvec(127 80)
\lvec(127 81)
\lvec(125 81)
\ifill f:0
\move(128 80)
\lvec(129 80)
\lvec(129 81)
\lvec(128 81)
\ifill f:0
\move(130 80)
\lvec(131 80)
\lvec(131 81)
\lvec(130 81)
\ifill f:0
\move(132 80)
\lvec(133 80)
\lvec(133 81)
\lvec(132 81)
\ifill f:0
\move(134 80)
\lvec(135 80)
\lvec(135 81)
\lvec(134 81)
\ifill f:0
\move(136 80)
\lvec(140 80)
\lvec(140 81)
\lvec(136 81)
\ifill f:0
\move(141 80)
\lvec(145 80)
\lvec(145 81)
\lvec(141 81)
\ifill f:0
\move(146 80)
\lvec(152 80)
\lvec(152 81)
\lvec(146 81)
\ifill f:0
\move(153 80)
\lvec(161 80)
\lvec(161 81)
\lvec(153 81)
\ifill f:0
\move(162 80)
\lvec(168 80)
\lvec(168 81)
\lvec(162 81)
\ifill f:0
\move(169 80)
\lvec(170 80)
\lvec(170 81)
\lvec(169 81)
\ifill f:0
\move(171 80)
\lvec(192 80)
\lvec(192 81)
\lvec(171 81)
\ifill f:0
\move(193 80)
\lvec(194 80)
\lvec(194 81)
\lvec(193 81)
\ifill f:0
\move(195 80)
\lvec(197 80)
\lvec(197 81)
\lvec(195 81)
\ifill f:0
\move(198 80)
\lvec(200 80)
\lvec(200 81)
\lvec(198 81)
\ifill f:0
\move(202 80)
\lvec(211 80)
\lvec(211 81)
\lvec(202 81)
\ifill f:0
\move(212 80)
\lvec(226 80)
\lvec(226 81)
\lvec(212 81)
\ifill f:0
\move(227 80)
\lvec(231 80)
\lvec(231 81)
\lvec(227 81)
\ifill f:0
\move(232 80)
\lvec(234 80)
\lvec(234 81)
\lvec(232 81)
\ifill f:0
\move(235 80)
\lvec(239 80)
\lvec(239 81)
\lvec(235 81)
\ifill f:0
\move(240 80)
\lvec(242 80)
\lvec(242 81)
\lvec(240 81)
\ifill f:0
\move(243 80)
\lvec(244 80)
\lvec(244 81)
\lvec(243 81)
\ifill f:0
\move(245 80)
\lvec(251 80)
\lvec(251 81)
\lvec(245 81)
\ifill f:0
\move(252 80)
\lvec(253 80)
\lvec(253 81)
\lvec(252 81)
\ifill f:0
\move(254 80)
\lvec(255 80)
\lvec(255 81)
\lvec(254 81)
\ifill f:0
\move(256 80)
\lvec(257 80)
\lvec(257 81)
\lvec(256 81)
\ifill f:0
\move(258 80)
\lvec(261 80)
\lvec(261 81)
\lvec(258 81)
\ifill f:0
\move(262 80)
\lvec(263 80)
\lvec(263 81)
\lvec(262 81)
\ifill f:0
\move(264 80)
\lvec(274 80)
\lvec(274 81)
\lvec(264 81)
\ifill f:0
\move(275 80)
\lvec(281 80)
\lvec(281 81)
\lvec(275 81)
\ifill f:0
\move(282 80)
\lvec(290 80)
\lvec(290 81)
\lvec(282 81)
\ifill f:0
\move(291 80)
\lvec(307 80)
\lvec(307 81)
\lvec(291 81)
\ifill f:0
\move(308 80)
\lvec(310 80)
\lvec(310 81)
\lvec(308 81)
\ifill f:0
\move(311 80)
\lvec(314 80)
\lvec(314 81)
\lvec(311 81)
\ifill f:0
\move(315 80)
\lvec(318 80)
\lvec(318 81)
\lvec(315 81)
\ifill f:0
\move(319 80)
\lvec(322 80)
\lvec(322 81)
\lvec(319 81)
\ifill f:0
\move(323 80)
\lvec(325 80)
\lvec(325 81)
\lvec(323 81)
\ifill f:0
\move(326 80)
\lvec(335 80)
\lvec(335 81)
\lvec(326 81)
\ifill f:0
\move(336 80)
\lvec(341 80)
\lvec(341 81)
\lvec(336 81)
\ifill f:0
\move(342 80)
\lvec(362 80)
\lvec(362 81)
\lvec(342 81)
\ifill f:0
\move(363 80)
\lvec(370 80)
\lvec(370 81)
\lvec(363 81)
\ifill f:0
\move(371 80)
\lvec(383 80)
\lvec(383 81)
\lvec(371 81)
\ifill f:0
\move(386 80)
\lvec(401 80)
\lvec(401 81)
\lvec(386 81)
\ifill f:0
\move(402 80)
\lvec(426 80)
\lvec(426 81)
\lvec(402 81)
\ifill f:0
\move(428 80)
\lvec(442 80)
\lvec(442 81)
\lvec(428 81)
\ifill f:0
\move(444 80)
\lvec(451 80)
\lvec(451 81)
\lvec(444 81)
\ifill f:0
\move(15 81)
\lvec(17 81)
\lvec(17 82)
\lvec(15 82)
\ifill f:0
\move(18 81)
\lvec(19 81)
\lvec(19 82)
\lvec(18 82)
\ifill f:0
\move(23 81)
\lvec(24 81)
\lvec(24 82)
\lvec(23 82)
\ifill f:0
\move(25 81)
\lvec(26 81)
\lvec(26 82)
\lvec(25 82)
\ifill f:0
\move(36 81)
\lvec(37 81)
\lvec(37 82)
\lvec(36 82)
\ifill f:0
\move(41 81)
\lvec(43 81)
\lvec(43 82)
\lvec(41 82)
\ifill f:0
\move(44 81)
\lvec(46 81)
\lvec(46 82)
\lvec(44 82)
\ifill f:0
\move(47 81)
\lvec(50 81)
\lvec(50 82)
\lvec(47 82)
\ifill f:0
\move(52 81)
\lvec(53 81)
\lvec(53 82)
\lvec(52 82)
\ifill f:0
\move(54 81)
\lvec(55 81)
\lvec(55 82)
\lvec(54 82)
\ifill f:0
\move(58 81)
\lvec(60 81)
\lvec(60 82)
\lvec(58 82)
\ifill f:0
\move(61 81)
\lvec(65 81)
\lvec(65 82)
\lvec(61 82)
\ifill f:0
\move(66 81)
\lvec(73 81)
\lvec(73 82)
\lvec(66 82)
\ifill f:0
\move(75 81)
\lvec(76 81)
\lvec(76 82)
\lvec(75 82)
\ifill f:0
\move(77 81)
\lvec(79 81)
\lvec(79 82)
\lvec(77 82)
\ifill f:0
\move(80 81)
\lvec(82 81)
\lvec(82 82)
\lvec(80 82)
\ifill f:0
\move(84 81)
\lvec(85 81)
\lvec(85 82)
\lvec(84 82)
\ifill f:0
\move(86 81)
\lvec(87 81)
\lvec(87 82)
\lvec(86 82)
\ifill f:0
\move(88 81)
\lvec(90 81)
\lvec(90 82)
\lvec(88 82)
\ifill f:0
\move(92 81)
\lvec(93 81)
\lvec(93 82)
\lvec(92 82)
\ifill f:0
\move(96 81)
\lvec(101 81)
\lvec(101 82)
\lvec(96 82)
\ifill f:0
\move(102 81)
\lvec(112 81)
\lvec(112 82)
\lvec(102 82)
\ifill f:0
\move(114 81)
\lvec(122 81)
\lvec(122 82)
\lvec(114 82)
\ifill f:0
\move(123 81)
\lvec(125 81)
\lvec(125 82)
\lvec(123 82)
\ifill f:0
\move(126 81)
\lvec(127 81)
\lvec(127 82)
\lvec(126 82)
\ifill f:0
\move(128 81)
\lvec(130 81)
\lvec(130 82)
\lvec(128 82)
\ifill f:0
\move(131 81)
\lvec(132 81)
\lvec(132 82)
\lvec(131 82)
\ifill f:0
\move(133 81)
\lvec(134 81)
\lvec(134 82)
\lvec(133 82)
\ifill f:0
\move(135 81)
\lvec(136 81)
\lvec(136 82)
\lvec(135 82)
\ifill f:0
\move(137 81)
\lvec(138 81)
\lvec(138 82)
\lvec(137 82)
\ifill f:0
\move(139 81)
\lvec(145 81)
\lvec(145 82)
\lvec(139 82)
\ifill f:0
\move(146 81)
\lvec(148 81)
\lvec(148 82)
\lvec(146 82)
\ifill f:0
\move(149 81)
\lvec(151 81)
\lvec(151 82)
\lvec(149 82)
\ifill f:0
\move(152 81)
\lvec(168 81)
\lvec(168 82)
\lvec(152 82)
\ifill f:0
\move(169 81)
\lvec(170 81)
\lvec(170 82)
\lvec(169 82)
\ifill f:0
\move(171 81)
\lvec(176 81)
\lvec(176 82)
\lvec(171 82)
\ifill f:0
\move(177 81)
\lvec(179 81)
\lvec(179 82)
\lvec(177 82)
\ifill f:0
\move(181 81)
\lvec(182 81)
\lvec(182 82)
\lvec(181 82)
\ifill f:0
\move(193 81)
\lvec(194 81)
\lvec(194 82)
\lvec(193 82)
\ifill f:0
\move(195 81)
\lvec(197 81)
\lvec(197 82)
\lvec(195 82)
\ifill f:0
\move(198 81)
\lvec(201 81)
\lvec(201 82)
\lvec(198 82)
\ifill f:0
\move(202 81)
\lvec(226 81)
\lvec(226 82)
\lvec(202 82)
\ifill f:0
\move(227 81)
\lvec(229 81)
\lvec(229 82)
\lvec(227 82)
\ifill f:0
\move(230 81)
\lvec(232 81)
\lvec(232 82)
\lvec(230 82)
\ifill f:0
\move(233 81)
\lvec(235 81)
\lvec(235 82)
\lvec(233 82)
\ifill f:0
\move(236 81)
\lvec(243 81)
\lvec(243 82)
\lvec(236 82)
\ifill f:0
\move(244 81)
\lvec(248 81)
\lvec(248 82)
\lvec(244 82)
\ifill f:0
\move(249 81)
\lvec(255 81)
\lvec(255 82)
\lvec(249 82)
\ifill f:0
\move(256 81)
\lvec(257 81)
\lvec(257 82)
\lvec(256 82)
\ifill f:0
\move(258 81)
\lvec(275 81)
\lvec(275 82)
\lvec(258 82)
\ifill f:0
\move(276 81)
\lvec(277 81)
\lvec(277 82)
\lvec(276 82)
\ifill f:0
\move(278 81)
\lvec(279 81)
\lvec(279 82)
\lvec(278 82)
\ifill f:0
\move(280 81)
\lvec(284 81)
\lvec(284 82)
\lvec(280 82)
\ifill f:0
\move(285 81)
\lvec(290 81)
\lvec(290 82)
\lvec(285 82)
\ifill f:0
\move(292 81)
\lvec(294 81)
\lvec(294 82)
\lvec(292 82)
\ifill f:0
\move(295 81)
\lvec(315 81)
\lvec(315 82)
\lvec(295 82)
\ifill f:0
\move(316 81)
\lvec(322 81)
\lvec(322 82)
\lvec(316 82)
\ifill f:0
\move(323 81)
\lvec(325 81)
\lvec(325 82)
\lvec(323 82)
\ifill f:0
\move(327 81)
\lvec(330 81)
\lvec(330 82)
\lvec(327 82)
\ifill f:0
\move(331 81)
\lvec(362 81)
\lvec(362 82)
\lvec(331 82)
\ifill f:0
\move(363 81)
\lvec(368 81)
\lvec(368 82)
\lvec(363 82)
\ifill f:0
\move(369 81)
\lvec(377 81)
\lvec(377 82)
\lvec(369 82)
\ifill f:0
\move(378 81)
\lvec(389 81)
\lvec(389 82)
\lvec(378 82)
\ifill f:0
\move(391 81)
\lvec(401 81)
\lvec(401 82)
\lvec(391 82)
\ifill f:0
\move(402 81)
\lvec(442 81)
\lvec(442 82)
\lvec(402 82)
\ifill f:0
\move(444 81)
\lvec(451 81)
\lvec(451 82)
\lvec(444 82)
\ifill f:0
\move(15 82)
\lvec(17 82)
\lvec(17 83)
\lvec(15 83)
\ifill f:0
\move(19 82)
\lvec(21 82)
\lvec(21 83)
\lvec(19 83)
\ifill f:0
\move(23 82)
\lvec(24 82)
\lvec(24 83)
\lvec(23 83)
\ifill f:0
\move(25 82)
\lvec(26 82)
\lvec(26 83)
\lvec(25 83)
\ifill f:0
\move(36 82)
\lvec(37 82)
\lvec(37 83)
\lvec(36 83)
\ifill f:0
\move(38 82)
\lvec(39 82)
\lvec(39 83)
\lvec(38 83)
\ifill f:0
\move(40 82)
\lvec(41 82)
\lvec(41 83)
\lvec(40 83)
\ifill f:0
\move(42 82)
\lvec(45 82)
\lvec(45 83)
\lvec(42 83)
\ifill f:0
\move(46 82)
\lvec(50 82)
\lvec(50 83)
\lvec(46 83)
\ifill f:0
\move(57 82)
\lvec(58 82)
\lvec(58 83)
\lvec(57 83)
\ifill f:0
\move(59 82)
\lvec(61 82)
\lvec(61 83)
\lvec(59 83)
\ifill f:0
\move(62 82)
\lvec(65 82)
\lvec(65 83)
\lvec(62 83)
\ifill f:0
\move(66 82)
\lvec(71 82)
\lvec(71 83)
\lvec(66 83)
\ifill f:0
\move(72 82)
\lvec(75 82)
\lvec(75 83)
\lvec(72 83)
\ifill f:0
\move(76 82)
\lvec(79 82)
\lvec(79 83)
\lvec(76 83)
\ifill f:0
\move(80 82)
\lvec(82 82)
\lvec(82 83)
\lvec(80 83)
\ifill f:0
\move(83 82)
\lvec(84 82)
\lvec(84 83)
\lvec(83 83)
\ifill f:0
\move(85 82)
\lvec(86 82)
\lvec(86 83)
\lvec(85 83)
\ifill f:0
\move(88 82)
\lvec(89 82)
\lvec(89 83)
\lvec(88 83)
\ifill f:0
\move(90 82)
\lvec(92 82)
\lvec(92 83)
\lvec(90 83)
\ifill f:0
\move(93 82)
\lvec(94 82)
\lvec(94 83)
\lvec(93 83)
\ifill f:0
\move(97 82)
\lvec(101 82)
\lvec(101 83)
\lvec(97 83)
\ifill f:0
\move(102 82)
\lvec(106 82)
\lvec(106 83)
\lvec(102 83)
\ifill f:0
\move(108 82)
\lvec(109 82)
\lvec(109 83)
\lvec(108 83)
\ifill f:0
\move(110 82)
\lvec(111 82)
\lvec(111 83)
\lvec(110 83)
\ifill f:0
\move(112 82)
\lvec(117 82)
\lvec(117 83)
\lvec(112 83)
\ifill f:0
\move(118 82)
\lvec(122 82)
\lvec(122 83)
\lvec(118 83)
\ifill f:0
\move(123 82)
\lvec(125 82)
\lvec(125 83)
\lvec(123 83)
\ifill f:0
\move(126 82)
\lvec(128 82)
\lvec(128 83)
\lvec(126 83)
\ifill f:0
\move(129 82)
\lvec(131 82)
\lvec(131 83)
\lvec(129 83)
\ifill f:0
\move(132 82)
\lvec(133 82)
\lvec(133 83)
\lvec(132 83)
\ifill f:0
\move(134 82)
\lvec(138 82)
\lvec(138 83)
\lvec(134 83)
\ifill f:0
\move(139 82)
\lvec(145 82)
\lvec(145 83)
\lvec(139 83)
\ifill f:0
\move(146 82)
\lvec(150 82)
\lvec(150 83)
\lvec(146 83)
\ifill f:0
\move(151 82)
\lvec(153 82)
\lvec(153 83)
\lvec(151 83)
\ifill f:0
\move(154 82)
\lvec(159 82)
\lvec(159 83)
\lvec(154 83)
\ifill f:0
\move(160 82)
\lvec(163 82)
\lvec(163 83)
\lvec(160 83)
\ifill f:0
\move(165 82)
\lvec(170 82)
\lvec(170 83)
\lvec(165 83)
\ifill f:0
\move(171 82)
\lvec(174 82)
\lvec(174 83)
\lvec(171 83)
\ifill f:0
\move(177 82)
\lvec(191 82)
\lvec(191 83)
\lvec(177 83)
\ifill f:0
\move(192 82)
\lvec(197 82)
\lvec(197 83)
\lvec(192 83)
\ifill f:0
\move(198 82)
\lvec(203 82)
\lvec(203 83)
\lvec(198 83)
\ifill f:0
\move(204 82)
\lvec(211 82)
\lvec(211 83)
\lvec(204 83)
\ifill f:0
\move(212 82)
\lvec(216 82)
\lvec(216 83)
\lvec(212 83)
\ifill f:0
\move(217 82)
\lvec(226 82)
\lvec(226 83)
\lvec(217 83)
\ifill f:0
\move(227 82)
\lvec(229 82)
\lvec(229 83)
\lvec(227 83)
\ifill f:0
\move(230 82)
\lvec(233 82)
\lvec(233 83)
\lvec(230 83)
\ifill f:0
\move(234 82)
\lvec(239 82)
\lvec(239 83)
\lvec(234 83)
\ifill f:0
\move(240 82)
\lvec(242 82)
\lvec(242 83)
\lvec(240 83)
\ifill f:0
\move(243 82)
\lvec(250 82)
\lvec(250 83)
\lvec(243 83)
\ifill f:0
\move(251 82)
\lvec(255 82)
\lvec(255 83)
\lvec(251 83)
\ifill f:0
\move(256 82)
\lvec(257 82)
\lvec(257 83)
\lvec(256 83)
\ifill f:0
\move(258 82)
\lvec(264 82)
\lvec(264 83)
\lvec(258 83)
\ifill f:0
\move(265 82)
\lvec(280 82)
\lvec(280 83)
\lvec(265 83)
\ifill f:0
\move(281 82)
\lvec(290 82)
\lvec(290 83)
\lvec(281 83)
\ifill f:0
\move(292 82)
\lvec(296 82)
\lvec(296 83)
\lvec(292 83)
\ifill f:0
\move(297 82)
\lvec(298 82)
\lvec(298 83)
\lvec(297 83)
\ifill f:0
\move(299 82)
\lvec(307 82)
\lvec(307 83)
\lvec(299 83)
\ifill f:0
\move(308 82)
\lvec(319 82)
\lvec(319 83)
\lvec(308 83)
\ifill f:0
\move(320 82)
\lvec(322 82)
\lvec(322 83)
\lvec(320 83)
\ifill f:0
\move(323 82)
\lvec(325 82)
\lvec(325 83)
\lvec(323 83)
\ifill f:0
\move(326 82)
\lvec(329 82)
\lvec(329 83)
\lvec(326 83)
\ifill f:0
\move(330 82)
\lvec(333 82)
\lvec(333 83)
\lvec(330 83)
\ifill f:0
\move(334 82)
\lvec(337 82)
\lvec(337 83)
\lvec(334 83)
\ifill f:0
\move(338 82)
\lvec(341 82)
\lvec(341 83)
\lvec(338 83)
\ifill f:0
\move(342 82)
\lvec(362 82)
\lvec(362 83)
\lvec(342 83)
\ifill f:0
\move(363 82)
\lvec(367 82)
\lvec(367 83)
\lvec(363 83)
\ifill f:0
\move(368 82)
\lvec(382 82)
\lvec(382 83)
\lvec(368 83)
\ifill f:0
\move(383 82)
\lvec(392 82)
\lvec(392 83)
\lvec(383 83)
\ifill f:0
\move(394 82)
\lvec(401 82)
\lvec(401 83)
\lvec(394 83)
\ifill f:0
\move(402 82)
\lvec(410 82)
\lvec(410 83)
\lvec(402 83)
\ifill f:0
\move(411 82)
\lvec(442 82)
\lvec(442 83)
\lvec(411 83)
\ifill f:0
\move(446 82)
\lvec(451 82)
\lvec(451 83)
\lvec(446 83)
\ifill f:0
\move(16 83)
\lvec(17 83)
\lvec(17 84)
\lvec(16 84)
\ifill f:0
\move(22 83)
\lvec(23 83)
\lvec(23 84)
\lvec(22 84)
\ifill f:0
\move(24 83)
\lvec(26 83)
\lvec(26 84)
\lvec(24 84)
\ifill f:0
\move(36 83)
\lvec(37 83)
\lvec(37 84)
\lvec(36 84)
\ifill f:0
\move(38 83)
\lvec(39 83)
\lvec(39 84)
\lvec(38 84)
\ifill f:0
\move(40 83)
\lvec(45 83)
\lvec(45 84)
\lvec(40 84)
\ifill f:0
\move(49 83)
\lvec(50 83)
\lvec(50 84)
\lvec(49 84)
\ifill f:0
\move(54 83)
\lvec(55 83)
\lvec(55 84)
\lvec(54 84)
\ifill f:0
\move(60 83)
\lvec(62 83)
\lvec(62 84)
\lvec(60 84)
\ifill f:0
\move(63 83)
\lvec(65 83)
\lvec(65 84)
\lvec(63 84)
\ifill f:0
\move(66 83)
\lvec(67 83)
\lvec(67 84)
\lvec(66 84)
\ifill f:0
\move(68 83)
\lvec(73 83)
\lvec(73 84)
\lvec(68 84)
\ifill f:0
\move(76 83)
\lvec(78 83)
\lvec(78 84)
\lvec(76 84)
\ifill f:0
\move(79 83)
\lvec(82 83)
\lvec(82 84)
\lvec(79 84)
\ifill f:0
\move(83 83)
\lvec(84 83)
\lvec(84 84)
\lvec(83 84)
\ifill f:0
\move(87 83)
\lvec(88 83)
\lvec(88 84)
\lvec(87 84)
\ifill f:0
\move(89 83)
\lvec(90 83)
\lvec(90 84)
\lvec(89 84)
\ifill f:0
\move(91 83)
\lvec(93 83)
\lvec(93 84)
\lvec(91 84)
\ifill f:0
\move(94 83)
\lvec(96 83)
\lvec(96 84)
\lvec(94 84)
\ifill f:0
\move(97 83)
\lvec(101 83)
\lvec(101 84)
\lvec(97 84)
\ifill f:0
\move(102 83)
\lvec(103 83)
\lvec(103 84)
\lvec(102 84)
\ifill f:0
\move(104 83)
\lvec(115 83)
\lvec(115 84)
\lvec(104 84)
\ifill f:0
\move(116 83)
\lvec(122 83)
\lvec(122 84)
\lvec(116 84)
\ifill f:0
\move(123 83)
\lvec(126 83)
\lvec(126 84)
\lvec(123 84)
\ifill f:0
\move(127 83)
\lvec(129 83)
\lvec(129 84)
\lvec(127 84)
\ifill f:0
\move(130 83)
\lvec(132 83)
\lvec(132 84)
\lvec(130 84)
\ifill f:0
\move(133 83)
\lvec(135 83)
\lvec(135 84)
\lvec(133 84)
\ifill f:0
\move(136 83)
\lvec(137 83)
\lvec(137 84)
\lvec(136 84)
\ifill f:0
\move(138 83)
\lvec(139 83)
\lvec(139 84)
\lvec(138 84)
\ifill f:0
\move(140 83)
\lvec(141 83)
\lvec(141 84)
\lvec(140 84)
\ifill f:0
\move(142 83)
\lvec(143 83)
\lvec(143 84)
\lvec(142 84)
\ifill f:0
\move(144 83)
\lvec(145 83)
\lvec(145 84)
\lvec(144 84)
\ifill f:0
\move(146 83)
\lvec(147 83)
\lvec(147 84)
\lvec(146 84)
\ifill f:0
\move(148 83)
\lvec(152 83)
\lvec(152 84)
\lvec(148 84)
\ifill f:0
\move(153 83)
\lvec(161 83)
\lvec(161 84)
\lvec(153 84)
\ifill f:0
\move(162 83)
\lvec(170 83)
\lvec(170 84)
\lvec(162 84)
\ifill f:0
\move(171 83)
\lvec(173 83)
\lvec(173 84)
\lvec(171 84)
\ifill f:0
\move(175 83)
\lvec(180 83)
\lvec(180 84)
\lvec(175 84)
\ifill f:0
\move(181 83)
\lvec(182 83)
\lvec(182 84)
\lvec(181 84)
\ifill f:0
\move(183 83)
\lvec(197 83)
\lvec(197 84)
\lvec(183 84)
\ifill f:0
\move(198 83)
\lvec(214 83)
\lvec(214 84)
\lvec(198 84)
\ifill f:0
\move(215 83)
\lvec(220 83)
\lvec(220 84)
\lvec(215 84)
\ifill f:0
\move(221 83)
\lvec(226 83)
\lvec(226 84)
\lvec(221 84)
\ifill f:0
\move(227 83)
\lvec(234 83)
\lvec(234 84)
\lvec(227 84)
\ifill f:0
\move(235 83)
\lvec(237 83)
\lvec(237 84)
\lvec(235 84)
\ifill f:0
\move(238 83)
\lvec(250 83)
\lvec(250 84)
\lvec(238 84)
\ifill f:0
\move(251 83)
\lvec(252 83)
\lvec(252 84)
\lvec(251 84)
\ifill f:0
\move(253 83)
\lvec(255 83)
\lvec(255 84)
\lvec(253 84)
\ifill f:0
\move(256 83)
\lvec(257 83)
\lvec(257 84)
\lvec(256 84)
\ifill f:0
\move(258 83)
\lvec(269 83)
\lvec(269 84)
\lvec(258 84)
\ifill f:0
\move(270 83)
\lvec(290 83)
\lvec(290 84)
\lvec(270 84)
\ifill f:0
\move(292 83)
\lvec(293 83)
\lvec(293 84)
\lvec(292 84)
\ifill f:0
\move(294 83)
\lvec(298 83)
\lvec(298 84)
\lvec(294 84)
\ifill f:0
\move(299 83)
\lvec(303 83)
\lvec(303 84)
\lvec(299 84)
\ifill f:0
\move(304 83)
\lvec(308 83)
\lvec(308 84)
\lvec(304 84)
\ifill f:0
\move(309 83)
\lvec(319 83)
\lvec(319 84)
\lvec(309 84)
\ifill f:0
\move(320 83)
\lvec(322 83)
\lvec(322 84)
\lvec(320 84)
\ifill f:0
\move(323 83)
\lvec(325 83)
\lvec(325 84)
\lvec(323 84)
\ifill f:0
\move(326 83)
\lvec(356 83)
\lvec(356 84)
\lvec(326 84)
\ifill f:0
\move(357 83)
\lvec(362 83)
\lvec(362 84)
\lvec(357 84)
\ifill f:0
\move(363 83)
\lvec(366 83)
\lvec(366 84)
\lvec(363 84)
\ifill f:0
\move(367 83)
\lvec(372 83)
\lvec(372 84)
\lvec(367 84)
\ifill f:0
\move(373 83)
\lvec(386 83)
\lvec(386 84)
\lvec(373 84)
\ifill f:0
\move(387 83)
\lvec(394 83)
\lvec(394 84)
\lvec(387 84)
\ifill f:0
\move(395 83)
\lvec(401 83)
\lvec(401 84)
\lvec(395 84)
\ifill f:0
\move(402 83)
\lvec(427 83)
\lvec(427 84)
\lvec(402 84)
\ifill f:0
\move(428 83)
\lvec(430 83)
\lvec(430 84)
\lvec(428 84)
\ifill f:0
\move(436 83)
\lvec(437 83)
\lvec(437 84)
\lvec(436 84)
\ifill f:0
\move(441 83)
\lvec(442 83)
\lvec(442 84)
\lvec(441 84)
\ifill f:0
\move(446 83)
\lvec(451 83)
\lvec(451 84)
\lvec(446 84)
\ifill f:0
\move(16 84)
\lvec(17 84)
\lvec(17 85)
\lvec(16 85)
\ifill f:0
\move(20 84)
\lvec(21 84)
\lvec(21 85)
\lvec(20 85)
\ifill f:0
\move(24 84)
\lvec(26 84)
\lvec(26 85)
\lvec(24 85)
\ifill f:0
\move(36 84)
\lvec(37 84)
\lvec(37 85)
\lvec(36 85)
\ifill f:0
\move(38 84)
\lvec(40 84)
\lvec(40 85)
\lvec(38 85)
\ifill f:0
\move(41 84)
\lvec(46 84)
\lvec(46 85)
\lvec(41 85)
\ifill f:0
\move(49 84)
\lvec(50 84)
\lvec(50 85)
\lvec(49 85)
\ifill f:0
\move(54 84)
\lvec(55 84)
\lvec(55 85)
\lvec(54 85)
\ifill f:0
\move(57 84)
\lvec(62 84)
\lvec(62 85)
\lvec(57 85)
\ifill f:0
\move(63 84)
\lvec(65 84)
\lvec(65 85)
\lvec(63 85)
\ifill f:0
\move(67 84)
\lvec(77 84)
\lvec(77 85)
\lvec(67 85)
\ifill f:0
\move(79 84)
\lvec(82 84)
\lvec(82 85)
\lvec(79 85)
\ifill f:0
\move(83 84)
\lvec(84 84)
\lvec(84 85)
\lvec(83 85)
\ifill f:0
\move(86 84)
\lvec(87 84)
\lvec(87 85)
\lvec(86 85)
\ifill f:0
\move(88 84)
\lvec(91 84)
\lvec(91 85)
\lvec(88 85)
\ifill f:0
\move(95 84)
\lvec(98 84)
\lvec(98 85)
\lvec(95 85)
\ifill f:0
\move(99 84)
\lvec(101 84)
\lvec(101 85)
\lvec(99 85)
\ifill f:0
\move(102 84)
\lvec(111 84)
\lvec(111 85)
\lvec(102 85)
\ifill f:0
\move(112 84)
\lvec(122 84)
\lvec(122 85)
\lvec(112 85)
\ifill f:0
\move(123 84)
\lvec(126 84)
\lvec(126 85)
\lvec(123 85)
\ifill f:0
\move(128 84)
\lvec(130 84)
\lvec(130 85)
\lvec(128 85)
\ifill f:0
\move(132 84)
\lvec(133 84)
\lvec(133 85)
\lvec(132 85)
\ifill f:0
\move(135 84)
\lvec(136 84)
\lvec(136 85)
\lvec(135 85)
\ifill f:0
\move(137 84)
\lvec(143 84)
\lvec(143 85)
\lvec(137 85)
\ifill f:0
\move(144 84)
\lvec(145 84)
\lvec(145 85)
\lvec(144 85)
\ifill f:0
\move(146 84)
\lvec(147 84)
\lvec(147 85)
\lvec(146 85)
\ifill f:0
\move(148 84)
\lvec(149 84)
\lvec(149 85)
\lvec(148 85)
\ifill f:0
\move(150 84)
\lvec(151 84)
\lvec(151 85)
\lvec(150 85)
\ifill f:0
\move(152 84)
\lvec(159 84)
\lvec(159 85)
\lvec(152 85)
\ifill f:0
\move(160 84)
\lvec(162 84)
\lvec(162 85)
\lvec(160 85)
\ifill f:0
\move(163 84)
\lvec(165 84)
\lvec(165 85)
\lvec(163 85)
\ifill f:0
\move(166 84)
\lvec(170 84)
\lvec(170 85)
\lvec(166 85)
\ifill f:0
\move(171 84)
\lvec(172 84)
\lvec(172 85)
\lvec(171 85)
\ifill f:0
\move(174 84)
\lvec(186 84)
\lvec(186 85)
\lvec(174 85)
\ifill f:0
\move(187 84)
\lvec(197 84)
\lvec(197 85)
\lvec(187 85)
\ifill f:0
\move(198 84)
\lvec(211 84)
\lvec(211 85)
\lvec(198 85)
\ifill f:0
\move(212 84)
\lvec(219 84)
\lvec(219 85)
\lvec(212 85)
\ifill f:0
\move(220 84)
\lvec(226 84)
\lvec(226 85)
\lvec(220 85)
\ifill f:0
\move(227 84)
\lvec(230 84)
\lvec(230 85)
\lvec(227 85)
\ifill f:0
\move(231 84)
\lvec(235 84)
\lvec(235 85)
\lvec(231 85)
\ifill f:0
\move(236 84)
\lvec(242 84)
\lvec(242 85)
\lvec(236 85)
\ifill f:0
\move(243 84)
\lvec(255 84)
\lvec(255 85)
\lvec(243 85)
\ifill f:0
\move(256 84)
\lvec(257 84)
\lvec(257 85)
\lvec(256 85)
\ifill f:0
\move(258 84)
\lvec(260 84)
\lvec(260 85)
\lvec(258 85)
\ifill f:0
\move(261 84)
\lvec(279 84)
\lvec(279 85)
\lvec(261 85)
\ifill f:0
\move(280 84)
\lvec(281 84)
\lvec(281 85)
\lvec(280 85)
\ifill f:0
\move(282 84)
\lvec(283 84)
\lvec(283 85)
\lvec(282 85)
\ifill f:0
\move(284 84)
\lvec(285 84)
\lvec(285 85)
\lvec(284 85)
\ifill f:0
\move(286 84)
\lvec(287 84)
\lvec(287 85)
\lvec(286 85)
\ifill f:0
\move(288 84)
\lvec(290 84)
\lvec(290 85)
\lvec(288 85)
\ifill f:0
\move(292 84)
\lvec(293 84)
\lvec(293 85)
\lvec(292 85)
\ifill f:0
\move(294 84)
\lvec(295 84)
\lvec(295 85)
\lvec(294 85)
\ifill f:0
\move(296 84)
\lvec(302 84)
\lvec(302 85)
\lvec(296 85)
\ifill f:0
\move(303 84)
\lvec(309 84)
\lvec(309 85)
\lvec(303 85)
\ifill f:0
\move(310 84)
\lvec(317 84)
\lvec(317 85)
\lvec(310 85)
\ifill f:0
\move(318 84)
\lvec(322 84)
\lvec(322 85)
\lvec(318 85)
\ifill f:0
\move(323 84)
\lvec(325 84)
\lvec(325 85)
\lvec(323 85)
\ifill f:0
\move(326 84)
\lvec(338 84)
\lvec(338 85)
\lvec(326 85)
\ifill f:0
\move(339 84)
\lvec(345 84)
\lvec(345 85)
\lvec(339 85)
\ifill f:0
\move(346 84)
\lvec(362 84)
\lvec(362 85)
\lvec(346 85)
\ifill f:0
\move(363 84)
\lvec(365 84)
\lvec(365 85)
\lvec(363 85)
\ifill f:0
\move(366 84)
\lvec(370 84)
\lvec(370 85)
\lvec(366 85)
\ifill f:0
\move(371 84)
\lvec(376 84)
\lvec(376 85)
\lvec(371 85)
\ifill f:0
\move(377 84)
\lvec(388 84)
\lvec(388 85)
\lvec(377 85)
\ifill f:0
\move(389 84)
\lvec(395 84)
\lvec(395 85)
\lvec(389 85)
\ifill f:0
\move(397 84)
\lvec(401 84)
\lvec(401 85)
\lvec(397 85)
\ifill f:0
\move(402 84)
\lvec(416 84)
\lvec(416 85)
\lvec(402 85)
\ifill f:0
\move(417 84)
\lvec(436 84)
\lvec(436 85)
\lvec(417 85)
\ifill f:0
\move(437 84)
\lvec(439 84)
\lvec(439 85)
\lvec(437 85)
\ifill f:0
\move(441 84)
\lvec(442 84)
\lvec(442 85)
\lvec(441 85)
\ifill f:0
\move(446 84)
\lvec(447 84)
\lvec(447 85)
\lvec(446 85)
\ifill f:0
\move(16 85)
\lvec(17 85)
\lvec(17 86)
\lvec(16 86)
\ifill f:0
\move(18 85)
\lvec(20 85)
\lvec(20 86)
\lvec(18 86)
\ifill f:0
\move(24 85)
\lvec(26 85)
\lvec(26 86)
\lvec(24 86)
\ifill f:0
\move(36 85)
\lvec(37 85)
\lvec(37 86)
\lvec(36 86)
\ifill f:0
\move(38 85)
\lvec(39 85)
\lvec(39 86)
\lvec(38 86)
\ifill f:0
\move(40 85)
\lvec(41 85)
\lvec(41 86)
\lvec(40 86)
\ifill f:0
\move(42 85)
\lvec(46 85)
\lvec(46 86)
\lvec(42 86)
\ifill f:0
\move(47 85)
\lvec(48 85)
\lvec(48 86)
\lvec(47 86)
\ifill f:0
\move(49 85)
\lvec(50 85)
\lvec(50 86)
\lvec(49 86)
\ifill f:0
\move(51 85)
\lvec(53 85)
\lvec(53 86)
\lvec(51 86)
\ifill f:0
\move(56 85)
\lvec(58 85)
\lvec(58 86)
\lvec(56 86)
\ifill f:0
\move(59 85)
\lvec(60 85)
\lvec(60 86)
\lvec(59 86)
\ifill f:0
\move(61 85)
\lvec(65 85)
\lvec(65 86)
\lvec(61 86)
\ifill f:0
\move(66 85)
\lvec(70 85)
\lvec(70 86)
\lvec(66 86)
\ifill f:0
\move(73 85)
\lvec(74 85)
\lvec(74 86)
\lvec(73 86)
\ifill f:0
\move(78 85)
\lvec(82 85)
\lvec(82 86)
\lvec(78 86)
\ifill f:0
\move(83 85)
\lvec(85 85)
\lvec(85 86)
\lvec(83 86)
\ifill f:0
\move(86 85)
\lvec(87 85)
\lvec(87 86)
\lvec(86 86)
\ifill f:0
\move(88 85)
\lvec(90 85)
\lvec(90 86)
\lvec(88 86)
\ifill f:0
\move(91 85)
\lvec(94 85)
\lvec(94 86)
\lvec(91 86)
\ifill f:0
\move(96 85)
\lvec(98 85)
\lvec(98 86)
\lvec(96 86)
\ifill f:0
\move(99 85)
\lvec(101 85)
\lvec(101 86)
\lvec(99 86)
\ifill f:0
\move(103 85)
\lvec(106 85)
\lvec(106 86)
\lvec(103 86)
\ifill f:0
\move(108 85)
\lvec(122 85)
\lvec(122 86)
\lvec(108 86)
\ifill f:0
\move(125 85)
\lvec(128 85)
\lvec(128 86)
\lvec(125 86)
\ifill f:0
\move(129 85)
\lvec(132 85)
\lvec(132 86)
\lvec(129 86)
\ifill f:0
\move(133 85)
\lvec(135 85)
\lvec(135 86)
\lvec(133 86)
\ifill f:0
\move(136 85)
\lvec(138 85)
\lvec(138 86)
\lvec(136 86)
\ifill f:0
\move(139 85)
\lvec(143 85)
\lvec(143 86)
\lvec(139 86)
\ifill f:0
\move(144 85)
\lvec(145 85)
\lvec(145 86)
\lvec(144 86)
\ifill f:0
\move(146 85)
\lvec(155 85)
\lvec(155 86)
\lvec(146 86)
\ifill f:0
\move(156 85)
\lvec(160 85)
\lvec(160 86)
\lvec(156 86)
\ifill f:0
\move(161 85)
\lvec(170 85)
\lvec(170 86)
\lvec(161 86)
\ifill f:0
\move(171 85)
\lvec(172 85)
\lvec(172 86)
\lvec(171 86)
\ifill f:0
\move(173 85)
\lvec(176 85)
\lvec(176 86)
\lvec(173 86)
\ifill f:0
\move(178 85)
\lvec(182 85)
\lvec(182 86)
\lvec(178 86)
\ifill f:0
\move(183 85)
\lvec(189 85)
\lvec(189 86)
\lvec(183 86)
\ifill f:0
\move(190 85)
\lvec(197 85)
\lvec(197 86)
\lvec(190 86)
\ifill f:0
\move(198 85)
\lvec(217 85)
\lvec(217 86)
\lvec(198 86)
\ifill f:0
\move(218 85)
\lvec(226 85)
\lvec(226 86)
\lvec(218 86)
\ifill f:0
\move(227 85)
\lvec(231 85)
\lvec(231 86)
\lvec(227 86)
\ifill f:0
\move(232 85)
\lvec(241 85)
\lvec(241 86)
\lvec(232 86)
\ifill f:0
\move(242 85)
\lvec(248 85)
\lvec(248 86)
\lvec(242 86)
\ifill f:0
\move(249 85)
\lvec(257 85)
\lvec(257 86)
\lvec(249 86)
\ifill f:0
\move(258 85)
\lvec(266 85)
\lvec(266 86)
\lvec(258 86)
\ifill f:0
\move(267 85)
\lvec(271 85)
\lvec(271 86)
\lvec(267 86)
\ifill f:0
\move(272 85)
\lvec(283 85)
\lvec(283 86)
\lvec(272 86)
\ifill f:0
\move(284 85)
\lvec(287 85)
\lvec(287 86)
\lvec(284 86)
\ifill f:0
\move(288 85)
\lvec(290 85)
\lvec(290 86)
\lvec(288 86)
\ifill f:0
\move(292 85)
\lvec(293 85)
\lvec(293 86)
\lvec(292 86)
\ifill f:0
\move(294 85)
\lvec(295 85)
\lvec(295 86)
\lvec(294 86)
\ifill f:0
\move(296 85)
\lvec(297 85)
\lvec(297 86)
\lvec(296 86)
\ifill f:0
\move(298 85)
\lvec(299 85)
\lvec(299 86)
\lvec(298 86)
\ifill f:0
\move(300 85)
\lvec(308 85)
\lvec(308 86)
\lvec(300 86)
\ifill f:0
\move(309 85)
\lvec(310 85)
\lvec(310 86)
\lvec(309 86)
\ifill f:0
\move(311 85)
\lvec(315 85)
\lvec(315 86)
\lvec(311 86)
\ifill f:0
\move(316 85)
\lvec(320 85)
\lvec(320 86)
\lvec(316 86)
\ifill f:0
\move(321 85)
\lvec(325 85)
\lvec(325 86)
\lvec(321 86)
\ifill f:0
\move(326 85)
\lvec(328 85)
\lvec(328 86)
\lvec(326 86)
\ifill f:0
\move(329 85)
\lvec(343 85)
\lvec(343 86)
\lvec(329 86)
\ifill f:0
\move(344 85)
\lvec(362 85)
\lvec(362 86)
\lvec(344 86)
\ifill f:0
\move(363 85)
\lvec(365 85)
\lvec(365 86)
\lvec(363 86)
\ifill f:0
\move(366 85)
\lvec(369 85)
\lvec(369 86)
\lvec(366 86)
\ifill f:0
\move(370 85)
\lvec(379 85)
\lvec(379 86)
\lvec(370 86)
\ifill f:0
\move(380 85)
\lvec(384 85)
\lvec(384 86)
\lvec(380 86)
\ifill f:0
\move(385 85)
\lvec(390 85)
\lvec(390 86)
\lvec(385 86)
\ifill f:0
\move(391 85)
\lvec(401 85)
\lvec(401 86)
\lvec(391 86)
\ifill f:0
\move(402 85)
\lvec(404 85)
\lvec(404 86)
\lvec(402 86)
\ifill f:0
\move(405 85)
\lvec(412 85)
\lvec(412 86)
\lvec(405 86)
\ifill f:0
\move(413 85)
\lvec(423 85)
\lvec(423 86)
\lvec(413 86)
\ifill f:0
\move(424 85)
\lvec(439 85)
\lvec(439 86)
\lvec(424 86)
\ifill f:0
\move(441 85)
\lvec(442 85)
\lvec(442 86)
\lvec(441 86)
\ifill f:0
\move(443 85)
\lvec(451 85)
\lvec(451 86)
\lvec(443 86)
\ifill f:0
\move(16 86)
\lvec(17 86)
\lvec(17 87)
\lvec(16 87)
\ifill f:0
\move(20 86)
\lvec(21 86)
\lvec(21 87)
\lvec(20 87)
\ifill f:0
\move(23 86)
\lvec(26 86)
\lvec(26 87)
\lvec(23 87)
\ifill f:0
\move(36 86)
\lvec(37 86)
\lvec(37 87)
\lvec(36 87)
\ifill f:0
\move(38 86)
\lvec(39 86)
\lvec(39 87)
\lvec(38 87)
\ifill f:0
\move(40 86)
\lvec(45 86)
\lvec(45 87)
\lvec(40 87)
\ifill f:0
\move(46 86)
\lvec(48 86)
\lvec(48 87)
\lvec(46 87)
\ifill f:0
\move(49 86)
\lvec(50 86)
\lvec(50 87)
\lvec(49 87)
\ifill f:0
\move(51 86)
\lvec(53 86)
\lvec(53 87)
\lvec(51 87)
\ifill f:0
\move(54 86)
\lvec(55 86)
\lvec(55 87)
\lvec(54 87)
\ifill f:0
\move(58 86)
\lvec(59 86)
\lvec(59 87)
\lvec(58 87)
\ifill f:0
\move(60 86)
\lvec(65 86)
\lvec(65 87)
\lvec(60 87)
\ifill f:0
\move(66 86)
\lvec(74 86)
\lvec(74 87)
\lvec(66 87)
\ifill f:0
\move(76 86)
\lvec(82 86)
\lvec(82 87)
\lvec(76 87)
\ifill f:0
\move(83 86)
\lvec(85 86)
\lvec(85 87)
\lvec(83 87)
\ifill f:0
\move(87 86)
\lvec(88 86)
\lvec(88 87)
\lvec(87 87)
\ifill f:0
\move(89 86)
\lvec(91 86)
\lvec(91 87)
\lvec(89 87)
\ifill f:0
\move(92 86)
\lvec(95 86)
\lvec(95 87)
\lvec(92 87)
\ifill f:0
\move(96 86)
\lvec(98 86)
\lvec(98 87)
\lvec(96 87)
\ifill f:0
\move(99 86)
\lvec(101 86)
\lvec(101 87)
\lvec(99 87)
\ifill f:0
\move(102 86)
\lvec(111 86)
\lvec(111 87)
\lvec(102 87)
\ifill f:0
\move(112 86)
\lvec(113 86)
\lvec(113 87)
\lvec(112 87)
\ifill f:0
\move(114 86)
\lvec(122 86)
\lvec(122 87)
\lvec(114 87)
\ifill f:0
\move(126 86)
\lvec(129 86)
\lvec(129 87)
\lvec(126 87)
\ifill f:0
\move(131 86)
\lvec(134 86)
\lvec(134 87)
\lvec(131 87)
\ifill f:0
\move(135 86)
\lvec(137 86)
\lvec(137 87)
\lvec(135 87)
\ifill f:0
\move(139 86)
\lvec(140 86)
\lvec(140 87)
\lvec(139 87)
\ifill f:0
\move(141 86)
\lvec(143 86)
\lvec(143 87)
\lvec(141 87)
\ifill f:0
\move(144 86)
\lvec(145 86)
\lvec(145 87)
\lvec(144 87)
\ifill f:0
\move(146 86)
\lvec(150 86)
\lvec(150 87)
\lvec(146 87)
\ifill f:0
\move(151 86)
\lvec(154 86)
\lvec(154 87)
\lvec(151 87)
\ifill f:0
\move(155 86)
\lvec(163 86)
\lvec(163 87)
\lvec(155 87)
\ifill f:0
\move(164 86)
\lvec(166 86)
\lvec(166 87)
\lvec(164 87)
\ifill f:0
\move(167 86)
\lvec(170 86)
\lvec(170 87)
\lvec(167 87)
\ifill f:0
\move(171 86)
\lvec(172 86)
\lvec(172 87)
\lvec(171 87)
\ifill f:0
\move(173 86)
\lvec(175 86)
\lvec(175 87)
\lvec(173 87)
\ifill f:0
\move(176 86)
\lvec(179 86)
\lvec(179 87)
\lvec(176 87)
\ifill f:0
\move(180 86)
\lvec(191 86)
\lvec(191 87)
\lvec(180 87)
\ifill f:0
\move(192 86)
\lvec(197 86)
\lvec(197 87)
\lvec(192 87)
\ifill f:0
\move(198 86)
\lvec(201 86)
\lvec(201 87)
\lvec(198 87)
\ifill f:0
\move(202 86)
\lvec(203 86)
\lvec(203 87)
\lvec(202 87)
\ifill f:0
\move(205 86)
\lvec(206 86)
\lvec(206 87)
\lvec(205 87)
\ifill f:0
\move(208 86)
\lvec(209 86)
\lvec(209 87)
\lvec(208 87)
\ifill f:0
\move(212 86)
\lvec(226 86)
\lvec(226 87)
\lvec(212 87)
\ifill f:0
\move(227 86)
\lvec(233 86)
\lvec(233 87)
\lvec(227 87)
\ifill f:0
\move(234 86)
\lvec(238 86)
\lvec(238 87)
\lvec(234 87)
\ifill f:0
\move(239 86)
\lvec(243 86)
\lvec(243 87)
\lvec(239 87)
\ifill f:0
\move(244 86)
\lvec(247 86)
\lvec(247 87)
\lvec(244 87)
\ifill f:0
\move(248 86)
\lvec(251 86)
\lvec(251 87)
\lvec(248 87)
\ifill f:0
\move(252 86)
\lvec(257 86)
\lvec(257 87)
\lvec(252 87)
\ifill f:0
\move(258 86)
\lvec(275 86)
\lvec(275 87)
\lvec(258 87)
\ifill f:0
\move(276 86)
\lvec(280 86)
\lvec(280 87)
\lvec(276 87)
\ifill f:0
\move(281 86)
\lvec(287 86)
\lvec(287 87)
\lvec(281 87)
\ifill f:0
\move(288 86)
\lvec(290 86)
\lvec(290 87)
\lvec(288 87)
\ifill f:0
\move(291 86)
\lvec(298 86)
\lvec(298 87)
\lvec(291 87)
\ifill f:0
\move(299 86)
\lvec(311 86)
\lvec(311 87)
\lvec(299 87)
\ifill f:0
\move(312 86)
\lvec(318 86)
\lvec(318 87)
\lvec(312 87)
\ifill f:0
\move(319 86)
\lvec(320 86)
\lvec(320 87)
\lvec(319 87)
\ifill f:0
\move(321 86)
\lvec(325 86)
\lvec(325 87)
\lvec(321 87)
\ifill f:0
\move(326 86)
\lvec(333 86)
\lvec(333 87)
\lvec(326 87)
\ifill f:0
\move(334 86)
\lvec(362 86)
\lvec(362 87)
\lvec(334 87)
\ifill f:0
\move(363 86)
\lvec(377 86)
\lvec(377 87)
\lvec(363 87)
\ifill f:0
\move(378 86)
\lvec(386 86)
\lvec(386 87)
\lvec(378 87)
\ifill f:0
\move(387 86)
\lvec(391 86)
\lvec(391 87)
\lvec(387 87)
\ifill f:0
\move(392 86)
\lvec(397 86)
\lvec(397 87)
\lvec(392 87)
\ifill f:0
\move(398 86)
\lvec(401 86)
\lvec(401 87)
\lvec(398 87)
\ifill f:0
\move(402 86)
\lvec(403 86)
\lvec(403 87)
\lvec(402 87)
\ifill f:0
\move(404 86)
\lvec(410 86)
\lvec(410 87)
\lvec(404 87)
\ifill f:0
\move(411 86)
\lvec(418 86)
\lvec(418 87)
\lvec(411 87)
\ifill f:0
\move(419 86)
\lvec(427 86)
\lvec(427 87)
\lvec(419 87)
\ifill f:0
\move(429 86)
\lvec(440 86)
\lvec(440 87)
\lvec(429 87)
\ifill f:0
\move(441 86)
\lvec(442 86)
\lvec(442 87)
\lvec(441 87)
\ifill f:0
\move(443 86)
\lvec(451 86)
\lvec(451 87)
\lvec(443 87)
\ifill f:0
\move(15 87)
\lvec(17 87)
\lvec(17 88)
\lvec(15 88)
\ifill f:0
\move(20 87)
\lvec(21 87)
\lvec(21 88)
\lvec(20 88)
\ifill f:0
\move(23 87)
\lvec(26 87)
\lvec(26 88)
\lvec(23 88)
\ifill f:0
\move(36 87)
\lvec(37 87)
\lvec(37 88)
\lvec(36 88)
\ifill f:0
\move(40 87)
\lvec(41 87)
\lvec(41 88)
\lvec(40 88)
\ifill f:0
\move(42 87)
\lvec(45 87)
\lvec(45 88)
\lvec(42 88)
\ifill f:0
\move(47 87)
\lvec(50 87)
\lvec(50 88)
\lvec(47 88)
\ifill f:0
\move(51 87)
\lvec(52 87)
\lvec(52 88)
\lvec(51 88)
\ifill f:0
\move(56 87)
\lvec(58 87)
\lvec(58 88)
\lvec(56 88)
\ifill f:0
\move(59 87)
\lvec(61 87)
\lvec(61 88)
\lvec(59 88)
\ifill f:0
\move(62 87)
\lvec(63 87)
\lvec(63 88)
\lvec(62 88)
\ifill f:0
\move(64 87)
\lvec(65 87)
\lvec(65 88)
\lvec(64 88)
\ifill f:0
\move(66 87)
\lvec(71 87)
\lvec(71 88)
\lvec(66 88)
\ifill f:0
\move(73 87)
\lvec(82 87)
\lvec(82 88)
\lvec(73 88)
\ifill f:0
\move(84 87)
\lvec(86 87)
\lvec(86 88)
\lvec(84 88)
\ifill f:0
\move(88 87)
\lvec(89 87)
\lvec(89 88)
\lvec(88 88)
\ifill f:0
\move(91 87)
\lvec(93 87)
\lvec(93 88)
\lvec(91 88)
\ifill f:0
\move(95 87)
\lvec(96 87)
\lvec(96 88)
\lvec(95 88)
\ifill f:0
\move(97 87)
\lvec(98 87)
\lvec(98 88)
\lvec(97 88)
\ifill f:0
\move(99 87)
\lvec(101 87)
\lvec(101 88)
\lvec(99 88)
\ifill f:0
\move(102 87)
\lvec(116 87)
\lvec(116 88)
\lvec(102 88)
\ifill f:0
\move(117 87)
\lvec(118 87)
\lvec(118 88)
\lvec(117 88)
\ifill f:0
\move(121 87)
\lvec(122 87)
\lvec(122 88)
\lvec(121 88)
\ifill f:0
\move(127 87)
\lvec(131 87)
\lvec(131 88)
\lvec(127 88)
\ifill f:0
\move(133 87)
\lvec(136 87)
\lvec(136 88)
\lvec(133 88)
\ifill f:0
\move(137 87)
\lvec(143 87)
\lvec(143 88)
\lvec(137 88)
\ifill f:0
\move(144 87)
\lvec(145 87)
\lvec(145 88)
\lvec(144 88)
\ifill f:0
\move(146 87)
\lvec(148 87)
\lvec(148 88)
\lvec(146 88)
\ifill f:0
\move(149 87)
\lvec(155 87)
\lvec(155 88)
\lvec(149 88)
\ifill f:0
\move(156 87)
\lvec(159 87)
\lvec(159 88)
\lvec(156 88)
\ifill f:0
\move(160 87)
\lvec(170 87)
\lvec(170 88)
\lvec(160 88)
\ifill f:0
\move(171 87)
\lvec(174 87)
\lvec(174 88)
\lvec(171 88)
\ifill f:0
\move(176 87)
\lvec(186 87)
\lvec(186 88)
\lvec(176 88)
\ifill f:0
\move(187 87)
\lvec(191 87)
\lvec(191 88)
\lvec(187 88)
\ifill f:0
\move(193 87)
\lvec(197 87)
\lvec(197 88)
\lvec(193 88)
\ifill f:0
\move(198 87)
\lvec(226 87)
\lvec(226 88)
\lvec(198 88)
\ifill f:0
\move(228 87)
\lvec(234 87)
\lvec(234 88)
\lvec(228 88)
\ifill f:0
\move(236 87)
\lvec(241 87)
\lvec(241 88)
\lvec(236 88)
\ifill f:0
\move(243 87)
\lvec(250 87)
\lvec(250 88)
\lvec(243 88)
\ifill f:0
\move(251 87)
\lvec(254 87)
\lvec(254 88)
\lvec(251 88)
\ifill f:0
\move(255 87)
\lvec(257 87)
\lvec(257 88)
\lvec(255 88)
\ifill f:0
\move(258 87)
\lvec(271 87)
\lvec(271 88)
\lvec(258 88)
\ifill f:0
\move(272 87)
\lvec(290 87)
\lvec(290 88)
\lvec(272 88)
\ifill f:0
\move(291 87)
\lvec(296 87)
\lvec(296 88)
\lvec(291 88)
\ifill f:0
\move(297 87)
\lvec(298 87)
\lvec(298 88)
\lvec(297 88)
\ifill f:0
\move(299 87)
\lvec(314 87)
\lvec(314 88)
\lvec(299 88)
\ifill f:0
\move(315 87)
\lvec(325 87)
\lvec(325 88)
\lvec(315 88)
\ifill f:0
\move(326 87)
\lvec(330 87)
\lvec(330 88)
\lvec(326 88)
\ifill f:0
\move(331 87)
\lvec(335 87)
\lvec(335 88)
\lvec(331 88)
\ifill f:0
\move(336 87)
\lvec(343 87)
\lvec(343 88)
\lvec(336 88)
\ifill f:0
\move(344 87)
\lvec(362 87)
\lvec(362 88)
\lvec(344 88)
\ifill f:0
\move(363 87)
\lvec(364 87)
\lvec(364 88)
\lvec(363 88)
\ifill f:0
\move(365 87)
\lvec(375 87)
\lvec(375 88)
\lvec(365 88)
\ifill f:0
\move(376 87)
\lvec(379 87)
\lvec(379 88)
\lvec(376 88)
\ifill f:0
\move(380 87)
\lvec(383 87)
\lvec(383 88)
\lvec(380 88)
\ifill f:0
\move(384 87)
\lvec(392 87)
\lvec(392 88)
\lvec(384 88)
\ifill f:0
\move(393 87)
\lvec(397 87)
\lvec(397 88)
\lvec(393 88)
\ifill f:0
\move(398 87)
\lvec(401 87)
\lvec(401 88)
\lvec(398 88)
\ifill f:0
\move(402 87)
\lvec(415 87)
\lvec(415 88)
\lvec(402 88)
\ifill f:0
\move(416 87)
\lvec(430 87)
\lvec(430 88)
\lvec(416 88)
\ifill f:0
\move(431 87)
\lvec(440 87)
\lvec(440 88)
\lvec(431 88)
\ifill f:0
\move(441 87)
\lvec(442 87)
\lvec(442 88)
\lvec(441 88)
\ifill f:0
\move(443 87)
\lvec(451 87)
\lvec(451 88)
\lvec(443 88)
\ifill f:0
\move(15 88)
\lvec(17 88)
\lvec(17 89)
\lvec(15 89)
\ifill f:0
\move(19 88)
\lvec(21 88)
\lvec(21 89)
\lvec(19 89)
\ifill f:0
\move(24 88)
\lvec(26 88)
\lvec(26 89)
\lvec(24 89)
\ifill f:0
\move(36 88)
\lvec(37 88)
\lvec(37 89)
\lvec(36 89)
\ifill f:0
\move(38 88)
\lvec(43 88)
\lvec(43 89)
\lvec(38 89)
\ifill f:0
\move(44 88)
\lvec(46 88)
\lvec(46 89)
\lvec(44 89)
\ifill f:0
\move(47 88)
\lvec(50 88)
\lvec(50 89)
\lvec(47 89)
\ifill f:0
\move(59 88)
\lvec(63 88)
\lvec(63 89)
\lvec(59 89)
\ifill f:0
\move(64 88)
\lvec(65 88)
\lvec(65 89)
\lvec(64 89)
\ifill f:0
\move(66 88)
\lvec(67 88)
\lvec(67 89)
\lvec(66 89)
\ifill f:0
\move(68 88)
\lvec(70 88)
\lvec(70 89)
\lvec(68 89)
\ifill f:0
\move(71 88)
\lvec(74 88)
\lvec(74 89)
\lvec(71 89)
\ifill f:0
\move(76 88)
\lvec(82 88)
\lvec(82 89)
\lvec(76 89)
\ifill f:0
\move(86 88)
\lvec(87 88)
\lvec(87 89)
\lvec(86 89)
\ifill f:0
\move(89 88)
\lvec(91 88)
\lvec(91 89)
\lvec(89 89)
\ifill f:0
\move(92 88)
\lvec(93 88)
\lvec(93 89)
\lvec(92 89)
\ifill f:0
\move(95 88)
\lvec(96 88)
\lvec(96 89)
\lvec(95 89)
\ifill f:0
\move(97 88)
\lvec(98 88)
\lvec(98 89)
\lvec(97 89)
\ifill f:0
\move(99 88)
\lvec(101 88)
\lvec(101 89)
\lvec(99 89)
\ifill f:0
\move(102 88)
\lvec(106 88)
\lvec(106 89)
\lvec(102 89)
\ifill f:0
\move(108 88)
\lvec(111 88)
\lvec(111 89)
\lvec(108 89)
\ifill f:0
\move(112 88)
\lvec(119 88)
\lvec(119 89)
\lvec(112 89)
\ifill f:0
\move(121 88)
\lvec(122 88)
\lvec(122 89)
\lvec(121 89)
\ifill f:0
\move(128 88)
\lvec(134 88)
\lvec(134 89)
\lvec(128 89)
\ifill f:0
\move(136 88)
\lvec(139 88)
\lvec(139 89)
\lvec(136 89)
\ifill f:0
\move(140 88)
\lvec(143 88)
\lvec(143 89)
\lvec(140 89)
\ifill f:0
\move(144 88)
\lvec(145 88)
\lvec(145 89)
\lvec(144 89)
\ifill f:0
\move(146 88)
\lvec(151 88)
\lvec(151 89)
\lvec(146 89)
\ifill f:0
\move(152 88)
\lvec(160 88)
\lvec(160 89)
\lvec(152 89)
\ifill f:0
\move(161 88)
\lvec(162 88)
\lvec(162 89)
\lvec(161 89)
\ifill f:0
\move(163 88)
\lvec(170 88)
\lvec(170 89)
\lvec(163 89)
\ifill f:0
\move(172 88)
\lvec(174 88)
\lvec(174 89)
\lvec(172 89)
\ifill f:0
\move(175 88)
\lvec(177 88)
\lvec(177 89)
\lvec(175 89)
\ifill f:0
\move(178 88)
\lvec(180 88)
\lvec(180 89)
\lvec(178 89)
\ifill f:0
\move(181 88)
\lvec(184 88)
\lvec(184 89)
\lvec(181 89)
\ifill f:0
\move(185 88)
\lvec(188 88)
\lvec(188 89)
\lvec(185 89)
\ifill f:0
\move(189 88)
\lvec(192 88)
\lvec(192 89)
\lvec(189 89)
\ifill f:0
\move(194 88)
\lvec(197 88)
\lvec(197 89)
\lvec(194 89)
\ifill f:0
\move(198 88)
\lvec(199 88)
\lvec(199 89)
\lvec(198 89)
\ifill f:0
\move(200 88)
\lvec(210 88)
\lvec(210 89)
\lvec(200 89)
\ifill f:0
\move(211 88)
\lvec(226 88)
\lvec(226 89)
\lvec(211 89)
\ifill f:0
\move(229 88)
\lvec(237 88)
\lvec(237 89)
\lvec(229 89)
\ifill f:0
\move(238 88)
\lvec(244 88)
\lvec(244 89)
\lvec(238 89)
\ifill f:0
\move(245 88)
\lvec(249 88)
\lvec(249 89)
\lvec(245 89)
\ifill f:0
\move(251 88)
\lvec(254 88)
\lvec(254 89)
\lvec(251 89)
\ifill f:0
\move(255 88)
\lvec(257 88)
\lvec(257 89)
\lvec(255 89)
\ifill f:0
\move(259 88)
\lvec(262 88)
\lvec(262 89)
\lvec(259 89)
\ifill f:0
\move(263 88)
\lvec(266 88)
\lvec(266 89)
\lvec(263 89)
\ifill f:0
\move(267 88)
\lvec(269 88)
\lvec(269 89)
\lvec(267 89)
\ifill f:0
\move(270 88)
\lvec(284 88)
\lvec(284 89)
\lvec(270 89)
\ifill f:0
\move(285 88)
\lvec(290 88)
\lvec(290 89)
\lvec(285 89)
\ifill f:0
\move(291 88)
\lvec(294 88)
\lvec(294 89)
\lvec(291 89)
\ifill f:0
\move(295 88)
\lvec(299 88)
\lvec(299 89)
\lvec(295 89)
\ifill f:0
\move(300 88)
\lvec(303 88)
\lvec(303 89)
\lvec(300 89)
\ifill f:0
\move(304 88)
\lvec(325 88)
\lvec(325 89)
\lvec(304 89)
\ifill f:0
\move(326 88)
\lvec(339 88)
\lvec(339 89)
\lvec(326 89)
\ifill f:0
\move(340 88)
\lvec(344 88)
\lvec(344 89)
\lvec(340 89)
\ifill f:0
\move(345 88)
\lvec(352 88)
\lvec(352 89)
\lvec(345 89)
\ifill f:0
\move(353 88)
\lvec(355 88)
\lvec(355 89)
\lvec(353 89)
\ifill f:0
\move(356 88)
\lvec(358 88)
\lvec(358 89)
\lvec(356 89)
\ifill f:0
\move(359 88)
\lvec(362 88)
\lvec(362 89)
\lvec(359 89)
\ifill f:0
\move(363 88)
\lvec(364 88)
\lvec(364 89)
\lvec(363 89)
\ifill f:0
\move(365 88)
\lvec(367 88)
\lvec(367 89)
\lvec(365 89)
\ifill f:0
\move(368 88)
\lvec(370 88)
\lvec(370 89)
\lvec(368 89)
\ifill f:0
\move(371 88)
\lvec(374 88)
\lvec(374 89)
\lvec(371 89)
\ifill f:0
\move(375 88)
\lvec(381 88)
\lvec(381 89)
\lvec(375 89)
\ifill f:0
\move(382 88)
\lvec(385 88)
\lvec(385 89)
\lvec(382 89)
\ifill f:0
\move(386 88)
\lvec(389 88)
\lvec(389 89)
\lvec(386 89)
\ifill f:0
\move(390 88)
\lvec(401 88)
\lvec(401 89)
\lvec(390 89)
\ifill f:0
\move(402 88)
\lvec(425 88)
\lvec(425 89)
\lvec(402 89)
\ifill f:0
\move(426 88)
\lvec(432 88)
\lvec(432 89)
\lvec(426 89)
\ifill f:0
\move(433 88)
\lvec(440 88)
\lvec(440 89)
\lvec(433 89)
\ifill f:0
\move(441 88)
\lvec(442 88)
\lvec(442 89)
\lvec(441 89)
\ifill f:0
\move(443 88)
\lvec(450 88)
\lvec(450 89)
\lvec(443 89)
\ifill f:0
\move(15 89)
\lvec(17 89)
\lvec(17 90)
\lvec(15 90)
\ifill f:0
\move(18 89)
\lvec(19 89)
\lvec(19 90)
\lvec(18 90)
\ifill f:0
\move(20 89)
\lvec(21 89)
\lvec(21 90)
\lvec(20 90)
\ifill f:0
\move(25 89)
\lvec(26 89)
\lvec(26 90)
\lvec(25 90)
\ifill f:0
\move(36 89)
\lvec(37 89)
\lvec(37 90)
\lvec(36 90)
\ifill f:0
\move(38 89)
\lvec(39 89)
\lvec(39 90)
\lvec(38 90)
\ifill f:0
\move(40 89)
\lvec(45 89)
\lvec(45 90)
\lvec(40 90)
\ifill f:0
\move(47 89)
\lvec(50 89)
\lvec(50 90)
\lvec(47 90)
\ifill f:0
\move(54 89)
\lvec(55 89)
\lvec(55 90)
\lvec(54 90)
\ifill f:0
\move(56 89)
\lvec(59 89)
\lvec(59 90)
\lvec(56 90)
\ifill f:0
\move(61 89)
\lvec(63 89)
\lvec(63 90)
\lvec(61 90)
\ifill f:0
\move(64 89)
\lvec(65 89)
\lvec(65 90)
\lvec(64 90)
\ifill f:0
\move(66 89)
\lvec(71 89)
\lvec(71 90)
\lvec(66 90)
\ifill f:0
\move(73 89)
\lvec(75 89)
\lvec(75 90)
\lvec(73 90)
\ifill f:0
\move(76 89)
\lvec(78 89)
\lvec(78 90)
\lvec(76 90)
\ifill f:0
\move(80 89)
\lvec(82 89)
\lvec(82 90)
\lvec(80 90)
\ifill f:0
\move(87 89)
\lvec(89 89)
\lvec(89 90)
\lvec(87 90)
\ifill f:0
\move(91 89)
\lvec(92 89)
\lvec(92 90)
\lvec(91 90)
\ifill f:0
\move(96 89)
\lvec(97 89)
\lvec(97 90)
\lvec(96 90)
\ifill f:0
\move(99 89)
\lvec(101 89)
\lvec(101 90)
\lvec(99 90)
\ifill f:0
\move(102 89)
\lvec(103 89)
\lvec(103 90)
\lvec(102 90)
\ifill f:0
\move(104 89)
\lvec(106 89)
\lvec(106 90)
\lvec(104 90)
\ifill f:0
\move(107 89)
\lvec(109 89)
\lvec(109 90)
\lvec(107 90)
\ifill f:0
\move(110 89)
\lvec(113 89)
\lvec(113 90)
\lvec(110 90)
\ifill f:0
\move(114 89)
\lvec(120 89)
\lvec(120 90)
\lvec(114 90)
\ifill f:0
\move(121 89)
\lvec(122 89)
\lvec(122 90)
\lvec(121 90)
\ifill f:0
\move(123 89)
\lvec(130 89)
\lvec(130 90)
\lvec(123 90)
\ifill f:0
\move(134 89)
\lvec(138 89)
\lvec(138 90)
\lvec(134 90)
\ifill f:0
\move(139 89)
\lvec(142 89)
\lvec(142 90)
\lvec(139 90)
\ifill f:0
\move(144 89)
\lvec(145 89)
\lvec(145 90)
\lvec(144 90)
\ifill f:0
\move(146 89)
\lvec(149 89)
\lvec(149 90)
\lvec(146 90)
\ifill f:0
\move(150 89)
\lvec(157 89)
\lvec(157 90)
\lvec(150 90)
\ifill f:0
\move(158 89)
\lvec(170 89)
\lvec(170 90)
\lvec(158 90)
\ifill f:0
\move(172 89)
\lvec(173 89)
\lvec(173 90)
\lvec(172 90)
\ifill f:0
\move(174 89)
\lvec(176 89)
\lvec(176 90)
\lvec(174 90)
\ifill f:0
\move(177 89)
\lvec(179 89)
\lvec(179 90)
\lvec(177 90)
\ifill f:0
\move(180 89)
\lvec(185 89)
\lvec(185 90)
\lvec(180 90)
\ifill f:0
\move(186 89)
\lvec(189 89)
\lvec(189 90)
\lvec(186 90)
\ifill f:0
\move(190 89)
\lvec(193 89)
\lvec(193 90)
\lvec(190 90)
\ifill f:0
\move(194 89)
\lvec(197 89)
\lvec(197 90)
\lvec(194 90)
\ifill f:0
\move(198 89)
\lvec(219 89)
\lvec(219 90)
\lvec(198 90)
\ifill f:0
\move(225 89)
\lvec(226 89)
\lvec(226 90)
\lvec(225 90)
\ifill f:0
\move(230 89)
\lvec(241 89)
\lvec(241 90)
\lvec(230 90)
\ifill f:0
\move(243 89)
\lvec(248 89)
\lvec(248 90)
\lvec(243 90)
\ifill f:0
\move(249 89)
\lvec(257 89)
\lvec(257 90)
\lvec(249 90)
\ifill f:0
\move(258 89)
\lvec(267 89)
\lvec(267 90)
\lvec(258 90)
\ifill f:0
\move(268 89)
\lvec(271 89)
\lvec(271 90)
\lvec(268 90)
\ifill f:0
\move(272 89)
\lvec(274 89)
\lvec(274 90)
\lvec(272 90)
\ifill f:0
\move(275 89)
\lvec(290 89)
\lvec(290 90)
\lvec(275 90)
\ifill f:0
\move(291 89)
\lvec(292 89)
\lvec(292 90)
\lvec(291 90)
\ifill f:0
\move(293 89)
\lvec(297 89)
\lvec(297 90)
\lvec(293 90)
\ifill f:0
\move(298 89)
\lvec(302 89)
\lvec(302 90)
\lvec(298 90)
\ifill f:0
\move(303 89)
\lvec(304 89)
\lvec(304 90)
\lvec(303 90)
\ifill f:0
\move(305 89)
\lvec(311 89)
\lvec(311 90)
\lvec(305 90)
\ifill f:0
\move(312 89)
\lvec(313 89)
\lvec(313 90)
\lvec(312 90)
\ifill f:0
\move(314 89)
\lvec(315 89)
\lvec(315 90)
\lvec(314 90)
\ifill f:0
\move(316 89)
\lvec(317 89)
\lvec(317 90)
\lvec(316 90)
\ifill f:0
\move(318 89)
\lvec(319 89)
\lvec(319 90)
\lvec(318 90)
\ifill f:0
\move(320 89)
\lvec(321 89)
\lvec(321 90)
\lvec(320 90)
\ifill f:0
\move(322 89)
\lvec(323 89)
\lvec(323 90)
\lvec(322 90)
\ifill f:0
\move(324 89)
\lvec(325 89)
\lvec(325 90)
\lvec(324 90)
\ifill f:0
\move(326 89)
\lvec(329 89)
\lvec(329 90)
\lvec(326 90)
\ifill f:0
\move(330 89)
\lvec(338 89)
\lvec(338 90)
\lvec(330 90)
\ifill f:0
\move(339 89)
\lvec(362 89)
\lvec(362 90)
\lvec(339 90)
\ifill f:0
\move(363 89)
\lvec(376 89)
\lvec(376 90)
\lvec(363 90)
\ifill f:0
\move(377 89)
\lvec(379 89)
\lvec(379 90)
\lvec(377 90)
\ifill f:0
\move(380 89)
\lvec(390 89)
\lvec(390 90)
\lvec(380 90)
\ifill f:0
\move(391 89)
\lvec(394 89)
\lvec(394 90)
\lvec(391 90)
\ifill f:0
\move(395 89)
\lvec(398 89)
\lvec(398 90)
\lvec(395 90)
\ifill f:0
\move(399 89)
\lvec(401 89)
\lvec(401 90)
\lvec(399 90)
\ifill f:0
\move(402 89)
\lvec(416 89)
\lvec(416 90)
\lvec(402 90)
\ifill f:0
\move(417 89)
\lvec(427 89)
\lvec(427 90)
\lvec(417 90)
\ifill f:0
\move(428 89)
\lvec(434 89)
\lvec(434 90)
\lvec(428 90)
\ifill f:0
\move(435 89)
\lvec(442 89)
\lvec(442 90)
\lvec(435 90)
\ifill f:0
\move(443 89)
\lvec(449 89)
\lvec(449 90)
\lvec(443 90)
\ifill f:0
\move(450 89)
\lvec(451 89)
\lvec(451 90)
\lvec(450 90)
\ifill f:0
\move(15 90)
\lvec(17 90)
\lvec(17 91)
\lvec(15 91)
\ifill f:0
\move(20 90)
\lvec(22 90)
\lvec(22 91)
\lvec(20 91)
\ifill f:0
\move(25 90)
\lvec(26 90)
\lvec(26 91)
\lvec(25 91)
\ifill f:0
\move(36 90)
\lvec(37 90)
\lvec(37 91)
\lvec(36 91)
\ifill f:0
\move(38 90)
\lvec(39 90)
\lvec(39 91)
\lvec(38 91)
\ifill f:0
\move(40 90)
\lvec(41 90)
\lvec(41 91)
\lvec(40 91)
\ifill f:0
\move(44 90)
\lvec(45 90)
\lvec(45 91)
\lvec(44 91)
\ifill f:0
\move(48 90)
\lvec(50 90)
\lvec(50 91)
\lvec(48 91)
\ifill f:0
\move(52 90)
\lvec(53 90)
\lvec(53 91)
\lvec(52 91)
\ifill f:0
\move(60 90)
\lvec(63 90)
\lvec(63 91)
\lvec(60 91)
\ifill f:0
\move(64 90)
\lvec(65 90)
\lvec(65 91)
\lvec(64 91)
\ifill f:0
\move(66 90)
\lvec(73 90)
\lvec(73 91)
\lvec(66 91)
\ifill f:0
\move(75 90)
\lvec(78 90)
\lvec(78 91)
\lvec(75 91)
\ifill f:0
\move(81 90)
\lvec(82 90)
\lvec(82 91)
\lvec(81 91)
\ifill f:0
\move(88 90)
\lvec(91 90)
\lvec(91 91)
\lvec(88 91)
\ifill f:0
\move(92 90)
\lvec(93 90)
\lvec(93 91)
\lvec(92 91)
\ifill f:0
\move(96 90)
\lvec(98 90)
\lvec(98 91)
\lvec(96 91)
\ifill f:0
\move(100 90)
\lvec(101 90)
\lvec(101 91)
\lvec(100 91)
\ifill f:0
\move(102 90)
\lvec(103 90)
\lvec(103 91)
\lvec(102 91)
\ifill f:0
\move(104 90)
\lvec(108 90)
\lvec(108 91)
\lvec(104 91)
\ifill f:0
\move(109 90)
\lvec(111 90)
\lvec(111 91)
\lvec(109 91)
\ifill f:0
\move(112 90)
\lvec(115 90)
\lvec(115 91)
\lvec(112 91)
\ifill f:0
\move(116 90)
\lvec(120 90)
\lvec(120 91)
\lvec(116 91)
\ifill f:0
\move(121 90)
\lvec(122 90)
\lvec(122 91)
\lvec(121 91)
\ifill f:0
\move(123 90)
\lvec(135 90)
\lvec(135 91)
\lvec(123 91)
\ifill f:0
\move(137 90)
\lvec(138 90)
\lvec(138 91)
\lvec(137 91)
\ifill f:0
\move(139 90)
\lvec(142 90)
\lvec(142 91)
\lvec(139 91)
\ifill f:0
\move(143 90)
\lvec(145 90)
\lvec(145 91)
\lvec(143 91)
\ifill f:0
\move(147 90)
\lvec(150 90)
\lvec(150 91)
\lvec(147 91)
\ifill f:0
\move(151 90)
\lvec(153 90)
\lvec(153 91)
\lvec(151 91)
\ifill f:0
\move(154 90)
\lvec(163 90)
\lvec(163 91)
\lvec(154 91)
\ifill f:0
\move(164 90)
\lvec(165 90)
\lvec(165 91)
\lvec(164 91)
\ifill f:0
\move(166 90)
\lvec(170 90)
\lvec(170 91)
\lvec(166 91)
\ifill f:0
\move(172 90)
\lvec(173 90)
\lvec(173 91)
\lvec(172 91)
\ifill f:0
\move(174 90)
\lvec(186 90)
\lvec(186 91)
\lvec(174 91)
\ifill f:0
\move(187 90)
\lvec(190 90)
\lvec(190 91)
\lvec(187 91)
\ifill f:0
\move(191 90)
\lvec(197 90)
\lvec(197 91)
\lvec(191 91)
\ifill f:0
\move(198 90)
\lvec(203 90)
\lvec(203 91)
\lvec(198 91)
\ifill f:0
\move(204 90)
\lvec(210 90)
\lvec(210 91)
\lvec(204 91)
\ifill f:0
\move(211 90)
\lvec(223 90)
\lvec(223 91)
\lvec(211 91)
\ifill f:0
\move(225 90)
\lvec(226 90)
\lvec(226 91)
\lvec(225 91)
\ifill f:0
\move(236 90)
\lvec(253 90)
\lvec(253 91)
\lvec(236 91)
\ifill f:0
\move(254 90)
\lvec(257 90)
\lvec(257 91)
\lvec(254 91)
\ifill f:0
\move(258 90)
\lvec(259 90)
\lvec(259 91)
\lvec(258 91)
\ifill f:0
\move(260 90)
\lvec(276 90)
\lvec(276 91)
\lvec(260 91)
\ifill f:0
\move(277 90)
\lvec(280 90)
\lvec(280 91)
\lvec(277 91)
\ifill f:0
\move(281 90)
\lvec(283 90)
\lvec(283 91)
\lvec(281 91)
\ifill f:0
\move(284 90)
\lvec(286 90)
\lvec(286 91)
\lvec(284 91)
\ifill f:0
\move(287 90)
\lvec(290 90)
\lvec(290 91)
\lvec(287 91)
\ifill f:0
\move(291 90)
\lvec(292 90)
\lvec(292 91)
\lvec(291 91)
\ifill f:0
\move(293 90)
\lvec(295 90)
\lvec(295 91)
\lvec(293 91)
\ifill f:0
\move(296 90)
\lvec(298 90)
\lvec(298 91)
\lvec(296 91)
\ifill f:0
\move(299 90)
\lvec(300 90)
\lvec(300 91)
\lvec(299 91)
\ifill f:0
\move(301 90)
\lvec(305 90)
\lvec(305 91)
\lvec(301 91)
\ifill f:0
\move(306 90)
\lvec(310 90)
\lvec(310 91)
\lvec(306 91)
\ifill f:0
\move(311 90)
\lvec(319 90)
\lvec(319 91)
\lvec(311 91)
\ifill f:0
\move(320 90)
\lvec(321 90)
\lvec(321 91)
\lvec(320 91)
\ifill f:0
\move(322 90)
\lvec(323 90)
\lvec(323 91)
\lvec(322 91)
\ifill f:0
\move(324 90)
\lvec(325 90)
\lvec(325 91)
\lvec(324 91)
\ifill f:0
\move(326 90)
\lvec(327 90)
\lvec(327 91)
\lvec(326 91)
\ifill f:0
\move(328 90)
\lvec(329 90)
\lvec(329 91)
\lvec(328 91)
\ifill f:0
\move(330 90)
\lvec(331 90)
\lvec(331 91)
\lvec(330 91)
\ifill f:0
\move(332 90)
\lvec(333 90)
\lvec(333 91)
\lvec(332 91)
\ifill f:0
\move(334 90)
\lvec(335 90)
\lvec(335 91)
\lvec(334 91)
\ifill f:0
\move(336 90)
\lvec(344 90)
\lvec(344 91)
\lvec(336 91)
\ifill f:0
\move(345 90)
\lvec(351 90)
\lvec(351 91)
\lvec(345 91)
\ifill f:0
\move(352 90)
\lvec(362 90)
\lvec(362 91)
\lvec(352 91)
\ifill f:0
\move(363 90)
\lvec(369 90)
\lvec(369 91)
\lvec(363 91)
\ifill f:0
\move(370 90)
\lvec(372 90)
\lvec(372 91)
\lvec(370 91)
\ifill f:0
\move(373 90)
\lvec(375 90)
\lvec(375 91)
\lvec(373 91)
\ifill f:0
\move(376 90)
\lvec(378 90)
\lvec(378 91)
\lvec(376 91)
\ifill f:0
\move(379 90)
\lvec(381 90)
\lvec(381 91)
\lvec(379 91)
\ifill f:0
\move(382 90)
\lvec(384 90)
\lvec(384 91)
\lvec(382 91)
\ifill f:0
\move(385 90)
\lvec(391 90)
\lvec(391 91)
\lvec(385 91)
\ifill f:0
\move(392 90)
\lvec(394 90)
\lvec(394 91)
\lvec(392 91)
\ifill f:0
\move(395 90)
\lvec(398 90)
\lvec(398 91)
\lvec(395 91)
\ifill f:0
\move(399 90)
\lvec(401 90)
\lvec(401 91)
\lvec(399 91)
\ifill f:0
\move(403 90)
\lvec(406 90)
\lvec(406 91)
\lvec(403 91)
\ifill f:0
\move(407 90)
\lvec(410 90)
\lvec(410 91)
\lvec(407 91)
\ifill f:0
\move(411 90)
\lvec(419 90)
\lvec(419 91)
\lvec(411 91)
\ifill f:0
\move(420 90)
\lvec(424 90)
\lvec(424 91)
\lvec(420 91)
\ifill f:0
\move(425 90)
\lvec(435 90)
\lvec(435 91)
\lvec(425 91)
\ifill f:0
\move(436 90)
\lvec(442 90)
\lvec(442 91)
\lvec(436 91)
\ifill f:0
\move(443 90)
\lvec(448 90)
\lvec(448 91)
\lvec(443 91)
\ifill f:0
\move(449 90)
\lvec(451 90)
\lvec(451 91)
\lvec(449 91)
\ifill f:0
\move(16 91)
\lvec(17 91)
\lvec(17 92)
\lvec(16 92)
\ifill f:0
\move(20 91)
\lvec(21 91)
\lvec(21 92)
\lvec(20 92)
\ifill f:0
\move(23 91)
\lvec(24 91)
\lvec(24 92)
\lvec(23 92)
\ifill f:0
\move(25 91)
\lvec(26 91)
\lvec(26 92)
\lvec(25 92)
\ifill f:0
\move(36 91)
\lvec(37 91)
\lvec(37 92)
\lvec(36 92)
\ifill f:0
\move(38 91)
\lvec(39 91)
\lvec(39 92)
\lvec(38 92)
\ifill f:0
\move(40 91)
\lvec(47 91)
\lvec(47 92)
\lvec(40 92)
\ifill f:0
\move(48 91)
\lvec(50 91)
\lvec(50 92)
\lvec(48 92)
\ifill f:0
\move(55 91)
\lvec(58 91)
\lvec(58 92)
\lvec(55 92)
\ifill f:0
\move(59 91)
\lvec(62 91)
\lvec(62 92)
\lvec(59 92)
\ifill f:0
\move(64 91)
\lvec(65 91)
\lvec(65 92)
\lvec(64 92)
\ifill f:0
\move(66 91)
\lvec(72 91)
\lvec(72 92)
\lvec(66 92)
\ifill f:0
\move(73 91)
\lvec(74 91)
\lvec(74 92)
\lvec(73 92)
\ifill f:0
\move(76 91)
\lvec(79 91)
\lvec(79 92)
\lvec(76 92)
\ifill f:0
\move(81 91)
\lvec(82 91)
\lvec(82 92)
\lvec(81 92)
\ifill f:0
\move(83 91)
\lvec(87 91)
\lvec(87 92)
\lvec(83 92)
\ifill f:0
\move(89 91)
\lvec(90 91)
\lvec(90 92)
\lvec(89 92)
\ifill f:0
\move(91 91)
\lvec(93 91)
\lvec(93 92)
\lvec(91 92)
\ifill f:0
\move(95 91)
\lvec(99 91)
\lvec(99 92)
\lvec(95 92)
\ifill f:0
\move(100 91)
\lvec(101 91)
\lvec(101 92)
\lvec(100 92)
\ifill f:0
\move(102 91)
\lvec(107 91)
\lvec(107 92)
\lvec(102 92)
\ifill f:0
\move(108 91)
\lvec(109 91)
\lvec(109 92)
\lvec(108 92)
\ifill f:0
\move(110 91)
\lvec(116 91)
\lvec(116 92)
\lvec(110 92)
\ifill f:0
\move(117 91)
\lvec(122 91)
\lvec(122 92)
\lvec(117 92)
\ifill f:0
\move(123 91)
\lvec(127 91)
\lvec(127 92)
\lvec(123 92)
\ifill f:0
\move(128 91)
\lvec(130 91)
\lvec(130 92)
\lvec(128 92)
\ifill f:0
\move(135 91)
\lvec(141 91)
\lvec(141 92)
\lvec(135 92)
\ifill f:0
\move(142 91)
\lvec(145 91)
\lvec(145 92)
\lvec(142 92)
\ifill f:0
\move(146 91)
\lvec(147 91)
\lvec(147 92)
\lvec(146 92)
\ifill f:0
\move(148 91)
\lvec(151 91)
\lvec(151 92)
\lvec(148 92)
\ifill f:0
\move(152 91)
\lvec(154 91)
\lvec(154 92)
\lvec(152 92)
\ifill f:0
\move(155 91)
\lvec(157 91)
\lvec(157 92)
\lvec(155 92)
\ifill f:0
\move(158 91)
\lvec(160 91)
\lvec(160 92)
\lvec(158 92)
\ifill f:0
\move(161 91)
\lvec(163 91)
\lvec(163 92)
\lvec(161 92)
\ifill f:0
\move(164 91)
\lvec(167 91)
\lvec(167 92)
\lvec(164 92)
\ifill f:0
\move(168 91)
\lvec(170 91)
\lvec(170 92)
\lvec(168 92)
\ifill f:0
\move(171 91)
\lvec(173 91)
\lvec(173 92)
\lvec(171 92)
\ifill f:0
\move(174 91)
\lvec(175 91)
\lvec(175 92)
\lvec(174 92)
\ifill f:0
\move(176 91)
\lvec(177 91)
\lvec(177 92)
\lvec(176 92)
\ifill f:0
\move(178 91)
\lvec(194 91)
\lvec(194 92)
\lvec(178 92)
\ifill f:0
\move(195 91)
\lvec(197 91)
\lvec(197 92)
\lvec(195 92)
\ifill f:0
\move(198 91)
\lvec(202 91)
\lvec(202 92)
\lvec(198 92)
\ifill f:0
\move(203 91)
\lvec(207 91)
\lvec(207 92)
\lvec(203 92)
\ifill f:0
\move(208 91)
\lvec(223 91)
\lvec(223 92)
\lvec(208 92)
\ifill f:0
\move(225 91)
\lvec(226 91)
\lvec(226 92)
\lvec(225 92)
\ifill f:0
\move(227 91)
\lvec(240 91)
\lvec(240 92)
\lvec(227 92)
\ifill f:0
\move(243 91)
\lvec(252 91)
\lvec(252 92)
\lvec(243 92)
\ifill f:0
\move(254 91)
\lvec(257 91)
\lvec(257 92)
\lvec(254 92)
\ifill f:0
\move(258 91)
\lvec(265 91)
\lvec(265 92)
\lvec(258 92)
\ifill f:0
\move(267 91)
\lvec(270 91)
\lvec(270 92)
\lvec(267 92)
\ifill f:0
\move(271 91)
\lvec(279 91)
\lvec(279 92)
\lvec(271 92)
\ifill f:0
\move(280 91)
\lvec(290 91)
\lvec(290 92)
\lvec(280 92)
\ifill f:0
\move(291 91)
\lvec(293 91)
\lvec(293 92)
\lvec(291 92)
\ifill f:0
\move(294 91)
\lvec(309 91)
\lvec(309 92)
\lvec(294 92)
\ifill f:0
\move(310 91)
\lvec(314 91)
\lvec(314 92)
\lvec(310 92)
\ifill f:0
\move(315 91)
\lvec(321 91)
\lvec(321 92)
\lvec(315 92)
\ifill f:0
\move(322 91)
\lvec(323 91)
\lvec(323 92)
\lvec(322 92)
\ifill f:0
\move(324 91)
\lvec(325 91)
\lvec(325 92)
\lvec(324 92)
\ifill f:0
\move(326 91)
\lvec(347 91)
\lvec(347 92)
\lvec(326 92)
\ifill f:0
\move(348 91)
\lvec(356 91)
\lvec(356 92)
\lvec(348 92)
\ifill f:0
\move(357 91)
\lvec(362 91)
\lvec(362 92)
\lvec(357 92)
\ifill f:0
\move(363 91)
\lvec(366 91)
\lvec(366 92)
\lvec(363 92)
\ifill f:0
\move(367 91)
\lvec(371 91)
\lvec(371 92)
\lvec(367 92)
\ifill f:0
\move(372 91)
\lvec(374 91)
\lvec(374 92)
\lvec(372 92)
\ifill f:0
\move(375 91)
\lvec(395 91)
\lvec(395 92)
\lvec(375 92)
\ifill f:0
\move(396 91)
\lvec(401 91)
\lvec(401 92)
\lvec(396 92)
\ifill f:0
\move(402 91)
\lvec(409 91)
\lvec(409 92)
\lvec(402 92)
\ifill f:0
\move(410 91)
\lvec(417 91)
\lvec(417 92)
\lvec(410 92)
\ifill f:0
\move(418 91)
\lvec(442 91)
\lvec(442 92)
\lvec(418 92)
\ifill f:0
\move(443 91)
\lvec(447 91)
\lvec(447 92)
\lvec(443 92)
\ifill f:0
\move(448 91)
\lvec(451 91)
\lvec(451 92)
\lvec(448 92)
\ifill f:0
\move(16 92)
\lvec(17 92)
\lvec(17 93)
\lvec(16 93)
\ifill f:0
\move(20 92)
\lvec(21 92)
\lvec(21 93)
\lvec(20 93)
\ifill f:0
\move(23 92)
\lvec(24 92)
\lvec(24 93)
\lvec(23 93)
\ifill f:0
\move(25 92)
\lvec(26 92)
\lvec(26 93)
\lvec(25 93)
\ifill f:0
\move(36 92)
\lvec(37 92)
\lvec(37 93)
\lvec(36 93)
\ifill f:0
\move(38 92)
\lvec(41 92)
\lvec(41 93)
\lvec(38 93)
\ifill f:0
\move(43 92)
\lvec(46 92)
\lvec(46 93)
\lvec(43 93)
\ifill f:0
\move(47 92)
\lvec(50 92)
\lvec(50 93)
\lvec(47 93)
\ifill f:0
\move(51 92)
\lvec(52 92)
\lvec(52 93)
\lvec(51 93)
\ifill f:0
\move(54 92)
\lvec(55 92)
\lvec(55 93)
\lvec(54 93)
\ifill f:0
\move(64 92)
\lvec(65 92)
\lvec(65 93)
\lvec(64 93)
\ifill f:0
\move(67 92)
\lvec(71 92)
\lvec(71 93)
\lvec(67 93)
\ifill f:0
\move(72 92)
\lvec(73 92)
\lvec(73 93)
\lvec(72 93)
\ifill f:0
\move(74 92)
\lvec(75 92)
\lvec(75 93)
\lvec(74 93)
\ifill f:0
\move(77 92)
\lvec(80 92)
\lvec(80 93)
\lvec(77 93)
\ifill f:0
\move(81 92)
\lvec(82 92)
\lvec(82 93)
\lvec(81 93)
\ifill f:0
\move(83 92)
\lvec(85 92)
\lvec(85 93)
\lvec(83 93)
\ifill f:0
\move(86 92)
\lvec(91 92)
\lvec(91 93)
\lvec(86 93)
\ifill f:0
\move(92 92)
\lvec(93 92)
\lvec(93 93)
\lvec(92 93)
\ifill f:0
\move(95 92)
\lvec(96 92)
\lvec(96 93)
\lvec(95 93)
\ifill f:0
\move(97 92)
\lvec(99 92)
\lvec(99 93)
\lvec(97 93)
\ifill f:0
\move(100 92)
\lvec(101 92)
\lvec(101 93)
\lvec(100 93)
\ifill f:0
\move(102 92)
\lvec(106 92)
\lvec(106 93)
\lvec(102 93)
\ifill f:0
\move(107 92)
\lvec(122 92)
\lvec(122 93)
\lvec(107 93)
\ifill f:0
\move(123 92)
\lvec(125 92)
\lvec(125 93)
\lvec(123 93)
\ifill f:0
\move(128 92)
\lvec(140 92)
\lvec(140 93)
\lvec(128 93)
\ifill f:0
\move(141 92)
\lvec(145 92)
\lvec(145 93)
\lvec(141 93)
\ifill f:0
\move(146 92)
\lvec(147 92)
\lvec(147 93)
\lvec(146 93)
\ifill f:0
\move(148 92)
\lvec(159 92)
\lvec(159 93)
\lvec(148 93)
\ifill f:0
\move(160 92)
\lvec(167 92)
\lvec(167 93)
\lvec(160 93)
\ifill f:0
\move(168 92)
\lvec(170 92)
\lvec(170 93)
\lvec(168 93)
\ifill f:0
\move(171 92)
\lvec(174 92)
\lvec(174 93)
\lvec(171 93)
\ifill f:0
\move(175 92)
\lvec(176 92)
\lvec(176 93)
\lvec(175 93)
\ifill f:0
\move(177 92)
\lvec(191 92)
\lvec(191 93)
\lvec(177 93)
\ifill f:0
\move(192 92)
\lvec(194 92)
\lvec(194 93)
\lvec(192 93)
\ifill f:0
\move(195 92)
\lvec(197 92)
\lvec(197 93)
\lvec(195 93)
\ifill f:0
\move(198 92)
\lvec(201 92)
\lvec(201 93)
\lvec(198 93)
\ifill f:0
\move(202 92)
\lvec(210 92)
\lvec(210 93)
\lvec(202 93)
\ifill f:0
\move(211 92)
\lvec(216 92)
\lvec(216 93)
\lvec(211 93)
\ifill f:0
\move(217 92)
\lvec(224 92)
\lvec(224 93)
\lvec(217 93)
\ifill f:0
\move(225 92)
\lvec(226 92)
\lvec(226 93)
\lvec(225 93)
\ifill f:0
\move(227 92)
\lvec(251 92)
\lvec(251 93)
\lvec(227 93)
\ifill f:0
\move(252 92)
\lvec(257 92)
\lvec(257 93)
\lvec(252 93)
\ifill f:0
\move(258 92)
\lvec(260 92)
\lvec(260 93)
\lvec(258 93)
\ifill f:0
\move(261 92)
\lvec(267 92)
\lvec(267 93)
\lvec(261 93)
\ifill f:0
\move(268 92)
\lvec(277 92)
\lvec(277 93)
\lvec(268 93)
\ifill f:0
\move(278 92)
\lvec(290 92)
\lvec(290 93)
\lvec(278 93)
\ifill f:0
\move(291 92)
\lvec(293 92)
\lvec(293 93)
\lvec(291 93)
\ifill f:0
\move(294 92)
\lvec(296 92)
\lvec(296 93)
\lvec(294 93)
\ifill f:0
\move(297 92)
\lvec(299 92)
\lvec(299 93)
\lvec(297 93)
\ifill f:0
\move(300 92)
\lvec(302 92)
\lvec(302 93)
\lvec(300 93)
\ifill f:0
\move(303 92)
\lvec(305 92)
\lvec(305 93)
\lvec(303 93)
\ifill f:0
\move(306 92)
\lvec(313 92)
\lvec(313 93)
\lvec(306 93)
\ifill f:0
\move(314 92)
\lvec(318 92)
\lvec(318 93)
\lvec(314 93)
\ifill f:0
\move(319 92)
\lvec(323 92)
\lvec(323 93)
\lvec(319 93)
\ifill f:0
\move(324 92)
\lvec(325 92)
\lvec(325 93)
\lvec(324 93)
\ifill f:0
\move(326 92)
\lvec(334 92)
\lvec(334 93)
\lvec(326 93)
\ifill f:0
\move(335 92)
\lvec(350 92)
\lvec(350 93)
\lvec(335 93)
\ifill f:0
\move(351 92)
\lvec(352 92)
\lvec(352 93)
\lvec(351 93)
\ifill f:0
\move(353 92)
\lvec(362 92)
\lvec(362 93)
\lvec(353 93)
\ifill f:0
\move(364 92)
\lvec(368 92)
\lvec(368 93)
\lvec(364 93)
\ifill f:0
\move(369 92)
\lvec(373 92)
\lvec(373 93)
\lvec(369 93)
\ifill f:0
\move(374 92)
\lvec(392 92)
\lvec(392 93)
\lvec(374 93)
\ifill f:0
\move(393 92)
\lvec(395 92)
\lvec(395 93)
\lvec(393 93)
\ifill f:0
\move(396 92)
\lvec(401 92)
\lvec(401 93)
\lvec(396 93)
\ifill f:0
\move(402 92)
\lvec(405 92)
\lvec(405 93)
\lvec(402 93)
\ifill f:0
\move(406 92)
\lvec(419 92)
\lvec(419 93)
\lvec(406 93)
\ifill f:0
\move(420 92)
\lvec(423 92)
\lvec(423 93)
\lvec(420 93)
\ifill f:0
\move(424 92)
\lvec(436 92)
\lvec(436 93)
\lvec(424 93)
\ifill f:0
\move(437 92)
\lvec(442 92)
\lvec(442 93)
\lvec(437 93)
\ifill f:0
\move(443 92)
\lvec(445 92)
\lvec(445 93)
\lvec(443 93)
\ifill f:0
\move(447 92)
\lvec(451 92)
\lvec(451 93)
\lvec(447 93)
\ifill f:0
\move(16 93)
\lvec(17 93)
\lvec(17 94)
\lvec(16 94)
\ifill f:0
\move(18 93)
\lvec(21 93)
\lvec(21 94)
\lvec(18 94)
\ifill f:0
\move(22 93)
\lvec(23 93)
\lvec(23 94)
\lvec(22 94)
\ifill f:0
\move(24 93)
\lvec(26 93)
\lvec(26 94)
\lvec(24 94)
\ifill f:0
\move(36 93)
\lvec(37 93)
\lvec(37 94)
\lvec(36 94)
\ifill f:0
\move(40 93)
\lvec(45 93)
\lvec(45 94)
\lvec(40 94)
\ifill f:0
\move(47 93)
\lvec(50 93)
\lvec(50 94)
\lvec(47 94)
\ifill f:0
\move(51 93)
\lvec(52 93)
\lvec(52 94)
\lvec(51 94)
\ifill f:0
\move(56 93)
\lvec(58 93)
\lvec(58 94)
\lvec(56 94)
\ifill f:0
\move(62 93)
\lvec(63 93)
\lvec(63 94)
\lvec(62 94)
\ifill f:0
\move(64 93)
\lvec(65 93)
\lvec(65 94)
\lvec(64 94)
\ifill f:0
\move(66 93)
\lvec(67 93)
\lvec(67 94)
\lvec(66 94)
\ifill f:0
\move(68 93)
\lvec(74 93)
\lvec(74 94)
\lvec(68 94)
\ifill f:0
\move(75 93)
\lvec(77 93)
\lvec(77 94)
\lvec(75 94)
\ifill f:0
\move(78 93)
\lvec(80 93)
\lvec(80 94)
\lvec(78 94)
\ifill f:0
\move(81 93)
\lvec(82 93)
\lvec(82 94)
\lvec(81 94)
\ifill f:0
\move(83 93)
\lvec(85 93)
\lvec(85 94)
\lvec(83 94)
\ifill f:0
\move(88 93)
\lvec(93 93)
\lvec(93 94)
\lvec(88 94)
\ifill f:0
\move(94 93)
\lvec(95 93)
\lvec(95 94)
\lvec(94 94)
\ifill f:0
\move(97 93)
\lvec(99 93)
\lvec(99 94)
\lvec(97 94)
\ifill f:0
\move(100 93)
\lvec(101 93)
\lvec(101 94)
\lvec(100 94)
\ifill f:0
\move(102 93)
\lvec(112 93)
\lvec(112 94)
\lvec(102 94)
\ifill f:0
\move(113 93)
\lvec(117 93)
\lvec(117 94)
\lvec(113 94)
\ifill f:0
\move(118 93)
\lvec(122 93)
\lvec(122 94)
\lvec(118 94)
\ifill f:0
\move(123 93)
\lvec(124 93)
\lvec(124 94)
\lvec(123 94)
\ifill f:0
\move(126 93)
\lvec(131 93)
\lvec(131 94)
\lvec(126 94)
\ifill f:0
\move(139 93)
\lvec(145 93)
\lvec(145 94)
\lvec(139 94)
\ifill f:0
\move(146 93)
\lvec(148 93)
\lvec(148 94)
\lvec(146 94)
\ifill f:0
\move(149 93)
\lvec(170 93)
\lvec(170 94)
\lvec(149 94)
\ifill f:0
\move(171 93)
\lvec(172 93)
\lvec(172 94)
\lvec(171 94)
\ifill f:0
\move(173 93)
\lvec(174 93)
\lvec(174 94)
\lvec(173 94)
\ifill f:0
\move(175 93)
\lvec(176 93)
\lvec(176 94)
\lvec(175 94)
\ifill f:0
\move(177 93)
\lvec(180 93)
\lvec(180 94)
\lvec(177 94)
\ifill f:0
\move(181 93)
\lvec(186 93)
\lvec(186 94)
\lvec(181 94)
\ifill f:0
\move(187 93)
\lvec(189 93)
\lvec(189 94)
\lvec(187 94)
\ifill f:0
\move(190 93)
\lvec(191 93)
\lvec(191 94)
\lvec(190 94)
\ifill f:0
\move(192 93)
\lvec(194 93)
\lvec(194 94)
\lvec(192 94)
\ifill f:0
\move(195 93)
\lvec(197 93)
\lvec(197 94)
\lvec(195 94)
\ifill f:0
\move(198 93)
\lvec(204 93)
\lvec(204 94)
\lvec(198 94)
\ifill f:0
\move(205 93)
\lvec(208 93)
\lvec(208 94)
\lvec(205 94)
\ifill f:0
\move(209 93)
\lvec(218 93)
\lvec(218 94)
\lvec(209 94)
\ifill f:0
\move(219 93)
\lvec(224 93)
\lvec(224 94)
\lvec(219 94)
\ifill f:0
\move(225 93)
\lvec(226 93)
\lvec(226 94)
\lvec(225 94)
\ifill f:0
\move(227 93)
\lvec(234 93)
\lvec(234 94)
\lvec(227 94)
\ifill f:0
\move(236 93)
\lvec(238 93)
\lvec(238 94)
\lvec(236 94)
\ifill f:0
\move(243 93)
\lvec(244 93)
\lvec(244 94)
\lvec(243 94)
\ifill f:0
\move(252 93)
\lvec(257 93)
\lvec(257 94)
\lvec(252 94)
\ifill f:0
\move(258 93)
\lvec(269 93)
\lvec(269 94)
\lvec(258 94)
\ifill f:0
\move(270 93)
\lvec(290 93)
\lvec(290 94)
\lvec(270 94)
\ifill f:0
\move(291 93)
\lvec(293 93)
\lvec(293 94)
\lvec(291 94)
\ifill f:0
\move(294 93)
\lvec(297 93)
\lvec(297 94)
\lvec(294 94)
\ifill f:0
\move(298 93)
\lvec(315 93)
\lvec(315 94)
\lvec(298 94)
\ifill f:0
\move(316 93)
\lvec(323 93)
\lvec(323 94)
\lvec(316 94)
\ifill f:0
\move(324 93)
\lvec(325 93)
\lvec(325 94)
\lvec(324 94)
\ifill f:0
\move(326 93)
\lvec(330 93)
\lvec(330 94)
\lvec(326 94)
\ifill f:0
\move(331 93)
\lvec(337 93)
\lvec(337 94)
\lvec(331 94)
\ifill f:0
\move(338 93)
\lvec(339 93)
\lvec(339 94)
\lvec(338 94)
\ifill f:0
\move(340 93)
\lvec(362 93)
\lvec(362 94)
\lvec(340 94)
\ifill f:0
\move(364 93)
\lvec(365 93)
\lvec(365 94)
\lvec(364 94)
\ifill f:0
\move(366 93)
\lvec(370 93)
\lvec(370 94)
\lvec(366 94)
\ifill f:0
\move(371 93)
\lvec(377 93)
\lvec(377 94)
\lvec(371 94)
\ifill f:0
\move(378 93)
\lvec(382 93)
\lvec(382 94)
\lvec(378 94)
\ifill f:0
\move(383 93)
\lvec(387 93)
\lvec(387 94)
\lvec(383 94)
\ifill f:0
\move(388 93)
\lvec(401 93)
\lvec(401 94)
\lvec(388 94)
\ifill f:0
\move(402 93)
\lvec(442 93)
\lvec(442 94)
\lvec(402 94)
\ifill f:0
\move(443 93)
\lvec(445 93)
\lvec(445 94)
\lvec(443 94)
\ifill f:0
\move(446 93)
\lvec(450 93)
\lvec(450 94)
\lvec(446 94)
\ifill f:0
\move(16 94)
\lvec(17 94)
\lvec(17 95)
\lvec(16 95)
\ifill f:0
\move(20 94)
\lvec(21 94)
\lvec(21 95)
\lvec(20 95)
\ifill f:0
\move(24 94)
\lvec(26 94)
\lvec(26 95)
\lvec(24 95)
\ifill f:0
\move(36 94)
\lvec(37 94)
\lvec(37 95)
\lvec(36 95)
\ifill f:0
\move(38 94)
\lvec(39 94)
\lvec(39 95)
\lvec(38 95)
\ifill f:0
\move(40 94)
\lvec(41 94)
\lvec(41 95)
\lvec(40 95)
\ifill f:0
\move(43 94)
\lvec(45 94)
\lvec(45 95)
\lvec(43 95)
\ifill f:0
\move(47 94)
\lvec(50 94)
\lvec(50 95)
\lvec(47 95)
\ifill f:0
\move(51 94)
\lvec(52 94)
\lvec(52 95)
\lvec(51 95)
\ifill f:0
\move(58 94)
\lvec(63 94)
\lvec(63 95)
\lvec(58 95)
\ifill f:0
\move(64 94)
\lvec(65 94)
\lvec(65 95)
\lvec(64 95)
\ifill f:0
\move(66 94)
\lvec(71 94)
\lvec(71 95)
\lvec(66 95)
\ifill f:0
\move(72 94)
\lvec(73 94)
\lvec(73 95)
\lvec(72 95)
\ifill f:0
\move(74 94)
\lvec(75 94)
\lvec(75 95)
\lvec(74 95)
\ifill f:0
\move(76 94)
\lvec(77 94)
\lvec(77 95)
\lvec(76 95)
\ifill f:0
\move(78 94)
\lvec(80 94)
\lvec(80 95)
\lvec(78 95)
\ifill f:0
\move(81 94)
\lvec(82 94)
\lvec(82 95)
\lvec(81 95)
\ifill f:0
\move(83 94)
\lvec(84 94)
\lvec(84 95)
\lvec(83 95)
\ifill f:0
\move(96 94)
\lvec(98 94)
\lvec(98 95)
\lvec(96 95)
\ifill f:0
\move(100 94)
\lvec(101 94)
\lvec(101 95)
\lvec(100 95)
\ifill f:0
\move(102 94)
\lvec(115 94)
\lvec(115 95)
\lvec(102 95)
\ifill f:0
\move(116 94)
\lvec(122 94)
\lvec(122 95)
\lvec(116 95)
\ifill f:0
\move(123 94)
\lvec(124 94)
\lvec(124 95)
\lvec(123 95)
\ifill f:0
\move(125 94)
\lvec(129 94)
\lvec(129 95)
\lvec(125 95)
\ifill f:0
\move(131 94)
\lvec(138 94)
\lvec(138 95)
\lvec(131 95)
\ifill f:0
\move(139 94)
\lvec(145 94)
\lvec(145 95)
\lvec(139 95)
\ifill f:0
\move(146 94)
\lvec(150 94)
\lvec(150 95)
\lvec(146 95)
\ifill f:0
\move(151 94)
\lvec(155 94)
\lvec(155 95)
\lvec(151 95)
\ifill f:0
\move(156 94)
\lvec(160 94)
\lvec(160 95)
\lvec(156 95)
\ifill f:0
\move(161 94)
\lvec(163 94)
\lvec(163 95)
\lvec(161 95)
\ifill f:0
\move(164 94)
\lvec(170 94)
\lvec(170 95)
\lvec(164 95)
\ifill f:0
\move(171 94)
\lvec(172 94)
\lvec(172 95)
\lvec(171 95)
\ifill f:0
\move(173 94)
\lvec(179 94)
\lvec(179 95)
\lvec(173 95)
\ifill f:0
\move(180 94)
\lvec(183 94)
\lvec(183 95)
\lvec(180 95)
\ifill f:0
\move(184 94)
\lvec(185 94)
\lvec(185 95)
\lvec(184 95)
\ifill f:0
\move(186 94)
\lvec(192 94)
\lvec(192 95)
\lvec(186 95)
\ifill f:0
\move(193 94)
\lvec(194 94)
\lvec(194 95)
\lvec(193 95)
\ifill f:0
\move(195 94)
\lvec(197 94)
\lvec(197 95)
\lvec(195 95)
\ifill f:0
\move(198 94)
\lvec(203 94)
\lvec(203 95)
\lvec(198 95)
\ifill f:0
\move(204 94)
\lvec(214 94)
\lvec(214 95)
\lvec(204 95)
\ifill f:0
\move(215 94)
\lvec(219 94)
\lvec(219 95)
\lvec(215 95)
\ifill f:0
\move(220 94)
\lvec(226 94)
\lvec(226 95)
\lvec(220 95)
\ifill f:0
\move(227 94)
\lvec(232 94)
\lvec(232 95)
\lvec(227 95)
\ifill f:0
\move(234 94)
\lvec(257 94)
\lvec(257 95)
\lvec(234 95)
\ifill f:0
\move(258 94)
\lvec(279 94)
\lvec(279 95)
\lvec(258 95)
\ifill f:0
\move(281 94)
\lvec(290 94)
\lvec(290 95)
\lvec(281 95)
\ifill f:0
\move(291 94)
\lvec(298 94)
\lvec(298 95)
\lvec(291 95)
\ifill f:0
\move(299 94)
\lvec(301 94)
\lvec(301 95)
\lvec(299 95)
\ifill f:0
\move(302 94)
\lvec(305 94)
\lvec(305 95)
\lvec(302 95)
\ifill f:0
\move(306 94)
\lvec(308 94)
\lvec(308 95)
\lvec(306 95)
\ifill f:0
\move(309 94)
\lvec(320 94)
\lvec(320 95)
\lvec(309 95)
\ifill f:0
\move(321 94)
\lvec(323 94)
\lvec(323 95)
\lvec(321 95)
\ifill f:0
\move(324 94)
\lvec(325 94)
\lvec(325 95)
\lvec(324 95)
\ifill f:0
\move(326 94)
\lvec(328 94)
\lvec(328 95)
\lvec(326 95)
\ifill f:0
\move(329 94)
\lvec(340 94)
\lvec(340 95)
\lvec(329 95)
\ifill f:0
\move(341 94)
\lvec(353 94)
\lvec(353 95)
\lvec(341 95)
\ifill f:0
\move(354 94)
\lvec(355 94)
\lvec(355 95)
\lvec(354 95)
\ifill f:0
\move(356 94)
\lvec(359 94)
\lvec(359 95)
\lvec(356 95)
\ifill f:0
\move(360 94)
\lvec(362 94)
\lvec(362 95)
\lvec(360 95)
\ifill f:0
\move(364 94)
\lvec(365 94)
\lvec(365 95)
\lvec(364 95)
\ifill f:0
\move(366 94)
\lvec(367 94)
\lvec(367 95)
\lvec(366 95)
\ifill f:0
\move(368 94)
\lvec(383 94)
\lvec(383 95)
\lvec(368 95)
\ifill f:0
\move(384 94)
\lvec(388 94)
\lvec(388 95)
\lvec(384 95)
\ifill f:0
\move(389 94)
\lvec(401 94)
\lvec(401 95)
\lvec(389 95)
\ifill f:0
\move(402 94)
\lvec(410 94)
\lvec(410 95)
\lvec(402 95)
\ifill f:0
\move(411 94)
\lvec(416 94)
\lvec(416 95)
\lvec(411 95)
\ifill f:0
\move(417 94)
\lvec(426 94)
\lvec(426 95)
\lvec(417 95)
\ifill f:0
\move(427 94)
\lvec(442 94)
\lvec(442 95)
\lvec(427 95)
\ifill f:0
\move(443 94)
\lvec(445 94)
\lvec(445 95)
\lvec(443 95)
\ifill f:0
\move(446 94)
\lvec(449 94)
\lvec(449 95)
\lvec(446 95)
\ifill f:0
\move(450 94)
\lvec(451 94)
\lvec(451 95)
\lvec(450 95)
\ifill f:0
\move(15 95)
\lvec(17 95)
\lvec(17 96)
\lvec(15 96)
\ifill f:0
\move(20 95)
\lvec(21 95)
\lvec(21 96)
\lvec(20 96)
\ifill f:0
\move(24 95)
\lvec(26 95)
\lvec(26 96)
\lvec(24 96)
\ifill f:0
\move(36 95)
\lvec(37 95)
\lvec(37 96)
\lvec(36 96)
\ifill f:0
\move(38 95)
\lvec(39 95)
\lvec(39 96)
\lvec(38 96)
\ifill f:0
\move(40 95)
\lvec(41 95)
\lvec(41 96)
\lvec(40 96)
\ifill f:0
\move(42 95)
\lvec(46 95)
\lvec(46 96)
\lvec(42 96)
\ifill f:0
\move(47 95)
\lvec(50 95)
\lvec(50 96)
\lvec(47 96)
\ifill f:0
\move(52 95)
\lvec(53 95)
\lvec(53 96)
\lvec(52 96)
\ifill f:0
\move(54 95)
\lvec(55 95)
\lvec(55 96)
\lvec(54 96)
\ifill f:0
\move(56 95)
\lvec(57 95)
\lvec(57 96)
\lvec(56 96)
\ifill f:0
\move(59 95)
\lvec(63 95)
\lvec(63 96)
\lvec(59 96)
\ifill f:0
\move(64 95)
\lvec(65 95)
\lvec(65 96)
\lvec(64 96)
\ifill f:0
\move(66 95)
\lvec(75 95)
\lvec(75 96)
\lvec(66 96)
\ifill f:0
\move(76 95)
\lvec(78 95)
\lvec(78 96)
\lvec(76 96)
\ifill f:0
\move(79 95)
\lvec(82 95)
\lvec(82 96)
\lvec(79 96)
\ifill f:0
\move(83 95)
\lvec(84 95)
\lvec(84 96)
\lvec(83 96)
\ifill f:0
\move(86 95)
\lvec(93 95)
\lvec(93 96)
\lvec(86 96)
\ifill f:0
\move(95 95)
\lvec(98 95)
\lvec(98 96)
\lvec(95 96)
\ifill f:0
\move(100 95)
\lvec(101 95)
\lvec(101 96)
\lvec(100 96)
\ifill f:0
\move(103 95)
\lvec(106 95)
\lvec(106 96)
\lvec(103 96)
\ifill f:0
\move(107 95)
\lvec(108 95)
\lvec(108 96)
\lvec(107 96)
\ifill f:0
\move(109 95)
\lvec(111 95)
\lvec(111 96)
\lvec(109 96)
\ifill f:0
\move(112 95)
\lvec(122 95)
\lvec(122 96)
\lvec(112 96)
\ifill f:0
\move(123 95)
\lvec(124 95)
\lvec(124 96)
\lvec(123 96)
\ifill f:0
\move(125 95)
\lvec(127 95)
\lvec(127 96)
\lvec(125 96)
\ifill f:0
\move(129 95)
\lvec(131 95)
\lvec(131 96)
\lvec(129 96)
\ifill f:0
\move(135 95)
\lvec(145 95)
\lvec(145 96)
\lvec(135 96)
\ifill f:0
\move(146 95)
\lvec(163 95)
\lvec(163 96)
\lvec(146 96)
\ifill f:0
\move(164 95)
\lvec(170 95)
\lvec(170 96)
\lvec(164 96)
\ifill f:0
\move(171 95)
\lvec(172 95)
\lvec(172 96)
\lvec(171 96)
\ifill f:0
\move(173 95)
\lvec(175 95)
\lvec(175 96)
\lvec(173 96)
\ifill f:0
\move(176 95)
\lvec(180 95)
\lvec(180 96)
\lvec(176 96)
\ifill f:0
\move(181 95)
\lvec(184 95)
\lvec(184 96)
\lvec(181 96)
\ifill f:0
\move(185 95)
\lvec(186 95)
\lvec(186 96)
\lvec(185 96)
\ifill f:0
\move(187 95)
\lvec(188 95)
\lvec(188 96)
\lvec(187 96)
\ifill f:0
\move(189 95)
\lvec(190 95)
\lvec(190 96)
\lvec(189 96)
\ifill f:0
\move(191 95)
\lvec(192 95)
\lvec(192 96)
\lvec(191 96)
\ifill f:0
\move(193 95)
\lvec(194 95)
\lvec(194 96)
\lvec(193 96)
\ifill f:0
\move(195 95)
\lvec(197 95)
\lvec(197 96)
\lvec(195 96)
\ifill f:0
\move(198 95)
\lvec(212 95)
\lvec(212 96)
\lvec(198 96)
\ifill f:0
\move(213 95)
\lvec(220 95)
\lvec(220 96)
\lvec(213 96)
\ifill f:0
\move(221 95)
\lvec(226 95)
\lvec(226 96)
\lvec(221 96)
\ifill f:0
\move(227 95)
\lvec(231 95)
\lvec(231 96)
\lvec(227 96)
\ifill f:0
\move(232 95)
\lvec(240 95)
\lvec(240 96)
\lvec(232 96)
\ifill f:0
\move(241 95)
\lvec(257 95)
\lvec(257 96)
\lvec(241 96)
\ifill f:0
\move(258 95)
\lvec(267 95)
\lvec(267 96)
\lvec(258 96)
\ifill f:0
\move(268 95)
\lvec(277 95)
\lvec(277 96)
\lvec(268 96)
\ifill f:0
\move(278 95)
\lvec(284 95)
\lvec(284 96)
\lvec(278 96)
\ifill f:0
\move(285 95)
\lvec(290 95)
\lvec(290 96)
\lvec(285 96)
\ifill f:0
\move(291 95)
\lvec(294 95)
\lvec(294 96)
\lvec(291 96)
\ifill f:0
\move(295 95)
\lvec(299 95)
\lvec(299 96)
\lvec(295 96)
\ifill f:0
\move(300 95)
\lvec(310 95)
\lvec(310 96)
\lvec(300 96)
\ifill f:0
\move(311 95)
\lvec(320 95)
\lvec(320 96)
\lvec(311 96)
\ifill f:0
\move(321 95)
\lvec(323 95)
\lvec(323 96)
\lvec(321 96)
\ifill f:0
\move(324 95)
\lvec(325 95)
\lvec(325 96)
\lvec(324 96)
\ifill f:0
\move(326 95)
\lvec(341 95)
\lvec(341 96)
\lvec(326 96)
\ifill f:0
\move(342 95)
\lvec(359 95)
\lvec(359 96)
\lvec(342 96)
\ifill f:0
\move(360 95)
\lvec(362 95)
\lvec(362 96)
\lvec(360 96)
\ifill f:0
\move(364 95)
\lvec(367 95)
\lvec(367 96)
\lvec(364 96)
\ifill f:0
\move(368 95)
\lvec(369 95)
\lvec(369 96)
\lvec(368 96)
\ifill f:0
\move(370 95)
\lvec(371 95)
\lvec(371 96)
\lvec(370 96)
\ifill f:0
\move(372 95)
\lvec(384 95)
\lvec(384 96)
\lvec(372 96)
\ifill f:0
\move(385 95)
\lvec(401 95)
\lvec(401 96)
\lvec(385 96)
\ifill f:0
\move(402 95)
\lvec(415 95)
\lvec(415 96)
\lvec(402 96)
\ifill f:0
\move(416 95)
\lvec(418 95)
\lvec(418 96)
\lvec(416 96)
\ifill f:0
\move(419 95)
\lvec(421 95)
\lvec(421 96)
\lvec(419 96)
\ifill f:0
\move(422 95)
\lvec(424 95)
\lvec(424 96)
\lvec(422 96)
\ifill f:0
\move(425 95)
\lvec(434 95)
\lvec(434 96)
\lvec(425 96)
\ifill f:0
\move(435 95)
\lvec(442 95)
\lvec(442 96)
\lvec(435 96)
\ifill f:0
\move(443 95)
\lvec(451 95)
\lvec(451 96)
\lvec(443 96)
\ifill f:0
\move(15 96)
\lvec(17 96)
\lvec(17 97)
\lvec(15 97)
\ifill f:0
\move(19 96)
\lvec(21 96)
\lvec(21 97)
\lvec(19 97)
\ifill f:0
\move(22 96)
\lvec(26 96)
\lvec(26 97)
\lvec(22 97)
\ifill f:0
\move(36 96)
\lvec(37 96)
\lvec(37 97)
\lvec(36 97)
\ifill f:0
\move(38 96)
\lvec(42 96)
\lvec(42 97)
\lvec(38 97)
\ifill f:0
\move(43 96)
\lvec(45 96)
\lvec(45 97)
\lvec(43 97)
\ifill f:0
\move(47 96)
\lvec(50 96)
\lvec(50 97)
\lvec(47 97)
\ifill f:0
\move(57 96)
\lvec(58 96)
\lvec(58 97)
\lvec(57 97)
\ifill f:0
\move(61 96)
\lvec(65 96)
\lvec(65 97)
\lvec(61 97)
\ifill f:0
\move(66 96)
\lvec(70 96)
\lvec(70 97)
\lvec(66 97)
\ifill f:0
\move(72 96)
\lvec(73 96)
\lvec(73 97)
\lvec(72 97)
\ifill f:0
\move(75 96)
\lvec(78 96)
\lvec(78 97)
\lvec(75 97)
\ifill f:0
\move(79 96)
\lvec(82 96)
\lvec(82 97)
\lvec(79 97)
\ifill f:0
\move(85 96)
\lvec(87 96)
\lvec(87 97)
\lvec(85 97)
\ifill f:0
\move(91 96)
\lvec(93 96)
\lvec(93 97)
\lvec(91 97)
\ifill f:0
\move(94 96)
\lvec(96 96)
\lvec(96 97)
\lvec(94 97)
\ifill f:0
\move(97 96)
\lvec(98 96)
\lvec(98 97)
\lvec(97 97)
\ifill f:0
\move(100 96)
\lvec(101 96)
\lvec(101 97)
\lvec(100 97)
\ifill f:0
\move(102 96)
\lvec(103 96)
\lvec(103 97)
\lvec(102 97)
\ifill f:0
\move(104 96)
\lvec(109 96)
\lvec(109 97)
\lvec(104 97)
\ifill f:0
\move(110 96)
\lvec(116 96)
\lvec(116 97)
\lvec(110 97)
\ifill f:0
\move(117 96)
\lvec(122 96)
\lvec(122 97)
\lvec(117 97)
\ifill f:0
\move(124 96)
\lvec(126 96)
\lvec(126 97)
\lvec(124 97)
\ifill f:0
\move(128 96)
\lvec(130 96)
\lvec(130 97)
\lvec(128 97)
\ifill f:0
\move(132 96)
\lvec(135 96)
\lvec(135 97)
\lvec(132 97)
\ifill f:0
\move(138 96)
\lvec(145 96)
\lvec(145 97)
\lvec(138 97)
\ifill f:0
\move(146 96)
\lvec(155 96)
\lvec(155 97)
\lvec(146 97)
\ifill f:0
\move(156 96)
\lvec(161 96)
\lvec(161 97)
\lvec(156 97)
\ifill f:0
\move(162 96)
\lvec(170 96)
\lvec(170 97)
\lvec(162 97)
\ifill f:0
\move(171 96)
\lvec(173 96)
\lvec(173 97)
\lvec(171 97)
\ifill f:0
\move(174 96)
\lvec(176 96)
\lvec(176 97)
\lvec(174 97)
\ifill f:0
\move(177 96)
\lvec(183 96)
\lvec(183 97)
\lvec(177 97)
\ifill f:0
\move(184 96)
\lvec(197 96)
\lvec(197 97)
\lvec(184 97)
\ifill f:0
\move(198 96)
\lvec(226 96)
\lvec(226 97)
\lvec(198 97)
\ifill f:0
\move(227 96)
\lvec(229 96)
\lvec(229 97)
\lvec(227 97)
\ifill f:0
\move(231 96)
\lvec(235 96)
\lvec(235 97)
\lvec(231 97)
\ifill f:0
\move(236 96)
\lvec(242 96)
\lvec(242 97)
\lvec(236 97)
\ifill f:0
\move(243 96)
\lvec(257 96)
\lvec(257 97)
\lvec(243 97)
\ifill f:0
\move(258 96)
\lvec(272 96)
\lvec(272 97)
\lvec(258 97)
\ifill f:0
\move(273 96)
\lvec(282 96)
\lvec(282 97)
\lvec(273 97)
\ifill f:0
\move(284 96)
\lvec(290 96)
\lvec(290 97)
\lvec(284 97)
\ifill f:0
\move(291 96)
\lvec(295 96)
\lvec(295 97)
\lvec(291 97)
\ifill f:0
\move(296 96)
\lvec(323 96)
\lvec(323 97)
\lvec(296 97)
\ifill f:0
\move(324 96)
\lvec(325 96)
\lvec(325 97)
\lvec(324 97)
\ifill f:0
\move(326 96)
\lvec(337 96)
\lvec(337 97)
\lvec(326 97)
\ifill f:0
\move(338 96)
\lvec(352 96)
\lvec(352 97)
\lvec(338 97)
\ifill f:0
\move(353 96)
\lvec(359 96)
\lvec(359 97)
\lvec(353 97)
\ifill f:0
\move(360 96)
\lvec(362 96)
\lvec(362 97)
\lvec(360 97)
\ifill f:0
\move(363 96)
\lvec(385 96)
\lvec(385 97)
\lvec(363 97)
\ifill f:0
\move(386 96)
\lvec(387 96)
\lvec(387 97)
\lvec(386 97)
\ifill f:0
\move(388 96)
\lvec(394 96)
\lvec(394 97)
\lvec(388 97)
\ifill f:0
\move(395 96)
\lvec(401 96)
\lvec(401 97)
\lvec(395 97)
\ifill f:0
\move(402 96)
\lvec(414 96)
\lvec(414 97)
\lvec(402 97)
\ifill f:0
\move(415 96)
\lvec(442 96)
\lvec(442 97)
\lvec(415 97)
\ifill f:0
\move(443 96)
\lvec(444 96)
\lvec(444 97)
\lvec(443 97)
\ifill f:0
\move(445 96)
\lvec(448 96)
\lvec(448 97)
\lvec(445 97)
\ifill f:0
\move(449 96)
\lvec(451 96)
\lvec(451 97)
\lvec(449 97)
\ifill f:0
\move(15 97)
\lvec(17 97)
\lvec(17 98)
\lvec(15 98)
\ifill f:0
\move(18 97)
\lvec(19 97)
\lvec(19 98)
\lvec(18 98)
\ifill f:0
\move(20 97)
\lvec(21 97)
\lvec(21 98)
\lvec(20 98)
\ifill f:0
\move(23 97)
\lvec(26 97)
\lvec(26 98)
\lvec(23 98)
\ifill f:0
\move(36 97)
\lvec(37 97)
\lvec(37 98)
\lvec(36 98)
\ifill f:0
\move(38 97)
\lvec(39 97)
\lvec(39 98)
\lvec(38 98)
\ifill f:0
\move(40 97)
\lvec(41 97)
\lvec(41 98)
\lvec(40 98)
\ifill f:0
\move(42 97)
\lvec(43 97)
\lvec(43 98)
\lvec(42 98)
\ifill f:0
\move(44 97)
\lvec(45 97)
\lvec(45 98)
\lvec(44 98)
\ifill f:0
\move(49 97)
\lvec(50 97)
\lvec(50 98)
\lvec(49 98)
\ifill f:0
\move(56 97)
\lvec(57 97)
\lvec(57 98)
\lvec(56 98)
\ifill f:0
\move(58 97)
\lvec(60 97)
\lvec(60 98)
\lvec(58 98)
\ifill f:0
\move(62 97)
\lvec(65 97)
\lvec(65 98)
\lvec(62 98)
\ifill f:0
\move(66 97)
\lvec(75 97)
\lvec(75 98)
\lvec(66 98)
\ifill f:0
\move(76 97)
\lvec(82 97)
\lvec(82 98)
\lvec(76 98)
\ifill f:0
\move(84 97)
\lvec(85 97)
\lvec(85 98)
\lvec(84 98)
\ifill f:0
\move(88 97)
\lvec(91 97)
\lvec(91 98)
\lvec(88 98)
\ifill f:0
\move(92 97)
\lvec(93 97)
\lvec(93 98)
\lvec(92 98)
\ifill f:0
\move(97 97)
\lvec(98 97)
\lvec(98 98)
\lvec(97 98)
\ifill f:0
\move(100 97)
\lvec(101 97)
\lvec(101 98)
\lvec(100 98)
\ifill f:0
\move(102 97)
\lvec(111 97)
\lvec(111 98)
\lvec(102 98)
\ifill f:0
\move(112 97)
\lvec(113 97)
\lvec(113 98)
\lvec(112 98)
\ifill f:0
\move(114 97)
\lvec(122 97)
\lvec(122 98)
\lvec(114 98)
\ifill f:0
\move(124 97)
\lvec(125 97)
\lvec(125 98)
\lvec(124 98)
\ifill f:0
\move(127 97)
\lvec(129 97)
\lvec(129 98)
\lvec(127 98)
\ifill f:0
\move(130 97)
\lvec(132 97)
\lvec(132 98)
\lvec(130 98)
\ifill f:0
\move(134 97)
\lvec(138 97)
\lvec(138 98)
\lvec(134 98)
\ifill f:0
\move(140 97)
\lvec(145 97)
\lvec(145 98)
\lvec(140 98)
\ifill f:0
\move(146 97)
\lvec(159 97)
\lvec(159 98)
\lvec(146 98)
\ifill f:0
\move(160 97)
\lvec(165 97)
\lvec(165 98)
\lvec(160 98)
\ifill f:0
\move(166 97)
\lvec(170 97)
\lvec(170 98)
\lvec(166 98)
\ifill f:0
\move(171 97)
\lvec(173 97)
\lvec(173 98)
\lvec(171 98)
\ifill f:0
\move(175 97)
\lvec(176 97)
\lvec(176 98)
\lvec(175 98)
\ifill f:0
\move(177 97)
\lvec(179 97)
\lvec(179 98)
\lvec(177 98)
\ifill f:0
\move(180 97)
\lvec(187 97)
\lvec(187 98)
\lvec(180 98)
\ifill f:0
\move(188 97)
\lvec(189 97)
\lvec(189 98)
\lvec(188 98)
\ifill f:0
\move(190 97)
\lvec(191 97)
\lvec(191 98)
\lvec(190 98)
\ifill f:0
\move(192 97)
\lvec(193 97)
\lvec(193 98)
\lvec(192 98)
\ifill f:0
\move(194 97)
\lvec(195 97)
\lvec(195 98)
\lvec(194 98)
\ifill f:0
\move(196 97)
\lvec(197 97)
\lvec(197 98)
\lvec(196 98)
\ifill f:0
\move(198 97)
\lvec(201 97)
\lvec(201 98)
\lvec(198 98)
\ifill f:0
\move(202 97)
\lvec(221 97)
\lvec(221 98)
\lvec(202 98)
\ifill f:0
\move(222 97)
\lvec(226 97)
\lvec(226 98)
\lvec(222 98)
\ifill f:0
\move(227 97)
\lvec(229 97)
\lvec(229 98)
\lvec(227 98)
\ifill f:0
\move(230 97)
\lvec(234 97)
\lvec(234 98)
\lvec(230 98)
\ifill f:0
\move(235 97)
\lvec(240 97)
\lvec(240 98)
\lvec(235 98)
\ifill f:0
\move(241 97)
\lvec(248 97)
\lvec(248 98)
\lvec(241 98)
\ifill f:0
\move(249 97)
\lvec(257 97)
\lvec(257 98)
\lvec(249 98)
\ifill f:0
\move(258 97)
\lvec(280 97)
\lvec(280 98)
\lvec(258 98)
\ifill f:0
\move(281 97)
\lvec(290 97)
\lvec(290 98)
\lvec(281 98)
\ifill f:0
\move(291 97)
\lvec(296 97)
\lvec(296 98)
\lvec(291 98)
\ifill f:0
\move(297 97)
\lvec(302 97)
\lvec(302 98)
\lvec(297 98)
\ifill f:0
\move(303 97)
\lvec(311 97)
\lvec(311 98)
\lvec(303 98)
\ifill f:0
\move(312 97)
\lvec(315 97)
\lvec(315 98)
\lvec(312 98)
\ifill f:0
\move(316 97)
\lvec(319 97)
\lvec(319 98)
\lvec(316 98)
\ifill f:0
\move(320 97)
\lvec(323 97)
\lvec(323 98)
\lvec(320 98)
\ifill f:0
\move(324 97)
\lvec(325 97)
\lvec(325 98)
\lvec(324 98)
\ifill f:0
\move(326 97)
\lvec(329 97)
\lvec(329 98)
\lvec(326 98)
\ifill f:0
\move(330 97)
\lvec(338 97)
\lvec(338 98)
\lvec(330 98)
\ifill f:0
\move(339 97)
\lvec(341 97)
\lvec(341 98)
\lvec(339 98)
\ifill f:0
\move(342 97)
\lvec(362 97)
\lvec(362 98)
\lvec(342 98)
\ifill f:0
\move(363 97)
\lvec(370 97)
\lvec(370 98)
\lvec(363 98)
\ifill f:0
\move(371 97)
\lvec(372 97)
\lvec(372 98)
\lvec(371 98)
\ifill f:0
\move(373 97)
\lvec(392 97)
\lvec(392 98)
\lvec(373 98)
\ifill f:0
\move(393 97)
\lvec(401 97)
\lvec(401 98)
\lvec(393 98)
\ifill f:0
\move(402 97)
\lvec(432 97)
\lvec(432 98)
\lvec(402 98)
\ifill f:0
\move(433 97)
\lvec(435 97)
\lvec(435 98)
\lvec(433 98)
\ifill f:0
\move(436 97)
\lvec(438 97)
\lvec(438 98)
\lvec(436 98)
\ifill f:0
\move(439 97)
\lvec(442 97)
\lvec(442 98)
\lvec(439 98)
\ifill f:0
\move(443 97)
\lvec(444 97)
\lvec(444 98)
\lvec(443 98)
\ifill f:0
\move(445 97)
\lvec(451 97)
\lvec(451 98)
\lvec(445 98)
\ifill f:0
\move(16 98)
\lvec(17 98)
\lvec(17 99)
\lvec(16 99)
\ifill f:0
\move(20 98)
\lvec(21 98)
\lvec(21 99)
\lvec(20 99)
\ifill f:0
\move(24 98)
\lvec(26 98)
\lvec(26 99)
\lvec(24 99)
\ifill f:0
\move(36 98)
\lvec(37 98)
\lvec(37 99)
\lvec(36 99)
\ifill f:0
\move(38 98)
\lvec(39 98)
\lvec(39 99)
\lvec(38 99)
\ifill f:0
\move(40 98)
\lvec(42 98)
\lvec(42 99)
\lvec(40 99)
\ifill f:0
\move(43 98)
\lvec(46 98)
\lvec(46 99)
\lvec(43 99)
\ifill f:0
\move(49 98)
\lvec(50 98)
\lvec(50 99)
\lvec(49 99)
\ifill f:0
\move(54 98)
\lvec(55 98)
\lvec(55 99)
\lvec(54 99)
\ifill f:0
\move(57 98)
\lvec(58 98)
\lvec(58 99)
\lvec(57 99)
\ifill f:0
\move(59 98)
\lvec(61 98)
\lvec(61 99)
\lvec(59 99)
\ifill f:0
\move(62 98)
\lvec(65 98)
\lvec(65 99)
\lvec(62 99)
\ifill f:0
\move(66 98)
\lvec(74 98)
\lvec(74 99)
\lvec(66 99)
\ifill f:0
\move(75 98)
\lvec(79 98)
\lvec(79 99)
\lvec(75 99)
\ifill f:0
\move(80 98)
\lvec(82 98)
\lvec(82 99)
\lvec(80 99)
\ifill f:0
\move(84 98)
\lvec(85 98)
\lvec(85 99)
\lvec(84 99)
\ifill f:0
\move(87 98)
\lvec(89 98)
\lvec(89 99)
\lvec(87 99)
\ifill f:0
\move(91 98)
\lvec(93 98)
\lvec(93 99)
\lvec(91 99)
\ifill f:0
\move(97 98)
\lvec(98 98)
\lvec(98 99)
\lvec(97 99)
\ifill f:0
\move(100 98)
\lvec(101 98)
\lvec(101 99)
\lvec(100 99)
\ifill f:0
\move(102 98)
\lvec(104 98)
\lvec(104 99)
\lvec(102 99)
\ifill f:0
\move(105 98)
\lvec(109 98)
\lvec(109 99)
\lvec(105 99)
\ifill f:0
\move(110 98)
\lvec(119 98)
\lvec(119 99)
\lvec(110 99)
\ifill f:0
\move(120 98)
\lvec(122 98)
\lvec(122 99)
\lvec(120 99)
\ifill f:0
\move(124 98)
\lvec(125 98)
\lvec(125 99)
\lvec(124 99)
\ifill f:0
\move(126 98)
\lvec(128 98)
\lvec(128 99)
\lvec(126 99)
\ifill f:0
\move(129 98)
\lvec(131 98)
\lvec(131 99)
\lvec(129 99)
\ifill f:0
\move(132 98)
\lvec(134 98)
\lvec(134 99)
\lvec(132 99)
\ifill f:0
\move(136 98)
\lvec(138 98)
\lvec(138 99)
\lvec(136 99)
\ifill f:0
\move(139 98)
\lvec(140 98)
\lvec(140 99)
\lvec(139 99)
\ifill f:0
\move(141 98)
\lvec(145 98)
\lvec(145 99)
\lvec(141 99)
\ifill f:0
\move(146 98)
\lvec(149 98)
\lvec(149 99)
\lvec(146 99)
\ifill f:0
\move(152 98)
\lvec(153 98)
\lvec(153 99)
\lvec(152 99)
\ifill f:0
\move(154 98)
\lvec(155 98)
\lvec(155 99)
\lvec(154 99)
\ifill f:0
\move(156 98)
\lvec(164 98)
\lvec(164 99)
\lvec(156 99)
\ifill f:0
\move(165 98)
\lvec(170 98)
\lvec(170 99)
\lvec(165 99)
\ifill f:0
\move(171 98)
\lvec(174 98)
\lvec(174 99)
\lvec(171 99)
\ifill f:0
\move(175 98)
\lvec(177 98)
\lvec(177 99)
\lvec(175 99)
\ifill f:0
\move(178 98)
\lvec(180 98)
\lvec(180 99)
\lvec(178 99)
\ifill f:0
\move(181 98)
\lvec(183 98)
\lvec(183 99)
\lvec(181 99)
\ifill f:0
\move(184 98)
\lvec(186 98)
\lvec(186 99)
\lvec(184 99)
\ifill f:0
\move(187 98)
\lvec(191 98)
\lvec(191 99)
\lvec(187 99)
\ifill f:0
\move(192 98)
\lvec(193 98)
\lvec(193 99)
\lvec(192 99)
\ifill f:0
\move(194 98)
\lvec(195 98)
\lvec(195 99)
\lvec(194 99)
\ifill f:0
\move(196 98)
\lvec(197 98)
\lvec(197 99)
\lvec(196 99)
\ifill f:0
\move(198 98)
\lvec(199 98)
\lvec(199 99)
\lvec(198 99)
\ifill f:0
\move(200 98)
\lvec(201 98)
\lvec(201 99)
\lvec(200 99)
\ifill f:0
\move(202 98)
\lvec(203 98)
\lvec(203 99)
\lvec(202 99)
\ifill f:0
\move(204 98)
\lvec(226 98)
\lvec(226 99)
\lvec(204 99)
\ifill f:0
\move(227 98)
\lvec(228 98)
\lvec(228 99)
\lvec(227 99)
\ifill f:0
\move(229 98)
\lvec(242 98)
\lvec(242 99)
\lvec(229 99)
\ifill f:0
\move(244 98)
\lvec(250 98)
\lvec(250 99)
\lvec(244 99)
\ifill f:0
\move(251 98)
\lvec(257 98)
\lvec(257 99)
\lvec(251 99)
\ifill f:0
\move(258 98)
\lvec(266 98)
\lvec(266 99)
\lvec(258 99)
\ifill f:0
\move(276 98)
\lvec(279 98)
\lvec(279 99)
\lvec(276 99)
\ifill f:0
\move(280 98)
\lvec(290 98)
\lvec(290 99)
\lvec(280 99)
\ifill f:0
\move(292 98)
\lvec(297 98)
\lvec(297 99)
\lvec(292 99)
\ifill f:0
\move(299 98)
\lvec(309 98)
\lvec(309 99)
\lvec(299 99)
\ifill f:0
\move(311 98)
\lvec(314 98)
\lvec(314 99)
\lvec(311 99)
\ifill f:0
\move(315 98)
\lvec(318 98)
\lvec(318 99)
\lvec(315 99)
\ifill f:0
\move(319 98)
\lvec(322 98)
\lvec(322 99)
\lvec(319 99)
\ifill f:0
\move(323 98)
\lvec(325 98)
\lvec(325 99)
\lvec(323 99)
\ifill f:0
\move(326 98)
\lvec(356 98)
\lvec(356 99)
\lvec(326 99)
\ifill f:0
\move(357 98)
\lvec(362 98)
\lvec(362 99)
\lvec(357 99)
\ifill f:0
\move(363 98)
\lvec(366 98)
\lvec(366 99)
\lvec(363 99)
\ifill f:0
\move(367 98)
\lvec(375 98)
\lvec(375 99)
\lvec(367 99)
\ifill f:0
\move(376 98)
\lvec(377 98)
\lvec(377 99)
\lvec(376 99)
\ifill f:0
\move(378 98)
\lvec(401 98)
\lvec(401 99)
\lvec(378 99)
\ifill f:0
\move(402 98)
\lvec(410 98)
\lvec(410 99)
\lvec(402 99)
\ifill f:0
\move(411 98)
\lvec(442 98)
\lvec(442 99)
\lvec(411 99)
\ifill f:0
\move(443 98)
\lvec(444 98)
\lvec(444 99)
\lvec(443 99)
\ifill f:0
\move(445 98)
\lvec(447 98)
\lvec(447 99)
\lvec(445 99)
\ifill f:0
\move(448 98)
\lvec(450 98)
\lvec(450 99)
\lvec(448 99)
\ifill f:0
\move(16 99)
\lvec(17 99)
\lvec(17 100)
\lvec(16 100)
\ifill f:0
\move(19 99)
\lvec(21 99)
\lvec(21 100)
\lvec(19 100)
\ifill f:0
\move(25 99)
\lvec(26 99)
\lvec(26 100)
\lvec(25 100)
\ifill f:0
\move(36 99)
\lvec(37 99)
\lvec(37 100)
\lvec(36 100)
\ifill f:0
\move(40 99)
\lvec(41 99)
\lvec(41 100)
\lvec(40 100)
\ifill f:0
\move(42 99)
\lvec(43 99)
\lvec(43 100)
\lvec(42 100)
\ifill f:0
\move(44 99)
\lvec(45 99)
\lvec(45 100)
\lvec(44 100)
\ifill f:0
\move(47 99)
\lvec(48 99)
\lvec(48 100)
\lvec(47 100)
\ifill f:0
\move(49 99)
\lvec(50 99)
\lvec(50 100)
\lvec(49 100)
\ifill f:0
\move(51 99)
\lvec(52 99)
\lvec(52 100)
\lvec(51 100)
\ifill f:0
\move(56 99)
\lvec(57 99)
\lvec(57 100)
\lvec(56 100)
\ifill f:0
\move(60 99)
\lvec(62 99)
\lvec(62 100)
\lvec(60 100)
\ifill f:0
\move(63 99)
\lvec(65 99)
\lvec(65 100)
\lvec(63 100)
\ifill f:0
\move(66 99)
\lvec(67 99)
\lvec(67 100)
\lvec(66 100)
\ifill f:0
\move(72 99)
\lvec(74 99)
\lvec(74 100)
\lvec(72 100)
\ifill f:0
\move(75 99)
\lvec(76 99)
\lvec(76 100)
\lvec(75 100)
\ifill f:0
\move(77 99)
\lvec(79 99)
\lvec(79 100)
\lvec(77 100)
\ifill f:0
\move(80 99)
\lvec(82 99)
\lvec(82 100)
\lvec(80 100)
\ifill f:0
\move(84 99)
\lvec(85 99)
\lvec(85 100)
\lvec(84 100)
\ifill f:0
\move(86 99)
\lvec(87 99)
\lvec(87 100)
\lvec(86 100)
\ifill f:0
\move(89 99)
\lvec(91 99)
\lvec(91 100)
\lvec(89 100)
\ifill f:0
\move(95 99)
\lvec(98 99)
\lvec(98 100)
\lvec(95 100)
\ifill f:0
\move(100 99)
\lvec(101 99)
\lvec(101 100)
\lvec(100 100)
\ifill f:0
\move(102 99)
\lvec(106 99)
\lvec(106 100)
\lvec(102 100)
\ifill f:0
\move(107 99)
\lvec(111 99)
\lvec(111 100)
\lvec(107 100)
\ifill f:0
\move(112 99)
\lvec(119 99)
\lvec(119 100)
\lvec(112 100)
\ifill f:0
\move(120 99)
\lvec(122 99)
\lvec(122 100)
\lvec(120 100)
\ifill f:0
\move(124 99)
\lvec(125 99)
\lvec(125 100)
\lvec(124 100)
\ifill f:0
\move(126 99)
\lvec(127 99)
\lvec(127 100)
\lvec(126 100)
\ifill f:0
\move(128 99)
\lvec(130 99)
\lvec(130 100)
\lvec(128 100)
\ifill f:0
\move(131 99)
\lvec(132 99)
\lvec(132 100)
\lvec(131 100)
\ifill f:0
\move(134 99)
\lvec(136 99)
\lvec(136 100)
\lvec(134 100)
\ifill f:0
\move(137 99)
\lvec(141 99)
\lvec(141 100)
\lvec(137 100)
\ifill f:0
\move(142 99)
\lvec(145 99)
\lvec(145 100)
\lvec(142 100)
\ifill f:0
\move(146 99)
\lvec(147 99)
\lvec(147 100)
\lvec(146 100)
\ifill f:0
\move(149 99)
\lvec(161 99)
\lvec(161 100)
\lvec(149 100)
\ifill f:0
\move(164 99)
\lvec(170 99)
\lvec(170 100)
\lvec(164 100)
\ifill f:0
\move(171 99)
\lvec(174 99)
\lvec(174 100)
\lvec(171 100)
\ifill f:0
\move(176 99)
\lvec(185 99)
\lvec(185 100)
\lvec(176 100)
\ifill f:0
\move(186 99)
\lvec(190 99)
\lvec(190 100)
\lvec(186 100)
\ifill f:0
\move(191 99)
\lvec(195 99)
\lvec(195 100)
\lvec(191 100)
\ifill f:0
\move(196 99)
\lvec(197 99)
\lvec(197 100)
\lvec(196 100)
\ifill f:0
\move(198 99)
\lvec(211 99)
\lvec(211 100)
\lvec(198 100)
\ifill f:0
\move(212 99)
\lvec(219 99)
\lvec(219 100)
\lvec(212 100)
\ifill f:0
\move(220 99)
\lvec(222 99)
\lvec(222 100)
\lvec(220 100)
\ifill f:0
\move(223 99)
\lvec(226 99)
\lvec(226 100)
\lvec(223 100)
\ifill f:0
\move(227 99)
\lvec(228 99)
\lvec(228 100)
\lvec(227 100)
\ifill f:0
\move(229 99)
\lvec(235 99)
\lvec(235 100)
\lvec(229 100)
\ifill f:0
\move(236 99)
\lvec(240 99)
\lvec(240 100)
\lvec(236 100)
\ifill f:0
\move(241 99)
\lvec(250 99)
\lvec(250 100)
\lvec(241 100)
\ifill f:0
\move(253 99)
\lvec(257 99)
\lvec(257 100)
\lvec(253 100)
\ifill f:0
\move(258 99)
\lvec(261 99)
\lvec(261 100)
\lvec(258 100)
\ifill f:0
\move(262 99)
\lvec(290 99)
\lvec(290 100)
\lvec(262 100)
\ifill f:0
\move(292 99)
\lvec(299 99)
\lvec(299 100)
\lvec(292 100)
\ifill f:0
\move(300 99)
\lvec(307 99)
\lvec(307 100)
\lvec(300 100)
\ifill f:0
\move(308 99)
\lvec(312 99)
\lvec(312 100)
\lvec(308 100)
\ifill f:0
\move(313 99)
\lvec(322 99)
\lvec(322 100)
\lvec(313 100)
\ifill f:0
\move(323 99)
\lvec(325 99)
\lvec(325 100)
\lvec(323 100)
\ifill f:0
\move(327 99)
\lvec(330 99)
\lvec(330 100)
\lvec(327 100)
\ifill f:0
\move(331 99)
\lvec(344 99)
\lvec(344 100)
\lvec(331 100)
\ifill f:0
\move(345 99)
\lvec(347 99)
\lvec(347 100)
\lvec(345 100)
\ifill f:0
\move(348 99)
\lvec(356 99)
\lvec(356 100)
\lvec(348 100)
\ifill f:0
\move(357 99)
\lvec(362 99)
\lvec(362 100)
\lvec(357 100)
\ifill f:0
\move(363 99)
\lvec(364 99)
\lvec(364 100)
\lvec(363 100)
\ifill f:0
\move(365 99)
\lvec(371 99)
\lvec(371 100)
\lvec(365 100)
\ifill f:0
\move(372 99)
\lvec(378 99)
\lvec(378 100)
\lvec(372 100)
\ifill f:0
\move(379 99)
\lvec(387 99)
\lvec(387 100)
\lvec(379 100)
\ifill f:0
\move(388 99)
\lvec(389 99)
\lvec(389 100)
\lvec(388 100)
\ifill f:0
\move(390 99)
\lvec(391 99)
\lvec(391 100)
\lvec(390 100)
\ifill f:0
\move(392 99)
\lvec(393 99)
\lvec(393 100)
\lvec(392 100)
\ifill f:0
\move(394 99)
\lvec(395 99)
\lvec(395 100)
\lvec(394 100)
\ifill f:0
\move(396 99)
\lvec(397 99)
\lvec(397 100)
\lvec(396 100)
\ifill f:0
\move(398 99)
\lvec(399 99)
\lvec(399 100)
\lvec(398 100)
\ifill f:0
\move(400 99)
\lvec(401 99)
\lvec(401 100)
\lvec(400 100)
\ifill f:0
\move(402 99)
\lvec(403 99)
\lvec(403 100)
\lvec(402 100)
\ifill f:0
\move(404 99)
\lvec(423 99)
\lvec(423 100)
\lvec(404 100)
\ifill f:0
\move(424 99)
\lvec(433 99)
\lvec(433 100)
\lvec(424 100)
\ifill f:0
\move(434 99)
\lvec(442 99)
\lvec(442 100)
\lvec(434 100)
\ifill f:0
\move(443 99)
\lvec(451 99)
\lvec(451 100)
\lvec(443 100)
\ifill f:0
\move(16 100)
\lvec(17 100)
\lvec(17 101)
\lvec(16 101)
\ifill f:0
\move(25 100)
\lvec(26 100)
\lvec(26 101)
\lvec(25 101)
\ifill f:0
\move(36 100)
\lvec(37 100)
\lvec(37 101)
\lvec(36 101)
\ifill f:0
\move(38 100)
\lvec(42 100)
\lvec(42 101)
\lvec(38 101)
\ifill f:0
\move(43 100)
\lvec(45 100)
\lvec(45 101)
\lvec(43 101)
\ifill f:0
\move(46 100)
\lvec(48 100)
\lvec(48 101)
\lvec(46 101)
\ifill f:0
\move(49 100)
\lvec(50 100)
\lvec(50 101)
\lvec(49 101)
\ifill f:0
\move(51 100)
\lvec(53 100)
\lvec(53 101)
\lvec(51 101)
\ifill f:0
\move(54 100)
\lvec(55 100)
\lvec(55 101)
\lvec(54 101)
\ifill f:0
\move(57 100)
\lvec(58 100)
\lvec(58 101)
\lvec(57 101)
\ifill f:0
\move(59 100)
\lvec(60 100)
\lvec(60 101)
\lvec(59 101)
\ifill f:0
\move(61 100)
\lvec(62 100)
\lvec(62 101)
\lvec(61 101)
\ifill f:0
\move(63 100)
\lvec(65 100)
\lvec(65 101)
\lvec(63 101)
\ifill f:0
\move(67 100)
\lvec(75 100)
\lvec(75 101)
\lvec(67 101)
\ifill f:0
\move(76 100)
\lvec(79 100)
\lvec(79 101)
\lvec(76 101)
\ifill f:0
\move(80 100)
\lvec(82 100)
\lvec(82 101)
\lvec(80 101)
\ifill f:0
\move(83 100)
\lvec(84 100)
\lvec(84 101)
\lvec(83 101)
\ifill f:0
\move(85 100)
\lvec(86 100)
\lvec(86 101)
\lvec(85 101)
\ifill f:0
\move(88 100)
\lvec(90 100)
\lvec(90 101)
\lvec(88 101)
\ifill f:0
\move(91 100)
\lvec(93 100)
\lvec(93 101)
\lvec(91 101)
\ifill f:0
\move(95 100)
\lvec(98 100)
\lvec(98 101)
\lvec(95 101)
\ifill f:0
\move(100 100)
\lvec(101 100)
\lvec(101 101)
\lvec(100 101)
\ifill f:0
\move(102 100)
\lvec(109 100)
\lvec(109 101)
\lvec(102 101)
\ifill f:0
\move(110 100)
\lvec(116 100)
\lvec(116 101)
\lvec(110 101)
\ifill f:0
\move(117 100)
\lvec(122 100)
\lvec(122 101)
\lvec(117 101)
\ifill f:0
\move(123 100)
\lvec(124 100)
\lvec(124 101)
\lvec(123 101)
\ifill f:0
\move(125 100)
\lvec(126 100)
\lvec(126 101)
\lvec(125 101)
\ifill f:0
\move(127 100)
\lvec(129 100)
\lvec(129 101)
\lvec(127 101)
\ifill f:0
\move(130 100)
\lvec(131 100)
\lvec(131 101)
\lvec(130 101)
\ifill f:0
\move(132 100)
\lvec(134 100)
\lvec(134 101)
\lvec(132 101)
\ifill f:0
\move(135 100)
\lvec(145 100)
\lvec(145 101)
\lvec(135 101)
\ifill f:0
\move(146 100)
\lvec(147 100)
\lvec(147 101)
\lvec(146 101)
\ifill f:0
\move(148 100)
\lvec(159 100)
\lvec(159 101)
\lvec(148 101)
\ifill f:0
\move(160 100)
\lvec(170 100)
\lvec(170 101)
\lvec(160 101)
\ifill f:0
\move(171 100)
\lvec(175 100)
\lvec(175 101)
\lvec(171 101)
\ifill f:0
\move(177 100)
\lvec(180 100)
\lvec(180 101)
\lvec(177 101)
\ifill f:0
\move(181 100)
\lvec(187 100)
\lvec(187 101)
\lvec(181 101)
\ifill f:0
\move(188 100)
\lvec(195 100)
\lvec(195 101)
\lvec(188 101)
\ifill f:0
\move(196 100)
\lvec(197 100)
\lvec(197 101)
\lvec(196 101)
\ifill f:0
\move(198 100)
\lvec(217 100)
\lvec(217 101)
\lvec(198 101)
\ifill f:0
\move(218 100)
\lvec(226 100)
\lvec(226 101)
\lvec(218 101)
\ifill f:0
\move(227 100)
\lvec(228 100)
\lvec(228 101)
\lvec(227 101)
\ifill f:0
\move(229 100)
\lvec(231 100)
\lvec(231 101)
\lvec(229 101)
\ifill f:0
\move(232 100)
\lvec(234 100)
\lvec(234 101)
\lvec(232 101)
\ifill f:0
\move(235 100)
\lvec(238 100)
\lvec(238 101)
\lvec(235 101)
\ifill f:0
\move(239 100)
\lvec(242 100)
\lvec(242 101)
\lvec(239 101)
\ifill f:0
\move(243 100)
\lvec(247 100)
\lvec(247 101)
\lvec(243 101)
\ifill f:0
\move(248 100)
\lvec(257 100)
\lvec(257 101)
\lvec(248 101)
\ifill f:0
\move(258 100)
\lvec(260 100)
\lvec(260 101)
\lvec(258 101)
\ifill f:0
\move(261 100)
\lvec(271 100)
\lvec(271 101)
\lvec(261 101)
\ifill f:0
\move(273 100)
\lvec(290 100)
\lvec(290 101)
\lvec(273 101)
\ifill f:0
\move(294 100)
\lvec(302 100)
\lvec(302 101)
\lvec(294 101)
\ifill f:0
\move(303 100)
\lvec(310 100)
\lvec(310 101)
\lvec(303 101)
\ifill f:0
\move(312 100)
\lvec(322 100)
\lvec(322 101)
\lvec(312 101)
\ifill f:0
\move(323 100)
\lvec(325 100)
\lvec(325 101)
\lvec(323 101)
\ifill f:0
\move(326 100)
\lvec(335 100)
\lvec(335 101)
\lvec(326 101)
\ifill f:0
\move(336 100)
\lvec(349 100)
\lvec(349 101)
\lvec(336 101)
\ifill f:0
\move(350 100)
\lvec(362 100)
\lvec(362 101)
\lvec(350 101)
\ifill f:0
\move(363 100)
\lvec(364 100)
\lvec(364 101)
\lvec(363 101)
\ifill f:0
\move(365 100)
\lvec(384 100)
\lvec(384 101)
\lvec(365 101)
\ifill f:0
\move(385 100)
\lvec(386 100)
\lvec(386 101)
\lvec(385 101)
\ifill f:0
\move(387 100)
\lvec(395 100)
\lvec(395 101)
\lvec(387 101)
\ifill f:0
\move(396 100)
\lvec(397 100)
\lvec(397 101)
\lvec(396 101)
\ifill f:0
\move(398 100)
\lvec(399 100)
\lvec(399 101)
\lvec(398 101)
\ifill f:0
\move(400 100)
\lvec(401 100)
\lvec(401 101)
\lvec(400 101)
\ifill f:0
\move(402 100)
\lvec(403 100)
\lvec(403 101)
\lvec(402 101)
\ifill f:0
\move(404 100)
\lvec(405 100)
\lvec(405 101)
\lvec(404 101)
\ifill f:0
\move(406 100)
\lvec(407 100)
\lvec(407 101)
\lvec(406 101)
\ifill f:0
\move(408 100)
\lvec(409 100)
\lvec(409 101)
\lvec(408 101)
\ifill f:0
\move(410 100)
\lvec(436 100)
\lvec(436 101)
\lvec(410 101)
\ifill f:0
\move(437 100)
\lvec(442 100)
\lvec(442 101)
\lvec(437 101)
\ifill f:0
\move(443 100)
\lvec(446 100)
\lvec(446 101)
\lvec(443 101)
\ifill f:0
\move(447 100)
\lvec(451 100)
\lvec(451 101)
\lvec(447 101)
\ifill f:0
\move(16 101)
\lvec(17 101)
\lvec(17 102)
\lvec(16 102)
\ifill f:0
\move(18 101)
\lvec(19 101)
\lvec(19 102)
\lvec(18 102)
\ifill f:0
\move(20 101)
\lvec(21 101)
\lvec(21 102)
\lvec(20 102)
\ifill f:0
\move(23 101)
\lvec(24 101)
\lvec(24 102)
\lvec(23 102)
\ifill f:0
\move(25 101)
\lvec(26 101)
\lvec(26 102)
\lvec(25 102)
\ifill f:0
\move(36 101)
\lvec(37 101)
\lvec(37 102)
\lvec(36 102)
\ifill f:0
\move(38 101)
\lvec(39 101)
\lvec(39 102)
\lvec(38 102)
\ifill f:0
\move(40 101)
\lvec(45 101)
\lvec(45 102)
\lvec(40 102)
\ifill f:0
\move(47 101)
\lvec(50 101)
\lvec(50 102)
\lvec(47 102)
\ifill f:0
\move(51 101)
\lvec(52 101)
\lvec(52 102)
\lvec(51 102)
\ifill f:0
\move(54 101)
\lvec(55 101)
\lvec(55 102)
\lvec(54 102)
\ifill f:0
\move(56 101)
\lvec(57 101)
\lvec(57 102)
\lvec(56 102)
\ifill f:0
\move(59 101)
\lvec(60 101)
\lvec(60 102)
\lvec(59 102)
\ifill f:0
\move(61 101)
\lvec(65 101)
\lvec(65 102)
\lvec(61 102)
\ifill f:0
\move(66 101)
\lvec(71 101)
\lvec(71 102)
\lvec(66 102)
\ifill f:0
\move(72 101)
\lvec(73 101)
\lvec(73 102)
\lvec(72 102)
\ifill f:0
\move(75 101)
\lvec(78 101)
\lvec(78 102)
\lvec(75 102)
\ifill f:0
\move(79 101)
\lvec(82 101)
\lvec(82 102)
\lvec(79 102)
\ifill f:0
\move(83 101)
\lvec(84 101)
\lvec(84 102)
\lvec(83 102)
\ifill f:0
\move(87 101)
\lvec(91 101)
\lvec(91 102)
\lvec(87 102)
\ifill f:0
\move(92 101)
\lvec(93 101)
\lvec(93 102)
\lvec(92 102)
\ifill f:0
\move(97 101)
\lvec(101 101)
\lvec(101 102)
\lvec(97 102)
\ifill f:0
\move(102 101)
\lvec(111 101)
\lvec(111 102)
\lvec(102 102)
\ifill f:0
\move(112 101)
\lvec(115 101)
\lvec(115 102)
\lvec(112 102)
\ifill f:0
\move(116 101)
\lvec(122 101)
\lvec(122 102)
\lvec(116 102)
\ifill f:0
\move(123 101)
\lvec(124 101)
\lvec(124 102)
\lvec(123 102)
\ifill f:0
\move(125 101)
\lvec(126 101)
\lvec(126 102)
\lvec(125 102)
\ifill f:0
\move(127 101)
\lvec(128 101)
\lvec(128 102)
\lvec(127 102)
\ifill f:0
\move(129 101)
\lvec(130 101)
\lvec(130 102)
\lvec(129 102)
\ifill f:0
\move(131 101)
\lvec(132 101)
\lvec(132 102)
\lvec(131 102)
\ifill f:0
\move(133 101)
\lvec(135 101)
\lvec(135 102)
\lvec(133 102)
\ifill f:0
\move(136 101)
\lvec(138 101)
\lvec(138 102)
\lvec(136 102)
\ifill f:0
\move(139 101)
\lvec(142 101)
\lvec(142 102)
\lvec(139 102)
\ifill f:0
\move(143 101)
\lvec(145 101)
\lvec(145 102)
\lvec(143 102)
\ifill f:0
\move(146 101)
\lvec(153 101)
\lvec(153 102)
\lvec(146 102)
\ifill f:0
\move(154 101)
\lvec(170 101)
\lvec(170 102)
\lvec(154 102)
\ifill f:0
\move(172 101)
\lvec(177 101)
\lvec(177 102)
\lvec(172 102)
\ifill f:0
\move(178 101)
\lvec(186 101)
\lvec(186 102)
\lvec(178 102)
\ifill f:0
\move(187 101)
\lvec(189 101)
\lvec(189 102)
\lvec(187 102)
\ifill f:0
\move(190 101)
\lvec(192 101)
\lvec(192 102)
\lvec(190 102)
\ifill f:0
\move(193 101)
\lvec(195 101)
\lvec(195 102)
\lvec(193 102)
\ifill f:0
\move(196 101)
\lvec(197 101)
\lvec(197 102)
\lvec(196 102)
\ifill f:0
\move(198 101)
\lvec(211 101)
\lvec(211 102)
\lvec(198 102)
\ifill f:0
\move(212 101)
\lvec(220 101)
\lvec(220 102)
\lvec(212 102)
\ifill f:0
\move(221 101)
\lvec(226 101)
\lvec(226 102)
\lvec(221 102)
\ifill f:0
\move(227 101)
\lvec(233 101)
\lvec(233 102)
\lvec(227 102)
\ifill f:0
\move(234 101)
\lvec(240 101)
\lvec(240 102)
\lvec(234 102)
\ifill f:0
\move(241 101)
\lvec(244 101)
\lvec(244 102)
\lvec(241 102)
\ifill f:0
\move(245 101)
\lvec(248 101)
\lvec(248 102)
\lvec(245 102)
\ifill f:0
\move(249 101)
\lvec(253 101)
\lvec(253 102)
\lvec(249 102)
\ifill f:0
\move(254 101)
\lvec(257 101)
\lvec(257 102)
\lvec(254 102)
\ifill f:0
\move(258 101)
\lvec(259 101)
\lvec(259 102)
\lvec(258 102)
\ifill f:0
\move(260 101)
\lvec(266 101)
\lvec(266 102)
\lvec(260 102)
\ifill f:0
\move(267 101)
\lvec(279 101)
\lvec(279 102)
\lvec(267 102)
\ifill f:0
\move(286 101)
\lvec(287 101)
\lvec(287 102)
\lvec(286 102)
\ifill f:0
\move(289 101)
\lvec(290 101)
\lvec(290 102)
\lvec(289 102)
\ifill f:0
\move(294 101)
\lvec(295 101)
\lvec(295 102)
\lvec(294 102)
\ifill f:0
\move(296 101)
\lvec(298 101)
\lvec(298 102)
\lvec(296 102)
\ifill f:0
\move(299 101)
\lvec(307 101)
\lvec(307 102)
\lvec(299 102)
\ifill f:0
\move(308 101)
\lvec(315 101)
\lvec(315 102)
\lvec(308 102)
\ifill f:0
\move(316 101)
\lvec(322 101)
\lvec(322 102)
\lvec(316 102)
\ifill f:0
\move(323 101)
\lvec(325 101)
\lvec(325 102)
\lvec(323 102)
\ifill f:0
\move(326 101)
\lvec(358 101)
\lvec(358 102)
\lvec(326 102)
\ifill f:0
\move(359 101)
\lvec(362 101)
\lvec(362 102)
\lvec(359 102)
\ifill f:0
\move(363 101)
\lvec(364 101)
\lvec(364 102)
\lvec(363 102)
\ifill f:0
\move(365 101)
\lvec(367 101)
\lvec(367 102)
\lvec(365 102)
\ifill f:0
\move(368 101)
\lvec(370 101)
\lvec(370 102)
\lvec(368 102)
\ifill f:0
\move(371 101)
\lvec(383 101)
\lvec(383 102)
\lvec(371 102)
\ifill f:0
\move(384 101)
\lvec(399 101)
\lvec(399 102)
\lvec(384 102)
\ifill f:0
\move(400 101)
\lvec(401 101)
\lvec(401 102)
\lvec(400 102)
\ifill f:0
\move(402 101)
\lvec(425 101)
\lvec(425 102)
\lvec(402 102)
\ifill f:0
\move(426 101)
\lvec(434 101)
\lvec(434 102)
\lvec(426 102)
\ifill f:0
\move(435 101)
\lvec(436 101)
\lvec(436 102)
\lvec(435 102)
\ifill f:0
\move(437 101)
\lvec(442 101)
\lvec(442 102)
\lvec(437 102)
\ifill f:0
\move(444 101)
\lvec(451 101)
\lvec(451 102)
\lvec(444 102)
\ifill f:0
\move(16 102)
\lvec(17 102)
\lvec(17 103)
\lvec(16 103)
\ifill f:0
\move(19 102)
\lvec(20 102)
\lvec(20 103)
\lvec(19 103)
\ifill f:0
\move(21 102)
\lvec(22 102)
\lvec(22 103)
\lvec(21 103)
\ifill f:0
\move(23 102)
\lvec(24 102)
\lvec(24 103)
\lvec(23 103)
\ifill f:0
\move(25 102)
\lvec(26 102)
\lvec(26 103)
\lvec(25 103)
\ifill f:0
\move(36 102)
\lvec(37 102)
\lvec(37 103)
\lvec(36 103)
\ifill f:0
\move(38 102)
\lvec(39 102)
\lvec(39 103)
\lvec(38 103)
\ifill f:0
\move(40 102)
\lvec(41 102)
\lvec(41 103)
\lvec(40 103)
\ifill f:0
\move(43 102)
\lvec(46 102)
\lvec(46 103)
\lvec(43 103)
\ifill f:0
\move(47 102)
\lvec(50 102)
\lvec(50 103)
\lvec(47 103)
\ifill f:0
\move(57 102)
\lvec(59 102)
\lvec(59 103)
\lvec(57 103)
\ifill f:0
\move(60 102)
\lvec(65 102)
\lvec(65 103)
\lvec(60 103)
\ifill f:0
\move(66 102)
\lvec(78 102)
\lvec(78 103)
\lvec(66 103)
\ifill f:0
\move(79 102)
\lvec(82 102)
\lvec(82 103)
\lvec(79 103)
\ifill f:0
\move(83 102)
\lvec(84 102)
\lvec(84 103)
\lvec(83 103)
\ifill f:0
\move(86 102)
\lvec(87 102)
\lvec(87 103)
\lvec(86 103)
\ifill f:0
\move(88 102)
\lvec(90 102)
\lvec(90 103)
\lvec(88 103)
\ifill f:0
\move(91 102)
\lvec(93 102)
\lvec(93 103)
\lvec(91 103)
\ifill f:0
\move(98 102)
\lvec(101 102)
\lvec(101 103)
\lvec(98 103)
\ifill f:0
\move(102 102)
\lvec(106 102)
\lvec(106 103)
\lvec(102 103)
\ifill f:0
\move(107 102)
\lvec(118 102)
\lvec(118 103)
\lvec(107 103)
\ifill f:0
\move(119 102)
\lvec(122 102)
\lvec(122 103)
\lvec(119 103)
\ifill f:0
\move(123 102)
\lvec(124 102)
\lvec(124 103)
\lvec(123 103)
\ifill f:0
\move(125 102)
\lvec(127 102)
\lvec(127 103)
\lvec(125 103)
\ifill f:0
\move(128 102)
\lvec(131 102)
\lvec(131 103)
\lvec(128 103)
\ifill f:0
\move(132 102)
\lvec(136 102)
\lvec(136 103)
\lvec(132 103)
\ifill f:0
\move(137 102)
\lvec(138 102)
\lvec(138 103)
\lvec(137 103)
\ifill f:0
\move(140 102)
\lvec(142 102)
\lvec(142 103)
\lvec(140 103)
\ifill f:0
\move(143 102)
\lvec(145 102)
\lvec(145 103)
\lvec(143 103)
\ifill f:0
\move(147 102)
\lvec(150 102)
\lvec(150 103)
\lvec(147 103)
\ifill f:0
\move(151 102)
\lvec(159 102)
\lvec(159 103)
\lvec(151 103)
\ifill f:0
\move(160 102)
\lvec(170 102)
\lvec(170 103)
\lvec(160 103)
\ifill f:0
\move(175 102)
\lvec(179 102)
\lvec(179 103)
\lvec(175 103)
\ifill f:0
\move(180 102)
\lvec(184 102)
\lvec(184 103)
\lvec(180 103)
\ifill f:0
\move(185 102)
\lvec(192 102)
\lvec(192 103)
\lvec(185 103)
\ifill f:0
\move(193 102)
\lvec(195 102)
\lvec(195 103)
\lvec(193 103)
\ifill f:0
\move(196 102)
\lvec(197 102)
\lvec(197 103)
\lvec(196 103)
\ifill f:0
\move(198 102)
\lvec(203 102)
\lvec(203 103)
\lvec(198 103)
\ifill f:0
\move(204 102)
\lvec(214 102)
\lvec(214 103)
\lvec(204 103)
\ifill f:0
\move(215 102)
\lvec(226 102)
\lvec(226 103)
\lvec(215 103)
\ifill f:0
\move(228 102)
\lvec(230 102)
\lvec(230 103)
\lvec(228 103)
\ifill f:0
\move(231 102)
\lvec(235 102)
\lvec(235 103)
\lvec(231 103)
\ifill f:0
\move(236 102)
\lvec(249 102)
\lvec(249 103)
\lvec(236 103)
\ifill f:0
\move(250 102)
\lvec(257 102)
\lvec(257 103)
\lvec(250 103)
\ifill f:0
\move(258 102)
\lvec(271 102)
\lvec(271 103)
\lvec(258 103)
\ifill f:0
\move(272 102)
\lvec(283 102)
\lvec(283 103)
\lvec(272 103)
\ifill f:0
\move(285 102)
\lvec(287 102)
\lvec(287 103)
\lvec(285 103)
\ifill f:0
\move(289 102)
\lvec(290 102)
\lvec(290 103)
\lvec(289 103)
\ifill f:0
\move(292 102)
\lvec(293 102)
\lvec(293 103)
\lvec(292 103)
\ifill f:0
\move(300 102)
\lvec(312 102)
\lvec(312 103)
\lvec(300 103)
\ifill f:0
\move(314 102)
\lvec(321 102)
\lvec(321 103)
\lvec(314 103)
\ifill f:0
\move(323 102)
\lvec(325 102)
\lvec(325 103)
\lvec(323 103)
\ifill f:0
\move(326 102)
\lvec(338 102)
\lvec(338 103)
\lvec(326 103)
\ifill f:0
\move(339 102)
\lvec(362 102)
\lvec(362 103)
\lvec(339 103)
\ifill f:0
\move(363 102)
\lvec(379 102)
\lvec(379 103)
\lvec(363 103)
\ifill f:0
\move(380 102)
\lvec(387 102)
\lvec(387 103)
\lvec(380 103)
\ifill f:0
\move(388 102)
\lvec(392 102)
\lvec(392 103)
\lvec(388 103)
\ifill f:0
\move(393 102)
\lvec(399 102)
\lvec(399 103)
\lvec(393 103)
\ifill f:0
\move(400 102)
\lvec(401 102)
\lvec(401 103)
\lvec(400 103)
\ifill f:0
\move(402 102)
\lvec(410 102)
\lvec(410 103)
\lvec(402 103)
\ifill f:0
\move(411 102)
\lvec(412 102)
\lvec(412 103)
\lvec(411 103)
\ifill f:0
\move(413 102)
\lvec(430 102)
\lvec(430 103)
\lvec(413 103)
\ifill f:0
\move(431 102)
\lvec(432 102)
\lvec(432 103)
\lvec(431 103)
\ifill f:0
\move(433 102)
\lvec(442 102)
\lvec(442 103)
\lvec(433 103)
\ifill f:0
\move(444 102)
\lvec(448 102)
\lvec(448 103)
\lvec(444 103)
\ifill f:0
\move(449 102)
\lvec(450 102)
\lvec(450 103)
\lvec(449 103)
\ifill f:0
\move(15 103)
\lvec(17 103)
\lvec(17 104)
\lvec(15 104)
\ifill f:0
\move(20 103)
\lvec(21 103)
\lvec(21 104)
\lvec(20 104)
\ifill f:0
\move(24 103)
\lvec(26 103)
\lvec(26 104)
\lvec(24 104)
\ifill f:0
\move(36 103)
\lvec(37 103)
\lvec(37 104)
\lvec(36 104)
\ifill f:0
\move(38 103)
\lvec(39 103)
\lvec(39 104)
\lvec(38 104)
\ifill f:0
\move(40 103)
\lvec(45 103)
\lvec(45 104)
\lvec(40 104)
\ifill f:0
\move(47 103)
\lvec(50 103)
\lvec(50 104)
\lvec(47 104)
\ifill f:0
\move(54 103)
\lvec(57 103)
\lvec(57 104)
\lvec(54 104)
\ifill f:0
\move(59 103)
\lvec(61 103)
\lvec(61 104)
\lvec(59 104)
\ifill f:0
\move(62 103)
\lvec(63 103)
\lvec(63 104)
\lvec(62 104)
\ifill f:0
\move(64 103)
\lvec(65 103)
\lvec(65 104)
\lvec(64 104)
\ifill f:0
\move(66 103)
\lvec(68 103)
\lvec(68 104)
\lvec(66 104)
\ifill f:0
\move(69 103)
\lvec(71 103)
\lvec(71 104)
\lvec(69 104)
\ifill f:0
\move(72 103)
\lvec(73 103)
\lvec(73 104)
\lvec(72 104)
\ifill f:0
\move(76 103)
\lvec(77 103)
\lvec(77 104)
\lvec(76 104)
\ifill f:0
\move(78 103)
\lvec(82 103)
\lvec(82 104)
\lvec(78 104)
\ifill f:0
\move(83 103)
\lvec(85 103)
\lvec(85 104)
\lvec(83 104)
\ifill f:0
\move(86 103)
\lvec(87 103)
\lvec(87 104)
\lvec(86 104)
\ifill f:0
\move(88 103)
\lvec(89 103)
\lvec(89 104)
\lvec(88 104)
\ifill f:0
\move(90 103)
\lvec(91 103)
\lvec(91 104)
\lvec(90 104)
\ifill f:0
\move(92 103)
\lvec(93 103)
\lvec(93 104)
\lvec(92 104)
\ifill f:0
\move(95 103)
\lvec(96 103)
\lvec(96 104)
\lvec(95 104)
\ifill f:0
\move(97 103)
\lvec(101 103)
\lvec(101 104)
\lvec(97 104)
\ifill f:0
\move(102 103)
\lvec(103 103)
\lvec(103 104)
\lvec(102 104)
\ifill f:0
\move(105 103)
\lvec(111 103)
\lvec(111 104)
\lvec(105 104)
\ifill f:0
\move(112 103)
\lvec(122 103)
\lvec(122 104)
\lvec(112 104)
\ifill f:0
\move(123 103)
\lvec(125 103)
\lvec(125 104)
\lvec(123 104)
\ifill f:0
\move(126 103)
\lvec(127 103)
\lvec(127 104)
\lvec(126 104)
\ifill f:0
\move(128 103)
\lvec(130 103)
\lvec(130 104)
\lvec(128 104)
\ifill f:0
\move(131 103)
\lvec(132 103)
\lvec(132 104)
\lvec(131 104)
\ifill f:0
\move(133 103)
\lvec(134 103)
\lvec(134 104)
\lvec(133 104)
\ifill f:0
\move(135 103)
\lvec(137 103)
\lvec(137 104)
\lvec(135 104)
\ifill f:0
\move(138 103)
\lvec(145 103)
\lvec(145 104)
\lvec(138 104)
\ifill f:0
\move(146 103)
\lvec(149 103)
\lvec(149 104)
\lvec(146 104)
\ifill f:0
\move(150 103)
\lvec(154 103)
\lvec(154 104)
\lvec(150 104)
\ifill f:0
\move(155 103)
\lvec(163 103)
\lvec(163 104)
\lvec(155 104)
\ifill f:0
\move(165 103)
\lvec(166 103)
\lvec(166 104)
\lvec(165 104)
\ifill f:0
\move(169 103)
\lvec(170 103)
\lvec(170 104)
\lvec(169 104)
\ifill f:0
\move(175 103)
\lvec(176 103)
\lvec(176 104)
\lvec(175 104)
\ifill f:0
\move(177 103)
\lvec(187 103)
\lvec(187 104)
\lvec(177 104)
\ifill f:0
\move(188 103)
\lvec(191 103)
\lvec(191 104)
\lvec(188 104)
\ifill f:0
\move(192 103)
\lvec(195 103)
\lvec(195 104)
\lvec(192 104)
\ifill f:0
\move(196 103)
\lvec(197 103)
\lvec(197 104)
\lvec(196 104)
\ifill f:0
\move(198 103)
\lvec(201 103)
\lvec(201 104)
\lvec(198 104)
\ifill f:0
\move(202 103)
\lvec(206 103)
\lvec(206 104)
\lvec(202 104)
\ifill f:0
\move(207 103)
\lvec(211 103)
\lvec(211 104)
\lvec(207 104)
\ifill f:0
\move(212 103)
\lvec(226 103)
\lvec(226 104)
\lvec(212 104)
\ifill f:0
\move(228 103)
\lvec(232 103)
\lvec(232 104)
\lvec(228 104)
\ifill f:0
\move(233 103)
\lvec(234 103)
\lvec(234 104)
\lvec(233 104)
\ifill f:0
\move(235 103)
\lvec(240 103)
\lvec(240 104)
\lvec(235 104)
\ifill f:0
\move(241 103)
\lvec(243 103)
\lvec(243 104)
\lvec(241 104)
\ifill f:0
\move(244 103)
\lvec(250 103)
\lvec(250 104)
\lvec(244 104)
\ifill f:0
\move(251 103)
\lvec(254 103)
\lvec(254 104)
\lvec(251 104)
\ifill f:0
\move(255 103)
\lvec(257 103)
\lvec(257 104)
\lvec(255 104)
\ifill f:0
\move(258 103)
\lvec(263 103)
\lvec(263 104)
\lvec(258 104)
\ifill f:0
\move(264 103)
\lvec(287 103)
\lvec(287 104)
\lvec(264 104)
\ifill f:0
\move(289 103)
\lvec(290 103)
\lvec(290 104)
\lvec(289 104)
\ifill f:0
\move(291 103)
\lvec(306 103)
\lvec(306 104)
\lvec(291 104)
\ifill f:0
\move(309 103)
\lvec(310 103)
\lvec(310 104)
\lvec(309 104)
\ifill f:0
\move(311 103)
\lvec(320 103)
\lvec(320 104)
\lvec(311 104)
\ifill f:0
\move(323 103)
\lvec(325 103)
\lvec(325 104)
\lvec(323 104)
\ifill f:0
\move(326 103)
\lvec(328 103)
\lvec(328 104)
\lvec(326 104)
\ifill f:0
\move(329 103)
\lvec(362 103)
\lvec(362 104)
\lvec(329 104)
\ifill f:0
\move(363 103)
\lvec(365 103)
\lvec(365 104)
\lvec(363 104)
\ifill f:0
\move(366 103)
\lvec(368 103)
\lvec(368 104)
\lvec(366 104)
\ifill f:0
\move(369 103)
\lvec(386 103)
\lvec(386 104)
\lvec(369 104)
\ifill f:0
\move(387 103)
\lvec(394 103)
\lvec(394 104)
\lvec(387 104)
\ifill f:0
\move(395 103)
\lvec(399 103)
\lvec(399 104)
\lvec(395 104)
\ifill f:0
\move(400 103)
\lvec(401 103)
\lvec(401 104)
\lvec(400 104)
\ifill f:0
\move(402 103)
\lvec(415 103)
\lvec(415 104)
\lvec(402 104)
\ifill f:0
\move(416 103)
\lvec(442 103)
\lvec(442 104)
\lvec(416 104)
\ifill f:0
\move(444 103)
\lvec(445 103)
\lvec(445 104)
\lvec(444 104)
\ifill f:0
\move(446 103)
\lvec(451 103)
\lvec(451 104)
\lvec(446 104)
\ifill f:0
\move(15 104)
\lvec(17 104)
\lvec(17 105)
\lvec(15 105)
\ifill f:0
\move(19 104)
\lvec(20 104)
\lvec(20 105)
\lvec(19 105)
\ifill f:0
\move(24 104)
\lvec(26 104)
\lvec(26 105)
\lvec(24 105)
\ifill f:0
\move(36 104)
\lvec(37 104)
\lvec(37 105)
\lvec(36 105)
\ifill f:0
\move(38 104)
\lvec(41 104)
\lvec(41 105)
\lvec(38 105)
\ifill f:0
\move(42 104)
\lvec(43 104)
\lvec(43 105)
\lvec(42 105)
\ifill f:0
\move(44 104)
\lvec(45 104)
\lvec(45 105)
\lvec(44 105)
\ifill f:0
\move(48 104)
\lvec(50 104)
\lvec(50 105)
\lvec(48 105)
\ifill f:0
\move(52 104)
\lvec(53 104)
\lvec(53 105)
\lvec(52 105)
\ifill f:0
\move(54 104)
\lvec(55 104)
\lvec(55 105)
\lvec(54 105)
\ifill f:0
\move(57 104)
\lvec(58 104)
\lvec(58 105)
\lvec(57 105)
\ifill f:0
\move(59 104)
\lvec(60 104)
\lvec(60 105)
\lvec(59 105)
\ifill f:0
\move(61 104)
\lvec(63 104)
\lvec(63 105)
\lvec(61 105)
\ifill f:0
\move(64 104)
\lvec(65 104)
\lvec(65 105)
\lvec(64 105)
\ifill f:0
\move(66 104)
\lvec(67 104)
\lvec(67 105)
\lvec(66 105)
\ifill f:0
\move(68 104)
\lvec(74 104)
\lvec(74 105)
\lvec(68 105)
\ifill f:0
\move(75 104)
\lvec(82 104)
\lvec(82 105)
\lvec(75 105)
\ifill f:0
\move(83 104)
\lvec(85 104)
\lvec(85 105)
\lvec(83 105)
\ifill f:0
\move(87 104)
\lvec(88 104)
\lvec(88 105)
\lvec(87 105)
\ifill f:0
\move(89 104)
\lvec(90 104)
\lvec(90 105)
\lvec(89 105)
\ifill f:0
\move(91 104)
\lvec(92 104)
\lvec(92 105)
\lvec(91 105)
\ifill f:0
\move(95 104)
\lvec(98 104)
\lvec(98 105)
\lvec(95 105)
\ifill f:0
\move(99 104)
\lvec(101 104)
\lvec(101 105)
\lvec(99 105)
\ifill f:0
\move(102 104)
\lvec(116 104)
\lvec(116 105)
\lvec(102 105)
\ifill f:0
\move(118 104)
\lvec(122 104)
\lvec(122 105)
\lvec(118 105)
\ifill f:0
\move(123 104)
\lvec(125 104)
\lvec(125 105)
\lvec(123 105)
\ifill f:0
\move(126 104)
\lvec(128 104)
\lvec(128 105)
\lvec(126 105)
\ifill f:0
\move(129 104)
\lvec(130 104)
\lvec(130 105)
\lvec(129 105)
\ifill f:0
\move(132 104)
\lvec(133 104)
\lvec(133 105)
\lvec(132 105)
\ifill f:0
\move(134 104)
\lvec(135 104)
\lvec(135 105)
\lvec(134 105)
\ifill f:0
\move(136 104)
\lvec(140 104)
\lvec(140 105)
\lvec(136 105)
\ifill f:0
\move(141 104)
\lvec(145 104)
\lvec(145 105)
\lvec(141 105)
\ifill f:0
\move(146 104)
\lvec(163 104)
\lvec(163 105)
\lvec(146 105)
\ifill f:0
\move(164 104)
\lvec(167 104)
\lvec(167 105)
\lvec(164 105)
\ifill f:0
\move(169 104)
\lvec(170 104)
\lvec(170 105)
\lvec(169 105)
\ifill f:0
\move(177 104)
\lvec(184 104)
\lvec(184 105)
\lvec(177 105)
\ifill f:0
\move(187 104)
\lvec(190 104)
\lvec(190 105)
\lvec(187 105)
\ifill f:0
\move(192 104)
\lvec(195 104)
\lvec(195 105)
\lvec(192 105)
\ifill f:0
\move(196 104)
\lvec(197 104)
\lvec(197 105)
\lvec(196 105)
\ifill f:0
\move(198 104)
\lvec(207 104)
\lvec(207 105)
\lvec(198 105)
\ifill f:0
\move(208 104)
\lvec(219 104)
\lvec(219 105)
\lvec(208 105)
\ifill f:0
\move(220 104)
\lvec(221 104)
\lvec(221 105)
\lvec(220 105)
\ifill f:0
\move(222 104)
\lvec(226 104)
\lvec(226 105)
\lvec(222 105)
\ifill f:0
\move(228 104)
\lvec(229 104)
\lvec(229 105)
\lvec(228 105)
\ifill f:0
\move(230 104)
\lvec(250 104)
\lvec(250 105)
\lvec(230 105)
\ifill f:0
\move(251 104)
\lvec(254 104)
\lvec(254 105)
\lvec(251 105)
\ifill f:0
\move(255 104)
\lvec(257 104)
\lvec(257 105)
\lvec(255 105)
\ifill f:0
\move(259 104)
\lvec(262 104)
\lvec(262 105)
\lvec(259 105)
\ifill f:0
\move(263 104)
\lvec(266 104)
\lvec(266 105)
\lvec(263 105)
\ifill f:0
\move(267 104)
\lvec(279 104)
\lvec(279 105)
\lvec(267 105)
\ifill f:0
\move(280 104)
\lvec(288 104)
\lvec(288 105)
\lvec(280 105)
\ifill f:0
\move(289 104)
\lvec(290 104)
\lvec(290 105)
\lvec(289 105)
\ifill f:0
\move(291 104)
\lvec(318 104)
\lvec(318 105)
\lvec(291 105)
\ifill f:0
\move(321 104)
\lvec(325 104)
\lvec(325 105)
\lvec(321 105)
\ifill f:0
\move(326 104)
\lvec(329 104)
\lvec(329 105)
\lvec(326 105)
\ifill f:0
\move(330 104)
\lvec(336 104)
\lvec(336 105)
\lvec(330 105)
\ifill f:0
\move(337 104)
\lvec(362 104)
\lvec(362 105)
\lvec(337 105)
\ifill f:0
\move(363 104)
\lvec(365 104)
\lvec(365 105)
\lvec(363 105)
\ifill f:0
\move(366 104)
\lvec(379 104)
\lvec(379 105)
\lvec(366 105)
\ifill f:0
\move(380 104)
\lvec(382 104)
\lvec(382 105)
\lvec(380 105)
\ifill f:0
\move(383 104)
\lvec(399 104)
\lvec(399 105)
\lvec(383 105)
\ifill f:0
\move(400 104)
\lvec(401 104)
\lvec(401 105)
\lvec(400 105)
\ifill f:0
\move(402 104)
\lvec(418 104)
\lvec(418 105)
\lvec(402 105)
\ifill f:0
\move(419 104)
\lvec(431 104)
\lvec(431 105)
\lvec(419 105)
\ifill f:0
\move(432 104)
\lvec(433 104)
\lvec(433 105)
\lvec(432 105)
\ifill f:0
\move(434 104)
\lvec(435 104)
\lvec(435 105)
\lvec(434 105)
\ifill f:0
\move(436 104)
\lvec(439 104)
\lvec(439 105)
\lvec(436 105)
\ifill f:0
\move(440 104)
\lvec(442 104)
\lvec(442 105)
\lvec(440 105)
\ifill f:0
\move(444 104)
\lvec(445 104)
\lvec(445 105)
\lvec(444 105)
\ifill f:0
\move(446 104)
\lvec(451 104)
\lvec(451 105)
\lvec(446 105)
\ifill f:0
\move(15 105)
\lvec(17 105)
\lvec(17 106)
\lvec(15 106)
\ifill f:0
\move(18 105)
\lvec(19 105)
\lvec(19 106)
\lvec(18 106)
\ifill f:0
\move(20 105)
\lvec(21 105)
\lvec(21 106)
\lvec(20 106)
\ifill f:0
\move(24 105)
\lvec(26 105)
\lvec(26 106)
\lvec(24 106)
\ifill f:0
\move(36 105)
\lvec(37 105)
\lvec(37 106)
\lvec(36 106)
\ifill f:0
\move(40 105)
\lvec(47 105)
\lvec(47 106)
\lvec(40 106)
\ifill f:0
\move(48 105)
\lvec(50 105)
\lvec(50 106)
\lvec(48 106)
\ifill f:0
\move(61 105)
\lvec(63 105)
\lvec(63 106)
\lvec(61 106)
\ifill f:0
\move(64 105)
\lvec(65 105)
\lvec(65 106)
\lvec(64 106)
\ifill f:0
\move(66 105)
\lvec(73 105)
\lvec(73 106)
\lvec(66 106)
\ifill f:0
\move(74 105)
\lvec(82 105)
\lvec(82 106)
\lvec(74 106)
\ifill f:0
\move(84 105)
\lvec(86 105)
\lvec(86 106)
\lvec(84 106)
\ifill f:0
\move(88 105)
\lvec(89 105)
\lvec(89 106)
\lvec(88 106)
\ifill f:0
\move(90 105)
\lvec(91 105)
\lvec(91 106)
\lvec(90 106)
\ifill f:0
\move(92 105)
\lvec(93 105)
\lvec(93 106)
\lvec(92 106)
\ifill f:0
\move(96 105)
\lvec(98 105)
\lvec(98 106)
\lvec(96 106)
\ifill f:0
\move(99 105)
\lvec(101 105)
\lvec(101 106)
\lvec(99 106)
\ifill f:0
\move(103 105)
\lvec(106 105)
\lvec(106 106)
\lvec(103 106)
\ifill f:0
\move(107 105)
\lvec(109 105)
\lvec(109 106)
\lvec(107 106)
\ifill f:0
\move(114 105)
\lvec(122 105)
\lvec(122 106)
\lvec(114 106)
\ifill f:0
\move(123 105)
\lvec(126 105)
\lvec(126 106)
\lvec(123 106)
\ifill f:0
\move(127 105)
\lvec(129 105)
\lvec(129 106)
\lvec(127 106)
\ifill f:0
\move(130 105)
\lvec(131 105)
\lvec(131 106)
\lvec(130 106)
\ifill f:0
\move(133 105)
\lvec(134 105)
\lvec(134 106)
\lvec(133 106)
\ifill f:0
\move(135 105)
\lvec(136 105)
\lvec(136 106)
\lvec(135 106)
\ifill f:0
\move(137 105)
\lvec(138 105)
\lvec(138 106)
\lvec(137 106)
\ifill f:0
\move(139 105)
\lvec(145 105)
\lvec(145 106)
\lvec(139 106)
\ifill f:0
\move(146 105)
\lvec(148 105)
\lvec(148 106)
\lvec(146 106)
\ifill f:0
\move(149 105)
\lvec(151 105)
\lvec(151 106)
\lvec(149 106)
\ifill f:0
\move(152 105)
\lvec(155 105)
\lvec(155 106)
\lvec(152 106)
\ifill f:0
\move(156 105)
\lvec(159 105)
\lvec(159 106)
\lvec(156 106)
\ifill f:0
\move(161 105)
\lvec(167 105)
\lvec(167 106)
\lvec(161 106)
\ifill f:0
\move(169 105)
\lvec(170 105)
\lvec(170 106)
\lvec(169 106)
\ifill f:0
\move(171 105)
\lvec(179 105)
\lvec(179 106)
\lvec(171 106)
\ifill f:0
\move(181 105)
\lvec(189 105)
\lvec(189 106)
\lvec(181 106)
\ifill f:0
\move(190 105)
\lvec(194 105)
\lvec(194 106)
\lvec(190 106)
\ifill f:0
\move(196 105)
\lvec(197 105)
\lvec(197 106)
\lvec(196 106)
\ifill f:0
\move(198 105)
\lvec(202 105)
\lvec(202 106)
\lvec(198 106)
\ifill f:0
\move(203 105)
\lvec(211 105)
\lvec(211 106)
\lvec(203 106)
\ifill f:0
\move(212 105)
\lvec(223 105)
\lvec(223 106)
\lvec(212 106)
\ifill f:0
\move(224 105)
\lvec(226 105)
\lvec(226 106)
\lvec(224 106)
\ifill f:0
\move(227 105)
\lvec(229 105)
\lvec(229 106)
\lvec(227 106)
\ifill f:0
\move(230 105)
\lvec(231 105)
\lvec(231 106)
\lvec(230 106)
\ifill f:0
\move(232 105)
\lvec(233 105)
\lvec(233 106)
\lvec(232 106)
\ifill f:0
\move(234 105)
\lvec(235 105)
\lvec(235 106)
\lvec(234 106)
\ifill f:0
\move(236 105)
\lvec(240 105)
\lvec(240 106)
\lvec(236 106)
\ifill f:0
\move(241 105)
\lvec(242 105)
\lvec(242 106)
\lvec(241 106)
\ifill f:0
\move(243 105)
\lvec(248 105)
\lvec(248 106)
\lvec(243 106)
\ifill f:0
\move(249 105)
\lvec(251 105)
\lvec(251 106)
\lvec(249 106)
\ifill f:0
\move(252 105)
\lvec(257 105)
\lvec(257 106)
\lvec(252 106)
\ifill f:0
\move(258 105)
\lvec(269 105)
\lvec(269 106)
\lvec(258 106)
\ifill f:0
\move(270 105)
\lvec(280 105)
\lvec(280 106)
\lvec(270 106)
\ifill f:0
\move(282 105)
\lvec(288 105)
\lvec(288 106)
\lvec(282 106)
\ifill f:0
\move(289 105)
\lvec(290 105)
\lvec(290 106)
\lvec(289 106)
\ifill f:0
\move(291 105)
\lvec(298 105)
\lvec(298 106)
\lvec(291 106)
\ifill f:0
\move(299 105)
\lvec(302 105)
\lvec(302 106)
\lvec(299 106)
\ifill f:0
\move(303 105)
\lvec(304 105)
\lvec(304 106)
\lvec(303 106)
\ifill f:0
\move(309 105)
\lvec(310 105)
\lvec(310 106)
\lvec(309 106)
\ifill f:0
\move(316 105)
\lvec(325 105)
\lvec(325 106)
\lvec(316 106)
\ifill f:0
\move(326 105)
\lvec(330 105)
\lvec(330 106)
\lvec(326 106)
\ifill f:0
\move(331 105)
\lvec(339 105)
\lvec(339 106)
\lvec(331 106)
\ifill f:0
\move(340 105)
\lvec(346 105)
\lvec(346 106)
\lvec(340 106)
\ifill f:0
\move(347 105)
\lvec(362 105)
\lvec(362 106)
\lvec(347 106)
\ifill f:0
\move(363 105)
\lvec(370 105)
\lvec(370 106)
\lvec(363 106)
\ifill f:0
\move(371 105)
\lvec(377 105)
\lvec(377 106)
\lvec(371 106)
\ifill f:0
\move(378 105)
\lvec(384 105)
\lvec(384 106)
\lvec(378 106)
\ifill f:0
\move(385 105)
\lvec(387 105)
\lvec(387 106)
\lvec(385 106)
\ifill f:0
\move(388 105)
\lvec(393 105)
\lvec(393 106)
\lvec(388 106)
\ifill f:0
\move(394 105)
\lvec(401 105)
\lvec(401 106)
\lvec(394 106)
\ifill f:0
\move(402 105)
\lvec(426 105)
\lvec(426 106)
\lvec(402 106)
\ifill f:0
\move(427 105)
\lvec(439 105)
\lvec(439 106)
\lvec(427 106)
\ifill f:0
\move(440 105)
\lvec(442 105)
\lvec(442 106)
\lvec(440 106)
\ifill f:0
\move(444 105)
\lvec(447 105)
\lvec(447 106)
\lvec(444 106)
\ifill f:0
\move(448 105)
\lvec(449 105)
\lvec(449 106)
\lvec(448 106)
\ifill f:0
\move(450 105)
\lvec(451 105)
\lvec(451 106)
\lvec(450 106)
\ifill f:0
\move(16 106)
\lvec(17 106)
\lvec(17 107)
\lvec(16 107)
\ifill f:0
\move(20 106)
\lvec(21 106)
\lvec(21 107)
\lvec(20 107)
\ifill f:0
\move(22 106)
\lvec(26 106)
\lvec(26 107)
\lvec(22 107)
\ifill f:0
\move(36 106)
\lvec(37 106)
\lvec(37 107)
\lvec(36 107)
\ifill f:0
\move(38 106)
\lvec(39 106)
\lvec(39 107)
\lvec(38 107)
\ifill f:0
\move(40 106)
\lvec(43 106)
\lvec(43 107)
\lvec(40 107)
\ifill f:0
\move(44 106)
\lvec(45 106)
\lvec(45 107)
\lvec(44 107)
\ifill f:0
\move(47 106)
\lvec(50 106)
\lvec(50 107)
\lvec(47 107)
\ifill f:0
\move(51 106)
\lvec(52 106)
\lvec(52 107)
\lvec(51 107)
\ifill f:0
\move(54 106)
\lvec(55 106)
\lvec(55 107)
\lvec(54 107)
\ifill f:0
\move(59 106)
\lvec(63 106)
\lvec(63 107)
\lvec(59 107)
\ifill f:0
\move(64 106)
\lvec(65 106)
\lvec(65 107)
\lvec(64 107)
\ifill f:0
\move(66 106)
\lvec(71 106)
\lvec(71 107)
\lvec(66 107)
\ifill f:0
\move(72 106)
\lvec(74 106)
\lvec(74 107)
\lvec(72 107)
\ifill f:0
\move(77 106)
\lvec(82 106)
\lvec(82 107)
\lvec(77 107)
\ifill f:0
\move(86 106)
\lvec(87 106)
\lvec(87 107)
\lvec(86 107)
\ifill f:0
\move(88 106)
\lvec(90 106)
\lvec(90 107)
\lvec(88 107)
\ifill f:0
\move(91 106)
\lvec(93 106)
\lvec(93 107)
\lvec(91 107)
\ifill f:0
\move(94 106)
\lvec(95 106)
\lvec(95 107)
\lvec(94 107)
\ifill f:0
\move(97 106)
\lvec(98 106)
\lvec(98 107)
\lvec(97 107)
\ifill f:0
\move(99 106)
\lvec(101 106)
\lvec(101 107)
\lvec(99 107)
\ifill f:0
\move(102 106)
\lvec(106 106)
\lvec(106 107)
\lvec(102 107)
\ifill f:0
\move(107 106)
\lvec(122 106)
\lvec(122 107)
\lvec(107 107)
\ifill f:0
\move(123 106)
\lvec(126 106)
\lvec(126 107)
\lvec(123 107)
\ifill f:0
\move(128 106)
\lvec(130 106)
\lvec(130 107)
\lvec(128 107)
\ifill f:0
\move(131 106)
\lvec(133 106)
\lvec(133 107)
\lvec(131 107)
\ifill f:0
\move(134 106)
\lvec(135 106)
\lvec(135 107)
\lvec(134 107)
\ifill f:0
\move(136 106)
\lvec(138 106)
\lvec(138 107)
\lvec(136 107)
\ifill f:0
\move(139 106)
\lvec(145 106)
\lvec(145 107)
\lvec(139 107)
\ifill f:0
\move(146 106)
\lvec(150 106)
\lvec(150 107)
\lvec(146 107)
\ifill f:0
\move(151 106)
\lvec(157 106)
\lvec(157 107)
\lvec(151 107)
\ifill f:0
\move(158 106)
\lvec(168 106)
\lvec(168 107)
\lvec(158 107)
\ifill f:0
\move(169 106)
\lvec(170 106)
\lvec(170 107)
\lvec(169 107)
\ifill f:0
\move(171 106)
\lvec(186 106)
\lvec(186 107)
\lvec(171 107)
\ifill f:0
\move(188 106)
\lvec(194 106)
\lvec(194 107)
\lvec(188 107)
\ifill f:0
\move(196 106)
\lvec(197 106)
\lvec(197 107)
\lvec(196 107)
\ifill f:0
\move(198 106)
\lvec(203 106)
\lvec(203 107)
\lvec(198 107)
\ifill f:0
\move(204 106)
\lvec(223 106)
\lvec(223 107)
\lvec(204 107)
\ifill f:0
\move(224 106)
\lvec(226 106)
\lvec(226 107)
\lvec(224 107)
\ifill f:0
\move(227 106)
\lvec(239 106)
\lvec(239 107)
\lvec(227 107)
\ifill f:0
\move(240 106)
\lvec(241 106)
\lvec(241 107)
\lvec(240 107)
\ifill f:0
\move(242 106)
\lvec(251 106)
\lvec(251 107)
\lvec(242 107)
\ifill f:0
\move(252 106)
\lvec(257 106)
\lvec(257 107)
\lvec(252 107)
\ifill f:0
\move(258 106)
\lvec(277 106)
\lvec(277 107)
\lvec(258 107)
\ifill f:0
\move(278 106)
\lvec(282 106)
\lvec(282 107)
\lvec(278 107)
\ifill f:0
\move(283 106)
\lvec(290 106)
\lvec(290 107)
\lvec(283 107)
\ifill f:0
\move(291 106)
\lvec(297 106)
\lvec(297 107)
\lvec(291 107)
\ifill f:0
\move(299 106)
\lvec(325 106)
\lvec(325 107)
\lvec(299 107)
\ifill f:0
\move(326 106)
\lvec(356 106)
\lvec(356 107)
\lvec(326 107)
\ifill f:0
\move(357 106)
\lvec(362 106)
\lvec(362 107)
\lvec(357 107)
\ifill f:0
\move(363 106)
\lvec(366 106)
\lvec(366 107)
\lvec(363 107)
\ifill f:0
\move(367 106)
\lvec(375 106)
\lvec(375 107)
\lvec(367 107)
\ifill f:0
\move(376 106)
\lvec(379 106)
\lvec(379 107)
\lvec(376 107)
\ifill f:0
\move(380 106)
\lvec(389 106)
\lvec(389 107)
\lvec(380 107)
\ifill f:0
\move(390 106)
\lvec(401 106)
\lvec(401 107)
\lvec(390 107)
\ifill f:0
\move(402 106)
\lvec(410 106)
\lvec(410 107)
\lvec(402 107)
\ifill f:0
\move(411 106)
\lvec(423 106)
\lvec(423 107)
\lvec(411 107)
\ifill f:0
\move(424 106)
\lvec(430 106)
\lvec(430 107)
\lvec(424 107)
\ifill f:0
\move(431 106)
\lvec(439 106)
\lvec(439 107)
\lvec(431 107)
\ifill f:0
\move(440 106)
\lvec(442 106)
\lvec(442 107)
\lvec(440 107)
\ifill f:0
\move(443 106)
\lvec(451 106)
\lvec(451 107)
\lvec(443 107)
\ifill f:0
\move(16 107)
\lvec(17 107)
\lvec(17 108)
\lvec(16 108)
\ifill f:0
\move(19 107)
\lvec(22 107)
\lvec(22 108)
\lvec(19 108)
\ifill f:0
\move(23 107)
\lvec(26 107)
\lvec(26 108)
\lvec(23 108)
\ifill f:0
\move(36 107)
\lvec(37 107)
\lvec(37 108)
\lvec(36 108)
\ifill f:0
\move(38 107)
\lvec(39 107)
\lvec(39 108)
\lvec(38 108)
\ifill f:0
\move(40 107)
\lvec(41 107)
\lvec(41 108)
\lvec(40 108)
\ifill f:0
\move(42 107)
\lvec(45 107)
\lvec(45 108)
\lvec(42 108)
\ifill f:0
\move(47 107)
\lvec(50 107)
\lvec(50 108)
\lvec(47 108)
\ifill f:0
\move(51 107)
\lvec(52 107)
\lvec(52 108)
\lvec(51 108)
\ifill f:0
\move(57 107)
\lvec(58 107)
\lvec(58 108)
\lvec(57 108)
\ifill f:0
\move(59 107)
\lvec(62 107)
\lvec(62 108)
\lvec(59 108)
\ifill f:0
\move(64 107)
\lvec(65 107)
\lvec(65 108)
\lvec(64 108)
\ifill f:0
\move(66 107)
\lvec(68 107)
\lvec(68 108)
\lvec(66 108)
\ifill f:0
\move(69 107)
\lvec(70 107)
\lvec(70 108)
\lvec(69 108)
\ifill f:0
\move(71 107)
\lvec(73 107)
\lvec(73 108)
\lvec(71 108)
\ifill f:0
\move(74 107)
\lvec(78 107)
\lvec(78 108)
\lvec(74 108)
\ifill f:0
\move(81 107)
\lvec(82 107)
\lvec(82 108)
\lvec(81 108)
\ifill f:0
\move(86 107)
\lvec(88 107)
\lvec(88 108)
\lvec(86 108)
\ifill f:0
\move(89 107)
\lvec(91 107)
\lvec(91 108)
\lvec(89 108)
\ifill f:0
\move(92 107)
\lvec(93 107)
\lvec(93 108)
\lvec(92 108)
\ifill f:0
\move(95 107)
\lvec(96 107)
\lvec(96 108)
\lvec(95 108)
\ifill f:0
\move(97 107)
\lvec(98 107)
\lvec(98 108)
\lvec(97 108)
\ifill f:0
\move(99 107)
\lvec(101 107)
\lvec(101 108)
\lvec(99 108)
\ifill f:0
\move(102 107)
\lvec(109 107)
\lvec(109 108)
\lvec(102 108)
\ifill f:0
\move(111 107)
\lvec(122 107)
\lvec(122 108)
\lvec(111 108)
\ifill f:0
\move(125 107)
\lvec(127 107)
\lvec(127 108)
\lvec(125 108)
\ifill f:0
\move(128 107)
\lvec(131 107)
\lvec(131 108)
\lvec(128 108)
\ifill f:0
\move(133 107)
\lvec(134 107)
\lvec(134 108)
\lvec(133 108)
\ifill f:0
\move(135 107)
\lvec(139 107)
\lvec(139 108)
\lvec(135 108)
\ifill f:0
\move(140 107)
\lvec(141 107)
\lvec(141 108)
\lvec(140 108)
\ifill f:0
\move(142 107)
\lvec(143 107)
\lvec(143 108)
\lvec(142 108)
\ifill f:0
\move(144 107)
\lvec(145 107)
\lvec(145 108)
\lvec(144 108)
\ifill f:0
\move(146 107)
\lvec(147 107)
\lvec(147 108)
\lvec(146 108)
\ifill f:0
\move(148 107)
\lvec(155 107)
\lvec(155 108)
\lvec(148 108)
\ifill f:0
\move(156 107)
\lvec(163 107)
\lvec(163 108)
\lvec(156 108)
\ifill f:0
\move(164 107)
\lvec(168 107)
\lvec(168 108)
\lvec(164 108)
\ifill f:0
\move(169 107)
\lvec(170 107)
\lvec(170 108)
\lvec(169 108)
\ifill f:0
\move(171 107)
\lvec(176 107)
\lvec(176 108)
\lvec(171 108)
\ifill f:0
\move(177 107)
\lvec(179 107)
\lvec(179 108)
\lvec(177 108)
\ifill f:0
\move(181 107)
\lvec(186 107)
\lvec(186 108)
\lvec(181 108)
\ifill f:0
\move(187 107)
\lvec(193 107)
\lvec(193 108)
\lvec(187 108)
\ifill f:0
\move(195 107)
\lvec(197 107)
\lvec(197 108)
\lvec(195 108)
\ifill f:0
\move(198 107)
\lvec(199 107)
\lvec(199 108)
\lvec(198 108)
\ifill f:0
\move(200 107)
\lvec(220 107)
\lvec(220 108)
\lvec(200 108)
\ifill f:0
\move(221 107)
\lvec(226 107)
\lvec(226 108)
\lvec(221 108)
\ifill f:0
\move(227 107)
\lvec(232 107)
\lvec(232 108)
\lvec(227 108)
\ifill f:0
\move(233 107)
\lvec(234 107)
\lvec(234 108)
\lvec(233 108)
\ifill f:0
\move(235 107)
\lvec(238 107)
\lvec(238 108)
\lvec(235 108)
\ifill f:0
\move(239 107)
\lvec(240 107)
\lvec(240 108)
\lvec(239 108)
\ifill f:0
\move(241 107)
\lvec(242 107)
\lvec(242 108)
\lvec(241 108)
\ifill f:0
\move(243 107)
\lvec(247 107)
\lvec(247 108)
\lvec(243 108)
\ifill f:0
\move(248 107)
\lvec(252 107)
\lvec(252 108)
\lvec(248 108)
\ifill f:0
\move(253 107)
\lvec(257 107)
\lvec(257 108)
\lvec(253 108)
\ifill f:0
\move(258 107)
\lvec(274 107)
\lvec(274 108)
\lvec(258 108)
\ifill f:0
\move(275 107)
\lvec(283 107)
\lvec(283 108)
\lvec(275 108)
\ifill f:0
\move(284 107)
\lvec(290 107)
\lvec(290 108)
\lvec(284 108)
\ifill f:0
\move(291 107)
\lvec(295 107)
\lvec(295 108)
\lvec(291 108)
\ifill f:0
\move(297 107)
\lvec(305 107)
\lvec(305 108)
\lvec(297 108)
\ifill f:0
\move(307 107)
\lvec(325 107)
\lvec(325 108)
\lvec(307 108)
\ifill f:0
\move(326 107)
\lvec(336 107)
\lvec(336 108)
\lvec(326 108)
\ifill f:0
\move(339 107)
\lvec(347 107)
\lvec(347 108)
\lvec(339 108)
\ifill f:0
\move(348 107)
\lvec(355 107)
\lvec(355 108)
\lvec(348 108)
\ifill f:0
\move(356 107)
\lvec(362 107)
\lvec(362 108)
\lvec(356 108)
\ifill f:0
\move(363 107)
\lvec(367 107)
\lvec(367 108)
\lvec(363 108)
\ifill f:0
\move(368 107)
\lvec(376 107)
\lvec(376 108)
\lvec(368 108)
\ifill f:0
\move(377 107)
\lvec(388 107)
\lvec(388 108)
\lvec(377 108)
\ifill f:0
\move(389 107)
\lvec(392 107)
\lvec(392 108)
\lvec(389 108)
\ifill f:0
\move(393 107)
\lvec(395 107)
\lvec(395 108)
\lvec(393 108)
\ifill f:0
\move(396 107)
\lvec(401 107)
\lvec(401 108)
\lvec(396 108)
\ifill f:0
\move(402 107)
\lvec(419 107)
\lvec(419 108)
\lvec(402 108)
\ifill f:0
\move(420 107)
\lvec(434 107)
\lvec(434 108)
\lvec(420 108)
\ifill f:0
\move(435 107)
\lvec(442 107)
\lvec(442 108)
\lvec(435 108)
\ifill f:0
\move(443 107)
\lvec(450 107)
\lvec(450 108)
\lvec(443 108)
\ifill f:0
\move(16 108)
\lvec(17 108)
\lvec(17 109)
\lvec(16 109)
\ifill f:0
\move(20 108)
\lvec(21 108)
\lvec(21 109)
\lvec(20 109)
\ifill f:0
\move(24 108)
\lvec(26 108)
\lvec(26 109)
\lvec(24 109)
\ifill f:0
\move(36 108)
\lvec(37 108)
\lvec(37 109)
\lvec(36 109)
\ifill f:0
\move(38 108)
\lvec(43 108)
\lvec(43 109)
\lvec(38 109)
\ifill f:0
\move(44 108)
\lvec(45 108)
\lvec(45 109)
\lvec(44 109)
\ifill f:0
\move(47 108)
\lvec(50 108)
\lvec(50 109)
\lvec(47 109)
\ifill f:0
\move(51 108)
\lvec(52 108)
\lvec(52 109)
\lvec(51 109)
\ifill f:0
\move(64 108)
\lvec(65 108)
\lvec(65 109)
\lvec(64 109)
\ifill f:0
\move(67 108)
\lvec(71 108)
\lvec(71 109)
\lvec(67 109)
\ifill f:0
\move(72 108)
\lvec(74 108)
\lvec(74 109)
\lvec(72 109)
\ifill f:0
\move(75 108)
\lvec(78 108)
\lvec(78 109)
\lvec(75 109)
\ifill f:0
\move(81 108)
\lvec(82 108)
\lvec(82 109)
\lvec(81 109)
\ifill f:0
\move(88 108)
\lvec(90 108)
\lvec(90 109)
\lvec(88 109)
\ifill f:0
\move(91 108)
\lvec(93 108)
\lvec(93 109)
\lvec(91 109)
\ifill f:0
\move(95 108)
\lvec(96 108)
\lvec(96 109)
\lvec(95 109)
\ifill f:0
\move(97 108)
\lvec(98 108)
\lvec(98 109)
\lvec(97 109)
\ifill f:0
\move(99 108)
\lvec(101 108)
\lvec(101 109)
\lvec(99 109)
\ifill f:0
\move(102 108)
\lvec(104 108)
\lvec(104 109)
\lvec(102 109)
\ifill f:0
\move(105 108)
\lvec(113 108)
\lvec(113 109)
\lvec(105 109)
\ifill f:0
\move(114 108)
\lvec(122 108)
\lvec(122 109)
\lvec(114 109)
\ifill f:0
\move(126 108)
\lvec(129 108)
\lvec(129 109)
\lvec(126 109)
\ifill f:0
\move(131 108)
\lvec(133 108)
\lvec(133 109)
\lvec(131 109)
\ifill f:0
\move(134 108)
\lvec(136 108)
\lvec(136 109)
\lvec(134 109)
\ifill f:0
\move(137 108)
\lvec(143 108)
\lvec(143 109)
\lvec(137 109)
\ifill f:0
\move(144 108)
\lvec(145 108)
\lvec(145 109)
\lvec(144 109)
\ifill f:0
\move(146 108)
\lvec(147 108)
\lvec(147 109)
\lvec(146 109)
\ifill f:0
\move(148 108)
\lvec(149 108)
\lvec(149 109)
\lvec(148 109)
\ifill f:0
\move(150 108)
\lvec(154 108)
\lvec(154 109)
\lvec(150 109)
\ifill f:0
\move(155 108)
\lvec(157 108)
\lvec(157 109)
\lvec(155 109)
\ifill f:0
\move(158 108)
\lvec(160 108)
\lvec(160 109)
\lvec(158 109)
\ifill f:0
\move(161 108)
\lvec(163 108)
\lvec(163 109)
\lvec(161 109)
\ifill f:0
\move(165 108)
\lvec(168 108)
\lvec(168 109)
\lvec(165 109)
\ifill f:0
\move(169 108)
\lvec(170 108)
\lvec(170 109)
\lvec(169 109)
\ifill f:0
\move(171 108)
\lvec(174 108)
\lvec(174 109)
\lvec(171 109)
\ifill f:0
\move(178 108)
\lvec(191 108)
\lvec(191 109)
\lvec(178 109)
\ifill f:0
\move(195 108)
\lvec(197 108)
\lvec(197 109)
\lvec(195 109)
\ifill f:0
\move(198 108)
\lvec(213 108)
\lvec(213 109)
\lvec(198 109)
\ifill f:0
\move(214 108)
\lvec(226 108)
\lvec(226 109)
\lvec(214 109)
\ifill f:0
\move(227 108)
\lvec(228 108)
\lvec(228 109)
\lvec(227 109)
\ifill f:0
\move(229 108)
\lvec(230 108)
\lvec(230 109)
\lvec(229 109)
\ifill f:0
\move(231 108)
\lvec(235 108)
\lvec(235 109)
\lvec(231 109)
\ifill f:0
\move(236 108)
\lvec(241 108)
\lvec(241 109)
\lvec(236 109)
\ifill f:0
\move(242 108)
\lvec(257 108)
\lvec(257 109)
\lvec(242 109)
\ifill f:0
\move(258 108)
\lvec(260 108)
\lvec(260 109)
\lvec(258 109)
\ifill f:0
\move(261 108)
\lvec(279 108)
\lvec(279 109)
\lvec(261 109)
\ifill f:0
\move(280 108)
\lvec(284 108)
\lvec(284 109)
\lvec(280 109)
\ifill f:0
\move(285 108)
\lvec(290 108)
\lvec(290 109)
\lvec(285 109)
\ifill f:0
\move(291 108)
\lvec(293 108)
\lvec(293 109)
\lvec(291 109)
\ifill f:0
\move(294 108)
\lvec(310 108)
\lvec(310 109)
\lvec(294 109)
\ifill f:0
\move(313 108)
\lvec(325 108)
\lvec(325 109)
\lvec(313 109)
\ifill f:0
\move(326 108)
\lvec(342 108)
\lvec(342 109)
\lvec(326 109)
\ifill f:0
\move(343 108)
\lvec(362 108)
\lvec(362 109)
\lvec(343 109)
\ifill f:0
\move(363 108)
\lvec(368 108)
\lvec(368 109)
\lvec(363 109)
\ifill f:0
\move(369 108)
\lvec(373 108)
\lvec(373 109)
\lvec(369 109)
\ifill f:0
\move(374 108)
\lvec(378 108)
\lvec(378 109)
\lvec(374 109)
\ifill f:0
\move(379 108)
\lvec(383 108)
\lvec(383 109)
\lvec(379 109)
\ifill f:0
\move(384 108)
\lvec(387 108)
\lvec(387 109)
\lvec(384 109)
\ifill f:0
\move(388 108)
\lvec(391 108)
\lvec(391 109)
\lvec(388 109)
\ifill f:0
\move(392 108)
\lvec(395 108)
\lvec(395 109)
\lvec(392 109)
\ifill f:0
\move(396 108)
\lvec(401 108)
\lvec(401 109)
\lvec(396 109)
\ifill f:0
\move(402 108)
\lvec(405 108)
\lvec(405 109)
\lvec(402 109)
\ifill f:0
\move(406 108)
\lvec(423 108)
\lvec(423 109)
\lvec(406 109)
\ifill f:0
\move(424 108)
\lvec(431 108)
\lvec(431 109)
\lvec(424 109)
\ifill f:0
\move(432 108)
\lvec(442 108)
\lvec(442 109)
\lvec(432 109)
\ifill f:0
\move(443 108)
\lvec(448 108)
\lvec(448 109)
\lvec(443 109)
\ifill f:0
\move(449 108)
\lvec(451 108)
\lvec(451 109)
\lvec(449 109)
\ifill f:0
\move(16 109)
\lvec(17 109)
\lvec(17 110)
\lvec(16 110)
\ifill f:0
\move(18 109)
\lvec(19 109)
\lvec(19 110)
\lvec(18 110)
\ifill f:0
\move(20 109)
\lvec(21 109)
\lvec(21 110)
\lvec(20 110)
\ifill f:0
\move(22 109)
\lvec(23 109)
\lvec(23 110)
\lvec(22 110)
\ifill f:0
\move(25 109)
\lvec(26 109)
\lvec(26 110)
\lvec(25 110)
\ifill f:0
\move(36 109)
\lvec(37 109)
\lvec(37 110)
\lvec(36 110)
\ifill f:0
\move(38 109)
\lvec(39 109)
\lvec(39 110)
\lvec(38 110)
\ifill f:0
\move(40 109)
\lvec(41 109)
\lvec(41 110)
\lvec(40 110)
\ifill f:0
\move(42 109)
\lvec(46 109)
\lvec(46 110)
\lvec(42 110)
\ifill f:0
\move(47 109)
\lvec(50 109)
\lvec(50 110)
\lvec(47 110)
\ifill f:0
\move(52 109)
\lvec(53 109)
\lvec(53 110)
\lvec(52 110)
\ifill f:0
\move(54 109)
\lvec(55 109)
\lvec(55 110)
\lvec(54 110)
\ifill f:0
\move(57 109)
\lvec(58 109)
\lvec(58 110)
\lvec(57 110)
\ifill f:0
\move(59 109)
\lvec(60 109)
\lvec(60 110)
\lvec(59 110)
\ifill f:0
\move(62 109)
\lvec(63 109)
\lvec(63 110)
\lvec(62 110)
\ifill f:0
\move(64 109)
\lvec(65 109)
\lvec(65 110)
\lvec(64 110)
\ifill f:0
\move(66 109)
\lvec(67 109)
\lvec(67 110)
\lvec(66 110)
\ifill f:0
\move(68 109)
\lvec(74 109)
\lvec(74 110)
\lvec(68 110)
\ifill f:0
\move(76 109)
\lvec(79 109)
\lvec(79 110)
\lvec(76 110)
\ifill f:0
\move(81 109)
\lvec(82 109)
\lvec(82 110)
\lvec(81 110)
\ifill f:0
\move(83 109)
\lvec(86 109)
\lvec(86 110)
\lvec(83 110)
\ifill f:0
\move(89 109)
\lvec(92 109)
\lvec(92 110)
\lvec(89 110)
\ifill f:0
\move(96 109)
\lvec(97 109)
\lvec(97 110)
\lvec(96 110)
\ifill f:0
\move(99 109)
\lvec(101 109)
\lvec(101 110)
\lvec(99 110)
\ifill f:0
\move(102 109)
\lvec(103 109)
\lvec(103 110)
\lvec(102 110)
\ifill f:0
\move(104 109)
\lvec(106 109)
\lvec(106 110)
\lvec(104 110)
\ifill f:0
\move(107 109)
\lvec(116 109)
\lvec(116 110)
\lvec(107 110)
\ifill f:0
\move(117 109)
\lvec(118 109)
\lvec(118 110)
\lvec(117 110)
\ifill f:0
\move(120 109)
\lvec(122 109)
\lvec(122 110)
\lvec(120 110)
\ifill f:0
\move(127 109)
\lvec(130 109)
\lvec(130 110)
\lvec(127 110)
\ifill f:0
\move(133 109)
\lvec(135 109)
\lvec(135 110)
\lvec(133 110)
\ifill f:0
\move(136 109)
\lvec(138 109)
\lvec(138 110)
\lvec(136 110)
\ifill f:0
\move(139 109)
\lvec(143 109)
\lvec(143 110)
\lvec(139 110)
\ifill f:0
\move(144 109)
\lvec(145 109)
\lvec(145 110)
\lvec(144 110)
\ifill f:0
\move(146 109)
\lvec(153 109)
\lvec(153 110)
\lvec(146 110)
\ifill f:0
\move(154 109)
\lvec(155 109)
\lvec(155 110)
\lvec(154 110)
\ifill f:0
\move(156 109)
\lvec(161 109)
\lvec(161 110)
\lvec(156 110)
\ifill f:0
\move(162 109)
\lvec(170 109)
\lvec(170 110)
\lvec(162 110)
\ifill f:0
\move(171 109)
\lvec(173 109)
\lvec(173 110)
\lvec(171 110)
\ifill f:0
\move(175 109)
\lvec(185 109)
\lvec(185 110)
\lvec(175 110)
\ifill f:0
\move(188 109)
\lvec(189 109)
\lvec(189 110)
\lvec(188 110)
\ifill f:0
\move(193 109)
\lvec(194 109)
\lvec(194 110)
\lvec(193 110)
\ifill f:0
\move(195 109)
\lvec(197 109)
\lvec(197 110)
\lvec(195 110)
\ifill f:0
\move(198 109)
\lvec(201 109)
\lvec(201 110)
\lvec(198 110)
\ifill f:0
\move(202 109)
\lvec(207 109)
\lvec(207 110)
\lvec(202 110)
\ifill f:0
\move(208 109)
\lvec(219 109)
\lvec(219 110)
\lvec(208 110)
\ifill f:0
\move(220 109)
\lvec(222 109)
\lvec(222 110)
\lvec(220 110)
\ifill f:0
\move(223 109)
\lvec(226 109)
\lvec(226 110)
\lvec(223 110)
\ifill f:0
\move(227 109)
\lvec(228 109)
\lvec(228 110)
\lvec(227 110)
\ifill f:0
\move(229 109)
\lvec(233 109)
\lvec(233 110)
\lvec(229 110)
\ifill f:0
\move(234 109)
\lvec(240 109)
\lvec(240 110)
\lvec(234 110)
\ifill f:0
\move(241 109)
\lvec(242 109)
\lvec(242 110)
\lvec(241 110)
\ifill f:0
\move(243 109)
\lvec(244 109)
\lvec(244 110)
\lvec(243 110)
\ifill f:0
\move(245 109)
\lvec(248 109)
\lvec(248 110)
\lvec(245 110)
\ifill f:0
\move(249 109)
\lvec(250 109)
\lvec(250 110)
\lvec(249 110)
\ifill f:0
\move(251 109)
\lvec(257 109)
\lvec(257 110)
\lvec(251 110)
\ifill f:0
\move(258 109)
\lvec(277 109)
\lvec(277 110)
\lvec(258 110)
\ifill f:0
\move(278 109)
\lvec(280 109)
\lvec(280 110)
\lvec(278 110)
\ifill f:0
\move(281 109)
\lvec(290 109)
\lvec(290 110)
\lvec(281 110)
\ifill f:0
\move(291 109)
\lvec(293 109)
\lvec(293 110)
\lvec(291 110)
\ifill f:0
\move(294 109)
\lvec(298 109)
\lvec(298 110)
\lvec(294 110)
\ifill f:0
\move(300 109)
\lvec(305 109)
\lvec(305 110)
\lvec(300 110)
\ifill f:0
\move(307 109)
\lvec(314 109)
\lvec(314 110)
\lvec(307 110)
\ifill f:0
\move(317 109)
\lvec(325 109)
\lvec(325 110)
\lvec(317 110)
\ifill f:0
\move(326 109)
\lvec(351 109)
\lvec(351 110)
\lvec(326 110)
\ifill f:0
\move(353 109)
\lvec(362 109)
\lvec(362 110)
\lvec(353 110)
\ifill f:0
\move(363 109)
\lvec(369 109)
\lvec(369 110)
\lvec(363 110)
\ifill f:0
\move(370 109)
\lvec(375 109)
\lvec(375 110)
\lvec(370 110)
\ifill f:0
\move(376 109)
\lvec(394 109)
\lvec(394 110)
\lvec(376 110)
\ifill f:0
\move(395 109)
\lvec(398 109)
\lvec(398 110)
\lvec(395 110)
\ifill f:0
\move(399 109)
\lvec(401 109)
\lvec(401 110)
\lvec(399 110)
\ifill f:0
\move(402 109)
\lvec(409 109)
\lvec(409 110)
\lvec(402 110)
\ifill f:0
\move(410 109)
\lvec(419 109)
\lvec(419 110)
\lvec(410 110)
\ifill f:0
\move(420 109)
\lvec(436 109)
\lvec(436 110)
\lvec(420 110)
\ifill f:0
\move(437 109)
\lvec(442 109)
\lvec(442 110)
\lvec(437 110)
\ifill f:0
\move(443 109)
\lvec(446 109)
\lvec(446 110)
\lvec(443 110)
\ifill f:0
\move(447 109)
\lvec(451 109)
\lvec(451 110)
\lvec(447 110)
\ifill f:0
\move(16 110)
\lvec(17 110)
\lvec(17 111)
\lvec(16 111)
\ifill f:0
\move(19 110)
\lvec(21 110)
\lvec(21 111)
\lvec(19 111)
\ifill f:0
\move(25 110)
\lvec(26 110)
\lvec(26 111)
\lvec(25 111)
\ifill f:0
\move(36 110)
\lvec(37 110)
\lvec(37 111)
\lvec(36 111)
\ifill f:0
\move(38 110)
\lvec(39 110)
\lvec(39 111)
\lvec(38 111)
\ifill f:0
\move(40 110)
\lvec(45 110)
\lvec(45 111)
\lvec(40 111)
\ifill f:0
\move(47 110)
\lvec(50 110)
\lvec(50 111)
\lvec(47 111)
\ifill f:0
\move(58 110)
\lvec(63 110)
\lvec(63 111)
\lvec(58 111)
\ifill f:0
\move(64 110)
\lvec(65 110)
\lvec(65 111)
\lvec(64 111)
\ifill f:0
\move(66 110)
\lvec(67 110)
\lvec(67 111)
\lvec(66 111)
\ifill f:0
\move(68 110)
\lvec(71 110)
\lvec(71 111)
\lvec(68 111)
\ifill f:0
\move(72 110)
\lvec(76 110)
\lvec(76 111)
\lvec(72 111)
\ifill f:0
\move(77 110)
\lvec(80 110)
\lvec(80 111)
\lvec(77 111)
\ifill f:0
\move(81 110)
\lvec(82 110)
\lvec(82 111)
\lvec(81 111)
\ifill f:0
\move(83 110)
\lvec(85 110)
\lvec(85 111)
\lvec(83 111)
\ifill f:0
\move(86 110)
\lvec(89 110)
\lvec(89 111)
\lvec(86 111)
\ifill f:0
\move(92 110)
\lvec(93 110)
\lvec(93 111)
\lvec(92 111)
\ifill f:0
\move(96 110)
\lvec(98 110)
\lvec(98 111)
\lvec(96 111)
\ifill f:0
\move(100 110)
\lvec(101 110)
\lvec(101 111)
\lvec(100 111)
\ifill f:0
\move(102 110)
\lvec(103 110)
\lvec(103 111)
\lvec(102 111)
\ifill f:0
\move(104 110)
\lvec(111 110)
\lvec(111 111)
\lvec(104 111)
\ifill f:0
\move(112 110)
\lvec(119 110)
\lvec(119 111)
\lvec(112 111)
\ifill f:0
\move(121 110)
\lvec(122 110)
\lvec(122 111)
\lvec(121 111)
\ifill f:0
\move(128 110)
\lvec(133 110)
\lvec(133 111)
\lvec(128 111)
\ifill f:0
\move(135 110)
\lvec(137 110)
\lvec(137 111)
\lvec(135 111)
\ifill f:0
\move(139 110)
\lvec(140 110)
\lvec(140 111)
\lvec(139 111)
\ifill f:0
\move(141 110)
\lvec(143 110)
\lvec(143 111)
\lvec(141 111)
\ifill f:0
\move(144 110)
\lvec(145 110)
\lvec(145 111)
\lvec(144 111)
\ifill f:0
\move(146 110)
\lvec(159 110)
\lvec(159 111)
\lvec(146 111)
\ifill f:0
\move(160 110)
\lvec(162 110)
\lvec(162 111)
\lvec(160 111)
\ifill f:0
\move(163 110)
\lvec(165 110)
\lvec(165 111)
\lvec(163 111)
\ifill f:0
\move(166 110)
\lvec(170 110)
\lvec(170 111)
\lvec(166 111)
\ifill f:0
\move(171 110)
\lvec(172 110)
\lvec(172 111)
\lvec(171 111)
\ifill f:0
\move(174 110)
\lvec(179 110)
\lvec(179 111)
\lvec(174 111)
\ifill f:0
\move(181 110)
\lvec(191 110)
\lvec(191 111)
\lvec(181 111)
\ifill f:0
\move(192 110)
\lvec(197 110)
\lvec(197 111)
\lvec(192 111)
\ifill f:0
\move(198 110)
\lvec(202 110)
\lvec(202 111)
\lvec(198 111)
\ifill f:0
\move(203 110)
\lvec(214 110)
\lvec(214 111)
\lvec(203 111)
\ifill f:0
\move(215 110)
\lvec(218 110)
\lvec(218 111)
\lvec(215 111)
\ifill f:0
\move(219 110)
\lvec(226 110)
\lvec(226 111)
\lvec(219 111)
\ifill f:0
\move(227 110)
\lvec(228 110)
\lvec(228 111)
\lvec(227 111)
\ifill f:0
\move(229 110)
\lvec(231 110)
\lvec(231 111)
\lvec(229 111)
\ifill f:0
\move(232 110)
\lvec(234 110)
\lvec(234 111)
\lvec(232 111)
\ifill f:0
\move(235 110)
\lvec(241 110)
\lvec(241 111)
\lvec(235 111)
\ifill f:0
\move(242 110)
\lvec(257 110)
\lvec(257 111)
\lvec(242 111)
\ifill f:0
\move(258 110)
\lvec(269 110)
\lvec(269 111)
\lvec(258 111)
\ifill f:0
\move(270 110)
\lvec(275 110)
\lvec(275 111)
\lvec(270 111)
\ifill f:0
\move(276 110)
\lvec(290 110)
\lvec(290 111)
\lvec(276 111)
\ifill f:0
\move(291 110)
\lvec(293 110)
\lvec(293 111)
\lvec(291 111)
\ifill f:0
\move(294 110)
\lvec(309 110)
\lvec(309 111)
\lvec(294 111)
\ifill f:0
\move(310 110)
\lvec(318 110)
\lvec(318 111)
\lvec(310 111)
\ifill f:0
\move(319 110)
\lvec(325 110)
\lvec(325 111)
\lvec(319 111)
\ifill f:0
\move(326 110)
\lvec(335 110)
\lvec(335 111)
\lvec(326 111)
\ifill f:0
\move(339 110)
\lvec(340 110)
\lvec(340 111)
\lvec(339 111)
\ifill f:0
\move(345 110)
\lvec(346 110)
\lvec(346 111)
\lvec(345 111)
\ifill f:0
\move(347 110)
\lvec(362 110)
\lvec(362 111)
\lvec(347 111)
\ifill f:0
\move(363 110)
\lvec(370 110)
\lvec(370 111)
\lvec(363 111)
\ifill f:0
\move(372 110)
\lvec(389 110)
\lvec(389 111)
\lvec(372 111)
\ifill f:0
\move(390 110)
\lvec(394 110)
\lvec(394 111)
\lvec(390 111)
\ifill f:0
\move(395 110)
\lvec(398 110)
\lvec(398 111)
\lvec(395 111)
\ifill f:0
\move(399 110)
\lvec(401 110)
\lvec(401 111)
\lvec(399 111)
\ifill f:0
\move(403 110)
\lvec(406 110)
\lvec(406 111)
\lvec(403 111)
\ifill f:0
\move(407 110)
\lvec(417 110)
\lvec(417 111)
\lvec(407 111)
\ifill f:0
\move(418 110)
\lvec(442 110)
\lvec(442 111)
\lvec(418 111)
\ifill f:0
\move(443 110)
\lvec(444 110)
\lvec(444 111)
\lvec(443 111)
\ifill f:0
\move(445 110)
\lvec(449 110)
\lvec(449 111)
\lvec(445 111)
\ifill f:0
\move(450 110)
\lvec(451 110)
\lvec(451 111)
\lvec(450 111)
\ifill f:0
\move(15 111)
\lvec(17 111)
\lvec(17 112)
\lvec(15 112)
\ifill f:0
\move(20 111)
\lvec(21 111)
\lvec(21 112)
\lvec(20 112)
\ifill f:0
\move(23 111)
\lvec(24 111)
\lvec(24 112)
\lvec(23 112)
\ifill f:0
\move(25 111)
\lvec(26 111)
\lvec(26 112)
\lvec(25 112)
\ifill f:0
\move(36 111)
\lvec(37 111)
\lvec(37 112)
\lvec(36 112)
\ifill f:0
\move(40 111)
\lvec(45 111)
\lvec(45 112)
\lvec(40 112)
\ifill f:0
\move(49 111)
\lvec(50 111)
\lvec(50 112)
\lvec(49 112)
\ifill f:0
\move(57 111)
\lvec(58 111)
\lvec(58 112)
\lvec(57 112)
\ifill f:0
\move(60 111)
\lvec(63 111)
\lvec(63 112)
\lvec(60 112)
\ifill f:0
\move(64 111)
\lvec(65 111)
\lvec(65 112)
\lvec(64 112)
\ifill f:0
\move(66 111)
\lvec(68 111)
\lvec(68 112)
\lvec(66 112)
\ifill f:0
\move(69 111)
\lvec(74 111)
\lvec(74 112)
\lvec(69 112)
\ifill f:0
\move(75 111)
\lvec(80 111)
\lvec(80 112)
\lvec(75 112)
\ifill f:0
\move(81 111)
\lvec(82 111)
\lvec(82 112)
\lvec(81 112)
\ifill f:0
\move(83 111)
\lvec(85 111)
\lvec(85 112)
\lvec(83 112)
\ifill f:0
\move(88 111)
\lvec(92 111)
\lvec(92 112)
\lvec(88 112)
\ifill f:0
\move(95 111)
\lvec(96 111)
\lvec(96 112)
\lvec(95 112)
\ifill f:0
\move(97 111)
\lvec(99 111)
\lvec(99 112)
\lvec(97 112)
\ifill f:0
\move(100 111)
\lvec(101 111)
\lvec(101 112)
\lvec(100 112)
\ifill f:0
\move(102 111)
\lvec(107 111)
\lvec(107 112)
\lvec(102 112)
\ifill f:0
\move(108 111)
\lvec(120 111)
\lvec(120 112)
\lvec(108 112)
\ifill f:0
\move(121 111)
\lvec(122 111)
\lvec(122 112)
\lvec(121 112)
\ifill f:0
\move(123 111)
\lvec(128 111)
\lvec(128 112)
\lvec(123 112)
\ifill f:0
\move(133 111)
\lvec(136 111)
\lvec(136 112)
\lvec(133 112)
\ifill f:0
\move(137 111)
\lvec(143 111)
\lvec(143 112)
\lvec(137 112)
\ifill f:0
\move(144 111)
\lvec(145 111)
\lvec(145 112)
\lvec(144 112)
\ifill f:0
\move(146 111)
\lvec(148 111)
\lvec(148 112)
\lvec(146 112)
\ifill f:0
\move(149 111)
\lvec(160 111)
\lvec(160 112)
\lvec(149 112)
\ifill f:0
\move(161 111)
\lvec(170 111)
\lvec(170 112)
\lvec(161 112)
\ifill f:0
\move(171 111)
\lvec(172 111)
\lvec(172 112)
\lvec(171 112)
\ifill f:0
\move(173 111)
\lvec(176 111)
\lvec(176 112)
\lvec(173 112)
\ifill f:0
\move(178 111)
\lvec(184 111)
\lvec(184 112)
\lvec(178 112)
\ifill f:0
\move(185 111)
\lvec(197 111)
\lvec(197 112)
\lvec(185 112)
\ifill f:0
\move(198 111)
\lvec(217 111)
\lvec(217 112)
\lvec(198 112)
\ifill f:0
\move(218 111)
\lvec(226 111)
\lvec(226 112)
\lvec(218 112)
\ifill f:0
\move(227 111)
\lvec(229 111)
\lvec(229 112)
\lvec(227 112)
\ifill f:0
\move(230 111)
\lvec(235 111)
\lvec(235 112)
\lvec(230 112)
\ifill f:0
\move(236 111)
\lvec(240 111)
\lvec(240 112)
\lvec(236 112)
\ifill f:0
\move(241 111)
\lvec(249 111)
\lvec(249 112)
\lvec(241 112)
\ifill f:0
\move(250 111)
\lvec(251 111)
\lvec(251 112)
\lvec(250 112)
\ifill f:0
\move(252 111)
\lvec(253 111)
\lvec(253 112)
\lvec(252 112)
\ifill f:0
\move(254 111)
\lvec(255 111)
\lvec(255 112)
\lvec(254 112)
\ifill f:0
\move(256 111)
\lvec(257 111)
\lvec(257 112)
\lvec(256 112)
\ifill f:0
\move(258 111)
\lvec(259 111)
\lvec(259 112)
\lvec(258 112)
\ifill f:0
\move(260 111)
\lvec(266 111)
\lvec(266 112)
\lvec(260 112)
\ifill f:0
\move(267 111)
\lvec(279 111)
\lvec(279 112)
\lvec(267 112)
\ifill f:0
\move(280 111)
\lvec(282 111)
\lvec(282 112)
\lvec(280 112)
\ifill f:0
\move(283 111)
\lvec(290 111)
\lvec(290 112)
\lvec(283 112)
\ifill f:0
\move(291 111)
\lvec(292 111)
\lvec(292 112)
\lvec(291 112)
\ifill f:0
\move(293 111)
\lvec(306 111)
\lvec(306 112)
\lvec(293 112)
\ifill f:0
\move(307 111)
\lvec(311 111)
\lvec(311 112)
\lvec(307 112)
\ifill f:0
\move(313 111)
\lvec(319 111)
\lvec(319 112)
\lvec(313 112)
\ifill f:0
\move(320 111)
\lvec(325 111)
\lvec(325 112)
\lvec(320 112)
\ifill f:0
\move(326 111)
\lvec(329 111)
\lvec(329 112)
\lvec(326 112)
\ifill f:0
\move(331 111)
\lvec(362 111)
\lvec(362 112)
\lvec(331 112)
\ifill f:0
\move(364 111)
\lvec(387 111)
\lvec(387 112)
\lvec(364 112)
\ifill f:0
\move(388 111)
\lvec(401 111)
\lvec(401 112)
\lvec(388 112)
\ifill f:0
\move(402 111)
\lvec(422 111)
\lvec(422 112)
\lvec(402 112)
\ifill f:0
\move(423 111)
\lvec(432 111)
\lvec(432 112)
\lvec(423 112)
\ifill f:0
\move(433 111)
\lvec(435 111)
\lvec(435 112)
\lvec(433 112)
\ifill f:0
\move(436 111)
\lvec(438 111)
\lvec(438 112)
\lvec(436 112)
\ifill f:0
\move(439 111)
\lvec(442 111)
\lvec(442 112)
\lvec(439 112)
\ifill f:0
\move(443 111)
\lvec(444 111)
\lvec(444 112)
\lvec(443 112)
\ifill f:0
\move(445 111)
\lvec(451 111)
\lvec(451 112)
\lvec(445 112)
\ifill f:0
\move(15 112)
\lvec(17 112)
\lvec(17 113)
\lvec(15 113)
\ifill f:0
\move(20 112)
\lvec(21 112)
\lvec(21 113)
\lvec(20 113)
\ifill f:0
\move(23 112)
\lvec(24 112)
\lvec(24 113)
\lvec(23 113)
\ifill f:0
\move(25 112)
\lvec(26 112)
\lvec(26 113)
\lvec(25 113)
\ifill f:0
\move(36 112)
\lvec(37 112)
\lvec(37 113)
\lvec(36 113)
\ifill f:0
\move(39 112)
\lvec(41 112)
\lvec(41 113)
\lvec(39 113)
\ifill f:0
\move(42 112)
\lvec(46 112)
\lvec(46 113)
\lvec(42 113)
\ifill f:0
\move(49 112)
\lvec(50 112)
\lvec(50 113)
\lvec(49 113)
\ifill f:0
\move(54 112)
\lvec(55 112)
\lvec(55 113)
\lvec(54 113)
\ifill f:0
\move(61 112)
\lvec(63 112)
\lvec(63 113)
\lvec(61 113)
\ifill f:0
\move(64 112)
\lvec(65 112)
\lvec(65 113)
\lvec(64 113)
\ifill f:0
\move(66 112)
\lvec(69 112)
\lvec(69 113)
\lvec(66 113)
\ifill f:0
\move(70 112)
\lvec(73 112)
\lvec(73 113)
\lvec(70 113)
\ifill f:0
\move(74 112)
\lvec(75 112)
\lvec(75 113)
\lvec(74 113)
\ifill f:0
\move(76 112)
\lvec(77 112)
\lvec(77 113)
\lvec(76 113)
\ifill f:0
\move(78 112)
\lvec(80 112)
\lvec(80 113)
\lvec(78 113)
\ifill f:0
\move(81 112)
\lvec(82 112)
\lvec(82 113)
\lvec(81 113)
\ifill f:0
\move(83 112)
\lvec(84 112)
\lvec(84 113)
\lvec(83 113)
\ifill f:0
\move(92 112)
\lvec(93 112)
\lvec(93 113)
\lvec(92 113)
\ifill f:0
\move(95 112)
\lvec(96 112)
\lvec(96 113)
\lvec(95 113)
\ifill f:0
\move(97 112)
\lvec(99 112)
\lvec(99 113)
\lvec(97 113)
\ifill f:0
\move(100 112)
\lvec(101 112)
\lvec(101 113)
\lvec(100 113)
\ifill f:0
\move(102 112)
\lvec(106 112)
\lvec(106 113)
\lvec(102 113)
\ifill f:0
\move(107 112)
\lvec(111 112)
\lvec(111 113)
\lvec(107 113)
\ifill f:0
\move(112 112)
\lvec(115 112)
\lvec(115 113)
\lvec(112 113)
\ifill f:0
\move(116 112)
\lvec(120 112)
\lvec(120 113)
\lvec(116 113)
\ifill f:0
\move(121 112)
\lvec(122 112)
\lvec(122 113)
\lvec(121 113)
\ifill f:0
\move(123 112)
\lvec(133 112)
\lvec(133 113)
\lvec(123 113)
\ifill f:0
\move(135 112)
\lvec(139 112)
\lvec(139 113)
\lvec(135 113)
\ifill f:0
\move(140 112)
\lvec(143 112)
\lvec(143 113)
\lvec(140 113)
\ifill f:0
\move(144 112)
\lvec(145 112)
\lvec(145 113)
\lvec(144 113)
\ifill f:0
\move(146 112)
\lvec(151 112)
\lvec(151 113)
\lvec(146 113)
\ifill f:0
\move(152 112)
\lvec(161 112)
\lvec(161 113)
\lvec(152 113)
\ifill f:0
\move(162 112)
\lvec(163 112)
\lvec(163 113)
\lvec(162 113)
\ifill f:0
\move(164 112)
\lvec(166 112)
\lvec(166 113)
\lvec(164 113)
\ifill f:0
\move(167 112)
\lvec(170 112)
\lvec(170 113)
\lvec(167 113)
\ifill f:0
\move(171 112)
\lvec(172 112)
\lvec(172 113)
\lvec(171 113)
\ifill f:0
\move(173 112)
\lvec(175 112)
\lvec(175 113)
\lvec(173 113)
\ifill f:0
\move(177 112)
\lvec(180 112)
\lvec(180 113)
\lvec(177 113)
\ifill f:0
\move(181 112)
\lvec(186 112)
\lvec(186 113)
\lvec(181 113)
\ifill f:0
\move(188 112)
\lvec(197 112)
\lvec(197 113)
\lvec(188 113)
\ifill f:0
\move(198 112)
\lvec(208 112)
\lvec(208 113)
\lvec(198 113)
\ifill f:0
\move(209 112)
\lvec(221 112)
\lvec(221 113)
\lvec(209 113)
\ifill f:0
\move(222 112)
\lvec(226 112)
\lvec(226 113)
\lvec(222 113)
\ifill f:0
\move(227 112)
\lvec(229 112)
\lvec(229 113)
\lvec(227 113)
\ifill f:0
\move(230 112)
\lvec(241 112)
\lvec(241 113)
\lvec(230 113)
\ifill f:0
\move(242 112)
\lvec(246 112)
\lvec(246 113)
\lvec(242 113)
\ifill f:0
\move(247 112)
\lvec(251 112)
\lvec(251 113)
\lvec(247 113)
\ifill f:0
\move(252 112)
\lvec(253 112)
\lvec(253 113)
\lvec(252 113)
\ifill f:0
\move(254 112)
\lvec(255 112)
\lvec(255 113)
\lvec(254 113)
\ifill f:0
\move(256 112)
\lvec(257 112)
\lvec(257 113)
\lvec(256 113)
\ifill f:0
\move(258 112)
\lvec(259 112)
\lvec(259 113)
\lvec(258 113)
\ifill f:0
\move(260 112)
\lvec(261 112)
\lvec(261 113)
\lvec(260 113)
\ifill f:0
\move(262 112)
\lvec(277 112)
\lvec(277 113)
\lvec(262 113)
\ifill f:0
\move(278 112)
\lvec(280 112)
\lvec(280 113)
\lvec(278 113)
\ifill f:0
\move(281 112)
\lvec(290 112)
\lvec(290 113)
\lvec(281 113)
\ifill f:0
\move(291 112)
\lvec(292 112)
\lvec(292 113)
\lvec(291 113)
\ifill f:0
\move(293 112)
\lvec(299 112)
\lvec(299 113)
\lvec(293 113)
\ifill f:0
\move(300 112)
\lvec(304 112)
\lvec(304 113)
\lvec(300 113)
\ifill f:0
\move(305 112)
\lvec(314 112)
\lvec(314 113)
\lvec(305 113)
\ifill f:0
\move(315 112)
\lvec(319 112)
\lvec(319 113)
\lvec(315 113)
\ifill f:0
\move(321 112)
\lvec(325 112)
\lvec(325 113)
\lvec(321 113)
\ifill f:0
\move(326 112)
\lvec(328 112)
\lvec(328 113)
\lvec(326 113)
\ifill f:0
\move(329 112)
\lvec(341 112)
\lvec(341 113)
\lvec(329 113)
\ifill f:0
\move(343 112)
\lvec(362 112)
\lvec(362 113)
\lvec(343 113)
\ifill f:0
\move(366 112)
\lvec(376 112)
\lvec(376 113)
\lvec(366 113)
\ifill f:0
\move(377 112)
\lvec(392 112)
\lvec(392 113)
\lvec(377 113)
\ifill f:0
\move(393 112)
\lvec(401 112)
\lvec(401 113)
\lvec(393 113)
\ifill f:0
\move(402 112)
\lvec(412 112)
\lvec(412 113)
\lvec(402 113)
\ifill f:0
\move(413 112)
\lvec(424 112)
\lvec(424 113)
\lvec(413 113)
\ifill f:0
\move(425 112)
\lvec(442 112)
\lvec(442 113)
\lvec(425 113)
\ifill f:0
\move(443 112)
\lvec(444 112)
\lvec(444 113)
\lvec(443 113)
\ifill f:0
\move(445 112)
\lvec(447 112)
\lvec(447 113)
\lvec(445 113)
\ifill f:0
\move(448 112)
\lvec(450 112)
\lvec(450 113)
\lvec(448 113)
\ifill f:0
\move(15 113)
\lvec(17 113)
\lvec(17 114)
\lvec(15 114)
\ifill f:0
\move(18 113)
\lvec(21 113)
\lvec(21 114)
\lvec(18 114)
\ifill f:0
\move(24 113)
\lvec(26 113)
\lvec(26 114)
\lvec(24 114)
\ifill f:0
\move(36 113)
\lvec(37 113)
\lvec(37 114)
\lvec(36 114)
\ifill f:0
\move(38 113)
\lvec(39 113)
\lvec(39 114)
\lvec(38 114)
\ifill f:0
\move(40 113)
\lvec(43 113)
\lvec(43 114)
\lvec(40 114)
\ifill f:0
\move(44 113)
\lvec(45 113)
\lvec(45 114)
\lvec(44 114)
\ifill f:0
\move(47 113)
\lvec(48 113)
\lvec(48 114)
\lvec(47 114)
\ifill f:0
\move(49 113)
\lvec(50 113)
\lvec(50 114)
\lvec(49 114)
\ifill f:0
\move(51 113)
\lvec(52 113)
\lvec(52 114)
\lvec(51 114)
\ifill f:0
\move(57 113)
\lvec(58 113)
\lvec(58 114)
\lvec(57 114)
\ifill f:0
\move(59 113)
\lvec(60 113)
\lvec(60 114)
\lvec(59 114)
\ifill f:0
\move(62 113)
\lvec(65 113)
\lvec(65 114)
\lvec(62 114)
\ifill f:0
\move(66 113)
\lvec(71 113)
\lvec(71 114)
\lvec(66 114)
\ifill f:0
\move(72 113)
\lvec(74 113)
\lvec(74 114)
\lvec(72 114)
\ifill f:0
\move(75 113)
\lvec(76 113)
\lvec(76 114)
\lvec(75 114)
\ifill f:0
\move(77 113)
\lvec(78 113)
\lvec(78 114)
\lvec(77 114)
\ifill f:0
\move(79 113)
\lvec(82 113)
\lvec(82 114)
\lvec(79 114)
\ifill f:0
\move(83 113)
\lvec(84 113)
\lvec(84 114)
\lvec(83 114)
\ifill f:0
\move(86 113)
\lvec(93 113)
\lvec(93 114)
\lvec(86 114)
\ifill f:0
\move(97 113)
\lvec(99 113)
\lvec(99 114)
\lvec(97 114)
\ifill f:0
\move(100 113)
\lvec(101 113)
\lvec(101 114)
\lvec(100 114)
\ifill f:0
\move(102 113)
\lvec(104 113)
\lvec(104 114)
\lvec(102 114)
\ifill f:0
\move(105 113)
\lvec(116 113)
\lvec(116 114)
\lvec(105 114)
\ifill f:0
\move(117 113)
\lvec(122 113)
\lvec(122 114)
\lvec(117 114)
\ifill f:0
\move(123 113)
\lvec(127 113)
\lvec(127 114)
\lvec(123 114)
\ifill f:0
\move(128 113)
\lvec(130 113)
\lvec(130 114)
\lvec(128 114)
\ifill f:0
\move(133 113)
\lvec(137 113)
\lvec(137 114)
\lvec(133 114)
\ifill f:0
\move(139 113)
\lvec(142 113)
\lvec(142 114)
\lvec(139 114)
\ifill f:0
\move(144 113)
\lvec(145 113)
\lvec(145 114)
\lvec(144 114)
\ifill f:0
\move(146 113)
\lvec(149 113)
\lvec(149 114)
\lvec(146 114)
\ifill f:0
\move(150 113)
\lvec(154 113)
\lvec(154 114)
\lvec(150 114)
\ifill f:0
\move(155 113)
\lvec(164 113)
\lvec(164 114)
\lvec(155 114)
\ifill f:0
\move(165 113)
\lvec(170 113)
\lvec(170 114)
\lvec(165 114)
\ifill f:0
\move(171 113)
\lvec(174 113)
\lvec(174 114)
\lvec(171 114)
\ifill f:0
\move(176 113)
\lvec(183 113)
\lvec(183 114)
\lvec(176 114)
\ifill f:0
\move(184 113)
\lvec(189 113)
\lvec(189 114)
\lvec(184 114)
\ifill f:0
\move(191 113)
\lvec(197 113)
\lvec(197 114)
\lvec(191 114)
\ifill f:0
\move(198 113)
\lvec(220 113)
\lvec(220 114)
\lvec(198 114)
\ifill f:0
\move(221 113)
\lvec(226 113)
\lvec(226 114)
\lvec(221 114)
\ifill f:0
\move(227 113)
\lvec(233 113)
\lvec(233 114)
\lvec(227 114)
\ifill f:0
\move(234 113)
\lvec(240 113)
\lvec(240 114)
\lvec(234 114)
\ifill f:0
\move(241 113)
\lvec(248 113)
\lvec(248 114)
\lvec(241 114)
\ifill f:0
\move(249 113)
\lvec(255 113)
\lvec(255 114)
\lvec(249 114)
\ifill f:0
\move(256 113)
\lvec(257 113)
\lvec(257 114)
\lvec(256 114)
\ifill f:0
\move(258 113)
\lvec(271 113)
\lvec(271 114)
\lvec(258 114)
\ifill f:0
\move(272 113)
\lvec(283 113)
\lvec(283 114)
\lvec(272 114)
\ifill f:0
\move(284 113)
\lvec(286 113)
\lvec(286 114)
\lvec(284 114)
\ifill f:0
\move(287 113)
\lvec(290 113)
\lvec(290 114)
\lvec(287 114)
\ifill f:0
\move(291 113)
\lvec(292 113)
\lvec(292 114)
\lvec(291 114)
\ifill f:0
\move(293 113)
\lvec(295 113)
\lvec(295 114)
\lvec(293 114)
\ifill f:0
\move(296 113)
\lvec(298 113)
\lvec(298 114)
\lvec(296 114)
\ifill f:0
\move(299 113)
\lvec(302 113)
\lvec(302 114)
\lvec(299 114)
\ifill f:0
\move(303 113)
\lvec(306 113)
\lvec(306 114)
\lvec(303 114)
\ifill f:0
\move(307 113)
\lvec(310 113)
\lvec(310 114)
\lvec(307 114)
\ifill f:0
\move(311 113)
\lvec(315 113)
\lvec(315 114)
\lvec(311 114)
\ifill f:0
\move(316 113)
\lvec(320 113)
\lvec(320 114)
\lvec(316 114)
\ifill f:0
\move(321 113)
\lvec(325 113)
\lvec(325 114)
\lvec(321 114)
\ifill f:0
\move(326 113)
\lvec(335 113)
\lvec(335 114)
\lvec(326 114)
\ifill f:0
\move(337 113)
\lvec(352 113)
\lvec(352 114)
\lvec(337 114)
\ifill f:0
\move(353 113)
\lvec(355 113)
\lvec(355 114)
\lvec(353 114)
\ifill f:0
\move(358 113)
\lvec(359 113)
\lvec(359 114)
\lvec(358 114)
\ifill f:0
\move(361 113)
\lvec(362 113)
\lvec(362 114)
\lvec(361 114)
\ifill f:0
\move(366 113)
\lvec(397 113)
\lvec(397 114)
\lvec(366 114)
\ifill f:0
\move(398 113)
\lvec(401 113)
\lvec(401 114)
\lvec(398 114)
\ifill f:0
\move(402 113)
\lvec(403 113)
\lvec(403 114)
\lvec(402 114)
\ifill f:0
\move(404 113)
\lvec(418 113)
\lvec(418 114)
\lvec(404 114)
\ifill f:0
\move(419 113)
\lvec(434 113)
\lvec(434 114)
\lvec(419 114)
\ifill f:0
\move(435 113)
\lvec(442 113)
\lvec(442 114)
\lvec(435 114)
\ifill f:0
\move(443 113)
\lvec(445 113)
\lvec(445 114)
\lvec(443 114)
\ifill f:0
\move(446 113)
\lvec(451 113)
\lvec(451 114)
\lvec(446 114)
\ifill f:0
\move(16 114)
\lvec(17 114)
\lvec(17 115)
\lvec(16 115)
\ifill f:0
\move(20 114)
\lvec(22 114)
\lvec(22 115)
\lvec(20 115)
\ifill f:0
\move(24 114)
\lvec(26 114)
\lvec(26 115)
\lvec(24 115)
\ifill f:0
\move(36 114)
\lvec(37 114)
\lvec(37 115)
\lvec(36 115)
\ifill f:0
\move(38 114)
\lvec(39 114)
\lvec(39 115)
\lvec(38 115)
\ifill f:0
\move(40 114)
\lvec(41 114)
\lvec(41 115)
\lvec(40 115)
\ifill f:0
\move(42 114)
\lvec(45 114)
\lvec(45 115)
\lvec(42 115)
\ifill f:0
\move(46 114)
\lvec(48 114)
\lvec(48 115)
\lvec(46 115)
\ifill f:0
\move(49 114)
\lvec(50 114)
\lvec(50 115)
\lvec(49 115)
\ifill f:0
\move(51 114)
\lvec(53 114)
\lvec(53 115)
\lvec(51 115)
\ifill f:0
\move(60 114)
\lvec(61 114)
\lvec(61 115)
\lvec(60 115)
\ifill f:0
\move(62 114)
\lvec(65 114)
\lvec(65 115)
\lvec(62 115)
\ifill f:0
\move(66 114)
\lvec(73 114)
\lvec(73 115)
\lvec(66 115)
\ifill f:0
\move(75 114)
\lvec(78 114)
\lvec(78 115)
\lvec(75 115)
\ifill f:0
\move(79 114)
\lvec(82 114)
\lvec(82 115)
\lvec(79 115)
\ifill f:0
\move(85 114)
\lvec(87 114)
\lvec(87 115)
\lvec(85 115)
\ifill f:0
\move(88 114)
\lvec(90 114)
\lvec(90 115)
\lvec(88 115)
\ifill f:0
\move(96 114)
\lvec(98 114)
\lvec(98 115)
\lvec(96 115)
\ifill f:0
\move(100 114)
\lvec(101 114)
\lvec(101 115)
\lvec(100 115)
\ifill f:0
\move(102 114)
\lvec(107 114)
\lvec(107 115)
\lvec(102 115)
\ifill f:0
\move(108 114)
\lvec(109 114)
\lvec(109 115)
\lvec(108 115)
\ifill f:0
\move(110 114)
\lvec(111 114)
\lvec(111 115)
\lvec(110 115)
\ifill f:0
\move(112 114)
\lvec(117 114)
\lvec(117 115)
\lvec(112 115)
\ifill f:0
\move(118 114)
\lvec(122 114)
\lvec(122 115)
\lvec(118 115)
\ifill f:0
\move(123 114)
\lvec(125 114)
\lvec(125 115)
\lvec(123 115)
\ifill f:0
\move(131 114)
\lvec(132 114)
\lvec(132 115)
\lvec(131 115)
\ifill f:0
\move(137 114)
\lvec(138 114)
\lvec(138 115)
\lvec(137 115)
\ifill f:0
\move(139 114)
\lvec(142 114)
\lvec(142 115)
\lvec(139 115)
\ifill f:0
\move(144 114)
\lvec(145 114)
\lvec(145 115)
\lvec(144 115)
\ifill f:0
\move(147 114)
\lvec(170 114)
\lvec(170 115)
\lvec(147 115)
\ifill f:0
\move(172 114)
\lvec(174 114)
\lvec(174 115)
\lvec(172 115)
\ifill f:0
\move(175 114)
\lvec(177 114)
\lvec(177 115)
\lvec(175 115)
\ifill f:0
\move(178 114)
\lvec(185 114)
\lvec(185 115)
\lvec(178 115)
\ifill f:0
\move(186 114)
\lvec(191 114)
\lvec(191 115)
\lvec(186 115)
\ifill f:0
\move(192 114)
\lvec(197 114)
\lvec(197 115)
\lvec(192 115)
\ifill f:0
\move(198 114)
\lvec(201 114)
\lvec(201 115)
\lvec(198 115)
\ifill f:0
\move(202 114)
\lvec(219 114)
\lvec(219 115)
\lvec(202 115)
\ifill f:0
\move(220 114)
\lvec(226 114)
\lvec(226 115)
\lvec(220 115)
\ifill f:0
\move(227 114)
\lvec(230 114)
\lvec(230 115)
\lvec(227 115)
\ifill f:0
\move(231 114)
\lvec(234 114)
\lvec(234 115)
\lvec(231 115)
\ifill f:0
\move(236 114)
\lvec(238 114)
\lvec(238 115)
\lvec(236 115)
\ifill f:0
\move(239 114)
\lvec(242 114)
\lvec(242 115)
\lvec(239 115)
\ifill f:0
\move(243 114)
\lvec(247 114)
\lvec(247 115)
\lvec(243 115)
\ifill f:0
\move(248 114)
\lvec(250 114)
\lvec(250 115)
\lvec(248 115)
\ifill f:0
\move(251 114)
\lvec(255 114)
\lvec(255 115)
\lvec(251 115)
\ifill f:0
\move(256 114)
\lvec(257 114)
\lvec(257 115)
\lvec(256 115)
\ifill f:0
\move(258 114)
\lvec(266 114)
\lvec(266 115)
\lvec(258 115)
\ifill f:0
\move(267 114)
\lvec(290 114)
\lvec(290 115)
\lvec(267 115)
\ifill f:0
\move(291 114)
\lvec(301 114)
\lvec(301 115)
\lvec(291 115)
\ifill f:0
\move(302 114)
\lvec(308 114)
\lvec(308 115)
\lvec(302 115)
\ifill f:0
\move(309 114)
\lvec(312 114)
\lvec(312 115)
\lvec(309 115)
\ifill f:0
\move(313 114)
\lvec(321 114)
\lvec(321 115)
\lvec(313 115)
\ifill f:0
\move(322 114)
\lvec(325 114)
\lvec(325 115)
\lvec(322 115)
\ifill f:0
\move(326 114)
\lvec(342 114)
\lvec(342 115)
\lvec(326 115)
\ifill f:0
\move(343 114)
\lvec(356 114)
\lvec(356 115)
\lvec(343 115)
\ifill f:0
\move(357 114)
\lvec(359 114)
\lvec(359 115)
\lvec(357 115)
\ifill f:0
\move(361 114)
\lvec(362 114)
\lvec(362 115)
\lvec(361 115)
\ifill f:0
\move(364 114)
\lvec(365 114)
\lvec(365 115)
\lvec(364 115)
\ifill f:0
\move(372 114)
\lvec(375 114)
\lvec(375 115)
\lvec(372 115)
\ifill f:0
\move(376 114)
\lvec(387 114)
\lvec(387 115)
\lvec(376 115)
\ifill f:0
\move(388 114)
\lvec(396 114)
\lvec(396 115)
\lvec(388 115)
\ifill f:0
\move(398 114)
\lvec(401 114)
\lvec(401 115)
\lvec(398 115)
\ifill f:0
\move(402 114)
\lvec(425 114)
\lvec(425 115)
\lvec(402 115)
\ifill f:0
\move(426 114)
\lvec(442 114)
\lvec(442 115)
\lvec(426 115)
\ifill f:0
\move(443 114)
\lvec(445 114)
\lvec(445 115)
\lvec(443 115)
\ifill f:0
\move(446 114)
\lvec(448 114)
\lvec(448 115)
\lvec(446 115)
\ifill f:0
\move(449 114)
\lvec(451 114)
\lvec(451 115)
\lvec(449 115)
\ifill f:0
\move(16 115)
\lvec(17 115)
\lvec(17 116)
\lvec(16 116)
\ifill f:0
\move(20 115)
\lvec(21 115)
\lvec(21 116)
\lvec(20 116)
\ifill f:0
\move(24 115)
\lvec(26 115)
\lvec(26 116)
\lvec(24 116)
\ifill f:0
\move(36 115)
\lvec(37 115)
\lvec(37 116)
\lvec(36 116)
\ifill f:0
\move(38 115)
\lvec(39 115)
\lvec(39 116)
\lvec(38 116)
\ifill f:0
\move(41 115)
\lvec(43 115)
\lvec(43 116)
\lvec(41 116)
\ifill f:0
\move(44 115)
\lvec(45 115)
\lvec(45 116)
\lvec(44 116)
\ifill f:0
\move(47 115)
\lvec(50 115)
\lvec(50 116)
\lvec(47 116)
\ifill f:0
\move(51 115)
\lvec(52 115)
\lvec(52 116)
\lvec(51 116)
\ifill f:0
\move(54 115)
\lvec(56 115)
\lvec(56 116)
\lvec(54 116)
\ifill f:0
\move(60 115)
\lvec(62 115)
\lvec(62 116)
\lvec(60 116)
\ifill f:0
\move(63 115)
\lvec(65 115)
\lvec(65 116)
\lvec(63 116)
\ifill f:0
\move(66 115)
\lvec(67 115)
\lvec(67 116)
\lvec(66 116)
\ifill f:0
\move(69 115)
\lvec(70 115)
\lvec(70 116)
\lvec(69 116)
\ifill f:0
\move(72 115)
\lvec(75 115)
\lvec(75 116)
\lvec(72 116)
\ifill f:0
\move(76 115)
\lvec(77 115)
\lvec(77 116)
\lvec(76 116)
\ifill f:0
\move(78 115)
\lvec(82 115)
\lvec(82 116)
\lvec(78 116)
\ifill f:0
\move(84 115)
\lvec(85 115)
\lvec(85 116)
\lvec(84 116)
\ifill f:0
\move(89 115)
\lvec(93 115)
\lvec(93 116)
\lvec(89 116)
\ifill f:0
\move(95 115)
\lvec(98 115)
\lvec(98 116)
\lvec(95 116)
\ifill f:0
\move(100 115)
\lvec(101 115)
\lvec(101 116)
\lvec(100 116)
\ifill f:0
\move(103 115)
\lvec(106 115)
\lvec(106 116)
\lvec(103 116)
\ifill f:0
\move(107 115)
\lvec(117 115)
\lvec(117 116)
\lvec(107 116)
\ifill f:0
\move(118 115)
\lvec(122 115)
\lvec(122 116)
\lvec(118 116)
\ifill f:0
\move(123 115)
\lvec(124 115)
\lvec(124 116)
\lvec(123 116)
\ifill f:0
\move(127 115)
\lvec(132 115)
\lvec(132 116)
\lvec(127 116)
\ifill f:0
\move(133 115)
\lvec(141 115)
\lvec(141 116)
\lvec(133 116)
\ifill f:0
\move(142 115)
\lvec(145 115)
\lvec(145 116)
\lvec(142 116)
\ifill f:0
\move(146 115)
\lvec(147 115)
\lvec(147 116)
\lvec(146 116)
\ifill f:0
\move(148 115)
\lvec(159 115)
\lvec(159 116)
\lvec(148 116)
\ifill f:0
\move(160 115)
\lvec(170 115)
\lvec(170 116)
\lvec(160 116)
\ifill f:0
\move(172 115)
\lvec(173 115)
\lvec(173 116)
\lvec(172 116)
\ifill f:0
\move(174 115)
\lvec(176 115)
\lvec(176 116)
\lvec(174 116)
\ifill f:0
\move(177 115)
\lvec(179 115)
\lvec(179 116)
\lvec(177 116)
\ifill f:0
\move(180 115)
\lvec(183 115)
\lvec(183 116)
\lvec(180 116)
\ifill f:0
\move(184 115)
\lvec(186 115)
\lvec(186 116)
\lvec(184 116)
\ifill f:0
\move(188 115)
\lvec(191 115)
\lvec(191 116)
\lvec(188 116)
\ifill f:0
\move(193 115)
\lvec(197 115)
\lvec(197 116)
\lvec(193 116)
\ifill f:0
\move(198 115)
\lvec(201 115)
\lvec(201 116)
\lvec(198 116)
\ifill f:0
\move(202 115)
\lvec(216 115)
\lvec(216 116)
\lvec(202 116)
\ifill f:0
\move(218 115)
\lvec(226 115)
\lvec(226 116)
\lvec(218 116)
\ifill f:0
\move(227 115)
\lvec(231 115)
\lvec(231 116)
\lvec(227 116)
\ifill f:0
\move(232 115)
\lvec(243 115)
\lvec(243 116)
\lvec(232 116)
\ifill f:0
\move(244 115)
\lvec(250 115)
\lvec(250 116)
\lvec(244 116)
\ifill f:0
\move(251 115)
\lvec(252 115)
\lvec(252 116)
\lvec(251 116)
\ifill f:0
\move(253 115)
\lvec(255 115)
\lvec(255 116)
\lvec(253 116)
\ifill f:0
\move(256 115)
\lvec(257 115)
\lvec(257 116)
\lvec(256 116)
\ifill f:0
\move(258 115)
\lvec(269 115)
\lvec(269 116)
\lvec(258 116)
\ifill f:0
\move(270 115)
\lvec(271 115)
\lvec(271 116)
\lvec(270 116)
\ifill f:0
\move(272 115)
\lvec(284 115)
\lvec(284 116)
\lvec(272 116)
\ifill f:0
\move(285 115)
\lvec(290 115)
\lvec(290 116)
\lvec(285 116)
\ifill f:0
\move(292 115)
\lvec(297 115)
\lvec(297 116)
\lvec(292 116)
\ifill f:0
\move(298 115)
\lvec(303 115)
\lvec(303 116)
\lvec(298 116)
\ifill f:0
\move(304 115)
\lvec(306 115)
\lvec(306 116)
\lvec(304 116)
\ifill f:0
\move(307 115)
\lvec(313 115)
\lvec(313 116)
\lvec(307 116)
\ifill f:0
\move(314 115)
\lvec(317 115)
\lvec(317 116)
\lvec(314 116)
\ifill f:0
\move(318 115)
\lvec(321 115)
\lvec(321 116)
\lvec(318 116)
\ifill f:0
\move(322 115)
\lvec(325 115)
\lvec(325 116)
\lvec(322 116)
\ifill f:0
\move(326 115)
\lvec(332 115)
\lvec(332 116)
\lvec(326 116)
\ifill f:0
\move(333 115)
\lvec(338 115)
\lvec(338 116)
\lvec(333 116)
\ifill f:0
\move(339 115)
\lvec(346 115)
\lvec(346 116)
\lvec(339 116)
\ifill f:0
\move(347 115)
\lvec(359 115)
\lvec(359 116)
\lvec(347 116)
\ifill f:0
\move(361 115)
\lvec(362 115)
\lvec(362 116)
\lvec(361 116)
\ifill f:0
\move(363 115)
\lvec(380 115)
\lvec(380 116)
\lvec(363 116)
\ifill f:0
\move(384 115)
\lvec(395 115)
\lvec(395 116)
\lvec(384 116)
\ifill f:0
\move(396 115)
\lvec(401 115)
\lvec(401 116)
\lvec(396 116)
\ifill f:0
\move(402 115)
\lvec(404 115)
\lvec(404 116)
\lvec(402 116)
\ifill f:0
\move(405 115)
\lvec(423 115)
\lvec(423 116)
\lvec(405 116)
\ifill f:0
\move(424 115)
\lvec(442 115)
\lvec(442 116)
\lvec(424 116)
\ifill f:0
\move(443 115)
\lvec(445 115)
\lvec(445 116)
\lvec(443 116)
\ifill f:0
\move(446 115)
\lvec(449 115)
\lvec(449 116)
\lvec(446 116)
\ifill f:0
\move(450 115)
\lvec(451 115)
\lvec(451 116)
\lvec(450 116)
\ifill f:0
\move(16 116)
\lvec(17 116)
\lvec(17 117)
\lvec(16 117)
\ifill f:0
\move(19 116)
\lvec(21 116)
\lvec(21 117)
\lvec(19 117)
\ifill f:0
\move(23 116)
\lvec(26 116)
\lvec(26 117)
\lvec(23 117)
\ifill f:0
\move(36 116)
\lvec(37 116)
\lvec(37 117)
\lvec(36 117)
\ifill f:0
\move(38 116)
\lvec(45 116)
\lvec(45 117)
\lvec(38 117)
\ifill f:0
\move(47 116)
\lvec(50 116)
\lvec(50 117)
\lvec(47 117)
\ifill f:0
\move(56 116)
\lvec(58 116)
\lvec(58 117)
\lvec(56 117)
\ifill f:0
\move(59 116)
\lvec(60 116)
\lvec(60 117)
\lvec(59 117)
\ifill f:0
\move(61 116)
\lvec(62 116)
\lvec(62 117)
\lvec(61 117)
\ifill f:0
\move(63 116)
\lvec(65 116)
\lvec(65 117)
\lvec(63 117)
\ifill f:0
\move(67 116)
\lvec(73 116)
\lvec(73 117)
\lvec(67 117)
\ifill f:0
\move(75 116)
\lvec(77 116)
\lvec(77 117)
\lvec(75 117)
\ifill f:0
\move(78 116)
\lvec(79 116)
\lvec(79 117)
\lvec(78 117)
\ifill f:0
\move(80 116)
\lvec(82 116)
\lvec(82 117)
\lvec(80 117)
\ifill f:0
\move(84 116)
\lvec(85 116)
\lvec(85 117)
\lvec(84 117)
\ifill f:0
\move(87 116)
\lvec(90 116)
\lvec(90 117)
\lvec(87 117)
\ifill f:0
\move(97 116)
\lvec(98 116)
\lvec(98 117)
\lvec(97 117)
\ifill f:0
\move(100 116)
\lvec(101 116)
\lvec(101 117)
\lvec(100 117)
\ifill f:0
\move(102 116)
\lvec(103 116)
\lvec(103 117)
\lvec(102 117)
\ifill f:0
\move(104 116)
\lvec(106 116)
\lvec(106 117)
\lvec(104 117)
\ifill f:0
\move(107 116)
\lvec(109 116)
\lvec(109 117)
\lvec(107 117)
\ifill f:0
\move(110 116)
\lvec(111 116)
\lvec(111 117)
\lvec(110 117)
\ifill f:0
\move(112 116)
\lvec(113 116)
\lvec(113 117)
\lvec(112 117)
\ifill f:0
\move(114 116)
\lvec(122 116)
\lvec(122 117)
\lvec(114 117)
\ifill f:0
\move(123 116)
\lvec(124 116)
\lvec(124 117)
\lvec(123 117)
\ifill f:0
\move(126 116)
\lvec(130 116)
\lvec(130 117)
\lvec(126 117)
\ifill f:0
\move(133 116)
\lvec(140 116)
\lvec(140 117)
\lvec(133 117)
\ifill f:0
\move(141 116)
\lvec(145 116)
\lvec(145 117)
\lvec(141 117)
\ifill f:0
\move(146 116)
\lvec(147 116)
\lvec(147 117)
\lvec(146 117)
\ifill f:0
\move(148 116)
\lvec(155 116)
\lvec(155 117)
\lvec(148 117)
\ifill f:0
\move(156 116)
\lvec(163 116)
\lvec(163 117)
\lvec(156 117)
\ifill f:0
\move(164 116)
\lvec(167 116)
\lvec(167 117)
\lvec(164 117)
\ifill f:0
\move(168 116)
\lvec(170 116)
\lvec(170 117)
\lvec(168 117)
\ifill f:0
\move(172 116)
\lvec(173 116)
\lvec(173 117)
\lvec(172 117)
\ifill f:0
\move(174 116)
\lvec(192 116)
\lvec(192 117)
\lvec(174 117)
\ifill f:0
\move(194 116)
\lvec(197 116)
\lvec(197 117)
\lvec(194 117)
\ifill f:0
\move(198 116)
\lvec(199 116)
\lvec(199 117)
\lvec(198 117)
\ifill f:0
\move(200 116)
\lvec(226 116)
\lvec(226 117)
\lvec(200 117)
\ifill f:0
\move(227 116)
\lvec(232 116)
\lvec(232 117)
\lvec(227 117)
\ifill f:0
\move(234 116)
\lvec(237 116)
\lvec(237 117)
\lvec(234 117)
\ifill f:0
\move(238 116)
\lvec(242 116)
\lvec(242 117)
\lvec(238 117)
\ifill f:0
\move(243 116)
\lvec(245 116)
\lvec(245 117)
\lvec(243 117)
\ifill f:0
\move(246 116)
\lvec(255 116)
\lvec(255 117)
\lvec(246 117)
\ifill f:0
\move(256 116)
\lvec(257 116)
\lvec(257 117)
\lvec(256 117)
\ifill f:0
\move(258 116)
\lvec(260 116)
\lvec(260 117)
\lvec(258 117)
\ifill f:0
\move(261 116)
\lvec(265 116)
\lvec(265 117)
\lvec(261 117)
\ifill f:0
\move(266 116)
\lvec(274 116)
\lvec(274 117)
\lvec(266 117)
\ifill f:0
\move(275 116)
\lvec(280 116)
\lvec(280 117)
\lvec(275 117)
\ifill f:0
\move(281 116)
\lvec(290 116)
\lvec(290 117)
\lvec(281 117)
\ifill f:0
\move(292 116)
\lvec(299 116)
\lvec(299 117)
\lvec(292 117)
\ifill f:0
\move(300 116)
\lvec(311 116)
\lvec(311 117)
\lvec(300 117)
\ifill f:0
\move(312 116)
\lvec(314 116)
\lvec(314 117)
\lvec(312 117)
\ifill f:0
\move(315 116)
\lvec(318 116)
\lvec(318 117)
\lvec(315 117)
\ifill f:0
\move(319 116)
\lvec(322 116)
\lvec(322 117)
\lvec(319 117)
\ifill f:0
\move(323 116)
\lvec(325 116)
\lvec(325 117)
\lvec(323 117)
\ifill f:0
\move(326 116)
\lvec(336 116)
\lvec(336 117)
\lvec(326 117)
\ifill f:0
\move(337 116)
\lvec(342 116)
\lvec(342 117)
\lvec(337 117)
\ifill f:0
\move(343 116)
\lvec(350 116)
\lvec(350 117)
\lvec(343 117)
\ifill f:0
\move(351 116)
\lvec(360 116)
\lvec(360 117)
\lvec(351 117)
\ifill f:0
\move(361 116)
\lvec(362 116)
\lvec(362 117)
\lvec(361 117)
\ifill f:0
\move(363 116)
\lvec(394 116)
\lvec(394 117)
\lvec(363 117)
\ifill f:0
\move(396 116)
\lvec(401 116)
\lvec(401 117)
\lvec(396 117)
\ifill f:0
\move(402 116)
\lvec(405 116)
\lvec(405 117)
\lvec(402 117)
\ifill f:0
\move(406 116)
\lvec(442 116)
\lvec(442 117)
\lvec(406 117)
\ifill f:0
\move(443 116)
\lvec(450 116)
\lvec(450 117)
\lvec(443 117)
\ifill f:0
\move(16 117)
\lvec(17 117)
\lvec(17 118)
\lvec(16 118)
\ifill f:0
\move(18 117)
\lvec(19 117)
\lvec(19 118)
\lvec(18 118)
\ifill f:0
\move(23 117)
\lvec(26 117)
\lvec(26 118)
\lvec(23 118)
\ifill f:0
\move(36 117)
\lvec(37 117)
\lvec(37 118)
\lvec(36 118)
\ifill f:0
\move(40 117)
\lvec(41 117)
\lvec(41 118)
\lvec(40 118)
\ifill f:0
\move(44 117)
\lvec(45 117)
\lvec(45 118)
\lvec(44 118)
\ifill f:0
\move(47 117)
\lvec(50 117)
\lvec(50 118)
\lvec(47 118)
\ifill f:0
\move(59 117)
\lvec(60 117)
\lvec(60 118)
\lvec(59 118)
\ifill f:0
\move(61 117)
\lvec(65 117)
\lvec(65 118)
\lvec(61 118)
\ifill f:0
\move(66 117)
\lvec(71 117)
\lvec(71 118)
\lvec(66 118)
\ifill f:0
\move(72 117)
\lvec(76 117)
\lvec(76 118)
\lvec(72 118)
\ifill f:0
\move(77 117)
\lvec(79 117)
\lvec(79 118)
\lvec(77 118)
\ifill f:0
\move(80 117)
\lvec(82 117)
\lvec(82 118)
\lvec(80 118)
\ifill f:0
\move(84 117)
\lvec(85 117)
\lvec(85 118)
\lvec(84 118)
\ifill f:0
\move(86 117)
\lvec(87 117)
\lvec(87 118)
\lvec(86 118)
\ifill f:0
\move(90 117)
\lvec(93 117)
\lvec(93 118)
\lvec(90 118)
\ifill f:0
\move(97 117)
\lvec(98 117)
\lvec(98 118)
\lvec(97 118)
\ifill f:0
\move(100 117)
\lvec(101 117)
\lvec(101 118)
\lvec(100 118)
\ifill f:0
\move(102 117)
\lvec(103 117)
\lvec(103 118)
\lvec(102 118)
\ifill f:0
\move(104 117)
\lvec(107 117)
\lvec(107 118)
\lvec(104 118)
\ifill f:0
\move(108 117)
\lvec(112 117)
\lvec(112 118)
\lvec(108 118)
\ifill f:0
\move(113 117)
\lvec(122 117)
\lvec(122 118)
\lvec(113 118)
\ifill f:0
\move(123 117)
\lvec(124 117)
\lvec(124 118)
\lvec(123 118)
\ifill f:0
\move(125 117)
\lvec(127 117)
\lvec(127 118)
\lvec(125 118)
\ifill f:0
\move(129 117)
\lvec(132 117)
\lvec(132 118)
\lvec(129 118)
\ifill f:0
\move(133 117)
\lvec(134 117)
\lvec(134 118)
\lvec(133 118)
\ifill f:0
\move(137 117)
\lvec(138 117)
\lvec(138 118)
\lvec(137 118)
\ifill f:0
\move(139 117)
\lvec(145 117)
\lvec(145 118)
\lvec(139 118)
\ifill f:0
\move(146 117)
\lvec(148 117)
\lvec(148 118)
\lvec(146 118)
\ifill f:0
\move(149 117)
\lvec(153 117)
\lvec(153 118)
\lvec(149 118)
\ifill f:0
\move(154 117)
\lvec(163 117)
\lvec(163 118)
\lvec(154 118)
\ifill f:0
\move(164 117)
\lvec(167 117)
\lvec(167 118)
\lvec(164 118)
\ifill f:0
\move(168 117)
\lvec(170 117)
\lvec(170 118)
\lvec(168 118)
\ifill f:0
\move(172 117)
\lvec(173 117)
\lvec(173 118)
\lvec(172 118)
\ifill f:0
\move(174 117)
\lvec(175 117)
\lvec(175 118)
\lvec(174 118)
\ifill f:0
\move(176 117)
\lvec(180 117)
\lvec(180 118)
\lvec(176 118)
\ifill f:0
\move(181 117)
\lvec(183 117)
\lvec(183 118)
\lvec(181 118)
\ifill f:0
\move(184 117)
\lvec(186 117)
\lvec(186 118)
\lvec(184 118)
\ifill f:0
\move(187 117)
\lvec(189 117)
\lvec(189 118)
\lvec(187 118)
\ifill f:0
\move(190 117)
\lvec(193 117)
\lvec(193 118)
\lvec(190 118)
\ifill f:0
\move(194 117)
\lvec(197 117)
\lvec(197 118)
\lvec(194 118)
\ifill f:0
\move(198 117)
\lvec(226 117)
\lvec(226 118)
\lvec(198 118)
\ifill f:0
\move(228 117)
\lvec(234 117)
\lvec(234 118)
\lvec(228 118)
\ifill f:0
\move(236 117)
\lvec(239 117)
\lvec(239 118)
\lvec(236 118)
\ifill f:0
\move(240 117)
\lvec(244 117)
\lvec(244 118)
\lvec(240 118)
\ifill f:0
\move(245 117)
\lvec(248 117)
\lvec(248 118)
\lvec(245 118)
\ifill f:0
\move(249 117)
\lvec(255 117)
\lvec(255 118)
\lvec(249 118)
\ifill f:0
\move(256 117)
\lvec(257 117)
\lvec(257 118)
\lvec(256 118)
\ifill f:0
\move(258 117)
\lvec(266 117)
\lvec(266 118)
\lvec(258 118)
\ifill f:0
\move(267 117)
\lvec(271 117)
\lvec(271 118)
\lvec(267 118)
\ifill f:0
\move(272 117)
\lvec(275 117)
\lvec(275 118)
\lvec(272 118)
\ifill f:0
\move(276 117)
\lvec(290 117)
\lvec(290 118)
\lvec(276 118)
\ifill f:0
\move(292 117)
\lvec(293 117)
\lvec(293 118)
\lvec(292 118)
\ifill f:0
\move(294 117)
\lvec(298 117)
\lvec(298 118)
\lvec(294 118)
\ifill f:0
\move(299 117)
\lvec(306 117)
\lvec(306 118)
\lvec(299 118)
\ifill f:0
\move(307 117)
\lvec(309 117)
\lvec(309 118)
\lvec(307 118)
\ifill f:0
\move(310 117)
\lvec(315 117)
\lvec(315 118)
\lvec(310 118)
\ifill f:0
\move(316 117)
\lvec(322 117)
\lvec(322 118)
\lvec(316 118)
\ifill f:0
\move(323 117)
\lvec(325 117)
\lvec(325 118)
\lvec(323 118)
\ifill f:0
\move(327 117)
\lvec(330 117)
\lvec(330 118)
\lvec(327 118)
\ifill f:0
\move(331 117)
\lvec(352 117)
\lvec(352 118)
\lvec(331 118)
\ifill f:0
\move(353 117)
\lvec(362 117)
\lvec(362 118)
\lvec(353 118)
\ifill f:0
\move(363 117)
\lvec(377 117)
\lvec(377 118)
\lvec(363 118)
\ifill f:0
\move(390 117)
\lvec(391 117)
\lvec(391 118)
\lvec(390 118)
\ifill f:0
\move(393 117)
\lvec(401 117)
\lvec(401 118)
\lvec(393 118)
\ifill f:0
\move(402 117)
\lvec(424 117)
\lvec(424 118)
\lvec(402 118)
\ifill f:0
\move(425 117)
\lvec(436 117)
\lvec(436 118)
\lvec(425 118)
\ifill f:0
\move(437 117)
\lvec(442 117)
\lvec(442 118)
\lvec(437 118)
\ifill f:0
\move(443 117)
\lvec(446 117)
\lvec(446 118)
\lvec(443 118)
\ifill f:0
\move(447 117)
\lvec(451 117)
\lvec(451 118)
\lvec(447 118)
\ifill f:0
\move(16 118)
\lvec(17 118)
\lvec(17 119)
\lvec(16 119)
\ifill f:0
\move(19 118)
\lvec(21 118)
\lvec(21 119)
\lvec(19 119)
\ifill f:0
\move(24 118)
\lvec(26 118)
\lvec(26 119)
\lvec(24 119)
\ifill f:0
\move(36 118)
\lvec(37 118)
\lvec(37 119)
\lvec(36 119)
\ifill f:0
\move(41 118)
\lvec(45 118)
\lvec(45 119)
\lvec(41 119)
\ifill f:0
\move(48 118)
\lvec(50 118)
\lvec(50 119)
\lvec(48 119)
\ifill f:0
\move(52 118)
\lvec(53 118)
\lvec(53 119)
\lvec(52 119)
\ifill f:0
\move(54 118)
\lvec(55 118)
\lvec(55 119)
\lvec(54 119)
\ifill f:0
\move(56 118)
\lvec(58 118)
\lvec(58 119)
\lvec(56 119)
\ifill f:0
\move(60 118)
\lvec(65 118)
\lvec(65 119)
\lvec(60 119)
\ifill f:0
\move(66 118)
\lvec(70 118)
\lvec(70 119)
\lvec(66 119)
\ifill f:0
\move(76 118)
\lvec(79 118)
\lvec(79 119)
\lvec(76 119)
\ifill f:0
\move(80 118)
\lvec(82 118)
\lvec(82 119)
\lvec(80 119)
\ifill f:0
\move(83 118)
\lvec(84 118)
\lvec(84 119)
\lvec(83 119)
\ifill f:0
\move(85 118)
\lvec(87 118)
\lvec(87 119)
\lvec(85 119)
\ifill f:0
\move(88 118)
\lvec(90 118)
\lvec(90 119)
\lvec(88 119)
\ifill f:0
\move(92 118)
\lvec(93 118)
\lvec(93 119)
\lvec(92 119)
\ifill f:0
\move(97 118)
\lvec(98 118)
\lvec(98 119)
\lvec(97 119)
\ifill f:0
\move(100 118)
\lvec(101 118)
\lvec(101 119)
\lvec(100 119)
\ifill f:0
\move(102 118)
\lvec(104 118)
\lvec(104 119)
\lvec(102 119)
\ifill f:0
\move(105 118)
\lvec(108 118)
\lvec(108 119)
\lvec(105 119)
\ifill f:0
\move(109 118)
\lvec(111 118)
\lvec(111 119)
\lvec(109 119)
\ifill f:0
\move(112 118)
\lvec(116 118)
\lvec(116 119)
\lvec(112 119)
\ifill f:0
\move(117 118)
\lvec(122 118)
\lvec(122 119)
\lvec(117 119)
\ifill f:0
\move(124 118)
\lvec(126 118)
\lvec(126 119)
\lvec(124 119)
\ifill f:0
\move(128 118)
\lvec(130 118)
\lvec(130 119)
\lvec(128 119)
\ifill f:0
\move(133 118)
\lvec(138 118)
\lvec(138 119)
\lvec(133 119)
\ifill f:0
\move(139 118)
\lvec(145 118)
\lvec(145 119)
\lvec(139 119)
\ifill f:0
\move(146 118)
\lvec(149 118)
\lvec(149 119)
\lvec(146 119)
\ifill f:0
\move(150 118)
\lvec(155 118)
\lvec(155 119)
\lvec(150 119)
\ifill f:0
\move(156 118)
\lvec(159 118)
\lvec(159 119)
\lvec(156 119)
\ifill f:0
\move(160 118)
\lvec(167 118)
\lvec(167 119)
\lvec(160 119)
\ifill f:0
\move(168 118)
\lvec(170 118)
\lvec(170 119)
\lvec(168 119)
\ifill f:0
\move(171 118)
\lvec(174 118)
\lvec(174 119)
\lvec(171 119)
\ifill f:0
\move(175 118)
\lvec(176 118)
\lvec(176 119)
\lvec(175 119)
\ifill f:0
\move(177 118)
\lvec(179 118)
\lvec(179 119)
\lvec(177 119)
\ifill f:0
\move(180 118)
\lvec(184 118)
\lvec(184 119)
\lvec(180 119)
\ifill f:0
\move(185 118)
\lvec(186 118)
\lvec(186 119)
\lvec(185 119)
\ifill f:0
\move(188 118)
\lvec(190 118)
\lvec(190 119)
\lvec(188 119)
\ifill f:0
\move(191 118)
\lvec(194 118)
\lvec(194 119)
\lvec(191 119)
\ifill f:0
\move(195 118)
\lvec(197 118)
\lvec(197 119)
\lvec(195 119)
\ifill f:0
\move(198 118)
\lvec(211 118)
\lvec(211 119)
\lvec(198 119)
\ifill f:0
\move(212 118)
\lvec(214 118)
\lvec(214 119)
\lvec(212 119)
\ifill f:0
\move(215 118)
\lvec(226 118)
\lvec(226 119)
\lvec(215 119)
\ifill f:0
\move(229 118)
\lvec(242 118)
\lvec(242 119)
\lvec(229 119)
\ifill f:0
\move(243 118)
\lvec(247 118)
\lvec(247 119)
\lvec(243 119)
\ifill f:0
\move(248 118)
\lvec(251 118)
\lvec(251 119)
\lvec(248 119)
\ifill f:0
\move(252 118)
\lvec(257 118)
\lvec(257 119)
\lvec(252 119)
\ifill f:0
\move(258 118)
\lvec(279 118)
\lvec(279 119)
\lvec(258 119)
\ifill f:0
\move(280 118)
\lvec(281 118)
\lvec(281 119)
\lvec(280 119)
\ifill f:0
\move(282 118)
\lvec(283 118)
\lvec(283 119)
\lvec(282 119)
\ifill f:0
\move(284 118)
\lvec(285 118)
\lvec(285 119)
\lvec(284 119)
\ifill f:0
\move(286 118)
\lvec(287 118)
\lvec(287 119)
\lvec(286 119)
\ifill f:0
\move(288 118)
\lvec(290 118)
\lvec(290 119)
\lvec(288 119)
\ifill f:0
\move(292 118)
\lvec(293 118)
\lvec(293 119)
\lvec(292 119)
\ifill f:0
\move(294 118)
\lvec(295 118)
\lvec(295 119)
\lvec(294 119)
\ifill f:0
\move(296 118)
\lvec(305 118)
\lvec(305 119)
\lvec(296 119)
\ifill f:0
\move(306 118)
\lvec(310 118)
\lvec(310 119)
\lvec(306 119)
\ifill f:0
\move(311 118)
\lvec(313 118)
\lvec(313 119)
\lvec(311 119)
\ifill f:0
\move(314 118)
\lvec(316 118)
\lvec(316 119)
\lvec(314 119)
\ifill f:0
\move(317 118)
\lvec(319 118)
\lvec(319 119)
\lvec(317 119)
\ifill f:0
\move(320 118)
\lvec(322 118)
\lvec(322 119)
\lvec(320 119)
\ifill f:0
\move(323 118)
\lvec(325 118)
\lvec(325 119)
\lvec(323 119)
\ifill f:0
\move(326 118)
\lvec(329 118)
\lvec(329 119)
\lvec(326 119)
\ifill f:0
\move(330 118)
\lvec(342 118)
\lvec(342 119)
\lvec(330 119)
\ifill f:0
\move(343 118)
\lvec(362 118)
\lvec(362 119)
\lvec(343 119)
\ifill f:0
\move(363 118)
\lvec(370 118)
\lvec(370 119)
\lvec(363 119)
\ifill f:0
\move(372 118)
\lvec(401 118)
\lvec(401 119)
\lvec(372 119)
\ifill f:0
\move(402 118)
\lvec(409 118)
\lvec(409 119)
\lvec(402 119)
\ifill f:0
\move(410 118)
\lvec(435 118)
\lvec(435 119)
\lvec(410 119)
\ifill f:0
\move(436 118)
\lvec(442 118)
\lvec(442 119)
\lvec(436 119)
\ifill f:0
\move(443 118)
\lvec(451 118)
\lvec(451 119)
\lvec(443 119)
\ifill f:0
\move(15 119)
\lvec(17 119)
\lvec(17 120)
\lvec(15 120)
\ifill f:0
\move(22 119)
\lvec(23 119)
\lvec(23 120)
\lvec(22 120)
\ifill f:0
\move(25 119)
\lvec(26 119)
\lvec(26 120)
\lvec(25 120)
\ifill f:0
\move(36 119)
\lvec(37 119)
\lvec(37 120)
\lvec(36 120)
\ifill f:0
\move(38 119)
\lvec(39 119)
\lvec(39 120)
\lvec(38 120)
\ifill f:0
\move(40 119)
\lvec(41 119)
\lvec(41 120)
\lvec(40 120)
\ifill f:0
\move(43 119)
\lvec(46 119)
\lvec(46 120)
\lvec(43 120)
\ifill f:0
\move(48 119)
\lvec(50 119)
\lvec(50 120)
\lvec(48 120)
\ifill f:0
\move(52 119)
\lvec(53 119)
\lvec(53 120)
\lvec(52 120)
\ifill f:0
\move(59 119)
\lvec(61 119)
\lvec(61 120)
\lvec(59 120)
\ifill f:0
\move(62 119)
\lvec(63 119)
\lvec(63 120)
\lvec(62 120)
\ifill f:0
\move(64 119)
\lvec(65 119)
\lvec(65 120)
\lvec(64 120)
\ifill f:0
\move(66 119)
\lvec(68 119)
\lvec(68 120)
\lvec(66 120)
\ifill f:0
\move(69 119)
\lvec(74 119)
\lvec(74 120)
\lvec(69 120)
\ifill f:0
\move(75 119)
\lvec(78 119)
\lvec(78 120)
\lvec(75 120)
\ifill f:0
\move(79 119)
\lvec(82 119)
\lvec(82 120)
\lvec(79 120)
\ifill f:0
\move(83 119)
\lvec(84 119)
\lvec(84 120)
\lvec(83 120)
\ifill f:0
\move(87 119)
\lvec(89 119)
\lvec(89 120)
\lvec(87 120)
\ifill f:0
\move(90 119)
\lvec(93 119)
\lvec(93 120)
\lvec(90 120)
\ifill f:0
\move(95 119)
\lvec(98 119)
\lvec(98 120)
\lvec(95 120)
\ifill f:0
\move(100 119)
\lvec(101 119)
\lvec(101 120)
\lvec(100 120)
\ifill f:0
\move(102 119)
\lvec(106 119)
\lvec(106 120)
\lvec(102 120)
\ifill f:0
\move(107 119)
\lvec(115 119)
\lvec(115 120)
\lvec(107 120)
\ifill f:0
\move(116 119)
\lvec(122 119)
\lvec(122 120)
\lvec(116 120)
\ifill f:0
\move(124 119)
\lvec(125 119)
\lvec(125 120)
\lvec(124 120)
\ifill f:0
\move(127 119)
\lvec(129 119)
\lvec(129 120)
\lvec(127 120)
\ifill f:0
\move(131 119)
\lvec(132 119)
\lvec(132 120)
\lvec(131 120)
\ifill f:0
\move(136 119)
\lvec(145 119)
\lvec(145 120)
\lvec(136 120)
\ifill f:0
\move(146 119)
\lvec(151 119)
\lvec(151 120)
\lvec(146 120)
\ifill f:0
\move(152 119)
\lvec(170 119)
\lvec(170 120)
\lvec(152 120)
\ifill f:0
\move(171 119)
\lvec(172 119)
\lvec(172 120)
\lvec(171 120)
\ifill f:0
\move(173 119)
\lvec(174 119)
\lvec(174 120)
\lvec(173 120)
\ifill f:0
\move(175 119)
\lvec(176 119)
\lvec(176 120)
\lvec(175 120)
\ifill f:0
\move(177 119)
\lvec(180 119)
\lvec(180 120)
\lvec(177 120)
\ifill f:0
\move(181 119)
\lvec(185 119)
\lvec(185 120)
\lvec(181 120)
\ifill f:0
\move(186 119)
\lvec(188 119)
\lvec(188 120)
\lvec(186 120)
\ifill f:0
\move(189 119)
\lvec(194 119)
\lvec(194 120)
\lvec(189 120)
\ifill f:0
\move(195 119)
\lvec(197 119)
\lvec(197 120)
\lvec(195 120)
\ifill f:0
\move(199 119)
\lvec(202 119)
\lvec(202 120)
\lvec(199 120)
\ifill f:0
\move(203 119)
\lvec(208 119)
\lvec(208 120)
\lvec(203 120)
\ifill f:0
\move(209 119)
\lvec(219 119)
\lvec(219 120)
\lvec(209 120)
\ifill f:0
\move(223 119)
\lvec(224 119)
\lvec(224 120)
\lvec(223 120)
\ifill f:0
\move(225 119)
\lvec(226 119)
\lvec(226 120)
\lvec(225 120)
\ifill f:0
\move(230 119)
\lvec(239 119)
\lvec(239 120)
\lvec(230 120)
\ifill f:0
\move(240 119)
\lvec(250 119)
\lvec(250 120)
\lvec(240 120)
\ifill f:0
\move(252 119)
\lvec(257 119)
\lvec(257 120)
\lvec(252 120)
\ifill f:0
\move(258 119)
\lvec(273 119)
\lvec(273 120)
\lvec(258 120)
\ifill f:0
\move(274 119)
\lvec(287 119)
\lvec(287 120)
\lvec(274 120)
\ifill f:0
\move(288 119)
\lvec(290 119)
\lvec(290 120)
\lvec(288 120)
\ifill f:0
\move(292 119)
\lvec(293 119)
\lvec(293 120)
\lvec(292 120)
\ifill f:0
\move(294 119)
\lvec(295 119)
\lvec(295 120)
\lvec(294 120)
\ifill f:0
\move(296 119)
\lvec(297 119)
\lvec(297 120)
\lvec(296 120)
\ifill f:0
\move(298 119)
\lvec(299 119)
\lvec(299 120)
\lvec(298 120)
\ifill f:0
\move(300 119)
\lvec(311 119)
\lvec(311 120)
\lvec(300 120)
\ifill f:0
\move(312 119)
\lvec(319 119)
\lvec(319 120)
\lvec(312 120)
\ifill f:0
\move(320 119)
\lvec(322 119)
\lvec(322 120)
\lvec(320 120)
\ifill f:0
\move(323 119)
\lvec(325 119)
\lvec(325 120)
\lvec(323 120)
\ifill f:0
\move(326 119)
\lvec(336 119)
\lvec(336 120)
\lvec(326 120)
\ifill f:0
\move(337 119)
\lvec(340 119)
\lvec(340 120)
\lvec(337 120)
\ifill f:0
\move(341 119)
\lvec(344 119)
\lvec(344 120)
\lvec(341 120)
\ifill f:0
\move(345 119)
\lvec(362 119)
\lvec(362 120)
\lvec(345 120)
\ifill f:0
\move(363 119)
\lvec(368 119)
\lvec(368 120)
\lvec(363 120)
\ifill f:0
\move(369 119)
\lvec(379 119)
\lvec(379 120)
\lvec(369 120)
\ifill f:0
\move(381 119)
\lvec(401 119)
\lvec(401 120)
\lvec(381 120)
\ifill f:0
\move(402 119)
\lvec(426 119)
\lvec(426 120)
\lvec(402 120)
\ifill f:0
\move(427 119)
\lvec(434 119)
\lvec(434 120)
\lvec(427 120)
\ifill f:0
\move(435 119)
\lvec(442 119)
\lvec(442 120)
\lvec(435 120)
\ifill f:0
\move(443 119)
\lvec(447 119)
\lvec(447 120)
\lvec(443 120)
\ifill f:0
\move(448 119)
\lvec(451 119)
\lvec(451 120)
\lvec(448 120)
\ifill f:0
\move(15 120)
\lvec(17 120)
\lvec(17 121)
\lvec(15 121)
\ifill f:0
\move(20 120)
\lvec(21 120)
\lvec(21 121)
\lvec(20 121)
\ifill f:0
\move(25 120)
\lvec(26 120)
\lvec(26 121)
\lvec(25 121)
\ifill f:0
\move(36 120)
\lvec(37 120)
\lvec(37 121)
\lvec(36 121)
\ifill f:0
\move(38 120)
\lvec(45 120)
\lvec(45 121)
\lvec(38 121)
\ifill f:0
\move(47 120)
\lvec(50 120)
\lvec(50 121)
\lvec(47 121)
\ifill f:0
\move(51 120)
\lvec(52 120)
\lvec(52 121)
\lvec(51 121)
\ifill f:0
\move(56 120)
\lvec(58 120)
\lvec(58 121)
\lvec(56 121)
\ifill f:0
\move(59 120)
\lvec(60 120)
\lvec(60 121)
\lvec(59 121)
\ifill f:0
\move(61 120)
\lvec(63 120)
\lvec(63 121)
\lvec(61 121)
\ifill f:0
\move(64 120)
\lvec(65 120)
\lvec(65 121)
\lvec(64 121)
\ifill f:0
\move(66 120)
\lvec(67 120)
\lvec(67 121)
\lvec(66 121)
\ifill f:0
\move(68 120)
\lvec(71 120)
\lvec(71 121)
\lvec(68 121)
\ifill f:0
\move(73 120)
\lvec(76 120)
\lvec(76 121)
\lvec(73 121)
\ifill f:0
\move(79 120)
\lvec(82 120)
\lvec(82 121)
\lvec(79 121)
\ifill f:0
\move(83 120)
\lvec(84 120)
\lvec(84 121)
\lvec(83 121)
\ifill f:0
\move(89 120)
\lvec(90 120)
\lvec(90 121)
\lvec(89 121)
\ifill f:0
\move(92 120)
\lvec(93 120)
\lvec(93 121)
\lvec(92 121)
\ifill f:0
\move(96 120)
\lvec(98 120)
\lvec(98 121)
\lvec(96 121)
\ifill f:0
\move(100 120)
\lvec(101 120)
\lvec(101 121)
\lvec(100 121)
\ifill f:0
\move(102 120)
\lvec(107 120)
\lvec(107 121)
\lvec(102 121)
\ifill f:0
\move(108 120)
\lvec(119 120)
\lvec(119 121)
\lvec(108 121)
\ifill f:0
\move(120 120)
\lvec(122 120)
\lvec(122 121)
\lvec(120 121)
\ifill f:0
\move(124 120)
\lvec(125 120)
\lvec(125 121)
\lvec(124 121)
\ifill f:0
\move(126 120)
\lvec(128 120)
\lvec(128 121)
\lvec(126 121)
\ifill f:0
\move(129 120)
\lvec(131 120)
\lvec(131 121)
\lvec(129 121)
\ifill f:0
\move(133 120)
\lvec(136 120)
\lvec(136 121)
\lvec(133 121)
\ifill f:0
\move(138 120)
\lvec(145 120)
\lvec(145 121)
\lvec(138 121)
\ifill f:0
\move(146 120)
\lvec(153 120)
\lvec(153 121)
\lvec(146 121)
\ifill f:0
\move(154 120)
\lvec(159 120)
\lvec(159 121)
\lvec(154 121)
\ifill f:0
\move(160 120)
\lvec(163 120)
\lvec(163 121)
\lvec(160 121)
\ifill f:0
\move(164 120)
\lvec(170 120)
\lvec(170 121)
\lvec(164 121)
\ifill f:0
\move(171 120)
\lvec(172 120)
\lvec(172 121)
\lvec(171 121)
\ifill f:0
\move(173 120)
\lvec(174 120)
\lvec(174 121)
\lvec(173 121)
\ifill f:0
\move(175 120)
\lvec(186 120)
\lvec(186 121)
\lvec(175 121)
\ifill f:0
\move(187 120)
\lvec(191 120)
\lvec(191 121)
\lvec(187 121)
\ifill f:0
\move(192 120)
\lvec(194 120)
\lvec(194 121)
\lvec(192 121)
\ifill f:0
\move(195 120)
\lvec(197 120)
\lvec(197 121)
\lvec(195 121)
\ifill f:0
\move(198 120)
\lvec(201 120)
\lvec(201 121)
\lvec(198 121)
\ifill f:0
\move(202 120)
\lvec(206 120)
\lvec(206 121)
\lvec(202 121)
\ifill f:0
\move(207 120)
\lvec(211 120)
\lvec(211 121)
\lvec(207 121)
\ifill f:0
\move(213 120)
\lvec(223 120)
\lvec(223 121)
\lvec(213 121)
\ifill f:0
\move(225 120)
\lvec(226 120)
\lvec(226 121)
\lvec(225 121)
\ifill f:0
\move(227 120)
\lvec(228 120)
\lvec(228 121)
\lvec(227 121)
\ifill f:0
\move(236 120)
\lvec(243 120)
\lvec(243 121)
\lvec(236 121)
\ifill f:0
\move(244 120)
\lvec(249 120)
\lvec(249 121)
\lvec(244 121)
\ifill f:0
\move(251 120)
\lvec(254 120)
\lvec(254 121)
\lvec(251 121)
\ifill f:0
\move(255 120)
\lvec(257 120)
\lvec(257 121)
\lvec(255 121)
\ifill f:0
\move(258 120)
\lvec(262 120)
\lvec(262 121)
\lvec(258 121)
\ifill f:0
\move(263 120)
\lvec(266 120)
\lvec(266 121)
\lvec(263 121)
\ifill f:0
\move(267 120)
\lvec(269 120)
\lvec(269 121)
\lvec(267 121)
\ifill f:0
\move(270 120)
\lvec(280 120)
\lvec(280 121)
\lvec(270 121)
\ifill f:0
\move(281 120)
\lvec(287 120)
\lvec(287 121)
\lvec(281 121)
\ifill f:0
\move(288 120)
\lvec(290 120)
\lvec(290 121)
\lvec(288 121)
\ifill f:0
\move(291 120)
\lvec(298 120)
\lvec(298 121)
\lvec(291 121)
\ifill f:0
\move(299 120)
\lvec(305 120)
\lvec(305 121)
\lvec(299 121)
\ifill f:0
\move(306 120)
\lvec(317 120)
\lvec(317 121)
\lvec(306 121)
\ifill f:0
\move(318 120)
\lvec(322 120)
\lvec(322 121)
\lvec(318 121)
\ifill f:0
\move(323 120)
\lvec(325 120)
\lvec(325 121)
\lvec(323 121)
\ifill f:0
\move(326 120)
\lvec(335 120)
\lvec(335 121)
\lvec(326 121)
\ifill f:0
\move(336 120)
\lvec(338 120)
\lvec(338 121)
\lvec(336 121)
\ifill f:0
\move(339 120)
\lvec(342 120)
\lvec(342 121)
\lvec(339 121)
\ifill f:0
\move(343 120)
\lvec(346 120)
\lvec(346 121)
\lvec(343 121)
\ifill f:0
\move(347 120)
\lvec(362 120)
\lvec(362 121)
\lvec(347 121)
\ifill f:0
\move(363 120)
\lvec(367 120)
\lvec(367 121)
\lvec(363 121)
\ifill f:0
\move(368 120)
\lvec(375 120)
\lvec(375 121)
\lvec(368 121)
\ifill f:0
\move(376 120)
\lvec(386 120)
\lvec(386 121)
\lvec(376 121)
\ifill f:0
\move(388 120)
\lvec(401 120)
\lvec(401 121)
\lvec(388 121)
\ifill f:0
\move(402 120)
\lvec(420 120)
\lvec(420 121)
\lvec(402 121)
\ifill f:0
\move(421 120)
\lvec(442 120)
\lvec(442 121)
\lvec(421 121)
\ifill f:0
\move(443 120)
\lvec(448 120)
\lvec(448 121)
\lvec(443 121)
\ifill f:0
\move(449 120)
\lvec(451 120)
\lvec(451 121)
\lvec(449 121)
\ifill f:0
\move(15 121)
\lvec(17 121)
\lvec(17 122)
\lvec(15 122)
\ifill f:0
\move(19 121)
\lvec(20 121)
\lvec(20 122)
\lvec(19 122)
\ifill f:0
\move(23 121)
\lvec(24 121)
\lvec(24 122)
\lvec(23 122)
\ifill f:0
\move(25 121)
\lvec(26 121)
\lvec(26 122)
\lvec(25 122)
\ifill f:0
\move(36 121)
\lvec(37 121)
\lvec(37 122)
\lvec(36 122)
\ifill f:0
\move(38 121)
\lvec(39 121)
\lvec(39 122)
\lvec(38 122)
\ifill f:0
\move(41 121)
\lvec(42 121)
\lvec(42 122)
\lvec(41 122)
\ifill f:0
\move(43 121)
\lvec(45 121)
\lvec(45 122)
\lvec(43 122)
\ifill f:0
\move(47 121)
\lvec(50 121)
\lvec(50 122)
\lvec(47 122)
\ifill f:0
\move(51 121)
\lvec(52 121)
\lvec(52 122)
\lvec(51 122)
\ifill f:0
\move(54 121)
\lvec(55 121)
\lvec(55 122)
\lvec(54 122)
\ifill f:0
\move(61 121)
\lvec(63 121)
\lvec(63 122)
\lvec(61 122)
\ifill f:0
\move(64 121)
\lvec(65 121)
\lvec(65 122)
\lvec(64 122)
\ifill f:0
\move(66 121)
\lvec(74 121)
\lvec(74 122)
\lvec(66 122)
\ifill f:0
\move(76 121)
\lvec(77 121)
\lvec(77 122)
\lvec(76 122)
\ifill f:0
\move(78 121)
\lvec(82 121)
\lvec(82 122)
\lvec(78 122)
\ifill f:0
\move(83 121)
\lvec(85 121)
\lvec(85 122)
\lvec(83 122)
\ifill f:0
\move(86 121)
\lvec(87 121)
\lvec(87 122)
\lvec(86 122)
\ifill f:0
\move(88 121)
\lvec(92 121)
\lvec(92 122)
\lvec(88 122)
\ifill f:0
\move(97 121)
\lvec(99 121)
\lvec(99 122)
\lvec(97 122)
\ifill f:0
\move(100 121)
\lvec(101 121)
\lvec(101 122)
\lvec(100 122)
\ifill f:0
\move(102 121)
\lvec(109 121)
\lvec(109 122)
\lvec(102 122)
\ifill f:0
\move(110 121)
\lvec(119 121)
\lvec(119 122)
\lvec(110 122)
\ifill f:0
\move(120 121)
\lvec(122 121)
\lvec(122 122)
\lvec(120 122)
\ifill f:0
\move(124 121)
\lvec(125 121)
\lvec(125 122)
\lvec(124 122)
\ifill f:0
\move(126 121)
\lvec(127 121)
\lvec(127 122)
\lvec(126 122)
\ifill f:0
\move(128 121)
\lvec(130 121)
\lvec(130 122)
\lvec(128 122)
\ifill f:0
\move(131 121)
\lvec(132 121)
\lvec(132 122)
\lvec(131 122)
\ifill f:0
\move(135 121)
\lvec(138 121)
\lvec(138 122)
\lvec(135 122)
\ifill f:0
\move(140 121)
\lvec(145 121)
\lvec(145 122)
\lvec(140 122)
\ifill f:0
\move(146 121)
\lvec(162 121)
\lvec(162 122)
\lvec(146 122)
\ifill f:0
\move(164 121)
\lvec(170 121)
\lvec(170 122)
\lvec(164 122)
\ifill f:0
\move(171 121)
\lvec(172 121)
\lvec(172 122)
\lvec(171 122)
\ifill f:0
\move(174 121)
\lvec(175 121)
\lvec(175 122)
\lvec(174 122)
\ifill f:0
\move(176 121)
\lvec(177 121)
\lvec(177 122)
\lvec(176 122)
\ifill f:0
\move(178 121)
\lvec(186 121)
\lvec(186 122)
\lvec(178 122)
\ifill f:0
\move(188 121)
\lvec(189 121)
\lvec(189 122)
\lvec(188 122)
\ifill f:0
\move(190 121)
\lvec(191 121)
\lvec(191 122)
\lvec(190 122)
\ifill f:0
\move(192 121)
\lvec(194 121)
\lvec(194 122)
\lvec(192 122)
\ifill f:0
\move(195 121)
\lvec(197 121)
\lvec(197 122)
\lvec(195 122)
\ifill f:0
\move(198 121)
\lvec(209 121)
\lvec(209 122)
\lvec(198 122)
\ifill f:0
\move(210 121)
\lvec(223 121)
\lvec(223 122)
\lvec(210 122)
\ifill f:0
\move(225 121)
\lvec(226 121)
\lvec(226 122)
\lvec(225 122)
\ifill f:0
\move(227 121)
\lvec(237 121)
\lvec(237 122)
\lvec(227 122)
\ifill f:0
\move(239 121)
\lvec(248 121)
\lvec(248 122)
\lvec(239 122)
\ifill f:0
\move(249 121)
\lvec(257 121)
\lvec(257 122)
\lvec(249 122)
\ifill f:0
\move(258 121)
\lvec(273 121)
\lvec(273 122)
\lvec(258 122)
\ifill f:0
\move(274 121)
\lvec(276 121)
\lvec(276 122)
\lvec(274 122)
\ifill f:0
\move(277 121)
\lvec(279 121)
\lvec(279 122)
\lvec(277 122)
\ifill f:0
\move(280 121)
\lvec(290 121)
\lvec(290 122)
\lvec(280 122)
\ifill f:0
\move(291 121)
\lvec(298 121)
\lvec(298 122)
\lvec(291 122)
\ifill f:0
\move(299 121)
\lvec(302 121)
\lvec(302 122)
\lvec(299 122)
\ifill f:0
\move(303 121)
\lvec(315 121)
\lvec(315 122)
\lvec(303 122)
\ifill f:0
\move(316 121)
\lvec(320 121)
\lvec(320 122)
\lvec(316 122)
\ifill f:0
\move(321 121)
\lvec(322 121)
\lvec(322 122)
\lvec(321 122)
\ifill f:0
\move(323 121)
\lvec(325 121)
\lvec(325 122)
\lvec(323 122)
\ifill f:0
\move(326 121)
\lvec(328 121)
\lvec(328 122)
\lvec(326 122)
\ifill f:0
\move(329 121)
\lvec(331 121)
\lvec(331 122)
\lvec(329 122)
\ifill f:0
\move(332 121)
\lvec(334 121)
\lvec(334 122)
\lvec(332 122)
\ifill f:0
\move(335 121)
\lvec(337 121)
\lvec(337 122)
\lvec(335 122)
\ifill f:0
\move(338 121)
\lvec(344 121)
\lvec(344 122)
\lvec(338 122)
\ifill f:0
\move(345 121)
\lvec(356 121)
\lvec(356 122)
\lvec(345 122)
\ifill f:0
\move(357 121)
\lvec(362 121)
\lvec(362 122)
\lvec(357 122)
\ifill f:0
\move(363 121)
\lvec(366 121)
\lvec(366 122)
\lvec(363 122)
\ifill f:0
\move(367 121)
\lvec(379 121)
\lvec(379 122)
\lvec(367 122)
\ifill f:0
\move(381 121)
\lvec(391 121)
\lvec(391 122)
\lvec(381 122)
\ifill f:0
\move(392 121)
\lvec(401 121)
\lvec(401 122)
\lvec(392 122)
\ifill f:0
\move(402 121)
\lvec(430 121)
\lvec(430 122)
\lvec(402 122)
\ifill f:0
\move(431 121)
\lvec(442 121)
\lvec(442 122)
\lvec(431 122)
\ifill f:0
\move(443 121)
\lvec(450 121)
\lvec(450 122)
\lvec(443 122)
\ifill f:0
\move(16 122)
\lvec(17 122)
\lvec(17 123)
\lvec(16 123)
\ifill f:0
\move(20 122)
\lvec(21 122)
\lvec(21 123)
\lvec(20 123)
\ifill f:0
\move(22 122)
\lvec(23 122)
\lvec(23 123)
\lvec(22 123)
\ifill f:0
\move(25 122)
\lvec(26 122)
\lvec(26 123)
\lvec(25 123)
\ifill f:0
\move(36 122)
\lvec(37 122)
\lvec(37 123)
\lvec(36 123)
\ifill f:0
\move(38 122)
\lvec(39 122)
\lvec(39 123)
\lvec(38 123)
\ifill f:0
\move(40 122)
\lvec(41 122)
\lvec(41 123)
\lvec(40 123)
\ifill f:0
\move(42 122)
\lvec(43 122)
\lvec(43 123)
\lvec(42 123)
\ifill f:0
\move(44 122)
\lvec(45 122)
\lvec(45 123)
\lvec(44 123)
\ifill f:0
\move(47 122)
\lvec(50 122)
\lvec(50 123)
\lvec(47 123)
\ifill f:0
\move(51 122)
\lvec(52 122)
\lvec(52 123)
\lvec(51 123)
\ifill f:0
\move(56 122)
\lvec(58 122)
\lvec(58 123)
\lvec(56 123)
\ifill f:0
\move(59 122)
\lvec(63 122)
\lvec(63 123)
\lvec(59 123)
\ifill f:0
\move(64 122)
\lvec(65 122)
\lvec(65 123)
\lvec(64 123)
\ifill f:0
\move(66 122)
\lvec(71 122)
\lvec(71 123)
\lvec(66 123)
\ifill f:0
\move(73 122)
\lvec(82 122)
\lvec(82 123)
\lvec(73 123)
\ifill f:0
\move(83 122)
\lvec(85 122)
\lvec(85 123)
\lvec(83 123)
\ifill f:0
\move(87 122)
\lvec(88 122)
\lvec(88 123)
\lvec(87 123)
\ifill f:0
\move(89 122)
\lvec(90 122)
\lvec(90 123)
\lvec(89 123)
\ifill f:0
\move(91 122)
\lvec(93 122)
\lvec(93 123)
\lvec(91 123)
\ifill f:0
\move(98 122)
\lvec(101 122)
\lvec(101 123)
\lvec(98 123)
\ifill f:0
\move(102 122)
\lvec(106 122)
\lvec(106 123)
\lvec(102 123)
\ifill f:0
\move(107 122)
\lvec(116 122)
\lvec(116 123)
\lvec(107 123)
\ifill f:0
\move(117 122)
\lvec(122 122)
\lvec(122 123)
\lvec(117 123)
\ifill f:0
\move(123 122)
\lvec(124 122)
\lvec(124 123)
\lvec(123 123)
\ifill f:0
\move(125 122)
\lvec(127 122)
\lvec(127 123)
\lvec(125 123)
\ifill f:0
\move(128 122)
\lvec(129 122)
\lvec(129 123)
\lvec(128 123)
\ifill f:0
\move(130 122)
\lvec(131 122)
\lvec(131 123)
\lvec(130 123)
\ifill f:0
\move(133 122)
\lvec(135 122)
\lvec(135 123)
\lvec(133 123)
\ifill f:0
\move(136 122)
\lvec(138 122)
\lvec(138 123)
\lvec(136 123)
\ifill f:0
\move(139 122)
\lvec(140 122)
\lvec(140 123)
\lvec(139 123)
\ifill f:0
\move(141 122)
\lvec(145 122)
\lvec(145 123)
\lvec(141 123)
\ifill f:0
\move(146 122)
\lvec(170 122)
\lvec(170 123)
\lvec(146 123)
\ifill f:0
\move(171 122)
\lvec(173 122)
\lvec(173 123)
\lvec(171 123)
\ifill f:0
\move(174 122)
\lvec(176 122)
\lvec(176 123)
\lvec(174 123)
\ifill f:0
\move(177 122)
\lvec(180 122)
\lvec(180 123)
\lvec(177 123)
\ifill f:0
\move(181 122)
\lvec(192 122)
\lvec(192 123)
\lvec(181 123)
\ifill f:0
\move(193 122)
\lvec(194 122)
\lvec(194 123)
\lvec(193 123)
\ifill f:0
\move(195 122)
\lvec(197 122)
\lvec(197 123)
\lvec(195 123)
\ifill f:0
\move(198 122)
\lvec(207 122)
\lvec(207 123)
\lvec(198 123)
\ifill f:0
\move(208 122)
\lvec(211 122)
\lvec(211 123)
\lvec(208 123)
\ifill f:0
\move(212 122)
\lvec(224 122)
\lvec(224 123)
\lvec(212 123)
\ifill f:0
\move(225 122)
\lvec(226 122)
\lvec(226 123)
\lvec(225 123)
\ifill f:0
\move(227 122)
\lvec(244 122)
\lvec(244 123)
\lvec(227 123)
\ifill f:0
\move(245 122)
\lvec(253 122)
\lvec(253 123)
\lvec(245 123)
\ifill f:0
\move(254 122)
\lvec(257 122)
\lvec(257 123)
\lvec(254 123)
\ifill f:0
\move(258 122)
\lvec(259 122)
\lvec(259 123)
\lvec(258 123)
\ifill f:0
\move(260 122)
\lvec(275 122)
\lvec(275 123)
\lvec(260 123)
\ifill f:0
\move(276 122)
\lvec(281 122)
\lvec(281 123)
\lvec(276 123)
\ifill f:0
\move(282 122)
\lvec(284 122)
\lvec(284 123)
\lvec(282 123)
\ifill f:0
\move(285 122)
\lvec(290 122)
\lvec(290 123)
\lvec(285 123)
\ifill f:0
\move(291 122)
\lvec(294 122)
\lvec(294 123)
\lvec(291 123)
\ifill f:0
\move(295 122)
\lvec(305 122)
\lvec(305 123)
\lvec(295 123)
\ifill f:0
\move(306 122)
\lvec(307 122)
\lvec(307 123)
\lvec(306 123)
\ifill f:0
\move(308 122)
\lvec(318 122)
\lvec(318 123)
\lvec(308 123)
\ifill f:0
\move(319 122)
\lvec(320 122)
\lvec(320 123)
\lvec(319 123)
\ifill f:0
\move(321 122)
\lvec(325 122)
\lvec(325 123)
\lvec(321 123)
\ifill f:0
\move(326 122)
\lvec(339 122)
\lvec(339 123)
\lvec(326 123)
\ifill f:0
\move(340 122)
\lvec(349 122)
\lvec(349 123)
\lvec(340 123)
\ifill f:0
\move(350 122)
\lvec(362 122)
\lvec(362 123)
\lvec(350 123)
\ifill f:0
\move(363 122)
\lvec(365 122)
\lvec(365 123)
\lvec(363 123)
\ifill f:0
\move(366 122)
\lvec(370 122)
\lvec(370 123)
\lvec(366 123)
\ifill f:0
\move(372 122)
\lvec(383 122)
\lvec(383 123)
\lvec(372 123)
\ifill f:0
\move(385 122)
\lvec(392 122)
\lvec(392 123)
\lvec(385 123)
\ifill f:0
\move(394 122)
\lvec(401 122)
\lvec(401 123)
\lvec(394 123)
\ifill f:0
\move(402 122)
\lvec(412 122)
\lvec(412 123)
\lvec(402 123)
\ifill f:0
\move(413 122)
\lvec(414 122)
\lvec(414 123)
\lvec(413 123)
\ifill f:0
\move(417 122)
\lvec(418 122)
\lvec(418 123)
\lvec(417 123)
\ifill f:0
\move(422 122)
\lvec(423 122)
\lvec(423 123)
\lvec(422 123)
\ifill f:0
\move(424 122)
\lvec(442 122)
\lvec(442 123)
\lvec(424 123)
\ifill f:0
\move(444 122)
\lvec(451 122)
\lvec(451 123)
\lvec(444 123)
\ifill f:0
\move(16 123)
\lvec(17 123)
\lvec(17 124)
\lvec(16 124)
\ifill f:0
\move(24 123)
\lvec(26 123)
\lvec(26 124)
\lvec(24 124)
\ifill f:0
\move(36 123)
\lvec(37 123)
\lvec(37 124)
\lvec(36 124)
\ifill f:0
\move(40 123)
\lvec(42 123)
\lvec(42 124)
\lvec(40 124)
\ifill f:0
\move(43 123)
\lvec(45 123)
\lvec(45 124)
\lvec(43 124)
\ifill f:0
\move(47 123)
\lvec(50 123)
\lvec(50 124)
\lvec(47 124)
\ifill f:0
\move(52 123)
\lvec(53 123)
\lvec(53 124)
\lvec(52 124)
\ifill f:0
\move(54 123)
\lvec(55 123)
\lvec(55 124)
\lvec(54 124)
\ifill f:0
\move(59 123)
\lvec(62 123)
\lvec(62 124)
\lvec(59 124)
\ifill f:0
\move(64 123)
\lvec(65 123)
\lvec(65 124)
\lvec(64 124)
\ifill f:0
\move(66 123)
\lvec(73 123)
\lvec(73 124)
\lvec(66 124)
\ifill f:0
\move(75 123)
\lvec(82 123)
\lvec(82 124)
\lvec(75 124)
\ifill f:0
\move(84 123)
\lvec(86 123)
\lvec(86 124)
\lvec(84 124)
\ifill f:0
\move(88 123)
\lvec(91 123)
\lvec(91 124)
\lvec(88 124)
\ifill f:0
\move(92 123)
\lvec(93 123)
\lvec(93 124)
\lvec(92 124)
\ifill f:0
\move(95 123)
\lvec(96 123)
\lvec(96 124)
\lvec(95 124)
\ifill f:0
\move(97 123)
\lvec(101 123)
\lvec(101 124)
\lvec(97 124)
\ifill f:0
\move(102 123)
\lvec(103 123)
\lvec(103 124)
\lvec(102 124)
\ifill f:0
\move(106 123)
\lvec(107 123)
\lvec(107 124)
\lvec(106 124)
\ifill f:0
\move(110 123)
\lvec(111 123)
\lvec(111 124)
\lvec(110 124)
\ifill f:0
\move(112 123)
\lvec(115 123)
\lvec(115 124)
\lvec(112 124)
\ifill f:0
\move(116 123)
\lvec(122 123)
\lvec(122 124)
\lvec(116 124)
\ifill f:0
\move(123 123)
\lvec(124 123)
\lvec(124 124)
\lvec(123 124)
\ifill f:0
\move(125 123)
\lvec(126 123)
\lvec(126 124)
\lvec(125 124)
\ifill f:0
\move(127 123)
\lvec(130 123)
\lvec(130 124)
\lvec(127 124)
\ifill f:0
\move(131 123)
\lvec(132 123)
\lvec(132 124)
\lvec(131 124)
\ifill f:0
\move(134 123)
\lvec(136 123)
\lvec(136 124)
\lvec(134 124)
\ifill f:0
\move(137 123)
\lvec(141 123)
\lvec(141 124)
\lvec(137 124)
\ifill f:0
\move(142 123)
\lvec(145 123)
\lvec(145 124)
\lvec(142 124)
\ifill f:0
\move(146 123)
\lvec(147 123)
\lvec(147 124)
\lvec(146 124)
\ifill f:0
\move(150 123)
\lvec(157 123)
\lvec(157 124)
\lvec(150 124)
\ifill f:0
\move(158 123)
\lvec(159 123)
\lvec(159 124)
\lvec(158 124)
\ifill f:0
\move(160 123)
\lvec(165 123)
\lvec(165 124)
\lvec(160 124)
\ifill f:0
\move(166 123)
\lvec(170 123)
\lvec(170 124)
\lvec(166 124)
\ifill f:0
\move(171 123)
\lvec(173 123)
\lvec(173 124)
\lvec(171 124)
\ifill f:0
\move(175 123)
\lvec(176 123)
\lvec(176 124)
\lvec(175 124)
\ifill f:0
\move(177 123)
\lvec(179 123)
\lvec(179 124)
\lvec(177 124)
\ifill f:0
\move(180 123)
\lvec(186 123)
\lvec(186 124)
\lvec(180 124)
\ifill f:0
\move(187 123)
\lvec(190 123)
\lvec(190 124)
\lvec(187 124)
\ifill f:0
\move(191 123)
\lvec(192 123)
\lvec(192 124)
\lvec(191 124)
\ifill f:0
\move(193 123)
\lvec(194 123)
\lvec(194 124)
\lvec(193 124)
\ifill f:0
\move(195 123)
\lvec(197 123)
\lvec(197 124)
\lvec(195 124)
\ifill f:0
\move(198 123)
\lvec(224 123)
\lvec(224 124)
\lvec(198 124)
\ifill f:0
\move(225 123)
\lvec(226 123)
\lvec(226 124)
\lvec(225 124)
\ifill f:0
\move(227 123)
\lvec(235 123)
\lvec(235 124)
\lvec(227 124)
\ifill f:0
\move(236 123)
\lvec(242 123)
\lvec(242 124)
\lvec(236 124)
\ifill f:0
\move(243 123)
\lvec(252 123)
\lvec(252 124)
\lvec(243 124)
\ifill f:0
\move(254 123)
\lvec(257 123)
\lvec(257 124)
\lvec(254 124)
\ifill f:0
\move(258 123)
\lvec(277 123)
\lvec(277 124)
\lvec(258 124)
\ifill f:0
\move(278 123)
\lvec(290 123)
\lvec(290 124)
\lvec(278 124)
\ifill f:0
\move(291 123)
\lvec(292 123)
\lvec(292 124)
\lvec(291 124)
\ifill f:0
\move(293 123)
\lvec(297 123)
\lvec(297 124)
\lvec(293 124)
\ifill f:0
\move(298 123)
\lvec(299 123)
\lvec(299 124)
\lvec(298 124)
\ifill f:0
\move(300 123)
\lvec(308 123)
\lvec(308 124)
\lvec(300 124)
\ifill f:0
\move(309 123)
\lvec(310 123)
\lvec(310 124)
\lvec(309 124)
\ifill f:0
\move(311 123)
\lvec(314 123)
\lvec(314 124)
\lvec(311 124)
\ifill f:0
\move(315 123)
\lvec(325 123)
\lvec(325 124)
\lvec(315 124)
\ifill f:0
\move(326 123)
\lvec(330 123)
\lvec(330 124)
\lvec(326 124)
\ifill f:0
\move(331 123)
\lvec(335 123)
\lvec(335 124)
\lvec(331 124)
\ifill f:0
\move(336 123)
\lvec(362 123)
\lvec(362 124)
\lvec(336 124)
\ifill f:0
\move(363 123)
\lvec(365 123)
\lvec(365 124)
\lvec(363 124)
\ifill f:0
\move(366 123)
\lvec(380 123)
\lvec(380 124)
\lvec(366 124)
\ifill f:0
\move(381 123)
\lvec(394 123)
\lvec(394 124)
\lvec(381 124)
\ifill f:0
\move(396 123)
\lvec(401 123)
\lvec(401 124)
\lvec(396 124)
\ifill f:0
\move(402 123)
\lvec(407 123)
\lvec(407 124)
\lvec(402 124)
\ifill f:0
\move(408 123)
\lvec(442 123)
\lvec(442 124)
\lvec(408 124)
\ifill f:0
\move(444 123)
\lvec(451 123)
\lvec(451 124)
\lvec(444 124)
\ifill f:0
\move(16 124)
\lvec(17 124)
\lvec(17 125)
\lvec(16 125)
\ifill f:0
\move(19 124)
\lvec(21 124)
\lvec(21 125)
\lvec(19 125)
\ifill f:0
\move(24 124)
\lvec(26 124)
\lvec(26 125)
\lvec(24 125)
\ifill f:0
\move(36 124)
\lvec(37 124)
\lvec(37 125)
\lvec(36 125)
\ifill f:0
\move(40 124)
\lvec(41 124)
\lvec(41 125)
\lvec(40 125)
\ifill f:0
\move(42 124)
\lvec(43 124)
\lvec(43 125)
\lvec(42 125)
\ifill f:0
\move(44 124)
\lvec(45 124)
\lvec(45 125)
\lvec(44 125)
\ifill f:0
\move(47 124)
\lvec(50 124)
\lvec(50 125)
\lvec(47 125)
\ifill f:0
\move(54 124)
\lvec(55 124)
\lvec(55 125)
\lvec(54 125)
\ifill f:0
\move(56 124)
\lvec(57 124)
\lvec(57 125)
\lvec(56 125)
\ifill f:0
\move(64 124)
\lvec(65 124)
\lvec(65 125)
\lvec(64 125)
\ifill f:0
\move(67 124)
\lvec(71 124)
\lvec(71 125)
\lvec(67 125)
\ifill f:0
\move(73 124)
\lvec(74 124)
\lvec(74 125)
\lvec(73 125)
\ifill f:0
\move(78 124)
\lvec(82 124)
\lvec(82 125)
\lvec(78 125)
\ifill f:0
\move(86 124)
\lvec(87 124)
\lvec(87 125)
\lvec(86 125)
\ifill f:0
\move(88 124)
\lvec(90 124)
\lvec(90 125)
\lvec(88 125)
\ifill f:0
\move(91 124)
\lvec(93 124)
\lvec(93 125)
\lvec(91 125)
\ifill f:0
\move(96 124)
\lvec(98 124)
\lvec(98 125)
\lvec(96 125)
\ifill f:0
\move(99 124)
\lvec(101 124)
\lvec(101 125)
\lvec(99 125)
\ifill f:0
\move(102 124)
\lvec(113 124)
\lvec(113 125)
\lvec(102 125)
\ifill f:0
\move(114 124)
\lvec(118 124)
\lvec(118 125)
\lvec(114 125)
\ifill f:0
\move(119 124)
\lvec(122 124)
\lvec(122 125)
\lvec(119 125)
\ifill f:0
\move(123 124)
\lvec(124 124)
\lvec(124 125)
\lvec(123 125)
\ifill f:0
\move(125 124)
\lvec(126 124)
\lvec(126 125)
\lvec(125 125)
\ifill f:0
\move(128 124)
\lvec(132 124)
\lvec(132 125)
\lvec(128 125)
\ifill f:0
\move(133 124)
\lvec(134 124)
\lvec(134 125)
\lvec(133 125)
\ifill f:0
\move(135 124)
\lvec(145 124)
\lvec(145 125)
\lvec(135 125)
\ifill f:0
\move(146 124)
\lvec(147 124)
\lvec(147 125)
\lvec(146 125)
\ifill f:0
\move(148 124)
\lvec(163 124)
\lvec(163 125)
\lvec(148 125)
\ifill f:0
\move(165 124)
\lvec(170 124)
\lvec(170 125)
\lvec(165 125)
\ifill f:0
\move(171 124)
\lvec(174 124)
\lvec(174 125)
\lvec(171 125)
\ifill f:0
\move(175 124)
\lvec(177 124)
\lvec(177 125)
\lvec(175 125)
\ifill f:0
\move(178 124)
\lvec(180 124)
\lvec(180 125)
\lvec(178 125)
\ifill f:0
\move(181 124)
\lvec(185 124)
\lvec(185 125)
\lvec(181 125)
\ifill f:0
\move(186 124)
\lvec(187 124)
\lvec(187 125)
\lvec(186 125)
\ifill f:0
\move(188 124)
\lvec(197 124)
\lvec(197 125)
\lvec(188 125)
\ifill f:0
\move(198 124)
\lvec(202 124)
\lvec(202 125)
\lvec(198 125)
\ifill f:0
\move(203 124)
\lvec(211 124)
\lvec(211 125)
\lvec(203 125)
\ifill f:0
\move(212 124)
\lvec(219 124)
\lvec(219 125)
\lvec(212 125)
\ifill f:0
\move(220 124)
\lvec(226 124)
\lvec(226 125)
\lvec(220 125)
\ifill f:0
\move(227 124)
\lvec(233 124)
\lvec(233 125)
\lvec(227 125)
\ifill f:0
\move(235 124)
\lvec(251 124)
\lvec(251 125)
\lvec(235 125)
\ifill f:0
\move(252 124)
\lvec(257 124)
\lvec(257 125)
\lvec(252 125)
\ifill f:0
\move(258 124)
\lvec(260 124)
\lvec(260 125)
\lvec(258 125)
\ifill f:0
\move(261 124)
\lvec(280 124)
\lvec(280 125)
\lvec(261 125)
\ifill f:0
\move(281 124)
\lvec(283 124)
\lvec(283 125)
\lvec(281 125)
\ifill f:0
\move(284 124)
\lvec(286 124)
\lvec(286 125)
\lvec(284 125)
\ifill f:0
\move(287 124)
\lvec(290 124)
\lvec(290 125)
\lvec(287 125)
\ifill f:0
\move(291 124)
\lvec(292 124)
\lvec(292 125)
\lvec(291 125)
\ifill f:0
\move(293 124)
\lvec(295 124)
\lvec(295 125)
\lvec(293 125)
\ifill f:0
\move(296 124)
\lvec(298 124)
\lvec(298 125)
\lvec(296 125)
\ifill f:0
\move(299 124)
\lvec(307 124)
\lvec(307 125)
\lvec(299 125)
\ifill f:0
\move(308 124)
\lvec(309 124)
\lvec(309 125)
\lvec(308 125)
\ifill f:0
\move(310 124)
\lvec(325 124)
\lvec(325 125)
\lvec(310 125)
\ifill f:0
\move(326 124)
\lvec(332 124)
\lvec(332 125)
\lvec(326 125)
\ifill f:0
\move(333 124)
\lvec(337 124)
\lvec(337 125)
\lvec(333 125)
\ifill f:0
\move(338 124)
\lvec(348 124)
\lvec(348 125)
\lvec(338 125)
\ifill f:0
\move(349 124)
\lvec(351 124)
\lvec(351 125)
\lvec(349 125)
\ifill f:0
\move(352 124)
\lvec(362 124)
\lvec(362 125)
\lvec(352 125)
\ifill f:0
\move(363 124)
\lvec(373 124)
\lvec(373 125)
\lvec(363 125)
\ifill f:0
\move(374 124)
\lvec(383 124)
\lvec(383 125)
\lvec(374 125)
\ifill f:0
\move(384 124)
\lvec(395 124)
\lvec(395 125)
\lvec(384 125)
\ifill f:0
\move(397 124)
\lvec(401 124)
\lvec(401 125)
\lvec(397 125)
\ifill f:0
\move(402 124)
\lvec(405 124)
\lvec(405 125)
\lvec(402 125)
\ifill f:0
\move(406 124)
\lvec(419 124)
\lvec(419 125)
\lvec(406 125)
\ifill f:0
\move(421 124)
\lvec(442 124)
\lvec(442 125)
\lvec(421 125)
\ifill f:0
\move(446 124)
\lvec(451 124)
\lvec(451 125)
\lvec(446 125)
\ifill f:0
\move(16 125)
\lvec(17 125)
\lvec(17 126)
\lvec(16 126)
\ifill f:0
\move(20 125)
\lvec(21 125)
\lvec(21 126)
\lvec(20 126)
\ifill f:0
\move(22 125)
\lvec(23 125)
\lvec(23 126)
\lvec(22 126)
\ifill f:0
\move(24 125)
\lvec(26 125)
\lvec(26 126)
\lvec(24 126)
\ifill f:0
\move(36 125)
\lvec(37 125)
\lvec(37 126)
\lvec(36 126)
\ifill f:0
\move(38 125)
\lvec(39 125)
\lvec(39 126)
\lvec(38 126)
\ifill f:0
\move(40 125)
\lvec(42 125)
\lvec(42 126)
\lvec(40 126)
\ifill f:0
\move(43 125)
\lvec(45 125)
\lvec(45 126)
\lvec(43 126)
\ifill f:0
\move(49 125)
\lvec(50 125)
\lvec(50 126)
\lvec(49 126)
\ifill f:0
\move(59 125)
\lvec(60 125)
\lvec(60 126)
\lvec(59 126)
\ifill f:0
\move(62 125)
\lvec(63 125)
\lvec(63 126)
\lvec(62 126)
\ifill f:0
\move(64 125)
\lvec(65 125)
\lvec(65 126)
\lvec(64 126)
\ifill f:0
\move(66 125)
\lvec(67 125)
\lvec(67 126)
\lvec(66 126)
\ifill f:0
\move(68 125)
\lvec(71 125)
\lvec(71 126)
\lvec(68 126)
\ifill f:0
\move(72 125)
\lvec(73 125)
\lvec(73 126)
\lvec(72 126)
\ifill f:0
\move(74 125)
\lvec(77 125)
\lvec(77 126)
\lvec(74 126)
\ifill f:0
\move(81 125)
\lvec(82 125)
\lvec(82 126)
\lvec(81 126)
\ifill f:0
\move(86 125)
\lvec(88 125)
\lvec(88 126)
\lvec(86 126)
\ifill f:0
\move(89 125)
\lvec(91 125)
\lvec(91 126)
\lvec(89 126)
\ifill f:0
\move(92 125)
\lvec(93 125)
\lvec(93 126)
\lvec(92 126)
\ifill f:0
\move(96 125)
\lvec(98 125)
\lvec(98 126)
\lvec(96 126)
\ifill f:0
\move(99 125)
\lvec(101 125)
\lvec(101 126)
\lvec(99 126)
\ifill f:0
\move(103 125)
\lvec(106 125)
\lvec(106 126)
\lvec(103 126)
\ifill f:0
\move(107 125)
\lvec(109 125)
\lvec(109 126)
\lvec(107 126)
\ifill f:0
\move(110 125)
\lvec(111 125)
\lvec(111 126)
\lvec(110 126)
\ifill f:0
\move(112 125)
\lvec(117 125)
\lvec(117 126)
\lvec(112 126)
\ifill f:0
\move(118 125)
\lvec(122 125)
\lvec(122 126)
\lvec(118 126)
\ifill f:0
\move(123 125)
\lvec(124 125)
\lvec(124 126)
\lvec(123 126)
\ifill f:0
\move(126 125)
\lvec(127 125)
\lvec(127 126)
\lvec(126 126)
\ifill f:0
\move(128 125)
\lvec(129 125)
\lvec(129 126)
\lvec(128 126)
\ifill f:0
\move(130 125)
\lvec(131 125)
\lvec(131 126)
\lvec(130 126)
\ifill f:0
\move(134 125)
\lvec(135 125)
\lvec(135 126)
\lvec(134 126)
\ifill f:0
\move(136 125)
\lvec(138 125)
\lvec(138 126)
\lvec(136 126)
\ifill f:0
\move(139 125)
\lvec(142 125)
\lvec(142 126)
\lvec(139 126)
\ifill f:0
\move(143 125)
\lvec(145 125)
\lvec(145 126)
\lvec(143 126)
\ifill f:0
\move(146 125)
\lvec(154 125)
\lvec(154 126)
\lvec(146 126)
\ifill f:0
\move(162 125)
\lvec(163 125)
\lvec(163 126)
\lvec(162 126)
\ifill f:0
\move(164 125)
\lvec(170 125)
\lvec(170 126)
\lvec(164 126)
\ifill f:0
\move(171 125)
\lvec(174 125)
\lvec(174 126)
\lvec(171 126)
\ifill f:0
\move(176 125)
\lvec(184 125)
\lvec(184 126)
\lvec(176 126)
\ifill f:0
\move(185 125)
\lvec(189 125)
\lvec(189 126)
\lvec(185 126)
\ifill f:0
\move(190 125)
\lvec(191 125)
\lvec(191 126)
\lvec(190 126)
\ifill f:0
\move(192 125)
\lvec(193 125)
\lvec(193 126)
\lvec(192 126)
\ifill f:0
\move(194 125)
\lvec(195 125)
\lvec(195 126)
\lvec(194 126)
\ifill f:0
\move(196 125)
\lvec(197 125)
\lvec(197 126)
\lvec(196 126)
\ifill f:0
\move(198 125)
\lvec(204 125)
\lvec(204 126)
\lvec(198 126)
\ifill f:0
\move(205 125)
\lvec(216 125)
\lvec(216 126)
\lvec(205 126)
\ifill f:0
\move(217 125)
\lvec(220 125)
\lvec(220 126)
\lvec(217 126)
\ifill f:0
\move(221 125)
\lvec(226 125)
\lvec(226 126)
\lvec(221 126)
\ifill f:0
\move(227 125)
\lvec(231 125)
\lvec(231 126)
\lvec(227 126)
\ifill f:0
\move(232 125)
\lvec(242 125)
\lvec(242 126)
\lvec(232 126)
\ifill f:0
\move(243 125)
\lvec(247 125)
\lvec(247 126)
\lvec(243 126)
\ifill f:0
\move(252 125)
\lvec(257 125)
\lvec(257 126)
\lvec(252 126)
\ifill f:0
\move(258 125)
\lvec(261 125)
\lvec(261 126)
\lvec(258 126)
\ifill f:0
\move(262 125)
\lvec(269 125)
\lvec(269 126)
\lvec(262 126)
\ifill f:0
\move(270 125)
\lvec(278 125)
\lvec(278 126)
\lvec(270 126)
\ifill f:0
\move(280 125)
\lvec(282 125)
\lvec(282 126)
\lvec(280 126)
\ifill f:0
\move(283 125)
\lvec(290 125)
\lvec(290 126)
\lvec(283 126)
\ifill f:0
\move(291 125)
\lvec(298 125)
\lvec(298 126)
\lvec(291 126)
\ifill f:0
\move(299 125)
\lvec(308 125)
\lvec(308 126)
\lvec(299 126)
\ifill f:0
\move(309 125)
\lvec(315 125)
\lvec(315 126)
\lvec(309 126)
\ifill f:0
\move(316 125)
\lvec(317 125)
\lvec(317 126)
\lvec(316 126)
\ifill f:0
\move(318 125)
\lvec(319 125)
\lvec(319 126)
\lvec(318 126)
\ifill f:0
\move(320 125)
\lvec(321 125)
\lvec(321 126)
\lvec(320 126)
\ifill f:0
\move(322 125)
\lvec(323 125)
\lvec(323 126)
\lvec(322 126)
\ifill f:0
\move(324 125)
\lvec(325 125)
\lvec(325 126)
\lvec(324 126)
\ifill f:0
\move(326 125)
\lvec(341 125)
\lvec(341 126)
\lvec(326 126)
\ifill f:0
\move(342 125)
\lvec(346 125)
\lvec(346 126)
\lvec(342 126)
\ifill f:0
\move(347 125)
\lvec(362 125)
\lvec(362 126)
\lvec(347 126)
\ifill f:0
\move(363 125)
\lvec(364 125)
\lvec(364 126)
\lvec(363 126)
\ifill f:0
\move(365 125)
\lvec(368 125)
\lvec(368 126)
\lvec(365 126)
\ifill f:0
\move(369 125)
\lvec(380 125)
\lvec(380 126)
\lvec(369 126)
\ifill f:0
\move(381 125)
\lvec(385 125)
\lvec(385 126)
\lvec(381 126)
\ifill f:0
\move(386 125)
\lvec(401 125)
\lvec(401 126)
\lvec(386 126)
\ifill f:0
\move(402 125)
\lvec(404 125)
\lvec(404 126)
\lvec(402 126)
\ifill f:0
\move(405 125)
\lvec(430 125)
\lvec(430 126)
\lvec(405 126)
\ifill f:0
\move(431 125)
\lvec(432 125)
\lvec(432 126)
\lvec(431 126)
\ifill f:0
\move(441 125)
\lvec(442 125)
\lvec(442 126)
\lvec(441 126)
\ifill f:0
\move(446 125)
\lvec(451 125)
\lvec(451 126)
\lvec(446 126)
\ifill f:0
\move(16 126)
\lvec(17 126)
\lvec(17 127)
\lvec(16 127)
\ifill f:0
\move(20 126)
\lvec(21 126)
\lvec(21 127)
\lvec(20 127)
\ifill f:0
\move(23 126)
\lvec(26 126)
\lvec(26 127)
\lvec(23 127)
\ifill f:0
\move(36 126)
\lvec(37 126)
\lvec(37 127)
\lvec(36 127)
\ifill f:0
\move(38 126)
\lvec(39 126)
\lvec(39 127)
\lvec(38 127)
\ifill f:0
\move(40 126)
\lvec(46 126)
\lvec(46 127)
\lvec(40 127)
\ifill f:0
\move(49 126)
\lvec(50 126)
\lvec(50 127)
\lvec(49 127)
\ifill f:0
\move(54 126)
\lvec(55 126)
\lvec(55 127)
\lvec(54 127)
\ifill f:0
\move(56 126)
\lvec(57 126)
\lvec(57 127)
\lvec(56 127)
\ifill f:0
\move(59 126)
\lvec(63 126)
\lvec(63 127)
\lvec(59 127)
\ifill f:0
\move(64 126)
\lvec(65 126)
\lvec(65 127)
\lvec(64 127)
\ifill f:0
\move(66 126)
\lvec(67 126)
\lvec(67 127)
\lvec(66 127)
\ifill f:0
\move(68 126)
\lvec(70 126)
\lvec(70 127)
\lvec(68 127)
\ifill f:0
\move(71 126)
\lvec(72 126)
\lvec(72 127)
\lvec(71 127)
\ifill f:0
\move(73 126)
\lvec(74 126)
\lvec(74 127)
\lvec(73 127)
\ifill f:0
\move(75 126)
\lvec(78 126)
\lvec(78 127)
\lvec(75 127)
\ifill f:0
\move(81 126)
\lvec(82 126)
\lvec(82 127)
\lvec(81 127)
\ifill f:0
\move(88 126)
\lvec(89 126)
\lvec(89 127)
\lvec(88 127)
\ifill f:0
\move(91 126)
\lvec(92 126)
\lvec(92 127)
\lvec(91 127)
\ifill f:0
\move(97 126)
\lvec(98 126)
\lvec(98 127)
\lvec(97 127)
\ifill f:0
\move(99 126)
\lvec(101 126)
\lvec(101 127)
\lvec(99 127)
\ifill f:0
\move(102 126)
\lvec(106 126)
\lvec(106 127)
\lvec(102 127)
\ifill f:0
\move(107 126)
\lvec(116 126)
\lvec(116 127)
\lvec(107 127)
\ifill f:0
\move(118 126)
\lvec(122 126)
\lvec(122 127)
\lvec(118 127)
\ifill f:0
\move(123 126)
\lvec(125 126)
\lvec(125 127)
\lvec(123 127)
\ifill f:0
\move(126 126)
\lvec(130 126)
\lvec(130 127)
\lvec(126 127)
\ifill f:0
\move(131 126)
\lvec(132 126)
\lvec(132 127)
\lvec(131 127)
\ifill f:0
\move(133 126)
\lvec(134 126)
\lvec(134 127)
\lvec(133 127)
\ifill f:0
\move(135 126)
\lvec(136 126)
\lvec(136 127)
\lvec(135 127)
\ifill f:0
\move(137 126)
\lvec(138 126)
\lvec(138 127)
\lvec(137 127)
\ifill f:0
\move(140 126)
\lvec(142 126)
\lvec(142 127)
\lvec(140 127)
\ifill f:0
\move(143 126)
\lvec(145 126)
\lvec(145 127)
\lvec(143 127)
\ifill f:0
\move(147 126)
\lvec(151 126)
\lvec(151 127)
\lvec(147 127)
\ifill f:0
\move(152 126)
\lvec(170 126)
\lvec(170 127)
\lvec(152 127)
\ifill f:0
\move(171 126)
\lvec(175 126)
\lvec(175 127)
\lvec(171 127)
\ifill f:0
\move(177 126)
\lvec(179 126)
\lvec(179 127)
\lvec(177 127)
\ifill f:0
\move(180 126)
\lvec(186 126)
\lvec(186 127)
\lvec(180 127)
\ifill f:0
\move(187 126)
\lvec(188 126)
\lvec(188 127)
\lvec(187 127)
\ifill f:0
\move(189 126)
\lvec(191 126)
\lvec(191 127)
\lvec(189 127)
\ifill f:0
\move(192 126)
\lvec(193 126)
\lvec(193 127)
\lvec(192 127)
\ifill f:0
\move(194 126)
\lvec(195 126)
\lvec(195 127)
\lvec(194 127)
\ifill f:0
\move(196 126)
\lvec(197 126)
\lvec(197 127)
\lvec(196 127)
\ifill f:0
\move(198 126)
\lvec(199 126)
\lvec(199 127)
\lvec(198 127)
\ifill f:0
\move(200 126)
\lvec(201 126)
\lvec(201 127)
\lvec(200 127)
\ifill f:0
\move(202 126)
\lvec(217 126)
\lvec(217 127)
\lvec(202 127)
\ifill f:0
\move(218 126)
\lvec(226 126)
\lvec(226 127)
\lvec(218 127)
\ifill f:0
\move(227 126)
\lvec(229 126)
\lvec(229 127)
\lvec(227 127)
\ifill f:0
\move(231 126)
\lvec(235 126)
\lvec(235 127)
\lvec(231 127)
\ifill f:0
\move(236 126)
\lvec(237 126)
\lvec(237 127)
\lvec(236 127)
\ifill f:0
\move(238 126)
\lvec(257 126)
\lvec(257 127)
\lvec(238 127)
\ifill f:0
\move(258 126)
\lvec(263 126)
\lvec(263 127)
\lvec(258 127)
\ifill f:0
\move(264 126)
\lvec(271 126)
\lvec(271 127)
\lvec(264 127)
\ifill f:0
\move(272 126)
\lvec(277 126)
\lvec(277 127)
\lvec(272 127)
\ifill f:0
\move(278 126)
\lvec(290 126)
\lvec(290 127)
\lvec(278 127)
\ifill f:0
\move(291 126)
\lvec(293 126)
\lvec(293 127)
\lvec(291 127)
\ifill f:0
\move(294 126)
\lvec(296 126)
\lvec(296 127)
\lvec(294 127)
\ifill f:0
\move(297 126)
\lvec(299 126)
\lvec(299 127)
\lvec(297 127)
\ifill f:0
\move(300 126)
\lvec(302 126)
\lvec(302 127)
\lvec(300 127)
\ifill f:0
\move(303 126)
\lvec(307 126)
\lvec(307 127)
\lvec(303 127)
\ifill f:0
\move(308 126)
\lvec(310 126)
\lvec(310 127)
\lvec(308 127)
\ifill f:0
\move(311 126)
\lvec(319 126)
\lvec(319 127)
\lvec(311 127)
\ifill f:0
\move(320 126)
\lvec(321 126)
\lvec(321 127)
\lvec(320 127)
\ifill f:0
\move(322 126)
\lvec(323 126)
\lvec(323 127)
\lvec(322 127)
\ifill f:0
\move(324 126)
\lvec(325 126)
\lvec(325 127)
\lvec(324 127)
\ifill f:0
\move(326 126)
\lvec(327 126)
\lvec(327 127)
\lvec(326 127)
\ifill f:0
\move(328 126)
\lvec(329 126)
\lvec(329 127)
\lvec(328 127)
\ifill f:0
\move(330 126)
\lvec(355 126)
\lvec(355 127)
\lvec(330 127)
\ifill f:0
\move(356 126)
\lvec(358 126)
\lvec(358 127)
\lvec(356 127)
\ifill f:0
\move(359 126)
\lvec(362 126)
\lvec(362 127)
\lvec(359 127)
\ifill f:0
\move(363 126)
\lvec(364 126)
\lvec(364 127)
\lvec(363 127)
\ifill f:0
\move(365 126)
\lvec(367 126)
\lvec(367 127)
\lvec(365 127)
\ifill f:0
\move(368 126)
\lvec(378 126)
\lvec(378 127)
\lvec(368 127)
\ifill f:0
\move(379 126)
\lvec(382 126)
\lvec(382 127)
\lvec(379 127)
\ifill f:0
\move(383 126)
\lvec(397 126)
\lvec(397 127)
\lvec(383 127)
\ifill f:0
\move(398 126)
\lvec(401 126)
\lvec(401 127)
\lvec(398 127)
\ifill f:0
\move(402 126)
\lvec(403 126)
\lvec(403 127)
\lvec(402 127)
\ifill f:0
\move(404 126)
\lvec(436 126)
\lvec(436 127)
\lvec(404 127)
\ifill f:0
\move(437 126)
\lvec(439 126)
\lvec(439 127)
\lvec(437 127)
\ifill f:0
\move(441 126)
\lvec(442 126)
\lvec(442 127)
\lvec(441 127)
\ifill f:0
\move(15 127)
\lvec(17 127)
\lvec(17 128)
\lvec(15 128)
\ifill f:0
\move(19 127)
\lvec(21 127)
\lvec(21 128)
\lvec(19 128)
\ifill f:0
\move(24 127)
\lvec(26 127)
\lvec(26 128)
\lvec(24 128)
\ifill f:0
\move(36 127)
\lvec(37 127)
\lvec(37 128)
\lvec(36 128)
\ifill f:0
\move(38 127)
\lvec(39 127)
\lvec(39 128)
\lvec(38 128)
\ifill f:0
\move(40 127)
\lvec(41 127)
\lvec(41 128)
\lvec(40 128)
\ifill f:0
\move(43 127)
\lvec(45 127)
\lvec(45 128)
\lvec(43 128)
\ifill f:0
\move(47 127)
\lvec(48 127)
\lvec(48 128)
\lvec(47 128)
\ifill f:0
\move(49 127)
\lvec(50 127)
\lvec(50 128)
\lvec(49 128)
\ifill f:0
\move(51 127)
\lvec(52 127)
\lvec(52 128)
\lvec(51 128)
\ifill f:0
\move(54 127)
\lvec(56 127)
\lvec(56 128)
\lvec(54 128)
\ifill f:0
\move(57 127)
\lvec(58 127)
\lvec(58 128)
\lvec(57 128)
\ifill f:0
\move(60 127)
\lvec(63 127)
\lvec(63 128)
\lvec(60 128)
\ifill f:0
\move(64 127)
\lvec(65 127)
\lvec(65 128)
\lvec(64 128)
\ifill f:0
\move(66 127)
\lvec(71 127)
\lvec(71 128)
\lvec(66 128)
\ifill f:0
\move(72 127)
\lvec(73 127)
\lvec(73 128)
\lvec(72 128)
\ifill f:0
\move(74 127)
\lvec(75 127)
\lvec(75 128)
\lvec(74 128)
\ifill f:0
\move(76 127)
\lvec(79 127)
\lvec(79 128)
\lvec(76 128)
\ifill f:0
\move(81 127)
\lvec(82 127)
\lvec(82 128)
\lvec(81 128)
\ifill f:0
\move(83 127)
\lvec(85 127)
\lvec(85 128)
\lvec(83 128)
\ifill f:0
\move(88 127)
\lvec(91 127)
\lvec(91 128)
\lvec(88 128)
\ifill f:0
\move(92 127)
\lvec(93 127)
\lvec(93 128)
\lvec(92 128)
\ifill f:0
\move(95 127)
\lvec(96 127)
\lvec(96 128)
\lvec(95 128)
\ifill f:0
\move(97 127)
\lvec(98 127)
\lvec(98 128)
\lvec(97 128)
\ifill f:0
\move(99 127)
\lvec(101 127)
\lvec(101 128)
\lvec(99 128)
\ifill f:0
\move(102 127)
\lvec(105 127)
\lvec(105 128)
\lvec(102 128)
\ifill f:0
\move(106 127)
\lvec(111 127)
\lvec(111 128)
\lvec(106 128)
\ifill f:0
\move(112 127)
\lvec(113 127)
\lvec(113 128)
\lvec(112 128)
\ifill f:0
\move(114 127)
\lvec(122 127)
\lvec(122 128)
\lvec(114 128)
\ifill f:0
\move(123 127)
\lvec(125 127)
\lvec(125 128)
\lvec(123 128)
\ifill f:0
\move(127 127)
\lvec(128 127)
\lvec(128 128)
\lvec(127 128)
\ifill f:0
\move(129 127)
\lvec(131 127)
\lvec(131 128)
\lvec(129 128)
\ifill f:0
\move(134 127)
\lvec(135 127)
\lvec(135 128)
\lvec(134 128)
\ifill f:0
\move(136 127)
\lvec(137 127)
\lvec(137 128)
\lvec(136 128)
\ifill f:0
\move(138 127)
\lvec(145 127)
\lvec(145 128)
\lvec(138 128)
\ifill f:0
\move(146 127)
\lvec(149 127)
\lvec(149 128)
\lvec(146 128)
\ifill f:0
\move(150 127)
\lvec(155 127)
\lvec(155 128)
\lvec(150 128)
\ifill f:0
\move(157 127)
\lvec(170 127)
\lvec(170 128)
\lvec(157 128)
\ifill f:0
\move(172 127)
\lvec(176 127)
\lvec(176 128)
\lvec(172 128)
\ifill f:0
\move(177 127)
\lvec(184 127)
\lvec(184 128)
\lvec(177 128)
\ifill f:0
\move(185 127)
\lvec(190 127)
\lvec(190 128)
\lvec(185 128)
\ifill f:0
\move(191 127)
\lvec(193 127)
\lvec(193 128)
\lvec(191 128)
\ifill f:0
\move(194 127)
\lvec(195 127)
\lvec(195 128)
\lvec(194 128)
\ifill f:0
\move(196 127)
\lvec(197 127)
\lvec(197 128)
\lvec(196 128)
\ifill f:0
\move(198 127)
\lvec(221 127)
\lvec(221 128)
\lvec(198 128)
\ifill f:0
\move(222 127)
\lvec(226 127)
\lvec(226 128)
\lvec(222 128)
\ifill f:0
\move(227 127)
\lvec(229 127)
\lvec(229 128)
\lvec(227 128)
\ifill f:0
\move(230 127)
\lvec(234 127)
\lvec(234 128)
\lvec(230 128)
\ifill f:0
\move(235 127)
\lvec(241 127)
\lvec(241 128)
\lvec(235 128)
\ifill f:0
\move(243 127)
\lvec(257 127)
\lvec(257 128)
\lvec(243 128)
\ifill f:0
\move(258 127)
\lvec(266 127)
\lvec(266 128)
\lvec(258 128)
\ifill f:0
\move(267 127)
\lvec(274 127)
\lvec(274 128)
\lvec(267 128)
\ifill f:0
\move(275 127)
\lvec(280 127)
\lvec(280 128)
\lvec(275 128)
\ifill f:0
\move(281 127)
\lvec(290 127)
\lvec(290 128)
\lvec(281 128)
\ifill f:0
\move(291 127)
\lvec(293 127)
\lvec(293 128)
\lvec(291 128)
\ifill f:0
\move(294 127)
\lvec(311 127)
\lvec(311 128)
\lvec(294 128)
\ifill f:0
\move(312 127)
\lvec(316 127)
\lvec(316 128)
\lvec(312 128)
\ifill f:0
\move(317 127)
\lvec(321 127)
\lvec(321 128)
\lvec(317 128)
\ifill f:0
\move(322 127)
\lvec(323 127)
\lvec(323 128)
\lvec(322 128)
\ifill f:0
\move(324 127)
\lvec(325 127)
\lvec(325 128)
\lvec(324 128)
\ifill f:0
\move(326 127)
\lvec(341 127)
\lvec(341 128)
\lvec(326 128)
\ifill f:0
\move(342 127)
\lvec(353 127)
\lvec(353 128)
\lvec(342 128)
\ifill f:0
\move(354 127)
\lvec(362 127)
\lvec(362 128)
\lvec(354 128)
\ifill f:0
\move(363 127)
\lvec(370 127)
\lvec(370 128)
\lvec(363 128)
\ifill f:0
\move(371 127)
\lvec(380 127)
\lvec(380 128)
\lvec(371 128)
\ifill f:0
\move(381 127)
\lvec(384 127)
\lvec(384 128)
\lvec(381 128)
\ifill f:0
\move(385 127)
\lvec(388 127)
\lvec(388 128)
\lvec(385 128)
\ifill f:0
\move(389 127)
\lvec(392 127)
\lvec(392 128)
\lvec(389 128)
\ifill f:0
\move(393 127)
\lvec(397 127)
\lvec(397 128)
\lvec(393 128)
\ifill f:0
\move(398 127)
\lvec(401 127)
\lvec(401 128)
\lvec(398 128)
\ifill f:0
\move(402 127)
\lvec(416 127)
\lvec(416 128)
\lvec(402 128)
\ifill f:0
\move(417 127)
\lvec(425 127)
\lvec(425 128)
\lvec(417 128)
\ifill f:0
\move(426 127)
\lvec(439 127)
\lvec(439 128)
\lvec(426 128)
\ifill f:0
\move(441 127)
\lvec(442 127)
\lvec(442 128)
\lvec(441 128)
\ifill f:0
\move(443 127)
\lvec(451 127)
\lvec(451 128)
\lvec(443 128)
\ifill f:0
\move(15 128)
\lvec(17 128)
\lvec(17 129)
\lvec(15 129)
\ifill f:0
\move(20 128)
\lvec(21 128)
\lvec(21 129)
\lvec(20 129)
\ifill f:0
\move(24 128)
\lvec(26 128)
\lvec(26 129)
\lvec(24 129)
\ifill f:0
\move(36 128)
\lvec(37 128)
\lvec(37 129)
\lvec(36 129)
\ifill f:0
\move(38 128)
\lvec(39 128)
\lvec(39 129)
\lvec(38 129)
\ifill f:0
\move(40 128)
\lvec(45 128)
\lvec(45 129)
\lvec(40 129)
\ifill f:0
\move(46 128)
\lvec(48 128)
\lvec(48 129)
\lvec(46 129)
\ifill f:0
\move(49 128)
\lvec(50 128)
\lvec(50 129)
\lvec(49 129)
\ifill f:0
\move(51 128)
\lvec(53 128)
\lvec(53 129)
\lvec(51 129)
\ifill f:0
\move(56 128)
\lvec(57 128)
\lvec(57 129)
\lvec(56 129)
\ifill f:0
\move(59 128)
\lvec(60 128)
\lvec(60 129)
\lvec(59 129)
\ifill f:0
\move(61 128)
\lvec(63 128)
\lvec(63 129)
\lvec(61 129)
\ifill f:0
\move(64 128)
\lvec(65 128)
\lvec(65 129)
\lvec(64 129)
\ifill f:0
\move(66 128)
\lvec(74 128)
\lvec(74 129)
\lvec(66 129)
\ifill f:0
\move(75 128)
\lvec(80 128)
\lvec(80 129)
\lvec(75 129)
\ifill f:0
\move(81 128)
\lvec(82 128)
\lvec(82 129)
\lvec(81 129)
\ifill f:0
\move(83 128)
\lvec(85 128)
\lvec(85 129)
\lvec(83 129)
\ifill f:0
\move(86 128)
\lvec(88 128)
\lvec(88 129)
\lvec(86 129)
\ifill f:0
\move(91 128)
\lvec(93 128)
\lvec(93 129)
\lvec(91 129)
\ifill f:0
\move(97 128)
\lvec(98 128)
\lvec(98 129)
\lvec(97 129)
\ifill f:0
\move(99 128)
\lvec(101 128)
\lvec(101 129)
\lvec(99 129)
\ifill f:0
\move(102 128)
\lvec(104 128)
\lvec(104 129)
\lvec(102 129)
\ifill f:0
\move(105 128)
\lvec(122 128)
\lvec(122 129)
\lvec(105 129)
\ifill f:0
\move(123 128)
\lvec(126 128)
\lvec(126 129)
\lvec(123 129)
\ifill f:0
\move(128 128)
\lvec(129 128)
\lvec(129 129)
\lvec(128 129)
\ifill f:0
\move(131 128)
\lvec(132 128)
\lvec(132 129)
\lvec(131 129)
\ifill f:0
\move(133 128)
\lvec(134 128)
\lvec(134 129)
\lvec(133 129)
\ifill f:0
\move(135 128)
\lvec(140 128)
\lvec(140 129)
\lvec(135 129)
\ifill f:0
\move(141 128)
\lvec(145 128)
\lvec(145 129)
\lvec(141 129)
\ifill f:0
\move(146 128)
\lvec(153 128)
\lvec(153 129)
\lvec(146 129)
\ifill f:0
\move(154 128)
\lvec(159 128)
\lvec(159 129)
\lvec(154 129)
\ifill f:0
\move(162 128)
\lvec(170 128)
\lvec(170 129)
\lvec(162 129)
\ifill f:0
\move(175 128)
\lvec(186 128)
\lvec(186 129)
\lvec(175 129)
\ifill f:0
\move(187 128)
\lvec(195 128)
\lvec(195 129)
\lvec(187 129)
\ifill f:0
\move(196 128)
\lvec(197 128)
\lvec(197 129)
\lvec(196 129)
\ifill f:0
\move(198 128)
\lvec(226 128)
\lvec(226 129)
\lvec(198 129)
\ifill f:0
\move(227 128)
\lvec(228 128)
\lvec(228 129)
\lvec(227 129)
\ifill f:0
\move(229 128)
\lvec(233 128)
\lvec(233 129)
\lvec(229 129)
\ifill f:0
\move(234 128)
\lvec(238 128)
\lvec(238 129)
\lvec(234 129)
\ifill f:0
\move(239 128)
\lvec(257 128)
\lvec(257 129)
\lvec(239 129)
\ifill f:0
\move(258 128)
\lvec(271 128)
\lvec(271 129)
\lvec(258 129)
\ifill f:0
\move(272 128)
\lvec(279 128)
\lvec(279 129)
\lvec(272 129)
\ifill f:0
\move(280 128)
\lvec(290 128)
\lvec(290 129)
\lvec(280 129)
\ifill f:0
\move(291 128)
\lvec(298 128)
\lvec(298 129)
\lvec(291 129)
\ifill f:0
\move(299 128)
\lvec(304 128)
\lvec(304 129)
\lvec(299 129)
\ifill f:0
\move(305 128)
\lvec(307 128)
\lvec(307 129)
\lvec(305 129)
\ifill f:0
\move(308 128)
\lvec(310 128)
\lvec(310 129)
\lvec(308 129)
\ifill f:0
\move(311 128)
\lvec(313 128)
\lvec(313 129)
\lvec(311 129)
\ifill f:0
\move(314 128)
\lvec(318 128)
\lvec(318 129)
\lvec(314 129)
\ifill f:0
\move(319 128)
\lvec(323 128)
\lvec(323 129)
\lvec(319 129)
\ifill f:0
\move(324 128)
\lvec(325 128)
\lvec(325 129)
\lvec(324 129)
\ifill f:0
\move(326 128)
\lvec(338 128)
\lvec(338 129)
\lvec(326 129)
\ifill f:0
\move(339 128)
\lvec(362 128)
\lvec(362 129)
\lvec(339 129)
\ifill f:0
\move(363 128)
\lvec(369 128)
\lvec(369 129)
\lvec(363 129)
\ifill f:0
\move(370 128)
\lvec(375 128)
\lvec(375 129)
\lvec(370 129)
\ifill f:0
\move(376 128)
\lvec(382 128)
\lvec(382 129)
\lvec(376 129)
\ifill f:0
\move(383 128)
\lvec(389 128)
\lvec(389 129)
\lvec(383 129)
\ifill f:0
\move(390 128)
\lvec(401 128)
\lvec(401 129)
\lvec(390 129)
\ifill f:0
\move(402 128)
\lvec(440 128)
\lvec(440 129)
\lvec(402 129)
\ifill f:0
\move(441 128)
\lvec(442 128)
\lvec(442 129)
\lvec(441 129)
\ifill f:0
\move(443 128)
\lvec(451 128)
\lvec(451 129)
\lvec(443 129)
\ifill f:0
\move(15 129)
\lvec(17 129)
\lvec(17 130)
\lvec(15 130)
\ifill f:0
\move(20 129)
\lvec(21 129)
\lvec(21 130)
\lvec(20 130)
\ifill f:0
\move(25 129)
\lvec(26 129)
\lvec(26 130)
\lvec(25 130)
\ifill f:0
\move(36 129)
\lvec(37 129)
\lvec(37 130)
\lvec(36 130)
\ifill f:0
\move(40 129)
\lvec(41 129)
\lvec(41 130)
\lvec(40 130)
\ifill f:0
\move(43 129)
\lvec(45 129)
\lvec(45 130)
\lvec(43 130)
\ifill f:0
\move(47 129)
\lvec(50 129)
\lvec(50 130)
\lvec(47 130)
\ifill f:0
\move(51 129)
\lvec(52 129)
\lvec(52 130)
\lvec(51 130)
\ifill f:0
\move(54 129)
\lvec(55 129)
\lvec(55 130)
\lvec(54 130)
\ifill f:0
\move(57 129)
\lvec(58 129)
\lvec(58 130)
\lvec(57 130)
\ifill f:0
\move(59 129)
\lvec(60 129)
\lvec(60 130)
\lvec(59 130)
\ifill f:0
\move(62 129)
\lvec(63 129)
\lvec(63 130)
\lvec(62 130)
\ifill f:0
\move(64 129)
\lvec(65 129)
\lvec(65 130)
\lvec(64 130)
\ifill f:0
\move(66 129)
\lvec(70 129)
\lvec(70 130)
\lvec(66 130)
\ifill f:0
\move(72 129)
\lvec(73 129)
\lvec(73 130)
\lvec(72 130)
\ifill f:0
\move(76 129)
\lvec(80 129)
\lvec(80 130)
\lvec(76 130)
\ifill f:0
\move(81 129)
\lvec(82 129)
\lvec(82 130)
\lvec(81 130)
\ifill f:0
\move(83 129)
\lvec(85 129)
\lvec(85 130)
\lvec(83 130)
\ifill f:0
\move(86 129)
\lvec(87 129)
\lvec(87 130)
\lvec(86 130)
\ifill f:0
\move(88 129)
\lvec(91 129)
\lvec(91 130)
\lvec(88 130)
\ifill f:0
\move(92 129)
\lvec(93 129)
\lvec(93 130)
\lvec(92 130)
\ifill f:0
\move(96 129)
\lvec(97 129)
\lvec(97 130)
\lvec(96 130)
\ifill f:0
\move(100 129)
\lvec(101 129)
\lvec(101 130)
\lvec(100 130)
\ifill f:0
\move(102 129)
\lvec(103 129)
\lvec(103 130)
\lvec(102 130)
\ifill f:0
\move(104 129)
\lvec(106 129)
\lvec(106 130)
\lvec(104 130)
\ifill f:0
\move(107 129)
\lvec(111 129)
\lvec(111 130)
\lvec(107 130)
\ifill f:0
\move(113 129)
\lvec(122 129)
\lvec(122 130)
\lvec(113 130)
\ifill f:0
\move(125 129)
\lvec(127 129)
\lvec(127 130)
\lvec(125 130)
\ifill f:0
\move(128 129)
\lvec(131 129)
\lvec(131 130)
\lvec(128 130)
\ifill f:0
\move(135 129)
\lvec(136 129)
\lvec(136 130)
\lvec(135 130)
\ifill f:0
\move(137 129)
\lvec(138 129)
\lvec(138 130)
\lvec(137 130)
\ifill f:0
\move(139 129)
\lvec(145 129)
\lvec(145 130)
\lvec(139 130)
\ifill f:0
\move(146 129)
\lvec(155 129)
\lvec(155 130)
\lvec(146 130)
\ifill f:0
\move(156 129)
\lvec(163 129)
\lvec(163 130)
\lvec(156 130)
\ifill f:0
\move(165 129)
\lvec(166 129)
\lvec(166 130)
\lvec(165 130)
\ifill f:0
\move(169 129)
\lvec(170 129)
\lvec(170 130)
\lvec(169 130)
\ifill f:0
\move(175 129)
\lvec(180 129)
\lvec(180 130)
\lvec(175 130)
\ifill f:0
\move(181 129)
\lvec(185 129)
\lvec(185 130)
\lvec(181 130)
\ifill f:0
\move(187 129)
\lvec(189 129)
\lvec(189 130)
\lvec(187 130)
\ifill f:0
\move(190 129)
\lvec(192 129)
\lvec(192 130)
\lvec(190 130)
\ifill f:0
\move(193 129)
\lvec(195 129)
\lvec(195 130)
\lvec(193 130)
\ifill f:0
\move(196 129)
\lvec(197 129)
\lvec(197 130)
\lvec(196 130)
\ifill f:0
\move(198 129)
\lvec(214 129)
\lvec(214 130)
\lvec(198 130)
\ifill f:0
\move(215 129)
\lvec(219 129)
\lvec(219 130)
\lvec(215 130)
\ifill f:0
\move(220 129)
\lvec(222 129)
\lvec(222 130)
\lvec(220 130)
\ifill f:0
\move(223 129)
\lvec(226 129)
\lvec(226 130)
\lvec(223 130)
\ifill f:0
\move(227 129)
\lvec(228 129)
\lvec(228 130)
\lvec(227 130)
\ifill f:0
\move(229 129)
\lvec(232 129)
\lvec(232 130)
\lvec(229 130)
\ifill f:0
\move(233 129)
\lvec(235 129)
\lvec(235 130)
\lvec(233 130)
\ifill f:0
\move(236 129)
\lvec(241 129)
\lvec(241 130)
\lvec(236 130)
\ifill f:0
\move(242 129)
\lvec(248 129)
\lvec(248 130)
\lvec(242 130)
\ifill f:0
\move(250 129)
\lvec(257 129)
\lvec(257 130)
\lvec(250 130)
\ifill f:0
\move(258 129)
\lvec(275 129)
\lvec(275 130)
\lvec(258 130)
\ifill f:0
\move(276 129)
\lvec(284 129)
\lvec(284 130)
\lvec(276 130)
\ifill f:0
\move(285 129)
\lvec(290 129)
\lvec(290 130)
\lvec(285 130)
\ifill f:0
\move(291 129)
\lvec(294 129)
\lvec(294 130)
\lvec(291 130)
\ifill f:0
\move(295 129)
\lvec(302 129)
\lvec(302 130)
\lvec(295 130)
\ifill f:0
\move(303 129)
\lvec(309 129)
\lvec(309 130)
\lvec(303 130)
\ifill f:0
\move(310 129)
\lvec(315 129)
\lvec(315 130)
\lvec(310 130)
\ifill f:0
\move(316 129)
\lvec(323 129)
\lvec(323 130)
\lvec(316 130)
\ifill f:0
\move(324 129)
\lvec(325 129)
\lvec(325 130)
\lvec(324 130)
\ifill f:0
\move(326 129)
\lvec(330 129)
\lvec(330 130)
\lvec(326 130)
\ifill f:0
\move(331 129)
\lvec(356 129)
\lvec(356 130)
\lvec(331 130)
\ifill f:0
\move(357 129)
\lvec(362 129)
\lvec(362 130)
\lvec(357 130)
\ifill f:0
\move(363 129)
\lvec(366 129)
\lvec(366 130)
\lvec(363 130)
\ifill f:0
\move(367 129)
\lvec(371 129)
\lvec(371 130)
\lvec(367 130)
\ifill f:0
\move(372 129)
\lvec(377 129)
\lvec(377 130)
\lvec(372 130)
\ifill f:0
\move(378 129)
\lvec(380 129)
\lvec(380 130)
\lvec(378 130)
\ifill f:0
\move(381 129)
\lvec(383 129)
\lvec(383 130)
\lvec(381 130)
\ifill f:0
\move(384 129)
\lvec(394 129)
\lvec(394 130)
\lvec(384 130)
\ifill f:0
\move(395 129)
\lvec(398 129)
\lvec(398 130)
\lvec(395 130)
\ifill f:0
\move(399 129)
\lvec(401 129)
\lvec(401 130)
\lvec(399 130)
\ifill f:0
\move(402 129)
\lvec(423 129)
\lvec(423 130)
\lvec(402 130)
\ifill f:0
\move(424 129)
\lvec(431 129)
\lvec(431 130)
\lvec(424 130)
\ifill f:0
\move(432 129)
\lvec(440 129)
\lvec(440 130)
\lvec(432 130)
\ifill f:0
\move(441 129)
\lvec(442 129)
\lvec(442 130)
\lvec(441 130)
\ifill f:0
\move(443 129)
\lvec(451 129)
\lvec(451 130)
\lvec(443 130)
\ifill f:0
\move(16 130)
\lvec(17 130)
\lvec(17 131)
\lvec(16 131)
\ifill f:0
\move(20 130)
\lvec(21 130)
\lvec(21 131)
\lvec(20 131)
\ifill f:0
\move(25 130)
\lvec(26 130)
\lvec(26 131)
\lvec(25 131)
\ifill f:0
\move(36 130)
\lvec(37 130)
\lvec(37 131)
\lvec(36 131)
\ifill f:0
\move(40 130)
\lvec(45 130)
\lvec(45 131)
\lvec(40 131)
\ifill f:0
\move(47 130)
\lvec(50 130)
\lvec(50 131)
\lvec(47 131)
\ifill f:0
\move(56 130)
\lvec(57 130)
\lvec(57 131)
\lvec(56 131)
\ifill f:0
\move(60 130)
\lvec(61 130)
\lvec(61 131)
\lvec(60 131)
\ifill f:0
\move(62 130)
\lvec(65 130)
\lvec(65 131)
\lvec(62 131)
\ifill f:0
\move(66 130)
\lvec(75 130)
\lvec(75 131)
\lvec(66 131)
\ifill f:0
\move(76 130)
\lvec(77 130)
\lvec(77 131)
\lvec(76 131)
\ifill f:0
\move(78 130)
\lvec(80 130)
\lvec(80 131)
\lvec(78 131)
\ifill f:0
\move(81 130)
\lvec(82 130)
\lvec(82 131)
\lvec(81 131)
\ifill f:0
\move(83 130)
\lvec(84 130)
\lvec(84 131)
\lvec(83 131)
\ifill f:0
\move(89 130)
\lvec(90 130)
\lvec(90 131)
\lvec(89 131)
\ifill f:0
\move(91 130)
\lvec(93 130)
\lvec(93 131)
\lvec(91 131)
\ifill f:0
\move(96 130)
\lvec(98 130)
\lvec(98 131)
\lvec(96 131)
\ifill f:0
\move(100 130)
\lvec(101 130)
\lvec(101 131)
\lvec(100 131)
\ifill f:0
\move(102 130)
\lvec(103 130)
\lvec(103 131)
\lvec(102 131)
\ifill f:0
\move(104 130)
\lvec(109 130)
\lvec(109 131)
\lvec(104 131)
\ifill f:0
\move(110 130)
\lvec(115 130)
\lvec(115 131)
\lvec(110 131)
\ifill f:0
\move(116 130)
\lvec(122 130)
\lvec(122 131)
\lvec(116 131)
\ifill f:0
\move(126 130)
\lvec(128 130)
\lvec(128 131)
\lvec(126 131)
\ifill f:0
\move(130 130)
\lvec(132 130)
\lvec(132 131)
\lvec(130 131)
\ifill f:0
\move(133 130)
\lvec(135 130)
\lvec(135 131)
\lvec(133 131)
\ifill f:0
\move(136 130)
\lvec(137 130)
\lvec(137 131)
\lvec(136 131)
\ifill f:0
\move(139 130)
\lvec(145 130)
\lvec(145 131)
\lvec(139 131)
\ifill f:0
\move(146 130)
\lvec(154 130)
\lvec(154 131)
\lvec(146 131)
\ifill f:0
\move(155 130)
\lvec(163 130)
\lvec(163 131)
\lvec(155 131)
\ifill f:0
\move(164 130)
\lvec(167 130)
\lvec(167 131)
\lvec(164 131)
\ifill f:0
\move(169 130)
\lvec(170 130)
\lvec(170 131)
\lvec(169 131)
\ifill f:0
\move(177 130)
\lvec(183 130)
\lvec(183 131)
\lvec(177 131)
\ifill f:0
\move(184 130)
\lvec(192 130)
\lvec(192 131)
\lvec(184 131)
\ifill f:0
\move(193 130)
\lvec(195 130)
\lvec(195 131)
\lvec(193 131)
\ifill f:0
\move(196 130)
\lvec(197 130)
\lvec(197 131)
\lvec(196 131)
\ifill f:0
\move(198 130)
\lvec(215 130)
\lvec(215 131)
\lvec(198 131)
\ifill f:0
\move(216 130)
\lvec(226 130)
\lvec(226 131)
\lvec(216 131)
\ifill f:0
\move(227 130)
\lvec(228 130)
\lvec(228 131)
\lvec(227 131)
\ifill f:0
\move(229 130)
\lvec(231 130)
\lvec(231 131)
\lvec(229 131)
\ifill f:0
\move(232 130)
\lvec(234 130)
\lvec(234 131)
\lvec(232 131)
\ifill f:0
\move(235 130)
\lvec(239 130)
\lvec(239 131)
\lvec(235 131)
\ifill f:0
\move(240 130)
\lvec(244 130)
\lvec(244 131)
\lvec(240 131)
\ifill f:0
\move(245 130)
\lvec(250 130)
\lvec(250 131)
\lvec(245 131)
\ifill f:0
\move(252 130)
\lvec(257 130)
\lvec(257 131)
\lvec(252 131)
\ifill f:0
\move(258 130)
\lvec(266 130)
\lvec(266 131)
\lvec(258 131)
\ifill f:0
\move(267 130)
\lvec(269 130)
\lvec(269 131)
\lvec(267 131)
\ifill f:0
\move(270 130)
\lvec(271 130)
\lvec(271 131)
\lvec(270 131)
\ifill f:0
\move(272 130)
\lvec(282 130)
\lvec(282 131)
\lvec(272 131)
\ifill f:0
\move(283 130)
\lvec(290 130)
\lvec(290 131)
\lvec(283 131)
\ifill f:0
\move(291 130)
\lvec(295 130)
\lvec(295 131)
\lvec(291 131)
\ifill f:0
\move(296 130)
\lvec(311 130)
\lvec(311 131)
\lvec(296 131)
\ifill f:0
\move(312 130)
\lvec(314 130)
\lvec(314 131)
\lvec(312 131)
\ifill f:0
\move(315 130)
\lvec(317 130)
\lvec(317 131)
\lvec(315 131)
\ifill f:0
\move(318 130)
\lvec(320 130)
\lvec(320 131)
\lvec(318 131)
\ifill f:0
\move(321 130)
\lvec(323 130)
\lvec(323 131)
\lvec(321 131)
\ifill f:0
\move(324 130)
\lvec(325 130)
\lvec(325 131)
\lvec(324 131)
\ifill f:0
\move(326 130)
\lvec(328 130)
\lvec(328 131)
\lvec(326 131)
\ifill f:0
\move(329 130)
\lvec(344 130)
\lvec(344 131)
\lvec(329 131)
\ifill f:0
\move(345 130)
\lvec(346 130)
\lvec(346 131)
\lvec(345 131)
\ifill f:0
\move(347 130)
\lvec(362 130)
\lvec(362 131)
\lvec(347 131)
\ifill f:0
\move(364 130)
\lvec(368 130)
\lvec(368 131)
\lvec(364 131)
\ifill f:0
\move(369 130)
\lvec(376 130)
\lvec(376 131)
\lvec(369 131)
\ifill f:0
\move(377 130)
\lvec(391 130)
\lvec(391 131)
\lvec(377 131)
\ifill f:0
\move(392 130)
\lvec(394 130)
\lvec(394 131)
\lvec(392 131)
\ifill f:0
\move(395 130)
\lvec(398 130)
\lvec(398 131)
\lvec(395 131)
\ifill f:0
\move(399 130)
\lvec(401 130)
\lvec(401 131)
\lvec(399 131)
\ifill f:0
\move(403 130)
\lvec(406 130)
\lvec(406 131)
\lvec(403 131)
\ifill f:0
\move(407 130)
\lvec(415 130)
\lvec(415 131)
\lvec(407 131)
\ifill f:0
\move(416 130)
\lvec(426 130)
\lvec(426 131)
\lvec(416 131)
\ifill f:0
\move(427 130)
\lvec(440 130)
\lvec(440 131)
\lvec(427 131)
\ifill f:0
\move(441 130)
\lvec(442 130)
\lvec(442 131)
\lvec(441 131)
\ifill f:0
\move(443 130)
\lvec(450 130)
\lvec(450 131)
\lvec(443 131)
\ifill f:0
\move(16 131)
\lvec(17 131)
\lvec(17 132)
\lvec(16 132)
\ifill f:0
\move(20 131)
\lvec(21 131)
\lvec(21 132)
\lvec(20 132)
\ifill f:0
\move(23 131)
\lvec(24 131)
\lvec(24 132)
\lvec(23 132)
\ifill f:0
\move(25 131)
\lvec(26 131)
\lvec(26 132)
\lvec(25 132)
\ifill f:0
\move(36 131)
\lvec(37 131)
\lvec(37 132)
\lvec(36 132)
\ifill f:0
\move(38 131)
\lvec(39 131)
\lvec(39 132)
\lvec(38 132)
\ifill f:0
\move(40 131)
\lvec(41 131)
\lvec(41 132)
\lvec(40 132)
\ifill f:0
\move(42 131)
\lvec(43 131)
\lvec(43 132)
\lvec(42 132)
\ifill f:0
\move(44 131)
\lvec(45 131)
\lvec(45 132)
\lvec(44 132)
\ifill f:0
\move(47 131)
\lvec(50 131)
\lvec(50 132)
\lvec(47 132)
\ifill f:0
\move(57 131)
\lvec(58 131)
\lvec(58 132)
\lvec(57 132)
\ifill f:0
\move(59 131)
\lvec(62 131)
\lvec(62 132)
\lvec(59 132)
\ifill f:0
\move(63 131)
\lvec(65 131)
\lvec(65 132)
\lvec(63 132)
\ifill f:0
\move(66 131)
\lvec(67 131)
\lvec(67 132)
\lvec(66 132)
\ifill f:0
\move(68 131)
\lvec(71 131)
\lvec(71 132)
\lvec(68 132)
\ifill f:0
\move(72 131)
\lvec(74 131)
\lvec(74 132)
\lvec(72 132)
\ifill f:0
\move(75 131)
\lvec(76 131)
\lvec(76 132)
\lvec(75 132)
\ifill f:0
\move(77 131)
\lvec(78 131)
\lvec(78 132)
\lvec(77 132)
\ifill f:0
\move(79 131)
\lvec(82 131)
\lvec(82 132)
\lvec(79 132)
\ifill f:0
\move(83 131)
\lvec(84 131)
\lvec(84 132)
\lvec(83 132)
\ifill f:0
\move(87 131)
\lvec(91 131)
\lvec(91 132)
\lvec(87 132)
\ifill f:0
\move(95 131)
\lvec(96 131)
\lvec(96 132)
\lvec(95 132)
\ifill f:0
\move(97 131)
\lvec(99 131)
\lvec(99 132)
\lvec(97 132)
\ifill f:0
\move(100 131)
\lvec(101 131)
\lvec(101 132)
\lvec(100 132)
\ifill f:0
\move(102 131)
\lvec(105 131)
\lvec(105 132)
\lvec(102 132)
\ifill f:0
\move(106 131)
\lvec(111 131)
\lvec(111 132)
\lvec(106 132)
\ifill f:0
\move(112 131)
\lvec(116 131)
\lvec(116 132)
\lvec(112 132)
\ifill f:0
\move(117 131)
\lvec(118 131)
\lvec(118 132)
\lvec(117 132)
\ifill f:0
\move(120 131)
\lvec(122 131)
\lvec(122 132)
\lvec(120 132)
\ifill f:0
\move(127 131)
\lvec(130 131)
\lvec(130 132)
\lvec(127 132)
\ifill f:0
\move(133 131)
\lvec(134 131)
\lvec(134 132)
\lvec(133 132)
\ifill f:0
\move(135 131)
\lvec(139 131)
\lvec(139 132)
\lvec(135 132)
\ifill f:0
\move(140 131)
\lvec(141 131)
\lvec(141 132)
\lvec(140 132)
\ifill f:0
\move(142 131)
\lvec(145 131)
\lvec(145 132)
\lvec(142 132)
\ifill f:0
\move(146 131)
\lvec(147 131)
\lvec(147 132)
\lvec(146 132)
\ifill f:0
\move(148 131)
\lvec(160 131)
\lvec(160 132)
\lvec(148 132)
\ifill f:0
\move(161 131)
\lvec(167 131)
\lvec(167 132)
\lvec(161 132)
\ifill f:0
\move(169 131)
\lvec(170 131)
\lvec(170 132)
\lvec(169 132)
\ifill f:0
\move(171 131)
\lvec(177 131)
\lvec(177 132)
\lvec(171 132)
\ifill f:0
\move(178 131)
\lvec(179 131)
\lvec(179 132)
\lvec(178 132)
\ifill f:0
\move(181 131)
\lvec(186 131)
\lvec(186 132)
\lvec(181 132)
\ifill f:0
\move(187 131)
\lvec(191 131)
\lvec(191 132)
\lvec(187 132)
\ifill f:0
\move(192 131)
\lvec(195 131)
\lvec(195 132)
\lvec(192 132)
\ifill f:0
\move(196 131)
\lvec(197 131)
\lvec(197 132)
\lvec(196 132)
\ifill f:0
\move(198 131)
\lvec(201 131)
\lvec(201 132)
\lvec(198 132)
\ifill f:0
\move(202 131)
\lvec(220 131)
\lvec(220 132)
\lvec(202 132)
\ifill f:0
\move(221 131)
\lvec(226 131)
\lvec(226 132)
\lvec(221 132)
\ifill f:0
\move(227 131)
\lvec(237 131)
\lvec(237 132)
\lvec(227 132)
\ifill f:0
\move(238 131)
\lvec(241 131)
\lvec(241 132)
\lvec(238 132)
\ifill f:0
\move(242 131)
\lvec(246 131)
\lvec(246 132)
\lvec(242 132)
\ifill f:0
\move(247 131)
\lvec(250 131)
\lvec(250 132)
\lvec(247 132)
\ifill f:0
\move(253 131)
\lvec(257 131)
\lvec(257 132)
\lvec(253 132)
\ifill f:0
\move(258 131)
\lvec(261 131)
\lvec(261 132)
\lvec(258 132)
\ifill f:0
\move(263 131)
\lvec(279 131)
\lvec(279 132)
\lvec(263 132)
\ifill f:0
\move(281 131)
\lvec(290 131)
\lvec(290 132)
\lvec(281 132)
\ifill f:0
\move(291 131)
\lvec(296 131)
\lvec(296 132)
\lvec(291 132)
\ifill f:0
\move(297 131)
\lvec(310 131)
\lvec(310 132)
\lvec(297 132)
\ifill f:0
\move(311 131)
\lvec(313 131)
\lvec(313 132)
\lvec(311 132)
\ifill f:0
\move(314 131)
\lvec(320 131)
\lvec(320 132)
\lvec(314 132)
\ifill f:0
\move(321 131)
\lvec(323 131)
\lvec(323 132)
\lvec(321 132)
\ifill f:0
\move(324 131)
\lvec(325 131)
\lvec(325 132)
\lvec(324 132)
\ifill f:0
\move(326 131)
\lvec(345 131)
\lvec(345 132)
\lvec(326 132)
\ifill f:0
\move(346 131)
\lvec(362 131)
\lvec(362 132)
\lvec(346 132)
\ifill f:0
\move(364 131)
\lvec(365 131)
\lvec(365 132)
\lvec(364 132)
\ifill f:0
\move(366 131)
\lvec(370 131)
\lvec(370 132)
\lvec(366 132)
\ifill f:0
\move(371 131)
\lvec(375 131)
\lvec(375 132)
\lvec(371 132)
\ifill f:0
\move(376 131)
\lvec(380 131)
\lvec(380 132)
\lvec(376 132)
\ifill f:0
\move(381 131)
\lvec(383 131)
\lvec(383 132)
\lvec(381 132)
\ifill f:0
\move(384 131)
\lvec(395 131)
\lvec(395 132)
\lvec(384 132)
\ifill f:0
\move(396 131)
\lvec(401 131)
\lvec(401 132)
\lvec(396 132)
\ifill f:0
\move(402 131)
\lvec(434 131)
\lvec(434 132)
\lvec(402 132)
\ifill f:0
\move(435 131)
\lvec(442 131)
\lvec(442 132)
\lvec(435 132)
\ifill f:0
\move(443 131)
\lvec(449 131)
\lvec(449 132)
\lvec(443 132)
\ifill f:0
\move(450 131)
\lvec(451 131)
\lvec(451 132)
\lvec(450 132)
\ifill f:0
\move(16 132)
\lvec(17 132)
\lvec(17 133)
\lvec(16 133)
\ifill f:0
\move(19 132)
\lvec(21 132)
\lvec(21 133)
\lvec(19 133)
\ifill f:0
\move(25 132)
\lvec(26 132)
\lvec(26 133)
\lvec(25 133)
\ifill f:0
\move(36 132)
\lvec(37 132)
\lvec(37 133)
\lvec(36 133)
\ifill f:0
\move(38 132)
\lvec(39 132)
\lvec(39 133)
\lvec(38 133)
\ifill f:0
\move(40 132)
\lvec(41 132)
\lvec(41 133)
\lvec(40 133)
\ifill f:0
\move(42 132)
\lvec(45 132)
\lvec(45 133)
\lvec(42 133)
\ifill f:0
\move(48 132)
\lvec(50 132)
\lvec(50 133)
\lvec(48 133)
\ifill f:0
\move(54 132)
\lvec(55 132)
\lvec(55 133)
\lvec(54 133)
\ifill f:0
\move(56 132)
\lvec(57 132)
\lvec(57 133)
\lvec(56 133)
\ifill f:0
\move(59 132)
\lvec(60 132)
\lvec(60 133)
\lvec(59 133)
\ifill f:0
\move(61 132)
\lvec(62 132)
\lvec(62 133)
\lvec(61 133)
\ifill f:0
\move(63 132)
\lvec(65 132)
\lvec(65 133)
\lvec(63 133)
\ifill f:0
\move(67 132)
\lvec(71 132)
\lvec(71 133)
\lvec(67 133)
\ifill f:0
\move(72 132)
\lvec(75 132)
\lvec(75 133)
\lvec(72 133)
\ifill f:0
\move(76 132)
\lvec(78 132)
\lvec(78 133)
\lvec(76 133)
\ifill f:0
\move(79 132)
\lvec(82 132)
\lvec(82 133)
\lvec(79 133)
\ifill f:0
\move(85 132)
\lvec(87 132)
\lvec(87 133)
\lvec(85 133)
\ifill f:0
\move(88 132)
\lvec(90 132)
\lvec(90 133)
\lvec(88 133)
\ifill f:0
\move(91 132)
\lvec(93 132)
\lvec(93 133)
\lvec(91 133)
\ifill f:0
\move(97 132)
\lvec(99 132)
\lvec(99 133)
\lvec(97 133)
\ifill f:0
\move(100 132)
\lvec(101 132)
\lvec(101 133)
\lvec(100 133)
\ifill f:0
\move(102 132)
\lvec(106 132)
\lvec(106 133)
\lvec(102 133)
\ifill f:0
\move(107 132)
\lvec(109 132)
\lvec(109 133)
\lvec(107 133)
\ifill f:0
\move(110 132)
\lvec(113 132)
\lvec(113 133)
\lvec(110 133)
\ifill f:0
\move(114 132)
\lvec(119 132)
\lvec(119 133)
\lvec(114 133)
\ifill f:0
\move(121 132)
\lvec(122 132)
\lvec(122 133)
\lvec(121 133)
\ifill f:0
\move(128 132)
\lvec(132 132)
\lvec(132 133)
\lvec(128 133)
\ifill f:0
\move(133 132)
\lvec(136 132)
\lvec(136 133)
\lvec(133 133)
\ifill f:0
\move(137 132)
\lvec(143 132)
\lvec(143 133)
\lvec(137 133)
\ifill f:0
\move(144 132)
\lvec(145 132)
\lvec(145 133)
\lvec(144 133)
\ifill f:0
\move(146 132)
\lvec(147 132)
\lvec(147 133)
\lvec(146 133)
\ifill f:0
\move(148 132)
\lvec(149 132)
\lvec(149 133)
\lvec(148 133)
\ifill f:0
\move(150 132)
\lvec(158 132)
\lvec(158 133)
\lvec(150 133)
\ifill f:0
\move(159 132)
\lvec(162 132)
\lvec(162 133)
\lvec(159 133)
\ifill f:0
\move(163 132)
\lvec(168 132)
\lvec(168 133)
\lvec(163 133)
\ifill f:0
\move(169 132)
\lvec(170 132)
\lvec(170 133)
\lvec(169 133)
\ifill f:0
\move(171 132)
\lvec(183 132)
\lvec(183 133)
\lvec(171 133)
\ifill f:0
\move(184 132)
\lvec(186 132)
\lvec(186 133)
\lvec(184 133)
\ifill f:0
\move(187 132)
\lvec(190 132)
\lvec(190 133)
\lvec(187 133)
\ifill f:0
\move(192 132)
\lvec(194 132)
\lvec(194 133)
\lvec(192 133)
\ifill f:0
\move(196 132)
\lvec(197 132)
\lvec(197 133)
\lvec(196 133)
\ifill f:0
\move(198 132)
\lvec(201 132)
\lvec(201 133)
\lvec(198 133)
\ifill f:0
\move(202 132)
\lvec(204 132)
\lvec(204 133)
\lvec(202 133)
\ifill f:0
\move(205 132)
\lvec(226 132)
\lvec(226 133)
\lvec(205 133)
\ifill f:0
\move(228 132)
\lvec(230 132)
\lvec(230 133)
\lvec(228 133)
\ifill f:0
\move(231 132)
\lvec(235 132)
\lvec(235 133)
\lvec(231 133)
\ifill f:0
\move(237 132)
\lvec(239 132)
\lvec(239 133)
\lvec(237 133)
\ifill f:0
\move(240 132)
\lvec(242 132)
\lvec(242 133)
\lvec(240 133)
\ifill f:0
\move(244 132)
\lvec(247 132)
\lvec(247 133)
\lvec(244 133)
\ifill f:0
\move(248 132)
\lvec(251 132)
\lvec(251 133)
\lvec(248 133)
\ifill f:0
\move(252 132)
\lvec(257 132)
\lvec(257 133)
\lvec(252 133)
\ifill f:0
\move(258 132)
\lvec(260 132)
\lvec(260 133)
\lvec(258 133)
\ifill f:0
\move(261 132)
\lvec(290 132)
\lvec(290 133)
\lvec(261 133)
\ifill f:0
\move(292 132)
\lvec(297 132)
\lvec(297 133)
\lvec(292 133)
\ifill f:0
\move(299 132)
\lvec(323 132)
\lvec(323 133)
\lvec(299 133)
\ifill f:0
\move(324 132)
\lvec(325 132)
\lvec(325 133)
\lvec(324 133)
\ifill f:0
\move(326 132)
\lvec(344 132)
\lvec(344 133)
\lvec(326 133)
\ifill f:0
\move(345 132)
\lvec(355 132)
\lvec(355 133)
\lvec(345 133)
\ifill f:0
\move(356 132)
\lvec(357 132)
\lvec(357 133)
\lvec(356 133)
\ifill f:0
\move(358 132)
\lvec(359 132)
\lvec(359 133)
\lvec(358 133)
\ifill f:0
\move(360 132)
\lvec(362 132)
\lvec(362 133)
\lvec(360 133)
\ifill f:0
\move(364 132)
\lvec(365 132)
\lvec(365 133)
\lvec(364 133)
\ifill f:0
\move(366 132)
\lvec(367 132)
\lvec(367 133)
\lvec(366 133)
\ifill f:0
\move(368 132)
\lvec(379 132)
\lvec(379 133)
\lvec(368 133)
\ifill f:0
\move(380 132)
\lvec(384 132)
\lvec(384 133)
\lvec(380 133)
\ifill f:0
\move(385 132)
\lvec(392 132)
\lvec(392 133)
\lvec(385 133)
\ifill f:0
\move(393 132)
\lvec(395 132)
\lvec(395 133)
\lvec(393 133)
\ifill f:0
\move(396 132)
\lvec(401 132)
\lvec(401 133)
\lvec(396 133)
\ifill f:0
\move(402 132)
\lvec(405 132)
\lvec(405 133)
\lvec(402 133)
\ifill f:0
\move(406 132)
\lvec(412 132)
\lvec(412 133)
\lvec(406 133)
\ifill f:0
\move(413 132)
\lvec(416 132)
\lvec(416 133)
\lvec(413 133)
\ifill f:0
\move(417 132)
\lvec(435 132)
\lvec(435 133)
\lvec(417 133)
\ifill f:0
\move(436 132)
\lvec(442 132)
\lvec(442 133)
\lvec(436 133)
\ifill f:0
\move(443 132)
\lvec(448 132)
\lvec(448 133)
\lvec(443 133)
\ifill f:0
\move(449 132)
\lvec(451 132)
\lvec(451 133)
\lvec(449 133)
\ifill f:0
\move(16 133)
\lvec(17 133)
\lvec(17 134)
\lvec(16 134)
\ifill f:0
\move(20 133)
\lvec(21 133)
\lvec(21 134)
\lvec(20 134)
\ifill f:0
\move(24 133)
\lvec(26 133)
\lvec(26 134)
\lvec(24 134)
\ifill f:0
\move(36 133)
\lvec(37 133)
\lvec(37 134)
\lvec(36 134)
\ifill f:0
\move(38 133)
\lvec(39 133)
\lvec(39 134)
\lvec(38 134)
\ifill f:0
\move(40 133)
\lvec(43 133)
\lvec(43 134)
\lvec(40 134)
\ifill f:0
\move(44 133)
\lvec(46 133)
\lvec(46 134)
\lvec(44 134)
\ifill f:0
\move(48 133)
\lvec(50 133)
\lvec(50 134)
\lvec(48 134)
\ifill f:0
\move(52 133)
\lvec(53 133)
\lvec(53 134)
\lvec(52 134)
\ifill f:0
\move(61 133)
\lvec(65 133)
\lvec(65 134)
\lvec(61 134)
\ifill f:0
\move(66 133)
\lvec(74 133)
\lvec(74 134)
\lvec(66 134)
\ifill f:0
\move(76 133)
\lvec(77 133)
\lvec(77 134)
\lvec(76 134)
\ifill f:0
\move(78 133)
\lvec(82 133)
\lvec(82 134)
\lvec(78 134)
\ifill f:0
\move(84 133)
\lvec(85 133)
\lvec(85 134)
\lvec(84 134)
\ifill f:0
\move(91 133)
\lvec(92 133)
\lvec(92 134)
\lvec(91 134)
\ifill f:0
\move(97 133)
\lvec(99 133)
\lvec(99 134)
\lvec(97 134)
\ifill f:0
\move(100 133)
\lvec(101 133)
\lvec(101 134)
\lvec(100 134)
\ifill f:0
\move(102 133)
\lvec(104 133)
\lvec(104 134)
\lvec(102 134)
\ifill f:0
\move(105 133)
\lvec(106 133)
\lvec(106 134)
\lvec(105 134)
\ifill f:0
\move(107 133)
\lvec(120 133)
\lvec(120 134)
\lvec(107 134)
\ifill f:0
\move(121 133)
\lvec(122 133)
\lvec(122 134)
\lvec(121 134)
\ifill f:0
\move(123 133)
\lvec(127 133)
\lvec(127 134)
\lvec(123 134)
\ifill f:0
\move(133 133)
\lvec(134 133)
\lvec(134 134)
\lvec(133 134)
\ifill f:0
\move(135 133)
\lvec(138 133)
\lvec(138 134)
\lvec(135 134)
\ifill f:0
\move(139 133)
\lvec(143 133)
\lvec(143 134)
\lvec(139 134)
\ifill f:0
\move(144 133)
\lvec(145 133)
\lvec(145 134)
\lvec(144 134)
\ifill f:0
\move(146 133)
\lvec(151 133)
\lvec(151 134)
\lvec(146 134)
\ifill f:0
\move(152 133)
\lvec(159 133)
\lvec(159 134)
\lvec(152 134)
\ifill f:0
\move(160 133)
\lvec(163 133)
\lvec(163 134)
\lvec(160 134)
\ifill f:0
\move(164 133)
\lvec(168 133)
\lvec(168 134)
\lvec(164 134)
\ifill f:0
\move(169 133)
\lvec(170 133)
\lvec(170 134)
\lvec(169 134)
\ifill f:0
\move(171 133)
\lvec(176 133)
\lvec(176 134)
\lvec(171 134)
\ifill f:0
\move(177 133)
\lvec(188 133)
\lvec(188 134)
\lvec(177 134)
\ifill f:0
\move(190 133)
\lvec(194 133)
\lvec(194 134)
\lvec(190 134)
\ifill f:0
\move(196 133)
\lvec(197 133)
\lvec(197 134)
\lvec(196 134)
\ifill f:0
\move(199 133)
\lvec(202 133)
\lvec(202 134)
\lvec(199 134)
\ifill f:0
\move(203 133)
\lvec(215 133)
\lvec(215 134)
\lvec(203 134)
\ifill f:0
\move(216 133)
\lvec(226 133)
\lvec(226 134)
\lvec(216 134)
\ifill f:0
\move(228 133)
\lvec(232 133)
\lvec(232 134)
\lvec(228 134)
\ifill f:0
\move(233 133)
\lvec(241 133)
\lvec(241 134)
\lvec(233 134)
\ifill f:0
\move(242 133)
\lvec(248 133)
\lvec(248 134)
\lvec(242 134)
\ifill f:0
\move(249 133)
\lvec(253 133)
\lvec(253 134)
\lvec(249 134)
\ifill f:0
\move(254 133)
\lvec(257 133)
\lvec(257 134)
\lvec(254 134)
\ifill f:0
\move(258 133)
\lvec(259 133)
\lvec(259 134)
\lvec(258 134)
\ifill f:0
\move(260 133)
\lvec(266 133)
\lvec(266 134)
\lvec(260 134)
\ifill f:0
\move(267 133)
\lvec(268 133)
\lvec(268 134)
\lvec(267 134)
\ifill f:0
\move(269 133)
\lvec(290 133)
\lvec(290 134)
\lvec(269 134)
\ifill f:0
\move(292 133)
\lvec(299 133)
\lvec(299 134)
\lvec(292 134)
\ifill f:0
\move(300 133)
\lvec(305 133)
\lvec(305 134)
\lvec(300 134)
\ifill f:0
\move(306 133)
\lvec(310 133)
\lvec(310 134)
\lvec(306 134)
\ifill f:0
\move(311 133)
\lvec(315 133)
\lvec(315 134)
\lvec(311 134)
\ifill f:0
\move(316 133)
\lvec(319 133)
\lvec(319 134)
\lvec(316 134)
\ifill f:0
\move(320 133)
\lvec(323 133)
\lvec(323 134)
\lvec(320 134)
\ifill f:0
\move(324 133)
\lvec(325 133)
\lvec(325 134)
\lvec(324 134)
\ifill f:0
\move(326 133)
\lvec(329 133)
\lvec(329 134)
\lvec(326 134)
\ifill f:0
\move(330 133)
\lvec(332 133)
\lvec(332 134)
\lvec(330 134)
\ifill f:0
\move(333 133)
\lvec(335 133)
\lvec(335 134)
\lvec(333 134)
\ifill f:0
\move(336 133)
\lvec(338 133)
\lvec(338 134)
\lvec(336 134)
\ifill f:0
\move(339 133)
\lvec(343 133)
\lvec(343 134)
\lvec(339 134)
\ifill f:0
\move(344 133)
\lvec(359 133)
\lvec(359 134)
\lvec(344 134)
\ifill f:0
\move(360 133)
\lvec(362 133)
\lvec(362 134)
\lvec(360 134)
\ifill f:0
\move(364 133)
\lvec(365 133)
\lvec(365 134)
\lvec(364 134)
\ifill f:0
\move(366 133)
\lvec(367 133)
\lvec(367 134)
\lvec(366 134)
\ifill f:0
\move(368 133)
\lvec(369 133)
\lvec(369 134)
\lvec(368 134)
\ifill f:0
\move(370 133)
\lvec(371 133)
\lvec(371 134)
\lvec(370 134)
\ifill f:0
\move(372 133)
\lvec(380 133)
\lvec(380 134)
\lvec(372 134)
\ifill f:0
\move(381 133)
\lvec(385 133)
\lvec(385 134)
\lvec(381 134)
\ifill f:0
\move(386 133)
\lvec(401 133)
\lvec(401 134)
\lvec(386 134)
\ifill f:0
\move(402 133)
\lvec(411 133)
\lvec(411 134)
\lvec(402 134)
\ifill f:0
\move(412 133)
\lvec(442 133)
\lvec(442 134)
\lvec(412 134)
\ifill f:0
\move(443 133)
\lvec(447 133)
\lvec(447 134)
\lvec(443 134)
\ifill f:0
\move(448 133)
\lvec(451 133)
\lvec(451 134)
\lvec(448 134)
\ifill f:0
\move(16 134)
\lvec(17 134)
\lvec(17 135)
\lvec(16 135)
\ifill f:0
\move(24 134)
\lvec(26 134)
\lvec(26 135)
\lvec(24 135)
\ifill f:0
\move(36 134)
\lvec(37 134)
\lvec(37 135)
\lvec(36 135)
\ifill f:0
\move(38 134)
\lvec(39 134)
\lvec(39 135)
\lvec(38 135)
\ifill f:0
\move(40 134)
\lvec(41 134)
\lvec(41 135)
\lvec(40 135)
\ifill f:0
\move(42 134)
\lvec(45 134)
\lvec(45 135)
\lvec(42 135)
\ifill f:0
\move(47 134)
\lvec(50 134)
\lvec(50 135)
\lvec(47 135)
\ifill f:0
\move(51 134)
\lvec(52 134)
\lvec(52 135)
\lvec(51 135)
\ifill f:0
\move(60 134)
\lvec(61 134)
\lvec(61 135)
\lvec(60 135)
\ifill f:0
\move(62 134)
\lvec(65 134)
\lvec(65 135)
\lvec(62 135)
\ifill f:0
\move(66 134)
\lvec(71 134)
\lvec(71 135)
\lvec(66 135)
\ifill f:0
\move(72 134)
\lvec(73 134)
\lvec(73 135)
\lvec(72 135)
\ifill f:0
\move(75 134)
\lvec(77 134)
\lvec(77 135)
\lvec(75 135)
\ifill f:0
\move(78 134)
\lvec(79 134)
\lvec(79 135)
\lvec(78 135)
\ifill f:0
\move(80 134)
\lvec(82 134)
\lvec(82 135)
\lvec(80 135)
\ifill f:0
\move(84 134)
\lvec(85 134)
\lvec(85 135)
\lvec(84 135)
\ifill f:0
\move(87 134)
\lvec(93 134)
\lvec(93 135)
\lvec(87 135)
\ifill f:0
\move(96 134)
\lvec(98 134)
\lvec(98 135)
\lvec(96 135)
\ifill f:0
\move(100 134)
\lvec(101 134)
\lvec(101 135)
\lvec(100 135)
\ifill f:0
\move(102 134)
\lvec(109 134)
\lvec(109 135)
\lvec(102 135)
\ifill f:0
\move(110 134)
\lvec(115 134)
\lvec(115 135)
\lvec(110 135)
\ifill f:0
\move(116 134)
\lvec(120 134)
\lvec(120 135)
\lvec(116 135)
\ifill f:0
\move(121 134)
\lvec(122 134)
\lvec(122 135)
\lvec(121 135)
\ifill f:0
\move(123 134)
\lvec(131 134)
\lvec(131 135)
\lvec(123 135)
\ifill f:0
\move(134 134)
\lvec(138 134)
\lvec(138 135)
\lvec(134 135)
\ifill f:0
\move(139 134)
\lvec(140 134)
\lvec(140 135)
\lvec(139 135)
\ifill f:0
\move(141 134)
\lvec(143 134)
\lvec(143 135)
\lvec(141 135)
\ifill f:0
\move(144 134)
\lvec(145 134)
\lvec(145 135)
\lvec(144 135)
\ifill f:0
\move(146 134)
\lvec(163 134)
\lvec(163 135)
\lvec(146 135)
\ifill f:0
\move(165 134)
\lvec(168 134)
\lvec(168 135)
\lvec(165 135)
\ifill f:0
\move(169 134)
\lvec(170 134)
\lvec(170 135)
\lvec(169 135)
\ifill f:0
\move(171 134)
\lvec(173 134)
\lvec(173 135)
\lvec(171 135)
\ifill f:0
\move(178 134)
\lvec(179 134)
\lvec(179 135)
\lvec(178 135)
\ifill f:0
\move(181 134)
\lvec(182 134)
\lvec(182 135)
\lvec(181 135)
\ifill f:0
\move(187 134)
\lvec(194 134)
\lvec(194 135)
\lvec(187 135)
\ifill f:0
\move(196 134)
\lvec(197 134)
\lvec(197 135)
\lvec(196 135)
\ifill f:0
\move(198 134)
\lvec(219 134)
\lvec(219 135)
\lvec(198 135)
\ifill f:0
\move(220 134)
\lvec(223 134)
\lvec(223 135)
\lvec(220 135)
\ifill f:0
\move(224 134)
\lvec(226 134)
\lvec(226 135)
\lvec(224 135)
\ifill f:0
\move(228 134)
\lvec(229 134)
\lvec(229 135)
\lvec(228 135)
\ifill f:0
\move(230 134)
\lvec(234 134)
\lvec(234 135)
\lvec(230 135)
\ifill f:0
\move(235 134)
\lvec(242 134)
\lvec(242 135)
\lvec(235 135)
\ifill f:0
\move(243 134)
\lvec(249 134)
\lvec(249 135)
\lvec(243 135)
\ifill f:0
\move(250 134)
\lvec(253 134)
\lvec(253 135)
\lvec(250 135)
\ifill f:0
\move(254 134)
\lvec(257 134)
\lvec(257 135)
\lvec(254 135)
\ifill f:0
\move(258 134)
\lvec(265 134)
\lvec(265 135)
\lvec(258 135)
\ifill f:0
\move(266 134)
\lvec(275 134)
\lvec(275 135)
\lvec(266 135)
\ifill f:0
\move(276 134)
\lvec(290 134)
\lvec(290 135)
\lvec(276 135)
\ifill f:0
\move(294 134)
\lvec(301 134)
\lvec(301 135)
\lvec(294 135)
\ifill f:0
\move(303 134)
\lvec(308 134)
\lvec(308 135)
\lvec(303 135)
\ifill f:0
\move(309 134)
\lvec(314 134)
\lvec(314 135)
\lvec(309 135)
\ifill f:0
\move(315 134)
\lvec(318 134)
\lvec(318 135)
\lvec(315 135)
\ifill f:0
\move(319 134)
\lvec(322 134)
\lvec(322 135)
\lvec(319 135)
\ifill f:0
\move(323 134)
\lvec(325 134)
\lvec(325 135)
\lvec(323 135)
\ifill f:0
\move(326 134)
\lvec(347 134)
\lvec(347 135)
\lvec(326 135)
\ifill f:0
\move(348 134)
\lvec(352 134)
\lvec(352 135)
\lvec(348 135)
\ifill f:0
\move(353 134)
\lvec(359 134)
\lvec(359 135)
\lvec(353 135)
\ifill f:0
\move(360 134)
\lvec(362 134)
\lvec(362 135)
\lvec(360 135)
\ifill f:0
\move(363 134)
\lvec(377 134)
\lvec(377 135)
\lvec(363 135)
\ifill f:0
\move(378 134)
\lvec(379 134)
\lvec(379 135)
\lvec(378 135)
\ifill f:0
\move(380 134)
\lvec(386 134)
\lvec(386 135)
\lvec(380 135)
\ifill f:0
\move(387 134)
\lvec(396 134)
\lvec(396 135)
\lvec(387 135)
\ifill f:0
\move(397 134)
\lvec(401 134)
\lvec(401 135)
\lvec(397 135)
\ifill f:0
\move(402 134)
\lvec(436 134)
\lvec(436 135)
\lvec(402 135)
\ifill f:0
\move(437 134)
\lvec(442 134)
\lvec(442 135)
\lvec(437 135)
\ifill f:0
\move(443 134)
\lvec(445 134)
\lvec(445 135)
\lvec(443 135)
\ifill f:0
\move(447 134)
\lvec(451 134)
\lvec(451 135)
\lvec(447 135)
\ifill f:0
\move(15 135)
\lvec(17 135)
\lvec(17 136)
\lvec(15 136)
\ifill f:0
\move(19 135)
\lvec(21 135)
\lvec(21 136)
\lvec(19 136)
\ifill f:0
\move(22 135)
\lvec(23 135)
\lvec(23 136)
\lvec(22 136)
\ifill f:0
\move(24 135)
\lvec(26 135)
\lvec(26 136)
\lvec(24 136)
\ifill f:0
\move(36 135)
\lvec(37 135)
\lvec(37 136)
\lvec(36 136)
\ifill f:0
\move(40 135)
\lvec(45 135)
\lvec(45 136)
\lvec(40 136)
\ifill f:0
\move(47 135)
\lvec(50 135)
\lvec(50 136)
\lvec(47 136)
\ifill f:0
\move(51 135)
\lvec(52 135)
\lvec(52 136)
\lvec(51 136)
\ifill f:0
\move(54 135)
\lvec(55 135)
\lvec(55 136)
\lvec(54 136)
\ifill f:0
\move(59 135)
\lvec(61 135)
\lvec(61 136)
\lvec(59 136)
\ifill f:0
\move(62 135)
\lvec(63 135)
\lvec(63 136)
\lvec(62 136)
\ifill f:0
\move(64 135)
\lvec(65 135)
\lvec(65 136)
\lvec(64 136)
\ifill f:0
\move(66 135)
\lvec(75 135)
\lvec(75 136)
\lvec(66 136)
\ifill f:0
\move(76 135)
\lvec(79 135)
\lvec(79 136)
\lvec(76 136)
\ifill f:0
\move(80 135)
\lvec(82 135)
\lvec(82 136)
\lvec(80 136)
\ifill f:0
\move(84 135)
\lvec(85 135)
\lvec(85 136)
\lvec(84 136)
\ifill f:0
\move(86 135)
\lvec(87 135)
\lvec(87 136)
\lvec(86 136)
\ifill f:0
\move(88 135)
\lvec(89 135)
\lvec(89 136)
\lvec(88 136)
\ifill f:0
\move(91 135)
\lvec(93 135)
\lvec(93 136)
\lvec(91 136)
\ifill f:0
\move(95 135)
\lvec(98 135)
\lvec(98 136)
\lvec(95 136)
\ifill f:0
\move(100 135)
\lvec(101 135)
\lvec(101 136)
\lvec(100 136)
\ifill f:0
\move(103 135)
\lvec(105 135)
\lvec(105 136)
\lvec(103 136)
\ifill f:0
\move(107 135)
\lvec(113 135)
\lvec(113 136)
\lvec(107 136)
\ifill f:0
\move(114 135)
\lvec(116 135)
\lvec(116 136)
\lvec(114 136)
\ifill f:0
\move(117 135)
\lvec(122 135)
\lvec(122 136)
\lvec(117 136)
\ifill f:0
\move(123 135)
\lvec(127 135)
\lvec(127 136)
\lvec(123 136)
\ifill f:0
\move(128 135)
\lvec(132 135)
\lvec(132 136)
\lvec(128 136)
\ifill f:0
\move(133 135)
\lvec(135 135)
\lvec(135 136)
\lvec(133 136)
\ifill f:0
\move(137 135)
\lvec(143 135)
\lvec(143 136)
\lvec(137 136)
\ifill f:0
\move(144 135)
\lvec(145 135)
\lvec(145 136)
\lvec(144 136)
\ifill f:0
\move(146 135)
\lvec(154 135)
\lvec(154 136)
\lvec(146 136)
\ifill f:0
\move(155 135)
\lvec(170 135)
\lvec(170 136)
\lvec(155 136)
\ifill f:0
\move(171 135)
\lvec(173 135)
\lvec(173 136)
\lvec(171 136)
\ifill f:0
\move(176 135)
\lvec(193 135)
\lvec(193 136)
\lvec(176 136)
\ifill f:0
\move(195 135)
\lvec(197 135)
\lvec(197 136)
\lvec(195 136)
\ifill f:0
\move(198 135)
\lvec(199 135)
\lvec(199 136)
\lvec(198 136)
\ifill f:0
\move(200 135)
\lvec(207 135)
\lvec(207 136)
\lvec(200 136)
\ifill f:0
\move(208 135)
\lvec(223 135)
\lvec(223 136)
\lvec(208 136)
\ifill f:0
\move(224 135)
\lvec(226 135)
\lvec(226 136)
\lvec(224 136)
\ifill f:0
\move(228 135)
\lvec(229 135)
\lvec(229 136)
\lvec(228 136)
\ifill f:0
\move(230 135)
\lvec(231 135)
\lvec(231 136)
\lvec(230 136)
\ifill f:0
\move(232 135)
\lvec(233 135)
\lvec(233 136)
\lvec(232 136)
\ifill f:0
\move(234 135)
\lvec(235 135)
\lvec(235 136)
\lvec(234 136)
\ifill f:0
\move(236 135)
\lvec(238 135)
\lvec(238 136)
\lvec(236 136)
\ifill f:0
\move(239 135)
\lvec(250 135)
\lvec(250 136)
\lvec(239 136)
\ifill f:0
\move(251 135)
\lvec(257 135)
\lvec(257 136)
\lvec(251 136)
\ifill f:0
\move(258 135)
\lvec(269 135)
\lvec(269 136)
\lvec(258 136)
\ifill f:0
\move(271 135)
\lvec(279 135)
\lvec(279 136)
\lvec(271 136)
\ifill f:0
\move(281 135)
\lvec(282 135)
\lvec(282 136)
\lvec(281 136)
\ifill f:0
\move(287 135)
\lvec(288 135)
\lvec(288 136)
\lvec(287 136)
\ifill f:0
\move(289 135)
\lvec(290 135)
\lvec(290 136)
\lvec(289 136)
\ifill f:0
\move(294 135)
\lvec(295 135)
\lvec(295 136)
\lvec(294 136)
\ifill f:0
\move(296 135)
\lvec(305 135)
\lvec(305 136)
\lvec(296 136)
\ifill f:0
\move(306 135)
\lvec(317 135)
\lvec(317 136)
\lvec(306 136)
\ifill f:0
\move(318 135)
\lvec(322 135)
\lvec(322 136)
\lvec(318 136)
\ifill f:0
\move(323 135)
\lvec(325 135)
\lvec(325 136)
\lvec(323 136)
\ifill f:0
\move(327 135)
\lvec(330 135)
\lvec(330 136)
\lvec(327 136)
\ifill f:0
\move(331 135)
\lvec(337 135)
\lvec(337 136)
\lvec(331 136)
\ifill f:0
\move(338 135)
\lvec(340 135)
\lvec(340 136)
\lvec(338 136)
\ifill f:0
\move(341 135)
\lvec(343 135)
\lvec(343 136)
\lvec(341 136)
\ifill f:0
\move(344 135)
\lvec(346 135)
\lvec(346 136)
\lvec(344 136)
\ifill f:0
\move(347 135)
\lvec(362 135)
\lvec(362 136)
\lvec(347 136)
\ifill f:0
\move(363 135)
\lvec(370 135)
\lvec(370 136)
\lvec(363 136)
\ifill f:0
\move(371 135)
\lvec(376 135)
\lvec(376 136)
\lvec(371 136)
\ifill f:0
\move(377 135)
\lvec(380 135)
\lvec(380 136)
\lvec(377 136)
\ifill f:0
\move(381 135)
\lvec(389 135)
\lvec(389 136)
\lvec(381 136)
\ifill f:0
\move(390 135)
\lvec(401 135)
\lvec(401 136)
\lvec(390 136)
\ifill f:0
\move(402 135)
\lvec(404 135)
\lvec(404 136)
\lvec(402 136)
\ifill f:0
\move(405 135)
\lvec(429 135)
\lvec(429 136)
\lvec(405 136)
\ifill f:0
\move(430 135)
\lvec(442 135)
\lvec(442 136)
\lvec(430 136)
\ifill f:0
\move(443 135)
\lvec(445 135)
\lvec(445 136)
\lvec(443 136)
\ifill f:0
\move(446 135)
\lvec(450 135)
\lvec(450 136)
\lvec(446 136)
\ifill f:0
\move(15 136)
\lvec(17 136)
\lvec(17 137)
\lvec(15 137)
\ifill f:0
\move(23 136)
\lvec(26 136)
\lvec(26 137)
\lvec(23 137)
\ifill f:0
\move(36 136)
\lvec(37 136)
\lvec(37 137)
\lvec(36 137)
\ifill f:0
\move(40 136)
\lvec(41 136)
\lvec(41 137)
\lvec(40 137)
\ifill f:0
\move(42 136)
\lvec(45 136)
\lvec(45 137)
\lvec(42 137)
\ifill f:0
\move(47 136)
\lvec(50 136)
\lvec(50 137)
\lvec(47 137)
\ifill f:0
\move(51 136)
\lvec(52 136)
\lvec(52 137)
\lvec(51 137)
\ifill f:0
\move(57 136)
\lvec(58 136)
\lvec(58 137)
\lvec(57 137)
\ifill f:0
\move(61 136)
\lvec(63 136)
\lvec(63 137)
\lvec(61 137)
\ifill f:0
\move(64 136)
\lvec(65 136)
\lvec(65 137)
\lvec(64 137)
\ifill f:0
\move(66 136)
\lvec(67 136)
\lvec(67 137)
\lvec(66 137)
\ifill f:0
\move(68 136)
\lvec(71 136)
\lvec(71 137)
\lvec(68 137)
\ifill f:0
\move(72 136)
\lvec(73 136)
\lvec(73 137)
\lvec(72 137)
\ifill f:0
\move(76 136)
\lvec(79 136)
\lvec(79 137)
\lvec(76 137)
\ifill f:0
\move(80 136)
\lvec(82 136)
\lvec(82 137)
\lvec(80 137)
\ifill f:0
\move(83 136)
\lvec(84 136)
\lvec(84 137)
\lvec(83 137)
\ifill f:0
\move(86 136)
\lvec(87 136)
\lvec(87 137)
\lvec(86 137)
\ifill f:0
\move(89 136)
\lvec(91 136)
\lvec(91 137)
\lvec(89 137)
\ifill f:0
\move(97 136)
\lvec(98 136)
\lvec(98 137)
\lvec(97 137)
\ifill f:0
\move(100 136)
\lvec(101 136)
\lvec(101 137)
\lvec(100 137)
\ifill f:0
\move(102 136)
\lvec(103 136)
\lvec(103 137)
\lvec(102 137)
\ifill f:0
\move(104 136)
\lvec(106 136)
\lvec(106 137)
\lvec(104 137)
\ifill f:0
\move(107 136)
\lvec(117 136)
\lvec(117 137)
\lvec(107 137)
\ifill f:0
\move(118 136)
\lvec(122 136)
\lvec(122 137)
\lvec(118 137)
\ifill f:0
\move(123 136)
\lvec(125 136)
\lvec(125 137)
\lvec(123 137)
\ifill f:0
\move(135 136)
\lvec(139 136)
\lvec(139 137)
\lvec(135 137)
\ifill f:0
\move(140 136)
\lvec(143 136)
\lvec(143 137)
\lvec(140 137)
\ifill f:0
\move(144 136)
\lvec(145 136)
\lvec(145 137)
\lvec(144 137)
\ifill f:0
\move(146 136)
\lvec(151 136)
\lvec(151 137)
\lvec(146 137)
\ifill f:0
\move(152 136)
\lvec(153 136)
\lvec(153 137)
\lvec(152 137)
\ifill f:0
\move(154 136)
\lvec(155 136)
\lvec(155 137)
\lvec(154 137)
\ifill f:0
\move(156 136)
\lvec(162 136)
\lvec(162 137)
\lvec(156 137)
\ifill f:0
\move(163 136)
\lvec(165 136)
\lvec(165 137)
\lvec(163 137)
\ifill f:0
\move(166 136)
\lvec(170 136)
\lvec(170 137)
\lvec(166 137)
\ifill f:0
\move(171 136)
\lvec(172 136)
\lvec(172 137)
\lvec(171 137)
\ifill f:0
\move(174 136)
\lvec(179 136)
\lvec(179 137)
\lvec(174 137)
\ifill f:0
\move(181 136)
\lvec(182 136)
\lvec(182 137)
\lvec(181 137)
\ifill f:0
\move(183 136)
\lvec(191 136)
\lvec(191 137)
\lvec(183 137)
\ifill f:0
\move(195 136)
\lvec(197 136)
\lvec(197 137)
\lvec(195 137)
\ifill f:0
\move(198 136)
\lvec(215 136)
\lvec(215 137)
\lvec(198 137)
\ifill f:0
\move(216 136)
\lvec(223 136)
\lvec(223 137)
\lvec(216 137)
\ifill f:0
\move(224 136)
\lvec(226 136)
\lvec(226 137)
\lvec(224 137)
\ifill f:0
\move(227 136)
\lvec(235 136)
\lvec(235 137)
\lvec(227 137)
\ifill f:0
\move(236 136)
\lvec(242 136)
\lvec(242 137)
\lvec(236 137)
\ifill f:0
\move(243 136)
\lvec(250 136)
\lvec(250 137)
\lvec(243 137)
\ifill f:0
\move(251 136)
\lvec(254 136)
\lvec(254 137)
\lvec(251 137)
\ifill f:0
\move(255 136)
\lvec(257 136)
\lvec(257 137)
\lvec(255 137)
\ifill f:0
\move(258 136)
\lvec(262 136)
\lvec(262 137)
\lvec(258 137)
\ifill f:0
\move(263 136)
\lvec(267 136)
\lvec(267 137)
\lvec(263 137)
\ifill f:0
\move(268 136)
\lvec(274 136)
\lvec(274 137)
\lvec(268 137)
\ifill f:0
\move(275 136)
\lvec(284 136)
\lvec(284 137)
\lvec(275 137)
\ifill f:0
\move(285 136)
\lvec(287 136)
\lvec(287 137)
\lvec(285 137)
\ifill f:0
\move(289 136)
\lvec(290 136)
\lvec(290 137)
\lvec(289 137)
\ifill f:0
\move(291 136)
\lvec(292 136)
\lvec(292 137)
\lvec(291 137)
\ifill f:0
\move(300 136)
\lvec(309 136)
\lvec(309 137)
\lvec(300 137)
\ifill f:0
\move(311 136)
\lvec(316 136)
\lvec(316 137)
\lvec(311 137)
\ifill f:0
\move(317 136)
\lvec(322 136)
\lvec(322 137)
\lvec(317 137)
\ifill f:0
\move(323 136)
\lvec(325 136)
\lvec(325 137)
\lvec(323 137)
\ifill f:0
\move(326 136)
\lvec(356 136)
\lvec(356 137)
\lvec(326 137)
\ifill f:0
\move(357 136)
\lvec(362 136)
\lvec(362 137)
\lvec(357 137)
\ifill f:0
\move(363 136)
\lvec(366 136)
\lvec(366 137)
\lvec(363 137)
\ifill f:0
\move(367 136)
\lvec(379 136)
\lvec(379 137)
\lvec(367 137)
\ifill f:0
\move(380 136)
\lvec(383 136)
\lvec(383 137)
\lvec(380 137)
\ifill f:0
\move(384 136)
\lvec(394 136)
\lvec(394 137)
\lvec(384 137)
\ifill f:0
\move(395 136)
\lvec(401 136)
\lvec(401 137)
\lvec(395 137)
\ifill f:0
\move(402 136)
\lvec(423 136)
\lvec(423 137)
\lvec(402 137)
\ifill f:0
\move(424 136)
\lvec(442 136)
\lvec(442 137)
\lvec(424 137)
\ifill f:0
\move(443 136)
\lvec(445 136)
\lvec(445 137)
\lvec(443 137)
\ifill f:0
\move(446 136)
\lvec(451 136)
\lvec(451 137)
\lvec(446 137)
\ifill f:0
\move(15 137)
\lvec(17 137)
\lvec(17 138)
\lvec(15 138)
\ifill f:0
\move(20 137)
\lvec(21 137)
\lvec(21 138)
\lvec(20 138)
\ifill f:0
\move(24 137)
\lvec(26 137)
\lvec(26 138)
\lvec(24 138)
\ifill f:0
\move(36 137)
\lvec(37 137)
\lvec(37 138)
\lvec(36 138)
\ifill f:0
\move(38 137)
\lvec(39 137)
\lvec(39 138)
\lvec(38 138)
\ifill f:0
\move(40 137)
\lvec(41 137)
\lvec(41 138)
\lvec(40 138)
\ifill f:0
\move(42 137)
\lvec(45 137)
\lvec(45 138)
\lvec(42 138)
\ifill f:0
\move(47 137)
\lvec(50 137)
\lvec(50 138)
\lvec(47 138)
\ifill f:0
\move(52 137)
\lvec(53 137)
\lvec(53 138)
\lvec(52 138)
\ifill f:0
\move(61 137)
\lvec(63 137)
\lvec(63 138)
\lvec(61 138)
\ifill f:0
\move(64 137)
\lvec(65 137)
\lvec(65 138)
\lvec(64 138)
\ifill f:0
\move(66 137)
\lvec(67 137)
\lvec(67 138)
\lvec(66 138)
\ifill f:0
\move(68 137)
\lvec(70 137)
\lvec(70 138)
\lvec(68 138)
\ifill f:0
\move(71 137)
\lvec(78 137)
\lvec(78 138)
\lvec(71 138)
\ifill f:0
\move(79 137)
\lvec(82 137)
\lvec(82 138)
\lvec(79 138)
\ifill f:0
\move(83 137)
\lvec(84 137)
\lvec(84 138)
\lvec(83 138)
\ifill f:0
\move(87 137)
\lvec(90 137)
\lvec(90 138)
\lvec(87 138)
\ifill f:0
\move(91 137)
\lvec(93 137)
\lvec(93 138)
\lvec(91 138)
\ifill f:0
\move(97 137)
\lvec(98 137)
\lvec(98 138)
\lvec(97 138)
\ifill f:0
\move(100 137)
\lvec(101 137)
\lvec(101 138)
\lvec(100 138)
\ifill f:0
\move(102 137)
\lvec(103 137)
\lvec(103 138)
\lvec(102 138)
\ifill f:0
\move(104 137)
\lvec(117 137)
\lvec(117 138)
\lvec(104 138)
\ifill f:0
\move(118 137)
\lvec(122 137)
\lvec(122 138)
\lvec(118 138)
\ifill f:0
\move(123 137)
\lvec(124 137)
\lvec(124 138)
\lvec(123 138)
\ifill f:0
\move(127 137)
\lvec(132 137)
\lvec(132 138)
\lvec(127 138)
\ifill f:0
\move(133 137)
\lvec(137 137)
\lvec(137 138)
\lvec(133 138)
\ifill f:0
\move(139 137)
\lvec(142 137)
\lvec(142 138)
\lvec(139 138)
\ifill f:0
\move(144 137)
\lvec(145 137)
\lvec(145 138)
\lvec(144 138)
\ifill f:0
\move(146 137)
\lvec(149 137)
\lvec(149 138)
\lvec(146 138)
\ifill f:0
\move(150 137)
\lvec(163 137)
\lvec(163 138)
\lvec(150 138)
\ifill f:0
\move(164 137)
\lvec(170 137)
\lvec(170 138)
\lvec(164 138)
\ifill f:0
\move(171 137)
\lvec(172 137)
\lvec(172 138)
\lvec(171 138)
\ifill f:0
\move(173 137)
\lvec(176 137)
\lvec(176 138)
\lvec(173 138)
\ifill f:0
\move(178 137)
\lvec(186 137)
\lvec(186 138)
\lvec(178 138)
\ifill f:0
\move(188 137)
\lvec(189 137)
\lvec(189 138)
\lvec(188 138)
\ifill f:0
\move(193 137)
\lvec(194 137)
\lvec(194 138)
\lvec(193 138)
\ifill f:0
\move(195 137)
\lvec(197 137)
\lvec(197 138)
\lvec(195 138)
\ifill f:0
\move(198 137)
\lvec(201 137)
\lvec(201 138)
\lvec(198 138)
\ifill f:0
\move(202 137)
\lvec(206 137)
\lvec(206 138)
\lvec(202 138)
\ifill f:0
\move(207 137)
\lvec(214 137)
\lvec(214 138)
\lvec(207 138)
\ifill f:0
\move(215 137)
\lvec(217 137)
\lvec(217 138)
\lvec(215 138)
\ifill f:0
\move(218 137)
\lvec(220 137)
\lvec(220 138)
\lvec(218 138)
\ifill f:0
\move(221 137)
\lvec(226 137)
\lvec(226 138)
\lvec(221 138)
\ifill f:0
\move(227 137)
\lvec(234 137)
\lvec(234 138)
\lvec(227 138)
\ifill f:0
\move(235 137)
\lvec(243 137)
\lvec(243 138)
\lvec(235 138)
\ifill f:0
\move(244 137)
\lvec(248 137)
\lvec(248 138)
\lvec(244 138)
\ifill f:0
\move(249 137)
\lvec(251 137)
\lvec(251 138)
\lvec(249 138)
\ifill f:0
\move(252 137)
\lvec(257 137)
\lvec(257 138)
\lvec(252 138)
\ifill f:0
\move(258 137)
\lvec(277 137)
\lvec(277 138)
\lvec(258 138)
\ifill f:0
\move(278 137)
\lvec(287 137)
\lvec(287 138)
\lvec(278 138)
\ifill f:0
\move(289 137)
\lvec(290 137)
\lvec(290 138)
\lvec(289 138)
\ifill f:0
\move(291 137)
\lvec(302 137)
\lvec(302 138)
\lvec(291 138)
\ifill f:0
\move(303 137)
\lvec(305 137)
\lvec(305 138)
\lvec(303 138)
\ifill f:0
\move(306 137)
\lvec(314 137)
\lvec(314 138)
\lvec(306 138)
\ifill f:0
\move(316 137)
\lvec(321 137)
\lvec(321 138)
\lvec(316 138)
\ifill f:0
\move(323 137)
\lvec(325 137)
\lvec(325 138)
\lvec(323 138)
\ifill f:0
\move(326 137)
\lvec(362 137)
\lvec(362 138)
\lvec(326 138)
\ifill f:0
\move(363 137)
\lvec(364 137)
\lvec(364 138)
\lvec(363 138)
\ifill f:0
\move(365 137)
\lvec(371 137)
\lvec(371 138)
\lvec(365 138)
\ifill f:0
\move(372 137)
\lvec(380 137)
\lvec(380 138)
\lvec(372 138)
\ifill f:0
\move(381 137)
\lvec(382 137)
\lvec(382 138)
\lvec(381 138)
\ifill f:0
\move(383 137)
\lvec(384 137)
\lvec(384 138)
\lvec(383 138)
\ifill f:0
\move(385 137)
\lvec(386 137)
\lvec(386 138)
\lvec(385 138)
\ifill f:0
\move(387 137)
\lvec(388 137)
\lvec(388 138)
\lvec(387 138)
\ifill f:0
\move(389 137)
\lvec(392 137)
\lvec(392 138)
\lvec(389 138)
\ifill f:0
\move(393 137)
\lvec(401 137)
\lvec(401 138)
\lvec(393 138)
\ifill f:0
\move(402 137)
\lvec(416 137)
\lvec(416 138)
\lvec(402 138)
\ifill f:0
\move(417 137)
\lvec(434 137)
\lvec(434 138)
\lvec(417 138)
\ifill f:0
\move(435 137)
\lvec(442 137)
\lvec(442 138)
\lvec(435 138)
\ifill f:0
\move(443 137)
\lvec(451 137)
\lvec(451 138)
\lvec(443 138)
\ifill f:0
\move(16 138)
\lvec(17 138)
\lvec(17 139)
\lvec(16 139)
\ifill f:0
\move(19 138)
\lvec(20 138)
\lvec(20 139)
\lvec(19 139)
\ifill f:0
\move(22 138)
\lvec(23 138)
\lvec(23 139)
\lvec(22 139)
\ifill f:0
\move(24 138)
\lvec(26 138)
\lvec(26 139)
\lvec(24 139)
\ifill f:0
\move(36 138)
\lvec(37 138)
\lvec(37 139)
\lvec(36 139)
\ifill f:0
\move(38 138)
\lvec(39 138)
\lvec(39 139)
\lvec(38 139)
\ifill f:0
\move(40 138)
\lvec(45 138)
\lvec(45 139)
\lvec(40 139)
\ifill f:0
\move(47 138)
\lvec(50 138)
\lvec(50 139)
\lvec(47 139)
\ifill f:0
\move(54 138)
\lvec(55 138)
\lvec(55 139)
\lvec(54 139)
\ifill f:0
\move(57 138)
\lvec(58 138)
\lvec(58 139)
\lvec(57 139)
\ifill f:0
\move(59 138)
\lvec(63 138)
\lvec(63 139)
\lvec(59 139)
\ifill f:0
\move(64 138)
\lvec(65 138)
\lvec(65 139)
\lvec(64 139)
\ifill f:0
\move(66 138)
\lvec(71 138)
\lvec(71 139)
\lvec(66 139)
\ifill f:0
\move(72 138)
\lvec(73 138)
\lvec(73 139)
\lvec(72 139)
\ifill f:0
\move(77 138)
\lvec(78 138)
\lvec(78 139)
\lvec(77 139)
\ifill f:0
\move(79 138)
\lvec(82 138)
\lvec(82 139)
\lvec(79 139)
\ifill f:0
\move(83 138)
\lvec(84 138)
\lvec(84 139)
\lvec(83 139)
\ifill f:0
\move(89 138)
\lvec(91 138)
\lvec(91 139)
\lvec(89 139)
\ifill f:0
\move(97 138)
\lvec(98 138)
\lvec(98 139)
\lvec(97 139)
\ifill f:0
\move(100 138)
\lvec(101 138)
\lvec(101 139)
\lvec(100 139)
\ifill f:0
\move(102 138)
\lvec(111 138)
\lvec(111 139)
\lvec(102 139)
\ifill f:0
\move(112 138)
\lvec(118 138)
\lvec(118 139)
\lvec(112 139)
\ifill f:0
\move(119 138)
\lvec(122 138)
\lvec(122 139)
\lvec(119 139)
\ifill f:0
\move(123 138)
\lvec(124 138)
\lvec(124 139)
\lvec(123 139)
\ifill f:0
\move(126 138)
\lvec(130 138)
\lvec(130 139)
\lvec(126 139)
\ifill f:0
\move(137 138)
\lvec(138 138)
\lvec(138 139)
\lvec(137 139)
\ifill f:0
\move(139 138)
\lvec(142 138)
\lvec(142 139)
\lvec(139 139)
\ifill f:0
\move(144 138)
\lvec(145 138)
\lvec(145 139)
\lvec(144 139)
\ifill f:0
\move(147 138)
\lvec(155 138)
\lvec(155 139)
\lvec(147 139)
\ifill f:0
\move(156 138)
\lvec(159 138)
\lvec(159 139)
\lvec(156 139)
\ifill f:0
\move(160 138)
\lvec(163 138)
\lvec(163 139)
\lvec(160 139)
\ifill f:0
\move(164 138)
\lvec(166 138)
\lvec(166 139)
\lvec(164 139)
\ifill f:0
\move(167 138)
\lvec(170 138)
\lvec(170 139)
\lvec(167 139)
\ifill f:0
\move(171 138)
\lvec(172 138)
\lvec(172 139)
\lvec(171 139)
\ifill f:0
\move(173 138)
\lvec(175 138)
\lvec(175 139)
\lvec(173 139)
\ifill f:0
\move(177 138)
\lvec(182 138)
\lvec(182 139)
\lvec(177 139)
\ifill f:0
\move(183 138)
\lvec(191 138)
\lvec(191 139)
\lvec(183 139)
\ifill f:0
\move(192 138)
\lvec(197 138)
\lvec(197 139)
\lvec(192 139)
\ifill f:0
\move(198 138)
\lvec(202 138)
\lvec(202 139)
\lvec(198 139)
\ifill f:0
\move(203 138)
\lvec(208 138)
\lvec(208 139)
\lvec(203 139)
\ifill f:0
\move(209 138)
\lvec(226 138)
\lvec(226 139)
\lvec(209 139)
\ifill f:0
\move(227 138)
\lvec(228 138)
\lvec(228 139)
\lvec(227 139)
\ifill f:0
\move(229 138)
\lvec(230 138)
\lvec(230 139)
\lvec(229 139)
\ifill f:0
\move(231 138)
\lvec(235 138)
\lvec(235 139)
\lvec(231 139)
\ifill f:0
\move(236 138)
\lvec(244 138)
\lvec(244 139)
\lvec(236 139)
\ifill f:0
\move(245 138)
\lvec(249 138)
\lvec(249 139)
\lvec(245 139)
\ifill f:0
\move(250 138)
\lvec(251 138)
\lvec(251 139)
\lvec(250 139)
\ifill f:0
\move(252 138)
\lvec(257 138)
\lvec(257 139)
\lvec(252 139)
\ifill f:0
\move(258 138)
\lvec(279 138)
\lvec(279 139)
\lvec(258 139)
\ifill f:0
\move(280 138)
\lvec(288 138)
\lvec(288 139)
\lvec(280 139)
\ifill f:0
\move(289 138)
\lvec(290 138)
\lvec(290 139)
\lvec(289 139)
\ifill f:0
\move(291 138)
\lvec(311 138)
\lvec(311 139)
\lvec(291 139)
\ifill f:0
\move(312 138)
\lvec(321 138)
\lvec(321 139)
\lvec(312 139)
\ifill f:0
\move(323 138)
\lvec(325 138)
\lvec(325 139)
\lvec(323 139)
\ifill f:0
\move(326 138)
\lvec(327 138)
\lvec(327 139)
\lvec(326 139)
\ifill f:0
\move(328 138)
\lvec(352 138)
\lvec(352 139)
\lvec(328 139)
\ifill f:0
\move(353 138)
\lvec(362 138)
\lvec(362 139)
\lvec(353 139)
\ifill f:0
\move(363 138)
\lvec(364 138)
\lvec(364 139)
\lvec(363 139)
\ifill f:0
\move(365 138)
\lvec(369 138)
\lvec(369 139)
\lvec(365 139)
\ifill f:0
\move(370 138)
\lvec(383 138)
\lvec(383 139)
\lvec(370 139)
\ifill f:0
\move(384 138)
\lvec(401 138)
\lvec(401 139)
\lvec(384 139)
\ifill f:0
\move(402 138)
\lvec(415 138)
\lvec(415 139)
\lvec(402 139)
\ifill f:0
\move(416 138)
\lvec(442 138)
\lvec(442 139)
\lvec(416 139)
\ifill f:0
\move(443 138)
\lvec(444 138)
\lvec(444 139)
\lvec(443 139)
\ifill f:0
\move(445 138)
\lvec(448 138)
\lvec(448 139)
\lvec(445 139)
\ifill f:0
\move(449 138)
\lvec(451 138)
\lvec(451 139)
\lvec(449 139)
\ifill f:0
\move(16 139)
\lvec(17 139)
\lvec(17 140)
\lvec(16 140)
\ifill f:0
\move(20 139)
\lvec(21 139)
\lvec(21 140)
\lvec(20 140)
\ifill f:0
\move(25 139)
\lvec(26 139)
\lvec(26 140)
\lvec(25 140)
\ifill f:0
\move(36 139)
\lvec(37 139)
\lvec(37 140)
\lvec(36 140)
\ifill f:0
\move(38 139)
\lvec(39 139)
\lvec(39 140)
\lvec(38 140)
\ifill f:0
\move(40 139)
\lvec(41 139)
\lvec(41 140)
\lvec(40 140)
\ifill f:0
\move(42 139)
\lvec(45 139)
\lvec(45 140)
\lvec(42 140)
\ifill f:0
\move(49 139)
\lvec(50 139)
\lvec(50 140)
\lvec(49 140)
\ifill f:0
\move(55 139)
\lvec(56 139)
\lvec(56 140)
\lvec(55 140)
\ifill f:0
\move(59 139)
\lvec(62 139)
\lvec(62 140)
\lvec(59 140)
\ifill f:0
\move(64 139)
\lvec(65 139)
\lvec(65 140)
\lvec(64 140)
\ifill f:0
\move(66 139)
\lvec(71 139)
\lvec(71 140)
\lvec(66 140)
\ifill f:0
\move(72 139)
\lvec(74 139)
\lvec(74 140)
\lvec(72 140)
\ifill f:0
\move(76 139)
\lvec(77 139)
\lvec(77 140)
\lvec(76 140)
\ifill f:0
\move(78 139)
\lvec(82 139)
\lvec(82 140)
\lvec(78 140)
\ifill f:0
\move(83 139)
\lvec(85 139)
\lvec(85 140)
\lvec(83 140)
\ifill f:0
\move(86 139)
\lvec(87 139)
\lvec(87 140)
\lvec(86 140)
\ifill f:0
\move(88 139)
\lvec(90 139)
\lvec(90 140)
\lvec(88 140)
\ifill f:0
\move(91 139)
\lvec(93 139)
\lvec(93 140)
\lvec(91 140)
\ifill f:0
\move(95 139)
\lvec(98 139)
\lvec(98 140)
\lvec(95 140)
\ifill f:0
\move(100 139)
\lvec(101 139)
\lvec(101 140)
\lvec(100 140)
\ifill f:0
\move(102 139)
\lvec(105 139)
\lvec(105 140)
\lvec(102 140)
\ifill f:0
\move(107 139)
\lvec(109 139)
\lvec(109 140)
\lvec(107 140)
\ifill f:0
\move(110 139)
\lvec(116 139)
\lvec(116 140)
\lvec(110 140)
\ifill f:0
\move(117 139)
\lvec(122 139)
\lvec(122 140)
\lvec(117 140)
\ifill f:0
\move(123 139)
\lvec(124 139)
\lvec(124 140)
\lvec(123 140)
\ifill f:0
\move(125 139)
\lvec(127 139)
\lvec(127 140)
\lvec(125 140)
\ifill f:0
\move(130 139)
\lvec(132 139)
\lvec(132 140)
\lvec(130 140)
\ifill f:0
\move(133 139)
\lvec(141 139)
\lvec(141 140)
\lvec(133 140)
\ifill f:0
\move(142 139)
\lvec(143 139)
\lvec(143 140)
\lvec(142 140)
\ifill f:0
\move(144 139)
\lvec(145 139)
\lvec(145 140)
\lvec(144 140)
\ifill f:0
\move(146 139)
\lvec(147 139)
\lvec(147 140)
\lvec(146 140)
\ifill f:0
\move(148 139)
\lvec(164 139)
\lvec(164 140)
\lvec(148 140)
\ifill f:0
\move(165 139)
\lvec(170 139)
\lvec(170 140)
\lvec(165 140)
\ifill f:0
\move(171 139)
\lvec(174 139)
\lvec(174 140)
\lvec(171 140)
\ifill f:0
\move(176 139)
\lvec(179 139)
\lvec(179 140)
\lvec(176 140)
\ifill f:0
\move(180 139)
\lvec(185 139)
\lvec(185 140)
\lvec(180 140)
\ifill f:0
\move(186 139)
\lvec(197 139)
\lvec(197 140)
\lvec(186 140)
\ifill f:0
\move(198 139)
\lvec(204 139)
\lvec(204 140)
\lvec(198 140)
\ifill f:0
\move(205 139)
\lvec(211 139)
\lvec(211 140)
\lvec(205 140)
\ifill f:0
\move(212 139)
\lvec(215 139)
\lvec(215 140)
\lvec(212 140)
\ifill f:0
\move(216 139)
\lvec(219 139)
\lvec(219 140)
\lvec(216 140)
\ifill f:0
\move(220 139)
\lvec(222 139)
\lvec(222 140)
\lvec(220 140)
\ifill f:0
\move(223 139)
\lvec(226 139)
\lvec(226 140)
\lvec(223 140)
\ifill f:0
\move(227 139)
\lvec(228 139)
\lvec(228 140)
\lvec(227 140)
\ifill f:0
\move(229 139)
\lvec(233 139)
\lvec(233 140)
\lvec(229 140)
\ifill f:0
\move(234 139)
\lvec(235 139)
\lvec(235 140)
\lvec(234 140)
\ifill f:0
\move(236 139)
\lvec(247 139)
\lvec(247 140)
\lvec(236 140)
\ifill f:0
\move(248 139)
\lvec(252 139)
\lvec(252 140)
\lvec(248 140)
\ifill f:0
\move(253 139)
\lvec(257 139)
\lvec(257 140)
\lvec(253 140)
\ifill f:0
\move(258 139)
\lvec(267 139)
\lvec(267 140)
\lvec(258 140)
\ifill f:0
\move(268 139)
\lvec(271 139)
\lvec(271 140)
\lvec(268 140)
\ifill f:0
\move(272 139)
\lvec(280 139)
\lvec(280 140)
\lvec(272 140)
\ifill f:0
\move(282 139)
\lvec(288 139)
\lvec(288 140)
\lvec(282 140)
\ifill f:0
\move(289 139)
\lvec(290 139)
\lvec(290 140)
\lvec(289 140)
\ifill f:0
\move(291 139)
\lvec(298 139)
\lvec(298 140)
\lvec(291 140)
\ifill f:0
\move(299 139)
\lvec(320 139)
\lvec(320 140)
\lvec(299 140)
\ifill f:0
\move(323 139)
\lvec(325 139)
\lvec(325 140)
\lvec(323 140)
\ifill f:0
\move(326 139)
\lvec(339 139)
\lvec(339 140)
\lvec(326 140)
\ifill f:0
\move(340 139)
\lvec(362 139)
\lvec(362 140)
\lvec(340 140)
\ifill f:0
\move(363 139)
\lvec(364 139)
\lvec(364 140)
\lvec(363 140)
\ifill f:0
\move(365 139)
\lvec(367 139)
\lvec(367 140)
\lvec(365 140)
\ifill f:0
\move(368 139)
\lvec(370 139)
\lvec(370 140)
\lvec(368 140)
\ifill f:0
\move(371 139)
\lvec(375 139)
\lvec(375 140)
\lvec(371 140)
\ifill f:0
\move(376 139)
\lvec(382 139)
\lvec(382 140)
\lvec(376 140)
\ifill f:0
\move(383 139)
\lvec(391 139)
\lvec(391 140)
\lvec(383 140)
\ifill f:0
\move(392 139)
\lvec(393 139)
\lvec(393 140)
\lvec(392 140)
\ifill f:0
\move(394 139)
\lvec(395 139)
\lvec(395 140)
\lvec(394 140)
\ifill f:0
\move(396 139)
\lvec(399 139)
\lvec(399 140)
\lvec(396 140)
\ifill f:0
\move(400 139)
\lvec(401 139)
\lvec(401 140)
\lvec(400 140)
\ifill f:0
\move(402 139)
\lvec(403 139)
\lvec(403 140)
\lvec(402 140)
\ifill f:0
\move(404 139)
\lvec(419 139)
\lvec(419 140)
\lvec(404 140)
\ifill f:0
\move(420 139)
\lvec(424 139)
\lvec(424 140)
\lvec(420 140)
\ifill f:0
\move(425 139)
\lvec(432 139)
\lvec(432 140)
\lvec(425 140)
\ifill f:0
\move(433 139)
\lvec(435 139)
\lvec(435 140)
\lvec(433 140)
\ifill f:0
\move(436 139)
\lvec(442 139)
\lvec(442 140)
\lvec(436 140)
\ifill f:0
\move(443 139)
\lvec(444 139)
\lvec(444 140)
\lvec(443 140)
\ifill f:0
\move(445 139)
\lvec(451 139)
\lvec(451 140)
\lvec(445 140)
\ifill f:0
\move(16 140)
\lvec(17 140)
\lvec(17 141)
\lvec(16 141)
\ifill f:0
\move(25 140)
\lvec(26 140)
\lvec(26 141)
\lvec(25 141)
\ifill f:0
\move(36 140)
\lvec(37 140)
\lvec(37 141)
\lvec(36 141)
\ifill f:0
\move(38 140)
\lvec(39 140)
\lvec(39 141)
\lvec(38 141)
\ifill f:0
\move(40 140)
\lvec(43 140)
\lvec(43 141)
\lvec(40 141)
\ifill f:0
\move(44 140)
\lvec(45 140)
\lvec(45 141)
\lvec(44 141)
\ifill f:0
\move(49 140)
\lvec(50 140)
\lvec(50 141)
\lvec(49 141)
\ifill f:0
\move(57 140)
\lvec(58 140)
\lvec(58 141)
\lvec(57 141)
\ifill f:0
\move(64 140)
\lvec(65 140)
\lvec(65 141)
\lvec(64 141)
\ifill f:0
\move(67 140)
\lvec(73 140)
\lvec(73 141)
\lvec(67 141)
\ifill f:0
\move(74 140)
\lvec(82 140)
\lvec(82 141)
\lvec(74 141)
\ifill f:0
\move(83 140)
\lvec(85 140)
\lvec(85 141)
\lvec(83 141)
\ifill f:0
\move(86 140)
\lvec(87 140)
\lvec(87 141)
\lvec(86 141)
\ifill f:0
\move(88 140)
\lvec(89 140)
\lvec(89 141)
\lvec(88 141)
\ifill f:0
\move(90 140)
\lvec(91 140)
\lvec(91 141)
\lvec(90 141)
\ifill f:0
\move(92 140)
\lvec(93 140)
\lvec(93 141)
\lvec(92 141)
\ifill f:0
\move(96 140)
\lvec(98 140)
\lvec(98 141)
\lvec(96 141)
\ifill f:0
\move(100 140)
\lvec(101 140)
\lvec(101 141)
\lvec(100 141)
\ifill f:0
\move(102 140)
\lvec(107 140)
\lvec(107 141)
\lvec(102 141)
\ifill f:0
\move(108 140)
\lvec(111 140)
\lvec(111 141)
\lvec(108 141)
\ifill f:0
\move(112 140)
\lvec(116 140)
\lvec(116 141)
\lvec(112 141)
\ifill f:0
\move(117 140)
\lvec(122 140)
\lvec(122 141)
\lvec(117 141)
\ifill f:0
\move(124 140)
\lvec(126 140)
\lvec(126 141)
\lvec(124 141)
\ifill f:0
\move(128 140)
\lvec(131 140)
\lvec(131 141)
\lvec(128 141)
\ifill f:0
\move(137 140)
\lvec(138 140)
\lvec(138 141)
\lvec(137 141)
\ifill f:0
\move(141 140)
\lvec(145 140)
\lvec(145 141)
\lvec(141 141)
\ifill f:0
\move(146 140)
\lvec(147 140)
\lvec(147 141)
\lvec(146 141)
\ifill f:0
\move(148 140)
\lvec(151 140)
\lvec(151 141)
\lvec(148 141)
\ifill f:0
\move(152 140)
\lvec(157 140)
\lvec(157 141)
\lvec(152 141)
\ifill f:0
\move(158 140)
\lvec(170 140)
\lvec(170 141)
\lvec(158 141)
\ifill f:0
\move(172 140)
\lvec(174 140)
\lvec(174 141)
\lvec(172 141)
\ifill f:0
\move(175 140)
\lvec(178 140)
\lvec(178 141)
\lvec(175 141)
\ifill f:0
\move(179 140)
\lvec(182 140)
\lvec(182 141)
\lvec(179 141)
\ifill f:0
\move(183 140)
\lvec(186 140)
\lvec(186 141)
\lvec(183 141)
\ifill f:0
\move(189 140)
\lvec(197 140)
\lvec(197 141)
\lvec(189 141)
\ifill f:0
\move(198 140)
\lvec(207 140)
\lvec(207 141)
\lvec(198 141)
\ifill f:0
\move(208 140)
\lvec(218 140)
\lvec(218 141)
\lvec(208 141)
\ifill f:0
\move(219 140)
\lvec(226 140)
\lvec(226 141)
\lvec(219 141)
\ifill f:0
\move(227 140)
\lvec(228 140)
\lvec(228 141)
\lvec(227 141)
\ifill f:0
\move(229 140)
\lvec(231 140)
\lvec(231 141)
\lvec(229 141)
\ifill f:0
\move(232 140)
\lvec(234 140)
\lvec(234 141)
\lvec(232 141)
\ifill f:0
\move(235 140)
\lvec(250 140)
\lvec(250 141)
\lvec(235 141)
\ifill f:0
\move(251 140)
\lvec(257 140)
\lvec(257 141)
\lvec(251 141)
\ifill f:0
\move(258 140)
\lvec(260 140)
\lvec(260 141)
\lvec(258 141)
\ifill f:0
\move(261 140)
\lvec(266 140)
\lvec(266 141)
\lvec(261 141)
\ifill f:0
\move(267 140)
\lvec(269 140)
\lvec(269 141)
\lvec(267 141)
\ifill f:0
\move(270 140)
\lvec(277 140)
\lvec(277 141)
\lvec(270 141)
\ifill f:0
\move(278 140)
\lvec(282 140)
\lvec(282 141)
\lvec(278 141)
\ifill f:0
\move(283 140)
\lvec(288 140)
\lvec(288 141)
\lvec(283 141)
\ifill f:0
\move(289 140)
\lvec(290 140)
\lvec(290 141)
\lvec(289 141)
\ifill f:0
\move(291 140)
\lvec(298 140)
\lvec(298 141)
\lvec(291 141)
\ifill f:0
\move(300 140)
\lvec(318 140)
\lvec(318 141)
\lvec(300 141)
\ifill f:0
\move(321 140)
\lvec(325 140)
\lvec(325 141)
\lvec(321 141)
\ifill f:0
\move(326 140)
\lvec(329 140)
\lvec(329 141)
\lvec(326 141)
\ifill f:0
\move(330 140)
\lvec(341 140)
\lvec(341 141)
\lvec(330 141)
\ifill f:0
\move(342 140)
\lvec(346 140)
\lvec(346 141)
\lvec(342 141)
\ifill f:0
\move(347 140)
\lvec(350 140)
\lvec(350 141)
\lvec(347 141)
\ifill f:0
\move(351 140)
\lvec(362 140)
\lvec(362 141)
\lvec(351 141)
\ifill f:0
\move(363 140)
\lvec(373 140)
\lvec(373 141)
\lvec(363 141)
\ifill f:0
\move(374 140)
\lvec(381 140)
\lvec(381 141)
\lvec(374 141)
\ifill f:0
\move(382 140)
\lvec(395 140)
\lvec(395 141)
\lvec(382 141)
\ifill f:0
\move(396 140)
\lvec(397 140)
\lvec(397 141)
\lvec(396 141)
\ifill f:0
\move(398 140)
\lvec(399 140)
\lvec(399 141)
\lvec(398 141)
\ifill f:0
\move(400 140)
\lvec(401 140)
\lvec(401 141)
\lvec(400 141)
\ifill f:0
\move(402 140)
\lvec(405 140)
\lvec(405 141)
\lvec(402 141)
\ifill f:0
\move(406 140)
\lvec(425 140)
\lvec(425 141)
\lvec(406 141)
\ifill f:0
\move(426 140)
\lvec(430 140)
\lvec(430 141)
\lvec(426 141)
\ifill f:0
\move(431 140)
\lvec(442 140)
\lvec(442 141)
\lvec(431 141)
\ifill f:0
\move(443 140)
\lvec(447 140)
\lvec(447 141)
\lvec(443 141)
\ifill f:0
\move(448 140)
\lvec(450 140)
\lvec(450 141)
\lvec(448 141)
\ifill f:0
\move(16 141)
\lvec(17 141)
\lvec(17 142)
\lvec(16 142)
\ifill f:0
\move(19 141)
\lvec(21 141)
\lvec(21 142)
\lvec(19 142)
\ifill f:0
\move(23 141)
\lvec(24 141)
\lvec(24 142)
\lvec(23 142)
\ifill f:0
\move(25 141)
\lvec(26 141)
\lvec(26 142)
\lvec(25 142)
\ifill f:0
\move(36 141)
\lvec(37 141)
\lvec(37 142)
\lvec(36 142)
\ifill f:0
\move(40 141)
\lvec(41 141)
\lvec(41 142)
\lvec(40 142)
\ifill f:0
\move(42 141)
\lvec(45 141)
\lvec(45 142)
\lvec(42 142)
\ifill f:0
\move(47 141)
\lvec(48 141)
\lvec(48 142)
\lvec(47 142)
\ifill f:0
\move(49 141)
\lvec(50 141)
\lvec(50 142)
\lvec(49 142)
\ifill f:0
\move(51 141)
\lvec(52 141)
\lvec(52 142)
\lvec(51 142)
\ifill f:0
\move(54 141)
\lvec(55 141)
\lvec(55 142)
\lvec(54 142)
\ifill f:0
\move(59 141)
\lvec(60 141)
\lvec(60 142)
\lvec(59 142)
\ifill f:0
\move(62 141)
\lvec(63 141)
\lvec(63 142)
\lvec(62 142)
\ifill f:0
\move(64 141)
\lvec(65 141)
\lvec(65 142)
\lvec(64 142)
\ifill f:0
\move(66 141)
\lvec(67 141)
\lvec(67 142)
\lvec(66 142)
\ifill f:0
\move(68 141)
\lvec(71 141)
\lvec(71 142)
\lvec(68 142)
\ifill f:0
\move(72 141)
\lvec(74 141)
\lvec(74 142)
\lvec(72 142)
\ifill f:0
\move(76 141)
\lvec(82 141)
\lvec(82 142)
\lvec(76 142)
\ifill f:0
\move(84 141)
\lvec(86 141)
\lvec(86 142)
\lvec(84 142)
\ifill f:0
\move(89 141)
\lvec(90 141)
\lvec(90 142)
\lvec(89 142)
\ifill f:0
\move(91 141)
\lvec(92 141)
\lvec(92 142)
\lvec(91 142)
\ifill f:0
\move(97 141)
\lvec(99 141)
\lvec(99 142)
\lvec(97 142)
\ifill f:0
\move(100 141)
\lvec(101 141)
\lvec(101 142)
\lvec(100 142)
\ifill f:0
\move(102 141)
\lvec(109 141)
\lvec(109 142)
\lvec(102 142)
\ifill f:0
\move(110 141)
\lvec(112 141)
\lvec(112 142)
\lvec(110 142)
\ifill f:0
\move(113 141)
\lvec(115 141)
\lvec(115 142)
\lvec(113 142)
\ifill f:0
\move(116 141)
\lvec(122 141)
\lvec(122 142)
\lvec(116 142)
\ifill f:0
\move(124 141)
\lvec(125 141)
\lvec(125 142)
\lvec(124 142)
\ifill f:0
\move(127 141)
\lvec(130 141)
\lvec(130 142)
\lvec(127 142)
\ifill f:0
\move(131 141)
\lvec(136 141)
\lvec(136 142)
\lvec(131 142)
\ifill f:0
\move(137 141)
\lvec(138 141)
\lvec(138 142)
\lvec(137 142)
\ifill f:0
\move(139 141)
\lvec(145 141)
\lvec(145 142)
\lvec(139 142)
\ifill f:0
\move(146 141)
\lvec(159 141)
\lvec(159 142)
\lvec(146 142)
\ifill f:0
\move(160 141)
\lvec(161 141)
\lvec(161 142)
\lvec(160 142)
\ifill f:0
\move(162 141)
\lvec(170 141)
\lvec(170 142)
\lvec(162 142)
\ifill f:0
\move(172 141)
\lvec(173 141)
\lvec(173 142)
\lvec(172 142)
\ifill f:0
\move(174 141)
\lvec(176 141)
\lvec(176 142)
\lvec(174 142)
\ifill f:0
\move(177 141)
\lvec(180 141)
\lvec(180 142)
\lvec(177 142)
\ifill f:0
\move(181 141)
\lvec(184 141)
\lvec(184 142)
\lvec(181 142)
\ifill f:0
\move(185 141)
\lvec(189 141)
\lvec(189 142)
\lvec(185 142)
\ifill f:0
\move(191 141)
\lvec(197 141)
\lvec(197 142)
\lvec(191 142)
\ifill f:0
\move(198 141)
\lvec(211 141)
\lvec(211 142)
\lvec(198 142)
\ifill f:0
\move(212 141)
\lvec(217 141)
\lvec(217 142)
\lvec(212 142)
\ifill f:0
\move(218 141)
\lvec(226 141)
\lvec(226 142)
\lvec(218 142)
\ifill f:0
\move(227 141)
\lvec(229 141)
\lvec(229 142)
\lvec(227 142)
\ifill f:0
\move(230 141)
\lvec(235 141)
\lvec(235 142)
\lvec(230 142)
\ifill f:0
\move(236 141)
\lvec(242 141)
\lvec(242 142)
\lvec(236 142)
\ifill f:0
\move(243 141)
\lvec(248 141)
\lvec(248 142)
\lvec(243 142)
\ifill f:0
\move(249 141)
\lvec(250 141)
\lvec(250 142)
\lvec(249 142)
\ifill f:0
\move(251 141)
\lvec(257 141)
\lvec(257 142)
\lvec(251 142)
\ifill f:0
\move(258 141)
\lvec(271 141)
\lvec(271 142)
\lvec(258 142)
\ifill f:0
\move(272 141)
\lvec(275 141)
\lvec(275 142)
\lvec(272 142)
\ifill f:0
\move(276 141)
\lvec(283 141)
\lvec(283 142)
\lvec(276 142)
\ifill f:0
\move(284 141)
\lvec(290 141)
\lvec(290 142)
\lvec(284 142)
\ifill f:0
\move(291 141)
\lvec(295 141)
\lvec(295 142)
\lvec(291 142)
\ifill f:0
\move(297 141)
\lvec(310 141)
\lvec(310 142)
\lvec(297 142)
\ifill f:0
\move(312 141)
\lvec(313 141)
\lvec(313 142)
\lvec(312 142)
\ifill f:0
\move(316 141)
\lvec(325 141)
\lvec(325 142)
\lvec(316 142)
\ifill f:0
\move(326 141)
\lvec(330 141)
\lvec(330 142)
\lvec(326 142)
\ifill f:0
\move(331 141)
\lvec(338 141)
\lvec(338 142)
\lvec(331 142)
\ifill f:0
\move(339 141)
\lvec(344 141)
\lvec(344 142)
\lvec(339 142)
\ifill f:0
\move(345 141)
\lvec(362 141)
\lvec(362 142)
\lvec(345 142)
\ifill f:0
\move(363 141)
\lvec(365 141)
\lvec(365 142)
\lvec(363 142)
\ifill f:0
\move(366 141)
\lvec(368 141)
\lvec(368 142)
\lvec(366 142)
\ifill f:0
\move(369 141)
\lvec(371 141)
\lvec(371 142)
\lvec(369 142)
\ifill f:0
\move(372 141)
\lvec(377 141)
\lvec(377 142)
\lvec(372 142)
\ifill f:0
\move(378 141)
\lvec(397 141)
\lvec(397 142)
\lvec(378 142)
\ifill f:0
\move(398 141)
\lvec(399 141)
\lvec(399 142)
\lvec(398 142)
\ifill f:0
\move(400 141)
\lvec(401 141)
\lvec(401 142)
\lvec(400 142)
\ifill f:0
\move(402 141)
\lvec(417 141)
\lvec(417 142)
\lvec(402 142)
\ifill f:0
\move(418 141)
\lvec(419 141)
\lvec(419 142)
\lvec(418 142)
\ifill f:0
\move(420 141)
\lvec(442 141)
\lvec(442 142)
\lvec(420 142)
\ifill f:0
\move(443 141)
\lvec(451 141)
\lvec(451 142)
\lvec(443 142)
\ifill f:0
\move(16 142)
\lvec(17 142)
\lvec(17 143)
\lvec(16 143)
\ifill f:0
\move(25 142)
\lvec(26 142)
\lvec(26 143)
\lvec(25 143)
\ifill f:0
\move(36 142)
\lvec(37 142)
\lvec(37 143)
\lvec(36 143)
\ifill f:0
\move(40 142)
\lvec(41 142)
\lvec(41 143)
\lvec(40 143)
\ifill f:0
\move(42 142)
\lvec(43 142)
\lvec(43 143)
\lvec(42 143)
\ifill f:0
\move(44 142)
\lvec(45 142)
\lvec(45 143)
\lvec(44 143)
\ifill f:0
\move(47 142)
\lvec(48 142)
\lvec(48 143)
\lvec(47 143)
\ifill f:0
\move(49 142)
\lvec(50 142)
\lvec(50 143)
\lvec(49 143)
\ifill f:0
\move(51 142)
\lvec(53 142)
\lvec(53 143)
\lvec(51 143)
\ifill f:0
\move(59 142)
\lvec(63 142)
\lvec(63 143)
\lvec(59 143)
\ifill f:0
\move(64 142)
\lvec(65 142)
\lvec(65 143)
\lvec(64 143)
\ifill f:0
\move(66 142)
\lvec(67 142)
\lvec(67 143)
\lvec(66 143)
\ifill f:0
\move(68 142)
\lvec(74 142)
\lvec(74 143)
\lvec(68 143)
\ifill f:0
\move(78 142)
\lvec(82 142)
\lvec(82 143)
\lvec(78 143)
\ifill f:0
\move(88 142)
\lvec(89 142)
\lvec(89 143)
\lvec(88 143)
\ifill f:0
\move(90 142)
\lvec(91 142)
\lvec(91 143)
\lvec(90 143)
\ifill f:0
\move(92 142)
\lvec(93 142)
\lvec(93 143)
\lvec(92 143)
\ifill f:0
\move(98 142)
\lvec(99 142)
\lvec(99 143)
\lvec(98 143)
\ifill f:0
\move(100 142)
\lvec(101 142)
\lvec(101 143)
\lvec(100 143)
\ifill f:0
\move(102 142)
\lvec(106 142)
\lvec(106 143)
\lvec(102 143)
\ifill f:0
\move(107 142)
\lvec(111 142)
\lvec(111 143)
\lvec(107 143)
\ifill f:0
\move(112 142)
\lvec(119 142)
\lvec(119 143)
\lvec(112 143)
\ifill f:0
\move(120 142)
\lvec(122 142)
\lvec(122 143)
\lvec(120 143)
\ifill f:0
\move(124 142)
\lvec(125 142)
\lvec(125 143)
\lvec(124 143)
\ifill f:0
\move(126 142)
\lvec(129 142)
\lvec(129 143)
\lvec(126 143)
\ifill f:0
\move(130 142)
\lvec(131 142)
\lvec(131 143)
\lvec(130 143)
\ifill f:0
\move(134 142)
\lvec(138 142)
\lvec(138 143)
\lvec(134 143)
\ifill f:0
\move(139 142)
\lvec(145 142)
\lvec(145 143)
\lvec(139 143)
\ifill f:0
\move(146 142)
\lvec(149 142)
\lvec(149 143)
\lvec(146 143)
\ifill f:0
\move(150 142)
\lvec(154 142)
\lvec(154 143)
\lvec(150 143)
\ifill f:0
\move(155 142)
\lvec(160 142)
\lvec(160 143)
\lvec(155 143)
\ifill f:0
\move(161 142)
\lvec(163 142)
\lvec(163 143)
\lvec(161 143)
\ifill f:0
\move(164 142)
\lvec(167 142)
\lvec(167 143)
\lvec(164 143)
\ifill f:0
\move(168 142)
\lvec(170 142)
\lvec(170 143)
\lvec(168 143)
\ifill f:0
\move(172 142)
\lvec(173 142)
\lvec(173 143)
\lvec(172 143)
\ifill f:0
\move(174 142)
\lvec(176 142)
\lvec(176 143)
\lvec(174 143)
\ifill f:0
\move(177 142)
\lvec(186 142)
\lvec(186 143)
\lvec(177 143)
\ifill f:0
\move(187 142)
\lvec(191 142)
\lvec(191 143)
\lvec(187 143)
\ifill f:0
\move(192 142)
\lvec(197 142)
\lvec(197 143)
\lvec(192 143)
\ifill f:0
\move(198 142)
\lvec(201 142)
\lvec(201 143)
\lvec(198 143)
\ifill f:0
\move(202 142)
\lvec(215 142)
\lvec(215 143)
\lvec(202 143)
\ifill f:0
\move(216 142)
\lvec(221 142)
\lvec(221 143)
\lvec(216 143)
\ifill f:0
\move(222 142)
\lvec(226 142)
\lvec(226 143)
\lvec(222 143)
\ifill f:0
\move(227 142)
\lvec(229 142)
\lvec(229 143)
\lvec(227 143)
\ifill f:0
\move(230 142)
\lvec(232 142)
\lvec(232 143)
\lvec(230 143)
\ifill f:0
\move(233 142)
\lvec(235 142)
\lvec(235 143)
\lvec(233 143)
\ifill f:0
\move(236 142)
\lvec(238 142)
\lvec(238 143)
\lvec(236 143)
\ifill f:0
\move(239 142)
\lvec(257 142)
\lvec(257 143)
\lvec(239 143)
\ifill f:0
\move(258 142)
\lvec(267 142)
\lvec(267 143)
\lvec(258 143)
\ifill f:0
\move(268 142)
\lvec(270 142)
\lvec(270 143)
\lvec(268 143)
\ifill f:0
\move(271 142)
\lvec(273 142)
\lvec(273 143)
\lvec(271 143)
\ifill f:0
\move(274 142)
\lvec(280 142)
\lvec(280 143)
\lvec(274 143)
\ifill f:0
\move(281 142)
\lvec(284 142)
\lvec(284 143)
\lvec(281 143)
\ifill f:0
\move(285 142)
\lvec(290 142)
\lvec(290 143)
\lvec(285 143)
\ifill f:0
\move(291 142)
\lvec(302 142)
\lvec(302 143)
\lvec(291 143)
\ifill f:0
\move(304 142)
\lvec(325 142)
\lvec(325 143)
\lvec(304 143)
\ifill f:0
\move(326 142)
\lvec(332 142)
\lvec(332 143)
\lvec(326 143)
\ifill f:0
\move(333 142)
\lvec(341 142)
\lvec(341 143)
\lvec(333 143)
\ifill f:0
\move(342 142)
\lvec(347 142)
\lvec(347 143)
\lvec(342 143)
\ifill f:0
\move(348 142)
\lvec(362 142)
\lvec(362 143)
\lvec(348 143)
\ifill f:0
\move(363 142)
\lvec(365 142)
\lvec(365 143)
\lvec(363 143)
\ifill f:0
\move(366 142)
\lvec(372 142)
\lvec(372 143)
\lvec(366 143)
\ifill f:0
\move(373 142)
\lvec(375 142)
\lvec(375 143)
\lvec(373 143)
\ifill f:0
\move(376 142)
\lvec(378 142)
\lvec(378 143)
\lvec(376 143)
\ifill f:0
\move(379 142)
\lvec(381 142)
\lvec(381 143)
\lvec(379 143)
\ifill f:0
\move(382 142)
\lvec(384 142)
\lvec(384 143)
\lvec(382 143)
\ifill f:0
\move(385 142)
\lvec(392 142)
\lvec(392 143)
\lvec(385 143)
\ifill f:0
\move(393 142)
\lvec(399 142)
\lvec(399 143)
\lvec(393 143)
\ifill f:0
\move(400 142)
\lvec(401 142)
\lvec(401 143)
\lvec(400 143)
\ifill f:0
\move(402 142)
\lvec(416 142)
\lvec(416 143)
\lvec(402 143)
\ifill f:0
\move(417 142)
\lvec(429 142)
\lvec(429 143)
\lvec(417 143)
\ifill f:0
\move(430 142)
\lvec(436 142)
\lvec(436 143)
\lvec(430 143)
\ifill f:0
\move(437 142)
\lvec(442 142)
\lvec(442 143)
\lvec(437 143)
\ifill f:0
\move(443 142)
\lvec(446 142)
\lvec(446 143)
\lvec(443 143)
\ifill f:0
\move(447 142)
\lvec(451 142)
\lvec(451 143)
\lvec(447 143)
\ifill f:0
\move(15 143)
\lvec(17 143)
\lvec(17 144)
\lvec(15 144)
\ifill f:0
\move(20 143)
\lvec(21 143)
\lvec(21 144)
\lvec(20 144)
\ifill f:0
\move(24 143)
\lvec(26 143)
\lvec(26 144)
\lvec(24 144)
\ifill f:0
\move(36 143)
\lvec(37 143)
\lvec(37 144)
\lvec(36 144)
\ifill f:0
\move(38 143)
\lvec(39 143)
\lvec(39 144)
\lvec(38 144)
\ifill f:0
\move(40 143)
\lvec(45 143)
\lvec(45 144)
\lvec(40 144)
\ifill f:0
\move(47 143)
\lvec(50 143)
\lvec(50 144)
\lvec(47 144)
\ifill f:0
\move(51 143)
\lvec(52 143)
\lvec(52 144)
\lvec(51 144)
\ifill f:0
\move(60 143)
\lvec(63 143)
\lvec(63 144)
\lvec(60 144)
\ifill f:0
\move(64 143)
\lvec(65 143)
\lvec(65 144)
\lvec(64 144)
\ifill f:0
\move(66 143)
\lvec(71 143)
\lvec(71 144)
\lvec(66 144)
\ifill f:0
\move(72 143)
\lvec(74 143)
\lvec(74 144)
\lvec(72 144)
\ifill f:0
\move(75 143)
\lvec(77 143)
\lvec(77 144)
\lvec(75 144)
\ifill f:0
\move(81 143)
\lvec(82 143)
\lvec(82 144)
\lvec(81 144)
\ifill f:0
\move(86 143)
\lvec(87 143)
\lvec(87 144)
\lvec(86 144)
\ifill f:0
\move(89 143)
\lvec(90 143)
\lvec(90 144)
\lvec(89 144)
\ifill f:0
\move(91 143)
\lvec(92 143)
\lvec(92 144)
\lvec(91 144)
\ifill f:0
\move(95 143)
\lvec(96 143)
\lvec(96 144)
\lvec(95 144)
\ifill f:0
\move(97 143)
\lvec(101 143)
\lvec(101 144)
\lvec(97 144)
\ifill f:0
\move(102 143)
\lvec(103 143)
\lvec(103 144)
\lvec(102 144)
\ifill f:0
\move(104 143)
\lvec(105 143)
\lvec(105 144)
\lvec(104 144)
\ifill f:0
\move(108 143)
\lvec(113 143)
\lvec(113 144)
\lvec(108 144)
\ifill f:0
\move(114 143)
\lvec(119 143)
\lvec(119 144)
\lvec(114 144)
\ifill f:0
\move(120 143)
\lvec(122 143)
\lvec(122 144)
\lvec(120 144)
\ifill f:0
\move(124 143)
\lvec(125 143)
\lvec(125 144)
\lvec(124 144)
\ifill f:0
\move(126 143)
\lvec(127 143)
\lvec(127 144)
\lvec(126 144)
\ifill f:0
\move(128 143)
\lvec(130 143)
\lvec(130 144)
\lvec(128 144)
\ifill f:0
\move(132 143)
\lvec(134 143)
\lvec(134 144)
\lvec(132 144)
\ifill f:0
\move(137 143)
\lvec(145 143)
\lvec(145 144)
\lvec(137 144)
\ifill f:0
\move(146 143)
\lvec(150 143)
\lvec(150 144)
\lvec(146 144)
\ifill f:0
\move(152 143)
\lvec(159 143)
\lvec(159 144)
\lvec(152 144)
\ifill f:0
\move(160 143)
\lvec(162 143)
\lvec(162 144)
\lvec(160 144)
\ifill f:0
\move(164 143)
\lvec(167 143)
\lvec(167 144)
\lvec(164 144)
\ifill f:0
\move(168 143)
\lvec(170 143)
\lvec(170 144)
\lvec(168 144)
\ifill f:0
\move(172 143)
\lvec(173 143)
\lvec(173 144)
\lvec(172 144)
\ifill f:0
\move(174 143)
\lvec(175 143)
\lvec(175 144)
\lvec(174 144)
\ifill f:0
\move(176 143)
\lvec(180 143)
\lvec(180 144)
\lvec(176 144)
\ifill f:0
\move(181 143)
\lvec(186 143)
\lvec(186 144)
\lvec(181 144)
\ifill f:0
\move(188 143)
\lvec(191 143)
\lvec(191 144)
\lvec(188 144)
\ifill f:0
\move(193 143)
\lvec(197 143)
\lvec(197 144)
\lvec(193 144)
\ifill f:0
\move(198 143)
\lvec(201 143)
\lvec(201 144)
\lvec(198 144)
\ifill f:0
\move(202 143)
\lvec(211 143)
\lvec(211 144)
\lvec(202 144)
\ifill f:0
\move(212 143)
\lvec(220 143)
\lvec(220 144)
\lvec(212 144)
\ifill f:0
\move(221 143)
\lvec(226 143)
\lvec(226 144)
\lvec(221 144)
\ifill f:0
\move(227 143)
\lvec(233 143)
\lvec(233 144)
\lvec(227 144)
\ifill f:0
\move(234 143)
\lvec(236 143)
\lvec(236 144)
\lvec(234 144)
\ifill f:0
\move(237 143)
\lvec(239 143)
\lvec(239 144)
\lvec(237 144)
\ifill f:0
\move(240 143)
\lvec(242 143)
\lvec(242 144)
\lvec(240 144)
\ifill f:0
\move(243 143)
\lvec(249 143)
\lvec(249 144)
\lvec(243 144)
\ifill f:0
\move(250 143)
\lvec(251 143)
\lvec(251 144)
\lvec(250 144)
\ifill f:0
\move(252 143)
\lvec(253 143)
\lvec(253 144)
\lvec(252 144)
\ifill f:0
\move(254 143)
\lvec(255 143)
\lvec(255 144)
\lvec(254 144)
\ifill f:0
\move(256 143)
\lvec(257 143)
\lvec(257 144)
\lvec(256 144)
\ifill f:0
\move(258 143)
\lvec(259 143)
\lvec(259 144)
\lvec(258 144)
\ifill f:0
\move(260 143)
\lvec(266 143)
\lvec(266 144)
\lvec(260 144)
\ifill f:0
\move(267 143)
\lvec(271 143)
\lvec(271 144)
\lvec(267 144)
\ifill f:0
\move(272 143)
\lvec(281 143)
\lvec(281 144)
\lvec(272 144)
\ifill f:0
\move(282 143)
\lvec(290 143)
\lvec(290 144)
\lvec(282 144)
\ifill f:0
\move(291 143)
\lvec(293 143)
\lvec(293 144)
\lvec(291 144)
\ifill f:0
\move(295 143)
\lvec(298 143)
\lvec(298 144)
\lvec(295 144)
\ifill f:0
\move(300 143)
\lvec(308 143)
\lvec(308 144)
\lvec(300 144)
\ifill f:0
\move(310 143)
\lvec(325 143)
\lvec(325 144)
\lvec(310 144)
\ifill f:0
\move(326 143)
\lvec(335 143)
\lvec(335 144)
\lvec(326 144)
\ifill f:0
\move(336 143)
\lvec(351 143)
\lvec(351 144)
\lvec(336 144)
\ifill f:0
\move(352 143)
\lvec(362 143)
\lvec(362 144)
\lvec(352 144)
\ifill f:0
\move(363 143)
\lvec(370 143)
\lvec(370 144)
\lvec(363 144)
\ifill f:0
\move(371 143)
\lvec(383 143)
\lvec(383 144)
\lvec(371 144)
\ifill f:0
\move(384 143)
\lvec(394 143)
\lvec(394 144)
\lvec(384 144)
\ifill f:0
\move(395 143)
\lvec(399 143)
\lvec(399 144)
\lvec(395 144)
\ifill f:0
\move(400 143)
\lvec(401 143)
\lvec(401 144)
\lvec(400 144)
\ifill f:0
\move(402 143)
\lvec(423 143)
\lvec(423 144)
\lvec(402 144)
\ifill f:0
\move(424 143)
\lvec(434 143)
\lvec(434 144)
\lvec(424 144)
\ifill f:0
\move(435 143)
\lvec(442 143)
\lvec(442 144)
\lvec(435 144)
\ifill f:0
\move(444 143)
\lvec(451 143)
\lvec(451 144)
\lvec(444 144)
\ifill f:0
\move(15 144)
\lvec(17 144)
\lvec(17 145)
\lvec(15 145)
\ifill f:0
\move(20 144)
\lvec(21 144)
\lvec(21 145)
\lvec(20 145)
\ifill f:0
\move(24 144)
\lvec(26 144)
\lvec(26 145)
\lvec(24 145)
\ifill f:0
\move(36 144)
\lvec(37 144)
\lvec(37 145)
\lvec(36 145)
\ifill f:0
\move(38 144)
\lvec(39 144)
\lvec(39 145)
\lvec(38 145)
\ifill f:0
\move(40 144)
\lvec(41 144)
\lvec(41 145)
\lvec(40 145)
\ifill f:0
\move(43 144)
\lvec(44 144)
\lvec(44 145)
\lvec(43 145)
\ifill f:0
\move(47 144)
\lvec(50 144)
\lvec(50 145)
\lvec(47 145)
\ifill f:0
\move(54 144)
\lvec(55 144)
\lvec(55 145)
\lvec(54 145)
\ifill f:0
\move(59 144)
\lvec(60 144)
\lvec(60 145)
\lvec(59 145)
\ifill f:0
\move(61 144)
\lvec(63 144)
\lvec(63 145)
\lvec(61 145)
\ifill f:0
\move(64 144)
\lvec(65 144)
\lvec(65 145)
\lvec(64 145)
\ifill f:0
\move(66 144)
\lvec(74 144)
\lvec(74 145)
\lvec(66 145)
\ifill f:0
\move(76 144)
\lvec(78 144)
\lvec(78 145)
\lvec(76 145)
\ifill f:0
\move(81 144)
\lvec(82 144)
\lvec(82 145)
\lvec(81 145)
\ifill f:0
\move(88 144)
\lvec(89 144)
\lvec(89 145)
\lvec(88 145)
\ifill f:0
\move(90 144)
\lvec(91 144)
\lvec(91 145)
\lvec(90 145)
\ifill f:0
\move(92 144)
\lvec(93 144)
\lvec(93 145)
\lvec(92 145)
\ifill f:0
\move(96 144)
\lvec(98 144)
\lvec(98 145)
\lvec(96 145)
\ifill f:0
\move(99 144)
\lvec(101 144)
\lvec(101 145)
\lvec(99 145)
\ifill f:0
\move(102 144)
\lvec(111 144)
\lvec(111 145)
\lvec(102 145)
\ifill f:0
\move(112 144)
\lvec(116 144)
\lvec(116 145)
\lvec(112 145)
\ifill f:0
\move(117 144)
\lvec(122 144)
\lvec(122 145)
\lvec(117 145)
\ifill f:0
\move(123 144)
\lvec(124 144)
\lvec(124 145)
\lvec(123 145)
\ifill f:0
\move(125 144)
\lvec(127 144)
\lvec(127 145)
\lvec(125 145)
\ifill f:0
\move(128 144)
\lvec(129 144)
\lvec(129 145)
\lvec(128 145)
\ifill f:0
\move(130 144)
\lvec(132 144)
\lvec(132 145)
\lvec(130 145)
\ifill f:0
\move(134 144)
\lvec(136 144)
\lvec(136 145)
\lvec(134 145)
\ifill f:0
\move(137 144)
\lvec(138 144)
\lvec(138 145)
\lvec(137 145)
\ifill f:0
\move(139 144)
\lvec(145 144)
\lvec(145 145)
\lvec(139 145)
\ifill f:0
\move(146 144)
\lvec(153 144)
\lvec(153 145)
\lvec(146 145)
\ifill f:0
\move(154 144)
\lvec(161 144)
\lvec(161 145)
\lvec(154 145)
\ifill f:0
\move(162 144)
\lvec(164 144)
\lvec(164 145)
\lvec(162 145)
\ifill f:0
\move(165 144)
\lvec(167 144)
\lvec(167 145)
\lvec(165 145)
\ifill f:0
\move(168 144)
\lvec(170 144)
\lvec(170 145)
\lvec(168 145)
\ifill f:0
\move(171 144)
\lvec(174 144)
\lvec(174 145)
\lvec(171 145)
\ifill f:0
\move(175 144)
\lvec(177 144)
\lvec(177 145)
\lvec(175 145)
\ifill f:0
\move(178 144)
\lvec(179 144)
\lvec(179 145)
\lvec(178 145)
\ifill f:0
\move(180 144)
\lvec(185 144)
\lvec(185 145)
\lvec(180 145)
\ifill f:0
\move(186 144)
\lvec(192 144)
\lvec(192 145)
\lvec(186 145)
\ifill f:0
\move(194 144)
\lvec(197 144)
\lvec(197 145)
\lvec(194 145)
\ifill f:0
\move(198 144)
\lvec(199 144)
\lvec(199 145)
\lvec(198 145)
\ifill f:0
\move(200 144)
\lvec(219 144)
\lvec(219 145)
\lvec(200 145)
\ifill f:0
\move(220 144)
\lvec(226 144)
\lvec(226 145)
\lvec(220 145)
\ifill f:0
\move(227 144)
\lvec(230 144)
\lvec(230 145)
\lvec(227 145)
\ifill f:0
\move(231 144)
\lvec(234 144)
\lvec(234 145)
\lvec(231 145)
\ifill f:0
\move(236 144)
\lvec(251 144)
\lvec(251 145)
\lvec(236 145)
\ifill f:0
\move(252 144)
\lvec(253 144)
\lvec(253 145)
\lvec(252 145)
\ifill f:0
\move(254 144)
\lvec(255 144)
\lvec(255 145)
\lvec(254 145)
\ifill f:0
\move(256 144)
\lvec(257 144)
\lvec(257 145)
\lvec(256 145)
\ifill f:0
\move(258 144)
\lvec(259 144)
\lvec(259 145)
\lvec(258 145)
\ifill f:0
\move(260 144)
\lvec(261 144)
\lvec(261 145)
\lvec(260 145)
\ifill f:0
\move(262 144)
\lvec(290 144)
\lvec(290 145)
\lvec(262 145)
\ifill f:0
\move(291 144)
\lvec(293 144)
\lvec(293 145)
\lvec(291 145)
\ifill f:0
\move(294 144)
\lvec(298 144)
\lvec(298 145)
\lvec(294 145)
\ifill f:0
\move(299 144)
\lvec(313 144)
\lvec(313 145)
\lvec(299 145)
\ifill f:0
\move(314 144)
\lvec(325 144)
\lvec(325 145)
\lvec(314 145)
\ifill f:0
\move(326 144)
\lvec(340 144)
\lvec(340 145)
\lvec(326 145)
\ifill f:0
\move(341 144)
\lvec(349 144)
\lvec(349 145)
\lvec(341 145)
\ifill f:0
\move(350 144)
\lvec(356 144)
\lvec(356 145)
\lvec(350 145)
\ifill f:0
\move(357 144)
\lvec(362 144)
\lvec(362 145)
\lvec(357 145)
\ifill f:0
\move(363 144)
\lvec(366 144)
\lvec(366 145)
\lvec(363 145)
\ifill f:0
\move(367 144)
\lvec(370 144)
\lvec(370 145)
\lvec(367 145)
\ifill f:0
\move(371 144)
\lvec(396 144)
\lvec(396 145)
\lvec(371 145)
\ifill f:0
\move(397 144)
\lvec(399 144)
\lvec(399 145)
\lvec(397 145)
\ifill f:0
\move(400 144)
\lvec(401 144)
\lvec(401 145)
\lvec(400 145)
\ifill f:0
\move(402 144)
\lvec(411 144)
\lvec(411 145)
\lvec(402 145)
\ifill f:0
\move(412 144)
\lvec(430 144)
\lvec(430 145)
\lvec(412 145)
\ifill f:0
\move(431 144)
\lvec(442 144)
\lvec(442 145)
\lvec(431 145)
\ifill f:0
\move(444 144)
\lvec(448 144)
\lvec(448 145)
\lvec(444 145)
\ifill f:0
\move(449 144)
\lvec(451 144)
\lvec(451 145)
\lvec(449 145)
\ifill f:0
\move(15 145)
\lvec(17 145)
\lvec(17 146)
\lvec(15 146)
\ifill f:0
\move(20 145)
\lvec(21 145)
\lvec(21 146)
\lvec(20 146)
\ifill f:0
\move(24 145)
\lvec(26 145)
\lvec(26 146)
\lvec(24 146)
\ifill f:0
\move(36 145)
\lvec(37 145)
\lvec(37 146)
\lvec(36 146)
\ifill f:0
\move(38 145)
\lvec(39 145)
\lvec(39 146)
\lvec(38 146)
\ifill f:0
\move(40 145)
\lvec(45 145)
\lvec(45 146)
\lvec(40 146)
\ifill f:0
\move(48 145)
\lvec(50 145)
\lvec(50 146)
\lvec(48 146)
\ifill f:0
\move(59 145)
\lvec(60 145)
\lvec(60 146)
\lvec(59 146)
\ifill f:0
\move(62 145)
\lvec(63 145)
\lvec(63 146)
\lvec(62 146)
\ifill f:0
\move(64 145)
\lvec(65 145)
\lvec(65 146)
\lvec(64 146)
\ifill f:0
\move(66 145)
\lvec(70 145)
\lvec(70 146)
\lvec(66 146)
\ifill f:0
\move(72 145)
\lvec(73 145)
\lvec(73 146)
\lvec(72 146)
\ifill f:0
\move(74 145)
\lvec(75 145)
\lvec(75 146)
\lvec(74 146)
\ifill f:0
\move(77 145)
\lvec(79 145)
\lvec(79 146)
\lvec(77 146)
\ifill f:0
\move(81 145)
\lvec(82 145)
\lvec(82 146)
\lvec(81 146)
\ifill f:0
\move(83 145)
\lvec(85 145)
\lvec(85 146)
\lvec(83 146)
\ifill f:0
\move(88 145)
\lvec(90 145)
\lvec(90 146)
\lvec(88 146)
\ifill f:0
\move(91 145)
\lvec(93 145)
\lvec(93 146)
\lvec(91 146)
\ifill f:0
\move(96 145)
\lvec(98 145)
\lvec(98 146)
\lvec(96 146)
\ifill f:0
\move(99 145)
\lvec(101 145)
\lvec(101 146)
\lvec(99 146)
\ifill f:0
\move(103 145)
\lvec(106 145)
\lvec(106 146)
\lvec(103 146)
\ifill f:0
\move(107 145)
\lvec(122 145)
\lvec(122 146)
\lvec(107 146)
\ifill f:0
\move(123 145)
\lvec(124 145)
\lvec(124 146)
\lvec(123 146)
\ifill f:0
\move(125 145)
\lvec(126 145)
\lvec(126 146)
\lvec(125 146)
\ifill f:0
\move(127 145)
\lvec(131 145)
\lvec(131 146)
\lvec(127 146)
\ifill f:0
\move(132 145)
\lvec(134 145)
\lvec(134 146)
\lvec(132 146)
\ifill f:0
\move(135 145)
\lvec(138 145)
\lvec(138 146)
\lvec(135 146)
\ifill f:0
\move(140 145)
\lvec(145 145)
\lvec(145 146)
\lvec(140 146)
\ifill f:0
\move(146 145)
\lvec(155 145)
\lvec(155 146)
\lvec(146 146)
\ifill f:0
\move(156 145)
\lvec(160 145)
\lvec(160 146)
\lvec(156 146)
\ifill f:0
\move(161 145)
\lvec(170 145)
\lvec(170 146)
\lvec(161 146)
\ifill f:0
\move(171 145)
\lvec(172 145)
\lvec(172 146)
\lvec(171 146)
\ifill f:0
\move(173 145)
\lvec(174 145)
\lvec(174 146)
\lvec(173 146)
\ifill f:0
\move(175 145)
\lvec(176 145)
\lvec(176 146)
\lvec(175 146)
\ifill f:0
\move(177 145)
\lvec(181 145)
\lvec(181 146)
\lvec(177 146)
\ifill f:0
\move(182 145)
\lvec(186 145)
\lvec(186 146)
\lvec(182 146)
\ifill f:0
\move(187 145)
\lvec(189 145)
\lvec(189 146)
\lvec(187 146)
\ifill f:0
\move(190 145)
\lvec(193 145)
\lvec(193 146)
\lvec(190 146)
\ifill f:0
\move(194 145)
\lvec(197 145)
\lvec(197 146)
\lvec(194 146)
\ifill f:0
\move(198 145)
\lvec(208 145)
\lvec(208 146)
\lvec(198 146)
\ifill f:0
\move(209 145)
\lvec(210 145)
\lvec(210 146)
\lvec(209 146)
\ifill f:0
\move(216 145)
\lvec(226 145)
\lvec(226 146)
\lvec(216 146)
\ifill f:0
\move(227 145)
\lvec(231 145)
\lvec(231 146)
\lvec(227 146)
\ifill f:0
\move(232 145)
\lvec(235 145)
\lvec(235 146)
\lvec(232 146)
\ifill f:0
\move(236 145)
\lvec(239 145)
\lvec(239 146)
\lvec(236 146)
\ifill f:0
\move(240 145)
\lvec(242 145)
\lvec(242 146)
\lvec(240 146)
\ifill f:0
\move(243 145)
\lvec(248 145)
\lvec(248 146)
\lvec(243 146)
\ifill f:0
\move(249 145)
\lvec(255 145)
\lvec(255 146)
\lvec(249 146)
\ifill f:0
\move(256 145)
\lvec(257 145)
\lvec(257 146)
\lvec(256 146)
\ifill f:0
\move(258 145)
\lvec(267 145)
\lvec(267 146)
\lvec(258 146)
\ifill f:0
\move(268 145)
\lvec(269 145)
\lvec(269 146)
\lvec(268 146)
\ifill f:0
\move(270 145)
\lvec(274 145)
\lvec(274 146)
\lvec(270 146)
\ifill f:0
\move(275 145)
\lvec(279 145)
\lvec(279 146)
\lvec(275 146)
\ifill f:0
\move(280 145)
\lvec(290 145)
\lvec(290 146)
\lvec(280 146)
\ifill f:0
\move(291 145)
\lvec(292 145)
\lvec(292 146)
\lvec(291 146)
\ifill f:0
\move(293 145)
\lvec(297 145)
\lvec(297 146)
\lvec(293 146)
\ifill f:0
\move(298 145)
\lvec(307 145)
\lvec(307 146)
\lvec(298 146)
\ifill f:0
\move(309 145)
\lvec(315 145)
\lvec(315 146)
\lvec(309 146)
\ifill f:0
\move(317 145)
\lvec(325 145)
\lvec(325 146)
\lvec(317 146)
\ifill f:0
\move(326 145)
\lvec(346 145)
\lvec(346 146)
\lvec(326 146)
\ifill f:0
\move(347 145)
\lvec(355 145)
\lvec(355 146)
\lvec(347 146)
\ifill f:0
\move(356 145)
\lvec(362 145)
\lvec(362 146)
\lvec(356 146)
\ifill f:0
\move(363 145)
\lvec(367 145)
\lvec(367 146)
\lvec(363 146)
\ifill f:0
\move(368 145)
\lvec(383 145)
\lvec(383 146)
\lvec(368 146)
\ifill f:0
\move(384 145)
\lvec(393 145)
\lvec(393 146)
\lvec(384 146)
\ifill f:0
\move(394 145)
\lvec(399 145)
\lvec(399 146)
\lvec(394 146)
\ifill f:0
\move(400 145)
\lvec(401 145)
\lvec(401 146)
\lvec(400 146)
\ifill f:0
\move(402 145)
\lvec(404 145)
\lvec(404 146)
\lvec(402 146)
\ifill f:0
\move(405 145)
\lvec(414 145)
\lvec(414 146)
\lvec(405 146)
\ifill f:0
\move(415 145)
\lvec(423 145)
\lvec(423 146)
\lvec(415 146)
\ifill f:0
\move(424 145)
\lvec(442 145)
\lvec(442 146)
\lvec(424 146)
\ifill f:0
\move(444 145)
\lvec(445 145)
\lvec(445 146)
\lvec(444 146)
\ifill f:0
\move(446 145)
\lvec(451 145)
\lvec(451 146)
\lvec(446 146)
\ifill f:0
\move(16 146)
\lvec(17 146)
\lvec(17 147)
\lvec(16 147)
\ifill f:0
\move(19 146)
\lvec(21 146)
\lvec(21 147)
\lvec(19 147)
\ifill f:0
\move(23 146)
\lvec(26 146)
\lvec(26 147)
\lvec(23 147)
\ifill f:0
\move(36 146)
\lvec(37 146)
\lvec(37 147)
\lvec(36 147)
\ifill f:0
\move(38 146)
\lvec(39 146)
\lvec(39 147)
\lvec(38 147)
\ifill f:0
\move(40 146)
\lvec(41 146)
\lvec(41 147)
\lvec(40 147)
\ifill f:0
\move(43 146)
\lvec(45 146)
\lvec(45 147)
\lvec(43 147)
\ifill f:0
\move(48 146)
\lvec(50 146)
\lvec(50 147)
\lvec(48 147)
\ifill f:0
\move(54 146)
\lvec(55 146)
\lvec(55 147)
\lvec(54 147)
\ifill f:0
\move(60 146)
\lvec(61 146)
\lvec(61 147)
\lvec(60 147)
\ifill f:0
\move(62 146)
\lvec(63 146)
\lvec(63 147)
\lvec(62 147)
\ifill f:0
\move(64 146)
\lvec(65 146)
\lvec(65 147)
\lvec(64 147)
\ifill f:0
\move(66 146)
\lvec(71 146)
\lvec(71 147)
\lvec(66 147)
\ifill f:0
\move(72 146)
\lvec(74 146)
\lvec(74 147)
\lvec(72 147)
\ifill f:0
\move(75 146)
\lvec(80 146)
\lvec(80 147)
\lvec(75 147)
\ifill f:0
\move(81 146)
\lvec(82 146)
\lvec(82 147)
\lvec(81 147)
\ifill f:0
\move(83 146)
\lvec(85 146)
\lvec(85 147)
\lvec(83 147)
\ifill f:0
\move(86 146)
\lvec(87 146)
\lvec(87 147)
\lvec(86 147)
\ifill f:0
\move(89 146)
\lvec(92 146)
\lvec(92 147)
\lvec(89 147)
\ifill f:0
\move(97 146)
\lvec(98 146)
\lvec(98 147)
\lvec(97 147)
\ifill f:0
\move(99 146)
\lvec(101 146)
\lvec(101 147)
\lvec(99 147)
\ifill f:0
\move(102 146)
\lvec(106 146)
\lvec(106 147)
\lvec(102 147)
\ifill f:0
\move(112 146)
\lvec(118 146)
\lvec(118 147)
\lvec(112 147)
\ifill f:0
\move(119 146)
\lvec(122 146)
\lvec(122 147)
\lvec(119 147)
\ifill f:0
\move(123 146)
\lvec(124 146)
\lvec(124 147)
\lvec(123 147)
\ifill f:0
\move(125 146)
\lvec(126 146)
\lvec(126 147)
\lvec(125 147)
\ifill f:0
\move(129 146)
\lvec(130 146)
\lvec(130 147)
\lvec(129 147)
\ifill f:0
\move(131 146)
\lvec(132 146)
\lvec(132 147)
\lvec(131 147)
\ifill f:0
\move(133 146)
\lvec(135 146)
\lvec(135 147)
\lvec(133 147)
\ifill f:0
\move(137 146)
\lvec(138 146)
\lvec(138 147)
\lvec(137 147)
\ifill f:0
\move(139 146)
\lvec(140 146)
\lvec(140 147)
\lvec(139 147)
\ifill f:0
\move(141 146)
\lvec(145 146)
\lvec(145 147)
\lvec(141 147)
\ifill f:0
\move(146 146)
\lvec(159 146)
\lvec(159 147)
\lvec(146 147)
\ifill f:0
\move(160 146)
\lvec(163 146)
\lvec(163 147)
\lvec(160 147)
\ifill f:0
\move(164 146)
\lvec(170 146)
\lvec(170 147)
\lvec(164 147)
\ifill f:0
\move(171 146)
\lvec(172 146)
\lvec(172 147)
\lvec(171 147)
\ifill f:0
\move(173 146)
\lvec(174 146)
\lvec(174 147)
\lvec(173 147)
\ifill f:0
\move(175 146)
\lvec(180 146)
\lvec(180 147)
\lvec(175 147)
\ifill f:0
\move(181 146)
\lvec(186 146)
\lvec(186 147)
\lvec(181 147)
\ifill f:0
\move(188 146)
\lvec(190 146)
\lvec(190 147)
\lvec(188 147)
\ifill f:0
\move(191 146)
\lvec(193 146)
\lvec(193 147)
\lvec(191 147)
\ifill f:0
\move(195 146)
\lvec(197 146)
\lvec(197 147)
\lvec(195 147)
\ifill f:0
\move(198 146)
\lvec(204 146)
\lvec(204 147)
\lvec(198 147)
\ifill f:0
\move(205 146)
\lvec(226 146)
\lvec(226 147)
\lvec(205 147)
\ifill f:0
\move(227 146)
\lvec(232 146)
\lvec(232 147)
\lvec(227 147)
\ifill f:0
\move(233 146)
\lvec(244 146)
\lvec(244 147)
\lvec(233 147)
\ifill f:0
\move(245 146)
\lvec(247 146)
\lvec(247 147)
\lvec(245 147)
\ifill f:0
\move(248 146)
\lvec(250 146)
\lvec(250 147)
\lvec(248 147)
\ifill f:0
\move(251 146)
\lvec(255 146)
\lvec(255 147)
\lvec(251 147)
\ifill f:0
\move(256 146)
\lvec(257 146)
\lvec(257 147)
\lvec(256 147)
\ifill f:0
\move(258 146)
\lvec(266 146)
\lvec(266 147)
\lvec(258 147)
\ifill f:0
\move(267 146)
\lvec(275 146)
\lvec(275 147)
\lvec(267 147)
\ifill f:0
\move(276 146)
\lvec(280 146)
\lvec(280 147)
\lvec(276 147)
\ifill f:0
\move(281 146)
\lvec(283 146)
\lvec(283 147)
\lvec(281 147)
\ifill f:0
\move(284 146)
\lvec(290 146)
\lvec(290 147)
\lvec(284 147)
\ifill f:0
\move(291 146)
\lvec(292 146)
\lvec(292 147)
\lvec(291 147)
\ifill f:0
\move(293 146)
\lvec(296 146)
\lvec(296 147)
\lvec(293 147)
\ifill f:0
\move(297 146)
\lvec(305 146)
\lvec(305 147)
\lvec(297 147)
\ifill f:0
\move(306 146)
\lvec(310 146)
\lvec(310 147)
\lvec(306 147)
\ifill f:0
\move(311 146)
\lvec(318 146)
\lvec(318 147)
\lvec(311 147)
\ifill f:0
\move(319 146)
\lvec(325 146)
\lvec(325 147)
\lvec(319 147)
\ifill f:0
\move(326 146)
\lvec(335 146)
\lvec(335 147)
\lvec(326 147)
\ifill f:0
\move(336 146)
\lvec(338 146)
\lvec(338 147)
\lvec(336 147)
\ifill f:0
\move(339 146)
\lvec(353 146)
\lvec(353 147)
\lvec(339 147)
\ifill f:0
\move(354 146)
\lvec(362 146)
\lvec(362 147)
\lvec(354 147)
\ifill f:0
\move(363 146)
\lvec(377 146)
\lvec(377 147)
\lvec(363 147)
\ifill f:0
\move(378 146)
\lvec(389 146)
\lvec(389 147)
\lvec(378 147)
\ifill f:0
\move(390 146)
\lvec(401 146)
\lvec(401 147)
\lvec(390 147)
\ifill f:0
\move(402 146)
\lvec(422 146)
\lvec(422 147)
\lvec(402 147)
\ifill f:0
\move(423 146)
\lvec(433 146)
\lvec(433 147)
\lvec(423 147)
\ifill f:0
\move(434 146)
\lvec(435 146)
\lvec(435 147)
\lvec(434 147)
\ifill f:0
\move(436 146)
\lvec(437 146)
\lvec(437 147)
\lvec(436 147)
\ifill f:0
\move(438 146)
\lvec(439 146)
\lvec(439 147)
\lvec(438 147)
\ifill f:0
\move(440 146)
\lvec(442 146)
\lvec(442 147)
\lvec(440 147)
\ifill f:0
\move(444 146)
\lvec(445 146)
\lvec(445 147)
\lvec(444 147)
\ifill f:0
\move(446 146)
\lvec(451 146)
\lvec(451 147)
\lvec(446 147)
\ifill f:0
\move(16 147)
\lvec(17 147)
\lvec(17 148)
\lvec(16 148)
\ifill f:0
\move(20 147)
\lvec(21 147)
\lvec(21 148)
\lvec(20 148)
\ifill f:0
\move(24 147)
\lvec(26 147)
\lvec(26 148)
\lvec(24 148)
\ifill f:0
\move(36 147)
\lvec(37 147)
\lvec(37 148)
\lvec(36 148)
\ifill f:0
\move(40 147)
\lvec(41 147)
\lvec(41 148)
\lvec(40 148)
\ifill f:0
\move(42 147)
\lvec(45 147)
\lvec(45 148)
\lvec(42 148)
\ifill f:0
\move(48 147)
\lvec(50 147)
\lvec(50 148)
\lvec(48 148)
\ifill f:0
\move(52 147)
\lvec(53 147)
\lvec(53 148)
\lvec(52 148)
\ifill f:0
\move(54 147)
\lvec(55 147)
\lvec(55 148)
\lvec(54 148)
\ifill f:0
\move(56 147)
\lvec(58 147)
\lvec(58 148)
\lvec(56 148)
\ifill f:0
\move(59 147)
\lvec(62 147)
\lvec(62 148)
\lvec(59 148)
\ifill f:0
\move(63 147)
\lvec(65 147)
\lvec(65 148)
\lvec(63 148)
\ifill f:0
\move(66 147)
\lvec(67 147)
\lvec(67 148)
\lvec(66 148)
\ifill f:0
\move(68 147)
\lvec(73 147)
\lvec(73 148)
\lvec(68 148)
\ifill f:0
\move(76 147)
\lvec(80 147)
\lvec(80 148)
\lvec(76 148)
\ifill f:0
\move(81 147)
\lvec(82 147)
\lvec(82 148)
\lvec(81 148)
\ifill f:0
\move(83 147)
\lvec(85 147)
\lvec(85 148)
\lvec(83 148)
\ifill f:0
\move(86 147)
\lvec(87 147)
\lvec(87 148)
\lvec(86 148)
\ifill f:0
\move(88 147)
\lvec(90 147)
\lvec(90 148)
\lvec(88 148)
\ifill f:0
\move(92 147)
\lvec(93 147)
\lvec(93 148)
\lvec(92 148)
\ifill f:0
\move(95 147)
\lvec(96 147)
\lvec(96 148)
\lvec(95 148)
\ifill f:0
\move(97 147)
\lvec(98 147)
\lvec(98 148)
\lvec(97 148)
\ifill f:0
\move(99 147)
\lvec(101 147)
\lvec(101 148)
\lvec(99 148)
\ifill f:0
\move(102 147)
\lvec(105 147)
\lvec(105 148)
\lvec(102 148)
\ifill f:0
\move(106 147)
\lvec(122 147)
\lvec(122 148)
\lvec(106 148)
\ifill f:0
\move(123 147)
\lvec(124 147)
\lvec(124 148)
\lvec(123 148)
\ifill f:0
\move(126 147)
\lvec(127 147)
\lvec(127 148)
\lvec(126 148)
\ifill f:0
\move(128 147)
\lvec(129 147)
\lvec(129 148)
\lvec(128 148)
\ifill f:0
\move(130 147)
\lvec(131 147)
\lvec(131 148)
\lvec(130 148)
\ifill f:0
\move(132 147)
\lvec(134 147)
\lvec(134 148)
\lvec(132 148)
\ifill f:0
\move(135 147)
\lvec(136 147)
\lvec(136 148)
\lvec(135 148)
\ifill f:0
\move(138 147)
\lvec(141 147)
\lvec(141 148)
\lvec(138 148)
\ifill f:0
\move(142 147)
\lvec(145 147)
\lvec(145 148)
\lvec(142 148)
\ifill f:0
\move(146 147)
\lvec(147 147)
\lvec(147 148)
\lvec(146 148)
\ifill f:0
\move(148 147)
\lvec(149 147)
\lvec(149 148)
\lvec(148 148)
\ifill f:0
\move(151 147)
\lvec(152 147)
\lvec(152 148)
\lvec(151 148)
\ifill f:0
\move(156 147)
\lvec(162 147)
\lvec(162 148)
\lvec(156 148)
\ifill f:0
\move(164 147)
\lvec(170 147)
\lvec(170 148)
\lvec(164 148)
\ifill f:0
\move(171 147)
\lvec(172 147)
\lvec(172 148)
\lvec(171 148)
\ifill f:0
\move(174 147)
\lvec(175 147)
\lvec(175 148)
\lvec(174 148)
\ifill f:0
\move(176 147)
\lvec(177 147)
\lvec(177 148)
\lvec(176 148)
\ifill f:0
\move(178 147)
\lvec(179 147)
\lvec(179 148)
\lvec(178 148)
\ifill f:0
\move(180 147)
\lvec(183 147)
\lvec(183 148)
\lvec(180 148)
\ifill f:0
\move(184 147)
\lvec(188 147)
\lvec(188 148)
\lvec(184 148)
\ifill f:0
\move(189 147)
\lvec(191 147)
\lvec(191 148)
\lvec(189 148)
\ifill f:0
\move(192 147)
\lvec(194 147)
\lvec(194 148)
\lvec(192 148)
\ifill f:0
\move(195 147)
\lvec(197 147)
\lvec(197 148)
\lvec(195 148)
\ifill f:0
\move(199 147)
\lvec(202 147)
\lvec(202 148)
\lvec(199 148)
\ifill f:0
\move(203 147)
\lvec(210 147)
\lvec(210 148)
\lvec(203 148)
\ifill f:0
\move(211 147)
\lvec(226 147)
\lvec(226 148)
\lvec(211 148)
\ifill f:0
\move(228 147)
\lvec(233 147)
\lvec(233 148)
\lvec(228 148)
\ifill f:0
\move(236 147)
\lvec(246 147)
\lvec(246 148)
\lvec(236 148)
\ifill f:0
\move(247 147)
\lvec(249 147)
\lvec(249 148)
\lvec(247 148)
\ifill f:0
\move(251 147)
\lvec(252 147)
\lvec(252 148)
\lvec(251 148)
\ifill f:0
\move(253 147)
\lvec(255 147)
\lvec(255 148)
\lvec(253 148)
\ifill f:0
\move(256 147)
\lvec(257 147)
\lvec(257 148)
\lvec(256 148)
\ifill f:0
\move(258 147)
\lvec(262 147)
\lvec(262 148)
\lvec(258 148)
\ifill f:0
\move(263 147)
\lvec(286 147)
\lvec(286 148)
\lvec(263 148)
\ifill f:0
\move(287 147)
\lvec(290 147)
\lvec(290 148)
\lvec(287 148)
\ifill f:0
\move(291 147)
\lvec(292 147)
\lvec(292 148)
\lvec(291 148)
\ifill f:0
\move(293 147)
\lvec(295 147)
\lvec(295 148)
\lvec(293 148)
\ifill f:0
\move(296 147)
\lvec(298 147)
\lvec(298 148)
\lvec(296 148)
\ifill f:0
\move(300 147)
\lvec(307 147)
\lvec(307 148)
\lvec(300 148)
\ifill f:0
\move(308 147)
\lvec(319 147)
\lvec(319 148)
\lvec(308 148)
\ifill f:0
\move(320 147)
\lvec(325 147)
\lvec(325 148)
\lvec(320 148)
\ifill f:0
\move(326 147)
\lvec(329 147)
\lvec(329 148)
\lvec(326 148)
\ifill f:0
\move(331 147)
\lvec(362 147)
\lvec(362 148)
\lvec(331 148)
\ifill f:0
\move(363 147)
\lvec(374 147)
\lvec(374 148)
\lvec(363 148)
\ifill f:0
\move(376 147)
\lvec(384 147)
\lvec(384 148)
\lvec(376 148)
\ifill f:0
\move(385 147)
\lvec(388 147)
\lvec(388 148)
\lvec(385 148)
\ifill f:0
\move(389 147)
\lvec(392 147)
\lvec(392 148)
\lvec(389 148)
\ifill f:0
\move(393 147)
\lvec(395 147)
\lvec(395 148)
\lvec(393 148)
\ifill f:0
\move(396 147)
\lvec(401 147)
\lvec(401 148)
\lvec(396 148)
\ifill f:0
\move(402 147)
\lvec(416 147)
\lvec(416 148)
\lvec(402 148)
\ifill f:0
\move(417 147)
\lvec(428 147)
\lvec(428 148)
\lvec(417 148)
\ifill f:0
\move(429 147)
\lvec(439 147)
\lvec(439 148)
\lvec(429 148)
\ifill f:0
\move(440 147)
\lvec(442 147)
\lvec(442 148)
\lvec(440 148)
\ifill f:0
\move(444 147)
\lvec(445 147)
\lvec(445 148)
\lvec(444 148)
\ifill f:0
\move(446 147)
\lvec(447 147)
\lvec(447 148)
\lvec(446 148)
\ifill f:0
\move(448 147)
\lvec(449 147)
\lvec(449 148)
\lvec(448 148)
\ifill f:0
\move(450 147)
\lvec(451 147)
\lvec(451 148)
\lvec(450 148)
\ifill f:0
\move(16 148)
\lvec(17 148)
\lvec(17 149)
\lvec(16 149)
\ifill f:0
\move(20 148)
\lvec(21 148)
\lvec(21 149)
\lvec(20 149)
\ifill f:0
\move(22 148)
\lvec(23 148)
\lvec(23 149)
\lvec(22 149)
\ifill f:0
\move(25 148)
\lvec(26 148)
\lvec(26 149)
\lvec(25 149)
\ifill f:0
\move(36 148)
\lvec(37 148)
\lvec(37 149)
\lvec(36 149)
\ifill f:0
\move(40 148)
\lvec(42 148)
\lvec(42 149)
\lvec(40 149)
\ifill f:0
\move(43 148)
\lvec(45 148)
\lvec(45 149)
\lvec(43 149)
\ifill f:0
\move(47 148)
\lvec(50 148)
\lvec(50 149)
\lvec(47 149)
\ifill f:0
\move(59 148)
\lvec(60 148)
\lvec(60 149)
\lvec(59 149)
\ifill f:0
\move(61 148)
\lvec(62 148)
\lvec(62 149)
\lvec(61 149)
\ifill f:0
\move(63 148)
\lvec(65 148)
\lvec(65 149)
\lvec(63 149)
\ifill f:0
\move(67 148)
\lvec(70 148)
\lvec(70 149)
\lvec(67 149)
\ifill f:0
\move(72 148)
\lvec(74 148)
\lvec(74 149)
\lvec(72 149)
\ifill f:0
\move(76 148)
\lvec(77 148)
\lvec(77 149)
\lvec(76 149)
\ifill f:0
\move(78 148)
\lvec(80 148)
\lvec(80 149)
\lvec(78 149)
\ifill f:0
\move(81 148)
\lvec(82 148)
\lvec(82 149)
\lvec(81 149)
\ifill f:0
\move(83 148)
\lvec(84 148)
\lvec(84 149)
\lvec(83 149)
\ifill f:0
\move(88 148)
\lvec(92 148)
\lvec(92 149)
\lvec(88 149)
\ifill f:0
\move(97 148)
\lvec(98 148)
\lvec(98 149)
\lvec(97 149)
\ifill f:0
\move(99 148)
\lvec(101 148)
\lvec(101 149)
\lvec(99 149)
\ifill f:0
\move(102 148)
\lvec(109 148)
\lvec(109 149)
\lvec(102 149)
\ifill f:0
\move(111 148)
\lvec(115 148)
\lvec(115 149)
\lvec(111 149)
\ifill f:0
\move(118 148)
\lvec(122 148)
\lvec(122 149)
\lvec(118 149)
\ifill f:0
\move(123 148)
\lvec(125 148)
\lvec(125 149)
\lvec(123 149)
\ifill f:0
\move(126 148)
\lvec(130 148)
\lvec(130 149)
\lvec(126 149)
\ifill f:0
\move(131 148)
\lvec(132 148)
\lvec(132 149)
\lvec(131 149)
\ifill f:0
\move(133 148)
\lvec(135 148)
\lvec(135 149)
\lvec(133 149)
\ifill f:0
\move(136 148)
\lvec(138 148)
\lvec(138 149)
\lvec(136 149)
\ifill f:0
\move(139 148)
\lvec(145 148)
\lvec(145 149)
\lvec(139 149)
\ifill f:0
\move(146 148)
\lvec(147 148)
\lvec(147 149)
\lvec(146 149)
\ifill f:0
\move(148 148)
\lvec(160 148)
\lvec(160 149)
\lvec(148 149)
\ifill f:0
\move(161 148)
\lvec(170 148)
\lvec(170 149)
\lvec(161 149)
\ifill f:0
\move(171 148)
\lvec(173 148)
\lvec(173 149)
\lvec(171 149)
\ifill f:0
\move(174 148)
\lvec(175 148)
\lvec(175 149)
\lvec(174 149)
\ifill f:0
\move(177 148)
\lvec(180 148)
\lvec(180 149)
\lvec(177 149)
\ifill f:0
\move(181 148)
\lvec(184 148)
\lvec(184 149)
\lvec(181 149)
\ifill f:0
\move(185 148)
\lvec(186 148)
\lvec(186 149)
\lvec(185 149)
\ifill f:0
\move(187 148)
\lvec(191 148)
\lvec(191 149)
\lvec(187 149)
\ifill f:0
\move(192 148)
\lvec(194 148)
\lvec(194 149)
\lvec(192 149)
\ifill f:0
\move(195 148)
\lvec(197 148)
\lvec(197 149)
\lvec(195 149)
\ifill f:0
\move(198 148)
\lvec(201 148)
\lvec(201 149)
\lvec(198 149)
\ifill f:0
\move(202 148)
\lvec(216 148)
\lvec(216 149)
\lvec(202 149)
\ifill f:0
\move(217 148)
\lvec(226 148)
\lvec(226 149)
\lvec(217 149)
\ifill f:0
\move(229 148)
\lvec(235 148)
\lvec(235 149)
\lvec(229 149)
\ifill f:0
\move(236 148)
\lvec(241 148)
\lvec(241 149)
\lvec(236 149)
\ifill f:0
\move(242 148)
\lvec(255 148)
\lvec(255 149)
\lvec(242 149)
\ifill f:0
\move(256 148)
\lvec(257 148)
\lvec(257 149)
\lvec(256 149)
\ifill f:0
\move(258 148)
\lvec(260 148)
\lvec(260 149)
\lvec(258 149)
\ifill f:0
\move(261 148)
\lvec(267 148)
\lvec(267 149)
\lvec(261 149)
\ifill f:0
\move(268 148)
\lvec(279 148)
\lvec(279 149)
\lvec(268 149)
\ifill f:0
\move(280 148)
\lvec(290 148)
\lvec(290 149)
\lvec(280 149)
\ifill f:0
\move(291 148)
\lvec(298 148)
\lvec(298 149)
\lvec(291 149)
\ifill f:0
\move(299 148)
\lvec(305 148)
\lvec(305 149)
\lvec(299 149)
\ifill f:0
\move(306 148)
\lvec(309 148)
\lvec(309 149)
\lvec(306 149)
\ifill f:0
\move(310 148)
\lvec(314 148)
\lvec(314 149)
\lvec(310 149)
\ifill f:0
\move(315 148)
\lvec(319 148)
\lvec(319 149)
\lvec(315 149)
\ifill f:0
\move(321 148)
\lvec(325 148)
\lvec(325 149)
\lvec(321 149)
\ifill f:0
\move(326 148)
\lvec(362 148)
\lvec(362 149)
\lvec(326 149)
\ifill f:0
\move(363 148)
\lvec(370 148)
\lvec(370 149)
\lvec(363 149)
\ifill f:0
\move(372 148)
\lvec(382 148)
\lvec(382 149)
\lvec(372 149)
\ifill f:0
\move(383 148)
\lvec(391 148)
\lvec(391 149)
\lvec(383 149)
\ifill f:0
\move(392 148)
\lvec(395 148)
\lvec(395 149)
\lvec(392 149)
\ifill f:0
\move(396 148)
\lvec(401 148)
\lvec(401 149)
\lvec(396 149)
\ifill f:0
\move(402 148)
\lvec(405 148)
\lvec(405 149)
\lvec(402 149)
\ifill f:0
\move(406 148)
\lvec(425 148)
\lvec(425 149)
\lvec(406 149)
\ifill f:0
\move(426 148)
\lvec(432 148)
\lvec(432 149)
\lvec(426 149)
\ifill f:0
\move(433 148)
\lvec(439 148)
\lvec(439 149)
\lvec(433 149)
\ifill f:0
\move(440 148)
\lvec(442 148)
\lvec(442 149)
\lvec(440 149)
\ifill f:0
\move(443 148)
\lvec(451 148)
\lvec(451 149)
\lvec(443 149)
\ifill f:0
\move(16 149)
\lvec(17 149)
\lvec(17 150)
\lvec(16 150)
\ifill f:0
\move(19 149)
\lvec(21 149)
\lvec(21 150)
\lvec(19 150)
\ifill f:0
\move(25 149)
\lvec(26 149)
\lvec(26 150)
\lvec(25 150)
\ifill f:0
\move(36 149)
\lvec(37 149)
\lvec(37 150)
\lvec(36 150)
\ifill f:0
\move(38 149)
\lvec(39 149)
\lvec(39 150)
\lvec(38 150)
\ifill f:0
\move(40 149)
\lvec(41 149)
\lvec(41 150)
\lvec(40 150)
\ifill f:0
\move(42 149)
\lvec(43 149)
\lvec(43 150)
\lvec(42 150)
\ifill f:0
\move(44 149)
\lvec(45 149)
\lvec(45 150)
\lvec(44 150)
\ifill f:0
\move(47 149)
\lvec(50 149)
\lvec(50 150)
\lvec(47 150)
\ifill f:0
\move(51 149)
\lvec(52 149)
\lvec(52 150)
\lvec(51 150)
\ifill f:0
\move(54 149)
\lvec(55 149)
\lvec(55 150)
\lvec(54 150)
\ifill f:0
\move(56 149)
\lvec(58 149)
\lvec(58 150)
\lvec(56 150)
\ifill f:0
\move(61 149)
\lvec(65 149)
\lvec(65 150)
\lvec(61 150)
\ifill f:0
\move(66 149)
\lvec(73 149)
\lvec(73 150)
\lvec(66 150)
\ifill f:0
\move(75 149)
\lvec(76 149)
\lvec(76 150)
\lvec(75 150)
\ifill f:0
\move(77 149)
\lvec(78 149)
\lvec(78 150)
\lvec(77 150)
\ifill f:0
\move(79 149)
\lvec(82 149)
\lvec(82 150)
\lvec(79 150)
\ifill f:0
\move(83 149)
\lvec(84 149)
\lvec(84 150)
\lvec(83 150)
\ifill f:0
\move(92 149)
\lvec(93 149)
\lvec(93 150)
\lvec(92 150)
\ifill f:0
\move(96 149)
\lvec(97 149)
\lvec(97 150)
\lvec(96 150)
\ifill f:0
\move(100 149)
\lvec(101 149)
\lvec(101 150)
\lvec(100 150)
\ifill f:0
\move(102 149)
\lvec(103 149)
\lvec(103 150)
\lvec(102 150)
\ifill f:0
\move(104 149)
\lvec(106 149)
\lvec(106 150)
\lvec(104 150)
\ifill f:0
\move(108 149)
\lvec(122 149)
\lvec(122 150)
\lvec(108 150)
\ifill f:0
\move(123 149)
\lvec(125 149)
\lvec(125 150)
\lvec(123 150)
\ifill f:0
\move(127 149)
\lvec(128 149)
\lvec(128 150)
\lvec(127 150)
\ifill f:0
\move(129 149)
\lvec(131 149)
\lvec(131 150)
\lvec(129 150)
\ifill f:0
\move(132 149)
\lvec(133 149)
\lvec(133 150)
\lvec(132 150)
\ifill f:0
\move(134 149)
\lvec(136 149)
\lvec(136 150)
\lvec(134 150)
\ifill f:0
\move(137 149)
\lvec(138 149)
\lvec(138 150)
\lvec(137 150)
\ifill f:0
\move(139 149)
\lvec(142 149)
\lvec(142 150)
\lvec(139 150)
\ifill f:0
\move(143 149)
\lvec(145 149)
\lvec(145 150)
\lvec(143 150)
\ifill f:0
\move(146 149)
\lvec(155 149)
\lvec(155 150)
\lvec(146 150)
\ifill f:0
\move(156 149)
\lvec(165 149)
\lvec(165 150)
\lvec(156 150)
\ifill f:0
\move(166 149)
\lvec(170 149)
\lvec(170 150)
\lvec(166 150)
\ifill f:0
\move(171 149)
\lvec(173 149)
\lvec(173 150)
\lvec(171 150)
\ifill f:0
\move(175 149)
\lvec(176 149)
\lvec(176 150)
\lvec(175 150)
\ifill f:0
\move(177 149)
\lvec(183 149)
\lvec(183 150)
\lvec(177 150)
\ifill f:0
\move(184 149)
\lvec(185 149)
\lvec(185 150)
\lvec(184 150)
\ifill f:0
\move(186 149)
\lvec(187 149)
\lvec(187 150)
\lvec(186 150)
\ifill f:0
\move(188 149)
\lvec(189 149)
\lvec(189 150)
\lvec(188 150)
\ifill f:0
\move(190 149)
\lvec(191 149)
\lvec(191 150)
\lvec(190 150)
\ifill f:0
\move(192 149)
\lvec(194 149)
\lvec(194 150)
\lvec(192 150)
\ifill f:0
\move(195 149)
\lvec(197 149)
\lvec(197 150)
\lvec(195 150)
\ifill f:0
\move(198 149)
\lvec(210 149)
\lvec(210 150)
\lvec(198 150)
\ifill f:0
\move(211 149)
\lvec(219 149)
\lvec(219 150)
\lvec(211 150)
\ifill f:0
\move(221 149)
\lvec(222 149)
\lvec(222 150)
\lvec(221 150)
\ifill f:0
\move(225 149)
\lvec(226 149)
\lvec(226 150)
\lvec(225 150)
\ifill f:0
\move(230 149)
\lvec(243 149)
\lvec(243 150)
\lvec(230 150)
\ifill f:0
\move(244 149)
\lvec(248 149)
\lvec(248 150)
\lvec(244 150)
\ifill f:0
\move(249 149)
\lvec(255 149)
\lvec(255 150)
\lvec(249 150)
\ifill f:0
\move(256 149)
\lvec(257 149)
\lvec(257 150)
\lvec(256 150)
\ifill f:0
\move(258 149)
\lvec(263 149)
\lvec(263 150)
\lvec(258 150)
\ifill f:0
\move(264 149)
\lvec(266 149)
\lvec(266 150)
\lvec(264 150)
\ifill f:0
\move(267 149)
\lvec(284 149)
\lvec(284 150)
\lvec(267 150)
\ifill f:0
\move(285 149)
\lvec(290 149)
\lvec(290 150)
\lvec(285 150)
\ifill f:0
\move(292 149)
\lvec(297 149)
\lvec(297 150)
\lvec(292 150)
\ifill f:0
\move(298 149)
\lvec(307 149)
\lvec(307 150)
\lvec(298 150)
\ifill f:0
\move(308 149)
\lvec(311 149)
\lvec(311 150)
\lvec(308 150)
\ifill f:0
\move(312 149)
\lvec(315 149)
\lvec(315 150)
\lvec(312 150)
\ifill f:0
\move(316 149)
\lvec(320 149)
\lvec(320 150)
\lvec(316 150)
\ifill f:0
\move(322 149)
\lvec(325 149)
\lvec(325 150)
\lvec(322 150)
\ifill f:0
\move(326 149)
\lvec(337 149)
\lvec(337 150)
\lvec(326 150)
\ifill f:0
\move(338 149)
\lvec(362 149)
\lvec(362 150)
\lvec(338 150)
\ifill f:0
\move(364 149)
\lvec(372 149)
\lvec(372 150)
\lvec(364 150)
\ifill f:0
\move(373 149)
\lvec(379 149)
\lvec(379 150)
\lvec(373 150)
\ifill f:0
\move(380 149)
\lvec(394 149)
\lvec(394 150)
\lvec(380 150)
\ifill f:0
\move(395 149)
\lvec(398 149)
\lvec(398 150)
\lvec(395 150)
\ifill f:0
\move(399 149)
\lvec(401 149)
\lvec(401 150)
\lvec(399 150)
\ifill f:0
\move(402 149)
\lvec(421 149)
\lvec(421 150)
\lvec(402 150)
\ifill f:0
\move(422 149)
\lvec(434 149)
\lvec(434 150)
\lvec(422 150)
\ifill f:0
\move(435 149)
\lvec(442 149)
\lvec(442 150)
\lvec(435 150)
\ifill f:0
\move(443 149)
\lvec(450 149)
\lvec(450 150)
\lvec(443 150)
\ifill f:0
\move(16 150)
\lvec(17 150)
\lvec(17 151)
\lvec(16 151)
\ifill f:0
\move(20 150)
\lvec(21 150)
\lvec(21 151)
\lvec(20 151)
\ifill f:0
\move(25 150)
\lvec(26 150)
\lvec(26 151)
\lvec(25 151)
\ifill f:0
\move(36 150)
\lvec(37 150)
\lvec(37 151)
\lvec(36 151)
\ifill f:0
\move(38 150)
\lvec(39 150)
\lvec(39 151)
\lvec(38 151)
\ifill f:0
\move(40 150)
\lvec(42 150)
\lvec(42 151)
\lvec(40 151)
\ifill f:0
\move(43 150)
\lvec(45 150)
\lvec(45 151)
\lvec(43 151)
\ifill f:0
\move(47 150)
\lvec(50 150)
\lvec(50 151)
\lvec(47 151)
\ifill f:0
\move(51 150)
\lvec(52 150)
\lvec(52 151)
\lvec(51 151)
\ifill f:0
\move(54 150)
\lvec(55 150)
\lvec(55 151)
\lvec(54 151)
\ifill f:0
\move(59 150)
\lvec(61 150)
\lvec(61 151)
\lvec(59 151)
\ifill f:0
\move(62 150)
\lvec(65 150)
\lvec(65 151)
\lvec(62 151)
\ifill f:0
\move(66 150)
\lvec(71 150)
\lvec(71 151)
\lvec(66 151)
\ifill f:0
\move(72 150)
\lvec(75 150)
\lvec(75 151)
\lvec(72 151)
\ifill f:0
\move(76 150)
\lvec(78 150)
\lvec(78 151)
\lvec(76 151)
\ifill f:0
\move(79 150)
\lvec(82 150)
\lvec(82 151)
\lvec(79 151)
\ifill f:0
\move(85 150)
\lvec(87 150)
\lvec(87 151)
\lvec(85 151)
\ifill f:0
\move(88 150)
\lvec(93 150)
\lvec(93 151)
\lvec(88 151)
\ifill f:0
\move(96 150)
\lvec(98 150)
\lvec(98 151)
\lvec(96 151)
\ifill f:0
\move(100 150)
\lvec(101 150)
\lvec(101 151)
\lvec(100 151)
\ifill f:0
\move(102 150)
\lvec(103 150)
\lvec(103 151)
\lvec(102 151)
\ifill f:0
\move(104 150)
\lvec(106 150)
\lvec(106 151)
\lvec(104 151)
\ifill f:0
\move(107 150)
\lvec(109 150)
\lvec(109 151)
\lvec(107 151)
\ifill f:0
\move(111 150)
\lvec(122 150)
\lvec(122 151)
\lvec(111 151)
\ifill f:0
\move(123 150)
\lvec(126 150)
\lvec(126 151)
\lvec(123 151)
\ifill f:0
\move(128 150)
\lvec(129 150)
\lvec(129 151)
\lvec(128 151)
\ifill f:0
\move(130 150)
\lvec(132 150)
\lvec(132 151)
\lvec(130 151)
\ifill f:0
\move(133 150)
\lvec(134 150)
\lvec(134 151)
\lvec(133 151)
\ifill f:0
\move(135 150)
\lvec(136 150)
\lvec(136 151)
\lvec(135 151)
\ifill f:0
\move(137 150)
\lvec(138 150)
\lvec(138 151)
\lvec(137 151)
\ifill f:0
\move(140 150)
\lvec(142 150)
\lvec(142 151)
\lvec(140 151)
\ifill f:0
\move(143 150)
\lvec(145 150)
\lvec(145 151)
\lvec(143 151)
\ifill f:0
\move(147 150)
\lvec(151 150)
\lvec(151 151)
\lvec(147 151)
\ifill f:0
\move(153 150)
\lvec(163 150)
\lvec(163 151)
\lvec(153 151)
\ifill f:0
\move(165 150)
\lvec(170 150)
\lvec(170 151)
\lvec(165 151)
\ifill f:0
\move(171 150)
\lvec(174 150)
\lvec(174 151)
\lvec(171 151)
\ifill f:0
\move(175 150)
\lvec(184 150)
\lvec(184 151)
\lvec(175 151)
\ifill f:0
\move(185 150)
\lvec(190 150)
\lvec(190 151)
\lvec(185 151)
\ifill f:0
\move(191 150)
\lvec(192 150)
\lvec(192 151)
\lvec(191 151)
\ifill f:0
\move(193 150)
\lvec(194 150)
\lvec(194 151)
\lvec(193 151)
\ifill f:0
\move(195 150)
\lvec(197 150)
\lvec(197 151)
\lvec(195 151)
\ifill f:0
\move(198 150)
\lvec(223 150)
\lvec(223 151)
\lvec(198 151)
\ifill f:0
\move(225 150)
\lvec(226 150)
\lvec(226 151)
\lvec(225 151)
\ifill f:0
\move(232 150)
\lvec(233 150)
\lvec(233 151)
\lvec(232 151)
\ifill f:0
\move(236 150)
\lvec(241 150)
\lvec(241 151)
\lvec(236 151)
\ifill f:0
\move(243 150)
\lvec(246 150)
\lvec(246 151)
\lvec(243 151)
\ifill f:0
\move(247 150)
\lvec(251 150)
\lvec(251 151)
\lvec(247 151)
\ifill f:0
\move(252 150)
\lvec(255 150)
\lvec(255 151)
\lvec(252 151)
\ifill f:0
\move(256 150)
\lvec(257 150)
\lvec(257 151)
\lvec(256 151)
\ifill f:0
\move(258 150)
\lvec(269 150)
\lvec(269 151)
\lvec(258 151)
\ifill f:0
\move(270 150)
\lvec(280 150)
\lvec(280 151)
\lvec(270 151)
\ifill f:0
\move(281 150)
\lvec(290 150)
\lvec(290 151)
\lvec(281 151)
\ifill f:0
\move(292 150)
\lvec(299 150)
\lvec(299 151)
\lvec(292 151)
\ifill f:0
\move(300 150)
\lvec(302 150)
\lvec(302 151)
\lvec(300 151)
\ifill f:0
\move(303 150)
\lvec(305 150)
\lvec(305 151)
\lvec(303 151)
\ifill f:0
\move(306 150)
\lvec(321 150)
\lvec(321 151)
\lvec(306 151)
\ifill f:0
\move(322 150)
\lvec(325 150)
\lvec(325 151)
\lvec(322 151)
\ifill f:0
\move(326 150)
\lvec(327 150)
\lvec(327 151)
\lvec(326 151)
\ifill f:0
\move(328 150)
\lvec(334 150)
\lvec(334 151)
\lvec(328 151)
\ifill f:0
\move(335 150)
\lvec(346 150)
\lvec(346 151)
\lvec(335 151)
\ifill f:0
\move(348 150)
\lvec(362 150)
\lvec(362 151)
\lvec(348 151)
\ifill f:0
\move(365 150)
\lvec(375 150)
\lvec(375 151)
\lvec(365 151)
\ifill f:0
\move(376 150)
\lvec(382 150)
\lvec(382 151)
\lvec(376 151)
\ifill f:0
\move(384 150)
\lvec(388 150)
\lvec(388 151)
\lvec(384 151)
\ifill f:0
\move(389 150)
\lvec(394 150)
\lvec(394 151)
\lvec(389 151)
\ifill f:0
\move(395 150)
\lvec(398 150)
\lvec(398 151)
\lvec(395 151)
\ifill f:0
\move(399 150)
\lvec(401 150)
\lvec(401 151)
\lvec(399 151)
\ifill f:0
\move(403 150)
\lvec(406 150)
\lvec(406 151)
\lvec(403 151)
\ifill f:0
\move(407 150)
\lvec(431 150)
\lvec(431 151)
\lvec(407 151)
\ifill f:0
\move(432 150)
\lvec(442 150)
\lvec(442 151)
\lvec(432 151)
\ifill f:0
\move(443 150)
\lvec(448 150)
\lvec(448 151)
\lvec(443 151)
\ifill f:0
\move(449 150)
\lvec(451 150)
\lvec(451 151)
\lvec(449 151)
\ifill f:0
\move(30 -20)
\htext{number of points $r \in [10,450]$}
\move(-2 20)
\vtext{multiplicity $m \in [1,150]$}
\arrowheadtype t:V
\linewd 1
\move(10 -5)
\avec(470 -5)
\move(5 0)
\avec(5 170)
}
\caption{Graph showing when {\sc NSglue} gives conjectured value of $\tau_2$}\label{pic1}
\end{figure}

\begin{figure}[ht!]
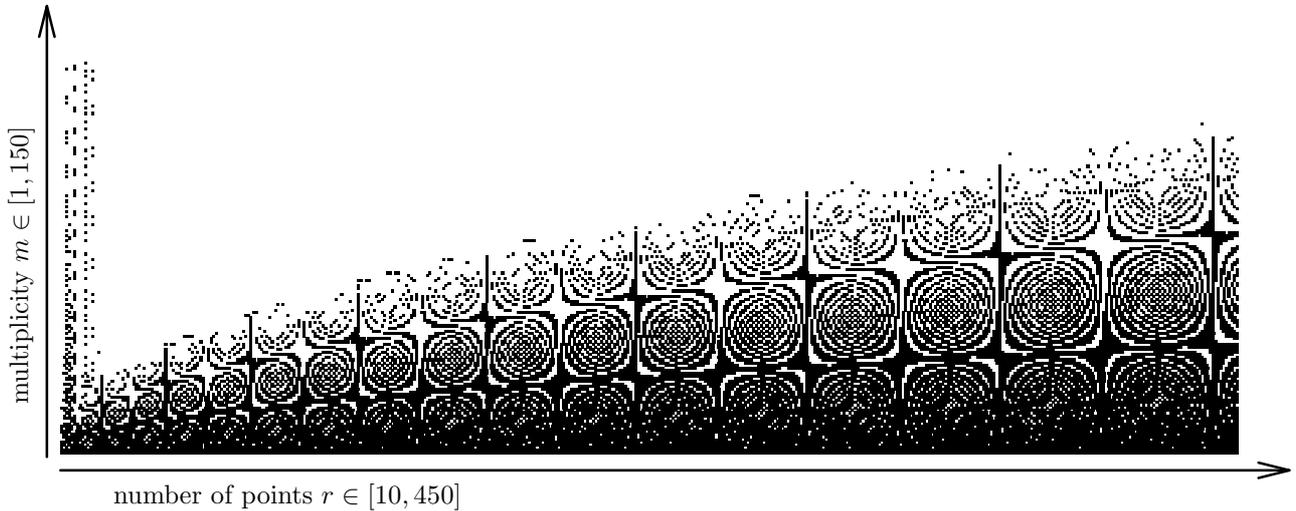

\centertexdraw{
\drawdim pt
\move(10 1)
\lvec(451 1)
\lvec(451 2)
\lvec(10 2)
\ifill f:0
\move(10 2)
\lvec(451 2)
\lvec(451 3)
\lvec(10 3)
\ifill f:0
\move(10 3)
\lvec(11 3)
\lvec(11 4)
\lvec(10 4)
\ifill f:0
\move(12 3)
\lvec(13 3)
\lvec(13 4)
\lvec(12 4)
\ifill f:0
\move(14 3)
\lvec(20 3)
\lvec(20 4)
\lvec(14 4)
\ifill f:0
\move(21 3)
\lvec(35 3)
\lvec(35 4)
\lvec(21 4)
\ifill f:0
\move(36 3)
\lvec(46 3)
\lvec(46 4)
\lvec(36 4)
\ifill f:0
\move(47 3)
\lvec(50 3)
\lvec(50 4)
\lvec(47 4)
\ifill f:0
\move(51 3)
\lvec(63 3)
\lvec(63 4)
\lvec(51 4)
\ifill f:0
\move(64 3)
\lvec(88 3)
\lvec(88 4)
\lvec(64 4)
\ifill f:0
\move(89 3)
\lvec(105 3)
\lvec(105 4)
\lvec(89 4)
\ifill f:0
\move(106 3)
\lvec(111 3)
\lvec(111 4)
\lvec(106 4)
\ifill f:0
\move(112 3)
\lvec(130 3)
\lvec(130 4)
\lvec(112 4)
\ifill f:0
\move(131 3)
\lvec(165 3)
\lvec(165 4)
\lvec(131 4)
\ifill f:0
\move(166 3)
\lvec(188 3)
\lvec(188 4)
\lvec(166 4)
\ifill f:0
\move(189 3)
\lvec(196 3)
\lvec(196 4)
\lvec(189 4)
\ifill f:0
\move(197 3)
\lvec(221 3)
\lvec(221 4)
\lvec(197 4)
\ifill f:0
\move(222 3)
\lvec(266 3)
\lvec(266 4)
\lvec(222 4)
\ifill f:0
\move(267 3)
\lvec(295 3)
\lvec(295 4)
\lvec(267 4)
\ifill f:0
\move(296 3)
\lvec(305 3)
\lvec(305 4)
\lvec(296 4)
\ifill f:0
\move(306 3)
\lvec(336 3)
\lvec(336 4)
\lvec(306 4)
\ifill f:0
\move(337 3)
\lvec(391 3)
\lvec(391 4)
\lvec(337 4)
\ifill f:0
\move(392 3)
\lvec(426 3)
\lvec(426 4)
\lvec(392 4)
\ifill f:0
\move(427 3)
\lvec(438 3)
\lvec(438 4)
\lvec(427 4)
\ifill f:0
\move(439 3)
\lvec(451 3)
\lvec(451 4)
\lvec(439 4)
\ifill f:0
\move(10 4)
\lvec(12 4)
\lvec(12 5)
\lvec(10 5)
\ifill f:0
\move(13 4)
\lvec(19 4)
\lvec(19 5)
\lvec(13 5)
\ifill f:0
\move(20 4)
\lvec(21 4)
\lvec(21 5)
\lvec(20 5)
\ifill f:0
\move(22 4)
\lvec(30 4)
\lvec(30 5)
\lvec(22 5)
\ifill f:0
\move(31 4)
\lvec(63 4)
\lvec(63 5)
\lvec(31 5)
\ifill f:0
\move(64 4)
\lvec(78 4)
\lvec(78 5)
\lvec(64 5)
\ifill f:0
\move(79 4)
\lvec(82 4)
\lvec(82 5)
\lvec(79 5)
\ifill f:0
\move(83 4)
\lvec(99 4)
\lvec(99 5)
\lvec(83 5)
\ifill f:0
\move(100 4)
\lvec(154 4)
\lvec(154 5)
\lvec(100 5)
\ifill f:0
\move(155 4)
\lvec(177 4)
\lvec(177 5)
\lvec(155 5)
\ifill f:0
\move(178 4)
\lvec(183 4)
\lvec(183 5)
\lvec(178 5)
\ifill f:0
\move(184 4)
\lvec(208 4)
\lvec(208 5)
\lvec(184 5)
\ifill f:0
\move(209 4)
\lvec(285 4)
\lvec(285 5)
\lvec(209 5)
\ifill f:0
\move(286 4)
\lvec(316 4)
\lvec(316 5)
\lvec(286 5)
\ifill f:0
\move(317 4)
\lvec(324 4)
\lvec(324 5)
\lvec(317 5)
\ifill f:0
\move(325 4)
\lvec(357 4)
\lvec(357 5)
\lvec(325 5)
\ifill f:0
\move(358 4)
\lvec(451 4)
\lvec(451 5)
\lvec(358 5)
\ifill f:0
\move(10 5)
\lvec(14 5)
\lvec(14 6)
\lvec(10 6)
\ifill f:0
\move(15 5)
\lvec(17 5)
\lvec(17 6)
\lvec(15 6)
\ifill f:0
\move(18 5)
\lvec(20 5)
\lvec(20 6)
\lvec(18 6)
\ifill f:0
\move(21 5)
\lvec(29 5)
\lvec(29 6)
\lvec(21 6)
\ifill f:0
\move(30 5)
\lvec(31 5)
\lvec(31 6)
\lvec(30 6)
\ifill f:0
\move(32 5)
\lvec(42 5)
\lvec(42 6)
\lvec(32 6)
\ifill f:0
\move(43 5)
\lvec(47 5)
\lvec(47 6)
\lvec(43 6)
\ifill f:0
\move(48 5)
\lvec(52 5)
\lvec(52 6)
\lvec(48 6)
\ifill f:0
\move(53 5)
\lvec(66 5)
\lvec(66 6)
\lvec(53 6)
\ifill f:0
\move(67 5)
\lvec(69 5)
\lvec(69 6)
\lvec(67 6)
\ifill f:0
\move(70 5)
\lvec(85 5)
\lvec(85 6)
\lvec(70 6)
\ifill f:0
\move(86 5)
\lvec(92 5)
\lvec(92 6)
\lvec(86 6)
\ifill f:0
\move(93 5)
\lvec(99 5)
\lvec(99 6)
\lvec(93 6)
\ifill f:0
\move(100 5)
\lvec(118 5)
\lvec(118 6)
\lvec(100 6)
\ifill f:0
\move(119 5)
\lvec(122 5)
\lvec(122 6)
\lvec(119 6)
\ifill f:0
\move(123 5)
\lvec(143 5)
\lvec(143 6)
\lvec(123 6)
\ifill f:0
\move(144 5)
\lvec(152 5)
\lvec(152 6)
\lvec(144 6)
\ifill f:0
\move(153 5)
\lvec(161 5)
\lvec(161 6)
\lvec(153 6)
\ifill f:0
\move(162 5)
\lvec(185 5)
\lvec(185 6)
\lvec(162 6)
\ifill f:0
\move(186 5)
\lvec(190 5)
\lvec(190 6)
\lvec(186 6)
\ifill f:0
\move(191 5)
\lvec(216 5)
\lvec(216 6)
\lvec(191 6)
\ifill f:0
\move(217 5)
\lvec(227 5)
\lvec(227 6)
\lvec(217 6)
\ifill f:0
\move(228 5)
\lvec(238 5)
\lvec(238 6)
\lvec(228 6)
\ifill f:0
\move(239 5)
\lvec(267 5)
\lvec(267 6)
\lvec(239 6)
\ifill f:0
\move(268 5)
\lvec(273 5)
\lvec(273 6)
\lvec(268 6)
\ifill f:0
\move(274 5)
\lvec(304 5)
\lvec(304 6)
\lvec(274 6)
\ifill f:0
\move(305 5)
\lvec(317 5)
\lvec(317 6)
\lvec(305 6)
\ifill f:0
\move(318 5)
\lvec(330 5)
\lvec(330 6)
\lvec(318 6)
\ifill f:0
\move(331 5)
\lvec(364 5)
\lvec(364 6)
\lvec(331 6)
\ifill f:0
\move(365 5)
\lvec(371 5)
\lvec(371 6)
\lvec(365 6)
\ifill f:0
\move(372 5)
\lvec(407 5)
\lvec(407 6)
\lvec(372 6)
\ifill f:0
\move(408 5)
\lvec(422 5)
\lvec(422 6)
\lvec(408 6)
\ifill f:0
\move(423 5)
\lvec(437 5)
\lvec(437 6)
\lvec(423 6)
\ifill f:0
\move(438 5)
\lvec(451 5)
\lvec(451 6)
\lvec(438 6)
\ifill f:0
\move(12 6)
\lvec(18 6)
\lvec(18 7)
\lvec(12 7)
\ifill f:0
\move(19 6)
\lvec(21 6)
\lvec(21 7)
\lvec(19 7)
\ifill f:0
\move(22 6)
\lvec(30 6)
\lvec(30 7)
\lvec(22 7)
\ifill f:0
\move(31 6)
\lvec(41 6)
\lvec(41 7)
\lvec(31 7)
\ifill f:0
\move(42 6)
\lvec(43 6)
\lvec(43 7)
\lvec(42 7)
\ifill f:0
\move(44 6)
\lvec(56 6)
\lvec(56 7)
\lvec(44 7)
\ifill f:0
\move(57 6)
\lvec(76 6)
\lvec(76 7)
\lvec(57 7)
\ifill f:0
\move(77 6)
\lvec(93 6)
\lvec(93 7)
\lvec(77 7)
\ifill f:0
\move(94 6)
\lvec(96 6)
\lvec(96 7)
\lvec(94 7)
\ifill f:0
\move(97 6)
\lvec(115 6)
\lvec(115 7)
\lvec(97 7)
\ifill f:0
\move(116 6)
\lvec(143 6)
\lvec(143 7)
\lvec(116 7)
\ifill f:0
\move(144 6)
\lvec(166 6)
\lvec(166 7)
\lvec(144 7)
\ifill f:0
\move(167 6)
\lvec(170 6)
\lvec(170 7)
\lvec(167 7)
\ifill f:0
\move(171 6)
\lvec(195 6)
\lvec(195 7)
\lvec(171 7)
\ifill f:0
\move(196 6)
\lvec(231 6)
\lvec(231 7)
\lvec(196 7)
\ifill f:0
\move(232 6)
\lvec(260 6)
\lvec(260 7)
\lvec(232 7)
\ifill f:0
\move(261 6)
\lvec(265 6)
\lvec(265 7)
\lvec(261 7)
\ifill f:0
\move(266 6)
\lvec(296 6)
\lvec(296 7)
\lvec(266 7)
\ifill f:0
\move(297 6)
\lvec(340 6)
\lvec(340 7)
\lvec(297 7)
\ifill f:0
\move(341 6)
\lvec(375 6)
\lvec(375 7)
\lvec(341 7)
\ifill f:0
\move(376 6)
\lvec(381 6)
\lvec(381 7)
\lvec(376 7)
\ifill f:0
\move(382 6)
\lvec(418 6)
\lvec(418 7)
\lvec(382 7)
\ifill f:0
\move(419 6)
\lvec(451 6)
\lvec(451 7)
\lvec(419 7)
\ifill f:0
\move(11 7)
\lvec(19 7)
\lvec(19 8)
\lvec(11 8)
\ifill f:0
\move(20 7)
\lvec(28 7)
\lvec(28 8)
\lvec(20 8)
\ifill f:0
\move(29 7)
\lvec(37 7)
\lvec(37 8)
\lvec(29 8)
\ifill f:0
\move(38 7)
\lvec(42 7)
\lvec(42 8)
\lvec(38 8)
\ifill f:0
\move(43 7)
\lvec(55 7)
\lvec(55 8)
\lvec(43 8)
\ifill f:0
\move(56 7)
\lvec(57 7)
\lvec(57 8)
\lvec(56 8)
\ifill f:0
\move(58 7)
\lvec(72 7)
\lvec(72 8)
\lvec(58 8)
\ifill f:0
\move(73 7)
\lvec(79 7)
\lvec(79 8)
\lvec(73 8)
\ifill f:0
\move(80 7)
\lvec(94 7)
\lvec(94 8)
\lvec(80 8)
\ifill f:0
\move(95 7)
\lvec(113 7)
\lvec(113 8)
\lvec(95 8)
\ifill f:0
\move(114 7)
\lvec(140 7)
\lvec(140 8)
\lvec(114 8)
\ifill f:0
\move(141 7)
\lvec(163 7)
\lvec(163 8)
\lvec(141 8)
\ifill f:0
\move(164 7)
\lvec(184 7)
\lvec(184 8)
\lvec(164 8)
\ifill f:0
\move(185 7)
\lvec(195 7)
\lvec(195 8)
\lvec(185 8)
\ifill f:0
\move(196 7)
\lvec(222 7)
\lvec(222 8)
\lvec(196 8)
\ifill f:0
\move(223 7)
\lvec(226 7)
\lvec(226 8)
\lvec(223 8)
\ifill f:0
\move(227 7)
\lvec(255 7)
\lvec(255 8)
\lvec(227 8)
\ifill f:0
\move(256 7)
\lvec(268 7)
\lvec(268 8)
\lvec(256 8)
\ifill f:0
\move(269 7)
\lvec(295 7)
\lvec(295 8)
\lvec(269 8)
\ifill f:0
\move(296 7)
\lvec(328 7)
\lvec(328 8)
\lvec(296 8)
\ifill f:0
\move(329 7)
\lvec(373 7)
\lvec(373 8)
\lvec(329 8)
\ifill f:0
\move(374 7)
\lvec(410 7)
\lvec(410 8)
\lvec(374 8)
\ifill f:0
\move(411 7)
\lvec(443 7)
\lvec(443 8)
\lvec(411 8)
\ifill f:0
\move(444 7)
\lvec(451 7)
\lvec(451 8)
\lvec(444 8)
\ifill f:0
\move(10 9)
\lvec(14 9)
\lvec(14 10)
\lvec(10 10)
\ifill f:0
\move(16 9)
\lvec(22 9)
\lvec(22 10)
\lvec(16 10)
\ifill f:0
\move(24 9)
\lvec(32 9)
\lvec(32 10)
\lvec(24 10)
\ifill f:0
\move(34 9)
\lvec(45 9)
\lvec(45 10)
\lvec(34 10)
\ifill f:0
\move(46 9)
\lvec(57 9)
\lvec(57 10)
\lvec(46 10)
\ifill f:0
\move(58 9)
\lvec(67 9)
\lvec(67 10)
\lvec(58 10)
\ifill f:0
\move(68 9)
\lvec(72 9)
\lvec(72 10)
\lvec(68 10)
\ifill f:0
\move(73 9)
\lvec(74 9)
\lvec(74 10)
\lvec(73 10)
\ifill f:0
\move(75 9)
\lvec(89 9)
\lvec(89 10)
\lvec(75 10)
\ifill f:0
\move(90 9)
\lvec(91 9)
\lvec(91 10)
\lvec(90 10)
\ifill f:0
\move(92 9)
\lvec(108 9)
\lvec(108 10)
\lvec(92 10)
\ifill f:0
\move(109 9)
\lvec(110 9)
\lvec(110 10)
\lvec(109 10)
\ifill f:0
\move(111 9)
\lvec(117 9)
\lvec(117 10)
\lvec(111 10)
\ifill f:0
\move(118 9)
\lvec(131 9)
\lvec(131 10)
\lvec(118 10)
\ifill f:0
\move(132 9)
\lvec(151 9)
\lvec(151 10)
\lvec(132 10)
\ifill f:0
\move(152 9)
\lvec(167 9)
\lvec(167 10)
\lvec(152 10)
\ifill f:0
\move(168 9)
\lvec(175 9)
\lvec(175 10)
\lvec(168 10)
\ifill f:0
\move(176 9)
\lvec(178 9)
\lvec(178 10)
\lvec(176 10)
\ifill f:0
\move(179 9)
\lvec(201 9)
\lvec(201 10)
\lvec(179 10)
\ifill f:0
\move(202 9)
\lvec(204 9)
\lvec(204 10)
\lvec(202 10)
\ifill f:0
\move(205 9)
\lvec(229 9)
\lvec(229 10)
\lvec(205 10)
\ifill f:0
\move(230 9)
\lvec(232 9)
\lvec(232 10)
\lvec(230 10)
\ifill f:0
\move(233 9)
\lvec(242 9)
\lvec(242 10)
\lvec(233 10)
\ifill f:0
\move(243 9)
\lvec(262 9)
\lvec(262 10)
\lvec(243 10)
\ifill f:0
\move(263 9)
\lvec(290 9)
\lvec(290 10)
\lvec(263 10)
\ifill f:0
\move(291 9)
\lvec(312 9)
\lvec(312 10)
\lvec(291 10)
\ifill f:0
\move(313 9)
\lvec(323 9)
\lvec(323 10)
\lvec(313 10)
\ifill f:0
\move(324 9)
\lvec(327 9)
\lvec(327 10)
\lvec(324 10)
\ifill f:0
\move(328 9)
\lvec(358 9)
\lvec(358 10)
\lvec(328 10)
\ifill f:0
\move(359 9)
\lvec(362 9)
\lvec(362 10)
\lvec(359 10)
\ifill f:0
\move(363 9)
\lvec(395 9)
\lvec(395 10)
\lvec(363 10)
\ifill f:0
\move(396 9)
\lvec(399 9)
\lvec(399 10)
\lvec(396 10)
\ifill f:0
\move(400 9)
\lvec(412 9)
\lvec(412 10)
\lvec(400 10)
\ifill f:0
\move(413 9)
\lvec(438 9)
\lvec(438 10)
\lvec(413 10)
\ifill f:0
\move(439 9)
\lvec(451 9)
\lvec(451 10)
\lvec(439 10)
\ifill f:0
\move(10 10)
\lvec(11 10)
\lvec(11 11)
\lvec(10 11)
\ifill f:0
\move(12 10)
\lvec(15 10)
\lvec(15 11)
\lvec(12 11)
\ifill f:0
\move(16 10)
\lvec(18 10)
\lvec(18 11)
\lvec(16 11)
\ifill f:0
\move(20 10)
\lvec(27 10)
\lvec(27 11)
\lvec(20 11)
\ifill f:0
\move(29 10)
\lvec(38 10)
\lvec(38 11)
\lvec(29 11)
\ifill f:0
\move(40 10)
\lvec(44 10)
\lvec(44 11)
\lvec(40 11)
\ifill f:0
\move(45 10)
\lvec(48 10)
\lvec(48 11)
\lvec(45 11)
\ifill f:0
\move(49 10)
\lvec(52 10)
\lvec(52 11)
\lvec(49 11)
\ifill f:0
\move(53 10)
\lvec(59 10)
\lvec(59 11)
\lvec(53 11)
\ifill f:0
\move(60 10)
\lvec(62 10)
\lvec(62 11)
\lvec(60 11)
\ifill f:0
\move(63 10)
\lvec(65 10)
\lvec(65 11)
\lvec(63 11)
\ifill f:0
\move(66 10)
\lvec(73 10)
\lvec(73 11)
\lvec(66 11)
\ifill f:0
\move(74 10)
\lvec(78 10)
\lvec(78 11)
\lvec(74 11)
\ifill f:0
\move(79 10)
\lvec(83 10)
\lvec(83 11)
\lvec(79 11)
\ifill f:0
\move(84 10)
\lvec(90 10)
\lvec(90 11)
\lvec(84 11)
\ifill f:0
\move(91 10)
\lvec(92 10)
\lvec(92 11)
\lvec(91 11)
\ifill f:0
\move(93 10)
\lvec(109 10)
\lvec(109 11)
\lvec(93 11)
\ifill f:0
\move(110 10)
\lvec(111 10)
\lvec(111 11)
\lvec(110 11)
\ifill f:0
\move(112 10)
\lvec(130 10)
\lvec(130 11)
\lvec(112 11)
\ifill f:0
\move(131 10)
\lvec(132 10)
\lvec(132 11)
\lvec(131 11)
\ifill f:0
\move(133 10)
\lvec(141 10)
\lvec(141 11)
\lvec(133 11)
\ifill f:0
\move(142 10)
\lvec(148 10)
\lvec(148 11)
\lvec(142 11)
\ifill f:0
\move(149 10)
\lvec(155 10)
\lvec(155 11)
\lvec(149 11)
\ifill f:0
\move(156 10)
\lvec(167 10)
\lvec(167 11)
\lvec(156 11)
\ifill f:0
\move(168 10)
\lvec(172 10)
\lvec(172 11)
\lvec(168 11)
\ifill f:0
\move(173 10)
\lvec(177 10)
\lvec(177 11)
\lvec(173 11)
\ifill f:0
\move(178 10)
\lvec(190 10)
\lvec(190 11)
\lvec(178 11)
\ifill f:0
\move(191 10)
\lvec(198 10)
\lvec(198 11)
\lvec(191 11)
\ifill f:0
\move(199 10)
\lvec(206 10)
\lvec(206 11)
\lvec(199 11)
\ifill f:0
\move(207 10)
\lvec(217 10)
\lvec(217 11)
\lvec(207 11)
\ifill f:0
\move(218 10)
\lvec(220 10)
\lvec(220 11)
\lvec(218 11)
\ifill f:0
\move(221 10)
\lvec(246 10)
\lvec(246 11)
\lvec(221 11)
\ifill f:0
\move(247 10)
\lvec(249 10)
\lvec(249 11)
\lvec(247 11)
\ifill f:0
\move(250 10)
\lvec(277 10)
\lvec(277 11)
\lvec(250 11)
\ifill f:0
\move(278 10)
\lvec(280 10)
\lvec(280 11)
\lvec(278 11)
\ifill f:0
\move(281 10)
\lvec(293 10)
\lvec(293 11)
\lvec(281 11)
\ifill f:0
\move(294 10)
\lvec(303 10)
\lvec(303 11)
\lvec(294 11)
\ifill f:0
\move(304 10)
\lvec(313 10)
\lvec(313 11)
\lvec(304 11)
\ifill f:0
\move(314 10)
\lvec(330 10)
\lvec(330 11)
\lvec(314 11)
\ifill f:0
\move(331 10)
\lvec(337 10)
\lvec(337 11)
\lvec(331 11)
\ifill f:0
\move(338 10)
\lvec(344 10)
\lvec(344 11)
\lvec(338 11)
\ifill f:0
\move(345 10)
\lvec(362 10)
\lvec(362 11)
\lvec(345 11)
\ifill f:0
\move(363 10)
\lvec(373 10)
\lvec(373 11)
\lvec(363 11)
\ifill f:0
\move(374 10)
\lvec(384 10)
\lvec(384 11)
\lvec(374 11)
\ifill f:0
\move(385 10)
\lvec(399 10)
\lvec(399 11)
\lvec(385 11)
\ifill f:0
\move(400 10)
\lvec(403 10)
\lvec(403 11)
\lvec(400 11)
\ifill f:0
\move(404 10)
\lvec(438 10)
\lvec(438 11)
\lvec(404 11)
\ifill f:0
\move(439 10)
\lvec(442 10)
\lvec(442 11)
\lvec(439 11)
\ifill f:0
\move(443 10)
\lvec(451 10)
\lvec(451 11)
\lvec(443 11)
\ifill f:0
\move(10 11)
\lvec(15 11)
\lvec(15 12)
\lvec(10 12)
\ifill f:0
\move(16 11)
\lvec(18 11)
\lvec(18 12)
\lvec(16 12)
\ifill f:0
\move(19 11)
\lvec(21 11)
\lvec(21 12)
\lvec(19 12)
\ifill f:0
\move(22 11)
\lvec(26 11)
\lvec(26 12)
\lvec(22 12)
\ifill f:0
\move(29 11)
\lvec(39 11)
\lvec(39 12)
\lvec(29 12)
\ifill f:0
\move(42 11)
\lvec(48 11)
\lvec(48 12)
\lvec(42 12)
\ifill f:0
\move(49 11)
\lvec(53 11)
\lvec(53 12)
\lvec(49 12)
\ifill f:0
\move(54 11)
\lvec(58 11)
\lvec(58 12)
\lvec(54 12)
\ifill f:0
\move(59 11)
\lvec(65 11)
\lvec(65 12)
\lvec(59 12)
\ifill f:0
\move(66 11)
\lvec(75 11)
\lvec(75 12)
\lvec(66 12)
\ifill f:0
\move(76 11)
\lvec(83 11)
\lvec(83 12)
\lvec(76 12)
\ifill f:0
\move(84 11)
\lvec(86 11)
\lvec(86 12)
\lvec(84 12)
\ifill f:0
\move(87 11)
\lvec(91 11)
\lvec(91 12)
\lvec(87 12)
\ifill f:0
\move(92 11)
\lvec(96 11)
\lvec(96 12)
\lvec(92 12)
\ifill f:0
\move(97 11)
\lvec(103 11)
\lvec(103 12)
\lvec(97 12)
\ifill f:0
\move(104 11)
\lvec(110 11)
\lvec(110 12)
\lvec(104 12)
\ifill f:0
\move(111 11)
\lvec(112 11)
\lvec(112 12)
\lvec(111 12)
\ifill f:0
\move(113 11)
\lvec(131 11)
\lvec(131 12)
\lvec(113 12)
\ifill f:0
\move(132 11)
\lvec(133 11)
\lvec(133 12)
\lvec(132 12)
\ifill f:0
\move(134 11)
\lvec(154 11)
\lvec(154 12)
\lvec(134 12)
\ifill f:0
\move(155 11)
\lvec(156 11)
\lvec(156 12)
\lvec(155 12)
\ifill f:0
\move(157 11)
\lvec(165 11)
\lvec(165 12)
\lvec(157 12)
\ifill f:0
\move(166 11)
\lvec(174 11)
\lvec(174 12)
\lvec(166 12)
\ifill f:0
\move(175 11)
\lvec(181 11)
\lvec(181 12)
\lvec(175 12)
\ifill f:0
\move(182 11)
\lvec(188 11)
\lvec(188 12)
\lvec(182 12)
\ifill f:0
\move(189 11)
\lvec(193 11)
\lvec(193 12)
\lvec(189 12)
\ifill f:0
\move(194 11)
\lvec(205 11)
\lvec(205 12)
\lvec(194 12)
\ifill f:0
\move(206 11)
\lvec(223 11)
\lvec(223 12)
\lvec(206 12)
\ifill f:0
\move(224 11)
\lvec(236 11)
\lvec(236 12)
\lvec(224 12)
\ifill f:0
\move(237 11)
\lvec(247 11)
\lvec(247 12)
\lvec(237 12)
\ifill f:0
\move(248 11)
\lvec(258 11)
\lvec(258 12)
\lvec(248 12)
\ifill f:0
\move(259 11)
\lvec(275 11)
\lvec(275 12)
\lvec(259 12)
\ifill f:0
\move(276 11)
\lvec(278 11)
\lvec(278 12)
\lvec(276 12)
\ifill f:0
\move(279 11)
\lvec(281 11)
\lvec(281 12)
\lvec(279 12)
\ifill f:0
\move(282 11)
\lvec(314 11)
\lvec(314 12)
\lvec(282 12)
\ifill f:0
\move(315 11)
\lvec(317 11)
\lvec(317 12)
\lvec(315 12)
\ifill f:0
\move(318 11)
\lvec(320 11)
\lvec(320 12)
\lvec(318 12)
\ifill f:0
\move(321 11)
\lvec(339 11)
\lvec(339 12)
\lvec(321 12)
\ifill f:0
\move(340 11)
\lvec(352 11)
\lvec(352 12)
\lvec(340 12)
\ifill f:0
\move(353 11)
\lvec(365 11)
\lvec(365 12)
\lvec(353 12)
\ifill f:0
\move(366 11)
\lvec(382 11)
\lvec(382 12)
\lvec(366 12)
\ifill f:0
\move(383 11)
\lvec(406 11)
\lvec(406 12)
\lvec(383 12)
\ifill f:0
\move(407 11)
\lvec(424 11)
\lvec(424 12)
\lvec(407 12)
\ifill f:0
\move(425 11)
\lvec(431 11)
\lvec(431 12)
\lvec(425 12)
\ifill f:0
\move(432 11)
\lvec(442 11)
\lvec(442 12)
\lvec(432 12)
\ifill f:0
\move(443 11)
\lvec(451 11)
\lvec(451 12)
\lvec(443 12)
\ifill f:0
\move(11 12)
\lvec(14 12)
\lvec(14 13)
\lvec(11 13)
\ifill f:0
\move(15 12)
\lvec(17 12)
\lvec(17 13)
\lvec(15 13)
\ifill f:0
\move(18 12)
\lvec(22 12)
\lvec(22 13)
\lvec(18 13)
\ifill f:0
\move(23 12)
\lvec(26 12)
\lvec(26 13)
\lvec(23 13)
\ifill f:0
\move(27 12)
\lvec(31 12)
\lvec(31 13)
\lvec(27 13)
\ifill f:0
\move(34 12)
\lvec(46 12)
\lvec(46 13)
\lvec(34 13)
\ifill f:0
\move(49 12)
\lvec(55 12)
\lvec(55 13)
\lvec(49 13)
\ifill f:0
\move(56 12)
\lvec(61 12)
\lvec(61 13)
\lvec(56 13)
\ifill f:0
\move(62 12)
\lvec(70 12)
\lvec(70 13)
\lvec(62 13)
\ifill f:0
\move(71 12)
\lvec(77 12)
\lvec(77 13)
\lvec(71 13)
\ifill f:0
\move(78 12)
\lvec(87 12)
\lvec(87 13)
\lvec(78 13)
\ifill f:0
\move(88 12)
\lvec(98 12)
\lvec(98 13)
\lvec(88 13)
\ifill f:0
\move(99 12)
\lvec(101 12)
\lvec(101 13)
\lvec(99 13)
\ifill f:0
\move(102 12)
\lvec(106 12)
\lvec(106 13)
\lvec(102 13)
\ifill f:0
\move(107 12)
\lvec(111 12)
\lvec(111 13)
\lvec(107 13)
\ifill f:0
\move(112 12)
\lvec(116 12)
\lvec(116 13)
\lvec(112 13)
\ifill f:0
\move(117 12)
\lvec(123 12)
\lvec(123 13)
\lvec(117 13)
\ifill f:0
\move(124 12)
\lvec(125 12)
\lvec(125 13)
\lvec(124 13)
\ifill f:0
\move(126 12)
\lvec(132 12)
\lvec(132 13)
\lvec(126 13)
\ifill f:0
\move(133 12)
\lvec(134 12)
\lvec(134 13)
\lvec(133 13)
\ifill f:0
\move(135 12)
\lvec(155 12)
\lvec(155 13)
\lvec(135 13)
\ifill f:0
\move(156 12)
\lvec(157 12)
\lvec(157 13)
\lvec(156 13)
\ifill f:0
\move(158 12)
\lvec(180 12)
\lvec(180 13)
\lvec(158 13)
\ifill f:0
\move(181 12)
\lvec(182 12)
\lvec(182 13)
\lvec(181 13)
\ifill f:0
\move(183 12)
\lvec(191 12)
\lvec(191 13)
\lvec(183 13)
\ifill f:0
\move(192 12)
\lvec(193 12)
\lvec(193 13)
\lvec(192 13)
\ifill f:0
\move(194 12)
\lvec(202 12)
\lvec(202 13)
\lvec(194 13)
\ifill f:0
\move(203 12)
\lvec(209 12)
\lvec(209 13)
\lvec(203 13)
\ifill f:0
\move(210 12)
\lvec(216 12)
\lvec(216 13)
\lvec(210 13)
\ifill f:0
\move(217 12)
\lvec(223 12)
\lvec(223 13)
\lvec(217 13)
\ifill f:0
\move(224 12)
\lvec(228 12)
\lvec(228 13)
\lvec(224 13)
\ifill f:0
\move(229 12)
\lvec(245 12)
\lvec(245 13)
\lvec(229 13)
\ifill f:0
\move(246 12)
\lvec(263 12)
\lvec(263 13)
\lvec(246 13)
\ifill f:0
\move(264 12)
\lvec(276 12)
\lvec(276 13)
\lvec(264 13)
\ifill f:0
\move(277 12)
\lvec(287 12)
\lvec(287 13)
\lvec(277 13)
\ifill f:0
\move(288 12)
\lvec(295 12)
\lvec(295 13)
\lvec(288 13)
\ifill f:0
\move(296 12)
\lvec(309 12)
\lvec(309 13)
\lvec(296 13)
\ifill f:0
\move(310 12)
\lvec(326 12)
\lvec(326 13)
\lvec(310 13)
\ifill f:0
\move(327 12)
\lvec(329 12)
\lvec(329 13)
\lvec(327 13)
\ifill f:0
\move(330 12)
\lvec(332 12)
\lvec(332 13)
\lvec(330 13)
\ifill f:0
\move(333 12)
\lvec(371 12)
\lvec(371 13)
\lvec(333 13)
\ifill f:0
\move(372 12)
\lvec(374 12)
\lvec(374 13)
\lvec(372 13)
\ifill f:0
\move(375 12)
\lvec(377 12)
\lvec(377 13)
\lvec(375 13)
\ifill f:0
\move(378 12)
\lvec(396 12)
\lvec(396 13)
\lvec(378 13)
\ifill f:0
\move(397 12)
\lvec(412 12)
\lvec(412 13)
\lvec(397 13)
\ifill f:0
\move(413 12)
\lvec(422 12)
\lvec(422 13)
\lvec(413 13)
\ifill f:0
\move(423 12)
\lvec(435 12)
\lvec(435 13)
\lvec(423 13)
\ifill f:0
\move(436 12)
\lvec(451 12)
\lvec(451 13)
\lvec(436 13)
\ifill f:0
\move(14 13)
\lvec(16 13)
\lvec(16 14)
\lvec(14 14)
\ifill f:0
\move(18 13)
\lvec(19 13)
\lvec(19 14)
\lvec(18 14)
\ifill f:0
\move(20 13)
\lvec(21 13)
\lvec(21 14)
\lvec(20 14)
\ifill f:0
\move(22 13)
\lvec(23 13)
\lvec(23 14)
\lvec(22 14)
\ifill f:0
\move(24 13)
\lvec(26 13)
\lvec(26 14)
\lvec(24 14)
\ifill f:0
\move(27 13)
\lvec(29 13)
\lvec(29 14)
\lvec(27 14)
\ifill f:0
\move(30 13)
\lvec(33 13)
\lvec(33 14)
\lvec(30 14)
\ifill f:0
\move(36 13)
\lvec(45 13)
\lvec(45 14)
\lvec(36 14)
\ifill f:0
\move(47 13)
\lvec(58 13)
\lvec(58 14)
\lvec(47 14)
\ifill f:0
\move(61 13)
\lvec(66 13)
\lvec(66 14)
\lvec(61 14)
\ifill f:0
\move(67 13)
\lvec(71 13)
\lvec(71 14)
\lvec(67 14)
\ifill f:0
\move(72 13)
\lvec(76 13)
\lvec(76 14)
\lvec(72 14)
\ifill f:0
\move(77 13)
\lvec(80 13)
\lvec(80 14)
\lvec(77 14)
\ifill f:0
\move(81 13)
\lvec(84 13)
\lvec(84 14)
\lvec(81 14)
\ifill f:0
\move(85 13)
\lvec(88 13)
\lvec(88 14)
\lvec(85 14)
\ifill f:0
\move(89 13)
\lvec(91 13)
\lvec(91 14)
\lvec(89 14)
\ifill f:0
\move(92 13)
\lvec(98 13)
\lvec(98 14)
\lvec(92 14)
\ifill f:0
\move(99 13)
\lvec(101 13)
\lvec(101 14)
\lvec(99 14)
\ifill f:0
\move(102 13)
\lvec(104 13)
\lvec(104 14)
\lvec(102 14)
\ifill f:0
\move(105 13)
\lvec(107 13)
\lvec(107 14)
\lvec(105 14)
\ifill f:0
\move(108 13)
\lvec(115 13)
\lvec(115 14)
\lvec(108 14)
\ifill f:0
\move(116 13)
\lvec(118 13)
\lvec(118 14)
\lvec(116 14)
\ifill f:0
\move(119 13)
\lvec(123 13)
\lvec(123 14)
\lvec(119 14)
\ifill f:0
\move(124 13)
\lvec(128 13)
\lvec(128 14)
\lvec(124 14)
\ifill f:0
\move(129 13)
\lvec(133 13)
\lvec(133 14)
\lvec(129 14)
\ifill f:0
\move(134 13)
\lvec(140 13)
\lvec(140 14)
\lvec(134 14)
\ifill f:0
\move(141 13)
\lvec(147 13)
\lvec(147 14)
\lvec(141 14)
\ifill f:0
\move(148 13)
\lvec(156 13)
\lvec(156 14)
\lvec(148 14)
\ifill f:0
\move(157 13)
\lvec(158 13)
\lvec(158 14)
\lvec(157 14)
\ifill f:0
\move(159 13)
\lvec(160 13)
\lvec(160 14)
\lvec(159 14)
\ifill f:0
\move(161 13)
\lvec(181 13)
\lvec(181 14)
\lvec(161 14)
\ifill f:0
\move(182 13)
\lvec(183 13)
\lvec(183 14)
\lvec(182 14)
\ifill f:0
\move(184 13)
\lvec(206 13)
\lvec(206 14)
\lvec(184 14)
\ifill f:0
\move(207 13)
\lvec(208 13)
\lvec(208 14)
\lvec(207 14)
\ifill f:0
\move(209 13)
\lvec(210 13)
\lvec(210 14)
\lvec(209 14)
\ifill f:0
\move(211 13)
\lvec(221 13)
\lvec(221 14)
\lvec(211 14)
\ifill f:0
\move(222 13)
\lvec(230 13)
\lvec(230 14)
\lvec(222 14)
\ifill f:0
\move(231 13)
\lvec(239 13)
\lvec(239 14)
\lvec(231 14)
\ifill f:0
\move(240 13)
\lvec(246 13)
\lvec(246 14)
\lvec(240 14)
\ifill f:0
\move(247 13)
\lvec(253 13)
\lvec(253 14)
\lvec(247 14)
\ifill f:0
\move(254 13)
\lvec(260 13)
\lvec(260 14)
\lvec(254 14)
\ifill f:0
\move(261 13)
\lvec(265 13)
\lvec(265 14)
\lvec(261 14)
\ifill f:0
\move(266 13)
\lvec(277 13)
\lvec(277 14)
\lvec(266 14)
\ifill f:0
\move(278 13)
\lvec(282 13)
\lvec(282 14)
\lvec(278 14)
\ifill f:0
\move(283 13)
\lvec(287 13)
\lvec(287 14)
\lvec(283 14)
\ifill f:0
\move(288 13)
\lvec(292 13)
\lvec(292 14)
\lvec(288 14)
\ifill f:0
\move(293 13)
\lvec(305 13)
\lvec(305 14)
\lvec(293 14)
\ifill f:0
\move(306 13)
\lvec(310 13)
\lvec(310 14)
\lvec(306 14)
\ifill f:0
\move(311 13)
\lvec(318 13)
\lvec(318 14)
\lvec(311 14)
\ifill f:0
\move(319 13)
\lvec(326 13)
\lvec(326 14)
\lvec(319 14)
\ifill f:0
\move(327 13)
\lvec(334 13)
\lvec(334 14)
\lvec(327 14)
\ifill f:0
\move(335 13)
\lvec(345 13)
\lvec(345 14)
\lvec(335 14)
\ifill f:0
\move(346 13)
\lvec(356 13)
\lvec(356 14)
\lvec(346 14)
\ifill f:0
\move(357 13)
\lvec(370 13)
\lvec(370 14)
\lvec(357 14)
\ifill f:0
\move(371 13)
\lvec(373 13)
\lvec(373 14)
\lvec(371 14)
\ifill f:0
\move(374 13)
\lvec(376 13)
\lvec(376 14)
\lvec(374 14)
\ifill f:0
\move(377 13)
\lvec(408 13)
\lvec(408 14)
\lvec(377 14)
\ifill f:0
\move(409 13)
\lvec(411 13)
\lvec(411 14)
\lvec(409 14)
\ifill f:0
\move(412 13)
\lvec(445 13)
\lvec(445 14)
\lvec(412 14)
\ifill f:0
\move(446 13)
\lvec(448 13)
\lvec(448 14)
\lvec(446 14)
\ifill f:0
\move(449 13)
\lvec(451 13)
\lvec(451 14)
\lvec(449 14)
\ifill f:0
\move(11 14)
\lvec(16 14)
\lvec(16 15)
\lvec(11 15)
\ifill f:0
\move(18 14)
\lvec(20 14)
\lvec(20 15)
\lvec(18 15)
\ifill f:0
\move(21 14)
\lvec(22 14)
\lvec(22 15)
\lvec(21 15)
\ifill f:0
\move(24 14)
\lvec(26 14)
\lvec(26 15)
\lvec(24 15)
\ifill f:0
\move(27 14)
\lvec(28 14)
\lvec(28 15)
\lvec(27 15)
\ifill f:0
\move(29 14)
\lvec(31 14)
\lvec(31 15)
\lvec(29 15)
\ifill f:0
\move(32 14)
\lvec(34 14)
\lvec(34 15)
\lvec(32 15)
\ifill f:0
\move(36 14)
\lvec(39 14)
\lvec(39 15)
\lvec(36 15)
\ifill f:0
\move(42 14)
\lvec(52 14)
\lvec(52 15)
\lvec(42 15)
\ifill f:0
\move(54 14)
\lvec(66 14)
\lvec(66 15)
\lvec(54 15)
\ifill f:0
\move(69 14)
\lvec(74 14)
\lvec(74 15)
\lvec(69 15)
\ifill f:0
\move(76 14)
\lvec(80 14)
\lvec(80 15)
\lvec(76 15)
\ifill f:0
\move(81 14)
\lvec(85 14)
\lvec(85 15)
\lvec(81 15)
\ifill f:0
\move(86 14)
\lvec(89 14)
\lvec(89 15)
\lvec(86 15)
\ifill f:0
\move(90 14)
\lvec(94 14)
\lvec(94 15)
\lvec(90 15)
\ifill f:0
\move(95 14)
\lvec(101 14)
\lvec(101 15)
\lvec(95 15)
\ifill f:0
\move(102 14)
\lvec(108 14)
\lvec(108 15)
\lvec(102 15)
\ifill f:0
\move(109 14)
\lvec(111 14)
\lvec(111 15)
\lvec(109 15)
\ifill f:0
\move(112 14)
\lvec(126 14)
\lvec(126 15)
\lvec(112 15)
\ifill f:0
\move(127 14)
\lvec(129 14)
\lvec(129 15)
\lvec(127 15)
\ifill f:0
\move(130 14)
\lvec(137 14)
\lvec(137 15)
\lvec(130 15)
\ifill f:0
\move(138 14)
\lvec(145 14)
\lvec(145 15)
\lvec(138 15)
\ifill f:0
\move(146 14)
\lvec(152 14)
\lvec(152 15)
\lvec(146 15)
\ifill f:0
\move(153 14)
\lvec(157 14)
\lvec(157 15)
\lvec(153 15)
\ifill f:0
\move(158 14)
\lvec(164 14)
\lvec(164 15)
\lvec(158 15)
\ifill f:0
\move(165 14)
\lvec(171 14)
\lvec(171 15)
\lvec(165 15)
\ifill f:0
\move(172 14)
\lvec(173 14)
\lvec(173 15)
\lvec(172 15)
\ifill f:0
\move(174 14)
\lvec(182 14)
\lvec(182 15)
\lvec(174 15)
\ifill f:0
\move(183 14)
\lvec(184 14)
\lvec(184 15)
\lvec(183 15)
\ifill f:0
\move(185 14)
\lvec(186 14)
\lvec(186 15)
\lvec(185 15)
\ifill f:0
\move(187 14)
\lvec(209 14)
\lvec(209 15)
\lvec(187 15)
\ifill f:0
\move(210 14)
\lvec(211 14)
\lvec(211 15)
\lvec(210 15)
\ifill f:0
\move(212 14)
\lvec(236 14)
\lvec(236 15)
\lvec(212 15)
\ifill f:0
\move(237 14)
\lvec(238 14)
\lvec(238 15)
\lvec(237 15)
\ifill f:0
\move(239 14)
\lvec(240 14)
\lvec(240 15)
\lvec(239 15)
\ifill f:0
\move(241 14)
\lvec(251 14)
\lvec(251 15)
\lvec(241 15)
\ifill f:0
\move(252 14)
\lvec(253 14)
\lvec(253 15)
\lvec(252 15)
\ifill f:0
\move(254 14)
\lvec(262 14)
\lvec(262 15)
\lvec(254 15)
\ifill f:0
\move(263 14)
\lvec(271 14)
\lvec(271 15)
\lvec(263 15)
\ifill f:0
\move(272 14)
\lvec(278 14)
\lvec(278 15)
\lvec(272 15)
\ifill f:0
\move(279 14)
\lvec(287 14)
\lvec(287 15)
\lvec(279 15)
\ifill f:0
\move(288 14)
\lvec(299 14)
\lvec(299 15)
\lvec(288 15)
\ifill f:0
\move(300 14)
\lvec(311 14)
\lvec(311 15)
\lvec(300 15)
\ifill f:0
\move(312 14)
\lvec(316 14)
\lvec(316 15)
\lvec(312 15)
\ifill f:0
\move(317 14)
\lvec(341 14)
\lvec(341 15)
\lvec(317 15)
\ifill f:0
\move(342 14)
\lvec(346 14)
\lvec(346 15)
\lvec(342 15)
\ifill f:0
\move(347 14)
\lvec(359 14)
\lvec(359 15)
\lvec(347 15)
\ifill f:0
\move(360 14)
\lvec(372 14)
\lvec(372 15)
\lvec(360 15)
\ifill f:0
\move(373 14)
\lvec(383 14)
\lvec(383 15)
\lvec(373 15)
\ifill f:0
\move(384 14)
\lvec(391 14)
\lvec(391 15)
\lvec(384 15)
\ifill f:0
\move(392 14)
\lvec(402 14)
\lvec(402 15)
\lvec(392 15)
\ifill f:0
\move(403 14)
\lvec(413 14)
\lvec(413 15)
\lvec(403 15)
\ifill f:0
\move(414 14)
\lvec(416 14)
\lvec(416 15)
\lvec(414 15)
\ifill f:0
\move(417 14)
\lvec(430 14)
\lvec(430 15)
\lvec(417 15)
\ifill f:0
\move(431 14)
\lvec(433 14)
\lvec(433 15)
\lvec(431 15)
\ifill f:0
\move(434 14)
\lvec(436 14)
\lvec(436 15)
\lvec(434 15)
\ifill f:0
\move(437 14)
\lvec(451 14)
\lvec(451 15)
\lvec(437 15)
\ifill f:0
\move(11 15)
\lvec(12 15)
\lvec(12 16)
\lvec(11 16)
\ifill f:0
\move(13 15)
\lvec(14 15)
\lvec(14 16)
\lvec(13 16)
\ifill f:0
\move(15 15)
\lvec(16 15)
\lvec(16 16)
\lvec(15 16)
\ifill f:0
\move(20 15)
\lvec(22 15)
\lvec(22 16)
\lvec(20 16)
\ifill f:0
\move(23 15)
\lvec(24 15)
\lvec(24 16)
\lvec(23 16)
\ifill f:0
\move(25 15)
\lvec(27 15)
\lvec(27 16)
\lvec(25 16)
\ifill f:0
\move(28 15)
\lvec(30 15)
\lvec(30 16)
\lvec(28 16)
\ifill f:0
\move(31 15)
\lvec(32 15)
\lvec(32 16)
\lvec(31 16)
\ifill f:0
\move(33 15)
\lvec(35 15)
\lvec(35 16)
\lvec(33 16)
\ifill f:0
\move(36 15)
\lvec(38 15)
\lvec(38 16)
\lvec(36 16)
\ifill f:0
\move(40 15)
\lvec(43 15)
\lvec(43 16)
\lvec(40 16)
\ifill f:0
\move(45 15)
\lvec(50 15)
\lvec(50 16)
\lvec(45 16)
\ifill f:0
\move(54 15)
\lvec(68 15)
\lvec(68 16)
\lvec(54 16)
\ifill f:0
\move(72 15)
\lvec(79 15)
\lvec(79 16)
\lvec(72 16)
\ifill f:0
\move(81 15)
\lvec(86 15)
\lvec(86 16)
\lvec(81 16)
\ifill f:0
\move(88 15)
\lvec(92 15)
\lvec(92 16)
\lvec(88 16)
\ifill f:0
\move(93 15)
\lvec(97 15)
\lvec(97 16)
\lvec(93 16)
\ifill f:0
\move(98 15)
\lvec(101 15)
\lvec(101 16)
\lvec(98 16)
\ifill f:0
\move(102 15)
\lvec(106 15)
\lvec(106 16)
\lvec(102 16)
\ifill f:0
\move(107 15)
\lvec(113 15)
\lvec(113 16)
\lvec(107 16)
\ifill f:0
\move(114 15)
\lvec(117 15)
\lvec(117 16)
\lvec(114 16)
\ifill f:0
\move(118 15)
\lvec(124 15)
\lvec(124 16)
\lvec(118 16)
\ifill f:0
\move(125 15)
\lvec(127 15)
\lvec(127 16)
\lvec(125 16)
\ifill f:0
\move(128 15)
\lvec(130 15)
\lvec(130 16)
\lvec(128 16)
\ifill f:0
\move(131 15)
\lvec(133 15)
\lvec(133 16)
\lvec(131 16)
\ifill f:0
\move(134 15)
\lvec(136 15)
\lvec(136 16)
\lvec(134 16)
\ifill f:0
\move(137 15)
\lvec(139 15)
\lvec(139 16)
\lvec(137 16)
\ifill f:0
\move(140 15)
\lvec(142 15)
\lvec(142 16)
\lvec(140 16)
\ifill f:0
\move(143 15)
\lvec(145 15)
\lvec(145 16)
\lvec(143 16)
\ifill f:0
\move(146 15)
\lvec(153 15)
\lvec(153 16)
\lvec(146 16)
\ifill f:0
\move(154 15)
\lvec(161 15)
\lvec(161 16)
\lvec(154 16)
\ifill f:0
\move(162 15)
\lvec(166 15)
\lvec(166 16)
\lvec(162 16)
\ifill f:0
\move(167 15)
\lvec(171 15)
\lvec(171 16)
\lvec(167 16)
\ifill f:0
\move(172 15)
\lvec(176 15)
\lvec(176 16)
\lvec(172 16)
\ifill f:0
\move(177 15)
\lvec(183 15)
\lvec(183 16)
\lvec(177 16)
\ifill f:0
\move(184 15)
\lvec(190 15)
\lvec(190 16)
\lvec(184 16)
\ifill f:0
\move(191 15)
\lvec(192 15)
\lvec(192 16)
\lvec(191 16)
\ifill f:0
\move(193 15)
\lvec(199 15)
\lvec(199 16)
\lvec(193 16)
\ifill f:0
\move(200 15)
\lvec(201 15)
\lvec(201 16)
\lvec(200 16)
\ifill f:0
\move(202 15)
\lvec(210 15)
\lvec(210 16)
\lvec(202 16)
\ifill f:0
\move(211 15)
\lvec(212 15)
\lvec(212 16)
\lvec(211 16)
\ifill f:0
\move(213 15)
\lvec(214 15)
\lvec(214 16)
\lvec(213 16)
\ifill f:0
\move(215 15)
\lvec(239 15)
\lvec(239 16)
\lvec(215 16)
\ifill f:0
\move(240 15)
\lvec(241 15)
\lvec(241 16)
\lvec(240 16)
\ifill f:0
\move(242 15)
\lvec(268 15)
\lvec(268 16)
\lvec(242 16)
\ifill f:0
\move(269 15)
\lvec(270 15)
\lvec(270 16)
\lvec(269 16)
\ifill f:0
\move(271 15)
\lvec(272 15)
\lvec(272 16)
\lvec(271 16)
\ifill f:0
\move(273 15)
\lvec(283 15)
\lvec(283 16)
\lvec(273 16)
\ifill f:0
\move(284 15)
\lvec(285 15)
\lvec(285 16)
\lvec(284 16)
\ifill f:0
\move(286 15)
\lvec(296 15)
\lvec(296 16)
\lvec(286 16)
\ifill f:0
\move(297 15)
\lvec(305 15)
\lvec(305 16)
\lvec(297 16)
\ifill f:0
\move(306 15)
\lvec(314 15)
\lvec(314 16)
\lvec(306 16)
\ifill f:0
\move(315 15)
\lvec(321 15)
\lvec(321 16)
\lvec(315 16)
\ifill f:0
\move(322 15)
\lvec(328 15)
\lvec(328 16)
\lvec(322 16)
\ifill f:0
\move(329 15)
\lvec(335 15)
\lvec(335 16)
\lvec(329 16)
\ifill f:0
\move(336 15)
\lvec(347 15)
\lvec(347 16)
\lvec(336 16)
\ifill f:0
\move(348 15)
\lvec(359 15)
\lvec(359 16)
\lvec(348 16)
\ifill f:0
\move(360 15)
\lvec(364 15)
\lvec(364 16)
\lvec(360 16)
\ifill f:0
\move(365 15)
\lvec(369 15)
\lvec(369 16)
\lvec(365 16)
\ifill f:0
\move(370 15)
\lvec(374 15)
\lvec(374 16)
\lvec(370 16)
\ifill f:0
\move(375 15)
\lvec(379 15)
\lvec(379 16)
\lvec(375 16)
\ifill f:0
\move(380 15)
\lvec(384 15)
\lvec(384 16)
\lvec(380 16)
\ifill f:0
\move(385 15)
\lvec(389 15)
\lvec(389 16)
\lvec(385 16)
\ifill f:0
\move(390 15)
\lvec(394 15)
\lvec(394 16)
\lvec(390 16)
\ifill f:0
\move(395 15)
\lvec(407 15)
\lvec(407 16)
\lvec(395 16)
\ifill f:0
\move(408 15)
\lvec(415 15)
\lvec(415 16)
\lvec(408 16)
\ifill f:0
\move(416 15)
\lvec(428 15)
\lvec(428 16)
\lvec(416 16)
\ifill f:0
\move(429 15)
\lvec(439 15)
\lvec(439 16)
\lvec(429 16)
\ifill f:0
\move(440 15)
\lvec(447 15)
\lvec(447 16)
\lvec(440 16)
\ifill f:0
\move(448 15)
\lvec(451 15)
\lvec(451 16)
\lvec(448 16)
\ifill f:0
\move(12 16)
\lvec(14 16)
\lvec(14 17)
\lvec(12 17)
\ifill f:0
\move(22 16)
\lvec(24 16)
\lvec(24 17)
\lvec(22 17)
\ifill f:0
\move(25 16)
\lvec(27 16)
\lvec(27 17)
\lvec(25 17)
\ifill f:0
\move(28 16)
\lvec(29 16)
\lvec(29 17)
\lvec(28 17)
\ifill f:0
\move(30 16)
\lvec(31 16)
\lvec(31 17)
\lvec(30 17)
\ifill f:0
\move(32 16)
\lvec(33 16)
\lvec(33 17)
\lvec(32 17)
\ifill f:0
\move(34 16)
\lvec(35 16)
\lvec(35 17)
\lvec(34 17)
\ifill f:0
\move(36 16)
\lvec(38 16)
\lvec(38 17)
\lvec(36 17)
\ifill f:0
\move(39 16)
\lvec(41 16)
\lvec(41 17)
\lvec(39 17)
\ifill f:0
\move(42 16)
\lvec(45 16)
\lvec(45 17)
\lvec(42 17)
\ifill f:0
\move(46 16)
\lvec(50 16)
\lvec(50 17)
\lvec(46 17)
\ifill f:0
\move(52 16)
\lvec(57 16)
\lvec(57 17)
\lvec(52 17)
\ifill f:0
\move(62 16)
\lvec(76 16)
\lvec(76 17)
\lvec(62 17)
\ifill f:0
\move(81 16)
\lvec(88 16)
\lvec(88 17)
\lvec(81 17)
\ifill f:0
\move(90 16)
\lvec(95 16)
\lvec(95 17)
\lvec(90 17)
\ifill f:0
\move(97 16)
\lvec(102 16)
\lvec(102 17)
\lvec(97 17)
\ifill f:0
\move(103 16)
\lvec(107 16)
\lvec(107 17)
\lvec(103 17)
\ifill f:0
\move(108 16)
\lvec(112 16)
\lvec(112 17)
\lvec(108 17)
\ifill f:0
\move(113 16)
\lvec(116 16)
\lvec(116 17)
\lvec(113 17)
\ifill f:0
\move(117 16)
\lvec(120 16)
\lvec(120 17)
\lvec(117 17)
\ifill f:0
\move(121 16)
\lvec(124 16)
\lvec(124 17)
\lvec(121 17)
\ifill f:0
\move(125 16)
\lvec(128 16)
\lvec(128 17)
\lvec(125 17)
\ifill f:0
\move(129 16)
\lvec(135 16)
\lvec(135 17)
\lvec(129 17)
\ifill f:0
\move(136 16)
\lvec(142 16)
\lvec(142 17)
\lvec(136 17)
\ifill f:0
\move(143 16)
\lvec(145 16)
\lvec(145 17)
\lvec(143 17)
\ifill f:0
\move(146 16)
\lvec(148 16)
\lvec(148 17)
\lvec(146 17)
\ifill f:0
\move(149 16)
\lvec(151 16)
\lvec(151 17)
\lvec(149 17)
\ifill f:0
\move(152 16)
\lvec(154 16)
\lvec(154 17)
\lvec(152 17)
\ifill f:0
\move(155 16)
\lvec(157 16)
\lvec(157 17)
\lvec(155 17)
\ifill f:0
\move(158 16)
\lvec(160 16)
\lvec(160 17)
\lvec(158 17)
\ifill f:0
\move(161 16)
\lvec(163 16)
\lvec(163 17)
\lvec(161 17)
\ifill f:0
\move(164 16)
\lvec(171 16)
\lvec(171 17)
\lvec(164 17)
\ifill f:0
\move(172 16)
\lvec(174 16)
\lvec(174 17)
\lvec(172 17)
\ifill f:0
\move(175 16)
\lvec(179 16)
\lvec(179 17)
\lvec(175 17)
\ifill f:0
\move(180 16)
\lvec(187 16)
\lvec(187 17)
\lvec(180 17)
\ifill f:0
\move(188 16)
\lvec(192 16)
\lvec(192 17)
\lvec(188 17)
\ifill f:0
\move(193 16)
\lvec(199 16)
\lvec(199 17)
\lvec(193 17)
\ifill f:0
\move(200 16)
\lvec(204 16)
\lvec(204 17)
\lvec(200 17)
\ifill f:0
\move(205 16)
\lvec(211 16)
\lvec(211 17)
\lvec(205 17)
\ifill f:0
\move(212 16)
\lvec(218 16)
\lvec(218 17)
\lvec(212 17)
\ifill f:0
\move(219 16)
\lvec(220 16)
\lvec(220 17)
\lvec(219 17)
\ifill f:0
\move(221 16)
\lvec(229 16)
\lvec(229 17)
\lvec(221 17)
\ifill f:0
\move(230 16)
\lvec(231 16)
\lvec(231 17)
\lvec(230 17)
\ifill f:0
\move(232 16)
\lvec(240 16)
\lvec(240 17)
\lvec(232 17)
\ifill f:0
\move(241 16)
\lvec(242 16)
\lvec(242 17)
\lvec(241 17)
\ifill f:0
\move(243 16)
\lvec(244 16)
\lvec(244 17)
\lvec(243 17)
\ifill f:0
\move(245 16)
\lvec(271 16)
\lvec(271 17)
\lvec(245 17)
\ifill f:0
\move(272 16)
\lvec(273 16)
\lvec(273 17)
\lvec(272 17)
\ifill f:0
\move(274 16)
\lvec(302 16)
\lvec(302 17)
\lvec(274 17)
\ifill f:0
\move(303 16)
\lvec(304 16)
\lvec(304 17)
\lvec(303 17)
\ifill f:0
\move(305 16)
\lvec(306 16)
\lvec(306 17)
\lvec(305 17)
\ifill f:0
\move(307 16)
\lvec(317 16)
\lvec(317 17)
\lvec(307 17)
\ifill f:0
\move(318 16)
\lvec(319 16)
\lvec(319 17)
\lvec(318 17)
\ifill f:0
\move(320 16)
\lvec(330 16)
\lvec(330 17)
\lvec(320 17)
\ifill f:0
\move(331 16)
\lvec(332 16)
\lvec(332 17)
\lvec(331 17)
\ifill f:0
\move(333 16)
\lvec(341 16)
\lvec(341 17)
\lvec(333 17)
\ifill f:0
\move(342 16)
\lvec(350 16)
\lvec(350 17)
\lvec(342 17)
\ifill f:0
\move(351 16)
\lvec(357 16)
\lvec(357 17)
\lvec(351 17)
\ifill f:0
\move(358 16)
\lvec(366 16)
\lvec(366 17)
\lvec(358 17)
\ifill f:0
\move(367 16)
\lvec(373 16)
\lvec(373 17)
\lvec(367 17)
\ifill f:0
\move(374 16)
\lvec(385 16)
\lvec(385 17)
\lvec(374 17)
\ifill f:0
\move(386 16)
\lvec(392 16)
\lvec(392 17)
\lvec(386 17)
\ifill f:0
\move(393 16)
\lvec(397 16)
\lvec(397 17)
\lvec(393 17)
\ifill f:0
\move(398 16)
\lvec(409 16)
\lvec(409 17)
\lvec(398 17)
\ifill f:0
\move(410 16)
\lvec(414 16)
\lvec(414 17)
\lvec(410 17)
\ifill f:0
\move(415 16)
\lvec(419 16)
\lvec(419 17)
\lvec(415 17)
\ifill f:0
\move(420 16)
\lvec(424 16)
\lvec(424 17)
\lvec(420 17)
\ifill f:0
\move(425 16)
\lvec(429 16)
\lvec(429 17)
\lvec(425 17)
\ifill f:0
\move(430 16)
\lvec(434 16)
\lvec(434 17)
\lvec(430 17)
\ifill f:0
\move(435 16)
\lvec(439 16)
\lvec(439 17)
\lvec(435 17)
\ifill f:0
\move(440 16)
\lvec(444 16)
\lvec(444 17)
\lvec(440 17)
\ifill f:0
\move(445 16)
\lvec(451 16)
\lvec(451 17)
\lvec(445 17)
\ifill f:0
\move(11 17)
\lvec(12 17)
\lvec(12 18)
\lvec(11 18)
\ifill f:0
\move(13 17)
\lvec(15 17)
\lvec(15 18)
\lvec(13 18)
\ifill f:0
\move(18 17)
\lvec(23 17)
\lvec(23 18)
\lvec(18 18)
\ifill f:0
\move(25 17)
\lvec(27 17)
\lvec(27 18)
\lvec(25 18)
\ifill f:0
\move(29 17)
\lvec(30 17)
\lvec(30 18)
\lvec(29 18)
\ifill f:0
\move(31 17)
\lvec(38 17)
\lvec(38 18)
\lvec(31 18)
\ifill f:0
\move(39 17)
\lvec(40 17)
\lvec(40 18)
\lvec(39 18)
\ifill f:0
\move(41 17)
\lvec(43 17)
\lvec(43 18)
\lvec(41 18)
\ifill f:0
\move(44 17)
\lvec(46 17)
\lvec(46 18)
\lvec(44 18)
\ifill f:0
\move(47 17)
\lvec(50 17)
\lvec(50 18)
\lvec(47 18)
\ifill f:0
\move(51 17)
\lvec(54 17)
\lvec(54 18)
\lvec(51 18)
\ifill f:0
\move(56 17)
\lvec(60 17)
\lvec(60 18)
\lvec(56 18)
\ifill f:0
\move(63 17)
\lvec(76 17)
\lvec(76 18)
\lvec(63 18)
\ifill f:0
\move(78 17)
\lvec(93 17)
\lvec(93 18)
\lvec(78 18)
\ifill f:0
\move(96 17)
\lvec(102 17)
\lvec(102 18)
\lvec(96 18)
\ifill f:0
\move(104 17)
\lvec(109 17)
\lvec(109 18)
\lvec(104 18)
\ifill f:0
\move(110 17)
\lvec(115 17)
\lvec(115 18)
\lvec(110 18)
\ifill f:0
\move(116 17)
\lvec(120 17)
\lvec(120 18)
\lvec(116 18)
\ifill f:0
\move(121 17)
\lvec(125 17)
\lvec(125 18)
\lvec(121 18)
\ifill f:0
\move(126 17)
\lvec(129 17)
\lvec(129 18)
\lvec(126 18)
\ifill f:0
\move(130 17)
\lvec(145 17)
\lvec(145 18)
\lvec(130 18)
\ifill f:0
\move(146 17)
\lvec(149 17)
\lvec(149 18)
\lvec(146 18)
\ifill f:0
\move(150 17)
\lvec(152 17)
\lvec(152 18)
\lvec(150 18)
\ifill f:0
\move(153 17)
\lvec(159 17)
\lvec(159 18)
\lvec(153 18)
\ifill f:0
\move(160 17)
\lvec(162 17)
\lvec(162 18)
\lvec(160 18)
\ifill f:0
\move(163 17)
\lvec(183 17)
\lvec(183 18)
\lvec(163 18)
\ifill f:0
\move(184 17)
\lvec(186 17)
\lvec(186 18)
\lvec(184 18)
\ifill f:0
\move(187 17)
\lvec(194 17)
\lvec(194 18)
\lvec(187 18)
\ifill f:0
\move(195 17)
\lvec(197 17)
\lvec(197 18)
\lvec(195 18)
\ifill f:0
\move(198 17)
\lvec(202 17)
\lvec(202 18)
\lvec(198 18)
\ifill f:0
\move(203 17)
\lvec(222 17)
\lvec(222 18)
\lvec(203 18)
\ifill f:0
\move(223 17)
\lvec(227 17)
\lvec(227 18)
\lvec(223 18)
\ifill f:0
\move(228 17)
\lvec(234 17)
\lvec(234 18)
\lvec(228 18)
\ifill f:0
\move(235 17)
\lvec(241 17)
\lvec(241 18)
\lvec(235 18)
\ifill f:0
\move(242 17)
\lvec(250 17)
\lvec(250 18)
\lvec(242 18)
\ifill f:0
\move(251 17)
\lvec(259 17)
\lvec(259 18)
\lvec(251 18)
\ifill f:0
\move(260 17)
\lvec(261 17)
\lvec(261 18)
\lvec(260 18)
\ifill f:0
\move(262 17)
\lvec(272 17)
\lvec(272 18)
\lvec(262 18)
\ifill f:0
\move(273 17)
\lvec(274 17)
\lvec(274 18)
\lvec(273 18)
\ifill f:0
\move(275 17)
\lvec(276 17)
\lvec(276 18)
\lvec(275 18)
\ifill f:0
\move(277 17)
\lvec(305 17)
\lvec(305 18)
\lvec(277 18)
\ifill f:0
\move(306 17)
\lvec(307 17)
\lvec(307 18)
\lvec(306 18)
\ifill f:0
\move(308 17)
\lvec(338 17)
\lvec(338 18)
\lvec(308 18)
\ifill f:0
\move(339 17)
\lvec(340 17)
\lvec(340 18)
\lvec(339 18)
\ifill f:0
\move(341 17)
\lvec(342 17)
\lvec(342 18)
\lvec(341 18)
\ifill f:0
\move(343 17)
\lvec(355 17)
\lvec(355 18)
\lvec(343 18)
\ifill f:0
\move(356 17)
\lvec(357 17)
\lvec(357 18)
\lvec(356 18)
\ifill f:0
\move(358 17)
\lvec(368 17)
\lvec(368 18)
\lvec(358 18)
\ifill f:0
\move(369 17)
\lvec(379 17)
\lvec(379 18)
\lvec(369 18)
\ifill f:0
\move(380 17)
\lvec(388 17)
\lvec(388 18)
\lvec(380 18)
\ifill f:0
\move(389 17)
\lvec(397 17)
\lvec(397 18)
\lvec(389 18)
\ifill f:0
\move(398 17)
\lvec(404 17)
\lvec(404 18)
\lvec(398 18)
\ifill f:0
\move(405 17)
\lvec(432 17)
\lvec(432 18)
\lvec(405 18)
\ifill f:0
\move(433 17)
\lvec(439 17)
\lvec(439 18)
\lvec(433 18)
\ifill f:0
\move(440 17)
\lvec(444 17)
\lvec(444 18)
\lvec(440 18)
\ifill f:0
\move(445 17)
\lvec(451 17)
\lvec(451 18)
\lvec(445 18)
\ifill f:0
\move(10 18)
\lvec(13 18)
\lvec(13 19)
\lvec(10 19)
\ifill f:0
\move(14 18)
\lvec(15 18)
\lvec(15 19)
\lvec(14 19)
\ifill f:0
\move(17 18)
\lvec(19 18)
\lvec(19 19)
\lvec(17 19)
\ifill f:0
\move(25 18)
\lvec(28 18)
\lvec(28 19)
\lvec(25 19)
\ifill f:0
\move(30 18)
\lvec(31 18)
\lvec(31 19)
\lvec(30 19)
\ifill f:0
\move(33 18)
\lvec(34 18)
\lvec(34 19)
\lvec(33 19)
\ifill f:0
\move(35 18)
\lvec(36 18)
\lvec(36 19)
\lvec(35 19)
\ifill f:0
\move(37 18)
\lvec(42 18)
\lvec(42 19)
\lvec(37 19)
\ifill f:0
\move(43 18)
\lvec(44 18)
\lvec(44 19)
\lvec(43 19)
\ifill f:0
\move(45 18)
\lvec(47 18)
\lvec(47 19)
\lvec(45 19)
\ifill f:0
\move(48 18)
\lvec(50 18)
\lvec(50 19)
\lvec(48 19)
\ifill f:0
\move(51 18)
\lvec(53 18)
\lvec(53 19)
\lvec(51 19)
\ifill f:0
\move(55 18)
\lvec(57 18)
\lvec(57 19)
\lvec(55 19)
\ifill f:0
\move(59 18)
\lvec(62 18)
\lvec(62 19)
\lvec(59 19)
\ifill f:0
\move(64 18)
\lvec(68 18)
\lvec(68 19)
\lvec(64 19)
\ifill f:0
\move(72 18)
\lvec(85 18)
\lvec(85 19)
\lvec(72 19)
\ifill f:0
\move(87 18)
\lvec(102 18)
\lvec(102 19)
\lvec(87 19)
\ifill f:0
\move(106 18)
\lvec(112 18)
\lvec(112 19)
\lvec(106 19)
\ifill f:0
\move(114 18)
\lvec(119 18)
\lvec(119 19)
\lvec(114 19)
\ifill f:0
\move(121 18)
\lvec(125 18)
\lvec(125 19)
\lvec(121 19)
\ifill f:0
\move(127 18)
\lvec(131 18)
\lvec(131 19)
\lvec(127 19)
\ifill f:0
\move(132 18)
\lvec(136 18)
\lvec(136 19)
\lvec(132 19)
\ifill f:0
\move(137 18)
\lvec(141 18)
\lvec(141 19)
\lvec(137 19)
\ifill f:0
\move(142 18)
\lvec(145 18)
\lvec(145 19)
\lvec(142 19)
\ifill f:0
\move(146 18)
\lvec(149 18)
\lvec(149 19)
\lvec(146 19)
\ifill f:0
\move(150 18)
\lvec(157 18)
\lvec(157 19)
\lvec(150 19)
\ifill f:0
\move(158 18)
\lvec(161 18)
\lvec(161 19)
\lvec(158 19)
\ifill f:0
\move(162 18)
\lvec(165 18)
\lvec(165 19)
\lvec(162 19)
\ifill f:0
\move(166 18)
\lvec(168 18)
\lvec(168 19)
\lvec(166 19)
\ifill f:0
\move(169 18)
\lvec(172 18)
\lvec(172 19)
\lvec(169 19)
\ifill f:0
\move(173 18)
\lvec(175 18)
\lvec(175 19)
\lvec(173 19)
\ifill f:0
\move(176 18)
\lvec(185 18)
\lvec(185 19)
\lvec(176 19)
\ifill f:0
\move(186 18)
\lvec(188 18)
\lvec(188 19)
\lvec(186 19)
\ifill f:0
\move(189 18)
\lvec(191 18)
\lvec(191 19)
\lvec(189 19)
\ifill f:0
\move(192 18)
\lvec(194 18)
\lvec(194 19)
\lvec(192 19)
\ifill f:0
\move(195 18)
\lvec(197 18)
\lvec(197 19)
\lvec(195 19)
\ifill f:0
\move(198 18)
\lvec(200 18)
\lvec(200 19)
\lvec(198 19)
\ifill f:0
\move(201 18)
\lvec(211 18)
\lvec(211 19)
\lvec(201 19)
\ifill f:0
\move(212 18)
\lvec(214 18)
\lvec(214 19)
\lvec(212 19)
\ifill f:0
\move(215 18)
\lvec(219 18)
\lvec(219 19)
\lvec(215 19)
\ifill f:0
\move(220 18)
\lvec(222 18)
\lvec(222 19)
\lvec(220 19)
\ifill f:0
\move(223 18)
\lvec(227 18)
\lvec(227 19)
\lvec(223 19)
\ifill f:0
\move(228 18)
\lvec(232 18)
\lvec(232 19)
\lvec(228 19)
\ifill f:0
\move(233 18)
\lvec(247 18)
\lvec(247 19)
\lvec(233 19)
\ifill f:0
\move(248 18)
\lvec(252 18)
\lvec(252 19)
\lvec(248 19)
\ifill f:0
\move(253 18)
\lvec(259 18)
\lvec(259 19)
\lvec(253 19)
\ifill f:0
\move(260 18)
\lvec(266 18)
\lvec(266 19)
\lvec(260 19)
\ifill f:0
\move(267 18)
\lvec(273 18)
\lvec(273 19)
\lvec(267 19)
\ifill f:0
\move(274 18)
\lvec(282 18)
\lvec(282 19)
\lvec(274 19)
\ifill f:0
\move(283 18)
\lvec(284 18)
\lvec(284 19)
\lvec(283 19)
\ifill f:0
\move(285 18)
\lvec(293 18)
\lvec(293 19)
\lvec(285 19)
\ifill f:0
\move(294 18)
\lvec(295 18)
\lvec(295 19)
\lvec(294 19)
\ifill f:0
\move(296 18)
\lvec(306 18)
\lvec(306 19)
\lvec(296 19)
\ifill f:0
\move(307 18)
\lvec(308 18)
\lvec(308 19)
\lvec(307 19)
\ifill f:0
\move(309 18)
\lvec(310 18)
\lvec(310 19)
\lvec(309 19)
\ifill f:0
\move(311 18)
\lvec(312 18)
\lvec(312 19)
\lvec(311 19)
\ifill f:0
\move(313 18)
\lvec(341 18)
\lvec(341 19)
\lvec(313 19)
\ifill f:0
\move(342 18)
\lvec(343 18)
\lvec(343 19)
\lvec(342 19)
\ifill f:0
\move(344 18)
\lvec(374 18)
\lvec(374 19)
\lvec(344 19)
\ifill f:0
\move(375 18)
\lvec(376 18)
\lvec(376 19)
\lvec(375 19)
\ifill f:0
\move(377 18)
\lvec(378 18)
\lvec(378 19)
\lvec(377 19)
\ifill f:0
\move(379 18)
\lvec(380 18)
\lvec(380 19)
\lvec(379 19)
\ifill f:0
\move(381 18)
\lvec(393 18)
\lvec(393 19)
\lvec(381 19)
\ifill f:0
\move(394 18)
\lvec(395 18)
\lvec(395 19)
\lvec(394 19)
\ifill f:0
\move(396 18)
\lvec(406 18)
\lvec(406 19)
\lvec(396 19)
\ifill f:0
\move(407 18)
\lvec(408 18)
\lvec(408 19)
\lvec(407 19)
\ifill f:0
\move(409 18)
\lvec(419 18)
\lvec(419 19)
\lvec(409 19)
\ifill f:0
\move(420 18)
\lvec(428 18)
\lvec(428 19)
\lvec(420 19)
\ifill f:0
\move(429 18)
\lvec(437 18)
\lvec(437 19)
\lvec(429 19)
\ifill f:0
\move(438 18)
\lvec(446 18)
\lvec(446 19)
\lvec(438 19)
\ifill f:0
\move(447 18)
\lvec(451 18)
\lvec(451 19)
\lvec(447 19)
\ifill f:0
\move(13 19)
\lvec(14 19)
\lvec(14 20)
\lvec(13 20)
\ifill f:0
\move(17 19)
\lvec(18 19)
\lvec(18 20)
\lvec(17 20)
\ifill f:0
\move(20 19)
\lvec(24 19)
\lvec(24 20)
\lvec(20 20)
\ifill f:0
\move(25 19)
\lvec(29 19)
\lvec(29 20)
\lvec(25 20)
\ifill f:0
\move(32 19)
\lvec(33 19)
\lvec(33 20)
\lvec(32 20)
\ifill f:0
\move(35 19)
\lvec(36 19)
\lvec(36 20)
\lvec(35 20)
\ifill f:0
\move(38 19)
\lvec(39 19)
\lvec(39 20)
\lvec(38 20)
\ifill f:0
\move(40 19)
\lvec(41 19)
\lvec(41 20)
\lvec(40 20)
\ifill f:0
\move(42 19)
\lvec(43 19)
\lvec(43 20)
\lvec(42 20)
\ifill f:0
\move(44 19)
\lvec(45 19)
\lvec(45 20)
\lvec(44 20)
\ifill f:0
\move(46 19)
\lvec(47 19)
\lvec(47 20)
\lvec(46 20)
\ifill f:0
\move(48 19)
\lvec(50 19)
\lvec(50 20)
\lvec(48 20)
\ifill f:0
\move(51 19)
\lvec(52 19)
\lvec(52 20)
\lvec(51 20)
\ifill f:0
\move(53 19)
\lvec(55 19)
\lvec(55 20)
\lvec(53 20)
\ifill f:0
\move(57 19)
\lvec(59 19)
\lvec(59 20)
\lvec(57 20)
\ifill f:0
\move(60 19)
\lvec(63 19)
\lvec(63 20)
\lvec(60 20)
\ifill f:0
\move(64 19)
\lvec(67 19)
\lvec(67 20)
\lvec(64 20)
\ifill f:0
\move(69 19)
\lvec(73 19)
\lvec(73 20)
\lvec(69 20)
\ifill f:0
\move(76 19)
\lvec(82 19)
\lvec(82 20)
\lvec(76 20)
\ifill f:0
\move(88 19)
\lvec(104 19)
\lvec(104 20)
\lvec(88 20)
\ifill f:0
\move(110 19)
\lvec(118 19)
\lvec(118 20)
\lvec(110 20)
\ifill f:0
\move(121 19)
\lvec(127 19)
\lvec(127 20)
\lvec(121 20)
\ifill f:0
\move(129 19)
\lvec(134 19)
\lvec(134 20)
\lvec(129 20)
\ifill f:0
\move(135 19)
\lvec(140 19)
\lvec(140 20)
\lvec(135 20)
\ifill f:0
\move(141 19)
\lvec(145 19)
\lvec(145 20)
\lvec(141 20)
\ifill f:0
\move(147 19)
\lvec(151 19)
\lvec(151 20)
\lvec(147 20)
\ifill f:0
\move(152 19)
\lvec(155 19)
\lvec(155 20)
\lvec(152 20)
\ifill f:0
\move(156 19)
\lvec(160 19)
\lvec(160 20)
\lvec(156 20)
\ifill f:0
\move(161 19)
\lvec(164 19)
\lvec(164 20)
\lvec(161 20)
\ifill f:0
\move(165 19)
\lvec(168 19)
\lvec(168 20)
\lvec(165 20)
\ifill f:0
\move(169 19)
\lvec(172 19)
\lvec(172 20)
\lvec(169 20)
\ifill f:0
\move(173 19)
\lvec(176 19)
\lvec(176 20)
\lvec(173 20)
\ifill f:0
\move(177 19)
\lvec(180 19)
\lvec(180 20)
\lvec(177 20)
\ifill f:0
\move(181 19)
\lvec(183 19)
\lvec(183 20)
\lvec(181 20)
\ifill f:0
\move(184 19)
\lvec(187 19)
\lvec(187 20)
\lvec(184 20)
\ifill f:0
\move(188 19)
\lvec(190 19)
\lvec(190 20)
\lvec(188 20)
\ifill f:0
\move(191 19)
\lvec(194 19)
\lvec(194 20)
\lvec(191 20)
\ifill f:0
\move(195 19)
\lvec(197 19)
\lvec(197 20)
\lvec(195 20)
\ifill f:0
\move(198 19)
\lvec(200 19)
\lvec(200 20)
\lvec(198 20)
\ifill f:0
\move(201 19)
\lvec(213 19)
\lvec(213 20)
\lvec(201 20)
\ifill f:0
\move(214 19)
\lvec(227 19)
\lvec(227 20)
\lvec(214 20)
\ifill f:0
\move(228 19)
\lvec(230 19)
\lvec(230 20)
\lvec(228 20)
\ifill f:0
\move(231 19)
\lvec(233 19)
\lvec(233 20)
\lvec(231 20)
\ifill f:0
\move(234 19)
\lvec(241 19)
\lvec(241 20)
\lvec(234 20)
\ifill f:0
\move(242 19)
\lvec(244 19)
\lvec(244 20)
\lvec(242 20)
\ifill f:0
\move(245 19)
\lvec(249 19)
\lvec(249 20)
\lvec(245 20)
\ifill f:0
\move(250 19)
\lvec(254 19)
\lvec(254 20)
\lvec(250 20)
\ifill f:0
\move(255 19)
\lvec(257 19)
\lvec(257 20)
\lvec(255 20)
\ifill f:0
\move(258 19)
\lvec(262 19)
\lvec(262 20)
\lvec(258 20)
\ifill f:0
\move(263 19)
\lvec(267 19)
\lvec(267 20)
\lvec(263 20)
\ifill f:0
\move(268 19)
\lvec(274 19)
\lvec(274 20)
\lvec(268 20)
\ifill f:0
\move(275 19)
\lvec(279 19)
\lvec(279 20)
\lvec(275 20)
\ifill f:0
\move(280 19)
\lvec(284 19)
\lvec(284 20)
\lvec(280 20)
\ifill f:0
\move(285 19)
\lvec(286 19)
\lvec(286 20)
\lvec(285 20)
\ifill f:0
\move(287 19)
\lvec(291 19)
\lvec(291 20)
\lvec(287 20)
\ifill f:0
\move(292 19)
\lvec(298 19)
\lvec(298 20)
\lvec(292 20)
\ifill f:0
\move(299 19)
\lvec(300 19)
\lvec(300 20)
\lvec(299 20)
\ifill f:0
\move(301 19)
\lvec(307 19)
\lvec(307 20)
\lvec(301 20)
\ifill f:0
\move(308 19)
\lvec(309 19)
\lvec(309 20)
\lvec(308 20)
\ifill f:0
\move(310 19)
\lvec(316 19)
\lvec(316 20)
\lvec(310 20)
\ifill f:0
\move(317 19)
\lvec(318 19)
\lvec(318 20)
\lvec(317 20)
\ifill f:0
\move(319 19)
\lvec(327 19)
\lvec(327 20)
\lvec(319 20)
\ifill f:0
\move(328 19)
\lvec(329 19)
\lvec(329 20)
\lvec(328 20)
\ifill f:0
\move(330 19)
\lvec(331 19)
\lvec(331 20)
\lvec(330 20)
\ifill f:0
\move(332 19)
\lvec(342 19)
\lvec(342 20)
\lvec(332 20)
\ifill f:0
\move(343 19)
\lvec(344 19)
\lvec(344 20)
\lvec(343 20)
\ifill f:0
\move(345 19)
\lvec(346 19)
\lvec(346 20)
\lvec(345 20)
\ifill f:0
\move(347 19)
\lvec(348 19)
\lvec(348 20)
\lvec(347 20)
\ifill f:0
\move(349 19)
\lvec(379 19)
\lvec(379 20)
\lvec(349 20)
\ifill f:0
\move(380 19)
\lvec(381 19)
\lvec(381 20)
\lvec(380 20)
\ifill f:0
\move(382 19)
\lvec(414 19)
\lvec(414 20)
\lvec(382 20)
\ifill f:0
\move(415 19)
\lvec(416 19)
\lvec(416 20)
\lvec(415 20)
\ifill f:0
\move(417 19)
\lvec(418 19)
\lvec(418 20)
\lvec(417 20)
\ifill f:0
\move(419 19)
\lvec(420 19)
\lvec(420 20)
\lvec(419 20)
\ifill f:0
\move(421 19)
\lvec(433 19)
\lvec(433 20)
\lvec(421 20)
\ifill f:0
\move(434 19)
\lvec(435 19)
\lvec(435 20)
\lvec(434 20)
\ifill f:0
\move(436 19)
\lvec(437 19)
\lvec(437 20)
\lvec(436 20)
\ifill f:0
\move(438 19)
\lvec(448 19)
\lvec(448 20)
\lvec(438 20)
\ifill f:0
\move(449 19)
\lvec(450 19)
\lvec(450 20)
\lvec(449 20)
\ifill f:0
\move(12 20)
\lvec(14 20)
\lvec(14 21)
\lvec(12 21)
\ifill f:0
\move(15 20)
\lvec(16 20)
\lvec(16 21)
\lvec(15 21)
\ifill f:0
\move(19 20)
\lvec(20 20)
\lvec(20 21)
\lvec(19 21)
\ifill f:0
\move(23 20)
\lvec(26 20)
\lvec(26 21)
\lvec(23 21)
\ifill f:0
\move(27 20)
\lvec(31 20)
\lvec(31 21)
\lvec(27 21)
\ifill f:0
\move(35 20)
\lvec(36 20)
\lvec(36 21)
\lvec(35 21)
\ifill f:0
\move(38 20)
\lvec(40 20)
\lvec(40 21)
\lvec(38 21)
\ifill f:0
\move(41 20)
\lvec(42 20)
\lvec(42 21)
\lvec(41 21)
\ifill f:0
\move(43 20)
\lvec(44 20)
\lvec(44 21)
\lvec(43 21)
\ifill f:0
\move(45 20)
\lvec(47 20)
\lvec(47 21)
\lvec(45 21)
\ifill f:0
\move(48 20)
\lvec(52 20)
\lvec(52 21)
\lvec(48 21)
\ifill f:0
\move(53 20)
\lvec(54 20)
\lvec(54 21)
\lvec(53 21)
\ifill f:0
\move(55 20)
\lvec(57 20)
\lvec(57 21)
\lvec(55 21)
\ifill f:0
\move(58 20)
\lvec(60 20)
\lvec(60 21)
\lvec(58 21)
\ifill f:0
\move(61 20)
\lvec(63 20)
\lvec(63 21)
\lvec(61 21)
\ifill f:0
\move(64 20)
\lvec(67 20)
\lvec(67 21)
\lvec(64 21)
\ifill f:0
\move(68 20)
\lvec(71 20)
\lvec(71 21)
\lvec(68 21)
\ifill f:0
\move(72 20)
\lvec(76 20)
\lvec(76 21)
\lvec(72 21)
\ifill f:0
\move(78 20)
\lvec(82 20)
\lvec(82 21)
\lvec(78 21)
\ifill f:0
\move(85 20)
\lvec(91 20)
\lvec(91 21)
\lvec(85 21)
\ifill f:0
\move(97 20)
\lvec(115 20)
\lvec(115 21)
\lvec(97 21)
\ifill f:0
\move(121 20)
\lvec(129 20)
\lvec(129 21)
\lvec(121 21)
\ifill f:0
\move(132 20)
\lvec(138 20)
\lvec(138 21)
\lvec(132 21)
\ifill f:0
\move(140 20)
\lvec(146 20)
\lvec(146 21)
\lvec(140 21)
\ifill f:0
\move(147 20)
\lvec(152 20)
\lvec(152 21)
\lvec(147 21)
\ifill f:0
\move(153 20)
\lvec(158 20)
\lvec(158 21)
\lvec(153 21)
\ifill f:0
\move(159 20)
\lvec(163 20)
\lvec(163 21)
\lvec(159 21)
\ifill f:0
\move(164 20)
\lvec(168 20)
\lvec(168 21)
\lvec(164 21)
\ifill f:0
\move(169 20)
\lvec(173 20)
\lvec(173 21)
\lvec(169 21)
\ifill f:0
\move(174 20)
\lvec(177 20)
\lvec(177 21)
\lvec(174 21)
\ifill f:0
\move(178 20)
\lvec(181 20)
\lvec(181 21)
\lvec(178 21)
\ifill f:0
\move(182 20)
\lvec(186 20)
\lvec(186 21)
\lvec(182 21)
\ifill f:0
\move(187 20)
\lvec(193 20)
\lvec(193 21)
\lvec(187 21)
\ifill f:0
\move(194 20)
\lvec(197 20)
\lvec(197 21)
\lvec(194 21)
\ifill f:0
\move(198 20)
\lvec(201 20)
\lvec(201 21)
\lvec(198 21)
\ifill f:0
\move(202 20)
\lvec(204 20)
\lvec(204 21)
\lvec(202 21)
\ifill f:0
\move(205 20)
\lvec(208 20)
\lvec(208 21)
\lvec(205 21)
\ifill f:0
\move(209 20)
\lvec(211 20)
\lvec(211 21)
\lvec(209 21)
\ifill f:0
\move(212 20)
\lvec(215 20)
\lvec(215 21)
\lvec(212 21)
\ifill f:0
\move(216 20)
\lvec(218 20)
\lvec(218 21)
\lvec(216 21)
\ifill f:0
\move(219 20)
\lvec(221 20)
\lvec(221 21)
\lvec(219 21)
\ifill f:0
\move(222 20)
\lvec(237 20)
\lvec(237 21)
\lvec(222 21)
\ifill f:0
\move(238 20)
\lvec(251 20)
\lvec(251 21)
\lvec(238 21)
\ifill f:0
\move(252 20)
\lvec(254 20)
\lvec(254 21)
\lvec(252 21)
\ifill f:0
\move(255 20)
\lvec(257 20)
\lvec(257 21)
\lvec(255 21)
\ifill f:0
\move(258 20)
\lvec(265 20)
\lvec(265 21)
\lvec(258 21)
\ifill f:0
\move(266 20)
\lvec(268 20)
\lvec(268 21)
\lvec(266 21)
\ifill f:0
\move(269 20)
\lvec(273 20)
\lvec(273 21)
\lvec(269 21)
\ifill f:0
\move(274 20)
\lvec(281 20)
\lvec(281 21)
\lvec(274 21)
\ifill f:0
\move(282 20)
\lvec(286 20)
\lvec(286 21)
\lvec(282 21)
\ifill f:0
\move(287 20)
\lvec(291 20)
\lvec(291 21)
\lvec(287 21)
\ifill f:0
\move(292 20)
\lvec(296 20)
\lvec(296 21)
\lvec(292 21)
\ifill f:0
\move(297 20)
\lvec(301 20)
\lvec(301 21)
\lvec(297 21)
\ifill f:0
\move(302 20)
\lvec(308 20)
\lvec(308 21)
\lvec(302 21)
\ifill f:0
\move(309 20)
\lvec(313 20)
\lvec(313 21)
\lvec(309 21)
\ifill f:0
\move(314 20)
\lvec(320 20)
\lvec(320 21)
\lvec(314 21)
\ifill f:0
\move(321 20)
\lvec(327 20)
\lvec(327 21)
\lvec(321 21)
\ifill f:0
\move(328 20)
\lvec(334 20)
\lvec(334 21)
\lvec(328 21)
\ifill f:0
\move(335 20)
\lvec(336 20)
\lvec(336 21)
\lvec(335 21)
\ifill f:0
\move(337 20)
\lvec(343 20)
\lvec(343 21)
\lvec(337 21)
\ifill f:0
\move(344 20)
\lvec(345 20)
\lvec(345 21)
\lvec(344 21)
\ifill f:0
\move(346 20)
\lvec(352 20)
\lvec(352 21)
\lvec(346 21)
\ifill f:0
\move(353 20)
\lvec(354 20)
\lvec(354 21)
\lvec(353 21)
\ifill f:0
\move(355 20)
\lvec(365 20)
\lvec(365 21)
\lvec(355 21)
\ifill f:0
\move(366 20)
\lvec(367 20)
\lvec(367 21)
\lvec(366 21)
\ifill f:0
\move(368 20)
\lvec(380 20)
\lvec(380 21)
\lvec(368 21)
\ifill f:0
\move(381 20)
\lvec(382 20)
\lvec(382 21)
\lvec(381 21)
\ifill f:0
\move(383 20)
\lvec(384 20)
\lvec(384 21)
\lvec(383 21)
\ifill f:0
\move(385 20)
\lvec(386 20)
\lvec(386 21)
\lvec(385 21)
\ifill f:0
\move(387 20)
\lvec(419 20)
\lvec(419 21)
\lvec(387 21)
\ifill f:0
\move(420 20)
\lvec(421 20)
\lvec(421 21)
\lvec(420 21)
\ifill f:0
\move(422 20)
\lvec(451 20)
\lvec(451 21)
\lvec(422 21)
\ifill f:0
\move(11 21)
\lvec(13 21)
\lvec(13 22)
\lvec(11 22)
\ifill f:0
\move(14 21)
\lvec(16 21)
\lvec(16 22)
\lvec(14 22)
\ifill f:0
\move(24 21)
\lvec(26 21)
\lvec(26 22)
\lvec(24 22)
\ifill f:0
\move(33 21)
\lvec(36 21)
\lvec(36 22)
\lvec(33 22)
\ifill f:0
\move(39 21)
\lvec(40 21)
\lvec(40 22)
\lvec(39 22)
\ifill f:0
\move(42 21)
\lvec(43 21)
\lvec(43 22)
\lvec(42 22)
\ifill f:0
\move(45 21)
\lvec(46 21)
\lvec(46 22)
\lvec(45 22)
\ifill f:0
\move(47 21)
\lvec(48 21)
\lvec(48 22)
\lvec(47 22)
\ifill f:0
\move(49 21)
\lvec(51 21)
\lvec(51 22)
\lvec(49 22)
\ifill f:0
\move(52 21)
\lvec(56 21)
\lvec(56 22)
\lvec(52 22)
\ifill f:0
\move(57 21)
\lvec(58 21)
\lvec(58 22)
\lvec(57 22)
\ifill f:0
\move(59 21)
\lvec(61 21)
\lvec(61 22)
\lvec(59 22)
\ifill f:0
\move(62 21)
\lvec(63 21)
\lvec(63 22)
\lvec(62 22)
\ifill f:0
\move(64 21)
\lvec(66 21)
\lvec(66 22)
\lvec(64 22)
\ifill f:0
\move(67 21)
\lvec(69 21)
\lvec(69 22)
\lvec(67 22)
\ifill f:0
\move(71 21)
\lvec(73 21)
\lvec(73 22)
\lvec(71 22)
\ifill f:0
\move(75 21)
\lvec(77 21)
\lvec(77 22)
\lvec(75 22)
\ifill f:0
\move(79 21)
\lvec(82 21)
\lvec(82 22)
\lvec(79 22)
\ifill f:0
\move(84 21)
\lvec(88 21)
\lvec(88 22)
\lvec(84 22)
\ifill f:0
\move(90 21)
\lvec(95 21)
\lvec(95 22)
\lvec(90 22)
\ifill f:0
\move(99 21)
\lvec(115 21)
\lvec(115 22)
\lvec(99 22)
\ifill f:0
\move(117 21)
\lvec(135 21)
\lvec(135 22)
\lvec(117 22)
\ifill f:0
\move(139 21)
\lvec(146 21)
\lvec(146 22)
\lvec(139 22)
\ifill f:0
\move(148 21)
\lvec(154 21)
\lvec(154 22)
\lvec(148 22)
\ifill f:0
\move(156 21)
\lvec(161 21)
\lvec(161 22)
\lvec(156 22)
\ifill f:0
\move(163 21)
\lvec(168 21)
\lvec(168 22)
\lvec(163 22)
\ifill f:0
\move(169 21)
\lvec(173 21)
\lvec(173 22)
\lvec(169 22)
\ifill f:0
\move(175 21)
\lvec(179 21)
\lvec(179 22)
\lvec(175 22)
\ifill f:0
\move(180 21)
\lvec(184 21)
\lvec(184 22)
\lvec(180 22)
\ifill f:0
\move(185 21)
\lvec(188 21)
\lvec(188 22)
\lvec(185 22)
\ifill f:0
\move(189 21)
\lvec(193 21)
\lvec(193 22)
\lvec(189 22)
\ifill f:0
\move(194 21)
\lvec(197 21)
\lvec(197 22)
\lvec(194 22)
\ifill f:0
\move(198 21)
\lvec(206 21)
\lvec(206 22)
\lvec(198 22)
\ifill f:0
\move(207 21)
\lvec(213 21)
\lvec(213 22)
\lvec(207 22)
\ifill f:0
\move(214 21)
\lvec(217 21)
\lvec(217 22)
\lvec(214 22)
\ifill f:0
\move(218 21)
\lvec(221 21)
\lvec(221 22)
\lvec(218 22)
\ifill f:0
\move(222 21)
\lvec(224 21)
\lvec(224 22)
\lvec(222 22)
\ifill f:0
\move(225 21)
\lvec(228 21)
\lvec(228 22)
\lvec(225 22)
\ifill f:0
\move(229 21)
\lvec(231 21)
\lvec(231 22)
\lvec(229 22)
\ifill f:0
\move(232 21)
\lvec(235 21)
\lvec(235 22)
\lvec(232 22)
\ifill f:0
\move(236 21)
\lvec(238 21)
\lvec(238 22)
\lvec(236 22)
\ifill f:0
\move(239 21)
\lvec(241 21)
\lvec(241 22)
\lvec(239 22)
\ifill f:0
\move(242 21)
\lvec(251 21)
\lvec(251 22)
\lvec(242 22)
\ifill f:0
\move(252 21)
\lvec(254 21)
\lvec(254 22)
\lvec(252 22)
\ifill f:0
\move(255 21)
\lvec(257 21)
\lvec(257 22)
\lvec(255 22)
\ifill f:0
\move(258 21)
\lvec(260 21)
\lvec(260 22)
\lvec(258 22)
\ifill f:0
\move(261 21)
\lvec(263 21)
\lvec(263 22)
\lvec(261 22)
\ifill f:0
\move(264 21)
\lvec(266 21)
\lvec(266 22)
\lvec(264 22)
\ifill f:0
\move(267 21)
\lvec(269 21)
\lvec(269 22)
\lvec(267 22)
\ifill f:0
\move(270 21)
\lvec(280 21)
\lvec(280 22)
\lvec(270 22)
\ifill f:0
\move(281 21)
\lvec(283 21)
\lvec(283 22)
\lvec(281 22)
\ifill f:0
\move(284 21)
\lvec(286 21)
\lvec(286 22)
\lvec(284 22)
\ifill f:0
\move(287 21)
\lvec(291 21)
\lvec(291 22)
\lvec(287 22)
\ifill f:0
\move(292 21)
\lvec(294 21)
\lvec(294 22)
\lvec(292 22)
\ifill f:0
\move(295 21)
\lvec(299 21)
\lvec(299 22)
\lvec(295 22)
\ifill f:0
\move(300 21)
\lvec(302 21)
\lvec(302 22)
\lvec(300 22)
\ifill f:0
\move(303 21)
\lvec(307 21)
\lvec(307 22)
\lvec(303 22)
\ifill f:0
\move(308 21)
\lvec(312 21)
\lvec(312 22)
\lvec(308 22)
\ifill f:0
\move(313 21)
\lvec(320 21)
\lvec(320 22)
\lvec(313 22)
\ifill f:0
\move(321 21)
\lvec(332 21)
\lvec(332 22)
\lvec(321 22)
\ifill f:0
\move(333 21)
\lvec(337 21)
\lvec(337 22)
\lvec(333 22)
\ifill f:0
\move(338 21)
\lvec(344 21)
\lvec(344 22)
\lvec(338 22)
\ifill f:0
\move(345 21)
\lvec(349 21)
\lvec(349 22)
\lvec(345 22)
\ifill f:0
\move(350 21)
\lvec(356 21)
\lvec(356 22)
\lvec(350 22)
\ifill f:0
\move(357 21)
\lvec(363 21)
\lvec(363 22)
\lvec(357 22)
\ifill f:0
\move(364 21)
\lvec(365 21)
\lvec(365 22)
\lvec(364 22)
\ifill f:0
\move(366 21)
\lvec(372 21)
\lvec(372 22)
\lvec(366 22)
\ifill f:0
\move(373 21)
\lvec(381 21)
\lvec(381 22)
\lvec(373 22)
\ifill f:0
\move(382 21)
\lvec(383 21)
\lvec(383 22)
\lvec(382 22)
\ifill f:0
\move(384 21)
\lvec(392 21)
\lvec(392 22)
\lvec(384 22)
\ifill f:0
\move(393 21)
\lvec(394 21)
\lvec(394 22)
\lvec(393 22)
\ifill f:0
\move(395 21)
\lvec(405 21)
\lvec(405 22)
\lvec(395 22)
\ifill f:0
\move(406 21)
\lvec(407 21)
\lvec(407 22)
\lvec(406 22)
\ifill f:0
\move(408 21)
\lvec(420 21)
\lvec(420 22)
\lvec(408 22)
\ifill f:0
\move(421 21)
\lvec(422 21)
\lvec(422 22)
\lvec(421 22)
\ifill f:0
\move(423 21)
\lvec(424 21)
\lvec(424 22)
\lvec(423 22)
\ifill f:0
\move(425 21)
\lvec(426 21)
\lvec(426 22)
\lvec(425 22)
\ifill f:0
\move(427 21)
\lvec(451 21)
\lvec(451 22)
\lvec(427 22)
\ifill f:0
\move(12 22)
\lvec(16 22)
\lvec(16 23)
\lvec(12 23)
\ifill f:0
\move(18 22)
\lvec(19 22)
\lvec(19 23)
\lvec(18 23)
\ifill f:0
\move(24 22)
\lvec(26 22)
\lvec(26 23)
\lvec(24 23)
\ifill f:0
\move(29 22)
\lvec(36 22)
\lvec(36 23)
\lvec(29 23)
\ifill f:0
\move(39 22)
\lvec(42 22)
\lvec(42 23)
\lvec(39 23)
\ifill f:0
\move(44 22)
\lvec(45 22)
\lvec(45 23)
\lvec(44 23)
\ifill f:0
\move(47 22)
\lvec(48 22)
\lvec(48 23)
\lvec(47 23)
\ifill f:0
\move(49 22)
\lvec(51 22)
\lvec(51 23)
\lvec(49 23)
\ifill f:0
\move(52 22)
\lvec(53 22)
\lvec(53 23)
\lvec(52 23)
\ifill f:0
\move(54 22)
\lvec(55 22)
\lvec(55 23)
\lvec(54 23)
\ifill f:0
\move(56 22)
\lvec(57 22)
\lvec(57 23)
\lvec(56 23)
\ifill f:0
\move(58 22)
\lvec(59 22)
\lvec(59 23)
\lvec(58 23)
\ifill f:0
\move(60 22)
\lvec(61 22)
\lvec(61 23)
\lvec(60 23)
\ifill f:0
\move(62 22)
\lvec(63 22)
\lvec(63 23)
\lvec(62 23)
\ifill f:0
\move(64 22)
\lvec(66 22)
\lvec(66 23)
\lvec(64 23)
\ifill f:0
\move(67 22)
\lvec(69 22)
\lvec(69 23)
\lvec(67 23)
\ifill f:0
\move(70 22)
\lvec(71 22)
\lvec(71 23)
\lvec(70 23)
\ifill f:0
\move(72 22)
\lvec(74 22)
\lvec(74 23)
\lvec(72 23)
\ifill f:0
\move(76 22)
\lvec(78 22)
\lvec(78 23)
\lvec(76 23)
\ifill f:0
\move(79 22)
\lvec(82 22)
\lvec(82 23)
\lvec(79 23)
\ifill f:0
\move(83 22)
\lvec(86 22)
\lvec(86 23)
\lvec(83 23)
\ifill f:0
\move(88 22)
\lvec(91 22)
\lvec(91 23)
\lvec(88 23)
\ifill f:0
\move(93 22)
\lvec(97 22)
\lvec(97 23)
\lvec(93 23)
\ifill f:0
\move(99 22)
\lvec(105 22)
\lvec(105 23)
\lvec(99 23)
\ifill f:0
\move(109 22)
\lvec(126 22)
\lvec(126 23)
\lvec(109 23)
\ifill f:0
\move(128 22)
\lvec(147 22)
\lvec(147 23)
\lvec(128 23)
\ifill f:0
\move(151 22)
\lvec(159 22)
\lvec(159 23)
\lvec(151 23)
\ifill f:0
\move(161 22)
\lvec(167 22)
\lvec(167 23)
\lvec(161 23)
\ifill f:0
\move(169 22)
\lvec(174 22)
\lvec(174 23)
\lvec(169 23)
\ifill f:0
\move(176 22)
\lvec(181 22)
\lvec(181 23)
\lvec(176 23)
\ifill f:0
\move(182 22)
\lvec(187 22)
\lvec(187 23)
\lvec(182 23)
\ifill f:0
\move(188 22)
\lvec(192 22)
\lvec(192 23)
\lvec(188 23)
\ifill f:0
\move(194 22)
\lvec(198 22)
\lvec(198 23)
\lvec(194 23)
\ifill f:0
\move(199 22)
\lvec(202 22)
\lvec(202 23)
\lvec(199 23)
\ifill f:0
\move(203 22)
\lvec(207 22)
\lvec(207 23)
\lvec(203 23)
\ifill f:0
\move(208 22)
\lvec(212 22)
\lvec(212 23)
\lvec(208 23)
\ifill f:0
\move(213 22)
\lvec(216 22)
\lvec(216 23)
\lvec(213 23)
\ifill f:0
\move(217 22)
\lvec(220 22)
\lvec(220 23)
\lvec(217 23)
\ifill f:0
\move(221 22)
\lvec(224 22)
\lvec(224 23)
\lvec(221 23)
\ifill f:0
\move(225 22)
\lvec(228 22)
\lvec(228 23)
\lvec(225 23)
\ifill f:0
\move(229 22)
\lvec(232 22)
\lvec(232 23)
\lvec(229 23)
\ifill f:0
\move(233 22)
\lvec(236 22)
\lvec(236 23)
\lvec(233 23)
\ifill f:0
\move(237 22)
\lvec(243 22)
\lvec(243 23)
\lvec(237 23)
\ifill f:0
\move(244 22)
\lvec(247 22)
\lvec(247 23)
\lvec(244 23)
\ifill f:0
\move(248 22)
\lvec(250 22)
\lvec(250 23)
\lvec(248 23)
\ifill f:0
\move(251 22)
\lvec(254 22)
\lvec(254 23)
\lvec(251 23)
\ifill f:0
\move(255 22)
\lvec(257 22)
\lvec(257 23)
\lvec(255 23)
\ifill f:0
\move(258 22)
\lvec(267 22)
\lvec(267 23)
\lvec(258 23)
\ifill f:0
\move(268 22)
\lvec(270 22)
\lvec(270 23)
\lvec(268 23)
\ifill f:0
\move(271 22)
\lvec(273 22)
\lvec(273 23)
\lvec(271 23)
\ifill f:0
\move(274 22)
\lvec(297 22)
\lvec(297 23)
\lvec(274 23)
\ifill f:0
\move(298 22)
\lvec(300 22)
\lvec(300 23)
\lvec(298 23)
\ifill f:0
\move(301 22)
\lvec(303 22)
\lvec(303 23)
\lvec(301 23)
\ifill f:0
\move(304 22)
\lvec(314 22)
\lvec(314 23)
\lvec(304 23)
\ifill f:0
\move(315 22)
\lvec(317 22)
\lvec(317 23)
\lvec(315 23)
\ifill f:0
\move(318 22)
\lvec(322 22)
\lvec(322 23)
\lvec(318 23)
\ifill f:0
\move(323 22)
\lvec(325 22)
\lvec(325 23)
\lvec(323 23)
\ifill f:0
\move(326 22)
\lvec(330 22)
\lvec(330 23)
\lvec(326 23)
\ifill f:0
\move(331 22)
\lvec(338 22)
\lvec(338 23)
\lvec(331 23)
\ifill f:0
\move(339 22)
\lvec(343 22)
\lvec(343 23)
\lvec(339 23)
\ifill f:0
\move(344 22)
\lvec(348 22)
\lvec(348 23)
\lvec(344 23)
\ifill f:0
\move(349 22)
\lvec(353 22)
\lvec(353 23)
\lvec(349 23)
\ifill f:0
\move(354 22)
\lvec(358 22)
\lvec(358 23)
\lvec(354 23)
\ifill f:0
\move(359 22)
\lvec(363 22)
\lvec(363 23)
\lvec(359 23)
\ifill f:0
\move(364 22)
\lvec(368 22)
\lvec(368 23)
\lvec(364 23)
\ifill f:0
\move(369 22)
\lvec(375 22)
\lvec(375 23)
\lvec(369 23)
\ifill f:0
\move(376 22)
\lvec(382 22)
\lvec(382 23)
\lvec(376 23)
\ifill f:0
\move(383 22)
\lvec(387 22)
\lvec(387 23)
\lvec(383 23)
\ifill f:0
\move(388 22)
\lvec(394 22)
\lvec(394 23)
\lvec(388 23)
\ifill f:0
\move(395 22)
\lvec(396 22)
\lvec(396 23)
\lvec(395 23)
\ifill f:0
\move(397 22)
\lvec(403 22)
\lvec(403 23)
\lvec(397 23)
\ifill f:0
\move(404 22)
\lvec(412 22)
\lvec(412 23)
\lvec(404 23)
\ifill f:0
\move(413 22)
\lvec(421 22)
\lvec(421 23)
\lvec(413 23)
\ifill f:0
\move(422 22)
\lvec(423 22)
\lvec(423 23)
\lvec(422 23)
\ifill f:0
\move(424 22)
\lvec(432 22)
\lvec(432 23)
\lvec(424 23)
\ifill f:0
\move(433 22)
\lvec(434 22)
\lvec(434 23)
\lvec(433 23)
\ifill f:0
\move(435 22)
\lvec(445 22)
\lvec(445 23)
\lvec(435 23)
\ifill f:0
\move(446 22)
\lvec(447 22)
\lvec(447 23)
\lvec(446 23)
\ifill f:0
\move(448 22)
\lvec(451 22)
\lvec(451 23)
\lvec(448 23)
\ifill f:0
\move(12 23)
\lvec(13 23)
\lvec(13 24)
\lvec(12 24)
\ifill f:0
\move(15 23)
\lvec(16 23)
\lvec(16 24)
\lvec(15 24)
\ifill f:0
\move(19 23)
\lvec(20 23)
\lvec(20 24)
\lvec(19 24)
\ifill f:0
\move(25 23)
\lvec(26 23)
\lvec(26 24)
\lvec(25 24)
\ifill f:0
\move(28 23)
\lvec(29 23)
\lvec(29 24)
\lvec(28 24)
\ifill f:0
\move(41 23)
\lvec(44 23)
\lvec(44 24)
\lvec(41 24)
\ifill f:0
\move(46 23)
\lvec(48 23)
\lvec(48 24)
\lvec(46 24)
\ifill f:0
\move(49 23)
\lvec(51 23)
\lvec(51 24)
\lvec(49 24)
\ifill f:0
\move(52 23)
\lvec(54 23)
\lvec(54 24)
\lvec(52 24)
\ifill f:0
\move(55 23)
\lvec(56 23)
\lvec(56 24)
\lvec(55 24)
\ifill f:0
\move(57 23)
\lvec(58 23)
\lvec(58 24)
\lvec(57 24)
\ifill f:0
\move(59 23)
\lvec(60 23)
\lvec(60 24)
\lvec(59 24)
\ifill f:0
\move(62 23)
\lvec(64 23)
\lvec(64 24)
\lvec(62 24)
\ifill f:0
\move(65 23)
\lvec(66 23)
\lvec(66 24)
\lvec(65 24)
\ifill f:0
\move(67 23)
\lvec(68 23)
\lvec(68 24)
\lvec(67 24)
\ifill f:0
\move(69 23)
\lvec(70 23)
\lvec(70 24)
\lvec(69 24)
\ifill f:0
\move(71 23)
\lvec(73 23)
\lvec(73 24)
\lvec(71 24)
\ifill f:0
\move(74 23)
\lvec(76 23)
\lvec(76 24)
\lvec(74 24)
\ifill f:0
\move(77 23)
\lvec(78 23)
\lvec(78 24)
\lvec(77 24)
\ifill f:0
\move(80 23)
\lvec(82 23)
\lvec(82 24)
\lvec(80 24)
\ifill f:0
\move(83 23)
\lvec(85 23)
\lvec(85 24)
\lvec(83 24)
\ifill f:0
\move(87 23)
\lvec(89 23)
\lvec(89 24)
\lvec(87 24)
\ifill f:0
\move(91 23)
\lvec(93 23)
\lvec(93 24)
\lvec(91 24)
\ifill f:0
\move(95 23)
\lvec(98 23)
\lvec(98 24)
\lvec(95 24)
\ifill f:0
\move(100 23)
\lvec(104 23)
\lvec(104 24)
\lvec(100 24)
\ifill f:0
\move(106 23)
\lvec(111 23)
\lvec(111 24)
\lvec(106 24)
\ifill f:0
\move(114 23)
\lvec(122 23)
\lvec(122 24)
\lvec(114 24)
\ifill f:0
\move(129 23)
\lvec(149 23)
\lvec(149 24)
\lvec(129 24)
\ifill f:0
\move(156 23)
\lvec(166 23)
\lvec(166 24)
\lvec(156 24)
\ifill f:0
\move(169 23)
\lvec(176 23)
\lvec(176 24)
\lvec(169 24)
\ifill f:0
\move(178 23)
\lvec(184 23)
\lvec(184 24)
\lvec(178 24)
\ifill f:0
\move(186 23)
\lvec(191 23)
\lvec(191 24)
\lvec(186 24)
\ifill f:0
\move(193 23)
\lvec(198 23)
\lvec(198 24)
\lvec(193 24)
\ifill f:0
\move(199 23)
\lvec(203 23)
\lvec(203 24)
\lvec(199 24)
\ifill f:0
\move(205 23)
\lvec(209 23)
\lvec(209 24)
\lvec(205 24)
\ifill f:0
\move(210 23)
\lvec(214 23)
\lvec(214 24)
\lvec(210 24)
\ifill f:0
\move(216 23)
\lvec(219 23)
\lvec(219 24)
\lvec(216 24)
\ifill f:0
\move(220 23)
\lvec(224 23)
\lvec(224 24)
\lvec(220 24)
\ifill f:0
\move(225 23)
\lvec(229 23)
\lvec(229 24)
\lvec(225 24)
\ifill f:0
\move(230 23)
\lvec(233 23)
\lvec(233 24)
\lvec(230 24)
\ifill f:0
\move(234 23)
\lvec(237 23)
\lvec(237 24)
\lvec(234 24)
\ifill f:0
\move(238 23)
\lvec(246 23)
\lvec(246 24)
\lvec(238 24)
\ifill f:0
\move(247 23)
\lvec(253 23)
\lvec(253 24)
\lvec(247 24)
\ifill f:0
\move(254 23)
\lvec(257 23)
\lvec(257 24)
\lvec(254 24)
\ifill f:0
\move(258 23)
\lvec(261 23)
\lvec(261 24)
\lvec(258 24)
\ifill f:0
\move(262 23)
\lvec(268 23)
\lvec(268 24)
\lvec(262 24)
\ifill f:0
\move(269 23)
\lvec(275 23)
\lvec(275 24)
\lvec(269 24)
\ifill f:0
\move(276 23)
\lvec(282 23)
\lvec(282 24)
\lvec(276 24)
\ifill f:0
\move(283 23)
\lvec(285 23)
\lvec(285 24)
\lvec(283 24)
\ifill f:0
\move(286 23)
\lvec(295 23)
\lvec(295 24)
\lvec(286 24)
\ifill f:0
\move(296 23)
\lvec(298 23)
\lvec(298 24)
\lvec(296 24)
\ifill f:0
\move(299 23)
\lvec(301 23)
\lvec(301 24)
\lvec(299 24)
\ifill f:0
\move(302 23)
\lvec(304 23)
\lvec(304 24)
\lvec(302 24)
\ifill f:0
\move(305 23)
\lvec(307 23)
\lvec(307 24)
\lvec(305 24)
\ifill f:0
\move(308 23)
\lvec(310 23)
\lvec(310 24)
\lvec(308 24)
\ifill f:0
\move(311 23)
\lvec(313 23)
\lvec(313 24)
\lvec(311 24)
\ifill f:0
\move(314 23)
\lvec(316 23)
\lvec(316 24)
\lvec(314 24)
\ifill f:0
\move(317 23)
\lvec(319 23)
\lvec(319 24)
\lvec(317 24)
\ifill f:0
\move(320 23)
\lvec(322 23)
\lvec(322 24)
\lvec(320 24)
\ifill f:0
\move(323 23)
\lvec(325 23)
\lvec(325 24)
\lvec(323 24)
\ifill f:0
\move(326 23)
\lvec(328 23)
\lvec(328 24)
\lvec(326 24)
\ifill f:0
\move(329 23)
\lvec(336 23)
\lvec(336 24)
\lvec(329 24)
\ifill f:0
\move(337 23)
\lvec(339 23)
\lvec(339 24)
\lvec(337 24)
\ifill f:0
\move(340 23)
\lvec(347 23)
\lvec(347 24)
\lvec(340 24)
\ifill f:0
\move(348 23)
\lvec(350 23)
\lvec(350 24)
\lvec(348 24)
\ifill f:0
\move(351 23)
\lvec(355 23)
\lvec(355 24)
\lvec(351 24)
\ifill f:0
\move(356 23)
\lvec(358 23)
\lvec(358 24)
\lvec(356 24)
\ifill f:0
\move(359 23)
\lvec(363 23)
\lvec(363 24)
\lvec(359 24)
\ifill f:0
\move(364 23)
\lvec(371 23)
\lvec(371 24)
\lvec(364 24)
\ifill f:0
\move(372 23)
\lvec(376 23)
\lvec(376 24)
\lvec(372 24)
\ifill f:0
\move(377 23)
\lvec(381 23)
\lvec(381 24)
\lvec(377 24)
\ifill f:0
\move(382 23)
\lvec(386 23)
\lvec(386 24)
\lvec(382 24)
\ifill f:0
\move(387 23)
\lvec(391 23)
\lvec(391 24)
\lvec(387 24)
\ifill f:0
\move(392 23)
\lvec(396 23)
\lvec(396 24)
\lvec(392 24)
\ifill f:0
\move(397 23)
\lvec(403 23)
\lvec(403 24)
\lvec(397 24)
\ifill f:0
\move(404 23)
\lvec(408 23)
\lvec(408 24)
\lvec(404 24)
\ifill f:0
\move(409 23)
\lvec(415 23)
\lvec(415 24)
\lvec(409 24)
\ifill f:0
\move(416 23)
\lvec(422 23)
\lvec(422 24)
\lvec(416 24)
\ifill f:0
\move(423 23)
\lvec(429 23)
\lvec(429 24)
\lvec(423 24)
\ifill f:0
\move(430 23)
\lvec(436 23)
\lvec(436 24)
\lvec(430 24)
\ifill f:0
\move(437 23)
\lvec(443 23)
\lvec(443 24)
\lvec(437 24)
\ifill f:0
\move(444 23)
\lvec(445 23)
\lvec(445 24)
\lvec(444 24)
\ifill f:0
\move(446 23)
\lvec(451 23)
\lvec(451 24)
\lvec(446 24)
\ifill f:0
\move(11 24)
\lvec(14 24)
\lvec(14 25)
\lvec(11 25)
\ifill f:0
\move(23 24)
\lvec(24 24)
\lvec(24 25)
\lvec(23 25)
\ifill f:0
\move(25 24)
\lvec(26 24)
\lvec(26 25)
\lvec(25 25)
\ifill f:0
\move(31 24)
\lvec(32 24)
\lvec(32 25)
\lvec(31 25)
\ifill f:0
\move(44 24)
\lvec(47 24)
\lvec(47 25)
\lvec(44 25)
\ifill f:0
\move(49 24)
\lvec(51 24)
\lvec(51 25)
\lvec(49 25)
\ifill f:0
\move(53 24)
\lvec(54 24)
\lvec(54 25)
\lvec(53 25)
\ifill f:0
\move(56 24)
\lvec(57 24)
\lvec(57 25)
\lvec(56 25)
\ifill f:0
\move(59 24)
\lvec(60 24)
\lvec(60 25)
\lvec(59 25)
\ifill f:0
\move(61 24)
\lvec(62 24)
\lvec(62 25)
\lvec(61 25)
\ifill f:0
\move(63 24)
\lvec(64 24)
\lvec(64 25)
\lvec(63 25)
\ifill f:0
\move(65 24)
\lvec(66 24)
\lvec(66 25)
\lvec(65 25)
\ifill f:0
\move(69 24)
\lvec(70 24)
\lvec(70 25)
\lvec(69 25)
\ifill f:0
\move(71 24)
\lvec(72 24)
\lvec(72 25)
\lvec(71 25)
\ifill f:0
\move(73 24)
\lvec(74 24)
\lvec(74 25)
\lvec(73 25)
\ifill f:0
\move(75 24)
\lvec(76 24)
\lvec(76 25)
\lvec(75 25)
\ifill f:0
\move(78 24)
\lvec(79 24)
\lvec(79 25)
\lvec(78 25)
\ifill f:0
\move(80 24)
\lvec(82 24)
\lvec(82 25)
\lvec(80 25)
\ifill f:0
\move(83 24)
\lvec(84 24)
\lvec(84 25)
\lvec(83 25)
\ifill f:0
\move(86 24)
\lvec(88 24)
\lvec(88 25)
\lvec(86 25)
\ifill f:0
\move(89 24)
\lvec(91 24)
\lvec(91 25)
\lvec(89 25)
\ifill f:0
\move(92 24)
\lvec(95 24)
\lvec(95 25)
\lvec(92 25)
\ifill f:0
\move(96 24)
\lvec(99 24)
\lvec(99 25)
\lvec(96 25)
\ifill f:0
\move(100 24)
\lvec(103 24)
\lvec(103 25)
\lvec(100 25)
\ifill f:0
\move(105 24)
\lvec(108 24)
\lvec(108 25)
\lvec(105 25)
\ifill f:0
\move(110 24)
\lvec(114 24)
\lvec(114 25)
\lvec(110 25)
\ifill f:0
\move(117 24)
\lvec(122 24)
\lvec(122 25)
\lvec(117 25)
\ifill f:0
\move(126 24)
\lvec(134 24)
\lvec(134 25)
\lvec(126 25)
\ifill f:0
\move(141 24)
\lvec(161 24)
\lvec(161 25)
\lvec(141 25)
\ifill f:0
\move(168 24)
\lvec(179 24)
\lvec(179 25)
\lvec(168 25)
\ifill f:0
\move(182 24)
\lvec(189 24)
\lvec(189 25)
\lvec(182 25)
\ifill f:0
\move(192 24)
\lvec(198 24)
\lvec(198 25)
\lvec(192 25)
\ifill f:0
\move(200 24)
\lvec(205 24)
\lvec(205 25)
\lvec(200 25)
\ifill f:0
\move(207 24)
\lvec(212 24)
\lvec(212 25)
\lvec(207 25)
\ifill f:0
\move(213 24)
\lvec(218 24)
\lvec(218 25)
\lvec(213 25)
\ifill f:0
\move(219 24)
\lvec(224 24)
\lvec(224 25)
\lvec(219 25)
\ifill f:0
\move(225 24)
\lvec(229 24)
\lvec(229 25)
\lvec(225 25)
\ifill f:0
\move(230 24)
\lvec(234 24)
\lvec(234 25)
\lvec(230 25)
\ifill f:0
\move(236 24)
\lvec(239 24)
\lvec(239 25)
\lvec(236 25)
\ifill f:0
\move(240 24)
\lvec(244 24)
\lvec(244 25)
\lvec(240 25)
\ifill f:0
\move(245 24)
\lvec(248 24)
\lvec(248 25)
\lvec(245 25)
\ifill f:0
\move(250 24)
\lvec(253 24)
\lvec(253 25)
\lvec(250 25)
\ifill f:0
\move(254 24)
\lvec(257 24)
\lvec(257 25)
\lvec(254 25)
\ifill f:0
\move(258 24)
\lvec(261 24)
\lvec(261 25)
\lvec(258 25)
\ifill f:0
\move(262 24)
\lvec(265 24)
\lvec(265 25)
\lvec(262 25)
\ifill f:0
\move(267 24)
\lvec(269 24)
\lvec(269 25)
\lvec(267 25)
\ifill f:0
\move(270 24)
\lvec(273 24)
\lvec(273 25)
\lvec(270 25)
\ifill f:0
\move(274 24)
\lvec(277 24)
\lvec(277 25)
\lvec(274 25)
\ifill f:0
\move(278 24)
\lvec(281 24)
\lvec(281 25)
\lvec(278 25)
\ifill f:0
\move(282 24)
\lvec(285 24)
\lvec(285 25)
\lvec(282 25)
\ifill f:0
\move(286 24)
\lvec(288 24)
\lvec(288 25)
\lvec(286 25)
\ifill f:0
\move(289 24)
\lvec(292 24)
\lvec(292 25)
\lvec(289 25)
\ifill f:0
\move(293 24)
\lvec(295 24)
\lvec(295 25)
\lvec(293 25)
\ifill f:0
\move(296 24)
\lvec(299 24)
\lvec(299 25)
\lvec(296 25)
\ifill f:0
\move(300 24)
\lvec(302 24)
\lvec(302 25)
\lvec(300 25)
\ifill f:0
\move(303 24)
\lvec(309 24)
\lvec(309 25)
\lvec(303 25)
\ifill f:0
\move(310 24)
\lvec(312 24)
\lvec(312 25)
\lvec(310 25)
\ifill f:0
\move(313 24)
\lvec(322 24)
\lvec(322 25)
\lvec(313 25)
\ifill f:0
\move(323 24)
\lvec(325 24)
\lvec(325 25)
\lvec(323 25)
\ifill f:0
\move(326 24)
\lvec(328 24)
\lvec(328 25)
\lvec(326 25)
\ifill f:0
\move(329 24)
\lvec(331 24)
\lvec(331 25)
\lvec(329 25)
\ifill f:0
\move(332 24)
\lvec(334 24)
\lvec(334 25)
\lvec(332 25)
\ifill f:0
\move(335 24)
\lvec(337 24)
\lvec(337 25)
\lvec(335 25)
\ifill f:0
\move(338 24)
\lvec(340 24)
\lvec(340 25)
\lvec(338 25)
\ifill f:0
\move(341 24)
\lvec(343 24)
\lvec(343 25)
\lvec(341 25)
\ifill f:0
\move(344 24)
\lvec(346 24)
\lvec(346 25)
\lvec(344 25)
\ifill f:0
\move(347 24)
\lvec(349 24)
\lvec(349 25)
\lvec(347 25)
\ifill f:0
\move(350 24)
\lvec(352 24)
\lvec(352 25)
\lvec(350 25)
\ifill f:0
\move(353 24)
\lvec(355 24)
\lvec(355 25)
\lvec(353 25)
\ifill f:0
\move(356 24)
\lvec(363 24)
\lvec(363 25)
\lvec(356 25)
\ifill f:0
\move(364 24)
\lvec(366 24)
\lvec(366 25)
\lvec(364 25)
\ifill f:0
\move(367 24)
\lvec(369 24)
\lvec(369 25)
\lvec(367 25)
\ifill f:0
\move(370 24)
\lvec(374 24)
\lvec(374 25)
\lvec(370 25)
\ifill f:0
\move(375 24)
\lvec(377 24)
\lvec(377 25)
\lvec(375 25)
\ifill f:0
\move(378 24)
\lvec(385 24)
\lvec(385 25)
\lvec(378 25)
\ifill f:0
\move(386 24)
\lvec(393 24)
\lvec(393 25)
\lvec(386 25)
\ifill f:0
\move(394 24)
\lvec(398 24)
\lvec(398 25)
\lvec(394 25)
\ifill f:0
\move(399 24)
\lvec(401 24)
\lvec(401 25)
\lvec(399 25)
\ifill f:0
\move(402 24)
\lvec(406 24)
\lvec(406 25)
\lvec(402 25)
\ifill f:0
\move(407 24)
\lvec(411 24)
\lvec(411 25)
\lvec(407 25)
\ifill f:0
\move(412 24)
\lvec(416 24)
\lvec(416 25)
\lvec(412 25)
\ifill f:0
\move(417 24)
\lvec(421 24)
\lvec(421 25)
\lvec(417 25)
\ifill f:0
\move(422 24)
\lvec(426 24)
\lvec(426 25)
\lvec(422 25)
\ifill f:0
\move(427 24)
\lvec(431 24)
\lvec(431 25)
\lvec(427 25)
\ifill f:0
\move(432 24)
\lvec(438 24)
\lvec(438 25)
\lvec(432 25)
\ifill f:0
\move(439 24)
\lvec(443 24)
\lvec(443 25)
\lvec(439 25)
\ifill f:0
\move(444 24)
\lvec(450 24)
\lvec(450 25)
\lvec(444 25)
\ifill f:0
\move(12 25)
\lvec(15 25)
\lvec(15 26)
\lvec(12 26)
\ifill f:0
\move(23 25)
\lvec(24 25)
\lvec(24 26)
\lvec(23 26)
\ifill f:0
\move(25 25)
\lvec(26 25)
\lvec(26 26)
\lvec(25 26)
\ifill f:0
\move(27 25)
\lvec(28 25)
\lvec(28 26)
\lvec(27 26)
\ifill f:0
\move(29 25)
\lvec(30 25)
\lvec(30 26)
\lvec(29 26)
\ifill f:0
\move(32 25)
\lvec(34 25)
\lvec(34 26)
\lvec(32 26)
\ifill f:0
\move(38 25)
\lvec(45 25)
\lvec(45 26)
\lvec(38 26)
\ifill f:0
\move(49 25)
\lvec(52 25)
\lvec(52 26)
\lvec(49 26)
\ifill f:0
\move(54 25)
\lvec(56 25)
\lvec(56 26)
\lvec(54 26)
\ifill f:0
\move(57 25)
\lvec(59 25)
\lvec(59 26)
\lvec(57 26)
\ifill f:0
\move(60 25)
\lvec(62 25)
\lvec(62 26)
\lvec(60 26)
\ifill f:0
\move(63 25)
\lvec(64 25)
\lvec(64 26)
\lvec(63 26)
\ifill f:0
\move(66 25)
\lvec(67 25)
\lvec(67 26)
\lvec(66 26)
\ifill f:0
\move(68 25)
\lvec(69 25)
\lvec(69 26)
\lvec(68 26)
\ifill f:0
\move(70 25)
\lvec(71 25)
\lvec(71 26)
\lvec(70 26)
\ifill f:0
\move(72 25)
\lvec(73 25)
\lvec(73 26)
\lvec(72 26)
\ifill f:0
\move(74 25)
\lvec(75 25)
\lvec(75 26)
\lvec(74 26)
\ifill f:0
\move(76 25)
\lvec(77 25)
\lvec(77 26)
\lvec(76 26)
\ifill f:0
\move(78 25)
\lvec(79 25)
\lvec(79 26)
\lvec(78 26)
\ifill f:0
\move(80 25)
\lvec(82 25)
\lvec(82 26)
\lvec(80 26)
\ifill f:0
\move(83 25)
\lvec(84 25)
\lvec(84 26)
\lvec(83 26)
\ifill f:0
\move(85 25)
\lvec(87 25)
\lvec(87 26)
\lvec(85 26)
\ifill f:0
\move(88 25)
\lvec(89 25)
\lvec(89 26)
\lvec(88 26)
\ifill f:0
\move(91 25)
\lvec(92 25)
\lvec(92 26)
\lvec(91 26)
\ifill f:0
\move(94 25)
\lvec(96 25)
\lvec(96 26)
\lvec(94 26)
\ifill f:0
\move(97 25)
\lvec(99 25)
\lvec(99 26)
\lvec(97 26)
\ifill f:0
\move(100 25)
\lvec(103 25)
\lvec(103 26)
\lvec(100 26)
\ifill f:0
\move(104 25)
\lvec(107 25)
\lvec(107 26)
\lvec(104 26)
\ifill f:0
\move(108 25)
\lvec(111 25)
\lvec(111 26)
\lvec(108 26)
\ifill f:0
\move(113 25)
\lvec(116 25)
\lvec(116 26)
\lvec(113 26)
\ifill f:0
\move(118 25)
\lvec(122 25)
\lvec(122 26)
\lvec(118 26)
\ifill f:0
\move(124 25)
\lvec(129 25)
\lvec(129 26)
\lvec(124 26)
\ifill f:0
\move(132 25)
\lvec(138 25)
\lvec(138 26)
\lvec(132 26)
\ifill f:0
\move(143 25)
\lvec(162 25)
\lvec(162 26)
\lvec(143 26)
\ifill f:0
\move(164 25)
\lvec(185 25)
\lvec(185 26)
\lvec(164 26)
\ifill f:0
\move(190 25)
\lvec(198 25)
\lvec(198 26)
\lvec(190 26)
\ifill f:0
\move(201 25)
\lvec(208 25)
\lvec(208 26)
\lvec(201 26)
\ifill f:0
\move(210 25)
\lvec(216 25)
\lvec(216 26)
\lvec(210 26)
\ifill f:0
\move(218 25)
\lvec(223 25)
\lvec(223 26)
\lvec(218 26)
\ifill f:0
\move(225 25)
\lvec(230 25)
\lvec(230 26)
\lvec(225 26)
\ifill f:0
\move(231 25)
\lvec(236 25)
\lvec(236 26)
\lvec(231 26)
\ifill f:0
\move(237 25)
\lvec(242 25)
\lvec(242 26)
\lvec(237 26)
\ifill f:0
\move(243 25)
\lvec(247 25)
\lvec(247 26)
\lvec(243 26)
\ifill f:0
\move(248 25)
\lvec(252 25)
\lvec(252 26)
\lvec(248 26)
\ifill f:0
\move(254 25)
\lvec(257 25)
\lvec(257 26)
\lvec(254 26)
\ifill f:0
\move(259 25)
\lvec(262 25)
\lvec(262 26)
\lvec(259 26)
\ifill f:0
\move(263 25)
\lvec(267 25)
\lvec(267 26)
\lvec(263 26)
\ifill f:0
\move(268 25)
\lvec(271 25)
\lvec(271 26)
\lvec(268 26)
\ifill f:0
\move(272 25)
\lvec(276 25)
\lvec(276 26)
\lvec(272 26)
\ifill f:0
\move(277 25)
\lvec(280 25)
\lvec(280 26)
\lvec(277 26)
\ifill f:0
\move(281 25)
\lvec(284 25)
\lvec(284 26)
\lvec(281 26)
\ifill f:0
\move(285 25)
\lvec(288 25)
\lvec(288 26)
\lvec(285 26)
\ifill f:0
\move(289 25)
\lvec(292 25)
\lvec(292 26)
\lvec(289 26)
\ifill f:0
\move(293 25)
\lvec(296 25)
\lvec(296 26)
\lvec(293 26)
\ifill f:0
\move(297 25)
\lvec(300 25)
\lvec(300 26)
\lvec(297 26)
\ifill f:0
\move(301 25)
\lvec(304 25)
\lvec(304 26)
\lvec(301 26)
\ifill f:0
\move(305 25)
\lvec(307 25)
\lvec(307 26)
\lvec(305 26)
\ifill f:0
\move(308 25)
\lvec(311 25)
\lvec(311 26)
\lvec(308 26)
\ifill f:0
\move(312 25)
\lvec(318 25)
\lvec(318 26)
\lvec(312 26)
\ifill f:0
\move(319 25)
\lvec(325 25)
\lvec(325 26)
\lvec(319 26)
\ifill f:0
\move(326 25)
\lvec(332 25)
\lvec(332 26)
\lvec(326 26)
\ifill f:0
\move(333 25)
\lvec(335 25)
\lvec(335 26)
\lvec(333 26)
\ifill f:0
\move(336 25)
\lvec(345 25)
\lvec(345 26)
\lvec(336 26)
\ifill f:0
\move(346 25)
\lvec(348 25)
\lvec(348 26)
\lvec(346 26)
\ifill f:0
\move(349 25)
\lvec(351 25)
\lvec(351 26)
\lvec(349 26)
\ifill f:0
\move(352 25)
\lvec(354 25)
\lvec(354 26)
\lvec(352 26)
\ifill f:0
\move(355 25)
\lvec(378 25)
\lvec(378 26)
\lvec(355 26)
\ifill f:0
\move(379 25)
\lvec(381 25)
\lvec(381 26)
\lvec(379 26)
\ifill f:0
\move(382 25)
\lvec(384 25)
\lvec(384 26)
\lvec(382 26)
\ifill f:0
\move(385 25)
\lvec(387 25)
\lvec(387 26)
\lvec(385 26)
\ifill f:0
\move(388 25)
\lvec(398 25)
\lvec(398 26)
\lvec(388 26)
\ifill f:0
\move(399 25)
\lvec(401 25)
\lvec(401 26)
\lvec(399 26)
\ifill f:0
\move(402 25)
\lvec(409 25)
\lvec(409 26)
\lvec(402 26)
\ifill f:0
\move(410 25)
\lvec(417 25)
\lvec(417 26)
\lvec(410 26)
\ifill f:0
\move(418 25)
\lvec(425 25)
\lvec(425 26)
\lvec(418 26)
\ifill f:0
\move(426 25)
\lvec(430 25)
\lvec(430 26)
\lvec(426 26)
\ifill f:0
\move(431 25)
\lvec(433 25)
\lvec(433 26)
\lvec(431 26)
\ifill f:0
\move(434 25)
\lvec(438 25)
\lvec(438 26)
\lvec(434 26)
\ifill f:0
\move(439 25)
\lvec(443 25)
\lvec(443 26)
\lvec(439 26)
\ifill f:0
\move(444 25)
\lvec(448 25)
\lvec(448 26)
\lvec(444 26)
\ifill f:0
\move(449 25)
\lvec(451 25)
\lvec(451 26)
\lvec(449 26)
\ifill f:0
\move(14 26)
\lvec(15 26)
\lvec(15 27)
\lvec(14 27)
\ifill f:0
\move(18 26)
\lvec(22 26)
\lvec(22 27)
\lvec(18 27)
\ifill f:0
\move(25 26)
\lvec(26 26)
\lvec(26 27)
\lvec(25 27)
\ifill f:0
\move(34 26)
\lvec(35 26)
\lvec(35 27)
\lvec(34 27)
\ifill f:0
\move(38 26)
\lvec(40 26)
\lvec(40 27)
\lvec(38 27)
\ifill f:0
\move(49 26)
\lvec(52 26)
\lvec(52 27)
\lvec(49 27)
\ifill f:0
\move(55 26)
\lvec(57 26)
\lvec(57 27)
\lvec(55 27)
\ifill f:0
\move(60 26)
\lvec(61 26)
\lvec(61 27)
\lvec(60 27)
\ifill f:0
\move(63 26)
\lvec(64 26)
\lvec(64 27)
\lvec(63 27)
\ifill f:0
\move(66 26)
\lvec(67 26)
\lvec(67 27)
\lvec(66 27)
\ifill f:0
\move(68 26)
\lvec(70 26)
\lvec(70 27)
\lvec(68 27)
\ifill f:0
\move(71 26)
\lvec(72 26)
\lvec(72 27)
\lvec(71 27)
\ifill f:0
\move(73 26)
\lvec(74 26)
\lvec(74 27)
\lvec(73 27)
\ifill f:0
\move(75 26)
\lvec(82 26)
\lvec(82 27)
\lvec(75 27)
\ifill f:0
\move(83 26)
\lvec(84 26)
\lvec(84 27)
\lvec(83 27)
\ifill f:0
\move(85 26)
\lvec(86 26)
\lvec(86 27)
\lvec(85 27)
\ifill f:0
\move(87 26)
\lvec(88 26)
\lvec(88 27)
\lvec(87 27)
\ifill f:0
\move(89 26)
\lvec(91 26)
\lvec(91 27)
\lvec(89 27)
\ifill f:0
\move(92 26)
\lvec(93 26)
\lvec(93 27)
\lvec(92 27)
\ifill f:0
\move(95 26)
\lvec(96 26)
\lvec(96 27)
\lvec(95 27)
\ifill f:0
\move(97 26)
\lvec(99 26)
\lvec(99 27)
\lvec(97 27)
\ifill f:0
\move(100 26)
\lvec(102 26)
\lvec(102 27)
\lvec(100 27)
\ifill f:0
\move(104 26)
\lvec(106 26)
\lvec(106 27)
\lvec(104 27)
\ifill f:0
\move(107 26)
\lvec(109 26)
\lvec(109 27)
\lvec(107 27)
\ifill f:0
\move(111 26)
\lvec(113 26)
\lvec(113 27)
\lvec(111 27)
\ifill f:0
\move(114 26)
\lvec(117 26)
\lvec(117 27)
\lvec(114 27)
\ifill f:0
\move(119 26)
\lvec(122 26)
\lvec(122 27)
\lvec(119 27)
\ifill f:0
\move(124 26)
\lvec(127 26)
\lvec(127 27)
\lvec(124 27)
\ifill f:0
\move(129 26)
\lvec(133 26)
\lvec(133 27)
\lvec(129 27)
\ifill f:0
\move(136 26)
\lvec(141 26)
\lvec(141 27)
\lvec(136 27)
\ifill f:0
\move(144 26)
\lvec(150 26)
\lvec(150 27)
\lvec(144 27)
\ifill f:0
\move(155 26)
\lvec(175 26)
\lvec(175 27)
\lvec(155 27)
\ifill f:0
\move(177 26)
\lvec(199 26)
\lvec(199 27)
\lvec(177 27)
\ifill f:0
\move(204 26)
\lvec(212 26)
\lvec(212 27)
\lvec(204 27)
\ifill f:0
\move(215 26)
\lvec(222 26)
\lvec(222 27)
\lvec(215 27)
\ifill f:0
\move(225 26)
\lvec(231 26)
\lvec(231 27)
\lvec(225 27)
\ifill f:0
\move(233 26)
\lvec(238 26)
\lvec(238 27)
\lvec(233 27)
\ifill f:0
\move(240 26)
\lvec(245 26)
\lvec(245 27)
\lvec(240 27)
\ifill f:0
\move(247 26)
\lvec(252 26)
\lvec(252 27)
\lvec(247 27)
\ifill f:0
\move(253 26)
\lvec(257 26)
\lvec(257 27)
\lvec(253 27)
\ifill f:0
\move(259 26)
\lvec(263 26)
\lvec(263 27)
\lvec(259 27)
\ifill f:0
\move(264 26)
\lvec(268 26)
\lvec(268 27)
\lvec(264 27)
\ifill f:0
\move(270 26)
\lvec(274 26)
\lvec(274 27)
\lvec(270 27)
\ifill f:0
\move(275 26)
\lvec(279 26)
\lvec(279 27)
\lvec(275 27)
\ifill f:0
\move(280 26)
\lvec(283 26)
\lvec(283 27)
\lvec(280 27)
\ifill f:0
\move(285 26)
\lvec(288 26)
\lvec(288 27)
\lvec(285 27)
\ifill f:0
\move(289 26)
\lvec(293 26)
\lvec(293 27)
\lvec(289 27)
\ifill f:0
\move(294 26)
\lvec(297 26)
\lvec(297 27)
\lvec(294 27)
\ifill f:0
\move(298 26)
\lvec(301 26)
\lvec(301 27)
\lvec(298 27)
\ifill f:0
\move(302 26)
\lvec(305 26)
\lvec(305 27)
\lvec(302 27)
\ifill f:0
\move(306 26)
\lvec(321 26)
\lvec(321 27)
\lvec(306 27)
\ifill f:0
\move(322 26)
\lvec(325 26)
\lvec(325 27)
\lvec(322 27)
\ifill f:0
\move(326 26)
\lvec(329 26)
\lvec(329 27)
\lvec(326 27)
\ifill f:0
\move(330 26)
\lvec(336 26)
\lvec(336 27)
\lvec(330 27)
\ifill f:0
\move(337 26)
\lvec(340 26)
\lvec(340 27)
\lvec(337 27)
\ifill f:0
\move(341 26)
\lvec(343 26)
\lvec(343 27)
\lvec(341 27)
\ifill f:0
\move(344 26)
\lvec(347 26)
\lvec(347 27)
\lvec(344 27)
\ifill f:0
\move(348 26)
\lvec(350 26)
\lvec(350 27)
\lvec(348 27)
\ifill f:0
\move(351 26)
\lvec(354 26)
\lvec(354 27)
\lvec(351 27)
\ifill f:0
\move(355 26)
\lvec(357 26)
\lvec(357 27)
\lvec(355 27)
\ifill f:0
\move(358 26)
\lvec(360 26)
\lvec(360 27)
\lvec(358 27)
\ifill f:0
\move(361 26)
\lvec(367 26)
\lvec(367 27)
\lvec(361 27)
\ifill f:0
\move(368 26)
\lvec(370 26)
\lvec(370 27)
\lvec(368 27)
\ifill f:0
\move(371 26)
\lvec(373 26)
\lvec(373 27)
\lvec(371 27)
\ifill f:0
\move(374 26)
\lvec(383 26)
\lvec(383 27)
\lvec(374 27)
\ifill f:0
\move(384 26)
\lvec(386 26)
\lvec(386 27)
\lvec(384 27)
\ifill f:0
\move(387 26)
\lvec(389 26)
\lvec(389 27)
\lvec(387 27)
\ifill f:0
\move(390 26)
\lvec(392 26)
\lvec(392 27)
\lvec(390 27)
\ifill f:0
\move(393 26)
\lvec(395 26)
\lvec(395 27)
\lvec(393 27)
\ifill f:0
\move(396 26)
\lvec(398 26)
\lvec(398 27)
\lvec(396 27)
\ifill f:0
\move(399 26)
\lvec(401 26)
\lvec(401 27)
\lvec(399 27)
\ifill f:0
\move(402 26)
\lvec(404 26)
\lvec(404 27)
\lvec(402 27)
\ifill f:0
\move(405 26)
\lvec(407 26)
\lvec(407 27)
\lvec(405 27)
\ifill f:0
\move(408 26)
\lvec(418 26)
\lvec(418 27)
\lvec(408 27)
\ifill f:0
\move(419 26)
\lvec(421 26)
\lvec(421 27)
\lvec(419 27)
\ifill f:0
\move(422 26)
\lvec(424 26)
\lvec(424 27)
\lvec(422 27)
\ifill f:0
\move(425 26)
\lvec(432 26)
\lvec(432 27)
\lvec(425 27)
\ifill f:0
\move(433 26)
\lvec(435 26)
\lvec(435 27)
\lvec(433 27)
\ifill f:0
\move(436 26)
\lvec(438 26)
\lvec(438 27)
\lvec(436 27)
\ifill f:0
\move(439 26)
\lvec(443 26)
\lvec(443 27)
\lvec(439 27)
\ifill f:0
\move(444 26)
\lvec(446 26)
\lvec(446 27)
\lvec(444 27)
\ifill f:0
\move(447 26)
\lvec(451 26)
\lvec(451 27)
\lvec(447 27)
\ifill f:0
\move(11 27)
\lvec(14 27)
\lvec(14 28)
\lvec(11 28)
\ifill f:0
\move(15 27)
\lvec(16 27)
\lvec(16 28)
\lvec(15 28)
\ifill f:0
\move(25 27)
\lvec(26 27)
\lvec(26 28)
\lvec(25 28)
\ifill f:0
\move(28 27)
\lvec(29 27)
\lvec(29 28)
\lvec(28 28)
\ifill f:0
\move(34 27)
\lvec(35 27)
\lvec(35 28)
\lvec(34 28)
\ifill f:0
\move(37 27)
\lvec(39 27)
\lvec(39 28)
\lvec(37 28)
\ifill f:0
\move(43 27)
\lvec(47 27)
\lvec(47 28)
\lvec(43 28)
\ifill f:0
\move(48 27)
\lvec(54 27)
\lvec(54 28)
\lvec(48 28)
\ifill f:0
\move(58 27)
\lvec(60 27)
\lvec(60 28)
\lvec(58 28)
\ifill f:0
\move(63 27)
\lvec(64 27)
\lvec(64 28)
\lvec(63 28)
\ifill f:0
\move(66 27)
\lvec(68 27)
\lvec(68 28)
\lvec(66 28)
\ifill f:0
\move(69 27)
\lvec(70 27)
\lvec(70 28)
\lvec(69 28)
\ifill f:0
\move(72 27)
\lvec(73 27)
\lvec(73 28)
\lvec(72 28)
\ifill f:0
\move(75 27)
\lvec(76 27)
\lvec(76 28)
\lvec(75 28)
\ifill f:0
\move(77 27)
\lvec(78 27)
\lvec(78 28)
\lvec(77 28)
\ifill f:0
\move(79 27)
\lvec(80 27)
\lvec(80 28)
\lvec(79 28)
\ifill f:0
\move(81 27)
\lvec(82 27)
\lvec(82 28)
\lvec(81 28)
\ifill f:0
\move(83 27)
\lvec(84 27)
\lvec(84 28)
\lvec(83 28)
\ifill f:0
\move(85 27)
\lvec(86 27)
\lvec(86 28)
\lvec(85 28)
\ifill f:0
\move(87 27)
\lvec(88 27)
\lvec(88 28)
\lvec(87 28)
\ifill f:0
\move(89 27)
\lvec(90 27)
\lvec(90 28)
\lvec(89 28)
\ifill f:0
\move(91 27)
\lvec(92 27)
\lvec(92 28)
\lvec(91 28)
\ifill f:0
\move(93 27)
\lvec(94 27)
\lvec(94 28)
\lvec(93 28)
\ifill f:0
\move(95 27)
\lvec(97 27)
\lvec(97 28)
\lvec(95 28)
\ifill f:0
\move(98 27)
\lvec(99 27)
\lvec(99 28)
\lvec(98 28)
\ifill f:0
\move(100 27)
\lvec(102 27)
\lvec(102 28)
\lvec(100 28)
\ifill f:0
\move(103 27)
\lvec(105 27)
\lvec(105 28)
\lvec(103 28)
\ifill f:0
\move(106 27)
\lvec(108 27)
\lvec(108 28)
\lvec(106 28)
\ifill f:0
\move(109 27)
\lvec(111 27)
\lvec(111 28)
\lvec(109 28)
\ifill f:0
\move(112 27)
\lvec(114 27)
\lvec(114 28)
\lvec(112 28)
\ifill f:0
\move(116 27)
\lvec(118 27)
\lvec(118 28)
\lvec(116 28)
\ifill f:0
\move(119 27)
\lvec(122 27)
\lvec(122 28)
\lvec(119 28)
\ifill f:0
\move(123 27)
\lvec(126 27)
\lvec(126 28)
\lvec(123 28)
\ifill f:0
\move(128 27)
\lvec(131 27)
\lvec(131 28)
\lvec(128 28)
\ifill f:0
\move(132 27)
\lvec(136 27)
\lvec(136 28)
\lvec(132 28)
\ifill f:0
\move(138 27)
\lvec(142 27)
\lvec(142 28)
\lvec(138 28)
\ifill f:0
\move(144 27)
\lvec(149 27)
\lvec(149 28)
\lvec(144 28)
\ifill f:0
\move(152 27)
\lvec(157 27)
\lvec(157 28)
\lvec(152 28)
\ifill f:0
\move(161 27)
\lvec(171 27)
\lvec(171 28)
\lvec(161 28)
\ifill f:0
\move(178 27)
\lvec(201 27)
\lvec(201 28)
\lvec(178 28)
\ifill f:0
\move(209 27)
\lvec(221 27)
\lvec(221 28)
\lvec(209 28)
\ifill f:0
\move(225 27)
\lvec(233 27)
\lvec(233 28)
\lvec(225 28)
\ifill f:0
\move(235 27)
\lvec(242 27)
\lvec(242 28)
\lvec(235 28)
\ifill f:0
\move(244 27)
\lvec(250 27)
\lvec(250 28)
\lvec(244 28)
\ifill f:0
\move(252 27)
\lvec(257 27)
\lvec(257 28)
\lvec(252 28)
\ifill f:0
\move(259 27)
\lvec(264 27)
\lvec(264 28)
\lvec(259 28)
\ifill f:0
\move(266 27)
\lvec(271 27)
\lvec(271 28)
\lvec(266 28)
\ifill f:0
\move(272 27)
\lvec(277 27)
\lvec(277 28)
\lvec(272 28)
\ifill f:0
\move(278 27)
\lvec(282 27)
\lvec(282 28)
\lvec(278 28)
\ifill f:0
\move(284 27)
\lvec(288 27)
\lvec(288 28)
\lvec(284 28)
\ifill f:0
\move(289 27)
\lvec(293 27)
\lvec(293 28)
\lvec(289 28)
\ifill f:0
\move(294 27)
\lvec(298 27)
\lvec(298 28)
\lvec(294 28)
\ifill f:0
\move(299 27)
\lvec(303 27)
\lvec(303 28)
\lvec(299 28)
\ifill f:0
\move(304 27)
\lvec(307 27)
\lvec(307 28)
\lvec(304 28)
\ifill f:0
\move(309 27)
\lvec(312 27)
\lvec(312 28)
\lvec(309 28)
\ifill f:0
\move(313 27)
\lvec(317 27)
\lvec(317 28)
\lvec(313 28)
\ifill f:0
\move(318 27)
\lvec(321 27)
\lvec(321 28)
\lvec(318 28)
\ifill f:0
\move(322 27)
\lvec(325 27)
\lvec(325 28)
\lvec(322 28)
\ifill f:0
\move(326 27)
\lvec(329 27)
\lvec(329 28)
\lvec(326 28)
\ifill f:0
\move(330 27)
\lvec(333 27)
\lvec(333 28)
\lvec(330 28)
\ifill f:0
\move(334 27)
\lvec(337 27)
\lvec(337 28)
\lvec(334 28)
\ifill f:0
\move(338 27)
\lvec(341 27)
\lvec(341 28)
\lvec(338 28)
\ifill f:0
\move(342 27)
\lvec(345 27)
\lvec(345 28)
\lvec(342 28)
\ifill f:0
\move(346 27)
\lvec(349 27)
\lvec(349 28)
\lvec(346 28)
\ifill f:0
\move(350 27)
\lvec(353 27)
\lvec(353 28)
\lvec(350 28)
\ifill f:0
\move(354 27)
\lvec(357 27)
\lvec(357 28)
\lvec(354 28)
\ifill f:0
\move(358 27)
\lvec(360 27)
\lvec(360 28)
\lvec(358 28)
\ifill f:0
\move(361 27)
\lvec(364 27)
\lvec(364 28)
\lvec(361 28)
\ifill f:0
\move(365 27)
\lvec(367 27)
\lvec(367 28)
\lvec(365 28)
\ifill f:0
\move(368 27)
\lvec(371 27)
\lvec(371 28)
\lvec(368 28)
\ifill f:0
\move(372 27)
\lvec(374 27)
\lvec(374 28)
\lvec(372 28)
\ifill f:0
\move(375 27)
\lvec(378 27)
\lvec(378 28)
\lvec(375 28)
\ifill f:0
\move(379 27)
\lvec(381 27)
\lvec(381 28)
\lvec(379 28)
\ifill f:0
\move(382 27)
\lvec(385 27)
\lvec(385 28)
\lvec(382 28)
\ifill f:0
\move(386 27)
\lvec(388 27)
\lvec(388 28)
\lvec(386 28)
\ifill f:0
\move(389 27)
\lvec(391 27)
\lvec(391 28)
\lvec(389 28)
\ifill f:0
\move(392 27)
\lvec(398 27)
\lvec(398 28)
\lvec(392 28)
\ifill f:0
\move(399 27)
\lvec(401 27)
\lvec(401 28)
\lvec(399 28)
\ifill f:0
\move(402 27)
\lvec(404 27)
\lvec(404 28)
\lvec(402 28)
\ifill f:0
\move(405 27)
\lvec(407 27)
\lvec(407 28)
\lvec(405 28)
\ifill f:0
\move(408 27)
\lvec(426 27)
\lvec(426 28)
\lvec(408 28)
\ifill f:0
\move(427 27)
\lvec(443 27)
\lvec(443 28)
\lvec(427 28)
\ifill f:0
\move(444 27)
\lvec(446 27)
\lvec(446 28)
\lvec(444 28)
\ifill f:0
\move(447 27)
\lvec(449 27)
\lvec(449 28)
\lvec(447 28)
\ifill f:0
\move(450 27)
\lvec(451 27)
\lvec(451 28)
\lvec(450 28)
\ifill f:0
\move(15 28)
\lvec(16 28)
\lvec(16 29)
\lvec(15 29)
\ifill f:0
\move(19 28)
\lvec(21 28)
\lvec(21 29)
\lvec(19 29)
\ifill f:0
\move(25 28)
\lvec(26 28)
\lvec(26 29)
\lvec(25 29)
\ifill f:0
\move(29 28)
\lvec(30 28)
\lvec(30 29)
\lvec(29 29)
\ifill f:0
\move(37 28)
\lvec(38 28)
\lvec(38 29)
\lvec(37 29)
\ifill f:0
\move(40 28)
\lvec(42 28)
\lvec(42 29)
\lvec(40 29)
\ifill f:0
\move(45 28)
\lvec(51 28)
\lvec(51 29)
\lvec(45 29)
\ifill f:0
\move(52 28)
\lvec(58 28)
\lvec(58 29)
\lvec(52 29)
\ifill f:0
\move(62 28)
\lvec(64 28)
\lvec(64 29)
\lvec(62 29)
\ifill f:0
\move(67 28)
\lvec(68 28)
\lvec(68 29)
\lvec(67 29)
\ifill f:0
\move(70 28)
\lvec(72 28)
\lvec(72 29)
\lvec(70 29)
\ifill f:0
\move(73 28)
\lvec(75 28)
\lvec(75 29)
\lvec(73 29)
\ifill f:0
\move(76 28)
\lvec(78 28)
\lvec(78 29)
\lvec(76 29)
\ifill f:0
\move(79 28)
\lvec(80 28)
\lvec(80 29)
\lvec(79 29)
\ifill f:0
\move(81 28)
\lvec(83 28)
\lvec(83 29)
\lvec(81 29)
\ifill f:0
\move(84 28)
\lvec(85 28)
\lvec(85 29)
\lvec(84 29)
\ifill f:0
\move(86 28)
\lvec(87 28)
\lvec(87 29)
\lvec(86 29)
\ifill f:0
\move(88 28)
\lvec(89 28)
\lvec(89 29)
\lvec(88 29)
\ifill f:0
\move(90 28)
\lvec(91 28)
\lvec(91 29)
\lvec(90 29)
\ifill f:0
\move(92 28)
\lvec(93 28)
\lvec(93 29)
\lvec(92 29)
\ifill f:0
\move(94 28)
\lvec(95 28)
\lvec(95 29)
\lvec(94 29)
\ifill f:0
\move(96 28)
\lvec(97 28)
\lvec(97 29)
\lvec(96 29)
\ifill f:0
\move(98 28)
\lvec(99 28)
\lvec(99 29)
\lvec(98 29)
\ifill f:0
\move(100 28)
\lvec(102 28)
\lvec(102 29)
\lvec(100 29)
\ifill f:0
\move(103 28)
\lvec(104 28)
\lvec(104 29)
\lvec(103 29)
\ifill f:0
\move(105 28)
\lvec(107 28)
\lvec(107 29)
\lvec(105 29)
\ifill f:0
\move(108 28)
\lvec(109 28)
\lvec(109 29)
\lvec(108 29)
\ifill f:0
\move(111 28)
\lvec(112 28)
\lvec(112 29)
\lvec(111 29)
\ifill f:0
\move(113 28)
\lvec(115 28)
\lvec(115 29)
\lvec(113 29)
\ifill f:0
\move(117 28)
\lvec(118 28)
\lvec(118 29)
\lvec(117 29)
\ifill f:0
\move(120 28)
\lvec(122 28)
\lvec(122 29)
\lvec(120 29)
\ifill f:0
\move(123 28)
\lvec(125 28)
\lvec(125 29)
\lvec(123 29)
\ifill f:0
\move(127 28)
\lvec(129 28)
\lvec(129 29)
\lvec(127 29)
\ifill f:0
\move(130 28)
\lvec(133 28)
\lvec(133 29)
\lvec(130 29)
\ifill f:0
\move(135 28)
\lvec(137 28)
\lvec(137 29)
\lvec(135 29)
\ifill f:0
\move(139 28)
\lvec(142 28)
\lvec(142 29)
\lvec(139 29)
\ifill f:0
\move(144 28)
\lvec(148 28)
\lvec(148 29)
\lvec(144 29)
\ifill f:0
\move(150 28)
\lvec(154 28)
\lvec(154 29)
\lvec(150 29)
\ifill f:0
\move(156 28)
\lvec(161 28)
\lvec(161 29)
\lvec(156 29)
\ifill f:0
\move(164 28)
\lvec(170 28)
\lvec(170 29)
\lvec(164 29)
\ifill f:0
\move(174 28)
\lvec(184 28)
\lvec(184 29)
\lvec(174 29)
\ifill f:0
\move(191 28)
\lvec(216 28)
\lvec(216 29)
\lvec(191 29)
\ifill f:0
\move(224 28)
\lvec(236 28)
\lvec(236 29)
\lvec(224 29)
\ifill f:0
\move(240 28)
\lvec(248 28)
\lvec(248 29)
\lvec(240 29)
\ifill f:0
\move(251 28)
\lvec(258 28)
\lvec(258 29)
\lvec(251 29)
\ifill f:0
\move(260 28)
\lvec(266 28)
\lvec(266 29)
\lvec(260 29)
\ifill f:0
\move(268 28)
\lvec(274 28)
\lvec(274 29)
\lvec(268 29)
\ifill f:0
\move(276 28)
\lvec(281 28)
\lvec(281 29)
\lvec(276 29)
\ifill f:0
\move(283 28)
\lvec(287 28)
\lvec(287 29)
\lvec(283 29)
\ifill f:0
\move(289 28)
\lvec(294 28)
\lvec(294 29)
\lvec(289 29)
\ifill f:0
\move(295 28)
\lvec(299 28)
\lvec(299 29)
\lvec(295 29)
\ifill f:0
\move(301 28)
\lvec(305 28)
\lvec(305 29)
\lvec(301 29)
\ifill f:0
\move(306 28)
\lvec(310 28)
\lvec(310 29)
\lvec(306 29)
\ifill f:0
\move(312 28)
\lvec(316 28)
\lvec(316 29)
\lvec(312 29)
\ifill f:0
\move(317 28)
\lvec(321 28)
\lvec(321 29)
\lvec(317 29)
\ifill f:0
\move(322 28)
\lvec(325 28)
\lvec(325 29)
\lvec(322 29)
\ifill f:0
\move(327 28)
\lvec(330 28)
\lvec(330 29)
\lvec(327 29)
\ifill f:0
\move(331 28)
\lvec(335 28)
\lvec(335 29)
\lvec(331 29)
\ifill f:0
\move(336 28)
\lvec(339 28)
\lvec(339 29)
\lvec(336 29)
\ifill f:0
\move(340 28)
\lvec(344 28)
\lvec(344 29)
\lvec(340 29)
\ifill f:0
\move(345 28)
\lvec(348 28)
\lvec(348 29)
\lvec(345 29)
\ifill f:0
\move(349 28)
\lvec(352 28)
\lvec(352 29)
\lvec(349 29)
\ifill f:0
\move(353 28)
\lvec(356 28)
\lvec(356 29)
\lvec(353 29)
\ifill f:0
\move(357 28)
\lvec(360 28)
\lvec(360 29)
\lvec(357 29)
\ifill f:0
\move(361 28)
\lvec(364 28)
\lvec(364 29)
\lvec(361 29)
\ifill f:0
\move(365 28)
\lvec(368 28)
\lvec(368 29)
\lvec(365 29)
\ifill f:0
\move(369 28)
\lvec(372 28)
\lvec(372 29)
\lvec(369 29)
\ifill f:0
\move(373 28)
\lvec(376 28)
\lvec(376 29)
\lvec(373 29)
\ifill f:0
\move(377 28)
\lvec(383 28)
\lvec(383 29)
\lvec(377 29)
\ifill f:0
\move(384 28)
\lvec(387 28)
\lvec(387 29)
\lvec(384 29)
\ifill f:0
\move(388 28)
\lvec(394 28)
\lvec(394 29)
\lvec(388 29)
\ifill f:0
\move(395 28)
\lvec(401 28)
\lvec(401 29)
\lvec(395 29)
\ifill f:0
\move(402 28)
\lvec(408 28)
\lvec(408 29)
\lvec(402 29)
\ifill f:0
\move(409 28)
\lvec(411 28)
\lvec(411 29)
\lvec(409 29)
\ifill f:0
\move(412 28)
\lvec(415 28)
\lvec(415 29)
\lvec(412 29)
\ifill f:0
\move(416 28)
\lvec(418 28)
\lvec(418 29)
\lvec(416 29)
\ifill f:0
\move(419 28)
\lvec(421 28)
\lvec(421 29)
\lvec(419 29)
\ifill f:0
\move(422 28)
\lvec(428 28)
\lvec(428 29)
\lvec(422 29)
\ifill f:0
\move(429 28)
\lvec(431 28)
\lvec(431 29)
\lvec(429 29)
\ifill f:0
\move(432 28)
\lvec(434 28)
\lvec(434 29)
\lvec(432 29)
\ifill f:0
\move(435 28)
\lvec(437 28)
\lvec(437 29)
\lvec(435 29)
\ifill f:0
\move(438 28)
\lvec(451 28)
\lvec(451 29)
\lvec(438 29)
\ifill f:0
\move(12 29)
\lvec(16 29)
\lvec(16 30)
\lvec(12 30)
\ifill f:0
\move(22 29)
\lvec(23 29)
\lvec(23 30)
\lvec(22 30)
\ifill f:0
\move(25 29)
\lvec(26 29)
\lvec(26 30)
\lvec(25 30)
\ifill f:0
\move(33 29)
\lvec(34 29)
\lvec(34 30)
\lvec(33 30)
\ifill f:0
\move(37 29)
\lvec(38 29)
\lvec(38 30)
\lvec(37 30)
\ifill f:0
\move(43 29)
\lvec(44 29)
\lvec(44 30)
\lvec(43 30)
\ifill f:0
\move(47 29)
\lvec(50 29)
\lvec(50 30)
\lvec(47 30)
\ifill f:0
\move(60 29)
\lvec(64 29)
\lvec(64 30)
\lvec(60 30)
\ifill f:0
\move(67 29)
\lvec(69 29)
\lvec(69 30)
\lvec(67 30)
\ifill f:0
\move(72 29)
\lvec(73 29)
\lvec(73 30)
\lvec(72 30)
\ifill f:0
\move(75 29)
\lvec(77 29)
\lvec(77 30)
\lvec(75 30)
\ifill f:0
\move(79 29)
\lvec(80 29)
\lvec(80 30)
\lvec(79 30)
\ifill f:0
\move(81 29)
\lvec(83 29)
\lvec(83 30)
\lvec(81 30)
\ifill f:0
\move(84 29)
\lvec(85 29)
\lvec(85 30)
\lvec(84 30)
\ifill f:0
\move(87 29)
\lvec(88 29)
\lvec(88 30)
\lvec(87 30)
\ifill f:0
\move(89 29)
\lvec(90 29)
\lvec(90 30)
\lvec(89 30)
\ifill f:0
\move(91 29)
\lvec(92 29)
\lvec(92 30)
\lvec(91 30)
\ifill f:0
\move(93 29)
\lvec(94 29)
\lvec(94 30)
\lvec(93 30)
\ifill f:0
\move(95 29)
\lvec(97 29)
\lvec(97 30)
\lvec(95 30)
\ifill f:0
\move(98 29)
\lvec(102 29)
\lvec(102 30)
\lvec(98 30)
\ifill f:0
\move(103 29)
\lvec(104 29)
\lvec(104 30)
\lvec(103 30)
\ifill f:0
\move(105 29)
\lvec(106 29)
\lvec(106 30)
\lvec(105 30)
\ifill f:0
\move(107 29)
\lvec(108 29)
\lvec(108 30)
\lvec(107 30)
\ifill f:0
\move(109 29)
\lvec(111 29)
\lvec(111 30)
\lvec(109 30)
\ifill f:0
\move(112 29)
\lvec(113 29)
\lvec(113 30)
\lvec(112 30)
\ifill f:0
\move(115 29)
\lvec(116 29)
\lvec(116 30)
\lvec(115 30)
\ifill f:0
\move(117 29)
\lvec(119 29)
\lvec(119 30)
\lvec(117 30)
\ifill f:0
\move(120 29)
\lvec(122 29)
\lvec(122 30)
\lvec(120 30)
\ifill f:0
\move(123 29)
\lvec(125 29)
\lvec(125 30)
\lvec(123 30)
\ifill f:0
\move(126 29)
\lvec(128 29)
\lvec(128 30)
\lvec(126 30)
\ifill f:0
\move(129 29)
\lvec(131 29)
\lvec(131 30)
\lvec(129 30)
\ifill f:0
\move(133 29)
\lvec(135 29)
\lvec(135 30)
\lvec(133 30)
\ifill f:0
\move(136 29)
\lvec(139 29)
\lvec(139 30)
\lvec(136 30)
\ifill f:0
\move(140 29)
\lvec(143 29)
\lvec(143 30)
\lvec(140 30)
\ifill f:0
\move(144 29)
\lvec(147 29)
\lvec(147 30)
\lvec(144 30)
\ifill f:0
\move(149 29)
\lvec(152 29)
\lvec(152 30)
\lvec(149 30)
\ifill f:0
\move(154 29)
\lvec(157 29)
\lvec(157 30)
\lvec(154 30)
\ifill f:0
\move(159 29)
\lvec(163 29)
\lvec(163 30)
\lvec(159 30)
\ifill f:0
\move(165 29)
\lvec(170 29)
\lvec(170 30)
\lvec(165 30)
\ifill f:0
\move(173 29)
\lvec(178 29)
\lvec(178 30)
\lvec(173 30)
\ifill f:0
\move(182 29)
\lvec(189 29)
\lvec(189 30)
\lvec(182 30)
\ifill f:0
\move(194 29)
\lvec(217 29)
\lvec(217 30)
\lvec(194 30)
\ifill f:0
\move(219 29)
\lvec(243 29)
\lvec(243 30)
\lvec(219 30)
\ifill f:0
\move(249 29)
\lvec(258 29)
\lvec(258 30)
\lvec(249 30)
\ifill f:0
\move(262 29)
\lvec(269 29)
\lvec(269 30)
\lvec(262 30)
\ifill f:0
\move(272 29)
\lvec(279 29)
\lvec(279 30)
\lvec(272 30)
\ifill f:0
\move(281 29)
\lvec(287 29)
\lvec(287 30)
\lvec(281 30)
\ifill f:0
\move(289 29)
\lvec(294 29)
\lvec(294 30)
\lvec(289 30)
\ifill f:0
\move(296 29)
\lvec(301 29)
\lvec(301 30)
\lvec(296 30)
\ifill f:0
\move(303 29)
\lvec(308 29)
\lvec(308 30)
\lvec(303 30)
\ifill f:0
\move(309 29)
\lvec(314 29)
\lvec(314 30)
\lvec(309 30)
\ifill f:0
\move(315 29)
\lvec(320 29)
\lvec(320 30)
\lvec(315 30)
\ifill f:0
\move(321 29)
\lvec(325 29)
\lvec(325 30)
\lvec(321 30)
\ifill f:0
\move(327 29)
\lvec(331 29)
\lvec(331 30)
\lvec(327 30)
\ifill f:0
\move(332 29)
\lvec(336 29)
\lvec(336 30)
\lvec(332 30)
\ifill f:0
\move(337 29)
\lvec(341 29)
\lvec(341 30)
\lvec(337 30)
\ifill f:0
\move(342 29)
\lvec(346 29)
\lvec(346 30)
\lvec(342 30)
\ifill f:0
\move(347 29)
\lvec(351 29)
\lvec(351 30)
\lvec(347 30)
\ifill f:0
\move(352 29)
\lvec(355 29)
\lvec(355 30)
\lvec(352 30)
\ifill f:0
\move(357 29)
\lvec(360 29)
\lvec(360 30)
\lvec(357 30)
\ifill f:0
\move(361 29)
\lvec(364 29)
\lvec(364 30)
\lvec(361 30)
\ifill f:0
\move(366 29)
\lvec(369 29)
\lvec(369 30)
\lvec(366 30)
\ifill f:0
\move(370 29)
\lvec(373 29)
\lvec(373 30)
\lvec(370 30)
\ifill f:0
\move(374 29)
\lvec(377 29)
\lvec(377 30)
\lvec(374 30)
\ifill f:0
\move(378 29)
\lvec(381 29)
\lvec(381 30)
\lvec(378 30)
\ifill f:0
\move(382 29)
\lvec(386 29)
\lvec(386 30)
\lvec(382 30)
\ifill f:0
\move(387 29)
\lvec(393 29)
\lvec(393 30)
\lvec(387 30)
\ifill f:0
\move(394 29)
\lvec(397 29)
\lvec(397 30)
\lvec(394 30)
\ifill f:0
\move(398 29)
\lvec(401 29)
\lvec(401 30)
\lvec(398 30)
\ifill f:0
\move(402 29)
\lvec(405 29)
\lvec(405 30)
\lvec(402 30)
\ifill f:0
\move(406 29)
\lvec(409 29)
\lvec(409 30)
\lvec(406 30)
\ifill f:0
\move(410 29)
\lvec(412 29)
\lvec(412 30)
\lvec(410 30)
\ifill f:0
\move(413 29)
\lvec(416 29)
\lvec(416 30)
\lvec(413 30)
\ifill f:0
\move(417 29)
\lvec(423 29)
\lvec(423 30)
\lvec(417 30)
\ifill f:0
\move(424 29)
\lvec(427 29)
\lvec(427 30)
\lvec(424 30)
\ifill f:0
\move(428 29)
\lvec(430 29)
\lvec(430 30)
\lvec(428 30)
\ifill f:0
\move(431 29)
\lvec(437 29)
\lvec(437 30)
\lvec(431 30)
\ifill f:0
\move(438 29)
\lvec(440 29)
\lvec(440 30)
\lvec(438 30)
\ifill f:0
\move(441 29)
\lvec(444 29)
\lvec(444 30)
\lvec(441 30)
\ifill f:0
\move(445 29)
\lvec(447 29)
\lvec(447 30)
\lvec(445 30)
\ifill f:0
\move(448 29)
\lvec(450 29)
\lvec(450 30)
\lvec(448 30)
\ifill f:0
\move(13 30)
\lvec(14 30)
\lvec(14 31)
\lvec(13 31)
\ifill f:0
\move(15 30)
\lvec(16 30)
\lvec(16 31)
\lvec(15 31)
\ifill f:0
\move(18 30)
\lvec(19 30)
\lvec(19 31)
\lvec(18 31)
\ifill f:0
\move(20 30)
\lvec(21 30)
\lvec(21 31)
\lvec(20 31)
\ifill f:0
\move(25 30)
\lvec(26 30)
\lvec(26 31)
\lvec(25 31)
\ifill f:0
\move(35 30)
\lvec(36 30)
\lvec(36 31)
\lvec(35 31)
\ifill f:0
\move(39 30)
\lvec(40 30)
\lvec(40 31)
\lvec(39 31)
\ifill f:0
\move(41 30)
\lvec(42 30)
\lvec(42 31)
\lvec(41 31)
\ifill f:0
\move(44 30)
\lvec(45 30)
\lvec(45 31)
\lvec(44 31)
\ifill f:0
\move(48 30)
\lvec(50 30)
\lvec(50 31)
\lvec(48 31)
\ifill f:0
\move(55 30)
\lvec(63 30)
\lvec(63 31)
\lvec(55 31)
\ifill f:0
\move(68 30)
\lvec(71 30)
\lvec(71 31)
\lvec(68 31)
\ifill f:0
\move(74 30)
\lvec(76 30)
\lvec(76 31)
\lvec(74 31)
\ifill f:0
\move(78 30)
\lvec(80 30)
\lvec(80 31)
\lvec(78 31)
\ifill f:0
\move(81 30)
\lvec(83 30)
\lvec(83 31)
\lvec(81 31)
\ifill f:0
\move(85 30)
\lvec(86 30)
\lvec(86 31)
\lvec(85 31)
\ifill f:0
\move(87 30)
\lvec(89 30)
\lvec(89 31)
\lvec(87 31)
\ifill f:0
\move(90 30)
\lvec(91 30)
\lvec(91 31)
\lvec(90 31)
\ifill f:0
\move(93 30)
\lvec(94 30)
\lvec(94 31)
\lvec(93 31)
\ifill f:0
\move(95 30)
\lvec(96 30)
\lvec(96 31)
\lvec(95 31)
\ifill f:0
\move(97 30)
\lvec(98 30)
\lvec(98 31)
\lvec(97 31)
\ifill f:0
\move(99 30)
\lvec(100 30)
\lvec(100 31)
\lvec(99 31)
\ifill f:0
\move(101 30)
\lvec(103 30)
\lvec(103 31)
\lvec(101 31)
\ifill f:0
\move(104 30)
\lvec(108 30)
\lvec(108 31)
\lvec(104 31)
\ifill f:0
\move(109 30)
\lvec(110 30)
\lvec(110 31)
\lvec(109 31)
\ifill f:0
\move(111 30)
\lvec(112 30)
\lvec(112 31)
\lvec(111 31)
\ifill f:0
\move(113 30)
\lvec(114 30)
\lvec(114 31)
\lvec(113 31)
\ifill f:0
\move(115 30)
\lvec(117 30)
\lvec(117 31)
\lvec(115 31)
\ifill f:0
\move(118 30)
\lvec(119 30)
\lvec(119 31)
\lvec(118 31)
\ifill f:0
\move(120 30)
\lvec(122 30)
\lvec(122 31)
\lvec(120 31)
\ifill f:0
\move(123 30)
\lvec(124 30)
\lvec(124 31)
\lvec(123 31)
\ifill f:0
\move(125 30)
\lvec(127 30)
\lvec(127 31)
\lvec(125 31)
\ifill f:0
\move(128 30)
\lvec(130 30)
\lvec(130 31)
\lvec(128 31)
\ifill f:0
\move(131 30)
\lvec(133 30)
\lvec(133 31)
\lvec(131 31)
\ifill f:0
\move(134 30)
\lvec(136 30)
\lvec(136 31)
\lvec(134 31)
\ifill f:0
\move(137 30)
\lvec(139 30)
\lvec(139 31)
\lvec(137 31)
\ifill f:0
\move(141 30)
\lvec(143 30)
\lvec(143 31)
\lvec(141 31)
\ifill f:0
\move(144 30)
\lvec(147 30)
\lvec(147 31)
\lvec(144 31)
\ifill f:0
\move(148 30)
\lvec(151 30)
\lvec(151 31)
\lvec(148 31)
\ifill f:0
\move(152 30)
\lvec(155 30)
\lvec(155 31)
\lvec(152 31)
\ifill f:0
\move(156 30)
\lvec(159 30)
\lvec(159 31)
\lvec(156 31)
\ifill f:0
\move(161 30)
\lvec(164 30)
\lvec(164 31)
\lvec(161 31)
\ifill f:0
\move(166 30)
\lvec(170 30)
\lvec(170 31)
\lvec(166 31)
\ifill f:0
\move(172 30)
\lvec(176 30)
\lvec(176 31)
\lvec(172 31)
\ifill f:0
\move(179 30)
\lvec(183 30)
\lvec(183 31)
\lvec(179 31)
\ifill f:0
\move(186 30)
\lvec(192 30)
\lvec(192 31)
\lvec(186 31)
\ifill f:0
\move(195 30)
\lvec(203 30)
\lvec(203 31)
\lvec(195 31)
\ifill f:0
\move(209 30)
\lvec(232 30)
\lvec(232 31)
\lvec(209 31)
\ifill f:0
\move(234 30)
\lvec(259 30)
\lvec(259 31)
\lvec(234 31)
\ifill f:0
\move(265 30)
\lvec(274 30)
\lvec(274 31)
\lvec(265 31)
\ifill f:0
\move(278 30)
\lvec(286 30)
\lvec(286 31)
\lvec(278 31)
\ifill f:0
\move(289 30)
\lvec(295 30)
\lvec(295 31)
\lvec(289 31)
\ifill f:0
\move(298 30)
\lvec(304 30)
\lvec(304 31)
\lvec(298 31)
\ifill f:0
\move(306 30)
\lvec(312 30)
\lvec(312 31)
\lvec(306 31)
\ifill f:0
\move(314 30)
\lvec(319 30)
\lvec(319 31)
\lvec(314 31)
\ifill f:0
\move(321 30)
\lvec(326 30)
\lvec(326 31)
\lvec(321 31)
\ifill f:0
\move(327 30)
\lvec(332 30)
\lvec(332 31)
\lvec(327 31)
\ifill f:0
\move(333 30)
\lvec(338 30)
\lvec(338 31)
\lvec(333 31)
\ifill f:0
\move(339 30)
\lvec(344 30)
\lvec(344 31)
\lvec(339 31)
\ifill f:0
\move(345 30)
\lvec(349 30)
\lvec(349 31)
\lvec(345 31)
\ifill f:0
\move(351 30)
\lvec(355 30)
\lvec(355 31)
\lvec(351 31)
\ifill f:0
\move(356 30)
\lvec(360 30)
\lvec(360 31)
\lvec(356 31)
\ifill f:0
\move(361 30)
\lvec(365 30)
\lvec(365 31)
\lvec(361 31)
\ifill f:0
\move(366 30)
\lvec(370 30)
\lvec(370 31)
\lvec(366 31)
\ifill f:0
\move(371 30)
\lvec(375 30)
\lvec(375 31)
\lvec(371 31)
\ifill f:0
\move(376 30)
\lvec(379 30)
\lvec(379 31)
\lvec(376 31)
\ifill f:0
\move(380 30)
\lvec(384 30)
\lvec(384 31)
\lvec(380 31)
\ifill f:0
\move(385 30)
\lvec(388 30)
\lvec(388 31)
\lvec(385 31)
\ifill f:0
\move(389 30)
\lvec(393 30)
\lvec(393 31)
\lvec(389 31)
\ifill f:0
\move(394 30)
\lvec(397 30)
\lvec(397 31)
\lvec(394 31)
\ifill f:0
\move(398 30)
\lvec(401 30)
\lvec(401 31)
\lvec(398 31)
\ifill f:0
\move(402 30)
\lvec(405 30)
\lvec(405 31)
\lvec(402 31)
\ifill f:0
\move(406 30)
\lvec(409 30)
\lvec(409 31)
\lvec(406 31)
\ifill f:0
\move(410 30)
\lvec(414 30)
\lvec(414 31)
\lvec(410 31)
\ifill f:0
\move(415 30)
\lvec(421 30)
\lvec(421 31)
\lvec(415 31)
\ifill f:0
\move(422 30)
\lvec(425 30)
\lvec(425 31)
\lvec(422 31)
\ifill f:0
\move(426 30)
\lvec(429 30)
\lvec(429 31)
\lvec(426 31)
\ifill f:0
\move(430 30)
\lvec(433 30)
\lvec(433 31)
\lvec(430 31)
\ifill f:0
\move(434 30)
\lvec(437 30)
\lvec(437 31)
\lvec(434 31)
\ifill f:0
\move(438 30)
\lvec(440 30)
\lvec(440 31)
\lvec(438 31)
\ifill f:0
\move(441 30)
\lvec(444 30)
\lvec(444 31)
\lvec(441 31)
\ifill f:0
\move(445 30)
\lvec(451 30)
\lvec(451 31)
\lvec(445 31)
\ifill f:0
\move(11 31)
\lvec(13 31)
\lvec(13 32)
\lvec(11 32)
\ifill f:0
\move(19 31)
\lvec(20 31)
\lvec(20 32)
\lvec(19 32)
\ifill f:0
\move(31 31)
\lvec(32 31)
\lvec(32 32)
\lvec(31 32)
\ifill f:0
\move(35 31)
\lvec(36 31)
\lvec(36 32)
\lvec(35 32)
\ifill f:0
\move(45 31)
\lvec(46 31)
\lvec(46 32)
\lvec(45 32)
\ifill f:0
\move(48 31)
\lvec(50 31)
\lvec(50 32)
\lvec(48 32)
\ifill f:0
\move(53 31)
\lvec(55 31)
\lvec(55 32)
\lvec(53 32)
\ifill f:0
\move(71 31)
\lvec(74 31)
\lvec(74 32)
\lvec(71 32)
\ifill f:0
\move(77 31)
\lvec(79 31)
\lvec(79 32)
\lvec(77 32)
\ifill f:0
\move(81 31)
\lvec(83 31)
\lvec(83 32)
\lvec(81 32)
\ifill f:0
\move(85 31)
\lvec(87 31)
\lvec(87 32)
\lvec(85 32)
\ifill f:0
\move(88 31)
\lvec(90 31)
\lvec(90 32)
\lvec(88 32)
\ifill f:0
\move(91 31)
\lvec(93 31)
\lvec(93 32)
\lvec(91 32)
\ifill f:0
\move(94 31)
\lvec(96 31)
\lvec(96 32)
\lvec(94 32)
\ifill f:0
\move(97 31)
\lvec(98 31)
\lvec(98 32)
\lvec(97 32)
\ifill f:0
\move(99 31)
\lvec(100 31)
\lvec(100 32)
\lvec(99 32)
\ifill f:0
\move(102 31)
\lvec(103 31)
\lvec(103 32)
\lvec(102 32)
\ifill f:0
\move(104 31)
\lvec(105 31)
\lvec(105 32)
\lvec(104 32)
\ifill f:0
\move(106 31)
\lvec(107 31)
\lvec(107 32)
\lvec(106 32)
\ifill f:0
\move(108 31)
\lvec(109 31)
\lvec(109 32)
\lvec(108 32)
\ifill f:0
\move(110 31)
\lvec(111 31)
\lvec(111 32)
\lvec(110 32)
\ifill f:0
\move(112 31)
\lvec(113 31)
\lvec(113 32)
\lvec(112 32)
\ifill f:0
\move(114 31)
\lvec(115 31)
\lvec(115 32)
\lvec(114 32)
\ifill f:0
\move(116 31)
\lvec(117 31)
\lvec(117 32)
\lvec(116 32)
\ifill f:0
\move(118 31)
\lvec(119 31)
\lvec(119 32)
\lvec(118 32)
\ifill f:0
\move(120 31)
\lvec(122 31)
\lvec(122 32)
\lvec(120 32)
\ifill f:0
\move(123 31)
\lvec(124 31)
\lvec(124 32)
\lvec(123 32)
\ifill f:0
\move(125 31)
\lvec(126 31)
\lvec(126 32)
\lvec(125 32)
\ifill f:0
\move(127 31)
\lvec(129 31)
\lvec(129 32)
\lvec(127 32)
\ifill f:0
\move(130 31)
\lvec(132 31)
\lvec(132 32)
\lvec(130 32)
\ifill f:0
\move(133 31)
\lvec(134 31)
\lvec(134 32)
\lvec(133 32)
\ifill f:0
\move(135 31)
\lvec(137 31)
\lvec(137 32)
\lvec(135 32)
\ifill f:0
\move(138 31)
\lvec(140 31)
\lvec(140 32)
\lvec(138 32)
\ifill f:0
\move(141 31)
\lvec(143 31)
\lvec(143 32)
\lvec(141 32)
\ifill f:0
\move(144 31)
\lvec(146 31)
\lvec(146 32)
\lvec(144 32)
\ifill f:0
\move(148 31)
\lvec(150 31)
\lvec(150 32)
\lvec(148 32)
\ifill f:0
\move(151 31)
\lvec(153 31)
\lvec(153 32)
\lvec(151 32)
\ifill f:0
\move(155 31)
\lvec(157 31)
\lvec(157 32)
\lvec(155 32)
\ifill f:0
\move(158 31)
\lvec(161 31)
\lvec(161 32)
\lvec(158 32)
\ifill f:0
\move(163 31)
\lvec(165 31)
\lvec(165 32)
\lvec(163 32)
\ifill f:0
\move(167 31)
\lvec(170 31)
\lvec(170 32)
\lvec(167 32)
\ifill f:0
\move(171 31)
\lvec(175 31)
\lvec(175 32)
\lvec(171 32)
\ifill f:0
\move(177 31)
\lvec(180 31)
\lvec(180 32)
\lvec(177 32)
\ifill f:0
\move(182 31)
\lvec(186 31)
\lvec(186 32)
\lvec(182 32)
\ifill f:0
\move(188 31)
\lvec(193 31)
\lvec(193 32)
\lvec(188 32)
\ifill f:0
\move(196 31)
\lvec(201 31)
\lvec(201 32)
\lvec(196 32)
\ifill f:0
\move(204 31)
\lvec(211 31)
\lvec(211 32)
\lvec(204 32)
\ifill f:0
\move(215 31)
\lvec(227 31)
\lvec(227 32)
\lvec(215 32)
\ifill f:0
\move(235 31)
\lvec(263 31)
\lvec(263 32)
\lvec(235 32)
\ifill f:0
\move(271 31)
\lvec(285 31)
\lvec(285 32)
\lvec(271 32)
\ifill f:0
\move(289 31)
\lvec(298 31)
\lvec(298 32)
\lvec(289 32)
\ifill f:0
\move(301 31)
\lvec(308 31)
\lvec(308 32)
\lvec(301 32)
\ifill f:0
\move(311 31)
\lvec(318 31)
\lvec(318 32)
\lvec(311 32)
\ifill f:0
\move(320 31)
\lvec(326 31)
\lvec(326 32)
\lvec(320 32)
\ifill f:0
\move(328 31)
\lvec(333 31)
\lvec(333 32)
\lvec(328 32)
\ifill f:0
\move(335 31)
\lvec(341 31)
\lvec(341 32)
\lvec(335 32)
\ifill f:0
\move(342 31)
\lvec(347 31)
\lvec(347 32)
\lvec(342 32)
\ifill f:0
\move(349 31)
\lvec(354 31)
\lvec(354 32)
\lvec(349 32)
\ifill f:0
\move(355 31)
\lvec(360 31)
\lvec(360 32)
\lvec(355 32)
\ifill f:0
\move(361 31)
\lvec(365 31)
\lvec(365 32)
\lvec(361 32)
\ifill f:0
\move(367 31)
\lvec(371 31)
\lvec(371 32)
\lvec(367 32)
\ifill f:0
\move(372 31)
\lvec(376 31)
\lvec(376 32)
\lvec(372 32)
\ifill f:0
\move(378 31)
\lvec(382 31)
\lvec(382 32)
\lvec(378 32)
\ifill f:0
\move(383 31)
\lvec(387 31)
\lvec(387 32)
\lvec(383 32)
\ifill f:0
\move(388 31)
\lvec(392 31)
\lvec(392 32)
\lvec(388 32)
\ifill f:0
\move(393 31)
\lvec(397 31)
\lvec(397 32)
\lvec(393 32)
\ifill f:0
\move(398 31)
\lvec(401 31)
\lvec(401 32)
\lvec(398 32)
\ifill f:0
\move(402 31)
\lvec(406 31)
\lvec(406 32)
\lvec(402 32)
\ifill f:0
\move(407 31)
\lvec(411 31)
\lvec(411 32)
\lvec(407 32)
\ifill f:0
\move(412 31)
\lvec(415 31)
\lvec(415 32)
\lvec(412 32)
\ifill f:0
\move(416 31)
\lvec(419 31)
\lvec(419 32)
\lvec(416 32)
\ifill f:0
\move(420 31)
\lvec(424 31)
\lvec(424 32)
\lvec(420 32)
\ifill f:0
\move(425 31)
\lvec(428 31)
\lvec(428 32)
\lvec(425 32)
\ifill f:0
\move(429 31)
\lvec(432 31)
\lvec(432 32)
\lvec(429 32)
\ifill f:0
\move(433 31)
\lvec(436 31)
\lvec(436 32)
\lvec(433 32)
\ifill f:0
\move(437 31)
\lvec(440 31)
\lvec(440 32)
\lvec(437 32)
\ifill f:0
\move(441 31)
\lvec(444 31)
\lvec(444 32)
\lvec(441 32)
\ifill f:0
\move(445 31)
\lvec(448 31)
\lvec(448 32)
\lvec(445 32)
\ifill f:0
\move(449 31)
\lvec(451 31)
\lvec(451 32)
\lvec(449 32)
\ifill f:0
\move(13 32)
\lvec(14 32)
\lvec(14 33)
\lvec(13 33)
\ifill f:0
\move(22 32)
\lvec(23 32)
\lvec(23 33)
\lvec(22 33)
\ifill f:0
\move(34 32)
\lvec(35 32)
\lvec(35 33)
\lvec(34 33)
\ifill f:0
\move(46 32)
\lvec(47 32)
\lvec(47 33)
\lvec(46 33)
\ifill f:0
\move(49 32)
\lvec(50 32)
\lvec(50 33)
\lvec(49 33)
\ifill f:0
\move(52 32)
\lvec(53 32)
\lvec(53 33)
\lvec(52 33)
\ifill f:0
\move(57 32)
\lvec(59 32)
\lvec(59 33)
\lvec(57 33)
\ifill f:0
\move(75 32)
\lvec(78 32)
\lvec(78 33)
\lvec(75 33)
\ifill f:0
\move(81 32)
\lvec(83 32)
\lvec(83 33)
\lvec(81 33)
\ifill f:0
\move(86 32)
\lvec(88 32)
\lvec(88 33)
\lvec(86 33)
\ifill f:0
\move(90 32)
\lvec(91 32)
\lvec(91 33)
\lvec(90 33)
\ifill f:0
\move(93 32)
\lvec(95 32)
\lvec(95 33)
\lvec(93 33)
\ifill f:0
\move(96 32)
\lvec(98 32)
\lvec(98 33)
\lvec(96 33)
\ifill f:0
\move(99 32)
\lvec(100 32)
\lvec(100 33)
\lvec(99 33)
\ifill f:0
\move(102 32)
\lvec(103 32)
\lvec(103 33)
\lvec(102 33)
\ifill f:0
\move(104 32)
\lvec(105 32)
\lvec(105 33)
\lvec(104 33)
\ifill f:0
\move(107 32)
\lvec(108 32)
\lvec(108 33)
\lvec(107 33)
\ifill f:0
\move(109 32)
\lvec(110 32)
\lvec(110 33)
\lvec(109 33)
\ifill f:0
\move(111 32)
\lvec(112 32)
\lvec(112 33)
\lvec(111 33)
\ifill f:0
\move(113 32)
\lvec(114 32)
\lvec(114 33)
\lvec(113 33)
\ifill f:0
\move(115 32)
\lvec(117 32)
\lvec(117 33)
\lvec(115 33)
\ifill f:0
\move(118 32)
\lvec(120 32)
\lvec(120 33)
\lvec(118 33)
\ifill f:0
\move(121 32)
\lvec(122 32)
\lvec(122 33)
\lvec(121 33)
\ifill f:0
\move(123 32)
\lvec(124 32)
\lvec(124 33)
\lvec(123 33)
\ifill f:0
\move(125 32)
\lvec(126 32)
\lvec(126 33)
\lvec(125 33)
\ifill f:0
\move(127 32)
\lvec(128 32)
\lvec(128 33)
\lvec(127 33)
\ifill f:0
\move(129 32)
\lvec(130 32)
\lvec(130 33)
\lvec(129 33)
\ifill f:0
\move(132 32)
\lvec(133 32)
\lvec(133 33)
\lvec(132 33)
\ifill f:0
\move(134 32)
\lvec(135 32)
\lvec(135 33)
\lvec(134 33)
\ifill f:0
\move(136 32)
\lvec(138 32)
\lvec(138 33)
\lvec(136 33)
\ifill f:0
\move(139 32)
\lvec(140 32)
\lvec(140 33)
\lvec(139 33)
\ifill f:0
\move(142 32)
\lvec(143 32)
\lvec(143 33)
\lvec(142 33)
\ifill f:0
\move(144 32)
\lvec(146 32)
\lvec(146 33)
\lvec(144 33)
\ifill f:0
\move(147 32)
\lvec(149 32)
\lvec(149 33)
\lvec(147 33)
\ifill f:0
\move(150 32)
\lvec(152 32)
\lvec(152 33)
\lvec(150 33)
\ifill f:0
\move(153 32)
\lvec(155 32)
\lvec(155 33)
\lvec(153 33)
\ifill f:0
\move(157 32)
\lvec(159 32)
\lvec(159 33)
\lvec(157 33)
\ifill f:0
\move(160 32)
\lvec(162 32)
\lvec(162 33)
\lvec(160 33)
\ifill f:0
\move(164 32)
\lvec(166 32)
\lvec(166 33)
\lvec(164 33)
\ifill f:0
\move(167 32)
\lvec(170 32)
\lvec(170 33)
\lvec(167 33)
\ifill f:0
\move(171 32)
\lvec(174 32)
\lvec(174 33)
\lvec(171 33)
\ifill f:0
\move(176 32)
\lvec(178 32)
\lvec(178 33)
\lvec(176 33)
\ifill f:0
\move(180 32)
\lvec(183 32)
\lvec(183 33)
\lvec(180 33)
\ifill f:0
\move(185 32)
\lvec(188 32)
\lvec(188 33)
\lvec(185 33)
\ifill f:0
\move(190 32)
\lvec(194 32)
\lvec(194 33)
\lvec(190 33)
\ifill f:0
\move(196 32)
\lvec(200 32)
\lvec(200 33)
\lvec(196 33)
\ifill f:0
\move(203 32)
\lvec(207 32)
\lvec(207 33)
\lvec(203 33)
\ifill f:0
\move(210 32)
\lvec(216 32)
\lvec(216 33)
\lvec(210 33)
\ifill f:0
\move(219 32)
\lvec(226 32)
\lvec(226 33)
\lvec(219 33)
\ifill f:0
\move(231 32)
\lvec(242 32)
\lvec(242 33)
\lvec(231 33)
\ifill f:0
\move(251 32)
\lvec(278 32)
\lvec(278 33)
\lvec(251 33)
\ifill f:0
\move(288 32)
\lvec(301 32)
\lvec(301 33)
\lvec(288 33)
\ifill f:0
\move(306 32)
\lvec(315 32)
\lvec(315 33)
\lvec(306 33)
\ifill f:0
\move(318 32)
\lvec(326 32)
\lvec(326 33)
\lvec(318 33)
\ifill f:0
\move(329 32)
\lvec(335 32)
\lvec(335 33)
\lvec(329 33)
\ifill f:0
\move(338 32)
\lvec(344 32)
\lvec(344 33)
\lvec(338 33)
\ifill f:0
\move(346 32)
\lvec(352 32)
\lvec(352 33)
\lvec(346 33)
\ifill f:0
\move(354 32)
\lvec(359 32)
\lvec(359 33)
\lvec(354 33)
\ifill f:0
\move(361 32)
\lvec(366 32)
\lvec(366 33)
\lvec(361 33)
\ifill f:0
\move(368 32)
\lvec(372 32)
\lvec(372 33)
\lvec(368 33)
\ifill f:0
\move(374 32)
\lvec(379 32)
\lvec(379 33)
\lvec(374 33)
\ifill f:0
\move(380 32)
\lvec(385 32)
\lvec(385 33)
\lvec(380 33)
\ifill f:0
\move(386 32)
\lvec(390 32)
\lvec(390 33)
\lvec(386 33)
\ifill f:0
\move(392 32)
\lvec(396 32)
\lvec(396 33)
\lvec(392 33)
\ifill f:0
\move(397 32)
\lvec(401 32)
\lvec(401 33)
\lvec(397 33)
\ifill f:0
\move(403 32)
\lvec(407 32)
\lvec(407 33)
\lvec(403 33)
\ifill f:0
\move(408 32)
\lvec(412 32)
\lvec(412 33)
\lvec(408 33)
\ifill f:0
\move(413 32)
\lvec(417 32)
\lvec(417 33)
\lvec(413 33)
\ifill f:0
\move(418 32)
\lvec(422 32)
\lvec(422 33)
\lvec(418 33)
\ifill f:0
\move(423 32)
\lvec(426 32)
\lvec(426 33)
\lvec(423 33)
\ifill f:0
\move(428 32)
\lvec(431 32)
\lvec(431 33)
\lvec(428 33)
\ifill f:0
\move(432 32)
\lvec(436 32)
\lvec(436 33)
\lvec(432 33)
\ifill f:0
\move(437 32)
\lvec(440 32)
\lvec(440 33)
\lvec(437 33)
\ifill f:0
\move(441 32)
\lvec(444 32)
\lvec(444 33)
\lvec(441 33)
\ifill f:0
\move(446 32)
\lvec(449 32)
\lvec(449 33)
\lvec(446 33)
\ifill f:0
\move(450 32)
\lvec(451 32)
\lvec(451 33)
\lvec(450 33)
\ifill f:0
\move(12 33)
\lvec(13 33)
\lvec(13 34)
\lvec(12 34)
\ifill f:0
\move(14 33)
\lvec(15 33)
\lvec(15 34)
\lvec(14 34)
\ifill f:0
\move(38 33)
\lvec(39 33)
\lvec(39 34)
\lvec(38 34)
\ifill f:0
\move(41 33)
\lvec(42 33)
\lvec(42 34)
\lvec(41 34)
\ifill f:0
\move(49 33)
\lvec(50 33)
\lvec(50 34)
\lvec(49 34)
\ifill f:0
\move(55 33)
\lvec(56 33)
\lvec(56 34)
\lvec(55 34)
\ifill f:0
\move(59 33)
\lvec(61 33)
\lvec(61 34)
\lvec(59 34)
\ifill f:0
\move(68 33)
\lvec(74 33)
\lvec(74 34)
\lvec(68 34)
\ifill f:0
\move(81 33)
\lvec(84 33)
\lvec(84 34)
\lvec(81 34)
\ifill f:0
\move(87 33)
\lvec(89 33)
\lvec(89 34)
\lvec(87 34)
\ifill f:0
\move(92 33)
\lvec(93 33)
\lvec(93 34)
\lvec(92 34)
\ifill f:0
\move(95 33)
\lvec(97 33)
\lvec(97 34)
\lvec(95 34)
\ifill f:0
\move(99 33)
\lvec(100 33)
\lvec(100 34)
\lvec(99 34)
\ifill f:0
\move(102 33)
\lvec(103 33)
\lvec(103 34)
\lvec(102 34)
\ifill f:0
\move(105 33)
\lvec(106 33)
\lvec(106 34)
\lvec(105 34)
\ifill f:0
\move(107 33)
\lvec(109 33)
\lvec(109 34)
\lvec(107 34)
\ifill f:0
\move(110 33)
\lvec(111 33)
\lvec(111 34)
\lvec(110 34)
\ifill f:0
\move(112 33)
\lvec(114 33)
\lvec(114 34)
\lvec(112 34)
\ifill f:0
\move(115 33)
\lvec(116 33)
\lvec(116 34)
\lvec(115 34)
\ifill f:0
\move(117 33)
\lvec(118 33)
\lvec(118 34)
\lvec(117 34)
\ifill f:0
\move(119 33)
\lvec(120 33)
\lvec(120 34)
\lvec(119 34)
\ifill f:0
\move(121 33)
\lvec(122 33)
\lvec(122 34)
\lvec(121 34)
\ifill f:0
\move(123 33)
\lvec(124 33)
\lvec(124 34)
\lvec(123 34)
\ifill f:0
\move(127 33)
\lvec(128 33)
\lvec(128 34)
\lvec(127 34)
\ifill f:0
\move(129 33)
\lvec(130 33)
\lvec(130 34)
\lvec(129 34)
\ifill f:0
\move(131 33)
\lvec(132 33)
\lvec(132 34)
\lvec(131 34)
\ifill f:0
\move(133 33)
\lvec(134 33)
\lvec(134 34)
\lvec(133 34)
\ifill f:0
\move(135 33)
\lvec(136 33)
\lvec(136 34)
\lvec(135 34)
\ifill f:0
\move(137 33)
\lvec(139 33)
\lvec(139 34)
\lvec(137 34)
\ifill f:0
\move(140 33)
\lvec(141 33)
\lvec(141 34)
\lvec(140 34)
\ifill f:0
\move(142 33)
\lvec(143 33)
\lvec(143 34)
\lvec(142 34)
\ifill f:0
\move(145 33)
\lvec(146 33)
\lvec(146 34)
\lvec(145 34)
\ifill f:0
\move(147 33)
\lvec(148 33)
\lvec(148 34)
\lvec(147 34)
\ifill f:0
\move(150 33)
\lvec(151 33)
\lvec(151 34)
\lvec(150 34)
\ifill f:0
\move(153 33)
\lvec(154 33)
\lvec(154 34)
\lvec(153 34)
\ifill f:0
\move(155 33)
\lvec(157 33)
\lvec(157 34)
\lvec(155 34)
\ifill f:0
\move(158 33)
\lvec(160 33)
\lvec(160 34)
\lvec(158 34)
\ifill f:0
\move(161 33)
\lvec(163 33)
\lvec(163 34)
\lvec(161 34)
\ifill f:0
\move(165 33)
\lvec(166 33)
\lvec(166 34)
\lvec(165 34)
\ifill f:0
\move(168 33)
\lvec(170 33)
\lvec(170 34)
\lvec(168 34)
\ifill f:0
\move(171 33)
\lvec(173 33)
\lvec(173 34)
\lvec(171 34)
\ifill f:0
\move(175 33)
\lvec(177 33)
\lvec(177 34)
\lvec(175 34)
\ifill f:0
\move(179 33)
\lvec(181 33)
\lvec(181 34)
\lvec(179 34)
\ifill f:0
\move(183 33)
\lvec(185 33)
\lvec(185 34)
\lvec(183 34)
\ifill f:0
\move(187 33)
\lvec(190 33)
\lvec(190 34)
\lvec(187 34)
\ifill f:0
\move(192 33)
\lvec(194 33)
\lvec(194 34)
\lvec(192 34)
\ifill f:0
\move(196 33)
\lvec(199 33)
\lvec(199 34)
\lvec(196 34)
\ifill f:0
\move(202 33)
\lvec(205 33)
\lvec(205 34)
\lvec(202 34)
\ifill f:0
\move(207 33)
\lvec(211 33)
\lvec(211 34)
\lvec(207 34)
\ifill f:0
\move(214 33)
\lvec(218 33)
\lvec(218 34)
\lvec(214 34)
\ifill f:0
\move(221 33)
\lvec(226 33)
\lvec(226 34)
\lvec(221 34)
\ifill f:0
\move(230 33)
\lvec(236 33)
\lvec(236 34)
\lvec(230 34)
\ifill f:0
\move(240 33)
\lvec(248 33)
\lvec(248 34)
\lvec(240 34)
\ifill f:0
\move(255 33)
\lvec(280 33)
\lvec(280 34)
\lvec(255 34)
\ifill f:0
\move(282 33)
\lvec(309 33)
\lvec(309 34)
\lvec(282 34)
\ifill f:0
\move(316 33)
\lvec(326 33)
\lvec(326 34)
\lvec(316 34)
\ifill f:0
\move(330 33)
\lvec(338 33)
\lvec(338 34)
\lvec(330 34)
\ifill f:0
\move(342 33)
\lvec(349 33)
\lvec(349 34)
\lvec(342 34)
\ifill f:0
\move(352 33)
\lvec(358 33)
\lvec(358 34)
\lvec(352 34)
\ifill f:0
\move(361 33)
\lvec(367 33)
\lvec(367 34)
\lvec(361 34)
\ifill f:0
\move(369 33)
\lvec(374 33)
\lvec(374 34)
\lvec(369 34)
\ifill f:0
\move(377 33)
\lvec(382 33)
\lvec(382 34)
\lvec(377 34)
\ifill f:0
\move(384 33)
\lvec(388 33)
\lvec(388 34)
\lvec(384 34)
\ifill f:0
\move(390 33)
\lvec(395 33)
\lvec(395 34)
\lvec(390 34)
\ifill f:0
\move(397 33)
\lvec(401 33)
\lvec(401 34)
\lvec(397 34)
\ifill f:0
\move(403 33)
\lvec(407 33)
\lvec(407 34)
\lvec(403 34)
\ifill f:0
\move(409 33)
\lvec(413 33)
\lvec(413 34)
\lvec(409 34)
\ifill f:0
\move(415 33)
\lvec(419 33)
\lvec(419 34)
\lvec(415 34)
\ifill f:0
\move(420 33)
\lvec(424 33)
\lvec(424 34)
\lvec(420 34)
\ifill f:0
\move(426 33)
\lvec(429 33)
\lvec(429 34)
\lvec(426 34)
\ifill f:0
\move(431 33)
\lvec(435 33)
\lvec(435 34)
\lvec(431 34)
\ifill f:0
\move(436 33)
\lvec(440 33)
\lvec(440 34)
\lvec(436 34)
\ifill f:0
\move(441 33)
\lvec(445 33)
\lvec(445 34)
\lvec(441 34)
\ifill f:0
\move(446 33)
\lvec(449 33)
\lvec(449 34)
\lvec(446 34)
\ifill f:0
\move(11 34)
\lvec(12 34)
\lvec(12 35)
\lvec(11 35)
\ifill f:0
\move(13 34)
\lvec(14 34)
\lvec(14 35)
\lvec(13 35)
\ifill f:0
\move(19 34)
\lvec(20 34)
\lvec(20 35)
\lvec(19 35)
\ifill f:0
\move(39 34)
\lvec(40 34)
\lvec(40 35)
\lvec(39 35)
\ifill f:0
\move(45 34)
\lvec(46 34)
\lvec(46 35)
\lvec(45 35)
\ifill f:0
\move(47 34)
\lvec(48 34)
\lvec(48 35)
\lvec(47 35)
\ifill f:0
\move(49 34)
\lvec(50 34)
\lvec(50 35)
\lvec(49 35)
\ifill f:0
\move(51 34)
\lvec(52 34)
\lvec(52 35)
\lvec(51 35)
\ifill f:0
\move(54 34)
\lvec(55 34)
\lvec(55 35)
\lvec(54 35)
\ifill f:0
\move(57 34)
\lvec(58 34)
\lvec(58 35)
\lvec(57 35)
\ifill f:0
\move(60 34)
\lvec(62 34)
\lvec(62 35)
\lvec(60 35)
\ifill f:0
\move(66 34)
\lvec(69 34)
\lvec(69 35)
\lvec(66 35)
\ifill f:0
\move(81 34)
\lvec(85 34)
\lvec(85 35)
\lvec(81 35)
\ifill f:0
\move(89 34)
\lvec(91 34)
\lvec(91 35)
\lvec(89 35)
\ifill f:0
\move(94 34)
\lvec(96 34)
\lvec(96 35)
\lvec(94 35)
\ifill f:0
\move(98 34)
\lvec(100 34)
\lvec(100 35)
\lvec(98 35)
\ifill f:0
\move(102 34)
\lvec(104 34)
\lvec(104 35)
\lvec(102 35)
\ifill f:0
\move(105 34)
\lvec(107 34)
\lvec(107 35)
\lvec(105 35)
\ifill f:0
\move(108 34)
\lvec(110 34)
\lvec(110 35)
\lvec(108 35)
\ifill f:0
\move(111 34)
\lvec(113 34)
\lvec(113 35)
\lvec(111 35)
\ifill f:0
\move(114 34)
\lvec(115 34)
\lvec(115 35)
\lvec(114 35)
\ifill f:0
\move(117 34)
\lvec(118 34)
\lvec(118 35)
\lvec(117 35)
\ifill f:0
\move(119 34)
\lvec(120 34)
\lvec(120 35)
\lvec(119 35)
\ifill f:0
\move(121 34)
\lvec(123 34)
\lvec(123 35)
\lvec(121 35)
\ifill f:0
\move(124 34)
\lvec(125 34)
\lvec(125 35)
\lvec(124 35)
\ifill f:0
\move(126 34)
\lvec(127 34)
\lvec(127 35)
\lvec(126 35)
\ifill f:0
\move(128 34)
\lvec(129 34)
\lvec(129 35)
\lvec(128 35)
\ifill f:0
\move(130 34)
\lvec(131 34)
\lvec(131 35)
\lvec(130 35)
\ifill f:0
\move(132 34)
\lvec(133 34)
\lvec(133 35)
\lvec(132 35)
\ifill f:0
\move(134 34)
\lvec(135 34)
\lvec(135 35)
\lvec(134 35)
\ifill f:0
\move(136 34)
\lvec(137 34)
\lvec(137 35)
\lvec(136 35)
\ifill f:0
\move(138 34)
\lvec(139 34)
\lvec(139 35)
\lvec(138 35)
\ifill f:0
\move(140 34)
\lvec(141 34)
\lvec(141 35)
\lvec(140 35)
\ifill f:0
\move(142 34)
\lvec(143 34)
\lvec(143 35)
\lvec(142 35)
\ifill f:0
\move(145 34)
\lvec(146 34)
\lvec(146 35)
\lvec(145 35)
\ifill f:0
\move(147 34)
\lvec(148 34)
\lvec(148 35)
\lvec(147 35)
\ifill f:0
\move(149 34)
\lvec(151 34)
\lvec(151 35)
\lvec(149 35)
\ifill f:0
\move(152 34)
\lvec(153 34)
\lvec(153 35)
\lvec(152 35)
\ifill f:0
\move(154 34)
\lvec(156 34)
\lvec(156 35)
\lvec(154 35)
\ifill f:0
\move(157 34)
\lvec(158 34)
\lvec(158 35)
\lvec(157 35)
\ifill f:0
\move(159 34)
\lvec(161 34)
\lvec(161 35)
\lvec(159 35)
\ifill f:0
\move(162 34)
\lvec(164 34)
\lvec(164 35)
\lvec(162 35)
\ifill f:0
\move(165 34)
\lvec(167 34)
\lvec(167 35)
\lvec(165 35)
\ifill f:0
\move(168 34)
\lvec(170 34)
\lvec(170 35)
\lvec(168 35)
\ifill f:0
\move(171 34)
\lvec(173 34)
\lvec(173 35)
\lvec(171 35)
\ifill f:0
\move(174 34)
\lvec(176 34)
\lvec(176 35)
\lvec(174 35)
\ifill f:0
\move(177 34)
\lvec(179 34)
\lvec(179 35)
\lvec(177 35)
\ifill f:0
\move(181 34)
\lvec(183 34)
\lvec(183 35)
\lvec(181 35)
\ifill f:0
\move(184 34)
\lvec(187 34)
\lvec(187 35)
\lvec(184 35)
\ifill f:0
\move(188 34)
\lvec(190 34)
\lvec(190 35)
\lvec(188 35)
\ifill f:0
\move(192 34)
\lvec(195 34)
\lvec(195 35)
\lvec(192 35)
\ifill f:0
\move(196 34)
\lvec(199 34)
\lvec(199 35)
\lvec(196 35)
\ifill f:0
\move(201 34)
\lvec(204 34)
\lvec(204 35)
\lvec(201 35)
\ifill f:0
\move(205 34)
\lvec(208 34)
\lvec(208 35)
\lvec(205 35)
\ifill f:0
\move(210 34)
\lvec(214 34)
\lvec(214 35)
\lvec(210 35)
\ifill f:0
\move(216 34)
\lvec(219 34)
\lvec(219 35)
\lvec(216 35)
\ifill f:0
\move(222 34)
\lvec(226 34)
\lvec(226 35)
\lvec(222 35)
\ifill f:0
\move(228 34)
\lvec(233 34)
\lvec(233 35)
\lvec(228 35)
\ifill f:0
\move(236 34)
\lvec(241 34)
\lvec(241 35)
\lvec(236 35)
\ifill f:0
\move(245 34)
\lvec(251 34)
\lvec(251 35)
\lvec(245 35)
\ifill f:0
\move(255 34)
\lvec(264 34)
\lvec(264 35)
\lvec(255 35)
\ifill f:0
\move(271 34)
\lvec(297 34)
\lvec(297 35)
\lvec(271 35)
\ifill f:0
\move(299 34)
\lvec(327 34)
\lvec(327 35)
\lvec(299 35)
\ifill f:0
\move(334 34)
\lvec(345 34)
\lvec(345 35)
\lvec(334 35)
\ifill f:0
\move(349 34)
\lvec(357 34)
\lvec(357 35)
\lvec(349 35)
\ifill f:0
\move(361 34)
\lvec(368 34)
\lvec(368 35)
\lvec(361 35)
\ifill f:0
\move(371 34)
\lvec(378 34)
\lvec(378 35)
\lvec(371 35)
\ifill f:0
\move(380 34)
\lvec(386 34)
\lvec(386 35)
\lvec(380 35)
\ifill f:0
\move(389 34)
\lvec(394 34)
\lvec(394 35)
\lvec(389 35)
\ifill f:0
\move(396 34)
\lvec(402 34)
\lvec(402 35)
\lvec(396 35)
\ifill f:0
\move(404 34)
\lvec(409 34)
\lvec(409 35)
\lvec(404 35)
\ifill f:0
\move(410 34)
\lvec(415 34)
\lvec(415 35)
\lvec(410 35)
\ifill f:0
\move(417 34)
\lvec(422 34)
\lvec(422 35)
\lvec(417 35)
\ifill f:0
\move(423 34)
\lvec(428 34)
\lvec(428 35)
\lvec(423 35)
\ifill f:0
\move(430 34)
\lvec(434 34)
\lvec(434 35)
\lvec(430 35)
\ifill f:0
\move(435 34)
\lvec(440 34)
\lvec(440 35)
\lvec(435 35)
\ifill f:0
\move(441 34)
\lvec(445 34)
\lvec(445 35)
\lvec(441 35)
\ifill f:0
\move(447 34)
\lvec(451 34)
\lvec(451 35)
\lvec(447 35)
\ifill f:0
\move(12 35)
\lvec(14 35)
\lvec(14 36)
\lvec(12 36)
\ifill f:0
\move(15 35)
\lvec(16 35)
\lvec(16 36)
\lvec(15 36)
\ifill f:0
\move(22 35)
\lvec(23 35)
\lvec(23 36)
\lvec(22 36)
\ifill f:0
\move(44 35)
\lvec(45 35)
\lvec(45 36)
\lvec(44 36)
\ifill f:0
\move(47 35)
\lvec(48 35)
\lvec(48 36)
\lvec(47 36)
\ifill f:0
\move(49 35)
\lvec(50 35)
\lvec(50 36)
\lvec(49 36)
\ifill f:0
\move(51 35)
\lvec(52 35)
\lvec(52 36)
\lvec(51 36)
\ifill f:0
\move(53 35)
\lvec(54 35)
\lvec(54 36)
\lvec(53 36)
\ifill f:0
\move(55 35)
\lvec(56 35)
\lvec(56 36)
\lvec(55 36)
\ifill f:0
\move(58 35)
\lvec(59 35)
\lvec(59 36)
\lvec(58 36)
\ifill f:0
\move(61 35)
\lvec(62 35)
\lvec(62 36)
\lvec(61 36)
\ifill f:0
\move(66 35)
\lvec(67 35)
\lvec(67 36)
\lvec(66 36)
\ifill f:0
\move(72 35)
\lvec(79 35)
\lvec(79 36)
\lvec(72 36)
\ifill f:0
\move(80 35)
\lvec(87 35)
\lvec(87 36)
\lvec(80 36)
\ifill f:0
\move(92 35)
\lvec(95 35)
\lvec(95 36)
\lvec(92 36)
\ifill f:0
\move(98 35)
\lvec(100 35)
\lvec(100 36)
\lvec(98 36)
\ifill f:0
\move(102 35)
\lvec(104 35)
\lvec(104 36)
\lvec(102 36)
\ifill f:0
\move(106 35)
\lvec(108 35)
\lvec(108 36)
\lvec(106 36)
\ifill f:0
\move(110 35)
\lvec(111 35)
\lvec(111 36)
\lvec(110 36)
\ifill f:0
\move(113 35)
\lvec(114 35)
\lvec(114 36)
\lvec(113 36)
\ifill f:0
\move(116 35)
\lvec(117 35)
\lvec(117 36)
\lvec(116 36)
\ifill f:0
\move(119 35)
\lvec(120 35)
\lvec(120 36)
\lvec(119 36)
\ifill f:0
\move(121 35)
\lvec(123 35)
\lvec(123 36)
\lvec(121 36)
\ifill f:0
\move(124 35)
\lvec(125 35)
\lvec(125 36)
\lvec(124 36)
\ifill f:0
\move(126 35)
\lvec(128 35)
\lvec(128 36)
\lvec(126 36)
\ifill f:0
\move(129 35)
\lvec(130 35)
\lvec(130 36)
\lvec(129 36)
\ifill f:0
\move(131 35)
\lvec(132 35)
\lvec(132 36)
\lvec(131 36)
\ifill f:0
\move(133 35)
\lvec(134 35)
\lvec(134 36)
\lvec(133 36)
\ifill f:0
\move(135 35)
\lvec(136 35)
\lvec(136 36)
\lvec(135 36)
\ifill f:0
\move(137 35)
\lvec(144 35)
\lvec(144 36)
\lvec(137 36)
\ifill f:0
\move(145 35)
\lvec(146 35)
\lvec(146 36)
\lvec(145 36)
\ifill f:0
\move(147 35)
\lvec(148 35)
\lvec(148 36)
\lvec(147 36)
\ifill f:0
\move(149 35)
\lvec(150 35)
\lvec(150 36)
\lvec(149 36)
\ifill f:0
\move(151 35)
\lvec(152 35)
\lvec(152 36)
\lvec(151 36)
\ifill f:0
\move(153 35)
\lvec(154 35)
\lvec(154 36)
\lvec(153 36)
\ifill f:0
\move(156 35)
\lvec(157 35)
\lvec(157 36)
\lvec(156 36)
\ifill f:0
\move(158 35)
\lvec(159 35)
\lvec(159 36)
\lvec(158 36)
\ifill f:0
\move(160 35)
\lvec(162 35)
\lvec(162 36)
\lvec(160 36)
\ifill f:0
\move(163 35)
\lvec(164 35)
\lvec(164 36)
\lvec(163 36)
\ifill f:0
\move(165 35)
\lvec(167 35)
\lvec(167 36)
\lvec(165 36)
\ifill f:0
\move(168 35)
\lvec(170 35)
\lvec(170 36)
\lvec(168 36)
\ifill f:0
\move(171 35)
\lvec(172 35)
\lvec(172 36)
\lvec(171 36)
\ifill f:0
\move(174 35)
\lvec(175 35)
\lvec(175 36)
\lvec(174 36)
\ifill f:0
\move(177 35)
\lvec(178 35)
\lvec(178 36)
\lvec(177 36)
\ifill f:0
\move(180 35)
\lvec(181 35)
\lvec(181 36)
\lvec(180 36)
\ifill f:0
\move(183 35)
\lvec(185 35)
\lvec(185 36)
\lvec(183 36)
\ifill f:0
\move(186 35)
\lvec(188 35)
\lvec(188 36)
\lvec(186 36)
\ifill f:0
\move(189 35)
\lvec(191 35)
\lvec(191 36)
\lvec(189 36)
\ifill f:0
\move(193 35)
\lvec(195 35)
\lvec(195 36)
\lvec(193 36)
\ifill f:0
\move(196 35)
\lvec(199 35)
\lvec(199 36)
\lvec(196 36)
\ifill f:0
\move(200 35)
\lvec(202 35)
\lvec(202 36)
\lvec(200 36)
\ifill f:0
\move(204 35)
\lvec(207 35)
\lvec(207 36)
\lvec(204 36)
\ifill f:0
\move(208 35)
\lvec(211 35)
\lvec(211 36)
\lvec(208 36)
\ifill f:0
\move(213 35)
\lvec(216 35)
\lvec(216 36)
\lvec(213 36)
\ifill f:0
\move(217 35)
\lvec(221 35)
\lvec(221 36)
\lvec(217 36)
\ifill f:0
\move(223 35)
\lvec(226 35)
\lvec(226 36)
\lvec(223 36)
\ifill f:0
\move(228 35)
\lvec(231 35)
\lvec(231 36)
\lvec(228 36)
\ifill f:0
\move(234 35)
\lvec(238 35)
\lvec(238 36)
\lvec(234 36)
\ifill f:0
\move(240 35)
\lvec(245 35)
\lvec(245 36)
\lvec(240 36)
\ifill f:0
\move(247 35)
\lvec(253 35)
\lvec(253 36)
\lvec(247 36)
\ifill f:0
\move(256 35)
\lvec(262 35)
\lvec(262 36)
\lvec(256 36)
\ifill f:0
\move(266 35)
\lvec(274 35)
\lvec(274 36)
\lvec(266 36)
\ifill f:0
\move(279 35)
\lvec(291 35)
\lvec(291 36)
\lvec(279 36)
\ifill f:0
\move(301 35)
\lvec(331 35)
\lvec(331 36)
\lvec(301 36)
\ifill f:0
\move(341 35)
\lvec(355 35)
\lvec(355 36)
\lvec(341 36)
\ifill f:0
\move(360 35)
\lvec(370 35)
\lvec(370 36)
\lvec(360 36)
\ifill f:0
\move(374 35)
\lvec(382 35)
\lvec(382 36)
\lvec(374 36)
\ifill f:0
\move(385 35)
\lvec(392 35)
\lvec(392 36)
\lvec(385 36)
\ifill f:0
\move(395 35)
\lvec(402 35)
\lvec(402 36)
\lvec(395 36)
\ifill f:0
\move(404 35)
\lvec(410 35)
\lvec(410 36)
\lvec(404 36)
\ifill f:0
\move(413 35)
\lvec(418 35)
\lvec(418 36)
\lvec(413 36)
\ifill f:0
\move(420 35)
\lvec(425 35)
\lvec(425 36)
\lvec(420 36)
\ifill f:0
\move(427 35)
\lvec(432 35)
\lvec(432 36)
\lvec(427 36)
\ifill f:0
\move(434 35)
\lvec(439 35)
\lvec(439 36)
\lvec(434 36)
\ifill f:0
\move(441 35)
\lvec(446 35)
\lvec(446 36)
\lvec(441 36)
\ifill f:0
\move(447 35)
\lvec(451 35)
\lvec(451 36)
\lvec(447 36)
\ifill f:0
\move(12 36)
\lvec(13 36)
\lvec(13 37)
\lvec(12 37)
\ifill f:0
\move(15 36)
\lvec(16 36)
\lvec(16 37)
\lvec(15 37)
\ifill f:0
\move(49 36)
\lvec(50 36)
\lvec(50 37)
\lvec(49 37)
\ifill f:0
\move(59 36)
\lvec(60 36)
\lvec(60 37)
\lvec(59 37)
\ifill f:0
\move(62 36)
\lvec(63 36)
\lvec(63 37)
\lvec(62 37)
\ifill f:0
\move(65 36)
\lvec(66 36)
\lvec(66 37)
\lvec(65 37)
\ifill f:0
\move(70 36)
\lvec(71 36)
\lvec(71 37)
\lvec(70 37)
\ifill f:0
\move(76 36)
\lvec(83 36)
\lvec(83 37)
\lvec(76 37)
\ifill f:0
\move(84 36)
\lvec(92 36)
\lvec(92 37)
\lvec(84 37)
\ifill f:0
\move(97 36)
\lvec(100 36)
\lvec(100 37)
\lvec(97 37)
\ifill f:0
\move(103 36)
\lvec(105 36)
\lvec(105 37)
\lvec(103 37)
\ifill f:0
\move(108 36)
\lvec(109 36)
\lvec(109 37)
\lvec(108 37)
\ifill f:0
\move(112 36)
\lvec(113 36)
\lvec(113 37)
\lvec(112 37)
\ifill f:0
\move(115 36)
\lvec(117 36)
\lvec(117 37)
\lvec(115 37)
\ifill f:0
\move(118 36)
\lvec(120 36)
\lvec(120 37)
\lvec(118 37)
\ifill f:0
\move(121 36)
\lvec(123 36)
\lvec(123 37)
\lvec(121 37)
\ifill f:0
\move(124 36)
\lvec(126 36)
\lvec(126 37)
\lvec(124 37)
\ifill f:0
\move(127 36)
\lvec(128 36)
\lvec(128 37)
\lvec(127 37)
\ifill f:0
\move(130 36)
\lvec(131 36)
\lvec(131 37)
\lvec(130 37)
\ifill f:0
\move(132 36)
\lvec(133 36)
\lvec(133 37)
\lvec(132 37)
\ifill f:0
\move(134 36)
\lvec(136 36)
\lvec(136 37)
\lvec(134 37)
\ifill f:0
\move(137 36)
\lvec(138 36)
\lvec(138 37)
\lvec(137 37)
\ifill f:0
\move(139 36)
\lvec(140 36)
\lvec(140 37)
\lvec(139 37)
\ifill f:0
\move(141 36)
\lvec(142 36)
\lvec(142 37)
\lvec(141 37)
\ifill f:0
\move(143 36)
\lvec(144 36)
\lvec(144 37)
\lvec(143 37)
\ifill f:0
\move(145 36)
\lvec(154 36)
\lvec(154 37)
\lvec(145 37)
\ifill f:0
\move(155 36)
\lvec(156 36)
\lvec(156 37)
\lvec(155 37)
\ifill f:0
\move(157 36)
\lvec(158 36)
\lvec(158 37)
\lvec(157 37)
\ifill f:0
\move(159 36)
\lvec(160 36)
\lvec(160 37)
\lvec(159 37)
\ifill f:0
\move(161 36)
\lvec(162 36)
\lvec(162 37)
\lvec(161 37)
\ifill f:0
\move(163 36)
\lvec(165 36)
\lvec(165 37)
\lvec(163 37)
\ifill f:0
\move(166 36)
\lvec(167 36)
\lvec(167 37)
\lvec(166 37)
\ifill f:0
\move(168 36)
\lvec(170 36)
\lvec(170 37)
\lvec(168 37)
\ifill f:0
\move(171 36)
\lvec(172 36)
\lvec(172 37)
\lvec(171 37)
\ifill f:0
\move(173 36)
\lvec(175 36)
\lvec(175 37)
\lvec(173 37)
\ifill f:0
\move(176 36)
\lvec(177 36)
\lvec(177 37)
\lvec(176 37)
\ifill f:0
\move(178 36)
\lvec(180 36)
\lvec(180 37)
\lvec(178 37)
\ifill f:0
\move(181 36)
\lvec(183 36)
\lvec(183 37)
\lvec(181 37)
\ifill f:0
\move(184 36)
\lvec(186 36)
\lvec(186 37)
\lvec(184 37)
\ifill f:0
\move(187 36)
\lvec(189 36)
\lvec(189 37)
\lvec(187 37)
\ifill f:0
\move(190 36)
\lvec(192 36)
\lvec(192 37)
\lvec(190 37)
\ifill f:0
\move(193 36)
\lvec(195 36)
\lvec(195 37)
\lvec(193 37)
\ifill f:0
\move(196 36)
\lvec(198 36)
\lvec(198 37)
\lvec(196 37)
\ifill f:0
\move(200 36)
\lvec(202 36)
\lvec(202 37)
\lvec(200 37)
\ifill f:0
\move(203 36)
\lvec(205 36)
\lvec(205 37)
\lvec(203 37)
\ifill f:0
\move(207 36)
\lvec(209 36)
\lvec(209 37)
\lvec(207 37)
\ifill f:0
\move(211 36)
\lvec(213 36)
\lvec(213 37)
\lvec(211 37)
\ifill f:0
\move(214 36)
\lvec(217 36)
\lvec(217 37)
\lvec(214 37)
\ifill f:0
\move(219 36)
\lvec(221 36)
\lvec(221 37)
\lvec(219 37)
\ifill f:0
\move(223 36)
\lvec(226 36)
\lvec(226 37)
\lvec(223 37)
\ifill f:0
\move(227 36)
\lvec(231 36)
\lvec(231 37)
\lvec(227 37)
\ifill f:0
\move(232 36)
\lvec(236 36)
\lvec(236 37)
\lvec(232 37)
\ifill f:0
\move(238 36)
\lvec(241 36)
\lvec(241 37)
\lvec(238 37)
\ifill f:0
\move(243 36)
\lvec(247 36)
\lvec(247 37)
\lvec(243 37)
\ifill f:0
\move(249 36)
\lvec(253 36)
\lvec(253 37)
\lvec(249 37)
\ifill f:0
\move(256 36)
\lvec(261 36)
\lvec(261 37)
\lvec(256 37)
\ifill f:0
\move(263 36)
\lvec(269 36)
\lvec(269 37)
\lvec(263 37)
\ifill f:0
\move(272 36)
\lvec(278 36)
\lvec(278 37)
\lvec(272 37)
\ifill f:0
\move(282 36)
\lvec(290 36)
\lvec(290 37)
\lvec(282 37)
\ifill f:0
\move(295 36)
\lvec(308 36)
\lvec(308 37)
\lvec(295 37)
\ifill f:0
\move(318 36)
\lvec(350 36)
\lvec(350 37)
\lvec(318 37)
\ifill f:0
\move(360 36)
\lvec(375 36)
\lvec(375 37)
\lvec(360 37)
\ifill f:0
\move(380 36)
\lvec(390 36)
\lvec(390 37)
\lvec(380 37)
\ifill f:0
\move(394 36)
\lvec(402 36)
\lvec(402 37)
\lvec(394 37)
\ifill f:0
\move(405 36)
\lvec(413 36)
\lvec(413 37)
\lvec(405 37)
\ifill f:0
\move(415 36)
\lvec(422 36)
\lvec(422 37)
\lvec(415 37)
\ifill f:0
\move(425 36)
\lvec(431 36)
\lvec(431 37)
\lvec(425 37)
\ifill f:0
\move(433 36)
\lvec(439 36)
\lvec(439 37)
\lvec(433 37)
\ifill f:0
\move(441 36)
\lvec(446 36)
\lvec(446 37)
\lvec(441 37)
\ifill f:0
\move(448 36)
\lvec(451 36)
\lvec(451 37)
\lvec(448 37)
\ifill f:0
\move(11 37)
\lvec(12 37)
\lvec(12 38)
\lvec(11 38)
\ifill f:0
\move(13 37)
\lvec(16 37)
\lvec(16 38)
\lvec(13 38)
\ifill f:0
\move(19 37)
\lvec(20 37)
\lvec(20 38)
\lvec(19 38)
\ifill f:0
\move(49 37)
\lvec(50 37)
\lvec(50 38)
\lvec(49 38)
\ifill f:0
\move(52 37)
\lvec(53 37)
\lvec(53 38)
\lvec(52 38)
\ifill f:0
\move(60 37)
\lvec(61 37)
\lvec(61 38)
\lvec(60 38)
\ifill f:0
\move(62 37)
\lvec(63 37)
\lvec(63 38)
\lvec(62 38)
\ifill f:0
\move(65 37)
\lvec(66 37)
\lvec(66 38)
\lvec(65 38)
\ifill f:0
\move(69 37)
\lvec(70 37)
\lvec(70 38)
\lvec(69 38)
\ifill f:0
\move(73 37)
\lvec(74 37)
\lvec(74 38)
\lvec(73 38)
\ifill f:0
\move(78 37)
\lvec(82 37)
\lvec(82 38)
\lvec(78 38)
\ifill f:0
\move(95 37)
\lvec(99 37)
\lvec(99 38)
\lvec(95 38)
\ifill f:0
\move(104 37)
\lvec(106 37)
\lvec(106 38)
\lvec(104 38)
\ifill f:0
\move(109 37)
\lvec(111 37)
\lvec(111 38)
\lvec(109 38)
\ifill f:0
\move(114 37)
\lvec(116 37)
\lvec(116 38)
\lvec(114 38)
\ifill f:0
\move(118 37)
\lvec(119 37)
\lvec(119 38)
\lvec(118 38)
\ifill f:0
\move(121 37)
\lvec(123 37)
\lvec(123 38)
\lvec(121 38)
\ifill f:0
\move(125 37)
\lvec(126 37)
\lvec(126 38)
\lvec(125 38)
\ifill f:0
\move(128 37)
\lvec(129 37)
\lvec(129 38)
\lvec(128 38)
\ifill f:0
\move(131 37)
\lvec(132 37)
\lvec(132 38)
\lvec(131 38)
\ifill f:0
\move(133 37)
\lvec(135 37)
\lvec(135 38)
\lvec(133 38)
\ifill f:0
\move(136 37)
\lvec(137 37)
\lvec(137 38)
\lvec(136 38)
\ifill f:0
\move(138 37)
\lvec(140 37)
\lvec(140 38)
\lvec(138 38)
\ifill f:0
\move(141 37)
\lvec(142 37)
\lvec(142 38)
\lvec(141 38)
\ifill f:0
\move(143 37)
\lvec(144 37)
\lvec(144 38)
\lvec(143 38)
\ifill f:0
\move(146 37)
\lvec(147 37)
\lvec(147 38)
\lvec(146 38)
\ifill f:0
\move(148 37)
\lvec(149 37)
\lvec(149 38)
\lvec(148 38)
\ifill f:0
\move(150 37)
\lvec(151 37)
\lvec(151 38)
\lvec(150 38)
\ifill f:0
\move(152 37)
\lvec(153 37)
\lvec(153 38)
\lvec(152 38)
\ifill f:0
\move(154 37)
\lvec(155 37)
\lvec(155 38)
\lvec(154 38)
\ifill f:0
\move(156 37)
\lvec(157 37)
\lvec(157 38)
\lvec(156 38)
\ifill f:0
\move(158 37)
\lvec(159 37)
\lvec(159 38)
\lvec(158 38)
\ifill f:0
\move(160 37)
\lvec(161 37)
\lvec(161 38)
\lvec(160 38)
\ifill f:0
\move(162 37)
\lvec(163 37)
\lvec(163 38)
\lvec(162 38)
\ifill f:0
\move(164 37)
\lvec(165 37)
\lvec(165 38)
\lvec(164 38)
\ifill f:0
\move(166 37)
\lvec(167 37)
\lvec(167 38)
\lvec(166 38)
\ifill f:0
\move(168 37)
\lvec(170 37)
\lvec(170 38)
\lvec(168 38)
\ifill f:0
\move(171 37)
\lvec(172 37)
\lvec(172 38)
\lvec(171 38)
\ifill f:0
\move(173 37)
\lvec(174 37)
\lvec(174 38)
\lvec(173 38)
\ifill f:0
\move(175 37)
\lvec(177 37)
\lvec(177 38)
\lvec(175 38)
\ifill f:0
\move(178 37)
\lvec(179 37)
\lvec(179 38)
\lvec(178 38)
\ifill f:0
\move(180 37)
\lvec(182 37)
\lvec(182 38)
\lvec(180 38)
\ifill f:0
\move(183 37)
\lvec(184 37)
\lvec(184 38)
\lvec(183 38)
\ifill f:0
\move(185 37)
\lvec(187 37)
\lvec(187 38)
\lvec(185 38)
\ifill f:0
\move(188 37)
\lvec(190 37)
\lvec(190 38)
\lvec(188 38)
\ifill f:0
\move(191 37)
\lvec(192 37)
\lvec(192 38)
\lvec(191 38)
\ifill f:0
\move(194 37)
\lvec(195 37)
\lvec(195 38)
\lvec(194 38)
\ifill f:0
\move(196 37)
\lvec(198 37)
\lvec(198 38)
\lvec(196 38)
\ifill f:0
\move(199 37)
\lvec(201 37)
\lvec(201 38)
\lvec(199 38)
\ifill f:0
\move(202 37)
\lvec(204 37)
\lvec(204 38)
\lvec(202 38)
\ifill f:0
\move(206 37)
\lvec(208 37)
\lvec(208 38)
\lvec(206 38)
\ifill f:0
\move(209 37)
\lvec(211 37)
\lvec(211 38)
\lvec(209 38)
\ifill f:0
\move(212 37)
\lvec(214 37)
\lvec(214 38)
\lvec(212 38)
\ifill f:0
\move(216 37)
\lvec(218 37)
\lvec(218 38)
\lvec(216 38)
\ifill f:0
\move(220 37)
\lvec(222 37)
\lvec(222 38)
\lvec(220 38)
\ifill f:0
\move(223 37)
\lvec(226 37)
\lvec(226 38)
\lvec(223 38)
\ifill f:0
\move(227 37)
\lvec(230 37)
\lvec(230 38)
\lvec(227 38)
\ifill f:0
\move(231 37)
\lvec(234 37)
\lvec(234 38)
\lvec(231 38)
\ifill f:0
\move(236 37)
\lvec(239 37)
\lvec(239 38)
\lvec(236 38)
\ifill f:0
\move(240 37)
\lvec(243 37)
\lvec(243 38)
\lvec(240 38)
\ifill f:0
\move(245 37)
\lvec(249 37)
\lvec(249 38)
\lvec(245 38)
\ifill f:0
\move(251 37)
\lvec(254 37)
\lvec(254 38)
\lvec(251 38)
\ifill f:0
\move(256 37)
\lvec(260 37)
\lvec(260 38)
\lvec(256 38)
\ifill f:0
\move(262 37)
\lvec(266 37)
\lvec(266 38)
\lvec(262 38)
\ifill f:0
\move(269 37)
\lvec(273 37)
\lvec(273 38)
\lvec(269 38)
\ifill f:0
\move(276 37)
\lvec(281 37)
\lvec(281 38)
\lvec(276 38)
\ifill f:0
\move(284 37)
\lvec(290 37)
\lvec(290 38)
\lvec(284 38)
\ifill f:0
\move(294 37)
\lvec(301 37)
\lvec(301 38)
\lvec(294 38)
\ifill f:0
\move(305 37)
\lvec(315 37)
\lvec(315 38)
\lvec(305 38)
\ifill f:0
\move(322 37)
\lvec(351 37)
\lvec(351 38)
\lvec(322 38)
\ifill f:0
\move(353 37)
\lvec(384 37)
\lvec(384 38)
\lvec(353 38)
\ifill f:0
\move(391 37)
\lvec(403 37)
\lvec(403 38)
\lvec(391 38)
\ifill f:0
\move(407 37)
\lvec(416 37)
\lvec(416 38)
\lvec(407 38)
\ifill f:0
\move(420 37)
\lvec(428 37)
\lvec(428 38)
\lvec(420 38)
\ifill f:0
\move(431 37)
\lvec(438 37)
\lvec(438 38)
\lvec(431 38)
\ifill f:0
\move(441 37)
\lvec(447 37)
\lvec(447 38)
\lvec(441 38)
\ifill f:0
\move(450 37)
\lvec(451 37)
\lvec(451 38)
\lvec(450 38)
\ifill f:0
\move(12 38)
\lvec(13 38)
\lvec(13 39)
\lvec(12 39)
\ifill f:0
\move(15 38)
\lvec(16 38)
\lvec(16 39)
\lvec(15 39)
\ifill f:0
\move(49 38)
\lvec(50 38)
\lvec(50 39)
\lvec(49 39)
\ifill f:0
\move(65 38)
\lvec(66 38)
\lvec(66 39)
\lvec(65 39)
\ifill f:0
\move(68 38)
\lvec(69 38)
\lvec(69 39)
\lvec(68 39)
\ifill f:0
\move(71 38)
\lvec(72 38)
\lvec(72 39)
\lvec(71 39)
\ifill f:0
\move(75 38)
\lvec(76 38)
\lvec(76 39)
\lvec(75 39)
\ifill f:0
\move(80 38)
\lvec(82 38)
\lvec(82 39)
\lvec(80 39)
\ifill f:0
\move(89 38)
\lvec(98 38)
\lvec(98 39)
\lvec(89 39)
\ifill f:0
\move(105 38)
\lvec(108 38)
\lvec(108 39)
\lvec(105 39)
\ifill f:0
\move(112 38)
\lvec(114 38)
\lvec(114 39)
\lvec(112 39)
\ifill f:0
\move(117 38)
\lvec(119 38)
\lvec(119 39)
\lvec(117 39)
\ifill f:0
\move(121 38)
\lvec(123 38)
\lvec(123 39)
\lvec(121 39)
\ifill f:0
\move(125 38)
\lvec(127 38)
\lvec(127 39)
\lvec(125 39)
\ifill f:0
\move(129 38)
\lvec(130 38)
\lvec(130 39)
\lvec(129 39)
\ifill f:0
\move(132 38)
\lvec(133 38)
\lvec(133 39)
\lvec(132 39)
\ifill f:0
\move(135 38)
\lvec(136 38)
\lvec(136 39)
\lvec(135 39)
\ifill f:0
\move(138 38)
\lvec(139 38)
\lvec(139 39)
\lvec(138 39)
\ifill f:0
\move(141 38)
\lvec(142 38)
\lvec(142 39)
\lvec(141 39)
\ifill f:0
\move(143 38)
\lvec(144 38)
\lvec(144 39)
\lvec(143 39)
\ifill f:0
\move(146 38)
\lvec(147 38)
\lvec(147 39)
\lvec(146 39)
\ifill f:0
\move(148 38)
\lvec(149 38)
\lvec(149 39)
\lvec(148 39)
\ifill f:0
\move(150 38)
\lvec(152 38)
\lvec(152 39)
\lvec(150 39)
\ifill f:0
\move(153 38)
\lvec(154 38)
\lvec(154 39)
\lvec(153 39)
\ifill f:0
\move(155 38)
\lvec(156 38)
\lvec(156 39)
\lvec(155 39)
\ifill f:0
\move(157 38)
\lvec(158 38)
\lvec(158 39)
\lvec(157 39)
\ifill f:0
\move(159 38)
\lvec(160 38)
\lvec(160 39)
\lvec(159 39)
\ifill f:0
\move(161 38)
\lvec(162 38)
\lvec(162 39)
\lvec(161 39)
\ifill f:0
\move(163 38)
\lvec(165 38)
\lvec(165 39)
\lvec(163 39)
\ifill f:0
\move(166 38)
\lvec(168 38)
\lvec(168 39)
\lvec(166 39)
\ifill f:0
\move(169 38)
\lvec(170 38)
\lvec(170 39)
\lvec(169 39)
\ifill f:0
\move(171 38)
\lvec(172 38)
\lvec(172 39)
\lvec(171 39)
\ifill f:0
\move(173 38)
\lvec(174 38)
\lvec(174 39)
\lvec(173 39)
\ifill f:0
\move(175 38)
\lvec(176 38)
\lvec(176 39)
\lvec(175 39)
\ifill f:0
\move(177 38)
\lvec(178 38)
\lvec(178 39)
\lvec(177 39)
\ifill f:0
\move(179 38)
\lvec(181 38)
\lvec(181 39)
\lvec(179 39)
\ifill f:0
\move(182 38)
\lvec(183 38)
\lvec(183 39)
\lvec(182 39)
\ifill f:0
\move(184 38)
\lvec(185 38)
\lvec(185 39)
\lvec(184 39)
\ifill f:0
\move(186 38)
\lvec(188 38)
\lvec(188 39)
\lvec(186 39)
\ifill f:0
\move(189 38)
\lvec(190 38)
\lvec(190 39)
\lvec(189 39)
\ifill f:0
\move(191 38)
\lvec(193 38)
\lvec(193 39)
\lvec(191 39)
\ifill f:0
\move(194 38)
\lvec(195 38)
\lvec(195 39)
\lvec(194 39)
\ifill f:0
\move(196 38)
\lvec(198 38)
\lvec(198 39)
\lvec(196 39)
\ifill f:0
\move(199 38)
\lvec(201 38)
\lvec(201 39)
\lvec(199 39)
\ifill f:0
\move(202 38)
\lvec(204 38)
\lvec(204 39)
\lvec(202 39)
\ifill f:0
\move(205 38)
\lvec(206 38)
\lvec(206 39)
\lvec(205 39)
\ifill f:0
\move(208 38)
\lvec(209 38)
\lvec(209 39)
\lvec(208 39)
\ifill f:0
\move(211 38)
\lvec(212 38)
\lvec(212 39)
\lvec(211 39)
\ifill f:0
\move(214 38)
\lvec(216 38)
\lvec(216 39)
\lvec(214 39)
\ifill f:0
\move(217 38)
\lvec(219 38)
\lvec(219 39)
\lvec(217 39)
\ifill f:0
\move(220 38)
\lvec(222 38)
\lvec(222 39)
\lvec(220 39)
\ifill f:0
\move(224 38)
\lvec(226 38)
\lvec(226 39)
\lvec(224 39)
\ifill f:0
\move(227 38)
\lvec(229 38)
\lvec(229 39)
\lvec(227 39)
\ifill f:0
\move(231 38)
\lvec(233 38)
\lvec(233 39)
\lvec(231 39)
\ifill f:0
\move(235 38)
\lvec(237 38)
\lvec(237 39)
\lvec(235 39)
\ifill f:0
\move(239 38)
\lvec(241 38)
\lvec(241 39)
\lvec(239 39)
\ifill f:0
\move(243 38)
\lvec(245 38)
\lvec(245 39)
\lvec(243 39)
\ifill f:0
\move(247 38)
\lvec(250 38)
\lvec(250 39)
\lvec(247 39)
\ifill f:0
\move(251 38)
\lvec(254 38)
\lvec(254 39)
\lvec(251 39)
\ifill f:0
\move(256 38)
\lvec(259 38)
\lvec(259 39)
\lvec(256 39)
\ifill f:0
\move(261 38)
\lvec(265 38)
\lvec(265 39)
\lvec(261 39)
\ifill f:0
\move(267 38)
\lvec(270 38)
\lvec(270 39)
\lvec(267 39)
\ifill f:0
\move(272 38)
\lvec(276 38)
\lvec(276 39)
\lvec(272 39)
\ifill f:0
\move(279 38)
\lvec(283 38)
\lvec(283 39)
\lvec(279 39)
\ifill f:0
\move(285 38)
\lvec(290 38)
\lvec(290 39)
\lvec(285 39)
\ifill f:0
\move(293 38)
\lvec(298 38)
\lvec(298 39)
\lvec(293 39)
\ifill f:0
\move(301 38)
\lvec(307 38)
\lvec(307 39)
\lvec(301 39)
\ifill f:0
\move(311 38)
\lvec(319 38)
\lvec(319 39)
\lvec(311 39)
\ifill f:0
\move(323 38)
\lvec(333 38)
\lvec(333 39)
\lvec(323 39)
\ifill f:0
\move(340 38)
\lvec(370 38)
\lvec(370 39)
\lvec(340 39)
\ifill f:0
\move(372 38)
\lvec(404 38)
\lvec(404 39)
\lvec(372 39)
\ifill f:0
\move(411 38)
\lvec(423 38)
\lvec(423 39)
\lvec(411 39)
\ifill f:0
\move(427 38)
\lvec(437 38)
\lvec(437 39)
\lvec(427 39)
\ifill f:0
\move(441 38)
\lvec(449 38)
\lvec(449 39)
\lvec(441 39)
\ifill f:0
\move(13 39)
\lvec(14 39)
\lvec(14 40)
\lvec(13 40)
\ifill f:0
\move(24 39)
\lvec(25 39)
\lvec(25 40)
\lvec(24 40)
\ifill f:0
\move(49 39)
\lvec(50 39)
\lvec(50 40)
\lvec(49 40)
\ifill f:0
\move(59 39)
\lvec(60 39)
\lvec(60 40)
\lvec(59 40)
\ifill f:0
\move(61 39)
\lvec(62 39)
\lvec(62 40)
\lvec(61 40)
\ifill f:0
\move(63 39)
\lvec(64 39)
\lvec(64 40)
\lvec(63 40)
\ifill f:0
\move(65 39)
\lvec(66 39)
\lvec(66 40)
\lvec(65 40)
\ifill f:0
\move(67 39)
\lvec(68 39)
\lvec(68 40)
\lvec(67 40)
\ifill f:0
\move(70 39)
\lvec(71 39)
\lvec(71 40)
\lvec(70 40)
\ifill f:0
\move(73 39)
\lvec(74 39)
\lvec(74 40)
\lvec(73 40)
\ifill f:0
\move(76 39)
\lvec(77 39)
\lvec(77 40)
\lvec(76 40)
\ifill f:0
\move(80 39)
\lvec(82 39)
\lvec(82 40)
\lvec(80 40)
\ifill f:0
\move(86 39)
\lvec(88 39)
\lvec(88 40)
\lvec(86 40)
\ifill f:0
\move(108 39)
\lvec(112 39)
\lvec(112 40)
\lvec(108 40)
\ifill f:0
\move(116 39)
\lvec(118 39)
\lvec(118 40)
\lvec(116 40)
\ifill f:0
\move(121 39)
\lvec(123 39)
\lvec(123 40)
\lvec(121 40)
\ifill f:0
\move(126 39)
\lvec(128 39)
\lvec(128 40)
\lvec(126 40)
\ifill f:0
\move(130 39)
\lvec(132 39)
\lvec(132 40)
\lvec(130 40)
\ifill f:0
\move(134 39)
\lvec(135 39)
\lvec(135 40)
\lvec(134 40)
\ifill f:0
\move(137 39)
\lvec(138 39)
\lvec(138 40)
\lvec(137 40)
\ifill f:0
\move(140 39)
\lvec(141 39)
\lvec(141 40)
\lvec(140 40)
\ifill f:0
\move(143 39)
\lvec(144 39)
\lvec(144 40)
\lvec(143 40)
\ifill f:0
\move(146 39)
\lvec(147 39)
\lvec(147 40)
\lvec(146 40)
\ifill f:0
\move(149 39)
\lvec(150 39)
\lvec(150 40)
\lvec(149 40)
\ifill f:0
\move(151 39)
\lvec(152 39)
\lvec(152 40)
\lvec(151 40)
\ifill f:0
\move(154 39)
\lvec(155 39)
\lvec(155 40)
\lvec(154 40)
\ifill f:0
\move(156 39)
\lvec(157 39)
\lvec(157 40)
\lvec(156 40)
\ifill f:0
\move(158 39)
\lvec(160 39)
\lvec(160 40)
\lvec(158 40)
\ifill f:0
\move(161 39)
\lvec(162 39)
\lvec(162 40)
\lvec(161 40)
\ifill f:0
\move(163 39)
\lvec(164 39)
\lvec(164 40)
\lvec(163 40)
\ifill f:0
\move(165 39)
\lvec(166 39)
\lvec(166 40)
\lvec(165 40)
\ifill f:0
\move(167 39)
\lvec(168 39)
\lvec(168 40)
\lvec(167 40)
\ifill f:0
\move(169 39)
\lvec(170 39)
\lvec(170 40)
\lvec(169 40)
\ifill f:0
\move(171 39)
\lvec(173 39)
\lvec(173 40)
\lvec(171 40)
\ifill f:0
\move(174 39)
\lvec(176 39)
\lvec(176 40)
\lvec(174 40)
\ifill f:0
\move(177 39)
\lvec(178 39)
\lvec(178 40)
\lvec(177 40)
\ifill f:0
\move(179 39)
\lvec(180 39)
\lvec(180 40)
\lvec(179 40)
\ifill f:0
\move(181 39)
\lvec(182 39)
\lvec(182 40)
\lvec(181 40)
\ifill f:0
\move(183 39)
\lvec(184 39)
\lvec(184 40)
\lvec(183 40)
\ifill f:0
\move(185 39)
\lvec(186 39)
\lvec(186 40)
\lvec(185 40)
\ifill f:0
\move(187 39)
\lvec(188 39)
\lvec(188 40)
\lvec(187 40)
\ifill f:0
\move(189 39)
\lvec(191 39)
\lvec(191 40)
\lvec(189 40)
\ifill f:0
\move(192 39)
\lvec(193 39)
\lvec(193 40)
\lvec(192 40)
\ifill f:0
\move(194 39)
\lvec(195 39)
\lvec(195 40)
\lvec(194 40)
\ifill f:0
\move(197 39)
\lvec(198 39)
\lvec(198 40)
\lvec(197 40)
\ifill f:0
\move(199 39)
\lvec(200 39)
\lvec(200 40)
\lvec(199 40)
\ifill f:0
\move(201 39)
\lvec(203 39)
\lvec(203 40)
\lvec(201 40)
\ifill f:0
\move(204 39)
\lvec(205 39)
\lvec(205 40)
\lvec(204 40)
\ifill f:0
\move(207 39)
\lvec(208 39)
\lvec(208 40)
\lvec(207 40)
\ifill f:0
\move(209 39)
\lvec(211 39)
\lvec(211 40)
\lvec(209 40)
\ifill f:0
\move(212 39)
\lvec(214 39)
\lvec(214 40)
\lvec(212 40)
\ifill f:0
\move(215 39)
\lvec(217 39)
\lvec(217 40)
\lvec(215 40)
\ifill f:0
\move(218 39)
\lvec(219 39)
\lvec(219 40)
\lvec(218 40)
\ifill f:0
\move(221 39)
\lvec(223 39)
\lvec(223 40)
\lvec(221 40)
\ifill f:0
\move(224 39)
\lvec(226 39)
\lvec(226 40)
\lvec(224 40)
\ifill f:0
\move(227 39)
\lvec(229 39)
\lvec(229 40)
\lvec(227 40)
\ifill f:0
\move(230 39)
\lvec(232 39)
\lvec(232 40)
\lvec(230 40)
\ifill f:0
\move(234 39)
\lvec(236 39)
\lvec(236 40)
\lvec(234 40)
\ifill f:0
\move(237 39)
\lvec(239 39)
\lvec(239 40)
\lvec(237 40)
\ifill f:0
\move(241 39)
\lvec(243 39)
\lvec(243 40)
\lvec(241 40)
\ifill f:0
\move(244 39)
\lvec(247 39)
\lvec(247 40)
\lvec(244 40)
\ifill f:0
\move(248 39)
\lvec(250 39)
\lvec(250 40)
\lvec(248 40)
\ifill f:0
\move(252 39)
\lvec(255 39)
\lvec(255 40)
\lvec(252 40)
\ifill f:0
\move(256 39)
\lvec(259 39)
\lvec(259 40)
\lvec(256 40)
\ifill f:0
\move(261 39)
\lvec(263 39)
\lvec(263 40)
\lvec(261 40)
\ifill f:0
\move(265 39)
\lvec(268 39)
\lvec(268 40)
\lvec(265 40)
\ifill f:0
\move(270 39)
\lvec(273 39)
\lvec(273 40)
\lvec(270 40)
\ifill f:0
\move(275 39)
\lvec(278 39)
\lvec(278 40)
\lvec(275 40)
\ifill f:0
\move(280 39)
\lvec(284 39)
\lvec(284 40)
\lvec(280 40)
\ifill f:0
\move(286 39)
\lvec(290 39)
\lvec(290 40)
\lvec(286 40)
\ifill f:0
\move(292 39)
\lvec(296 39)
\lvec(296 40)
\lvec(292 40)
\ifill f:0
\move(299 39)
\lvec(303 39)
\lvec(303 40)
\lvec(299 40)
\ifill f:0
\move(306 39)
\lvec(311 39)
\lvec(311 40)
\lvec(306 40)
\ifill f:0
\move(314 39)
\lvec(320 39)
\lvec(320 40)
\lvec(314 40)
\ifill f:0
\move(324 39)
\lvec(331 39)
\lvec(331 40)
\lvec(324 40)
\ifill f:0
\move(335 39)
\lvec(344 39)
\lvec(344 40)
\lvec(335 40)
\ifill f:0
\move(349 39)
\lvec(363 39)
\lvec(363 40)
\lvec(349 40)
\ifill f:0
\move(374 39)
\lvec(408 39)
\lvec(408 40)
\lvec(374 40)
\ifill f:0
\move(419 39)
\lvec(435 39)
\lvec(435 40)
\lvec(419 40)
\ifill f:0
\move(440 39)
\lvec(451 39)
\lvec(451 40)
\lvec(440 40)
\ifill f:0
\move(12 40)
\lvec(13 40)
\lvec(13 41)
\lvec(12 41)
\ifill f:0
\move(49 40)
\lvec(50 40)
\lvec(50 41)
\lvec(49 41)
\ifill f:0
\move(58 40)
\lvec(59 40)
\lvec(59 41)
\lvec(58 41)
\ifill f:0
\move(63 40)
\lvec(64 40)
\lvec(64 41)
\lvec(63 41)
\ifill f:0
\move(65 40)
\lvec(66 40)
\lvec(66 41)
\lvec(65 41)
\ifill f:0
\move(67 40)
\lvec(68 40)
\lvec(68 41)
\lvec(67 41)
\ifill f:0
\move(69 40)
\lvec(70 40)
\lvec(70 41)
\lvec(69 41)
\ifill f:0
\move(74 40)
\lvec(75 40)
\lvec(75 41)
\lvec(74 41)
\ifill f:0
\move(77 40)
\lvec(78 40)
\lvec(78 41)
\lvec(77 41)
\ifill f:0
\move(81 40)
\lvec(82 40)
\lvec(82 41)
\lvec(81 41)
\ifill f:0
\move(85 40)
\lvec(86 40)
\lvec(86 41)
\lvec(85 41)
\ifill f:0
\move(91 40)
\lvec(93 40)
\lvec(93 41)
\lvec(91 41)
\ifill f:0
\move(113 40)
\lvec(117 40)
\lvec(117 41)
\lvec(113 41)
\ifill f:0
\move(121 40)
\lvec(124 40)
\lvec(124 41)
\lvec(121 41)
\ifill f:0
\move(127 40)
\lvec(129 40)
\lvec(129 41)
\lvec(127 41)
\ifill f:0
\move(132 40)
\lvec(133 40)
\lvec(133 41)
\lvec(132 41)
\ifill f:0
\move(136 40)
\lvec(137 40)
\lvec(137 41)
\lvec(136 41)
\ifill f:0
\move(139 40)
\lvec(141 40)
\lvec(141 41)
\lvec(139 41)
\ifill f:0
\move(143 40)
\lvec(144 40)
\lvec(144 41)
\lvec(143 41)
\ifill f:0
\move(146 40)
\lvec(147 40)
\lvec(147 41)
\lvec(146 41)
\ifill f:0
\move(149 40)
\lvec(150 40)
\lvec(150 41)
\lvec(149 41)
\ifill f:0
\move(152 40)
\lvec(153 40)
\lvec(153 41)
\lvec(152 41)
\ifill f:0
\move(155 40)
\lvec(156 40)
\lvec(156 41)
\lvec(155 41)
\ifill f:0
\move(157 40)
\lvec(159 40)
\lvec(159 41)
\lvec(157 41)
\ifill f:0
\move(160 40)
\lvec(161 40)
\lvec(161 41)
\lvec(160 41)
\ifill f:0
\move(162 40)
\lvec(163 40)
\lvec(163 41)
\lvec(162 41)
\ifill f:0
\move(165 40)
\lvec(166 40)
\lvec(166 41)
\lvec(165 41)
\ifill f:0
\move(167 40)
\lvec(168 40)
\lvec(168 41)
\lvec(167 41)
\ifill f:0
\move(169 40)
\lvec(170 40)
\lvec(170 41)
\lvec(169 41)
\ifill f:0
\move(172 40)
\lvec(173 40)
\lvec(173 41)
\lvec(172 41)
\ifill f:0
\move(174 40)
\lvec(175 40)
\lvec(175 41)
\lvec(174 41)
\ifill f:0
\move(176 40)
\lvec(177 40)
\lvec(177 41)
\lvec(176 41)
\ifill f:0
\move(178 40)
\lvec(179 40)
\lvec(179 41)
\lvec(178 41)
\ifill f:0
\move(180 40)
\lvec(181 40)
\lvec(181 41)
\lvec(180 41)
\ifill f:0
\move(182 40)
\lvec(183 40)
\lvec(183 41)
\lvec(182 41)
\ifill f:0
\move(184 40)
\lvec(185 40)
\lvec(185 41)
\lvec(184 41)
\ifill f:0
\move(186 40)
\lvec(187 40)
\lvec(187 41)
\lvec(186 41)
\ifill f:0
\move(188 40)
\lvec(189 40)
\lvec(189 41)
\lvec(188 41)
\ifill f:0
\move(190 40)
\lvec(191 40)
\lvec(191 41)
\lvec(190 41)
\ifill f:0
\move(192 40)
\lvec(193 40)
\lvec(193 41)
\lvec(192 41)
\ifill f:0
\move(194 40)
\lvec(195 40)
\lvec(195 41)
\lvec(194 41)
\ifill f:0
\move(197 40)
\lvec(198 40)
\lvec(198 41)
\lvec(197 41)
\ifill f:0
\move(199 40)
\lvec(200 40)
\lvec(200 41)
\lvec(199 41)
\ifill f:0
\move(201 40)
\lvec(202 40)
\lvec(202 41)
\lvec(201 41)
\ifill f:0
\move(203 40)
\lvec(205 40)
\lvec(205 41)
\lvec(203 41)
\ifill f:0
\move(206 40)
\lvec(207 40)
\lvec(207 41)
\lvec(206 41)
\ifill f:0
\move(208 40)
\lvec(210 40)
\lvec(210 41)
\lvec(208 41)
\ifill f:0
\move(211 40)
\lvec(212 40)
\lvec(212 41)
\lvec(211 41)
\ifill f:0
\move(213 40)
\lvec(215 40)
\lvec(215 41)
\lvec(213 41)
\ifill f:0
\move(216 40)
\lvec(217 40)
\lvec(217 41)
\lvec(216 41)
\ifill f:0
\move(219 40)
\lvec(220 40)
\lvec(220 41)
\lvec(219 41)
\ifill f:0
\move(221 40)
\lvec(223 40)
\lvec(223 41)
\lvec(221 41)
\ifill f:0
\move(224 40)
\lvec(226 40)
\lvec(226 41)
\lvec(224 41)
\ifill f:0
\move(227 40)
\lvec(229 40)
\lvec(229 41)
\lvec(227 41)
\ifill f:0
\move(230 40)
\lvec(231 40)
\lvec(231 41)
\lvec(230 41)
\ifill f:0
\move(233 40)
\lvec(235 40)
\lvec(235 41)
\lvec(233 41)
\ifill f:0
\move(236 40)
\lvec(238 40)
\lvec(238 41)
\lvec(236 41)
\ifill f:0
\move(239 40)
\lvec(241 40)
\lvec(241 41)
\lvec(239 41)
\ifill f:0
\move(242 40)
\lvec(244 40)
\lvec(244 41)
\lvec(242 41)
\ifill f:0
\move(246 40)
\lvec(248 40)
\lvec(248 41)
\lvec(246 41)
\ifill f:0
\move(249 40)
\lvec(251 40)
\lvec(251 41)
\lvec(249 41)
\ifill f:0
\move(253 40)
\lvec(255 40)
\lvec(255 41)
\lvec(253 41)
\ifill f:0
\move(256 40)
\lvec(259 40)
\lvec(259 41)
\lvec(256 41)
\ifill f:0
\move(260 40)
\lvec(262 40)
\lvec(262 41)
\lvec(260 41)
\ifill f:0
\move(264 40)
\lvec(267 40)
\lvec(267 41)
\lvec(264 41)
\ifill f:0
\move(268 40)
\lvec(271 40)
\lvec(271 41)
\lvec(268 41)
\ifill f:0
\move(273 40)
\lvec(275 40)
\lvec(275 41)
\lvec(273 41)
\ifill f:0
\move(277 40)
\lvec(280 40)
\lvec(280 41)
\lvec(277 41)
\ifill f:0
\move(282 40)
\lvec(285 40)
\lvec(285 41)
\lvec(282 41)
\ifill f:0
\move(287 40)
\lvec(290 40)
\lvec(290 41)
\lvec(287 41)
\ifill f:0
\move(292 40)
\lvec(295 40)
\lvec(295 41)
\lvec(292 41)
\ifill f:0
\move(297 40)
\lvec(301 40)
\lvec(301 41)
\lvec(297 41)
\ifill f:0
\move(303 40)
\lvec(307 40)
\lvec(307 41)
\lvec(303 41)
\ifill f:0
\move(310 40)
\lvec(314 40)
\lvec(314 41)
\lvec(310 41)
\ifill f:0
\move(316 40)
\lvec(321 40)
\lvec(321 41)
\lvec(316 41)
\ifill f:0
\move(324 40)
\lvec(329 40)
\lvec(329 41)
\lvec(324 41)
\ifill f:0
\move(332 40)
\lvec(338 40)
\lvec(338 41)
\lvec(332 41)
\ifill f:0
\move(342 40)
\lvec(349 40)
\lvec(349 41)
\lvec(342 41)
\ifill f:0
\move(353 40)
\lvec(362 40)
\lvec(362 41)
\lvec(353 41)
\ifill f:0
\move(368 40)
\lvec(382 40)
\lvec(382 41)
\lvec(368 41)
\ifill f:0
\move(394 40)
\lvec(428 40)
\lvec(428 41)
\lvec(394 41)
\ifill f:0
\move(440 40)
\lvec(451 40)
\lvec(451 41)
\lvec(440 41)
\ifill f:0
\move(14 41)
\lvec(15 41)
\lvec(15 42)
\lvec(14 42)
\ifill f:0
\move(60 41)
\lvec(61 41)
\lvec(61 42)
\lvec(60 42)
\ifill f:0
\move(63 41)
\lvec(64 41)
\lvec(64 42)
\lvec(63 42)
\ifill f:0
\move(75 41)
\lvec(76 41)
\lvec(76 42)
\lvec(75 42)
\ifill f:0
\move(78 41)
\lvec(79 41)
\lvec(79 42)
\lvec(78 42)
\ifill f:0
\move(81 41)
\lvec(82 41)
\lvec(82 42)
\lvec(81 42)
\ifill f:0
\move(84 41)
\lvec(85 41)
\lvec(85 42)
\lvec(84 42)
\ifill f:0
\move(88 41)
\lvec(90 41)
\lvec(90 42)
\lvec(88 42)
\ifill f:0
\move(94 41)
\lvec(96 41)
\lvec(96 42)
\lvec(94 42)
\ifill f:0
\move(104 41)
\lvec(113 41)
\lvec(113 42)
\lvec(104 42)
\ifill f:0
\move(121 41)
\lvec(124 41)
\lvec(124 42)
\lvec(121 42)
\ifill f:0
\move(128 41)
\lvec(131 41)
\lvec(131 42)
\lvec(128 42)
\ifill f:0
\move(134 41)
\lvec(136 41)
\lvec(136 42)
\lvec(134 42)
\ifill f:0
\move(138 41)
\lvec(140 41)
\lvec(140 42)
\lvec(138 42)
\ifill f:0
\move(142 41)
\lvec(144 41)
\lvec(144 42)
\lvec(142 42)
\ifill f:0
\move(146 41)
\lvec(148 41)
\lvec(148 42)
\lvec(146 42)
\ifill f:0
\move(150 41)
\lvec(151 41)
\lvec(151 42)
\lvec(150 42)
\ifill f:0
\move(153 41)
\lvec(154 41)
\lvec(154 42)
\lvec(153 42)
\ifill f:0
\move(156 41)
\lvec(157 41)
\lvec(157 42)
\lvec(156 42)
\ifill f:0
\move(159 41)
\lvec(160 41)
\lvec(160 42)
\lvec(159 42)
\ifill f:0
\move(162 41)
\lvec(163 41)
\lvec(163 42)
\lvec(162 42)
\ifill f:0
\move(164 41)
\lvec(166 41)
\lvec(166 42)
\lvec(164 42)
\ifill f:0
\move(167 41)
\lvec(168 41)
\lvec(168 42)
\lvec(167 42)
\ifill f:0
\move(169 41)
\lvec(171 41)
\lvec(171 42)
\lvec(169 42)
\ifill f:0
\move(172 41)
\lvec(173 41)
\lvec(173 42)
\lvec(172 42)
\ifill f:0
\move(174 41)
\lvec(175 41)
\lvec(175 42)
\lvec(174 42)
\ifill f:0
\move(177 41)
\lvec(178 41)
\lvec(178 42)
\lvec(177 42)
\ifill f:0
\move(179 41)
\lvec(180 41)
\lvec(180 42)
\lvec(179 42)
\ifill f:0
\move(181 41)
\lvec(182 41)
\lvec(182 42)
\lvec(181 42)
\ifill f:0
\move(183 41)
\lvec(184 41)
\lvec(184 42)
\lvec(183 42)
\ifill f:0
\move(185 41)
\lvec(186 41)
\lvec(186 42)
\lvec(185 42)
\ifill f:0
\move(187 41)
\lvec(188 41)
\lvec(188 42)
\lvec(187 42)
\ifill f:0
\move(189 41)
\lvec(190 41)
\lvec(190 42)
\lvec(189 42)
\ifill f:0
\move(193 41)
\lvec(194 41)
\lvec(194 42)
\lvec(193 42)
\ifill f:0
\move(195 41)
\lvec(196 41)
\lvec(196 42)
\lvec(195 42)
\ifill f:0
\move(197 41)
\lvec(198 41)
\lvec(198 42)
\lvec(197 42)
\ifill f:0
\move(199 41)
\lvec(200 41)
\lvec(200 42)
\lvec(199 42)
\ifill f:0
\move(201 41)
\lvec(202 41)
\lvec(202 42)
\lvec(201 42)
\ifill f:0
\move(203 41)
\lvec(204 41)
\lvec(204 42)
\lvec(203 42)
\ifill f:0
\move(205 41)
\lvec(206 41)
\lvec(206 42)
\lvec(205 42)
\ifill f:0
\move(207 41)
\lvec(209 41)
\lvec(209 42)
\lvec(207 42)
\ifill f:0
\move(210 41)
\lvec(211 41)
\lvec(211 42)
\lvec(210 42)
\ifill f:0
\move(212 41)
\lvec(213 41)
\lvec(213 42)
\lvec(212 42)
\ifill f:0
\move(214 41)
\lvec(216 41)
\lvec(216 42)
\lvec(214 42)
\ifill f:0
\move(217 41)
\lvec(218 41)
\lvec(218 42)
\lvec(217 42)
\ifill f:0
\move(219 41)
\lvec(220 41)
\lvec(220 42)
\lvec(219 42)
\ifill f:0
\move(222 41)
\lvec(223 41)
\lvec(223 42)
\lvec(222 42)
\ifill f:0
\move(224 41)
\lvec(226 41)
\lvec(226 42)
\lvec(224 42)
\ifill f:0
\move(227 41)
\lvec(228 41)
\lvec(228 42)
\lvec(227 42)
\ifill f:0
\move(229 41)
\lvec(231 41)
\lvec(231 42)
\lvec(229 42)
\ifill f:0
\move(232 41)
\lvec(234 41)
\lvec(234 42)
\lvec(232 42)
\ifill f:0
\move(235 41)
\lvec(237 41)
\lvec(237 42)
\lvec(235 42)
\ifill f:0
\move(238 41)
\lvec(239 41)
\lvec(239 42)
\lvec(238 42)
\ifill f:0
\move(241 41)
\lvec(242 41)
\lvec(242 42)
\lvec(241 42)
\ifill f:0
\move(244 41)
\lvec(245 41)
\lvec(245 42)
\lvec(244 42)
\ifill f:0
\move(247 41)
\lvec(249 41)
\lvec(249 42)
\lvec(247 42)
\ifill f:0
\move(250 41)
\lvec(252 41)
\lvec(252 42)
\lvec(250 42)
\ifill f:0
\move(253 41)
\lvec(255 41)
\lvec(255 42)
\lvec(253 42)
\ifill f:0
\move(256 41)
\lvec(258 41)
\lvec(258 42)
\lvec(256 42)
\ifill f:0
\move(260 41)
\lvec(262 41)
\lvec(262 42)
\lvec(260 42)
\ifill f:0
\move(263 41)
\lvec(265 41)
\lvec(265 42)
\lvec(263 42)
\ifill f:0
\move(267 41)
\lvec(269 41)
\lvec(269 42)
\lvec(267 42)
\ifill f:0
\move(271 41)
\lvec(273 41)
\lvec(273 42)
\lvec(271 42)
\ifill f:0
\move(275 41)
\lvec(277 41)
\lvec(277 42)
\lvec(275 42)
\ifill f:0
\move(279 41)
\lvec(281 41)
\lvec(281 42)
\lvec(279 42)
\ifill f:0
\move(283 41)
\lvec(285 41)
\lvec(285 42)
\lvec(283 42)
\ifill f:0
\move(287 41)
\lvec(290 41)
\lvec(290 42)
\lvec(287 42)
\ifill f:0
\move(292 41)
\lvec(294 41)
\lvec(294 42)
\lvec(292 42)
\ifill f:0
\move(296 41)
\lvec(299 41)
\lvec(299 42)
\lvec(296 42)
\ifill f:0
\move(301 41)
\lvec(304 41)
\lvec(304 42)
\lvec(301 42)
\ifill f:0
\move(306 41)
\lvec(310 41)
\lvec(310 42)
\lvec(306 42)
\ifill f:0
\move(312 41)
\lvec(316 41)
\lvec(316 42)
\lvec(312 42)
\ifill f:0
\move(318 41)
\lvec(322 41)
\lvec(322 42)
\lvec(318 42)
\ifill f:0
\move(324 41)
\lvec(328 41)
\lvec(328 42)
\lvec(324 42)
\ifill f:0
\move(331 41)
\lvec(335 41)
\lvec(335 42)
\lvec(331 42)
\ifill f:0
\move(338 41)
\lvec(343 41)
\lvec(343 42)
\lvec(338 42)
\ifill f:0
\move(346 41)
\lvec(352 41)
\lvec(352 42)
\lvec(346 42)
\ifill f:0
\move(356 41)
\lvec(362 41)
\lvec(362 42)
\lvec(356 42)
\ifill f:0
\move(366 41)
\lvec(374 41)
\lvec(374 42)
\lvec(366 42)
\ifill f:0
\move(379 41)
\lvec(390 41)
\lvec(390 42)
\lvec(379 42)
\ifill f:0
\move(398 41)
\lvec(430 41)
\lvec(430 42)
\lvec(398 42)
\ifill f:0
\move(432 41)
\lvec(451 41)
\lvec(451 42)
\lvec(432 42)
\ifill f:0
\move(12 42)
\lvec(14 42)
\lvec(14 43)
\lvec(12 43)
\ifill f:0
\move(19 42)
\lvec(20 42)
\lvec(20 43)
\lvec(19 43)
\ifill f:0
\move(49 42)
\lvec(50 42)
\lvec(50 43)
\lvec(49 43)
\ifill f:0
\move(51 42)
\lvec(53 42)
\lvec(53 43)
\lvec(51 43)
\ifill f:0
\move(59 42)
\lvec(60 42)
\lvec(60 43)
\lvec(59 43)
\ifill f:0
\move(68 42)
\lvec(69 42)
\lvec(69 43)
\lvec(68 43)
\ifill f:0
\move(76 42)
\lvec(77 42)
\lvec(77 43)
\lvec(76 43)
\ifill f:0
\move(78 42)
\lvec(79 42)
\lvec(79 43)
\lvec(78 43)
\ifill f:0
\move(81 42)
\lvec(82 42)
\lvec(82 43)
\lvec(81 43)
\ifill f:0
\move(84 42)
\lvec(85 42)
\lvec(85 43)
\lvec(84 43)
\ifill f:0
\move(87 42)
\lvec(88 42)
\lvec(88 43)
\lvec(87 43)
\ifill f:0
\move(91 42)
\lvec(92 42)
\lvec(92 43)
\lvec(91 43)
\ifill f:0
\move(95 42)
\lvec(97 42)
\lvec(97 43)
\lvec(95 43)
\ifill f:0
\move(102 42)
\lvec(106 42)
\lvec(106 43)
\lvec(102 43)
\ifill f:0
\move(121 42)
\lvec(126 42)
\lvec(126 43)
\lvec(121 43)
\ifill f:0
\move(131 42)
\lvec(134 42)
\lvec(134 43)
\lvec(131 43)
\ifill f:0
\move(137 42)
\lvec(139 42)
\lvec(139 43)
\lvec(137 43)
\ifill f:0
\move(142 42)
\lvec(144 42)
\lvec(144 43)
\lvec(142 43)
\ifill f:0
\move(146 42)
\lvec(148 42)
\lvec(148 43)
\lvec(146 43)
\ifill f:0
\move(150 42)
\lvec(152 42)
\lvec(152 43)
\lvec(150 43)
\ifill f:0
\move(154 42)
\lvec(156 42)
\lvec(156 43)
\lvec(154 43)
\ifill f:0
\move(157 42)
\lvec(159 42)
\lvec(159 43)
\lvec(157 43)
\ifill f:0
\move(161 42)
\lvec(162 42)
\lvec(162 43)
\lvec(161 43)
\ifill f:0
\move(164 42)
\lvec(165 42)
\lvec(165 43)
\lvec(164 43)
\ifill f:0
\move(167 42)
\lvec(168 42)
\lvec(168 43)
\lvec(167 43)
\ifill f:0
\move(169 42)
\lvec(171 42)
\lvec(171 43)
\lvec(169 43)
\ifill f:0
\move(172 42)
\lvec(173 42)
\lvec(173 43)
\lvec(172 43)
\ifill f:0
\move(175 42)
\lvec(176 42)
\lvec(176 43)
\lvec(175 43)
\ifill f:0
\move(177 42)
\lvec(178 42)
\lvec(178 43)
\lvec(177 43)
\ifill f:0
\move(180 42)
\lvec(181 42)
\lvec(181 43)
\lvec(180 43)
\ifill f:0
\move(182 42)
\lvec(183 42)
\lvec(183 43)
\lvec(182 43)
\ifill f:0
\move(184 42)
\lvec(186 42)
\lvec(186 43)
\lvec(184 43)
\ifill f:0
\move(187 42)
\lvec(188 42)
\lvec(188 43)
\lvec(187 43)
\ifill f:0
\move(189 42)
\lvec(190 42)
\lvec(190 43)
\lvec(189 43)
\ifill f:0
\move(191 42)
\lvec(192 42)
\lvec(192 43)
\lvec(191 43)
\ifill f:0
\move(193 42)
\lvec(194 42)
\lvec(194 43)
\lvec(193 43)
\ifill f:0
\move(195 42)
\lvec(196 42)
\lvec(196 43)
\lvec(195 43)
\ifill f:0
\move(197 42)
\lvec(201 42)
\lvec(201 43)
\lvec(197 43)
\ifill f:0
\move(202 42)
\lvec(206 42)
\lvec(206 43)
\lvec(202 43)
\ifill f:0
\move(207 42)
\lvec(208 42)
\lvec(208 43)
\lvec(207 43)
\ifill f:0
\move(209 42)
\lvec(210 42)
\lvec(210 43)
\lvec(209 43)
\ifill f:0
\move(211 42)
\lvec(212 42)
\lvec(212 43)
\lvec(211 43)
\ifill f:0
\move(213 42)
\lvec(214 42)
\lvec(214 43)
\lvec(213 43)
\ifill f:0
\move(215 42)
\lvec(216 42)
\lvec(216 43)
\lvec(215 43)
\ifill f:0
\move(217 42)
\lvec(219 42)
\lvec(219 43)
\lvec(217 43)
\ifill f:0
\move(220 42)
\lvec(221 42)
\lvec(221 43)
\lvec(220 43)
\ifill f:0
\move(222 42)
\lvec(223 42)
\lvec(223 43)
\lvec(222 43)
\ifill f:0
\move(224 42)
\lvec(226 42)
\lvec(226 43)
\lvec(224 43)
\ifill f:0
\move(227 42)
\lvec(228 42)
\lvec(228 43)
\lvec(227 43)
\ifill f:0
\move(229 42)
\lvec(230 42)
\lvec(230 43)
\lvec(229 43)
\ifill f:0
\move(232 42)
\lvec(233 42)
\lvec(233 43)
\lvec(232 43)
\ifill f:0
\move(234 42)
\lvec(236 42)
\lvec(236 43)
\lvec(234 43)
\ifill f:0
\move(237 42)
\lvec(238 42)
\lvec(238 43)
\lvec(237 43)
\ifill f:0
\move(239 42)
\lvec(241 42)
\lvec(241 43)
\lvec(239 43)
\ifill f:0
\move(242 42)
\lvec(244 42)
\lvec(244 43)
\lvec(242 43)
\ifill f:0
\move(245 42)
\lvec(246 42)
\lvec(246 43)
\lvec(245 43)
\ifill f:0
\move(248 42)
\lvec(249 42)
\lvec(249 43)
\lvec(248 43)
\ifill f:0
\move(250 42)
\lvec(252 42)
\lvec(252 43)
\lvec(250 43)
\ifill f:0
\move(253 42)
\lvec(255 42)
\lvec(255 43)
\lvec(253 43)
\ifill f:0
\move(256 42)
\lvec(258 42)
\lvec(258 43)
\lvec(256 43)
\ifill f:0
\move(259 42)
\lvec(261 42)
\lvec(261 43)
\lvec(259 43)
\ifill f:0
\move(263 42)
\lvec(265 42)
\lvec(265 43)
\lvec(263 43)
\ifill f:0
\move(266 42)
\lvec(268 42)
\lvec(268 43)
\lvec(266 43)
\ifill f:0
\move(269 42)
\lvec(271 42)
\lvec(271 43)
\lvec(269 43)
\ifill f:0
\move(273 42)
\lvec(275 42)
\lvec(275 43)
\lvec(273 43)
\ifill f:0
\move(276 42)
\lvec(278 42)
\lvec(278 43)
\lvec(276 43)
\ifill f:0
\move(280 42)
\lvec(282 42)
\lvec(282 43)
\lvec(280 43)
\ifill f:0
\move(283 42)
\lvec(286 42)
\lvec(286 43)
\lvec(283 43)
\ifill f:0
\move(287 42)
\lvec(290 42)
\lvec(290 43)
\lvec(287 43)
\ifill f:0
\move(291 42)
\lvec(294 42)
\lvec(294 43)
\lvec(291 43)
\ifill f:0
\move(295 42)
\lvec(298 42)
\lvec(298 43)
\lvec(295 43)
\ifill f:0
\move(300 42)
\lvec(302 42)
\lvec(302 43)
\lvec(300 43)
\ifill f:0
\move(304 42)
\lvec(307 42)
\lvec(307 43)
\lvec(304 43)
\ifill f:0
\move(309 42)
\lvec(312 42)
\lvec(312 43)
\lvec(309 43)
\ifill f:0
\move(314 42)
\lvec(317 42)
\lvec(317 43)
\lvec(314 43)
\ifill f:0
\move(319 42)
\lvec(322 42)
\lvec(322 43)
\lvec(319 43)
\ifill f:0
\move(324 42)
\lvec(328 42)
\lvec(328 43)
\lvec(324 43)
\ifill f:0
\move(330 42)
\lvec(334 42)
\lvec(334 43)
\lvec(330 43)
\ifill f:0
\move(336 42)
\lvec(340 42)
\lvec(340 43)
\lvec(336 43)
\ifill f:0
\move(342 42)
\lvec(347 42)
\lvec(347 43)
\lvec(342 43)
\ifill f:0
\move(349 42)
\lvec(354 42)
\lvec(354 43)
\lvec(349 43)
\ifill f:0
\move(357 42)
\lvec(362 42)
\lvec(362 43)
\lvec(357 43)
\ifill f:0
\move(365 42)
\lvec(371 42)
\lvec(371 43)
\lvec(365 43)
\ifill f:0
\move(374 42)
\lvec(381 42)
\lvec(381 43)
\lvec(374 43)
\ifill f:0
\move(385 42)
\lvec(394 42)
\lvec(394 43)
\lvec(385 43)
\ifill f:0
\move(399 42)
\lvec(410 42)
\lvec(410 43)
\lvec(399 43)
\ifill f:0
\move(418 42)
\lvec(451 42)
\lvec(451 43)
\lvec(418 43)
\ifill f:0
\move(15 43)
\lvec(16 43)
\lvec(16 44)
\lvec(15 44)
\ifill f:0
\move(62 43)
\lvec(63 43)
\lvec(63 44)
\lvec(62 44)
\ifill f:0
\move(66 43)
\lvec(67 43)
\lvec(67 44)
\lvec(66 44)
\ifill f:0
\move(81 43)
\lvec(82 43)
\lvec(82 44)
\lvec(81 44)
\ifill f:0
\move(83 43)
\lvec(84 43)
\lvec(84 44)
\lvec(83 44)
\ifill f:0
\move(86 43)
\lvec(87 43)
\lvec(87 44)
\lvec(86 44)
\ifill f:0
\move(89 43)
\lvec(90 43)
\lvec(90 44)
\lvec(89 44)
\ifill f:0
\move(93 43)
\lvec(94 43)
\lvec(94 44)
\lvec(93 44)
\ifill f:0
\move(97 43)
\lvec(98 43)
\lvec(98 44)
\lvec(97 44)
\ifill f:0
\move(102 43)
\lvec(104 43)
\lvec(104 44)
\lvec(102 44)
\ifill f:0
\move(110 43)
\lvec(118 43)
\lvec(118 44)
\lvec(110 44)
\ifill f:0
\move(119 43)
\lvec(128 43)
\lvec(128 44)
\lvec(119 44)
\ifill f:0
\move(135 43)
\lvec(138 43)
\lvec(138 44)
\lvec(135 44)
\ifill f:0
\move(142 43)
\lvec(144 43)
\lvec(144 44)
\lvec(142 44)
\ifill f:0
\move(147 43)
\lvec(149 43)
\lvec(149 44)
\lvec(147 44)
\ifill f:0
\move(152 43)
\lvec(153 43)
\lvec(153 44)
\lvec(152 44)
\ifill f:0
\move(156 43)
\lvec(157 43)
\lvec(157 44)
\lvec(156 44)
\ifill f:0
\move(159 43)
\lvec(161 43)
\lvec(161 44)
\lvec(159 44)
\ifill f:0
\move(163 43)
\lvec(164 43)
\lvec(164 44)
\lvec(163 44)
\ifill f:0
\move(166 43)
\lvec(168 43)
\lvec(168 44)
\lvec(166 44)
\ifill f:0
\move(169 43)
\lvec(171 43)
\lvec(171 44)
\lvec(169 44)
\ifill f:0
\move(172 43)
\lvec(174 43)
\lvec(174 44)
\lvec(172 44)
\ifill f:0
\move(175 43)
\lvec(177 43)
\lvec(177 44)
\lvec(175 44)
\ifill f:0
\move(178 43)
\lvec(179 43)
\lvec(179 44)
\lvec(178 44)
\ifill f:0
\move(181 43)
\lvec(182 43)
\lvec(182 44)
\lvec(181 44)
\ifill f:0
\move(183 43)
\lvec(185 43)
\lvec(185 44)
\lvec(183 44)
\ifill f:0
\move(186 43)
\lvec(187 43)
\lvec(187 44)
\lvec(186 44)
\ifill f:0
\move(188 43)
\lvec(189 43)
\lvec(189 44)
\lvec(188 44)
\ifill f:0
\move(191 43)
\lvec(192 43)
\lvec(192 44)
\lvec(191 44)
\ifill f:0
\move(193 43)
\lvec(194 43)
\lvec(194 44)
\lvec(193 44)
\ifill f:0
\move(195 43)
\lvec(196 43)
\lvec(196 44)
\lvec(195 44)
\ifill f:0
\move(198 43)
\lvec(199 43)
\lvec(199 44)
\lvec(198 44)
\ifill f:0
\move(200 43)
\lvec(201 43)
\lvec(201 44)
\lvec(200 44)
\ifill f:0
\move(202 43)
\lvec(203 43)
\lvec(203 44)
\lvec(202 44)
\ifill f:0
\move(204 43)
\lvec(205 43)
\lvec(205 44)
\lvec(204 44)
\ifill f:0
\move(206 43)
\lvec(207 43)
\lvec(207 44)
\lvec(206 44)
\ifill f:0
\move(208 43)
\lvec(209 43)
\lvec(209 44)
\lvec(208 44)
\ifill f:0
\move(210 43)
\lvec(211 43)
\lvec(211 44)
\lvec(210 44)
\ifill f:0
\move(212 43)
\lvec(213 43)
\lvec(213 44)
\lvec(212 44)
\ifill f:0
\move(214 43)
\lvec(215 43)
\lvec(215 44)
\lvec(214 44)
\ifill f:0
\move(216 43)
\lvec(217 43)
\lvec(217 44)
\lvec(216 44)
\ifill f:0
\move(218 43)
\lvec(219 43)
\lvec(219 44)
\lvec(218 44)
\ifill f:0
\move(220 43)
\lvec(221 43)
\lvec(221 44)
\lvec(220 44)
\ifill f:0
\move(222 43)
\lvec(223 43)
\lvec(223 44)
\lvec(222 44)
\ifill f:0
\move(224 43)
\lvec(226 43)
\lvec(226 44)
\lvec(224 44)
\ifill f:0
\move(227 43)
\lvec(228 43)
\lvec(228 44)
\lvec(227 44)
\ifill f:0
\move(229 43)
\lvec(230 43)
\lvec(230 44)
\lvec(229 44)
\ifill f:0
\move(231 43)
\lvec(232 43)
\lvec(232 44)
\lvec(231 44)
\ifill f:0
\move(234 43)
\lvec(235 43)
\lvec(235 44)
\lvec(234 44)
\ifill f:0
\move(236 43)
\lvec(237 43)
\lvec(237 44)
\lvec(236 44)
\ifill f:0
\move(238 43)
\lvec(240 43)
\lvec(240 44)
\lvec(238 44)
\ifill f:0
\move(241 43)
\lvec(242 43)
\lvec(242 44)
\lvec(241 44)
\ifill f:0
\move(243 43)
\lvec(245 43)
\lvec(245 44)
\lvec(243 44)
\ifill f:0
\move(246 43)
\lvec(247 43)
\lvec(247 44)
\lvec(246 44)
\ifill f:0
\move(248 43)
\lvec(250 43)
\lvec(250 44)
\lvec(248 44)
\ifill f:0
\move(251 43)
\lvec(253 43)
\lvec(253 44)
\lvec(251 44)
\ifill f:0
\move(254 43)
\lvec(255 43)
\lvec(255 44)
\lvec(254 44)
\ifill f:0
\move(256 43)
\lvec(258 43)
\lvec(258 44)
\lvec(256 44)
\ifill f:0
\move(259 43)
\lvec(261 43)
\lvec(261 44)
\lvec(259 44)
\ifill f:0
\move(262 43)
\lvec(264 43)
\lvec(264 44)
\lvec(262 44)
\ifill f:0
\move(265 43)
\lvec(267 43)
\lvec(267 44)
\lvec(265 44)
\ifill f:0
\move(268 43)
\lvec(270 43)
\lvec(270 44)
\lvec(268 44)
\ifill f:0
\move(271 43)
\lvec(273 43)
\lvec(273 44)
\lvec(271 44)
\ifill f:0
\move(274 43)
\lvec(276 43)
\lvec(276 44)
\lvec(274 44)
\ifill f:0
\move(277 43)
\lvec(279 43)
\lvec(279 44)
\lvec(277 44)
\ifill f:0
\move(281 43)
\lvec(283 43)
\lvec(283 44)
\lvec(281 44)
\ifill f:0
\move(284 43)
\lvec(286 43)
\lvec(286 44)
\lvec(284 44)
\ifill f:0
\move(288 43)
\lvec(290 43)
\lvec(290 44)
\lvec(288 44)
\ifill f:0
\move(291 43)
\lvec(293 43)
\lvec(293 44)
\lvec(291 44)
\ifill f:0
\move(295 43)
\lvec(297 43)
\lvec(297 44)
\lvec(295 44)
\ifill f:0
\move(299 43)
\lvec(301 43)
\lvec(301 44)
\lvec(299 44)
\ifill f:0
\move(303 43)
\lvec(305 43)
\lvec(305 44)
\lvec(303 44)
\ifill f:0
\move(307 43)
\lvec(309 43)
\lvec(309 44)
\lvec(307 44)
\ifill f:0
\move(311 43)
\lvec(313 43)
\lvec(313 44)
\lvec(311 44)
\ifill f:0
\move(315 43)
\lvec(318 43)
\lvec(318 44)
\lvec(315 44)
\ifill f:0
\move(320 43)
\lvec(322 43)
\lvec(322 44)
\lvec(320 44)
\ifill f:0
\move(324 43)
\lvec(327 43)
\lvec(327 44)
\lvec(324 44)
\ifill f:0
\move(329 43)
\lvec(332 43)
\lvec(332 44)
\lvec(329 44)
\ifill f:0
\move(334 43)
\lvec(337 43)
\lvec(337 44)
\lvec(334 44)
\ifill f:0
\move(340 43)
\lvec(343 43)
\lvec(343 44)
\lvec(340 44)
\ifill f:0
\move(345 43)
\lvec(349 43)
\lvec(349 44)
\lvec(345 44)
\ifill f:0
\move(351 43)
\lvec(355 43)
\lvec(355 44)
\lvec(351 44)
\ifill f:0
\move(358 43)
\lvec(362 43)
\lvec(362 44)
\lvec(358 44)
\ifill f:0
\move(364 43)
\lvec(369 43)
\lvec(369 44)
\lvec(364 44)
\ifill f:0
\move(372 43)
\lvec(377 43)
\lvec(377 44)
\lvec(372 44)
\ifill f:0
\move(380 43)
\lvec(386 43)
\lvec(386 44)
\lvec(380 44)
\ifill f:0
\move(389 43)
\lvec(396 43)
\lvec(396 44)
\lvec(389 44)
\ifill f:0
\move(399 43)
\lvec(407 43)
\lvec(407 44)
\lvec(399 44)
\ifill f:0
\move(412 43)
\lvec(422 43)
\lvec(422 44)
\lvec(412 44)
\ifill f:0
\move(428 43)
\lvec(443 43)
\lvec(443 44)
\lvec(428 44)
\ifill f:0
\move(12 44)
\lvec(16 44)
\lvec(16 45)
\lvec(12 45)
\ifill f:0
\move(66 44)
\lvec(67 44)
\lvec(67 45)
\lvec(66 45)
\ifill f:0
\move(75 44)
\lvec(76 44)
\lvec(76 45)
\lvec(75 45)
\ifill f:0
\move(77 44)
\lvec(78 44)
\lvec(78 45)
\lvec(77 45)
\ifill f:0
\move(81 44)
\lvec(82 44)
\lvec(82 45)
\lvec(81 45)
\ifill f:0
\move(83 44)
\lvec(84 44)
\lvec(84 45)
\lvec(83 45)
\ifill f:0
\move(88 44)
\lvec(89 44)
\lvec(89 45)
\lvec(88 45)
\ifill f:0
\move(91 44)
\lvec(92 44)
\lvec(92 45)
\lvec(91 45)
\ifill f:0
\move(94 44)
\lvec(95 44)
\lvec(95 45)
\lvec(94 45)
\ifill f:0
\move(102 44)
\lvec(103 44)
\lvec(103 45)
\lvec(102 45)
\ifill f:0
\move(107 44)
\lvec(109 44)
\lvec(109 45)
\lvec(107 45)
\ifill f:0
\move(116 44)
\lvec(124 44)
\lvec(124 45)
\lvec(116 45)
\ifill f:0
\move(125 44)
\lvec(134 44)
\lvec(134 45)
\lvec(125 45)
\ifill f:0
\move(140 44)
\lvec(143 44)
\lvec(143 45)
\lvec(140 45)
\ifill f:0
\move(148 44)
\lvec(150 44)
\lvec(150 45)
\lvec(148 45)
\ifill f:0
\move(153 44)
\lvec(155 44)
\lvec(155 45)
\lvec(153 45)
\ifill f:0
\move(158 44)
\lvec(160 44)
\lvec(160 45)
\lvec(158 45)
\ifill f:0
\move(162 44)
\lvec(164 44)
\lvec(164 45)
\lvec(162 45)
\ifill f:0
\move(166 44)
\lvec(167 44)
\lvec(167 45)
\lvec(166 45)
\ifill f:0
\move(169 44)
\lvec(171 44)
\lvec(171 45)
\lvec(169 45)
\ifill f:0
\move(173 44)
\lvec(174 44)
\lvec(174 45)
\lvec(173 45)
\ifill f:0
\move(176 44)
\lvec(177 44)
\lvec(177 45)
\lvec(176 45)
\ifill f:0
\move(179 44)
\lvec(180 44)
\lvec(180 45)
\lvec(179 45)
\ifill f:0
\move(182 44)
\lvec(183 44)
\lvec(183 45)
\lvec(182 45)
\ifill f:0
\move(185 44)
\lvec(186 44)
\lvec(186 45)
\lvec(185 45)
\ifill f:0
\move(187 44)
\lvec(189 44)
\lvec(189 45)
\lvec(187 45)
\ifill f:0
\move(190 44)
\lvec(191 44)
\lvec(191 45)
\lvec(190 45)
\ifill f:0
\move(193 44)
\lvec(194 44)
\lvec(194 45)
\lvec(193 45)
\ifill f:0
\move(195 44)
\lvec(196 44)
\lvec(196 45)
\lvec(195 45)
\ifill f:0
\move(198 44)
\lvec(199 44)
\lvec(199 45)
\lvec(198 45)
\ifill f:0
\move(200 44)
\lvec(201 44)
\lvec(201 45)
\lvec(200 45)
\ifill f:0
\move(202 44)
\lvec(203 44)
\lvec(203 45)
\lvec(202 45)
\ifill f:0
\move(205 44)
\lvec(206 44)
\lvec(206 45)
\lvec(205 45)
\ifill f:0
\move(207 44)
\lvec(208 44)
\lvec(208 45)
\lvec(207 45)
\ifill f:0
\move(209 44)
\lvec(210 44)
\lvec(210 45)
\lvec(209 45)
\ifill f:0
\move(211 44)
\lvec(212 44)
\lvec(212 45)
\lvec(211 45)
\ifill f:0
\move(213 44)
\lvec(214 44)
\lvec(214 45)
\lvec(213 45)
\ifill f:0
\move(215 44)
\lvec(216 44)
\lvec(216 45)
\lvec(215 45)
\ifill f:0
\move(217 44)
\lvec(218 44)
\lvec(218 45)
\lvec(217 45)
\ifill f:0
\move(219 44)
\lvec(222 44)
\lvec(222 45)
\lvec(219 45)
\ifill f:0
\move(223 44)
\lvec(224 44)
\lvec(224 45)
\lvec(223 45)
\ifill f:0
\move(225 44)
\lvec(226 44)
\lvec(226 45)
\lvec(225 45)
\ifill f:0
\move(227 44)
\lvec(228 44)
\lvec(228 45)
\lvec(227 45)
\ifill f:0
\move(229 44)
\lvec(230 44)
\lvec(230 45)
\lvec(229 45)
\ifill f:0
\move(231 44)
\lvec(232 44)
\lvec(232 45)
\lvec(231 45)
\ifill f:0
\move(233 44)
\lvec(234 44)
\lvec(234 45)
\lvec(233 45)
\ifill f:0
\move(235 44)
\lvec(236 44)
\lvec(236 45)
\lvec(235 45)
\ifill f:0
\move(237 44)
\lvec(239 44)
\lvec(239 45)
\lvec(237 45)
\ifill f:0
\move(240 44)
\lvec(241 44)
\lvec(241 45)
\lvec(240 45)
\ifill f:0
\move(242 44)
\lvec(243 44)
\lvec(243 45)
\lvec(242 45)
\ifill f:0
\move(244 44)
\lvec(246 44)
\lvec(246 45)
\lvec(244 45)
\ifill f:0
\move(247 44)
\lvec(248 44)
\lvec(248 45)
\lvec(247 45)
\ifill f:0
\move(249 44)
\lvec(250 44)
\lvec(250 45)
\lvec(249 45)
\ifill f:0
\move(251 44)
\lvec(253 44)
\lvec(253 45)
\lvec(251 45)
\ifill f:0
\move(254 44)
\lvec(255 44)
\lvec(255 45)
\lvec(254 45)
\ifill f:0
\move(257 44)
\lvec(258 44)
\lvec(258 45)
\lvec(257 45)
\ifill f:0
\move(259 44)
\lvec(261 44)
\lvec(261 45)
\lvec(259 45)
\ifill f:0
\move(262 44)
\lvec(263 44)
\lvec(263 45)
\lvec(262 45)
\ifill f:0
\move(264 44)
\lvec(266 44)
\lvec(266 45)
\lvec(264 45)
\ifill f:0
\move(267 44)
\lvec(269 44)
\lvec(269 45)
\lvec(267 45)
\ifill f:0
\move(270 44)
\lvec(271 44)
\lvec(271 45)
\lvec(270 45)
\ifill f:0
\move(273 44)
\lvec(274 44)
\lvec(274 45)
\lvec(273 45)
\ifill f:0
\move(276 44)
\lvec(277 44)
\lvec(277 45)
\lvec(276 45)
\ifill f:0
\move(279 44)
\lvec(280 44)
\lvec(280 45)
\lvec(279 45)
\ifill f:0
\move(282 44)
\lvec(283 44)
\lvec(283 45)
\lvec(282 45)
\ifill f:0
\move(285 44)
\lvec(286 44)
\lvec(286 45)
\lvec(285 45)
\ifill f:0
\move(288 44)
\lvec(290 44)
\lvec(290 45)
\lvec(288 45)
\ifill f:0
\move(291 44)
\lvec(293 44)
\lvec(293 45)
\lvec(291 45)
\ifill f:0
\move(294 44)
\lvec(296 44)
\lvec(296 45)
\lvec(294 45)
\ifill f:0
\move(298 44)
\lvec(300 44)
\lvec(300 45)
\lvec(298 45)
\ifill f:0
\move(301 44)
\lvec(303 44)
\lvec(303 45)
\lvec(301 45)
\ifill f:0
\move(305 44)
\lvec(307 44)
\lvec(307 45)
\lvec(305 45)
\ifill f:0
\move(308 44)
\lvec(311 44)
\lvec(311 45)
\lvec(308 45)
\ifill f:0
\move(312 44)
\lvec(315 44)
\lvec(315 45)
\lvec(312 45)
\ifill f:0
\move(316 44)
\lvec(319 44)
\lvec(319 45)
\lvec(316 45)
\ifill f:0
\move(320 44)
\lvec(323 44)
\lvec(323 45)
\lvec(320 45)
\ifill f:0
\move(324 44)
\lvec(327 44)
\lvec(327 45)
\lvec(324 45)
\ifill f:0
\move(329 44)
\lvec(331 44)
\lvec(331 45)
\lvec(329 45)
\ifill f:0
\move(333 44)
\lvec(336 44)
\lvec(336 45)
\lvec(333 45)
\ifill f:0
\move(338 44)
\lvec(341 44)
\lvec(341 45)
\lvec(338 45)
\ifill f:0
\move(343 44)
\lvec(346 44)
\lvec(346 45)
\lvec(343 45)
\ifill f:0
\move(348 44)
\lvec(351 44)
\lvec(351 45)
\lvec(348 45)
\ifill f:0
\move(353 44)
\lvec(356 44)
\lvec(356 45)
\lvec(353 45)
\ifill f:0
\move(358 44)
\lvec(362 44)
\lvec(362 45)
\lvec(358 45)
\ifill f:0
\move(364 44)
\lvec(368 44)
\lvec(368 45)
\lvec(364 45)
\ifill f:0
\move(370 44)
\lvec(374 44)
\lvec(374 45)
\lvec(370 45)
\ifill f:0
\move(377 44)
\lvec(381 44)
\lvec(381 45)
\lvec(377 45)
\ifill f:0
\move(384 44)
\lvec(389 44)
\lvec(389 45)
\lvec(384 45)
\ifill f:0
\move(391 44)
\lvec(397 44)
\lvec(397 45)
\lvec(391 45)
\ifill f:0
\move(400 44)
\lvec(406 44)
\lvec(406 45)
\lvec(400 45)
\ifill f:0
\move(409 44)
\lvec(416 44)
\lvec(416 45)
\lvec(409 45)
\ifill f:0
\move(420 44)
\lvec(428 44)
\lvec(428 45)
\lvec(420 45)
\ifill f:0
\move(432 44)
\lvec(442 44)
\lvec(442 45)
\lvec(432 45)
\ifill f:0
\move(449 44)
\lvec(451 44)
\lvec(451 45)
\lvec(449 45)
\ifill f:0
\move(15 45)
\lvec(16 45)
\lvec(16 46)
\lvec(15 46)
\ifill f:0
\move(19 45)
\lvec(20 45)
\lvec(20 46)
\lvec(19 46)
\ifill f:0
\move(22 45)
\lvec(23 45)
\lvec(23 46)
\lvec(22 46)
\ifill f:0
\move(56 45)
\lvec(60 45)
\lvec(60 46)
\lvec(56 46)
\ifill f:0
\move(67 45)
\lvec(68 45)
\lvec(68 46)
\lvec(67 46)
\ifill f:0
\move(74 45)
\lvec(75 45)
\lvec(75 46)
\lvec(74 46)
\ifill f:0
\move(79 45)
\lvec(80 45)
\lvec(80 46)
\lvec(79 46)
\ifill f:0
\move(81 45)
\lvec(82 45)
\lvec(82 46)
\lvec(81 46)
\ifill f:0
\move(85 45)
\lvec(86 45)
\lvec(86 46)
\lvec(85 46)
\ifill f:0
\move(87 45)
\lvec(88 45)
\lvec(88 46)
\lvec(87 46)
\ifill f:0
\move(92 45)
\lvec(93 45)
\lvec(93 46)
\lvec(92 46)
\ifill f:0
\move(95 45)
\lvec(96 45)
\lvec(96 46)
\lvec(95 46)
\ifill f:0
\move(98 45)
\lvec(99 45)
\lvec(99 46)
\lvec(98 46)
\ifill f:0
\move(101 45)
\lvec(102 45)
\lvec(102 46)
\lvec(101 46)
\ifill f:0
\move(105 45)
\lvec(107 45)
\lvec(107 46)
\lvec(105 46)
\ifill f:0
\move(110 45)
\lvec(113 45)
\lvec(113 46)
\lvec(110 46)
\ifill f:0
\move(118 45)
\lvec(122 45)
\lvec(122 46)
\lvec(118 46)
\ifill f:0
\move(138 45)
\lvec(143 45)
\lvec(143 46)
\lvec(138 46)
\ifill f:0
\move(148 45)
\lvec(152 45)
\lvec(152 46)
\lvec(148 46)
\ifill f:0
\move(155 45)
\lvec(158 45)
\lvec(158 46)
\lvec(155 46)
\ifill f:0
\move(161 45)
\lvec(163 45)
\lvec(163 46)
\lvec(161 46)
\ifill f:0
\move(165 45)
\lvec(167 45)
\lvec(167 46)
\lvec(165 46)
\ifill f:0
\move(169 45)
\lvec(171 45)
\lvec(171 46)
\lvec(169 46)
\ifill f:0
\move(173 45)
\lvec(175 45)
\lvec(175 46)
\lvec(173 46)
\ifill f:0
\move(177 45)
\lvec(178 45)
\lvec(178 46)
\lvec(177 46)
\ifill f:0
\move(180 45)
\lvec(182 45)
\lvec(182 46)
\lvec(180 46)
\ifill f:0
\move(183 45)
\lvec(185 45)
\lvec(185 46)
\lvec(183 46)
\ifill f:0
\move(187 45)
\lvec(188 45)
\lvec(188 46)
\lvec(187 46)
\ifill f:0
\move(189 45)
\lvec(191 45)
\lvec(191 46)
\lvec(189 46)
\ifill f:0
\move(192 45)
\lvec(194 45)
\lvec(194 46)
\lvec(192 46)
\ifill f:0
\move(195 45)
\lvec(196 45)
\lvec(196 46)
\lvec(195 46)
\ifill f:0
\move(198 45)
\lvec(199 45)
\lvec(199 46)
\lvec(198 46)
\ifill f:0
\move(200 45)
\lvec(202 45)
\lvec(202 46)
\lvec(200 46)
\ifill f:0
\move(203 45)
\lvec(204 45)
\lvec(204 46)
\lvec(203 46)
\ifill f:0
\move(205 45)
\lvec(207 45)
\lvec(207 46)
\lvec(205 46)
\ifill f:0
\move(208 45)
\lvec(209 45)
\lvec(209 46)
\lvec(208 46)
\ifill f:0
\move(210 45)
\lvec(211 45)
\lvec(211 46)
\lvec(210 46)
\ifill f:0
\move(212 45)
\lvec(214 45)
\lvec(214 46)
\lvec(212 46)
\ifill f:0
\move(215 45)
\lvec(216 45)
\lvec(216 46)
\lvec(215 46)
\ifill f:0
\move(217 45)
\lvec(218 45)
\lvec(218 46)
\lvec(217 46)
\ifill f:0
\move(219 45)
\lvec(220 45)
\lvec(220 46)
\lvec(219 46)
\ifill f:0
\move(221 45)
\lvec(222 45)
\lvec(222 46)
\lvec(221 46)
\ifill f:0
\move(223 45)
\lvec(224 45)
\lvec(224 46)
\lvec(223 46)
\ifill f:0
\move(225 45)
\lvec(236 45)
\lvec(236 46)
\lvec(225 46)
\ifill f:0
\move(237 45)
\lvec(238 45)
\lvec(238 46)
\lvec(237 46)
\ifill f:0
\move(239 45)
\lvec(240 45)
\lvec(240 46)
\lvec(239 46)
\ifill f:0
\move(241 45)
\lvec(242 45)
\lvec(242 46)
\lvec(241 46)
\ifill f:0
\move(243 45)
\lvec(244 45)
\lvec(244 46)
\lvec(243 46)
\ifill f:0
\move(245 45)
\lvec(246 45)
\lvec(246 46)
\lvec(245 46)
\ifill f:0
\move(247 45)
\lvec(249 45)
\lvec(249 46)
\lvec(247 46)
\ifill f:0
\move(250 45)
\lvec(251 45)
\lvec(251 46)
\lvec(250 46)
\ifill f:0
\move(252 45)
\lvec(253 45)
\lvec(253 46)
\lvec(252 46)
\ifill f:0
\move(254 45)
\lvec(255 45)
\lvec(255 46)
\lvec(254 46)
\ifill f:0
\move(256 45)
\lvec(258 45)
\lvec(258 46)
\lvec(256 46)
\ifill f:0
\move(259 45)
\lvec(260 45)
\lvec(260 46)
\lvec(259 46)
\ifill f:0
\move(261 45)
\lvec(263 45)
\lvec(263 46)
\lvec(261 46)
\ifill f:0
\move(264 45)
\lvec(265 45)
\lvec(265 46)
\lvec(264 46)
\ifill f:0
\move(266 45)
\lvec(268 45)
\lvec(268 46)
\lvec(266 46)
\ifill f:0
\move(269 45)
\lvec(270 45)
\lvec(270 46)
\lvec(269 46)
\ifill f:0
\move(271 45)
\lvec(273 45)
\lvec(273 46)
\lvec(271 46)
\ifill f:0
\move(274 45)
\lvec(276 45)
\lvec(276 46)
\lvec(274 46)
\ifill f:0
\move(277 45)
\lvec(278 45)
\lvec(278 46)
\lvec(277 46)
\ifill f:0
\move(279 45)
\lvec(281 45)
\lvec(281 46)
\lvec(279 46)
\ifill f:0
\move(282 45)
\lvec(284 45)
\lvec(284 46)
\lvec(282 46)
\ifill f:0
\move(285 45)
\lvec(287 45)
\lvec(287 46)
\lvec(285 46)
\ifill f:0
\move(288 45)
\lvec(290 45)
\lvec(290 46)
\lvec(288 46)
\ifill f:0
\move(291 45)
\lvec(293 45)
\lvec(293 46)
\lvec(291 46)
\ifill f:0
\move(294 45)
\lvec(296 45)
\lvec(296 46)
\lvec(294 46)
\ifill f:0
\move(297 45)
\lvec(299 45)
\lvec(299 46)
\lvec(297 46)
\ifill f:0
\move(300 45)
\lvec(302 45)
\lvec(302 46)
\lvec(300 46)
\ifill f:0
\move(303 45)
\lvec(305 45)
\lvec(305 46)
\lvec(303 46)
\ifill f:0
\move(307 45)
\lvec(309 45)
\lvec(309 46)
\lvec(307 46)
\ifill f:0
\move(310 45)
\lvec(312 45)
\lvec(312 46)
\lvec(310 46)
\ifill f:0
\move(313 45)
\lvec(316 45)
\lvec(316 46)
\lvec(313 46)
\ifill f:0
\move(317 45)
\lvec(319 45)
\lvec(319 46)
\lvec(317 46)
\ifill f:0
\move(321 45)
\lvec(323 45)
\lvec(323 46)
\lvec(321 46)
\ifill f:0
\move(324 45)
\lvec(327 45)
\lvec(327 46)
\lvec(324 46)
\ifill f:0
\move(328 45)
\lvec(331 45)
\lvec(331 46)
\lvec(328 46)
\ifill f:0
\move(332 45)
\lvec(335 45)
\lvec(335 46)
\lvec(332 46)
\ifill f:0
\move(336 45)
\lvec(339 45)
\lvec(339 46)
\lvec(336 46)
\ifill f:0
\move(340 45)
\lvec(343 45)
\lvec(343 46)
\lvec(340 46)
\ifill f:0
\move(345 45)
\lvec(347 45)
\lvec(347 46)
\lvec(345 46)
\ifill f:0
\move(349 45)
\lvec(352 45)
\lvec(352 46)
\lvec(349 46)
\ifill f:0
\move(354 45)
\lvec(357 45)
\lvec(357 46)
\lvec(354 46)
\ifill f:0
\move(359 45)
\lvec(362 45)
\lvec(362 46)
\lvec(359 46)
\ifill f:0
\move(364 45)
\lvec(367 45)
\lvec(367 46)
\lvec(364 46)
\ifill f:0
\move(369 45)
\lvec(372 45)
\lvec(372 46)
\lvec(369 46)
\ifill f:0
\move(374 45)
\lvec(378 45)
\lvec(378 46)
\lvec(374 46)
\ifill f:0
\move(380 45)
\lvec(384 45)
\lvec(384 46)
\lvec(380 46)
\ifill f:0
\move(386 45)
\lvec(391 45)
\lvec(391 46)
\lvec(386 46)
\ifill f:0
\move(393 45)
\lvec(397 45)
\lvec(397 46)
\lvec(393 46)
\ifill f:0
\move(400 45)
\lvec(405 45)
\lvec(405 46)
\lvec(400 46)
\ifill f:0
\move(407 45)
\lvec(413 45)
\lvec(413 46)
\lvec(407 46)
\ifill f:0
\move(415 45)
\lvec(421 45)
\lvec(421 46)
\lvec(415 46)
\ifill f:0
\move(424 45)
\lvec(431 45)
\lvec(431 46)
\lvec(424 46)
\ifill f:0
\move(435 45)
\lvec(442 45)
\lvec(442 46)
\lvec(435 46)
\ifill f:0
\move(446 45)
\lvec(451 45)
\lvec(451 46)
\lvec(446 46)
\ifill f:0
\move(12 46)
\lvec(13 46)
\lvec(13 47)
\lvec(12 47)
\ifill f:0
\move(15 46)
\lvec(16 46)
\lvec(16 47)
\lvec(15 47)
\ifill f:0
\move(68 46)
\lvec(69 46)
\lvec(69 47)
\lvec(68 47)
\ifill f:0
\move(81 46)
\lvec(82 46)
\lvec(82 47)
\lvec(81 47)
\ifill f:0
\move(93 46)
\lvec(94 46)
\lvec(94 47)
\lvec(93 47)
\ifill f:0
\move(98 46)
\lvec(99 46)
\lvec(99 47)
\lvec(98 47)
\ifill f:0
\move(101 46)
\lvec(102 46)
\lvec(102 47)
\lvec(101 47)
\ifill f:0
\move(105 46)
\lvec(106 46)
\lvec(106 47)
\lvec(105 47)
\ifill f:0
\move(109 46)
\lvec(110 46)
\lvec(110 47)
\lvec(109 47)
\ifill f:0
\move(113 46)
\lvec(115 46)
\lvec(115 47)
\lvec(113 47)
\ifill f:0
\move(119 46)
\lvec(122 46)
\lvec(122 47)
\lvec(119 47)
\ifill f:0
\move(131 46)
\lvec(141 46)
\lvec(141 47)
\lvec(131 47)
\ifill f:0
\move(150 46)
\lvec(154 46)
\lvec(154 47)
\lvec(150 47)
\ifill f:0
\move(158 46)
\lvec(161 46)
\lvec(161 47)
\lvec(158 47)
\ifill f:0
\move(164 46)
\lvec(166 46)
\lvec(166 47)
\lvec(164 47)
\ifill f:0
\move(169 46)
\lvec(171 46)
\lvec(171 47)
\lvec(169 47)
\ifill f:0
\move(174 46)
\lvec(176 46)
\lvec(176 47)
\lvec(174 47)
\ifill f:0
\move(178 46)
\lvec(180 46)
\lvec(180 47)
\lvec(178 47)
\ifill f:0
\move(182 46)
\lvec(183 46)
\lvec(183 47)
\lvec(182 47)
\ifill f:0
\move(185 46)
\lvec(187 46)
\lvec(187 47)
\lvec(185 47)
\ifill f:0
\move(189 46)
\lvec(190 46)
\lvec(190 47)
\lvec(189 47)
\ifill f:0
\move(192 46)
\lvec(193 46)
\lvec(193 47)
\lvec(192 47)
\ifill f:0
\move(195 46)
\lvec(196 46)
\lvec(196 47)
\lvec(195 47)
\ifill f:0
\move(198 46)
\lvec(199 46)
\lvec(199 47)
\lvec(198 47)
\ifill f:0
\move(201 46)
\lvec(202 46)
\lvec(202 47)
\lvec(201 47)
\ifill f:0
\move(203 46)
\lvec(205 46)
\lvec(205 47)
\lvec(203 47)
\ifill f:0
\move(206 46)
\lvec(207 46)
\lvec(207 47)
\lvec(206 47)
\ifill f:0
\move(209 46)
\lvec(210 46)
\lvec(210 47)
\lvec(209 47)
\ifill f:0
\move(211 46)
\lvec(212 46)
\lvec(212 47)
\lvec(211 47)
\ifill f:0
\move(214 46)
\lvec(215 46)
\lvec(215 47)
\lvec(214 47)
\ifill f:0
\move(216 46)
\lvec(217 46)
\lvec(217 47)
\lvec(216 47)
\ifill f:0
\move(219 46)
\lvec(220 46)
\lvec(220 47)
\lvec(219 47)
\ifill f:0
\move(221 46)
\lvec(222 46)
\lvec(222 47)
\lvec(221 47)
\ifill f:0
\move(223 46)
\lvec(224 46)
\lvec(224 47)
\lvec(223 47)
\ifill f:0
\move(225 46)
\lvec(226 46)
\lvec(226 47)
\lvec(225 47)
\ifill f:0
\move(228 46)
\lvec(229 46)
\lvec(229 47)
\lvec(228 47)
\ifill f:0
\move(230 46)
\lvec(231 46)
\lvec(231 47)
\lvec(230 47)
\ifill f:0
\move(232 46)
\lvec(233 46)
\lvec(233 47)
\lvec(232 47)
\ifill f:0
\move(234 46)
\lvec(235 46)
\lvec(235 47)
\lvec(234 47)
\ifill f:0
\move(236 46)
\lvec(237 46)
\lvec(237 47)
\lvec(236 47)
\ifill f:0
\move(238 46)
\lvec(239 46)
\lvec(239 47)
\lvec(238 47)
\ifill f:0
\move(240 46)
\lvec(241 46)
\lvec(241 47)
\lvec(240 47)
\ifill f:0
\move(242 46)
\lvec(243 46)
\lvec(243 47)
\lvec(242 47)
\ifill f:0
\move(244 46)
\lvec(245 46)
\lvec(245 47)
\lvec(244 47)
\ifill f:0
\move(246 46)
\lvec(247 46)
\lvec(247 47)
\lvec(246 47)
\ifill f:0
\move(248 46)
\lvec(249 46)
\lvec(249 47)
\lvec(248 47)
\ifill f:0
\move(250 46)
\lvec(251 46)
\lvec(251 47)
\lvec(250 47)
\ifill f:0
\move(252 46)
\lvec(253 46)
\lvec(253 47)
\lvec(252 47)
\ifill f:0
\move(254 46)
\lvec(255 46)
\lvec(255 47)
\lvec(254 47)
\ifill f:0
\move(257 46)
\lvec(258 46)
\lvec(258 47)
\lvec(257 47)
\ifill f:0
\move(259 46)
\lvec(260 46)
\lvec(260 47)
\lvec(259 47)
\ifill f:0
\move(261 46)
\lvec(262 46)
\lvec(262 47)
\lvec(261 47)
\ifill f:0
\move(263 46)
\lvec(265 46)
\lvec(265 47)
\lvec(263 47)
\ifill f:0
\move(266 46)
\lvec(267 46)
\lvec(267 47)
\lvec(266 47)
\ifill f:0
\move(268 46)
\lvec(269 46)
\lvec(269 47)
\lvec(268 47)
\ifill f:0
\move(270 46)
\lvec(272 46)
\lvec(272 47)
\lvec(270 47)
\ifill f:0
\move(273 46)
\lvec(274 46)
\lvec(274 47)
\lvec(273 47)
\ifill f:0
\move(275 46)
\lvec(277 46)
\lvec(277 47)
\lvec(275 47)
\ifill f:0
\move(278 46)
\lvec(279 46)
\lvec(279 47)
\lvec(278 47)
\ifill f:0
\move(280 46)
\lvec(282 46)
\lvec(282 47)
\lvec(280 47)
\ifill f:0
\move(283 46)
\lvec(284 46)
\lvec(284 47)
\lvec(283 47)
\ifill f:0
\move(285 46)
\lvec(287 46)
\lvec(287 47)
\lvec(285 47)
\ifill f:0
\move(288 46)
\lvec(290 46)
\lvec(290 47)
\lvec(288 47)
\ifill f:0
\move(291 46)
\lvec(292 46)
\lvec(292 47)
\lvec(291 47)
\ifill f:0
\move(294 46)
\lvec(295 46)
\lvec(295 47)
\lvec(294 47)
\ifill f:0
\move(296 46)
\lvec(298 46)
\lvec(298 47)
\lvec(296 47)
\ifill f:0
\move(299 46)
\lvec(301 46)
\lvec(301 47)
\lvec(299 47)
\ifill f:0
\move(302 46)
\lvec(304 46)
\lvec(304 47)
\lvec(302 47)
\ifill f:0
\move(305 46)
\lvec(307 46)
\lvec(307 47)
\lvec(305 47)
\ifill f:0
\move(308 46)
\lvec(310 46)
\lvec(310 47)
\lvec(308 47)
\ifill f:0
\move(311 46)
\lvec(313 46)
\lvec(313 47)
\lvec(311 47)
\ifill f:0
\move(315 46)
\lvec(316 46)
\lvec(316 47)
\lvec(315 47)
\ifill f:0
\move(318 46)
\lvec(320 46)
\lvec(320 47)
\lvec(318 47)
\ifill f:0
\move(321 46)
\lvec(323 46)
\lvec(323 47)
\lvec(321 47)
\ifill f:0
\move(324 46)
\lvec(326 46)
\lvec(326 47)
\lvec(324 47)
\ifill f:0
\move(328 46)
\lvec(330 46)
\lvec(330 47)
\lvec(328 47)
\ifill f:0
\move(331 46)
\lvec(334 46)
\lvec(334 47)
\lvec(331 47)
\ifill f:0
\move(335 46)
\lvec(337 46)
\lvec(337 47)
\lvec(335 47)
\ifill f:0
\move(339 46)
\lvec(341 46)
\lvec(341 47)
\lvec(339 47)
\ifill f:0
\move(343 46)
\lvec(345 46)
\lvec(345 47)
\lvec(343 47)
\ifill f:0
\move(347 46)
\lvec(349 46)
\lvec(349 47)
\lvec(347 47)
\ifill f:0
\move(351 46)
\lvec(353 46)
\lvec(353 47)
\lvec(351 47)
\ifill f:0
\move(355 46)
\lvec(357 46)
\lvec(357 47)
\lvec(355 47)
\ifill f:0
\move(359 46)
\lvec(362 46)
\lvec(362 47)
\lvec(359 47)
\ifill f:0
\move(364 46)
\lvec(366 46)
\lvec(366 47)
\lvec(364 47)
\ifill f:0
\move(368 46)
\lvec(371 46)
\lvec(371 47)
\lvec(368 47)
\ifill f:0
\move(373 46)
\lvec(376 46)
\lvec(376 47)
\lvec(373 47)
\ifill f:0
\move(378 46)
\lvec(381 46)
\lvec(381 47)
\lvec(378 47)
\ifill f:0
\move(383 46)
\lvec(386 46)
\lvec(386 47)
\lvec(383 47)
\ifill f:0
\move(389 46)
\lvec(392 46)
\lvec(392 47)
\lvec(389 47)
\ifill f:0
\move(394 46)
\lvec(398 46)
\lvec(398 47)
\lvec(394 47)
\ifill f:0
\move(400 46)
\lvec(404 46)
\lvec(404 47)
\lvec(400 47)
\ifill f:0
\move(406 46)
\lvec(411 46)
\lvec(411 47)
\lvec(406 47)
\ifill f:0
\move(413 46)
\lvec(417 46)
\lvec(417 47)
\lvec(413 47)
\ifill f:0
\move(420 46)
\lvec(425 46)
\lvec(425 47)
\lvec(420 47)
\ifill f:0
\move(428 46)
\lvec(433 46)
\lvec(433 47)
\lvec(428 47)
\ifill f:0
\move(436 46)
\lvec(442 46)
\lvec(442 47)
\lvec(436 47)
\ifill f:0
\move(446 46)
\lvec(451 46)
\lvec(451 47)
\lvec(446 47)
\ifill f:0
\move(13 47)
\lvec(14 47)
\lvec(14 48)
\lvec(13 48)
\ifill f:0
\move(78 47)
\lvec(79 47)
\lvec(79 48)
\lvec(78 48)
\ifill f:0
\move(81 47)
\lvec(82 47)
\lvec(82 48)
\lvec(81 48)
\ifill f:0
\move(86 47)
\lvec(87 47)
\lvec(87 48)
\lvec(86 48)
\ifill f:0
\move(96 47)
\lvec(97 47)
\lvec(97 48)
\lvec(96 48)
\ifill f:0
\move(101 47)
\lvec(102 47)
\lvec(102 48)
\lvec(101 48)
\ifill f:0
\move(104 47)
\lvec(105 47)
\lvec(105 48)
\lvec(104 48)
\ifill f:0
\move(107 47)
\lvec(108 47)
\lvec(108 48)
\lvec(107 48)
\ifill f:0
\move(111 47)
\lvec(112 47)
\lvec(112 48)
\lvec(111 48)
\ifill f:0
\move(115 47)
\lvec(116 47)
\lvec(116 48)
\lvec(115 48)
\ifill f:0
\move(120 47)
\lvec(122 47)
\lvec(122 48)
\lvec(120 48)
\ifill f:0
\move(127 47)
\lvec(130 47)
\lvec(130 48)
\lvec(127 48)
\ifill f:0
\move(154 47)
\lvec(158 47)
\lvec(158 48)
\lvec(154 48)
\ifill f:0
\move(163 47)
\lvec(165 47)
\lvec(165 48)
\lvec(163 48)
\ifill f:0
\move(169 47)
\lvec(172 47)
\lvec(172 48)
\lvec(169 48)
\ifill f:0
\move(175 47)
\lvec(177 47)
\lvec(177 48)
\lvec(175 48)
\ifill f:0
\move(179 47)
\lvec(181 47)
\lvec(181 48)
\lvec(179 48)
\ifill f:0
\move(184 47)
\lvec(185 47)
\lvec(185 48)
\lvec(184 48)
\ifill f:0
\move(188 47)
\lvec(189 47)
\lvec(189 48)
\lvec(188 48)
\ifill f:0
\move(191 47)
\lvec(193 47)
\lvec(193 48)
\lvec(191 48)
\ifill f:0
\move(195 47)
\lvec(196 47)
\lvec(196 48)
\lvec(195 48)
\ifill f:0
\move(198 47)
\lvec(199 47)
\lvec(199 48)
\lvec(198 48)
\ifill f:0
\move(201 47)
\lvec(203 47)
\lvec(203 48)
\lvec(201 48)
\ifill f:0
\move(204 47)
\lvec(206 47)
\lvec(206 48)
\lvec(204 48)
\ifill f:0
\move(207 47)
\lvec(208 47)
\lvec(208 48)
\lvec(207 48)
\ifill f:0
\move(210 47)
\lvec(211 47)
\lvec(211 48)
\lvec(210 48)
\ifill f:0
\move(213 47)
\lvec(214 47)
\lvec(214 48)
\lvec(213 48)
\ifill f:0
\move(215 47)
\lvec(217 47)
\lvec(217 48)
\lvec(215 48)
\ifill f:0
\move(218 47)
\lvec(219 47)
\lvec(219 48)
\lvec(218 48)
\ifill f:0
\move(221 47)
\lvec(222 47)
\lvec(222 48)
\lvec(221 48)
\ifill f:0
\move(223 47)
\lvec(224 47)
\lvec(224 48)
\lvec(223 48)
\ifill f:0
\move(225 47)
\lvec(227 47)
\lvec(227 48)
\lvec(225 48)
\ifill f:0
\move(228 47)
\lvec(229 47)
\lvec(229 48)
\lvec(228 48)
\ifill f:0
\move(230 47)
\lvec(231 47)
\lvec(231 48)
\lvec(230 48)
\ifill f:0
\move(232 47)
\lvec(234 47)
\lvec(234 48)
\lvec(232 48)
\ifill f:0
\move(235 47)
\lvec(236 47)
\lvec(236 48)
\lvec(235 48)
\ifill f:0
\move(237 47)
\lvec(238 47)
\lvec(238 48)
\lvec(237 48)
\ifill f:0
\move(239 47)
\lvec(240 47)
\lvec(240 48)
\lvec(239 48)
\ifill f:0
\move(241 47)
\lvec(242 47)
\lvec(242 48)
\lvec(241 48)
\ifill f:0
\move(243 47)
\lvec(244 47)
\lvec(244 48)
\lvec(243 48)
\ifill f:0
\move(245 47)
\lvec(246 47)
\lvec(246 48)
\lvec(245 48)
\ifill f:0
\move(247 47)
\lvec(248 47)
\lvec(248 48)
\lvec(247 48)
\ifill f:0
\move(249 47)
\lvec(251 47)
\lvec(251 48)
\lvec(249 48)
\ifill f:0
\move(252 47)
\lvec(254 47)
\lvec(254 48)
\lvec(252 48)
\ifill f:0
\move(255 47)
\lvec(256 47)
\lvec(256 48)
\lvec(255 48)
\ifill f:0
\move(257 47)
\lvec(258 47)
\lvec(258 48)
\lvec(257 48)
\ifill f:0
\move(259 47)
\lvec(260 47)
\lvec(260 48)
\lvec(259 48)
\ifill f:0
\move(261 47)
\lvec(262 47)
\lvec(262 48)
\lvec(261 48)
\ifill f:0
\move(263 47)
\lvec(264 47)
\lvec(264 48)
\lvec(263 48)
\ifill f:0
\move(265 47)
\lvec(266 47)
\lvec(266 48)
\lvec(265 48)
\ifill f:0
\move(267 47)
\lvec(268 47)
\lvec(268 48)
\lvec(267 48)
\ifill f:0
\move(269 47)
\lvec(271 47)
\lvec(271 48)
\lvec(269 48)
\ifill f:0
\move(272 47)
\lvec(273 47)
\lvec(273 48)
\lvec(272 48)
\ifill f:0
\move(274 47)
\lvec(275 47)
\lvec(275 48)
\lvec(274 48)
\ifill f:0
\move(276 47)
\lvec(277 47)
\lvec(277 48)
\lvec(276 48)
\ifill f:0
\move(279 47)
\lvec(280 47)
\lvec(280 48)
\lvec(279 48)
\ifill f:0
\move(281 47)
\lvec(282 47)
\lvec(282 48)
\lvec(281 48)
\ifill f:0
\move(283 47)
\lvec(285 47)
\lvec(285 48)
\lvec(283 48)
\ifill f:0
\move(286 47)
\lvec(287 47)
\lvec(287 48)
\lvec(286 48)
\ifill f:0
\move(288 47)
\lvec(290 47)
\lvec(290 48)
\lvec(288 48)
\ifill f:0
\move(291 47)
\lvec(292 47)
\lvec(292 48)
\lvec(291 48)
\ifill f:0
\move(293 47)
\lvec(295 47)
\lvec(295 48)
\lvec(293 48)
\ifill f:0
\move(296 47)
\lvec(297 47)
\lvec(297 48)
\lvec(296 48)
\ifill f:0
\move(299 47)
\lvec(300 47)
\lvec(300 48)
\lvec(299 48)
\ifill f:0
\move(301 47)
\lvec(303 47)
\lvec(303 48)
\lvec(301 48)
\ifill f:0
\move(304 47)
\lvec(305 47)
\lvec(305 48)
\lvec(304 48)
\ifill f:0
\move(307 47)
\lvec(308 47)
\lvec(308 48)
\lvec(307 48)
\ifill f:0
\move(310 47)
\lvec(311 47)
\lvec(311 48)
\lvec(310 48)
\ifill f:0
\move(312 47)
\lvec(314 47)
\lvec(314 48)
\lvec(312 48)
\ifill f:0
\move(315 47)
\lvec(317 47)
\lvec(317 48)
\lvec(315 48)
\ifill f:0
\move(318 47)
\lvec(320 47)
\lvec(320 48)
\lvec(318 48)
\ifill f:0
\move(321 47)
\lvec(323 47)
\lvec(323 48)
\lvec(321 48)
\ifill f:0
\move(324 47)
\lvec(326 47)
\lvec(326 48)
\lvec(324 48)
\ifill f:0
\move(328 47)
\lvec(329 47)
\lvec(329 48)
\lvec(328 48)
\ifill f:0
\move(331 47)
\lvec(333 47)
\lvec(333 48)
\lvec(331 48)
\ifill f:0
\move(334 47)
\lvec(336 47)
\lvec(336 48)
\lvec(334 48)
\ifill f:0
\move(337 47)
\lvec(339 47)
\lvec(339 48)
\lvec(337 48)
\ifill f:0
\move(341 47)
\lvec(343 47)
\lvec(343 48)
\lvec(341 48)
\ifill f:0
\move(344 47)
\lvec(346 47)
\lvec(346 48)
\lvec(344 48)
\ifill f:0
\move(348 47)
\lvec(350 47)
\lvec(350 48)
\lvec(348 48)
\ifill f:0
\move(352 47)
\lvec(354 47)
\lvec(354 48)
\lvec(352 48)
\ifill f:0
\move(355 47)
\lvec(358 47)
\lvec(358 48)
\lvec(355 48)
\ifill f:0
\move(359 47)
\lvec(362 47)
\lvec(362 48)
\lvec(359 48)
\ifill f:0
\move(363 47)
\lvec(366 47)
\lvec(366 48)
\lvec(363 48)
\ifill f:0
\move(367 47)
\lvec(370 47)
\lvec(370 48)
\lvec(367 48)
\ifill f:0
\move(372 47)
\lvec(374 47)
\lvec(374 48)
\lvec(372 48)
\ifill f:0
\move(376 47)
\lvec(379 47)
\lvec(379 48)
\lvec(376 48)
\ifill f:0
\move(381 47)
\lvec(383 47)
\lvec(383 48)
\lvec(381 48)
\ifill f:0
\move(385 47)
\lvec(388 47)
\lvec(388 48)
\lvec(385 48)
\ifill f:0
\move(390 47)
\lvec(393 47)
\lvec(393 48)
\lvec(390 48)
\ifill f:0
\move(395 47)
\lvec(398 47)
\lvec(398 48)
\lvec(395 48)
\ifill f:0
\move(400 47)
\lvec(403 47)
\lvec(403 48)
\lvec(400 48)
\ifill f:0
\move(406 47)
\lvec(409 47)
\lvec(409 48)
\lvec(406 48)
\ifill f:0
\move(411 47)
\lvec(415 47)
\lvec(415 48)
\lvec(411 48)
\ifill f:0
\move(417 47)
\lvec(421 47)
\lvec(421 48)
\lvec(417 48)
\ifill f:0
\move(424 47)
\lvec(428 47)
\lvec(428 48)
\lvec(424 48)
\ifill f:0
\move(430 47)
\lvec(435 47)
\lvec(435 48)
\lvec(430 48)
\ifill f:0
\move(437 47)
\lvec(442 47)
\lvec(442 48)
\lvec(437 48)
\ifill f:0
\move(445 47)
\lvec(450 47)
\lvec(450 48)
\lvec(445 48)
\ifill f:0
\move(12 48)
\lvec(13 48)
\lvec(13 49)
\lvec(12 49)
\ifill f:0
\move(14 48)
\lvec(15 48)
\lvec(15 49)
\lvec(14 49)
\ifill f:0
\move(19 48)
\lvec(20 48)
\lvec(20 49)
\lvec(19 49)
\ifill f:0
\move(22 48)
\lvec(23 48)
\lvec(23 49)
\lvec(22 49)
\ifill f:0
\move(81 48)
\lvec(82 48)
\lvec(82 49)
\lvec(81 49)
\ifill f:0
\move(84 48)
\lvec(85 48)
\lvec(85 49)
\lvec(84 49)
\ifill f:0
\move(101 48)
\lvec(102 48)
\lvec(102 49)
\lvec(101 49)
\ifill f:0
\move(106 48)
\lvec(107 48)
\lvec(107 49)
\lvec(106 49)
\ifill f:0
\move(109 48)
\lvec(110 48)
\lvec(110 49)
\lvec(109 49)
\ifill f:0
\move(116 48)
\lvec(117 48)
\lvec(117 49)
\lvec(116 49)
\ifill f:0
\move(120 48)
\lvec(122 48)
\lvec(122 49)
\lvec(120 49)
\ifill f:0
\move(126 48)
\lvec(127 48)
\lvec(127 49)
\lvec(126 49)
\ifill f:0
\move(133 48)
\lvec(136 48)
\lvec(136 49)
\lvec(133 49)
\ifill f:0
\move(160 48)
\lvec(164 48)
\lvec(164 49)
\lvec(160 49)
\ifill f:0
\move(169 48)
\lvec(172 48)
\lvec(172 49)
\lvec(169 49)
\ifill f:0
\move(176 48)
\lvec(178 48)
\lvec(178 49)
\lvec(176 49)
\ifill f:0
\move(181 48)
\lvec(183 48)
\lvec(183 49)
\lvec(181 49)
\ifill f:0
\move(186 48)
\lvec(188 48)
\lvec(188 49)
\lvec(186 49)
\ifill f:0
\move(190 48)
\lvec(192 48)
\lvec(192 49)
\lvec(190 49)
\ifill f:0
\move(194 48)
\lvec(196 48)
\lvec(196 49)
\lvec(194 49)
\ifill f:0
\move(198 48)
\lvec(200 48)
\lvec(200 49)
\lvec(198 49)
\ifill f:0
\move(202 48)
\lvec(203 48)
\lvec(203 49)
\lvec(202 49)
\ifill f:0
\move(205 48)
\lvec(207 48)
\lvec(207 49)
\lvec(205 49)
\ifill f:0
\move(208 48)
\lvec(210 48)
\lvec(210 49)
\lvec(208 49)
\ifill f:0
\move(211 48)
\lvec(213 48)
\lvec(213 49)
\lvec(211 49)
\ifill f:0
\move(214 48)
\lvec(216 48)
\lvec(216 49)
\lvec(214 49)
\ifill f:0
\move(217 48)
\lvec(219 48)
\lvec(219 49)
\lvec(217 49)
\ifill f:0
\move(220 48)
\lvec(221 48)
\lvec(221 49)
\lvec(220 49)
\ifill f:0
\move(223 48)
\lvec(224 48)
\lvec(224 49)
\lvec(223 49)
\ifill f:0
\move(225 48)
\lvec(227 48)
\lvec(227 49)
\lvec(225 49)
\ifill f:0
\move(228 48)
\lvec(229 48)
\lvec(229 49)
\lvec(228 49)
\ifill f:0
\move(231 48)
\lvec(232 48)
\lvec(232 49)
\lvec(231 49)
\ifill f:0
\move(233 48)
\lvec(234 48)
\lvec(234 49)
\lvec(233 49)
\ifill f:0
\move(235 48)
\lvec(237 48)
\lvec(237 49)
\lvec(235 49)
\ifill f:0
\move(238 48)
\lvec(239 48)
\lvec(239 49)
\lvec(238 49)
\ifill f:0
\move(240 48)
\lvec(241 48)
\lvec(241 49)
\lvec(240 49)
\ifill f:0
\move(242 48)
\lvec(243 48)
\lvec(243 49)
\lvec(242 49)
\ifill f:0
\move(245 48)
\lvec(246 48)
\lvec(246 49)
\lvec(245 49)
\ifill f:0
\move(247 48)
\lvec(248 48)
\lvec(248 49)
\lvec(247 49)
\ifill f:0
\move(249 48)
\lvec(250 48)
\lvec(250 49)
\lvec(249 49)
\ifill f:0
\move(251 48)
\lvec(252 48)
\lvec(252 49)
\lvec(251 49)
\ifill f:0
\move(253 48)
\lvec(254 48)
\lvec(254 49)
\lvec(253 49)
\ifill f:0
\move(255 48)
\lvec(256 48)
\lvec(256 49)
\lvec(255 49)
\ifill f:0
\move(257 48)
\lvec(258 48)
\lvec(258 49)
\lvec(257 49)
\ifill f:0
\move(259 48)
\lvec(261 48)
\lvec(261 49)
\lvec(259 49)
\ifill f:0
\move(262 48)
\lvec(264 48)
\lvec(264 49)
\lvec(262 49)
\ifill f:0
\move(265 48)
\lvec(266 48)
\lvec(266 49)
\lvec(265 49)
\ifill f:0
\move(267 48)
\lvec(268 48)
\lvec(268 49)
\lvec(267 49)
\ifill f:0
\move(269 48)
\lvec(270 48)
\lvec(270 49)
\lvec(269 49)
\ifill f:0
\move(271 48)
\lvec(272 48)
\lvec(272 49)
\lvec(271 49)
\ifill f:0
\move(273 48)
\lvec(274 48)
\lvec(274 49)
\lvec(273 49)
\ifill f:0
\move(275 48)
\lvec(276 48)
\lvec(276 49)
\lvec(275 49)
\ifill f:0
\move(277 48)
\lvec(278 48)
\lvec(278 49)
\lvec(277 49)
\ifill f:0
\move(279 48)
\lvec(280 48)
\lvec(280 49)
\lvec(279 49)
\ifill f:0
\move(282 48)
\lvec(283 48)
\lvec(283 49)
\lvec(282 49)
\ifill f:0
\move(284 48)
\lvec(285 48)
\lvec(285 49)
\lvec(284 49)
\ifill f:0
\move(286 48)
\lvec(287 48)
\lvec(287 49)
\lvec(286 49)
\ifill f:0
\move(288 48)
\lvec(290 48)
\lvec(290 49)
\lvec(288 49)
\ifill f:0
\move(291 48)
\lvec(292 48)
\lvec(292 49)
\lvec(291 49)
\ifill f:0
\move(293 48)
\lvec(294 48)
\lvec(294 49)
\lvec(293 49)
\ifill f:0
\move(295 48)
\lvec(297 48)
\lvec(297 49)
\lvec(295 49)
\ifill f:0
\move(298 48)
\lvec(299 48)
\lvec(299 49)
\lvec(298 49)
\ifill f:0
\move(300 48)
\lvec(302 48)
\lvec(302 49)
\lvec(300 49)
\ifill f:0
\move(303 48)
\lvec(304 48)
\lvec(304 49)
\lvec(303 49)
\ifill f:0
\move(305 48)
\lvec(307 48)
\lvec(307 49)
\lvec(305 49)
\ifill f:0
\move(308 48)
\lvec(309 48)
\lvec(309 49)
\lvec(308 49)
\ifill f:0
\move(311 48)
\lvec(312 48)
\lvec(312 49)
\lvec(311 49)
\ifill f:0
\move(313 48)
\lvec(315 48)
\lvec(315 49)
\lvec(313 49)
\ifill f:0
\move(316 48)
\lvec(318 48)
\lvec(318 49)
\lvec(316 49)
\ifill f:0
\move(319 48)
\lvec(320 48)
\lvec(320 49)
\lvec(319 49)
\ifill f:0
\move(322 48)
\lvec(323 48)
\lvec(323 49)
\lvec(322 49)
\ifill f:0
\move(324 48)
\lvec(326 48)
\lvec(326 49)
\lvec(324 49)
\ifill f:0
\move(327 48)
\lvec(329 48)
\lvec(329 49)
\lvec(327 49)
\ifill f:0
\move(330 48)
\lvec(332 48)
\lvec(332 49)
\lvec(330 49)
\ifill f:0
\move(333 48)
\lvec(335 48)
\lvec(335 49)
\lvec(333 49)
\ifill f:0
\move(336 48)
\lvec(338 48)
\lvec(338 49)
\lvec(336 49)
\ifill f:0
\move(340 48)
\lvec(341 48)
\lvec(341 49)
\lvec(340 49)
\ifill f:0
\move(343 48)
\lvec(345 48)
\lvec(345 49)
\lvec(343 49)
\ifill f:0
\move(346 48)
\lvec(348 48)
\lvec(348 49)
\lvec(346 49)
\ifill f:0
\move(349 48)
\lvec(351 48)
\lvec(351 49)
\lvec(349 49)
\ifill f:0
\move(353 48)
\lvec(355 48)
\lvec(355 49)
\lvec(353 49)
\ifill f:0
\move(356 48)
\lvec(358 48)
\lvec(358 49)
\lvec(356 49)
\ifill f:0
\move(360 48)
\lvec(362 48)
\lvec(362 49)
\lvec(360 49)
\ifill f:0
\move(363 48)
\lvec(365 48)
\lvec(365 49)
\lvec(363 49)
\ifill f:0
\move(367 48)
\lvec(369 48)
\lvec(369 49)
\lvec(367 49)
\ifill f:0
\move(371 48)
\lvec(373 48)
\lvec(373 49)
\lvec(371 49)
\ifill f:0
\move(375 48)
\lvec(377 48)
\lvec(377 49)
\lvec(375 49)
\ifill f:0
\move(379 48)
\lvec(381 48)
\lvec(381 49)
\lvec(379 49)
\ifill f:0
\move(383 48)
\lvec(385 48)
\lvec(385 49)
\lvec(383 49)
\ifill f:0
\move(387 48)
\lvec(389 48)
\lvec(389 49)
\lvec(387 49)
\ifill f:0
\move(391 48)
\lvec(394 48)
\lvec(394 49)
\lvec(391 49)
\ifill f:0
\move(396 48)
\lvec(398 48)
\lvec(398 49)
\lvec(396 49)
\ifill f:0
\move(400 48)
\lvec(403 48)
\lvec(403 49)
\lvec(400 49)
\ifill f:0
\move(405 48)
\lvec(408 48)
\lvec(408 49)
\lvec(405 49)
\ifill f:0
\move(410 48)
\lvec(413 48)
\lvec(413 49)
\lvec(410 49)
\ifill f:0
\move(415 48)
\lvec(418 48)
\lvec(418 49)
\lvec(415 49)
\ifill f:0
\move(420 48)
\lvec(424 48)
\lvec(424 49)
\lvec(420 49)
\ifill f:0
\move(426 48)
\lvec(430 48)
\lvec(430 49)
\lvec(426 49)
\ifill f:0
\move(432 48)
\lvec(436 48)
\lvec(436 49)
\lvec(432 49)
\ifill f:0
\move(438 48)
\lvec(442 48)
\lvec(442 49)
\lvec(438 49)
\ifill f:0
\move(444 48)
\lvec(449 48)
\lvec(449 49)
\lvec(444 49)
\ifill f:0
\move(12 49)
\lvec(14 49)
\lvec(14 50)
\lvec(12 50)
\ifill f:0
\move(81 49)
\lvec(82 49)
\lvec(82 50)
\lvec(81 50)
\ifill f:0
\move(95 49)
\lvec(96 49)
\lvec(96 50)
\lvec(95 50)
\ifill f:0
\move(97 49)
\lvec(98 49)
\lvec(98 50)
\lvec(97 50)
\ifill f:0
\move(99 49)
\lvec(100 49)
\lvec(100 50)
\lvec(99 50)
\ifill f:0
\move(101 49)
\lvec(102 49)
\lvec(102 50)
\lvec(101 50)
\ifill f:0
\move(103 49)
\lvec(104 49)
\lvec(104 50)
\lvec(103 50)
\ifill f:0
\move(108 49)
\lvec(109 49)
\lvec(109 50)
\lvec(108 50)
\ifill f:0
\move(111 49)
\lvec(112 49)
\lvec(112 50)
\lvec(111 50)
\ifill f:0
\move(114 49)
\lvec(115 49)
\lvec(115 50)
\lvec(114 50)
\ifill f:0
\move(117 49)
\lvec(118 49)
\lvec(118 50)
\lvec(117 50)
\ifill f:0
\move(121 49)
\lvec(122 49)
\lvec(122 50)
\lvec(121 50)
\ifill f:0
\move(125 49)
\lvec(126 49)
\lvec(126 50)
\lvec(125 50)
\ifill f:0
\move(130 49)
\lvec(131 49)
\lvec(131 50)
\lvec(130 50)
\ifill f:0
\move(136 49)
\lvec(139 49)
\lvec(139 50)
\lvec(136 50)
\ifill f:0
\move(149 49)
\lvec(159 49)
\lvec(159 50)
\lvec(149 50)
\ifill f:0
\move(169 49)
\lvec(173 49)
\lvec(173 50)
\lvec(169 50)
\ifill f:0
\move(178 49)
\lvec(180 49)
\lvec(180 50)
\lvec(178 50)
\ifill f:0
\move(184 49)
\lvec(186 49)
\lvec(186 50)
\lvec(184 50)
\ifill f:0
\move(189 49)
\lvec(191 49)
\lvec(191 50)
\lvec(189 50)
\ifill f:0
\move(194 49)
\lvec(196 49)
\lvec(196 50)
\lvec(194 50)
\ifill f:0
\move(198 49)
\lvec(200 49)
\lvec(200 50)
\lvec(198 50)
\ifill f:0
\move(202 49)
\lvec(204 49)
\lvec(204 50)
\lvec(202 50)
\ifill f:0
\move(206 49)
\lvec(208 49)
\lvec(208 50)
\lvec(206 50)
\ifill f:0
\move(210 49)
\lvec(211 49)
\lvec(211 50)
\lvec(210 50)
\ifill f:0
\move(213 49)
\lvec(215 49)
\lvec(215 50)
\lvec(213 50)
\ifill f:0
\move(216 49)
\lvec(218 49)
\lvec(218 50)
\lvec(216 50)
\ifill f:0
\move(219 49)
\lvec(221 49)
\lvec(221 50)
\lvec(219 50)
\ifill f:0
\move(222 49)
\lvec(224 49)
\lvec(224 50)
\lvec(222 50)
\ifill f:0
\move(225 49)
\lvec(227 49)
\lvec(227 50)
\lvec(225 50)
\ifill f:0
\move(228 49)
\lvec(229 49)
\lvec(229 50)
\lvec(228 50)
\ifill f:0
\move(231 49)
\lvec(232 49)
\lvec(232 50)
\lvec(231 50)
\ifill f:0
\move(234 49)
\lvec(235 49)
\lvec(235 50)
\lvec(234 50)
\ifill f:0
\move(236 49)
\lvec(237 49)
\lvec(237 50)
\lvec(236 50)
\ifill f:0
\move(239 49)
\lvec(240 49)
\lvec(240 50)
\lvec(239 50)
\ifill f:0
\move(241 49)
\lvec(242 49)
\lvec(242 50)
\lvec(241 50)
\ifill f:0
\move(244 49)
\lvec(245 49)
\lvec(245 50)
\lvec(244 50)
\ifill f:0
\move(246 49)
\lvec(247 49)
\lvec(247 50)
\lvec(246 50)
\ifill f:0
\move(248 49)
\lvec(250 49)
\lvec(250 50)
\lvec(248 50)
\ifill f:0
\move(251 49)
\lvec(252 49)
\lvec(252 50)
\lvec(251 50)
\ifill f:0
\move(253 49)
\lvec(254 49)
\lvec(254 50)
\lvec(253 50)
\ifill f:0
\move(255 49)
\lvec(256 49)
\lvec(256 50)
\lvec(255 50)
\ifill f:0
\move(258 49)
\lvec(259 49)
\lvec(259 50)
\lvec(258 50)
\ifill f:0
\move(260 49)
\lvec(261 49)
\lvec(261 50)
\lvec(260 50)
\ifill f:0
\move(262 49)
\lvec(263 49)
\lvec(263 50)
\lvec(262 50)
\ifill f:0
\move(264 49)
\lvec(265 49)
\lvec(265 50)
\lvec(264 50)
\ifill f:0
\move(266 49)
\lvec(267 49)
\lvec(267 50)
\lvec(266 50)
\ifill f:0
\move(268 49)
\lvec(269 49)
\lvec(269 50)
\lvec(268 50)
\ifill f:0
\move(270 49)
\lvec(271 49)
\lvec(271 50)
\lvec(270 50)
\ifill f:0
\move(272 49)
\lvec(273 49)
\lvec(273 50)
\lvec(272 50)
\ifill f:0
\move(274 49)
\lvec(275 49)
\lvec(275 50)
\lvec(274 50)
\ifill f:0
\move(276 49)
\lvec(277 49)
\lvec(277 50)
\lvec(276 50)
\ifill f:0
\move(278 49)
\lvec(279 49)
\lvec(279 50)
\lvec(278 50)
\ifill f:0
\move(280 49)
\lvec(281 49)
\lvec(281 50)
\lvec(280 50)
\ifill f:0
\move(282 49)
\lvec(283 49)
\lvec(283 50)
\lvec(282 50)
\ifill f:0
\move(284 49)
\lvec(285 49)
\lvec(285 50)
\lvec(284 50)
\ifill f:0
\move(286 49)
\lvec(287 49)
\lvec(287 50)
\lvec(286 50)
\ifill f:0
\move(288 49)
\lvec(290 49)
\lvec(290 50)
\lvec(288 50)
\ifill f:0
\move(291 49)
\lvec(292 49)
\lvec(292 50)
\lvec(291 50)
\ifill f:0
\move(293 49)
\lvec(294 49)
\lvec(294 50)
\lvec(293 50)
\ifill f:0
\move(295 49)
\lvec(296 49)
\lvec(296 50)
\lvec(295 50)
\ifill f:0
\move(297 49)
\lvec(299 49)
\lvec(299 50)
\lvec(297 50)
\ifill f:0
\move(300 49)
\lvec(301 49)
\lvec(301 50)
\lvec(300 50)
\ifill f:0
\move(302 49)
\lvec(303 49)
\lvec(303 50)
\lvec(302 50)
\ifill f:0
\move(304 49)
\lvec(306 49)
\lvec(306 50)
\lvec(304 50)
\ifill f:0
\move(307 49)
\lvec(308 49)
\lvec(308 50)
\lvec(307 50)
\ifill f:0
\move(309 49)
\lvec(311 49)
\lvec(311 50)
\lvec(309 50)
\ifill f:0
\move(312 49)
\lvec(313 49)
\lvec(313 50)
\lvec(312 50)
\ifill f:0
\move(314 49)
\lvec(316 49)
\lvec(316 50)
\lvec(314 50)
\ifill f:0
\move(317 49)
\lvec(318 49)
\lvec(318 50)
\lvec(317 50)
\ifill f:0
\move(319 49)
\lvec(321 49)
\lvec(321 50)
\lvec(319 50)
\ifill f:0
\move(322 49)
\lvec(323 49)
\lvec(323 50)
\lvec(322 50)
\ifill f:0
\move(325 49)
\lvec(326 49)
\lvec(326 50)
\lvec(325 50)
\ifill f:0
\move(327 49)
\lvec(329 49)
\lvec(329 50)
\lvec(327 50)
\ifill f:0
\move(330 49)
\lvec(331 49)
\lvec(331 50)
\lvec(330 50)
\ifill f:0
\move(333 49)
\lvec(334 49)
\lvec(334 50)
\lvec(333 50)
\ifill f:0
\move(336 49)
\lvec(337 49)
\lvec(337 50)
\lvec(336 50)
\ifill f:0
\move(338 49)
\lvec(340 49)
\lvec(340 50)
\lvec(338 50)
\ifill f:0
\move(341 49)
\lvec(343 49)
\lvec(343 50)
\lvec(341 50)
\ifill f:0
\move(344 49)
\lvec(346 49)
\lvec(346 50)
\lvec(344 50)
\ifill f:0
\move(347 49)
\lvec(349 49)
\lvec(349 50)
\lvec(347 50)
\ifill f:0
\move(350 49)
\lvec(352 49)
\lvec(352 50)
\lvec(350 50)
\ifill f:0
\move(353 49)
\lvec(355 49)
\lvec(355 50)
\lvec(353 50)
\ifill f:0
\move(357 49)
\lvec(358 49)
\lvec(358 50)
\lvec(357 50)
\ifill f:0
\move(360 49)
\lvec(362 49)
\lvec(362 50)
\lvec(360 50)
\ifill f:0
\move(363 49)
\lvec(365 49)
\lvec(365 50)
\lvec(363 50)
\ifill f:0
\move(366 49)
\lvec(368 49)
\lvec(368 50)
\lvec(366 50)
\ifill f:0
\move(370 49)
\lvec(372 49)
\lvec(372 50)
\lvec(370 50)
\ifill f:0
\move(373 49)
\lvec(375 49)
\lvec(375 50)
\lvec(373 50)
\ifill f:0
\move(377 49)
\lvec(379 49)
\lvec(379 50)
\lvec(377 50)
\ifill f:0
\move(381 49)
\lvec(383 49)
\lvec(383 50)
\lvec(381 50)
\ifill f:0
\move(384 49)
\lvec(387 49)
\lvec(387 50)
\lvec(384 50)
\ifill f:0
\move(388 49)
\lvec(391 49)
\lvec(391 50)
\lvec(388 50)
\ifill f:0
\move(392 49)
\lvec(395 49)
\lvec(395 50)
\lvec(392 50)
\ifill f:0
\move(396 49)
\lvec(399 49)
\lvec(399 50)
\lvec(396 50)
\ifill f:0
\move(400 49)
\lvec(403 49)
\lvec(403 50)
\lvec(400 50)
\ifill f:0
\move(405 49)
\lvec(407 49)
\lvec(407 50)
\lvec(405 50)
\ifill f:0
\move(409 49)
\lvec(412 49)
\lvec(412 50)
\lvec(409 50)
\ifill f:0
\move(414 49)
\lvec(416 49)
\lvec(416 50)
\lvec(414 50)
\ifill f:0
\move(418 49)
\lvec(421 49)
\lvec(421 50)
\lvec(418 50)
\ifill f:0
\move(423 49)
\lvec(426 49)
\lvec(426 50)
\lvec(423 50)
\ifill f:0
\move(428 49)
\lvec(431 49)
\lvec(431 50)
\lvec(428 50)
\ifill f:0
\move(433 49)
\lvec(436 49)
\lvec(436 50)
\lvec(433 50)
\ifill f:0
\move(438 49)
\lvec(442 49)
\lvec(442 50)
\lvec(438 50)
\ifill f:0
\move(444 49)
\lvec(448 49)
\lvec(448 50)
\lvec(444 50)
\ifill f:0
\move(450 49)
\lvec(451 49)
\lvec(451 50)
\lvec(450 50)
\ifill f:0
\move(81 50)
\lvec(82 50)
\lvec(82 51)
\lvec(81 51)
\ifill f:0
\move(99 50)
\lvec(100 50)
\lvec(100 51)
\lvec(99 51)
\ifill f:0
\move(101 50)
\lvec(102 50)
\lvec(102 51)
\lvec(101 51)
\ifill f:0
\move(103 50)
\lvec(104 50)
\lvec(104 51)
\lvec(103 51)
\ifill f:0
\move(105 50)
\lvec(106 50)
\lvec(106 51)
\lvec(105 51)
\ifill f:0
\move(107 50)
\lvec(108 50)
\lvec(108 51)
\lvec(107 51)
\ifill f:0
\move(112 50)
\lvec(113 50)
\lvec(113 51)
\lvec(112 51)
\ifill f:0
\move(121 50)
\lvec(122 50)
\lvec(122 51)
\lvec(121 51)
\ifill f:0
\move(124 50)
\lvec(125 50)
\lvec(125 51)
\lvec(124 51)
\ifill f:0
\move(128 50)
\lvec(129 50)
\lvec(129 51)
\lvec(128 51)
\ifill f:0
\move(133 50)
\lvec(134 50)
\lvec(134 51)
\lvec(133 51)
\ifill f:0
\move(139 50)
\lvec(141 50)
\lvec(141 51)
\lvec(139 51)
\ifill f:0
\move(147 50)
\lvec(151 50)
\lvec(151 51)
\lvec(147 51)
\ifill f:0
\move(169 50)
\lvec(174 50)
\lvec(174 51)
\lvec(169 51)
\ifill f:0
\move(180 50)
\lvec(183 50)
\lvec(183 51)
\lvec(180 51)
\ifill f:0
\move(188 50)
\lvec(190 50)
\lvec(190 51)
\lvec(188 51)
\ifill f:0
\move(194 50)
\lvec(196 50)
\lvec(196 51)
\lvec(194 51)
\ifill f:0
\move(199 50)
\lvec(201 50)
\lvec(201 51)
\lvec(199 51)
\ifill f:0
\move(203 50)
\lvec(205 50)
\lvec(205 51)
\lvec(203 51)
\ifill f:0
\move(208 50)
\lvec(209 50)
\lvec(209 51)
\lvec(208 51)
\ifill f:0
\move(212 50)
\lvec(213 50)
\lvec(213 51)
\lvec(212 51)
\ifill f:0
\move(215 50)
\lvec(217 50)
\lvec(217 51)
\lvec(215 51)
\ifill f:0
\move(219 50)
\lvec(220 50)
\lvec(220 51)
\lvec(219 51)
\ifill f:0
\move(222 50)
\lvec(224 50)
\lvec(224 51)
\lvec(222 51)
\ifill f:0
\move(225 50)
\lvec(227 50)
\lvec(227 51)
\lvec(225 51)
\ifill f:0
\move(228 50)
\lvec(230 50)
\lvec(230 51)
\lvec(228 51)
\ifill f:0
\move(231 50)
\lvec(233 50)
\lvec(233 51)
\lvec(231 51)
\ifill f:0
\move(234 50)
\lvec(236 50)
\lvec(236 51)
\lvec(234 51)
\ifill f:0
\move(237 50)
\lvec(238 50)
\lvec(238 51)
\lvec(237 51)
\ifill f:0
\move(240 50)
\lvec(241 50)
\lvec(241 51)
\lvec(240 51)
\ifill f:0
\move(243 50)
\lvec(244 50)
\lvec(244 51)
\lvec(243 51)
\ifill f:0
\move(245 50)
\lvec(246 50)
\lvec(246 51)
\lvec(245 51)
\ifill f:0
\move(248 50)
\lvec(249 50)
\lvec(249 51)
\lvec(248 51)
\ifill f:0
\move(250 50)
\lvec(251 50)
\lvec(251 51)
\lvec(250 51)
\ifill f:0
\move(253 50)
\lvec(254 50)
\lvec(254 51)
\lvec(253 51)
\ifill f:0
\move(255 50)
\lvec(256 50)
\lvec(256 51)
\lvec(255 51)
\ifill f:0
\move(258 50)
\lvec(259 50)
\lvec(259 51)
\lvec(258 51)
\ifill f:0
\move(260 50)
\lvec(261 50)
\lvec(261 51)
\lvec(260 51)
\ifill f:0
\move(262 50)
\lvec(263 50)
\lvec(263 51)
\lvec(262 51)
\ifill f:0
\move(264 50)
\lvec(266 50)
\lvec(266 51)
\lvec(264 51)
\ifill f:0
\move(267 50)
\lvec(268 50)
\lvec(268 51)
\lvec(267 51)
\ifill f:0
\move(269 50)
\lvec(270 50)
\lvec(270 51)
\lvec(269 51)
\ifill f:0
\move(271 50)
\lvec(272 50)
\lvec(272 51)
\lvec(271 51)
\ifill f:0
\move(273 50)
\lvec(274 50)
\lvec(274 51)
\lvec(273 51)
\ifill f:0
\move(275 50)
\lvec(276 50)
\lvec(276 51)
\lvec(275 51)
\ifill f:0
\move(277 50)
\lvec(278 50)
\lvec(278 51)
\lvec(277 51)
\ifill f:0
\move(279 50)
\lvec(280 50)
\lvec(280 51)
\lvec(279 51)
\ifill f:0
\move(281 50)
\lvec(282 50)
\lvec(282 51)
\lvec(281 51)
\ifill f:0
\move(285 50)
\lvec(286 50)
\lvec(286 51)
\lvec(285 51)
\ifill f:0
\move(287 50)
\lvec(288 50)
\lvec(288 51)
\lvec(287 51)
\ifill f:0
\move(289 50)
\lvec(290 50)
\lvec(290 51)
\lvec(289 51)
\ifill f:0
\move(291 50)
\lvec(292 50)
\lvec(292 51)
\lvec(291 51)
\ifill f:0
\move(293 50)
\lvec(294 50)
\lvec(294 51)
\lvec(293 51)
\ifill f:0
\move(295 50)
\lvec(296 50)
\lvec(296 51)
\lvec(295 51)
\ifill f:0
\move(297 50)
\lvec(298 50)
\lvec(298 51)
\lvec(297 51)
\ifill f:0
\move(299 50)
\lvec(300 50)
\lvec(300 51)
\lvec(299 51)
\ifill f:0
\move(301 50)
\lvec(302 50)
\lvec(302 51)
\lvec(301 51)
\ifill f:0
\move(303 50)
\lvec(305 50)
\lvec(305 51)
\lvec(303 51)
\ifill f:0
\move(306 50)
\lvec(307 50)
\lvec(307 51)
\lvec(306 51)
\ifill f:0
\move(308 50)
\lvec(309 50)
\lvec(309 51)
\lvec(308 51)
\ifill f:0
\move(310 50)
\lvec(311 50)
\lvec(311 51)
\lvec(310 51)
\ifill f:0
\move(313 50)
\lvec(314 50)
\lvec(314 51)
\lvec(313 51)
\ifill f:0
\move(315 50)
\lvec(316 50)
\lvec(316 51)
\lvec(315 51)
\ifill f:0
\move(317 50)
\lvec(318 50)
\lvec(318 51)
\lvec(317 51)
\ifill f:0
\move(320 50)
\lvec(321 50)
\lvec(321 51)
\lvec(320 51)
\ifill f:0
\move(322 50)
\lvec(323 50)
\lvec(323 51)
\lvec(322 51)
\ifill f:0
\move(325 50)
\lvec(326 50)
\lvec(326 51)
\lvec(325 51)
\ifill f:0
\move(327 50)
\lvec(328 50)
\lvec(328 51)
\lvec(327 51)
\ifill f:0
\move(330 50)
\lvec(331 50)
\lvec(331 51)
\lvec(330 51)
\ifill f:0
\move(332 50)
\lvec(334 50)
\lvec(334 51)
\lvec(332 51)
\ifill f:0
\move(335 50)
\lvec(336 50)
\lvec(336 51)
\lvec(335 51)
\ifill f:0
\move(337 50)
\lvec(339 50)
\lvec(339 51)
\lvec(337 51)
\ifill f:0
\move(340 50)
\lvec(342 50)
\lvec(342 51)
\lvec(340 51)
\ifill f:0
\move(343 50)
\lvec(344 50)
\lvec(344 51)
\lvec(343 51)
\ifill f:0
\move(346 50)
\lvec(347 50)
\lvec(347 51)
\lvec(346 51)
\ifill f:0
\move(348 50)
\lvec(350 50)
\lvec(350 51)
\lvec(348 51)
\ifill f:0
\move(351 50)
\lvec(353 50)
\lvec(353 51)
\lvec(351 51)
\ifill f:0
\move(354 50)
\lvec(356 50)
\lvec(356 51)
\lvec(354 51)
\ifill f:0
\move(357 50)
\lvec(359 50)
\lvec(359 51)
\lvec(357 51)
\ifill f:0
\move(360 50)
\lvec(362 50)
\lvec(362 51)
\lvec(360 51)
\ifill f:0
\move(363 50)
\lvec(365 50)
\lvec(365 51)
\lvec(363 51)
\ifill f:0
\move(366 50)
\lvec(368 50)
\lvec(368 51)
\lvec(366 51)
\ifill f:0
\move(369 50)
\lvec(371 50)
\lvec(371 51)
\lvec(369 51)
\ifill f:0
\move(372 50)
\lvec(374 50)
\lvec(374 51)
\lvec(372 51)
\ifill f:0
\move(376 50)
\lvec(378 50)
\lvec(378 51)
\lvec(376 51)
\ifill f:0
\move(379 50)
\lvec(381 50)
\lvec(381 51)
\lvec(379 51)
\ifill f:0
\move(382 50)
\lvec(384 50)
\lvec(384 51)
\lvec(382 51)
\ifill f:0
\move(386 50)
\lvec(388 50)
\lvec(388 51)
\lvec(386 51)
\ifill f:0
\move(389 50)
\lvec(391 50)
\lvec(391 51)
\lvec(389 51)
\ifill f:0
\move(393 50)
\lvec(395 50)
\lvec(395 51)
\lvec(393 51)
\ifill f:0
\move(397 50)
\lvec(399 50)
\lvec(399 51)
\lvec(397 51)
\ifill f:0
\move(400 50)
\lvec(403 50)
\lvec(403 51)
\lvec(400 51)
\ifill f:0
\move(404 50)
\lvec(407 50)
\lvec(407 51)
\lvec(404 51)
\ifill f:0
\move(408 50)
\lvec(411 50)
\lvec(411 51)
\lvec(408 51)
\ifill f:0
\move(412 50)
\lvec(415 50)
\lvec(415 51)
\lvec(412 51)
\ifill f:0
\move(416 50)
\lvec(419 50)
\lvec(419 51)
\lvec(416 51)
\ifill f:0
\move(421 50)
\lvec(423 50)
\lvec(423 51)
\lvec(421 51)
\ifill f:0
\move(425 50)
\lvec(428 50)
\lvec(428 51)
\lvec(425 51)
\ifill f:0
\move(429 50)
\lvec(432 50)
\lvec(432 51)
\lvec(429 51)
\ifill f:0
\move(434 50)
\lvec(437 50)
\lvec(437 51)
\lvec(434 51)
\ifill f:0
\move(439 50)
\lvec(442 50)
\lvec(442 51)
\lvec(439 51)
\ifill f:0
\move(444 50)
\lvec(447 50)
\lvec(447 51)
\lvec(444 51)
\ifill f:0
\move(449 50)
\lvec(451 50)
\lvec(451 51)
\lvec(449 51)
\ifill f:0
\move(11 8)
\lvec(13 8)
\lvec(13 9)
\lvec(11 9)
\ifill f:0
\move(14 8)
\lvec(24 8)
\lvec(24 9)
\lvec(14 9)
\ifill f:0
\move(25 8)
\lvec(37 8)
\lvec(37 9)
\lvec(25 9)
\ifill f:0
\move(38 8)
\lvec(43 8)
\lvec(43 9)
\lvec(38 9)
\ifill f:0
\move(44 8)
\lvec(46 8)
\lvec(46 9)
\lvec(44 9)
\ifill f:0
\move(47 8)
\lvec(51 8)
\lvec(51 9)
\lvec(47 9)
\ifill f:0
\move(52 8)
\lvec(56 8)
\lvec(56 9)
\lvec(52 9)
\ifill f:0
\move(57 8)
\lvec(58 8)
\lvec(58 9)
\lvec(57 9)
\ifill f:0
\move(59 8)
\lvec(71 8)
\lvec(71 9)
\lvec(59 9)
\ifill f:0
\move(72 8)
\lvec(73 8)
\lvec(73 9)
\lvec(72 9)
\ifill f:0
\move(74 8)
\lvec(88 8)
\lvec(88 9)
\lvec(74 9)
\ifill f:0
\move(89 8)
\lvec(90 8)
\lvec(90 9)
\lvec(89 9)
\ifill f:0
\move(91 8)
\lvec(97 8)
\lvec(97 9)
\lvec(91 9)
\ifill f:0
\move(98 8)
\lvec(104 8)
\lvec(104 9)
\lvec(98 9)
\ifill f:0
\move(105 8)
\lvec(109 8)
\lvec(109 9)
\lvec(105 9)
\ifill f:0
\move(110 8)
\lvec(119 8)
\lvec(119 9)
\lvec(110 9)
\ifill f:0
\move(120 8)
\lvec(146 8)
\lvec(146 9)
\lvec(120 9)
\ifill f:0
\move(147 8)
\lvec(149 8)
\lvec(149 9)
\lvec(147 9)
\ifill f:0
\move(150 8)
\lvec(176 8)
\lvec(176 9)
\lvec(150 9)
\ifill f:0
\move(177 8)
\lvec(179 8)
\lvec(179 9)
\lvec(177 9)
\ifill f:0
\move(180 8)
\lvec(212 8)
\lvec(212 9)
\lvec(180 9)
\ifill f:0
\move(213 8)
\lvec(226 8)
\lvec(226 9)
\lvec(213 9)
\ifill f:0
\move(227 8)
\lvec(233 8)
\lvec(233 9)
\lvec(227 9)
\ifill f:0
\move(234 8)
\lvec(244 8)
\lvec(244 9)
\lvec(234 9)
\ifill f:0
\move(245 8)
\lvec(255 8)
\lvec(255 9)
\lvec(245 9)
\ifill f:0
\move(256 8)
\lvec(259 8)
\lvec(259 9)
\lvec(256 9)
\ifill f:0
\move(260 8)
\lvec(286 8)
\lvec(286 9)
\lvec(260 9)
\ifill f:0
\move(287 8)
\lvec(290 8)
\lvec(290 9)
\lvec(287 9)
\ifill f:0
\move(291 8)
\lvec(319 8)
\lvec(319 9)
\lvec(291 9)
\ifill f:0
\move(320 8)
\lvec(323 8)
\lvec(323 9)
\lvec(320 9)
\ifill f:0
\move(324 8)
\lvec(336 8)
\lvec(336 9)
\lvec(324 9)
\ifill f:0
\move(337 8)
\lvec(349 8)
\lvec(349 9)
\lvec(337 9)
\ifill f:0
\move(350 8)
\lvec(358 8)
\lvec(358 9)
\lvec(350 9)
\ifill f:0
\move(359 8)
\lvec(376 8)
\lvec(376 9)
\lvec(359 9)
\ifill f:0
\move(377 8)
\lvec(423 8)
\lvec(423 9)
\lvec(377 9)
\ifill f:0
\move(424 8)
\lvec(428 8)
\lvec(428 9)
\lvec(424 9)
\ifill f:0
\move(429 8)
\lvec(451 8)
\lvec(451 9)
\lvec(429 9)
\ifill f:0
\move(12 51)
\lvec(13 51)
\lvec(13 52)
\lvec(12 52)
\ifill f:0
\move(15 51)
\lvec(16 51)
\lvec(16 52)
\lvec(15 52)
\ifill f:0
\move(19 51)
\lvec(20 51)
\lvec(20 52)
\lvec(19 52)
\ifill f:0
\move(81 51)
\lvec(82 51)
\lvec(82 52)
\lvec(81 52)
\ifill f:0
\move(99 51)
\lvec(100 51)
\lvec(100 52)
\lvec(99 52)
\ifill f:0
\move(113 51)
\lvec(114 51)
\lvec(114 52)
\lvec(113 52)
\ifill f:0
\move(118 51)
\lvec(119 51)
\lvec(119 52)
\lvec(118 52)
\ifill f:0
\move(121 51)
\lvec(122 51)
\lvec(122 52)
\lvec(121 52)
\ifill f:0
\move(124 51)
\lvec(125 51)
\lvec(125 52)
\lvec(124 52)
\ifill f:0
\move(127 51)
\lvec(128 51)
\lvec(128 52)
\lvec(127 52)
\ifill f:0
\move(131 51)
\lvec(132 51)
\lvec(132 52)
\lvec(131 52)
\ifill f:0
\move(135 51)
\lvec(136 51)
\lvec(136 52)
\lvec(135 52)
\ifill f:0
\move(140 51)
\lvec(141 51)
\lvec(141 52)
\lvec(140 52)
\ifill f:0
\move(146 51)
\lvec(149 51)
\lvec(149 52)
\lvec(146 52)
\ifill f:0
\move(156 51)
\lvec(166 51)
\lvec(166 52)
\lvec(156 52)
\ifill f:0
\move(167 51)
\lvec(177 51)
\lvec(177 52)
\lvec(167 52)
\ifill f:0
\move(185 51)
\lvec(188 51)
\lvec(188 52)
\lvec(185 52)
\ifill f:0
\move(193 51)
\lvec(196 51)
\lvec(196 52)
\lvec(193 52)
\ifill f:0
\move(199 51)
\lvec(202 51)
\lvec(202 52)
\lvec(199 52)
\ifill f:0
\move(205 51)
\lvec(207 51)
\lvec(207 52)
\lvec(205 52)
\ifill f:0
\move(209 51)
\lvec(211 51)
\lvec(211 52)
\lvec(209 52)
\ifill f:0
\move(214 51)
\lvec(216 51)
\lvec(216 52)
\lvec(214 52)
\ifill f:0
\move(218 51)
\lvec(220 51)
\lvec(220 52)
\lvec(218 52)
\ifill f:0
\move(222 51)
\lvec(223 51)
\lvec(223 52)
\lvec(222 52)
\ifill f:0
\move(225 51)
\lvec(227 51)
\lvec(227 52)
\lvec(225 52)
\ifill f:0
\move(229 51)
\lvec(230 51)
\lvec(230 52)
\lvec(229 52)
\ifill f:0
\move(232 51)
\lvec(234 51)
\lvec(234 52)
\lvec(232 52)
\ifill f:0
\move(235 51)
\lvec(237 51)
\lvec(237 52)
\lvec(235 52)
\ifill f:0
\move(238 51)
\lvec(240 51)
\lvec(240 52)
\lvec(238 52)
\ifill f:0
\move(241 51)
\lvec(243 51)
\lvec(243 52)
\lvec(241 52)
\ifill f:0
\move(244 51)
\lvec(246 51)
\lvec(246 52)
\lvec(244 52)
\ifill f:0
\move(247 51)
\lvec(248 51)
\lvec(248 52)
\lvec(247 52)
\ifill f:0
\move(250 51)
\lvec(251 51)
\lvec(251 52)
\lvec(250 52)
\ifill f:0
\move(252 51)
\lvec(254 51)
\lvec(254 52)
\lvec(252 52)
\ifill f:0
\move(255 51)
\lvec(256 51)
\lvec(256 52)
\lvec(255 52)
\ifill f:0
\move(258 51)
\lvec(259 51)
\lvec(259 52)
\lvec(258 52)
\ifill f:0
\move(260 51)
\lvec(261 51)
\lvec(261 52)
\lvec(260 52)
\ifill f:0
\move(263 51)
\lvec(264 51)
\lvec(264 52)
\lvec(263 52)
\ifill f:0
\move(265 51)
\lvec(266 51)
\lvec(266 52)
\lvec(265 52)
\ifill f:0
\move(267 51)
\lvec(269 51)
\lvec(269 52)
\lvec(267 52)
\ifill f:0
\move(270 51)
\lvec(271 51)
\lvec(271 52)
\lvec(270 52)
\ifill f:0
\move(272 51)
\lvec(273 51)
\lvec(273 52)
\lvec(272 52)
\ifill f:0
\move(274 51)
\lvec(275 51)
\lvec(275 52)
\lvec(274 52)
\ifill f:0
\move(277 51)
\lvec(278 51)
\lvec(278 52)
\lvec(277 52)
\ifill f:0
\move(279 51)
\lvec(280 51)
\lvec(280 52)
\lvec(279 52)
\ifill f:0
\move(281 51)
\lvec(282 51)
\lvec(282 52)
\lvec(281 52)
\ifill f:0
\move(283 51)
\lvec(284 51)
\lvec(284 52)
\lvec(283 52)
\ifill f:0
\move(285 51)
\lvec(286 51)
\lvec(286 52)
\lvec(285 52)
\ifill f:0
\move(287 51)
\lvec(288 51)
\lvec(288 52)
\lvec(287 52)
\ifill f:0
\move(289 51)
\lvec(290 51)
\lvec(290 52)
\lvec(289 52)
\ifill f:0
\move(291 51)
\lvec(292 51)
\lvec(292 52)
\lvec(291 52)
\ifill f:0
\move(293 51)
\lvec(295 51)
\lvec(295 52)
\lvec(293 52)
\ifill f:0
\move(296 51)
\lvec(298 51)
\lvec(298 52)
\lvec(296 52)
\ifill f:0
\move(299 51)
\lvec(300 51)
\lvec(300 52)
\lvec(299 52)
\ifill f:0
\move(301 51)
\lvec(302 51)
\lvec(302 52)
\lvec(301 52)
\ifill f:0
\move(303 51)
\lvec(304 51)
\lvec(304 52)
\lvec(303 52)
\ifill f:0
\move(305 51)
\lvec(306 51)
\lvec(306 52)
\lvec(305 52)
\ifill f:0
\move(307 51)
\lvec(308 51)
\lvec(308 52)
\lvec(307 52)
\ifill f:0
\move(309 51)
\lvec(310 51)
\lvec(310 52)
\lvec(309 52)
\ifill f:0
\move(311 51)
\lvec(312 51)
\lvec(312 52)
\lvec(311 52)
\ifill f:0
\move(313 51)
\lvec(314 51)
\lvec(314 52)
\lvec(313 52)
\ifill f:0
\move(315 51)
\lvec(317 51)
\lvec(317 52)
\lvec(315 52)
\ifill f:0
\move(318 51)
\lvec(319 51)
\lvec(319 52)
\lvec(318 52)
\ifill f:0
\move(320 51)
\lvec(321 51)
\lvec(321 52)
\lvec(320 52)
\ifill f:0
\move(322 51)
\lvec(323 51)
\lvec(323 52)
\lvec(322 52)
\ifill f:0
\move(325 51)
\lvec(326 51)
\lvec(326 52)
\lvec(325 52)
\ifill f:0
\move(327 51)
\lvec(328 51)
\lvec(328 52)
\lvec(327 52)
\ifill f:0
\move(329 51)
\lvec(330 51)
\lvec(330 52)
\lvec(329 52)
\ifill f:0
\move(332 51)
\lvec(333 51)
\lvec(333 52)
\lvec(332 52)
\ifill f:0
\move(334 51)
\lvec(335 51)
\lvec(335 52)
\lvec(334 52)
\ifill f:0
\move(336 51)
\lvec(338 51)
\lvec(338 52)
\lvec(336 52)
\ifill f:0
\move(339 51)
\lvec(340 51)
\lvec(340 52)
\lvec(339 52)
\ifill f:0
\move(341 51)
\lvec(343 51)
\lvec(343 52)
\lvec(341 52)
\ifill f:0
\move(344 51)
\lvec(345 51)
\lvec(345 52)
\lvec(344 52)
\ifill f:0
\move(347 51)
\lvec(348 51)
\lvec(348 52)
\lvec(347 52)
\ifill f:0
\move(349 51)
\lvec(351 51)
\lvec(351 52)
\lvec(349 52)
\ifill f:0
\move(352 51)
\lvec(353 51)
\lvec(353 52)
\lvec(352 52)
\ifill f:0
\move(355 51)
\lvec(356 51)
\lvec(356 52)
\lvec(355 52)
\ifill f:0
\move(357 51)
\lvec(359 51)
\lvec(359 52)
\lvec(357 52)
\ifill f:0
\move(360 51)
\lvec(362 51)
\lvec(362 52)
\lvec(360 52)
\ifill f:0
\move(363 51)
\lvec(364 51)
\lvec(364 52)
\lvec(363 52)
\ifill f:0
\move(366 51)
\lvec(367 51)
\lvec(367 52)
\lvec(366 52)
\ifill f:0
\move(369 51)
\lvec(370 51)
\lvec(370 52)
\lvec(369 52)
\ifill f:0
\move(372 51)
\lvec(373 51)
\lvec(373 52)
\lvec(372 52)
\ifill f:0
\move(375 51)
\lvec(376 51)
\lvec(376 52)
\lvec(375 52)
\ifill f:0
\move(378 51)
\lvec(379 51)
\lvec(379 52)
\lvec(378 52)
\ifill f:0
\move(381 51)
\lvec(382 51)
\lvec(382 52)
\lvec(381 52)
\ifill f:0
\move(384 51)
\lvec(386 51)
\lvec(386 52)
\lvec(384 52)
\ifill f:0
\move(387 51)
\lvec(389 51)
\lvec(389 52)
\lvec(387 52)
\ifill f:0
\move(390 51)
\lvec(392 51)
\lvec(392 52)
\lvec(390 52)
\ifill f:0
\move(394 51)
\lvec(396 51)
\lvec(396 52)
\lvec(394 52)
\ifill f:0
\move(397 51)
\lvec(399 51)
\lvec(399 52)
\lvec(397 52)
\ifill f:0
\move(400 51)
\lvec(402 51)
\lvec(402 52)
\lvec(400 52)
\ifill f:0
\move(404 51)
\lvec(406 51)
\lvec(406 52)
\lvec(404 52)
\ifill f:0
\move(407 51)
\lvec(410 51)
\lvec(410 52)
\lvec(407 52)
\ifill f:0
\move(411 51)
\lvec(413 51)
\lvec(413 52)
\lvec(411 52)
\ifill f:0
\move(415 51)
\lvec(417 51)
\lvec(417 52)
\lvec(415 52)
\ifill f:0
\move(419 51)
\lvec(421 51)
\lvec(421 52)
\lvec(419 52)
\ifill f:0
\move(423 51)
\lvec(425 51)
\lvec(425 52)
\lvec(423 52)
\ifill f:0
\move(427 51)
\lvec(429 51)
\lvec(429 52)
\lvec(427 52)
\ifill f:0
\move(431 51)
\lvec(433 51)
\lvec(433 52)
\lvec(431 52)
\ifill f:0
\move(435 51)
\lvec(437 51)
\lvec(437 52)
\lvec(435 52)
\ifill f:0
\move(439 51)
\lvec(442 51)
\lvec(442 52)
\lvec(439 52)
\ifill f:0
\move(444 51)
\lvec(446 51)
\lvec(446 52)
\lvec(444 52)
\ifill f:0
\move(448 51)
\lvec(451 51)
\lvec(451 52)
\lvec(448 52)
\ifill f:0
\move(13 52)
\lvec(16 52)
\lvec(16 53)
\lvec(13 53)
\ifill f:0
\move(81 52)
\lvec(82 52)
\lvec(82 53)
\lvec(81 53)
\ifill f:0
\move(88 52)
\lvec(89 52)
\lvec(89 53)
\lvec(88 53)
\ifill f:0
\move(104 52)
\lvec(105 52)
\lvec(105 53)
\lvec(104 53)
\ifill f:0
\move(106 52)
\lvec(107 52)
\lvec(107 53)
\lvec(106 53)
\ifill f:0
\move(116 52)
\lvec(117 52)
\lvec(117 53)
\lvec(116 53)
\ifill f:0
\move(121 52)
\lvec(122 52)
\lvec(122 53)
\lvec(121 53)
\ifill f:0
\move(126 52)
\lvec(127 52)
\lvec(127 53)
\lvec(126 53)
\ifill f:0
\move(133 52)
\lvec(134 52)
\lvec(134 53)
\lvec(133 53)
\ifill f:0
\move(137 52)
\lvec(138 52)
\lvec(138 53)
\lvec(137 53)
\ifill f:0
\move(141 52)
\lvec(142 52)
\lvec(142 53)
\lvec(141 53)
\ifill f:0
\move(146 52)
\lvec(147 52)
\lvec(147 53)
\lvec(146 53)
\ifill f:0
\move(152 52)
\lvec(155 52)
\lvec(155 53)
\lvec(152 53)
\ifill f:0
\move(163 52)
\lvec(172 52)
\lvec(172 53)
\lvec(163 53)
\ifill f:0
\move(173 52)
\lvec(184 52)
\lvec(184 53)
\lvec(173 53)
\ifill f:0
\move(192 52)
\lvec(195 52)
\lvec(195 53)
\lvec(192 53)
\ifill f:0
\move(200 52)
\lvec(203 52)
\lvec(203 53)
\lvec(200 53)
\ifill f:0
\move(206 52)
\lvec(209 52)
\lvec(209 53)
\lvec(206 53)
\ifill f:0
\move(212 52)
\lvec(214 52)
\lvec(214 53)
\lvec(212 53)
\ifill f:0
\move(217 52)
\lvec(219 52)
\lvec(219 53)
\lvec(217 53)
\ifill f:0
\move(221 52)
\lvec(223 52)
\lvec(223 53)
\lvec(221 53)
\ifill f:0
\move(225 52)
\lvec(227 52)
\lvec(227 53)
\lvec(225 53)
\ifill f:0
\move(229 52)
\lvec(231 52)
\lvec(231 53)
\lvec(229 53)
\ifill f:0
\move(233 52)
\lvec(234 52)
\lvec(234 53)
\lvec(233 53)
\ifill f:0
\move(236 52)
\lvec(238 52)
\lvec(238 53)
\lvec(236 53)
\ifill f:0
\move(240 52)
\lvec(241 52)
\lvec(241 53)
\lvec(240 53)
\ifill f:0
\move(243 52)
\lvec(244 52)
\lvec(244 53)
\lvec(243 53)
\ifill f:0
\move(246 52)
\lvec(248 52)
\lvec(248 53)
\lvec(246 53)
\ifill f:0
\move(249 52)
\lvec(251 52)
\lvec(251 53)
\lvec(249 53)
\ifill f:0
\move(252 52)
\lvec(253 52)
\lvec(253 53)
\lvec(252 53)
\ifill f:0
\move(255 52)
\lvec(256 52)
\lvec(256 53)
\lvec(255 53)
\ifill f:0
\move(258 52)
\lvec(259 52)
\lvec(259 53)
\lvec(258 53)
\ifill f:0
\move(260 52)
\lvec(262 52)
\lvec(262 53)
\lvec(260 53)
\ifill f:0
\move(263 52)
\lvec(264 52)
\lvec(264 53)
\lvec(263 53)
\ifill f:0
\move(266 52)
\lvec(267 52)
\lvec(267 53)
\lvec(266 53)
\ifill f:0
\move(268 52)
\lvec(269 52)
\lvec(269 53)
\lvec(268 53)
\ifill f:0
\move(271 52)
\lvec(272 52)
\lvec(272 53)
\lvec(271 53)
\ifill f:0
\move(273 52)
\lvec(274 52)
\lvec(274 53)
\lvec(273 53)
\ifill f:0
\move(276 52)
\lvec(277 52)
\lvec(277 53)
\lvec(276 53)
\ifill f:0
\move(278 52)
\lvec(279 52)
\lvec(279 53)
\lvec(278 53)
\ifill f:0
\move(280 52)
\lvec(281 52)
\lvec(281 53)
\lvec(280 53)
\ifill f:0
\move(283 52)
\lvec(284 52)
\lvec(284 53)
\lvec(283 53)
\ifill f:0
\move(285 52)
\lvec(286 52)
\lvec(286 53)
\lvec(285 53)
\ifill f:0
\move(287 52)
\lvec(288 52)
\lvec(288 53)
\lvec(287 53)
\ifill f:0
\move(289 52)
\lvec(290 52)
\lvec(290 53)
\lvec(289 53)
\ifill f:0
\move(292 52)
\lvec(293 52)
\lvec(293 53)
\lvec(292 53)
\ifill f:0
\move(294 52)
\lvec(295 52)
\lvec(295 53)
\lvec(294 53)
\ifill f:0
\move(296 52)
\lvec(297 52)
\lvec(297 53)
\lvec(296 53)
\ifill f:0
\move(298 52)
\lvec(299 52)
\lvec(299 53)
\lvec(298 53)
\ifill f:0
\move(300 52)
\lvec(301 52)
\lvec(301 53)
\lvec(300 53)
\ifill f:0
\move(302 52)
\lvec(303 52)
\lvec(303 53)
\lvec(302 53)
\ifill f:0
\move(304 52)
\lvec(305 52)
\lvec(305 53)
\lvec(304 53)
\ifill f:0
\move(306 52)
\lvec(307 52)
\lvec(307 53)
\lvec(306 53)
\ifill f:0
\move(308 52)
\lvec(309 52)
\lvec(309 53)
\lvec(308 53)
\ifill f:0
\move(310 52)
\lvec(311 52)
\lvec(311 53)
\lvec(310 53)
\ifill f:0
\move(312 52)
\lvec(313 52)
\lvec(313 53)
\lvec(312 53)
\ifill f:0
\move(314 52)
\lvec(315 52)
\lvec(315 53)
\lvec(314 53)
\ifill f:0
\move(316 52)
\lvec(317 52)
\lvec(317 53)
\lvec(316 53)
\ifill f:0
\move(318 52)
\lvec(319 52)
\lvec(319 53)
\lvec(318 53)
\ifill f:0
\move(320 52)
\lvec(321 52)
\lvec(321 53)
\lvec(320 53)
\ifill f:0
\move(322 52)
\lvec(324 52)
\lvec(324 53)
\lvec(322 53)
\ifill f:0
\move(325 52)
\lvec(326 52)
\lvec(326 53)
\lvec(325 53)
\ifill f:0
\move(327 52)
\lvec(328 52)
\lvec(328 53)
\lvec(327 53)
\ifill f:0
\move(329 52)
\lvec(330 52)
\lvec(330 53)
\lvec(329 53)
\ifill f:0
\move(331 52)
\lvec(332 52)
\lvec(332 53)
\lvec(331 53)
\ifill f:0
\move(333 52)
\lvec(335 52)
\lvec(335 53)
\lvec(333 53)
\ifill f:0
\move(336 52)
\lvec(337 52)
\lvec(337 53)
\lvec(336 53)
\ifill f:0
\move(338 52)
\lvec(339 52)
\lvec(339 53)
\lvec(338 53)
\ifill f:0
\move(340 52)
\lvec(342 52)
\lvec(342 53)
\lvec(340 53)
\ifill f:0
\move(343 52)
\lvec(344 52)
\lvec(344 53)
\lvec(343 53)
\ifill f:0
\move(345 52)
\lvec(346 52)
\lvec(346 53)
\lvec(345 53)
\ifill f:0
\move(348 52)
\lvec(349 52)
\lvec(349 53)
\lvec(348 53)
\ifill f:0
\move(350 52)
\lvec(351 52)
\lvec(351 53)
\lvec(350 53)
\ifill f:0
\move(353 52)
\lvec(354 52)
\lvec(354 53)
\lvec(353 53)
\ifill f:0
\move(355 52)
\lvec(356 52)
\lvec(356 53)
\lvec(355 53)
\ifill f:0
\move(358 52)
\lvec(359 52)
\lvec(359 53)
\lvec(358 53)
\ifill f:0
\move(360 52)
\lvec(362 52)
\lvec(362 53)
\lvec(360 53)
\ifill f:0
\move(363 52)
\lvec(364 52)
\lvec(364 53)
\lvec(363 53)
\ifill f:0
\move(365 52)
\lvec(367 52)
\lvec(367 53)
\lvec(365 53)
\ifill f:0
\move(368 52)
\lvec(370 52)
\lvec(370 53)
\lvec(368 53)
\ifill f:0
\move(371 52)
\lvec(372 52)
\lvec(372 53)
\lvec(371 53)
\ifill f:0
\move(374 52)
\lvec(375 52)
\lvec(375 53)
\lvec(374 53)
\ifill f:0
\move(376 52)
\lvec(378 52)
\lvec(378 53)
\lvec(376 53)
\ifill f:0
\move(379 52)
\lvec(381 52)
\lvec(381 53)
\lvec(379 53)
\ifill f:0
\move(382 52)
\lvec(384 52)
\lvec(384 53)
\lvec(382 53)
\ifill f:0
\move(385 52)
\lvec(387 52)
\lvec(387 53)
\lvec(385 53)
\ifill f:0
\move(388 52)
\lvec(390 52)
\lvec(390 53)
\lvec(388 53)
\ifill f:0
\move(391 52)
\lvec(393 52)
\lvec(393 53)
\lvec(391 53)
\ifill f:0
\move(394 52)
\lvec(396 52)
\lvec(396 53)
\lvec(394 53)
\ifill f:0
\move(397 52)
\lvec(399 52)
\lvec(399 53)
\lvec(397 53)
\ifill f:0
\move(400 52)
\lvec(402 52)
\lvec(402 53)
\lvec(400 53)
\ifill f:0
\move(404 52)
\lvec(406 52)
\lvec(406 53)
\lvec(404 53)
\ifill f:0
\move(407 52)
\lvec(409 52)
\lvec(409 53)
\lvec(407 53)
\ifill f:0
\move(410 52)
\lvec(412 52)
\lvec(412 53)
\lvec(410 53)
\ifill f:0
\move(414 52)
\lvec(416 52)
\lvec(416 53)
\lvec(414 53)
\ifill f:0
\move(417 52)
\lvec(419 52)
\lvec(419 53)
\lvec(417 53)
\ifill f:0
\move(421 52)
\lvec(423 52)
\lvec(423 53)
\lvec(421 53)
\ifill f:0
\move(424 52)
\lvec(426 52)
\lvec(426 53)
\lvec(424 53)
\ifill f:0
\move(428 52)
\lvec(430 52)
\lvec(430 53)
\lvec(428 53)
\ifill f:0
\move(432 52)
\lvec(434 52)
\lvec(434 53)
\lvec(432 53)
\ifill f:0
\move(435 52)
\lvec(438 52)
\lvec(438 53)
\lvec(435 53)
\ifill f:0
\move(439 52)
\lvec(442 52)
\lvec(442 53)
\lvec(439 53)
\ifill f:0
\move(443 52)
\lvec(446 52)
\lvec(446 53)
\lvec(443 53)
\ifill f:0
\move(447 52)
\lvec(450 52)
\lvec(450 53)
\lvec(447 53)
\ifill f:0
\move(12 53)
\lvec(13 53)
\lvec(13 54)
\lvec(12 54)
\ifill f:0
\move(15 53)
\lvec(16 53)
\lvec(16 54)
\lvec(15 54)
\ifill f:0
\move(19 53)
\lvec(20 53)
\lvec(20 54)
\lvec(19 54)
\ifill f:0
\move(79 53)
\lvec(81 53)
\lvec(81 54)
\lvec(79 54)
\ifill f:0
\move(95 53)
\lvec(96 53)
\lvec(96 54)
\lvec(95 54)
\ifill f:0
\move(107 53)
\lvec(108 53)
\lvec(108 54)
\lvec(107 54)
\ifill f:0
\move(121 53)
\lvec(122 53)
\lvec(122 54)
\lvec(121 54)
\ifill f:0
\move(123 53)
\lvec(124 53)
\lvec(124 54)
\lvec(123 54)
\ifill f:0
\move(131 53)
\lvec(132 53)
\lvec(132 54)
\lvec(131 54)
\ifill f:0
\move(134 53)
\lvec(135 53)
\lvec(135 54)
\lvec(134 54)
\ifill f:0
\move(138 53)
\lvec(139 53)
\lvec(139 54)
\lvec(138 54)
\ifill f:0
\move(141 53)
\lvec(143 53)
\lvec(143 54)
\lvec(141 54)
\ifill f:0
\move(146 53)
\lvec(147 53)
\lvec(147 54)
\lvec(146 54)
\ifill f:0
\move(150 53)
\lvec(152 53)
\lvec(152 54)
\lvec(150 54)
\ifill f:0
\move(157 53)
\lvec(159 53)
\lvec(159 54)
\lvec(157 54)
\ifill f:0
\move(165 53)
\lvec(170 53)
\lvec(170 54)
\lvec(165 54)
\ifill f:0
\move(189 53)
\lvec(195 53)
\lvec(195 54)
\lvec(189 54)
\ifill f:0
\move(201 53)
\lvec(204 53)
\lvec(204 54)
\lvec(201 54)
\ifill f:0
\move(209 53)
\lvec(211 53)
\lvec(211 54)
\lvec(209 54)
\ifill f:0
\move(215 53)
\lvec(217 53)
\lvec(217 54)
\lvec(215 54)
\ifill f:0
\move(220 53)
\lvec(223 53)
\lvec(223 54)
\lvec(220 54)
\ifill f:0
\move(225 53)
\lvec(227 53)
\lvec(227 54)
\lvec(225 54)
\ifill f:0
\move(230 53)
\lvec(232 53)
\lvec(232 54)
\lvec(230 54)
\ifill f:0
\move(234 53)
\lvec(236 53)
\lvec(236 54)
\lvec(234 54)
\ifill f:0
\move(238 53)
\lvec(239 53)
\lvec(239 54)
\lvec(238 54)
\ifill f:0
\move(242 53)
\lvec(243 53)
\lvec(243 54)
\lvec(242 54)
\ifill f:0
\move(245 53)
\lvec(247 53)
\lvec(247 54)
\lvec(245 54)
\ifill f:0
\move(248 53)
\lvec(250 53)
\lvec(250 54)
\lvec(248 54)
\ifill f:0
\move(252 53)
\lvec(253 53)
\lvec(253 54)
\lvec(252 54)
\ifill f:0
\move(255 53)
\lvec(256 53)
\lvec(256 54)
\lvec(255 54)
\ifill f:0
\move(258 53)
\lvec(259 53)
\lvec(259 54)
\lvec(258 54)
\ifill f:0
\move(261 53)
\lvec(262 53)
\lvec(262 54)
\lvec(261 54)
\ifill f:0
\move(264 53)
\lvec(265 53)
\lvec(265 54)
\lvec(264 54)
\ifill f:0
\move(267 53)
\lvec(268 53)
\lvec(268 54)
\lvec(267 54)
\ifill f:0
\move(269 53)
\lvec(271 53)
\lvec(271 54)
\lvec(269 54)
\ifill f:0
\move(272 53)
\lvec(273 53)
\lvec(273 54)
\lvec(272 54)
\ifill f:0
\move(275 53)
\lvec(276 53)
\lvec(276 54)
\lvec(275 54)
\ifill f:0
\move(277 53)
\lvec(278 53)
\lvec(278 54)
\lvec(277 54)
\ifill f:0
\move(280 53)
\lvec(281 53)
\lvec(281 54)
\lvec(280 54)
\ifill f:0
\move(282 53)
\lvec(283 53)
\lvec(283 54)
\lvec(282 54)
\ifill f:0
\move(285 53)
\lvec(286 53)
\lvec(286 54)
\lvec(285 54)
\ifill f:0
\move(287 53)
\lvec(288 53)
\lvec(288 54)
\lvec(287 54)
\ifill f:0
\move(289 53)
\lvec(291 53)
\lvec(291 54)
\lvec(289 54)
\ifill f:0
\move(292 53)
\lvec(293 53)
\lvec(293 54)
\lvec(292 54)
\ifill f:0
\move(294 53)
\lvec(295 53)
\lvec(295 54)
\lvec(294 54)
\ifill f:0
\move(296 53)
\lvec(297 53)
\lvec(297 54)
\lvec(296 54)
\ifill f:0
\move(299 53)
\lvec(300 53)
\lvec(300 54)
\lvec(299 54)
\ifill f:0
\move(301 53)
\lvec(302 53)
\lvec(302 54)
\lvec(301 54)
\ifill f:0
\move(303 53)
\lvec(304 53)
\lvec(304 54)
\lvec(303 54)
\ifill f:0
\move(305 53)
\lvec(306 53)
\lvec(306 54)
\lvec(305 54)
\ifill f:0
\move(307 53)
\lvec(308 53)
\lvec(308 54)
\lvec(307 54)
\ifill f:0
\move(309 53)
\lvec(310 53)
\lvec(310 54)
\lvec(309 54)
\ifill f:0
\move(311 53)
\lvec(312 53)
\lvec(312 54)
\lvec(311 54)
\ifill f:0
\move(313 53)
\lvec(314 53)
\lvec(314 54)
\lvec(313 54)
\ifill f:0
\move(315 53)
\lvec(322 53)
\lvec(322 54)
\lvec(315 54)
\ifill f:0
\move(323 53)
\lvec(324 53)
\lvec(324 54)
\lvec(323 54)
\ifill f:0
\move(325 53)
\lvec(326 53)
\lvec(326 54)
\lvec(325 54)
\ifill f:0
\move(327 53)
\lvec(328 53)
\lvec(328 54)
\lvec(327 54)
\ifill f:0
\move(329 53)
\lvec(330 53)
\lvec(330 54)
\lvec(329 54)
\ifill f:0
\move(331 53)
\lvec(332 53)
\lvec(332 54)
\lvec(331 54)
\ifill f:0
\move(333 53)
\lvec(334 53)
\lvec(334 54)
\lvec(333 54)
\ifill f:0
\move(335 53)
\lvec(336 53)
\lvec(336 54)
\lvec(335 54)
\ifill f:0
\move(337 53)
\lvec(338 53)
\lvec(338 54)
\lvec(337 54)
\ifill f:0
\move(339 53)
\lvec(341 53)
\lvec(341 54)
\lvec(339 54)
\ifill f:0
\move(342 53)
\lvec(343 53)
\lvec(343 54)
\lvec(342 54)
\ifill f:0
\move(344 53)
\lvec(345 53)
\lvec(345 54)
\lvec(344 54)
\ifill f:0
\move(346 53)
\lvec(347 53)
\lvec(347 54)
\lvec(346 54)
\ifill f:0
\move(348 53)
\lvec(350 53)
\lvec(350 54)
\lvec(348 54)
\ifill f:0
\move(351 53)
\lvec(352 53)
\lvec(352 54)
\lvec(351 54)
\ifill f:0
\move(353 53)
\lvec(354 53)
\lvec(354 54)
\lvec(353 54)
\ifill f:0
\move(355 53)
\lvec(357 53)
\lvec(357 54)
\lvec(355 54)
\ifill f:0
\move(358 53)
\lvec(359 53)
\lvec(359 54)
\lvec(358 54)
\ifill f:0
\move(360 53)
\lvec(362 53)
\lvec(362 54)
\lvec(360 54)
\ifill f:0
\move(363 53)
\lvec(364 53)
\lvec(364 54)
\lvec(363 54)
\ifill f:0
\move(365 53)
\lvec(367 53)
\lvec(367 54)
\lvec(365 54)
\ifill f:0
\move(368 53)
\lvec(369 53)
\lvec(369 54)
\lvec(368 54)
\ifill f:0
\move(370 53)
\lvec(372 53)
\lvec(372 54)
\lvec(370 54)
\ifill f:0
\move(373 53)
\lvec(374 53)
\lvec(374 54)
\lvec(373 54)
\ifill f:0
\move(375 53)
\lvec(377 53)
\lvec(377 54)
\lvec(375 54)
\ifill f:0
\move(378 53)
\lvec(380 53)
\lvec(380 54)
\lvec(378 54)
\ifill f:0
\move(381 53)
\lvec(382 53)
\lvec(382 54)
\lvec(381 54)
\ifill f:0
\move(383 53)
\lvec(385 53)
\lvec(385 54)
\lvec(383 54)
\ifill f:0
\move(386 53)
\lvec(388 53)
\lvec(388 54)
\lvec(386 54)
\ifill f:0
\move(389 53)
\lvec(391 53)
\lvec(391 54)
\lvec(389 54)
\ifill f:0
\move(392 53)
\lvec(393 53)
\lvec(393 54)
\lvec(392 54)
\ifill f:0
\move(395 53)
\lvec(396 53)
\lvec(396 54)
\lvec(395 54)
\ifill f:0
\move(398 53)
\lvec(399 53)
\lvec(399 54)
\lvec(398 54)
\ifill f:0
\move(400 53)
\lvec(402 53)
\lvec(402 54)
\lvec(400 54)
\ifill f:0
\move(403 53)
\lvec(405 53)
\lvec(405 54)
\lvec(403 54)
\ifill f:0
\move(406 53)
\lvec(408 53)
\lvec(408 54)
\lvec(406 54)
\ifill f:0
\move(410 53)
\lvec(411 53)
\lvec(411 54)
\lvec(410 54)
\ifill f:0
\move(413 53)
\lvec(414 53)
\lvec(414 54)
\lvec(413 54)
\ifill f:0
\move(416 53)
\lvec(418 53)
\lvec(418 54)
\lvec(416 54)
\ifill f:0
\move(419 53)
\lvec(421 53)
\lvec(421 54)
\lvec(419 54)
\ifill f:0
\move(422 53)
\lvec(424 53)
\lvec(424 54)
\lvec(422 54)
\ifill f:0
\move(426 53)
\lvec(428 53)
\lvec(428 54)
\lvec(426 54)
\ifill f:0
\move(429 53)
\lvec(431 53)
\lvec(431 54)
\lvec(429 54)
\ifill f:0
\move(432 53)
\lvec(435 53)
\lvec(435 54)
\lvec(432 54)
\ifill f:0
\move(436 53)
\lvec(438 53)
\lvec(438 54)
\lvec(436 54)
\ifill f:0
\move(440 53)
\lvec(442 53)
\lvec(442 54)
\lvec(440 54)
\ifill f:0
\move(443 53)
\lvec(445 53)
\lvec(445 54)
\lvec(443 54)
\ifill f:0
\move(447 53)
\lvec(449 53)
\lvec(449 54)
\lvec(447 54)
\ifill f:0
\move(13 54)
\lvec(14 54)
\lvec(14 55)
\lvec(13 55)
\ifill f:0
\move(15 54)
\lvec(16 54)
\lvec(16 55)
\lvec(15 55)
\ifill f:0
\move(82 54)
\lvec(84 54)
\lvec(84 55)
\lvec(82 55)
\ifill f:0
\move(102 54)
\lvec(103 54)
\lvec(103 55)
\lvec(102 55)
\ifill f:0
\move(105 54)
\lvec(106 54)
\lvec(106 55)
\lvec(105 55)
\ifill f:0
\move(108 54)
\lvec(109 54)
\lvec(109 55)
\lvec(108 55)
\ifill f:0
\move(115 54)
\lvec(116 54)
\lvec(116 55)
\lvec(115 55)
\ifill f:0
\move(117 54)
\lvec(118 54)
\lvec(118 55)
\lvec(117 55)
\ifill f:0
\move(121 54)
\lvec(122 54)
\lvec(122 55)
\lvec(121 55)
\ifill f:0
\move(123 54)
\lvec(124 54)
\lvec(124 55)
\lvec(123 55)
\ifill f:0
\move(125 54)
\lvec(126 54)
\lvec(126 55)
\lvec(125 55)
\ifill f:0
\move(130 54)
\lvec(131 54)
\lvec(131 55)
\lvec(130 55)
\ifill f:0
\move(133 54)
\lvec(134 54)
\lvec(134 55)
\lvec(133 55)
\ifill f:0
\move(142 54)
\lvec(143 54)
\lvec(143 55)
\lvec(142 55)
\ifill f:0
\move(145 54)
\lvec(147 54)
\lvec(147 55)
\lvec(145 55)
\ifill f:0
\move(149 54)
\lvec(151 54)
\lvec(151 55)
\lvec(149 55)
\ifill f:0
\move(154 54)
\lvec(156 54)
\lvec(156 55)
\lvec(154 55)
\ifill f:0
\move(160 54)
\lvec(162 54)
\lvec(162 55)
\lvec(160 55)
\ifill f:0
\move(167 54)
\lvec(170 54)
\lvec(170 55)
\lvec(167 55)
\ifill f:0
\move(181 54)
\lvec(193 54)
\lvec(193 55)
\lvec(181 55)
\ifill f:0
\move(203 54)
\lvec(207 54)
\lvec(207 55)
\lvec(203 55)
\ifill f:0
\move(212 54)
\lvec(215 54)
\lvec(215 55)
\lvec(212 55)
\ifill f:0
\move(219 54)
\lvec(222 54)
\lvec(222 55)
\lvec(219 55)
\ifill f:0
\move(225 54)
\lvec(227 54)
\lvec(227 55)
\lvec(225 55)
\ifill f:0
\move(230 54)
\lvec(232 54)
\lvec(232 55)
\lvec(230 55)
\ifill f:0
\move(235 54)
\lvec(237 54)
\lvec(237 55)
\lvec(235 55)
\ifill f:0
\move(240 54)
\lvec(241 54)
\lvec(241 55)
\lvec(240 55)
\ifill f:0
\move(244 54)
\lvec(245 54)
\lvec(245 55)
\lvec(244 55)
\ifill f:0
\move(247 54)
\lvec(249 54)
\lvec(249 55)
\lvec(247 55)
\ifill f:0
\move(251 54)
\lvec(253 54)
\lvec(253 55)
\lvec(251 55)
\ifill f:0
\move(255 54)
\lvec(256 54)
\lvec(256 55)
\lvec(255 55)
\ifill f:0
\move(258 54)
\lvec(259 54)
\lvec(259 55)
\lvec(258 55)
\ifill f:0
\move(261 54)
\lvec(263 54)
\lvec(263 55)
\lvec(261 55)
\ifill f:0
\move(264 54)
\lvec(266 54)
\lvec(266 55)
\lvec(264 55)
\ifill f:0
\move(268 54)
\lvec(269 54)
\lvec(269 55)
\lvec(268 55)
\ifill f:0
\move(270 54)
\lvec(272 54)
\lvec(272 55)
\lvec(270 55)
\ifill f:0
\move(273 54)
\lvec(275 54)
\lvec(275 55)
\lvec(273 55)
\ifill f:0
\move(276 54)
\lvec(277 54)
\lvec(277 55)
\lvec(276 55)
\ifill f:0
\move(279 54)
\lvec(280 54)
\lvec(280 55)
\lvec(279 55)
\ifill f:0
\move(282 54)
\lvec(283 54)
\lvec(283 55)
\lvec(282 55)
\ifill f:0
\move(284 54)
\lvec(285 54)
\lvec(285 55)
\lvec(284 55)
\ifill f:0
\move(287 54)
\lvec(288 54)
\lvec(288 55)
\lvec(287 55)
\ifill f:0
\move(289 54)
\lvec(291 54)
\lvec(291 55)
\lvec(289 55)
\ifill f:0
\move(292 54)
\lvec(293 54)
\lvec(293 55)
\lvec(292 55)
\ifill f:0
\move(294 54)
\lvec(296 54)
\lvec(296 55)
\lvec(294 55)
\ifill f:0
\move(297 54)
\lvec(298 54)
\lvec(298 55)
\lvec(297 55)
\ifill f:0
\move(299 54)
\lvec(300 54)
\lvec(300 55)
\lvec(299 55)
\ifill f:0
\move(302 54)
\lvec(303 54)
\lvec(303 55)
\lvec(302 55)
\ifill f:0
\move(304 54)
\lvec(305 54)
\lvec(305 55)
\lvec(304 55)
\ifill f:0
\move(306 54)
\lvec(307 54)
\lvec(307 55)
\lvec(306 55)
\ifill f:0
\move(308 54)
\lvec(309 54)
\lvec(309 55)
\lvec(308 55)
\ifill f:0
\move(311 54)
\lvec(312 54)
\lvec(312 55)
\lvec(311 55)
\ifill f:0
\move(313 54)
\lvec(314 54)
\lvec(314 55)
\lvec(313 55)
\ifill f:0
\move(315 54)
\lvec(316 54)
\lvec(316 55)
\lvec(315 55)
\ifill f:0
\move(317 54)
\lvec(318 54)
\lvec(318 55)
\lvec(317 55)
\ifill f:0
\move(319 54)
\lvec(320 54)
\lvec(320 55)
\lvec(319 55)
\ifill f:0
\move(321 54)
\lvec(322 54)
\lvec(322 55)
\lvec(321 55)
\ifill f:0
\move(323 54)
\lvec(324 54)
\lvec(324 55)
\lvec(323 55)
\ifill f:0
\move(325 54)
\lvec(326 54)
\lvec(326 55)
\lvec(325 55)
\ifill f:0
\move(327 54)
\lvec(328 54)
\lvec(328 55)
\lvec(327 55)
\ifill f:0
\move(329 54)
\lvec(334 54)
\lvec(334 55)
\lvec(329 55)
\ifill f:0
\move(335 54)
\lvec(336 54)
\lvec(336 55)
\lvec(335 55)
\ifill f:0
\move(337 54)
\lvec(338 54)
\lvec(338 55)
\lvec(337 55)
\ifill f:0
\move(339 54)
\lvec(340 54)
\lvec(340 55)
\lvec(339 55)
\ifill f:0
\move(341 54)
\lvec(342 54)
\lvec(342 55)
\lvec(341 55)
\ifill f:0
\move(343 54)
\lvec(344 54)
\lvec(344 55)
\lvec(343 55)
\ifill f:0
\move(345 54)
\lvec(346 54)
\lvec(346 55)
\lvec(345 55)
\ifill f:0
\move(347 54)
\lvec(348 54)
\lvec(348 55)
\lvec(347 55)
\ifill f:0
\move(349 54)
\lvec(350 54)
\lvec(350 55)
\lvec(349 55)
\ifill f:0
\move(351 54)
\lvec(353 54)
\lvec(353 55)
\lvec(351 55)
\ifill f:0
\move(354 54)
\lvec(355 54)
\lvec(355 55)
\lvec(354 55)
\ifill f:0
\move(356 54)
\lvec(357 54)
\lvec(357 55)
\lvec(356 55)
\ifill f:0
\move(358 54)
\lvec(359 54)
\lvec(359 55)
\lvec(358 55)
\ifill f:0
\move(360 54)
\lvec(362 54)
\lvec(362 55)
\lvec(360 55)
\ifill f:0
\move(363 54)
\lvec(364 54)
\lvec(364 55)
\lvec(363 55)
\ifill f:0
\move(365 54)
\lvec(366 54)
\lvec(366 55)
\lvec(365 55)
\ifill f:0
\move(367 54)
\lvec(369 54)
\lvec(369 55)
\lvec(367 55)
\ifill f:0
\move(370 54)
\lvec(371 54)
\lvec(371 55)
\lvec(370 55)
\ifill f:0
\move(372 54)
\lvec(373 54)
\lvec(373 55)
\lvec(372 55)
\ifill f:0
\move(375 54)
\lvec(376 54)
\lvec(376 55)
\lvec(375 55)
\ifill f:0
\move(377 54)
\lvec(378 54)
\lvec(378 55)
\lvec(377 55)
\ifill f:0
\move(380 54)
\lvec(381 54)
\lvec(381 55)
\lvec(380 55)
\ifill f:0
\move(382 54)
\lvec(383 54)
\lvec(383 55)
\lvec(382 55)
\ifill f:0
\move(385 54)
\lvec(386 54)
\lvec(386 55)
\lvec(385 55)
\ifill f:0
\move(387 54)
\lvec(389 54)
\lvec(389 55)
\lvec(387 55)
\ifill f:0
\move(390 54)
\lvec(391 54)
\lvec(391 55)
\lvec(390 55)
\ifill f:0
\move(392 54)
\lvec(394 54)
\lvec(394 55)
\lvec(392 55)
\ifill f:0
\move(395 54)
\lvec(397 54)
\lvec(397 55)
\lvec(395 55)
\ifill f:0
\move(398 54)
\lvec(399 54)
\lvec(399 55)
\lvec(398 55)
\ifill f:0
\move(400 54)
\lvec(402 54)
\lvec(402 55)
\lvec(400 55)
\ifill f:0
\move(403 54)
\lvec(405 54)
\lvec(405 55)
\lvec(403 55)
\ifill f:0
\move(406 54)
\lvec(408 54)
\lvec(408 55)
\lvec(406 55)
\ifill f:0
\move(409 54)
\lvec(411 54)
\lvec(411 55)
\lvec(409 55)
\ifill f:0
\move(412 54)
\lvec(413 54)
\lvec(413 55)
\lvec(412 55)
\ifill f:0
\move(415 54)
\lvec(416 54)
\lvec(416 55)
\lvec(415 55)
\ifill f:0
\move(418 54)
\lvec(419 54)
\lvec(419 55)
\lvec(418 55)
\ifill f:0
\move(421 54)
\lvec(422 54)
\lvec(422 55)
\lvec(421 55)
\ifill f:0
\move(424 54)
\lvec(426 54)
\lvec(426 55)
\lvec(424 55)
\ifill f:0
\move(427 54)
\lvec(429 54)
\lvec(429 55)
\lvec(427 55)
\ifill f:0
\move(430 54)
\lvec(432 54)
\lvec(432 55)
\lvec(430 55)
\ifill f:0
\move(433 54)
\lvec(435 54)
\lvec(435 55)
\lvec(433 55)
\ifill f:0
\move(436 54)
\lvec(438 54)
\lvec(438 55)
\lvec(436 55)
\ifill f:0
\move(440 54)
\lvec(442 54)
\lvec(442 55)
\lvec(440 55)
\ifill f:0
\move(443 54)
\lvec(445 54)
\lvec(445 55)
\lvec(443 55)
\ifill f:0
\move(447 54)
\lvec(449 54)
\lvec(449 55)
\lvec(447 55)
\ifill f:0
\move(450 54)
\lvec(451 54)
\lvec(451 55)
\lvec(450 55)
\ifill f:0
\move(12 55)
\lvec(13 55)
\lvec(13 56)
\lvec(12 56)
\ifill f:0
\move(90 55)
\lvec(91 55)
\lvec(91 56)
\lvec(90 56)
\ifill f:0
\move(119 55)
\lvec(120 55)
\lvec(120 56)
\lvec(119 56)
\ifill f:0
\move(121 55)
\lvec(122 55)
\lvec(122 56)
\lvec(121 56)
\ifill f:0
\move(125 55)
\lvec(126 55)
\lvec(126 56)
\lvec(125 56)
\ifill f:0
\move(127 55)
\lvec(128 55)
\lvec(128 56)
\lvec(127 56)
\ifill f:0
\move(134 55)
\lvec(135 55)
\lvec(135 56)
\lvec(134 56)
\ifill f:0
\move(142 55)
\lvec(143 55)
\lvec(143 56)
\lvec(142 56)
\ifill f:0
\move(149 55)
\lvec(150 55)
\lvec(150 56)
\lvec(149 56)
\ifill f:0
\move(153 55)
\lvec(154 55)
\lvec(154 56)
\lvec(153 56)
\ifill f:0
\move(157 55)
\lvec(158 55)
\lvec(158 56)
\lvec(157 56)
\ifill f:0
\move(162 55)
\lvec(163 55)
\lvec(163 56)
\lvec(162 56)
\ifill f:0
\move(168 55)
\lvec(170 55)
\lvec(170 56)
\lvec(168 56)
\ifill f:0
\move(176 55)
\lvec(179 55)
\lvec(179 56)
\lvec(176 56)
\ifill f:0
\move(207 55)
\lvec(211 55)
\lvec(211 56)
\lvec(207 56)
\ifill f:0
\move(218 55)
\lvec(221 55)
\lvec(221 56)
\lvec(218 56)
\ifill f:0
\move(225 55)
\lvec(228 55)
\lvec(228 56)
\lvec(225 56)
\ifill f:0
\move(231 55)
\lvec(234 55)
\lvec(234 56)
\lvec(231 56)
\ifill f:0
\move(237 55)
\lvec(239 55)
\lvec(239 56)
\lvec(237 56)
\ifill f:0
\move(242 55)
\lvec(244 55)
\lvec(244 56)
\lvec(242 56)
\ifill f:0
\move(246 55)
\lvec(248 55)
\lvec(248 56)
\lvec(246 56)
\ifill f:0
\move(250 55)
\lvec(252 55)
\lvec(252 56)
\lvec(250 56)
\ifill f:0
\move(254 55)
\lvec(256 55)
\lvec(256 56)
\lvec(254 56)
\ifill f:0
\move(258 55)
\lvec(260 55)
\lvec(260 56)
\lvec(258 56)
\ifill f:0
\move(262 55)
\lvec(263 55)
\lvec(263 56)
\lvec(262 56)
\ifill f:0
\move(265 55)
\lvec(267 55)
\lvec(267 56)
\lvec(265 56)
\ifill f:0
\move(269 55)
\lvec(270 55)
\lvec(270 56)
\lvec(269 56)
\ifill f:0
\move(272 55)
\lvec(273 55)
\lvec(273 56)
\lvec(272 56)
\ifill f:0
\move(275 55)
\lvec(276 55)
\lvec(276 56)
\lvec(275 56)
\ifill f:0
\move(278 55)
\lvec(279 55)
\lvec(279 56)
\lvec(278 56)
\ifill f:0
\move(281 55)
\lvec(282 55)
\lvec(282 56)
\lvec(281 56)
\ifill f:0
\move(284 55)
\lvec(285 55)
\lvec(285 56)
\lvec(284 56)
\ifill f:0
\move(287 55)
\lvec(288 55)
\lvec(288 56)
\lvec(287 56)
\ifill f:0
\move(289 55)
\lvec(291 55)
\lvec(291 56)
\lvec(289 56)
\ifill f:0
\move(292 55)
\lvec(293 55)
\lvec(293 56)
\lvec(292 56)
\ifill f:0
\move(295 55)
\lvec(296 55)
\lvec(296 56)
\lvec(295 56)
\ifill f:0
\move(297 55)
\lvec(298 55)
\lvec(298 56)
\lvec(297 56)
\ifill f:0
\move(300 55)
\lvec(301 55)
\lvec(301 56)
\lvec(300 56)
\ifill f:0
\move(302 55)
\lvec(303 55)
\lvec(303 56)
\lvec(302 56)
\ifill f:0
\move(305 55)
\lvec(306 55)
\lvec(306 56)
\lvec(305 56)
\ifill f:0
\move(307 55)
\lvec(308 55)
\lvec(308 56)
\lvec(307 56)
\ifill f:0
\move(310 55)
\lvec(311 55)
\lvec(311 56)
\lvec(310 56)
\ifill f:0
\move(312 55)
\lvec(313 55)
\lvec(313 56)
\lvec(312 56)
\ifill f:0
\move(314 55)
\lvec(315 55)
\lvec(315 56)
\lvec(314 56)
\ifill f:0
\move(317 55)
\lvec(318 55)
\lvec(318 56)
\lvec(317 56)
\ifill f:0
\move(319 55)
\lvec(320 55)
\lvec(320 56)
\lvec(319 56)
\ifill f:0
\move(321 55)
\lvec(322 55)
\lvec(322 56)
\lvec(321 56)
\ifill f:0
\move(323 55)
\lvec(324 55)
\lvec(324 56)
\lvec(323 56)
\ifill f:0
\move(326 55)
\lvec(327 55)
\lvec(327 56)
\lvec(326 56)
\ifill f:0
\move(328 55)
\lvec(329 55)
\lvec(329 56)
\lvec(328 56)
\ifill f:0
\move(330 55)
\lvec(331 55)
\lvec(331 56)
\lvec(330 56)
\ifill f:0
\move(332 55)
\lvec(333 55)
\lvec(333 56)
\lvec(332 56)
\ifill f:0
\move(334 55)
\lvec(335 55)
\lvec(335 56)
\lvec(334 56)
\ifill f:0
\move(336 55)
\lvec(337 55)
\lvec(337 56)
\lvec(336 56)
\ifill f:0
\move(338 55)
\lvec(339 55)
\lvec(339 56)
\lvec(338 56)
\ifill f:0
\move(340 55)
\lvec(341 55)
\lvec(341 56)
\lvec(340 56)
\ifill f:0
\move(342 55)
\lvec(343 55)
\lvec(343 56)
\lvec(342 56)
\ifill f:0
\move(344 55)
\lvec(345 55)
\lvec(345 56)
\lvec(344 56)
\ifill f:0
\move(346 55)
\lvec(347 55)
\lvec(347 56)
\lvec(346 56)
\ifill f:0
\move(348 55)
\lvec(349 55)
\lvec(349 56)
\lvec(348 56)
\ifill f:0
\move(350 55)
\lvec(351 55)
\lvec(351 56)
\lvec(350 56)
\ifill f:0
\move(352 55)
\lvec(353 55)
\lvec(353 56)
\lvec(352 56)
\ifill f:0
\move(354 55)
\lvec(355 55)
\lvec(355 56)
\lvec(354 56)
\ifill f:0
\move(356 55)
\lvec(357 55)
\lvec(357 56)
\lvec(356 56)
\ifill f:0
\move(358 55)
\lvec(359 55)
\lvec(359 56)
\lvec(358 56)
\ifill f:0
\move(361 55)
\lvec(362 55)
\lvec(362 56)
\lvec(361 56)
\ifill f:0
\move(363 55)
\lvec(364 55)
\lvec(364 56)
\lvec(363 56)
\ifill f:0
\move(365 55)
\lvec(366 55)
\lvec(366 56)
\lvec(365 56)
\ifill f:0
\move(367 55)
\lvec(368 55)
\lvec(368 56)
\lvec(367 56)
\ifill f:0
\move(369 55)
\lvec(370 55)
\lvec(370 56)
\lvec(369 56)
\ifill f:0
\move(372 55)
\lvec(373 55)
\lvec(373 56)
\lvec(372 56)
\ifill f:0
\move(374 55)
\lvec(375 55)
\lvec(375 56)
\lvec(374 56)
\ifill f:0
\move(376 55)
\lvec(377 55)
\lvec(377 56)
\lvec(376 56)
\ifill f:0
\move(378 55)
\lvec(380 55)
\lvec(380 56)
\lvec(378 56)
\ifill f:0
\move(381 55)
\lvec(382 55)
\lvec(382 56)
\lvec(381 56)
\ifill f:0
\move(383 55)
\lvec(384 55)
\lvec(384 56)
\lvec(383 56)
\ifill f:0
\move(386 55)
\lvec(387 55)
\lvec(387 56)
\lvec(386 56)
\ifill f:0
\move(388 55)
\lvec(389 55)
\lvec(389 56)
\lvec(388 56)
\ifill f:0
\move(390 55)
\lvec(392 55)
\lvec(392 56)
\lvec(390 56)
\ifill f:0
\move(393 55)
\lvec(394 55)
\lvec(394 56)
\lvec(393 56)
\ifill f:0
\move(395 55)
\lvec(397 55)
\lvec(397 56)
\lvec(395 56)
\ifill f:0
\move(398 55)
\lvec(399 55)
\lvec(399 56)
\lvec(398 56)
\ifill f:0
\move(401 55)
\lvec(402 55)
\lvec(402 56)
\lvec(401 56)
\ifill f:0
\move(403 55)
\lvec(405 55)
\lvec(405 56)
\lvec(403 56)
\ifill f:0
\move(406 55)
\lvec(407 55)
\lvec(407 56)
\lvec(406 56)
\ifill f:0
\move(408 55)
\lvec(410 55)
\lvec(410 56)
\lvec(408 56)
\ifill f:0
\move(411 55)
\lvec(413 55)
\lvec(413 56)
\lvec(411 56)
\ifill f:0
\move(414 55)
\lvec(415 55)
\lvec(415 56)
\lvec(414 56)
\ifill f:0
\move(417 55)
\lvec(418 55)
\lvec(418 56)
\lvec(417 56)
\ifill f:0
\move(419 55)
\lvec(421 55)
\lvec(421 56)
\lvec(419 56)
\ifill f:0
\move(422 55)
\lvec(424 55)
\lvec(424 56)
\lvec(422 56)
\ifill f:0
\move(425 55)
\lvec(427 55)
\lvec(427 56)
\lvec(425 56)
\ifill f:0
\move(428 55)
\lvec(430 55)
\lvec(430 56)
\lvec(428 56)
\ifill f:0
\move(431 55)
\lvec(433 55)
\lvec(433 56)
\lvec(431 56)
\ifill f:0
\move(434 55)
\lvec(436 55)
\lvec(436 56)
\lvec(434 56)
\ifill f:0
\move(437 55)
\lvec(439 55)
\lvec(439 56)
\lvec(437 56)
\ifill f:0
\move(440 55)
\lvec(442 55)
\lvec(442 56)
\lvec(440 56)
\ifill f:0
\move(443 55)
\lvec(445 55)
\lvec(445 56)
\lvec(443 56)
\ifill f:0
\move(446 55)
\lvec(448 55)
\lvec(448 56)
\lvec(446 56)
\ifill f:0
\move(449 55)
\lvec(451 55)
\lvec(451 56)
\lvec(449 56)
\ifill f:0
\move(14 56)
\lvec(15 56)
\lvec(15 57)
\lvec(14 57)
\ifill f:0
\move(19 56)
\lvec(20 56)
\lvec(20 57)
\lvec(19 57)
\ifill f:0
\move(121 56)
\lvec(122 56)
\lvec(122 57)
\lvec(121 57)
\ifill f:0
\move(135 56)
\lvec(136 56)
\lvec(136 57)
\lvec(135 57)
\ifill f:0
\move(140 56)
\lvec(141 56)
\lvec(141 57)
\lvec(140 57)
\ifill f:0
\move(145 56)
\lvec(146 56)
\lvec(146 57)
\lvec(145 57)
\ifill f:0
\move(148 56)
\lvec(149 56)
\lvec(149 57)
\lvec(148 57)
\ifill f:0
\move(155 56)
\lvec(156 56)
\lvec(156 57)
\lvec(155 57)
\ifill f:0
\move(159 56)
\lvec(160 56)
\lvec(160 57)
\lvec(159 57)
\ifill f:0
\move(163 56)
\lvec(164 56)
\lvec(164 57)
\lvec(163 57)
\ifill f:0
\move(168 56)
\lvec(170 56)
\lvec(170 57)
\lvec(168 57)
\ifill f:0
\move(174 56)
\lvec(176 56)
\lvec(176 57)
\lvec(174 57)
\ifill f:0
\move(183 56)
\lvec(186 56)
\lvec(186 57)
\lvec(183 57)
\ifill f:0
\move(214 56)
\lvec(219 56)
\lvec(219 57)
\lvec(214 57)
\ifill f:0
\move(225 56)
\lvec(228 56)
\lvec(228 57)
\lvec(225 57)
\ifill f:0
\move(233 56)
\lvec(235 56)
\lvec(235 57)
\lvec(233 57)
\ifill f:0
\move(239 56)
\lvec(241 56)
\lvec(241 57)
\lvec(239 57)
\ifill f:0
\move(245 56)
\lvec(247 56)
\lvec(247 57)
\lvec(245 57)
\ifill f:0
\move(250 56)
\lvec(251 56)
\lvec(251 57)
\lvec(250 57)
\ifill f:0
\move(254 56)
\lvec(256 56)
\lvec(256 57)
\lvec(254 57)
\ifill f:0
\move(258 56)
\lvec(260 56)
\lvec(260 57)
\lvec(258 57)
\ifill f:0
\move(262 56)
\lvec(264 56)
\lvec(264 57)
\lvec(262 57)
\ifill f:0
\move(266 56)
\lvec(268 56)
\lvec(268 57)
\lvec(266 57)
\ifill f:0
\move(270 56)
\lvec(271 56)
\lvec(271 57)
\lvec(270 57)
\ifill f:0
\move(273 56)
\lvec(275 56)
\lvec(275 57)
\lvec(273 57)
\ifill f:0
\move(277 56)
\lvec(278 56)
\lvec(278 57)
\lvec(277 57)
\ifill f:0
\move(280 56)
\lvec(282 56)
\lvec(282 57)
\lvec(280 57)
\ifill f:0
\move(283 56)
\lvec(285 56)
\lvec(285 57)
\lvec(283 57)
\ifill f:0
\move(286 56)
\lvec(288 56)
\lvec(288 57)
\lvec(286 57)
\ifill f:0
\move(289 56)
\lvec(291 56)
\lvec(291 57)
\lvec(289 57)
\ifill f:0
\move(292 56)
\lvec(294 56)
\lvec(294 57)
\lvec(292 57)
\ifill f:0
\move(295 56)
\lvec(296 56)
\lvec(296 57)
\lvec(295 57)
\ifill f:0
\move(298 56)
\lvec(299 56)
\lvec(299 57)
\lvec(298 57)
\ifill f:0
\move(301 56)
\lvec(302 56)
\lvec(302 57)
\lvec(301 57)
\ifill f:0
\move(303 56)
\lvec(305 56)
\lvec(305 57)
\lvec(303 57)
\ifill f:0
\move(306 56)
\lvec(307 56)
\lvec(307 57)
\lvec(306 57)
\ifill f:0
\move(309 56)
\lvec(310 56)
\lvec(310 57)
\lvec(309 57)
\ifill f:0
\move(311 56)
\lvec(312 56)
\lvec(312 57)
\lvec(311 57)
\ifill f:0
\move(314 56)
\lvec(315 56)
\lvec(315 57)
\lvec(314 57)
\ifill f:0
\move(316 56)
\lvec(317 56)
\lvec(317 57)
\lvec(316 57)
\ifill f:0
\move(319 56)
\lvec(320 56)
\lvec(320 57)
\lvec(319 57)
\ifill f:0
\move(321 56)
\lvec(322 56)
\lvec(322 57)
\lvec(321 57)
\ifill f:0
\move(323 56)
\lvec(324 56)
\lvec(324 57)
\lvec(323 57)
\ifill f:0
\move(326 56)
\lvec(327 56)
\lvec(327 57)
\lvec(326 57)
\ifill f:0
\move(328 56)
\lvec(329 56)
\lvec(329 57)
\lvec(328 57)
\ifill f:0
\move(330 56)
\lvec(331 56)
\lvec(331 57)
\lvec(330 57)
\ifill f:0
\move(332 56)
\lvec(333 56)
\lvec(333 57)
\lvec(332 57)
\ifill f:0
\move(335 56)
\lvec(336 56)
\lvec(336 57)
\lvec(335 57)
\ifill f:0
\move(337 56)
\lvec(338 56)
\lvec(338 57)
\lvec(337 57)
\ifill f:0
\move(339 56)
\lvec(340 56)
\lvec(340 57)
\lvec(339 57)
\ifill f:0
\move(341 56)
\lvec(342 56)
\lvec(342 57)
\lvec(341 57)
\ifill f:0
\move(343 56)
\lvec(344 56)
\lvec(344 57)
\lvec(343 57)
\ifill f:0
\move(345 56)
\lvec(346 56)
\lvec(346 57)
\lvec(345 57)
\ifill f:0
\move(347 56)
\lvec(348 56)
\lvec(348 57)
\lvec(347 57)
\ifill f:0
\move(349 56)
\lvec(350 56)
\lvec(350 57)
\lvec(349 57)
\ifill f:0
\move(351 56)
\lvec(352 56)
\lvec(352 57)
\lvec(351 57)
\ifill f:0
\move(353 56)
\lvec(354 56)
\lvec(354 57)
\lvec(353 57)
\ifill f:0
\move(357 56)
\lvec(358 56)
\lvec(358 57)
\lvec(357 57)
\ifill f:0
\move(359 56)
\lvec(360 56)
\lvec(360 57)
\lvec(359 57)
\ifill f:0
\move(361 56)
\lvec(362 56)
\lvec(362 57)
\lvec(361 57)
\ifill f:0
\move(363 56)
\lvec(364 56)
\lvec(364 57)
\lvec(363 57)
\ifill f:0
\move(365 56)
\lvec(366 56)
\lvec(366 57)
\lvec(365 57)
\ifill f:0
\move(367 56)
\lvec(368 56)
\lvec(368 57)
\lvec(367 57)
\ifill f:0
\move(369 56)
\lvec(370 56)
\lvec(370 57)
\lvec(369 57)
\ifill f:0
\move(371 56)
\lvec(372 56)
\lvec(372 57)
\lvec(371 57)
\ifill f:0
\move(373 56)
\lvec(374 56)
\lvec(374 57)
\lvec(373 57)
\ifill f:0
\move(375 56)
\lvec(376 56)
\lvec(376 57)
\lvec(375 57)
\ifill f:0
\move(378 56)
\lvec(379 56)
\lvec(379 57)
\lvec(378 57)
\ifill f:0
\move(380 56)
\lvec(381 56)
\lvec(381 57)
\lvec(380 57)
\ifill f:0
\move(382 56)
\lvec(383 56)
\lvec(383 57)
\lvec(382 57)
\ifill f:0
\move(384 56)
\lvec(385 56)
\lvec(385 57)
\lvec(384 57)
\ifill f:0
\move(386 56)
\lvec(388 56)
\lvec(388 57)
\lvec(386 57)
\ifill f:0
\move(389 56)
\lvec(390 56)
\lvec(390 57)
\lvec(389 57)
\ifill f:0
\move(391 56)
\lvec(392 56)
\lvec(392 57)
\lvec(391 57)
\ifill f:0
\move(393 56)
\lvec(395 56)
\lvec(395 57)
\lvec(393 57)
\ifill f:0
\move(396 56)
\lvec(397 56)
\lvec(397 57)
\lvec(396 57)
\ifill f:0
\move(398 56)
\lvec(399 56)
\lvec(399 57)
\lvec(398 57)
\ifill f:0
\move(401 56)
\lvec(402 56)
\lvec(402 57)
\lvec(401 57)
\ifill f:0
\move(403 56)
\lvec(404 56)
\lvec(404 57)
\lvec(403 57)
\ifill f:0
\move(405 56)
\lvec(407 56)
\lvec(407 57)
\lvec(405 57)
\ifill f:0
\move(408 56)
\lvec(409 56)
\lvec(409 57)
\lvec(408 57)
\ifill f:0
\move(410 56)
\lvec(412 56)
\lvec(412 57)
\lvec(410 57)
\ifill f:0
\move(413 56)
\lvec(414 56)
\lvec(414 57)
\lvec(413 57)
\ifill f:0
\move(416 56)
\lvec(417 56)
\lvec(417 57)
\lvec(416 57)
\ifill f:0
\move(418 56)
\lvec(420 56)
\lvec(420 57)
\lvec(418 57)
\ifill f:0
\move(421 56)
\lvec(422 56)
\lvec(422 57)
\lvec(421 57)
\ifill f:0
\move(423 56)
\lvec(425 56)
\lvec(425 57)
\lvec(423 57)
\ifill f:0
\move(426 56)
\lvec(428 56)
\lvec(428 57)
\lvec(426 57)
\ifill f:0
\move(429 56)
\lvec(430 56)
\lvec(430 57)
\lvec(429 57)
\ifill f:0
\move(432 56)
\lvec(433 56)
\lvec(433 57)
\lvec(432 57)
\ifill f:0
\move(434 56)
\lvec(436 56)
\lvec(436 57)
\lvec(434 57)
\ifill f:0
\move(437 56)
\lvec(439 56)
\lvec(439 57)
\lvec(437 57)
\ifill f:0
\move(440 56)
\lvec(442 56)
\lvec(442 57)
\lvec(440 57)
\ifill f:0
\move(443 56)
\lvec(445 56)
\lvec(445 57)
\lvec(443 57)
\ifill f:0
\move(446 56)
\lvec(448 56)
\lvec(448 57)
\lvec(446 57)
\ifill f:0
\move(449 56)
\lvec(451 56)
\lvec(451 57)
\lvec(449 57)
\ifill f:0
\move(12 57)
\lvec(14 57)
\lvec(14 58)
\lvec(12 58)
\ifill f:0
\move(91 57)
\lvec(92 57)
\lvec(92 58)
\lvec(91 58)
\ifill f:0
\move(93 57)
\lvec(94 57)
\lvec(94 58)
\lvec(93 58)
\ifill f:0
\move(121 57)
\lvec(122 57)
\lvec(122 58)
\lvec(121 58)
\ifill f:0
\move(126 57)
\lvec(127 57)
\lvec(127 58)
\lvec(126 58)
\ifill f:0
\move(138 57)
\lvec(139 57)
\lvec(139 58)
\lvec(138 58)
\ifill f:0
\move(140 57)
\lvec(141 57)
\lvec(141 58)
\lvec(140 58)
\ifill f:0
\move(145 57)
\lvec(146 57)
\lvec(146 58)
\lvec(145 58)
\ifill f:0
\move(148 57)
\lvec(149 57)
\lvec(149 58)
\lvec(148 58)
\ifill f:0
\move(157 57)
\lvec(158 57)
\lvec(158 58)
\lvec(157 58)
\ifill f:0
\move(160 57)
\lvec(161 57)
\lvec(161 58)
\lvec(160 58)
\ifill f:0
\move(164 57)
\lvec(165 57)
\lvec(165 58)
\lvec(164 58)
\ifill f:0
\move(168 57)
\lvec(170 57)
\lvec(170 58)
\lvec(168 58)
\ifill f:0
\move(173 57)
\lvec(175 57)
\lvec(175 58)
\lvec(173 58)
\ifill f:0
\move(179 57)
\lvec(181 57)
\lvec(181 58)
\lvec(179 58)
\ifill f:0
\move(187 57)
\lvec(190 57)
\lvec(190 58)
\lvec(187 58)
\ifill f:0
\move(201 57)
\lvec(214 57)
\lvec(214 58)
\lvec(201 58)
\ifill f:0
\move(225 57)
\lvec(229 57)
\lvec(229 58)
\lvec(225 58)
\ifill f:0
\move(235 57)
\lvec(238 57)
\lvec(238 58)
\lvec(235 58)
\ifill f:0
\move(242 57)
\lvec(245 57)
\lvec(245 58)
\lvec(242 58)
\ifill f:0
\move(248 57)
\lvec(251 57)
\lvec(251 58)
\lvec(248 58)
\ifill f:0
\move(254 57)
\lvec(256 57)
\lvec(256 58)
\lvec(254 58)
\ifill f:0
\move(259 57)
\lvec(261 57)
\lvec(261 58)
\lvec(259 58)
\ifill f:0
\move(263 57)
\lvec(265 57)
\lvec(265 58)
\lvec(263 58)
\ifill f:0
\move(268 57)
\lvec(269 57)
\lvec(269 58)
\lvec(268 58)
\ifill f:0
\move(272 57)
\lvec(273 57)
\lvec(273 58)
\lvec(272 58)
\ifill f:0
\move(275 57)
\lvec(277 57)
\lvec(277 58)
\lvec(275 58)
\ifill f:0
\move(279 57)
\lvec(281 57)
\lvec(281 58)
\lvec(279 58)
\ifill f:0
\move(283 57)
\lvec(284 57)
\lvec(284 58)
\lvec(283 58)
\ifill f:0
\move(286 57)
\lvec(288 57)
\lvec(288 58)
\lvec(286 58)
\ifill f:0
\move(289 57)
\lvec(291 57)
\lvec(291 58)
\lvec(289 58)
\ifill f:0
\move(293 57)
\lvec(294 57)
\lvec(294 58)
\lvec(293 58)
\ifill f:0
\move(296 57)
\lvec(297 57)
\lvec(297 58)
\lvec(296 58)
\ifill f:0
\move(299 57)
\lvec(300 57)
\lvec(300 58)
\lvec(299 58)
\ifill f:0
\move(302 57)
\lvec(303 57)
\lvec(303 58)
\lvec(302 58)
\ifill f:0
\move(305 57)
\lvec(306 57)
\lvec(306 58)
\lvec(305 58)
\ifill f:0
\move(307 57)
\lvec(309 57)
\lvec(309 58)
\lvec(307 58)
\ifill f:0
\move(310 57)
\lvec(311 57)
\lvec(311 58)
\lvec(310 58)
\ifill f:0
\move(313 57)
\lvec(314 57)
\lvec(314 58)
\lvec(313 58)
\ifill f:0
\move(315 57)
\lvec(317 57)
\lvec(317 58)
\lvec(315 58)
\ifill f:0
\move(318 57)
\lvec(319 57)
\lvec(319 58)
\lvec(318 58)
\ifill f:0
\move(321 57)
\lvec(322 57)
\lvec(322 58)
\lvec(321 58)
\ifill f:0
\move(323 57)
\lvec(324 57)
\lvec(324 58)
\lvec(323 58)
\ifill f:0
\move(326 57)
\lvec(327 57)
\lvec(327 58)
\lvec(326 58)
\ifill f:0
\move(328 57)
\lvec(329 57)
\lvec(329 58)
\lvec(328 58)
\ifill f:0
\move(331 57)
\lvec(332 57)
\lvec(332 58)
\lvec(331 58)
\ifill f:0
\move(333 57)
\lvec(334 57)
\lvec(334 58)
\lvec(333 58)
\ifill f:0
\move(335 57)
\lvec(336 57)
\lvec(336 58)
\lvec(335 58)
\ifill f:0
\move(338 57)
\lvec(339 57)
\lvec(339 58)
\lvec(338 58)
\ifill f:0
\move(340 57)
\lvec(341 57)
\lvec(341 58)
\lvec(340 58)
\ifill f:0
\move(342 57)
\lvec(343 57)
\lvec(343 58)
\lvec(342 58)
\ifill f:0
\move(344 57)
\lvec(345 57)
\lvec(345 58)
\lvec(344 58)
\ifill f:0
\move(347 57)
\lvec(348 57)
\lvec(348 58)
\lvec(347 58)
\ifill f:0
\move(349 57)
\lvec(350 57)
\lvec(350 58)
\lvec(349 58)
\ifill f:0
\move(351 57)
\lvec(352 57)
\lvec(352 58)
\lvec(351 58)
\ifill f:0
\move(353 57)
\lvec(354 57)
\lvec(354 58)
\lvec(353 58)
\ifill f:0
\move(355 57)
\lvec(356 57)
\lvec(356 58)
\lvec(355 58)
\ifill f:0
\move(357 57)
\lvec(358 57)
\lvec(358 58)
\lvec(357 58)
\ifill f:0
\move(359 57)
\lvec(360 57)
\lvec(360 58)
\lvec(359 58)
\ifill f:0
\move(361 57)
\lvec(362 57)
\lvec(362 58)
\lvec(361 58)
\ifill f:0
\move(363 57)
\lvec(364 57)
\lvec(364 58)
\lvec(363 58)
\ifill f:0
\move(365 57)
\lvec(367 57)
\lvec(367 58)
\lvec(365 58)
\ifill f:0
\move(368 57)
\lvec(372 57)
\lvec(372 58)
\lvec(368 58)
\ifill f:0
\move(373 57)
\lvec(374 57)
\lvec(374 58)
\lvec(373 58)
\ifill f:0
\move(375 57)
\lvec(376 57)
\lvec(376 58)
\lvec(375 58)
\ifill f:0
\move(377 57)
\lvec(378 57)
\lvec(378 58)
\lvec(377 58)
\ifill f:0
\move(379 57)
\lvec(380 57)
\lvec(380 58)
\lvec(379 58)
\ifill f:0
\move(381 57)
\lvec(382 57)
\lvec(382 58)
\lvec(381 58)
\ifill f:0
\move(383 57)
\lvec(384 57)
\lvec(384 58)
\lvec(383 58)
\ifill f:0
\move(385 57)
\lvec(386 57)
\lvec(386 58)
\lvec(385 58)
\ifill f:0
\move(387 57)
\lvec(388 57)
\lvec(388 58)
\lvec(387 58)
\ifill f:0
\move(389 57)
\lvec(391 57)
\lvec(391 58)
\lvec(389 58)
\ifill f:0
\move(392 57)
\lvec(393 57)
\lvec(393 58)
\lvec(392 58)
\ifill f:0
\move(394 57)
\lvec(395 57)
\lvec(395 58)
\lvec(394 58)
\ifill f:0
\move(396 57)
\lvec(397 57)
\lvec(397 58)
\lvec(396 58)
\ifill f:0
\move(398 57)
\lvec(399 57)
\lvec(399 58)
\lvec(398 58)
\ifill f:0
\move(401 57)
\lvec(402 57)
\lvec(402 58)
\lvec(401 58)
\ifill f:0
\move(403 57)
\lvec(404 57)
\lvec(404 58)
\lvec(403 58)
\ifill f:0
\move(405 57)
\lvec(406 57)
\lvec(406 58)
\lvec(405 58)
\ifill f:0
\move(407 57)
\lvec(409 57)
\lvec(409 58)
\lvec(407 58)
\ifill f:0
\move(410 57)
\lvec(411 57)
\lvec(411 58)
\lvec(410 58)
\ifill f:0
\move(412 57)
\lvec(414 57)
\lvec(414 58)
\lvec(412 58)
\ifill f:0
\move(415 57)
\lvec(416 57)
\lvec(416 58)
\lvec(415 58)
\ifill f:0
\move(417 57)
\lvec(418 57)
\lvec(418 58)
\lvec(417 58)
\ifill f:0
\move(420 57)
\lvec(421 57)
\lvec(421 58)
\lvec(420 58)
\ifill f:0
\move(422 57)
\lvec(423 57)
\lvec(423 58)
\lvec(422 58)
\ifill f:0
\move(425 57)
\lvec(426 57)
\lvec(426 58)
\lvec(425 58)
\ifill f:0
\move(427 57)
\lvec(428 57)
\lvec(428 58)
\lvec(427 58)
\ifill f:0
\move(430 57)
\lvec(431 57)
\lvec(431 58)
\lvec(430 58)
\ifill f:0
\move(432 57)
\lvec(434 57)
\lvec(434 58)
\lvec(432 58)
\ifill f:0
\move(435 57)
\lvec(436 57)
\lvec(436 58)
\lvec(435 58)
\ifill f:0
\move(437 57)
\lvec(439 57)
\lvec(439 58)
\lvec(437 58)
\ifill f:0
\move(440 57)
\lvec(442 57)
\lvec(442 58)
\lvec(440 58)
\ifill f:0
\move(443 57)
\lvec(444 57)
\lvec(444 58)
\lvec(443 58)
\ifill f:0
\move(446 57)
\lvec(447 57)
\lvec(447 58)
\lvec(446 58)
\ifill f:0
\move(448 57)
\lvec(450 57)
\lvec(450 58)
\lvec(448 58)
\ifill f:0
\move(22 58)
\lvec(23 58)
\lvec(23 59)
\lvec(22 59)
\ifill f:0
\move(118 58)
\lvec(119 58)
\lvec(119 59)
\lvec(118 59)
\ifill f:0
\move(121 58)
\lvec(122 58)
\lvec(122 59)
\lvec(121 59)
\ifill f:0
\move(124 58)
\lvec(125 58)
\lvec(125 59)
\lvec(124 59)
\ifill f:0
\move(129 58)
\lvec(130 58)
\lvec(130 59)
\lvec(129 59)
\ifill f:0
\move(145 58)
\lvec(146 58)
\lvec(146 59)
\lvec(145 59)
\ifill f:0
\move(150 58)
\lvec(151 58)
\lvec(151 59)
\lvec(150 59)
\ifill f:0
\move(155 58)
\lvec(156 58)
\lvec(156 59)
\lvec(155 59)
\ifill f:0
\move(158 58)
\lvec(159 58)
\lvec(159 59)
\lvec(158 59)
\ifill f:0
\move(165 58)
\lvec(166 58)
\lvec(166 59)
\lvec(165 59)
\ifill f:0
\move(169 58)
\lvec(170 58)
\lvec(170 59)
\lvec(169 59)
\ifill f:0
\move(173 58)
\lvec(174 58)
\lvec(174 59)
\lvec(173 59)
\ifill f:0
\move(177 58)
\lvec(179 58)
\lvec(179 59)
\lvec(177 59)
\ifill f:0
\move(183 58)
\lvec(184 58)
\lvec(184 59)
\lvec(183 59)
\ifill f:0
\move(189 58)
\lvec(192 58)
\lvec(192 59)
\lvec(189 59)
\ifill f:0
\move(199 58)
\lvec(204 58)
\lvec(204 59)
\lvec(199 59)
\ifill f:0
\move(225 58)
\lvec(231 58)
\lvec(231 59)
\lvec(225 59)
\ifill f:0
\move(238 58)
\lvec(241 58)
\lvec(241 59)
\lvec(238 59)
\ifill f:0
\move(247 58)
\lvec(249 58)
\lvec(249 59)
\lvec(247 59)
\ifill f:0
\move(253 58)
\lvec(256 58)
\lvec(256 59)
\lvec(253 59)
\ifill f:0
\move(259 58)
\lvec(261 58)
\lvec(261 59)
\lvec(259 59)
\ifill f:0
\move(264 58)
\lvec(266 58)
\lvec(266 59)
\lvec(264 59)
\ifill f:0
\move(269 58)
\lvec(271 58)
\lvec(271 59)
\lvec(269 59)
\ifill f:0
\move(274 58)
\lvec(275 58)
\lvec(275 59)
\lvec(274 59)
\ifill f:0
\move(278 58)
\lvec(280 58)
\lvec(280 59)
\lvec(278 59)
\ifill f:0
\move(282 58)
\lvec(283 58)
\lvec(283 59)
\lvec(282 59)
\ifill f:0
\move(286 58)
\lvec(287 58)
\lvec(287 59)
\lvec(286 59)
\ifill f:0
\move(289 58)
\lvec(291 58)
\lvec(291 59)
\lvec(289 59)
\ifill f:0
\move(293 58)
\lvec(294 58)
\lvec(294 59)
\lvec(293 59)
\ifill f:0
\move(296 58)
\lvec(298 58)
\lvec(298 59)
\lvec(296 59)
\ifill f:0
\move(300 58)
\lvec(301 58)
\lvec(301 59)
\lvec(300 59)
\ifill f:0
\move(303 58)
\lvec(304 58)
\lvec(304 59)
\lvec(303 59)
\ifill f:0
\move(306 58)
\lvec(307 58)
\lvec(307 59)
\lvec(306 59)
\ifill f:0
\move(309 58)
\lvec(310 58)
\lvec(310 59)
\lvec(309 59)
\ifill f:0
\move(312 58)
\lvec(313 58)
\lvec(313 59)
\lvec(312 59)
\ifill f:0
\move(315 58)
\lvec(316 58)
\lvec(316 59)
\lvec(315 59)
\ifill f:0
\move(318 58)
\lvec(319 58)
\lvec(319 59)
\lvec(318 59)
\ifill f:0
\move(320 58)
\lvec(322 58)
\lvec(322 59)
\lvec(320 59)
\ifill f:0
\move(323 58)
\lvec(324 58)
\lvec(324 59)
\lvec(323 59)
\ifill f:0
\move(326 58)
\lvec(327 58)
\lvec(327 59)
\lvec(326 59)
\ifill f:0
\move(328 58)
\lvec(330 58)
\lvec(330 59)
\lvec(328 59)
\ifill f:0
\move(331 58)
\lvec(332 58)
\lvec(332 59)
\lvec(331 59)
\ifill f:0
\move(333 58)
\lvec(335 58)
\lvec(335 59)
\lvec(333 59)
\ifill f:0
\move(336 58)
\lvec(337 58)
\lvec(337 59)
\lvec(336 59)
\ifill f:0
\move(338 58)
\lvec(340 58)
\lvec(340 59)
\lvec(338 59)
\ifill f:0
\move(341 58)
\lvec(342 58)
\lvec(342 59)
\lvec(341 59)
\ifill f:0
\move(343 58)
\lvec(344 58)
\lvec(344 59)
\lvec(343 59)
\ifill f:0
\move(346 58)
\lvec(347 58)
\lvec(347 59)
\lvec(346 59)
\ifill f:0
\move(348 58)
\lvec(349 58)
\lvec(349 59)
\lvec(348 59)
\ifill f:0
\move(350 58)
\lvec(351 58)
\lvec(351 59)
\lvec(350 59)
\ifill f:0
\move(353 58)
\lvec(354 58)
\lvec(354 59)
\lvec(353 59)
\ifill f:0
\move(355 58)
\lvec(356 58)
\lvec(356 59)
\lvec(355 59)
\ifill f:0
\move(357 58)
\lvec(358 58)
\lvec(358 59)
\lvec(357 59)
\ifill f:0
\move(359 58)
\lvec(360 58)
\lvec(360 59)
\lvec(359 59)
\ifill f:0
\move(361 58)
\lvec(362 58)
\lvec(362 59)
\lvec(361 59)
\ifill f:0
\move(364 58)
\lvec(365 58)
\lvec(365 59)
\lvec(364 59)
\ifill f:0
\move(366 58)
\lvec(367 58)
\lvec(367 59)
\lvec(366 59)
\ifill f:0
\move(368 58)
\lvec(369 58)
\lvec(369 59)
\lvec(368 59)
\ifill f:0
\move(370 58)
\lvec(371 58)
\lvec(371 59)
\lvec(370 59)
\ifill f:0
\move(372 58)
\lvec(373 58)
\lvec(373 59)
\lvec(372 59)
\ifill f:0
\move(374 58)
\lvec(375 58)
\lvec(375 59)
\lvec(374 59)
\ifill f:0
\move(376 58)
\lvec(377 58)
\lvec(377 59)
\lvec(376 59)
\ifill f:0
\move(378 58)
\lvec(379 58)
\lvec(379 59)
\lvec(378 59)
\ifill f:0
\move(380 58)
\lvec(381 58)
\lvec(381 59)
\lvec(380 59)
\ifill f:0
\move(382 58)
\lvec(383 58)
\lvec(383 59)
\lvec(382 59)
\ifill f:0
\move(384 58)
\lvec(385 58)
\lvec(385 59)
\lvec(384 59)
\ifill f:0
\move(386 58)
\lvec(387 58)
\lvec(387 59)
\lvec(386 59)
\ifill f:0
\move(388 58)
\lvec(389 58)
\lvec(389 59)
\lvec(388 59)
\ifill f:0
\move(390 58)
\lvec(391 58)
\lvec(391 59)
\lvec(390 59)
\ifill f:0
\move(392 58)
\lvec(393 58)
\lvec(393 59)
\lvec(392 59)
\ifill f:0
\move(394 58)
\lvec(395 58)
\lvec(395 59)
\lvec(394 59)
\ifill f:0
\move(396 58)
\lvec(397 58)
\lvec(397 59)
\lvec(396 59)
\ifill f:0
\move(398 58)
\lvec(400 58)
\lvec(400 59)
\lvec(398 59)
\ifill f:0
\move(401 58)
\lvec(402 58)
\lvec(402 59)
\lvec(401 59)
\ifill f:0
\move(403 58)
\lvec(404 58)
\lvec(404 59)
\lvec(403 59)
\ifill f:0
\move(405 58)
\lvec(406 58)
\lvec(406 59)
\lvec(405 59)
\ifill f:0
\move(407 58)
\lvec(408 58)
\lvec(408 59)
\lvec(407 59)
\ifill f:0
\move(409 58)
\lvec(411 58)
\lvec(411 59)
\lvec(409 59)
\ifill f:0
\move(412 58)
\lvec(413 58)
\lvec(413 59)
\lvec(412 59)
\ifill f:0
\move(414 58)
\lvec(415 58)
\lvec(415 59)
\lvec(414 59)
\ifill f:0
\move(416 58)
\lvec(417 58)
\lvec(417 59)
\lvec(416 59)
\ifill f:0
\move(418 58)
\lvec(420 58)
\lvec(420 59)
\lvec(418 59)
\ifill f:0
\move(421 58)
\lvec(422 58)
\lvec(422 59)
\lvec(421 59)
\ifill f:0
\move(423 58)
\lvec(424 58)
\lvec(424 59)
\lvec(423 59)
\ifill f:0
\move(426 58)
\lvec(427 58)
\lvec(427 59)
\lvec(426 59)
\ifill f:0
\move(428 58)
\lvec(429 58)
\lvec(429 59)
\lvec(428 59)
\ifill f:0
\move(430 58)
\lvec(432 58)
\lvec(432 59)
\lvec(430 59)
\ifill f:0
\move(433 58)
\lvec(434 58)
\lvec(434 59)
\lvec(433 59)
\ifill f:0
\move(435 58)
\lvec(437 58)
\lvec(437 59)
\lvec(435 59)
\ifill f:0
\move(438 58)
\lvec(439 58)
\lvec(439 59)
\lvec(438 59)
\ifill f:0
\move(440 58)
\lvec(442 58)
\lvec(442 59)
\lvec(440 59)
\ifill f:0
\move(443 58)
\lvec(444 58)
\lvec(444 59)
\lvec(443 59)
\ifill f:0
\move(445 58)
\lvec(447 58)
\lvec(447 59)
\lvec(445 59)
\ifill f:0
\move(448 58)
\lvec(449 58)
\lvec(449 59)
\lvec(448 59)
\ifill f:0
\move(12 59)
\lvec(13 59)
\lvec(13 60)
\lvec(12 60)
\ifill f:0
\move(15 59)
\lvec(16 59)
\lvec(16 60)
\lvec(15 60)
\ifill f:0
\move(19 59)
\lvec(20 59)
\lvec(20 60)
\lvec(19 60)
\ifill f:0
\move(113 59)
\lvec(114 59)
\lvec(114 60)
\lvec(113 60)
\ifill f:0
\move(121 59)
\lvec(122 59)
\lvec(122 60)
\lvec(121 60)
\ifill f:0
\move(130 59)
\lvec(131 59)
\lvec(131 60)
\lvec(130 60)
\ifill f:0
\move(137 59)
\lvec(138 59)
\lvec(138 60)
\lvec(137 60)
\ifill f:0
\move(139 59)
\lvec(140 59)
\lvec(140 60)
\lvec(139 60)
\ifill f:0
\move(141 59)
\lvec(142 59)
\lvec(142 60)
\lvec(141 60)
\ifill f:0
\move(143 59)
\lvec(144 59)
\lvec(144 60)
\lvec(143 60)
\ifill f:0
\move(145 59)
\lvec(146 59)
\lvec(146 60)
\lvec(145 60)
\ifill f:0
\move(147 59)
\lvec(148 59)
\lvec(148 60)
\lvec(147 60)
\ifill f:0
\move(154 59)
\lvec(155 59)
\lvec(155 60)
\lvec(154 60)
\ifill f:0
\move(162 59)
\lvec(163 59)
\lvec(163 60)
\lvec(162 60)
\ifill f:0
\move(169 59)
\lvec(170 59)
\lvec(170 60)
\lvec(169 60)
\ifill f:0
\move(172 59)
\lvec(173 59)
\lvec(173 60)
\lvec(172 60)
\ifill f:0
\move(176 59)
\lvec(177 59)
\lvec(177 60)
\lvec(176 60)
\ifill f:0
\move(180 59)
\lvec(182 59)
\lvec(182 60)
\lvec(180 60)
\ifill f:0
\move(185 59)
\lvec(187 59)
\lvec(187 60)
\lvec(185 60)
\ifill f:0
\move(191 59)
\lvec(193 59)
\lvec(193 60)
\lvec(191 60)
\ifill f:0
\move(198 59)
\lvec(201 59)
\lvec(201 60)
\lvec(198 60)
\ifill f:0
\move(210 59)
\lvec(221 59)
\lvec(221 60)
\lvec(210 60)
\ifill f:0
\move(222 59)
\lvec(235 59)
\lvec(235 60)
\lvec(222 60)
\ifill f:0
\move(243 59)
\lvec(247 59)
\lvec(247 60)
\lvec(243 60)
\ifill f:0
\move(253 59)
\lvec(255 59)
\lvec(255 60)
\lvec(253 60)
\ifill f:0
\move(260 59)
\lvec(262 59)
\lvec(262 60)
\lvec(260 60)
\ifill f:0
\move(266 59)
\lvec(268 59)
\lvec(268 60)
\lvec(266 60)
\ifill f:0
\move(271 59)
\lvec(273 59)
\lvec(273 60)
\lvec(271 60)
\ifill f:0
\move(276 59)
\lvec(278 59)
\lvec(278 60)
\lvec(276 60)
\ifill f:0
\move(281 59)
\lvec(283 59)
\lvec(283 60)
\lvec(281 60)
\ifill f:0
\move(285 59)
\lvec(287 59)
\lvec(287 60)
\lvec(285 60)
\ifill f:0
\move(289 59)
\lvec(291 59)
\lvec(291 60)
\lvec(289 60)
\ifill f:0
\move(293 59)
\lvec(295 59)
\lvec(295 60)
\lvec(293 60)
\ifill f:0
\move(297 59)
\lvec(299 59)
\lvec(299 60)
\lvec(297 60)
\ifill f:0
\move(301 59)
\lvec(302 59)
\lvec(302 60)
\lvec(301 60)
\ifill f:0
\move(304 59)
\lvec(306 59)
\lvec(306 60)
\lvec(304 60)
\ifill f:0
\move(307 59)
\lvec(309 59)
\lvec(309 60)
\lvec(307 60)
\ifill f:0
\move(311 59)
\lvec(312 59)
\lvec(312 60)
\lvec(311 60)
\ifill f:0
\move(314 59)
\lvec(315 59)
\lvec(315 60)
\lvec(314 60)
\ifill f:0
\move(317 59)
\lvec(318 59)
\lvec(318 60)
\lvec(317 60)
\ifill f:0
\move(320 59)
\lvec(321 59)
\lvec(321 60)
\lvec(320 60)
\ifill f:0
\move(323 59)
\lvec(324 59)
\lvec(324 60)
\lvec(323 60)
\ifill f:0
\move(326 59)
\lvec(327 59)
\lvec(327 60)
\lvec(326 60)
\ifill f:0
\move(329 59)
\lvec(330 59)
\lvec(330 60)
\lvec(329 60)
\ifill f:0
\move(331 59)
\lvec(333 59)
\lvec(333 60)
\lvec(331 60)
\ifill f:0
\move(334 59)
\lvec(335 59)
\lvec(335 60)
\lvec(334 60)
\ifill f:0
\move(337 59)
\lvec(338 59)
\lvec(338 60)
\lvec(337 60)
\ifill f:0
\move(339 59)
\lvec(341 59)
\lvec(341 60)
\lvec(339 60)
\ifill f:0
\move(342 59)
\lvec(343 59)
\lvec(343 60)
\lvec(342 60)
\ifill f:0
\move(345 59)
\lvec(346 59)
\lvec(346 60)
\lvec(345 60)
\ifill f:0
\move(347 59)
\lvec(348 59)
\lvec(348 60)
\lvec(347 60)
\ifill f:0
\move(350 59)
\lvec(351 59)
\lvec(351 60)
\lvec(350 60)
\ifill f:0
\move(352 59)
\lvec(353 59)
\lvec(353 60)
\lvec(352 60)
\ifill f:0
\move(354 59)
\lvec(356 59)
\lvec(356 60)
\lvec(354 60)
\ifill f:0
\move(357 59)
\lvec(358 59)
\lvec(358 60)
\lvec(357 60)
\ifill f:0
\move(359 59)
\lvec(360 59)
\lvec(360 60)
\lvec(359 60)
\ifill f:0
\move(361 59)
\lvec(363 59)
\lvec(363 60)
\lvec(361 60)
\ifill f:0
\move(364 59)
\lvec(365 59)
\lvec(365 60)
\lvec(364 60)
\ifill f:0
\move(366 59)
\lvec(367 59)
\lvec(367 60)
\lvec(366 60)
\ifill f:0
\move(368 59)
\lvec(369 59)
\lvec(369 60)
\lvec(368 60)
\ifill f:0
\move(370 59)
\lvec(371 59)
\lvec(371 60)
\lvec(370 60)
\ifill f:0
\move(373 59)
\lvec(374 59)
\lvec(374 60)
\lvec(373 60)
\ifill f:0
\move(375 59)
\lvec(376 59)
\lvec(376 60)
\lvec(375 60)
\ifill f:0
\move(377 59)
\lvec(378 59)
\lvec(378 60)
\lvec(377 60)
\ifill f:0
\move(379 59)
\lvec(380 59)
\lvec(380 60)
\lvec(379 60)
\ifill f:0
\move(381 59)
\lvec(382 59)
\lvec(382 60)
\lvec(381 60)
\ifill f:0
\move(383 59)
\lvec(384 59)
\lvec(384 60)
\lvec(383 60)
\ifill f:0
\move(385 59)
\lvec(386 59)
\lvec(386 60)
\lvec(385 60)
\ifill f:0
\move(387 59)
\lvec(388 59)
\lvec(388 60)
\lvec(387 60)
\ifill f:0
\move(389 59)
\lvec(390 59)
\lvec(390 60)
\lvec(389 60)
\ifill f:0
\move(391 59)
\lvec(392 59)
\lvec(392 60)
\lvec(391 60)
\ifill f:0
\move(394 59)
\lvec(396 59)
\lvec(396 60)
\lvec(394 60)
\ifill f:0
\move(397 59)
\lvec(398 59)
\lvec(398 60)
\lvec(397 60)
\ifill f:0
\move(399 59)
\lvec(400 59)
\lvec(400 60)
\lvec(399 60)
\ifill f:0
\move(401 59)
\lvec(402 59)
\lvec(402 60)
\lvec(401 60)
\ifill f:0
\move(403 59)
\lvec(404 59)
\lvec(404 60)
\lvec(403 60)
\ifill f:0
\move(405 59)
\lvec(406 59)
\lvec(406 60)
\lvec(405 60)
\ifill f:0
\move(407 59)
\lvec(408 59)
\lvec(408 60)
\lvec(407 60)
\ifill f:0
\move(409 59)
\lvec(410 59)
\lvec(410 60)
\lvec(409 60)
\ifill f:0
\move(411 59)
\lvec(412 59)
\lvec(412 60)
\lvec(411 60)
\ifill f:0
\move(413 59)
\lvec(414 59)
\lvec(414 60)
\lvec(413 60)
\ifill f:0
\move(415 59)
\lvec(416 59)
\lvec(416 60)
\lvec(415 60)
\ifill f:0
\move(418 59)
\lvec(419 59)
\lvec(419 60)
\lvec(418 60)
\ifill f:0
\move(420 59)
\lvec(421 59)
\lvec(421 60)
\lvec(420 60)
\ifill f:0
\move(422 59)
\lvec(423 59)
\lvec(423 60)
\lvec(422 60)
\ifill f:0
\move(424 59)
\lvec(425 59)
\lvec(425 60)
\lvec(424 60)
\ifill f:0
\move(426 59)
\lvec(428 59)
\lvec(428 60)
\lvec(426 60)
\ifill f:0
\move(429 59)
\lvec(430 59)
\lvec(430 60)
\lvec(429 60)
\ifill f:0
\move(431 59)
\lvec(432 59)
\lvec(432 60)
\lvec(431 60)
\ifill f:0
\move(433 59)
\lvec(434 59)
\lvec(434 60)
\lvec(433 60)
\ifill f:0
\move(436 59)
\lvec(437 59)
\lvec(437 60)
\lvec(436 60)
\ifill f:0
\move(438 59)
\lvec(439 59)
\lvec(439 60)
\lvec(438 60)
\ifill f:0
\move(440 59)
\lvec(442 59)
\lvec(442 60)
\lvec(440 60)
\ifill f:0
\move(443 59)
\lvec(444 59)
\lvec(444 60)
\lvec(443 60)
\ifill f:0
\move(445 59)
\lvec(446 59)
\lvec(446 60)
\lvec(445 60)
\ifill f:0
\move(448 59)
\lvec(449 59)
\lvec(449 60)
\lvec(448 60)
\ifill f:0
\move(450 59)
\lvec(451 59)
\lvec(451 60)
\lvec(450 60)
\ifill f:0
\move(15 60)
\lvec(16 60)
\lvec(16 61)
\lvec(15 61)
\ifill f:0
\move(121 60)
\lvec(122 60)
\lvec(122 61)
\lvec(121 61)
\ifill f:0
\move(125 60)
\lvec(126 60)
\lvec(126 61)
\lvec(125 61)
\ifill f:0
\move(128 60)
\lvec(129 60)
\lvec(129 61)
\lvec(128 61)
\ifill f:0
\move(131 60)
\lvec(132 60)
\lvec(132 61)
\lvec(131 61)
\ifill f:0
\move(141 60)
\lvec(142 60)
\lvec(142 61)
\lvec(141 61)
\ifill f:0
\move(143 60)
\lvec(144 60)
\lvec(144 61)
\lvec(143 61)
\ifill f:0
\move(145 60)
\lvec(146 60)
\lvec(146 61)
\lvec(145 61)
\ifill f:0
\move(147 60)
\lvec(148 60)
\lvec(148 61)
\lvec(147 61)
\ifill f:0
\move(149 60)
\lvec(150 60)
\lvec(150 61)
\lvec(149 61)
\ifill f:0
\move(151 60)
\lvec(152 60)
\lvec(152 61)
\lvec(151 61)
\ifill f:0
\move(153 60)
\lvec(154 60)
\lvec(154 61)
\lvec(153 61)
\ifill f:0
\move(158 60)
\lvec(159 60)
\lvec(159 61)
\lvec(158 61)
\ifill f:0
\move(163 60)
\lvec(164 60)
\lvec(164 61)
\lvec(163 61)
\ifill f:0
\move(166 60)
\lvec(167 60)
\lvec(167 61)
\lvec(166 61)
\ifill f:0
\move(169 60)
\lvec(170 60)
\lvec(170 61)
\lvec(169 61)
\ifill f:0
\move(172 60)
\lvec(173 60)
\lvec(173 61)
\lvec(172 61)
\ifill f:0
\move(175 60)
\lvec(176 60)
\lvec(176 61)
\lvec(175 61)
\ifill f:0
\move(179 60)
\lvec(180 60)
\lvec(180 61)
\lvec(179 61)
\ifill f:0
\move(183 60)
\lvec(184 60)
\lvec(184 61)
\lvec(183 61)
\ifill f:0
\move(187 60)
\lvec(188 60)
\lvec(188 61)
\lvec(187 61)
\ifill f:0
\move(192 60)
\lvec(194 60)
\lvec(194 61)
\lvec(192 61)
\ifill f:0
\move(198 60)
\lvec(200 60)
\lvec(200 61)
\lvec(198 61)
\ifill f:0
\move(205 60)
\lvec(208 60)
\lvec(208 61)
\lvec(205 61)
\ifill f:0
\move(217 60)
\lvec(229 60)
\lvec(229 61)
\lvec(217 61)
\ifill f:0
\move(230 60)
\lvec(242 60)
\lvec(242 61)
\lvec(230 61)
\ifill f:0
\move(251 60)
\lvec(255 60)
\lvec(255 61)
\lvec(251 61)
\ifill f:0
\move(261 60)
\lvec(264 60)
\lvec(264 61)
\lvec(261 61)
\ifill f:0
\move(268 60)
\lvec(270 60)
\lvec(270 61)
\lvec(268 61)
\ifill f:0
\move(274 60)
\lvec(276 60)
\lvec(276 61)
\lvec(274 61)
\ifill f:0
\move(280 60)
\lvec(282 60)
\lvec(282 61)
\lvec(280 61)
\ifill f:0
\move(285 60)
\lvec(287 60)
\lvec(287 61)
\lvec(285 61)
\ifill f:0
\move(289 60)
\lvec(291 60)
\lvec(291 61)
\lvec(289 61)
\ifill f:0
\move(294 60)
\lvec(296 60)
\lvec(296 61)
\lvec(294 61)
\ifill f:0
\move(298 60)
\lvec(300 60)
\lvec(300 61)
\lvec(298 61)
\ifill f:0
\move(302 60)
\lvec(304 60)
\lvec(304 61)
\lvec(302 61)
\ifill f:0
\move(306 60)
\lvec(307 60)
\lvec(307 61)
\lvec(306 61)
\ifill f:0
\move(309 60)
\lvec(311 60)
\lvec(311 61)
\lvec(309 61)
\ifill f:0
\move(313 60)
\lvec(314 60)
\lvec(314 61)
\lvec(313 61)
\ifill f:0
\move(316 60)
\lvec(318 60)
\lvec(318 61)
\lvec(316 61)
\ifill f:0
\move(320 60)
\lvec(321 60)
\lvec(321 61)
\lvec(320 61)
\ifill f:0
\move(323 60)
\lvec(324 60)
\lvec(324 61)
\lvec(323 61)
\ifill f:0
\move(326 60)
\lvec(327 60)
\lvec(327 61)
\lvec(326 61)
\ifill f:0
\move(329 60)
\lvec(330 60)
\lvec(330 61)
\lvec(329 61)
\ifill f:0
\move(332 60)
\lvec(333 60)
\lvec(333 61)
\lvec(332 61)
\ifill f:0
\move(335 60)
\lvec(336 60)
\lvec(336 61)
\lvec(335 61)
\ifill f:0
\move(338 60)
\lvec(339 60)
\lvec(339 61)
\lvec(338 61)
\ifill f:0
\move(341 60)
\lvec(342 60)
\lvec(342 61)
\lvec(341 61)
\ifill f:0
\move(343 60)
\lvec(345 60)
\lvec(345 61)
\lvec(343 61)
\ifill f:0
\move(346 60)
\lvec(347 60)
\lvec(347 61)
\lvec(346 61)
\ifill f:0
\move(349 60)
\lvec(350 60)
\lvec(350 61)
\lvec(349 61)
\ifill f:0
\move(351 60)
\lvec(353 60)
\lvec(353 61)
\lvec(351 61)
\ifill f:0
\move(354 60)
\lvec(355 60)
\lvec(355 61)
\lvec(354 61)
\ifill f:0
\move(356 60)
\lvec(358 60)
\lvec(358 61)
\lvec(356 61)
\ifill f:0
\move(359 60)
\lvec(360 60)
\lvec(360 61)
\lvec(359 61)
\ifill f:0
\move(361 60)
\lvec(363 60)
\lvec(363 61)
\lvec(361 61)
\ifill f:0
\move(364 60)
\lvec(365 60)
\lvec(365 61)
\lvec(364 61)
\ifill f:0
\move(366 60)
\lvec(367 60)
\lvec(367 61)
\lvec(366 61)
\ifill f:0
\move(369 60)
\lvec(370 60)
\lvec(370 61)
\lvec(369 61)
\ifill f:0
\move(371 60)
\lvec(372 60)
\lvec(372 61)
\lvec(371 61)
\ifill f:0
\move(373 60)
\lvec(374 60)
\lvec(374 61)
\lvec(373 61)
\ifill f:0
\move(376 60)
\lvec(377 60)
\lvec(377 61)
\lvec(376 61)
\ifill f:0
\move(378 60)
\lvec(379 60)
\lvec(379 61)
\lvec(378 61)
\ifill f:0
\move(380 60)
\lvec(381 60)
\lvec(381 61)
\lvec(380 61)
\ifill f:0
\move(382 60)
\lvec(383 60)
\lvec(383 61)
\lvec(382 61)
\ifill f:0
\move(385 60)
\lvec(386 60)
\lvec(386 61)
\lvec(385 61)
\ifill f:0
\move(387 60)
\lvec(388 60)
\lvec(388 61)
\lvec(387 61)
\ifill f:0
\move(389 60)
\lvec(390 60)
\lvec(390 61)
\lvec(389 61)
\ifill f:0
\move(391 60)
\lvec(392 60)
\lvec(392 61)
\lvec(391 61)
\ifill f:0
\move(393 60)
\lvec(394 60)
\lvec(394 61)
\lvec(393 61)
\ifill f:0
\move(395 60)
\lvec(396 60)
\lvec(396 61)
\lvec(395 61)
\ifill f:0
\move(397 60)
\lvec(398 60)
\lvec(398 61)
\lvec(397 61)
\ifill f:0
\move(399 60)
\lvec(400 60)
\lvec(400 61)
\lvec(399 61)
\ifill f:0
\move(401 60)
\lvec(402 60)
\lvec(402 61)
\lvec(401 61)
\ifill f:0
\move(403 60)
\lvec(407 60)
\lvec(407 61)
\lvec(403 61)
\ifill f:0
\move(408 60)
\lvec(412 60)
\lvec(412 61)
\lvec(408 61)
\ifill f:0
\move(413 60)
\lvec(414 60)
\lvec(414 61)
\lvec(413 61)
\ifill f:0
\move(415 60)
\lvec(416 60)
\lvec(416 61)
\lvec(415 61)
\ifill f:0
\move(417 60)
\lvec(418 60)
\lvec(418 61)
\lvec(417 61)
\ifill f:0
\move(419 60)
\lvec(420 60)
\lvec(420 61)
\lvec(419 61)
\ifill f:0
\move(421 60)
\lvec(422 60)
\lvec(422 61)
\lvec(421 61)
\ifill f:0
\move(423 60)
\lvec(424 60)
\lvec(424 61)
\lvec(423 61)
\ifill f:0
\move(425 60)
\lvec(426 60)
\lvec(426 61)
\lvec(425 61)
\ifill f:0
\move(427 60)
\lvec(428 60)
\lvec(428 61)
\lvec(427 61)
\ifill f:0
\move(429 60)
\lvec(430 60)
\lvec(430 61)
\lvec(429 61)
\ifill f:0
\move(431 60)
\lvec(433 60)
\lvec(433 61)
\lvec(431 61)
\ifill f:0
\move(434 60)
\lvec(435 60)
\lvec(435 61)
\lvec(434 61)
\ifill f:0
\move(436 60)
\lvec(437 60)
\lvec(437 61)
\lvec(436 61)
\ifill f:0
\move(438 60)
\lvec(439 60)
\lvec(439 61)
\lvec(438 61)
\ifill f:0
\move(440 60)
\lvec(442 60)
\lvec(442 61)
\lvec(440 61)
\ifill f:0
\move(443 60)
\lvec(444 60)
\lvec(444 61)
\lvec(443 61)
\ifill f:0
\move(445 60)
\lvec(446 60)
\lvec(446 61)
\lvec(445 61)
\ifill f:0
\move(447 60)
\lvec(448 60)
\lvec(448 61)
\lvec(447 61)
\ifill f:0
\move(450 60)
\lvec(451 60)
\lvec(451 61)
\lvec(450 61)
\ifill f:0
\move(12 61)
\lvec(13 61)
\lvec(13 62)
\lvec(12 62)
\ifill f:0
\move(15 61)
\lvec(16 61)
\lvec(16 62)
\lvec(15 62)
\ifill f:0
\move(22 61)
\lvec(23 61)
\lvec(23 62)
\lvec(22 62)
\ifill f:0
\move(115 61)
\lvec(116 61)
\lvec(116 62)
\lvec(115 62)
\ifill f:0
\move(121 61)
\lvec(122 61)
\lvec(122 62)
\lvec(121 62)
\ifill f:0
\move(129 61)
\lvec(130 61)
\lvec(130 62)
\lvec(129 62)
\ifill f:0
\move(138 61)
\lvec(139 61)
\lvec(139 62)
\lvec(138 62)
\ifill f:0
\move(143 61)
\lvec(144 61)
\lvec(144 62)
\lvec(143 62)
\ifill f:0
\move(159 61)
\lvec(160 61)
\lvec(160 62)
\lvec(159 62)
\ifill f:0
\move(161 61)
\lvec(162 61)
\lvec(162 62)
\lvec(161 62)
\ifill f:0
\move(166 61)
\lvec(167 61)
\lvec(167 62)
\lvec(166 62)
\ifill f:0
\move(169 61)
\lvec(170 61)
\lvec(170 62)
\lvec(169 62)
\ifill f:0
\move(178 61)
\lvec(179 61)
\lvec(179 62)
\lvec(178 62)
\ifill f:0
\move(181 61)
\lvec(182 61)
\lvec(182 62)
\lvec(181 62)
\ifill f:0
\move(185 61)
\lvec(186 61)
\lvec(186 62)
\lvec(185 62)
\ifill f:0
\move(189 61)
\lvec(190 61)
\lvec(190 62)
\lvec(189 62)
\ifill f:0
\move(193 61)
\lvec(194 61)
\lvec(194 62)
\lvec(193 62)
\ifill f:0
\move(198 61)
\lvec(199 61)
\lvec(199 62)
\lvec(198 62)
\ifill f:0
\move(203 61)
\lvec(205 61)
\lvec(205 62)
\lvec(203 62)
\ifill f:0
\move(211 61)
\lvec(213 61)
\lvec(213 62)
\lvec(211 62)
\ifill f:0
\move(221 61)
\lvec(226 61)
\lvec(226 62)
\lvec(221 62)
\ifill f:0
\move(248 61)
\lvec(254 61)
\lvec(254 62)
\lvec(248 62)
\ifill f:0
\move(262 61)
\lvec(265 61)
\lvec(265 62)
\lvec(262 62)
\ifill f:0
\move(271 61)
\lvec(274 61)
\lvec(274 62)
\lvec(271 62)
\ifill f:0
\move(278 61)
\lvec(280 61)
\lvec(280 62)
\lvec(278 62)
\ifill f:0
\move(284 61)
\lvec(286 61)
\lvec(286 62)
\lvec(284 62)
\ifill f:0
\move(289 61)
\lvec(291 61)
\lvec(291 62)
\lvec(289 62)
\ifill f:0
\move(294 61)
\lvec(296 61)
\lvec(296 62)
\lvec(294 62)
\ifill f:0
\move(299 61)
\lvec(301 61)
\lvec(301 62)
\lvec(299 62)
\ifill f:0
\move(303 61)
\lvec(305 61)
\lvec(305 62)
\lvec(303 62)
\ifill f:0
\move(308 61)
\lvec(309 61)
\lvec(309 62)
\lvec(308 62)
\ifill f:0
\move(312 61)
\lvec(313 61)
\lvec(313 62)
\lvec(312 62)
\ifill f:0
\move(315 61)
\lvec(317 61)
\lvec(317 62)
\lvec(315 62)
\ifill f:0
\move(319 61)
\lvec(321 61)
\lvec(321 62)
\lvec(319 62)
\ifill f:0
\move(323 61)
\lvec(324 61)
\lvec(324 62)
\lvec(323 62)
\ifill f:0
\move(326 61)
\lvec(328 61)
\lvec(328 62)
\lvec(326 62)
\ifill f:0
\move(329 61)
\lvec(331 61)
\lvec(331 62)
\lvec(329 62)
\ifill f:0
\move(333 61)
\lvec(334 61)
\lvec(334 62)
\lvec(333 62)
\ifill f:0
\move(336 61)
\lvec(337 61)
\lvec(337 62)
\lvec(336 62)
\ifill f:0
\move(339 61)
\lvec(340 61)
\lvec(340 62)
\lvec(339 62)
\ifill f:0
\move(342 61)
\lvec(343 61)
\lvec(343 62)
\lvec(342 62)
\ifill f:0
\move(345 61)
\lvec(346 61)
\lvec(346 62)
\lvec(345 62)
\ifill f:0
\move(348 61)
\lvec(349 61)
\lvec(349 62)
\lvec(348 62)
\ifill f:0
\move(351 61)
\lvec(352 61)
\lvec(352 62)
\lvec(351 62)
\ifill f:0
\move(353 61)
\lvec(355 61)
\lvec(355 62)
\lvec(353 62)
\ifill f:0
\move(356 61)
\lvec(357 61)
\lvec(357 62)
\lvec(356 62)
\ifill f:0
\move(359 61)
\lvec(360 61)
\lvec(360 62)
\lvec(359 62)
\ifill f:0
\move(361 61)
\lvec(363 61)
\lvec(363 62)
\lvec(361 62)
\ifill f:0
\move(364 61)
\lvec(365 61)
\lvec(365 62)
\lvec(364 62)
\ifill f:0
\move(367 61)
\lvec(368 61)
\lvec(368 62)
\lvec(367 62)
\ifill f:0
\move(369 61)
\lvec(370 61)
\lvec(370 62)
\lvec(369 62)
\ifill f:0
\move(372 61)
\lvec(373 61)
\lvec(373 62)
\lvec(372 62)
\ifill f:0
\move(374 61)
\lvec(375 61)
\lvec(375 62)
\lvec(374 62)
\ifill f:0
\move(376 61)
\lvec(378 61)
\lvec(378 62)
\lvec(376 62)
\ifill f:0
\move(379 61)
\lvec(380 61)
\lvec(380 62)
\lvec(379 62)
\ifill f:0
\move(381 61)
\lvec(382 61)
\lvec(382 62)
\lvec(381 62)
\ifill f:0
\move(384 61)
\lvec(385 61)
\lvec(385 62)
\lvec(384 62)
\ifill f:0
\move(386 61)
\lvec(387 61)
\lvec(387 62)
\lvec(386 62)
\ifill f:0
\move(388 61)
\lvec(389 61)
\lvec(389 62)
\lvec(388 62)
\ifill f:0
\move(390 61)
\lvec(392 61)
\lvec(392 62)
\lvec(390 62)
\ifill f:0
\move(393 61)
\lvec(394 61)
\lvec(394 62)
\lvec(393 62)
\ifill f:0
\move(395 61)
\lvec(396 61)
\lvec(396 62)
\lvec(395 62)
\ifill f:0
\move(397 61)
\lvec(398 61)
\lvec(398 62)
\lvec(397 62)
\ifill f:0
\move(399 61)
\lvec(400 61)
\lvec(400 62)
\lvec(399 62)
\ifill f:0
\move(401 61)
\lvec(403 61)
\lvec(403 62)
\lvec(401 62)
\ifill f:0
\move(404 61)
\lvec(405 61)
\lvec(405 62)
\lvec(404 62)
\ifill f:0
\move(406 61)
\lvec(407 61)
\lvec(407 62)
\lvec(406 62)
\ifill f:0
\move(408 61)
\lvec(409 61)
\lvec(409 62)
\lvec(408 62)
\ifill f:0
\move(410 61)
\lvec(411 61)
\lvec(411 62)
\lvec(410 62)
\ifill f:0
\move(412 61)
\lvec(413 61)
\lvec(413 62)
\lvec(412 62)
\ifill f:0
\move(414 61)
\lvec(415 61)
\lvec(415 62)
\lvec(414 62)
\ifill f:0
\move(416 61)
\lvec(417 61)
\lvec(417 62)
\lvec(416 62)
\ifill f:0
\move(418 61)
\lvec(419 61)
\lvec(419 62)
\lvec(418 62)
\ifill f:0
\move(420 61)
\lvec(421 61)
\lvec(421 62)
\lvec(420 62)
\ifill f:0
\move(422 61)
\lvec(423 61)
\lvec(423 62)
\lvec(422 62)
\ifill f:0
\move(424 61)
\lvec(425 61)
\lvec(425 62)
\lvec(424 62)
\ifill f:0
\move(426 61)
\lvec(427 61)
\lvec(427 62)
\lvec(426 62)
\ifill f:0
\move(428 61)
\lvec(429 61)
\lvec(429 62)
\lvec(428 62)
\ifill f:0
\move(430 61)
\lvec(431 61)
\lvec(431 62)
\lvec(430 62)
\ifill f:0
\move(432 61)
\lvec(433 61)
\lvec(433 62)
\lvec(432 62)
\ifill f:0
\move(434 61)
\lvec(435 61)
\lvec(435 62)
\lvec(434 62)
\ifill f:0
\move(436 61)
\lvec(437 61)
\lvec(437 62)
\lvec(436 62)
\ifill f:0
\move(438 61)
\lvec(439 61)
\lvec(439 62)
\lvec(438 62)
\ifill f:0
\move(440 61)
\lvec(442 61)
\lvec(442 62)
\lvec(440 62)
\ifill f:0
\move(443 61)
\lvec(444 61)
\lvec(444 62)
\lvec(443 62)
\ifill f:0
\move(445 61)
\lvec(446 61)
\lvec(446 62)
\lvec(445 62)
\ifill f:0
\move(447 61)
\lvec(448 61)
\lvec(448 62)
\lvec(447 62)
\ifill f:0
\move(449 61)
\lvec(450 61)
\lvec(450 62)
\lvec(449 62)
\ifill f:0
\move(13 62)
\lvec(14 62)
\lvec(14 63)
\lvec(13 63)
\ifill f:0
\move(19 62)
\lvec(20 62)
\lvec(20 63)
\lvec(19 63)
\ifill f:0
\move(134 62)
\lvec(135 62)
\lvec(135 63)
\lvec(134 63)
\ifill f:0
\move(140 62)
\lvec(141 62)
\lvec(141 63)
\lvec(140 63)
\ifill f:0
\move(150 62)
\lvec(151 62)
\lvec(151 63)
\lvec(150 63)
\ifill f:0
\move(162 62)
\lvec(163 62)
\lvec(163 63)
\lvec(162 63)
\ifill f:0
\move(164 62)
\lvec(165 62)
\lvec(165 63)
\lvec(164 63)
\ifill f:0
\move(169 62)
\lvec(170 62)
\lvec(170 63)
\lvec(169 63)
\ifill f:0
\move(174 62)
\lvec(175 62)
\lvec(175 63)
\lvec(174 63)
\ifill f:0
\move(177 62)
\lvec(178 62)
\lvec(178 63)
\lvec(177 63)
\ifill f:0
\move(180 62)
\lvec(181 62)
\lvec(181 63)
\lvec(180 63)
\ifill f:0
\move(183 62)
\lvec(184 62)
\lvec(184 63)
\lvec(183 63)
\ifill f:0
\move(186 62)
\lvec(187 62)
\lvec(187 63)
\lvec(186 63)
\ifill f:0
\move(190 62)
\lvec(191 62)
\lvec(191 63)
\lvec(190 63)
\ifill f:0
\move(194 62)
\lvec(195 62)
\lvec(195 63)
\lvec(194 63)
\ifill f:0
\move(198 62)
\lvec(199 62)
\lvec(199 63)
\lvec(198 63)
\ifill f:0
\move(202 62)
\lvec(204 62)
\lvec(204 63)
\lvec(202 63)
\ifill f:0
\move(208 62)
\lvec(209 62)
\lvec(209 63)
\lvec(208 63)
\ifill f:0
\move(214 62)
\lvec(216 62)
\lvec(216 63)
\lvec(214 63)
\ifill f:0
\move(223 62)
\lvec(226 62)
\lvec(226 63)
\lvec(223 63)
\ifill f:0
\move(239 62)
\lvec(251 62)
\lvec(251 63)
\lvec(239 63)
\ifill f:0
\move(264 62)
\lvec(268 62)
\lvec(268 63)
\lvec(264 63)
\ifill f:0
\move(275 62)
\lvec(278 62)
\lvec(278 63)
\lvec(275 63)
\ifill f:0
\move(283 62)
\lvec(285 62)
\lvec(285 63)
\lvec(283 63)
\ifill f:0
\move(289 62)
\lvec(292 62)
\lvec(292 63)
\lvec(289 63)
\ifill f:0
\move(295 62)
\lvec(297 62)
\lvec(297 63)
\lvec(295 63)
\ifill f:0
\move(300 62)
\lvec(302 62)
\lvec(302 63)
\lvec(300 63)
\ifill f:0
\move(305 62)
\lvec(307 62)
\lvec(307 63)
\lvec(305 63)
\ifill f:0
\move(310 62)
\lvec(312 62)
\lvec(312 63)
\lvec(310 63)
\ifill f:0
\move(314 62)
\lvec(316 62)
\lvec(316 63)
\lvec(314 63)
\ifill f:0
\move(318 62)
\lvec(320 62)
\lvec(320 63)
\lvec(318 63)
\ifill f:0
\move(322 62)
\lvec(324 62)
\lvec(324 63)
\lvec(322 63)
\ifill f:0
\move(326 62)
\lvec(328 62)
\lvec(328 63)
\lvec(326 63)
\ifill f:0
\move(330 62)
\lvec(331 62)
\lvec(331 63)
\lvec(330 63)
\ifill f:0
\move(333 62)
\lvec(335 62)
\lvec(335 63)
\lvec(333 63)
\ifill f:0
\move(337 62)
\lvec(338 62)
\lvec(338 63)
\lvec(337 63)
\ifill f:0
\move(340 62)
\lvec(342 62)
\lvec(342 63)
\lvec(340 63)
\ifill f:0
\move(343 62)
\lvec(345 62)
\lvec(345 63)
\lvec(343 63)
\ifill f:0
\move(347 62)
\lvec(348 62)
\lvec(348 63)
\lvec(347 63)
\ifill f:0
\move(350 62)
\lvec(351 62)
\lvec(351 63)
\lvec(350 63)
\ifill f:0
\move(353 62)
\lvec(354 62)
\lvec(354 63)
\lvec(353 63)
\ifill f:0
\move(356 62)
\lvec(357 62)
\lvec(357 63)
\lvec(356 63)
\ifill f:0
\move(359 62)
\lvec(360 62)
\lvec(360 63)
\lvec(359 63)
\ifill f:0
\move(361 62)
\lvec(363 62)
\lvec(363 63)
\lvec(361 63)
\ifill f:0
\move(364 62)
\lvec(365 62)
\lvec(365 63)
\lvec(364 63)
\ifill f:0
\move(367 62)
\lvec(368 62)
\lvec(368 63)
\lvec(367 63)
\ifill f:0
\move(370 62)
\lvec(371 62)
\lvec(371 63)
\lvec(370 63)
\ifill f:0
\move(372 62)
\lvec(373 62)
\lvec(373 63)
\lvec(372 63)
\ifill f:0
\move(375 62)
\lvec(376 62)
\lvec(376 63)
\lvec(375 63)
\ifill f:0
\move(377 62)
\lvec(379 62)
\lvec(379 63)
\lvec(377 63)
\ifill f:0
\move(380 62)
\lvec(381 62)
\lvec(381 63)
\lvec(380 63)
\ifill f:0
\move(383 62)
\lvec(384 62)
\lvec(384 63)
\lvec(383 63)
\ifill f:0
\move(385 62)
\lvec(386 62)
\lvec(386 63)
\lvec(385 63)
\ifill f:0
\move(387 62)
\lvec(389 62)
\lvec(389 63)
\lvec(387 63)
\ifill f:0
\move(390 62)
\lvec(391 62)
\lvec(391 63)
\lvec(390 63)
\ifill f:0
\move(392 62)
\lvec(393 62)
\lvec(393 63)
\lvec(392 63)
\ifill f:0
\move(395 62)
\lvec(396 62)
\lvec(396 63)
\lvec(395 63)
\ifill f:0
\move(397 62)
\lvec(398 62)
\lvec(398 63)
\lvec(397 63)
\ifill f:0
\move(399 62)
\lvec(400 62)
\lvec(400 63)
\lvec(399 63)
\ifill f:0
\move(402 62)
\lvec(403 62)
\lvec(403 63)
\lvec(402 63)
\ifill f:0
\move(404 62)
\lvec(405 62)
\lvec(405 63)
\lvec(404 63)
\ifill f:0
\move(406 62)
\lvec(407 62)
\lvec(407 63)
\lvec(406 63)
\ifill f:0
\move(408 62)
\lvec(409 62)
\lvec(409 63)
\lvec(408 63)
\ifill f:0
\move(410 62)
\lvec(412 62)
\lvec(412 63)
\lvec(410 63)
\ifill f:0
\move(413 62)
\lvec(414 62)
\lvec(414 63)
\lvec(413 63)
\ifill f:0
\move(415 62)
\lvec(416 62)
\lvec(416 63)
\lvec(415 63)
\ifill f:0
\move(417 62)
\lvec(418 62)
\lvec(418 63)
\lvec(417 63)
\ifill f:0
\move(419 62)
\lvec(420 62)
\lvec(420 63)
\lvec(419 63)
\ifill f:0
\move(421 62)
\lvec(422 62)
\lvec(422 63)
\lvec(421 63)
\ifill f:0
\move(423 62)
\lvec(424 62)
\lvec(424 63)
\lvec(423 63)
\ifill f:0
\move(425 62)
\lvec(426 62)
\lvec(426 63)
\lvec(425 63)
\ifill f:0
\move(427 62)
\lvec(428 62)
\lvec(428 63)
\lvec(427 63)
\ifill f:0
\move(429 62)
\lvec(430 62)
\lvec(430 63)
\lvec(429 63)
\ifill f:0
\move(431 62)
\lvec(432 62)
\lvec(432 63)
\lvec(431 63)
\ifill f:0
\move(433 62)
\lvec(434 62)
\lvec(434 63)
\lvec(433 63)
\ifill f:0
\move(435 62)
\lvec(436 62)
\lvec(436 63)
\lvec(435 63)
\ifill f:0
\move(437 62)
\lvec(438 62)
\lvec(438 63)
\lvec(437 63)
\ifill f:0
\move(439 62)
\lvec(440 62)
\lvec(440 63)
\lvec(439 63)
\ifill f:0
\move(441 62)
\lvec(442 62)
\lvec(442 63)
\lvec(441 63)
\ifill f:0
\move(443 62)
\lvec(444 62)
\lvec(444 63)
\lvec(443 63)
\ifill f:0
\move(445 62)
\lvec(446 62)
\lvec(446 63)
\lvec(445 63)
\ifill f:0
\move(447 62)
\lvec(448 62)
\lvec(448 63)
\lvec(447 63)
\ifill f:0
\move(449 62)
\lvec(450 62)
\lvec(450 63)
\lvec(449 63)
\ifill f:0
\move(121 63)
\lvec(122 63)
\lvec(122 64)
\lvec(121 64)
\ifill f:0
\move(148 63)
\lvec(149 63)
\lvec(149 64)
\lvec(148 64)
\ifill f:0
\move(153 63)
\lvec(154 63)
\lvec(154 64)
\lvec(153 64)
\ifill f:0
\move(169 63)
\lvec(170 63)
\lvec(170 64)
\lvec(169 64)
\ifill f:0
\move(171 63)
\lvec(172 63)
\lvec(172 64)
\lvec(171 64)
\ifill f:0
\move(176 63)
\lvec(177 63)
\lvec(177 64)
\lvec(176 64)
\ifill f:0
\move(184 63)
\lvec(185 63)
\lvec(185 64)
\lvec(184 64)
\ifill f:0
\move(187 63)
\lvec(188 63)
\lvec(188 64)
\lvec(187 64)
\ifill f:0
\move(194 63)
\lvec(195 63)
\lvec(195 64)
\lvec(194 64)
\ifill f:0
\move(198 63)
\lvec(199 63)
\lvec(199 64)
\lvec(198 64)
\ifill f:0
\move(202 63)
\lvec(203 63)
\lvec(203 64)
\lvec(202 64)
\ifill f:0
\move(206 63)
\lvec(207 63)
\lvec(207 64)
\lvec(206 64)
\ifill f:0
\move(211 63)
\lvec(212 63)
\lvec(212 64)
\lvec(211 64)
\ifill f:0
\move(216 63)
\lvec(218 63)
\lvec(218 64)
\lvec(216 64)
\ifill f:0
\move(223 63)
\lvec(226 63)
\lvec(226 64)
\lvec(223 64)
\ifill f:0
\move(233 63)
\lvec(237 63)
\lvec(237 64)
\lvec(233 64)
\ifill f:0
\move(269 63)
\lvec(274 63)
\lvec(274 64)
\lvec(269 64)
\ifill f:0
\move(281 63)
\lvec(284 63)
\lvec(284 64)
\lvec(281 64)
\ifill f:0
\move(289 63)
\lvec(292 63)
\lvec(292 64)
\lvec(289 64)
\ifill f:0
\move(296 63)
\lvec(299 63)
\lvec(299 64)
\lvec(296 64)
\ifill f:0
\move(302 63)
\lvec(305 63)
\lvec(305 64)
\lvec(302 64)
\ifill f:0
\move(308 63)
\lvec(310 63)
\lvec(310 64)
\lvec(308 64)
\ifill f:0
\move(313 63)
\lvec(315 63)
\lvec(315 64)
\lvec(313 64)
\ifill f:0
\move(318 63)
\lvec(320 63)
\lvec(320 64)
\lvec(318 64)
\ifill f:0
\move(322 63)
\lvec(324 63)
\lvec(324 64)
\lvec(322 64)
\ifill f:0
\move(326 63)
\lvec(328 63)
\lvec(328 64)
\lvec(326 64)
\ifill f:0
\move(330 63)
\lvec(332 63)
\lvec(332 64)
\lvec(330 64)
\ifill f:0
\move(334 63)
\lvec(336 63)
\lvec(336 64)
\lvec(334 64)
\ifill f:0
\move(338 63)
\lvec(340 63)
\lvec(340 64)
\lvec(338 64)
\ifill f:0
\move(342 63)
\lvec(343 63)
\lvec(343 64)
\lvec(342 64)
\ifill f:0
\move(345 63)
\lvec(347 63)
\lvec(347 64)
\lvec(345 64)
\ifill f:0
\move(349 63)
\lvec(350 63)
\lvec(350 64)
\lvec(349 64)
\ifill f:0
\move(352 63)
\lvec(353 63)
\lvec(353 64)
\lvec(352 64)
\ifill f:0
\move(355 63)
\lvec(357 63)
\lvec(357 64)
\lvec(355 64)
\ifill f:0
\move(358 63)
\lvec(360 63)
\lvec(360 64)
\lvec(358 64)
\ifill f:0
\move(361 63)
\lvec(363 63)
\lvec(363 64)
\lvec(361 64)
\ifill f:0
\move(364 63)
\lvec(366 63)
\lvec(366 64)
\lvec(364 64)
\ifill f:0
\move(367 63)
\lvec(369 63)
\lvec(369 64)
\lvec(367 64)
\ifill f:0
\move(370 63)
\lvec(372 63)
\lvec(372 64)
\lvec(370 64)
\ifill f:0
\move(373 63)
\lvec(374 63)
\lvec(374 64)
\lvec(373 64)
\ifill f:0
\move(376 63)
\lvec(377 63)
\lvec(377 64)
\lvec(376 64)
\ifill f:0
\move(379 63)
\lvec(380 63)
\lvec(380 64)
\lvec(379 64)
\ifill f:0
\move(381 63)
\lvec(383 63)
\lvec(383 64)
\lvec(381 64)
\ifill f:0
\move(384 63)
\lvec(385 63)
\lvec(385 64)
\lvec(384 64)
\ifill f:0
\move(387 63)
\lvec(388 63)
\lvec(388 64)
\lvec(387 64)
\ifill f:0
\move(389 63)
\lvec(390 63)
\lvec(390 64)
\lvec(389 64)
\ifill f:0
\move(392 63)
\lvec(393 63)
\lvec(393 64)
\lvec(392 64)
\ifill f:0
\move(394 63)
\lvec(395 63)
\lvec(395 64)
\lvec(394 64)
\ifill f:0
\move(397 63)
\lvec(398 63)
\lvec(398 64)
\lvec(397 64)
\ifill f:0
\move(399 63)
\lvec(400 63)
\lvec(400 64)
\lvec(399 64)
\ifill f:0
\move(402 63)
\lvec(403 63)
\lvec(403 64)
\lvec(402 64)
\ifill f:0
\move(404 63)
\lvec(405 63)
\lvec(405 64)
\lvec(404 64)
\ifill f:0
\move(406 63)
\lvec(407 63)
\lvec(407 64)
\lvec(406 64)
\ifill f:0
\move(409 63)
\lvec(410 63)
\lvec(410 64)
\lvec(409 64)
\ifill f:0
\move(411 63)
\lvec(412 63)
\lvec(412 64)
\lvec(411 64)
\ifill f:0
\move(413 63)
\lvec(414 63)
\lvec(414 64)
\lvec(413 64)
\ifill f:0
\move(416 63)
\lvec(417 63)
\lvec(417 64)
\lvec(416 64)
\ifill f:0
\move(418 63)
\lvec(419 63)
\lvec(419 64)
\lvec(418 64)
\ifill f:0
\move(420 63)
\lvec(421 63)
\lvec(421 64)
\lvec(420 64)
\ifill f:0
\move(422 63)
\lvec(423 63)
\lvec(423 64)
\lvec(422 64)
\ifill f:0
\move(425 63)
\lvec(426 63)
\lvec(426 64)
\lvec(425 64)
\ifill f:0
\move(427 63)
\lvec(428 63)
\lvec(428 64)
\lvec(427 64)
\ifill f:0
\move(429 63)
\lvec(430 63)
\lvec(430 64)
\lvec(429 64)
\ifill f:0
\move(431 63)
\lvec(432 63)
\lvec(432 64)
\lvec(431 64)
\ifill f:0
\move(433 63)
\lvec(434 63)
\lvec(434 64)
\lvec(433 64)
\ifill f:0
\move(435 63)
\lvec(436 63)
\lvec(436 64)
\lvec(435 64)
\ifill f:0
\move(437 63)
\lvec(438 63)
\lvec(438 64)
\lvec(437 64)
\ifill f:0
\move(439 63)
\lvec(440 63)
\lvec(440 64)
\lvec(439 64)
\ifill f:0
\move(441 63)
\lvec(442 63)
\lvec(442 64)
\lvec(441 64)
\ifill f:0
\move(443 63)
\lvec(444 63)
\lvec(444 64)
\lvec(443 64)
\ifill f:0
\move(445 63)
\lvec(446 63)
\lvec(446 64)
\lvec(445 64)
\ifill f:0
\move(447 63)
\lvec(448 63)
\lvec(448 64)
\lvec(447 64)
\ifill f:0
\move(449 63)
\lvec(450 63)
\lvec(450 64)
\lvec(449 64)
\ifill f:0
\move(12 64)
\lvec(13 64)
\lvec(13 65)
\lvec(12 65)
\ifill f:0
\move(121 64)
\lvec(122 64)
\lvec(122 65)
\lvec(121 65)
\ifill f:0
\move(139 64)
\lvec(140 64)
\lvec(140 65)
\lvec(139 65)
\ifill f:0
\move(154 64)
\lvec(155 64)
\lvec(155 65)
\lvec(154 65)
\ifill f:0
\move(161 64)
\lvec(162 64)
\lvec(162 65)
\lvec(161 65)
\ifill f:0
\move(163 64)
\lvec(164 64)
\lvec(164 65)
\lvec(163 65)
\ifill f:0
\move(165 64)
\lvec(166 64)
\lvec(166 65)
\lvec(165 65)
\ifill f:0
\move(167 64)
\lvec(168 64)
\lvec(168 65)
\lvec(167 65)
\ifill f:0
\move(169 64)
\lvec(170 64)
\lvec(170 65)
\lvec(169 65)
\ifill f:0
\move(171 64)
\lvec(172 64)
\lvec(172 65)
\lvec(171 65)
\ifill f:0
\move(173 64)
\lvec(174 64)
\lvec(174 65)
\lvec(173 65)
\ifill f:0
\move(180 64)
\lvec(181 64)
\lvec(181 65)
\lvec(180 65)
\ifill f:0
\move(188 64)
\lvec(189 64)
\lvec(189 65)
\lvec(188 65)
\ifill f:0
\move(191 64)
\lvec(192 64)
\lvec(192 65)
\lvec(191 65)
\ifill f:0
\move(194 64)
\lvec(195 64)
\lvec(195 65)
\lvec(194 65)
\ifill f:0
\move(197 64)
\lvec(198 64)
\lvec(198 65)
\lvec(197 65)
\ifill f:0
\move(201 64)
\lvec(202 64)
\lvec(202 65)
\lvec(201 65)
\ifill f:0
\move(205 64)
\lvec(206 64)
\lvec(206 65)
\lvec(205 65)
\ifill f:0
\move(209 64)
\lvec(210 64)
\lvec(210 65)
\lvec(209 65)
\ifill f:0
\move(213 64)
\lvec(214 64)
\lvec(214 65)
\lvec(213 65)
\ifill f:0
\move(218 64)
\lvec(220 64)
\lvec(220 65)
\lvec(218 65)
\ifill f:0
\move(224 64)
\lvec(226 64)
\lvec(226 65)
\lvec(224 65)
\ifill f:0
\move(231 64)
\lvec(233 64)
\lvec(233 65)
\lvec(231 65)
\ifill f:0
\move(241 64)
\lvec(245 64)
\lvec(245 65)
\lvec(241 65)
\ifill f:0
\move(277 64)
\lvec(282 64)
\lvec(282 65)
\lvec(277 65)
\ifill f:0
\move(289 64)
\lvec(293 64)
\lvec(293 65)
\lvec(289 65)
\ifill f:0
\move(298 64)
\lvec(301 64)
\lvec(301 65)
\lvec(298 65)
\ifill f:0
\move(305 64)
\lvec(307 64)
\lvec(307 65)
\lvec(305 65)
\ifill f:0
\move(311 64)
\lvec(313 64)
\lvec(313 65)
\lvec(311 65)
\ifill f:0
\move(317 64)
\lvec(319 64)
\lvec(319 65)
\lvec(317 65)
\ifill f:0
\move(322 64)
\lvec(324 64)
\lvec(324 65)
\lvec(322 65)
\ifill f:0
\move(327 64)
\lvec(329 64)
\lvec(329 65)
\lvec(327 65)
\ifill f:0
\move(331 64)
\lvec(333 64)
\lvec(333 65)
\lvec(331 65)
\ifill f:0
\move(336 64)
\lvec(337 64)
\lvec(337 65)
\lvec(336 65)
\ifill f:0
\move(340 64)
\lvec(341 64)
\lvec(341 65)
\lvec(340 65)
\ifill f:0
\move(344 64)
\lvec(345 64)
\lvec(345 65)
\lvec(344 65)
\ifill f:0
\move(347 64)
\lvec(349 64)
\lvec(349 65)
\lvec(347 65)
\ifill f:0
\move(351 64)
\lvec(353 64)
\lvec(353 65)
\lvec(351 65)
\ifill f:0
\move(355 64)
\lvec(356 64)
\lvec(356 65)
\lvec(355 65)
\ifill f:0
\move(358 64)
\lvec(359 64)
\lvec(359 65)
\lvec(358 65)
\ifill f:0
\move(361 64)
\lvec(363 64)
\lvec(363 65)
\lvec(361 65)
\ifill f:0
\move(365 64)
\lvec(366 64)
\lvec(366 65)
\lvec(365 65)
\ifill f:0
\move(368 64)
\lvec(369 64)
\lvec(369 65)
\lvec(368 65)
\ifill f:0
\move(371 64)
\lvec(372 64)
\lvec(372 65)
\lvec(371 65)
\ifill f:0
\move(374 64)
\lvec(375 64)
\lvec(375 65)
\lvec(374 65)
\ifill f:0
\move(377 64)
\lvec(378 64)
\lvec(378 65)
\lvec(377 65)
\ifill f:0
\move(380 64)
\lvec(381 64)
\lvec(381 65)
\lvec(380 65)
\ifill f:0
\move(383 64)
\lvec(384 64)
\lvec(384 65)
\lvec(383 65)
\ifill f:0
\move(386 64)
\lvec(387 64)
\lvec(387 65)
\lvec(386 65)
\ifill f:0
\move(388 64)
\lvec(390 64)
\lvec(390 65)
\lvec(388 65)
\ifill f:0
\move(391 64)
\lvec(392 64)
\lvec(392 65)
\lvec(391 65)
\ifill f:0
\move(394 64)
\lvec(395 64)
\lvec(395 65)
\lvec(394 65)
\ifill f:0
\move(396 64)
\lvec(398 64)
\lvec(398 65)
\lvec(396 65)
\ifill f:0
\move(399 64)
\lvec(400 64)
\lvec(400 65)
\lvec(399 65)
\ifill f:0
\move(402 64)
\lvec(403 64)
\lvec(403 65)
\lvec(402 65)
\ifill f:0
\move(404 64)
\lvec(405 64)
\lvec(405 65)
\lvec(404 65)
\ifill f:0
\move(407 64)
\lvec(408 64)
\lvec(408 65)
\lvec(407 65)
\ifill f:0
\move(409 64)
\lvec(410 64)
\lvec(410 65)
\lvec(409 65)
\ifill f:0
\move(412 64)
\lvec(413 64)
\lvec(413 65)
\lvec(412 65)
\ifill f:0
\move(414 64)
\lvec(415 64)
\lvec(415 65)
\lvec(414 65)
\ifill f:0
\move(417 64)
\lvec(418 64)
\lvec(418 65)
\lvec(417 65)
\ifill f:0
\move(419 64)
\lvec(420 64)
\lvec(420 65)
\lvec(419 65)
\ifill f:0
\move(421 64)
\lvec(422 64)
\lvec(422 65)
\lvec(421 65)
\ifill f:0
\move(424 64)
\lvec(425 64)
\lvec(425 65)
\lvec(424 65)
\ifill f:0
\move(426 64)
\lvec(427 64)
\lvec(427 65)
\lvec(426 65)
\ifill f:0
\move(428 64)
\lvec(429 64)
\lvec(429 65)
\lvec(428 65)
\ifill f:0
\move(430 64)
\lvec(431 64)
\lvec(431 65)
\lvec(430 65)
\ifill f:0
\move(433 64)
\lvec(434 64)
\lvec(434 65)
\lvec(433 65)
\ifill f:0
\move(435 64)
\lvec(436 64)
\lvec(436 65)
\lvec(435 65)
\ifill f:0
\move(437 64)
\lvec(438 64)
\lvec(438 65)
\lvec(437 65)
\ifill f:0
\move(439 64)
\lvec(440 64)
\lvec(440 65)
\lvec(439 65)
\ifill f:0
\move(441 64)
\lvec(442 64)
\lvec(442 65)
\lvec(441 65)
\ifill f:0
\move(444 64)
\lvec(445 64)
\lvec(445 65)
\lvec(444 65)
\ifill f:0
\move(446 64)
\lvec(447 64)
\lvec(447 65)
\lvec(446 65)
\ifill f:0
\move(448 64)
\lvec(449 64)
\lvec(449 65)
\lvec(448 65)
\ifill f:0
\move(450 64)
\lvec(451 64)
\lvec(451 65)
\lvec(450 65)
\ifill f:0
\move(155 65)
\lvec(156 65)
\lvec(156 66)
\lvec(155 66)
\ifill f:0
\move(160 65)
\lvec(161 65)
\lvec(161 66)
\lvec(160 66)
\ifill f:0
\move(167 65)
\lvec(168 65)
\lvec(168 66)
\lvec(167 66)
\ifill f:0
\move(169 65)
\lvec(170 65)
\lvec(170 66)
\lvec(169 66)
\ifill f:0
\move(171 65)
\lvec(172 65)
\lvec(172 66)
\lvec(171 66)
\ifill f:0
\move(173 65)
\lvec(174 65)
\lvec(174 66)
\lvec(173 66)
\ifill f:0
\move(175 65)
\lvec(176 65)
\lvec(176 66)
\lvec(175 66)
\ifill f:0
\move(177 65)
\lvec(178 65)
\lvec(178 66)
\lvec(177 66)
\ifill f:0
\move(184 65)
\lvec(185 65)
\lvec(185 66)
\lvec(184 66)
\ifill f:0
\move(189 65)
\lvec(190 65)
\lvec(190 66)
\lvec(189 66)
\ifill f:0
\move(197 65)
\lvec(198 65)
\lvec(198 66)
\lvec(197 66)
\ifill f:0
\move(207 65)
\lvec(208 65)
\lvec(208 66)
\lvec(207 66)
\ifill f:0
\move(211 65)
\lvec(212 65)
\lvec(212 66)
\lvec(211 66)
\ifill f:0
\move(215 65)
\lvec(216 65)
\lvec(216 66)
\lvec(215 66)
\ifill f:0
\move(219 65)
\lvec(221 65)
\lvec(221 66)
\lvec(219 66)
\ifill f:0
\move(224 65)
\lvec(226 65)
\lvec(226 66)
\lvec(224 66)
\ifill f:0
\move(230 65)
\lvec(232 65)
\lvec(232 66)
\lvec(230 66)
\ifill f:0
\move(237 65)
\lvec(239 65)
\lvec(239 66)
\lvec(237 66)
\ifill f:0
\move(246 65)
\lvec(249 65)
\lvec(249 66)
\lvec(246 66)
\ifill f:0
\move(263 65)
\lvec(275 65)
\lvec(275 66)
\lvec(263 66)
\ifill f:0
\move(289 65)
\lvec(293 65)
\lvec(293 66)
\lvec(289 66)
\ifill f:0
\move(300 65)
\lvec(303 65)
\lvec(303 66)
\lvec(300 66)
\ifill f:0
\move(308 65)
\lvec(311 65)
\lvec(311 66)
\lvec(308 66)
\ifill f:0
\move(315 65)
\lvec(318 65)
\lvec(318 66)
\lvec(315 66)
\ifill f:0
\move(321 65)
\lvec(324 65)
\lvec(324 66)
\lvec(321 66)
\ifill f:0
\move(327 65)
\lvec(329 65)
\lvec(329 66)
\lvec(327 66)
\ifill f:0
\move(332 65)
\lvec(334 65)
\lvec(334 66)
\lvec(332 66)
\ifill f:0
\move(337 65)
\lvec(339 65)
\lvec(339 66)
\lvec(337 66)
\ifill f:0
\move(342 65)
\lvec(343 65)
\lvec(343 66)
\lvec(342 66)
\ifill f:0
\move(346 65)
\lvec(347 65)
\lvec(347 66)
\lvec(346 66)
\ifill f:0
\move(350 65)
\lvec(352 65)
\lvec(352 66)
\lvec(350 66)
\ifill f:0
\move(354 65)
\lvec(355 65)
\lvec(355 66)
\lvec(354 66)
\ifill f:0
\move(358 65)
\lvec(359 65)
\lvec(359 66)
\lvec(358 66)
\ifill f:0
\move(361 65)
\lvec(363 65)
\lvec(363 66)
\lvec(361 66)
\ifill f:0
\move(365 65)
\lvec(366 65)
\lvec(366 66)
\lvec(365 66)
\ifill f:0
\move(368 65)
\lvec(370 65)
\lvec(370 66)
\lvec(368 66)
\ifill f:0
\move(372 65)
\lvec(373 65)
\lvec(373 66)
\lvec(372 66)
\ifill f:0
\move(375 65)
\lvec(376 65)
\lvec(376 66)
\lvec(375 66)
\ifill f:0
\move(378 65)
\lvec(380 65)
\lvec(380 66)
\lvec(378 66)
\ifill f:0
\move(381 65)
\lvec(383 65)
\lvec(383 66)
\lvec(381 66)
\ifill f:0
\move(384 65)
\lvec(386 65)
\lvec(386 66)
\lvec(384 66)
\ifill f:0
\move(387 65)
\lvec(389 65)
\lvec(389 66)
\lvec(387 66)
\ifill f:0
\move(390 65)
\lvec(392 65)
\lvec(392 66)
\lvec(390 66)
\ifill f:0
\move(393 65)
\lvec(395 65)
\lvec(395 66)
\lvec(393 66)
\ifill f:0
\move(396 65)
\lvec(397 65)
\lvec(397 66)
\lvec(396 66)
\ifill f:0
\move(399 65)
\lvec(400 65)
\lvec(400 66)
\lvec(399 66)
\ifill f:0
\move(402 65)
\lvec(403 65)
\lvec(403 66)
\lvec(402 66)
\ifill f:0
\move(404 65)
\lvec(406 65)
\lvec(406 66)
\lvec(404 66)
\ifill f:0
\move(407 65)
\lvec(408 65)
\lvec(408 66)
\lvec(407 66)
\ifill f:0
\move(410 65)
\lvec(411 65)
\lvec(411 66)
\lvec(410 66)
\ifill f:0
\move(412 65)
\lvec(414 65)
\lvec(414 66)
\lvec(412 66)
\ifill f:0
\move(415 65)
\lvec(416 65)
\lvec(416 66)
\lvec(415 66)
\ifill f:0
\move(418 65)
\lvec(419 65)
\lvec(419 66)
\lvec(418 66)
\ifill f:0
\move(420 65)
\lvec(421 65)
\lvec(421 66)
\lvec(420 66)
\ifill f:0
\move(423 65)
\lvec(424 65)
\lvec(424 66)
\lvec(423 66)
\ifill f:0
\move(425 65)
\lvec(426 65)
\lvec(426 66)
\lvec(425 66)
\ifill f:0
\move(427 65)
\lvec(428 65)
\lvec(428 66)
\lvec(427 66)
\ifill f:0
\move(430 65)
\lvec(431 65)
\lvec(431 66)
\lvec(430 66)
\ifill f:0
\move(432 65)
\lvec(433 65)
\lvec(433 66)
\lvec(432 66)
\ifill f:0
\move(435 65)
\lvec(436 65)
\lvec(436 66)
\lvec(435 66)
\ifill f:0
\move(437 65)
\lvec(438 65)
\lvec(438 66)
\lvec(437 66)
\ifill f:0
\move(439 65)
\lvec(440 65)
\lvec(440 66)
\lvec(439 66)
\ifill f:0
\move(441 65)
\lvec(442 65)
\lvec(442 66)
\lvec(441 66)
\ifill f:0
\move(444 65)
\lvec(445 65)
\lvec(445 66)
\lvec(444 66)
\ifill f:0
\move(446 65)
\lvec(447 65)
\lvec(447 66)
\lvec(446 66)
\ifill f:0
\move(448 65)
\lvec(449 65)
\lvec(449 66)
\lvec(448 66)
\ifill f:0
\move(450 65)
\lvec(451 65)
\lvec(451 66)
\lvec(450 66)
\ifill f:0
\move(12 66)
\lvec(13 66)
\lvec(13 67)
\lvec(12 67)
\ifill f:0
\move(121 66)
\lvec(123 66)
\lvec(123 67)
\lvec(121 67)
\ifill f:0
\move(124 66)
\lvec(126 66)
\lvec(126 67)
\lvec(124 67)
\ifill f:0
\move(136 66)
\lvec(137 66)
\lvec(137 67)
\lvec(136 67)
\ifill f:0
\move(142 66)
\lvec(143 66)
\lvec(143 67)
\lvec(142 67)
\ifill f:0
\move(150 66)
\lvec(151 66)
\lvec(151 67)
\lvec(150 67)
\ifill f:0
\move(162 66)
\lvec(163 66)
\lvec(163 67)
\lvec(162 67)
\ifill f:0
\move(167 66)
\lvec(168 66)
\lvec(168 67)
\lvec(167 67)
\ifill f:0
\move(169 66)
\lvec(170 66)
\lvec(170 67)
\lvec(169 67)
\ifill f:0
\move(185 66)
\lvec(186 66)
\lvec(186 67)
\lvec(185 67)
\ifill f:0
\move(187 66)
\lvec(188 66)
\lvec(188 67)
\lvec(187 67)
\ifill f:0
\move(192 66)
\lvec(193 66)
\lvec(193 67)
\lvec(192 67)
\ifill f:0
\move(197 66)
\lvec(198 66)
\lvec(198 67)
\lvec(197 67)
\ifill f:0
\move(200 66)
\lvec(201 66)
\lvec(201 67)
\lvec(200 67)
\ifill f:0
\move(203 66)
\lvec(204 66)
\lvec(204 67)
\lvec(203 67)
\ifill f:0
\move(206 66)
\lvec(207 66)
\lvec(207 67)
\lvec(206 67)
\ifill f:0
\move(209 66)
\lvec(210 66)
\lvec(210 67)
\lvec(209 67)
\ifill f:0
\move(216 66)
\lvec(217 66)
\lvec(217 67)
\lvec(216 67)
\ifill f:0
\move(220 66)
\lvec(221 66)
\lvec(221 67)
\lvec(220 67)
\ifill f:0
\move(224 66)
\lvec(226 66)
\lvec(226 67)
\lvec(224 67)
\ifill f:0
\move(229 66)
\lvec(231 66)
\lvec(231 67)
\lvec(229 67)
\ifill f:0
\move(234 66)
\lvec(236 66)
\lvec(236 67)
\lvec(234 67)
\ifill f:0
\move(241 66)
\lvec(243 66)
\lvec(243 67)
\lvec(241 67)
\ifill f:0
\move(248 66)
\lvec(251 66)
\lvec(251 67)
\lvec(248 67)
\ifill f:0
\move(259 66)
\lvec(265 66)
\lvec(265 67)
\lvec(259 67)
\ifill f:0
\move(289 66)
\lvec(296 66)
\lvec(296 67)
\lvec(289 67)
\ifill f:0
\move(304 66)
\lvec(308 66)
\lvec(308 67)
\lvec(304 67)
\ifill f:0
\move(313 66)
\lvec(316 66)
\lvec(316 67)
\lvec(313 67)
\ifill f:0
\move(321 66)
\lvec(324 66)
\lvec(324 67)
\lvec(321 67)
\ifill f:0
\move(327 66)
\lvec(330 66)
\lvec(330 67)
\lvec(327 67)
\ifill f:0
\move(333 66)
\lvec(336 66)
\lvec(336 67)
\lvec(333 67)
\ifill f:0
\move(339 66)
\lvec(341 66)
\lvec(341 67)
\lvec(339 67)
\ifill f:0
\move(344 66)
\lvec(346 66)
\lvec(346 67)
\lvec(344 67)
\ifill f:0
\move(349 66)
\lvec(350 66)
\lvec(350 67)
\lvec(349 67)
\ifill f:0
\move(353 66)
\lvec(355 66)
\lvec(355 67)
\lvec(353 67)
\ifill f:0
\move(357 66)
\lvec(359 66)
\lvec(359 67)
\lvec(357 67)
\ifill f:0
\move(361 66)
\lvec(363 66)
\lvec(363 67)
\lvec(361 67)
\ifill f:0
\move(365 66)
\lvec(367 66)
\lvec(367 67)
\lvec(365 67)
\ifill f:0
\move(369 66)
\lvec(371 66)
\lvec(371 67)
\lvec(369 67)
\ifill f:0
\move(373 66)
\lvec(374 66)
\lvec(374 67)
\lvec(373 67)
\ifill f:0
\move(376 66)
\lvec(378 66)
\lvec(378 67)
\lvec(376 67)
\ifill f:0
\move(380 66)
\lvec(381 66)
\lvec(381 67)
\lvec(380 67)
\ifill f:0
\move(383 66)
\lvec(385 66)
\lvec(385 67)
\lvec(383 67)
\ifill f:0
\move(386 66)
\lvec(388 66)
\lvec(388 67)
\lvec(386 67)
\ifill f:0
\move(390 66)
\lvec(391 66)
\lvec(391 67)
\lvec(390 67)
\ifill f:0
\move(393 66)
\lvec(394 66)
\lvec(394 67)
\lvec(393 67)
\ifill f:0
\move(396 66)
\lvec(397 66)
\lvec(397 67)
\lvec(396 67)
\ifill f:0
\move(399 66)
\lvec(400 66)
\lvec(400 67)
\lvec(399 67)
\ifill f:0
\move(402 66)
\lvec(403 66)
\lvec(403 67)
\lvec(402 67)
\ifill f:0
\move(405 66)
\lvec(406 66)
\lvec(406 67)
\lvec(405 67)
\ifill f:0
\move(408 66)
\lvec(409 66)
\lvec(409 67)
\lvec(408 67)
\ifill f:0
\move(410 66)
\lvec(412 66)
\lvec(412 67)
\lvec(410 67)
\ifill f:0
\move(413 66)
\lvec(415 66)
\lvec(415 67)
\lvec(413 67)
\ifill f:0
\move(416 66)
\lvec(417 66)
\lvec(417 67)
\lvec(416 67)
\ifill f:0
\move(419 66)
\lvec(420 66)
\lvec(420 67)
\lvec(419 67)
\ifill f:0
\move(421 66)
\lvec(423 66)
\lvec(423 67)
\lvec(421 67)
\ifill f:0
\move(424 66)
\lvec(425 66)
\lvec(425 67)
\lvec(424 67)
\ifill f:0
\move(427 66)
\lvec(428 66)
\lvec(428 67)
\lvec(427 67)
\ifill f:0
\move(429 66)
\lvec(430 66)
\lvec(430 67)
\lvec(429 67)
\ifill f:0
\move(432 66)
\lvec(433 66)
\lvec(433 67)
\lvec(432 67)
\ifill f:0
\move(434 66)
\lvec(435 66)
\lvec(435 67)
\lvec(434 67)
\ifill f:0
\move(437 66)
\lvec(438 66)
\lvec(438 67)
\lvec(437 67)
\ifill f:0
\move(439 66)
\lvec(440 66)
\lvec(440 67)
\lvec(439 67)
\ifill f:0
\move(441 66)
\lvec(443 66)
\lvec(443 67)
\lvec(441 67)
\ifill f:0
\move(444 66)
\lvec(445 66)
\lvec(445 67)
\lvec(444 67)
\ifill f:0
\move(446 66)
\lvec(447 66)
\lvec(447 67)
\lvec(446 67)
\ifill f:0
\move(449 66)
\lvec(450 66)
\lvec(450 67)
\lvec(449 67)
\ifill f:0
\move(13 67)
\lvec(14 67)
\lvec(14 68)
\lvec(13 68)
\ifill f:0
\move(15 67)
\lvec(16 67)
\lvec(16 68)
\lvec(15 68)
\ifill f:0
\move(19 67)
\lvec(20 67)
\lvec(20 68)
\lvec(19 68)
\ifill f:0
\move(141 67)
\lvec(142 67)
\lvec(142 68)
\lvec(141 68)
\ifill f:0
\move(151 67)
\lvec(152 67)
\lvec(152 68)
\lvec(151 68)
\ifill f:0
\move(158 67)
\lvec(159 67)
\lvec(159 68)
\lvec(158 68)
\ifill f:0
\move(161 67)
\lvec(162 67)
\lvec(162 68)
\lvec(161 68)
\ifill f:0
\move(164 67)
\lvec(165 67)
\lvec(165 68)
\lvec(164 68)
\ifill f:0
\move(169 67)
\lvec(170 67)
\lvec(170 68)
\lvec(169 68)
\ifill f:0
\move(174 67)
\lvec(175 67)
\lvec(175 68)
\lvec(174 68)
\ifill f:0
\move(176 67)
\lvec(177 67)
\lvec(177 68)
\lvec(176 68)
\ifill f:0
\move(188 67)
\lvec(189 67)
\lvec(189 68)
\lvec(188 68)
\ifill f:0
\move(190 67)
\lvec(191 67)
\lvec(191 68)
\lvec(190 68)
\ifill f:0
\move(197 67)
\lvec(198 67)
\lvec(198 68)
\lvec(197 68)
\ifill f:0
\move(202 67)
\lvec(203 67)
\lvec(203 68)
\lvec(202 68)
\ifill f:0
\move(205 67)
\lvec(206 67)
\lvec(206 68)
\lvec(205 68)
\ifill f:0
\move(208 67)
\lvec(209 67)
\lvec(209 68)
\lvec(208 68)
\ifill f:0
\move(211 67)
\lvec(212 67)
\lvec(212 68)
\lvec(211 68)
\ifill f:0
\move(214 67)
\lvec(215 67)
\lvec(215 68)
\lvec(214 68)
\ifill f:0
\move(217 67)
\lvec(218 67)
\lvec(218 68)
\lvec(217 68)
\ifill f:0
\move(221 67)
\lvec(222 67)
\lvec(222 68)
\lvec(221 68)
\ifill f:0
\move(225 67)
\lvec(226 67)
\lvec(226 68)
\lvec(225 68)
\ifill f:0
\move(229 67)
\lvec(230 67)
\lvec(230 68)
\lvec(229 68)
\ifill f:0
\move(233 67)
\lvec(234 67)
\lvec(234 68)
\lvec(233 68)
\ifill f:0
\move(238 67)
\lvec(240 67)
\lvec(240 68)
\lvec(238 68)
\ifill f:0
\move(244 67)
\lvec(245 67)
\lvec(245 68)
\lvec(244 68)
\ifill f:0
\move(250 67)
\lvec(252 67)
\lvec(252 68)
\lvec(250 68)
\ifill f:0
\move(259 67)
\lvec(262 67)
\lvec(262 68)
\lvec(259 68)
\ifill f:0
\move(272 67)
\lvec(285 67)
\lvec(285 68)
\lvec(272 68)
\ifill f:0
\move(286 67)
\lvec(299 67)
\lvec(299 68)
\lvec(286 68)
\ifill f:0
\move(310 67)
\lvec(314 67)
\lvec(314 68)
\lvec(310 68)
\ifill f:0
\move(320 67)
\lvec(323 67)
\lvec(323 68)
\lvec(320 68)
\ifill f:0
\move(328 67)
\lvec(331 67)
\lvec(331 68)
\lvec(328 68)
\ifill f:0
\move(335 67)
\lvec(337 67)
\lvec(337 68)
\lvec(335 68)
\ifill f:0
\move(341 67)
\lvec(343 67)
\lvec(343 68)
\lvec(341 68)
\ifill f:0
\move(347 67)
\lvec(349 67)
\lvec(349 68)
\lvec(347 68)
\ifill f:0
\move(352 67)
\lvec(354 67)
\lvec(354 68)
\lvec(352 68)
\ifill f:0
\move(357 67)
\lvec(359 67)
\lvec(359 68)
\lvec(357 68)
\ifill f:0
\move(361 67)
\lvec(363 67)
\lvec(363 68)
\lvec(361 68)
\ifill f:0
\move(366 67)
\lvec(367 67)
\lvec(367 68)
\lvec(366 68)
\ifill f:0
\move(370 67)
\lvec(372 67)
\lvec(372 68)
\lvec(370 68)
\ifill f:0
\move(374 67)
\lvec(376 67)
\lvec(376 68)
\lvec(374 68)
\ifill f:0
\move(378 67)
\lvec(379 67)
\lvec(379 68)
\lvec(378 68)
\ifill f:0
\move(382 67)
\lvec(383 67)
\lvec(383 68)
\lvec(382 68)
\ifill f:0
\move(385 67)
\lvec(387 67)
\lvec(387 68)
\lvec(385 68)
\ifill f:0
\move(389 67)
\lvec(390 67)
\lvec(390 68)
\lvec(389 68)
\ifill f:0
\move(392 67)
\lvec(394 67)
\lvec(394 68)
\lvec(392 68)
\ifill f:0
\move(395 67)
\lvec(397 67)
\lvec(397 68)
\lvec(395 68)
\ifill f:0
\move(399 67)
\lvec(400 67)
\lvec(400 68)
\lvec(399 68)
\ifill f:0
\move(402 67)
\lvec(403 67)
\lvec(403 68)
\lvec(402 68)
\ifill f:0
\move(405 67)
\lvec(406 67)
\lvec(406 68)
\lvec(405 68)
\ifill f:0
\move(408 67)
\lvec(410 67)
\lvec(410 68)
\lvec(408 68)
\ifill f:0
\move(411 67)
\lvec(413 67)
\lvec(413 68)
\lvec(411 68)
\ifill f:0
\move(414 67)
\lvec(415 67)
\lvec(415 68)
\lvec(414 68)
\ifill f:0
\move(417 67)
\lvec(418 67)
\lvec(418 68)
\lvec(417 68)
\ifill f:0
\move(420 67)
\lvec(421 67)
\lvec(421 68)
\lvec(420 68)
\ifill f:0
\move(423 67)
\lvec(424 67)
\lvec(424 68)
\lvec(423 68)
\ifill f:0
\move(426 67)
\lvec(427 67)
\lvec(427 68)
\lvec(426 68)
\ifill f:0
\move(428 67)
\lvec(429 67)
\lvec(429 68)
\lvec(428 68)
\ifill f:0
\move(431 67)
\lvec(432 67)
\lvec(432 68)
\lvec(431 68)
\ifill f:0
\move(434 67)
\lvec(435 67)
\lvec(435 68)
\lvec(434 68)
\ifill f:0
\move(436 67)
\lvec(437 67)
\lvec(437 68)
\lvec(436 68)
\ifill f:0
\move(439 67)
\lvec(440 67)
\lvec(440 68)
\lvec(439 68)
\ifill f:0
\move(441 67)
\lvec(443 67)
\lvec(443 68)
\lvec(441 68)
\ifill f:0
\move(444 67)
\lvec(445 67)
\lvec(445 68)
\lvec(444 68)
\ifill f:0
\move(446 67)
\lvec(448 67)
\lvec(448 68)
\lvec(446 68)
\ifill f:0
\move(449 67)
\lvec(450 67)
\lvec(450 68)
\lvec(449 68)
\ifill f:0
\move(12 68)
\lvec(13 68)
\lvec(13 69)
\lvec(12 69)
\ifill f:0
\move(15 68)
\lvec(16 68)
\lvec(16 69)
\lvec(15 69)
\ifill f:0
\move(140 68)
\lvec(141 68)
\lvec(141 69)
\lvec(140 69)
\ifill f:0
\move(147 68)
\lvec(148 68)
\lvec(148 69)
\lvec(147 69)
\ifill f:0
\move(160 68)
\lvec(161 68)
\lvec(161 69)
\lvec(160 69)
\ifill f:0
\move(169 68)
\lvec(170 68)
\lvec(170 69)
\lvec(169 69)
\ifill f:0
\move(179 68)
\lvec(180 68)
\lvec(180 69)
\lvec(179 69)
\ifill f:0
\move(197 68)
\lvec(198 68)
\lvec(198 69)
\lvec(197 69)
\ifill f:0
\move(204 68)
\lvec(205 68)
\lvec(205 69)
\lvec(204 69)
\ifill f:0
\move(212 68)
\lvec(213 68)
\lvec(213 69)
\lvec(212 69)
\ifill f:0
\move(215 68)
\lvec(216 68)
\lvec(216 69)
\lvec(215 69)
\ifill f:0
\move(218 68)
\lvec(219 68)
\lvec(219 69)
\lvec(218 69)
\ifill f:0
\move(221 68)
\lvec(222 68)
\lvec(222 69)
\lvec(221 69)
\ifill f:0
\move(225 68)
\lvec(226 68)
\lvec(226 69)
\lvec(225 69)
\ifill f:0
\move(228 68)
\lvec(229 68)
\lvec(229 69)
\lvec(228 69)
\ifill f:0
\move(232 68)
\lvec(233 68)
\lvec(233 69)
\lvec(232 69)
\ifill f:0
\move(236 68)
\lvec(237 68)
\lvec(237 69)
\lvec(236 69)
\ifill f:0
\move(241 68)
\lvec(242 68)
\lvec(242 69)
\lvec(241 69)
\ifill f:0
\move(246 68)
\lvec(247 68)
\lvec(247 69)
\lvec(246 69)
\ifill f:0
\move(252 68)
\lvec(253 68)
\lvec(253 69)
\lvec(252 69)
\ifill f:0
\move(258 68)
\lvec(260 68)
\lvec(260 69)
\lvec(258 69)
\ifill f:0
\move(267 68)
\lvec(270 68)
\lvec(270 69)
\lvec(267 69)
\ifill f:0
\move(280 68)
\lvec(293 68)
\lvec(293 69)
\lvec(280 69)
\ifill f:0
\move(294 68)
\lvec(308 68)
\lvec(308 69)
\lvec(294 69)
\ifill f:0
\move(318 68)
\lvec(323 68)
\lvec(323 69)
\lvec(318 69)
\ifill f:0
\move(329 68)
\lvec(332 68)
\lvec(332 69)
\lvec(329 69)
\ifill f:0
\move(337 68)
\lvec(340 68)
\lvec(340 69)
\lvec(337 69)
\ifill f:0
\move(344 68)
\lvec(347 68)
\lvec(347 69)
\lvec(344 69)
\ifill f:0
\move(350 68)
\lvec(353 68)
\lvec(353 69)
\lvec(350 69)
\ifill f:0
\move(356 68)
\lvec(358 68)
\lvec(358 69)
\lvec(356 69)
\ifill f:0
\move(361 68)
\lvec(363 68)
\lvec(363 69)
\lvec(361 69)
\ifill f:0
\move(366 68)
\lvec(368 68)
\lvec(368 69)
\lvec(366 69)
\ifill f:0
\move(371 68)
\lvec(373 68)
\lvec(373 69)
\lvec(371 69)
\ifill f:0
\move(375 68)
\lvec(377 68)
\lvec(377 69)
\lvec(375 69)
\ifill f:0
\move(380 68)
\lvec(381 68)
\lvec(381 69)
\lvec(380 69)
\ifill f:0
\move(384 68)
\lvec(385 68)
\lvec(385 69)
\lvec(384 69)
\ifill f:0
\move(388 68)
\lvec(389 68)
\lvec(389 69)
\lvec(388 69)
\ifill f:0
\move(391 68)
\lvec(393 68)
\lvec(393 69)
\lvec(391 69)
\ifill f:0
\move(395 68)
\lvec(397 68)
\lvec(397 69)
\lvec(395 69)
\ifill f:0
\move(399 68)
\lvec(400 68)
\lvec(400 69)
\lvec(399 69)
\ifill f:0
\move(402 68)
\lvec(404 68)
\lvec(404 69)
\lvec(402 69)
\ifill f:0
\move(405 68)
\lvec(407 68)
\lvec(407 69)
\lvec(405 69)
\ifill f:0
\move(409 68)
\lvec(410 68)
\lvec(410 69)
\lvec(409 69)
\ifill f:0
\move(412 68)
\lvec(413 68)
\lvec(413 69)
\lvec(412 69)
\ifill f:0
\move(415 68)
\lvec(417 68)
\lvec(417 69)
\lvec(415 69)
\ifill f:0
\move(418 68)
\lvec(420 68)
\lvec(420 69)
\lvec(418 69)
\ifill f:0
\move(421 68)
\lvec(423 68)
\lvec(423 69)
\lvec(421 69)
\ifill f:0
\move(424 68)
\lvec(426 68)
\lvec(426 69)
\lvec(424 69)
\ifill f:0
\move(427 68)
\lvec(429 68)
\lvec(429 69)
\lvec(427 69)
\ifill f:0
\move(430 68)
\lvec(432 68)
\lvec(432 69)
\lvec(430 69)
\ifill f:0
\move(433 68)
\lvec(434 68)
\lvec(434 69)
\lvec(433 69)
\ifill f:0
\move(436 68)
\lvec(437 68)
\lvec(437 69)
\lvec(436 69)
\ifill f:0
\move(439 68)
\lvec(440 68)
\lvec(440 69)
\lvec(439 69)
\ifill f:0
\move(441 68)
\lvec(443 68)
\lvec(443 69)
\lvec(441 69)
\ifill f:0
\move(444 68)
\lvec(445 68)
\lvec(445 69)
\lvec(444 69)
\ifill f:0
\move(447 68)
\lvec(448 68)
\lvec(448 69)
\lvec(447 69)
\ifill f:0
\move(449 68)
\lvec(451 68)
\lvec(451 69)
\lvec(449 69)
\ifill f:0
\move(15 69)
\lvec(16 69)
\lvec(16 70)
\lvec(15 70)
\ifill f:0
\move(132 69)
\lvec(134 69)
\lvec(134 70)
\lvec(132 70)
\ifill f:0
\move(135 69)
\lvec(137 69)
\lvec(137 70)
\lvec(135 70)
\ifill f:0
\move(148 69)
\lvec(149 69)
\lvec(149 70)
\lvec(148 70)
\ifill f:0
\move(154 69)
\lvec(155 69)
\lvec(155 70)
\lvec(154 70)
\ifill f:0
\move(166 69)
\lvec(167 69)
\lvec(167 70)
\lvec(166 70)
\ifill f:0
\move(169 69)
\lvec(170 69)
\lvec(170 70)
\lvec(169 70)
\ifill f:0
\move(172 69)
\lvec(173 69)
\lvec(173 70)
\lvec(172 70)
\ifill f:0
\move(175 69)
\lvec(176 69)
\lvec(176 70)
\lvec(175 70)
\ifill f:0
\move(180 69)
\lvec(181 69)
\lvec(181 70)
\lvec(180 70)
\ifill f:0
\move(187 69)
\lvec(188 69)
\lvec(188 70)
\lvec(187 70)
\ifill f:0
\move(189 69)
\lvec(190 69)
\lvec(190 70)
\lvec(189 70)
\ifill f:0
\move(191 69)
\lvec(192 69)
\lvec(192 70)
\lvec(191 70)
\ifill f:0
\move(195 69)
\lvec(196 69)
\lvec(196 70)
\lvec(195 70)
\ifill f:0
\move(197 69)
\lvec(198 69)
\lvec(198 70)
\lvec(197 70)
\ifill f:0
\move(199 69)
\lvec(200 69)
\lvec(200 70)
\lvec(199 70)
\ifill f:0
\move(201 69)
\lvec(202 69)
\lvec(202 70)
\lvec(201 70)
\ifill f:0
\move(206 69)
\lvec(207 69)
\lvec(207 70)
\lvec(206 70)
\ifill f:0
\move(208 69)
\lvec(209 69)
\lvec(209 70)
\lvec(208 70)
\ifill f:0
\move(213 69)
\lvec(214 69)
\lvec(214 70)
\lvec(213 70)
\ifill f:0
\move(216 69)
\lvec(217 69)
\lvec(217 70)
\lvec(216 70)
\ifill f:0
\move(219 69)
\lvec(220 69)
\lvec(220 70)
\lvec(219 70)
\ifill f:0
\move(225 69)
\lvec(226 69)
\lvec(226 70)
\lvec(225 70)
\ifill f:0
\move(228 69)
\lvec(229 69)
\lvec(229 70)
\lvec(228 70)
\ifill f:0
\move(231 69)
\lvec(232 69)
\lvec(232 70)
\lvec(231 70)
\ifill f:0
\move(235 69)
\lvec(236 69)
\lvec(236 70)
\lvec(235 70)
\ifill f:0
\move(239 69)
\lvec(240 69)
\lvec(240 70)
\lvec(239 70)
\ifill f:0
\move(243 69)
\lvec(244 69)
\lvec(244 70)
\lvec(243 70)
\ifill f:0
\move(247 69)
\lvec(249 69)
\lvec(249 70)
\lvec(247 70)
\ifill f:0
\move(252 69)
\lvec(254 69)
\lvec(254 70)
\lvec(252 70)
\ifill f:0
\move(258 69)
\lvec(260 69)
\lvec(260 70)
\lvec(258 70)
\ifill f:0
\move(264 69)
\lvec(267 69)
\lvec(267 70)
\lvec(264 70)
\ifill f:0
\move(272 69)
\lvec(276 69)
\lvec(276 70)
\lvec(272 70)
\ifill f:0
\move(284 69)
\lvec(290 69)
\lvec(290 70)
\lvec(284 70)
\ifill f:0
\move(315 69)
\lvec(322 69)
\lvec(322 70)
\lvec(315 70)
\ifill f:0
\move(330 69)
\lvec(335 69)
\lvec(335 70)
\lvec(330 70)
\ifill f:0
\move(340 69)
\lvec(344 69)
\lvec(344 70)
\lvec(340 70)
\ifill f:0
\move(348 69)
\lvec(351 69)
\lvec(351 70)
\lvec(348 70)
\ifill f:0
\move(355 69)
\lvec(358 69)
\lvec(358 70)
\lvec(355 70)
\ifill f:0
\move(361 69)
\lvec(364 69)
\lvec(364 70)
\lvec(361 70)
\ifill f:0
\move(367 69)
\lvec(369 69)
\lvec(369 70)
\lvec(367 70)
\ifill f:0
\move(372 69)
\lvec(374 69)
\lvec(374 70)
\lvec(372 70)
\ifill f:0
\move(377 69)
\lvec(379 69)
\lvec(379 70)
\lvec(377 70)
\ifill f:0
\move(382 69)
\lvec(384 69)
\lvec(384 70)
\lvec(382 70)
\ifill f:0
\move(386 69)
\lvec(388 69)
\lvec(388 70)
\lvec(386 70)
\ifill f:0
\move(390 69)
\lvec(392 69)
\lvec(392 70)
\lvec(390 70)
\ifill f:0
\move(395 69)
\lvec(396 69)
\lvec(396 70)
\lvec(395 70)
\ifill f:0
\move(398 69)
\lvec(400 69)
\lvec(400 70)
\lvec(398 70)
\ifill f:0
\move(402 69)
\lvec(404 69)
\lvec(404 70)
\lvec(402 70)
\ifill f:0
\move(406 69)
\lvec(408 69)
\lvec(408 70)
\lvec(406 70)
\ifill f:0
\move(410 69)
\lvec(411 69)
\lvec(411 70)
\lvec(410 70)
\ifill f:0
\move(413 69)
\lvec(415 69)
\lvec(415 70)
\lvec(413 70)
\ifill f:0
\move(416 69)
\lvec(418 69)
\lvec(418 70)
\lvec(416 70)
\ifill f:0
\move(420 69)
\lvec(421 69)
\lvec(421 70)
\lvec(420 70)
\ifill f:0
\move(423 69)
\lvec(425 69)
\lvec(425 70)
\lvec(423 70)
\ifill f:0
\move(426 69)
\lvec(428 69)
\lvec(428 70)
\lvec(426 70)
\ifill f:0
\move(429 69)
\lvec(431 69)
\lvec(431 70)
\lvec(429 70)
\ifill f:0
\move(432 69)
\lvec(434 69)
\lvec(434 70)
\lvec(432 70)
\ifill f:0
\move(436 69)
\lvec(437 69)
\lvec(437 70)
\lvec(436 70)
\ifill f:0
\move(438 69)
\lvec(440 69)
\lvec(440 70)
\lvec(438 70)
\ifill f:0
\move(441 69)
\lvec(443 69)
\lvec(443 70)
\lvec(441 70)
\ifill f:0
\move(444 69)
\lvec(446 69)
\lvec(446 70)
\lvec(444 70)
\ifill f:0
\move(447 69)
\lvec(448 69)
\lvec(448 70)
\lvec(447 70)
\ifill f:0
\move(450 69)
\lvec(451 69)
\lvec(451 70)
\lvec(450 70)
\ifill f:0
\move(12 70)
\lvec(13 70)
\lvec(13 71)
\lvec(12 71)
\ifill f:0
\move(19 70)
\lvec(20 70)
\lvec(20 71)
\lvec(19 71)
\ifill f:0
\move(169 70)
\lvec(170 70)
\lvec(170 71)
\lvec(169 71)
\ifill f:0
\move(181 70)
\lvec(182 70)
\lvec(182 71)
\lvec(181 71)
\ifill f:0
\move(186 70)
\lvec(187 70)
\lvec(187 71)
\lvec(186 71)
\ifill f:0
\move(193 70)
\lvec(194 70)
\lvec(194 71)
\lvec(193 71)
\ifill f:0
\move(195 70)
\lvec(196 70)
\lvec(196 71)
\lvec(195 71)
\ifill f:0
\move(197 70)
\lvec(198 70)
\lvec(198 71)
\lvec(197 71)
\ifill f:0
\move(201 70)
\lvec(202 70)
\lvec(202 71)
\lvec(201 71)
\ifill f:0
\move(203 70)
\lvec(204 70)
\lvec(204 71)
\lvec(203 71)
\ifill f:0
\move(205 70)
\lvec(206 70)
\lvec(206 71)
\lvec(205 71)
\ifill f:0
\move(212 70)
\lvec(213 70)
\lvec(213 71)
\lvec(212 71)
\ifill f:0
\move(222 70)
\lvec(223 70)
\lvec(223 71)
\lvec(222 71)
\ifill f:0
\move(225 70)
\lvec(226 70)
\lvec(226 71)
\lvec(225 71)
\ifill f:0
\move(231 70)
\lvec(232 70)
\lvec(232 71)
\lvec(231 71)
\ifill f:0
\move(234 70)
\lvec(235 70)
\lvec(235 71)
\lvec(234 71)
\ifill f:0
\move(237 70)
\lvec(238 70)
\lvec(238 71)
\lvec(237 71)
\ifill f:0
\move(241 70)
\lvec(242 70)
\lvec(242 71)
\lvec(241 71)
\ifill f:0
\move(245 70)
\lvec(246 70)
\lvec(246 71)
\lvec(245 71)
\ifill f:0
\move(249 70)
\lvec(250 70)
\lvec(250 71)
\lvec(249 71)
\ifill f:0
\move(253 70)
\lvec(254 70)
\lvec(254 71)
\lvec(253 71)
\ifill f:0
\move(258 70)
\lvec(259 70)
\lvec(259 71)
\lvec(258 71)
\ifill f:0
\move(263 70)
\lvec(265 70)
\lvec(265 71)
\lvec(263 71)
\ifill f:0
\move(269 70)
\lvec(271 70)
\lvec(271 71)
\lvec(269 71)
\ifill f:0
\move(277 70)
\lvec(279 70)
\lvec(279 71)
\lvec(277 71)
\ifill f:0
\move(286 70)
\lvec(290 70)
\lvec(290 71)
\lvec(286 71)
\ifill f:0
\move(304 70)
\lvec(319 70)
\lvec(319 71)
\lvec(304 71)
\ifill f:0
\move(333 70)
\lvec(338 70)
\lvec(338 71)
\lvec(333 71)
\ifill f:0
\move(345 70)
\lvec(348 70)
\lvec(348 71)
\lvec(345 71)
\ifill f:0
\move(354 70)
\lvec(357 70)
\lvec(357 71)
\lvec(354 71)
\ifill f:0
\move(361 70)
\lvec(364 70)
\lvec(364 71)
\lvec(361 71)
\ifill f:0
\move(368 70)
\lvec(370 70)
\lvec(370 71)
\lvec(368 71)
\ifill f:0
\move(374 70)
\lvec(376 70)
\lvec(376 71)
\lvec(374 71)
\ifill f:0
\move(379 70)
\lvec(381 70)
\lvec(381 71)
\lvec(379 71)
\ifill f:0
\move(384 70)
\lvec(386 70)
\lvec(386 71)
\lvec(384 71)
\ifill f:0
\move(389 70)
\lvec(391 70)
\lvec(391 71)
\lvec(389 71)
\ifill f:0
\move(394 70)
\lvec(396 70)
\lvec(396 71)
\lvec(394 71)
\ifill f:0
\move(398 70)
\lvec(400 70)
\lvec(400 71)
\lvec(398 71)
\ifill f:0
\move(402 70)
\lvec(404 70)
\lvec(404 71)
\lvec(402 71)
\ifill f:0
\move(407 70)
\lvec(408 70)
\lvec(408 71)
\lvec(407 71)
\ifill f:0
\move(410 70)
\lvec(412 70)
\lvec(412 71)
\lvec(410 71)
\ifill f:0
\move(414 70)
\lvec(416 70)
\lvec(416 71)
\lvec(414 71)
\ifill f:0
\move(418 70)
\lvec(419 70)
\lvec(419 71)
\lvec(418 71)
\ifill f:0
\move(422 70)
\lvec(423 70)
\lvec(423 71)
\lvec(422 71)
\ifill f:0
\move(425 70)
\lvec(426 70)
\lvec(426 71)
\lvec(425 71)
\ifill f:0
\move(428 70)
\lvec(430 70)
\lvec(430 71)
\lvec(428 71)
\ifill f:0
\move(432 70)
\lvec(433 70)
\lvec(433 71)
\lvec(432 71)
\ifill f:0
\move(435 70)
\lvec(436 70)
\lvec(436 71)
\lvec(435 71)
\ifill f:0
\move(438 70)
\lvec(440 70)
\lvec(440 71)
\lvec(438 71)
\ifill f:0
\move(441 70)
\lvec(443 70)
\lvec(443 71)
\lvec(441 71)
\ifill f:0
\move(444 70)
\lvec(446 70)
\lvec(446 71)
\lvec(444 71)
\ifill f:0
\move(448 70)
\lvec(449 70)
\lvec(449 71)
\lvec(448 71)
\ifill f:0
\move(450 70)
\lvec(451 70)
\lvec(451 71)
\lvec(450 71)
\ifill f:0
\move(22 71)
\lvec(23 71)
\lvec(23 72)
\lvec(22 72)
\ifill f:0
\move(152 71)
\lvec(153 71)
\lvec(153 72)
\lvec(152 72)
\ifill f:0
\move(169 71)
\lvec(170 71)
\lvec(170 72)
\lvec(169 72)
\ifill f:0
\move(185 71)
\lvec(186 71)
\lvec(186 72)
\lvec(185 72)
\ifill f:0
\move(195 71)
\lvec(196 71)
\lvec(196 72)
\lvec(195 72)
\ifill f:0
\move(213 71)
\lvec(214 71)
\lvec(214 72)
\lvec(213 72)
\ifill f:0
\move(215 71)
\lvec(216 71)
\lvec(216 72)
\lvec(215 72)
\ifill f:0
\move(220 71)
\lvec(221 71)
\lvec(221 72)
\lvec(220 72)
\ifill f:0
\move(225 71)
\lvec(226 71)
\lvec(226 72)
\lvec(225 72)
\ifill f:0
\move(230 71)
\lvec(231 71)
\lvec(231 72)
\lvec(230 72)
\ifill f:0
\move(233 71)
\lvec(234 71)
\lvec(234 72)
\lvec(233 72)
\ifill f:0
\move(236 71)
\lvec(237 71)
\lvec(237 72)
\lvec(236 72)
\ifill f:0
\move(239 71)
\lvec(240 71)
\lvec(240 72)
\lvec(239 72)
\ifill f:0
\move(246 71)
\lvec(247 71)
\lvec(247 72)
\lvec(246 72)
\ifill f:0
\move(250 71)
\lvec(251 71)
\lvec(251 72)
\lvec(250 72)
\ifill f:0
\move(254 71)
\lvec(255 71)
\lvec(255 72)
\lvec(254 72)
\ifill f:0
\move(258 71)
\lvec(259 71)
\lvec(259 72)
\lvec(258 72)
\ifill f:0
\move(262 71)
\lvec(263 71)
\lvec(263 72)
\lvec(262 72)
\ifill f:0
\move(267 71)
\lvec(269 71)
\lvec(269 72)
\lvec(267 72)
\ifill f:0
\move(273 71)
\lvec(274 71)
\lvec(274 72)
\lvec(273 72)
\ifill f:0
\move(279 71)
\lvec(281 71)
\lvec(281 72)
\lvec(279 72)
\ifill f:0
\move(287 71)
\lvec(290 71)
\lvec(290 72)
\lvec(287 72)
\ifill f:0
\move(298 71)
\lvec(302 71)
\lvec(302 72)
\lvec(298 72)
\ifill f:0
\move(338 71)
\lvec(344 71)
\lvec(344 72)
\lvec(338 72)
\ifill f:0
\move(352 71)
\lvec(355 71)
\lvec(355 72)
\lvec(352 72)
\ifill f:0
\move(361 71)
\lvec(364 71)
\lvec(364 72)
\lvec(361 72)
\ifill f:0
\move(369 71)
\lvec(372 71)
\lvec(372 72)
\lvec(369 72)
\ifill f:0
\move(376 71)
\lvec(378 71)
\lvec(378 72)
\lvec(376 72)
\ifill f:0
\move(382 71)
\lvec(384 71)
\lvec(384 72)
\lvec(382 72)
\ifill f:0
\move(388 71)
\lvec(390 71)
\lvec(390 72)
\lvec(388 72)
\ifill f:0
\move(393 71)
\lvec(395 71)
\lvec(395 72)
\lvec(393 72)
\ifill f:0
\move(398 71)
\lvec(400 71)
\lvec(400 72)
\lvec(398 72)
\ifill f:0
\move(403 71)
\lvec(404 71)
\lvec(404 72)
\lvec(403 72)
\ifill f:0
\move(407 71)
\lvec(409 71)
\lvec(409 72)
\lvec(407 72)
\ifill f:0
\move(411 71)
\lvec(413 71)
\lvec(413 72)
\lvec(411 72)
\ifill f:0
\move(416 71)
\lvec(417 71)
\lvec(417 72)
\lvec(416 72)
\ifill f:0
\move(420 71)
\lvec(421 71)
\lvec(421 72)
\lvec(420 72)
\ifill f:0
\move(424 71)
\lvec(425 71)
\lvec(425 72)
\lvec(424 72)
\ifill f:0
\move(427 71)
\lvec(429 71)
\lvec(429 72)
\lvec(427 72)
\ifill f:0
\move(431 71)
\lvec(432 71)
\lvec(432 72)
\lvec(431 72)
\ifill f:0
\move(434 71)
\lvec(436 71)
\lvec(436 72)
\lvec(434 72)
\ifill f:0
\move(438 71)
\lvec(439 71)
\lvec(439 72)
\lvec(438 72)
\ifill f:0
\move(441 71)
\lvec(443 71)
\lvec(443 72)
\lvec(441 72)
\ifill f:0
\move(445 71)
\lvec(446 71)
\lvec(446 72)
\lvec(445 72)
\ifill f:0
\move(448 71)
\lvec(449 71)
\lvec(449 72)
\lvec(448 72)
\ifill f:0
\move(12 72)
\lvec(14 72)
\lvec(14 73)
\lvec(12 73)
\ifill f:0
\move(169 72)
\lvec(170 72)
\lvec(170 73)
\lvec(169 73)
\ifill f:0
\move(184 72)
\lvec(185 72)
\lvec(185 73)
\lvec(184 73)
\ifill f:0
\move(187 72)
\lvec(188 72)
\lvec(188 73)
\lvec(187 73)
\ifill f:0
\move(195 72)
\lvec(196 72)
\lvec(196 73)
\lvec(195 73)
\ifill f:0
\move(202 72)
\lvec(203 72)
\lvec(203 73)
\lvec(202 73)
\ifill f:0
\move(204 72)
\lvec(205 72)
\lvec(205 73)
\lvec(204 73)
\ifill f:0
\move(218 72)
\lvec(219 72)
\lvec(219 73)
\lvec(218 73)
\ifill f:0
\move(220 72)
\lvec(221 72)
\lvec(221 73)
\lvec(220 73)
\ifill f:0
\move(225 72)
\lvec(226 72)
\lvec(226 73)
\lvec(225 73)
\ifill f:0
\move(230 72)
\lvec(231 72)
\lvec(231 73)
\lvec(230 73)
\ifill f:0
\move(235 72)
\lvec(236 72)
\lvec(236 73)
\lvec(235 73)
\ifill f:0
\move(238 72)
\lvec(239 72)
\lvec(239 73)
\lvec(238 73)
\ifill f:0
\move(241 72)
\lvec(242 72)
\lvec(242 73)
\lvec(241 73)
\ifill f:0
\move(244 72)
\lvec(245 72)
\lvec(245 73)
\lvec(244 73)
\ifill f:0
\move(247 72)
\lvec(248 72)
\lvec(248 73)
\lvec(247 73)
\ifill f:0
\move(250 72)
\lvec(251 72)
\lvec(251 73)
\lvec(250 73)
\ifill f:0
\move(254 72)
\lvec(255 72)
\lvec(255 73)
\lvec(254 73)
\ifill f:0
\move(258 72)
\lvec(259 72)
\lvec(259 73)
\lvec(258 73)
\ifill f:0
\move(262 72)
\lvec(263 72)
\lvec(263 73)
\lvec(262 73)
\ifill f:0
\move(266 72)
\lvec(267 72)
\lvec(267 73)
\lvec(266 73)
\ifill f:0
\move(270 72)
\lvec(272 72)
\lvec(272 73)
\lvec(270 73)
\ifill f:0
\move(275 72)
\lvec(277 72)
\lvec(277 73)
\lvec(275 73)
\ifill f:0
\move(281 72)
\lvec(283 72)
\lvec(283 73)
\lvec(281 73)
\ifill f:0
\move(288 72)
\lvec(290 72)
\lvec(290 73)
\lvec(288 73)
\ifill f:0
\move(296 72)
\lvec(298 72)
\lvec(298 73)
\lvec(296 73)
\ifill f:0
\move(306 72)
\lvec(311 72)
\lvec(311 73)
\lvec(306 73)
\ifill f:0
\move(347 72)
\lvec(353 72)
\lvec(353 73)
\lvec(347 73)
\ifill f:0
\move(361 72)
\lvec(365 72)
\lvec(365 73)
\lvec(361 73)
\ifill f:0
\move(371 72)
\lvec(374 72)
\lvec(374 73)
\lvec(371 73)
\ifill f:0
\move(379 72)
\lvec(381 72)
\lvec(381 73)
\lvec(379 73)
\ifill f:0
\move(386 72)
\lvec(388 72)
\lvec(388 73)
\lvec(386 73)
\ifill f:0
\move(392 72)
\lvec(394 72)
\lvec(394 73)
\lvec(392 73)
\ifill f:0
\move(398 72)
\lvec(400 72)
\lvec(400 73)
\lvec(398 73)
\ifill f:0
\move(403 72)
\lvec(405 72)
\lvec(405 73)
\lvec(403 73)
\ifill f:0
\move(408 72)
\lvec(410 72)
\lvec(410 73)
\lvec(408 73)
\ifill f:0
\move(413 72)
\lvec(415 72)
\lvec(415 73)
\lvec(413 73)
\ifill f:0
\move(417 72)
\lvec(419 72)
\lvec(419 73)
\lvec(417 73)
\ifill f:0
\move(422 72)
\lvec(423 72)
\lvec(423 73)
\lvec(422 73)
\ifill f:0
\move(426 72)
\lvec(428 72)
\lvec(428 73)
\lvec(426 73)
\ifill f:0
\move(430 72)
\lvec(432 72)
\lvec(432 73)
\lvec(430 73)
\ifill f:0
\move(434 72)
\lvec(435 72)
\lvec(435 73)
\lvec(434 73)
\ifill f:0
\move(438 72)
\lvec(439 72)
\lvec(439 73)
\lvec(438 73)
\ifill f:0
\move(441 72)
\lvec(443 72)
\lvec(443 73)
\lvec(441 73)
\ifill f:0
\move(445 72)
\lvec(446 72)
\lvec(446 73)
\lvec(445 73)
\ifill f:0
\move(448 72)
\lvec(450 72)
\lvec(450 73)
\lvec(448 73)
\ifill f:0
\move(19 73)
\lvec(20 73)
\lvec(20 74)
\lvec(19 74)
\ifill f:0
\move(162 73)
\lvec(163 73)
\lvec(163 74)
\lvec(162 74)
\ifill f:0
\move(169 73)
\lvec(170 73)
\lvec(170 74)
\lvec(169 74)
\ifill f:0
\move(189 73)
\lvec(190 73)
\lvec(190 74)
\lvec(189 74)
\ifill f:0
\move(192 73)
\lvec(193 73)
\lvec(193 74)
\lvec(192 74)
\ifill f:0
\move(200 73)
\lvec(201 73)
\lvec(201 74)
\lvec(200 74)
\ifill f:0
\move(207 73)
\lvec(208 73)
\lvec(208 74)
\lvec(207 74)
\ifill f:0
\move(225 73)
\lvec(226 73)
\lvec(226 74)
\lvec(225 74)
\ifill f:0
\move(227 73)
\lvec(228 73)
\lvec(228 74)
\lvec(227 74)
\ifill f:0
\move(232 73)
\lvec(233 73)
\lvec(233 74)
\lvec(232 74)
\ifill f:0
\move(234 73)
\lvec(235 73)
\lvec(235 74)
\lvec(234 74)
\ifill f:0
\move(237 73)
\lvec(238 73)
\lvec(238 74)
\lvec(237 74)
\ifill f:0
\move(242 73)
\lvec(243 73)
\lvec(243 74)
\lvec(242 74)
\ifill f:0
\move(245 73)
\lvec(246 73)
\lvec(246 74)
\lvec(245 74)
\ifill f:0
\move(248 73)
\lvec(249 73)
\lvec(249 74)
\lvec(248 74)
\ifill f:0
\move(251 73)
\lvec(252 73)
\lvec(252 74)
\lvec(251 74)
\ifill f:0
\move(254 73)
\lvec(255 73)
\lvec(255 74)
\lvec(254 74)
\ifill f:0
\move(261 73)
\lvec(262 73)
\lvec(262 74)
\lvec(261 74)
\ifill f:0
\move(265 73)
\lvec(266 73)
\lvec(266 74)
\lvec(265 74)
\ifill f:0
\move(269 73)
\lvec(270 73)
\lvec(270 74)
\lvec(269 74)
\ifill f:0
\move(273 73)
\lvec(274 73)
\lvec(274 74)
\lvec(273 74)
\ifill f:0
\move(277 73)
\lvec(279 73)
\lvec(279 74)
\lvec(277 74)
\ifill f:0
\move(282 73)
\lvec(284 73)
\lvec(284 74)
\lvec(282 74)
\ifill f:0
\move(288 73)
\lvec(290 73)
\lvec(290 74)
\lvec(288 74)
\ifill f:0
\move(294 73)
\lvec(296 73)
\lvec(296 74)
\lvec(294 74)
\ifill f:0
\move(302 73)
\lvec(305 73)
\lvec(305 74)
\lvec(302 74)
\ifill f:0
\move(312 73)
\lvec(316 73)
\lvec(316 74)
\lvec(312 74)
\ifill f:0
\move(332 73)
\lvec(346 73)
\lvec(346 74)
\lvec(332 74)
\ifill f:0
\move(361 73)
\lvec(366 73)
\lvec(366 74)
\lvec(361 74)
\ifill f:0
\move(373 73)
\lvec(377 73)
\lvec(377 74)
\lvec(373 74)
\ifill f:0
\move(383 73)
\lvec(386 73)
\lvec(386 74)
\lvec(383 74)
\ifill f:0
\move(390 73)
\lvec(393 73)
\lvec(393 74)
\lvec(390 74)
\ifill f:0
\move(397 73)
\lvec(400 73)
\lvec(400 74)
\lvec(397 74)
\ifill f:0
\move(403 73)
\lvec(406 73)
\lvec(406 74)
\lvec(403 74)
\ifill f:0
\move(409 73)
\lvec(411 73)
\lvec(411 74)
\lvec(409 74)
\ifill f:0
\move(414 73)
\lvec(416 73)
\lvec(416 74)
\lvec(414 74)
\ifill f:0
\move(419 73)
\lvec(421 73)
\lvec(421 74)
\lvec(419 74)
\ifill f:0
\move(424 73)
\lvec(426 73)
\lvec(426 74)
\lvec(424 74)
\ifill f:0
\move(429 73)
\lvec(430 73)
\lvec(430 74)
\lvec(429 74)
\ifill f:0
\move(433 73)
\lvec(435 73)
\lvec(435 74)
\lvec(433 74)
\ifill f:0
\move(437 73)
\lvec(439 73)
\lvec(439 74)
\lvec(437 74)
\ifill f:0
\move(441 73)
\lvec(443 73)
\lvec(443 74)
\lvec(441 74)
\ifill f:0
\move(445 73)
\lvec(447 73)
\lvec(447 74)
\lvec(445 74)
\ifill f:0
\move(449 73)
\lvec(451 73)
\lvec(451 74)
\lvec(449 74)
\ifill f:0
\move(22 74)
\lvec(23 74)
\lvec(23 75)
\lvec(22 75)
\ifill f:0
\move(153 74)
\lvec(154 74)
\lvec(154 75)
\lvec(153 75)
\ifill f:0
\move(155 74)
\lvec(157 74)
\lvec(157 75)
\lvec(155 75)
\ifill f:0
\move(169 74)
\lvec(170 74)
\lvec(170 75)
\lvec(169 75)
\ifill f:0
\move(188 74)
\lvec(189 74)
\lvec(189 75)
\lvec(188 75)
\ifill f:0
\move(203 74)
\lvec(204 74)
\lvec(204 75)
\lvec(203 75)
\ifill f:0
\move(208 74)
\lvec(209 74)
\lvec(209 75)
\lvec(208 75)
\ifill f:0
\move(215 74)
\lvec(216 74)
\lvec(216 75)
\lvec(215 75)
\ifill f:0
\move(217 74)
\lvec(218 74)
\lvec(218 75)
\lvec(217 75)
\ifill f:0
\move(219 74)
\lvec(220 74)
\lvec(220 75)
\lvec(219 75)
\ifill f:0
\move(221 74)
\lvec(222 74)
\lvec(222 75)
\lvec(221 75)
\ifill f:0
\move(223 74)
\lvec(224 74)
\lvec(224 75)
\lvec(223 75)
\ifill f:0
\move(225 74)
\lvec(226 74)
\lvec(226 75)
\lvec(225 75)
\ifill f:0
\move(227 74)
\lvec(228 74)
\lvec(228 75)
\lvec(227 75)
\ifill f:0
\move(229 74)
\lvec(230 74)
\lvec(230 75)
\lvec(229 75)
\ifill f:0
\move(236 74)
\lvec(237 74)
\lvec(237 75)
\lvec(236 75)
\ifill f:0
\move(238 74)
\lvec(239 74)
\lvec(239 75)
\lvec(238 75)
\ifill f:0
\move(243 74)
\lvec(244 74)
\lvec(244 75)
\lvec(243 75)
\ifill f:0
\move(246 74)
\lvec(247 74)
\lvec(247 75)
\lvec(246 75)
\ifill f:0
\move(251 74)
\lvec(252 74)
\lvec(252 75)
\lvec(251 75)
\ifill f:0
\move(254 74)
\lvec(255 74)
\lvec(255 75)
\lvec(254 75)
\ifill f:0
\move(257 74)
\lvec(258 74)
\lvec(258 75)
\lvec(257 75)
\ifill f:0
\move(264 74)
\lvec(265 74)
\lvec(265 75)
\lvec(264 75)
\ifill f:0
\move(267 74)
\lvec(268 74)
\lvec(268 75)
\lvec(267 75)
\ifill f:0
\move(271 74)
\lvec(272 74)
\lvec(272 75)
\lvec(271 75)
\ifill f:0
\move(275 74)
\lvec(276 74)
\lvec(276 75)
\lvec(275 75)
\ifill f:0
\move(279 74)
\lvec(280 74)
\lvec(280 75)
\lvec(279 75)
\ifill f:0
\move(283 74)
\lvec(285 74)
\lvec(285 75)
\lvec(283 75)
\ifill f:0
\move(288 74)
\lvec(290 74)
\lvec(290 75)
\lvec(288 75)
\ifill f:0
\move(294 74)
\lvec(295 74)
\lvec(295 75)
\lvec(294 75)
\ifill f:0
\move(299 74)
\lvec(302 74)
\lvec(302 75)
\lvec(299 75)
\ifill f:0
\move(307 74)
\lvec(309 74)
\lvec(309 75)
\lvec(307 75)
\ifill f:0
\move(315 74)
\lvec(319 74)
\lvec(319 75)
\lvec(315 75)
\ifill f:0
\move(328 74)
\lvec(334 74)
\lvec(334 75)
\lvec(328 75)
\ifill f:0
\move(361 74)
\lvec(368 74)
\lvec(368 75)
\lvec(361 75)
\ifill f:0
\move(377 74)
\lvec(382 74)
\lvec(382 75)
\lvec(377 75)
\ifill f:0
\move(388 74)
\lvec(391 74)
\lvec(391 75)
\lvec(388 75)
\ifill f:0
\move(396 74)
\lvec(399 74)
\lvec(399 75)
\lvec(396 75)
\ifill f:0
\move(404 74)
\lvec(407 74)
\lvec(407 75)
\lvec(404 75)
\ifill f:0
\move(410 74)
\lvec(413 74)
\lvec(413 75)
\lvec(410 75)
\ifill f:0
\move(416 74)
\lvec(419 74)
\lvec(419 75)
\lvec(416 75)
\ifill f:0
\move(422 74)
\lvec(424 74)
\lvec(424 75)
\lvec(422 75)
\ifill f:0
\move(427 74)
\lvec(429 74)
\lvec(429 75)
\lvec(427 75)
\ifill f:0
\move(432 74)
\lvec(434 74)
\lvec(434 75)
\lvec(432 75)
\ifill f:0
\move(437 74)
\lvec(439 74)
\lvec(439 75)
\lvec(437 75)
\ifill f:0
\move(441 74)
\lvec(443 74)
\lvec(443 75)
\lvec(441 75)
\ifill f:0
\move(446 74)
\lvec(447 74)
\lvec(447 75)
\lvec(446 75)
\ifill f:0
\move(450 74)
\lvec(451 74)
\lvec(451 75)
\lvec(450 75)
\ifill f:0
\move(15 75)
\lvec(16 75)
\lvec(16 76)
\lvec(15 76)
\ifill f:0
\move(169 75)
\lvec(170 75)
\lvec(170 76)
\lvec(169 76)
\ifill f:0
\move(177 75)
\lvec(178 75)
\lvec(178 76)
\lvec(177 76)
\ifill f:0
\move(187 75)
\lvec(188 75)
\lvec(188 76)
\lvec(187 76)
\ifill f:0
\move(191 75)
\lvec(192 75)
\lvec(192 76)
\lvec(191 76)
\ifill f:0
\move(201 75)
\lvec(202 75)
\lvec(202 76)
\lvec(201 76)
\ifill f:0
\move(209 75)
\lvec(210 75)
\lvec(210 76)
\lvec(209 76)
\ifill f:0
\move(214 75)
\lvec(215 75)
\lvec(215 76)
\lvec(214 76)
\ifill f:0
\move(221 75)
\lvec(222 75)
\lvec(222 76)
\lvec(221 76)
\ifill f:0
\move(223 75)
\lvec(224 75)
\lvec(224 76)
\lvec(223 76)
\ifill f:0
\move(225 75)
\lvec(226 75)
\lvec(226 76)
\lvec(225 76)
\ifill f:0
\move(227 75)
\lvec(228 75)
\lvec(228 76)
\lvec(227 76)
\ifill f:0
\move(229 75)
\lvec(230 75)
\lvec(230 76)
\lvec(229 76)
\ifill f:0
\move(231 75)
\lvec(232 75)
\lvec(232 76)
\lvec(231 76)
\ifill f:0
\move(233 75)
\lvec(234 75)
\lvec(234 76)
\lvec(233 76)
\ifill f:0
\move(235 75)
\lvec(236 75)
\lvec(236 76)
\lvec(235 76)
\ifill f:0
\move(242 75)
\lvec(243 75)
\lvec(243 76)
\lvec(242 76)
\ifill f:0
\move(244 75)
\lvec(245 75)
\lvec(245 76)
\lvec(244 76)
\ifill f:0
\move(249 75)
\lvec(250 75)
\lvec(250 76)
\lvec(249 76)
\ifill f:0
\move(252 75)
\lvec(253 75)
\lvec(253 76)
\lvec(252 76)
\ifill f:0
\move(257 75)
\lvec(258 75)
\lvec(258 76)
\lvec(257 76)
\ifill f:0
\move(260 75)
\lvec(261 75)
\lvec(261 76)
\lvec(260 76)
\ifill f:0
\move(263 75)
\lvec(264 75)
\lvec(264 76)
\lvec(263 76)
\ifill f:0
\move(266 75)
\lvec(267 75)
\lvec(267 76)
\lvec(266 76)
\ifill f:0
\move(269 75)
\lvec(270 75)
\lvec(270 76)
\lvec(269 76)
\ifill f:0
\move(273 75)
\lvec(274 75)
\lvec(274 76)
\lvec(273 76)
\ifill f:0
\move(276 75)
\lvec(277 75)
\lvec(277 76)
\lvec(276 76)
\ifill f:0
\move(280 75)
\lvec(281 75)
\lvec(281 76)
\lvec(280 76)
\ifill f:0
\move(284 75)
\lvec(285 75)
\lvec(285 76)
\lvec(284 76)
\ifill f:0
\move(288 75)
\lvec(290 75)
\lvec(290 76)
\lvec(288 76)
\ifill f:0
\move(293 75)
\lvec(294 75)
\lvec(294 76)
\lvec(293 76)
\ifill f:0
\move(298 75)
\lvec(300 75)
\lvec(300 76)
\lvec(298 76)
\ifill f:0
\move(304 75)
\lvec(305 75)
\lvec(305 76)
\lvec(304 76)
\ifill f:0
\move(310 75)
\lvec(312 75)
\lvec(312 76)
\lvec(310 76)
\ifill f:0
\move(317 75)
\lvec(320 75)
\lvec(320 76)
\lvec(317 76)
\ifill f:0
\move(327 75)
\lvec(331 75)
\lvec(331 76)
\lvec(327 76)
\ifill f:0
\move(341 75)
\lvec(356 75)
\lvec(356 76)
\lvec(341 76)
\ifill f:0
\move(357 75)
\lvec(373 75)
\lvec(373 76)
\lvec(357 76)
\ifill f:0
\move(384 75)
\lvec(389 75)
\lvec(389 76)
\lvec(384 76)
\ifill f:0
\move(396 75)
\lvec(399 75)
\lvec(399 76)
\lvec(396 76)
\ifill f:0
\move(404 75)
\lvec(408 75)
\lvec(408 76)
\lvec(404 76)
\ifill f:0
\move(412 75)
\lvec(415 75)
\lvec(415 76)
\lvec(412 76)
\ifill f:0
\move(419 75)
\lvec(421 75)
\lvec(421 76)
\lvec(419 76)
\ifill f:0
\move(425 75)
\lvec(428 75)
\lvec(428 76)
\lvec(425 76)
\ifill f:0
\move(431 75)
\lvec(433 75)
\lvec(433 76)
\lvec(431 76)
\ifill f:0
\move(436 75)
\lvec(438 75)
\lvec(438 76)
\lvec(436 76)
\ifill f:0
\move(441 75)
\lvec(443 75)
\lvec(443 76)
\lvec(441 76)
\ifill f:0
\move(446 75)
\lvec(448 75)
\lvec(448 76)
\lvec(446 76)
\ifill f:0
\move(15 76)
\lvec(16 76)
\lvec(16 77)
\lvec(15 77)
\ifill f:0
\move(19 76)
\lvec(20 76)
\lvec(20 77)
\lvec(19 77)
\ifill f:0
\move(179 76)
\lvec(180 76)
\lvec(180 77)
\lvec(179 77)
\ifill f:0
\move(185 76)
\lvec(186 76)
\lvec(186 77)
\lvec(185 77)
\ifill f:0
\move(190 76)
\lvec(191 76)
\lvec(191 77)
\lvec(190 77)
\ifill f:0
\move(194 76)
\lvec(195 76)
\lvec(195 77)
\lvec(194 77)
\ifill f:0
\move(198 76)
\lvec(199 76)
\lvec(199 77)
\lvec(198 77)
\ifill f:0
\move(213 76)
\lvec(214 76)
\lvec(214 77)
\lvec(213 77)
\ifill f:0
\move(218 76)
\lvec(219 76)
\lvec(219 77)
\lvec(218 77)
\ifill f:0
\move(223 76)
\lvec(224 76)
\lvec(224 77)
\lvec(223 77)
\ifill f:0
\move(225 76)
\lvec(226 76)
\lvec(226 77)
\lvec(225 77)
\ifill f:0
\move(243 76)
\lvec(244 76)
\lvec(244 77)
\lvec(243 77)
\ifill f:0
\move(245 76)
\lvec(246 76)
\lvec(246 77)
\lvec(245 77)
\ifill f:0
\move(252 76)
\lvec(253 76)
\lvec(253 77)
\lvec(252 77)
\ifill f:0
\move(257 76)
\lvec(258 76)
\lvec(258 77)
\lvec(257 77)
\ifill f:0
\move(260 76)
\lvec(261 76)
\lvec(261 77)
\lvec(260 77)
\ifill f:0
\move(265 76)
\lvec(266 76)
\lvec(266 77)
\lvec(265 77)
\ifill f:0
\move(268 76)
\lvec(269 76)
\lvec(269 77)
\lvec(268 77)
\ifill f:0
\move(271 76)
\lvec(272 76)
\lvec(272 77)
\lvec(271 77)
\ifill f:0
\move(278 76)
\lvec(279 76)
\lvec(279 77)
\lvec(278 77)
\ifill f:0
\move(281 76)
\lvec(282 76)
\lvec(282 77)
\lvec(281 77)
\ifill f:0
\move(285 76)
\lvec(286 76)
\lvec(286 77)
\lvec(285 77)
\ifill f:0
\move(289 76)
\lvec(290 76)
\lvec(290 77)
\lvec(289 77)
\ifill f:0
\move(293 76)
\lvec(294 76)
\lvec(294 77)
\lvec(293 77)
\ifill f:0
\move(297 76)
\lvec(298 76)
\lvec(298 77)
\lvec(297 77)
\ifill f:0
\move(302 76)
\lvec(303 76)
\lvec(303 77)
\lvec(302 77)
\ifill f:0
\move(307 76)
\lvec(308 76)
\lvec(308 77)
\lvec(307 77)
\ifill f:0
\move(312 76)
\lvec(314 76)
\lvec(314 77)
\lvec(312 77)
\ifill f:0
\move(319 76)
\lvec(321 76)
\lvec(321 77)
\lvec(319 77)
\ifill f:0
\move(326 76)
\lvec(329 76)
\lvec(329 77)
\lvec(326 77)
\ifill f:0
\move(336 76)
\lvec(340 76)
\lvec(340 77)
\lvec(336 77)
\ifill f:0
\move(351 76)
\lvec(366 76)
\lvec(366 77)
\lvec(351 77)
\ifill f:0
\move(367 76)
\lvec(382 76)
\lvec(382 77)
\lvec(367 77)
\ifill f:0
\move(394 76)
\lvec(398 76)
\lvec(398 77)
\lvec(394 77)
\ifill f:0
\move(406 76)
\lvec(409 76)
\lvec(409 77)
\lvec(406 77)
\ifill f:0
\move(415 76)
\lvec(418 76)
\lvec(418 77)
\lvec(415 77)
\ifill f:0
\move(422 76)
\lvec(425 76)
\lvec(425 77)
\lvec(422 77)
\ifill f:0
\move(429 76)
\lvec(432 76)
\lvec(432 77)
\lvec(429 77)
\ifill f:0
\move(436 76)
\lvec(438 76)
\lvec(438 77)
\lvec(436 77)
\ifill f:0
\move(441 76)
\lvec(443 76)
\lvec(443 77)
\lvec(441 77)
\ifill f:0
\move(447 76)
\lvec(449 76)
\lvec(449 77)
\lvec(447 77)
\ifill f:0
\move(12 77)
\lvec(14 77)
\lvec(14 78)
\lvec(12 78)
\ifill f:0
\move(15 77)
\lvec(16 77)
\lvec(16 78)
\lvec(15 78)
\ifill f:0
\move(194 77)
\lvec(195 77)
\lvec(195 78)
\lvec(194 78)
\ifill f:0
\move(198 77)
\lvec(199 77)
\lvec(199 78)
\lvec(198 78)
\ifill f:0
\move(202 77)
\lvec(203 77)
\lvec(203 78)
\lvec(202 78)
\ifill f:0
\move(212 77)
\lvec(213 77)
\lvec(213 78)
\lvec(212 78)
\ifill f:0
\move(220 77)
\lvec(221 77)
\lvec(221 78)
\lvec(220 78)
\ifill f:0
\move(225 77)
\lvec(226 77)
\lvec(226 78)
\lvec(225 78)
\ifill f:0
\move(232 77)
\lvec(233 77)
\lvec(233 78)
\lvec(232 78)
\ifill f:0
\move(234 77)
\lvec(235 77)
\lvec(235 78)
\lvec(234 78)
\ifill f:0
\move(248 77)
\lvec(249 77)
\lvec(249 78)
\lvec(248 78)
\ifill f:0
\move(250 77)
\lvec(251 77)
\lvec(251 78)
\lvec(250 78)
\ifill f:0
\move(257 77)
\lvec(258 77)
\lvec(258 78)
\lvec(257 78)
\ifill f:0
\move(262 77)
\lvec(263 77)
\lvec(263 78)
\lvec(262 78)
\ifill f:0
\move(267 77)
\lvec(268 77)
\lvec(268 78)
\lvec(267 78)
\ifill f:0
\move(270 77)
\lvec(271 77)
\lvec(271 78)
\lvec(270 78)
\ifill f:0
\move(273 77)
\lvec(274 77)
\lvec(274 78)
\lvec(273 78)
\ifill f:0
\move(279 77)
\lvec(280 77)
\lvec(280 78)
\lvec(279 78)
\ifill f:0
\move(282 77)
\lvec(283 77)
\lvec(283 78)
\lvec(282 78)
\ifill f:0
\move(285 77)
\lvec(286 77)
\lvec(286 78)
\lvec(285 78)
\ifill f:0
\move(289 77)
\lvec(290 77)
\lvec(290 78)
\lvec(289 78)
\ifill f:0
\move(292 77)
\lvec(293 77)
\lvec(293 78)
\lvec(292 78)
\ifill f:0
\move(296 77)
\lvec(297 77)
\lvec(297 78)
\lvec(296 78)
\ifill f:0
\move(300 77)
\lvec(301 77)
\lvec(301 78)
\lvec(300 78)
\ifill f:0
\move(305 77)
\lvec(306 77)
\lvec(306 78)
\lvec(305 78)
\ifill f:0
\move(309 77)
\lvec(310 77)
\lvec(310 78)
\lvec(309 78)
\ifill f:0
\move(314 77)
\lvec(316 77)
\lvec(316 78)
\lvec(314 78)
\ifill f:0
\move(320 77)
\lvec(322 77)
\lvec(322 78)
\lvec(320 78)
\ifill f:0
\move(326 77)
\lvec(328 77)
\lvec(328 78)
\lvec(326 78)
\ifill f:0
\move(333 77)
\lvec(336 77)
\lvec(336 78)
\lvec(333 78)
\ifill f:0
\move(342 77)
\lvec(346 77)
\lvec(346 78)
\lvec(342 78)
\ifill f:0
\move(355 77)
\lvec(362 77)
\lvec(362 78)
\lvec(355 78)
\ifill f:0
\move(390 77)
\lvec(397 77)
\lvec(397 78)
\lvec(390 78)
\ifill f:0
\move(407 77)
\lvec(411 77)
\lvec(411 78)
\lvec(407 78)
\ifill f:0
\move(418 77)
\lvec(422 77)
\lvec(422 78)
\lvec(418 78)
\ifill f:0
\move(427 77)
\lvec(430 77)
\lvec(430 78)
\lvec(427 78)
\ifill f:0
\move(434 77)
\lvec(437 77)
\lvec(437 78)
\lvec(434 78)
\ifill f:0
\move(441 77)
\lvec(444 77)
\lvec(444 78)
\lvec(441 78)
\ifill f:0
\move(448 77)
\lvec(450 77)
\lvec(450 78)
\lvec(448 78)
\ifill f:0
\move(225 78)
\lvec(226 78)
\lvec(226 79)
\lvec(225 79)
\ifill f:0
\move(230 78)
\lvec(231 78)
\lvec(231 79)
\lvec(230 79)
\ifill f:0
\move(237 78)
\lvec(238 78)
\lvec(238 79)
\lvec(237 79)
\ifill f:0
\move(257 78)
\lvec(258 78)
\lvec(258 79)
\lvec(257 79)
\ifill f:0
\move(264 78)
\lvec(265 78)
\lvec(265 79)
\lvec(264 79)
\ifill f:0
\move(269 78)
\lvec(270 78)
\lvec(270 79)
\lvec(269 79)
\ifill f:0
\move(274 78)
\lvec(275 78)
\lvec(275 79)
\lvec(274 79)
\ifill f:0
\move(277 78)
\lvec(278 78)
\lvec(278 79)
\lvec(277 79)
\ifill f:0
\move(289 78)
\lvec(290 78)
\lvec(290 79)
\lvec(289 79)
\ifill f:0
\move(292 78)
\lvec(293 78)
\lvec(293 79)
\lvec(292 79)
\ifill f:0
\move(299 78)
\lvec(300 78)
\lvec(300 79)
\lvec(299 79)
\ifill f:0
\move(303 78)
\lvec(304 78)
\lvec(304 79)
\lvec(303 79)
\ifill f:0
\move(307 78)
\lvec(308 78)
\lvec(308 79)
\lvec(307 79)
\ifill f:0
\move(311 78)
\lvec(312 78)
\lvec(312 79)
\lvec(311 79)
\ifill f:0
\move(316 78)
\lvec(317 78)
\lvec(317 79)
\lvec(316 79)
\ifill f:0
\move(321 78)
\lvec(322 78)
\lvec(322 79)
\lvec(321 79)
\ifill f:0
\move(326 78)
\lvec(328 78)
\lvec(328 79)
\lvec(326 79)
\ifill f:0
\move(332 78)
\lvec(334 78)
\lvec(334 79)
\lvec(332 79)
\ifill f:0
\move(339 78)
\lvec(341 78)
\lvec(341 79)
\lvec(339 79)
\ifill f:0
\move(347 78)
\lvec(350 78)
\lvec(350 79)
\lvec(347 79)
\ifill f:0
\move(358 78)
\lvec(362 78)
\lvec(362 79)
\lvec(358 79)
\ifill f:0
\move(378 78)
\lvec(394 78)
\lvec(394 79)
\lvec(378 79)
\ifill f:0
\move(410 78)
\lvec(415 78)
\lvec(415 79)
\lvec(410 79)
\ifill f:0
\move(423 78)
\lvec(427 78)
\lvec(427 79)
\lvec(423 79)
\ifill f:0
\move(433 78)
\lvec(436 78)
\lvec(436 79)
\lvec(433 79)
\ifill f:0
\move(441 78)
\lvec(444 78)
\lvec(444 79)
\lvec(441 79)
\ifill f:0
\move(448 78)
\lvec(451 78)
\lvec(451 79)
\lvec(448 79)
\ifill f:0
\move(12 79)
\lvec(13 79)
\lvec(13 80)
\lvec(12 80)
\ifill f:0
\move(199 79)
\lvec(200 79)
\lvec(200 80)
\lvec(199 80)
\ifill f:0
\move(204 79)
\lvec(205 79)
\lvec(205 80)
\lvec(204 80)
\ifill f:0
\move(219 79)
\lvec(220 79)
\lvec(220 80)
\lvec(219 80)
\ifill f:0
\move(225 79)
\lvec(226 79)
\lvec(226 80)
\lvec(225 80)
\ifill f:0
\move(228 79)
\lvec(229 79)
\lvec(229 80)
\lvec(228 80)
\ifill f:0
\move(233 79)
\lvec(234 79)
\lvec(234 80)
\lvec(233 80)
\ifill f:0
\move(238 79)
\lvec(239 79)
\lvec(239 80)
\lvec(238 80)
\ifill f:0
\move(247 79)
\lvec(248 79)
\lvec(248 80)
\lvec(247 80)
\ifill f:0
\move(249 79)
\lvec(250 79)
\lvec(250 80)
\lvec(249 80)
\ifill f:0
\move(251 79)
\lvec(252 79)
\lvec(252 80)
\lvec(251 80)
\ifill f:0
\move(255 79)
\lvec(256 79)
\lvec(256 80)
\lvec(255 80)
\ifill f:0
\move(257 79)
\lvec(258 79)
\lvec(258 80)
\lvec(257 80)
\ifill f:0
\move(259 79)
\lvec(260 79)
\lvec(260 80)
\lvec(259 80)
\ifill f:0
\move(261 79)
\lvec(262 79)
\lvec(262 80)
\lvec(261 80)
\ifill f:0
\move(268 79)
\lvec(269 79)
\lvec(269 80)
\lvec(268 80)
\ifill f:0
\move(275 79)
\lvec(276 79)
\lvec(276 80)
\lvec(275 80)
\ifill f:0
\move(283 79)
\lvec(284 79)
\lvec(284 80)
\lvec(283 80)
\ifill f:0
\move(286 79)
\lvec(287 79)
\lvec(287 80)
\lvec(286 80)
\ifill f:0
\move(289 79)
\lvec(290 79)
\lvec(290 80)
\lvec(289 80)
\ifill f:0
\move(292 79)
\lvec(293 79)
\lvec(293 80)
\lvec(292 80)
\ifill f:0
\move(295 79)
\lvec(296 79)
\lvec(296 80)
\lvec(295 80)
\ifill f:0
\move(298 79)
\lvec(299 79)
\lvec(299 80)
\lvec(298 80)
\ifill f:0
\move(305 79)
\lvec(306 79)
\lvec(306 80)
\lvec(305 80)
\ifill f:0
\move(309 79)
\lvec(310 79)
\lvec(310 80)
\lvec(309 80)
\ifill f:0
\move(313 79)
\lvec(314 79)
\lvec(314 80)
\lvec(313 80)
\ifill f:0
\move(317 79)
\lvec(318 79)
\lvec(318 80)
\lvec(317 80)
\ifill f:0
\move(321 79)
\lvec(322 79)
\lvec(322 80)
\lvec(321 80)
\ifill f:0
\move(326 79)
\lvec(327 79)
\lvec(327 80)
\lvec(326 80)
\ifill f:0
\move(331 79)
\lvec(332 79)
\lvec(332 80)
\lvec(331 80)
\ifill f:0
\move(337 79)
\lvec(338 79)
\lvec(338 80)
\lvec(337 80)
\ifill f:0
\move(343 79)
\lvec(345 79)
\lvec(345 80)
\lvec(343 80)
\ifill f:0
\move(350 79)
\lvec(352 79)
\lvec(352 80)
\lvec(350 80)
\ifill f:0
\move(359 79)
\lvec(362 79)
\lvec(362 80)
\lvec(359 80)
\ifill f:0
\move(371 79)
\lvec(376 79)
\lvec(376 80)
\lvec(371 80)
\ifill f:0
\move(416 79)
\lvec(422 79)
\lvec(422 80)
\lvec(416 80)
\ifill f:0
\move(431 79)
\lvec(435 79)
\lvec(435 80)
\lvec(431 80)
\ifill f:0
\move(441 79)
\lvec(444 79)
\lvec(444 80)
\lvec(441 80)
\ifill f:0
\move(450 79)
\lvec(451 79)
\lvec(451 80)
\lvec(450 80)
\ifill f:0
\move(225 80)
\lvec(226 80)
\lvec(226 81)
\lvec(225 81)
\ifill f:0
\move(231 80)
\lvec(232 80)
\lvec(232 81)
\lvec(231 81)
\ifill f:0
\move(239 80)
\lvec(240 80)
\lvec(240 81)
\lvec(239 81)
\ifill f:0
\move(244 80)
\lvec(245 80)
\lvec(245 81)
\lvec(244 81)
\ifill f:0
\move(253 80)
\lvec(254 80)
\lvec(254 81)
\lvec(253 81)
\ifill f:0
\move(255 80)
\lvec(256 80)
\lvec(256 81)
\lvec(255 81)
\ifill f:0
\move(257 80)
\lvec(258 80)
\lvec(258 81)
\lvec(257 81)
\ifill f:0
\move(261 80)
\lvec(262 80)
\lvec(262 81)
\lvec(261 81)
\ifill f:0
\move(263 80)
\lvec(264 80)
\lvec(264 81)
\lvec(263 81)
\ifill f:0
\move(265 80)
\lvec(266 80)
\lvec(266 81)
\lvec(265 81)
\ifill f:0
\move(267 80)
\lvec(268 80)
\lvec(268 81)
\lvec(267 81)
\ifill f:0
\move(274 80)
\lvec(275 80)
\lvec(275 81)
\lvec(274 81)
\ifill f:0
\move(276 80)
\lvec(277 80)
\lvec(277 81)
\lvec(276 81)
\ifill f:0
\move(281 80)
\lvec(282 80)
\lvec(282 81)
\lvec(281 81)
\ifill f:0
\move(286 80)
\lvec(287 80)
\lvec(287 81)
\lvec(286 81)
\ifill f:0
\move(289 80)
\lvec(290 80)
\lvec(290 81)
\lvec(289 81)
\ifill f:0
\move(307 80)
\lvec(308 80)
\lvec(308 81)
\lvec(307 81)
\ifill f:0
\move(310 80)
\lvec(311 80)
\lvec(311 81)
\lvec(310 81)
\ifill f:0
\move(314 80)
\lvec(315 80)
\lvec(315 81)
\lvec(314 81)
\ifill f:0
\move(318 80)
\lvec(319 80)
\lvec(319 81)
\lvec(318 81)
\ifill f:0
\move(322 80)
\lvec(323 80)
\lvec(323 81)
\lvec(322 81)
\ifill f:0
\move(326 80)
\lvec(327 80)
\lvec(327 81)
\lvec(326 81)
\ifill f:0
\move(330 80)
\lvec(331 80)
\lvec(331 81)
\lvec(330 81)
\ifill f:0
\move(335 80)
\lvec(336 80)
\lvec(336 81)
\lvec(335 81)
\ifill f:0
\move(340 80)
\lvec(342 80)
\lvec(342 81)
\lvec(340 81)
\ifill f:0
\move(346 80)
\lvec(347 80)
\lvec(347 81)
\lvec(346 81)
\ifill f:0
\move(352 80)
\lvec(354 80)
\lvec(354 81)
\lvec(352 81)
\ifill f:0
\move(359 80)
\lvec(362 80)
\lvec(362 81)
\lvec(359 81)
\ifill f:0
\move(368 80)
\lvec(371 80)
\lvec(371 81)
\lvec(368 81)
\ifill f:0
\move(380 80)
\lvec(386 80)
\lvec(386 81)
\lvec(380 81)
\ifill f:0
\move(426 80)
\lvec(432 80)
\lvec(432 81)
\lvec(426 81)
\ifill f:0
\move(441 80)
\lvec(445 80)
\lvec(445 81)
\lvec(441 81)
\ifill f:0
\move(12 81)
\lvec(13 81)
\lvec(13 82)
\lvec(12 82)
\ifill f:0
\move(19 81)
\lvec(20 81)
\lvec(20 82)
\lvec(19 82)
\ifill f:0
\move(183 81)
\lvec(188 81)
\lvec(188 82)
\lvec(183 82)
\ifill f:0
\move(201 81)
\lvec(202 81)
\lvec(202 82)
\lvec(201 82)
\ifill f:0
\move(213 81)
\lvec(214 81)
\lvec(214 82)
\lvec(213 82)
\ifill f:0
\move(225 81)
\lvec(226 81)
\lvec(226 82)
\lvec(225 82)
\ifill f:0
\move(243 81)
\lvec(244 81)
\lvec(244 82)
\lvec(243 82)
\ifill f:0
\move(248 81)
\lvec(249 81)
\lvec(249 82)
\lvec(248 82)
\ifill f:0
\move(255 81)
\lvec(256 81)
\lvec(256 82)
\lvec(255 82)
\ifill f:0
\move(275 81)
\lvec(276 81)
\lvec(276 82)
\lvec(275 82)
\ifill f:0
\move(277 81)
\lvec(278 81)
\lvec(278 82)
\lvec(277 82)
\ifill f:0
\move(279 81)
\lvec(280 81)
\lvec(280 82)
\lvec(279 82)
\ifill f:0
\move(284 81)
\lvec(285 81)
\lvec(285 82)
\lvec(284 82)
\ifill f:0
\move(289 81)
\lvec(290 81)
\lvec(290 82)
\lvec(289 82)
\ifill f:0
\move(294 81)
\lvec(295 81)
\lvec(295 82)
\lvec(294 82)
\ifill f:0
\move(302 81)
\lvec(303 81)
\lvec(303 82)
\lvec(302 82)
\ifill f:0
\move(305 81)
\lvec(306 81)
\lvec(306 82)
\lvec(305 82)
\ifill f:0
\move(308 81)
\lvec(309 81)
\lvec(309 82)
\lvec(308 82)
\ifill f:0
\move(315 81)
\lvec(316 81)
\lvec(316 82)
\lvec(315 82)
\ifill f:0
\move(318 81)
\lvec(319 81)
\lvec(319 82)
\lvec(318 82)
\ifill f:0
\move(322 81)
\lvec(323 81)
\lvec(323 82)
\lvec(322 82)
\ifill f:0
\move(325 81)
\lvec(327 81)
\lvec(327 82)
\lvec(325 82)
\ifill f:0
\move(329 81)
\lvec(331 81)
\lvec(331 82)
\lvec(329 82)
\ifill f:0
\move(334 81)
\lvec(335 81)
\lvec(335 82)
\lvec(334 82)
\ifill f:0
\move(338 81)
\lvec(339 81)
\lvec(339 82)
\lvec(338 82)
\ifill f:0
\move(343 81)
\lvec(344 81)
\lvec(344 82)
\lvec(343 82)
\ifill f:0
\move(348 81)
\lvec(349 81)
\lvec(349 82)
\lvec(348 82)
\ifill f:0
\move(353 81)
\lvec(355 81)
\lvec(355 82)
\lvec(353 82)
\ifill f:0
\move(360 81)
\lvec(362 81)
\lvec(362 82)
\lvec(360 82)
\ifill f:0
\move(367 81)
\lvec(369 81)
\lvec(369 82)
\lvec(367 82)
\ifill f:0
\move(375 81)
\lvec(378 81)
\lvec(378 82)
\lvec(375 82)
\ifill f:0
\move(386 81)
\lvec(391 81)
\lvec(391 82)
\lvec(386 82)
\ifill f:0
\move(407 81)
\lvec(425 81)
\lvec(425 82)
\lvec(407 82)
\ifill f:0
\move(441 81)
\lvec(447 81)
\lvec(447 82)
\lvec(441 82)
\ifill f:0
\move(13 82)
\lvec(14 82)
\lvec(14 83)
\lvec(13 83)
\ifill f:0
\move(216 82)
\lvec(217 82)
\lvec(217 83)
\lvec(216 83)
\ifill f:0
\move(225 82)
\lvec(226 82)
\lvec(226 83)
\lvec(225 83)
\ifill f:0
\move(229 82)
\lvec(230 82)
\lvec(230 83)
\lvec(229 83)
\ifill f:0
\move(250 82)
\lvec(251 82)
\lvec(251 83)
\lvec(250 83)
\ifill f:0
\move(255 82)
\lvec(256 82)
\lvec(256 83)
\lvec(255 83)
\ifill f:0
\move(264 82)
\lvec(265 82)
\lvec(265 83)
\lvec(264 83)
\ifill f:0
\move(280 82)
\lvec(281 82)
\lvec(281 83)
\lvec(280 83)
\ifill f:0
\move(282 82)
\lvec(283 82)
\lvec(283 83)
\lvec(282 83)
\ifill f:0
\move(289 82)
\lvec(290 82)
\lvec(290 83)
\lvec(289 83)
\ifill f:0
\move(291 82)
\lvec(292 82)
\lvec(292 83)
\lvec(291 83)
\ifill f:0
\move(296 82)
\lvec(297 82)
\lvec(297 83)
\lvec(296 83)
\ifill f:0
\move(304 82)
\lvec(305 82)
\lvec(305 83)
\lvec(304 83)
\ifill f:0
\move(307 82)
\lvec(308 82)
\lvec(308 83)
\lvec(307 83)
\ifill f:0
\move(319 82)
\lvec(320 82)
\lvec(320 83)
\lvec(319 83)
\ifill f:0
\move(322 82)
\lvec(323 82)
\lvec(323 83)
\lvec(322 83)
\ifill f:0
\move(329 82)
\lvec(330 82)
\lvec(330 83)
\lvec(329 83)
\ifill f:0
\move(333 82)
\lvec(334 82)
\lvec(334 83)
\lvec(333 83)
\ifill f:0
\move(337 82)
\lvec(338 82)
\lvec(338 83)
\lvec(337 83)
\ifill f:0
\move(341 82)
\lvec(342 82)
\lvec(342 83)
\lvec(341 83)
\ifill f:0
\move(345 82)
\lvec(346 82)
\lvec(346 83)
\lvec(345 83)
\ifill f:0
\move(350 82)
\lvec(351 82)
\lvec(351 83)
\lvec(350 83)
\ifill f:0
\move(355 82)
\lvec(356 82)
\lvec(356 83)
\lvec(355 83)
\ifill f:0
\move(360 82)
\lvec(362 82)
\lvec(362 83)
\lvec(360 83)
\ifill f:0
\move(366 82)
\lvec(368 82)
\lvec(368 83)
\lvec(366 83)
\ifill f:0
\move(373 82)
\lvec(375 82)
\lvec(375 83)
\lvec(373 83)
\ifill f:0
\move(381 82)
\lvec(383 82)
\lvec(383 83)
\lvec(381 83)
\ifill f:0
\move(390 82)
\lvec(394 82)
\lvec(394 83)
\lvec(390 83)
\ifill f:0
\move(404 82)
\lvec(411 82)
\lvec(411 83)
\lvec(404 83)
\ifill f:0
\move(441 82)
\lvec(449 82)
\lvec(449 83)
\lvec(441 83)
\ifill f:0
\move(12 83)
\lvec(13 83)
\lvec(13 84)
\lvec(12 84)
\ifill f:0
\move(15 83)
\lvec(16 83)
\lvec(16 84)
\lvec(15 84)
\ifill f:0
\move(220 83)
\lvec(221 83)
\lvec(221 84)
\lvec(220 84)
\ifill f:0
\move(225 83)
\lvec(226 83)
\lvec(226 84)
\lvec(225 84)
\ifill f:0
\move(237 83)
\lvec(238 83)
\lvec(238 84)
\lvec(237 84)
\ifill f:0
\move(252 83)
\lvec(253 83)
\lvec(253 84)
\lvec(252 84)
\ifill f:0
\move(269 83)
\lvec(270 83)
\lvec(270 84)
\lvec(269 84)
\ifill f:0
\move(289 83)
\lvec(290 83)
\lvec(290 84)
\lvec(289 84)
\ifill f:0
\move(291 83)
\lvec(292 83)
\lvec(292 84)
\lvec(291 84)
\ifill f:0
\move(298 83)
\lvec(299 83)
\lvec(299 84)
\lvec(298 84)
\ifill f:0
\move(303 83)
\lvec(304 83)
\lvec(304 84)
\lvec(303 84)
\ifill f:0
\move(308 83)
\lvec(309 83)
\lvec(309 84)
\lvec(308 84)
\ifill f:0
\move(319 83)
\lvec(320 83)
\lvec(320 84)
\lvec(319 84)
\ifill f:0
\move(322 83)
\lvec(323 83)
\lvec(323 84)
\lvec(322 84)
\ifill f:0
\move(332 83)
\lvec(333 83)
\lvec(333 84)
\lvec(332 84)
\ifill f:0
\move(339 83)
\lvec(340 83)
\lvec(340 84)
\lvec(339 84)
\ifill f:0
\move(343 83)
\lvec(344 83)
\lvec(344 84)
\lvec(343 84)
\ifill f:0
\move(347 83)
\lvec(348 83)
\lvec(348 84)
\lvec(347 84)
\ifill f:0
\move(351 83)
\lvec(352 83)
\lvec(352 84)
\lvec(351 84)
\ifill f:0
\move(356 83)
\lvec(357 83)
\lvec(357 84)
\lvec(356 84)
\ifill f:0
\move(360 83)
\lvec(362 83)
\lvec(362 84)
\lvec(360 84)
\ifill f:0
\move(366 83)
\lvec(367 83)
\lvec(367 84)
\lvec(366 84)
\ifill f:0
\move(371 83)
\lvec(373 83)
\lvec(373 84)
\lvec(371 84)
\ifill f:0
\move(377 83)
\lvec(379 83)
\lvec(379 84)
\lvec(377 84)
\ifill f:0
\move(384 83)
\lvec(387 83)
\lvec(387 84)
\lvec(384 84)
\ifill f:0
\move(393 83)
\lvec(395 83)
\lvec(395 84)
\lvec(393 84)
\ifill f:0
\move(403 83)
\lvec(407 83)
\lvec(407 84)
\lvec(403 84)
\ifill f:0
\move(420 83)
\lvec(436 83)
\lvec(436 84)
\lvec(420 84)
\ifill f:0
\move(437 83)
\lvec(451 83)
\lvec(451 84)
\lvec(437 84)
\ifill f:0
\move(15 84)
\lvec(16 84)
\lvec(16 85)
\lvec(15 85)
\ifill f:0
\move(19 84)
\lvec(20 84)
\lvec(20 85)
\lvec(19 85)
\ifill f:0
\move(219 84)
\lvec(220 84)
\lvec(220 85)
\lvec(219 85)
\ifill f:0
\move(225 84)
\lvec(226 84)
\lvec(226 85)
\lvec(225 85)
\ifill f:0
\move(230 84)
\lvec(231 84)
\lvec(231 85)
\lvec(230 85)
\ifill f:0
\move(242 84)
\lvec(243 84)
\lvec(243 85)
\lvec(242 85)
\ifill f:0
\move(260 84)
\lvec(261 84)
\lvec(261 85)
\lvec(260 85)
\ifill f:0
\move(270 84)
\lvec(271 84)
\lvec(271 85)
\lvec(270 85)
\ifill f:0
\move(279 84)
\lvec(280 84)
\lvec(280 85)
\lvec(279 85)
\ifill f:0
\move(281 84)
\lvec(282 84)
\lvec(282 85)
\lvec(281 85)
\ifill f:0
\move(283 84)
\lvec(284 84)
\lvec(284 85)
\lvec(283 85)
\ifill f:0
\move(285 84)
\lvec(286 84)
\lvec(286 85)
\lvec(285 85)
\ifill f:0
\move(287 84)
\lvec(288 84)
\lvec(288 85)
\lvec(287 85)
\ifill f:0
\move(289 84)
\lvec(290 84)
\lvec(290 85)
\lvec(289 85)
\ifill f:0
\move(291 84)
\lvec(292 84)
\lvec(292 85)
\lvec(291 85)
\ifill f:0
\move(293 84)
\lvec(294 84)
\lvec(294 85)
\lvec(293 85)
\ifill f:0
\move(302 84)
\lvec(303 84)
\lvec(303 85)
\lvec(302 85)
\ifill f:0
\move(309 84)
\lvec(310 84)
\lvec(310 85)
\lvec(309 85)
\ifill f:0
\move(314 84)
\lvec(315 84)
\lvec(315 85)
\lvec(314 85)
\ifill f:0
\move(317 84)
\lvec(318 84)
\lvec(318 85)
\lvec(317 85)
\ifill f:0
\move(325 84)
\lvec(326 84)
\lvec(326 85)
\lvec(325 85)
\ifill f:0
\move(328 84)
\lvec(329 84)
\lvec(329 85)
\lvec(328 85)
\ifill f:0
\move(331 84)
\lvec(332 84)
\lvec(332 85)
\lvec(331 85)
\ifill f:0
\move(334 84)
\lvec(335 84)
\lvec(335 85)
\lvec(334 85)
\ifill f:0
\move(338 84)
\lvec(339 84)
\lvec(339 85)
\lvec(338 85)
\ifill f:0
\move(341 84)
\lvec(342 84)
\lvec(342 85)
\lvec(341 85)
\ifill f:0
\move(345 84)
\lvec(346 84)
\lvec(346 85)
\lvec(345 85)
\ifill f:0
\move(348 84)
\lvec(349 84)
\lvec(349 85)
\lvec(348 85)
\ifill f:0
\move(352 84)
\lvec(353 84)
\lvec(353 85)
\lvec(352 85)
\ifill f:0
\move(356 84)
\lvec(357 84)
\lvec(357 85)
\lvec(356 85)
\ifill f:0
\move(360 84)
\lvec(362 84)
\lvec(362 85)
\lvec(360 85)
\ifill f:0
\move(365 84)
\lvec(366 84)
\lvec(366 85)
\lvec(365 85)
\ifill f:0
\move(370 84)
\lvec(371 84)
\lvec(371 85)
\lvec(370 85)
\ifill f:0
\move(375 84)
\lvec(377 84)
\lvec(377 85)
\lvec(375 85)
\ifill f:0
\move(381 84)
\lvec(382 84)
\lvec(382 85)
\lvec(381 85)
\ifill f:0
\move(387 84)
\lvec(389 84)
\lvec(389 85)
\lvec(387 85)
\ifill f:0
\move(394 84)
\lvec(397 84)
\lvec(397 85)
\lvec(394 85)
\ifill f:0
\move(403 84)
\lvec(405 84)
\lvec(405 85)
\lvec(403 85)
\ifill f:0
\move(413 84)
\lvec(417 84)
\lvec(417 85)
\lvec(413 85)
\ifill f:0
\move(430 84)
\lvec(446 84)
\lvec(446 85)
\lvec(430 85)
\ifill f:0
\move(447 84)
\lvec(451 84)
\lvec(451 85)
\lvec(447 85)
\ifill f:0
\move(12 85)
\lvec(13 85)
\lvec(13 86)
\lvec(12 86)
\ifill f:0
\move(15 85)
\lvec(16 85)
\lvec(16 86)
\lvec(15 86)
\ifill f:0
\move(248 85)
\lvec(249 85)
\lvec(249 86)
\lvec(248 86)
\ifill f:0
\move(266 85)
\lvec(267 85)
\lvec(267 86)
\lvec(266 86)
\ifill f:0
\move(271 85)
\lvec(272 85)
\lvec(272 86)
\lvec(271 86)
\ifill f:0
\move(276 85)
\lvec(277 85)
\lvec(277 86)
\lvec(276 86)
\ifill f:0
\move(285 85)
\lvec(286 85)
\lvec(286 86)
\lvec(285 86)
\ifill f:0
\move(287 85)
\lvec(288 85)
\lvec(288 86)
\lvec(287 86)
\ifill f:0
\move(289 85)
\lvec(290 85)
\lvec(290 86)
\lvec(289 86)
\ifill f:0
\move(291 85)
\lvec(292 85)
\lvec(292 86)
\lvec(291 86)
\ifill f:0
\move(293 85)
\lvec(294 85)
\lvec(294 86)
\lvec(293 86)
\ifill f:0
\move(295 85)
\lvec(296 85)
\lvec(296 86)
\lvec(295 86)
\ifill f:0
\move(297 85)
\lvec(298 85)
\lvec(298 86)
\lvec(297 86)
\ifill f:0
\move(299 85)
\lvec(300 85)
\lvec(300 86)
\lvec(299 86)
\ifill f:0
\move(308 85)
\lvec(309 85)
\lvec(309 86)
\lvec(308 86)
\ifill f:0
\move(310 85)
\lvec(311 85)
\lvec(311 86)
\lvec(310 86)
\ifill f:0
\move(315 85)
\lvec(316 85)
\lvec(316 86)
\lvec(315 86)
\ifill f:0
\move(320 85)
\lvec(321 85)
\lvec(321 86)
\lvec(320 86)
\ifill f:0
\move(325 85)
\lvec(326 85)
\lvec(326 86)
\lvec(325 86)
\ifill f:0
\move(328 85)
\lvec(329 85)
\lvec(329 86)
\lvec(328 86)
\ifill f:0
\move(343 85)
\lvec(344 85)
\lvec(344 86)
\lvec(343 86)
\ifill f:0
\move(346 85)
\lvec(347 85)
\lvec(347 86)
\lvec(346 86)
\ifill f:0
\move(353 85)
\lvec(354 85)
\lvec(354 86)
\lvec(353 86)
\ifill f:0
\move(357 85)
\lvec(358 85)
\lvec(358 86)
\lvec(357 86)
\ifill f:0
\move(361 85)
\lvec(362 85)
\lvec(362 86)
\lvec(361 86)
\ifill f:0
\move(365 85)
\lvec(366 85)
\lvec(366 86)
\lvec(365 86)
\ifill f:0
\move(369 85)
\lvec(370 85)
\lvec(370 86)
\lvec(369 86)
\ifill f:0
\move(374 85)
\lvec(375 85)
\lvec(375 86)
\lvec(374 86)
\ifill f:0
\move(378 85)
\lvec(380 85)
\lvec(380 86)
\lvec(378 86)
\ifill f:0
\move(384 85)
\lvec(385 85)
\lvec(385 86)
\lvec(384 86)
\ifill f:0
\move(389 85)
\lvec(391 85)
\lvec(391 86)
\lvec(389 86)
\ifill f:0
\move(395 85)
\lvec(397 85)
\lvec(397 86)
\lvec(395 86)
\ifill f:0
\move(402 85)
\lvec(405 85)
\lvec(405 86)
\lvec(402 86)
\ifill f:0
\move(411 85)
\lvec(413 85)
\lvec(413 86)
\lvec(411 86)
\ifill f:0
\move(421 85)
\lvec(424 85)
\lvec(424 86)
\lvec(421 86)
\ifill f:0
\move(435 85)
\lvec(442 85)
\lvec(442 86)
\lvec(435 86)
\ifill f:0
\move(225 86)
\lvec(226 86)
\lvec(226 87)
\lvec(225 87)
\ifill f:0
\move(247 86)
\lvec(248 86)
\lvec(248 87)
\lvec(247 87)
\ifill f:0
\move(261 86)
\lvec(262 86)
\lvec(262 87)
\lvec(261 87)
\ifill f:0
\move(264 86)
\lvec(265 86)
\lvec(265 87)
\lvec(264 87)
\ifill f:0
\move(267 86)
\lvec(268 86)
\lvec(268 87)
\lvec(267 87)
\ifill f:0
\move(275 86)
\lvec(276 86)
\lvec(276 87)
\lvec(275 87)
\ifill f:0
\move(280 86)
\lvec(281 86)
\lvec(281 87)
\lvec(280 87)
\ifill f:0
\move(287 86)
\lvec(288 86)
\lvec(288 87)
\lvec(287 87)
\ifill f:0
\move(289 86)
\lvec(290 86)
\lvec(290 87)
\lvec(289 87)
\ifill f:0
\move(311 86)
\lvec(312 86)
\lvec(312 87)
\lvec(311 87)
\ifill f:0
\move(318 86)
\lvec(319 86)
\lvec(319 87)
\lvec(318 87)
\ifill f:0
\move(320 86)
\lvec(321 86)
\lvec(321 87)
\lvec(320 87)
\ifill f:0
\move(325 86)
\lvec(326 86)
\lvec(326 87)
\lvec(325 87)
\ifill f:0
\move(330 86)
\lvec(331 86)
\lvec(331 87)
\lvec(330 87)
\ifill f:0
\move(333 86)
\lvec(334 86)
\lvec(334 87)
\lvec(333 87)
\ifill f:0
\move(341 86)
\lvec(342 86)
\lvec(342 87)
\lvec(341 87)
\ifill f:0
\move(344 86)
\lvec(345 86)
\lvec(345 87)
\lvec(344 87)
\ifill f:0
\move(347 86)
\lvec(348 86)
\lvec(348 87)
\lvec(347 87)
\ifill f:0
\move(354 86)
\lvec(355 86)
\lvec(355 87)
\lvec(354 87)
\ifill f:0
\move(357 86)
\lvec(358 86)
\lvec(358 87)
\lvec(357 87)
\ifill f:0
\move(361 86)
\lvec(362 86)
\lvec(362 87)
\lvec(361 87)
\ifill f:0
\move(364 86)
\lvec(365 86)
\lvec(365 87)
\lvec(364 87)
\ifill f:0
\move(368 86)
\lvec(369 86)
\lvec(369 87)
\lvec(368 87)
\ifill f:0
\move(372 86)
\lvec(373 86)
\lvec(373 87)
\lvec(372 87)
\ifill f:0
\move(376 86)
\lvec(378 86)
\lvec(378 87)
\lvec(376 87)
\ifill f:0
\move(381 86)
\lvec(382 86)
\lvec(382 87)
\lvec(381 87)
\ifill f:0
\move(386 86)
\lvec(387 86)
\lvec(387 87)
\lvec(386 87)
\ifill f:0
\move(391 86)
\lvec(392 86)
\lvec(392 87)
\lvec(391 87)
\ifill f:0
\move(396 86)
\lvec(398 86)
\lvec(398 87)
\lvec(396 87)
\ifill f:0
\move(402 86)
\lvec(404 86)
\lvec(404 87)
\lvec(402 87)
\ifill f:0
\move(409 86)
\lvec(411 86)
\lvec(411 87)
\lvec(409 87)
\ifill f:0
\move(416 86)
\lvec(419 86)
\lvec(419 87)
\lvec(416 87)
\ifill f:0
\move(425 86)
\lvec(429 86)
\lvec(429 87)
\lvec(425 87)
\ifill f:0
\move(437 86)
\lvec(442 86)
\lvec(442 87)
\lvec(437 87)
\ifill f:0
\move(13 87)
\lvec(14 87)
\lvec(14 88)
\lvec(13 88)
\ifill f:0
\move(19 87)
\lvec(20 87)
\lvec(20 88)
\lvec(19 88)
\ifill f:0
\move(22 87)
\lvec(23 87)
\lvec(23 88)
\lvec(22 88)
\ifill f:0
\move(234 87)
\lvec(235 87)
\lvec(235 88)
\lvec(234 88)
\ifill f:0
\move(250 87)
\lvec(251 87)
\lvec(251 88)
\lvec(250 88)
\ifill f:0
\move(254 87)
\lvec(255 87)
\lvec(255 88)
\lvec(254 88)
\ifill f:0
\move(258 87)
\lvec(259 87)
\lvec(259 88)
\lvec(258 88)
\ifill f:0
\move(268 87)
\lvec(269 87)
\lvec(269 88)
\lvec(268 88)
\ifill f:0
\move(271 87)
\lvec(272 87)
\lvec(272 88)
\lvec(271 88)
\ifill f:0
\move(289 87)
\lvec(290 87)
\lvec(290 88)
\lvec(289 88)
\ifill f:0
\move(296 87)
\lvec(297 87)
\lvec(297 88)
\lvec(296 88)
\ifill f:0
\move(298 87)
\lvec(299 87)
\lvec(299 88)
\lvec(298 88)
\ifill f:0
\move(314 87)
\lvec(315 87)
\lvec(315 88)
\lvec(314 88)
\ifill f:0
\move(316 87)
\lvec(317 87)
\lvec(317 88)
\lvec(316 88)
\ifill f:0
\move(325 87)
\lvec(326 87)
\lvec(326 88)
\lvec(325 88)
\ifill f:0
\move(330 87)
\lvec(331 87)
\lvec(331 88)
\lvec(330 88)
\ifill f:0
\move(335 87)
\lvec(336 87)
\lvec(336 88)
\lvec(335 88)
\ifill f:0
\move(340 87)
\lvec(341 87)
\lvec(341 88)
\lvec(340 88)
\ifill f:0
\move(343 87)
\lvec(344 87)
\lvec(344 88)
\lvec(343 88)
\ifill f:0
\move(361 87)
\lvec(362 87)
\lvec(362 88)
\lvec(361 88)
\ifill f:0
\move(364 87)
\lvec(365 87)
\lvec(365 88)
\lvec(364 88)
\ifill f:0
\move(371 87)
\lvec(372 87)
\lvec(372 88)
\lvec(371 88)
\ifill f:0
\move(375 87)
\lvec(376 87)
\lvec(376 88)
\lvec(375 88)
\ifill f:0
\move(379 87)
\lvec(380 87)
\lvec(380 88)
\lvec(379 88)
\ifill f:0
\move(383 87)
\lvec(384 87)
\lvec(384 88)
\lvec(383 88)
\ifill f:0
\move(387 87)
\lvec(388 87)
\lvec(388 88)
\lvec(387 88)
\ifill f:0
\move(392 87)
\lvec(393 87)
\lvec(393 88)
\lvec(392 88)
\ifill f:0
\move(397 87)
\lvec(398 87)
\lvec(398 88)
\lvec(397 88)
\ifill f:0
\move(402 87)
\lvec(403 87)
\lvec(403 88)
\lvec(402 88)
\ifill f:0
\move(408 87)
\lvec(409 87)
\lvec(409 88)
\lvec(408 88)
\ifill f:0
\move(414 87)
\lvec(416 87)
\lvec(416 88)
\lvec(414 88)
\ifill f:0
\move(421 87)
\lvec(423 87)
\lvec(423 88)
\lvec(421 88)
\ifill f:0
\move(429 87)
\lvec(431 87)
\lvec(431 88)
\lvec(429 88)
\ifill f:0
\move(438 87)
\lvec(442 87)
\lvec(442 88)
\lvec(438 88)
\ifill f:0
\move(249 88)
\lvec(250 88)
\lvec(250 89)
\lvec(249 89)
\ifill f:0
\move(258 88)
\lvec(259 88)
\lvec(259 89)
\lvec(258 89)
\ifill f:0
\move(262 88)
\lvec(263 88)
\lvec(263 89)
\lvec(262 89)
\ifill f:0
\move(269 88)
\lvec(270 88)
\lvec(270 89)
\lvec(269 89)
\ifill f:0
\move(284 88)
\lvec(285 88)
\lvec(285 89)
\lvec(284 89)
\ifill f:0
\move(289 88)
\lvec(290 88)
\lvec(290 89)
\lvec(289 89)
\ifill f:0
\move(294 88)
\lvec(295 88)
\lvec(295 89)
\lvec(294 89)
\ifill f:0
\move(301 88)
\lvec(302 88)
\lvec(302 89)
\lvec(301 89)
\ifill f:0
\move(303 88)
\lvec(304 88)
\lvec(304 89)
\lvec(303 89)
\ifill f:0
\move(325 88)
\lvec(326 88)
\lvec(326 89)
\lvec(325 89)
\ifill f:0
\move(334 88)
\lvec(335 88)
\lvec(335 89)
\lvec(334 89)
\ifill f:0
\move(339 88)
\lvec(340 88)
\lvec(340 89)
\lvec(339 89)
\ifill f:0
\move(344 88)
\lvec(345 88)
\lvec(345 89)
\lvec(344 89)
\ifill f:0
\move(352 88)
\lvec(353 88)
\lvec(353 89)
\lvec(352 89)
\ifill f:0
\move(355 88)
\lvec(356 88)
\lvec(356 89)
\lvec(355 89)
\ifill f:0
\move(358 88)
\lvec(359 88)
\lvec(359 89)
\lvec(358 89)
\ifill f:0
\move(361 88)
\lvec(362 88)
\lvec(362 89)
\lvec(361 89)
\ifill f:0
\move(364 88)
\lvec(365 88)
\lvec(365 89)
\lvec(364 89)
\ifill f:0
\move(367 88)
\lvec(368 88)
\lvec(368 89)
\lvec(367 89)
\ifill f:0
\move(370 88)
\lvec(371 88)
\lvec(371 89)
\lvec(370 89)
\ifill f:0
\move(374 88)
\lvec(375 88)
\lvec(375 89)
\lvec(374 89)
\ifill f:0
\move(377 88)
\lvec(378 88)
\lvec(378 89)
\lvec(377 89)
\ifill f:0
\move(381 88)
\lvec(382 88)
\lvec(382 89)
\lvec(381 89)
\ifill f:0
\move(385 88)
\lvec(386 88)
\lvec(386 89)
\lvec(385 89)
\ifill f:0
\move(389 88)
\lvec(390 88)
\lvec(390 89)
\lvec(389 89)
\ifill f:0
\move(393 88)
\lvec(394 88)
\lvec(394 89)
\lvec(393 89)
\ifill f:0
\move(397 88)
\lvec(398 88)
\lvec(398 89)
\lvec(397 89)
\ifill f:0
\move(402 88)
\lvec(403 88)
\lvec(403 89)
\lvec(402 89)
\ifill f:0
\move(407 88)
\lvec(408 88)
\lvec(408 89)
\lvec(407 89)
\ifill f:0
\move(412 88)
\lvec(413 88)
\lvec(413 89)
\lvec(412 89)
\ifill f:0
\move(418 88)
\lvec(419 88)
\lvec(419 89)
\lvec(418 89)
\ifill f:0
\move(424 88)
\lvec(426 88)
\lvec(426 89)
\lvec(424 89)
\ifill f:0
\move(431 88)
\lvec(433 88)
\lvec(433 89)
\lvec(431 89)
\ifill f:0
\move(439 88)
\lvec(442 88)
\lvec(442 89)
\lvec(439 89)
\ifill f:0
\move(449 88)
\lvec(451 88)
\lvec(451 89)
\lvec(449 89)
\ifill f:0
\move(263 89)
\lvec(264 89)
\lvec(264 90)
\lvec(263 90)
\ifill f:0
\move(267 89)
\lvec(268 89)
\lvec(268 90)
\lvec(267 90)
\ifill f:0
\move(274 89)
\lvec(275 89)
\lvec(275 90)
\lvec(274 90)
\ifill f:0
\move(289 89)
\lvec(290 89)
\lvec(290 90)
\lvec(289 90)
\ifill f:0
\move(297 89)
\lvec(298 89)
\lvec(298 90)
\lvec(297 90)
\ifill f:0
\move(304 89)
\lvec(305 89)
\lvec(305 90)
\lvec(304 90)
\ifill f:0
\move(313 89)
\lvec(314 89)
\lvec(314 90)
\lvec(313 90)
\ifill f:0
\move(315 89)
\lvec(316 89)
\lvec(316 90)
\lvec(315 90)
\ifill f:0
\move(317 89)
\lvec(318 89)
\lvec(318 90)
\lvec(317 90)
\ifill f:0
\move(319 89)
\lvec(320 89)
\lvec(320 90)
\lvec(319 90)
\ifill f:0
\move(321 89)
\lvec(322 89)
\lvec(322 90)
\lvec(321 90)
\ifill f:0
\move(323 89)
\lvec(324 89)
\lvec(324 90)
\lvec(323 90)
\ifill f:0
\move(325 89)
\lvec(326 89)
\lvec(326 90)
\lvec(325 90)
\ifill f:0
\move(327 89)
\lvec(328 89)
\lvec(328 90)
\lvec(327 90)
\ifill f:0
\move(329 89)
\lvec(330 89)
\lvec(330 90)
\lvec(329 90)
\ifill f:0
\move(338 89)
\lvec(339 89)
\lvec(339 90)
\lvec(338 90)
\ifill f:0
\move(345 89)
\lvec(346 89)
\lvec(346 90)
\lvec(345 90)
\ifill f:0
\move(350 89)
\lvec(351 89)
\lvec(351 90)
\lvec(350 90)
\ifill f:0
\move(358 89)
\lvec(359 89)
\lvec(359 90)
\lvec(358 90)
\ifill f:0
\move(361 89)
\lvec(362 89)
\lvec(362 90)
\lvec(361 90)
\ifill f:0
\move(376 89)
\lvec(377 89)
\lvec(377 90)
\lvec(376 90)
\ifill f:0
\move(379 89)
\lvec(380 89)
\lvec(380 90)
\lvec(379 90)
\ifill f:0
\move(386 89)
\lvec(387 89)
\lvec(387 90)
\lvec(386 90)
\ifill f:0
\move(390 89)
\lvec(391 89)
\lvec(391 90)
\lvec(390 90)
\ifill f:0
\move(394 89)
\lvec(395 89)
\lvec(395 90)
\lvec(394 90)
\ifill f:0
\move(398 89)
\lvec(399 89)
\lvec(399 90)
\lvec(398 90)
\ifill f:0
\move(402 89)
\lvec(403 89)
\lvec(403 90)
\lvec(402 90)
\ifill f:0
\move(406 89)
\lvec(407 89)
\lvec(407 90)
\lvec(406 90)
\ifill f:0
\move(411 89)
\lvec(412 89)
\lvec(412 90)
\lvec(411 90)
\ifill f:0
\move(416 89)
\lvec(417 89)
\lvec(417 90)
\lvec(416 90)
\ifill f:0
\move(421 89)
\lvec(422 89)
\lvec(422 90)
\lvec(421 90)
\ifill f:0
\move(426 89)
\lvec(428 89)
\lvec(428 90)
\lvec(426 90)
\ifill f:0
\move(433 89)
\lvec(435 89)
\lvec(435 90)
\lvec(433 90)
\ifill f:0
\move(440 89)
\lvec(442 89)
\lvec(442 90)
\lvec(440 90)
\ifill f:0
\move(447 89)
\lvec(450 89)
\lvec(450 90)
\lvec(447 90)
\ifill f:0
\move(12 90)
\lvec(13 90)
\lvec(13 91)
\lvec(12 91)
\ifill f:0
\move(19 90)
\lvec(20 90)
\lvec(20 91)
\lvec(19 91)
\ifill f:0
\move(22 90)
\lvec(23 90)
\lvec(23 91)
\lvec(22 91)
\ifill f:0
\move(253 90)
\lvec(254 90)
\lvec(254 91)
\lvec(253 91)
\ifill f:0
\move(276 90)
\lvec(277 90)
\lvec(277 91)
\lvec(276 91)
\ifill f:0
\move(283 90)
\lvec(284 90)
\lvec(284 91)
\lvec(283 91)
\ifill f:0
\move(286 90)
\lvec(287 90)
\lvec(287 91)
\lvec(286 91)
\ifill f:0
\move(289 90)
\lvec(290 90)
\lvec(290 91)
\lvec(289 91)
\ifill f:0
\move(292 90)
\lvec(293 90)
\lvec(293 91)
\lvec(292 91)
\ifill f:0
\move(300 90)
\lvec(301 90)
\lvec(301 91)
\lvec(300 91)
\ifill f:0
\move(305 90)
\lvec(306 90)
\lvec(306 91)
\lvec(305 91)
\ifill f:0
\move(310 90)
\lvec(311 90)
\lvec(311 91)
\lvec(310 91)
\ifill f:0
\move(321 90)
\lvec(322 90)
\lvec(322 91)
\lvec(321 91)
\ifill f:0
\move(323 90)
\lvec(324 90)
\lvec(324 91)
\lvec(323 91)
\ifill f:0
\move(325 90)
\lvec(326 90)
\lvec(326 91)
\lvec(325 91)
\ifill f:0
\move(327 90)
\lvec(328 90)
\lvec(328 91)
\lvec(327 91)
\ifill f:0
\move(329 90)
\lvec(330 90)
\lvec(330 91)
\lvec(329 91)
\ifill f:0
\move(331 90)
\lvec(332 90)
\lvec(332 91)
\lvec(331 91)
\ifill f:0
\move(333 90)
\lvec(334 90)
\lvec(334 91)
\lvec(333 91)
\ifill f:0
\move(335 90)
\lvec(336 90)
\lvec(336 91)
\lvec(335 91)
\ifill f:0
\move(344 90)
\lvec(345 90)
\lvec(345 91)
\lvec(344 91)
\ifill f:0
\move(346 90)
\lvec(347 90)
\lvec(347 91)
\lvec(346 91)
\ifill f:0
\move(351 90)
\lvec(352 90)
\lvec(352 91)
\lvec(351 91)
\ifill f:0
\move(358 90)
\lvec(359 90)
\lvec(359 91)
\lvec(358 91)
\ifill f:0
\move(361 90)
\lvec(362 90)
\lvec(362 91)
\lvec(361 91)
\ifill f:0
\move(366 90)
\lvec(367 90)
\lvec(367 91)
\lvec(366 91)
\ifill f:0
\move(369 90)
\lvec(370 90)
\lvec(370 91)
\lvec(369 91)
\ifill f:0
\move(372 90)
\lvec(373 90)
\lvec(373 91)
\lvec(372 91)
\ifill f:0
\move(375 90)
\lvec(376 90)
\lvec(376 91)
\lvec(375 91)
\ifill f:0
\move(378 90)
\lvec(379 90)
\lvec(379 91)
\lvec(378 91)
\ifill f:0
\move(381 90)
\lvec(382 90)
\lvec(382 91)
\lvec(381 91)
\ifill f:0
\move(384 90)
\lvec(385 90)
\lvec(385 91)
\lvec(384 91)
\ifill f:0
\move(387 90)
\lvec(388 90)
\lvec(388 91)
\lvec(387 91)
\ifill f:0
\move(391 90)
\lvec(392 90)
\lvec(392 91)
\lvec(391 91)
\ifill f:0
\move(394 90)
\lvec(395 90)
\lvec(395 91)
\lvec(394 91)
\ifill f:0
\move(398 90)
\lvec(399 90)
\lvec(399 91)
\lvec(398 91)
\ifill f:0
\move(402 90)
\lvec(403 90)
\lvec(403 91)
\lvec(402 91)
\ifill f:0
\move(406 90)
\lvec(407 90)
\lvec(407 91)
\lvec(406 91)
\ifill f:0
\move(410 90)
\lvec(411 90)
\lvec(411 91)
\lvec(410 91)
\ifill f:0
\move(414 90)
\lvec(415 90)
\lvec(415 91)
\lvec(414 91)
\ifill f:0
\move(419 90)
\lvec(420 90)
\lvec(420 91)
\lvec(419 91)
\ifill f:0
\move(423 90)
\lvec(425 90)
\lvec(425 91)
\lvec(423 91)
\ifill f:0
\move(428 90)
\lvec(430 90)
\lvec(430 91)
\lvec(428 91)
\ifill f:0
\move(434 90)
\lvec(436 90)
\lvec(436 91)
\lvec(434 91)
\ifill f:0
\move(440 90)
\lvec(442 90)
\lvec(442 91)
\lvec(440 91)
\ifill f:0
\move(447 90)
\lvec(449 90)
\lvec(449 91)
\lvec(447 91)
\ifill f:0
\move(15 91)
\lvec(16 91)
\lvec(16 92)
\lvec(15 92)
\ifill f:0
\move(252 91)
\lvec(253 91)
\lvec(253 92)
\lvec(252 92)
\ifill f:0
\move(265 91)
\lvec(266 91)
\lvec(266 92)
\lvec(265 92)
\ifill f:0
\move(270 91)
\lvec(271 91)
\lvec(271 92)
\lvec(270 92)
\ifill f:0
\move(289 91)
\lvec(290 91)
\lvec(290 92)
\lvec(289 92)
\ifill f:0
\move(301 91)
\lvec(302 91)
\lvec(302 92)
\lvec(301 92)
\ifill f:0
\move(309 91)
\lvec(310 91)
\lvec(310 92)
\lvec(309 92)
\ifill f:0
\move(323 91)
\lvec(324 91)
\lvec(324 92)
\lvec(323 92)
\ifill f:0
\move(347 91)
\lvec(348 91)
\lvec(348 92)
\lvec(347 92)
\ifill f:0
\move(349 91)
\lvec(350 91)
\lvec(350 92)
\lvec(349 92)
\ifill f:0
\move(356 91)
\lvec(357 91)
\lvec(357 92)
\lvec(356 92)
\ifill f:0
\move(361 91)
\lvec(362 91)
\lvec(362 92)
\lvec(361 92)
\ifill f:0
\move(366 91)
\lvec(367 91)
\lvec(367 92)
\lvec(366 92)
\ifill f:0
\move(371 91)
\lvec(372 91)
\lvec(372 92)
\lvec(371 92)
\ifill f:0
\move(374 91)
\lvec(375 91)
\lvec(375 92)
\lvec(374 92)
\ifill f:0
\move(382 91)
\lvec(383 91)
\lvec(383 92)
\lvec(382 92)
\ifill f:0
\move(385 91)
\lvec(386 91)
\lvec(386 92)
\lvec(385 92)
\ifill f:0
\move(388 91)
\lvec(389 91)
\lvec(389 92)
\lvec(388 92)
\ifill f:0
\move(395 91)
\lvec(396 91)
\lvec(396 92)
\lvec(395 92)
\ifill f:0
\move(398 91)
\lvec(399 91)
\lvec(399 92)
\lvec(398 92)
\ifill f:0
\move(405 91)
\lvec(406 91)
\lvec(406 92)
\lvec(405 92)
\ifill f:0
\move(409 91)
\lvec(410 91)
\lvec(410 92)
\lvec(409 92)
\ifill f:0
\move(413 91)
\lvec(414 91)
\lvec(414 92)
\lvec(413 92)
\ifill f:0
\move(417 91)
\lvec(418 91)
\lvec(418 92)
\lvec(417 92)
\ifill f:0
\move(421 91)
\lvec(422 91)
\lvec(422 92)
\lvec(421 92)
\ifill f:0
\move(425 91)
\lvec(426 91)
\lvec(426 92)
\lvec(425 92)
\ifill f:0
\move(430 91)
\lvec(431 91)
\lvec(431 92)
\lvec(430 92)
\ifill f:0
\move(435 91)
\lvec(436 91)
\lvec(436 92)
\lvec(435 92)
\ifill f:0
\move(440 91)
\lvec(442 91)
\lvec(442 92)
\lvec(440 92)
\ifill f:0
\move(446 91)
\lvec(448 91)
\lvec(448 92)
\lvec(446 92)
\ifill f:0
\move(12 92)
\lvec(13 92)
\lvec(13 93)
\lvec(12 93)
\ifill f:0
\move(15 92)
\lvec(16 92)
\lvec(16 93)
\lvec(15 93)
\ifill f:0
\move(19 92)
\lvec(20 92)
\lvec(20 93)
\lvec(19 93)
\ifill f:0
\move(260 92)
\lvec(261 92)
\lvec(261 93)
\lvec(260 93)
\ifill f:0
\move(277 92)
\lvec(278 92)
\lvec(278 93)
\lvec(277 93)
\ifill f:0
\move(289 92)
\lvec(290 92)
\lvec(290 93)
\lvec(289 93)
\ifill f:0
\move(299 92)
\lvec(300 92)
\lvec(300 93)
\lvec(299 93)
\ifill f:0
\move(302 92)
\lvec(303 92)
\lvec(303 93)
\lvec(302 93)
\ifill f:0
\move(305 92)
\lvec(306 92)
\lvec(306 93)
\lvec(305 93)
\ifill f:0
\move(313 92)
\lvec(314 92)
\lvec(314 93)
\lvec(313 93)
\ifill f:0
\move(318 92)
\lvec(319 92)
\lvec(319 93)
\lvec(318 93)
\ifill f:0
\move(323 92)
\lvec(324 92)
\lvec(324 93)
\lvec(323 93)
\ifill f:0
\move(332 92)
\lvec(333 92)
\lvec(333 93)
\lvec(332 93)
\ifill f:0
\move(334 92)
\lvec(335 92)
\lvec(335 93)
\lvec(334 93)
\ifill f:0
\move(350 92)
\lvec(351 92)
\lvec(351 93)
\lvec(350 93)
\ifill f:0
\move(352 92)
\lvec(353 92)
\lvec(353 93)
\lvec(352 93)
\ifill f:0
\move(354 92)
\lvec(355 92)
\lvec(355 93)
\lvec(354 93)
\ifill f:0
\move(361 92)
\lvec(362 92)
\lvec(362 93)
\lvec(361 93)
\ifill f:0
\move(363 92)
\lvec(364 92)
\lvec(364 93)
\lvec(363 93)
\ifill f:0
\move(368 92)
\lvec(369 92)
\lvec(369 93)
\lvec(368 93)
\ifill f:0
\move(373 92)
\lvec(374 92)
\lvec(374 93)
\lvec(373 93)
\ifill f:0
\move(378 92)
\lvec(379 92)
\lvec(379 93)
\lvec(378 93)
\ifill f:0
\move(386 92)
\lvec(387 92)
\lvec(387 93)
\lvec(386 93)
\ifill f:0
\move(389 92)
\lvec(390 92)
\lvec(390 93)
\lvec(389 93)
\ifill f:0
\move(392 92)
\lvec(393 92)
\lvec(393 93)
\lvec(392 93)
\ifill f:0
\move(395 92)
\lvec(396 92)
\lvec(396 93)
\lvec(395 93)
\ifill f:0
\move(398 92)
\lvec(399 92)
\lvec(399 93)
\lvec(398 93)
\ifill f:0
\move(401 92)
\lvec(402 92)
\lvec(402 93)
\lvec(401 93)
\ifill f:0
\move(405 92)
\lvec(406 92)
\lvec(406 93)
\lvec(405 93)
\ifill f:0
\move(408 92)
\lvec(409 92)
\lvec(409 93)
\lvec(408 93)
\ifill f:0
\move(415 92)
\lvec(416 92)
\lvec(416 93)
\lvec(415 93)
\ifill f:0
\move(419 92)
\lvec(420 92)
\lvec(420 93)
\lvec(419 93)
\ifill f:0
\move(423 92)
\lvec(424 92)
\lvec(424 93)
\lvec(423 93)
\ifill f:0
\move(427 92)
\lvec(428 92)
\lvec(428 93)
\lvec(427 93)
\ifill f:0
\move(431 92)
\lvec(432 92)
\lvec(432 93)
\lvec(431 93)
\ifill f:0
\move(436 92)
\lvec(437 92)
\lvec(437 93)
\lvec(436 93)
\ifill f:0
\move(440 92)
\lvec(442 92)
\lvec(442 93)
\lvec(440 93)
\ifill f:0
\move(445 92)
\lvec(447 92)
\lvec(447 93)
\lvec(445 93)
\ifill f:0
\move(269 93)
\lvec(270 93)
\lvec(270 94)
\lvec(269 94)
\ifill f:0
\move(289 93)
\lvec(290 93)
\lvec(290 94)
\lvec(289 94)
\ifill f:0
\move(293 93)
\lvec(294 93)
\lvec(294 94)
\lvec(293 94)
\ifill f:0
\move(300 93)
\lvec(301 93)
\lvec(301 94)
\lvec(300 94)
\ifill f:0
\move(315 93)
\lvec(316 93)
\lvec(316 94)
\lvec(315 94)
\ifill f:0
\move(330 93)
\lvec(331 93)
\lvec(331 94)
\lvec(330 94)
\ifill f:0
\move(337 93)
\lvec(338 93)
\lvec(338 94)
\lvec(337 94)
\ifill f:0
\move(339 93)
\lvec(340 93)
\lvec(340 94)
\lvec(339 94)
\ifill f:0
\move(361 93)
\lvec(362 93)
\lvec(362 94)
\lvec(361 94)
\ifill f:0
\move(363 93)
\lvec(364 93)
\lvec(364 94)
\lvec(363 94)
\ifill f:0
\move(370 93)
\lvec(371 93)
\lvec(371 94)
\lvec(370 94)
\ifill f:0
\move(372 93)
\lvec(373 93)
\lvec(373 94)
\lvec(372 94)
\ifill f:0
\move(377 93)
\lvec(378 93)
\lvec(378 94)
\lvec(377 94)
\ifill f:0
\move(382 93)
\lvec(383 93)
\lvec(383 94)
\lvec(382 94)
\ifill f:0
\move(387 93)
\lvec(388 93)
\lvec(388 94)
\lvec(387 94)
\ifill f:0
\move(390 93)
\lvec(391 93)
\lvec(391 94)
\lvec(390 94)
\ifill f:0
\move(401 93)
\lvec(402 93)
\lvec(402 94)
\lvec(401 94)
\ifill f:0
\move(414 93)
\lvec(415 93)
\lvec(415 94)
\lvec(414 94)
\ifill f:0
\move(417 93)
\lvec(418 93)
\lvec(418 94)
\lvec(417 94)
\ifill f:0
\move(421 93)
\lvec(422 93)
\lvec(422 94)
\lvec(421 94)
\ifill f:0
\move(428 93)
\lvec(429 93)
\lvec(429 94)
\lvec(428 94)
\ifill f:0
\move(432 93)
\lvec(433 93)
\lvec(433 94)
\lvec(432 94)
\ifill f:0
\move(436 93)
\lvec(437 93)
\lvec(437 94)
\lvec(436 94)
\ifill f:0
\move(441 93)
\lvec(442 93)
\lvec(442 94)
\lvec(441 94)
\ifill f:0
\move(445 93)
\lvec(446 93)
\lvec(446 94)
\lvec(445 94)
\ifill f:0
\move(450 93)
\lvec(451 93)
\lvec(451 94)
\lvec(450 94)
\ifill f:0
\move(12 94)
\lvec(13 94)
\lvec(13 95)
\lvec(12 95)
\ifill f:0
\move(289 94)
\lvec(290 94)
\lvec(290 95)
\lvec(289 95)
\ifill f:0
\move(301 94)
\lvec(302 94)
\lvec(302 95)
\lvec(301 95)
\ifill f:0
\move(308 94)
\lvec(309 94)
\lvec(309 95)
\lvec(308 95)
\ifill f:0
\move(320 94)
\lvec(321 94)
\lvec(321 95)
\lvec(320 95)
\ifill f:0
\move(328 94)
\lvec(329 94)
\lvec(329 95)
\lvec(328 95)
\ifill f:0
\move(333 94)
\lvec(334 94)
\lvec(334 95)
\lvec(333 95)
\ifill f:0
\move(340 94)
\lvec(341 94)
\lvec(341 95)
\lvec(340 95)
\ifill f:0
\move(349 94)
\lvec(350 94)
\lvec(350 95)
\lvec(349 95)
\ifill f:0
\move(351 94)
\lvec(352 94)
\lvec(352 95)
\lvec(351 95)
\ifill f:0
\move(353 94)
\lvec(354 94)
\lvec(354 95)
\lvec(353 95)
\ifill f:0
\move(355 94)
\lvec(356 94)
\lvec(356 95)
\lvec(355 95)
\ifill f:0
\move(359 94)
\lvec(360 94)
\lvec(360 95)
\lvec(359 95)
\ifill f:0
\move(361 94)
\lvec(362 94)
\lvec(362 95)
\lvec(361 95)
\ifill f:0
\move(363 94)
\lvec(364 94)
\lvec(364 95)
\lvec(363 95)
\ifill f:0
\move(365 94)
\lvec(366 94)
\lvec(366 95)
\lvec(365 95)
\ifill f:0
\move(376 94)
\lvec(377 94)
\lvec(377 95)
\lvec(376 95)
\ifill f:0
\move(383 94)
\lvec(384 94)
\lvec(384 95)
\lvec(383 95)
\ifill f:0
\move(388 94)
\lvec(389 94)
\lvec(389 95)
\lvec(388 95)
\ifill f:0
\move(393 94)
\lvec(394 94)
\lvec(394 95)
\lvec(393 95)
\ifill f:0
\move(401 94)
\lvec(402 94)
\lvec(402 95)
\lvec(401 95)
\ifill f:0
\move(404 94)
\lvec(405 94)
\lvec(405 95)
\lvec(404 95)
\ifill f:0
\move(407 94)
\lvec(408 94)
\lvec(408 95)
\lvec(407 95)
\ifill f:0
\move(410 94)
\lvec(411 94)
\lvec(411 95)
\lvec(410 95)
\ifill f:0
\move(413 94)
\lvec(414 94)
\lvec(414 95)
\lvec(413 95)
\ifill f:0
\move(416 94)
\lvec(417 94)
\lvec(417 95)
\lvec(416 95)
\ifill f:0
\move(419 94)
\lvec(420 94)
\lvec(420 95)
\lvec(419 95)
\ifill f:0
\move(426 94)
\lvec(427 94)
\lvec(427 95)
\lvec(426 95)
\ifill f:0
\move(429 94)
\lvec(430 94)
\lvec(430 95)
\lvec(429 95)
\ifill f:0
\move(433 94)
\lvec(434 94)
\lvec(434 95)
\lvec(433 95)
\ifill f:0
\move(437 94)
\lvec(438 94)
\lvec(438 95)
\lvec(437 95)
\ifill f:0
\move(441 94)
\lvec(442 94)
\lvec(442 95)
\lvec(441 95)
\ifill f:0
\move(445 94)
\lvec(446 94)
\lvec(446 95)
\lvec(445 95)
\ifill f:0
\move(449 94)
\lvec(450 94)
\lvec(450 95)
\lvec(449 95)
\ifill f:0
\move(19 95)
\lvec(20 95)
\lvec(20 96)
\lvec(19 96)
\ifill f:0
\move(267 95)
\lvec(268 95)
\lvec(268 96)
\lvec(267 96)
\ifill f:0
\move(289 95)
\lvec(290 95)
\lvec(290 96)
\lvec(289 96)
\ifill f:0
\move(294 95)
\lvec(295 95)
\lvec(295 96)
\lvec(294 96)
\ifill f:0
\move(310 95)
\lvec(311 95)
\lvec(311 96)
\lvec(310 96)
\ifill f:0
\move(331 95)
\lvec(332 95)
\lvec(332 96)
\lvec(331 96)
\ifill f:0
\move(341 95)
\lvec(342 95)
\lvec(342 96)
\lvec(341 96)
\ifill f:0
\move(346 95)
\lvec(347 95)
\lvec(347 96)
\lvec(346 96)
\ifill f:0
\move(355 95)
\lvec(356 95)
\lvec(356 96)
\lvec(355 96)
\ifill f:0
\move(357 95)
\lvec(358 95)
\lvec(358 96)
\lvec(357 96)
\ifill f:0
\move(359 95)
\lvec(360 95)
\lvec(360 96)
\lvec(359 96)
\ifill f:0
\move(361 95)
\lvec(362 95)
\lvec(362 96)
\lvec(361 96)
\ifill f:0
\move(363 95)
\lvec(364 95)
\lvec(364 96)
\lvec(363 96)
\ifill f:0
\move(367 95)
\lvec(368 95)
\lvec(368 96)
\lvec(367 96)
\ifill f:0
\move(369 95)
\lvec(370 95)
\lvec(370 96)
\lvec(369 96)
\ifill f:0
\move(371 95)
\lvec(372 95)
\lvec(372 96)
\lvec(371 96)
\ifill f:0
\move(373 95)
\lvec(374 95)
\lvec(374 96)
\lvec(373 96)
\ifill f:0
\move(375 95)
\lvec(376 95)
\lvec(376 96)
\lvec(375 96)
\ifill f:0
\move(384 95)
\lvec(385 95)
\lvec(385 96)
\lvec(384 96)
\ifill f:0
\move(391 95)
\lvec(392 95)
\lvec(392 96)
\lvec(391 96)
\ifill f:0
\move(396 95)
\lvec(397 95)
\lvec(397 96)
\lvec(396 96)
\ifill f:0
\move(401 95)
\lvec(402 95)
\lvec(402 96)
\lvec(401 96)
\ifill f:0
\move(409 95)
\lvec(410 95)
\lvec(410 96)
\lvec(409 96)
\ifill f:0
\move(412 95)
\lvec(413 95)
\lvec(413 96)
\lvec(412 96)
\ifill f:0
\move(415 95)
\lvec(416 95)
\lvec(416 96)
\lvec(415 96)
\ifill f:0
\move(418 95)
\lvec(419 95)
\lvec(419 96)
\lvec(418 96)
\ifill f:0
\move(421 95)
\lvec(422 95)
\lvec(422 96)
\lvec(421 96)
\ifill f:0
\move(424 95)
\lvec(425 95)
\lvec(425 96)
\lvec(424 96)
\ifill f:0
\move(427 95)
\lvec(428 95)
\lvec(428 96)
\lvec(427 96)
\ifill f:0
\move(430 95)
\lvec(431 95)
\lvec(431 96)
\lvec(430 96)
\ifill f:0
\move(434 95)
\lvec(435 95)
\lvec(435 96)
\lvec(434 96)
\ifill f:0
\move(437 95)
\lvec(438 95)
\lvec(438 96)
\lvec(437 96)
\ifill f:0
\move(441 95)
\lvec(442 95)
\lvec(442 96)
\lvec(441 96)
\ifill f:0
\move(444 95)
\lvec(445 95)
\lvec(445 96)
\lvec(444 96)
\ifill f:0
\move(448 95)
\lvec(449 95)
\lvec(449 96)
\lvec(448 96)
\ifill f:0
\move(12 96)
\lvec(13 96)
\lvec(13 97)
\lvec(12 97)
\ifill f:0
\move(272 96)
\lvec(273 96)
\lvec(273 97)
\lvec(272 97)
\ifill f:0
\move(282 96)
\lvec(283 96)
\lvec(283 97)
\lvec(282 97)
\ifill f:0
\move(289 96)
\lvec(290 96)
\lvec(290 97)
\lvec(289 97)
\ifill f:0
\move(337 96)
\lvec(338 96)
\lvec(338 97)
\lvec(337 97)
\ifill f:0
\move(345 96)
\lvec(346 96)
\lvec(346 97)
\lvec(345 97)
\ifill f:0
\move(352 96)
\lvec(353 96)
\lvec(353 97)
\lvec(352 97)
\ifill f:0
\move(359 96)
\lvec(360 96)
\lvec(360 97)
\lvec(359 97)
\ifill f:0
\move(361 96)
\lvec(362 96)
\lvec(362 97)
\lvec(361 97)
\ifill f:0
\move(385 96)
\lvec(386 96)
\lvec(386 97)
\lvec(385 97)
\ifill f:0
\move(387 96)
\lvec(388 96)
\lvec(388 97)
\lvec(387 97)
\ifill f:0
\move(394 96)
\lvec(395 96)
\lvec(395 97)
\lvec(394 97)
\ifill f:0
\move(401 96)
\lvec(402 96)
\lvec(402 97)
\lvec(401 97)
\ifill f:0
\move(406 96)
\lvec(407 96)
\lvec(407 97)
\lvec(406 97)
\ifill f:0
\move(411 96)
\lvec(412 96)
\lvec(412 97)
\lvec(411 97)
\ifill f:0
\move(414 96)
\lvec(415 96)
\lvec(415 97)
\lvec(414 97)
\ifill f:0
\move(422 96)
\lvec(423 96)
\lvec(423 97)
\lvec(422 97)
\ifill f:0
\move(425 96)
\lvec(426 96)
\lvec(426 97)
\lvec(425 97)
\ifill f:0
\move(428 96)
\lvec(429 96)
\lvec(429 97)
\lvec(428 97)
\ifill f:0
\move(431 96)
\lvec(432 96)
\lvec(432 97)
\lvec(431 97)
\ifill f:0
\move(434 96)
\lvec(435 96)
\lvec(435 97)
\lvec(434 97)
\ifill f:0
\move(441 96)
\lvec(442 96)
\lvec(442 97)
\lvec(441 97)
\ifill f:0
\move(444 96)
\lvec(445 96)
\lvec(445 97)
\lvec(444 97)
\ifill f:0
\move(448 96)
\lvec(449 96)
\lvec(449 97)
\lvec(448 97)
\ifill f:0
\move(329 97)
\lvec(330 97)
\lvec(330 98)
\lvec(329 98)
\ifill f:0
\move(338 97)
\lvec(339 97)
\lvec(339 98)
\lvec(338 98)
\ifill f:0
\move(341 97)
\lvec(342 97)
\lvec(342 98)
\lvec(341 98)
\ifill f:0
\move(349 97)
\lvec(350 97)
\lvec(350 98)
\lvec(349 98)
\ifill f:0
\move(354 97)
\lvec(355 97)
\lvec(355 98)
\lvec(354 98)
\ifill f:0
\move(361 97)
\lvec(362 97)
\lvec(362 98)
\lvec(361 98)
\ifill f:0
\move(370 97)
\lvec(371 97)
\lvec(371 98)
\lvec(370 98)
\ifill f:0
\move(372 97)
\lvec(373 97)
\lvec(373 98)
\lvec(372 98)
\ifill f:0
\move(390 97)
\lvec(391 97)
\lvec(391 98)
\lvec(390 98)
\ifill f:0
\move(392 97)
\lvec(393 97)
\lvec(393 98)
\lvec(392 98)
\ifill f:0
\move(401 97)
\lvec(402 97)
\lvec(402 98)
\lvec(401 98)
\ifill f:0
\move(408 97)
\lvec(409 97)
\lvec(409 98)
\lvec(408 98)
\ifill f:0
\move(413 97)
\lvec(414 97)
\lvec(414 98)
\lvec(413 98)
\ifill f:0
\move(418 97)
\lvec(419 97)
\lvec(419 98)
\lvec(418 98)
\ifill f:0
\move(426 97)
\lvec(427 97)
\lvec(427 98)
\lvec(426 98)
\ifill f:0
\move(429 97)
\lvec(430 97)
\lvec(430 98)
\lvec(429 98)
\ifill f:0
\move(432 97)
\lvec(433 97)
\lvec(433 98)
\lvec(432 98)
\ifill f:0
\move(435 97)
\lvec(436 97)
\lvec(436 98)
\lvec(435 98)
\ifill f:0
\move(438 97)
\lvec(439 97)
\lvec(439 98)
\lvec(438 98)
\ifill f:0
\move(441 97)
\lvec(442 97)
\lvec(442 98)
\lvec(441 98)
\ifill f:0
\move(444 97)
\lvec(445 97)
\lvec(445 98)
\lvec(444 98)
\ifill f:0
\move(447 97)
\lvec(448 97)
\lvec(448 98)
\lvec(447 98)
\ifill f:0
\move(12 98)
\lvec(13 98)
\lvec(13 99)
\lvec(12 99)
\ifill f:0
\move(15 98)
\lvec(16 98)
\lvec(16 99)
\lvec(15 99)
\ifill f:0
\move(19 98)
\lvec(20 98)
\lvec(20 99)
\lvec(19 99)
\ifill f:0
\move(268 98)
\lvec(272 98)
\lvec(272 99)
\lvec(268 99)
\ifill f:0
\move(289 98)
\lvec(290 98)
\lvec(290 99)
\lvec(289 99)
\ifill f:0
\move(297 98)
\lvec(298 98)
\lvec(298 99)
\lvec(297 99)
\ifill f:0
\move(309 98)
\lvec(310 98)
\lvec(310 99)
\lvec(309 99)
\ifill f:0
\move(318 98)
\lvec(319 98)
\lvec(319 99)
\lvec(318 99)
\ifill f:0
\move(322 98)
\lvec(323 98)
\lvec(323 99)
\lvec(322 99)
\ifill f:0
\move(326 98)
\lvec(327 98)
\lvec(327 99)
\lvec(326 99)
\ifill f:0
\move(333 98)
\lvec(334 98)
\lvec(334 99)
\lvec(333 99)
\ifill f:0
\move(348 98)
\lvec(349 98)
\lvec(349 99)
\lvec(348 99)
\ifill f:0
\move(351 98)
\lvec(352 98)
\lvec(352 99)
\lvec(351 99)
\ifill f:0
\move(356 98)
\lvec(357 98)
\lvec(357 99)
\lvec(356 99)
\ifill f:0
\move(361 98)
\lvec(362 98)
\lvec(362 99)
\lvec(361 99)
\ifill f:0
\move(366 98)
\lvec(367 98)
\lvec(367 99)
\lvec(366 99)
\ifill f:0
\move(375 98)
\lvec(376 98)
\lvec(376 99)
\lvec(375 99)
\ifill f:0
\move(377 98)
\lvec(378 98)
\lvec(378 99)
\lvec(377 99)
\ifill f:0
\move(401 98)
\lvec(402 98)
\lvec(402 99)
\lvec(401 99)
\ifill f:0
\move(403 98)
\lvec(404 98)
\lvec(404 99)
\lvec(403 99)
\ifill f:0
\move(410 98)
\lvec(411 98)
\lvec(411 99)
\lvec(410 99)
\ifill f:0
\move(412 98)
\lvec(413 98)
\lvec(413 99)
\lvec(412 99)
\ifill f:0
\move(417 98)
\lvec(418 98)
\lvec(418 99)
\lvec(417 99)
\ifill f:0
\move(422 98)
\lvec(423 98)
\lvec(423 99)
\lvec(422 99)
\ifill f:0
\move(427 98)
\lvec(428 98)
\lvec(428 99)
\lvec(427 99)
\ifill f:0
\move(435 98)
\lvec(436 98)
\lvec(436 99)
\lvec(435 99)
\ifill f:0
\move(438 98)
\lvec(439 98)
\lvec(439 99)
\lvec(438 99)
\ifill f:0
\move(441 98)
\lvec(442 98)
\lvec(442 99)
\lvec(441 99)
\ifill f:0
\move(444 98)
\lvec(445 98)
\lvec(445 99)
\lvec(444 99)
\ifill f:0
\move(447 98)
\lvec(448 98)
\lvec(448 99)
\lvec(447 99)
\ifill f:0
\move(450 98)
\lvec(451 98)
\lvec(451 99)
\lvec(450 99)
\ifill f:0
\move(15 99)
\lvec(16 99)
\lvec(16 100)
\lvec(15 100)
\ifill f:0
\move(289 99)
\lvec(290 99)
\lvec(290 100)
\lvec(289 100)
\ifill f:0
\move(299 99)
\lvec(300 99)
\lvec(300 100)
\lvec(299 100)
\ifill f:0
\move(312 99)
\lvec(313 99)
\lvec(313 100)
\lvec(312 100)
\ifill f:0
\move(322 99)
\lvec(323 99)
\lvec(323 100)
\lvec(322 100)
\ifill f:0
\move(326 99)
\lvec(327 99)
\lvec(327 100)
\lvec(326 100)
\ifill f:0
\move(330 99)
\lvec(331 99)
\lvec(331 100)
\lvec(330 100)
\ifill f:0
\move(344 99)
\lvec(345 99)
\lvec(345 100)
\lvec(344 100)
\ifill f:0
\move(347 99)
\lvec(348 99)
\lvec(348 100)
\lvec(347 100)
\ifill f:0
\move(350 99)
\lvec(351 99)
\lvec(351 100)
\lvec(350 100)
\ifill f:0
\move(353 99)
\lvec(354 99)
\lvec(354 100)
\lvec(353 100)
\ifill f:0
\move(361 99)
\lvec(362 99)
\lvec(362 100)
\lvec(361 100)
\ifill f:0
\move(371 99)
\lvec(372 99)
\lvec(372 100)
\lvec(371 100)
\ifill f:0
\move(378 99)
\lvec(379 99)
\lvec(379 100)
\lvec(378 100)
\ifill f:0
\move(387 99)
\lvec(388 99)
\lvec(388 100)
\lvec(387 100)
\ifill f:0
\move(389 99)
\lvec(390 99)
\lvec(390 100)
\lvec(389 100)
\ifill f:0
\move(391 99)
\lvec(392 99)
\lvec(392 100)
\lvec(391 100)
\ifill f:0
\move(393 99)
\lvec(394 99)
\lvec(394 100)
\lvec(393 100)
\ifill f:0
\move(395 99)
\lvec(396 99)
\lvec(396 100)
\lvec(395 100)
\ifill f:0
\move(397 99)
\lvec(398 99)
\lvec(398 100)
\lvec(397 100)
\ifill f:0
\move(399 99)
\lvec(400 99)
\lvec(400 100)
\lvec(399 100)
\ifill f:0
\move(401 99)
\lvec(402 99)
\lvec(402 100)
\lvec(401 100)
\ifill f:0
\move(403 99)
\lvec(404 99)
\lvec(404 100)
\lvec(403 100)
\ifill f:0
\move(405 99)
\lvec(406 99)
\lvec(406 100)
\lvec(405 100)
\ifill f:0
\move(416 99)
\lvec(417 99)
\lvec(417 100)
\lvec(416 100)
\ifill f:0
\move(423 99)
\lvec(424 99)
\lvec(424 100)
\lvec(423 100)
\ifill f:0
\move(428 99)
\lvec(429 99)
\lvec(429 100)
\lvec(428 100)
\ifill f:0
\move(433 99)
\lvec(434 99)
\lvec(434 100)
\lvec(433 100)
\ifill f:0
\move(438 99)
\lvec(439 99)
\lvec(439 100)
\lvec(438 100)
\ifill f:0
\move(441 99)
\lvec(442 99)
\lvec(442 100)
\lvec(441 100)
\ifill f:0
\move(446 99)
\lvec(447 99)
\lvec(447 100)
\lvec(446 100)
\ifill f:0
\move(449 99)
\lvec(450 99)
\lvec(450 100)
\lvec(449 100)
\ifill f:0
\move(15 100)
\lvec(16 100)
\lvec(16 101)
\lvec(15 101)
\ifill f:0
\move(22 100)
\lvec(23 100)
\lvec(23 101)
\lvec(22 101)
\ifill f:0
\move(302 100)
\lvec(303 100)
\lvec(303 101)
\lvec(302 101)
\ifill f:0
\move(310 100)
\lvec(311 100)
\lvec(311 101)
\lvec(310 101)
\ifill f:0
\move(335 100)
\lvec(336 100)
\lvec(336 101)
\lvec(335 101)
\ifill f:0
\move(349 100)
\lvec(350 100)
\lvec(350 101)
\lvec(349 101)
\ifill f:0
\move(361 100)
\lvec(362 100)
\lvec(362 101)
\lvec(361 101)
\ifill f:0
\move(364 100)
\lvec(365 100)
\lvec(365 101)
\lvec(364 101)
\ifill f:0
\move(384 100)
\lvec(385 100)
\lvec(385 101)
\lvec(384 101)
\ifill f:0
\move(386 100)
\lvec(387 100)
\lvec(387 101)
\lvec(386 101)
\ifill f:0
\move(395 100)
\lvec(396 100)
\lvec(396 101)
\lvec(395 101)
\ifill f:0
\move(397 100)
\lvec(398 100)
\lvec(398 101)
\lvec(397 101)
\ifill f:0
\move(399 100)
\lvec(400 100)
\lvec(400 101)
\lvec(399 101)
\ifill f:0
\move(401 100)
\lvec(402 100)
\lvec(402 101)
\lvec(401 101)
\ifill f:0
\move(403 100)
\lvec(404 100)
\lvec(404 101)
\lvec(403 101)
\ifill f:0
\move(405 100)
\lvec(406 100)
\lvec(406 101)
\lvec(405 101)
\ifill f:0
\move(407 100)
\lvec(408 100)
\lvec(408 101)
\lvec(407 101)
\ifill f:0
\move(409 100)
\lvec(410 100)
\lvec(410 101)
\lvec(409 101)
\ifill f:0
\move(411 100)
\lvec(412 100)
\lvec(412 101)
\lvec(411 101)
\ifill f:0
\move(413 100)
\lvec(414 100)
\lvec(414 101)
\lvec(413 101)
\ifill f:0
\move(424 100)
\lvec(425 100)
\lvec(425 101)
\lvec(424 101)
\ifill f:0
\move(431 100)
\lvec(432 100)
\lvec(432 101)
\lvec(431 101)
\ifill f:0
\move(436 100)
\lvec(437 100)
\lvec(437 101)
\lvec(436 101)
\ifill f:0
\move(441 100)
\lvec(442 100)
\lvec(442 101)
\lvec(441 101)
\ifill f:0
\move(446 100)
\lvec(447 100)
\lvec(447 101)
\lvec(446 101)
\ifill f:0
\move(19 101)
\lvec(20 101)
\lvec(20 102)
\lvec(19 102)
\ifill f:0
\move(348 101)
\lvec(349 101)
\lvec(349 102)
\lvec(348 102)
\ifill f:0
\move(358 101)
\lvec(359 101)
\lvec(359 102)
\lvec(358 102)
\ifill f:0
\move(361 101)
\lvec(362 101)
\lvec(362 102)
\lvec(361 102)
\ifill f:0
\move(364 101)
\lvec(365 101)
\lvec(365 102)
\lvec(364 102)
\ifill f:0
\move(367 101)
\lvec(368 101)
\lvec(368 102)
\lvec(367 102)
\ifill f:0
\move(390 101)
\lvec(391 101)
\lvec(391 102)
\lvec(390 102)
\ifill f:0
\move(399 101)
\lvec(400 101)
\lvec(400 102)
\lvec(399 102)
\ifill f:0
\move(425 101)
\lvec(426 101)
\lvec(426 102)
\lvec(425 102)
\ifill f:0
\move(427 101)
\lvec(428 101)
\lvec(428 102)
\lvec(427 102)
\ifill f:0
\move(434 101)
\lvec(435 101)
\lvec(435 102)
\lvec(434 102)
\ifill f:0
\move(436 101)
\lvec(437 101)
\lvec(437 102)
\lvec(436 102)
\ifill f:0
\move(441 101)
\lvec(442 101)
\lvec(442 102)
\lvec(441 102)
\ifill f:0
\move(448 101)
\lvec(449 101)
\lvec(449 102)
\lvec(448 102)
\ifill f:0
\move(361 102)
\lvec(362 102)
\lvec(362 103)
\lvec(361 103)
\ifill f:0
\move(379 102)
\lvec(380 102)
\lvec(380 103)
\lvec(379 103)
\ifill f:0
\move(399 102)
\lvec(400 102)
\lvec(400 103)
\lvec(399 103)
\ifill f:0
\move(408 102)
\lvec(409 102)
\lvec(409 103)
\lvec(408 103)
\ifill f:0
\move(410 102)
\lvec(411 102)
\lvec(411 103)
\lvec(410 103)
\ifill f:0
\move(412 102)
\lvec(413 102)
\lvec(413 103)
\lvec(412 103)
\ifill f:0
\move(430 102)
\lvec(431 102)
\lvec(431 103)
\lvec(430 103)
\ifill f:0
\move(432 102)
\lvec(433 102)
\lvec(433 103)
\lvec(432 103)
\ifill f:0
\move(441 102)
\lvec(442 102)
\lvec(442 103)
\lvec(441 103)
\ifill f:0
\move(443 102)
\lvec(444 102)
\lvec(444 103)
\lvec(443 103)
\ifill f:0
\move(448 102)
\lvec(449 102)
\lvec(449 103)
\lvec(448 103)
\ifill f:0
\move(22 103)
\lvec(23 103)
\lvec(23 104)
\lvec(22 104)
\ifill f:0
\move(306 103)
\lvec(307 103)
\lvec(307 104)
\lvec(306 104)
\ifill f:0
\move(328 103)
\lvec(329 103)
\lvec(329 104)
\lvec(328 104)
\ifill f:0
\move(361 103)
\lvec(362 103)
\lvec(362 104)
\lvec(361 104)
\ifill f:0
\move(368 103)
\lvec(369 103)
\lvec(369 104)
\lvec(368 104)
\ifill f:0
\move(386 103)
\lvec(387 103)
\lvec(387 104)
\lvec(386 104)
\ifill f:0
\move(394 103)
\lvec(395 103)
\lvec(395 104)
\lvec(394 104)
\ifill f:0
\move(399 103)
\lvec(400 103)
\lvec(400 104)
\lvec(399 104)
\ifill f:0
\move(406 103)
\lvec(407 103)
\lvec(407 104)
\lvec(406 104)
\ifill f:0
\move(413 103)
\lvec(414 103)
\lvec(414 104)
\lvec(413 104)
\ifill f:0
\move(415 103)
\lvec(416 103)
\lvec(416 104)
\lvec(415 104)
\ifill f:0
\move(441 103)
\lvec(442 103)
\lvec(442 104)
\lvec(441 104)
\ifill f:0
\move(443 103)
\lvec(444 103)
\lvec(444 104)
\lvec(443 104)
\ifill f:0
\move(336 104)
\lvec(337 104)
\lvec(337 105)
\lvec(336 105)
\ifill f:0
\move(357 104)
\lvec(358 104)
\lvec(358 105)
\lvec(357 105)
\ifill f:0
\move(361 104)
\lvec(362 104)
\lvec(362 105)
\lvec(361 105)
\ifill f:0
\move(365 104)
\lvec(366 104)
\lvec(366 105)
\lvec(365 105)
\ifill f:0
\move(382 104)
\lvec(383 104)
\lvec(383 105)
\lvec(382 105)
\ifill f:0
\move(385 104)
\lvec(386 104)
\lvec(386 105)
\lvec(385 105)
\ifill f:0
\move(404 104)
\lvec(405 104)
\lvec(405 105)
\lvec(404 105)
\ifill f:0
\move(416 104)
\lvec(417 104)
\lvec(417 105)
\lvec(416 105)
\ifill f:0
\move(418 104)
\lvec(419 104)
\lvec(419 105)
\lvec(418 105)
\ifill f:0
\move(427 104)
\lvec(428 104)
\lvec(428 105)
\lvec(427 105)
\ifill f:0
\move(429 104)
\lvec(430 104)
\lvec(430 105)
\lvec(429 105)
\ifill f:0
\move(431 104)
\lvec(432 104)
\lvec(432 105)
\lvec(431 105)
\ifill f:0
\move(433 104)
\lvec(434 104)
\lvec(434 105)
\lvec(433 105)
\ifill f:0
\move(435 104)
\lvec(436 104)
\lvec(436 105)
\lvec(435 105)
\ifill f:0
\move(439 104)
\lvec(440 104)
\lvec(440 105)
\lvec(439 105)
\ifill f:0
\move(441 104)
\lvec(442 104)
\lvec(442 105)
\lvec(441 105)
\ifill f:0
\move(443 104)
\lvec(444 104)
\lvec(444 105)
\lvec(443 105)
\ifill f:0
\move(445 104)
\lvec(446 104)
\lvec(446 105)
\lvec(445 105)
\ifill f:0
\move(447 104)
\lvec(448 104)
\lvec(448 105)
\lvec(447 105)
\ifill f:0
\move(12 105)
\lvec(13 105)
\lvec(13 106)
\lvec(12 106)
\ifill f:0
\move(361 105)
\lvec(362 105)
\lvec(362 106)
\lvec(361 106)
\ifill f:0
\move(390 105)
\lvec(391 105)
\lvec(391 106)
\lvec(390 106)
\ifill f:0
\move(393 105)
\lvec(394 105)
\lvec(394 106)
\lvec(393 106)
\ifill f:0
\move(396 105)
\lvec(397 105)
\lvec(397 106)
\lvec(396 106)
\ifill f:0
\move(412 105)
\lvec(413 105)
\lvec(413 106)
\lvec(412 106)
\ifill f:0
\move(426 105)
\lvec(427 105)
\lvec(427 106)
\lvec(426 106)
\ifill f:0
\move(437 105)
\lvec(438 105)
\lvec(438 106)
\lvec(437 106)
\ifill f:0
\move(439 105)
\lvec(440 105)
\lvec(440 106)
\lvec(439 106)
\ifill f:0
\move(441 105)
\lvec(442 105)
\lvec(442 106)
\lvec(441 106)
\ifill f:0
\move(443 105)
\lvec(444 105)
\lvec(444 106)
\lvec(443 106)
\ifill f:0
\move(447 105)
\lvec(448 105)
\lvec(448 106)
\lvec(447 106)
\ifill f:0
\move(449 105)
\lvec(450 105)
\lvec(450 106)
\lvec(449 106)
\ifill f:0
\move(15 106)
\lvec(16 106)
\lvec(16 107)
\lvec(15 107)
\ifill f:0
\move(19 106)
\lvec(20 106)
\lvec(20 107)
\lvec(19 107)
\ifill f:0
\move(361 106)
\lvec(362 106)
\lvec(362 107)
\lvec(361 107)
\ifill f:0
\move(366 106)
\lvec(367 106)
\lvec(367 107)
\lvec(366 107)
\ifill f:0
\move(389 106)
\lvec(390 106)
\lvec(390 107)
\lvec(389 107)
\ifill f:0
\move(410 106)
\lvec(411 106)
\lvec(411 107)
\lvec(410 107)
\ifill f:0
\move(430 106)
\lvec(431 106)
\lvec(431 107)
\lvec(430 107)
\ifill f:0
\move(439 106)
\lvec(440 106)
\lvec(440 107)
\lvec(439 107)
\ifill f:0
\move(441 106)
\lvec(442 106)
\lvec(442 107)
\lvec(441 107)
\ifill f:0
\move(12 107)
\lvec(13 107)
\lvec(13 108)
\lvec(12 108)
\ifill f:0
\move(15 107)
\lvec(16 107)
\lvec(16 108)
\lvec(15 108)
\ifill f:0
\move(336 107)
\lvec(337 107)
\lvec(337 108)
\lvec(336 108)
\ifill f:0
\move(347 107)
\lvec(348 107)
\lvec(348 108)
\lvec(347 108)
\ifill f:0
\move(361 107)
\lvec(362 107)
\lvec(362 108)
\lvec(361 108)
\ifill f:0
\move(376 107)
\lvec(377 107)
\lvec(377 108)
\lvec(376 108)
\ifill f:0
\move(388 107)
\lvec(389 107)
\lvec(389 108)
\lvec(388 108)
\ifill f:0
\move(395 107)
\lvec(396 107)
\lvec(396 108)
\lvec(395 108)
\ifill f:0
\move(419 107)
\lvec(420 107)
\lvec(420 108)
\lvec(419 108)
\ifill f:0
\move(434 107)
\lvec(435 107)
\lvec(435 108)
\lvec(434 108)
\ifill f:0
\move(441 107)
\lvec(442 107)
\lvec(442 108)
\lvec(441 108)
\ifill f:0
\move(450 107)
\lvec(451 107)
\lvec(451 108)
\lvec(450 108)
\ifill f:0
\move(15 108)
\lvec(16 108)
\lvec(16 109)
\lvec(15 109)
\ifill f:0
\move(342 108)
\lvec(343 108)
\lvec(343 109)
\lvec(342 109)
\ifill f:0
\move(361 108)
\lvec(362 108)
\lvec(362 109)
\lvec(361 109)
\ifill f:0
\move(373 108)
\lvec(374 108)
\lvec(374 109)
\lvec(373 109)
\ifill f:0
\move(378 108)
\lvec(379 108)
\lvec(379 109)
\lvec(378 109)
\ifill f:0
\move(391 108)
\lvec(392 108)
\lvec(392 109)
\lvec(391 109)
\ifill f:0
\move(405 108)
\lvec(406 108)
\lvec(406 109)
\lvec(405 109)
\ifill f:0
\move(423 108)
\lvec(424 108)
\lvec(424 109)
\lvec(423 109)
\ifill f:0
\move(426 108)
\lvec(427 108)
\lvec(427 109)
\lvec(426 109)
\ifill f:0
\move(431 108)
\lvec(432 108)
\lvec(432 109)
\lvec(431 109)
\ifill f:0
\move(441 108)
\lvec(442 108)
\lvec(442 109)
\lvec(441 109)
\ifill f:0
\move(448 108)
\lvec(449 108)
\lvec(449 109)
\lvec(448 109)
\ifill f:0
\move(12 109)
\lvec(13 109)
\lvec(13 110)
\lvec(12 110)
\ifill f:0
\move(19 109)
\lvec(20 109)
\lvec(20 110)
\lvec(19 110)
\ifill f:0
\move(361 109)
\lvec(362 109)
\lvec(362 110)
\lvec(361 110)
\ifill f:0
\move(394 109)
\lvec(395 109)
\lvec(395 110)
\lvec(394 110)
\ifill f:0
\move(398 109)
\lvec(399 109)
\lvec(399 110)
\lvec(398 110)
\ifill f:0
\move(402 109)
\lvec(403 109)
\lvec(403 110)
\lvec(402 110)
\ifill f:0
\move(409 109)
\lvec(410 109)
\lvec(410 110)
\lvec(409 110)
\ifill f:0
\move(436 109)
\lvec(437 109)
\lvec(437 110)
\lvec(436 110)
\ifill f:0
\move(441 109)
\lvec(442 109)
\lvec(442 110)
\lvec(441 110)
\ifill f:0
\move(446 109)
\lvec(447 109)
\lvec(447 110)
\lvec(446 110)
\ifill f:0
\move(402 110)
\lvec(403 110)
\lvec(403 111)
\lvec(402 111)
\ifill f:0
\move(406 110)
\lvec(407 110)
\lvec(407 111)
\lvec(406 111)
\ifill f:0
\move(417 110)
\lvec(418 110)
\lvec(418 111)
\lvec(417 111)
\ifill f:0
\move(441 110)
\lvec(442 110)
\lvec(442 111)
\lvec(441 111)
\ifill f:0
\move(449 110)
\lvec(450 110)
\lvec(450 111)
\lvec(449 111)
\ifill f:0
\move(422 111)
\lvec(423 111)
\lvec(423 112)
\lvec(422 112)
\ifill f:0
\move(432 111)
\lvec(433 111)
\lvec(433 112)
\lvec(432 112)
\ifill f:0
\move(435 111)
\lvec(436 111)
\lvec(436 112)
\lvec(435 112)
\ifill f:0
\move(438 111)
\lvec(439 111)
\lvec(439 112)
\lvec(438 112)
\ifill f:0
\move(441 111)
\lvec(442 111)
\lvec(442 112)
\lvec(441 112)
\ifill f:0
\move(444 111)
\lvec(445 111)
\lvec(445 112)
\lvec(444 112)
\ifill f:0
\move(19 112)
\lvec(20 112)
\lvec(20 113)
\lvec(19 113)
\ifill f:0
\move(412 112)
\lvec(413 112)
\lvec(413 113)
\lvec(412 113)
\ifill f:0
\move(424 112)
\lvec(425 112)
\lvec(425 113)
\lvec(424 113)
\ifill f:0
\move(438 112)
\lvec(439 112)
\lvec(439 113)
\lvec(438 113)
\ifill f:0
\move(441 112)
\lvec(442 112)
\lvec(442 113)
\lvec(441 113)
\ifill f:0
\move(444 112)
\lvec(445 112)
\lvec(445 113)
\lvec(444 113)
\ifill f:0
\move(447 112)
\lvec(448 112)
\lvec(448 113)
\lvec(447 113)
\ifill f:0
\move(450 112)
\lvec(451 112)
\lvec(451 113)
\lvec(450 113)
\ifill f:0
\move(22 113)
\lvec(23 113)
\lvec(23 114)
\lvec(22 114)
\ifill f:0
\move(403 113)
\lvec(404 113)
\lvec(404 114)
\lvec(403 114)
\ifill f:0
\move(418 113)
\lvec(419 113)
\lvec(419 114)
\lvec(418 114)
\ifill f:0
\move(434 113)
\lvec(435 113)
\lvec(435 114)
\lvec(434 114)
\ifill f:0
\move(441 113)
\lvec(442 113)
\lvec(442 114)
\lvec(441 114)
\ifill f:0
\move(15 114)
\lvec(16 114)
\lvec(16 115)
\lvec(15 115)
\ifill f:0
\move(365 114)
\lvec(366 114)
\lvec(366 115)
\lvec(365 115)
\ifill f:0
\move(396 114)
\lvec(397 114)
\lvec(397 115)
\lvec(396 115)
\ifill f:0
\move(425 114)
\lvec(426 114)
\lvec(426 115)
\lvec(425 115)
\ifill f:0
\move(441 114)
\lvec(442 114)
\lvec(442 115)
\lvec(441 115)
\ifill f:0
\move(448 114)
\lvec(449 114)
\lvec(449 115)
\lvec(448 115)
\ifill f:0
\move(15 115)
\lvec(16 115)
\lvec(16 116)
\lvec(15 116)
\ifill f:0
\move(19 115)
\lvec(20 115)
\lvec(20 116)
\lvec(19 116)
\ifill f:0
\move(395 115)
\lvec(396 115)
\lvec(396 116)
\lvec(395 116)
\ifill f:0
\move(404 115)
\lvec(405 115)
\lvec(405 116)
\lvec(404 116)
\ifill f:0
\move(423 115)
\lvec(424 115)
\lvec(424 116)
\lvec(423 116)
\ifill f:0
\move(428 115)
\lvec(429 115)
\lvec(429 116)
\lvec(428 116)
\ifill f:0
\move(437 115)
\lvec(438 115)
\lvec(438 116)
\lvec(437 116)
\ifill f:0
\move(441 115)
\lvec(442 115)
\lvec(442 116)
\lvec(441 116)
\ifill f:0
\move(445 115)
\lvec(446 115)
\lvec(446 116)
\lvec(445 116)
\ifill f:0
\move(15 116)
\lvec(16 116)
\lvec(16 117)
\lvec(15 117)
\ifill f:0
\move(22 116)
\lvec(23 116)
\lvec(23 117)
\lvec(22 117)
\ifill f:0
\move(405 116)
\lvec(406 116)
\lvec(406 117)
\lvec(405 117)
\ifill f:0
\move(441 116)
\lvec(442 116)
\lvec(442 117)
\lvec(441 117)
\ifill f:0
\move(436 117)
\lvec(437 117)
\lvec(437 118)
\lvec(436 118)
\ifill f:0
\move(441 117)
\lvec(442 117)
\lvec(442 118)
\lvec(441 118)
\ifill f:0
\move(12 118)
\lvec(13 118)
\lvec(13 119)
\lvec(12 119)
\ifill f:0
\move(409 118)
\lvec(410 118)
\lvec(410 119)
\lvec(409 119)
\ifill f:0
\move(435 118)
\lvec(436 118)
\lvec(436 119)
\lvec(435 119)
\ifill f:0
\move(441 118)
\lvec(442 118)
\lvec(442 119)
\lvec(441 119)
\ifill f:0
\move(434 119)
\lvec(435 119)
\lvec(435 120)
\lvec(434 120)
\ifill f:0
\move(441 119)
\lvec(442 119)
\lvec(442 120)
\lvec(441 120)
\ifill f:0
\move(447 119)
\lvec(448 119)
\lvec(448 120)
\lvec(447 120)
\ifill f:0
\move(12 120)
\lvec(13 120)
\lvec(13 121)
\lvec(12 121)
\ifill f:0
\move(19 120)
\lvec(20 120)
\lvec(20 121)
\lvec(19 121)
\ifill f:0
\move(441 120)
\lvec(442 120)
\lvec(442 121)
\lvec(441 121)
\ifill f:0
\move(448 120)
\lvec(449 120)
\lvec(449 121)
\lvec(448 121)
\ifill f:0
\move(12 122)
\lvec(13 122)
\lvec(13 123)
\lvec(12 123)
\ifill f:0
\move(15 122)
\lvec(16 122)
\lvec(16 123)
\lvec(15 123)
\ifill f:0
\move(15 123)
\lvec(16 123)
\lvec(16 124)
\lvec(15 124)
\ifill f:0
\move(19 123)
\lvec(20 123)
\lvec(20 124)
\lvec(19 124)
\ifill f:0
\move(437 125)
\lvec(438 125)
\lvec(438 126)
\lvec(437 126)
\ifill f:0
\move(19 126)
\lvec(20 126)
\lvec(20 127)
\lvec(19 127)
\ifill f:0
\move(19 129)
\lvec(20 129)
\lvec(20 130)
\lvec(19 130)
\ifill f:0
\move(22 129)
\lvec(23 129)
\lvec(23 130)
\lvec(22 130)
\ifill f:0
\move(15 130)
\lvec(16 130)
\lvec(16 131)
\lvec(15 131)
\ifill f:0
\move(15 131)
\lvec(16 131)
\lvec(16 132)
\lvec(15 132)
\ifill f:0
\move(19 131)
\lvec(20 131)
\lvec(20 132)
\lvec(19 132)
\ifill f:0
\move(22 132)
\lvec(23 132)
\lvec(23 133)
\lvec(22 133)
\ifill f:0
\move(12 133)
\lvec(13 133)
\lvec(13 134)
\lvec(12 134)
\ifill f:0
\move(19 134)
\lvec(20 134)
\lvec(20 135)
\lvec(19 135)
\ifill f:0
\move(12 135)
\lvec(13 135)
\lvec(13 136)
\lvec(12 136)
\ifill f:0
\move(19 137)
\lvec(20 137)
\lvec(20 138)
\lvec(19 138)
\ifill f:0
\move(15 138)
\lvec(16 138)
\lvec(16 139)
\lvec(15 139)
\ifill f:0
\move(15 139)
\lvec(16 139)
\lvec(16 140)
\lvec(15 140)
\ifill f:0
\move(19 140)
\lvec(20 140)
\lvec(20 141)
\lvec(19 141)
\ifill f:0
\move(22 142)
\lvec(23 142)
\lvec(23 143)
\lvec(22 143)
\ifill f:0
\move(19 143)
\lvec(20 143)
\lvec(20 144)
\lvec(19 144)
\ifill f:0
\move(19 145)
\lvec(20 145)
\lvec(20 146)
\lvec(19 146)
\ifill f:0
\move(22 145)
\lvec(23 145)
\lvec(23 146)
\lvec(22 146)
\ifill f:0
\move(12 146)
\lvec(13 146)
\lvec(13 147)
\lvec(12 147)
\ifill f:0
\move(15 146)
\lvec(16 146)
\lvec(16 147)
\lvec(15 147)
\ifill f:0
\move(15 147)
\lvec(16 147)
\lvec(16 148)
\lvec(15 148)
\ifill f:0
\move(19 148)
\lvec(20 148)
\lvec(20 149)
\lvec(19 149)
\ifill f:0
\move(30 -20)
\htext{number of points $r \in [10,450]$}
\move(-2 20)
\vtext{multiplicity $m \in [1,150]$}
\arrowheadtype t:V
\linewd 1
\move(10 -5)
\avec(470 -5)
\move(5 0)
\avec(5 170)
}
\caption{Graph showing when {\sc NSglue} gives conjectured value of $\alpha_2$}\label{pic2}
\end{figure}

Bounding $\tau_2$ with {\sc NSglue} is much easier than bounding $\alpha_2$, 
since, for any $m$, there exists
non-special system with very low positive virtual dimension
($\vdim \sys_2(3m;m^{\times 9}) = 1$), which is not true for
systems based on $s^2$ points and negative virtual dimension.

\subsection{Case $n=3,4$}

For $n=3$ and $n=4$ very few is known. The systems with multiplicities
all equal to $2$ have been characterized in \cite{AH2}. For $r \leq 8$, $n=3$
the necessary and sufficient condition of non-specialty can be found in 
\cite{p3}.
For $r \leq n+1$ the methods presented in the proof of Thm. \ref{cremona}
are sufficient to decide about non-specialty.
For $r \geq 8$ there are,
so far, no nontrivial algorithm to bound $\alpha_3$ and $\alpha_4$.

For $n=3$ we will consider homogeneous systems and use {\sc NSglue} with
\begin{align*}
\text{{\sc Eclass}} & = \{ L : L \text{ is non-special by Crit. \ref{crendcr}}\},\\
\text{{\sc Gclass}} & = \{ L = \sys_3(d;m^{\times r}) : r \in \{4,\dots,8\} \text{ and $L$ is non-special by Crit. \ref{crendcr}}\}.
\end{align*}

\begin{example}
We will show that $\alpha_3(5^{\times 50}) = 20$. Since (by computation
of virtual dimension) we have $\alpha_3(5^{\times 50}) \leq 20$, it
is enough to show that $\sys_3(19;5^{\times 50})$ is empty.
Observe that (by Crit. \ref{crendcr}) systems
$\sys_3(9;5^{\times 7})$ and $\sys_3(8;5^{\times 5})$ are empty.
By ``glueing'' points (six times using $\sys_3(9;5^{\times 7})$ and then
$\sys_3(8;5^{\times 5})$ once) it is enough to show that
$\sys_3(19;10^{\times 6},9,5^{\times 3})$ is empty.
This is true by Crit. \ref{crendcr}.
\end{example}

Some bounds for $\alpha_3(m^{\times r})$ obtained
with {\sc NSglue} are presented in Tab. \ref{tab2}.
First we present the conjectured value, then (after comma) our bound.

\begin{table}[ht!]
$$
\begin{array}{|c|c|c|c|c|}
\hline
  n=3   & r=30 & r=50 & r=100 & r=150 \\
\hline
m=5 & 17,16 & 20,20 & 26,24 & 30,28 \\
\hline
m=10 & 33,31 & 39,37 & 49,46 & 57,52 \\
\hline
m=20 & 64,61 & 76,74 & 96,91 & 110,104 \\
\hline
m=50 & 157,151 & 186,182 & 235,226 & 269,258 \\
\hline
m=100 & 312,301 & 371,364 & 467,452 & 535,514 \\
\hline
\end{array}
$$
\caption{Some bounds for $\alpha_3$ obtained with {\sc NSglue}}\label{tab2}
\end{table}

The same can be done for $n=4$ (here, in Tab. \ref{tab3}, we present 
results bounding $\alpha_4$, then we give an example).

\begin{table}[ht!]
$$
\begin{array}{|c|c|c|c|}
\hline
n=4  & r=30 & r=50 & r=100 \\
\hline
m=5 & 13,12 & 15,14 & 18,16 \\
\hline
m=10 & 25,23 & 28,26 & 34,31 \\
\hline
m=20 & 48,44 & 55,51 & 66,59 \\
\hline
m=50 & 119,110 & 135,126 & 161,147 \\
\hline
\end{array}
$$
\caption{Some bounds for $\alpha_4$ obtained with {\sc NSglue}}\label{tab3}
\end{table}

\begin{example}
The conjectured value of $\alpha_4(5^{\times 50})$ is $15$.
We will show that $\alpha_4(5^{\times 50}) \geq 14$.
By Crit. \ref{crendcr} the system $\sys_4(8;5^{\times 8})$ is empty.
We begin with system $\sys_4(13;5^{\times 50})$, glue points
to obtain $\sys_4(13;9^{\times 6},5^{\times 2})$ and use Crit.
\ref{crendcr} to show that the last system is empty.
\end{example}

Our algorithm can also be used in
higher dimension. Observe that usually, for $n \geq 3$, if
$(n-1)d-\sum_{j=1}^{n+1}m_j < 0$ then
$\vdim \Cr(\sys_n(d;m_1,\dots,m_r)) > \vdim \sys_n(d;m_1,\dots,m_r)$.
It follows that the birational transformation is useless in proving
regularity for $n \geq 3$.

\section{Bounds for Seshadri constants}

Let $p_1,\dots,p_r \in \PP^2$ be points
in general position. Define the {\dff Seshadri constant}
of the line bundle $\stbun$ (see \cite{Dema} or \cite{szemb})
$$\varepsilon(\stbun,r) \rdf \inf \bigg\{ \frac{\deg C}{\sum_{j=1}^r \mult_{p_j}C} : C \text{ is a curve on } \PP^2 \bigg\}.$$
The famous Nagata Conjecture (see \cite{nagata}) states (in the language of Seshadri constants) that
$$\varepsilon(\stbun,r) = \frac{1}{\sqrt{r}}$$ 
for $r \geq 9$. The conjecture is still open. It is known that
$\varepsilon(\stbun,r) \leq \frac{1}{\sqrt{r}}$, so we may look
for lower bounds for Seshadri constant.
We will use the following fact:

\begin{theorem}
\label{quasi}
Any irreducible and reduced curve $C \subset \PP^2$ disproving Nagata
Conjecture for $r \geq 10$
is almost homogeneous, that is there exists $d,m,m_0 \in \N$ such that
$C \in \sys_2(d;m^{\times r-1},m_0)$.
\end{theorem}

\begin{proof}
See \cite{syzdek}, Prop. 2.10, where the above theorem is stated
for a polarized surface $(X,L)$ with Picard number equal one, two or three,
or \cite{harroe2}, Prop. 3.2, where additional conditions on $d,m,m_0$
are given.
\end{proof}

\begin{proposition}
\label{22}
$$
\varepsilon(\stbun,22) \geq \frac{42}{197} = \frac{1}{\sqrt{22+\frac{1}{1764}}}.
$$
\end{proposition}

\begin{proof}
Observe that it is enough to exclude the existence of
an irreducible and reduced curve $C \subset \PP^2$ such that
$$\frac{\deg C}{\sum_{j=1}^{22} \mult_{p_j}C} < \frac{42}{197}.$$
By Thm. \ref{quasi} such curve $C$ would be a member
of $\sys_2(d;m^{\times 21},h)$ for some $d,m,h \in \N$ such that
$\frac{d}{21m+h} < \frac{42}{197}$. The last system is empty (see Ex. \ref{197}).
\end{proof}

Using the same technique one can show that (see \cite{MYWWW})
$$
\varepsilon(\stbun,19) \geq \frac{39}{170} = \frac{1}{\sqrt{19+\frac{1}{1521}}}.
$$
For $r=10$ the above technique gives the lower bound equal to $\frac{3}{10}$
(glue $\sys_2(d;m^{\times 9},h)$ to $\sys_2(d;2m,m^{\times 5},h)$;
consider cases $m \leq h$ and $m > h$ separately), but
much better bounds are known (see \cite{harroe2}). We note here that for $r=19$ and $r=22$
our bounds are equal to bounds presented in \cite{harroe2} (given after Biran \cite{Biran}).

Since {\sc NSglue} works better for regularity (than for $\alpha_2$) we
will rather make use of the following Theorem 
(proposed and proved by Eckl in \cite{ecklc}):

\begin{theorem}[Eckl]
\label{eckl}
Let $r>9$, let $(d_j,m_j)$ be a sequence such that
\begin{itemize}
\item
$d_j \geq \tau_2((m_j+1)^{\times r})$,
\item
$\frac{d_j^2}{m_j^2 r} \stackrel{j \to \infty}\longrightarrow \frac{1}{a^2} \geq 1$.
\end{itemize}
Then
$$\varepsilon(\stbun,r) \geq a\frac{1}{\sqrt{r}}.$$
\end{theorem}

To estimate the Seshadri constant for $r$ points
we run {\sc NSglue} to find $d(r,m) \geq \tau_2(m^{\times r})$ for $m$ varying
from $1$ to
some $M \in \N$, depending on time of computations.
For each $m$ and $d(r,m)$ we consider $a=\frac{m-1}{d(r,m)}\sqrt{r}$ and 
constant sequence $(d_j,m_j)=(d(r,m),m-1)$
satisfying the assumptions of Eckl's Theorem.
Hence
$$\varepsilon(\stbun,r) \geq \frac{m-1}{d(r,m)}\sqrt{r}\frac{1}{\sqrt{r}} = \frac{m-1}{d(r,m)}.$$

The bound is given by the maximum of
$$\frac{m-1}{d(r,m)}, \quad m=1,\dots,M.$$
A suitable computer program can be found at \cite{MYWWW}.
Here (Tab. \ref{tab4}) we present our bounds for $r=10,\dots,24$ together
with non-special system giving
the above maximum. Note that for $r=15,20,24$ our bounds are slightly better than bounds
presented in \cite{harroe2}.

\begin{table}[ht!]
$$
\begin{array}{|c|c|c||c|c|c|}
\hline
r & \text{bound} & \text{system} & r & \text{bound} & \text{system} \\
\hline
10 & 100/331 & \sys_2(26480;8001^{\times 10}) & 18 & 8096/34413 & \sys_2(34413;8097^{\times 18})  \\
\hline
11 & 8098/26997 & \sys_2(26997;8099^{\times 11}) & 19 & 39/170 & - \\
\hline
12 & 8098/28347 & \sys_2(28347;8099^{\times 12}) & 20 & 1617/7235 & \sys_2(36175;8086^{\times 20}) \\
\hline
13 & 8098/29697 & \sys_2(29697;8099^{\times 13}) & 21 & 4000/18403 & \sys_2(36806;8001^{\times 21}) \\
\hline
14 & 4049/15186 & \sys_2(30372;8099^{\times 14}) & 22 & 42/197 & - \\
\hline
15 & 426/1651 & \sys_2(31369;8095^{\times 15}) & 23 & 8095/38862 & \sys_2(38862;8096^{\times 23}) \\
\hline
17 & 2000/8473 & \sys_2(33892;8001^{\times 17}) & 24 & 8092/39657 & \sys_2(39657;8093^{\times 24}) \\
\hline
\end{array}
$$
\caption{Lower bounds for $\varepsilon(\stbun,r)$ for $r=10,\dots,24$}\label{tab4}
\end{table}

Studying behavior of {\sc NSglue} for $n=2$ and homogeneous systems
we may observe that it works very well for $r=k^2-1$, on the other
hand the worst results are obtained for $r=k^2+1$ and, surprisigly, $r=13$.

One may ask if {\sc NSglue} can compute Seshadri constant,
i.e. if
$$\lim_{m \limto \infty} \frac{m-1}{d(r,m)} = \varepsilon(\stbun,r),$$
where $d(r,m)$ is the bound on $\tau_2(m^{\times r})$ obtained with
our algorithm. The answer is not known.


\begin{thebibliography}{Laf--Uga 03b}
%
%
\bibitem[Ale--Hir 95]{AH2}
J. Alexander, A. Hirschowitz,
\textit{Polynomial interpolation in several variables},
J. Algebraic Geometry 4,
201--222 (1995).
%
\bibitem[Bir 99]{Biran}
P. Biran,
\textit{Constructing new ample divisors out of old ones},
Duke Math. J. 98,
113--135 (1999).
%
\bibitem[Cha 01]{chandler}
K. Chandler,
\textit{A brief proof of a maximal rank theorem for generic double points in projective space},
Trans. Amer. Math. Soc. 353,
1907--1920 (2001).
%
%
%
\bibitem[Cil--Mir 98]{CMirdeg}
C. Ciliberto, R. Miranda,
\textit{Degeneration of planar linear systems},
J. Reine Angew. Math. 501,
191--220 (1998).
%
\bibitem[Cil--Mir 00]{CMir}
C. Ciliberto, R. Miranda,
\textit{Linear systems of plane curves with base points of equal multiplicity},
Trans. Amer. Math. Soc. 352,
4037--4050 (2000).
%
\bibitem[Dem 92]{Dema}
J.P. Demailly,
\textit{Singular Hermitian metrics on positive line bundles},
Complex Algebraic Varieties (K. Hulek at al., eds.), LNM, vol. 1507, 87--104,
Springer, 1992.
%
\bibitem[Vol--Laf 07]{p3}
C. De Volder, A. Laface,
\textit{On linear systems of $\PP^3$ through multiple points},
J. Algebra 310, 207--217 (2007)
%
\bibitem[Dum 08]{MYWWW}
M. Dumnicki,
\texttt{http://gamma.im.uj.edu.pl/dumnicki/interpol.htm}, February 7, 2008.
%
\bibitem[Dum--Jar 07]{mdwj}
M. Dumnicki, W.Jarnicki,
\textit{New effective bounds on the dimension of a linear system in $\PP^2$},
J. Symbolic Comput. 42,
621--635 (2007).
%
\bibitem[Dum 07]{md}
M. Dumnicki,
\textit{Cutting diagram method for systems of plane curves with base points},
Ann. Polon. Math. 90,
131--143 (2007).
%
\bibitem[Eckl 05]{ecklc}
T. Eckl,
\textit{Seshadri constants via Lelong numbers},
arXiv:math/0508561,
to be published in Math. Nachr. (2005).
%
\bibitem[Har 08]{harwww}
B. Harbourne,
\texttt{http://www.math.unl.edu/\~{}bharbourne1/FatPointAlgorithms.html}, February 7 2008.
%
\bibitem[Har 02]{surv}
B. Harbourne,
\textit{Problems and progress: A survey on fat points in $\PP^2$},
Zero-dimensional schemes and applications (Naples, 2000), Queen's Papers in Pure and Appl. Math., 123, 
Queen's Univ., Kingston, ON, 85--132 (2002).
%
\bibitem[Har 03]{harsesh}
B. Harbourne,
\textit{Seshadri constants and very ample divisors on algebraic surfaces},
J. Reine Angew. Math. 559, 115--122 (2003).
%
\bibitem[Har 05]{knowalpha}
B. Harbourne,
\textit{The (unexpected) importance of knowing $\alpha$},
Projective varieties with unexpected properties,
Walter de Gruyter GmbH \& Co. KG, Berlin, 267--272 (2005).
%
\bibitem[HHF 03]{tausquare}
B. Harbourne, S. Holay, S. Fitchett,
\textit{Resolution of ideals of quasiuniform fat points subschemes of $\PP^2$},
Trans. Amer. Math. Soc. 355, 593--608 (2003).
%
\bibitem[Har--Ro\'e 04]{harroe}
B. Harbourne, J. Ro\'e,
\textit{Linear systems with multiple base points in $\PP^2$},
Adv. Geom. 4, 41--59  (2004).
%
\bibitem[Har--Ro\'e 07]{harroe2}
B. Harbourne, J. Ro\'e,
\textit{Computing multi-point Seshadri constants on $\PP^2$},
arXiv:math/0309064v3 (2007).
%
\bibitem[Hir 85]{HirF}
A. Hirshowitz,
\textit{La m\'ethode d'Horace pour l'interpolation \`a plusieurs variables},
Manus. Math. 50, 337--388 (1985).
%
\bibitem[Hir 89]{HCON}
A. Hirshowitz,
\textit{Une conjecture pour la cohomologie des diviseurs sur les surfaces rationelles g\'en\'eriques},
J. Reine Angew. Math. 397, 208--213 (1989).
%
\bibitem[Laf--Uga 06]{conj3}
A. Laface, L. Ugaglia,
\textit{On a class of special linear systems of $\PP^3$},
Trans. Amer. Math. Soc. 358, 5485--5500 (2006).
%
\bibitem[Laf--Uga 03]{base5}
A. Laface, L. Ugaglia,
\textit{Quasi-homogeneous linear system on $\PP^2$ with base points
of multiplicity 5},
Canad. J. Math. 55,
561--575 (2003).
%
\bibitem[Mon 07]{M}
F. Monserrat,
\textit{Curves having one place at infinity and linear systems on rational surfaces},
J. Pure Appl. Algebra 211, 685--701  (2007).
%
\bibitem[Nag 59]{nagata}
M. Nagata,
\textit{On the 14-th problem of Hilbert},
Amer. J. Math. 81,
766--772 (1959).
%
\bibitem[Ro\'e 01a]{roe1}
J. Ro\'e,
\textit{On the existence of plane curves with imposed multiple points},
J. Pure Appl. Alg., 156, 115--126 (2001).
%
\bibitem[Ro\'e 01b]{roe2}
J. Ro\'e,
\textit{Linear systems of plane curves with imposed multiple points},
Illinois J. Math. 45, 895--906  (2001).
%
\bibitem[Syz 05]{syzdek}
V. Syzdek,
\textit{Seshadri constants and geometry of surfaces},
\texttt{http://www.uni-essen.de/\~{}mat903/alggeo.html}, February 7, 2008 (2005).
%
\bibitem[Sze 01]{szemb}
T. Szemberg,
\textit{Global and local positivity of line bundles},
Habilitation,
Essen (2001).
%
\end{thebibliography}
\end{document}